%% file: 447instances.tex
\documentclass{article}
\usepackage{array}
\usepackage{amsmath}
\usepackage{graphicx}
\usepackage{url}
\usepackage{maplestd2e}
\def\emptyline{\vspace{12pt}}
\DefineParaStyle{Maple Heading 1}
\DefineParaStyle{Maple Text Output}
\DefineParaStyle{Maple Dash Item}
\DefineParaStyle{Maple Bullet Item}
\DefineParaStyle{Maple Normal}
\DefineParaStyle{Maple Heading 4}
\DefineParaStyle{Maple Heading 3}
\DefineParaStyle{Maple Heading 2}
\DefineParaStyle{Maple Warning}
\DefineParaStyle{Maple Title}
\DefineParaStyle{Maple Error}
\DefineCharStyle{Maple Hyperlink}
\DefineCharStyle{Maple 2D Math}
\DefineCharStyle{Maple Maple Input}
\DefineCharStyle{Maple 2D Output}
\DefineCharStyle{Maple 2D Input}
\usepackage[top=.9in,bottom=.9in]{geometry}
 \numberwithin{equation}{section}
 \def\emptyline{\vspace{12pt}}
\DefineParaStyle{Maple Output}
\DefineCharStyle{2D Comment}
\DefineCharStyle{2D Math}
\DefineCharStyle{2D Output}
\begin{document}
\pagestyle{plain}
\begin{maplegroup}
\widowpenalty=600
\clubpenalty=600
\title{447 Instances of Hypergeometric 3F2(1)}
\author{Michael S. Milgram \footnote{mike@geometrics-unlimited.com} \\
\it Consulting Physicist, Geometrics Unlimited, Ltd. \\
\it Box 1484, Deep River, Ont. Canada. K0J 1P0 }
\date{May 12, 2011}
\maketitle

\setcounter{page}{1}
\pagenumbering{arabic}

\begin{flushleft}
\begin{abstract}
\vspace*{0.2cm}
259 new instances of hypergeometric 3F2$(...| 1)$ evaluations are obtained by systematically testing three-part 
transformations among these functions, against the previously known database of 133 such evaluations. A complete database of 447 evaluations is listed in a series of appendices, including one appendix where all cases are given in the form of packed pseudo-code, suitable for pasting as input instructions, to a computer algebra code.

\end{abstract}


\section {Introduction}

The need to evaluate instances of the function $_3F_2(...|1)$ in closed-form for particular combinations of its parameters is well-known and ubiquitous in mathematical physics, complicated by the fact that 
no universal closed-form exists, in contrast to the case of $_2F_1(...|1)$. As has been eloquently demonstrated by 
Petkov\u sek, Wilf and Zeilberger \cite{PWZ} in their book `A=B', not only does a complete 
database of such evaluations not exist, but it cannot exist because the WZ and 
related algorithms will always generate new and unexpected identities from any new (or old) $_3F_2(...|1)$ identity that comes to light. 
This was amply demonstrated by Gessel \cite{Gessel} who produced copious combinatorial identities utilizing the WZ algorithms, some of which involved $_3F_2(...|1)$. 
In other work, Chu \cite{Chu93}, Exton \cite{Exton}, Coffey \cite{Coffey}, Vidunas \cite{Vidunas} (and references therein), Krattenthaler and Rao \cite{KratRao} and Kim, Rakha and Rathie \cite{KRR} have shown alternate ways to generate hypergeometric identities based on contiguity relations among known results, inversion techniques, integral representations and various transformations and reductions of Kamp\' e de F\' eriet functions, some of which produce identities involving $_3F_2(...|1)$.\newline

Unfortunately, most attempts to apply any of these methods to obtain a `needed' new result will encounter a common inadequacy, notably the inability to predict the properties of any new identity 
that may emerge from the application of any of these algorithms. In fact, for the WZ algorithm, there is no assurance that the resulting identity will even involve a hypergeometric function $_pF_q(...|1)$.
Thus, although Gessel's listing of new identities contains many interesting elements alongside those presented in \cite{PWZ}, the form of none of them could have been predicted prior 
to the application of the WZ algorithm. 
The same can be said of the results presented in most of the other works cited above. Furthermore, it so happens that many of the ``new'' identities that reduce to a hypergeometric identity involving $_pF_q(...|1)$ also happen to share the property that they either terminate, and/or have a top parameter that exceeds a bottom parameter by a small positive integer and, as was first shown by Minton \cite{Minton} and elaborated upon by others \cite{Rosengren},\cite{GottMas},\cite{Karlsson},
such a hypergeometric function can always be reduced to a combination of $_{p-1}F_{q-1}(...|1)$ functions. In the case that $p=3$ and $q=2$, this means that any claimed ``new'' identity 
``discovered'' in these works can always be reduced to several $_2F_1(...|1)$, and their novelty challenged.\footnote{It is noted that both Kim and Rathie \cite{KimRat}, and Coffey \cite{Coffey} have invoked Minton's result to challenge the novelty of claimed ``new'' evaluations by both Miller \cite{Miller05} and Exton \cite{Exton} respectively.} In fact, I have found no new $_3F_2(...|1)$ evaluations in any of these works.\newline 

So, what is an analyst to do when faced with the need to evaluate a particular $_3F_2(...|1)$ that may have arisen in some physics calculation? Typically, (s)he will survey the literature, praying that
the particular case that has arisen happens to be included in a compendium of known identities, notwithstanding the comment in \cite{PWZ} that such a database cannot be complete. 
It is the purpose of this work to provide such an (incomplete) compendium.  

\section {Three-Part Relations}

\normalsize
In two previous papers \cite{Milgram1},\cite{Milgram2}, the two-part Thomae relations between hypergeometric functions $_3F_2(...|1)$ were systematically inspected in a search for new identities \footnote{ It is likely that many of the  Kamp\' e de F\' eriet reductions cited above consist of the 2-part Thomae transformations embedded into the Kamp\' e de F\' eriet notation.}. 
The initial  elements for the search were obtained by transcribing all relevant evaluations found in Prudnikov et. al. \cite{Prudnikov}, augmented by a collection of 
other results garnered from the literature. A database consisting of 133 $_3F_2(...|1)$ identies were cumulatively listed in those papers.\footnote{For completeness sake, 
these 133 cases are included here.}
The use of the 2-part Thomae identities to generate other $_3F_2(...|1)$ identities is well-established - Bailey long ago showed that the theorems of Watson and Whipple (\cite{Bailey}) 
could be obtained by 
applying one of the 10 Thomae relations to Dixon's Theorem (\cite{Bailey}), and in \cite{Milgram1} I observed that the 2-part Thomae identities closed this relationship among any of the contiguous Watson, Dixon or Whipple 
relations investigated by others (e.g. Lewanowicz \cite{Lewan}, and references in \cite{Milgram1}).\newline

However, other transformations between $_3F_2(...|1)$ sums are also well-known, these being the 3-part transformations introduced\footnote{ see \cite{KratRao} for complete historical references.} in 1879 by Thomae, categorized by Whipple in 1923, later summarized by Bailey (\cite{Bailey}, Section 3.5 to 3.7) and more recently reproduced by Luke (\cite{Luke}, Section 3.13.3).
The existence of such transformations gives rise to the possibility that additional $_3F_2(...|1)$ evaluations can be found, by applying each of these three-part transformations to the 
existing 133 member database and searching for cases where only one member of a multi-element transformation equation does not correspond to a previously known 
evaluation. This work is a report on such a search; it is worth noting that this method guarantees that any new identities uncovered will be of, at most, the same degree as the starting identity - $_3F_2(...|1)$ - in contrast to some of the methods cited in the previous section.\newline      

The collection of available three-part transformations is typified by six examples given identically in each of Bailey's and Luke's book, based on permutations among 120 different hypergeometric series, applied to a single three-part identity derived from a contour integral representation. The 120 different hypergeometric series involved are characterized by six parameters $r_{\it{i}}$ defined by a set of equations and the constraint $\sum_{i=0}^{5}r_{i}= 0$. \newline

Conveniently, Miller \cite{Miller03} has recently obtained two multi-part relations between hypergeometric $_4F_3(...|1)$. By equating 
Miller's parameters $g=d$, the first relation will reduce to the equivalent transformation among four $_3F_2(...|1)$ functions. It is probable that this transformation is already embedded somewhere in Whipple's results quoted in Bailey's book, but since it provides a reasonable starting point for a search for new instances of $_3F_2(...|1)$, this transformation was used here. It is referred 
to as ``Miller, Eq. 1'' in the results. Miller also comments that his second result reduces, in various limits, to one or another of the six instances of 3-part transformations presented by Bailey (and Luke). 
Since such a reduction provides a convenient characterization
of Bailey's results, various reductions of Miller's second transformation were used here to label different 3-part transformations of $_3F_2(...|1)$.\newline 

Thus the calculation proceeds as follows. Starting with the database of 133 known $_3F_2(...|1)$ identities, successively apply each of the multi-part transformations given by the reduction of 
one of Miller's two $_4F_3(...|1)$ transformations in various limits, and search for cases where only one of the elements of a multi-part transformation is not already part of the database\footnote{No attempt was made to search for simultaneous solutions to combinations of two (or more) equations.}. When such a case 
is found, add the new identity to the database, if it is viable 
({ \it `viable' - does not contain a factor characterized by $\Gamma (-n)$, n $\ge 0 $, and does not have a parametric excess equal to a
non-positive integer unless one of the top parameters is a negative integer}). Any new, viable, progeny generated by any of the ten 2-part Thomae relations belonging to each new 
case is also added to the database. 
As each new identity is added to the database, 
the search is repeated (recursively) until all transformations have been tested against all elements of the enlarged database.\newline
 
To check that all independent transformations were sampled after all of the reductions from Miller's transformations were searched, the same procedure was 
employed using the six specific examples cited by Bailey and Luke. All calculations were performed using the Maple 11 computer program.

\section {New Database}

Throughout this work, the variables $L$, $m$ and $n$ are positive integers, and the results are presented in a series of appendices. The starting database of 133 identities was described in previous work, with one significant modification. 
Recently, Chu (\cite{Chu11}) has presented new, closed expressions for elements contiguous to Whipple's, Watson's and Dixon's theorems. \newline

These are significantly improved from the results obtained by Lewanowicz \cite{Lewan}, and others ( see \cite{Milgram1}, \cite{Milgram2} and 
references therein), where equivalent relations are summarized and expressed recursively. Chu's new results were prepended to the list of known identities, appear in all the appendices that follow as case numbers 1 to 19, and are therefore counted among the 133 ``starting'' cases. 
Matching the first 19 cases in appendices A and B effectively defines symbols that appear elsewhere in the appendices; these symbols were chosen to correspond to each of Chu's cases, and are labelled $\it{Chu, Math.Comp}$ in the appendices.\newline

Thus WPPR(m,n,a,b,c) corresponds to Chu's result ${W_{m,n}}(a,b,c)$ for Watson's theorem (\cite{Chu11} Theorem 5) when 
the variables $m$ and $n$ are both associated with $\it{Plus}$   signs.
WMPR corresponds to the 
case when the variable $m$ is preceded by a $\it{Minus}$ and $n$ by a $\it{Plus}$ sign, and WMMR corresponds to the case that both parameters are associated with $\it{Minus}$ signs. 
The symbols DIXONPM(m,n,a,b,c) and WHIPPLEMP(m,n,a,b,c) embed corresponding nomenclature, corresponding to Chu's results for Dixon's and Whipple's theorems respectively (\cite{Chu11}). 
The equivalence between Chu's symbols as they are used here and Lewanowicz' results, 
denoted by ${F_{m,n}}(a, \,b, \,c)$ in \cite{Milgram2}, is

\footnotesize
\vspace*{.5 cm}
\maplemultiline{
{F_{2\,m, \,2\,n - 1}}(X, \,Y, \,Z)=\mathrm{WPMR}(2\,m + 2\,n - 1
, \,2\,n - 1, \,X, \,2\,n - 1 + 2\,Z - Y, \,2\,n - 1 + Z) \times \\
\Gamma ( - {\displaystyle \frac {X}{2}}  + 2\,n - {\displaystyle 
\frac {1}{2}}  + Z - {\displaystyle \frac {Y}{2}}  + m)\,\Gamma (
{\displaystyle \frac {X}{2}}  + {\displaystyle \frac {1}{2}}  + 
{\displaystyle \frac {Y}{2}}  + m) \left/ {\vrule 
height0.80em width0em depth0.80em} \right. \!  \! (\Gamma ( - 
{\displaystyle \frac {X}{2}}  + {\displaystyle \frac {1}{2}}  + 
{\displaystyle \frac {Y}{2}}  + m) 
\Gamma ({\displaystyle \frac {X}{2}}  + 2\,n - {\displaystyle 
\frac {1}{2}}  + Z - {\displaystyle \frac {Y}{2}}  + m)) }
\vspace*{.5 cm}
\normalsize

Relatives of Chu's fundamental cases survive in the new database because of a strict condition in the identification procedure requiring that the parameters $L$, $m$ and $n$ be positive integers. 
At times, a mapping found by the search algorithm will fail to associate a putative ``new'' result with a corresponding ``known'' result if 
a mapping of parameters between the two violates this requirement. For example, two cases of $_3F_2(...|1)$ that are related by a parametric mapping of the form $n \rightarrow n-1$ 
will be treated as different, to preclude the possibility that cases corresponding to $n=0$ will enter the mix. Such cases are identified in the results (see Appendix E, notes) and are maintained within 
the database to improve the chance that any search for a parametric mapping between a candidate and a member of the database, will find one or the other of such closely-related cases,
since, in general, it may not find both.\newline

The database results are given in the appendices. Each appendix lists the appropriate property, labelled by a case number. Missing labels in these appendices correspond 
to special results in two categories: 22 cases which escaped the viability test until the last run when all labels had already been established \footnote{ the following labels are therefore omitted from the appendices: 22, 32, 49, 58, 59, 60, 83, 251, 291, 292, 365, 369-374, 464-468}, and 116 cases that are special instances of more general cases included in the prior appendices. 
These 116 cases are retained in the tables and database because, in many instances, it is easier to evaluate a special case starting from a special result, compared to reducing a general result to a
special one (especially if the reduction involves a complicated limiting process), but these cases are not included in the count of ``new'' results. 
Thus the full database consists of 447 identities, 335 of which are general and viable, 259 of which are new, according to the search algorithm, and 112 0f which are special cases of one of the (335) general results.\newline

The 335 main results appear in appendices A and B. Appendix B lists the closed form for each of the corresponding sums appearing in Appendix A. The symbols corresponding to Chu's evaluation of Watson-Whipple-Dixon
cases are described above. As well, symbols $V1(a,n)$, $V7(a,n)$ and $V8(a,n)$ arise, corresponding to progeny of Maier's result\cite{Maier}. The recursive definition of each of these functions can be found
in \cite{Milgram2}.\newline

Many of the results in Appendix B have been $\it{simplified}$, from the form in which they were first generated, but the simplification is not always complete. This is a by-product of Maple's ``simplify'', ``combine'', ``factor'' and ``expand'' commands, which do not always capture the essence of a very long expression. Even in those cases where further reduction is obvious (e.g. extraneous multiples of $ 2\pi $), the reader should be assured that the results given are always shorter than the original. Each reduced version was tested numerically with several sets of parameters.\newline

Appendix C gives a traceback for each of the corresponding cases, from which it is possible to rederive the corresponding result. For example, to reproduce case 392: to Lewanowicz' case 82 \cite{Lewan} ( $m - odd, n-even$ ), apply Miller's \cite{Miller03} equation (1.2) with the reduction $e=a$, to yield case 368\footnote{confession: In each of the 3-part transformations, I neglected to retain the identity of the known term(s) that appear(s) in  the right-hand side of the equation. It will always be a case whose numerical label is smaller than the label of the case that is found (e.g. in this example, its label will be less than 368).}. To that result, apply the eighth Thomae 2-part transformation (T:8) listed in \cite{Milgram1}, from which case 383 will be found. Add this new case to the database and re-apply Miller's Eq.(1.2) with $e=a$. Case 392 will be found.\newline

Appendix D lists 112 cases analogous to Appendix A, except that all entries in this list are special results corresponding to more general cases appearing in Appendix A. Appendix E contains special notes that relate the contents of Appendix D to corresponding results presented in appendices A and B, and Appendix F gives the corresponding results for those entries appearing in Appendix D. Appendix G is a packed cut-and-paste set of pseudo-code that can be used as direct input to a computer algebra code.

\section {Comments on the Database}

The search program used throughout this work does not consider contiguity relations, although many of the results are closely related, but distinct. For example, compare case 142 after the mapping $a \rightarrow a-1/2$, $ n \rightarrow n+1 $ with case 140. The two only differ by unity in one of the bottom parameters. This suggests that a complete new family of results could be found by exploring a 3-part recursion in that parameter between those two cases.\newline

Cases 91 to 93 should be noted. In \cite{Milgram1} these cases were introduced in the form of conjectures; later, Chu \cite{Chu_private} provided a proof of these conjectures, after which it was belatedly noticed that Maier \cite{Maier06} had already identified these cases as 2-part Thomae progeny of a well-known result, misprinted in Slater's book \cite{Slater}. This is discussed in \cite{Milgram2}, and because of the sequence of awareness, these cases are labelled ``CHU Wenchang'' in the appendices.\newline

To cases 55, 56 and 57, all of which equate to zero, may be added the new case 252. This case has parametric excess $s=2$. There are also 5 other cases with $s=2$, these being cases 86, and 247 from Appendix A, and cases 53, 90, and 93 from Appendix D.\newline

Cases 29, 34 and 37 from the original set have $s=1$ and are therefore Saalsch\" utzian (\cite{Luke} , Eq. 3.13.3(4)). To these cases, add the following from Appendix A: cases 150, 164, 217, 223, 228, 235, 239, all of which are non-terminating and Saalsch\" utzian. Additionally, special cases labelled 40, 66, 146, 147, 180 and 189 from Appendix D share this property, although only one of the more general cases (217) for which these represent special values of the parameters, appears in the list from Appendix A.\newline

To the list of cases with constant values of $s$, add cases 28 and 248 with $s=n$, the terminating case 24 with $s=m$ and, from appendix A, cases 281, 314, 331, 340 and 362 with $s=1/2+n$. An interesting set of results corresponds to six cases with $s=1/2\,-\,n$ - cases 108, 275, 288, 392, 327 and 360. The numerical validity of these cases has been checked when $n=0$. Under the mapping $n\rightarrow -n, a\rightarrow 1/2-a $, each of these cases except one, maps into one of the cases with $s=n+1/2$. Specifically, case 108 maps to case 362, case 275 maps to case 340, case 288 maps to 331, case 327 maps to case 281 and case 360 maps to case 314. The exception is case 302 which, after setting $n \rightarrow -n$ maps into case 123 (where $s=n-1/2$) under the mapping ${a \rightarrow 1/6-a- n/2}$, $ {n \rightarrow 2\,n+1}$. From this, I conjecture that the right-hand sides of those cases with $s=1/2-n$ convey meaning to the left-hand sides for $n>0$, capturing the spirit (but not the embodiment!), of analytic continuation in the (discrete) variable $n$.
Finally, from Appendix A, cases 133, 423, 432, 443, 448 and 452 all share the property $s=1/2$, and from Appendix D, case 103 has $s=3/2$.\newline

Overall, 40 of the cases are nearly-poised, 16 appearing in Appendix A and 24 in Appendix D. These cases are marked by a $\dagger$ prepended to the case numbers in each of these two appendices. Because of the requirement that $m,n>0$, Dixon's classical result for a well-poised series corresponding to one of cases 17,18 or 19 (with $m,n=0$) are not so-distinguished.\newline

Finally, it should be noted that, although the emphasis on the results presented here has been directed toward the identification and reduction of infinite series, in various limits, cases involving terminating series can be analyzed. Recently Rao and Suresh (see \cite{RaoS} and references therein) have claimed a number of generalized summation theorems for terminating $_{3}F_{2}(1)$ series by transliterating $3-j$ relationships from quantum mechanics. As an example, Eq. (16) of \cite{RaoS} reads

\begin{equation}
\mapleinline{inert}{2d}{hypergeom([m-n, m-n, 1-2*n],[m-n+2, m-n+2],1) =
2*(-1)^(n-m)*m*GAMMA(2*n)*GAMMA(m-n+2)^2*GAMMA(n-m+1)/GAMMA(m+n+1);}{%
\maplemultiline{
\mathrm{_{3}F_{2}}([m - n, \,m - n, \,1 - 2\,n], \,[m - n + 2, \,
m - n + 2], \,1) \\ 
{\displaystyle = \frac {2\,(-1)^{(n - m)}\,m\,\Gamma (2\,n)\,\Gamma
 (m - n + 2)^{2}\,\Gamma (n - m + 1)}{\Gamma (m + n + 1)}}  }
}
\label{RaoSEq16}
\end{equation}

\vspace*{.5 cm}

Considered as a general result, both sides of (\ref{RaoSEq16}) are infinite for $n>m+1$; the right-hand side is finite only for $n=m$ and $n=m+1$, whereas the left-hand side is finite for $-3 < n \leq m+1$, well-defined (and finite) for $n<m$, non-terminating for $n < 0$ and infinite for $n<-2$. The left-hand side of this result can be identified with several cases established in the database. One possibility is a straightforward, ordered, substitutive mapping \newline
\begin{equation}
 \{a = -n+m,\, b = -n+m, \,c = m-n+2,\, n \rightarrow 2n \} 
\label{map88} 
\end{equation}
\emptyline
applied to case 88 of appendices A and B. The result is

\begin{equation}
\mapleinline{inert}{2d}{(m-n)*(-n+m+1)*(-2*n+(m-n)^2-(m-n)*(m-n+2))*GAMMA(m-n+2)*GAMMA(2*n)*S
um(GAMMA(-n+m+k)*GAMMA(2+k)/GAMMA(k+1)/GAMMA(m-n+2+k),k = 0 ..
2*n)/GAMMA(m+n+1)-1/2*(m-n)*(-n+m+1)*GAMMA(2*n+3)*GAMMA(m-n+2)/n/GAMMA
(m+n+2);}{%
\maplemultiline{
\mathrm{_{3}F_{2}}([m - n, \,m - n, \,1 - 2\,n], \,[m - n + 2, \,
m - n + 2], \,1)= \\
(m - n)\,( - n + m + 1)\,( - 2\,n + (m - n)^{2} - (m - n)\,(m - n
 + 2))\,\Gamma (m - n + 2)\,\Gamma (2\,n) \\
 \left(  \! {\displaystyle \sum _{k=0}^{2\,n}} \,{\displaystyle 
\frac {\Gamma ( - n + m + k)\,\Gamma (2 + k)}{\Gamma (k + 1)\,
\Gamma (m - n + 2 + k)}}  \!  \right) /\Gamma (m + n + 1) 
\mbox{} - {\displaystyle \frac {1}{2}} \,{\displaystyle \frac {(m
 - n)\,( - n + m + 1)\,\Gamma (2\,n + 3)\,\Gamma (m - n + 2)}{n\,
\Gamma (m + n + 2)}}  }
}
\label{subsEq88}
\end{equation}
\emptyline

The right-hand side of (\ref{subsEq88}) has the interesting property that it corresponds with the (infinite) right-hand side of (\ref{RaoSEq16}) for $n>m+1$, disagrees with the right-hand side of (\ref{RaoSEq16}), by vanishing, for $n=m$ and $n=m+1$, and agrees with the left-hand side of both (\ref{subsEq88}) and (\ref{RaoSEq16}) for $n<m$. The simple explanation for this behaviour has been introduced previously - straightforward, ordered, substitutive mapping between cases in the Appendices will not always suffice to obtain the correct answer in a reduction from a general to a specific case.\newline

Focussing once more on (\ref{map88}) for the case $m=n$, apply the substitution $a=0$ to case 88 of Appendix B, identify the resulting series as follows:
\emptyline

\begin{equation}
\mapleinline{inert}{2d}{Sum(GAMMA(b+k)/GAMMA(k+1), k = 0 .. n) = (n+1)*GAMMA(b+n+1)/(b*GAMMA(2+n))}{\[\displaystyle \sum _{k=0}^{n}{\frac {\Gamma  \left( b+k \right) }{\Gamma  \left( k+1 \right) }}={\frac { \left( n+1 \right) \Gamma  \left( b+n+1 \right) }{b \, \Gamma  \left( 2+n \right) }}\]} \, ,
\end{equation}
\emptyline

then apply the ordered substitutions\footnote{technically the limit $b \rightarrow 0$} $b=0$, $c=2$, $n \rightarrow 2n$ to find the well-known result
\emptyline

\begin{equation}
\mapleinline{inert}{2d}{F32([0, 0, 1-2*n], [2, 2], 1) = 1}{\[\displaystyle {\it _{3}F_{2} } \left( [0,0,1-2\,n],[2,2],1 \right) =1\]}
\end{equation}
\emptyline

in complete agreement with (\ref{RaoSEq16}) for the same case. A similar calculation for the case n=m+1, yields
\begin{equation}
\mapleinline{inert}{2d}{F32([-1, -1, 1-2*n], [1, 1], 1) = 2-2*n}{\[\displaystyle {\it _{3}F_{2}} \left( [-1,-1,1-2\,n],[1,1],1 \right) =2-2\,n\]}
\end{equation}
\emptyline

again in complete accord with the right-hand side of (\ref{RaoSEq16}) for the same case. A similar set of calculations with similar results can be applied to (17) of \cite{RaoS}; however (15) of that same paper confounds me because the left-hand side is real, while the right-hand side of that equation appears to be purely imaginary.

\section{Summary}

A large number of new results for the hypergeometric series $_{3}F_{2}(1)$ have been obtained by systematically exploring the 3-part Thomae/Whipple expressions that relate these series with different sets of parameters. A second 4-part expression was also explored. From this investigation, a database consisting of 469 closed forms for $_{3}F_{2}(1)$ was developed, 22 of which were non-viable, and 259 of which are believed to be new. The most general results are presented in Appendices A and B. Of  those cases that are not categorized as ``new'', 112 were identified as special cases and were isolated in a separate set of appendices. By means of an example, the importance of carefully evaluating a mapping of parameters between a candidate and a member of the database was demonstrated. Those cases with parametric excess equal to a numerical (or parametric) constant were identified, as were those cases that are nearly-poised. No results that are not included in this database have been found in the literature. Consequent to this work, I believe that no new Thomae-generated $_3F_2(...|1)$ identities will be found, unless a completely new starting identity someday comes to light. Maier's work \cite{Maier} was an example of one such new identity, prior to it's inclusion in \cite{Milgram2}. \newline 

\emptyline

\section{Acknowledgements}

I would like to thank Wenchang Chu for timely correspondance, and especially for providing a preprint copy of ref. \cite{Chu11}.

\emptyline
\normalsize
\begin{flushleft}

\newpage

\appendix
\section{Appendix A: List of viable $_3F_2(...|1)$ cases}

This Appendix lists the left-hand sides of all the viable and most general $_3F_2(...|1)$ cases that were found, labelled by a case number and the parametric identification.

\input{AppendixA}

\emptyline
$\dagger $ \,\, these cases are nearly-poised 
\emptyline

\hoffset = -50pt
\section{Appendix B: List of results}

The results in this appendix are the right-hand sides of each of the corresponding cases in Appendix A. The various special notations used are discussed in Section 3.
\emptyline
\footnotesize
\begin{flushleft}
\input{AppendixB1to30}

\input{AppendixB31to60}
\input{AppendixB61to100}
\input{AppendixB101to140}

\input{AppendixB141to150}

\input{AppendixB151to170}

\input{AppendixB171to200}

\input{AppendixB201to220}

\input{AppendixB221to240}
\input{AppendixB241to260}
\input{AppendixB261to290}
\input{AppendixB291to320}
\input{AppendixB321to350}
\input{AppendixB351to380}
\input{AppendixB381to410}
\input{AppendixB411to443}
\input{AppendixB444to469}

\newpage
\section{Appendix C: Traceback}

Each case listed here yields a derivation (traceback) for each of the corresponding cases in appendices A and B. See section 3 for an example.
\emptyline
\input{AppendixC}

\newpage
\section{Appendix D: Special Cases}
\normalsize
This Appendix lists the left-hand sides of viable and special cases that were identified, labelled by the case number together with the parametric identification of each case. The case labels are consistent with the labels used previously and represent some of the missing labels from the lists presented in previous appendices. The corresponding right-hand sides appear in Appendix F. For the mapping between each of these special cases and the more general cases appearing in Appendix A, see the corresponding entry in Appendix E.
\vspace*{1 cm}
\footnotesize
\input{AppendixD}

\emptyline
$\dagger $ \,\, these cases are nearly-poised 
\emptyline

\section{Appendix E: Relationship with previous cases}
\normalsize
This appendix lists special notes attached to particular cases.\newline

\footnotesize
\input{AppendixE}

\section{Appendix F: Results - Special cases}
\normalsize
This appendix lists the right-hand sides of the cases corresponding to those listed in Appendix D.
\emptyline
\footnotesize
\input{AppendixF}
\emptyline

\pagestyle{empty}
\thispagestyle{empty}
\section{Appendix G: Cut and Paste}
\normalsize
This appendix lists the results in packed pseudo-code for all the hypergeometric $_{3}F_{2}$ cases listed in appendices A and D. By judicious use of cut/paste and a minor amount of editing, (page numbers have - hopefully - been removed), the contents of the following can be copied, pasted and employed as input to a computer algebra code to allow searching the database.\newline

\tiny
\emptyline

\input{cut_paste}
\end{flushleft}
\end{flushleft}
\end{flushleft}
\end{maplegroup}
\end{document}

%% file: AppendixA.tex
\begin{maplegroup}
\begin{flushleft}
\mapleresult
\begin{maplelatex}
\mapleinline{inert}{2d}{1, _3F_2([a, b, c],[2*c+n, 1/2+1/2*a+1/2*b+1/2*m],1);}{%
\[
1: \,\mathrm{_3F_2}([a, \,b, \,c], \,[2\,c + n, \,{\displaystyle 
\frac {1}{2}}  + {\displaystyle \frac {a}{2}}  + {\displaystyle 
\frac {b}{2}}  + {\displaystyle \frac {m}{2}} ], \,1)
\]
}
\end{maplelatex}

\begin{maplelatex}
\mapleinline{inert}{2d}{2, _3F_2([a, b, c],[2*c+n, 1/2+1/2*a+1/2*b+m],1);}{%
\[
2: \,\mathrm{_3F_2}([a, \,b, \,c], \,[2\,c + n, \,{\displaystyle 
\frac {1}{2}}  + {\displaystyle \frac {a}{2}}  + {\displaystyle 
\frac {b}{2}}  + m], \,1)
\]
}
\end{maplelatex}

\begin{maplelatex}
\mapleinline{inert}{2d}{3, _3F_2([a, b, c],[2*c+n, 1+1/2*a+1/2*b+m],1);}{%
\[
3: \,\mathrm{_3F_2}([a, \,b, \,c], \,[2\,c + n, \,1 + 
{\displaystyle \frac {a}{2}}  + {\displaystyle \frac {b}{2}}  + m
], \,1)
\]
}
\end{maplelatex}

\begin{maplelatex}
\mapleinline{inert}{2d}{4, _3F_2([a, b, c],[2*c-n, 1/2+1/2*a+1/2*b+1/2*m],1);}{%
\[
4: \,\mathrm{_3F_2}([a, \,b, \,c], \,[2\,c - n, \,{\displaystyle 
\frac {1}{2}}  + {\displaystyle \frac {a}{2}}  + {\displaystyle 
\frac {b}{2}}  + {\displaystyle \frac {m}{2}} ], \,1)
\]
}
\end{maplelatex}

\begin{maplelatex}
\mapleinline{inert}{2d}{5, _3F_2([a, b, c],[2*c-n, 1/2+1/2*a+1/2*b+m],1);}{%
\[
5: \,\mathrm{_3F_2}([a, \,b, \,c], \,[2\,c - n, \,{\displaystyle 
\frac {1}{2}}  + {\displaystyle \frac {a}{2}}  + {\displaystyle 
\frac {b}{2}}  + m], \,1)
\]
}
\end{maplelatex}

\begin{maplelatex}
\mapleinline{inert}{2d}{6, _3F_2([a, b, c],[2*c-n, 1+1/2*a+1/2*b+m],1);}{%
\[
6: \,\mathrm{_3F_2}([a, \,b, \,c], \,[2\,c - n, \,1 + 
{\displaystyle \frac {a}{2}}  + {\displaystyle \frac {b}{2}}  + m
], \,1)
\]
}
\end{maplelatex}

\begin{maplelatex}
\mapleinline{inert}{2d}{7, _3F_2([a, b, c],[2*c+n, 1/2+1/2*a+1/2*b-1/2*m],1);}{%
\[
7: \,\mathrm{_3F_2}([a, \,b, \,c], \,[2\,c + n, \,{\displaystyle 
\frac {1}{2}}  + {\displaystyle \frac {a}{2}}  + {\displaystyle 
\frac {b}{2}}  - {\displaystyle \frac {m}{2}} ], \,1)
\]
}
\end{maplelatex}

\begin{maplelatex}
\mapleinline{inert}{2d}{8, _3F_2([a, b, c],[2*c+n, 1/2+1/2*a+1/2*b-m],1);}{%
\[
8: \,\mathrm{_3F_2}([a, \,b, \,c], \,[2\,c + n, \,{\displaystyle 
\frac {1}{2}}  + {\displaystyle \frac {a}{2}}  + {\displaystyle 
\frac {b}{2}}  - m], \,1)
\]
}
\end{maplelatex}

\begin{maplelatex}
\mapleinline{inert}{2d}{9, _3F_2([a, b, c],[2*c+n, 1/2*a+1/2*b-m],1);}{%
\[
9: \,\mathrm{_3F_2}([a, \,b, \,c], \,[2\,c + n, \,{\displaystyle 
\frac {a}{2}}  + {\displaystyle \frac {b}{2}}  - m], \,1)
\]
}
\end{maplelatex}

\begin{maplelatex}
\mapleinline{inert}{2d}{10, _3F_2([a, b, c],[2*c-n, 1/2+1/2*a+1/2*b-1/2*m],1);}{%
\[
10: \,\mathrm{_3F_2}([a, \,b, \,c], \,[2\,c - n, \,{\displaystyle 
\frac {1}{2}}  + {\displaystyle \frac {a}{2}}  + {\displaystyle 
\frac {b}{2}}  - {\displaystyle \frac {m}{2}} ], \,1)
\]
}
\end{maplelatex}

\begin{maplelatex}
\mapleinline{inert}{2d}{11, _3F_2([a, b, c],[2*c-n, 1/2+1/2*a+1/2*b-m],1);}{%
\[
11: \,\mathrm{_3F_2}([a, \,b, \,c], \,[2\,c - n, \,{\displaystyle 
\frac {1}{2}}  + {\displaystyle \frac {a}{2}}  + {\displaystyle 
\frac {b}{2}}  - m], \,1)
\]
}
\end{maplelatex}

\begin{maplelatex}
\mapleinline{inert}{2d}{12, _3F_2([a, b, c],[2*c-n, 1/2*a+1/2*b-m],1);}{%
\[
12: \,\mathrm{_3F_2}([a, \,b, \,c], \,[2\,c - n, \,{\displaystyle 
\frac {a}{2}}  + {\displaystyle \frac {b}{2}}  - m], \,1)
\]
}
\end{maplelatex}

\begin{maplelatex}
\mapleinline{inert}{2d}{13, _3F_2([a, b, 1-a+m],[c, 1+2*b-c+n],1);}{%
\[
13: \,\mathrm{_3F_2}([a, \,b, \,1 - a + m], \,[c, \,1 + 2\,b - c + 
n], \,1)
\]
}
\end{maplelatex}

\begin{maplelatex}
\mapleinline{inert}{2d}{14, _3F_2([a, b, 1-a+m],[c, 1+2*b-c-n],1);}{%
\[
14: \,\mathrm{_3F_2}([a, \,b, \,1 - a + m], \,[c, \,1 + 2\,b - c - 
n], \,1)
\]
}
\end{maplelatex}

\begin{maplelatex}
\mapleinline{inert}{2d}{15, _3F_2([a, b, 1-a-m],[c, 1+2*b-c+n],1);}{%
\[
15: \,\mathrm{_3F_2}([a, \,b, \,1 - a - m], \,[c, \,1 + 2\,b - c + 
n], \,1)
\]
}
\end{maplelatex}

\begin{maplelatex}
\mapleinline{inert}{2d}{16, _3F_2([a, b, 1-a-m],[c, 1+2*b-c-n],1);}{%
\[
16: \,\mathrm{_3F_2}([a, \,b, \,1 - a - m], \,[c, \,1 + 2\,b - c - 
n], \,1)
\]
}
\end{maplelatex}

\begin{maplelatex}
\mapleinline{inert}{2d}{17, _3F_2([a, b, c],[1+a-c+n, 1+a-b+m],1);}{%
\[
17: \,\mathrm{_3F_2}([a, \,b, \,c], \,[1 + a - c + n, \,1 + a - b
 + m], \,1)
\]
}
\end{maplelatex}

\begin{maplelatex}
\mapleinline{inert}{2d}{18, _3F_2([a, b, c],[1+a-c-n, 1+a-b+m],1);}{%
\[
18: \,\mathrm{_3F_2}([a, \,b, \,c], \,[1 + a - c - n, \,1 + a - b
 + m], \,1)
\]
}
\end{maplelatex}

\begin{maplelatex}
\mapleinline{inert}{2d}{19, _3F_2([a, b, c],[1+a-c-n, 1+a-b-m],1);}{%
\[
19: \,\mathrm{_3F_2}([a, \,b, \,c], \,[1 + a - c - n, \,1 + a - b
 - m], \,1)
\]
}
\end{maplelatex}

\begin{maplelatex}
\mapleinline{inert}{2d}{25, _3F_2([a, 1, b],[n+1, c],1);}{%
\[
25: \,\mathrm{_3F_2}([a, \,1, \,b], \,[n + 1, \,c], \,1)
\]
}
\end{maplelatex}

\begin{maplelatex}
\mapleinline{inert}{2d}{26, _3F_2([a, b, n],[n+1, c],1);}{%
\[
26: \,\mathrm{_3F_2}([a, \,b, \,n], \,[n + 1, \,c], \,1)
\]
}
\end{maplelatex}

\begin{maplelatex}
\mapleinline{inert}{2d}{27, _3F_2([a, b, c],[c+1, b+n],1);}{%
\[
27: \,\mathrm{_3F_2}([a, \,b, \,c], \,[c + 1, \,b + n], \,1)
\]
}
\end{maplelatex}

\begin{maplelatex}
\mapleinline{inert}{2d}{28, _3F_2([1, a, b],[c, -c+b+a+n+1],1);}{%
\[
28: \,\mathrm{_3F_2}([1, \,a, \,b], \,[c, \, - c + b + a + n + 1], 
\,1)
\]
}
\end{maplelatex}

\begin{maplelatex}
\mapleinline{inert}{2d}{29, _3F_2([n, a, b],[c, -c+b+a+n+1],1);}{%
\[
29: \,\mathrm{_3F_2}([n, \,a, \,b], \,[c, \, - c + b + a + n + 1], 
\,1)
\]
}
\end{maplelatex}

\begin{maplelatex}
\mapleinline{inert}{2d}{30, _3F_2([a, b, c],[b+1, -n+b+1],1);}{%
\[
30: \,\mathrm{_3F_2}([a, \,b, \,c], \,[b + 1, \, - n + b + 1], \,1)
\]
}
\end{maplelatex}

\begin{maplelatex}
\mapleinline{inert}{2d}{31, _3F_2([a, b, 1],[c, -n+b+1],1);}{%
\[
31: \,\mathrm{_3F_2}([a, \,b, \,1], \,[c, \, - n + b + 1], \,1)
\]
}
\end{maplelatex}

\begin{maplelatex}
\mapleinline{inert}{2d}{33, _3F_2([a, b, 1-n],[c, b+1],1);}{%
\[
33: \,\mathrm{_3F_2}([a, \,b, \,1 - n], \,[c, \,b + 1], \,1)
\]
}
\end{maplelatex}

\begin{maplelatex}
\mapleinline{inert}{2d}{34, _3F_2([a, b, c],[a-n+1, c+b+n],1);}{%
\[
34: \,\mathrm{_3F_2}([a, \,b, \,c], \,[a - n + 1, \,c + b + n], \,1
)
\]
}
\end{maplelatex}

\begin{maplelatex}
\mapleinline{inert}{2d}{35, _3F_2([a, a-n, b],[a-n+1, c],1);}{%
\[
35: \,\mathrm{_3F_2}([a, \,a - n, \,b], \,[a - n + 1, \,c], \,1)
\]
}
\end{maplelatex}

\begin{maplelatex}
\mapleinline{inert}{2d}{36, _3F_2([1, 1-n, a],[b, c],1);}{%
\[
36: \,\mathrm{_3F_2}([1, \,1 - n, \,a], \,[b, \,c], \,1)
\]
}
\end{maplelatex}

\begin{maplelatex}
\mapleinline{inert}{2d}{37, _3F_2([a, b, n+1/2*a-b-3/2],[1/2*a+b+1/2, n+a-b-1],1);}{%
\[
37: \,\mathrm{_3F_2}([a, \,b, \,n + {\displaystyle \frac {a}{2}} 
 - b - {\displaystyle \frac {3}{2}} ], \,[{\displaystyle \frac {a
}{2}}  + b + {\displaystyle \frac {1}{2}} , \,n + a - b - 1], \,1
)
\]
}
\end{maplelatex}

\begin{maplelatex}
\mapleinline{inert}{2d}{38, _3F_2([a, b, n-b-1],[n-b-2+2*a, b+1],1);}{%
\[
38: \,\mathrm{_3F_2}([a, \,b, \,n - b - 1], \,[n - b - 2 + 2\,a, \,
b + 1], \,1)
\]
}
\end{maplelatex}

\begin{maplelatex}
\mapleinline{inert}{2d}{39, _3F_2([1, a, b],[n-b, n-a],1);}{%
\[
39: \,\mathrm{_3F_2}([1, \,a, \,b], \,[n - b, \,n - a], \,1)
\]
}
\end{maplelatex}

\begin{maplelatex}
\mapleinline{inert}{2d}{63, _3F_2([a, b, 1/2*a+1/2],[1+a, n-1+1/2*a+1/2*b],1);}{%
\[
63: \,\mathrm{_3F_2}([a, \,b, \,{\displaystyle \frac {a}{2}}  + 
{\displaystyle \frac {1}{2}} ], \,[1 + a, \,n - 1 + 
{\displaystyle \frac {a}{2}}  + {\displaystyle \frac {b}{2}} ], 
\,1)
\]
}
\end{maplelatex}

\begin{maplelatex}
\mapleinline{inert}{2d}{64, _3F_2([a, b, 1],[2*b, n-1/2+1/2*a],1);}{%
\[
64: \,\mathrm{_3F_2}([a, \,b, \,1], \,[2\,b, \,n - {\displaystyle 
\frac {1}{2}}  + {\displaystyle \frac {a}{2}} ], \,1)
\]
}
\end{maplelatex}

\begin{maplelatex}
\mapleinline{inert}{2d}{65, _3F_2([a, n-3/2+1/2*a, b],[2*n-2+a-b, n-1/2+1/2*a],1);}{%
\[
65: \,\mathrm{_3F_2}([a, \,n - {\displaystyle \frac {3}{2}}  + 
{\displaystyle \frac {a}{2}} , \,b], \,[2\,n - 2 + a - b, \,n - 
{\displaystyle \frac {1}{2}}  + {\displaystyle \frac {a}{2}} ], 
\,1)
\]
}
\end{maplelatex}

\begin{maplelatex}
\mapleinline{inert}{2d}{80, _3F_2([a, b, c],[1/2*c+1/2*b+n+m-1/2, 2*a-2*n+1],1);}{%
\[
80: \,\mathrm{_3F_2}([a, \,b, \,c], \,[{\displaystyle \frac {c}{2}
}  + {\displaystyle \frac {b}{2}}  + n + m - {\displaystyle 
\frac {1}{2}} , \,2\,a - 2\,n + 1], \,1)
\]
}
\end{maplelatex}

\begin{maplelatex}
\mapleinline{inert}{2d}{81, _3F_2([a, b, c],[1/2*c+1/2*b+n+m, 2*a-2*n+1],1);}{%
\[
81: \,\mathrm{_3F_2}([a, \,b, \,c], \,[{\displaystyle \frac {c}{2}
}  + {\displaystyle \frac {b}{2}}  + n + m, \,2\,a - 2\,n + 1], 
\,1)
\]
}
\end{maplelatex}

\begin{maplelatex}
\mapleinline{inert}{2d}{82, _3F_2([a, b, c],[1/2*c+1/2*b+n+m, 2*a-2*n],1);}{%
\[
82: \,\mathrm{_3F_2}([a, \,b, \,c], \,[{\displaystyle \frac {c}{2}
}  + {\displaystyle \frac {b}{2}}  + n + m, \,2\,a - 2\,n], \,1)
\]
}
\end{maplelatex}

\begin{maplelatex}
\mapleinline{inert}{2d}{84, _3F_2([a, b, c],[2+a, a-n+1],1);}{%
\[
84: \,\mathrm{_3F_2}([a, \,b, \,c], \,[2 + a, \,a - n + 1], \,1)
\]
}
\end{maplelatex}

\begin{maplelatex}
\mapleinline{inert}{2d}{85, _3F_2([a, b, 2],[c, a-n+1],1);}{%
\[
85: \,\mathrm{_3F_2}([a, \,b, \,2], \,[c, \,a - n + 1], \,1)
\]
}
\end{maplelatex}

\begin{maplelatex}
\mapleinline{inert}{2d}{86, _3F_2([a, b, c],[n+1+c+b, a-n+1],1);}{%
\[
86: \,\mathrm{_3F_2}([a, \,b, \,c], \,[n + 1 + c + b, \,a - n + 1]
, \,1)
\]
}
\end{maplelatex}

\begin{maplelatex}
\mapleinline{inert}{2d}{87, _3F_2([a, a-n-1, b],[c, a-n+1],1);}{%
\[
87: \,\mathrm{_3F_2}([a, \,a - n - 1, \,b], \,[c, \,a - n + 1], \,1
)
\]
}
\end{maplelatex}

\begin{maplelatex}
\mapleinline{inert}{2d}{88, _3F_2([a, b, 1-n],[c, 2+b],1);}{%
\[
88: \,\mathrm{_3F_2}([a, \,b, \,1 - n], \,[c, \,2 + b], \,1)
\]
}
\end{maplelatex}

\begin{maplelatex}
\mapleinline{inert}{2d}{89, _3F_2([2, a, 1-n],[b, c],1);}{%
\[
89: \,\mathrm{_3F_2}([2, \,a, \,1 - n], \,[b, \,c], \,1)
\]
}
\end{maplelatex}

\begin{maplelatex}
\mapleinline{inert}{2d}{91, _3F_2([a, b, c],[2+a, (a*c+a*b-a-b*c+b+c-1)/a],1);}{%
\[
91: \,\mathrm{_3F_2}([a, \,b, \,c], \,[2 + a, \,{\displaystyle 
\frac {a\,c + a\,b - a - b\,c + b + c - 1}{a}} ], \,1)
\]
}
\end{maplelatex}

\begin{maplelatex}
\mapleinline{inert}{2d}{92, _3F_2([2, a, b],[c, (-3+a*b-b-a+2*c)/(c-2)],1);}{%
\[
92: \,\mathrm{_3F_2}([2, \,a, \,b], \,[c, \,{\displaystyle \frac {
 - 3 + a\,b - b - a + 2\,c}{c - 2}} ], \,1)
\]
}
\end{maplelatex}

\begin{maplelatex}
\mapleinline{inert}{2d}{104, _3F_2([2*a, a+1/3-1/2*n, a-1/3-1/2*n],[1/2+2*a-n,
3*a-1/2*n],1);}{%
\[
104: \,\mathrm{_3F_2}([2\,a, \,a + {\displaystyle \frac {1}{3}}  - 
{\displaystyle \frac {n}{2}} , \,a - {\displaystyle \frac {1}{3}
}  - {\displaystyle \frac {n}{2}} ], \,[{\displaystyle \frac {1}{
2}}  + 2\,a - n, \,3\,a - {\displaystyle \frac {n}{2}} ], \,1)
\]
}
\end{maplelatex}

\begin{maplelatex}
\mapleinline{inert}{2d}{105, _3F_2([2*a, a-1/2*n+5/6, a+1/6-1/2*n],[1/2+2*a-n,
3*a+1/2-1/2*n],1);}{%
\[
105: \,\mathrm{_3F_2}([2\,a, \,a - {\displaystyle \frac {n}{2}}  + 
{\displaystyle \frac {5}{6}} , \,a + {\displaystyle \frac {1}{6}
}  - {\displaystyle \frac {n}{2}} ], \,[{\displaystyle \frac {1}{
2}}  + 2\,a - n, \,3\,a + {\displaystyle \frac {1}{2}}  - 
{\displaystyle \frac {n}{2}} ], \,1)
\]
}
\end{maplelatex}

\begin{maplelatex}
\mapleinline{inert}{2d}{106, _3F_2([a+1/3-1/2*n, a-1/2*n+5/6, -n+1/2],[1/2+2*a-n,
2*a+5/6-n],1);}{%
\[
106: \,\mathrm{_3F_2}([a + {\displaystyle \frac {1}{3}}  - 
{\displaystyle \frac {n}{2}} , \,a - {\displaystyle \frac {n}{2}
}  + {\displaystyle \frac {5}{6}} , \, - n + {\displaystyle 
\frac {1}{2}} ], \,[{\displaystyle \frac {1}{2}}  + 2\,a - n, \,2
\,a + {\displaystyle \frac {5}{6}}  - n], \,1)
\]
}
\end{maplelatex}

\begin{maplelatex}
\mapleinline{inert}{2d}{107, _3F_2([a-1/3-1/2*n, a+1/6-1/2*n, -n+1/2],[1/2+2*a-n,
2*a+1/6-n],1);}{%
\[
107: \,\mathrm{_3F_2}([a - {\displaystyle \frac {1}{3}}  - 
{\displaystyle \frac {n}{2}} , \,a + {\displaystyle \frac {1}{6}
}  - {\displaystyle \frac {n}{2}} , \, - n + {\displaystyle 
\frac {1}{2}} ], \,[{\displaystyle \frac {1}{2}}  + 2\,a - n, \,2
\,a + {\displaystyle \frac {1}{6}}  - n], \,1)
\]
}
\end{maplelatex}

\begin{maplelatex}
\mapleinline{inert}{2d}{108, _3F_2([2*a, 2*a+1/3, 2*a-1/3],[3*a-1/2*n, 3*a+1/2-1/2*n],1);}{%
\[
108: \,\mathrm{_3F_2}([2\,a, \,2\,a + {\displaystyle \frac {1}{3}} 
, \,2\,a - {\displaystyle \frac {1}{3}} ], \,[3\,a - 
{\displaystyle \frac {n}{2}} , \,3\,a + {\displaystyle \frac {1}{
2}}  - {\displaystyle \frac {n}{2}} ], \,1)
\]
}
\end{maplelatex}

\begin{maplelatex}
\mapleinline{inert}{2d}{109, _3F_2([a+1/3-1/2*n, 2*a+1/3, a-1/2*n],[3*a-1/2*n,
2*a+5/6-n],1);}{%
\[
109: \,\mathrm{_3F_2}([a + {\displaystyle \frac {1}{3}}  - 
{\displaystyle \frac {n}{2}} , \,2\,a + {\displaystyle \frac {1}{
3}} , \,a - {\displaystyle \frac {n}{2}} ], \,[3\,a - 
{\displaystyle \frac {n}{2}} , \,2\,a + {\displaystyle \frac {5}{
6}}  - n], \,1)
\]
}
\end{maplelatex}

\begin{maplelatex}
\mapleinline{inert}{2d}{110, _3F_2([a-1/3-1/2*n, 2*a-1/3, a-1/2*n],[3*a-1/2*n,
2*a+1/6-n],1);}{%
\[
110: \,\mathrm{_3F_2}([a - {\displaystyle \frac {1}{3}}  - 
{\displaystyle \frac {n}{2}} , \,2\,a - {\displaystyle \frac {1}{
3}} , \,a - {\displaystyle \frac {n}{2}} ], \,[3\,a - 
{\displaystyle \frac {n}{2}} , \,2\,a + {\displaystyle \frac {1}{
6}}  - n], \,1)
\]
}
\end{maplelatex}

\begin{maplelatex}
\mapleinline{inert}{2d}{111, _3F_2([a-1/2*n+5/6, 2*a+1/3, 1/2+a-1/2*n],[3*a+1/2-1/2*n,
2*a+5/6-n],1);}{%
\[
111: \,\mathrm{_3F_2}([a - {\displaystyle \frac {n}{2}}  + 
{\displaystyle \frac {5}{6}} , \,2\,a + {\displaystyle \frac {1}{
3}} , \,{\displaystyle \frac {1}{2}}  + a - {\displaystyle 
\frac {n}{2}} ], \,[3\,a + {\displaystyle \frac {1}{2}}  - 
{\displaystyle \frac {n}{2}} , \,2\,a + {\displaystyle \frac {5}{
6}}  - n], \,1)
\]
}
\end{maplelatex}

\begin{maplelatex}
\mapleinline{inert}{2d}{112, _3F_2([a+1/6-1/2*n, 2*a-1/3, 1/2+a-1/2*n],[3*a+1/2-1/2*n,
2*a+1/6-n],1);}{%
\[
112: \,\mathrm{_3F_2}([a + {\displaystyle \frac {1}{6}}  - 
{\displaystyle \frac {n}{2}} , \,2\,a - {\displaystyle \frac {1}{
3}} , \,{\displaystyle \frac {1}{2}}  + a - {\displaystyle 
\frac {n}{2}} ], \,[3\,a + {\displaystyle \frac {1}{2}}  - 
{\displaystyle \frac {n}{2}} , \,2\,a + {\displaystyle \frac {1}{
6}}  - n], \,1)
\]
}
\end{maplelatex}

\begin{maplelatex}
\mapleinline{inert}{2d}{113, _3F_2([-n+1/2, a-1/2*n, 1/2+a-1/2*n],[2*a+5/6-n, 2*a+1/6-n],1);}{%
\[
113: \,\mathrm{_3F_2}([ - n + {\displaystyle \frac {1}{2}} , \,a - 
{\displaystyle \frac {n}{2}} , \,{\displaystyle \frac {1}{2}}  + 
a - {\displaystyle \frac {n}{2}} ], \,[2\,a + {\displaystyle 
\frac {5}{6}}  - n, \,2\,a + {\displaystyle \frac {1}{6}}  - n], 
\,1)
\]
}
\end{maplelatex}

\begin{maplelatex}
\mapleinline{inert}{2d}{114, _3F_2([-1/2+n, 2*a+n-2/3, -1/6-a],[n-1/2+a, -1/6+a+n],1);}{%
\[
114: \,\mathrm{_3F_2}([ - {\displaystyle \frac {1}{2}}  + n, \,2\,a
 + n - {\displaystyle \frac {2}{3}} , \, - {\displaystyle \frac {
1}{6}}  - a], \,[n - {\displaystyle \frac {1}{2}}  + a, \, - 
{\displaystyle \frac {1}{6}}  + a + n], \,1)
\]
}
\end{maplelatex}

\begin{maplelatex}
\mapleinline{inert}{2d}{115, _3F_2([-1/2+n, 2*a-1/3+n, -a+1/6],[n-1/2+a, 1/6+a+n],1);}{%
\[
115: \,\mathrm{_3F_2}([ - {\displaystyle \frac {1}{2}}  + n, \,2\,a
 - {\displaystyle \frac {1}{3}}  + n, \, - a + {\displaystyle 
\frac {1}{6}} ], \,[n - {\displaystyle \frac {1}{2}}  + a, \,
{\displaystyle \frac {1}{6}}  + a + n], \,1)
\]
}
\end{maplelatex}

\begin{maplelatex}
\mapleinline{inert}{2d}{116, _3F_2([2*a+n-2/3, 2*a-1/3+n, a],[n-1/2+a, 3*a+n],1);}{%
\[
116: \,\mathrm{_3F_2}([2\,a + n - {\displaystyle \frac {2}{3}} , \,
2\,a - {\displaystyle \frac {1}{3}}  + n, \,a], \,[n - 
{\displaystyle \frac {1}{2}}  + a, \,3\,a + n], \,1)
\]
}
\end{maplelatex}

\begin{maplelatex}
\mapleinline{inert}{2d}{117, _3F_2([-1/6-a, -a+1/6, a],[n-1/2+a, 1/2],1);}{%
\[
117: \,\mathrm{_3F_2}([ - {\displaystyle \frac {1}{6}}  - a, \, - a
 + {\displaystyle \frac {1}{6}} , \,a], \,[n - {\displaystyle 
\frac {1}{2}}  + a, \,{\displaystyle \frac {1}{2}} ], \,1)
\]
}
\end{maplelatex}

\begin{maplelatex}
\mapleinline{inert}{2d}{118, _3F_2([-1/2+n, 2*a+n, 1/2-a],[-1/6+a+n, 1/6+a+n],1);}{%
\[
118: \,\mathrm{_3F_2}([ - {\displaystyle \frac {1}{2}}  + n, \,2\,a
 + n, \,{\displaystyle \frac {1}{2}}  - a], \,[ - {\displaystyle 
\frac {1}{6}}  + a + n, \,{\displaystyle \frac {1}{6}}  + a + n]
, \,1)
\]
}
\end{maplelatex}

\begin{maplelatex}
\mapleinline{inert}{2d}{119, _3F_2([2*a+n-2/3, 2*a+n, a+1/3],[-1/6+a+n, 3*a+n],1);}{%
\[
119: \,\mathrm{_3F_2}([2\,a + n - {\displaystyle \frac {2}{3}} , \,
2\,a + n, \,a + {\displaystyle \frac {1}{3}} ], \,[ - 
{\displaystyle \frac {1}{6}}  + a + n, \,3\,a + n], \,1)
\]
}
\end{maplelatex}

\begin{maplelatex}
\mapleinline{inert}{2d}{120, _3F_2([-1/6-a, 1/2-a, a+1/3],[-1/6+a+n, 1/2],1);}{%
\[
120: \,\mathrm{_3F_2}([ - {\displaystyle \frac {1}{6}}  - a, \,
{\displaystyle \frac {1}{2}}  - a, \,a + {\displaystyle \frac {1
}{3}} ], \,[ - {\displaystyle \frac {1}{6}}  + a + n, \,
{\displaystyle \frac {1}{2}} ], \,1)
\]
}
\end{maplelatex}

\begin{maplelatex}
\mapleinline{inert}{2d}{121, _3F_2([2*a-1/3+n, 2*a+n, a+2/3],[1/6+a+n, 3*a+n],1);}{%
\[
121: \,\mathrm{_3F_2}([2\,a - {\displaystyle \frac {1}{3}}  + n, \,
2\,a + n, \,a + {\displaystyle \frac {2}{3}} ], \,[
{\displaystyle \frac {1}{6}}  + a + n, \,3\,a + n], \,1)
\]
}
\end{maplelatex}

\begin{maplelatex}
\mapleinline{inert}{2d}{122, _3F_2([-a+1/6, 1/2-a, a+2/3],[1/6+a+n, 1/2],1);}{%
\[
122: \,\mathrm{_3F_2}([ - a + {\displaystyle \frac {1}{6}} , \,
{\displaystyle \frac {1}{2}}  - a, \,a + {\displaystyle \frac {2
}{3}} ], \,[{\displaystyle \frac {1}{6}}  + a + n, \,
{\displaystyle \frac {1}{2}} ], \,1)
\]
}
\end{maplelatex}

\begin{maplelatex}
\mapleinline{inert}{2d}{123, _3F_2([a, a+1/3, a+2/3],[3*a+n, 1/2],1);}{%
\[
123: \,\mathrm{_3F_2}([a, \,a + {\displaystyle \frac {1}{3}} , \,a
 + {\displaystyle \frac {2}{3}} ], \,[3\,a + n, \,{\displaystyle 
\frac {1}{2}} ], \,1)
\]
}
\end{maplelatex}

\begin{maplelatex}
\mapleinline{inert}{2d}{124, _3F_2([1/2, 2*a+1/3-n, -1/6+n-a],[a+1/2, a+5/6],1);}{%
\[
124: \,\mathrm{_3F_2}([{\displaystyle \frac {1}{2}} , \,2\,a + 
{\displaystyle \frac {1}{3}}  - n, \, - {\displaystyle \frac {1}{
6}}  + n - a], \,[a + {\displaystyle \frac {1}{2}} , \,a + 
{\displaystyle \frac {5}{6}} ], \,1)
\]
}
\end{maplelatex}

\begin{maplelatex}
\mapleinline{inert}{2d}{125, _3F_2([1/2, 2*a+2/3-n, 1/6+n-a],[a+1/2, 7/6+a],1);}{%
\[
125: \,\mathrm{_3F_2}([{\displaystyle \frac {1}{2}} , \,2\,a + 
{\displaystyle \frac {2}{3}}  - n, \,{\displaystyle \frac {1}{6}
}  + n - a], \,[a + {\displaystyle \frac {1}{2}} , \,
{\displaystyle \frac {7}{6}}  + a], \,1)
\]
}
\end{maplelatex}

\begin{maplelatex}
\mapleinline{inert}{2d}{126, _3F_2([2*a+1/3-n, 2*a+2/3-n, a],[a+1/2, 3*a+1-n],1);}{%
\[
126: \,\mathrm{_3F_2}([2\,a + {\displaystyle \frac {1}{3}}  - n, \,
2\,a + {\displaystyle \frac {2}{3}}  - n, \,a], \,[a + 
{\displaystyle \frac {1}{2}} , \,3\,a + 1 - n], \,1)
\]
}
\end{maplelatex}

\begin{maplelatex}
\mapleinline{inert}{2d}{127, _3F_2([-1/6+n-a, 1/6+n-a, a],[a+1/2, 1/2+n],1);}{%
\[
127: \,\mathrm{_3F_2}([ - {\displaystyle \frac {1}{6}}  + n - a, \,
{\displaystyle \frac {1}{6}}  + n - a, \,a], \,[a + 
{\displaystyle \frac {1}{2}} , \,{\displaystyle \frac {1}{2}}  + 
n], \,1)
\]
}
\end{maplelatex}

\begin{maplelatex}
\mapleinline{inert}{2d}{128, _3F_2([1/2, 2*a+1-n, 1/2+n-a],[a+5/6, 7/6+a],1);}{%
\[
128: \,\mathrm{_3F_2}([{\displaystyle \frac {1}{2}} , \,2\,a + 1 - 
n, \,{\displaystyle \frac {1}{2}}  + n - a], \,[a + 
{\displaystyle \frac {5}{6}} , \,{\displaystyle \frac {7}{6}}  + 
a], \,1)
\]
}
\end{maplelatex}

\begin{maplelatex}
\mapleinline{inert}{2d}{129, _3F_2([2*a+1/3-n, 2*a+1-n, a+1/3],[a+5/6, 3*a+1-n],1);}{%
\[
129: \,\mathrm{_3F_2}([2\,a + {\displaystyle \frac {1}{3}}  - n, \,
2\,a + 1 - n, \,a + {\displaystyle \frac {1}{3}} ], \,[a + 
{\displaystyle \frac {5}{6}} , \,3\,a + 1 - n], \,1)
\]
}
\end{maplelatex}

\begin{maplelatex}
\mapleinline{inert}{2d}{130, _3F_2([-1/6+n-a, 1/2+n-a, a+1/3],[a+5/6, 1/2+n],1);}{%
\[
130: \,\mathrm{_3F_2}([ - {\displaystyle \frac {1}{6}}  + n - a, \,
{\displaystyle \frac {1}{2}}  + n - a, \,a + {\displaystyle 
\frac {1}{3}} ], \,[a + {\displaystyle \frac {5}{6}} , \,
{\displaystyle \frac {1}{2}}  + n], \,1)
\]
}
\end{maplelatex}

\begin{maplelatex}
\mapleinline{inert}{2d}{131, _3F_2([2*a+2/3-n, 2*a+1-n, a+2/3],[7/6+a, 3*a+1-n],1);}{%
\[
131: \,\mathrm{_3F_2}([2\,a + {\displaystyle \frac {2}{3}}  - n, \,
2\,a + 1 - n, \,a + {\displaystyle \frac {2}{3}} ], \,[
{\displaystyle \frac {7}{6}}  + a, \,3\,a + 1 - n], \,1)
\]
}
\end{maplelatex}

\begin{maplelatex}
\mapleinline{inert}{2d}{132, _3F_2([1/6+n-a, 1/2+n-a, a+2/3],[7/6+a, 1/2+n],1);}{%
\[
132: \,\mathrm{_3F_2}([{\displaystyle \frac {1}{6}}  + n - a, \,
{\displaystyle \frac {1}{2}}  + n - a, \,a + {\displaystyle 
\frac {2}{3}} ], \,[{\displaystyle \frac {7}{6}}  + a, \,
{\displaystyle \frac {1}{2}}  + n], \,1)
\]
}
\end{maplelatex}

\begin{maplelatex}
\mapleinline{inert}{2d}{133, _3F_2([a, a+1/3, a+2/3],[3*a+1-n, 1/2+n],1);}{%
\[
133: \,\mathrm{_3F_2}([a, \,a + {\displaystyle \frac {1}{3}} , \,a
 + {\displaystyle \frac {2}{3}} ], \,[3\,a + 1 - n, \,
{\displaystyle \frac {1}{2}}  + n], \,1)
\]
}
\end{maplelatex}

\begin{maplelatex}
\mapleinline{inert}{2d}{134, _3F_2([1, 2-n+a, b-n+2],[-a+2, -b+2],1);}{%
\[
134: \,\mathrm{_3F_2}([1, \,2 - n + a, \,b - n + 2], \,[ - a + 2, 
\, - b + 2], \,1)
\]
}
\end{maplelatex}

\begin{maplelatex}
\mapleinline{inert}{2d}{139, _3F_2([1, 2-2*b, 5/2-n-1/2*a],[-a+2, -b+2],1);}{%
\[
139: \,\mathrm{_3F_2}([1, \,2 - 2\,b, \,{\displaystyle \frac {5}{2}
}  - n - {\displaystyle \frac {a}{2}} ], \,[ - a + 2, \, - b + 2]
, \,1)
\]
}
\end{maplelatex}

\begin{maplelatex}
\mapleinline{inert}{2d}{148, _3F_2([a, b, 1/2*a+1/2],[1+a, 1/2*a+3/2+1/2*b-1/2*n],1);}{%
\[
148: \,\mathrm{_3F_2}([a, \,b, \,{\displaystyle \frac {a}{2}}  + 
{\displaystyle \frac {1}{2}} ], \,[1 + a, \,{\displaystyle 
\frac {a}{2}}  + {\displaystyle \frac {3}{2}}  + {\displaystyle 
\frac {b}{2}}  - {\displaystyle \frac {n}{2}} ], \,1)
\]
}
\end{maplelatex}

\begin{maplelatex}
\mapleinline{inert}{2d}{149, _3F_2([a, b, 1],[2*b, 1/2*a-1/2*n+2],1);}{%
\[
149: \,\mathrm{_3F_2}([a, \,b, \,1], \,[2\,b, \,{\displaystyle 
\frac {a}{2}}  - {\displaystyle \frac {n}{2}}  + 2], \,1)
\]
}
\end{maplelatex}

\begin{maplelatex}
\mapleinline{inert}{2d}{150, _3F_2([a, b, -b+5/2+1/2*a-n],[1/2*a+b+1/2, 3+a-n-b],1);}{%
\[
150: \,\mathrm{_3F_2}([a, \,b, \, - b + {\displaystyle \frac {5}{2}
}  + {\displaystyle \frac {a}{2}}  - n], \,[{\displaystyle 
\frac {a}{2}}  + b + {\displaystyle \frac {1}{2}} , \,3 + a - n
 - b], \,1)
\]
}
\end{maplelatex}

\begin{maplelatex}
\mapleinline{inert}{2d}{151, _3F_2([a, 1/2*a-1/2*n+1, b],[3+a-n-b, 1/2*a-1/2*n+2],1);}{%
\[
151: \,\mathrm{_3F_2}([a, \,{\displaystyle \frac {a}{2}}  - 
{\displaystyle \frac {n}{2}}  + 1, \,b], \,[3 + a - n - b, \,
{\displaystyle \frac {a}{2}}  - {\displaystyle \frac {n}{2}}  + 2
], \,1)
\]
}
\end{maplelatex}

\begin{maplelatex}
\mapleinline{inert}{2d}{152, _3F_2([a, b, -b+3-n],[-n-b+2+2*a, b+1],1);}{%
\[
152: \,\mathrm{_3F_2}([a, \,b, \, - b + 3 - n], \,[ - n - b + 2 + 2
\,a, \,b + 1], \,1)
\]
}
\end{maplelatex}

\begin{maplelatex}
\mapleinline{inert}{2d}{160, _3F_2([a, 2*n-2-a, b],[1+a, 2*b-a],1);}{%
\[
160: \,\mathrm{_3F_2}([a, \,2\,n - 2 - a, \,b], \,[1 + a, \,2\,b - 
a], \,1)
\]
}
\end{maplelatex}

\begin{maplelatex}
\mapleinline{inert}{2d}{161, _3F_2([a, b, 2*a-2*n+3],[1+a, 2*a-2*n+4-b],1);}{%
\[
\dagger 161: \,\mathrm{_3F_2}([a, \,b, \,2\,a - 2\,n + 3], \,[1 + a, \,2\,a
 - 2\,n + 4 - b], \,1)
\]
}
\end{maplelatex}

\begin{maplelatex}
\mapleinline{inert}{2d}{162, _3F_2([a, b, 1],[2*n-1-a, -b+2],1);}{%
\[
\dagger 162: \,\mathrm{_3F_2}([a, \,b, \,1], \,[2\,n - 1 - a, \, - b + 2], 
\,1)
\]
}
\end{maplelatex}

\begin{maplelatex}
\mapleinline{inert}{2d}{163, _3F_2([a, b, 1],[n-1/2+1/2*b, 2*a-2*n+3],1);}{%
\[
163: \,\mathrm{_3F_2}([a, \,b, \,1], \,[n - {\displaystyle \frac {1
}{2}}  + {\displaystyle \frac {b}{2}} , \,2\,a - 2\,n + 3], \,1)
\]
}
\end{maplelatex}

\begin{maplelatex}
\mapleinline{inert}{2d}{164, _3F_2([a, b, 2*a-2*n+2*b+2],[2*b+a, 2*a-2*n+b+3],1);}{%
\[
\dagger 164: \,\mathrm{_3F_2}([a, \,b, \,2\,a - 2\,n + 2\,b + 2], \,[2\,b
 + a, \,2\,a - 2\,n + b + 3], \,1)
\]
}
\end{maplelatex}

\begin{maplelatex}
\mapleinline{inert}{2d}{165, _3F_2([a, b, 2*b],[-a+2*n-2+2*b, b+1],1);}{%
\[
\dagger 165: \,\mathrm{_3F_2}([a, \,b, \,2\,b], \,[ - a + 2\,n - 2 + 2\,b, 
\,b + 1], \,1)
\]
}
\end{maplelatex}

\begin{maplelatex}
\mapleinline{inert}{2d}{166, _3F_2([a, 2-2*n+2*a, b],[a+1/2*b, 2*a-2*n+3],1);}{%
\[
166: \,\mathrm{_3F_2}([a, \,2 - 2\,n + 2\,a, \,b], \,[a + 
{\displaystyle \frac {b}{2}} , \,2\,a - 2\,n + 3], \,1)
\]
}
\end{maplelatex}

\begin{maplelatex}
\mapleinline{inert}{2d}{167, _3F_2([a, -a+1, b],[a+2*b-4+2*n, -a+2],1);}{%
\[
167: \,\mathrm{_3F_2}([a, \, - a + 1, \,b], \,[a + 2\,b - 4 + 2\,n
, \, - a + 2], \,1)
\]
}
\end{maplelatex}

\begin{maplelatex}
\mapleinline{inert}{2d}{168, _3F_2([a, b, 2-n+1/2*b],[1/2*a+1/2+1/2*b, b+1],1);}{%
\[
168: \,\mathrm{_3F_2}([a, \,b, \,2 - n + {\displaystyle \frac {b}{2
}} ], \,[{\displaystyle \frac {a}{2}}  + {\displaystyle \frac {1
}{2}}  + {\displaystyle \frac {b}{2}} , \,b + 1], \,1)
\]
}
\end{maplelatex}

\begin{maplelatex}
\mapleinline{inert}{2d}{194, _3F_2([-2*a+2, n-b-1, n-2*b-1],[n-2*b, n-a-b],1);}{%
\[
194: \,\mathrm{_3F_2}([ - 2\,a + 2, \,n - b - 1, \,n - 2\,b - 1], 
\,[n - 2\,b, \,n - a - b], \,1)
\]
}
\end{maplelatex}

\begin{maplelatex}
\mapleinline{inert}{2d}{195, _3F_2([a, 3-2*n+b, 1/2*a-n+3/2],[1+a-b, 5/2+1/2*a-n],1);}{%
\[
195: \,\mathrm{_3F_2}([a, \,3 - 2\,n + b, \,{\displaystyle \frac {a
}{2}}  - n + {\displaystyle \frac {3}{2}} ], \,[1 + a - b, \,
{\displaystyle \frac {5}{2}}  + {\displaystyle \frac {a}{2}}  - n
], \,1)
\]
}
\end{maplelatex}

\begin{maplelatex}
\mapleinline{inert}{2d}{196, _3F_2([1, 2-2*b, -1/2*a+1/2*n],[-b+2, -a+2],1);}{%
\[
196: \,\mathrm{_3F_2}([1, \,2 - 2\,b, \, - {\displaystyle \frac {a
}{2}}  + {\displaystyle \frac {n}{2}} ], \,[ - b + 2, \, - a + 2]
, \,1)
\]
}
\end{maplelatex}

\begin{maplelatex}
\mapleinline{inert}{2d}{197, _3F_2([-2*a+2, -b+3-n, -2*b+3-n],[-n-2*b+4, -n-a+4-b],1);}{%
\[
197: \,\mathrm{_3F_2}([ - 2\,a + 2, \, - b + 3 - n, \, - 2\,b + 3
 - n], \,[ - n - 2\,b + 4, \, - n - a + 4 - b], \,1)
\]
}
\end{maplelatex}

\begin{maplelatex}
\mapleinline{inert}{2d}{198, _3F_2([1, b, 3-2*n+a],[-b+2, -a+2],1);}{%
\[
\dagger 198: \,\mathrm{_3F_2}([1, \,b, \,3 - 2\,n + a], \,[ - b + 2, \, - a
 + 2], \,1)
\]
}
\end{maplelatex}

\begin{maplelatex}
\mapleinline{inert}{2d}{199, _3F_2([1, 2*n-2*a-1, 5/2-n-1/2*b],[-b+2, -a+2],1);}{%
\[
199: \,\mathrm{_3F_2}([1, \,2\,n - 2\,a - 1, \,{\displaystyle 
\frac {5}{2}}  - n - {\displaystyle \frac {b}{2}} ], \,[ - b + 2
, \, - a + 2], \,1)
\]
}
\end{maplelatex}

\begin{maplelatex}
\mapleinline{inert}{2d}{200, _3F_2([b, 2*b, 3-2*n+a],[b+1, 2*b+1-a],1);}{%
\[
\dagger 200: \,\mathrm{_3F_2}([b, \,2\,b, \,3 - 2\,n + a], \,[b + 1, \,2\,b
 + 1 - a], \,1)
\]
}
\end{maplelatex}

\begin{maplelatex}
\mapleinline{inert}{2d}{201, _3F_2([b, -1+a+b, 5-a-b-2*n],[a+b, b-a+1],1);}{%
\[
201: \,\mathrm{_3F_2}([b, \, - 1 + a + b, \,5 - a - b - 2\,n], \,[a
 + b, \,b - a + 1], \,1)
\]
}
\end{maplelatex}

\begin{maplelatex}
\mapleinline{inert}{2d}{210, _3F_2([a, b, 1],[-a+2, n-b],1);}{%
\[
\dagger 210: \,\mathrm{_3F_2}([a, \,b, \,1], \,[ - a + 2, \,n - b], \,1)
\]
}
\end{maplelatex}

\begin{maplelatex}
\mapleinline{inert}{2d}{211, _3F_2([a, -a+1, b],[-a+2, n+a+2*b-3],1);}{%
\[
211: \,\mathrm{_3F_2}([a, \, - a + 1, \,b], \,[ - a + 2, \,n + a + 
2\,b - 3], \,1)
\]
}
\end{maplelatex}

\begin{maplelatex}
\mapleinline{inert}{2d}{212, _3F_2([a, b, 2*b],[b+1, n+2*b-1-a],1);}{%
\[
\dagger 212: \,\mathrm{_3F_2}([a, \,b, \,2\,b], \,[b + 1, \,n + 2\,b - 1 - 
a], \,1)
\]
}
\end{maplelatex}

\begin{maplelatex}
\mapleinline{inert}{2d}{213, _3F_2([1, a, b],[1/2*b+1, n+2*a-2],1);}{%
\[
213: \,\mathrm{_3F_2}([1, \,a, \,b], \,[{\displaystyle \frac {b}{2}
}  + 1, \,n + 2\,a - 2], \,1)
\]
}
\end{maplelatex}

\begin{maplelatex}
\mapleinline{inert}{2d}{214, _3F_2([a, b, -n+2*b+2],[b+1, -n-a+2*b+3],1);}{%
\[
\dagger 214: \,\mathrm{_3F_2}([a, \,b, \, - n + 2\,b + 2], \,[b + 1, \, - n
 - a + 2\,b + 3], \,1)
\]
}
\end{maplelatex}

\begin{maplelatex}
\mapleinline{inert}{2d}{215, _3F_2([a, n-a-1, b],[n-a, a+2*b-n+1],1);}{%
\[
215: \,\mathrm{_3F_2}([a, \,n - a - 1, \,b], \,[n - a, \,a + 2\,b
 - n + 1], \,1)
\]
}
\end{maplelatex}

\begin{maplelatex}
\mapleinline{inert}{2d}{216, _3F_2([1, a, b],[1/2*n+1/2*a, -n+2*b+2],1);}{%
\[
216: \,\mathrm{_3F_2}([1, \,a, \,b], \,[{\displaystyle \frac {n}{2}
}  + {\displaystyle \frac {a}{2}} , \, - n + 2\,b + 2], \,1)
\]
}
\end{maplelatex}

\begin{maplelatex}
\mapleinline{inert}{2d}{217, _3F_2([a, b, -n+2*b+2*a+1],[-n+a+2*b+2, 2*a+b],1);}{%
\[
\dagger 217: \,\mathrm{_3F_2}([a, \,b, \, - n + 2\,b + 2\,a + 1], \,[ - n
 + a + 2\,b + 2, \,2\,a + b], \,1)
\]
}
\end{maplelatex}

\begin{maplelatex}
\mapleinline{inert}{2d}{218, _3F_2([a, b, 2*a+n-3],[1/2*n+1/2*b+a-1, n+2*a-2],1);}{%
\[
218: \,\mathrm{_3F_2}([a, \,b, \,2\,a + n - 3], \,[{\displaystyle 
\frac {n}{2}}  + {\displaystyle \frac {b}{2}}  + a - 1, \,n + 2\,
a - 2], \,1)
\]
}
\end{maplelatex}

\begin{maplelatex}
\mapleinline{inert}{2d}{219, _3F_2([a, 1, b],[2*b, 5/2+1/2*a-n],1);}{%
\[
219: \,\mathrm{_3F_2}([a, \,1, \,b], \,[2\,b, \,{\displaystyle 
\frac {5}{2}}  + {\displaystyle \frac {a}{2}}  - n], \,1)
\]
}
\end{maplelatex}

\begin{maplelatex}
\mapleinline{inert}{2d}{220, _3F_2([a, b, 2*b-1],[2*b, b+1/2*a-n+3/2],1);}{%
\[
220: \,\mathrm{_3F_2}([a, \,b, \,2\,b - 1], \,[2\,b, \,b + 
{\displaystyle \frac {a}{2}}  - n + {\displaystyle \frac {3}{2}} 
], \,1)
\]
}
\end{maplelatex}

\begin{maplelatex}
\mapleinline{inert}{2d}{221, _3F_2([1, a, b],[5-b-2*n, 5-a-2*n],1);}{%
\[
221: \,\mathrm{_3F_2}([1, \,a, \,b], \,[5 - b - 2\,n, \,5 - a - 2\,
n], \,1)
\]
}
\end{maplelatex}

\begin{maplelatex}
\mapleinline{inert}{2d}{222, _3F_2([a, b, -b+4-2*n],[b+1, 3+2*a-2*n-b],1);}{%
\[
222: \,\mathrm{_3F_2}([a, \,b, \, - b + 4 - 2\,n], \,[b + 1, \,3 + 
2\,a - 2\,n - b], \,1)
\]
}
\end{maplelatex}

\begin{maplelatex}
\mapleinline{inert}{2d}{223, _3F_2([a, b, 1/2*a-2*n+7/2-b],[1/2*a+b+1/2, -b+4+a-2*n],1);}{%
\[
223: \,\mathrm{_3F_2}([a, \,b, \,{\displaystyle \frac {a}{2}}  - 2
\,n + {\displaystyle \frac {7}{2}}  - b], \,[{\displaystyle 
\frac {a}{2}}  + b + {\displaystyle \frac {1}{2}} , \, - b + 4 + 
a - 2\,n], \,1)
\]
}
\end{maplelatex}

\begin{maplelatex}
\mapleinline{inert}{2d}{224, _3F_2([a, b, -n-a+3],[1+a, 2*b-a],1);}{%
\[
224: \,\mathrm{_3F_2}([a, \,b, \, - n - a + 3], \,[1 + a, \,2\,b - 
a], \,1)
\]
}
\end{maplelatex}

\begin{maplelatex}
\mapleinline{inert}{2d}{225, _3F_2([a, n+2*a-2, b],[1+a, -1-b+2*a+n],1);}{%
\[
\dagger 225: \,\mathrm{_3F_2}([a, \,n + 2\,a - 2, \,b], \,[1 + a, \, - 1 - 
b + 2\,a + n], \,1)
\]
}
\end{maplelatex}

\begin{maplelatex}
\mapleinline{inert}{2d}{226, _3F_2([a, b, 1],[1/2*b-1/2*n+2, n+2*a-2],1);}{%
\[
226: \,\mathrm{_3F_2}([a, \,b, \,1], \,[{\displaystyle \frac {b}{2}
}  - {\displaystyle \frac {n}{2}}  + 2, \,n + 2\,a - 2], \,1)
\]
}
\end{maplelatex}

\begin{maplelatex}
\mapleinline{inert}{2d}{227, _3F_2([a, b, 1],[-n-a+4, -b+2],1);}{%
\[
\dagger 227: \,\mathrm{_3F_2}([a, \,b, \,1], \,[ - n - a + 4, \, - b + 2], 
\,1)
\]
}
\end{maplelatex}

\begin{maplelatex}
\mapleinline{inert}{2d}{228, _3F_2([a, b, -1/2*n-a+3/2+1/2*b],[b+3-a-n,
a+1/2*b+1/2*n-1/2],1);}{%
\[
\dagger 228: \,\mathrm{_3F_2}([a, \,b, \, - {\displaystyle \frac {n}{2}} 
 - a + {\displaystyle \frac {3}{2}}  + {\displaystyle \frac {b}{2
}} ], \,[b + 3 - a - n, \,a + {\displaystyle \frac {b}{2}}  + 
{\displaystyle \frac {n}{2}}  - {\displaystyle \frac {1}{2}} ], 
\,1)
\]
}
\end{maplelatex}

\begin{maplelatex}
\mapleinline{inert}{2d}{229, _3F_2([a, b, 2*b],[-n-a+2*b+3, b+1],1);}{%
\[
\dagger 229: \,\mathrm{_3F_2}([a, \,b, \,2\,b], \,[ - n - a + 2\,b + 3, \,b
 + 1], \,1)
\]
}
\end{maplelatex}

\begin{maplelatex}
\mapleinline{inert}{2d}{230, _3F_2([a, b, 1/2*b-1/2+1/2*n],[1/2*a+1/2+1/2*b, b+1],1);}{%
\[
230: \,\mathrm{_3F_2}([a, \,b, \,{\displaystyle \frac {b}{2}}  - 
{\displaystyle \frac {1}{2}}  + {\displaystyle \frac {n}{2}} ], 
\,[{\displaystyle \frac {a}{2}}  + {\displaystyle \frac {1}{2}} 
 + {\displaystyle \frac {b}{2}} , \,b + 1], \,1)
\]
}
\end{maplelatex}

\begin{maplelatex}
\mapleinline{inert}{2d}{231, _3F_2([a, -a+1, b],[a+2*b-n+1, -a+2],1);}{%
\[
231: \,\mathrm{_3F_2}([a, \, - a + 1, \,b], \,[a + 2\,b - n + 1, \,
 - a + 2], \,1)
\]
}
\end{maplelatex}

\begin{maplelatex}
\mapleinline{inert}{2d}{232, _3F_2([a, b, -1+1/2*a+n],[1+a, 1/2*a+1/2+1/2*b],1);}{%
\[
232: \,\mathrm{_3F_2}([a, \,b, \, - 1 + {\displaystyle \frac {a}{2}
}  + n], \,[1 + a, \,{\displaystyle \frac {a}{2}}  + 
{\displaystyle \frac {1}{2}}  + {\displaystyle \frac {b}{2}} ], 
\,1)
\]
}
\end{maplelatex}

\begin{maplelatex}
\mapleinline{inert}{2d}{233, _3F_2([a, 2-n+1/2*a, b],[1+a, 1/2*a+1/2*b+2-n],1);}{%
\[
233: \,\mathrm{_3F_2}([a, \,2 - n + {\displaystyle \frac {a}{2}} , 
\,b], \,[1 + a, \,{\displaystyle \frac {a}{2}}  + {\displaystyle 
\frac {b}{2}}  + 2 - n], \,1)
\]
}
\end{maplelatex}

\begin{maplelatex}
\mapleinline{inert}{2d}{234, _3F_2([a, b, 1],[2*a-2*n+3, 1/2*b+1],1);}{%
\[
234: \,\mathrm{_3F_2}([a, \,b, \,1], \,[2\,a - 2\,n + 3, \,
{\displaystyle \frac {b}{2}}  + 1], \,1)
\]
}
\end{maplelatex}

\begin{maplelatex}
\mapleinline{inert}{2d}{235, _3F_2([a, b, 2-n+1/2*a-b],[1/2*a+n+b-1, -b+4+a-2*n],1);}{%
\[
\dagger 235: \,\mathrm{_3F_2}([a, \,b, \,2 - n + {\displaystyle \frac {a}{2
}}  - b], \,[{\displaystyle \frac {a}{2}}  + n + b - 1, \, - b + 
4 + a - 2\,n], \,1)
\]
}
\end{maplelatex}

\begin{maplelatex}
\mapleinline{inert}{2d}{236, _3F_2([a, 1/2*a-n+3/2, b],[1+a-b, 5/2+1/2*a-n],1);}{%
\[
\dagger 236: \,\mathrm{_3F_2}([a, \,{\displaystyle \frac {a}{2}}  - n + 
{\displaystyle \frac {3}{2}} , \,b], \,[1 + a - b, \,
{\displaystyle \frac {5}{2}}  + {\displaystyle \frac {a}{2}}  - n
], \,1)
\]
}
\end{maplelatex}

\begin{maplelatex}
\mapleinline{inert}{2d}{237, _3F_2([a, b, -b+1],[3+2*a-2*n-b, b+1],1);}{%
\[
237: \,\mathrm{_3F_2}([a, \,b, \, - b + 1], \,[3 + 2\,a - 2\,n - b
, \,b + 1], \,1)
\]
}
\end{maplelatex}

\begin{maplelatex}
\mapleinline{inert}{2d}{238, _3F_2([a, n+b-2, 1/2*n+1/2*a-1],[1+a-b, 1/2*n+1/2*a],1);}{%
\[
238: \,\mathrm{_3F_2}([a, \,n + b - 2, \,{\displaystyle \frac {n}{2
}}  + {\displaystyle \frac {a}{2}}  - 1], \,[1 + a - b, \,
{\displaystyle \frac {n}{2}}  + {\displaystyle \frac {a}{2}} ], 
\,1)
\]
}
\end{maplelatex}

\begin{maplelatex}
\mapleinline{inert}{2d}{239, _3F_2([a, -2*b+1, 2*n-2-a-b],[1+a-b, -2*b+2*n-a-1],1);}{%
\[
239: \,\mathrm{_3F_2}([a, \, - 2\,b + 1, \,2\,n - 2 - a - b], \,[1
 + a - b, \, - 2\,b + 2\,n - a - 1], \,1)
\]
}
\end{maplelatex}

\begin{maplelatex}
\mapleinline{inert}{2d}{240, _3F_2([a, 2*n-2-a, -1/2*b+1],[2*n-1-a, 1+a-b],1);}{%
\[
240: \,\mathrm{_3F_2}([a, \,2\,n - 2 - a, \, - {\displaystyle 
\frac {b}{2}}  + 1], \,[2\,n - 1 - a, \,1 + a - b], \,1)
\]
}
\end{maplelatex}

\begin{maplelatex}
\mapleinline{inert}{2d}{241, _3F_2([a, 2*a-1, -n-2*b+4],[2*a, 1+a-b],1);}{%
\[
241: \,\mathrm{_3F_2}([a, \,2\,a - 1, \, - n - 2\,b + 4], \,[2\,a, 
\,1 + a - b], \,1)
\]
}
\end{maplelatex}

\begin{maplelatex}
\mapleinline{inert}{2d}{242, _3F_2([1, -1/2*b+1, -n-2*a+4],[-a+2, -b+2],1);}{%
\[
242: \,\mathrm{_3F_2}([1, \, - {\displaystyle \frac {b}{2}}  + 1, 
\, - n - 2\,a + 4], \,[ - a + 2, \, - b + 2], \,1)
\]
}
\end{maplelatex}

\begin{maplelatex}
\mapleinline{inert}{2d}{243, _3F_2([1, 2*n-3+b, a-3+2*n],[-a+2, -b+2],1);}{%
\[
243: \,\mathrm{_3F_2}([1, \,2\,n - 3 + b, \,a - 3 + 2\,n], \,[ - a
 + 2, \, - b + 2], \,1)
\]
}
\end{maplelatex}

\begin{maplelatex}
\mapleinline{inert}{2d}{244, _3F_2([n, b-a+n, c-a+n],[2+n, n+1-a],1);}{%
\[
244: \,\mathrm{_3F_2}([n, \,b - a + n, \,c - a + n], \,[2 + n, \,n
 + 1 - a], \,1)
\]
}
\end{maplelatex}

\begin{maplelatex}
\mapleinline{inert}{2d}{245, _3F_2([2, 3-c, 3-b],[3-a, 2+n],1);}{%
\[
245: \,\mathrm{_3F_2}([2, \,3 - c, \,3 - b], \,[3 - a, \,2 + n], \,
1)
\]
}
\end{maplelatex}

\begin{maplelatex}
\mapleinline{inert}{2d}{246, _3F_2([a-c+3-b-n, -c+1, 2+a-c],[3-b+a-c, 3-c],1);}{%
\[
246: \,\mathrm{_3F_2}([a - c + 3 - b - n, \, - c + 1, \,2 + a - c]
, \,[3 - b + a - c, \,3 - c], \,1)
\]
}
\end{maplelatex}

\begin{maplelatex}
\mapleinline{inert}{2d}{247, _3F_2([a, b, n],[c, n+2-c+b+a],1);}{%
\[
247: \,\mathrm{_3F_2}([a, \,b, \,n], \,[c, \,n + 2 - c + b + a], \,
1)
\]
}
\end{maplelatex}

\begin{maplelatex}
\mapleinline{inert}{2d}{248, _3F_2([a, b, 2],[c, n+2-c+b+a],1);}{%
\[
248: \,\mathrm{_3F_2}([a, \,b, \,2], \,[c, \,n + 2 - c + b + a], \,
1)
\]
}
\end{maplelatex}

\begin{maplelatex}
\mapleinline{inert}{2d}{253, _3F_2([4/3-2*a, a-1/2*n+5/6, 4/3-a+1/2*n],[5/3,
-a-1/2*n+11/6],1);}{%
\[
253: \,\mathrm{_3F_2}([{\displaystyle \frac {4}{3}}  - 2\,a, \,a - 
{\displaystyle \frac {n}{2}}  + {\displaystyle \frac {5}{6}} , \,
{\displaystyle \frac {4}{3}}  - a + {\displaystyle \frac {n}{2}} 
], \,[{\displaystyle \frac {5}{3}} , \, - a - {\displaystyle 
\frac {n}{2}}  + {\displaystyle \frac {11}{6}} ], \,1)
\]
}
\end{maplelatex}

\begin{maplelatex}
\mapleinline{inert}{2d}{254, _3F_2([1-a+1/2*n, a-1/2*n+5/6, 4/3-a+1/2*n],[3/2,
4/3+a+1/2*n],1);}{%
\[
254: \,\mathrm{_3F_2}([1 - a + {\displaystyle \frac {n}{2}} , \,a
 - {\displaystyle \frac {n}{2}}  + {\displaystyle \frac {5}{6}} 
, \,{\displaystyle \frac {4}{3}}  - a + {\displaystyle \frac {n}{
2}} ], \,[{\displaystyle \frac {3}{2}} , \,{\displaystyle \frac {
4}{3}}  + a + {\displaystyle \frac {n}{2}} ], \,1)
\]
}
\end{maplelatex}

\begin{maplelatex}
\mapleinline{inert}{2d}{255, _3F_2([1-a+1/2*n, 2/3-a+1/2*n, a+1/6-1/2*n],[3/2,
2/3+a+1/2*n],1);}{%
\[
255: \,\mathrm{_3F_2}([1 - a + {\displaystyle \frac {n}{2}} , \,
{\displaystyle \frac {2}{3}}  - a + {\displaystyle \frac {n}{2}} 
, \,a + {\displaystyle \frac {1}{6}}  - {\displaystyle \frac {n}{
2}} ], \,[{\displaystyle \frac {3}{2}} , \,{\displaystyle \frac {
2}{3}}  + a + {\displaystyle \frac {n}{2}} ], \,1)
\]
}
\end{maplelatex}

\begin{maplelatex}
\mapleinline{inert}{2d}{256, _3F_2([4/3-2*a, 1-a+1/2*n, a-1/2*n+5/6],[4/3, 3/2-a-1/2*n],1);}{%
\[
256: \,\mathrm{_3F_2}([{\displaystyle \frac {4}{3}}  - 2\,a, \,1 - 
a + {\displaystyle \frac {n}{2}} , \,a - {\displaystyle \frac {n
}{2}}  + {\displaystyle \frac {5}{6}} ], \,[{\displaystyle 
\frac {4}{3}} , \,{\displaystyle \frac {3}{2}}  - a - 
{\displaystyle \frac {n}{2}} ], \,1)
\]
}
\end{maplelatex}

\begin{maplelatex}
\mapleinline{inert}{2d}{257, _3F_2([2/3-2*a, 1-a+1/2*n, a+1/6-1/2*n],[2/3, 3/2-a-1/2*n],1);}{%
\[
257: \,\mathrm{_3F_2}([{\displaystyle \frac {2}{3}}  - 2\,a, \,1 - 
a + {\displaystyle \frac {n}{2}} , \,a + {\displaystyle \frac {1
}{6}}  - {\displaystyle \frac {n}{2}} ], \,[{\displaystyle 
\frac {2}{3}} , \,{\displaystyle \frac {3}{2}}  - a - 
{\displaystyle \frac {n}{2}} ], \,1)
\]
}
\end{maplelatex}

\begin{maplelatex}
\mapleinline{inert}{2d}{258, _3F_2([2/3-2*a, -n+1/2, 4/3-2*a],[1-1/2*n-a, 3/2-a-1/2*n],1);}{%
\[
258: \,\mathrm{_3F_2}([{\displaystyle \frac {2}{3}}  - 2\,a, \, - n
 + {\displaystyle \frac {1}{2}} , \,{\displaystyle \frac {4}{3}} 
 - 2\,a], \,[1 - {\displaystyle \frac {n}{2}}  - a, \,
{\displaystyle \frac {3}{2}}  - a - {\displaystyle \frac {n}{2}} 
], \,1)
\]
}
\end{maplelatex}

\begin{maplelatex}
\mapleinline{inert}{2d}{259, _3F_2([a, c, b+n],[b, e],1);}{%
\[
259: \,\mathrm{_3F_2}([a, \,c, \,b + n], \,[b, \,e], \,1)
\]
}
\end{maplelatex}

\begin{maplelatex}
\mapleinline{inert}{2d}{260, _3F_2([a, b, -n],[c, e],1);}{%
\[
260: \,\mathrm{_3F_2}([a, \,b, \, - n], \,[c, \,e], \,1)
\]
}
\end{maplelatex}

\begin{maplelatex}
\mapleinline{inert}{2d}{267, _3F_2([a-1/2*n, -n+1/2, a+1/3-1/2*n],[4/3-a-1/2*n,
2*a+5/6-n],1);}{%
\[
267: \,\mathrm{_3F_2}([a - {\displaystyle \frac {n}{2}} , \, - n + 
{\displaystyle \frac {1}{2}} , \,a + {\displaystyle \frac {1}{3}
}  - {\displaystyle \frac {n}{2}} ], \,[{\displaystyle \frac {4}{
3}}  - a - {\displaystyle \frac {n}{2}} , \,2\,a + 
{\displaystyle \frac {5}{6}}  - n], \,1)
\]
}
\end{maplelatex}

\begin{maplelatex}
\mapleinline{inert}{2d}{268, _3F_2([a-1/2*n, 1-2*a, -a+5/6+1/2*n],[4/3-a-1/2*n, 4/3],1);}{%
\[
268: \,\mathrm{_3F_2}([a - {\displaystyle \frac {n}{2}} , \,1 - 2\,
a, \, - a + {\displaystyle \frac {5}{6}}  + {\displaystyle 
\frac {n}{2}} ], \,[{\displaystyle \frac {4}{3}}  - a - 
{\displaystyle \frac {n}{2}} , \,{\displaystyle \frac {4}{3}} ], 
\,1)
\]
}
\end{maplelatex}

\begin{maplelatex}
\mapleinline{inert}{2d}{269, _3F_2([-n+1/2, 1-2*a, 4/3-2*a],[4/3-a-1/2*n, -a-1/2*n+11/6],1);}{%
\[
269: \,\mathrm{_3F_2}([ - n + {\displaystyle \frac {1}{2}} , \,1 - 
2\,a, \,{\displaystyle \frac {4}{3}}  - 2\,a], \,[{\displaystyle 
\frac {4}{3}}  - a - {\displaystyle \frac {n}{2}} , \, - a - 
{\displaystyle \frac {n}{2}}  + {\displaystyle \frac {11}{6}} ], 
\,1)
\]
}
\end{maplelatex}

\begin{maplelatex}
\mapleinline{inert}{2d}{270, _3F_2([a+1/3-1/2*n, -a+5/6+1/2*n, 4/3-2*a],[4/3-a-1/2*n,
5/3],1);}{%
\[
270: \,\mathrm{_3F_2}([a + {\displaystyle \frac {1}{3}}  - 
{\displaystyle \frac {n}{2}} , \, - a + {\displaystyle \frac {5}{
6}}  + {\displaystyle \frac {n}{2}} , \,{\displaystyle \frac {4}{
3}}  - 2\,a], \,[{\displaystyle \frac {4}{3}}  - a - 
{\displaystyle \frac {n}{2}} , \,{\displaystyle \frac {5}{3}} ], 
\,1)
\]
}
\end{maplelatex}

\begin{maplelatex}
\mapleinline{inert}{2d}{271, _3F_2([a-1/2*n, 1/2+a-1/2*n, 2*a+1/3],[2*a+5/6-n, 4/3],1);}{%
\[
271: \,\mathrm{_3F_2}([a - {\displaystyle \frac {n}{2}} , \,
{\displaystyle \frac {1}{2}}  + a - {\displaystyle \frac {n}{2}} 
, \,2\,a + {\displaystyle \frac {1}{3}} ], \,[2\,a + 
{\displaystyle \frac {5}{6}}  - n, \,{\displaystyle \frac {4}{3}
} ], \,1)
\]
}
\end{maplelatex}

\begin{maplelatex}
\mapleinline{inert}{2d}{272, _3F_2([-n+1/2, 1/2+a-1/2*n, a-1/2*n+5/6],[2*a+5/6-n,
-a-1/2*n+11/6],1);}{%
\[
272: \,\mathrm{_3F_2}([ - n + {\displaystyle \frac {1}{2}} , \,
{\displaystyle \frac {1}{2}}  + a - {\displaystyle \frac {n}{2}} 
, \,a - {\displaystyle \frac {n}{2}}  + {\displaystyle \frac {5}{
6}} ], \,[2\,a + {\displaystyle \frac {5}{6}}  - n, \, - a - 
{\displaystyle \frac {n}{2}}  + {\displaystyle \frac {11}{6}} ], 
\,1)
\]
}
\end{maplelatex}

\begin{maplelatex}
\mapleinline{inert}{2d}{273, _3F_2([a+1/3-1/2*n, 2*a+1/3, a-1/2*n+5/6],[2*a+5/6-n, 5/3],1);}{%
\[
273: \,\mathrm{_3F_2}([a + {\displaystyle \frac {1}{3}}  - 
{\displaystyle \frac {n}{2}} , \,2\,a + {\displaystyle \frac {1}{
3}} , \,a - {\displaystyle \frac {n}{2}}  + {\displaystyle 
\frac {5}{6}} ], \,[2\,a + {\displaystyle \frac {5}{6}}  - n, \,
{\displaystyle \frac {5}{3}} ], \,1)
\]
}
\end{maplelatex}

\begin{maplelatex}
\mapleinline{inert}{2d}{274, _3F_2([1-2*a, 1/2+a-1/2*n, 4/3-a+1/2*n],[4/3,
-a-1/2*n+11/6],1);}{%
\[
274: \,\mathrm{_3F_2}([1 - 2\,a, \,{\displaystyle \frac {1}{2}}  + 
a - {\displaystyle \frac {n}{2}} , \,{\displaystyle \frac {4}{3}
}  - a + {\displaystyle \frac {n}{2}} ], \,[{\displaystyle 
\frac {4}{3}} , \, - a - {\displaystyle \frac {n}{2}}  + 
{\displaystyle \frac {11}{6}} ], \,1)
\]
}
\end{maplelatex}

\begin{maplelatex}
\mapleinline{inert}{2d}{275, _3F_2([-a+5/6+1/2*n, 2*a+1/3, 4/3-a+1/2*n],[4/3, 5/3],1);}{%
\[
275: \,\mathrm{_3F_2}([ - a + {\displaystyle \frac {5}{6}}  + 
{\displaystyle \frac {n}{2}} , \,2\,a + {\displaystyle \frac {1}{
3}} , \,{\displaystyle \frac {4}{3}}  - a + {\displaystyle 
\frac {n}{2}} ], \,[{\displaystyle \frac {4}{3}} , \,
{\displaystyle \frac {5}{3}} ], \,1)
\]
}
\end{maplelatex}

\begin{maplelatex}
\mapleinline{inert}{2d}{276, _3F_2([2*a-1/3, 2*a, a+1/6-1/2*n],[2/3+a+1/2*n,
3*a+1/2-1/2*n],1);}{%
\[
276: \,\mathrm{_3F_2}([2\,a - {\displaystyle \frac {1}{3}} , \,2\,a
, \,a + {\displaystyle \frac {1}{6}}  - {\displaystyle \frac {n}{
2}} ], \,[{\displaystyle \frac {2}{3}}  + a + {\displaystyle 
\frac {n}{2}} , \,3\,a + {\displaystyle \frac {1}{2}}  - 
{\displaystyle \frac {n}{2}} ], \,1)
\]
}
\end{maplelatex}

\begin{maplelatex}
\mapleinline{inert}{2d}{277, _3F_2([2*a-1/3, 1/2+n, 2/3-a+1/2*n],[2/3+a+1/2*n,
1+a+1/2*n],1);}{%
\[
277: \,\mathrm{_3F_2}([2\,a - {\displaystyle \frac {1}{3}} , \,
{\displaystyle \frac {1}{2}}  + n, \,{\displaystyle \frac {2}{3}
}  - a + {\displaystyle \frac {n}{2}} ], \,[{\displaystyle 
\frac {2}{3}}  + a + {\displaystyle \frac {n}{2}} , \,1 + a + 
{\displaystyle \frac {n}{2}} ], \,1)
\]
}
\end{maplelatex}

\begin{maplelatex}
\mapleinline{inert}{2d}{278, _3F_2([2*a, 1/2+n, 1-a+1/2*n],[2/3+a+1/2*n, 4/3+a+1/2*n],1);}{%
\[
278: \,\mathrm{_3F_2}([2\,a, \,{\displaystyle \frac {1}{2}}  + n, 
\,1 - a + {\displaystyle \frac {n}{2}} ], \,[{\displaystyle 
\frac {2}{3}}  + a + {\displaystyle \frac {n}{2}} , \,
{\displaystyle \frac {4}{3}}  + a + {\displaystyle \frac {n}{2}} 
], \,1)
\]
}
\end{maplelatex}

\begin{maplelatex}
\mapleinline{inert}{2d}{279, _3F_2([2*a-1/3, 2*a+1/3, 1/2+a-1/2*n],[3*a+1/2-1/2*n,
1+a+1/2*n],1);}{%
\[
279: \,\mathrm{_3F_2}([2\,a - {\displaystyle \frac {1}{3}} , \,2\,a
 + {\displaystyle \frac {1}{3}} , \,{\displaystyle \frac {1}{2}} 
 + a - {\displaystyle \frac {n}{2}} ], \,[3\,a + {\displaystyle 
\frac {1}{2}}  - {\displaystyle \frac {n}{2}} , \,1 + a + 
{\displaystyle \frac {n}{2}} ], \,1)
\]
}
\end{maplelatex}

\begin{maplelatex}
\mapleinline{inert}{2d}{280, _3F_2([2*a, 2*a+1/3, a-1/2*n+5/6],[3*a+1/2-1/2*n,
4/3+a+1/2*n],1);}{%
\[
280: \,\mathrm{_3F_2}([2\,a, \,2\,a + {\displaystyle \frac {1}{3}} 
, \,a - {\displaystyle \frac {n}{2}}  + {\displaystyle \frac {5}{
6}} ], \,[3\,a + {\displaystyle \frac {1}{2}}  - {\displaystyle 
\frac {n}{2}} , \,{\displaystyle \frac {4}{3}}  + a + 
{\displaystyle \frac {n}{2}} ], \,1)
\]
}
\end{maplelatex}

\begin{maplelatex}
\mapleinline{inert}{2d}{281, _3F_2([a+1/6-1/2*n, 1/2+a-1/2*n, a-1/2*n+5/6],[3*a+1/2-1/2*n,
3/2],1);}{%
\[
281: \,\mathrm{_3F_2}([a + {\displaystyle \frac {1}{6}}  - 
{\displaystyle \frac {n}{2}} , \,{\displaystyle \frac {1}{2}}  + 
a - {\displaystyle \frac {n}{2}} , \,a - {\displaystyle \frac {n
}{2}}  + {\displaystyle \frac {5}{6}} ], \,[3\,a + 
{\displaystyle \frac {1}{2}}  - {\displaystyle \frac {n}{2}} , \,
{\displaystyle \frac {3}{2}} ], \,1)
\]
}
\end{maplelatex}

\begin{maplelatex}
\mapleinline{inert}{2d}{282, _3F_2([1/2+n, 2*a+1/3, 4/3-a+1/2*n],[1+a+1/2*n,
4/3+a+1/2*n],1);}{%
\[
282: \,\mathrm{_3F_2}([{\displaystyle \frac {1}{2}}  + n, \,2\,a + 
{\displaystyle \frac {1}{3}} , \,{\displaystyle \frac {4}{3}}  - 
a + {\displaystyle \frac {n}{2}} ], \,[1 + a + {\displaystyle 
\frac {n}{2}} , \,{\displaystyle \frac {4}{3}}  + a + 
{\displaystyle \frac {n}{2}} ], \,1)
\]
}
\end{maplelatex}

\begin{maplelatex}
\mapleinline{inert}{2d}{283, _3F_2([2/3-a+1/2*n, 1/2+a-1/2*n, 4/3-a+1/2*n],[1+a+1/2*n,
3/2],1);}{%
\[
283: \,\mathrm{_3F_2}([{\displaystyle \frac {2}{3}}  - a + 
{\displaystyle \frac {n}{2}} , \,{\displaystyle \frac {1}{2}}  + 
a - {\displaystyle \frac {n}{2}} , \,{\displaystyle \frac {4}{3}
}  - a + {\displaystyle \frac {n}{2}} ], \,[1 + a + 
{\displaystyle \frac {n}{2}} , \,{\displaystyle \frac {3}{2}} ], 
\,1)
\]
}
\end{maplelatex}

\begin{maplelatex}
\mapleinline{inert}{2d}{284, _3F_2([a-1/3-1/2*n, 2/3-2*a, 1/2-a+1/2*n],[1-1/2*n-a, 2/3],1);}{%
\[
284: \,\mathrm{_3F_2}([a - {\displaystyle \frac {1}{3}}  - 
{\displaystyle \frac {n}{2}} , \,{\displaystyle \frac {2}{3}}  - 
2\,a, \,{\displaystyle \frac {1}{2}}  - a + {\displaystyle 
\frac {n}{2}} ], \,[1 - {\displaystyle \frac {n}{2}}  - a, \,
{\displaystyle \frac {2}{3}} ], \,1)
\]
}
\end{maplelatex}

\begin{maplelatex}
\mapleinline{inert}{2d}{285, _3F_2([a-1/3-1/2*n, -n+1/2, a+1/3-1/2*n],[1-1/2*n-a,
1/2+2*a-n],1);}{%
\[
285: \,\mathrm{_3F_2}([a - {\displaystyle \frac {1}{3}}  - 
{\displaystyle \frac {n}{2}} , \, - n + {\displaystyle \frac {1}{
2}} , \,a + {\displaystyle \frac {1}{3}}  - {\displaystyle 
\frac {n}{2}} ], \,[1 - {\displaystyle \frac {n}{2}}  - a, \,
{\displaystyle \frac {1}{2}}  + 2\,a - n], \,1)
\]
}
\end{maplelatex}

\begin{maplelatex}
\mapleinline{inert}{2d}{286, _3F_2([1/2-a+1/2*n, a+1/3-1/2*n, 4/3-2*a],[1-1/2*n-a, 4/3],1);}{%
\[
286: \,\mathrm{_3F_2}([{\displaystyle \frac {1}{2}}  - a + 
{\displaystyle \frac {n}{2}} , \,a + {\displaystyle \frac {1}{3}
}  - {\displaystyle \frac {n}{2}} , \,{\displaystyle \frac {4}{3}
}  - 2\,a], \,[1 - {\displaystyle \frac {n}{2}}  - a, \,
{\displaystyle \frac {4}{3}} ], \,1)
\]
}
\end{maplelatex}

\begin{maplelatex}
\mapleinline{inert}{2d}{287, _3F_2([a-1/3-1/2*n, a+1/6-1/2*n, 2*a],[2/3, 1/2+2*a-n],1);}{%
\[
287: \,\mathrm{_3F_2}([a - {\displaystyle \frac {1}{3}}  - 
{\displaystyle \frac {n}{2}} , \,a + {\displaystyle \frac {1}{6}
}  - {\displaystyle \frac {n}{2}} , \,2\,a], \,[{\displaystyle 
\frac {2}{3}} , \,{\displaystyle \frac {1}{2}}  + 2\,a - n], \,1)
\]
}
\end{maplelatex}

\begin{maplelatex}
\mapleinline{inert}{2d}{288, _3F_2([1/2-a+1/2*n, 2*a, 1-a+1/2*n],[2/3, 4/3],1);}{%
\[
288: \,\mathrm{_3F_2}([{\displaystyle \frac {1}{2}}  - a + 
{\displaystyle \frac {n}{2}} , \,2\,a, \,1 - a + {\displaystyle 
\frac {n}{2}} ], \,[{\displaystyle \frac {2}{3}} , \,
{\displaystyle \frac {4}{3}} ], \,1)
\]
}
\end{maplelatex}

\begin{maplelatex}
\mapleinline{inert}{2d}{289, _3F_2([-n+1/2, a+1/6-1/2*n, a-1/2*n+5/6],[1/2+2*a-n,
3/2-a-1/2*n],1);}{%
\[
289: \,\mathrm{_3F_2}([ - n + {\displaystyle \frac {1}{2}} , \,a + 
{\displaystyle \frac {1}{6}}  - {\displaystyle \frac {n}{2}} , \,
a - {\displaystyle \frac {n}{2}}  + {\displaystyle \frac {5}{6}} 
], \,[{\displaystyle \frac {1}{2}}  + 2\,a - n, \,{\displaystyle 
\frac {3}{2}}  - a - {\displaystyle \frac {n}{2}} ], \,1)
\]
}
\end{maplelatex}

\begin{maplelatex}
\mapleinline{inert}{2d}{290, _3F_2([a+1/3-1/2*n, 2*a, a-1/2*n+5/6],[1/2+2*a-n, 4/3],1);}{%
\[
290: \,\mathrm{_3F_2}([a + {\displaystyle \frac {1}{3}}  - 
{\displaystyle \frac {n}{2}} , \,2\,a, \,a - {\displaystyle 
\frac {n}{2}}  + {\displaystyle \frac {5}{6}} ], \,[
{\displaystyle \frac {1}{2}}  + 2\,a - n, \,{\displaystyle 
\frac {4}{3}} ], \,1)
\]
}
\end{maplelatex}

\begin{maplelatex}
\mapleinline{inert}{2d}{293, _3F_2([4/3-2*a, 2/3-2*a, 1/2-a+1/2*n],[1-1/2*n-a,
3/2-3*a+1/2*n],1);}{%
\[
293: \,\mathrm{_3F_2}([{\displaystyle \frac {4}{3}}  - 2\,a, \,
{\displaystyle \frac {2}{3}}  - 2\,a, \,{\displaystyle \frac {1}{
2}}  - a + {\displaystyle \frac {n}{2}} ], \,[1 - {\displaystyle 
\frac {n}{2}}  - a, \,{\displaystyle \frac {3}{2}}  - 3\,a + 
{\displaystyle \frac {n}{2}} ], \,1)
\]
}
\end{maplelatex}

\begin{maplelatex}
\mapleinline{inert}{2d}{294, _3F_2([2/3-2*a, 1-a+1/2*n, 1/2-a+1/2*n],[2/3, 7/6-2*a+n],1);}{%
\[
294: \,\mathrm{_3F_2}([{\displaystyle \frac {2}{3}}  - 2\,a, \,1 - 
a + {\displaystyle \frac {n}{2}} , \,{\displaystyle \frac {1}{2}
}  - a + {\displaystyle \frac {n}{2}} ], \,[{\displaystyle 
\frac {2}{3}} , \,{\displaystyle \frac {7}{6}}  - 2\,a + n], \,1)
\]
}
\end{maplelatex}

\begin{maplelatex}
\mapleinline{inert}{2d}{295, _3F_2([2/3-2*a, 2/3-a+1/2*n, 1/6-a+1/2*n],[1/3, 7/6-2*a+n],1);}{%
\[
295: \,\mathrm{_3F_2}([{\displaystyle \frac {2}{3}}  - 2\,a, \,
{\displaystyle \frac {2}{3}}  - a + {\displaystyle \frac {n}{2}} 
, \,{\displaystyle \frac {1}{6}}  - a + {\displaystyle \frac {n}{
2}} ], \,[{\displaystyle \frac {1}{3}} , \,{\displaystyle \frac {
7}{6}}  - 2\,a + n], \,1)
\]
}
\end{maplelatex}

\begin{maplelatex}
\mapleinline{inert}{2d}{298, _3F_2([1-2*a, 2/3-2*a, 2/3-a+1/2*n],[2-3*a+1/2*n,
7/6-a-1/2*n],1);}{%
\[
298: \,\mathrm{_3F_2}([1 - 2\,a, \,{\displaystyle \frac {2}{3}}  - 
2\,a, \,{\displaystyle \frac {2}{3}}  - a + {\displaystyle 
\frac {n}{2}} ], \,[2 - 3\,a + {\displaystyle \frac {n}{2}} , \,
{\displaystyle \frac {7}{6}}  - a - {\displaystyle \frac {n}{2}} 
], \,1)
\]
}
\end{maplelatex}

\begin{maplelatex}
\mapleinline{inert}{2d}{299, _3F_2([1-2*a, 2/3-a+1/2*n, 1/6-a+1/2*n],[2/3, 3/2-2*a+n],1);}{%
\[
299: \,\mathrm{_3F_2}([1 - 2\,a, \,{\displaystyle \frac {2}{3}}  - 
a + {\displaystyle \frac {n}{2}} , \,{\displaystyle \frac {1}{6}
}  - a + {\displaystyle \frac {n}{2}} ], \,[{\displaystyle 
\frac {2}{3}} , \,{\displaystyle \frac {3}{2}}  - 2\,a + n], \,1)
\]
}
\end{maplelatex}

\begin{maplelatex}
\mapleinline{inert}{2d}{300, _3F_2([1/2+n, -a+5/6+1/2*n, 1/6-a+1/2*n],[3/2-2*a+n,
a+1/2+1/2*n],1);}{%
\[
300: \,\mathrm{_3F_2}([{\displaystyle \frac {1}{2}}  + n, \, - a + 
{\displaystyle \frac {5}{6}}  + {\displaystyle \frac {n}{2}} , \,
{\displaystyle \frac {1}{6}}  - a + {\displaystyle \frac {n}{2}} 
], \,[{\displaystyle \frac {3}{2}}  - 2\,a + n, \,a + 
{\displaystyle \frac {1}{2}}  + {\displaystyle \frac {n}{2}} ], 
\,1)
\]
}
\end{maplelatex}

\begin{maplelatex}
\mapleinline{inert}{2d}{301, _3F_2([1-2*a, 2/3-2*a, 1/6-a+1/2*n],[3/2-3*a+1/2*n,
2/3-a-1/2*n],1);}{%
\[
301: \,\mathrm{_3F_2}([1 - 2\,a, \,{\displaystyle \frac {2}{3}}  - 
2\,a, \,{\displaystyle \frac {1}{6}}  - a + {\displaystyle 
\frac {n}{2}} ], \,[{\displaystyle \frac {3}{2}}  - 3\,a + 
{\displaystyle \frac {n}{2}} , \,{\displaystyle \frac {2}{3}}  - 
a - {\displaystyle \frac {n}{2}} ], \,1)
\]
}
\end{maplelatex}

\begin{maplelatex}
\mapleinline{inert}{2d}{302, _3F_2([-a+5/6+1/2*n, 1/2-a+1/2*n, 1/6-a+1/2*n],[1/2,
3/2-3*a+1/2*n],1);}{%
\[
302: \,\mathrm{_3F_2}([ - a + {\displaystyle \frac {5}{6}}  + 
{\displaystyle \frac {n}{2}} , \,{\displaystyle \frac {1}{2}}  - 
a + {\displaystyle \frac {n}{2}} , \,{\displaystyle \frac {1}{6}
}  - a + {\displaystyle \frac {n}{2}} ], \,[{\displaystyle 
\frac {1}{2}} , \,{\displaystyle \frac {3}{2}}  - 3\,a + 
{\displaystyle \frac {n}{2}} ], \,1)
\]
}
\end{maplelatex}

\begin{maplelatex}
\mapleinline{inert}{2d}{303, _3F_2([1/2+n, -a+5/6+1/2*n, 1/2-a+1/2*n],[11/6+n-2*a,
5/6+a+1/2*n],1);}{%
\[
303: \,\mathrm{_3F_2}([{\displaystyle \frac {1}{2}}  + n, \, - a + 
{\displaystyle \frac {5}{6}}  + {\displaystyle \frac {n}{2}} , \,
{\displaystyle \frac {1}{2}}  - a + {\displaystyle \frac {n}{2}} 
], \,[{\displaystyle \frac {11}{6}}  + n - 2\,a, \,
{\displaystyle \frac {5}{6}}  + a + {\displaystyle \frac {n}{2}} 
], \,1)
\]
}
\end{maplelatex}

\begin{maplelatex}
\mapleinline{inert}{2d}{304, _3F_2([1/2+n, 1/2-a+1/2*n, 1/6-a+1/2*n],[7/6-2*a+n,
1/6+a+1/2*n],1);}{%
\[
304: \,\mathrm{_3F_2}([{\displaystyle \frac {1}{2}}  + n, \,
{\displaystyle \frac {1}{2}}  - a + {\displaystyle \frac {n}{2}} 
, \,{\displaystyle \frac {1}{6}}  - a + {\displaystyle \frac {n}{
2}} ], \,[{\displaystyle \frac {7}{6}}  - 2\,a + n, \,
{\displaystyle \frac {1}{6}}  + a + {\displaystyle \frac {n}{2}} 
], \,1)
\]
}
\end{maplelatex}

\begin{maplelatex}
\mapleinline{inert}{2d}{305, _3F_2([a-1/2*n, 1-2*a, -n+1/2],[4/3-a-1/2*n, 2/3-a-1/2*n],1);}{%
\[
305: \,\mathrm{_3F_2}([a - {\displaystyle \frac {n}{2}} , \,1 - 2\,
a, \, - n + {\displaystyle \frac {1}{2}} ], \,[{\displaystyle 
\frac {4}{3}}  - a - {\displaystyle \frac {n}{2}} , \,
{\displaystyle \frac {2}{3}}  - a - {\displaystyle \frac {n}{2}} 
], \,1)
\]
}
\end{maplelatex}

\begin{maplelatex}
\mapleinline{inert}{2d}{306, _3F_2([a-1/2*n, -a+5/6+1/2*n, a+1/3-1/2*n],[4/3-a-1/2*n,
1/2],1);}{%
\[
306: \,\mathrm{_3F_2}([a - {\displaystyle \frac {n}{2}} , \, - a + 
{\displaystyle \frac {5}{6}}  + {\displaystyle \frac {n}{2}} , \,
a + {\displaystyle \frac {1}{3}}  - {\displaystyle \frac {n}{2}} 
], \,[{\displaystyle \frac {4}{3}}  - a - {\displaystyle \frac {n
}{2}} , \,{\displaystyle \frac {1}{2}} ], \,1)
\]
}
\end{maplelatex}

\begin{maplelatex}
\mapleinline{inert}{2d}{307, _3F_2([1-2*a, -a+5/6+1/2*n, 4/3-2*a],[4/3-a-1/2*n,
3/2-3*a+1/2*n],1);}{%
\[
307: \,\mathrm{_3F_2}([1 - 2\,a, \, - a + {\displaystyle \frac {5}{
6}}  + {\displaystyle \frac {n}{2}} , \,{\displaystyle \frac {4}{
3}}  - 2\,a], \,[{\displaystyle \frac {4}{3}}  - a - 
{\displaystyle \frac {n}{2}} , \,{\displaystyle \frac {3}{2}}  - 
3\,a + {\displaystyle \frac {n}{2}} ], \,1)
\]
}
\end{maplelatex}

\begin{maplelatex}
\mapleinline{inert}{2d}{308, _3F_2([-n+1/2, a+1/3-1/2*n, 4/3-2*a],[4/3-a-1/2*n,
1-1/2*n-a],1);}{%
\[
308: \,\mathrm{_3F_2}([ - n + {\displaystyle \frac {1}{2}} , \,a + 
{\displaystyle \frac {1}{3}}  - {\displaystyle \frac {n}{2}} , \,
{\displaystyle \frac {4}{3}}  - 2\,a], \,[{\displaystyle \frac {4
}{3}}  - a - {\displaystyle \frac {n}{2}} , \,1 - {\displaystyle 
\frac {n}{2}}  - a], \,1)
\]
}
\end{maplelatex}

\begin{maplelatex}
\mapleinline{inert}{2d}{309, _3F_2([a-1/2*n, 1/6-a+1/2*n, a-1/3-1/2*n],[2/3-a-1/2*n,
1/2],1);}{%
\[
309: \,\mathrm{_3F_2}([a - {\displaystyle \frac {n}{2}} , \,
{\displaystyle \frac {1}{6}}  - a + {\displaystyle \frac {n}{2}} 
, \,a - {\displaystyle \frac {1}{3}}  - {\displaystyle \frac {n}{
2}} ], \,[{\displaystyle \frac {2}{3}}  - a - {\displaystyle 
\frac {n}{2}} , \,{\displaystyle \frac {1}{2}} ], \,1)
\]
}
\end{maplelatex}

\begin{maplelatex}
\mapleinline{inert}{2d}{310, _3F_2([-n+1/2, a-1/3-1/2*n, 2/3-2*a],[2/3-a-1/2*n,
1-1/2*n-a],1);}{%
\[
310: \,\mathrm{_3F_2}([ - n + {\displaystyle \frac {1}{2}} , \,a - 
{\displaystyle \frac {1}{3}}  - {\displaystyle \frac {n}{2}} , \,
{\displaystyle \frac {2}{3}}  - 2\,a], \,[{\displaystyle \frac {2
}{3}}  - a - {\displaystyle \frac {n}{2}} , \,1 - {\displaystyle 
\frac {n}{2}}  - a], \,1)
\]
}
\end{maplelatex}

\begin{maplelatex}
\mapleinline{inert}{2d}{311, _3F_2([a+1/3-1/2*n, a-1/3-1/2*n, 1/2-a+1/2*n],[1/2,
1-1/2*n-a],1);}{%
\[
311: \,\mathrm{_3F_2}([a + {\displaystyle \frac {1}{3}}  - 
{\displaystyle \frac {n}{2}} , \,a - {\displaystyle \frac {1}{3}
}  - {\displaystyle \frac {n}{2}} , \,{\displaystyle \frac {1}{2}
}  - a + {\displaystyle \frac {n}{2}} ], \,[{\displaystyle 
\frac {1}{2}} , \,1 - {\displaystyle \frac {n}{2}}  - a], \,1)
\]
}
\end{maplelatex}

\begin{maplelatex}
\mapleinline{inert}{2d}{312, _3F_2([2*a-1/3, 2/3-a+1/2*n, a+1/6-1/2*n],[1/3,
2/3+a+1/2*n],1);}{%
\[
312: \,\mathrm{_3F_2}([2\,a - {\displaystyle \frac {1}{3}} , \,
{\displaystyle \frac {2}{3}}  - a + {\displaystyle \frac {n}{2}} 
, \,a + {\displaystyle \frac {1}{6}}  - {\displaystyle \frac {n}{
2}} ], \,[{\displaystyle \frac {1}{3}} , \,{\displaystyle \frac {
2}{3}}  + a + {\displaystyle \frac {n}{2}} ], \,1)
\]
}
\end{maplelatex}

\begin{maplelatex}
\mapleinline{inert}{2d}{313, _3F_2([2*a-1/3, 1/6-a+1/2*n, a-1/3-1/2*n],[1/3,
1/6+a+1/2*n],1);}{%
\[
313: \,\mathrm{_3F_2}([2\,a - {\displaystyle \frac {1}{3}} , \,
{\displaystyle \frac {1}{6}}  - a + {\displaystyle \frac {n}{2}} 
, \,a - {\displaystyle \frac {1}{3}}  - {\displaystyle \frac {n}{
2}} ], \,[{\displaystyle \frac {1}{3}} , \,{\displaystyle \frac {
1}{6}}  + a + {\displaystyle \frac {n}{2}} ], \,1)
\]
}
\end{maplelatex}

\begin{maplelatex}
\mapleinline{inert}{2d}{314, _3F_2([a+1/6-1/2*n, a-1/3-1/2*n, 2/3-2*a],[1/3, 2/3],1);}{%
\[
314: \,\mathrm{_3F_2}([a + {\displaystyle \frac {1}{6}}  - 
{\displaystyle \frac {n}{2}} , \,a - {\displaystyle \frac {1}{3}
}  - {\displaystyle \frac {n}{2}} , \,{\displaystyle \frac {2}{3}
}  - 2\,a], \,[{\displaystyle \frac {1}{3}} , \,{\displaystyle 
\frac {2}{3}} ], \,1)
\]
}
\end{maplelatex}

\begin{maplelatex}
\mapleinline{inert}{2d}{315, _3F_2([2*a-1/3, 1/2+n, 2*a],[2/3+a+1/2*n, 1/6+a+1/2*n],1);}{%
\[
315: \,\mathrm{_3F_2}([2\,a - {\displaystyle \frac {1}{3}} , \,
{\displaystyle \frac {1}{2}}  + n, \,2\,a], \,[{\displaystyle 
\frac {2}{3}}  + a + {\displaystyle \frac {n}{2}} , \,
{\displaystyle \frac {1}{6}}  + a + {\displaystyle \frac {n}{2}} 
], \,1)
\]
}
\end{maplelatex}

\begin{maplelatex}
\mapleinline{inert}{2d}{316, _3F_2([2/3-a+1/2*n, 1/2+n, 1-a+1/2*n],[2/3+a+1/2*n,
7/6-2*a+n],1);}{%
\[
316: \,\mathrm{_3F_2}([{\displaystyle \frac {2}{3}}  - a + 
{\displaystyle \frac {n}{2}} , \,{\displaystyle \frac {1}{2}}  + 
n, \,1 - a + {\displaystyle \frac {n}{2}} ], \,[{\displaystyle 
\frac {2}{3}}  + a + {\displaystyle \frac {n}{2}} , \,
{\displaystyle \frac {7}{6}}  - 2\,a + n], \,1)
\]
}
\end{maplelatex}

\begin{maplelatex}
\mapleinline{inert}{2d}{317, _3F_2([a+1/6-1/2*n, 2*a, 1-a+1/2*n],[2/3+a+1/2*n, 2/3],1);}{%
\[
317: \,\mathrm{_3F_2}([a + {\displaystyle \frac {1}{6}}  - 
{\displaystyle \frac {n}{2}} , \,2\,a, \,1 - a + {\displaystyle 
\frac {n}{2}} ], \,[{\displaystyle \frac {2}{3}}  + a + 
{\displaystyle \frac {n}{2}} , \,{\displaystyle \frac {2}{3}} ], 
\,1)
\]
}
\end{maplelatex}

\begin{maplelatex}
\mapleinline{inert}{2d}{318, _3F_2([a-1/3-1/2*n, 2*a, 1/2-a+1/2*n],[1/6+a+1/2*n, 2/3],1);}{%
\[
318: \,\mathrm{_3F_2}([a - {\displaystyle \frac {1}{3}}  - 
{\displaystyle \frac {n}{2}} , \,2\,a, \,{\displaystyle \frac {1
}{2}}  - a + {\displaystyle \frac {n}{2}} ], \,[{\displaystyle 
\frac {1}{6}}  + a + {\displaystyle \frac {n}{2}} , \,
{\displaystyle \frac {2}{3}} ], \,1)
\]
}
\end{maplelatex}

\begin{maplelatex}
\mapleinline{inert}{2d}{319, _3F_2([a-1/2*n+5/6, -n+1/2, 4/3-2*a],[-a-1/2*n+11/6,
3/2-a-1/2*n],1);}{%
\[
319: \,\mathrm{_3F_2}([a - {\displaystyle \frac {n}{2}}  + 
{\displaystyle \frac {5}{6}} , \, - n + {\displaystyle \frac {1}{
2}} , \,{\displaystyle \frac {4}{3}}  - 2\,a], \,[ - a - 
{\displaystyle \frac {n}{2}}  + {\displaystyle \frac {11}{6}} , 
\,{\displaystyle \frac {3}{2}}  - a - {\displaystyle \frac {n}{2}
} ], \,1)
\]
}
\end{maplelatex}

\begin{maplelatex}
\mapleinline{inert}{2d}{320, _3F_2([a-1/2*n+5/6, 1/2+a-1/2*n, 4/3-a+1/2*n],[-a-1/2*n+11/6,
3/2],1);}{%
\[
320: \,\mathrm{_3F_2}([a - {\displaystyle \frac {n}{2}}  + 
{\displaystyle \frac {5}{6}} , \,{\displaystyle \frac {1}{2}}  + 
a - {\displaystyle \frac {n}{2}} , \,{\displaystyle \frac {4}{3}
}  - a + {\displaystyle \frac {n}{2}} ], \,[ - a - 
{\displaystyle \frac {n}{2}}  + {\displaystyle \frac {11}{6}} , 
\,{\displaystyle \frac {3}{2}} ], \,1)
\]
}
\end{maplelatex}

\begin{maplelatex}
\mapleinline{inert}{2d}{321, _3F_2([-n+1/2, 1/2+a-1/2*n, 1-2*a],[-a-1/2*n+11/6,
7/6-a-1/2*n],1);}{%
\[
321: \,\mathrm{_3F_2}([ - n + {\displaystyle \frac {1}{2}} , \,
{\displaystyle \frac {1}{2}}  + a - {\displaystyle \frac {n}{2}} 
, \,1 - 2\,a], \,[ - a - {\displaystyle \frac {n}{2}}  + 
{\displaystyle \frac {11}{6}} , \,{\displaystyle \frac {7}{6}} 
 - a - {\displaystyle \frac {n}{2}} ], \,1)
\]
}
\end{maplelatex}

\begin{maplelatex}
\mapleinline{inert}{2d}{322, _3F_2([4/3-2*a, 4/3-a+1/2*n, 1-2*a],[-a-1/2*n+11/6,
2-3*a+1/2*n],1);}{%
\[
322: \,\mathrm{_3F_2}([{\displaystyle \frac {4}{3}}  - 2\,a, \,
{\displaystyle \frac {4}{3}}  - a + {\displaystyle \frac {n}{2}} 
, \,1 - 2\,a], \,[ - a - {\displaystyle \frac {n}{2}}  + 
{\displaystyle \frac {11}{6}} , \,2 - 3\,a + {\displaystyle 
\frac {n}{2}} ], \,1)
\]
}
\end{maplelatex}

\begin{maplelatex}
\mapleinline{inert}{2d}{323, _3F_2([a-1/2*n+5/6, a+1/6-1/2*n, 1-a+1/2*n],[3/2-a-1/2*n,
3/2],1);}{%
\[
323: \,\mathrm{_3F_2}([a - {\displaystyle \frac {n}{2}}  + 
{\displaystyle \frac {5}{6}} , \,a + {\displaystyle \frac {1}{6}
}  - {\displaystyle \frac {n}{2}} , \,1 - a + {\displaystyle 
\frac {n}{2}} ], \,[{\displaystyle \frac {3}{2}}  - a - 
{\displaystyle \frac {n}{2}} , \,{\displaystyle \frac {3}{2}} ], 
\,1)
\]
}
\end{maplelatex}

\begin{maplelatex}
\mapleinline{inert}{2d}{324, _3F_2([-n+1/2, a+1/6-1/2*n, 2/3-2*a],[3/2-a-1/2*n,
7/6-a-1/2*n],1);}{%
\[
324: \,\mathrm{_3F_2}([ - n + {\displaystyle \frac {1}{2}} , \,a + 
{\displaystyle \frac {1}{6}}  - {\displaystyle \frac {n}{2}} , \,
{\displaystyle \frac {2}{3}}  - 2\,a], \,[{\displaystyle \frac {3
}{2}}  - a - {\displaystyle \frac {n}{2}} , \,{\displaystyle 
\frac {7}{6}}  - a - {\displaystyle \frac {n}{2}} ], \,1)
\]
}
\end{maplelatex}

\begin{maplelatex}
\mapleinline{inert}{2d}{325, _3F_2([4/3-2*a, 1-a+1/2*n, 2/3-2*a],[3/2-a-1/2*n,
2-3*a+1/2*n],1);}{%
\[
325: \,\mathrm{_3F_2}([{\displaystyle \frac {4}{3}}  - 2\,a, \,1 - 
a + {\displaystyle \frac {n}{2}} , \,{\displaystyle \frac {2}{3}
}  - 2\,a], \,[{\displaystyle \frac {3}{2}}  - a - 
{\displaystyle \frac {n}{2}} , \,2 - 3\,a + {\displaystyle 
\frac {n}{2}} ], \,1)
\]
}
\end{maplelatex}

\begin{maplelatex}
\mapleinline{inert}{2d}{326, _3F_2([1/2+a-1/2*n, a+1/6-1/2*n, 2/3-a+1/2*n],[3/2,
7/6-a-1/2*n],1);}{%
\[
326: \,\mathrm{_3F_2}([{\displaystyle \frac {1}{2}}  + a - 
{\displaystyle \frac {n}{2}} , \,a + {\displaystyle \frac {1}{6}
}  - {\displaystyle \frac {n}{2}} , \,{\displaystyle \frac {2}{3}
}  - a + {\displaystyle \frac {n}{2}} ], \,[{\displaystyle 
\frac {3}{2}} , \,{\displaystyle \frac {7}{6}}  - a - 
{\displaystyle \frac {n}{2}} ], \,1)
\]
}
\end{maplelatex}

\begin{maplelatex}
\mapleinline{inert}{2d}{327, _3F_2([4/3-a+1/2*n, 1-a+1/2*n, 2/3-a+1/2*n],[3/2,
2-3*a+1/2*n],1);}{%
\[
327: \,\mathrm{_3F_2}([{\displaystyle \frac {4}{3}}  - a + 
{\displaystyle \frac {n}{2}} , \,1 - a + {\displaystyle \frac {n
}{2}} , \,{\displaystyle \frac {2}{3}}  - a + {\displaystyle 
\frac {n}{2}} ], \,[{\displaystyle \frac {3}{2}} , \,2 - 3\,a + 
{\displaystyle \frac {n}{2}} ], \,1)
\]
}
\end{maplelatex}

\begin{maplelatex}
\mapleinline{inert}{2d}{328, _3F_2([2*a+1/3, 4/3-a+1/2*n, 1/2+a-1/2*n],[4/3, 1+a+1/2*n],1);}{%
\[
328: \,\mathrm{_3F_2}([2\,a + {\displaystyle \frac {1}{3}} , \,
{\displaystyle \frac {4}{3}}  - a + {\displaystyle \frac {n}{2}} 
, \,{\displaystyle \frac {1}{2}}  + a - {\displaystyle \frac {n}{
2}} ], \,[{\displaystyle \frac {4}{3}} , \,1 + a + 
{\displaystyle \frac {n}{2}} ], \,1)
\]
}
\end{maplelatex}

\begin{maplelatex}
\mapleinline{inert}{2d}{329, _3F_2([2*a+1/3, -a+5/6+1/2*n, a-1/2*n],[4/3, a+1/2+1/2*n],1);}{%
\[
329: \,\mathrm{_3F_2}([2\,a + {\displaystyle \frac {1}{3}} , \, - a
 + {\displaystyle \frac {5}{6}}  + {\displaystyle \frac {n}{2}} 
, \,a - {\displaystyle \frac {n}{2}} ], \,[{\displaystyle \frac {
4}{3}} , \,a + {\displaystyle \frac {1}{2}}  + {\displaystyle 
\frac {n}{2}} ], \,1)
\]
}
\end{maplelatex}

\begin{maplelatex}
\mapleinline{inert}{2d}{330, _3F_2([4/3-a+1/2*n, -a+5/6+1/2*n, 1-2*a],[4/3, 3/2-2*a+n],1);}{%
\[
330: \,\mathrm{_3F_2}([{\displaystyle \frac {4}{3}}  - a + 
{\displaystyle \frac {n}{2}} , \, - a + {\displaystyle \frac {5}{
6}}  + {\displaystyle \frac {n}{2}} , \,1 - 2\,a], \,[
{\displaystyle \frac {4}{3}} , \,{\displaystyle \frac {3}{2}}  - 
2\,a + n], \,1)
\]
}
\end{maplelatex}

\begin{maplelatex}
\mapleinline{inert}{2d}{331, _3F_2([1/2+a-1/2*n, a-1/2*n, 1-2*a],[4/3, 2/3],1);}{%
\[
331: \,\mathrm{_3F_2}([{\displaystyle \frac {1}{2}}  + a - 
{\displaystyle \frac {n}{2}} , \,a - {\displaystyle \frac {n}{2}
} , \,1 - 2\,a], \,[{\displaystyle \frac {4}{3}} , \,
{\displaystyle \frac {2}{3}} ], \,1)
\]
}
\end{maplelatex}

\begin{maplelatex}
\mapleinline{inert}{2d}{332, _3F_2([2*a+1/3, 1/2+n, 2*a-1/3],[1+a+1/2*n, a+1/2+1/2*n],1);}{%
\[
332: \,\mathrm{_3F_2}([2\,a + {\displaystyle \frac {1}{3}} , \,
{\displaystyle \frac {1}{2}}  + n, \,2\,a - {\displaystyle 
\frac {1}{3}} ], \,[1 + a + {\displaystyle \frac {n}{2}} , \,a + 
{\displaystyle \frac {1}{2}}  + {\displaystyle \frac {n}{2}} ], 
\,1)
\]
}
\end{maplelatex}

\begin{maplelatex}
\mapleinline{inert}{2d}{333, _3F_2([4/3-a+1/2*n, 1/2+n, 2/3-a+1/2*n],[1+a+1/2*n,
3/2-2*a+n],1);}{%
\[
333: \,\mathrm{_3F_2}([{\displaystyle \frac {4}{3}}  - a + 
{\displaystyle \frac {n}{2}} , \,{\displaystyle \frac {1}{2}}  + 
n, \,{\displaystyle \frac {2}{3}}  - a + {\displaystyle \frac {n
}{2}} ], \,[1 + a + {\displaystyle \frac {n}{2}} , \,
{\displaystyle \frac {3}{2}}  - 2\,a + n], \,1)
\]
}
\end{maplelatex}

\begin{maplelatex}
\mapleinline{inert}{2d}{334, _3F_2([1/2+a-1/2*n, 2*a-1/3, 2/3-a+1/2*n],[1+a+1/2*n, 2/3],1);}{%
\[
334: \,\mathrm{_3F_2}([{\displaystyle \frac {1}{2}}  + a - 
{\displaystyle \frac {n}{2}} , \,2\,a - {\displaystyle \frac {1}{
3}} , \,{\displaystyle \frac {2}{3}}  - a + {\displaystyle 
\frac {n}{2}} ], \,[1 + a + {\displaystyle \frac {n}{2}} , \,
{\displaystyle \frac {2}{3}} ], \,1)
\]
}
\end{maplelatex}

\begin{maplelatex}
\mapleinline{inert}{2d}{335, _3F_2([a-1/2*n, 2*a-1/3, 1/6-a+1/2*n],[a+1/2+1/2*n, 2/3],1);}{%
\[
335: \,\mathrm{_3F_2}([a - {\displaystyle \frac {n}{2}} , \,2\,a - 
{\displaystyle \frac {1}{3}} , \,{\displaystyle \frac {1}{6}}  - 
a + {\displaystyle \frac {n}{2}} ], \,[a + {\displaystyle \frac {
1}{2}}  + {\displaystyle \frac {n}{2}} , \,{\displaystyle \frac {
2}{3}} ], \,1)
\]
}
\end{maplelatex}

\begin{maplelatex}
\mapleinline{inert}{2d}{336, _3F_2([a-1/2*n+5/6, 2*a+1/3, 4/3-a+1/2*n],[4/3+a+1/2*n,
5/3],1);}{%
\[
336: \,\mathrm{_3F_2}([a - {\displaystyle \frac {n}{2}}  + 
{\displaystyle \frac {5}{6}} , \,2\,a + {\displaystyle \frac {1}{
3}} , \,{\displaystyle \frac {4}{3}}  - a + {\displaystyle 
\frac {n}{2}} ], \,[{\displaystyle \frac {4}{3}}  + a + 
{\displaystyle \frac {n}{2}} , \,{\displaystyle \frac {5}{3}} ], 
\,1)
\]
}
\end{maplelatex}

\begin{maplelatex}
\mapleinline{inert}{2d}{337, _3F_2([a-1/2*n+5/6, 2*a, 1-a+1/2*n],[4/3+a+1/2*n, 4/3],1);}{%
\[
337: \,\mathrm{_3F_2}([a - {\displaystyle \frac {n}{2}}  + 
{\displaystyle \frac {5}{6}} , \,2\,a, \,1 - a + {\displaystyle 
\frac {n}{2}} ], \,[{\displaystyle \frac {4}{3}}  + a + 
{\displaystyle \frac {n}{2}} , \,{\displaystyle \frac {4}{3}} ], 
\,1)
\]
}
\end{maplelatex}

\begin{maplelatex}
\mapleinline{inert}{2d}{338, _3F_2([2*a+1/3, 2*a, 1/2+n],[4/3+a+1/2*n, 5/6+a+1/2*n],1);}{%
\[
338: \,\mathrm{_3F_2}([2\,a + {\displaystyle \frac {1}{3}} , \,2\,a
, \,{\displaystyle \frac {1}{2}}  + n], \,[{\displaystyle \frac {
4}{3}}  + a + {\displaystyle \frac {n}{2}} , \,{\displaystyle 
\frac {5}{6}}  + a + {\displaystyle \frac {n}{2}} ], \,1)
\]
}
\end{maplelatex}

\begin{maplelatex}
\mapleinline{inert}{2d}{339, _3F_2([4/3-a+1/2*n, 1-a+1/2*n, 1/2+n],[4/3+a+1/2*n,
11/6+n-2*a],1);}{%
\[
339: \,\mathrm{_3F_2}([{\displaystyle \frac {4}{3}}  - a + 
{\displaystyle \frac {n}{2}} , \,1 - a + {\displaystyle \frac {n
}{2}} , \,{\displaystyle \frac {1}{2}}  + n], \,[{\displaystyle 
\frac {4}{3}}  + a + {\displaystyle \frac {n}{2}} , \,
{\displaystyle \frac {11}{6}}  + n - 2\,a], \,1)
\]
}
\end{maplelatex}

\begin{maplelatex}
\mapleinline{inert}{2d}{340, _3F_2([a-1/2*n+5/6, a+1/3-1/2*n, 4/3-2*a],[5/3, 4/3],1);}{%
\[
340: \,\mathrm{_3F_2}([a - {\displaystyle \frac {n}{2}}  + 
{\displaystyle \frac {5}{6}} , \,a + {\displaystyle \frac {1}{3}
}  - {\displaystyle \frac {n}{2}} , \,{\displaystyle \frac {4}{3}
}  - 2\,a], \,[{\displaystyle \frac {5}{3}} , \,{\displaystyle 
\frac {4}{3}} ], \,1)
\]
}
\end{maplelatex}

\begin{maplelatex}
\mapleinline{inert}{2d}{341, _3F_2([2*a+1/3, a+1/3-1/2*n, -a+5/6+1/2*n],[5/3,
5/6+a+1/2*n],1);}{%
\[
341: \,\mathrm{_3F_2}([2\,a + {\displaystyle \frac {1}{3}} , \,a + 
{\displaystyle \frac {1}{3}}  - {\displaystyle \frac {n}{2}} , \,
 - a + {\displaystyle \frac {5}{6}}  + {\displaystyle \frac {n}{2
}} ], \,[{\displaystyle \frac {5}{3}} , \,{\displaystyle \frac {5
}{6}}  + a + {\displaystyle \frac {n}{2}} ], \,1)
\]
}
\end{maplelatex}

\begin{maplelatex}
\mapleinline{inert}{2d}{342, _3F_2([4/3-a+1/2*n, 4/3-2*a, -a+5/6+1/2*n],[5/3,
11/6+n-2*a],1);}{%
\[
342: \,\mathrm{_3F_2}([{\displaystyle \frac {4}{3}}  - a + 
{\displaystyle \frac {n}{2}} , \,{\displaystyle \frac {4}{3}}  - 
2\,a, \, - a + {\displaystyle \frac {5}{6}}  + {\displaystyle 
\frac {n}{2}} ], \,[{\displaystyle \frac {5}{3}} , \,
{\displaystyle \frac {11}{6}}  + n - 2\,a], \,1)
\]
}
\end{maplelatex}

\begin{maplelatex}
\mapleinline{inert}{2d}{343, _3F_2([2*a, a+1/3-1/2*n, 1/2-a+1/2*n],[4/3, 5/6+a+1/2*n],1);}{%
\[
343: \,\mathrm{_3F_2}([2\,a, \,a + {\displaystyle \frac {1}{3}}  - 
{\displaystyle \frac {n}{2}} , \,{\displaystyle \frac {1}{2}}  - 
a + {\displaystyle \frac {n}{2}} ], \,[{\displaystyle \frac {4}{3
}} , \,{\displaystyle \frac {5}{6}}  + a + {\displaystyle \frac {
n}{2}} ], \,1)
\]
}
\end{maplelatex}

\begin{maplelatex}
\mapleinline{inert}{2d}{344, _3F_2([1-a+1/2*n, 4/3-2*a, 1/2-a+1/2*n],[4/3, 11/6+n-2*a],1);}{%
\[
344: \,\mathrm{_3F_2}([1 - a + {\displaystyle \frac {n}{2}} , \,
{\displaystyle \frac {4}{3}}  - 2\,a, \,{\displaystyle \frac {1}{
2}}  - a + {\displaystyle \frac {n}{2}} ], \,[{\displaystyle 
\frac {4}{3}} , \,{\displaystyle \frac {11}{6}}  + n - 2\,a], \,1
)
\]
}
\end{maplelatex}

\begin{maplelatex}
\mapleinline{inert}{2d}{345, _3F_2([2/3-2*a, a+1/6-1/2*n, 2/3-a+1/2*n],[1/3,
7/6-a-1/2*n],1);}{%
\[
345: \,\mathrm{_3F_2}([{\displaystyle \frac {2}{3}}  - 2\,a, \,a + 
{\displaystyle \frac {1}{6}}  - {\displaystyle \frac {n}{2}} , \,
{\displaystyle \frac {2}{3}}  - a + {\displaystyle \frac {n}{2}} 
], \,[{\displaystyle \frac {1}{3}} , \,{\displaystyle \frac {7}{6
}}  - a - {\displaystyle \frac {n}{2}} ], \,1)
\]
}
\end{maplelatex}

\begin{maplelatex}
\mapleinline{inert}{2d}{346, _3F_2([1-2*a, 2/3-2*a, -n+1/2],[2/3-a-1/2*n, 7/6-a-1/2*n],1);}{%
\[
346: \,\mathrm{_3F_2}([1 - 2\,a, \,{\displaystyle \frac {2}{3}}  - 
2\,a, \, - n + {\displaystyle \frac {1}{2}} ], \,[{\displaystyle 
\frac {2}{3}}  - a - {\displaystyle \frac {n}{2}} , \,
{\displaystyle \frac {7}{6}}  - a - {\displaystyle \frac {n}{2}} 
], \,1)
\]
}
\end{maplelatex}

\begin{maplelatex}
\mapleinline{inert}{2d}{347, _3F_2([1-2*a, 2/3-a+1/2*n, 4/3-a+1/2*n],[2-3*a+1/2*n,
3/2-2*a+n],1);}{%
\[
347: \,\mathrm{_3F_2}([1 - 2\,a, \,{\displaystyle \frac {2}{3}}  - 
a + {\displaystyle \frac {n}{2}} , \,{\displaystyle \frac {4}{3}
}  - a + {\displaystyle \frac {n}{2}} ], \,[2 - 3\,a + 
{\displaystyle \frac {n}{2}} , \,{\displaystyle \frac {3}{2}}  - 
2\,a + n], \,1)
\]
}
\end{maplelatex}

\begin{maplelatex}
\mapleinline{inert}{2d}{348, _3F_2([1-2*a, 2/3-a+1/2*n, 1/2+a-1/2*n],[2/3, 7/6-a-1/2*n],1);}{%
\[
348: \,\mathrm{_3F_2}([1 - 2\,a, \,{\displaystyle \frac {2}{3}}  - 
a + {\displaystyle \frac {n}{2}} , \,{\displaystyle \frac {1}{2}
}  + a - {\displaystyle \frac {n}{2}} ], \,[{\displaystyle 
\frac {2}{3}} , \,{\displaystyle \frac {7}{6}}  - a - 
{\displaystyle \frac {n}{2}} ], \,1)
\]
}
\end{maplelatex}

\begin{maplelatex}
\mapleinline{inert}{2d}{349, _3F_2([-n+1/2, a+1/6-1/2*n, 1/2+a-1/2*n],[2*a+1/6-n,
7/6-a-1/2*n],1);}{%
\[
349: \,\mathrm{_3F_2}([ - n + {\displaystyle \frac {1}{2}} , \,a + 
{\displaystyle \frac {1}{6}}  - {\displaystyle \frac {n}{2}} , \,
{\displaystyle \frac {1}{2}}  + a - {\displaystyle \frac {n}{2}} 
], \,[2\,a + {\displaystyle \frac {1}{6}}  - n, \,{\displaystyle 
\frac {7}{6}}  - a - {\displaystyle \frac {n}{2}} ], \,1)
\]
}
\end{maplelatex}

\begin{maplelatex}
\mapleinline{inert}{2d}{350, _3F_2([1/2+n, 1/6-a+1/2*n, 2/3-a+1/2*n],[7/6-2*a+n,
3/2-2*a+n],1);}{%
\[
350: \,\mathrm{_3F_2}([{\displaystyle \frac {1}{2}}  + n, \,
{\displaystyle \frac {1}{6}}  - a + {\displaystyle \frac {n}{2}} 
, \,{\displaystyle \frac {2}{3}}  - a + {\displaystyle \frac {n}{
2}} ], \,[{\displaystyle \frac {7}{6}}  - 2\,a + n, \,
{\displaystyle \frac {3}{2}}  - 2\,a + n], \,1)
\]
}
\end{maplelatex}

\begin{maplelatex}
\mapleinline{inert}{2d}{351, _3F_2([2*a-1/3, a-1/3-1/2*n, a+1/6-1/2*n],[1/3, 2*a+1/6-n],1);}{%
\[
351: \,\mathrm{_3F_2}([2\,a - {\displaystyle \frac {1}{3}} , \,a - 
{\displaystyle \frac {1}{3}}  - {\displaystyle \frac {n}{2}} , \,
a + {\displaystyle \frac {1}{6}}  - {\displaystyle \frac {n}{2}} 
], \,[{\displaystyle \frac {1}{3}} , \,2\,a + {\displaystyle 
\frac {1}{6}}  - n], \,1)
\]
}
\end{maplelatex}

\begin{maplelatex}
\mapleinline{inert}{2d}{352, _3F_2([1-2*a, 1/6-a+1/2*n, -a+5/6+1/2*n],[3/2-2*a+n,
3/2-3*a+1/2*n],1);}{%
\[
352: \,\mathrm{_3F_2}([1 - 2\,a, \,{\displaystyle \frac {1}{6}}  - 
a + {\displaystyle \frac {n}{2}} , \, - a + {\displaystyle 
\frac {5}{6}}  + {\displaystyle \frac {n}{2}} ], \,[
{\displaystyle \frac {3}{2}}  - 2\,a + n, \,{\displaystyle 
\frac {3}{2}}  - 3\,a + {\displaystyle \frac {n}{2}} ], \,1)
\]
}
\end{maplelatex}

\begin{maplelatex}
\mapleinline{inert}{2d}{353, _3F_2([a-1/2*n, -n+1/2, a-1/3-1/2*n],[2*a+1/6-n,
2/3-a-1/2*n],1);}{%
\[
353: \,\mathrm{_3F_2}([a - {\displaystyle \frac {n}{2}} , \, - n + 
{\displaystyle \frac {1}{2}} , \,a - {\displaystyle \frac {1}{3}
}  - {\displaystyle \frac {n}{2}} ], \,[2\,a + {\displaystyle 
\frac {1}{6}}  - n, \,{\displaystyle \frac {2}{3}}  - a - 
{\displaystyle \frac {n}{2}} ], \,1)
\]
}
\end{maplelatex}

\begin{maplelatex}
\mapleinline{inert}{2d}{354, _3F_2([a-1/2*n, 2*a-1/3, 1/2+a-1/2*n],[2/3, 2*a+1/6-n],1);}{%
\[
354: \,\mathrm{_3F_2}([a - {\displaystyle \frac {n}{2}} , \,2\,a - 
{\displaystyle \frac {1}{3}} , \,{\displaystyle \frac {1}{2}}  + 
a - {\displaystyle \frac {n}{2}} ], \,[{\displaystyle \frac {2}{3
}} , \,2\,a + {\displaystyle \frac {1}{6}}  - n], \,1)
\]
}
\end{maplelatex}

\begin{maplelatex}
\mapleinline{inert}{2d}{355, _3F_2([4/3-2*a, 1-a+1/2*n, 4/3-a+1/2*n],[2-3*a+1/2*n,
11/6+n-2*a],1);}{%
\[
355: \,\mathrm{_3F_2}([{\displaystyle \frac {4}{3}}  - 2\,a, \,1 - 
a + {\displaystyle \frac {n}{2}} , \,{\displaystyle \frac {4}{3}
}  - a + {\displaystyle \frac {n}{2}} ], \,[2 - 3\,a + 
{\displaystyle \frac {n}{2}} , \,{\displaystyle \frac {11}{6}} 
 + n - 2\,a], \,1)
\]
}
\end{maplelatex}

\begin{maplelatex}
\mapleinline{inert}{2d}{356, _3F_2([4/3-2*a, 1/2-a+1/2*n, -a+5/6+1/2*n],[11/6+n-2*a,
3/2-3*a+1/2*n],1);}{%
\[
356: \,\mathrm{_3F_2}([{\displaystyle \frac {4}{3}}  - 2\,a, \,
{\displaystyle \frac {1}{2}}  - a + {\displaystyle \frac {n}{2}} 
, \, - a + {\displaystyle \frac {5}{6}}  + {\displaystyle \frac {
n}{2}} ], \,[{\displaystyle \frac {11}{6}}  + n - 2\,a, \,
{\displaystyle \frac {3}{2}}  - 3\,a + {\displaystyle \frac {n}{2
}} ], \,1)
\]
}
\end{maplelatex}

\begin{maplelatex}
\mapleinline{inert}{2d}{357, _3F_2([1/2+n, 1-a+1/2*n, 1/2-a+1/2*n],[7/6-2*a+n,
11/6+n-2*a],1);}{%
\[
357: \,\mathrm{_3F_2}([{\displaystyle \frac {1}{2}}  + n, \,1 - a
 + {\displaystyle \frac {n}{2}} , \,{\displaystyle \frac {1}{2}} 
 - a + {\displaystyle \frac {n}{2}} ], \,[{\displaystyle \frac {7
}{6}}  - 2\,a + n, \,{\displaystyle \frac {11}{6}}  + n - 2\,a], 
\,1)
\]
}
\end{maplelatex}

\begin{maplelatex}
\mapleinline{inert}{2d}{358, _3F_2([a-1/2*n, 1-2*a, 1/6-a+1/2*n],[2/3-a-1/2*n, 2/3],1);}{%
\[
358: \,\mathrm{_3F_2}([a - {\displaystyle \frac {n}{2}} , \,1 - 2\,
a, \,{\displaystyle \frac {1}{6}}  - a + {\displaystyle \frac {n
}{2}} ], \,[{\displaystyle \frac {2}{3}}  - a - {\displaystyle 
\frac {n}{2}} , \,{\displaystyle \frac {2}{3}} ], \,1)
\]
}
\end{maplelatex}

\begin{maplelatex}
\mapleinline{inert}{2d}{359, _3F_2([a-1/3-1/2*n, 1/6-a+1/2*n, 2/3-2*a],[2/3-a-1/2*n,
1/3],1);}{%
\[
359: \,\mathrm{_3F_2}([a - {\displaystyle \frac {1}{3}}  - 
{\displaystyle \frac {n}{2}} , \,{\displaystyle \frac {1}{6}}  - 
a + {\displaystyle \frac {n}{2}} , \,{\displaystyle \frac {2}{3}
}  - 2\,a], \,[{\displaystyle \frac {2}{3}}  - a - 
{\displaystyle \frac {n}{2}} , \,{\displaystyle \frac {1}{3}} ], 
\,1)
\]
}
\end{maplelatex}

\begin{maplelatex}
\mapleinline{inert}{2d}{360, _3F_2([1/6-a+1/2*n, 2*a-1/3, 2/3-a+1/2*n],[2/3, 1/3],1);}{%
\[
360: \,\mathrm{_3F_2}([{\displaystyle \frac {1}{6}}  - a + 
{\displaystyle \frac {n}{2}} , \,2\,a - {\displaystyle \frac {1}{
3}} , \,{\displaystyle \frac {2}{3}}  - a + {\displaystyle 
\frac {n}{2}} ], \,[{\displaystyle \frac {2}{3}} , \,
{\displaystyle \frac {1}{3}} ], \,1)
\]
}
\end{maplelatex}

\begin{maplelatex}
\mapleinline{inert}{2d}{361, _3F_2([1/2-a+1/2*n, 1/6-a+1/2*n, 2/3-2*a],[3/2-3*a+1/2*n,
7/6-2*a+n],1);}{%
\[
361: \,\mathrm{_3F_2}([{\displaystyle \frac {1}{2}}  - a + 
{\displaystyle \frac {n}{2}} , \,{\displaystyle \frac {1}{6}}  - 
a + {\displaystyle \frac {n}{2}} , \,{\displaystyle \frac {2}{3}
}  - 2\,a], \,[{\displaystyle \frac {3}{2}}  - 3\,a + 
{\displaystyle \frac {n}{2}} , \,{\displaystyle \frac {7}{6}}  - 
2\,a + n], \,1)
\]
}
\end{maplelatex}

\begin{maplelatex}
\mapleinline{inert}{2d}{362, _3F_2([2/3-2*a, 4/3-2*a, 1-2*a],[3/2-3*a+1/2*n,
2-3*a+1/2*n],1);}{%
\[
362: \,\mathrm{_3F_2}([{\displaystyle \frac {2}{3}}  - 2\,a, \,
{\displaystyle \frac {4}{3}}  - 2\,a, \,1 - 2\,a], \,[
{\displaystyle \frac {3}{2}}  - 3\,a + {\displaystyle \frac {n}{2
}} , \,2 - 3\,a + {\displaystyle \frac {n}{2}} ], \,1)
\]
}
\end{maplelatex}

\begin{maplelatex}
\mapleinline{inert}{2d}{363, _3F_2([2/3-2*a, 1-a+1/2*n, 2/3-a+1/2*n],[7/6-2*a+n,
2-3*a+1/2*n],1);}{%
\[
363: \,\mathrm{_3F_2}([{\displaystyle \frac {2}{3}}  - 2\,a, \,1 - 
a + {\displaystyle \frac {n}{2}} , \,{\displaystyle \frac {2}{3}
}  - a + {\displaystyle \frac {n}{2}} ], \,[{\displaystyle 
\frac {7}{6}}  - 2\,a + n, \,2 - 3\,a + {\displaystyle \frac {n}{
2}} ], \,1)
\]
}
\end{maplelatex}

\begin{maplelatex}
\mapleinline{inert}{2d}{364, _3F_2([-a+5/6+1/2*n, 1/2+n, 4/3-a+1/2*n],[11/6+n-2*a,
3/2-2*a+n],1);}{%
\[
364: \,\mathrm{_3F_2}([ - a + {\displaystyle \frac {5}{6}}  + 
{\displaystyle \frac {n}{2}} , \,{\displaystyle \frac {1}{2}}  + 
n, \,{\displaystyle \frac {4}{3}}  - a + {\displaystyle \frac {n
}{2}} ], \,[{\displaystyle \frac {11}{6}}  + n - 2\,a, \,
{\displaystyle \frac {3}{2}}  - 2\,a + n], \,1)
\]
}
\end{maplelatex}

\begin{maplelatex}
\mapleinline{inert}{2d}{366, _3F_2([a, -a+2*n, 3/2+a-1/2*c-1/2*b-n-m],[1+a-c, 1+a-b],1);}{%
\[
366: \,\mathrm{_3F_2}([a, \, - a + 2\,n, \,{\displaystyle \frac {3
}{2}}  + a - {\displaystyle \frac {c}{2}}  - {\displaystyle 
\frac {b}{2}}  - n - m], \,[1 + a - c, \,1 + a - b], \,1)
\]
}
\end{maplelatex}

\begin{maplelatex}
\mapleinline{inert}{2d}{367, _3F_2([a, -a+2*n, 1+a-1/2*c-1/2*b-n-m],[1+a-c, 1+a-b],1);}{%
\[
367: \,\mathrm{_3F_2}([a, \, - a + 2\,n, \,1 + a - {\displaystyle 
\frac {c}{2}}  - {\displaystyle \frac {b}{2}}  - n - m], \,[1 + a
 - c, \,1 + a - b], \,1)
\]
}
\end{maplelatex}

\begin{maplelatex}
\mapleinline{inert}{2d}{368, _3F_2([a, 1-a+2*n, 1+a-1/2*c-1/2*b-n-m],[1+a-c, 1+a-b],1);}{%
\[
368: \,\mathrm{_3F_2}([a, \,1 - a + 2\,n, \,1 + a - {\displaystyle 
\frac {c}{2}}  - {\displaystyle \frac {b}{2}}  - n - m], \,[1 + a
 - c, \,1 + a - b], \,1)
\]
}
\end{maplelatex}

\begin{maplelatex}
\mapleinline{inert}{2d}{375, _3F_2([a, b, -1+2*n-b+2*m],[c, -c+2*a+2*n],1);}{%
\[
375: \,\mathrm{_3F_2}([a, \,b, \, - 1 + 2\,n - b + 2\,m], \,[c, \,
 - c + 2\,a + 2\,n], \,1)
\]
}
\end{maplelatex}

\begin{maplelatex}
\mapleinline{inert}{2d}{376, _3F_2([a, b, c],[-a-1+b+2*n+2*m, -c+b+2*n],1);}{%
\[
376: \,\mathrm{_3F_2}([a, \,b, \,c], \,[ - a - 1 + b + 2\,n + 2\,m
, \, - c + b + 2\,n], \,1)
\]
}
\end{maplelatex}

\begin{maplelatex}
\mapleinline{inert}{2d}{377, _3F_2([a, b, c],[2*c+2*m-1, 1/2*a+1/2*b+n],1);}{%
\[
377: \,\mathrm{_3F_2}([a, \,b, \,c], \,[2\,c + 2\,m - 1, \,
{\displaystyle \frac {a}{2}}  + {\displaystyle \frac {b}{2}}  + n
], \,1)
\]
}
\end{maplelatex}

\begin{maplelatex}
\mapleinline{inert}{2d}{378, _3F_2([a, b, 2*n-b+2*m],[c, -c+2*a+2*n],1);}{%
\[
378: \,\mathrm{_3F_2}([a, \,b, \,2\,n - b + 2\,m], \,[c, \, - c + 2
\,a + 2\,n], \,1)
\]
}
\end{maplelatex}

\begin{maplelatex}
\mapleinline{inert}{2d}{379, _3F_2([a, b, c],[-a+b+2*n+2*m, -c+b+2*n],1);}{%
\[
379: \,\mathrm{_3F_2}([a, \,b, \,c], \,[ - a + b + 2\,n + 2\,m, \,
 - c + b + 2\,n], \,1)
\]
}
\end{maplelatex}

\begin{maplelatex}
\mapleinline{inert}{2d}{380, _3F_2([a, b, c],[2*c+2*m, 1/2*a+1/2*b+n],1);}{%
\[
380: \,\mathrm{_3F_2}([a, \,b, \,c], \,[2\,c + 2\,m, \,
{\displaystyle \frac {a}{2}}  + {\displaystyle \frac {b}{2}}  + n
], \,1)
\]
}
\end{maplelatex}

\begin{maplelatex}
\mapleinline{inert}{2d}{381, _3F_2([a, b, 2*n-b+2*m],[c, -c+2*a+2*n+1],1);}{%
\[
381: \,\mathrm{_3F_2}([a, \,b, \,2\,n - b + 2\,m], \,[c, \, - c + 2
\,a + 2\,n + 1], \,1)
\]
}
\end{maplelatex}

\begin{maplelatex}
\mapleinline{inert}{2d}{382, _3F_2([a, b, c],[-a+b+2*n+2*m, 1-c+b+2*n],1);}{%
\[
382: \,\mathrm{_3F_2}([a, \,b, \,c], \,[ - a + b + 2\,n + 2\,m, \,1
 - c + b + 2\,n], \,1)
\]
}
\end{maplelatex}

\begin{maplelatex}
\mapleinline{inert}{2d}{383, _3F_2([a, b, c],[2*c+2*m-1, 1/2+1/2*a+n+1/2*b],1);}{%
\[
383: \,\mathrm{_3F_2}([a, \,b, \,c], \,[2\,c + 2\,m - 1, \,
{\displaystyle \frac {1}{2}}  + {\displaystyle \frac {a}{2}}  + n
 + {\displaystyle \frac {b}{2}} ], \,1)
\]
}
\end{maplelatex}

\begin{maplelatex}
\mapleinline{inert}{2d}{384, _3F_2([b, -2*a+2*n+b, 3/2-1/2*c+1/2*b-n-m],[1+b-c, b-a+1],1);}{%
\[
384: \,\mathrm{_3F_2}([b, \, - 2\,a + 2\,n + b, \,{\displaystyle 
\frac {3}{2}}  - {\displaystyle \frac {c}{2}}  + {\displaystyle 
\frac {b}{2}}  - n - m], \,[1 + b - c, \,b - a + 1], \,1)
\]
}
\end{maplelatex}

\begin{maplelatex}
\mapleinline{inert}{2d}{385, _3F_2([b, -2*a+2*n+b, 1-1/2*c+1/2*b-n-m],[1+b-c, b-a+1],1);}{%
\[
385: \,\mathrm{_3F_2}([b, \, - 2\,a + 2\,n + b, \,1 - 
{\displaystyle \frac {c}{2}}  + {\displaystyle \frac {b}{2}}  - n
 - m], \,[1 + b - c, \,b - a + 1], \,1)
\]
}
\end{maplelatex}

\begin{maplelatex}
\mapleinline{inert}{2d}{386, _3F_2([b, 1+b-2*a+2*n, 1-1/2*c+1/2*b-n-m],[1+b-c, b-a+1],1);}{%
\[
386: \,\mathrm{_3F_2}([b, \,1 + b - 2\,a + 2\,n, \,1 - 
{\displaystyle \frac {c}{2}}  + {\displaystyle \frac {b}{2}}  - n
 - m], \,[1 + b - c, \,b - a + 1], \,1)
\]
}
\end{maplelatex}

\begin{maplelatex}
\mapleinline{inert}{2d}{387, _3F_2([c-2*n+1, 1+2*c-b-2*n, a+1-b+c-2*n],[-a+2*m+c,
2+c-b-2*n],1);}{%
\maplemultiline{
387: \mathrm{_3F_2}([c - 2\,n + 1, \,1 + 2\,c - b - 2\,n, \,a + 1
 - b + c - 2\,n],  \\
[ - a + 2\,m + c, \,2 + c - b - 2\,n], \,1) }
}
\end{maplelatex}

\begin{maplelatex}
\mapleinline{inert}{2d}{388, _3F_2([1/2*a-1/2*b-n+1, -1/2*a-n+1+1/2*b,
1+c-1/2*a-1/2*b-n],[2-1/2*a-1/2*b-n, 2*c+2*m-1/2*a-1/2*b-n],1);}{%
\maplemultiline{
388: \mathrm{_3F_2}([{\displaystyle \frac {a}{2}}  - 
{\displaystyle \frac {b}{2}}  - n + 1, \, - {\displaystyle 
\frac {a}{2}}  - n + 1 + {\displaystyle \frac {b}{2}} , \,1 + c
 - {\displaystyle \frac {a}{2}}  - {\displaystyle \frac {b}{2}} 
 - n],  \\
[2 - {\displaystyle \frac {a}{2}}  - {\displaystyle \frac {b}{2}
}  - n, \,2\,c + 2\,m - {\displaystyle \frac {a}{2}}  - 
{\displaystyle \frac {b}{2}}  - n], \,1) }
}
\end{maplelatex}

\begin{maplelatex}
\mapleinline{inert}{2d}{389, _3F_2([c-2*n+1, 1+2*c-b-2*n, a+1-b+c-2*n],[1-a+2*m+c,
2+c-b-2*n],1);}{%
\maplemultiline{
389: \mathrm{_3F_2}([c - 2\,n + 1, \,1 + 2\,c - b - 2\,n, \,a + 1
 - b + c - 2\,n],  \\
[1 - a + 2\,m + c, \,2 + c - b - 2\,n], \,1) }
}
\end{maplelatex}

\begin{maplelatex}
\mapleinline{inert}{2d}{390, _3F_2([1/2*a-1/2*b-n+1, -1/2*a-n+1+1/2*b,
1+c-1/2*a-1/2*b-n],[2-1/2*a-1/2*b-n, 1+2*c+2*m-1/2*a-1/2*b-n],1);}{%
\maplemultiline{
390: \mathrm{_3F_2}([{\displaystyle \frac {a}{2}}  - 
{\displaystyle \frac {b}{2}}  - n + 1, \, - {\displaystyle 
\frac {a}{2}}  - n + 1 + {\displaystyle \frac {b}{2}} , \,1 + c
 - {\displaystyle \frac {a}{2}}  - {\displaystyle \frac {b}{2}} 
 - n],  \\
[2 - {\displaystyle \frac {a}{2}}  - {\displaystyle \frac {b}{2}
}  - n, \,1 + 2\,c + 2\,m - {\displaystyle \frac {a}{2}}  - 
{\displaystyle \frac {b}{2}}  - n], \,1) }
}
\end{maplelatex}

\begin{maplelatex}
\mapleinline{inert}{2d}{391, _3F_2([c-2*n, 2*c-b-2*n, a-b+c-2*n],[-a+2*m+c, 1+c-b-2*n],1);}{%
\[
391: \,\mathrm{_3F_2}([c - 2\,n, \,2\,c - b - 2\,n, \,a - b + c - 2
\,n], \,[ - a + 2\,m + c, \,1 + c - b - 2\,n], \,1)
\]
}
\end{maplelatex}

\begin{maplelatex}
\mapleinline{inert}{2d}{392, _3F_2([1/2-1/2*a-n+1/2*b, 1/2+1/2*a-n-1/2*b,
1/2+c-1/2*a-n-1/2*b],[3/2-1/2*a-n-1/2*b,
2*c+2*m-1/2-1/2*a-n-1/2*b],1);}{%
\maplemultiline{
392: \mathrm{_3F_2}([{\displaystyle \frac {1}{2}}  - 
{\displaystyle \frac {a}{2}}  - n + {\displaystyle \frac {b}{2}} 
, \,{\displaystyle \frac {1}{2}}  + {\displaystyle \frac {a}{2}} 
 - n - {\displaystyle \frac {b}{2}} , \,{\displaystyle \frac {1}{
2}}  + c - {\displaystyle \frac {a}{2}}  - n - {\displaystyle 
\frac {b}{2}} ],  \\
[{\displaystyle \frac {3}{2}}  - {\displaystyle \frac {a}{2}}  - 
n - {\displaystyle \frac {b}{2}} , \,2\,c + 2\,m - 
{\displaystyle \frac {1}{2}}  - {\displaystyle \frac {a}{2}}  - n
 - {\displaystyle \frac {b}{2}} ], \,1) }
}
\end{maplelatex}

\begin{maplelatex}
\mapleinline{inert}{2d}{393, _3F_2([a, -a+2*m, b],[c, 2-c+2*b-2*n],1);}{%
\[
393: \,\mathrm{_3F_2}([a, \, - a + 2\,m, \,b], \,[c, \,2 - c + 2\,b
 - 2\,n], \,1)
\]
}
\end{maplelatex}

\begin{maplelatex}
\mapleinline{inert}{2d}{394, _3F_2([a, b, c],[1/2*c+m+1/2*b, 2+2*a-2*n-2*m],1);}{%
\[
394: \,\mathrm{_3F_2}([a, \,b, \,c], \,[{\displaystyle \frac {c}{2}
}  + m + {\displaystyle \frac {b}{2}} , \,2 + 2\,a - 2\,n - 2\,m]
, \,1)
\]
}
\end{maplelatex}

\begin{maplelatex}
\mapleinline{inert}{2d}{395, _3F_2([a, 1-a+2*m, b],[c, 2-c+2*b-2*n],1);}{%
\[
395: \,\mathrm{_3F_2}([a, \,1 - a + 2\,m, \,b], \,[c, \,2 - c + 2\,
b - 2\,n], \,1)
\]
}
\end{maplelatex}

\begin{maplelatex}
\mapleinline{inert}{2d}{396, _3F_2([a, b, c],[1/2*b+1/2+1/2*c+m, 1+2*a-2*n-2*m],1);}{%
\[
396: \,\mathrm{_3F_2}([a, \,b, \,c], \,[{\displaystyle \frac {b}{2}
}  + {\displaystyle \frac {1}{2}}  + {\displaystyle \frac {c}{2}
}  + m, \,1 + 2\,a - 2\,n - 2\,m], \,1)
\]
}
\end{maplelatex}

\begin{maplelatex}
\mapleinline{inert}{2d}{397, _3F_2([a, -a+2*m, b],[c, 1-c+2*b-2*n],1);}{%
\[
397: \,\mathrm{_3F_2}([a, \, - a + 2\,m, \,b], \,[c, \,1 - c + 2\,b
 - 2\,n], \,1)
\]
}
\end{maplelatex}

\begin{maplelatex}
\mapleinline{inert}{2d}{398, _3F_2([a, b, c],[1/2*c+m+1/2*b, 1+2*a-2*n-2*m],1);}{%
\[
398: \,\mathrm{_3F_2}([a, \,b, \,c], \,[{\displaystyle \frac {c}{2}
}  + m + {\displaystyle \frac {b}{2}} , \,1 + 2\,a - 2\,n - 2\,m]
, \,1)
\]
}
\end{maplelatex}

\begin{maplelatex}
\mapleinline{inert}{2d}{399, _3F_2([b, c-2*n+1, 2-2*n+a-2*m],[1+b-c, b-a+1],1);}{%
\[
399: \,\mathrm{_3F_2}([b, \,c - 2\,n + 1, \,2 - 2\,n + a - 2\,m], 
\,[1 + b - c, \,b - a + 1], \,1)
\]
}
\end{maplelatex}

\begin{maplelatex}
\mapleinline{inert}{2d}{400, _3F_2([b, c-2*n+1, 1-2*n+a-2*m],[1+b-c, b-a+1],1);}{%
\[
400: \,\mathrm{_3F_2}([b, \,c - 2\,n + 1, \,1 - 2\,n + a - 2\,m], 
\,[1 + b - c, \,b - a + 1], \,1)
\]
}
\end{maplelatex}

\begin{maplelatex}
\mapleinline{inert}{2d}{401, _3F_2([b, c-2*n, 1-2*n+a-2*m],[1+b-c, b-a+1],1);}{%
\[
401: \,\mathrm{_3F_2}([b, \,c - 2\,n, \,1 - 2\,n + a - 2\,m], \,[1
 + b - c, \,b - a + 1], \,1)
\]
}
\end{maplelatex}

\begin{maplelatex}
\mapleinline{inert}{2d}{402, _3F_2([a, 1+a-c, a+c-2*b+2*n],[1+a-b, 1+2*a-2*m],1);}{%
\[
402: \,\mathrm{_3F_2}([a, \,1 + a - c, \,a + c - 2\,b + 2\,n], \,[1
 + a - b, \,1 + 2\,a - 2\,m], \,1)
\]
}
\end{maplelatex}

\begin{maplelatex}
\mapleinline{inert}{2d}{403, _3F_2([a, b, c],[1+2*c-2*m, 1/2*a+1/2*b-n+1],1);}{%
\[
403: \,\mathrm{_3F_2}([a, \,b, \,c], \,[1 + 2\,c - 2\,m, \,
{\displaystyle \frac {a}{2}}  + {\displaystyle \frac {b}{2}}  - n
 + 1], \,1)
\]
}
\end{maplelatex}

\begin{maplelatex}
\mapleinline{inert}{2d}{404, _3F_2([a, b, c],[2*b-1+2*m, 1/2*a+1/2*c+3/2-n-m],1);}{%
\[
404: \,\mathrm{_3F_2}([a, \,b, \,c], \,[2\,b - 1 + 2\,m, \,
{\displaystyle \frac {a}{2}}  + {\displaystyle \frac {c}{2}}  + 
{\displaystyle \frac {3}{2}}  - n - m], \,1)
\]
}
\end{maplelatex}

\begin{maplelatex}
\mapleinline{inert}{2d}{405, _3F_2([a, b, -b+2-2*n],[c, 3-c+2*a-2*n-2*m],1);}{%
\[
405: \,\mathrm{_3F_2}([a, \,b, \, - b + 2 - 2\,n], \,[c, \,3 - c + 
2\,a - 2\,n - 2\,m], \,1)
\]
}
\end{maplelatex}

\begin{maplelatex}
\mapleinline{inert}{2d}{406, _3F_2([a, b, 3-2*n-b-2*m],[c, -c+2*a+2-2*n],1);}{%
\[
406: \,\mathrm{_3F_2}([a, \,b, \,3 - 2\,n - b - 2\,m], \,[c, \, - c
 + 2\,a + 2 - 2\,n], \,1)
\]
}
\end{maplelatex}

\begin{maplelatex}
\mapleinline{inert}{2d}{407, _3F_2([a, b, c],[-2*m+2*c, 1/2*a+1/2*b-n+1],1);}{%
\[
407: \,\mathrm{_3F_2}([a, \,b, \,c], \,[ - 2\,m + 2\,c, \,
{\displaystyle \frac {a}{2}}  + {\displaystyle \frac {b}{2}}  - n
 + 1], \,1)
\]
}
\end{maplelatex}

\begin{maplelatex}
\mapleinline{inert}{2d}{408, _3F_2([a, b, c],[2*m+2*b, 1/2*c+1/2*a-n-m+1],1);}{%
\[
408: \,\mathrm{_3F_2}([a, \,b, \,c], \,[2\,m + 2\,b, \,
{\displaystyle \frac {c}{2}}  + {\displaystyle \frac {a}{2}}  - n
 - m + 1], \,1)
\]
}
\end{maplelatex}

\begin{maplelatex}
\mapleinline{inert}{2d}{409, _3F_2([a, b, -b+2-2*n],[c, 2*a-c-2*n-2*m+2],1);}{%
\[
409: \,\mathrm{_3F_2}([a, \,b, \, - b + 2 - 2\,n], \,[c, \,2\,a - c
 - 2\,n - 2\,m + 2], \,1)
\]
}
\end{maplelatex}

\begin{maplelatex}
\mapleinline{inert}{2d}{410, _3F_2([a, b, 2-2*n-b-2*m],[c, -c+2*a+2-2*n],1);}{%
\[
410: \,\mathrm{_3F_2}([a, \,b, \,2 - 2\,n - b - 2\,m], \,[c, \, - c
 + 2\,a + 2 - 2\,n], \,1)
\]
}
\end{maplelatex}

\begin{maplelatex}
\mapleinline{inert}{2d}{411, _3F_2([a, b, c],[2*b-1+2*m, 1/2*c+1/2*a-n-m+1],1);}{%
\[
411: \,\mathrm{_3F_2}([a, \,b, \,c], \,[2\,b - 1 + 2\,m, \,
{\displaystyle \frac {c}{2}}  + {\displaystyle \frac {a}{2}}  - n
 - m + 1], \,1)
\]
}
\end{maplelatex}

\begin{maplelatex}
\mapleinline{inert}{2d}{412, _3F_2([a, b, -b+1-2*n],[c, 2*a-c-2*n-2*m+2],1);}{%
\[
412: \,\mathrm{_3F_2}([a, \,b, \, - b + 1 - 2\,n], \,[c, \,2\,a - c
 - 2\,n - 2\,m + 2], \,1)
\]
}
\end{maplelatex}

\begin{maplelatex}
\mapleinline{inert}{2d}{413, _3F_2([a, b, 2-2*n-b-2*m],[c, 2*a-2*n+1-c],1);}{%
\[
413: \,\mathrm{_3F_2}([a, \,b, \,2 - 2\,n - b - 2\,m], \,[c, \,2\,a
 - 2\,n + 1 - c], \,1)
\]
}
\end{maplelatex}

\begin{maplelatex}
\mapleinline{inert}{2d}{414, _3F_2([2*a+1/3-n, 5/6+a-n, a-n+1/2],[3*a-2*n+3/2,
2*a+5/6-n],1);}{%
\[
414: \,\mathrm{_3F_2}([2\,a + {\displaystyle \frac {1}{3}}  - n, \,
{\displaystyle \frac {5}{6}}  + a - n, \,a - n + {\displaystyle 
\frac {1}{2}} ], \,[3\,a - 2\,n + {\displaystyle \frac {3}{2}} , 
\,2\,a + {\displaystyle \frac {5}{6}}  - n], \,1)
\]
}
\end{maplelatex}

\begin{maplelatex}
\mapleinline{inert}{2d}{415, _3F_2([7/6+a-n, 2*a+2/3-n, a-n+1/2],[3*a-2*n+3/2,
7/6+2*a-n],1);}{%
\[
415: \,\mathrm{_3F_2}([{\displaystyle \frac {7}{6}}  + a - n, \,2\,
a + {\displaystyle \frac {2}{3}}  - n, \,a - n + {\displaystyle 
\frac {1}{2}} ], \,[3\,a - 2\,n + {\displaystyle \frac {3}{2}} , 
\,{\displaystyle \frac {7}{6}}  + 2\,a - n], \,1)
\]
}
\end{maplelatex}

\begin{maplelatex}
\mapleinline{inert}{2d}{416, _3F_2([2/3-a, 7/6+a-n, 2*a+2/3-n],[4/3, 5/3+a-n],1);}{%
\[
416: \,\mathrm{_3F_2}([{\displaystyle \frac {2}{3}}  - a, \,
{\displaystyle \frac {7}{6}}  + a - n, \,2\,a + {\displaystyle 
\frac {2}{3}}  - n], \,[{\displaystyle \frac {4}{3}} , \,
{\displaystyle \frac {5}{3}}  + a - n], \,1)
\]
}
\end{maplelatex}

\begin{maplelatex}
\mapleinline{inert}{2d}{417, _3F_2([2*a+1-n, 7/6+a-n, 5/6+a-n],[3*a-2*n+3/2, 2*a+3/2-n],1);}{%
\[
417: \,\mathrm{_3F_2}([2\,a + 1 - n, \,{\displaystyle \frac {7}{6}
}  + a - n, \,{\displaystyle \frac {5}{6}}  + a - n], \,[3\,a - 2
\,n + {\displaystyle \frac {3}{2}} , \,2\,a + {\displaystyle 
\frac {3}{2}}  - n], \,1)
\]
}
\end{maplelatex}

\begin{maplelatex}
\mapleinline{inert}{2d}{418, _3F_2([-a+1, 2*a+1-n, 7/6+a-n],[5/3, 5/3+a-n],1);}{%
\[
418: \,\mathrm{_3F_2}([ - a + 1, \,2\,a + 1 - n, \,{\displaystyle 
\frac {7}{6}}  + a - n], \,[{\displaystyle \frac {5}{3}} , \,
{\displaystyle \frac {5}{3}}  + a - n], \,1)
\]
}
\end{maplelatex}

\begin{maplelatex}
\mapleinline{inert}{2d}{419, _3F_2([-a+1, 2*a+1-n, 5/6+a-n],[4/3, 4/3+a-n],1);}{%
\[
419: \,\mathrm{_3F_2}([ - a + 1, \,2\,a + 1 - n, \,{\displaystyle 
\frac {5}{6}}  + a - n], \,[{\displaystyle \frac {4}{3}} , \,
{\displaystyle \frac {4}{3}}  + a - n], \,1)
\]
}
\end{maplelatex}

\begin{maplelatex}
\mapleinline{inert}{2d}{420, _3F_2([a+2/3, a+1/3, 2*a+1-n],[3*a+1-n, 2*a+3/2-n],1);}{%
\[
420: \,\mathrm{_3F_2}([a + {\displaystyle \frac {2}{3}} , \,a + 
{\displaystyle \frac {1}{3}} , \,2\,a + 1 - n], \,[3\,a + 1 - n, 
\,2\,a + {\displaystyle \frac {3}{2}}  - n], \,1)
\]
}
\end{maplelatex}

\begin{maplelatex}
\mapleinline{inert}{2d}{421, _3F_2([a+2/3, a, 2*a+2/3-n],[3*a+1-n, 7/6+2*a-n],1);}{%
\[
421: \,\mathrm{_3F_2}([a + {\displaystyle \frac {2}{3}} , \,a, \,2
\,a + {\displaystyle \frac {2}{3}}  - n], \,[3\,a + 1 - n, \,
{\displaystyle \frac {7}{6}}  + 2\,a - n], \,1)
\]
}
\end{maplelatex}

\begin{maplelatex}
\mapleinline{inert}{2d}{422, _3F_2([a+1/3, a, 2*a+1/3-n],[3*a+1-n, 2*a+5/6-n],1);}{%
\[
422: \,\mathrm{_3F_2}([a + {\displaystyle \frac {1}{3}} , \,a, \,2
\,a + {\displaystyle \frac {1}{3}}  - n], \,[3\,a + 1 - n, \,2\,a
 + {\displaystyle \frac {5}{6}}  - n], \,1)
\]
}
\end{maplelatex}

\begin{maplelatex}
\mapleinline{inert}{2d}{423, _3F_2([2*a+1-n, 2*a+2/3-n, 2*a+1/3-n],[3*a+1-n,
3*a-2*n+3/2],1);}{%
\[
423: \,\mathrm{_3F_2}([2\,a + 1 - n, \,2\,a + {\displaystyle 
\frac {2}{3}}  - n, \,2\,a + {\displaystyle \frac {1}{3}}  - n], 
\,[3\,a + 1 - n, \,3\,a - 2\,n + {\displaystyle \frac {3}{2}} ], 
\,1)
\]
}
\end{maplelatex}

\begin{maplelatex}
\mapleinline{inert}{2d}{424, _3F_2([a+2/3, 1/2, 7/6+a-n],[2*a+3/2-n, 7/6+2*a-n],1);}{%
\[
424: \,\mathrm{_3F_2}([a + {\displaystyle \frac {2}{3}} , \,
{\displaystyle \frac {1}{2}} , \,{\displaystyle \frac {7}{6}}  + 
a - n], \,[2\,a + {\displaystyle \frac {3}{2}}  - n, \,
{\displaystyle \frac {7}{6}}  + 2\,a - n], \,1)
\]
}
\end{maplelatex}

\begin{maplelatex}
\mapleinline{inert}{2d}{425, _3F_2([a+1/3, 1/2, 5/6+a-n],[2*a+3/2-n, 2*a+5/6-n],1);}{%
\[
425: \,\mathrm{_3F_2}([a + {\displaystyle \frac {1}{3}} , \,
{\displaystyle \frac {1}{2}} , \,{\displaystyle \frac {5}{6}}  + 
a - n], \,[2\,a + {\displaystyle \frac {3}{2}}  - n, \,2\,a + 
{\displaystyle \frac {5}{6}}  - n], \,1)
\]
}
\end{maplelatex}

\begin{maplelatex}
\mapleinline{inert}{2d}{426, _3F_2([a, 1/2, a-n+1/2],[7/6+2*a-n, 2*a+5/6-n],1);}{%
\[
426: \,\mathrm{_3F_2}([a, \,{\displaystyle \frac {1}{2}} , \,a - n
 + {\displaystyle \frac {1}{2}} ], \,[{\displaystyle \frac {7}{6}
}  + 2\,a - n, \,2\,a + {\displaystyle \frac {5}{6}}  - n], \,1)
\]
}
\end{maplelatex}

\begin{maplelatex}
\mapleinline{inert}{2d}{427, _3F_2([1/2+n-a, -a+1, 2/3+n-2*a],[7/6-2*a+n, 5/3],1);}{%
\[
427: \,\mathrm{_3F_2}([{\displaystyle \frac {1}{2}}  + n - a, \, - 
a + 1, \,{\displaystyle \frac {2}{3}}  + n - 2\,a], \,[
{\displaystyle \frac {7}{6}}  - 2\,a + n, \,{\displaystyle 
\frac {5}{3}} ], \,1)
\]
}
\end{maplelatex}

\begin{maplelatex}
\mapleinline{inert}{2d}{428, _3F_2([1/2+n-a, 1/2, 1/6+n-a],[7/6-2*a+n, 7/6+a],1);}{%
\[
428: \,\mathrm{_3F_2}([{\displaystyle \frac {1}{2}}  + n - a, \,
{\displaystyle \frac {1}{2}} , \,{\displaystyle \frac {1}{6}}  + 
n - a], \,[{\displaystyle \frac {7}{6}}  - 2\,a + n, \,
{\displaystyle \frac {7}{6}}  + a], \,1)
\]
}
\end{maplelatex}

\begin{maplelatex}
\mapleinline{inert}{2d}{429, _3F_2([-a+1, 1/2, 2/3-a],[7/6-2*a+n, 5/3+a-n],1);}{%
\[
429: \,\mathrm{_3F_2}([ - a + 1, \,{\displaystyle \frac {1}{2}} , 
\,{\displaystyle \frac {2}{3}}  - a], \,[{\displaystyle \frac {7
}{6}}  - 2\,a + n, \,{\displaystyle \frac {5}{3}}  + a - n], \,1)
\]
}
\end{maplelatex}

\begin{maplelatex}
\mapleinline{inert}{2d}{430, _3F_2([2/3+n-2*a, 1/6+n-a, 2/3-a],[7/6-2*a+n, 4/3],1);}{%
\[
430: \,\mathrm{_3F_2}([{\displaystyle \frac {2}{3}}  + n - 2\,a, \,
{\displaystyle \frac {1}{6}}  + n - a, \,{\displaystyle \frac {2
}{3}}  - a], \,[{\displaystyle \frac {7}{6}}  - 2\,a + n, \,
{\displaystyle \frac {4}{3}} ], \,1)
\]
}
\end{maplelatex}

\begin{maplelatex}
\mapleinline{inert}{2d}{431, _3F_2([1/2+n-a, 2*a+1-n, a+2/3],[5/3, 7/6+a],1);}{%
\[
431: \,\mathrm{_3F_2}([{\displaystyle \frac {1}{2}}  + n - a, \,2\,
a + 1 - n, \,a + {\displaystyle \frac {2}{3}} ], \,[
{\displaystyle \frac {5}{3}} , \,{\displaystyle \frac {7}{6}}  + 
a], \,1)
\]
}
\end{maplelatex}

\begin{maplelatex}
\mapleinline{inert}{2d}{432, _3F_2([2/3+n-2*a, a+2/3, 7/6+a-n],[5/3, 4/3],1);}{%
\[
432: \,\mathrm{_3F_2}([{\displaystyle \frac {2}{3}}  + n - 2\,a, \,
a + {\displaystyle \frac {2}{3}} , \,{\displaystyle \frac {7}{6}
}  + a - n], \,[{\displaystyle \frac {5}{3}} , \,{\displaystyle 
\frac {4}{3}} ], \,1)
\]
}
\end{maplelatex}

\begin{maplelatex}
\mapleinline{inert}{2d}{433, _3F_2([1/2, 2*a+1-n, 2*a+2/3-n],[7/6+a, 5/3+a-n],1);}{%
\[
433: \,\mathrm{_3F_2}([{\displaystyle \frac {1}{2}} , \,2\,a + 1 - 
n, \,2\,a + {\displaystyle \frac {2}{3}}  - n], \,[
{\displaystyle \frac {7}{6}}  + a, \,{\displaystyle \frac {5}{3}
}  + a - n], \,1)
\]
}
\end{maplelatex}

\begin{maplelatex}
\mapleinline{inert}{2d}{434, _3F_2([1/6+n-a, a+2/3, 2*a+2/3-n],[7/6+a, 4/3],1);}{%
\[
434: \,\mathrm{_3F_2}([{\displaystyle \frac {1}{6}}  + n - a, \,a
 + {\displaystyle \frac {2}{3}} , \,2\,a + {\displaystyle \frac {
2}{3}}  - n], \,[{\displaystyle \frac {7}{6}}  + a, \,
{\displaystyle \frac {4}{3}} ], \,1)
\]
}
\end{maplelatex}

\begin{maplelatex}
\mapleinline{inert}{2d}{435, _3F_2([-1/6+n-a, 1/2, 1/2+n-a],[5/6+n-2*a, a+5/6],1);}{%
\[
435: \,\mathrm{_3F_2}([ - {\displaystyle \frac {1}{6}}  + n - a, \,
{\displaystyle \frac {1}{2}} , \,{\displaystyle \frac {1}{2}}  + 
n - a], \,[{\displaystyle \frac {5}{6}}  + n - 2\,a, \,a + 
{\displaystyle \frac {5}{6}} ], \,1)
\]
}
\end{maplelatex}

\begin{maplelatex}
\mapleinline{inert}{2d}{436, _3F_2([-1/6+n-a, -a+1/3, 1/3+n-2*a],[5/6+n-2*a, 2/3],1);}{%
\[
436: \,\mathrm{_3F_2}([ - {\displaystyle \frac {1}{6}}  + n - a, \,
 - a + {\displaystyle \frac {1}{3}} , \,{\displaystyle \frac {1}{
3}}  + n - 2\,a], \,[{\displaystyle \frac {5}{6}}  + n - 2\,a, \,
{\displaystyle \frac {2}{3}} ], \,1)
\]
}
\end{maplelatex}

\begin{maplelatex}
\mapleinline{inert}{2d}{437, _3F_2([1/2, -a+1/3, -a+1],[5/6+n-2*a, 4/3+a-n],1);}{%
\[
437: \,\mathrm{_3F_2}([{\displaystyle \frac {1}{2}} , \, - a + 
{\displaystyle \frac {1}{3}} , \, - a + 1], \,[{\displaystyle 
\frac {5}{6}}  + n - 2\,a, \,{\displaystyle \frac {4}{3}}  + a - 
n], \,1)
\]
}
\end{maplelatex}

\begin{maplelatex}
\mapleinline{inert}{2d}{438, _3F_2([1/2+n-a, 1/3+n-2*a, -a+1],[5/6+n-2*a, 4/3],1);}{%
\[
438: \,\mathrm{_3F_2}([{\displaystyle \frac {1}{2}}  + n - a, \,
{\displaystyle \frac {1}{3}}  + n - 2\,a, \, - a + 1], \,[
{\displaystyle \frac {5}{6}}  + n - 2\,a, \,{\displaystyle 
\frac {4}{3}} ], \,1)
\]
}
\end{maplelatex}

\begin{maplelatex}
\mapleinline{inert}{2d}{439, _3F_2([-1/6+n-a, 2*a+1/3-n, a+1/3],[a+5/6, 2/3],1);}{%
\[
439: \,\mathrm{_3F_2}([ - {\displaystyle \frac {1}{6}}  + n - a, \,
2\,a + {\displaystyle \frac {1}{3}}  - n, \,a + {\displaystyle 
\frac {1}{3}} ], \,[a + {\displaystyle \frac {5}{6}} , \,
{\displaystyle \frac {2}{3}} ], \,1)
\]
}
\end{maplelatex}

\begin{maplelatex}
\mapleinline{inert}{2d}{440, _3F_2([1/2, 2*a+1/3-n, 2*a+1-n],[a+5/6, 4/3+a-n],1);}{%
\[
440: \,\mathrm{_3F_2}([{\displaystyle \frac {1}{2}} , \,2\,a + 
{\displaystyle \frac {1}{3}}  - n, \,2\,a + 1 - n], \,[a + 
{\displaystyle \frac {5}{6}} , \,{\displaystyle \frac {4}{3}}  + 
a - n], \,1)
\]
}
\end{maplelatex}

\begin{maplelatex}
\mapleinline{inert}{2d}{441, _3F_2([1/2+n-a, a+1/3, 2*a+1-n],[a+5/6, 4/3],1);}{%
\[
441: \,\mathrm{_3F_2}([{\displaystyle \frac {1}{2}}  + n - a, \,a
 + {\displaystyle \frac {1}{3}} , \,2\,a + 1 - n], \,[a + 
{\displaystyle \frac {5}{6}} , \,{\displaystyle \frac {4}{3}} ], 
\,1)
\]
}
\end{maplelatex}

\begin{maplelatex}
\mapleinline{inert}{2d}{442, _3F_2([-a+1/3, 2*a+1/3-n, 5/6+a-n],[2/3, 4/3+a-n],1);}{%
\[
442: \,\mathrm{_3F_2}([ - a + {\displaystyle \frac {1}{3}} , \,2\,a
 + {\displaystyle \frac {1}{3}}  - n, \,{\displaystyle \frac {5}{
6}}  + a - n], \,[{\displaystyle \frac {2}{3}} , \,
{\displaystyle \frac {4}{3}}  + a - n], \,1)
\]
}
\end{maplelatex}

\begin{maplelatex}
\mapleinline{inert}{2d}{443, _3F_2([1/3+n-2*a, a+1/3, 5/6+a-n],[2/3, 4/3],1);}{%
\[
443: \,\mathrm{_3F_2}([{\displaystyle \frac {1}{3}}  + n - 2\,a, \,
a + {\displaystyle \frac {1}{3}} , \,{\displaystyle \frac {5}{6}
}  + a - n], \,[{\displaystyle \frac {2}{3}} , \,{\displaystyle 
\frac {4}{3}} ], \,1)
\]
}
\end{maplelatex}

\begin{maplelatex}
\mapleinline{inert}{2d}{444, _3F_2([n-2*a, -a+1/3, a-n+1/2],[1/3, 5/6-a],1);}{%
\[
444: \,\mathrm{_3F_2}([n - 2\,a, \, - a + {\displaystyle \frac {1}{
3}} , \,a - n + {\displaystyle \frac {1}{2}} ], \,[
{\displaystyle \frac {1}{3}} , \,{\displaystyle \frac {5}{6}}  - 
a], \,1)
\]
}
\end{maplelatex}

\begin{maplelatex}
\mapleinline{inert}{2d}{445, _3F_2([-a+1/3, 5/6+a-n, 1/3+n-2*a],[2/3, 5/6-a],1);}{%
\[
445: \,\mathrm{_3F_2}([ - a + {\displaystyle \frac {1}{3}} , \,
{\displaystyle \frac {5}{6}}  + a - n, \,{\displaystyle \frac {1
}{3}}  + n - 2\,a], \,[{\displaystyle \frac {2}{3}} , \,
{\displaystyle \frac {5}{6}}  - a], \,1)
\]
}
\end{maplelatex}

\begin{maplelatex}
\mapleinline{inert}{2d}{446, _3F_2([1/2, 5/6+a-n, a-n+1/2],[5/6-a, 2*a+5/6-n],1);}{%
\[
446: \,\mathrm{_3F_2}([{\displaystyle \frac {1}{2}} , \,
{\displaystyle \frac {5}{6}}  + a - n, \,a - n + {\displaystyle 
\frac {1}{2}} ], \,[{\displaystyle \frac {5}{6}}  - a, \,2\,a + 
{\displaystyle \frac {5}{6}}  - n], \,1)
\]
}
\end{maplelatex}

\begin{maplelatex}
\mapleinline{inert}{2d}{447, _3F_2([1/2, -a+1/3, 2*a+1/3-n],[a-n+1, 4/3+a-n],1);}{%
\[
447: \,\mathrm{_3F_2}([{\displaystyle \frac {1}{2}} , \, - a + 
{\displaystyle \frac {1}{3}} , \,2\,a + {\displaystyle \frac {1}{
3}}  - n], \,[a - n + 1, \,{\displaystyle \frac {4}{3}}  + a - n]
, \,1)
\]
}
\end{maplelatex}

\begin{maplelatex}
\mapleinline{inert}{2d}{448, _3F_2([-a+1/3, 2*a+1/3-n, -1/6+n-a],[1/3, 2/3],1);}{%
\[
448: \,\mathrm{_3F_2}([ - a + {\displaystyle \frac {1}{3}} , \,2\,a
 + {\displaystyle \frac {1}{3}}  - n, \, - {\displaystyle \frac {
1}{6}}  + n - a], \,[{\displaystyle \frac {1}{3}} , \,
{\displaystyle \frac {2}{3}} ], \,1)
\]
}
\end{maplelatex}

\begin{maplelatex}
\mapleinline{inert}{2d}{449, _3F_2([-a+1, -a+1/3, 5/6+a-n],[3/2-n, 4/3+a-n],1);}{%
\[
449: \,\mathrm{_3F_2}([ - a + 1, \, - a + {\displaystyle \frac {1}{
3}} , \,{\displaystyle \frac {5}{6}}  + a - n], \,[
{\displaystyle \frac {3}{2}}  - n, \,{\displaystyle \frac {4}{3}
}  + a - n], \,1)
\]
}
\end{maplelatex}

\begin{maplelatex}
\mapleinline{inert}{2d}{450, _3F_2([2/3-a, -a+1/3, a-n+1/2],[a-n+1, 3/2-n],1);}{%
\[
450: \,\mathrm{_3F_2}([{\displaystyle \frac {2}{3}}  - a, \, - a + 
{\displaystyle \frac {1}{3}} , \,a - n + {\displaystyle \frac {1
}{2}} ], \,[a - n + 1, \,{\displaystyle \frac {3}{2}}  - n], \,1)
\]
}
\end{maplelatex}

\begin{maplelatex}
\mapleinline{inert}{2d}{451, _3F_2([2*a+1-n, 5/6+a-n, 2*a+1/3-n],[3*a-2*n+3/2, 4/3+a-n],1);}{%
\[
451: \,\mathrm{_3F_2}([2\,a + 1 - n, \,{\displaystyle \frac {5}{6}
}  + a - n, \,2\,a + {\displaystyle \frac {1}{3}}  - n], \,[3\,a
 - 2\,n + {\displaystyle \frac {3}{2}} , \,{\displaystyle \frac {
4}{3}}  + a - n], \,1)
\]
}
\end{maplelatex}

\begin{maplelatex}
\mapleinline{inert}{2d}{452, _3F_2([7/6+a-n, 5/6+a-n, a-n+1/2],[3/2-n, 3*a-2*n+3/2],1);}{%
\[
452: \,\mathrm{_3F_2}([{\displaystyle \frac {7}{6}}  + a - n, \,
{\displaystyle \frac {5}{6}}  + a - n, \,a - n + {\displaystyle 
\frac {1}{2}} ], \,[{\displaystyle \frac {3}{2}}  - n, \,3\,a - 2
\,n + {\displaystyle \frac {3}{2}} ], \,1)
\]
}
\end{maplelatex}

\begin{maplelatex}
\mapleinline{inert}{2d}{453, _3F_2([2*a+2/3-n, a-n+1/2, 2*a+1/3-n],[a-n+1, 3*a-2*n+3/2],1);}{%
\[
453: \,\mathrm{_3F_2}([2\,a + {\displaystyle \frac {2}{3}}  - n, \,
a - n + {\displaystyle \frac {1}{2}} , \,2\,a + {\displaystyle 
\frac {1}{3}}  - n], \,[a - n + 1, \,3\,a - 2\,n + 
{\displaystyle \frac {3}{2}} ], \,1)
\]
}
\end{maplelatex}

\begin{maplelatex}
\mapleinline{inert}{2d}{454, _3F_2([1/2, 2/3-a, 2*a+2/3-n],[a-n+1, 5/3+a-n],1);}{%
\[
454: \,\mathrm{_3F_2}([{\displaystyle \frac {1}{2}} , \,
{\displaystyle \frac {2}{3}}  - a, \,2\,a + {\displaystyle 
\frac {2}{3}}  - n], \,[a - n + 1, \,{\displaystyle \frac {5}{3}
}  + a - n], \,1)
\]
}
\end{maplelatex}

\begin{maplelatex}
\mapleinline{inert}{2d}{455, _3F_2([a, a-n+1/2, 2*a+1/3-n],[1/3, 2*a+5/6-n],1);}{%
\[
455: \,\mathrm{_3F_2}([a, \,a - n + {\displaystyle \frac {1}{2}} , 
\,2\,a + {\displaystyle \frac {1}{3}}  - n], \,[{\displaystyle 
\frac {1}{3}} , \,2\,a + {\displaystyle \frac {5}{6}}  - n], \,1)
\]
}
\end{maplelatex}

\begin{maplelatex}
\mapleinline{inert}{2d}{456, _3F_2([a+1/3, 1/3+n-2*a, -1/6+n-a],[1/3+n-a, 2/3],1);}{%
\[
456: \,\mathrm{_3F_2}([a + {\displaystyle \frac {1}{3}} , \,
{\displaystyle \frac {1}{3}}  + n - 2\,a, \, - {\displaystyle 
\frac {1}{6}}  + n - a], \,[{\displaystyle \frac {1}{3}}  + n - a
, \,{\displaystyle \frac {2}{3}} ], \,1)
\]
}
\end{maplelatex}

\begin{maplelatex}
\mapleinline{inert}{2d}{457, _3F_2([a+1/3, 1/2, a],[1/3+n-a, 2*a+5/6-n],1);}{%
\[
457: \,\mathrm{_3F_2}([a + {\displaystyle \frac {1}{3}} , \,
{\displaystyle \frac {1}{2}} , \,a], \,[{\displaystyle \frac {1}{
3}}  + n - a, \,2\,a + {\displaystyle \frac {5}{6}}  - n], \,1)
\]
}
\end{maplelatex}

\begin{maplelatex}
\mapleinline{inert}{2d}{458, _3F_2([1/3+n-2*a, 1/2, n-2*a],[1/3+n-a, 5/6-a],1);}{%
\[
458: \,\mathrm{_3F_2}([{\displaystyle \frac {1}{3}}  + n - 2\,a, \,
{\displaystyle \frac {1}{2}} , \,n - 2\,a], \,[{\displaystyle 
\frac {1}{3}}  + n - a, \,{\displaystyle \frac {5}{6}}  - a], \,1
)
\]
}
\end{maplelatex}

\begin{maplelatex}
\mapleinline{inert}{2d}{459, _3F_2([-1/6+n-a, a, n-2*a],[1/3+n-a, 1/3],1);}{%
\[
459: \,\mathrm{_3F_2}([ - {\displaystyle \frac {1}{6}}  + n - a, \,
a, \,n - 2\,a], \,[{\displaystyle \frac {1}{3}}  + n - a, \,
{\displaystyle \frac {1}{3}} ], \,1)
\]
}
\end{maplelatex}

\begin{maplelatex}
\mapleinline{inert}{2d}{460, _3F_2([a+1/3, 5/6+a-n, 2*a+1/3-n],[2/3, 2*a+5/6-n],1);}{%
\[
460: \,\mathrm{_3F_2}([a + {\displaystyle \frac {1}{3}} , \,
{\displaystyle \frac {5}{6}}  + a - n, \,2\,a + {\displaystyle 
\frac {1}{3}}  - n], \,[{\displaystyle \frac {2}{3}} , \,2\,a + 
{\displaystyle \frac {5}{6}}  - n], \,1)
\]
}
\end{maplelatex}

\begin{maplelatex}
\mapleinline{inert}{2d}{461, _3F_2([7/6+a-n, -a+1, 2/3-a],[5/3+a-n, 3/2-n],1);}{%
\[
461: \,\mathrm{_3F_2}([{\displaystyle \frac {7}{6}}  + a - n, \, - 
a + 1, \,{\displaystyle \frac {2}{3}}  - a], \,[{\displaystyle 
\frac {5}{3}}  + a - n, \,{\displaystyle \frac {3}{2}}  - n], \,1
)
\]
}
\end{maplelatex}

\begin{maplelatex}
\mapleinline{inert}{2d}{462, _3F_2([7/6+a-n, 2*a+1-n, 2*a+2/3-n],[5/3+a-n, 3*a-2*n+3/2],1);}{%
\[
462: \,\mathrm{_3F_2}([{\displaystyle \frac {7}{6}}  + a - n, \,2\,
a + 1 - n, \,2\,a + {\displaystyle \frac {2}{3}}  - n], \,[
{\displaystyle \frac {5}{3}}  + a - n, \,3\,a - 2\,n + 
{\displaystyle \frac {3}{2}} ], \,1)
\]
}
\end{maplelatex}

\begin{maplelatex}
\mapleinline{inert}{2d}{463, _3F_2([-a+1, 2*a+1-n, 1/2],[5/3+a-n, 4/3+a-n],1);}{%
\[
463: \,\mathrm{_3F_2}([ - a + 1, \,2\,a + 1 - n, \,{\displaystyle 
\frac {1}{2}} ], \,[{\displaystyle \frac {5}{3}}  + a - n, \,
{\displaystyle \frac {4}{3}}  + a - n], \,1)
\]
}
\end{maplelatex}

\begin{maplelatex}
\mapleinline{inert}{2d}{469, _3F_2([c, b, (a+c)*(a+b)/a],[1+b+a+c, (a+a*c+a*b+b*c)/a],1);}{%
\[
469: \,\mathrm{F32}([c, \,b, \,{\displaystyle \frac {(a + c)\,(a
 + b)}{a}} ], \,[1 + b + a + c, \,{\displaystyle \frac {a + a\,c
 + a\,b + b\,c}{a}} ], \,1)
\]
}
\end{maplelatex}
\end{flushleft}
\end{maplegroup}

%% file: AppendixB1to30.tex
\begin{maplegroup}
\begin{flushleft}
\mapleresult
\begin{maplelatex}
\mapleinline{inert}{2d}{1, ":   ", WPPR(m,n,a,b,c);}{%
\[
1\mbox{:~~~}     \mathrm{WPPR}(m, \,n, \,a, \,b, \,c)
\]
}
\end{maplelatex}

\begin{maplelatex}
\mapleinline{inert}{2d}{2, ":   ", WPPR(2*m,n,a,b,c);}{%
\[
2\mbox{:~~~}     \mathrm{WPPR}(2\,m, \,n, \,a, \,b, \,c)
\]
}
\end{maplelatex}

\begin{maplelatex}
\mapleinline{inert}{2d}{3, ":   ", WPPR(2*m+1,n,a,b,c);}{%
\[
3\mbox{:~~~}     \mathrm{WPPR}(2\,m + 1, \,n, \,a, \,b, \,
c)
\]
}
\end{maplelatex}

\begin{maplelatex}
\mapleinline{inert}{2d}{4, ":   ", WPMR(m,n,a,b,c);}{%
\[
4\mbox{:~~~}     \mathrm{WPMR}(m, \,n, \,a, \,b, \,c)
\]
}
\end{maplelatex}

\begin{maplelatex}
\mapleinline{inert}{2d}{5, ":   ", WPMR(2*m,n,a,b,c);}{%
\[
5\mbox{:~~~}     \mathrm{WPMR}(2\,m, \,n, \,a, \,b, \,c)
\]
}
\end{maplelatex}

\begin{maplelatex}
\mapleinline{inert}{2d}{6, ":   ", WPMR(2*m+1,n,a,b,c);}{%
\[
6\mbox{:~~~}     \mathrm{WPMR}(2\,m + 1, \,n, \,a, \,b, \,
c)
\]
}
\end{maplelatex}

\begin{maplelatex}
\mapleinline{inert}{2d}{7, ":   ", WMPR(m,n,a,b,c);}{%
\[
7\mbox{:~~~}     \mathrm{WMPR}(m, \,n, \,a, \,b, \,c)
\]
}
\end{maplelatex}

\begin{maplelatex}
\mapleinline{inert}{2d}{8, ":   ", WMPR(2*m,n,a,b,c);}{%
\[
8\mbox{:~~~}     \mathrm{WMPR}(2\,m, \,n, \,a, \,b, \,c)
\]
}
\end{maplelatex}

\begin{maplelatex}
\mapleinline{inert}{2d}{9, ":   ", WMPR(2*m+1,n,a,b,c);}{%
\[
9\mbox{:~~~}     \mathrm{WMPR}(2\,m + 1, \,n, \,a, \,b, \,
c)
\]
}
\end{maplelatex}

\begin{maplelatex}
\mapleinline{inert}{2d}{10, ":   ", WMMR(m,n,a,b,c);}{%
\[
10\mbox{:~~~}     \mathrm{WMMR}(m, \,n, \,a, \,b, \,c)
\]
}
\end{maplelatex}

\begin{maplelatex}
\mapleinline{inert}{2d}{11, ":   ", WMMR(2*m,n,a,b,c);}{%
\[
11\mbox{:~~~}     \mathrm{WMMR}(2\,m, \,n, \,a, \,b, \,c)
\]
}
\end{maplelatex}

\begin{maplelatex}
\mapleinline{inert}{2d}{12, ":   ", WMMR(2*m+1,n,a,b,c);}{%
\[
12\mbox{:~~~}     \mathrm{WMMR}(2\,m + 1, \,n, \,a, \,b, 
\,c)
\]
}
\end{maplelatex}

\begin{maplelatex}
\mapleinline{inert}{2d}{13, ":   ", WhipplePP(m,n,a,b,c);}{%
\[
13\mbox{:~~~}     \mathrm{WhipplePP}(m, \,n, \,a, \,b, \,c
)
\]
}
\end{maplelatex}

\begin{maplelatex}
\mapleinline{inert}{2d}{14, ":   ", WhipplePM(m,n,a,b,c);}{%
\[
14\mbox{:~~~}     \mathrm{WhipplePM}(m, \,n, \,a, \,b, \,c
)
\]
}
\end{maplelatex}

\begin{maplelatex}
\mapleinline{inert}{2d}{15, ":   ", WhippleMP(m,n,a,b,c);}{%
\[
15\mbox{:~~~}     \mathrm{WhippleMP}(m, \,n, \,a, \,b, \,c
)
\]
}
\end{maplelatex}

\begin{maplelatex}
\mapleinline{inert}{2d}{16, ":   ", WhippleMM(m,n,a,b,c);}{%
\[
16\mbox{:~~~}     \mathrm{WhippleMM}(m, \,n, \,a, \,b, \,c
)
\]
}
\end{maplelatex}

\begin{maplelatex}
\mapleinline{inert}{2d}{17, ":   ", DixonPP(m,n,a,b,c);}{%
\[
17\mbox{:~~~}     \mathrm{DixonPP}(m, \,n, \,a, \,b, \,c)
\]
}
\end{maplelatex}

\begin{maplelatex}
\mapleinline{inert}{2d}{18, ":   ", DixonPM(m,n,a,b,c);}{%
\[
18\mbox{:~~~}     \mathrm{DixonPM}(m, \,n, \,a, \,b, \,c)
\]
}
\end{maplelatex}

\begin{maplelatex}
\mapleinline{inert}{2d}{19, ":   ", DixonMM(m,n,a,b,c);}{%
\[
19\mbox{:~~~}     \mathrm{DixonMM}(m, \,n, \,a, \,b, \,c)
\]
}
\end{maplelatex}

\begin{maplelatex}
\mapleinline{inert}{2d}{25, ":   ",
(-b+c-a+n)*GAMMA(-b+1)*GAMMA(n+1)*GAMMA(c)/GAMMA(a)/GAMMA(-b+c-a+n+1)*
GAMMA(-b+c-a+n)*Sum(GAMMA(-k+a-1)*(-1)^k/GAMMA(c-1-k)/GAMMA(n-k)/GAMMA
(-b+2+k),k = 0 ..
n-1)+GAMMA(a-n)/GAMMA(b)*GAMMA(b-n)/GAMMA(c-b)*GAMMA(n+1)*GAMMA(c)/GAM
MA(a)/GAMMA(c-a)*GAMMA(-b+c-a+n);}{%
\maplemultiline{
25  \mbox{:~~~} ( - b + c - a + n)\,\Gamma ( - b + 1)\,
\Gamma (n + 1)\,\Gamma (c)\,\Gamma ( - b + c - a + n) \\
 \left(  \! {\displaystyle \sum _{k=0}^{n - 1}} \,{\displaystyle 
\frac {\Gamma ( - k + a - 1)\,(-1)^{k}}{\Gamma (c - 1 - k)\,
\Gamma (n - k)\,\Gamma ( - b + 2 + k)}}  \!  \right) /(\Gamma (a)
\,\Gamma ( - b + c - a + n + 1)) \\
\mbox{} + {\displaystyle \frac {\Gamma (a - n)\,\Gamma (b - n)\,
\Gamma (n + 1)\,\Gamma (c)\,\Gamma ( - b + c - a + n)}{\Gamma (b)
\,\Gamma (c - b)\,\Gamma (a)\,\Gamma (c - a)}}  }
}
\end{maplelatex}

\begin{maplelatex}
\mapleinline{inert}{2d}{26, ":   ",
(c-a)*GAMMA(b-n)*Sum(GAMMA(-k+a-1)*(-1)^k/GAMMA(-b+c-k)/GAMMA(n-k)/GAM
MA(b+1+k-n),k = 0 ..
n-1)*GAMMA(n+1)*GAMMA(c)/GAMMA(a)/GAMMA(c+1-a)*GAMMA(-b+c-a+1)+GAMMA(a
-n)/GAMMA(n+1-b)*GAMMA(-b+1)/GAMMA(c-n)*GAMMA(n+1)*GAMMA(c)/GAMMA(a);}
{%
\maplemultiline{
26  \mbox{:~~~} (c - a)\,\Gamma (b - n)\, \left(  \! 
{\displaystyle \sum _{k=0}^{n - 1}} \,{\displaystyle \frac {
\Gamma ( - k + a - 1)\,(-1)^{k}}{\Gamma ( - b + c - k)\,\Gamma (n
 - k)\,\Gamma (b + 1 + k - n)}}  \!  \right) \,\Gamma (n + 1)\,
\Gamma (c) \\
\Gamma ( - b + c - a + 1)/(\Gamma (a)\,\Gamma (c + 1 - a))\mbox{}
 + {\displaystyle \frac {\Gamma (a - n)\,\Gamma ( - b + 1)\,
\Gamma (n + 1)\,\Gamma (c)}{\Gamma (n + 1 - b)\,\Gamma (c - n)\,
\Gamma (a)}}  }
}
\end{maplelatex}

\begin{maplelatex}
\mapleinline{inert}{2d}{27, ":   ",
(b-a+n)*GAMMA(b-c)*Sum(GAMMA(-k+a-1)*(-1)^k/GAMMA(c-k)/GAMMA(n-k)/GAMM
A(b-c+1+k),k = 0 ..
n-1)*GAMMA(n+1-a)*GAMMA(c+1)*GAMMA(b+n)/GAMMA(a)/GAMMA(b-a+n+1)+GAMMA(
a-n)/GAMMA(c+1-b)*GAMMA(c-b+1-n)/GAMMA(b)*GAMMA(n+1-a)*GAMMA(c+1)*GAMM
A(b+n)/GAMMA(a)/GAMMA(c+1-a);}{%
\maplemultiline{
27  \mbox{:~~~} (b - a + n)\,\Gamma (b - c)\, \left(  \! 
{\displaystyle \sum _{k=0}^{n - 1}} \,{\displaystyle \frac {
\Gamma ( - k + a - 1)\,(-1)^{k}}{\Gamma (c - k)\,\Gamma (n - k)\,
\Gamma (b - c + 1 + k)}}  \!  \right) \,\Gamma (n + 1 - a) \\
\Gamma (c + 1)\,\Gamma (b + n)/(\Gamma (a)\,\Gamma (b - a + n + 1
)) \\
\mbox{} + {\displaystyle \frac {\Gamma (a - n)\,\Gamma (c - b + 1
 - n)\,\Gamma (n + 1 - a)\,\Gamma (c + 1)\,\Gamma (b + n)}{\Gamma
 (c + 1 - b)\,\Gamma (b)\,\Gamma (a)\,\Gamma (c + 1 - a)}}  }
}
\end{maplelatex}

\begin{maplelatex}
\mapleinline{inert}{2d}{28, ":   ",
(-c+b+a+n)*GAMMA(1+a-c)*Sum(GAMMA(-k+c-b-1)*(-1)^k/GAMMA(c-1-k)/GAMMA(
n-k)/GAMMA(-c+a+2+k),k = 0 ..
n-1)*GAMMA(n)*GAMMA(c)/GAMMA(c-b)+GAMMA(-c+b+a+n+1)*GAMMA(-b+c-n)/GAMM
A(c-a)*GAMMA(c-a-n)/GAMMA(a)*GAMMA(n)*GAMMA(c)/GAMMA(c-b)/GAMMA(b);}{%
\maplemultiline{
28  \mbox{:~~~} {\displaystyle \frac {( - c + b + a + n)\,
\Gamma (1 + a - c)\, \left(  \! {\displaystyle \sum _{k=0}^{n - 1
}} \,{\displaystyle \frac {\Gamma ( - k + c - b - 1)\,(-1)^{k}}{
\Gamma (c - 1 - k)\,\Gamma (n - k)\,\Gamma ( - c + a + 2 + k)}} 
 \!  \right) \,\Gamma (n)\,\Gamma (c)}{\Gamma (c - b)}}  \\
\mbox{} + {\displaystyle \frac {\Gamma ( - c + b + a + n + 1)\,
\Gamma ( - b + c - n)\,\Gamma (c - a - n)\,\Gamma (n)\,\Gamma (c)
}{\Gamma (c - a)\,\Gamma (a)\,\Gamma (c - b)\,\Gamma (b)}}  }
}
\end{maplelatex}

\begin{maplelatex}
\mapleinline{inert}{2d}{29, ":   ",
b*GAMMA(-c+b+a+n+1)*GAMMA(c-a-n)*Sum(GAMMA(-k+c-b-1)*(-1)^k/GAMMA(a-k)
/GAMMA(n-k)/GAMMA(c+1-a-n+k),k = 0 ..
n-1)*GAMMA(c)/GAMMA(c-b)/GAMMA(b+1)+GAMMA(-c+b+a+n+1)*GAMMA(-b+c-n)/GA
MMA(1+a-c+n)*GAMMA(1+a-c)/GAMMA(c-n)*GAMMA(c)/GAMMA(c-b)/GAMMA(-c+b+1+
a);}{%
\maplemultiline{
29  \mbox{:~~~} {\displaystyle \frac {b\,\Gamma ( - c + b
 + a + n + 1)\,\Gamma (c - a - n)\, \left(  \! {\displaystyle 
\sum _{k=0}^{n - 1}} \,{\displaystyle \frac {\Gamma ( - k + c - b
 - 1)\,(-1)^{k}}{\Gamma (a - k)\,\Gamma (n - k)\,\Gamma (c + 1 - 
a - n + k)}}  \!  \right) \,\Gamma (c)}{\Gamma (c - b)\,\Gamma (b
 + 1)}}  \\
\mbox{} + {\displaystyle \frac {\Gamma ( - c + b + a + n + 1)\,
\Gamma ( - b + c - n)\,\Gamma (1 + a - c)\,\Gamma (c)}{\Gamma (1
 + a - c + n)\,\Gamma (c - n)\,\Gamma (c - b)\,\Gamma ( - c + b
 + 1 + a)}}  }
}
\end{maplelatex}

\begin{maplelatex}
\mapleinline{inert}{2d}{30, ":   ",
(1+b-a-n)*GAMMA(-c+1)*GAMMA(n)/(c-1+n)/GAMMA(n-b+a)/GAMMA(-c-n+1)*GAMM
A(a-b)/GAMMA(1+b-c)/GAMMA(b)*Sum(GAMMA(-b+a+k)/GAMMA(a-b)*GAMMA(-c-n+b
+1+k)/GAMMA(b-c-n+1)/GAMMA(k+1)/GAMMA(-c-n+2+k)*GAMMA(-c-n+2),k = 0 ..
n-1)*GAMMA(b-c-n+2-a)*GAMMA(b+1)*GAMMA(-n+b+1)/GAMMA(2+b-a-n);}{%
\maplemultiline{
30  \mbox{:~~~} (1 + b - a - n)\,\Gamma ( - c + 1)\,\Gamma
 (n)\,\Gamma (a - b) \\
 \left(  \! {\displaystyle \sum _{k=0}^{n - 1}} \,{\displaystyle 
\frac {\Gamma ( - b + a + k)\,\Gamma ( - c - n + b + 1 + k)\,
\Gamma ( - c - n + 2)}{\Gamma (a - b)\,\Gamma (b - c - n + 1)\,
\Gamma (k + 1)\,\Gamma ( - c - n + 2 + k)}}  \!  \right) \,\Gamma
 (b - c - n + 2 - a) \\
\Gamma (b + 1)\,\Gamma ( - n + b + 1)/((c - 1 + n)\,\Gamma (n - b
 + a)\,\Gamma ( - c - n + 1)\,\Gamma (1 + b - c)\,\Gamma (b) \Gamma (2 + b - a - n)) }
}
\end{maplelatex}
\end{flushleft}
\end{maplegroup}

%% file: AppendixB31to60.tex
\begin{maplegroup}
\mapleresult
\begin{maplelatex}
\mapleinline{inert}{2d}{31, ":   ",
(c-a-n)*GAMMA(1+b-c)*GAMMA(n)/(c-1-b+n)/GAMMA(1+a-c+n)/GAMMA(b-c-n+1)*
GAMMA(1+a-c)/GAMMA(b)/GAMMA(c-1)*Sum(GAMMA(1+a-c+k)/GAMMA(1+a-c)*GAMMA
(b+k-n)/GAMMA(b-n)/GAMMA(k+1)/GAMMA(2+b-n-c+k)*GAMMA(2+b-n-c),k = 0 ..
n-1)*GAMMA(c-a-n)*GAMMA(c)*GAMMA(-n+b+1)/GAMMA(-a+1+c-n);}{%
\maplemultiline{
31\mbox{:~~~} (c - a - n)\,\Gamma (1 + b - c)\,\Gamma (n
)\,\Gamma (1 + a - c) \\
 \left(  \! {\displaystyle \sum _{k=0}^{n - 1}} \,{\displaystyle 
\frac {\Gamma (1 + a - c + k)\,\Gamma (b + k - n)\,\Gamma (2 + b
 - n - c)}{\Gamma (1 + a - c)\,\Gamma (b - n)\,\Gamma (k + 1)\,
\Gamma (2 + b - n - c + k)}}  \!  \right) \,\Gamma (c - a - n)\,
\Gamma (c) \\
\Gamma ( - n + b + 1)/((c - 1 - b + n)\,\Gamma (1 + a - c + n)\,
\Gamma (b - c - n + 1)\,\Gamma (b)\,\Gamma (c - 1) 
\Gamma ( - a + 1 + c - n)) }
}
\end{maplelatex}

\begin{maplelatex}
\mapleinline{inert}{2d}{33, ":   ",
(c-a)*GAMMA(1+b-c)*GAMMA(n)/(c-1-b+n)/GAMMA(1+a-c)/GAMMA(b-c-n+1)*GAMM
A(1+a-c-n)/GAMMA(b+n)/GAMMA(c-1+n)*Sum(GAMMA(a+k+1-n-c)/GAMMA(1+a-c-n)
*GAMMA(b+k)/GAMMA(b)/GAMMA(k+1)/GAMMA(2+b-n-c+k)*GAMMA(2+b-n-c),k = 0
.. n-1)*GAMMA(c-a+n)*GAMMA(c)*GAMMA(b+1)/GAMMA(c+1-a);}{%
\maplemultiline{
33\mbox{:~~~} (c - a)\,\Gamma (1 + b - c)\,\Gamma (n)\,
\Gamma (1 + a - c - n) \\
 \left(  \! {\displaystyle \sum _{k=0}^{n - 1}} \,{\displaystyle 
\frac {\Gamma (a + k + 1 - n - c)\,\Gamma (b + k)\,\Gamma (2 + b
 - n - c)}{\Gamma (1 + a - c - n)\,\Gamma (b)\,\Gamma (k + 1)\,
\Gamma (2 + b - n - c + k)}}  \!  \right) \,\Gamma (c - a + n)\,
\Gamma (c) \\
\Gamma (b + 1)/((c - 1 - b + n)\,\Gamma (1 + a - c)\,\Gamma (b - 
c - n + 1)\,\Gamma (b + n)\,\Gamma (c - 1 + n) 
\Gamma (c + 1 - a)) }
}
\end{maplelatex}

\begin{maplelatex}
\mapleinline{inert}{2d}{34, ":   ",
c*GAMMA(b-a+n)*GAMMA(n)/(a-b)/GAMMA(-c+1)/GAMMA(b-a)*GAMMA(-c-n+1)/GAM
MA(b+n)*GAMMA(c+b+n)/GAMMA(a)*Sum(GAMMA(1-c-n+k)/GAMMA(-c-n+1)*GAMMA(b
+k)/GAMMA(b)/GAMMA(k+1)/GAMMA(b-a+1+k)*GAMMA(b-a+1),k = 0 ..
n-1)*GAMMA(a-n+1)/GAMMA(c+1);}{%
\maplemultiline{
34\mbox{:~~~} c\,\Gamma (b - a + n)\,\Gamma (n)\,\Gamma 
( - c - n + 1)\,\Gamma (c + b + n) \\
 \left(  \! {\displaystyle \sum _{k=0}^{n - 1}} \,{\displaystyle 
\frac {\Gamma (1 - c - n + k)\,\Gamma (b + k)\,\Gamma (b - a + 1)
}{\Gamma ( - c - n + 1)\,\Gamma (b)\,\Gamma (k + 1)\,\Gamma (b - 
a + 1 + k)}}  \!  \right) \,\Gamma (a - n + 1)/((a - b) \\
\Gamma ( - c + 1)\,\Gamma (b - a)\,\Gamma (b + n)\,\Gamma (a)\,
\Gamma (c + 1)) }
}
\end{maplelatex}

\begin{maplelatex}
\mapleinline{inert}{2d}{35, ":   ",
(c-a)*GAMMA(b-c+a)*GAMMA(n)/(-b+c-a+n)/GAMMA(1+a-c)/GAMMA(a-c+b-n)*GAM
MA(1+a-c-n)/GAMMA(a)*GAMMA(c)/GAMMA(c-b)*Sum(GAMMA(a+k+1-n-c)/GAMMA(1+
a-c-n)*GAMMA(a+k-n)/GAMMA(a-n)/GAMMA(k+1)/GAMMA(b-c+a-n+1+k)*GAMMA(1+a
+b-c-n),k = 0 .. n-1)*GAMMA(-b+c-a+1)*GAMMA(a-n+1)/GAMMA(c+1-a);}{%
\maplemultiline{
35\mbox{:~~~} (c - a)\,\Gamma (b - c + a)\,\Gamma (n)\,
\Gamma (1 + a - c - n)\,\Gamma (c) \\
 \left(  \! {\displaystyle \sum _{k=0}^{n - 1}} \,{\displaystyle 
\frac {\Gamma (a + k + 1 - n - c)\,\Gamma (a + k - n)\,\Gamma (1
 + a + b - c - n)}{\Gamma (1 + a - c - n)\,\Gamma (a - n)\,\Gamma
 (k + 1)\,\Gamma (b - c + a - n + 1 + k)}}  \!  \right)  \\
\Gamma ( - b + c - a + 1)\,\Gamma (a - n + 1)/(( - b + c - a + n)
\,\Gamma (1 + a - c)\,\Gamma (a - c + b - n)\,\Gamma (a) \Gamma (c - b)\,\Gamma (c + 1 - a)) }
}
\end{maplelatex}

\begin{maplelatex}
\mapleinline{inert}{2d}{36, ":   ",
(c-1)*GAMMA(1+a-c)*GAMMA(n)/(c-a+n-1)/GAMMA(2-c)/GAMMA(1+a-c-n)*GAMMA(
-c-n+2)/GAMMA(b-1+n)*Sum(GAMMA(-c-n+2+k)/GAMMA(-c-n+2)*GAMMA(b-1+k)/GA
MMA(b-1)/GAMMA(k+1)/GAMMA(a-c+2-n+k)*GAMMA(a-c+2-n),k = 0 ..
n-1)*GAMMA(b);}{%
\maplemultiline{
36\mbox{:~~~} (c - 1)\,\Gamma (1 + a - c)\,\Gamma (n)\,
\Gamma ( - c - n + 2) 
 \left(  \! {\displaystyle \sum _{k=0}^{n - 1}} \,{\displaystyle 
\frac {\Gamma ( - c - n + 2 + k)\,\Gamma (b - 1 + k)\,\Gamma (a
 - c + 2 - n)}{\Gamma ( - c - n + 2)\,\Gamma (b - 1)\,\Gamma (k
 + 1)\,\Gamma (a - c + 2 - n + k)}}  \!  \right) \,\Gamma (b) \\ 
/((c - a + n - 1) 
\Gamma (2 - c)\,\Gamma (1 + a - c - n)\,\Gamma (b - 1 + n)) }
}
\end{maplelatex}

\begin{maplelatex}
\mapleinline{inert}{2d}{37, ":   ",
-1/2*GAMMA(1/2*a+b+1/2)*GAMMA(n+a-b-1)/GAMMA(a)/GAMMA(b+1)/GAMMA(n+1/2
*a-b-1/2)*Sum(GAMMA(3/2-1/2*a+b-n+k)/GAMMA(-1/2*a+b+3/2-n)*GAMMA(-b+k)
/GAMMA(-b)/GAMMA(1/2+1/2*a-b+k)*GAMMA(1/2*a-b+1/2)/GAMMA(2+b-n+k)*GAMM
A(b-n+2),k = 1 ..
n-2)+Pi^(1/2)*GAMMA(2*n+a-2*b-1)/GAMMA(n+a-2*b-1)*GAMMA(2*b+2)/GAMMA(-
n+2*b+2)/GAMMA(1/2*a+1/2)/(2^(2*n+1))*GAMMA(1/2*a)*GAMMA(1/2*a-b+1/2)*
GAMMA(b-n+2)/GAMMA(n+1/2*a-b)/GAMMA(3/2+b)*GAMMA(1/2*a+b+1/2)*GAMMA(n+
a-b-1)/GAMMA(a)/GAMMA(b+1)/GAMMA(n+1/2*a-b-1/2);}{%
\maplemultiline{
37\mbox{:~~~}  - {\displaystyle \frac {1}{2}} \Gamma (
{\displaystyle \frac {a}{2}}  + b + {\displaystyle \frac {1}{2}} 
)\,\Gamma (n + a - b - 1) 
 \left(  \! {\displaystyle \sum _{k=1}^{n - 2}} \,{\displaystyle 
\frac {\Gamma ({\displaystyle \frac {3}{2}}  - {\displaystyle 
\frac {a}{2}}  + b - n + k)\,\Gamma ( - b + k)\,\Gamma (
{\displaystyle \frac {a}{2}}  - b + {\displaystyle \frac {1}{2}} 
)\,\Gamma (b - n + 2)}{\Gamma ( - {\displaystyle \frac {a}{2}} 
 + b + {\displaystyle \frac {3}{2}}  - n)\,\Gamma ( - b)\,\Gamma 
({\displaystyle \frac {1}{2}}  + {\displaystyle \frac {a}{2}}  - 
b + k)\,\Gamma (2 + b - n + k)}}  \!  \right)  \left/ {\vrule 
height0.80em width0em depth0.80em} \right. \!  \!  (\Gamma (a) 
\Gamma (b + 1)\,  \\ \Gamma (n + {\displaystyle \frac {a}{2}}  - b - 
{\displaystyle \frac {1}{2}} ))\mbox{} +  \sqrt{\pi }\,\Gamma (2\,
n + a - 2\,b - 1)\,\Gamma (2\,b + 2)\,\Gamma ({\displaystyle 
\frac {a}{2}} )\,\Gamma ({\displaystyle \frac {a}{2}}  - b + 
{\displaystyle \frac {1}{2}} ) \Gamma (b - n + 2)\, 
\Gamma ({\displaystyle \frac {a}{2}}  + b + 
{\displaystyle \frac {1}{2}} )\, \Gamma (n + a - b - 1) \\ \left/ 
{\vrule height0.80em width0em depth0.80em} \right. \!  \! (\Gamma
 (n + a - 2\,b - 1) 
\Gamma ( - n + 2\,b + 2)\,\Gamma ({\displaystyle \frac {a}{2}} 
 + {\displaystyle \frac {1}{2}} )\,2^{(2\,n + 1)}\,\Gamma (n + 
{\displaystyle \frac {a}{2}}  - b)\,\Gamma ({\displaystyle 
\frac {3}{2}}  + b)\,\Gamma (a)\,\Gamma (b + 1) 
\Gamma (n + {\displaystyle \frac {a}{2}}  - b - {\displaystyle 
\frac {1}{2}} )) }
}
\end{maplelatex}

\begin{maplelatex}
\mapleinline{inert}{2d}{38, ":   ",
-1/2*GAMMA(a)/GAMMA(n+a-b-1)/GAMMA(2*a-1)*GAMMA(n-b-2+2*a)*Sum(GAMMA(-
b+k)/GAMMA(-b)*GAMMA(-a+b+2+k-n)/GAMMA(2+b-a-n)/GAMMA(2+b-n+k)*GAMMA(b
-n+2)/GAMMA(-b+a+k)*GAMMA(a-b),k = 1 ..
n-2)+Pi^(1/2)*GAMMA(2*b+2)/GAMMA(-n+2*b+2)*GAMMA(2*n+2*a-2*b-2)/GAMMA(
n+2*a-2*b-2)/(2^(2*n+1))*GAMMA(a-1/2)*GAMMA(b-n+2)/GAMMA(3/2+b)/GAMMA(
-1/2+n+a-b)*GAMMA(a-b)/GAMMA(n+a-b-1)/GAMMA(2*a-1)*GAMMA(n-b-2+2*a);}{
\maplemultiline{
38\mbox{:~~~}  - {\displaystyle \frac {1}{2}} \,
{\displaystyle \frac {\Gamma (a)\,\Gamma (n - b - 2 + 2\,a)\,
 \left(  \! {\displaystyle \sum _{k=1}^{n - 2}} \,{\displaystyle 
\frac {\Gamma ( - b + k)\,\Gamma ( - a + b + 2 + k - n)\,\Gamma (
b - n + 2)\,\Gamma (a - b)}{\Gamma ( - b)\,\Gamma (2 + b - a - n)
\,\Gamma (2 + b - n + k)\,\Gamma ( - b + a + k)}}  \!  \right) }{
\Gamma (n + a - b - 1)\,\Gamma (2\,a - 1)}}  \\
\mbox{} + \sqrt{\pi }\,\Gamma (2\,b + 2)\,\Gamma (2\,n + 2\,a - 2
\,b - 2)\,\Gamma (a - {\displaystyle \frac {1}{2}} )\,\Gamma (b
 - n + 2)\,\Gamma (a - b) 
\Gamma (n - b - 2 + 2\,a) \\ \left/ {\vrule 
height0.80em width0em depth0.80em} \right. \!  \! (\Gamma ( - n
 + 2\,b + 2)\,\Gamma (n + 2\,a - 2\,b - 2)\,2^{(2\,n + 1)}\,
\Gamma ({\displaystyle \frac {3}{2}}  + b) 
\Gamma ( - {\displaystyle \frac {1}{2}}  + n + a - b)\,\Gamma (n
 + a - b - 1)\,\Gamma (2\,a - 1)) }
}
\end{maplelatex}

\begin{maplelatex}
\mapleinline{inert}{2d}{39, ":   ",
-1/2*Sum(GAMMA(a+1+k-n)/GAMMA(a-n+1)*GAMMA(b+1+k-n)/GAMMA(-n+b+1)/GAMM
A(1-a+k)*GAMMA(-a+1)/GAMMA(1-b+k)*GAMMA(-b+1),k = 1 ..
n-2)+Pi^(1/2)*GAMMA(2*n-2*a)/GAMMA(n-2*a)*GAMMA(2*n-2*b)/GAMMA(n-2*b)/
GAMMA(n-a-b)/(2^(2*n+1))*GAMMA(n-a-b-1/2)*GAMMA(-a+1)*GAMMA(-b+1)/GAMM
A(1/2+n-a)/GAMMA(1/2+n-b);}{%
\maplemultiline{
39\mbox{:~~~}  - {\displaystyle \frac {1}{2}} \, \left( 
 \! {\displaystyle \sum _{k=1}^{n - 2}} \,{\displaystyle \frac {
\Gamma (a + 1 + k - n)\,\Gamma (b + 1 + k - n)\,\Gamma ( - a + 1)
\,\Gamma ( - b + 1)}{\Gamma (a - n + 1)\,\Gamma ( - n + b + 1)\,
\Gamma (1 - a + k)\,\Gamma (1 - b + k)}}  \!  \right)  \\
\mbox{} + {\displaystyle \frac {\sqrt{\pi }\,\Gamma (2\,n - 2\,a)
\,\Gamma (2\,n - 2\,b)\,\Gamma (n - a - b - {\displaystyle 
\frac {1}{2}} )\,\Gamma ( - a + 1)\,\Gamma ( - b + 1)}{\Gamma (n
 - 2\,a)\,\Gamma (n - 2\,b)\,\Gamma (n - a - b)\,2^{(2\,n + 1)}\,
\Gamma ({\displaystyle \frac {1}{2}}  + n - a)\,\Gamma (
{\displaystyle \frac {1}{2}}  + n - b)}}  }
}
\end{maplelatex}

\end{maplegroup}

%% file: AppendixB61to100.tex
\begin{maplegroup}
\mapleresult
\begin{maplelatex}
\mapleinline{inert}{2d}{63, ":   ",
1/4*GAMMA(n-1+1/2*a+1/2*b)*GAMMA(1/2*a-n+2-1/2*b)*(-a*sin(-Pi*b+1/2*Pi
*a)-a*sin(1/2*Pi*a))/cos(1/2*Pi*(b-2*n))/Pi*Sum(GAMMA(1-1/2*a-n+1/2*b+
k)*GAMMA(-1/2*b-n+3/2+k)/GAMMA(2+1/2*a-n-1/2*b+k)/GAMMA(3/2+1/2*b-n+k)
,k = 1 ..
2*n-3)-1/4*GAMMA(n-1+1/2*a+1/2*b)*GAMMA(1/2*a-n+2-1/2*b)*2^(1+a)*Pi^(1
/2)*GAMMA(1/2*a+1)/cos(1/2*Pi*(b-2*n))/GAMMA(1+a-b)/GAMMA(1/2*a+1/2)/G
AMMA(b);}{%
\maplemultiline{
63\mbox{:~~~} {\displaystyle \frac {1}{4}} \Gamma (n - 1
 + {\displaystyle \frac {a}{2}}  + {\displaystyle \frac {b}{2}} )
\,\Gamma ({\displaystyle \frac {a}{2}}  - n + 2 - {\displaystyle 
\frac {b}{2}} )\,( - a\,\mathrm{sin}( - \pi \,b + {\displaystyle 
\frac {1}{2}} \,\pi \,a) - a\,\mathrm{sin}({\displaystyle \frac {
\pi \,a}{2}} )) \\
 \left(  \! {\displaystyle \sum _{k=1}^{2\,n - 3}} \,
{\displaystyle \frac {\Gamma (1 - {\displaystyle \frac {a}{2}} 
 - n + {\displaystyle \frac {b}{2}}  + k)\,\Gamma ( - 
{\displaystyle \frac {b}{2}}  - n + {\displaystyle \frac {3}{2}} 
 + k)}{\Gamma (2 + {\displaystyle \frac {a}{2}}  - n - 
{\displaystyle \frac {b}{2}}  + k)\,\Gamma ({\displaystyle 
\frac {3}{2}}  + {\displaystyle \frac {b}{2}}  - n + k)}}  \! 
 \right)  \left/ {\vrule height0.80em width0em depth0.80em}
 \right. \!  \! (\mathrm{cos}({\displaystyle \frac {\pi \,(b - 2
\,n)}{2}} )\,\pi ) \\
\mbox{} - {\displaystyle \frac {1}{4}} \,{\displaystyle \frac {
\Gamma (n - 1 + {\displaystyle \frac {a}{2}}  + {\displaystyle 
\frac {b}{2}} )\,\Gamma ({\displaystyle \frac {a}{2}}  - n + 2 - 
{\displaystyle \frac {b}{2}} )\,2^{(1 + a)}\,\sqrt{\pi }\,\Gamma 
({\displaystyle \frac {a}{2}}  + 1)}{\mathrm{cos}({\displaystyle 
\frac {\pi \,(b - 2\,n)}{2}} )\,\Gamma (1 + a - b)\,\Gamma (
{\displaystyle \frac {a}{2}}  + {\displaystyle \frac {1}{2}} )\,
\Gamma (b)}}  }
}
\end{maplelatex}

\begin{maplelatex}
\mapleinline{inert}{2d}{64, ":   ",
-1/2*GAMMA(n-3/2+b-1/2*a)*GAMMA(2*b)/GAMMA(2*b-1)/GAMMA(n-1/2+b-1/2*a)
*Sum(GAMMA(3/2-b-n+1/2*a+k)/GAMMA(-b-n+3/2+1/2*a)*GAMMA(-1/2*a-n+3/2+k
)/GAMMA(-n+3/2-1/2*a)/GAMMA(3/2-1/2*a+b-n+k)*GAMMA(-1/2*a+b+3/2-n)/GAM
MA(3/2+1/2*a-n+k)*GAMMA(1/2*a-n+3/2),k = 1 ..
2*n-3)+Pi^(1/2)*GAMMA(2*n-1+2*b-a)/GAMMA(2*b-a)*GAMMA(a+2*n-1)/GAMMA(a
)/GAMMA(b)/(2^(4*n-1))*GAMMA(b-1/2)*GAMMA(-1/2*a+b+3/2-n)*GAMMA(1/2*a-
n+3/2)/GAMMA(n+b-1/2*a)/GAMMA(n+1/2*a)*GAMMA(n-3/2+b-1/2*a)*GAMMA(2*b)
/GAMMA(2*b-1)/GAMMA(n-1/2+b-1/2*a);}{%
\maplemultiline{
64\mbox{:~~~}  - {\displaystyle \frac {1}{2}} \Gamma (n
 - {\displaystyle \frac {3}{2}}  + b - {\displaystyle \frac {a}{2
}} )\,\Gamma (2\,b) \\
 \left(  \! {\displaystyle \sum _{k=1}^{2\,n - 3}} \,
{\displaystyle \frac {\Gamma ({\displaystyle \frac {3}{2}}  - b
 - n + {\displaystyle \frac {a}{2}}  + k)\,\Gamma ( - 
{\displaystyle \frac {a}{2}}  - n + {\displaystyle \frac {3}{2}} 
 + k)\,\Gamma ( - {\displaystyle \frac {a}{2}}  + b + 
{\displaystyle \frac {3}{2}}  - n)\,\Gamma ({\displaystyle 
\frac {a}{2}}  - n + {\displaystyle \frac {3}{2}} )}{\Gamma ( - b
 - n + {\displaystyle \frac {3}{2}}  + {\displaystyle \frac {a}{2
}} )\,\Gamma ( - n + {\displaystyle \frac {3}{2}}  - 
{\displaystyle \frac {a}{2}} )\,\Gamma ({\displaystyle \frac {3}{
2}}  - {\displaystyle \frac {a}{2}}  + b - n + k)\,\Gamma (
{\displaystyle \frac {3}{2}}  + {\displaystyle \frac {a}{2}}  - n
 + k)}}  \!  \right)  \\
 \left/ {\vrule height0.80em width0em depth0.80em} \right. \! 
 \! (\Gamma (2\,b - 1)\,\Gamma (n - {\displaystyle \frac {1}{2}} 
 + b - {\displaystyle \frac {a}{2}} ))\mbox{} + \sqrt{\pi }\,
\Gamma (2\,n - 1 + 2\,b - a)\,\Gamma (a + 2\,n - 1) \\
\Gamma (b - {\displaystyle \frac {1}{2}} )\,\Gamma ( - 
{\displaystyle \frac {a}{2}}  + b + {\displaystyle \frac {3}{2}} 
 - n)\,\Gamma ({\displaystyle \frac {a}{2}}  - n + 
{\displaystyle \frac {3}{2}} )\,\Gamma (n - {\displaystyle 
\frac {3}{2}}  + b - {\displaystyle \frac {a}{2}} )\,\Gamma (2\,b
) \left/ {\vrule height0.80em width0em depth0.80em} \right. \! 
 \! ( \\
\Gamma (2\,b - a)\,\Gamma (a)\,\Gamma (b)\,2^{(4\,n - 1)}\,\Gamma
 (n + b - {\displaystyle \frac {a}{2}} )\,\Gamma (n + 
{\displaystyle \frac {a}{2}} )\,\Gamma (2\,b - 1)\,\Gamma (n - 
{\displaystyle \frac {1}{2}}  + b - {\displaystyle \frac {a}{2}} 
)) }
}
\end{maplelatex}

\begin{maplelatex}
\mapleinline{inert}{2d}{65, ":   ",
-1/2*GAMMA(2*n-1-2*b)/GAMMA(2*n-2+a-2*b)*GAMMA(2*n-2+a-b)/GAMMA(2*n-1-
b)*Sum(GAMMA(b+2+k-2*n)/GAMMA(-2*n+2+b)*GAMMA(-1/2*a-n+3/2+k)/GAMMA(-n
+3/2-1/2*a)/GAMMA(1-b+k)*GAMMA(-b+1)/GAMMA(3/2+1/2*a-n+k)*GAMMA(1/2*a-
n+3/2),k = 1 ..
2*n-3)+Pi^(1/2)*GAMMA(4*n-2-2*b)*GAMMA(a+2*n-1)/GAMMA(a)/GAMMA(n+1/2*a
-b-1/2)/(2^(4*n-1))*GAMMA(n-1-b+1/2*a)*GAMMA(-b+1)/GAMMA(-1/2+2*n-b)/G
AMMA(n+1/2*a)*GAMMA(1/2*a-n+3/2)/GAMMA(2*n-2+a-2*b)*GAMMA(2*n-2+a-b)/G
AMMA(2*n-1-b);}{%
\maplemultiline{
65\mbox{:~~~}  - {\displaystyle \frac {1}{2}} \Gamma (2
\,n - 1 - 2\,b)\,\Gamma (2\,n - 2 + a - b) \\
 \left(  \! {\displaystyle \sum _{k=1}^{2\,n - 3}} \,
{\displaystyle \frac {\Gamma (b + 2 + k - 2\,n)\,\Gamma ( - 
{\displaystyle \frac {a}{2}}  - n + {\displaystyle \frac {3}{2}} 
 + k)\,\Gamma ( - b + 1)\,\Gamma ({\displaystyle \frac {a}{2}} 
 - n + {\displaystyle \frac {3}{2}} )}{\Gamma ( - 2\,n + 2 + b)\,
\Gamma ( - n + {\displaystyle \frac {3}{2}}  - {\displaystyle 
\frac {a}{2}} )\,\Gamma (1 - b + k)\,\Gamma ({\displaystyle 
\frac {3}{2}}  + {\displaystyle \frac {a}{2}}  - n + k)}}  \! 
 \right) /( \\
\Gamma (2\,n - 2 + a - 2\,b)\,\Gamma (2\,n - 1 - b))\mbox{} + 
\sqrt{\pi }\,\Gamma (4\,n - 2 - 2\,b)\,\Gamma (a + 2\,n - 1) \\
\Gamma (n - 1 - b + {\displaystyle \frac {a}{2}} )\,\Gamma ( - b
 + 1)\,\Gamma ({\displaystyle \frac {a}{2}}  - n + 
{\displaystyle \frac {3}{2}} )\,\Gamma (2\,n - 2 + a - b) \left/ 
{\vrule height0.80em width0em depth0.80em} \right. \!  \! (\Gamma
 (a) \\
\Gamma (n + {\displaystyle \frac {a}{2}}  - b - {\displaystyle 
\frac {1}{2}} )\,2^{(4\,n - 1)}\,\Gamma ( - {\displaystyle 
\frac {1}{2}}  + 2\,n - b)\,\Gamma (n + {\displaystyle \frac {a}{
2}} )\,\Gamma (2\,n - 2 + a - 2\,b) \Gamma (2\,n - 1 - b)) }
}
\end{maplelatex}

\begin{maplelatex}
\mapleinline{inert}{2d}{80, ":   ", WPMR(2*m+2*n-2,2*n-1,c,b,a);}{%
\[
80\mbox{:~~~} \,\mathrm{WPMR}(2\,m + 2\,n - 2, \,2\,n - 
1, \,c, \,b, \,a)
\]
}
\end{maplelatex}

\begin{maplelatex}
\mapleinline{inert}{2d}{81, ":   ", WPMR(2*m+2*n-1,2*n-1,c,b,a);}{%
\[
81\mbox{:~~~} \,\mathrm{WPMR}(2\,m + 2\,n - 1, \,2\,n - 
1, \,c, \,b, \,a)
\]
}
\end{maplelatex}

\begin{maplelatex}
\mapleinline{inert}{2d}{82, ":   ", WPMR(2*m+2*n-1,2*n,c,b,a);}{%
\[
82\mbox{:~~~} \,\mathrm{WPMR}(2\,m + 2\,n - 1, \,2\,n, 
\,c, \,b, \,a)
\]
}
\end{maplelatex}

\begin{maplelatex}
\mapleinline{inert}{2d}{84, ":   ",
(b*c-c-a*c-a*b+1+2*a-b-a*n+a^2)*GAMMA(a+2-b-c-n)/GAMMA(2+a-c)/GAMMA(2+
a-b)/GAMMA(-b-n+1+a)*GAMMA(2+a)/GAMMA(1+a-c-n)*GAMMA(a-n+1)*GAMMA(n)*S
um(GAMMA(a+k+1-n-c)*GAMMA(-b-n+1+a+k)/GAMMA(k+1)/GAMMA(a+1+k-n),k = 0
..
n)+GAMMA(a+2-b-c-n)*(-1-a)/GAMMA(-b-n+1+a)/n/GAMMA(1+a-c-n)*GAMMA(a-n+
1);}{%
\maplemultiline{
84\mbox{:~~~} (b\,c - c - a\,c - a\,b + 1 + 2\,a - b - a
\,n + a^{2})\,\Gamma (a + 2 - b - c - n)\,\Gamma (2 + a) \\
\Gamma (a - n + 1)\,\Gamma (n)\, \left(  \! {\displaystyle \sum 
_{k=0}^{n}} \,{\displaystyle \frac {\Gamma (a + k + 1 - n - c)\,
\Gamma ( - b - n + 1 + a + k)}{\Gamma (k + 1)\,\Gamma (a + 1 + k
 - n)}}  \!  \right) /(\Gamma (2 + a - c) \\
\Gamma (2 + a - b)\,\Gamma ( - b - n + 1 + a)\,\Gamma (1 + a - c
 - n)) 
\mbox{} + {\displaystyle \frac {\Gamma (a + 2 - b - c - n)\,( - 1
 - a)\,\Gamma (a - n + 1)}{\Gamma ( - b - n + 1 + a)\,n\,\Gamma (
1 + a - c - n)}}  }
}
\end{maplelatex}

\begin{maplelatex}
\mapleinline{inert}{2d}{85, ":   ",
-(a-n-1)*(a-n)*(a*b-b-1-2*n+a-a*c+n*c+c)/(a-n-b)*GAMMA(c)/GAMMA(a)*GAM
MA(n)/GAMMA(c-b)*Sum(GAMMA(a-1-n+k)*GAMMA(c-1-n-b+k)/GAMMA(k+1)/GAMMA(
c-n-1+k),k = 0 .. n)-(a-n-1)*(a-n)*(c-1)/n/(a-n-b);}{%
\maplemultiline{
85\mbox{:~~~}  - (a - n - 1)\,(a - n)\,(a\,b - b - 1 - 2
\,n + a - a\,c + n\,c + c)\,\Gamma (c)\,\Gamma (n) \\
 \left(  \! {\displaystyle \sum _{k=0}^{n}} \,{\displaystyle 
\frac {\Gamma (a - 1 - n + k)\,\Gamma (c - 1 - n - b + k)}{\Gamma
 (k + 1)\,\Gamma (c - n - 1 + k)}}  \!  \right) /((a - n - b)\,
\Gamma (a)\,\Gamma (c - b)) 
\mbox{} - {\displaystyle \frac {(a - n - 1)\,(a - n)\,(c - 1)}{n
\,(a - n - b)}}  }
}
\end{maplelatex}

\begin{maplelatex}
\mapleinline{inert}{2d}{86, ":   ",
-GAMMA(n+1+c+b)*GAMMA(a-n+1)*(n*c-a*n+n^2+b*n+b*c)/(a-c-b-n)/GAMMA(c)/
GAMMA(b)/GAMMA(c+n+1)/GAMMA(b+n+1)*GAMMA(n)*Sum(GAMMA(b+k)*GAMMA(c+k)/
GAMMA(k+1)/GAMMA(a+1+k-n),k = 0 ..
n)+1/n/(a-c-b-n)/GAMMA(1+a)/GAMMA(c)/GAMMA(b)*GAMMA(n+1+c+b)*GAMMA(a-n
+1);}{%
\maplemultiline{
86\mbox{:~~~}  - \Gamma (n + 1 + c + b)\,\Gamma (a - n
 + 1)\,(n\,c - a\,n + n^{2} + b\,n + b\,c)\,\Gamma (n) \\
 \left(  \! {\displaystyle \sum _{k=0}^{n}} \,{\displaystyle 
\frac {\Gamma (b + k)\,\Gamma (c + k)}{\Gamma (k + 1)\,\Gamma (a
 + 1 + k - n)}}  \!  \right) /((a - c - b - n)\,\Gamma (c)\,
\Gamma (b)\,\Gamma (c + n + 1) 
\Gamma (b + n + 1)) \\ \mbox{} + {\displaystyle \frac {\Gamma (n + 1
 + c + b)\,\Gamma (a - n + 1)}{n\,(a - c - b - n)\,\Gamma (1 + a)
\,\Gamma (c)\,\Gamma (b)}}  }
}
\end{maplelatex}

\begin{maplelatex}
\mapleinline{inert}{2d}{87, ":   ",
-(a-n)*(a-n-1)*(-b*n-a*c+n*c+a^2-a*n+c-a+n)/(b-1)/GAMMA(n+c+1-a)*GAMMA
(-b+2+c-a)/GAMMA(c-a)*GAMMA(c)/GAMMA(a)*GAMMA(n)*Sum(GAMMA(a-1-n+k)*GA
MMA(c-a+k)/GAMMA(k+1)/GAMMA(1-n+c-b+k),k = 0 ..
n)-(a-n)*(a-n-1)/n/(b-1)/GAMMA(c+1-b)/GAMMA(c-a)*GAMMA(-b+2+c-a)*GAMMA
(c);}{%
\maplemultiline{
87\mbox{:~~~}  - (a - n)\,(a - n - 1)\,( - b\,n - a\,c
 + n\,c + a^{2} - a\,n + c - a + n)\,\Gamma ( - b + 2 + c - a)
 \\
\Gamma (c)\,\Gamma (n)\, \left(  \! {\displaystyle \sum _{k=0}^{n
}} \,{\displaystyle \frac {\Gamma (a - 1 - n + k)\,\Gamma (c - a
 + k)}{\Gamma (k + 1)\,\Gamma (1 - n + c - b + k)}}  \!  \right) 
/((b - 1)\,\Gamma (n + c + 1 - a)\,\Gamma (c - a) 
\Gamma (a)) \\ \mbox{} - {\displaystyle \frac {(a - n)\,(a - n - 1)\,
\Gamma ( - b + 2 + c - a)\,\Gamma (c)}{n\,(b - 1)\,\Gamma (c + 1
 - b)\,\Gamma (c - a)}}  }
}
\end{maplelatex}

\begin{maplelatex}
\mapleinline{inert}{2d}{88, ":   ",
b*(b+1)*(a*n+a*b-b*c-b*n-n)/(-b+a-1)/GAMMA(b+n+1)/GAMMA(c-a)*GAMMA(c)*
GAMMA(n)*Sum(GAMMA(b+k)*GAMMA(c-a+k)/GAMMA(k+1)/GAMMA(c+k),k = 0 ..
n)+b*(b+1)/n/(-b+a-1)/GAMMA(c-a)*GAMMA(n+c+1-a)/GAMMA(c+n)*GAMMA(c);}{
\maplemultiline{
88\mbox{:~~~} {\displaystyle \frac {b\,(b + 1)\,(a\,n + 
a\,b - b\,c - b\,n - n)\,\Gamma (c)\,\Gamma (n)\, \left(  \! 
{\displaystyle \sum _{k=0}^{n}} \,{\displaystyle \frac {\Gamma (b
 + k)\,\Gamma (c - a + k)}{\Gamma (k + 1)\,\Gamma (c + k)}}  \! 
 \right) }{( - b + a - 1)\,\Gamma (b + n + 1)\,\Gamma (c - a)}} 
 \\
\mbox{} + {\displaystyle \frac {b\,(b + 1)\,\Gamma (n + c + 1 - a
)\,\Gamma (c)}{n\,( - b + a - 1)\,\Gamma (c - a)\,\Gamma (c + n)}
}  }
}
\end{maplelatex}

\begin{maplelatex}
\mapleinline{inert}{2d}{89, ":   ",
-(c-1)*(b-1)*(b-2)*(c-2)*(b-3+c+n-a)*(a*n+4-2*c-n-2*b+b*c)/(a-1)/(a-c+
3-b-n)*GAMMA(b-3+c+n-a)/GAMMA(b-1+n)/GAMMA(c-1+n)*GAMMA(n)*Sum(GAMMA(-
2+b+k)*GAMMA(c-2+k)/GAMMA(k+1)/GAMMA(c-2-a+k+b),k = 0 ..
n)+(c-1)*(b-1)*(b-2)*(c-2)/(a-1)/(a-c+3-b-n)/n;}{%
\maplemultiline{
89\mbox{:~~~}  - (c - 1)\,(b - 1)\,(b - 2)\,(c - 2)\,(b
 - 3 + c + n - a)\,(a\,n + 4 - 2\,c - n - 2\,b + b\,c) \\
\Gamma (b - 3 + c + n - a)\,\Gamma (n)\, \left(  \! 
{\displaystyle \sum _{k=0}^{n}} \,{\displaystyle \frac {\Gamma (
 - 2 + b + k)\,\Gamma (c - 2 + k)}{\Gamma (k + 1)\,\Gamma (c - 2
 - a + k + b)}}  \!  \right) /((a - 1) \\
(a - c + 3 - b - n)\,\Gamma (b - 1 + n)\,\Gamma (c - 1 + n))
\mbox{} + {\displaystyle \frac {(c - 1)\,(b - 1)\,(b - 2)\,(c - 2
)}{(a - 1)\,(a - c + 3 - b - n)\,n}}  }
}
\end{maplelatex}

\begin{maplelatex}
\mapleinline{inert}{2d}{91, ":   ",
(1+a)/GAMMA((-b*c+b+a*c+c-1)/a)/GAMMA((a*b-b*c+b+c-1)/a)*GAMMA((a*c+a*
b-a-b*c+b+c-1)/a)*GAMMA((-b*c+b+c-1+a)/a);}{%
\[
91\mbox{:~~~} \,{\displaystyle \frac {(1 + a)\,\Gamma (
{\displaystyle \frac {a\,c + a\,b - a - b\,c + b + c - 1}{a}} )\,
\Gamma ({\displaystyle \frac { - b\,c + b + c - 1 + a}{a}} )}{
\Gamma ({\displaystyle \frac { - b\,c + b + a\,c + c - 1}{a}} )\,
\Gamma ({\displaystyle \frac {a\,b - b\,c + b + c - 1}{a}} )}} 
\]
}
\end{maplelatex}

\begin{maplelatex}
\mapleinline{inert}{2d}{92, ":   ", (c-1)*(-1+a*b-b-a+c)/(1+b-c)/(1+a-c);}{%
\[
92\mbox{:~~~} \,{\displaystyle \frac {(c - 1)\,( - 1 + a
\,b - b - a + c)}{(1 + b - c)\,(1 + a - c)}} 
\]
}
\end{maplelatex}

\end{maplegroup}
\

%% file: AppendixB101to140.tex
\begin{maplegroup}
\mapleresult
\begin{maplelatex}
\mapleinline{inert}{2d}{104, ":   ",
V8(a,n)*GAMMA(1/2+a-1/2*n)*GAMMA(1/2+2*a-n)*GAMMA(3*a-1/2*n)/GAMMA(2*a
)/GAMMA(2*a+5/6-n)/GAMMA(2*a+1/6-n);}{%
\[
104\mbox{:~~~} {\displaystyle \frac {\mathrm{V8}(a, \,
n)\,\Gamma ({\displaystyle \frac {1}{2}}  + a - {\displaystyle 
\frac {n}{2}} )\,\Gamma ({\displaystyle \frac {1}{2}}  + 2\,a - n
)\,\Gamma (3\,a - {\displaystyle \frac {n}{2}} )}{\Gamma (2\,a)\,
\Gamma (2\,a + {\displaystyle \frac {5}{6}}  - n)\,\Gamma (2\,a
 + {\displaystyle \frac {1}{6}}  - n)}} 
\]
}
\end{maplelatex}

\begin{maplelatex}
\mapleinline{inert}{2d}{105, ":   ",
V8(a,n)*GAMMA(a-1/2*n)*GAMMA(1/2+2*a-n)*GAMMA(3*a+1/2-1/2*n)/GAMMA(2*a
)/GAMMA(2*a+5/6-n)/GAMMA(2*a+1/6-n);}{%
\[
105\mbox{:~~~} {\displaystyle \frac {\mathrm{V8}(a, \,
n)\,\Gamma (a - {\displaystyle \frac {n}{2}} )\,\Gamma (
{\displaystyle \frac {1}{2}}  + 2\,a - n)\,\Gamma (3\,a + 
{\displaystyle \frac {1}{2}}  - {\displaystyle \frac {n}{2}} )}{
\Gamma (2\,a)\,\Gamma (2\,a + {\displaystyle \frac {5}{6}}  - n)
\,\Gamma (2\,a + {\displaystyle \frac {1}{6}}  - n)}} 
\]
}
\end{maplelatex}

\begin{maplelatex}
\mapleinline{inert}{2d}{106, ":   ",
V8(a,n)*GAMMA(2*a-1/3)*GAMMA(1/2+2*a-n)/GAMMA(2*a)/GAMMA(2*a+1/6-n);}{
\[
106\mbox{:~~~} {\displaystyle \frac {\mathrm{V8}(a, \,
n)\,\Gamma (2\,a - {\displaystyle \frac {1}{3}} )\,\Gamma (
{\displaystyle \frac {1}{2}}  + 2\,a - n)}{\Gamma (2\,a)\,\Gamma 
(2\,a + {\displaystyle \frac {1}{6}}  - n)}} 
\]
}
\end{maplelatex}

\begin{maplelatex}
\mapleinline{inert}{2d}{107, ":   ",
V8(a,n)*GAMMA(2*a+1/3)*GAMMA(1/2+2*a-n)/GAMMA(2*a)/GAMMA(2*a+5/6-n);}{
\[
107\mbox{:~~~} {\displaystyle \frac {\mathrm{V8}(a, \,
n)\,\Gamma (2\,a + {\displaystyle \frac {1}{3}} )\,\Gamma (
{\displaystyle \frac {1}{2}}  + 2\,a - n)}{\Gamma (2\,a)\,\Gamma 
(2\,a + {\displaystyle \frac {5}{6}}  - n)}} 
\]
}
\end{maplelatex}

\begin{maplelatex}
\mapleinline{inert}{2d}{108, ":   ",
V8(a,n)*GAMMA(-n+1/2)*GAMMA(3*a-1/2*n)*GAMMA(3*a+1/2-1/2*n)/GAMMA(2*a)
/GAMMA(2*a+5/6-n)/GAMMA(2*a+1/6-n);}{%
\[
108\mbox{:~~~} {\displaystyle \frac {\mathrm{V8}(a, \,
n)\,\Gamma ( - n + {\displaystyle \frac {1}{2}} )\,\Gamma (3\,a
 - {\displaystyle \frac {n}{2}} )\,\Gamma (3\,a + {\displaystyle 
\frac {1}{2}}  - {\displaystyle \frac {n}{2}} )}{\Gamma (2\,a)\,
\Gamma (2\,a + {\displaystyle \frac {5}{6}}  - n)\,\Gamma (2\,a
 + {\displaystyle \frac {1}{6}}  - n)}} 
\]
}
\end{maplelatex}

\begin{maplelatex}
\mapleinline{inert}{2d}{109, ":   ",
V8(a,n)*GAMMA(a+1/6-1/2*n)*GAMMA(3*a-1/2*n)/GAMMA(2*a)/GAMMA(2*a+1/6-n
);}{%
\[
109\mbox{:~~~} {\displaystyle \frac {\mathrm{V8}(a, \,
n)\,\Gamma (a + {\displaystyle \frac {1}{6}}  - {\displaystyle 
\frac {n}{2}} )\,\Gamma (3\,a - {\displaystyle \frac {n}{2}} )}{
\Gamma (2\,a)\,\Gamma (2\,a + {\displaystyle \frac {1}{6}}  - n)}
} 
\]
}
\end{maplelatex}

\begin{maplelatex}
\mapleinline{inert}{2d}{110, ":   ",
V8(a,n)*GAMMA(a-1/2*n+5/6)*GAMMA(3*a-1/2*n)/GAMMA(2*a)/GAMMA(2*a+5/6-n
);}{%
\[
110\mbox{:~~~} {\displaystyle \frac {\mathrm{V8}(a, \,
n)\,\Gamma (a - {\displaystyle \frac {n}{2}}  + {\displaystyle 
\frac {5}{6}} )\,\Gamma (3\,a - {\displaystyle \frac {n}{2}} )}{
\Gamma (2\,a)\,\Gamma (2\,a + {\displaystyle \frac {5}{6}}  - n)}
} 
\]
}
\end{maplelatex}

\begin{maplelatex}
\mapleinline{inert}{2d}{111, ":   ",
V8(a,n)*GAMMA(a-1/3-1/2*n)*GAMMA(3*a+1/2-1/2*n)/GAMMA(2*a)/GAMMA(2*a+1
/6-n);}{%
\[
111\mbox{:~~~} {\displaystyle \frac {\mathrm{V8}(a, \,
n)\,\Gamma (a - {\displaystyle \frac {1}{3}}  - {\displaystyle 
\frac {n}{2}} )\,\Gamma (3\,a + {\displaystyle \frac {1}{2}}  - 
{\displaystyle \frac {n}{2}} )}{\Gamma (2\,a)\,\Gamma (2\,a + 
{\displaystyle \frac {1}{6}}  - n)}} 
\]
}
\end{maplelatex}

\begin{maplelatex}
\mapleinline{inert}{2d}{112, ":   ",
V8(a,n)*GAMMA(a+1/3-1/2*n)*GAMMA(3*a+1/2-1/2*n)/GAMMA(2*a)/GAMMA(2*a+5
/6-n);}{%
\[
112\mbox{:~~~} {\displaystyle \frac {\mathrm{V8}(a, \,
n)\,\Gamma (a + {\displaystyle \frac {1}{3}}  - {\displaystyle 
\frac {n}{2}} )\,\Gamma (3\,a + {\displaystyle \frac {1}{2}}  - 
{\displaystyle \frac {n}{2}} )}{\Gamma (2\,a)\,\Gamma (2\,a + 
{\displaystyle \frac {5}{6}}  - n)}} 
\]
}
\end{maplelatex}

\begin{maplelatex}
\mapleinline{inert}{2d}{113, ":   ", V8(a,n);}{%
\[
113\mbox{:~~~} \mathrm{V8}(a, \,n)
\]
}
\end{maplelatex}

\begin{maplelatex}
\mapleinline{inert}{2d}{114, ":   ",
V7(a,n)*GAMMA(a+2/3)*GAMMA(n-1/2+a)*GAMMA(-1/6+a+n)/GAMMA(-1/2+n)/GAMM
A(3*a+n)/Pi^(1/2);}{%
\[
114\mbox{:~~~} {\displaystyle \frac {\mathrm{V7}(a, \,
n)\,\Gamma (a + {\displaystyle \frac {2}{3}} )\,\Gamma (n - 
{\displaystyle \frac {1}{2}}  + a)\,\Gamma ( - {\displaystyle 
\frac {1}{6}}  + a + n)}{\Gamma ( - {\displaystyle \frac {1}{2}} 
 + n)\,\Gamma (3\,a + n)\,\sqrt{\pi }}} 
\]
}
\end{maplelatex}

\begin{maplelatex}
\mapleinline{inert}{2d}{115, ":   ",
V7(a,n)*GAMMA(a+1/3)*GAMMA(n-1/2+a)*GAMMA(1/6+a+n)/GAMMA(-1/2+n)/GAMMA
(3*a+n)/Pi^(1/2);}{%
\[
115\mbox{:~~~} {\displaystyle \frac {\mathrm{V7}(a, \,
n)\,\Gamma (a + {\displaystyle \frac {1}{3}} )\,\Gamma (n - 
{\displaystyle \frac {1}{2}}  + a)\,\Gamma ({\displaystyle 
\frac {1}{6}}  + a + n)}{\Gamma ( - {\displaystyle \frac {1}{2}} 
 + n)\,\Gamma (3\,a + n)\,\sqrt{\pi }}} 
\]
}
\end{maplelatex}

\begin{maplelatex}
\mapleinline{inert}{2d}{116, ":   ",
V7(a,n)*GAMMA(1/2-a)*GAMMA(n-1/2+a)/GAMMA(-1/2+n)/Pi^(1/2);}{%
\[
116\mbox{:~~~} {\displaystyle \frac {\mathrm{V7}(a, \,
n)\,\Gamma ({\displaystyle \frac {1}{2}}  - a)\,\Gamma (n - 
{\displaystyle \frac {1}{2}}  + a)}{\Gamma ( - {\displaystyle 
\frac {1}{2}}  + n)\,\sqrt{\pi }}} 
\]
}
\end{maplelatex}

\begin{maplelatex}
\mapleinline{inert}{2d}{117, ":   ",
V7(a,n)*GAMMA(2*a+n)*GAMMA(n-1/2+a)/GAMMA(-1/2+n)/GAMMA(3*a+n);}{%
\[
117\mbox{:~~~} {\displaystyle \frac {\mathrm{V7}(a, \,
n)\,\Gamma (2\,a + n)\,\Gamma (n - {\displaystyle \frac {1}{2}} 
 + a)}{\Gamma ( - {\displaystyle \frac {1}{2}}  + n)\,\Gamma (3\,
a + n)}} 
\]
}
\end{maplelatex}

\begin{maplelatex}
\mapleinline{inert}{2d}{118, ":   ",
V7(a,n)*GAMMA(a)*GAMMA(-1/6+a+n)*GAMMA(1/6+a+n)/GAMMA(-1/2+n)/GAMMA(3*
a+n)/Pi^(1/2);}{%
\[
118\mbox{:~~~} {\displaystyle \frac {\mathrm{V7}(a, \,
n)\,\Gamma (a)\,\Gamma ( - {\displaystyle \frac {1}{6}}  + a + n)
\,\Gamma ({\displaystyle \frac {1}{6}}  + a + n)}{\Gamma ( - 
{\displaystyle \frac {1}{2}}  + n)\,\Gamma (3\,a + n)\,\sqrt{\pi 
}}} 
\]
}
\end{maplelatex}

\begin{maplelatex}
\mapleinline{inert}{2d}{119, ":   ",
V7(a,n)*GAMMA(-a+1/6)*GAMMA(-1/6+a+n)/GAMMA(-1/2+n)/Pi^(1/2);}{%
\[
119\mbox{:~~~} {\displaystyle \frac {\mathrm{V7}(a, \,
n)\,\Gamma ( - a + {\displaystyle \frac {1}{6}} )\,\Gamma ( - 
{\displaystyle \frac {1}{6}}  + a + n)}{\Gamma ( - 
{\displaystyle \frac {1}{2}}  + n)\,\sqrt{\pi }}} 
\]
}
\end{maplelatex}

\begin{maplelatex}
\mapleinline{inert}{2d}{120, ":   ",
V7(a,n)*GAMMA(2*a-1/3+n)*GAMMA(-1/6+a+n)/GAMMA(-1/2+n)/GAMMA(3*a+n);}{
\[
120\mbox{:~~~} {\displaystyle \frac {\mathrm{V7}(a, \,
n)\,\Gamma (2\,a - {\displaystyle \frac {1}{3}}  + n)\,\Gamma (
 - {\displaystyle \frac {1}{6}}  + a + n)}{\Gamma ( - 
{\displaystyle \frac {1}{2}}  + n)\,\Gamma (3\,a + n)}} 
\]
}
\end{maplelatex}

\begin{maplelatex}
\mapleinline{inert}{2d}{121, ":   ",
V7(a,n)*GAMMA(-1/6-a)*GAMMA(1/6+a+n)/GAMMA(-1/2+n)/Pi^(1/2);}{%
\[
121\mbox{:~~~} {\displaystyle \frac {\mathrm{V7}(a, \,
n)\,\Gamma ( - {\displaystyle \frac {1}{6}}  - a)\,\Gamma (
{\displaystyle \frac {1}{6}}  + a + n)}{\Gamma ( - 
{\displaystyle \frac {1}{2}}  + n)\,\sqrt{\pi }}} 
\]
}
\end{maplelatex}

\begin{maplelatex}
\mapleinline{inert}{2d}{122, ":   ",
V7(a,n)*GAMMA(2*a+n-2/3)*GAMMA(1/6+a+n)/GAMMA(-1/2+n)/GAMMA(3*a+n);}{%
\[
122\mbox{:~~~} {\displaystyle \frac {\mathrm{V7}(a, \,
n)\,\Gamma (2\,a + n - {\displaystyle \frac {2}{3}} )\,\Gamma (
{\displaystyle \frac {1}{6}}  + a + n)}{\Gamma ( - 
{\displaystyle \frac {1}{2}}  + n)\,\Gamma (3\,a + n)}} 
\]
}
\end{maplelatex}

\begin{maplelatex}
\mapleinline{inert}{2d}{123, ":   ", V7(a,n);}{%
\[
123\mbox{:~~~} \mathrm{V7}(a, \,n)
\]
}
\end{maplelatex}

\begin{maplelatex}
\mapleinline{inert}{2d}{124, ":   ",
V1(a,n)*GAMMA(a+2/3)*GAMMA(a+1/2)*GAMMA(a+5/6)/Pi^(1/2)/GAMMA(3*a+1-n)
/GAMMA(1/2+n);}{%
\[
124\mbox{:~~~} {\displaystyle \frac {\mathrm{V1}(a, \,
n)\,\Gamma (a + {\displaystyle \frac {2}{3}} )\,\Gamma (a + 
{\displaystyle \frac {1}{2}} )\,\Gamma (a + {\displaystyle 
\frac {5}{6}} )}{\sqrt{\pi }\,\Gamma (3\,a + 1 - n)\,\Gamma (
{\displaystyle \frac {1}{2}}  + n)}} 
\]
}
\end{maplelatex}

\begin{maplelatex}
\mapleinline{inert}{2d}{125, ":   ",
V1(a,n)*GAMMA(a+1/3)*GAMMA(a+1/2)*GAMMA(7/6+a)/Pi^(1/2)/GAMMA(3*a+1-n)
/GAMMA(1/2+n);}{%
\[
125\mbox{:~~~} {\displaystyle \frac {\mathrm{V1}(a, \,
n)\,\Gamma (a + {\displaystyle \frac {1}{3}} )\,\Gamma (a + 
{\displaystyle \frac {1}{2}} )\,\Gamma ({\displaystyle \frac {7}{
6}}  + a)}{\sqrt{\pi }\,\Gamma (3\,a + 1 - n)\,\Gamma (
{\displaystyle \frac {1}{2}}  + n)}} 
\]
}
\end{maplelatex}

\begin{maplelatex}
\mapleinline{inert}{2d}{126, ":   ",
V1(a,n)*GAMMA(1/2+n-a)*GAMMA(a+1/2)/Pi^(1/2)/GAMMA(1/2+n);}{%
\[
126\mbox{:~~~} {\displaystyle \frac {\mathrm{V1}(a, \,
n)\,\Gamma ({\displaystyle \frac {1}{2}}  + n - a)\,\Gamma (a + 
{\displaystyle \frac {1}{2}} )}{\sqrt{\pi }\,\Gamma (
{\displaystyle \frac {1}{2}}  + n)}} 
\]
}
\end{maplelatex}

\begin{maplelatex}
\mapleinline{inert}{2d}{127, ":   ",
V1(a,n)*GAMMA(2*a+1-n)*GAMMA(a+1/2)/Pi^(1/2)/GAMMA(3*a+1-n);}{%
\[
127\mbox{:~~~} {\displaystyle \frac {\mathrm{V1}(a, \,
n)\,\Gamma (2\,a + 1 - n)\,\Gamma (a + {\displaystyle \frac {1}{2
}} )}{\sqrt{\pi }\,\Gamma (3\,a + 1 - n)}} 
\]
}
\end{maplelatex}

\begin{maplelatex}
\mapleinline{inert}{2d}{128, ":   ",
V1(a,n)*GAMMA(a)*GAMMA(a+5/6)*GAMMA(7/6+a)/Pi^(1/2)/GAMMA(3*a+1-n)/GAM
MA(1/2+n);}{%
\[
128\mbox{:~~~} {\displaystyle \frac {\mathrm{V1}(a, \,
n)\,\Gamma (a)\,\Gamma (a + {\displaystyle \frac {5}{6}} )\,
\Gamma ({\displaystyle \frac {7}{6}}  + a)}{\sqrt{\pi }\,\Gamma (
3\,a + 1 - n)\,\Gamma ({\displaystyle \frac {1}{2}}  + n)}} 
\]
}
\end{maplelatex}

\begin{maplelatex}
\mapleinline{inert}{2d}{129, ":   ",
V1(a,n)*GAMMA(1/6+n-a)*GAMMA(a+5/6)/Pi^(1/2)/GAMMA(1/2+n);}{%
\[
129\mbox{:~~~} {\displaystyle \frac {\mathrm{V1}(a, \,
n)\,\Gamma ({\displaystyle \frac {1}{6}}  + n - a)\,\Gamma (a + 
{\displaystyle \frac {5}{6}} )}{\sqrt{\pi }\,\Gamma (
{\displaystyle \frac {1}{2}}  + n)}} 
\]
}
\end{maplelatex}

\begin{maplelatex}
\mapleinline{inert}{2d}{130, ":   ",
V1(a,n)*GAMMA(2*a+2/3-n)*GAMMA(a+5/6)/Pi^(1/2)/GAMMA(3*a+1-n);}{%
\[
130\mbox{:~~~} {\displaystyle \frac {\mathrm{V1}(a, \,
n)\,\Gamma (2\,a + {\displaystyle \frac {2}{3}}  - n)\,\Gamma (a
 + {\displaystyle \frac {5}{6}} )}{\sqrt{\pi }\,\Gamma (3\,a + 1
 - n)}} 
\]
}
\end{maplelatex}

\begin{maplelatex}
\mapleinline{inert}{2d}{131, ":   ",
V1(a,n)*GAMMA(-1/6+n-a)*GAMMA(7/6+a)/Pi^(1/2)/GAMMA(1/2+n);}{%
\[
131\mbox{:~~~} {\displaystyle \frac {\mathrm{V1}(a, \,
n)\,\Gamma ( - {\displaystyle \frac {1}{6}}  + n - a)\,\Gamma (
{\displaystyle \frac {7}{6}}  + a)}{\sqrt{\pi }\,\Gamma (
{\displaystyle \frac {1}{2}}  + n)}} 
\]
}
\end{maplelatex}

\begin{maplelatex}
\mapleinline{inert}{2d}{132, ":   ",
V1(a,n)*GAMMA(2*a+1/3-n)*GAMMA(7/6+a)/Pi^(1/2)/GAMMA(3*a+1-n);}{%
\[
132\mbox{:~~~} {\displaystyle \frac {\mathrm{V1}(a, \,
n)\,\Gamma (2\,a + {\displaystyle \frac {1}{3}}  - n)\,\Gamma (
{\displaystyle \frac {7}{6}}  + a)}{\sqrt{\pi }\,\Gamma (3\,a + 1
 - n)}} 
\]
}
\end{maplelatex}

\begin{maplelatex}
\mapleinline{inert}{2d}{133, ":   ", V1(a,n);}{%
\[
133\mbox{:~~~} \mathrm{V1}(a, \,n)
\]
}
\end{maplelatex}

\begin{maplelatex}
\mapleinline{inert}{2d}{134, ":   ",
1/2*Pi^(3/2)*2^(2*n-2*a-2*b-2)/GAMMA(n-a-b)*GAMMA(n-b-1)*GAMMA(n-a-b-1
/2)*GAMMA(n-a-1)/sin(Pi*b)/sin(Pi*a)/GAMMA(n-2*b)/GAMMA(n-2*a)/GAMMA(b
-1)/GAMMA(a-1)+1/2*GAMMA(n-b-1)*GAMMA(n-a-1)*Sum(GAMMA(a+1+k-n)/GAMMA(
a-n+1)*GAMMA(b+1+k-n)/GAMMA(-n+b+1)/GAMMA(1-a+k)*GAMMA(-a+1)/GAMMA(1-b
+k)*GAMMA(-b+1),k = 1 ..
n-2)*GAMMA(a)*GAMMA(b)/GAMMA(n-b)/GAMMA(n-a)/GAMMA(a-1)/GAMMA(b-1);}{%
\[
\maplemultiline{
134\mbox{:~~~} {\displaystyle \frac {1}{2}} \,
{\displaystyle \frac {\pi ^{(3/2)}\,2^{(2\,n - 2\,a - 2\,b - 2)}
\,\Gamma (n - b - 1)\,\Gamma (n - a - b - {\displaystyle \frac {1
}{2}} )\,\Gamma (n - a - 1)}{\Gamma (n - a - b)\,\mathrm{sin}(\pi
 \,b)\,\mathrm{sin}(\pi \,a)\,\Gamma (n - 2\,b)\,\Gamma (n - 2\,a
)\,\Gamma (b - 1)\,\Gamma (a - 1)}} 
\\ + {\displaystyle \frac {1}{2
}} \,{\displaystyle \frac {\Gamma (n - b - 1)\,\Gamma (n - a - 1)
\, \left(  \! {\displaystyle \sum _{k=1}^{n - 2}} \,
{\displaystyle \frac {\Gamma (a + 1 + k - n)\,\Gamma (b + 1 + k
 - n)\,\Gamma ( - a + 1)\,\Gamma ( - b + 1)}{\Gamma (a - n + 1)\,
\Gamma ( - n + b + 1)\,\Gamma (1 - a + k)\,\Gamma (1 - b + k)}} 
 \!  \right) \,\Gamma (a)\,\Gamma (b)}{\Gamma (n - b)\,\Gamma (n
 - a)\,\Gamma (a - 1)\,\Gamma (b - 1)}} }
\]
}
\end{maplelatex}

\begin{maplelatex}
\mapleinline{inert}{2d}{139, ":   ",
-1/32*4^b*(sin(2*Pi*a-2*Pi*b)-sin(2*Pi*b)+sin(2*Pi*a)+2*sin(Pi*a-2*Pi*
b)+2*sin(Pi*a))*GAMMA(1/2*a-n+3/2)*GAMMA(-1/2*a+b+3/2-n)*GAMMA(n-3/2+b
-1/2*a)*GAMMA(n-3/2+1/2*a)*GAMMA(b-1/2)/Pi^(1/2)/sin(Pi*a)/sin(Pi*b)/s
in(Pi*(b-a))/GAMMA(a-1)/GAMMA(b-1)/GAMMA(2*b-a)+1/2*(a-1)*(b-1)*sin(1/
2*Pi*(2*b-a-1))*sin(1/2*Pi*(a-1))*Sum(GAMMA(3/2-b-n+1/2*a+k)*GAMMA(-1/
2*a-n+3/2+k)/GAMMA(3/2-1/2*a+b-n+k)/GAMMA(3/2+1/2*a-n+k),k = 1 ..
2*n-3)/Pi^2*GAMMA(n-3/2+1/2*a)*GAMMA(n-3/2+b-1/2*a)*GAMMA(-1/2*a+b+3/2
-n)*GAMMA(1/2*a-n+3/2);}{%
\maplemultiline{
139 \mbox{:~~~}  - {\displaystyle \frac {1}{32}} 4^{b}
( \mathrm{sin}(2\,\pi \,a - 2\,\pi \,b) - \mathrm{sin}(2\,\pi \,b)
 + \mathrm{sin}(2\,\pi \,a) + 2\,\mathrm{sin}(\pi \,a - 2\,\pi \,
b) + 2\,\mathrm{sin}(\pi \,a)) \\
\Gamma ({\displaystyle \frac {a}{2}}  - n + {\displaystyle 
\frac {3}{2}} )\,\Gamma ( - {\displaystyle \frac {a}{2}}  + b + 
{\displaystyle \frac {3}{2}}  - n)\,\Gamma (n - {\displaystyle 
\frac {3}{2}}  + b - {\displaystyle \frac {a}{2}} )\,\Gamma (n - 
{\displaystyle \frac {3}{2}}  + {\displaystyle \frac {a}{2}} )\,
\Gamma (b - {\displaystyle \frac {1}{2}} ) \left/ {\vrule 
height0.41em width0em depth0.41em} \right. \!  \! \\ (\sqrt{\pi }
\mathrm{sin}(\pi \,a)\,\mathrm{sin}(\pi \,b)\,\mathrm{sin}(\pi \,
(b - a))\,\Gamma (a - 1)\,\Gamma (b - 1)\,\Gamma (2\,b - a))
\mbox{} + {\displaystyle \frac {1}{2}} (a - 1)\,(b - 1)
\mathrm{sin}({\displaystyle \frac {\pi \,(2\,b - a - 1)}{2}} )\,
\mathrm{sin}({\displaystyle \frac {\pi \,(a - 1)}{2}} ) \\
 \left(  \! {\displaystyle \sum _{k=1}^{2\,n - 3}} \,
{\displaystyle \frac {\Gamma ({\displaystyle \frac {3}{2}}  - b
 - n + {\displaystyle \frac {a}{2}}  + k)\,\Gamma ( - 
{\displaystyle \frac {a}{2}}  - n + {\displaystyle \frac {3}{2}} 
 + k)}{\Gamma ({\displaystyle \frac {3}{2}}  - {\displaystyle 
\frac {a}{2}}  + b - n + k)\,\Gamma ({\displaystyle \frac {3}{2}
}  + {\displaystyle \frac {a}{2}}  - n + k)}}  \!  \right) \,
\Gamma (n - {\displaystyle \frac {3}{2}}  + {\displaystyle 
\frac {a}{2}} )\,\Gamma (n - {\displaystyle \frac {3}{2}}  + b - 
{\displaystyle \frac {a}{2}} )
\Gamma ( - {\displaystyle \frac {a}{2}}  + b + {\displaystyle 
\frac {3}{2}}  - n)\,\Gamma ({\displaystyle \frac {a}{2}}  - n + 
{\displaystyle \frac {3}{2}} ) \left/ {\vrule 
height0.43em width0em depth0.43em} \right. \!  \! \pi ^{2} }
}
\end{maplelatex}

\end{maplegroup}

%% file: AppendixB141to150.tex
\begin{maplegroup}
\mapleresult
\begin{maplelatex}
\mapleinline{inert}{2d}{148, ":   ",
-1/4*cos(-1/2*Pi*b+1/2*Pi*n+1/2*Pi*a)*(a-n+b+1)*GAMMA(1/2*n+1/2*a-1/2-
1/2*b)*GAMMA(1/2*a+1/2+1/2*b-1/2*n)*a*Sum(GAMMA(-1/2*b-1/2*n+1+k)*GAMM
A(-1/2*a+1/2+1/2*b-1/2*n+k)/Pi/GAMMA(3/2+1/2*a-1/2*b-1/2*n+k)/GAMMA(1+
1/2*b-1/2*n+k),k = 1 ..
n-2)-1/2*1/sin(1/2*Pi*b+1/2*Pi*n)*GAMMA(1/2*n+1/2*a-1/2-1/2*b)*GAMMA(1
/2*a+3/2+1/2*b-1/2*n)/GAMMA(1/2*a+1/2)/GAMMA(1+a-b)*GAMMA(1/2*a+1)/GAM
MA(b)*2^a*Pi^(1/2);}{%
\maplemultiline{
148\mbox{:~~~}   - {\displaystyle \frac {1}{4}} \mathrm{
cos}( - {\displaystyle \frac {\pi \,b}{2}}  + {\displaystyle 
\frac {\pi \,n}{2}}  + {\displaystyle \frac {\pi \,a}{2}} )\,(a
 - n + b + 1)\,\Gamma ({\displaystyle \frac {n}{2}}  + 
{\displaystyle \frac {a}{2}}  - {\displaystyle \frac {1}{2}}  - 
{\displaystyle \frac {b}{2}} )\,\Gamma ({\displaystyle \frac {a}{
2}}  + {\displaystyle \frac {1}{2}}  + {\displaystyle \frac {b}{2
}}  - {\displaystyle \frac {n}{2}} ) 
\, a \\  \left(  \! {\displaystyle \sum _{k=1}^{n - 2}} \,
{\displaystyle \frac {\Gamma ( - {\displaystyle \frac {b}{2}}  - 
{\displaystyle \frac {n}{2}}  + 1 + k)\,\Gamma ( - 
{\displaystyle \frac {a}{2}}  + {\displaystyle \frac {1}{2}}  + 
{\displaystyle \frac {b}{2}}  - {\displaystyle \frac {n}{2}}  + k
)}{\pi \,\Gamma ({\displaystyle \frac {3}{2}}  + {\displaystyle 
\frac {a}{2}}  - {\displaystyle \frac {b}{2}}  - {\displaystyle 
\frac {n}{2}}  + k)\,\Gamma (1 + {\displaystyle \frac {b}{2}}  - 
{\displaystyle \frac {n}{2}}  + k)}}  \!  \right)
\mbox{} - {\displaystyle \frac {1}{2}} \,{\displaystyle \frac {
\Gamma ({\displaystyle \frac {n}{2}}  + {\displaystyle \frac {a}{
2}}  - {\displaystyle \frac {1}{2}}  - {\displaystyle \frac {b}{2
}} )\,\Gamma ({\displaystyle \frac {a}{2}}  + {\displaystyle 
\frac {3}{2}}  + {\displaystyle \frac {b}{2}}  - {\displaystyle 
\frac {n}{2}} )\,\Gamma ({\displaystyle \frac {a}{2}}  + 1)\,2^{a
}\,\sqrt{\pi }}{\mathrm{sin}({\displaystyle \frac {1}{2}} \,\pi 
\,b + {\displaystyle \frac {1}{2}} \,\pi \,n)\,\Gamma (
{\displaystyle \frac {a}{2}}  + {\displaystyle \frac {1}{2}} )\,
\Gamma (1 + a - b)\,\Gamma (b)}}  }
}
\end{maplelatex}

\begin{maplelatex}
\mapleinline{inert}{2d}{149, ":   ",
1/4*(2*b-1)*(-cos(Pi*(b-a))+(-1)^n*cos(Pi*b))/sin(1/2*Pi*(n-a))/Pi/GAM
MA(-1/2*a+1/2*n-1)*GAMMA(-1/2*a-1/2*n+b+1)*GAMMA(1/2*n+b-1/2*a-1)*GAMM
A(1/2*n+1/2*a-1)*Sum(GAMMA(-1/2*a-1/2*n+1+k)*GAMMA(-b-1/2*n+1/2*a+1+k)
/GAMMA(1+1/2*a-1/2*n+k)/GAMMA(1+b-1/2*n-1/2*a+k),k = 1 ..
n-2)-1/4*1/sin(1/2*Pi*n-1/2*Pi*a)*GAMMA(-1/2*a-1/2*n+b+1)/GAMMA(-1/2*a
+1/2*n-1)*GAMMA(1/2*n+b-1/2*a-1)*GAMMA(1/2*n+1/2*a-1)*GAMMA(1/2+b)/GAM
MA(2*b-a)*2^(2*b)/GAMMA(b)/GAMMA(a)*Pi^(1/2);}{%
\maplemultiline{
149\mbox{:~~~}  {\displaystyle \frac {1}{4}} (2\,b - 1)\,
( - \mathrm{cos}(\pi \,(b - a)) + (-1)^{n}\,\mathrm{cos}(\pi \,b)
)\,\Gamma ( - {\displaystyle \frac {a}{2}}  - {\displaystyle 
\frac {n}{2}}  + b + 1)
\Gamma ({\displaystyle \frac {n}{2}}  + b - {\displaystyle 
\frac {a}{2}}  - 1)\,\Gamma ({\displaystyle \frac {n}{2}}  + 
{\displaystyle \frac {a}{2}}  - 1) \\
 \left(  \! {\displaystyle \sum _{k=1}^{n - 2}} \,{\displaystyle 
\frac {\Gamma ( - {\displaystyle \frac {a}{2}}  - {\displaystyle 
\frac {n}{2}}  + 1 + k)\,\Gamma ( - b - {\displaystyle \frac {n}{
2}}  + {\displaystyle \frac {a}{2}}  + 1 + k)}{\Gamma (1 + 
{\displaystyle \frac {a}{2}}  - {\displaystyle \frac {n}{2}}  + k
)\,\Gamma (1 + b - {\displaystyle \frac {n}{2}}  - 
{\displaystyle \frac {a}{2}}  + k)}}  \!  \right)  \left/ 
{\vrule height0.80em width0em depth0.80em} \right. \!  \! (
\mathrm{sin}({\displaystyle \frac {\pi \,(n - a)}{2}} )\,\pi
\Gamma ( - {\displaystyle \frac {a}{2}}  + {\displaystyle \frac {
n}{2}}  - 1)) \\
\mbox{} - {\displaystyle \frac {1}{4}} \,{\displaystyle \frac {
\Gamma ( - {\displaystyle \frac {a}{2}}  - {\displaystyle \frac {
n}{2}}  + b + 1)\,\Gamma ({\displaystyle \frac {n}{2}}  + b - 
{\displaystyle \frac {a}{2}}  - 1)\,\Gamma ({\displaystyle 
\frac {n}{2}}  + {\displaystyle \frac {a}{2}}  - 1)\,\Gamma (
{\displaystyle \frac {1}{2}}  + b)\,2^{(2\,b)}\,\sqrt{\pi }}{
\mathrm{sin}({\displaystyle \frac {1}{2}} \,\pi \,n - 
{\displaystyle \frac {1}{2}} \,\pi \,a)\,\Gamma ( - 
{\displaystyle \frac {a}{2}}  + {\displaystyle \frac {n}{2}}  - 1
)\,\Gamma (2\,b - a)\,\Gamma (b)\,\Gamma (a)}}  }
}
\end{maplelatex}

\mapleinline{inert}{2d}{150, ":   ",
-1/4*1/GAMMA(a)/Pi^2*GAMMA(1/2*a+b+1/2)*GAMMA(1/2*a-b+1/2)*GAMMA(3+a-n
-b)*GAMMA(n+b-2)*(-1)^n*Sum(1/GAMMA(5/2+1/2*a-b-n+k)/GAMMA(b+k)*GAMMA(
-b+2+k-n)*GAMMA(-1/2*a-1/2+b+k),k = 1 ..
n-2)*(-sin(1/2*Pi*a)+sin(1/2*Pi*(a-4*b)))+1/2*GAMMA(3+a-n-b)/GAMMA(3-2
*b-n+a)/GAMMA(n+2*b-2)*GAMMA(n+b-2)/GAMMA(1/2*a+1/2)^2*GAMMA(1/2*a+b+1
/2)*GAMMA(1/2*a-b+1/2);}{%
\maplemultiline{
150\mbox{:~~~}   - {\displaystyle \frac {1}{4}} \Gamma (
{\displaystyle \frac {a}{2}}  + b + {\displaystyle \frac {1}{2}} 
)\,\Gamma ({\displaystyle \frac {a}{2}}  - b + {\displaystyle 
\frac {1}{2}} )\,\Gamma (3 + a - n - b)\,\Gamma (n + b - 2)\,(-1)
^{n} \\
 \left(  \! {\displaystyle \sum _{k=1}^{n - 2}} \,{\displaystyle 
\frac {\Gamma ( - b + 2 + k - n)\,\Gamma ( - {\displaystyle 
\frac {a}{2}}  - {\displaystyle \frac {1}{2}}  + b + k)}{\Gamma (
{\displaystyle \frac {5}{2}}  + {\displaystyle \frac {a}{2}}  - b
 - n + k)\,\Gamma (b + k)}}  \!  \right) \,( - \mathrm{sin}(
{\displaystyle \frac {\pi \,a}{2}} ) + \mathrm{sin}(
{\displaystyle \frac {\pi \,(a - 4\,b)}{2}} )) \left/ {\vrule 
height0.43em width0em depth0.43em} \right. \!  \! (
\Gamma (a)\,\pi ^{2}) \\ \mbox{} + {\displaystyle \frac {1}{2}} \,
{\displaystyle \frac {\Gamma (3 + a - n - b)\,\Gamma (n + b - 2)
\,\Gamma ({\displaystyle \frac {a}{2}}  + b + {\displaystyle 
\frac {1}{2}} )\,\Gamma ({\displaystyle \frac {a}{2}}  - b + 
{\displaystyle \frac {1}{2}} )}{\Gamma (3 - 2\,b - n + a)\,\Gamma
 (n + 2\,b - 2)\,\Gamma ({\displaystyle \frac {a}{2}}  + 
{\displaystyle \frac {1}{2}} )^{2}}}  }
}

\end{maplegroup}

%% file: AppendixB151to170.tex
\begin{maplegroup}
\mapleresult
\begin{maplelatex}
\mapleinline{inert}{2d}{151, ":   ",
-1/2*sin(1/2*Pi*(a+n))/Pi*GAMMA(1/2*a-1/2*n+2)*GAMMA(1/2*n+1/2*a-1)*GA
MMA(3+a-n-b)*GAMMA(-n-2*b+4)*Sum(GAMMA(b-1+k)*GAMMA(-1/2*a-1/2*n+1+k)/
GAMMA(1+1/2*a-1/2*n+k)/GAMMA(3-b-n+k),k = 1 ..
n-2)/GAMMA(3-2*b-n+a)/GAMMA(b)-1/2*Pi/sin(1/2*Pi*n-1/2*Pi*a)*GAMMA(-b+
1)*GAMMA(1/2*n+1/2*a-1)/GAMMA(-1/2*a+1/2*n-1)/GAMMA(-1/2*n-b+2+1/2*a)^
2*GAMMA(3+a-n-b)/GAMMA(a);}{%
\maplemultiline{
151\mbox{:~~~}   - {\displaystyle \frac {1}{2}} \mathrm{
sin}({\displaystyle \frac {\pi \,(a + n)}{2}} )\,\Gamma (
{\displaystyle \frac {a}{2}}  - {\displaystyle \frac {n}{2}}  + 2
)\,\Gamma ({\displaystyle \frac {n}{2}}  + {\displaystyle \frac {
a}{2}}  - 1)\,\Gamma (3 + a - n - b)\,\Gamma ( - n - 2\,b + 4)
 \\
 \left(  \! {\displaystyle \sum _{k=1}^{n - 2}} \,{\displaystyle 
\frac {\Gamma (b - 1 + k)\,\Gamma ( - {\displaystyle \frac {a}{2}
}  - {\displaystyle \frac {n}{2}}  + 1 + k)}{\Gamma (1 + 
{\displaystyle \frac {a}{2}}  - {\displaystyle \frac {n}{2}}  + k
)\,\Gamma (3 - b - n + k)}}  \!  \right) /(\pi \,\Gamma (3 - 2\,b
 - n + a)\,\Gamma (b)) \\
\mbox{} - {\displaystyle \frac {1}{2}} \,{\displaystyle \frac {
\pi \,\Gamma ( - b + 1)\,\Gamma ({\displaystyle \frac {n}{2}}  + 
{\displaystyle \frac {a}{2}}  - 1)\,\Gamma (3 + a - n - b)}{
\mathrm{sin}({\displaystyle \frac {1}{2}} \,\pi \,n - 
{\displaystyle \frac {1}{2}} \,\pi \,a)\,\Gamma ( - 
{\displaystyle \frac {a}{2}}  + {\displaystyle \frac {n}{2}}  - 1
)\,\Gamma ( - {\displaystyle \frac {n}{2}}  - b + 2 + 
{\displaystyle \frac {a}{2}} )^{2}\,\Gamma (a)}}  }
}
\end{maplelatex}

\begin{maplelatex}
\mapleinline{inert}{2d}{152, ":   ",
4^(-a)*(-1)^n*(cos(Pi*(a-2*b))-cos(Pi*a))/Pi^(3/2)/GAMMA(a-1/2)*GAMMA(
-n-b+2+2*a)*GAMMA(n+b-2)*GAMMA(b+1)*GAMMA(a-b)*Sum(GAMMA(-b+2+k-n)*GAM
MA(-a+b+k)/GAMMA(2+a-b-n+k)/GAMMA(b+k),k = 1 ..
n-2)-1/2*Pi/sin(Pi*b)*GAMMA(-n-b+2+2*a)*GAMMA(n+b-2)/GAMMA(n+2*b-2)/GA
MMA(-n+2*a-2*b+2)*GAMMA(a-b)/GAMMA(a)/GAMMA(-b);}{%
\maplemultiline{
152\mbox{:~~~}  4^{( - a)}\,(-1)^{n}\,(\mathrm{cos}(\pi 
\,(a - 2\,b)) - \mathrm{cos}(\pi \,a))\,\Gamma ( - n - b + 2 + 2
\,a)\,\Gamma (n + b - 2)\,\Gamma (b + 1) \\
\Gamma (a - b)\, \left(  \! {\displaystyle \sum _{k=1}^{n - 2}} 
\,{\displaystyle \frac {\Gamma ( - b + 2 + k - n)\,\Gamma ( - a
 + b + k)}{\Gamma (2 + a - b - n + k)\,\Gamma (b + k)}}  \! 
 \right)  \left/ {\vrule height0.80em width0em depth0.80em}
 \right. \!  \! (\pi ^{(3/2)}\,\Gamma (a - {\displaystyle \frac {
1}{2}} )) \\
\mbox{} - {\displaystyle \frac {1}{2}} \,{\displaystyle \frac {
\pi \,\Gamma ( - n - b + 2 + 2\,a)\,\Gamma (n + b - 2)\,\Gamma (a
 - b)}{\mathrm{sin}(\pi \,b)\,\Gamma (n + 2\,b - 2)\,\Gamma ( - n
 + 2\,a - 2\,b + 2)\,\Gamma (a)\,\Gamma ( - b)}}  }
}
\end{maplelatex}

\begin{maplelatex}
\mapleinline{inert}{2d}{160, ":   ",
1/2*GAMMA(1+a)/GAMMA(2-2*n+2*b)*GAMMA(3-2*n+a)/GAMMA(b-a)*GAMMA(2*b-a)
*GAMMA(-b+1)*sin(Pi*a)/Pi*Sum(GAMMA(-a+k)*GAMMA(b+2+k-2*n)/GAMMA(3-2*n
+a+k)/GAMMA(1-b+k),k = 1 ..
2*n-3)-1/4*2^(2*a-2*b)*GAMMA(3-2*n+a)*GAMMA(a-b+1/2)*GAMMA(2*b-a)*GAMM
A(1+a)*GAMMA(-b+1)*(-2*sin(2*Pi*a)-2*sin(2*Pi*b)+sin(2*Pi*(a+b))+sin(2
*Pi*(a-b))+sin(4*Pi*b))/Pi^(3/2)/GAMMA(2*a-2*n+3)/(cos(Pi*a)-cos(Pi*(2
*b+a)));}{%
\maplemultiline{
160\mbox{:~~~}  {\displaystyle \frac {1}{2}} \Gamma (1 + 
a)\,\Gamma (3 - 2\,n + a)\,\Gamma (2\,b - a)\,\Gamma ( - b + 1)\,
\mathrm{sin}(\pi \,a) \\
 \left(  \! {\displaystyle \sum _{k=1}^{2\,n - 3}} \,
{\displaystyle \frac {\Gamma ( - a + k)\,\Gamma (b + 2 + k - 2\,n
)}{\Gamma (3 - 2\,n + a + k)\,\Gamma (1 - b + k)}}  \!  \right) /
(\Gamma (2 - 2\,n + 2\,b)\,\Gamma (b - a)\,\pi )\mbox{} \\ - 
{\displaystyle \frac {1}{4}} 2^{(2\,a - 2\,b)}
\Gamma (3 - 2\,n + a)\,\Gamma (a - b + {\displaystyle \frac {1}{2
}} )\,\Gamma (2\,b - a)\,\Gamma (1 + a)\,\Gamma ( - b + 1) \\
( - 2\,\mathrm{sin}(2\,\pi \,a) - 2\,\mathrm{sin}(2\,\pi \,b) + 
\mathrm{sin}(2\,\pi \,(a + b)) + \mathrm{sin}(2\,\pi \,(a - b))
 + \mathrm{sin}(4\,\pi \,b)) \left/ {\vrule 
height0.63em width0em depth0.63em} \right. \!  \!  \\
(\pi ^{(3/2)}\,\Gamma (2\,a - 2\,n + 3)\,(\mathrm{cos}(\pi \,a)
 - \mathrm{cos}(\pi \,(2\,b + a)))) }
}
\end{maplelatex}

\begin{maplelatex}
\mapleinline{inert}{2d}{161, ":   ",
1/2*4^(-b)*(cos(Pi*b)-cos(Pi*(2*a-b)))/Pi^(5/2)*GAMMA(3/2-b)*GAMMA(2*n
-3-a+b)*GAMMA(2*a-2*n+4-b)*GAMMA(3-2*n+a)*GAMMA(b-a)*GAMMA(1+a)*Sum(GA
MMA(-a+k)*GAMMA(3+a-2*n-b+k)/GAMMA(3-2*n+a+k)/GAMMA(-a+b+k),k = 1 ..
2*n-3)/GAMMA(2*a-2*n+4-2*b)-1/16*1/Pi^(5/2)/GAMMA(2*a-2*n+3)/sin(Pi*(2
*a-b))*4^b*GAMMA(2*n-3-a+b)*GAMMA(3-2*n+a)*GAMMA(2*a-2*n+4-b)*GAMMA(b+
1/2)*GAMMA(-2*b+1)*GAMMA(b-a)*GAMMA(1+a)*(-sin(2*Pi*b)+2*sin(2*Pi*a)+2
*sin(2*Pi*(a-b))-sin(2*Pi*(2*a-b))-sin(4*Pi*(a-b)));}{%
\maplemultiline{
161\mbox{:~~~}  {\displaystyle \frac {1}{2}} 4^{( - b)}\,
(\mathrm{cos}(\pi \,b) - \mathrm{cos}(\pi \,(2\,a - b)))\,\Gamma 
({\displaystyle \frac {3}{2}}  - b)\,\Gamma (2\,n - 3 - a + b)\,
\Gamma (2\,a - 2\,n + 4 - b) \\
\Gamma (3 - 2\,n + a)\,\Gamma (b - a)\,\Gamma (1 + a)\, \left( 
 \! {\displaystyle \sum _{k=1}^{2\,n - 3}} \,{\displaystyle 
\frac {\Gamma ( - a + k)\,\Gamma (3 + a - 2\,n - b + k)}{\Gamma (
3 - 2\,n + a + k)\,\Gamma ( - a + b + k)}}  \!  \right)  \left/ 
{\vrule height0.51em width0em depth0.51em} \right. \!  \! (\pi ^{
(5/2)} \\
\Gamma (2\,a - 2\,n + 4 - 2\,b))\mbox{} - {\displaystyle \frac {1
}{16}} 4^{b}\,\Gamma (2\,n - 3 - a + b)\,\Gamma (3 - 2\,n + a)\,
\Gamma (2\,a - 2\,n + 4 - b) \\
\Gamma (b + {\displaystyle \frac {1}{2}} )\,\Gamma ( - 2\,b + 1)
\,\Gamma (b - a)\,\Gamma (1 + a) \\
( - \mathrm{sin}(2\,\pi \,b) + 2\,\mathrm{sin}(2\,\pi \,a) + 2\,
\mathrm{sin}(2\,\pi \,(a - b)) - \mathrm{sin}(2\,\pi \,(2\,a - b)
) - \mathrm{sin}(4\,\pi \,(a - b))) \\
 \left/ {\vrule height0.51em width0em depth0.51em} \right. \! 
 \! (\pi ^{(5/2)}\,\Gamma (2\,a - 2\,n + 3)\,\mathrm{sin}(\pi \,(
2\,a - b))) }
}
\end{maplelatex}

\begin{maplelatex}
\mapleinline{inert}{2d}{162, ":   ",
-1/4*GAMMA(a+2+b-2*n)*GAMMA(2*n-1-a)*GAMMA(-1+a+b)*GAMMA(-a+1)*Sum(GAM
MA(a+2-2*n+k)*GAMMA(-a-b+1+k)/GAMMA(1-a+k)/GAMMA(2+a-2*n+b+k),k = 1 ..
2*n-3)*(cos(Pi*b)-cos(Pi*(2*a+b)))*(b-1)/Pi^2+1/8*4^b*GAMMA(2*n-1-a)*G
AMMA(a+2+b-2*n)*GAMMA(-a+1)*GAMMA(-1+a+b)*GAMMA(b-1/2)*(sin(4*Pi*(a+b)
)+sin(2*Pi*(2*a+b))-2*sin(2*Pi*(a+b))-2*sin(2*Pi*a)-sin(2*Pi*b))*(b-1)
/Pi^(1/2)/GAMMA(2*n-2*a-1)/GAMMA(2*a-2*n+2*b+1)/GAMMA(b)/(-sin(2*Pi*b)
+sin(4*Pi*(a+b))-sin(2*Pi*(2*a+b)));}{%
\maplemultiline{
162\mbox{:~~~}   - {\displaystyle \frac {1}{4}} \Gamma (a
 + 2 + b - 2\,n)\,\Gamma (2\,n - 1 - a)\,\Gamma ( - 1 + a + b)\,
\Gamma ( - a + 1) \\
 \left(  \! {\displaystyle \sum _{k=1}^{2\,n - 3}} \,
{\displaystyle \frac {\Gamma (a + 2 - 2\,n + k)\,\Gamma ( - a - b
 + 1 + k)}{\Gamma (1 - a + k)\,\Gamma (2 + a - 2\,n + b + k)}} 
 \!  \right) \,(\mathrm{cos}(\pi \,b) - \mathrm{cos}(\pi \,(2\,a
 + b)))\,(b - 1) \\
 \left/ {\vrule height0.43em width0em depth0.43em} \right. \! 
 \! \pi ^{2}\mbox{} + {\displaystyle \frac {1}{8}} 4^{b}\,\Gamma 
(2\,n - 1 - a)\,\Gamma (a + 2 + b - 2\,n)\,\Gamma ( - a + 1)\,
\Gamma ( - 1 + a + b)\,\Gamma (b - {\displaystyle \frac {1}{2}} )
( \\
\mathrm{sin}(4\,\pi \,(a + b)) + \mathrm{sin}(2\,\pi \,(2\,a + b)
) - 2\,\mathrm{sin}(2\,\pi \,(a + b)) - 2\,\mathrm{sin}(2\,\pi \,
a) - \mathrm{sin}(2\,\pi \,b)) \\
(b - 1) \left/ {\vrule height0.41em width0em depth0.41em}
 \right. \!  \! (\sqrt{\pi }\,\Gamma (2\,n - 2\,a - 1)\,\Gamma (2
\,a - 2\,n + 2\,b + 1)\,\Gamma (b) \\
( - \mathrm{sin}(2\,\pi \,b) + \mathrm{sin}(4\,\pi \,(a + b)) - 
\mathrm{sin}(2\,\pi \,(2\,a + b)))) }
}
\end{maplelatex}

\begin{maplelatex}
\mapleinline{inert}{2d}{163, ":   ",
1/2*GAMMA(1/2*b-n+3/2)*GAMMA(n-1/2+1/2*b)*GAMMA(2*n-2-a)*GAMMA(-a+1)*(
-1)^n*Sum(GAMMA(-1/2*b-n+3/2+k)*GAMMA(a+2-2*n+k)/GAMMA(3/2+1/2*b-n+k)/
GAMMA(1-a+k),k = 1 ..
2*n-3)*(sin(1/2*Pi*(2*a+b))+sin(1/2*Pi*(2*a-b)))*(n-a-1)/Pi^2+1/32*(-1
)^n*2^(-2*a+2*n+b)*(-2*sin(Pi*b)+sin(Pi*(2*a+b))+2*sin(2*Pi*a)-sin(Pi*
(2*a-b))-sin(4*Pi*a))/Pi^(3/2)/(sin(1/2*Pi*(2*a-b))+sin(1/2*Pi*(6*a+b)
))/GAMMA(-2*a-2+2*n)*GAMMA(1/2*b-n+3/2)*GAMMA(-1/2*b+a-n+3/2)*GAMMA(2*
n-2-a)*GAMMA(n-1/2+1/2*b)*GAMMA(1/2*b-a+n-1)*GAMMA(-a+1)/GAMMA(b);}{%
\maplemultiline{
163\mbox{:~~~}  {\displaystyle \frac {1}{2}} \Gamma (
{\displaystyle \frac {b}{2}}  - n + {\displaystyle \frac {3}{2}} 
)\,\Gamma (n - {\displaystyle \frac {1}{2}}  + {\displaystyle 
\frac {b}{2}} )\,\Gamma (2\,n - 2 - a)\,\Gamma ( - a + 1)\,(-1)^{
n}
 \left(  \! {\displaystyle \sum _{k=1}^{2\,n - 3}} \,
{\displaystyle \frac {\Gamma ( - {\displaystyle \frac {b}{2}}  - 
n + {\displaystyle \frac {3}{2}}  + k)\,\Gamma (a + 2 - 2\,n + k)
}{\Gamma ({\displaystyle \frac {3}{2}}  + {\displaystyle \frac {b
}{2}}  - n + k)\,\Gamma (1 - a + k)}}  \!  \right)  \\
(\mathrm{sin}({\displaystyle \frac {\pi \,(2\,a + b)}{2}} ) + 
\mathrm{sin}({\displaystyle \frac {\pi \,(2\,a - b)}{2}} ))\,(n
 - a - 1) \left/ {\vrule height0.43em width0em depth0.43em}
 \right. \!  \! \pi ^{2}\mbox{} + {\displaystyle \frac {1}{32}} (
-1)^{n}\,2^{( - 2\,a + 2\,n + b)} \\
( - 2\,\mathrm{sin}(\pi \,b) + \mathrm{sin}(\pi \,(2\,a + b)) + 2
\,\mathrm{sin}(2\,\pi \,a) - \mathrm{sin}(\pi \,(2\,a - b)) - 
\mathrm{sin}(4\,\pi \,a)) \\
\Gamma ({\displaystyle \frac {b}{2}}  - n + {\displaystyle 
\frac {3}{2}} )\,\Gamma ( - {\displaystyle \frac {b}{2}}  + a - n
 + {\displaystyle \frac {3}{2}} )\,\Gamma (2\,n - 2 - a)\,\Gamma 
(n - {\displaystyle \frac {1}{2}}  + {\displaystyle \frac {b}{2}
} )\,\Gamma ({\displaystyle \frac {b}{2}}  - a + n - 1) \\
\Gamma ( - a + 1) \left/ {\vrule 
height0.80em width0em depth0.80em} \right. \!  \! (\pi ^{(3/2)}\,
(\mathrm{sin}({\displaystyle \frac {\pi \,(2\,a - b)}{2}} ) + 
\mathrm{sin}({\displaystyle \frac {\pi \,(6\,a + b)}{2}} ))\,
\Gamma ( - 2\,a - 2 + 2\,n) \,
\Gamma (b)) }
}
\end{maplelatex}

\begin{maplelatex}
\mapleinline{inert}{2d}{164, ":   ",
-1/4*(cos(Pi*b)-cos(Pi*(2*a+b)))/Pi^2*GAMMA(-b+1-a)*GAMMA(2*n-2-a-b)*G
AMMA(2*a-2*n+b+3)*GAMMA(3-2*n+a)*GAMMA(2*b+a)*Sum(GAMMA(-a+k)*GAMMA(2+
a-2*n+b+k)/GAMMA(3-2*n+a+k)/GAMMA(-a-b+1+k),k = 1 ..
2*n-3)/GAMMA(2*a-2*n+2*b+2)/GAMMA(b)-1/8*4^(-b)*(sin(2*Pi*b)+2*sin(2*P
i*(a+b))+2*sin(2*Pi*a)-sin(4*Pi*(a+b))-sin(2*Pi*(2*a+b)))/Pi^(5/2)/sin
(Pi*(2*a+b))*GAMMA(-b+1-a)*GAMMA(2*n-2-a-b)*GAMMA(3-2*n+a)*GAMMA(2*a-2
*n+b+3)*GAMMA(2*b+a)*GAMMA(-b+1/2)/GAMMA(2*a-2*n+3);}{%
\maplemultiline{
164\mbox{:~~~}   - {\displaystyle \frac {1}{4}} (\mathrm{
cos}(\pi \,b) - \mathrm{cos}(\pi \,(2\,a + b)))\,\Gamma ( - b + 1
 - a)\,\Gamma (2\,n - 2 - a - b) 
\Gamma (2\,a - 2\,n + b + 3)\,\Gamma (3 - 2\,n + a) \\ \,\Gamma (2\,b
 + a) \,
 \left(  \! {\displaystyle \sum _{k=1}^{2\,n - 3}} \,
{\displaystyle \frac {\Gamma ( - a + k)\,\Gamma (2 + a - 2\,n + b
 + k)}{\Gamma (3 - 2\,n + a + k)\,\Gamma ( - a - b + 1 + k)}} 
 \!  \right)  \left/ {\vrule height0.43em width0em depth0.43em}
 \right. \!  \! (\pi ^{2}\,\Gamma (2\,a - 2\,n + 2\,b + 2)\,
\Gamma (b))\mbox{}  \\ -
{\displaystyle \frac {1}{8}} 4^{( - b)}(
\mathrm{sin}(2\,\pi \,b) + 2\,\mathrm{sin}(2\,\pi \,(a + b)) + 2
\,\mathrm{sin}(2\,\pi \,a) - \mathrm{sin}(4\,\pi \,(a + b)) - 
\mathrm{sin}(2\,\pi \,(2\,a + b))) \\
\Gamma ( - b + 1 - a)\,\Gamma (2\,n - 2 - a - b)\,\Gamma (3 - 2\,
n + a)\,\Gamma (2\,a - 2\,n + b + 3)\,\Gamma (2\,b + a) 
\Gamma ( - b + {\displaystyle \frac {1}{2}} ) \left/ {\vrule 
height0.63em width0em depth0.63em} \right. \!  \! (\pi ^{(5/2)}\,
\mathrm{sin}(\pi \,(2\,a + b))\,\Gamma (2\,a - 2\,n + 3)) }
}
\end{maplelatex}

\begin{maplelatex}
\mapleinline{inert}{2d}{165, ":   ",
-1/4*GAMMA(a-2*n+3-b)*GAMMA(-a+2*n-2+2*b)*GAMMA(2*n-2*a-1)*GAMMA(-a+1)
*GAMMA(a-b)*b*Sum(GAMMA(-a+b+k)*GAMMA(a+2-2*n+k)/GAMMA(3+a-2*n-b+k)/GA
MMA(1-a+k),k = 1 ..
2*n-3)*(cos(Pi*b)-cos(Pi*(2*a-b)))/GAMMA(2*n-2+2*b-2*a)/Pi^2+1/4*4^(-b
)*GAMMA(a-2*n+3-b)*GAMMA(-a+2*n-2+2*b)*GAMMA(-b+1/2)*GAMMA(-a+1)*GAMMA
(a-b)*(sin(2*Pi*b)-2*sin(2*Pi*a)-2*sin(2*Pi*(a-b))+sin(2*Pi*(2*a-b))+s
in(4*Pi*(a-b)))/Pi^(3/2)/GAMMA(-b)/(-cos(2*Pi*(a-b))+cos(2*Pi*a));}{%
\maplemultiline{
165\mbox{:~~~}   - {\displaystyle \frac {1}{4}} \Gamma (a
 - 2\,n + 3 - b)\,\Gamma ( - a + 2\,n - 2 + 2\,b)\,\Gamma (2\,n
 - 2\,a - 1)\,\Gamma ( - a + 1)\,\Gamma (a - b) \\
b\, \left(  \! {\displaystyle \sum _{k=1}^{2\,n - 3}} \,
{\displaystyle \frac {\Gamma ( - a + b + k)\,\Gamma (a + 2 - 2\,n
 + k)}{\Gamma (3 + a - 2\,n - b + k)\,\Gamma (1 - a + k)}}  \! 
 \right) \,(\mathrm{cos}(\pi \,b) - \mathrm{cos}(\pi \,(2\,a - b)
)) \left/ {\vrule height0.43em width0em depth0.43em} \right. \! 
 \! (
\Gamma (2\,n - 2 + 2\,b - 2\,a)\,\pi ^{2})\mbox{} \\ + 
{\displaystyle \frac {1}{4}} 4^{( - b)}\,\Gamma (a - 2\,n + 3 - b
)\,\Gamma ( - a + 2\,n - 2 + 2\,b)
\Gamma ( - b + {\displaystyle \frac {1}{2}} )\,\Gamma ( - a + 1)
\,\Gamma (a - b) \\ (
\mathrm{sin}(2\,\pi \,b) - 2\,\mathrm{sin}(2\,\pi \,a) - 2\,
\mathrm{sin}(2\,\pi \,(a - b)) + \mathrm{sin}(2\,\pi \,(2\,a - b)
) + \mathrm{sin}(4\,\pi \,(a - b))) \\
 \left/ {\vrule height0.63em width0em depth0.63em} \right. \! 
 \! (\pi ^{(3/2)}\,\Gamma ( - b)\,( - \mathrm{cos}(2\,\pi \,(a - 
b)) + \mathrm{cos}(2\,\pi \,a))) }
}
\end{maplelatex}

\begin{maplelatex}
\mapleinline{inert}{2d}{166, ":   ",
-1/2*GAMMA(a-2*n+3-1/2*b)*GAMMA(2*n-2-a)*GAMMA(a+1/2*b)*GAMMA(-a+1)*Su
m(GAMMA(1/2*b-a+k)*GAMMA(a+2-2*n+k)/GAMMA(3+a-2*n-1/2*b+k)/GAMMA(1-a+k
),k = 1 ..
2*n-3)*(-cos(1/2*Pi*b)+cos(1/2*Pi*(4*a-b)))*(n-a-1)/sin(1/2*Pi*b)/GAMM
A(1/2*b)^2/Pi+1/2*2^(-b)*GAMMA(2*n-2-a)*GAMMA(a-2*n+3-1/2*b)*GAMMA(-1/
2*b+1/2)*GAMMA(a+1/2*b)*GAMMA(-a+1)*(sin(Pi*(4*a-b))-2*sin(Pi*(2*a-b))
-sin(Pi*b)-2*sin(2*Pi*a)+sin(4*Pi*a))/Pi^(1/2)/GAMMA(1/2*b)/GAMMA(3+2*
a-2*n-b)/GAMMA(-2*a-2+2*n)/(-sin(Pi*(4*a-b))-sin(Pi*b)+sin(4*Pi*a));}{
\maplemultiline{
166\mbox{:~~~}   - {\displaystyle \frac {1}{2}} \Gamma (a
 - 2\,n + 3 - {\displaystyle \frac {b}{2}} )\,\Gamma (2\,n - 2 - 
a)\,\Gamma (a + {\displaystyle \frac {b}{2}} )\,\Gamma ( - a + 1)
 \\
 \left(  \! {\displaystyle \sum _{k=1}^{2\,n - 3}} \,
{\displaystyle \frac {\Gamma ({\displaystyle \frac {b}{2}}  - a
 + k)\,\Gamma (a + 2 - 2\,n + k)}{\Gamma (3 + a - 2\,n - 
{\displaystyle \frac {b}{2}}  + k)\,\Gamma (1 - a + k)}}  \! 
 \right) \,( - \mathrm{cos}({\displaystyle \frac {\pi \,b}{2}} )
 + \mathrm{cos}({\displaystyle \frac {\pi \,(4\,a - b)}{2}} ))
 \\
(n - a - 1) \left/ {\vrule height0.87em width0em depth0.87em}
 \right. \!  \! (\mathrm{sin}({\displaystyle \frac {\pi \,b}{2}} 
)\,\Gamma ({\displaystyle \frac {b}{2}} )^{2}\,\pi )\mbox{} + 
{\displaystyle \frac {1}{2}} 2^{( - b)}\,\Gamma (2\,n - 2 - a)\,
\Gamma (a - 2\,n + 3 - {\displaystyle \frac {b}{2}} ) \\
\Gamma ( - {\displaystyle \frac {b}{2}}  + {\displaystyle \frac {
1}{2}} )\,\Gamma (a + {\displaystyle \frac {b}{2}} )\,\Gamma ( - 
a + 1) 
(\mathrm{sin}(\pi \,(4\,a - b)) - 2\,\mathrm{sin}(\pi \,(2\,a - b
)) - \mathrm{sin}(\pi \,b) - 2\,\mathrm{sin}(2\,\pi \,a) + 
\mathrm{sin}(4\,\pi \,a)) \left/ {\vrule 
height0.80em width0em depth0.80em} \right. \!  \! ( \\
\sqrt{\pi }\,\Gamma ({\displaystyle \frac {b}{2}} )\,\Gamma (3 + 
2\,a - 2\,n - b)\,\Gamma ( - 2\,a - 2 + 2\,n) 
( - \mathrm{sin}(\pi \,(4\,a - b)) - \mathrm{sin}(\pi \,b) + 
\mathrm{sin}(4\,\pi \,a))) }
}
\end{maplelatex}

\begin{maplelatex}
\mapleinline{inert}{2d}{167, ":   ",
1/2*sin(Pi*(a+b))*(a-1)*GAMMA(a+2*b-4+2*n)*GAMMA(-1+a+b)*GAMMA(-b+1)*S
um(GAMMA(4-a-2*n-b+k)*GAMMA(b-1+k)/GAMMA(-1+b+a+k)/GAMMA(4-b-2*n+k),k
= 1 ..
2*n-3)/GAMMA(2*n+2*b-4)/Pi+1/8*2^(2*a)*GAMMA(a+2*b-4+2*n)*GAMMA(a-1/2)
*GAMMA(-1+a+b)*GAMMA(-b+1)*(2*cos(Pi*a)-cos(Pi*(a-2*b))-cos(Pi*(2*b+a)
))/Pi^(1/2)/GAMMA(2*a+2*b-5+2*n)/GAMMA(a-1)/(cos(Pi*(a-b))-cos(Pi*(a+b
)));}{%
\maplemultiline{
167\mbox{:~~~}  {\displaystyle \frac {1}{2}} \mathrm{sin}
(\pi \,(a + b))\,(a - 1)\,\Gamma (a + 2\,b - 4 + 2\,n)\,\Gamma (
 - 1 + a + b)\,\Gamma ( - b + 1) 
 \left(  \! {\displaystyle \sum _{k=1}^{2\,n - 3}} \,
{\displaystyle \frac {\Gamma (4 - a - 2\,n - b + k)\,\Gamma (b - 
1 + k)}{\Gamma ( - 1 + b + a + k)\,\Gamma (4 - b - 2\,n + k)}} 
 \!  \right) \\ /(\Gamma (2\,n + 2\,b - 4)\,\pi )\mbox{} + 
{\displaystyle \frac {1}{8}} 2^{(2\,a)} 
\Gamma (a + 2\,b - 4 + 2\,n)\,\Gamma (a - {\displaystyle \frac {1
}{2}} )\,\Gamma ( - 1 + a + b)\,\Gamma ( - b + 1) \\
(2\,\mathrm{cos}(\pi \,a) - \mathrm{cos}(\pi \,(a - 2\,b)) - 
\mathrm{cos}(\pi \,(2\,b + a))) \left/ {\vrule 
height0.41em width0em depth0.41em} \right. \!  \! (\sqrt{\pi }\,
\Gamma (2\,a + 2\,b - 5 + 2\,n)
\Gamma (a - 1)\,(\mathrm{cos}(\pi \,(a - b)) - \mathrm{cos}(\pi 
\,(a + b)))) }
}
\end{maplelatex}

\begin{maplelatex}
\mapleinline{inert}{2d}{168, ":   ",
-1/2*GAMMA(n-1-1/2*b)/GAMMA(1/2*b+1/2-1/2*a)*GAMMA(1/2*a+1/2+1/2*b)*GA
MMA(-1/2*b+2-n)*sin(1/2*Pi*b)*(-1)^n/Pi*b*Sum(GAMMA(-1/2*a-n+3/2+k)*GA
MMA(1/2*b+1-n+k)/GAMMA(3/2+1/2*a-n+k)/GAMMA(2-1/2*b-n+k),k = 1 ..
2*n-3)-1/4*(-1)^n*2^(a-b)*GAMMA(-1/2*b+2-n)*GAMMA(n-1-1/2*b)*GAMMA(1/2
*a-1/2*b)*GAMMA(1/2*a+1/2+1/2*b)*(2*sin(Pi*a)-2*sin(Pi*b)+sin(2*Pi*b)-
sin(Pi*(a-b))-sin(Pi*(a+b)))/Pi^(1/2)/GAMMA(a)/GAMMA(-b)/(sin(1/2*Pi*(
2*a+3*b))+sin(1/2*Pi*b)-sin(1/2*Pi*(2*a-b))+sin(3/2*Pi*b));}{%
\maplemultiline{
168\mbox{:~~~}   - {\displaystyle \frac {1}{2}} \Gamma (n
 - 1 - {\displaystyle \frac {b}{2}} )\,\Gamma ({\displaystyle 
\frac {a}{2}}  + {\displaystyle \frac {1}{2}}  + {\displaystyle 
\frac {b}{2}} )\,\Gamma ( - {\displaystyle \frac {b}{2}}  + 2 - n
)\,\mathrm{sin}({\displaystyle \frac {\pi \,b}{2}} )\,(-1)^{n}\,b
 \\
 \left(  \! {\displaystyle \sum _{k=1}^{2\,n - 3}} \,
{\displaystyle \frac {\Gamma ( - {\displaystyle \frac {a}{2}}  - 
n + {\displaystyle \frac {3}{2}}  + k)\,\Gamma ({\displaystyle 
\frac {b}{2}}  + 1 - n + k)}{\Gamma ({\displaystyle \frac {3}{2}
}  + {\displaystyle \frac {a}{2}}  - n + k)\,\Gamma (2 - 
{\displaystyle \frac {b}{2}}  - n + k)}}  \!  \right)  \left/ 
{\vrule height0.80em width0em depth0.80em} \right. \!  \! (\Gamma
 ({\displaystyle \frac {b}{2}}  + {\displaystyle \frac {1}{2}} 
 - {\displaystyle \frac {a}{2}} )\,\pi )\mbox{} \\ - {\displaystyle 
\frac {1}{4}} (-1)^{n}
2^{(a - b)}\,\Gamma ( - {\displaystyle \frac {b}{2}}  + 2 - n)\,
\Gamma (n - 1 - {\displaystyle \frac {b}{2}} )\,\Gamma (
{\displaystyle \frac {a}{2}}  - {\displaystyle \frac {b}{2}} )\,
\Gamma ({\displaystyle \frac {a}{2}}  + {\displaystyle \frac {1}{
2}}  + {\displaystyle \frac {b}{2}} ) \\
(2\,\mathrm{sin}(\pi \,a) - 2\,\mathrm{sin}(\pi \,b) + \mathrm{
sin}(2\,\pi \,b) - \mathrm{sin}(\pi \,(a - b)) - \mathrm{sin}(\pi
 \,(a + b))) \left/ {\vrule height0.80em width0em depth0.80em}
 \right. \!  \! (\sqrt{\pi } \\
\Gamma (a)\,\Gamma ( - b)\,(\mathrm{sin}({\displaystyle \frac {
\pi \,(2\,a + 3\,b)}{2}} ) + \mathrm{sin}({\displaystyle \frac {
\pi \,b}{2}} ) - \mathrm{sin}({\displaystyle \frac {\pi \,(2\,a
 - b)}{2}} ) + \mathrm{sin}({\displaystyle \frac {3\,\pi \,b}{2}
} ))) }
}
\end{maplelatex}

\end{maplegroup}

%% file: AppendixB171to200.tex
\begin{maplegroup}
\begin{mapleinput}
\end{mapleinput}

\mapleresult
\begin{maplelatex}
\mapleinline{inert}{2d}{195, ":   ",
-1/8*1/GAMMA(2*n-2+a-2*b)/Pi^2/sin(Pi*a)*GAMMA(n-3/2+1/2*a)*GAMMA(5/2+
1/2*a-n)*GAMMA(2*n-1-2*b)*GAMMA(2*n-2-b)*GAMMA(1+a-b)*(-1)^n*(cos(1/2*
Pi*(a-2*b))-cos(1/2*Pi*(3*a+2*b))-cos(1/2*Pi*(2*b+a))+cos(1/2*Pi*(3*a-
2*b)))*Sum(GAMMA(-1/2*a-n+3/2+k)*GAMMA(b+2+k-2*n)/GAMMA(3/2+1/2*a-n+k)
/GAMMA(1-b+k),k = 1 ..
2*n-3)+1/8*(2*n+a-2*b-1)^2*GAMMA(2*n-2-b)*GAMMA(1+a-b)*GAMMA(5/2+1/2*a
-n)*GAMMA(n-3/2+1/2*a)/GAMMA(a)/GAMMA(n+1/2*a-b+1/2)^2;}{%
\maplemultiline{
195\mbox{:~~~}  - {\displaystyle \frac {1}{8}} \Gamma (n
 - {\displaystyle \frac {3}{2}}  + {\displaystyle \frac {a}{2}} )
\,\Gamma ({\displaystyle \frac {5}{2}}  + {\displaystyle \frac {a
}{2}}  - n)\,\Gamma (2\,n - 1 - 2\,b)\,\Gamma (2\,n - 2 - b)\,
\Gamma (1 + a - b)\,(-1)^{n} \\
(\mathrm{cos}({\displaystyle \frac {\pi \,(a - 2\,b)}{2}} ) - 
\mathrm{cos}({\displaystyle \frac {\pi \,(3\,a + 2\,b)}{2}} ) - 
\mathrm{cos}({\displaystyle \frac {\pi \,(2\,b + a)}{2}} ) + 
\mathrm{cos}({\displaystyle \frac {\pi \,(3\,a - 2\,b)}{2}} ))
 \\
 \left(  \! {\displaystyle \sum _{k=1}^{2\,n - 3}} \,
{\displaystyle \frac {\Gamma ( - {\displaystyle \frac {a}{2}}  - 
n + {\displaystyle \frac {3}{2}}  + k)\,\Gamma (b + 2 + k - 2\,n)
}{\Gamma ({\displaystyle \frac {3}{2}}  + {\displaystyle \frac {a
}{2}}  - n + k)\,\Gamma (1 - b + k)}}  \!  \right)  \left/ 
{\vrule height0.43em width0em depth0.43em} \right. \!  \! (\Gamma
 (2\,n - 2 + a - 2\,b)\,\pi ^{2} \\
\mathrm{sin}(\pi \,a)) \\
\mbox{} + {\displaystyle \frac {1}{8}} \,{\displaystyle \frac {(2
\,n + a - 2\,b - 1)^{2}\,\Gamma (2\,n - 2 - b)\,\Gamma (1 + a - b
)\,\Gamma ({\displaystyle \frac {5}{2}}  + {\displaystyle \frac {
a}{2}}  - n)\,\Gamma (n - {\displaystyle \frac {3}{2}}  + 
{\displaystyle \frac {a}{2}} )}{\Gamma (a)\,\Gamma (n + 
{\displaystyle \frac {a}{2}}  - b + {\displaystyle \frac {1}{2}} 
)^{2}}}  }
}
\end{maplelatex}

\begin{maplelatex}
\mapleinline{inert}{2d}{196, ":   ",
1/2*sin(1/2*Pi*(a+n))*(a-1)*(b-1)/sin(1/2*Pi*(-2*b+n+a))*GAMMA(1/2*a-1
/2*n+1)*GAMMA(1/2*n+1/2*a-1)/GAMMA(-1/2*n-b+2+1/2*a)/GAMMA(1/2*a+1/2*n
-b)*Sum(GAMMA(-1/2*a-1/2*n+1+k)*GAMMA(-b-1/2*n+1/2*a+1+k)/GAMMA(1+1/2*
a-1/2*n+k)/GAMMA(1+b-1/2*n-1/2*a+k),k = 1 ..
n-2)-1/2*Pi^(3/2)*4^b*GAMMA(1/2*n+1/2*a-1)*GAMMA(1/2*a-1/2*n+1)*GAMMA(
b-1/2)*(-2*cos(Pi*(b+n))+cos(Pi*(a-b))+cos(Pi*(a+b)))/GAMMA(1/2*a+1/2*
n-b)/GAMMA(-1/2*n-b+2+1/2*a)/GAMMA(a-1)/GAMMA(b-1)/GAMMA(2*b-a)/(-2*co
s(Pi*(a-b))*(-1)^n+cos(Pi*b)+cos(Pi*(-3*b+2*a))+2*cos(Pi*(a+b))*(-1)^n
-cos(3*Pi*b)-cos(Pi*(2*a-b)));}{%
\maplemultiline{
196\mbox{:~~~} {\displaystyle \frac {1}{2}} \mathrm{sin}
({\displaystyle \frac {\pi \,(a + n)}{2}} )\,(a - 1)\,(b - 1)\,
\Gamma ({\displaystyle \frac {a}{2}}  - {\displaystyle \frac {n}{
2}}  + 1)\,\Gamma ({\displaystyle \frac {n}{2}}  + 
{\displaystyle \frac {a}{2}}  - 1) \\
 \left(  \! {\displaystyle \sum _{k=1}^{n - 2}} \,{\displaystyle 
\frac {\Gamma ( - {\displaystyle \frac {a}{2}}  - {\displaystyle 
\frac {n}{2}}  + 1 + k)\,\Gamma ( - b - {\displaystyle \frac {n}{
2}}  + {\displaystyle \frac {a}{2}}  + 1 + k)}{\Gamma (1 + 
{\displaystyle \frac {a}{2}}  - {\displaystyle \frac {n}{2}}  + k
)\,\Gamma (1 + b - {\displaystyle \frac {n}{2}}  - 
{\displaystyle \frac {a}{2}}  + k)}}  \!  \right)  \left/ 
{\vrule height0.80em width0em depth0.80em} \right. \!  \! (
\mathrm{sin}({\displaystyle \frac {\pi \,( - 2\,b + n + a)}{2}} )
 \\
\Gamma ( - {\displaystyle \frac {n}{2}}  - b + 2 + 
{\displaystyle \frac {a}{2}} )\,\Gamma ({\displaystyle \frac {a}{
2}}  + {\displaystyle \frac {n}{2}}  - b))\mbox{} - 
{\displaystyle \frac {1}{2}} \pi ^{(3/2)}\,4^{b}\,\Gamma (
{\displaystyle \frac {n}{2}}  + {\displaystyle \frac {a}{2}}  - 1
)\,\Gamma ({\displaystyle \frac {a}{2}}  - {\displaystyle \frac {
n}{2}}  + 1) \\
\Gamma (b - {\displaystyle \frac {1}{2}} )\,( - 2\,\mathrm{cos}(
\pi \,(b + n)) + \mathrm{cos}(\pi \,(a - b)) + \mathrm{cos}(\pi 
\,(a + b))) \left/ {\vrule height0.80em width0em depth0.80em}
 \right. \!  \! (\Gamma ({\displaystyle \frac {a}{2}}  + 
{\displaystyle \frac {n}{2}}  - b) \\
\Gamma ( - {\displaystyle \frac {n}{2}}  - b + 2 + 
{\displaystyle \frac {a}{2}} )\,\Gamma (a - 1)\,\Gamma (b - 1)\,
\Gamma (2\,b - a)( - 2\,\mathrm{cos}(\pi \,(a - b))\,(-1)^{n} + 
\mathrm{cos}(\pi \,b) \\
\mbox{} + \mathrm{cos}(\pi \,( - 3\,b + 2\,a)) + 2\,\mathrm{cos}(
\pi \,(a + b))\,(-1)^{n} - \mathrm{cos}(3\,\pi \,b) - \mathrm{cos
}(\pi \,(2\,a - b)))) }
}
\end{maplelatex}

\begin{maplelatex}
\mapleinline{inert}{2d}{197, ":   ",
1/2*(n+2*b-3)*(-sin(Pi*b)*sin(1/2*Pi*(2*a+1))+sin(Pi*(a+b)))*sin(Pi*b)
*sin(Pi*(a-b))*Sum(GAMMA(-b+2+k-n)*GAMMA(-a+b+k)/GAMMA(2+a-b-n+k)/GAMM
A(b+k),k = 1 ..
n-2)*GAMMA(n+b-2)*GAMMA(b)*GAMMA(a-b)*GAMMA(a)^2/Pi^2/sin(Pi*(a+b))/si
n(Pi*(b-1/2))/GAMMA(b+n+a-3)-1/2*GAMMA(n+b-2)*GAMMA(2*a-1)*GAMMA(b)*GA
MMA(a-b)*sin(Pi*b)*cos(Pi*a)/GAMMA(b+n+a-3)/GAMMA(n+2*b-3)/GAMMA(-n+2*
a-2*b+2)/cos(Pi*(b+n))/sin(Pi*(a+b));}{%
\maplemultiline{
197\mbox{:~~~} {\displaystyle \frac {1}{2}} (n + 2\,b - 
3)\,( - \mathrm{sin}(\pi \,b)\,\mathrm{sin}({\displaystyle 
\frac {\pi \,(2\,a + 1)}{2}} ) + \mathrm{sin}(\pi \,(a + b)))\,
\mathrm{sin}(\pi \,b) \\
\mathrm{sin}(\pi \,(a - b))\, \left(  \! {\displaystyle \sum _{k=
1}^{n - 2}} \,{\displaystyle \frac {\Gamma ( - b + 2 + k - n)\,
\Gamma ( - a + b + k)}{\Gamma (2 + a - b - n + k)\,\Gamma (b + k)
}}  \!  \right) \,\Gamma (n + b - 2)\,\Gamma (b)\,\Gamma (a - b)
 \\
\Gamma (a)^{2} \left/ {\vrule height0.80em width0em depth0.80em}
 \right. \!  \! (\pi ^{2}\,\mathrm{sin}(\pi \,(a + b))\,\mathrm{
sin}(\pi \,(b - {\displaystyle \frac {1}{2}} ))\,\Gamma (b + n + 
a - 3))\mbox{} - {\displaystyle \frac {1}{2}}  \\
{\displaystyle \frac {\Gamma (n + b - 2)\,\Gamma (2\,a - 1)\,
\Gamma (b)\,\Gamma (a - b)\,\mathrm{sin}(\pi \,b)\,\mathrm{cos}(
\pi \,a)}{\Gamma (b + n + a - 3)\,\Gamma (n + 2\,b - 3)\,\Gamma (
 - n + 2\,a - 2\,b + 2)\,\mathrm{cos}(\pi \,(b + n))\,\mathrm{sin
}(\pi \,(a + b))}}  }
}
\end{maplelatex}

\begin{maplelatex}
\mapleinline{inert}{2d}{198, ":   ",
-1/2*sin(Pi*(a+b))*(b-1)/Pi/GAMMA(a-1)*GAMMA(a+2+b-2*n)*GAMMA(2*n-2-a)
*GAMMA(-1+a+b)*Sum(GAMMA(a+2-2*n+k)*GAMMA(-a-b+1+k)/GAMMA(1-a+k)/GAMMA
(2+a-2*n+b+k),k = 1 ..
2*n-3)+Pi^2*b*GAMMA(a+2+b-2*n)*GAMMA(2*n-2-a)*GAMMA(-1+a+b)*(b-1)*(sin
(2*Pi*b+5*Pi*a)+sin(Pi*a+4*Pi*b)-sin(4*Pi*b+7*Pi*a)-sin(6*Pi*b+7*Pi*a)
+3*sin(4*Pi*b+5*Pi*a)-sin(6*Pi*b+5*Pi*a)-sin(8*Pi*b+5*Pi*a)+2*sin(Pi*a
-4*Pi*b)-sin(3*Pi*a+2*Pi*b)-sin(Pi*a+2*Pi*b)-sin(4*Pi*b+3*Pi*a)+sin(2*
Pi*b+7*Pi*a)+sin(6*Pi*b+3*Pi*a)+sin(8*Pi*b+7*Pi*a)+sin(Pi*a+6*Pi*b)-3*
sin(Pi*a)-2*sin(5*Pi*a)+sin(3*Pi*a))/GAMMA(2*a-2*n+2*b+1)/GAMMA(2*n-2*
a-1)/GAMMA(a-1)/GAMMA(b+1)^2/GAMMA(-2*b)/(2*b-1)/(-12*sin(2*Pi*b+2*Pi*
a)-2*sin(2*Pi*a)+11*sin(2*Pi*b)+2*sin(2*Pi*a+6*Pi*b)+10*sin(-2*Pi*b+2*
Pi*a)+2*sin(6*Pi*a)-6*sin(8*Pi*b+6*Pi*a)+sin(8*Pi*b)-sin(10*Pi*b+8*Pi*
a)-2*sin(6*Pi*b+8*Pi*a)+3*sin(8*Pi*b+8*Pi*a)+4*sin(4*Pi*b+6*Pi*a)-sin(
8*Pi*a)+2*sin(10*Pi*b+6*Pi*a)-2*sin(8*Pi*b+2*Pi*a)-sin(10*Pi*b+4*Pi*a)
-2*sin(8*Pi*a+4*Pi*b)+4*sin(8*Pi*b+4*Pi*a)-3*sin(4*Pi*a+6*Pi*b)+2*sin(
4*Pi*a-4*Pi*b)-sin(6*Pi*b)+8*sin(2*Pi*a+4*Pi*b)-5*sin(4*Pi*a-2*Pi*b)-6
*sin(4*Pi*b)-6*sin(4*Pi*a+4*Pi*b)+9*sin(4*Pi*a+2*Pi*b)-4*sin(2*Pi*a-4*
Pi*b)+3*sin(2*Pi*b+8*Pi*a)+4*sin(6*Pi*b+6*Pi*a)-6*sin(2*Pi*b+6*Pi*a));
}{%
\maplemultiline{
198\mbox{:~~~}  - {\displaystyle \frac {1}{2}} \mathrm{
sin}(\pi \,(a + b))\,(b - 1)\,\Gamma (a + 2 + b - 2\,n)\,\Gamma (
2\,n - 2 - a)\,\Gamma ( - 1 + a + b) \\
 \left(  \! {\displaystyle \sum _{k=1}^{2\,n - 3}} \,
{\displaystyle \frac {\Gamma (a + 2 - 2\,n + k)\,\Gamma ( - a - b
 + 1 + k)}{\Gamma (1 - a + k)\,\Gamma (2 + a - 2\,n + b + k)}} 
 \!  \right) /(\pi \,\Gamma (a - 1))\mbox{} + \pi ^{2}\,b\,\Gamma
 (a + 2 + b - 2\,n) \\
\Gamma (2\,n - 2 - a)\,\Gamma ( - 1 + a + b)\,(b - 1)(\mathrm{sin
}(2\,\pi \,b + 5\,\pi \,a) + \mathrm{sin}(\pi \,a + 4\,\pi \,b)
 \\
\mbox{} - \mathrm{sin}(4\,\pi \,b + 7\,\pi \,a) - \mathrm{sin}(6
\,\pi \,b + 7\,\pi \,a) + 3\,\mathrm{sin}(4\,\pi \,b + 5\,\pi \,a
) - \mathrm{sin}(6\,\pi \,b + 5\,\pi \,a) \\
\mbox{} - \mathrm{sin}(8\,\pi \,b + 5\,\pi \,a) + 2\,\mathrm{sin}
(\pi \,a - 4\,\pi \,b) - \mathrm{sin}(3\,\pi \,a + 2\,\pi \,b) - 
\mathrm{sin}(\pi \,a + 2\,\pi \,b) \\
\mbox{} - \mathrm{sin}(4\,\pi \,b + 3\,\pi \,a) + \mathrm{sin}(2
\,\pi \,b + 7\,\pi \,a) + \mathrm{sin}(6\,\pi \,b + 3\,\pi \,a)
 + \mathrm{sin}(8\,\pi \,b + 7\,\pi \,a) \\
\mbox{} + \mathrm{sin}(\pi \,a + 6\,\pi \,b) - 3\,\mathrm{sin}(
\pi \,a) - 2\,\mathrm{sin}(5\,\pi \,a) + \mathrm{sin}(3\,\pi \,a)
) \left/ {\vrule height0.43em width0em depth0.43em} \right. \! 
 \! ( \\
\Gamma (2\,a - 2\,n + 2\,b + 1)\,\Gamma (2\,n - 2\,a - 1)\,\Gamma
 (a - 1)\,\Gamma (b + 1)^{2}\,\Gamma ( - 2\,b)\,(2\,b - 1)( \\
 - 12\,\mathrm{sin}(2\,\pi \,b + 2\,\pi \,a) - 2\,\mathrm{sin}(2
\,\pi \,a) + 11\,\mathrm{sin}(2\,\pi \,b) + 2\,\mathrm{sin}(2\,
\pi \,a + 6\,\pi \,b) \\
\mbox{} + 10\,\mathrm{sin}( - 2\,\pi \,b + 2\,\pi \,a) + 2\,
\mathrm{sin}(6\,\pi \,a) - 6\,\mathrm{sin}(8\,\pi \,b + 6\,\pi \,
a) + \mathrm{sin}(8\,\pi \,b) \\
\mbox{} - \mathrm{sin}(10\,\pi \,b + 8\,\pi \,a) - 2\,\mathrm{sin
}(6\,\pi \,b + 8\,\pi \,a) + 3\,\mathrm{sin}(8\,\pi \,b + 8\,\pi 
\,a) \\
\mbox{} + 4\,\mathrm{sin}(4\,\pi \,b + 6\,\pi \,a) - \mathrm{sin}
(8\,\pi \,a) + 2\,\mathrm{sin}(10\,\pi \,b + 6\,\pi \,a) - 2\,
\mathrm{sin}(8\,\pi \,b + 2\,\pi \,a) \\
\mbox{} - \mathrm{sin}(10\,\pi \,b + 4\,\pi \,a) - 2\,\mathrm{sin
}(8\,\pi \,a + 4\,\pi \,b) + 4\,\mathrm{sin}(8\,\pi \,b + 4\,\pi 
\,a) \\
\mbox{} - 3\,\mathrm{sin}(4\,\pi \,a + 6\,\pi \,b) + 2\,\mathrm{
sin}(4\,\pi \,a - 4\,\pi \,b) - \mathrm{sin}(6\,\pi \,b) + 8\,
\mathrm{sin}(2\,\pi \,a + 4\,\pi \,b) \\
\mbox{} - 5\,\mathrm{sin}(4\,\pi \,a - 2\,\pi \,b) - 6\,\mathrm{
sin}(4\,\pi \,b) - 6\,\mathrm{sin}(4\,\pi \,a + 4\,\pi \,b) + 9\,
\mathrm{sin}(4\,\pi \,a + 2\,\pi \,b) \\
\mbox{} - 4\,\mathrm{sin}(2\,\pi \,a - 4\,\pi \,b) + 3\,\mathrm{
sin}(2\,\pi \,b + 8\,\pi \,a) + 4\,\mathrm{sin}(6\,\pi \,b + 6\,
\pi \,a) \\
\mbox{} - 6\,\mathrm{sin}(2\,\pi \,b + 6\,\pi \,a))) }
}
\end{maplelatex}

\begin{maplelatex}
\mapleinline{inert}{2d}{199, ":   ",
1/2*GAMMA(-3/2+1/2*b+n)*GAMMA(1/2*b-n+3/2)*cos(1/2*Pi*b)*GAMMA(2*n-2-a
)/GAMMA(a-1)*(b-1)*(-1)^n/Pi*Sum(GAMMA(-1/2*b-n+3/2+k)*GAMMA(a+2-2*n+k
)/GAMMA(3/2+1/2*b-n+k)/GAMMA(1-a+k),k = 1 ..
2*n-3)-1/8*Pi^(1/2)*2^(-2*a+2*n+b)*GAMMA(1/2*b-n+3/2)*GAMMA(2-2*n+2*a)
*GAMMA(-1/2*b+a-n+3/2)*GAMMA(-3/2+1/2*b+n)*(2+2*cos(Pi*(2*a+b))-cos(Pi
*(4*a+b))-2*cos(Pi*(2*a-b))+2*cos(Pi*b)-2*cos(4*Pi*a)-cos(2*Pi*(a-b))+
cos(2*Pi*(a+b))-cos(Pi*(4*a-b)))/GAMMA(a-1)/GAMMA(b-1)/GAMMA(-1/2*b+a-
n+2)/GAMMA(3-2*n+a)/(-2*cos(Pi*(2*a-b))+2*cos(Pi*(2*a+b))+2-2*cos(2*Pi
*b)-cos(Pi*(4*a+b))+cos(Pi*(4*a-b))-2*cos(2*Pi*a)+cos(2*Pi*(a-b))+cos(
2*Pi*(a+b)));}{%
\maplemultiline{
199\mbox{:~~~} {\displaystyle \frac {1}{2}} \Gamma ( - 
{\displaystyle \frac {3}{2}}  + {\displaystyle \frac {b}{2}}  + n
)\,\Gamma ({\displaystyle \frac {b}{2}}  - n + {\displaystyle 
\frac {3}{2}} )\,\mathrm{cos}({\displaystyle \frac {\pi \,b}{2}} 
)\,\Gamma (2\,n - 2 - a)\,(b - 1)\,(-1)^{n} \\
 \left(  \! {\displaystyle \sum _{k=1}^{2\,n - 3}} \,
{\displaystyle \frac {\Gamma ( - {\displaystyle \frac {b}{2}}  - 
n + {\displaystyle \frac {3}{2}}  + k)\,\Gamma (a + 2 - 2\,n + k)
}{\Gamma ({\displaystyle \frac {3}{2}}  + {\displaystyle \frac {b
}{2}}  - n + k)\,\Gamma (1 - a + k)}}  \!  \right) /(\Gamma (a - 
1)\,\pi )\mbox{} - {\displaystyle \frac {1}{8}} \sqrt{\pi } \\
2^{( - 2\,a + 2\,n + b)}\,\Gamma ({\displaystyle \frac {b}{2}} 
 - n + {\displaystyle \frac {3}{2}} )\,\Gamma (2 - 2\,n + 2\,a)\,
\Gamma ( - {\displaystyle \frac {b}{2}}  + a - n + 
{\displaystyle \frac {3}{2}} )\,\Gamma ( - {\displaystyle \frac {
3}{2}}  + {\displaystyle \frac {b}{2}}  + n)(2 \\
\mbox{} + 2\,\mathrm{cos}(\pi \,(2\,a + b)) - \mathrm{cos}(\pi \,
(4\,a + b)) - 2\,\mathrm{cos}(\pi \,(2\,a - b)) + 2\,\mathrm{cos}
(\pi \,b) \\
\mbox{} - 2\,\mathrm{cos}(4\,\pi \,a) - \mathrm{cos}(2\,\pi \,(a
 - b)) + \mathrm{cos}(2\,\pi \,(a + b)) - \mathrm{cos}(\pi \,(4\,
a - b))) \left/ {\vrule height0.80em width0em depth0.80em}
 \right. \!  \! ( \\
\Gamma (a - 1)\,\Gamma (b - 1)\,\Gamma ( - {\displaystyle \frac {
b}{2}}  + a - n + 2)\,\Gamma (3 - 2\,n + a)( - 2\,\mathrm{cos}(
\pi \,(2\,a - b)) \\
\mbox{} + 2\,\mathrm{cos}(\pi \,(2\,a + b)) + 2 - 2\,\mathrm{cos}
(2\,\pi \,b) - \mathrm{cos}(\pi \,(4\,a + b)) + \mathrm{cos}(\pi 
\,(4\,a - b)) \\
\mbox{} - 2\,\mathrm{cos}(2\,\pi \,a) + \mathrm{cos}(2\,\pi \,(a
 - b)) + \mathrm{cos}(2\,\pi \,(a + b)))) }
}
\end{maplelatex}

\mapleinline{inert}{2d}{200, ":   ",
1/2*GAMMA(2*n-2-a)*2^(-2*b-1)*GAMMA(-b+1/2)/Pi^(1/2)*GAMMA(1+2*b)*GAMM
A(b-a)*GAMMA(2*b-a)*GAMMA(1+a-2*b)*GAMMA(a)*GAMMA(2*n-2*a-1)*GAMMA(-b+
1)/GAMMA(a-2*b)/GAMMA(2*b)/GAMMA(2*n-2+2*b-2*a)/GAMMA(b-a+1)/GAMMA(-a)
/GAMMA(1+a)/GAMMA(-2*b+1)/GAMMA(2*n-1-a)*Sum(GAMMA(a+2-2*n+k)/GAMMA(-2
*n+2+a)*GAMMA(-a+b+k)/GAMMA(b-a)/GAMMA(1-a+k)*GAMMA(-a+1)/GAMMA(3+a-2*
n-b+k)*GAMMA(a-2*n+3-b),k = 1 ..
2*n-3)+GAMMA(b+1)*GAMMA(2*n-2-a)*GAMMA(a-2*n+3-b)*GAMMA(a-b)*4^(-b)*(-
1/2+1/2*cos(-2*Pi*b+2*Pi*a))/Pi^(1/2)/GAMMA(a-2*b)/GAMMA(b+1/2)/sin(Pi
*(a-2*b));}{%
\maplemultiline{
200\mbox{:~~~} {\displaystyle \frac {1}{2}} \Gamma (2\,n
 - 2 - a)\,2^{( - 2\,b - 1)}\,\Gamma ( - b + {\displaystyle 
\frac {1}{2}} )\,\Gamma (1 + 2\,b)\,\Gamma (b - a)\,\Gamma (2\,b
 - a) \\
\Gamma (1 + a - 2\,b)\,\Gamma (a)\,\Gamma (2\,n - 2\,a - 1)\,
\Gamma ( - b + 1) \\
 \left(  \! {\displaystyle \sum _{k=1}^{2\,n - 3}} \,
{\displaystyle \frac {\Gamma (a + 2 - 2\,n + k)\,\Gamma ( - a + b
 + k)\,\Gamma ( - a + 1)\,\Gamma (a - 2\,n + 3 - b)}{\Gamma ( - 2
\,n + 2 + a)\,\Gamma (b - a)\,\Gamma (1 - a + k)\,\Gamma (3 + a
 - 2\,n - b + k)}}  \!  \right)  \left/ {\vrule 
height0.41em width0em depth0.41em} \right. \!  \! (\sqrt{\pi }
 \\
\Gamma (a - 2\,b)\,\Gamma (2\,b)\,\Gamma (2\,n - 2 + 2\,b - 2\,a)
\,\Gamma (b - a + 1)\,\Gamma ( - a)\,\Gamma (1 + a)\,\Gamma ( - 2
\,b + 1)
\Gamma (2\,n - 1 - a)) \\ \mbox{} + {\displaystyle \frac {\Gamma (b
 + 1)\,\Gamma (2\,n - 2 - a)\,\Gamma (a - 2\,n + 3 - b)\,\Gamma (
a - b)\,4^{( - b)}\,( - {\displaystyle \frac {1}{2}}  + 
{\displaystyle \frac {1}{2}} \,\mathrm{cos}( - 2\,\pi \,b + 2\,
\pi \,a))}{\sqrt{\pi }\,\Gamma (a - 2\,b)\,\Gamma (b + 
{\displaystyle \frac {1}{2}} )\,\mathrm{sin}(\pi \,(a - 2\,b))}} 
 }
}

\end{maplegroup}

%% file: AppendixB201to220.tex
\begin{maplegroup}
\begin{mapleinput}
\end{mapleinput}
\mapleresult
\begin{maplelatex}
\mapleinline{inert}{2d}{201, ":   ",
(-1/2*cos(Pi*b)+1/2*cos(Pi*(2*a+b)))*GAMMA(b+a-4+2*n)*GAMMA(a+b)*GAMMA
(a)/(-cos(Pi*(a-2*b))+cos(Pi*a))/GAMMA(2*n+2*b-4)/GAMMA(a-b)/GAMMA(b)*
Sum(GAMMA(4-a-2*n-b+k)*GAMMA(b-1+k)/GAMMA(-1+b+a+k)/GAMMA(4-b-2*n+k),k
= 1 ..
2*n-3)+(-2^(2*a-2)*Pi^(1/2)*(2*cos(Pi*b)-cos(Pi*(2*a-b))-cos(Pi*(2*a+b
)))/(sin(2*Pi*(a-b))+sin(2*Pi*b)-sin(2*Pi*a))+1/4*(2^a)^2*Pi^(1/2)/sin
(Pi*b))*GAMMA(a+b)*GAMMA(b+a-4+2*n)/GAMMA(2*a+2*b-5+2*n)*GAMMA(a-1/2)/
GAMMA(b)/GAMMA(a-b);}{%
\maplemultiline{
201\mbox{:~~~} ( - {\displaystyle \frac {1}{2}} \,
\mathrm{cos}(\pi \,b) + {\displaystyle \frac {1}{2}} \,\mathrm{
cos}(\pi \,(2\,a + b)))\,\Gamma (b + a - 4 + 2\,n)\,\Gamma (a + b
)\,\Gamma (a) \\
 \left(  \! {\displaystyle \sum _{k=1}^{2\,n - 3}} \,
{\displaystyle \frac {\Gamma (4 - a - 2\,n - b + k)\,\Gamma (b - 
1 + k)}{\Gamma ( - 1 + b + a + k)\,\Gamma (4 - b - 2\,n + k)}} 
 \!  \right) /(( - \mathrm{cos}(\pi \,(a - 2\,b)) + \mathrm{cos}(
\pi \,a)) \\
\Gamma (2\,n + 2\,b - 4)\,\Gamma (a - b)\,\Gamma (b))\mbox{} + 
 \\
( - {\displaystyle \frac {2^{(2\,a - 2)}\,\sqrt{\pi }\,(2\,
\mathrm{cos}(\pi \,b) - \mathrm{cos}(\pi \,(2\,a - b)) - \mathrm{
cos}(\pi \,(2\,a + b)))}{\mathrm{sin}(2\,\pi \,(a - b)) + 
\mathrm{sin}(2\,\pi \,b) - \mathrm{sin}(2\,\pi \,a)}}  + 
{\displaystyle \frac {1}{4}} \,{\displaystyle \frac {(2^{a})^{2}
\,\sqrt{\pi }}{\mathrm{sin}(\pi \,b)}} ) \\
\Gamma (a + b)\,\Gamma (b + a - 4 + 2\,n)\,\Gamma (a - 
{\displaystyle \frac {1}{2}} )/(\Gamma (2\,a + 2\,b - 5 + 2\,n)\,
\Gamma (b)\,\Gamma (a - b)) }
}
\end{maplelatex}

\begin{maplelatex}
\mapleinline{inert}{2d}{210, ":   ",
(-1/4*sin(Pi*(2*b+a))-1/2*sin(Pi*a)+1/4*sin(Pi*(3*a+2*b)))*(a-1)*GAMMA
(a-n+b+1)*GAMMA(-1+a+b)*GAMMA(-b+1)/GAMMA(-n+b+1)/Pi/(cos(Pi*b)-cos(Pi
*(2*a+b)))*Sum(GAMMA(-a-b+1+k)*GAMMA(b+1+k-n)/GAMMA(1+b+a-n+k)/GAMMA(1
-b+k),k = 1 ..
n-2)+1/8*4^a*GAMMA(n-b)*GAMMA(a-n+b+1)*GAMMA(a-1/2)*GAMMA(-a+1)*GAMMA(
-1+a+b)*(cos(Pi*(2*a+b))-cos(Pi*b))*(a-1)/Pi^(1/2)/GAMMA(n-2*b)/GAMMA(
-n+2*b+2*a)/GAMMA(b)/(cos(Pi*a)-cos(Pi*(2*b+a)))-1/2*Pi^(1/2)*4^a*GAMM
A(-b+1)*GAMMA(a-1/2)*GAMMA(-1+a+b)*(-1)^n*GAMMA(n-b)*(cos(Pi*a)-cos(Pi
*(3*a+2*b)))*(a-1)/GAMMA(-n+2*b+2*a)/GAMMA(n-2*b)/GAMMA(n-a-b)/GAMMA(a
)/(sin(Pi*(2*a+b))-sin(Pi*b)-sin(Pi*(4*a+3*b))+sin(Pi*(3*b+2*a)));}{%
\maplemultiline{
210\mbox{:~~~} ( - {\displaystyle \frac {1}{4}} \,
\mathrm{sin}(\pi \,(2\,b + a)) - {\displaystyle \frac {1}{2}} \,
\mathrm{sin}(\pi \,a) + {\displaystyle \frac {1}{4}} \,\mathrm{
sin}(\pi \,(3\,a + 2\,b)))\,(a - 1)\,\Gamma (a - n + b + 1) \\
\Gamma ( - 1 + a + b)\,\Gamma ( - b + 1)\, \left(  \! 
{\displaystyle \sum _{k=1}^{n - 2}} \,{\displaystyle \frac {
\Gamma ( - a - b + 1 + k)\,\Gamma (b + 1 + k - n)}{\Gamma (1 + b
 + a - n + k)\,\Gamma (1 - b + k)}}  \!  \right) /(\Gamma ( - n
 + b + 1) \\
\pi \,(\mathrm{cos}(\pi \,b) - \mathrm{cos}(\pi \,(2\,a + b))))
\mbox{} + {\displaystyle \frac {1}{8}} 4^{a}\,\Gamma (n - b)\,
\Gamma (a - n + b + 1)\,\Gamma (a - {\displaystyle \frac {1}{2}} 
)\,\Gamma ( - a + 1) \\
\Gamma ( - 1 + a + b)\,(\mathrm{cos}(\pi \,(2\,a + b)) - \mathrm{
cos}(\pi \,b))\,(a - 1) \left/ {\vrule 
height0.41em width0em depth0.41em} \right. \!  \! (\sqrt{\pi }\,
\Gamma (n - 2\,b) \\
\Gamma ( - n + 2\,b + 2\,a)\,\Gamma (b)\,(\mathrm{cos}(\pi \,a)
 - \mathrm{cos}(\pi \,(2\,b + a))))\mbox{} - {\displaystyle 
\frac {1}{2}} \sqrt{\pi }\,4^{a}\,\Gamma ( - b + 1)\,\Gamma (a - 
{\displaystyle \frac {1}{2}} ) \\
\Gamma ( - 1 + a + b)\,(-1)^{n}\,\Gamma (n - b)\,(\mathrm{cos}(
\pi \,a) - \mathrm{cos}(\pi \,(3\,a + 2\,b)))\,(a - 1)/( \\
\Gamma ( - n + 2\,b + 2\,a)\,\Gamma (n - 2\,b)\,\Gamma (n - a - b
)\,\Gamma (a) \\
(\mathrm{sin}(\pi \,(2\,a + b)) - \mathrm{sin}(\pi \,b) - 
\mathrm{sin}(\pi \,(4\,a + 3\,b)) + \mathrm{sin}(\pi \,(3\,b + 2
\,a)))) }
}
\end{maplelatex}

\begin{maplelatex}
\mapleinline{inert}{2d}{211, ":   ",
1/2*Sum(GAMMA(b-1+k)*GAMMA(-a-b+3+k-n)/GAMMA(3-b-n+k)/GAMMA(-1+b+a+k),
k = 1 ..
n-2)*sin(Pi*(a+b))*(a-1)/sin(Pi*b)/GAMMA(b)/GAMMA(n+2*b-3)*GAMMA(n+a+2
*b-3)*GAMMA(-1+a+b)+2^(2*a-2)*Pi^(1/2)*GAMMA(a-1/2)*GAMMA(n+a+2*b-3)*G
AMMA(-1+a+b)/GAMMA(a-1)/GAMMA(b)/GAMMA(n+2*a+2*b-4)*cos(Pi*a)/sin(Pi*a
);}{%
\maplemultiline{
211\mbox{:~~~} {\displaystyle \frac {1}{2}}  \left(  \! 
{\displaystyle \sum _{k=1}^{n - 2}} \,{\displaystyle \frac {
\Gamma (b - 1 + k)\,\Gamma ( - a - b + 3 + k - n)}{\Gamma (3 - b
 - n + k)\,\Gamma ( - 1 + b + a + k)}}  \!  \right) \,\mathrm{sin
}(\pi \,(a + b))\,(a - 1) \\
\Gamma (n + a + 2\,b - 3)\,\Gamma ( - 1 + a + b)/(\mathrm{sin}(
\pi \,b)\,\Gamma (b)\,\Gamma (n + 2\,b - 3)) \\
\mbox{} + {\displaystyle \frac {2^{(2\,a - 2)}\,\sqrt{\pi }\,
\Gamma (a - {\displaystyle \frac {1}{2}} )\,\Gamma (n + a + 2\,b
 - 3)\,\Gamma ( - 1 + a + b)\,\mathrm{cos}(\pi \,a)}{\Gamma (a - 
1)\,\Gamma (b)\,\Gamma (n + 2\,a + 2\,b - 4)\,\mathrm{sin}(\pi \,
a)}}  }
}
\end{maplelatex}

\begin{maplelatex}
\mapleinline{inert}{2d}{212, ":   ",
1/2*b*sin(Pi*a)*Sum(GAMMA(-a+b+k)*GAMMA(a+1+k-n)/GAMMA(2+a-b-n+k)/GAMM
A(1-a+k),k = 1 ..
n-2)/Pi/GAMMA(n-a+b-1)/GAMMA(n-2*a+2*b-1)*GAMMA(n+2*b-1-a)*GAMMA(n-2*a
)*GAMMA(-a+1)*GAMMA(a-b)-1/4*4^(-b)*(cos(Pi*(2*a-3*b))+cos(Pi*(2*a-b))
-2*cos(Pi*b))/Pi^(1/2)/sin(Pi*(a-2*b))/sin(Pi*a)/GAMMA(-2*b+2-n+a)*GAM
MA(a-n-b+2)*GAMMA(b+1)*GAMMA(-b+1/2)*GAMMA(a-b)/GAMMA(a);}{%
\maplemultiline{
212\mbox{:~~~} {\displaystyle \frac {1}{2}} b\,\mathrm{
sin}(\pi \,a)\, \left(  \! {\displaystyle \sum _{k=1}^{n - 2}} \,
{\displaystyle \frac {\Gamma ( - a + b + k)\,\Gamma (a + 1 + k - 
n)}{\Gamma (2 + a - b - n + k)\,\Gamma (1 - a + k)}}  \! 
 \right) \,\Gamma (n + 2\,b - 1 - a) \\
\Gamma (n - 2\,a)\,\Gamma ( - a + 1)\,\Gamma (a - b)/(\pi \,
\Gamma (n - a + b - 1)\,\Gamma (n - 2\,a + 2\,b - 1))\mbox{} - 
{\displaystyle \frac {1}{4}} 4^{( - b)} \\
(\mathrm{cos}(\pi \,(2\,a - 3\,b)) + \mathrm{cos}(\pi \,(2\,a - b
)) - 2\,\mathrm{cos}(\pi \,b))\,\Gamma (a - n - b + 2)\,\Gamma (b
 + 1) \\
\Gamma ( - b + {\displaystyle \frac {1}{2}} )\,\Gamma (a - b)
 \left/ {\vrule height0.41em width0em depth0.41em} \right. \! 
 \! (\sqrt{\pi }\,\mathrm{sin}(\pi \,(a - 2\,b))\,\mathrm{sin}(
\pi \,a)\,\Gamma ( - 2\,b + 2 - n + a)\,\Gamma (a)) }
}
\end{maplelatex}

\begin{maplelatex}
\mapleinline{inert}{2d}{213, ":   ",
1/4*(2*a+n-3)*b*sin(Pi*(-1/2*b+a))*Sum(GAMMA(a-1+k)*GAMMA(1/2*b+2-a+k-
n)/GAMMA(3-n-a+k)/GAMMA(a-1/2*b+k),k = 1 ..
n-2)/sin(Pi*a)/GAMMA(a)/GAMMA(n+a-2)*GAMMA(a-2+n-1/2*b)*GAMMA(-1/2*b+a
)-1/2*(n+a-2)*GAMMA(1/2*b+1)*GAMMA(n+2*a-2)*GAMMA(a-2+n-1/2*b)*GAMMA(-
1/2*b+a)*GAMMA(-1/2*b+1/2)*GAMMA(2-a-n)*(-2^(-b+1)*sin(1/2*Pi*(b+2*n+2
*a))+2^(-b+1)*sin(1/2*Pi*(b+2*n+2*a))*cos(Pi*a)^2+2^(-b)*sin(2*Pi*a)*s
in(1/2*Pi*b)*sin(Pi*(a+n)))/Pi^(3/2)/GAMMA(a)/GAMMA(n-b-2+2*a)/(-1+cos
(Pi*a))/(cos(Pi*a)+1);}{%
\maplemultiline{
213\mbox{:~~~} {\displaystyle \frac {1}{4}} (2\,a + n - 
3)\,b\,\mathrm{sin}(\pi \,( - {\displaystyle \frac {b}{2}}  + a))
\, \left(  \! {\displaystyle \sum _{k=1}^{n - 2}} \,
{\displaystyle \frac {\Gamma (a - 1 + k)\,\Gamma ({\displaystyle 
\frac {b}{2}}  + 2 - a + k - n)}{\Gamma (3 - n - a + k)\,\Gamma (
a - {\displaystyle \frac {b}{2}}  + k)}}  \!  \right)  \\
\Gamma (a - 2 + n - {\displaystyle \frac {b}{2}} )\,\Gamma ( - 
{\displaystyle \frac {b}{2}}  + a)/(\mathrm{sin}(\pi \,a)\,\Gamma
 (a)\,\Gamma (n + a - 2))\mbox{} - {\displaystyle \frac {1}{2}} (
n + a - 2)\,\Gamma ({\displaystyle \frac {b}{2}}  + 1) \\
\Gamma (n + 2\,a - 2)\,\Gamma (a - 2 + n - {\displaystyle \frac {
b}{2}} )\,\Gamma ( - {\displaystyle \frac {b}{2}}  + a)\,\Gamma (
 - {\displaystyle \frac {b}{2}}  + {\displaystyle \frac {1}{2}} )
\,\Gamma (2 - a - n)( \\
 - 2^{( - b + 1)}\,\mathrm{sin}({\displaystyle \frac {\pi \,(b + 
2\,n + 2\,a)}{2}} ) + 2^{( - b + 1)}\,\mathrm{sin}(
{\displaystyle \frac {\pi \,(b + 2\,n + 2\,a)}{2}} )\,\mathrm{cos
}(\pi \,a)^{2} \\
\mbox{} + 2^{( - b)}\,\mathrm{sin}(2\,\pi \,a)\,\mathrm{sin}(
{\displaystyle \frac {\pi \,b}{2}} )\,\mathrm{sin}(\pi \,(a + n))
) \left/ {\vrule height0.63em width0em depth0.63em} \right. \! 
 \! (\pi ^{(3/2)}\,\Gamma (a)\,\Gamma (n - b - 2 + 2\,a) 
( - 1 + \mathrm{cos}(\pi \,a))\,(\mathrm{cos}(\pi \,a) + 1)) }
}
\end{maplelatex}

\begin{maplelatex}
\mapleinline{inert}{2d}{214, ":   ",
-1/2*Sum(GAMMA(-a+b+2+k-n)*GAMMA(-b+k)/GAMMA(2+b-n+k)/GAMMA(-b+a+k),k
= 1 ..
n-2)/Pi^3*GAMMA(n+2*a-2*b-2)*GAMMA(-n-a+2*b+3)*GAMMA(-b+n+a-2)*GAMMA(-
2*a+2)*GAMMA(b-n+2)*GAMMA(a-b)*GAMMA(a)/GAMMA(-b)*sin(Pi*a)*sin(-2*Pi*
b+2*Pi*a)*sin(Pi*a-Pi*b)+1/4*(-cos(Pi*(3*a-2*b))-cos(Pi*(a-2*b))+2*cos
(Pi*a))/Pi^2/sin(2*Pi*a)*GAMMA(-n-a+2*b+3)*GAMMA(-b+n+a-2)*sin(2*Pi*b)
*GAMMA(n-2*b-1)*GAMMA(b-n+2)*GAMMA(b+1)*GAMMA(a-b)/GAMMA(a);}{%
\maplemultiline{
214\mbox{:~~~}  - {\displaystyle \frac {1}{2}}  \left( 
 \! {\displaystyle \sum _{k=1}^{n - 2}} \,{\displaystyle \frac {
\Gamma ( - a + b + 2 + k - n)\,\Gamma ( - b + k)}{\Gamma (2 + b
 - n + k)\,\Gamma ( - b + a + k)}}  \!  \right) \,\Gamma (n + 2\,
a - 2\,b - 2) \\
\Gamma ( - n - a + 2\,b + 3)\,\Gamma ( - b + n + a - 2)\,\Gamma (
 - 2\,a + 2)\,\Gamma (b - n + 2)\,\Gamma (a - b)\,\Gamma (a) \\
\mathrm{sin}(\pi \,a)\,\mathrm{sin}( - 2\,\pi \,b + 2\,\pi \,a)\,
\mathrm{sin}(\pi \,a - \pi \,b) \left/ {\vrule 
height0.43em width0em depth0.43em} \right. \!  \! (\pi ^{3}\,
\Gamma ( - b))\mbox{} + {\displaystyle \frac {1}{4}}  \\
( - \mathrm{cos}(\pi \,(3\,a - 2\,b)) - \mathrm{cos}(\pi \,(a - 2
\,b)) + 2\,\mathrm{cos}(\pi \,a))\,\Gamma ( - n - a + 2\,b + 3)
 \\
\Gamma ( - b + n + a - 2)\,\mathrm{sin}(2\,\pi \,b)\,\Gamma (n - 
2\,b - 1)\,\Gamma (b - n + 2)\,\Gamma (b + 1)\,\Gamma (a - b)/(
\pi ^{2}
\mathrm{sin}(2\,\pi \,a)\,\Gamma (a)) }
}
\end{maplelatex}

\begin{maplelatex}
\mapleinline{inert}{2d}{215, ":   ",
1/2*Sum(GAMMA(a+1+k-n)*GAMMA(b+1+k-n)/GAMMA(1-a+k)/GAMMA(1-b+k),k = 1
..
n-2)*sin(2*Pi*b)*sin(Pi*(-a-b))/sin(Pi*b)/GAMMA(b)/Pi*GAMMA(a+2*b-n+1)
*GAMMA(n-2*b)*GAMMA(n-a-b)/GAMMA(a-n+1)*GAMMA(-a+1)-(1/2*Pi^(1/2)*2^(2
*n-2*a-2*b-2)*(-1)^n*sin(Pi*(a+b))/sin(Pi*b)*sin(2*Pi*b)-1/2*Pi^(1/2)*
2^(2*n-2*b-2*a-1)*(-1)^n*sin(Pi*(2*b+a)))/sin(Pi*b)/sin(Pi*a)*GAMMA(a+
2*b-n+1)*(-1)^n*Pi/sin(Pi*(a+b+1/2))/GAMMA(3/2-n+a+b)/GAMMA(n-2*a)*GAM
MA(n-a)/GAMMA(b)/GAMMA(a);}{%
\maplemultiline{
215\mbox{:~~~} {\displaystyle \frac {1}{2}}  \left(  \! 
{\displaystyle \sum _{k=1}^{n - 2}} \,{\displaystyle \frac {
\Gamma (a + 1 + k - n)\,\Gamma (b + 1 + k - n)}{\Gamma (1 - a + k
)\,\Gamma (1 - b + k)}}  \!  \right) \,\mathrm{sin}(2\,\pi \,b)\,
\mathrm{sin}(\pi \,( - a - b)) \\
\Gamma (a + 2\,b - n + 1)\,\Gamma (n - 2\,b)\,\Gamma (n - a - b)
\,\Gamma ( - a + 1)/(\mathrm{sin}(\pi \,b)\,\Gamma (b)\,\pi \,
\Gamma (a - n + 1)) \\
\mbox{} - ({\displaystyle \frac {1}{2}} \,{\displaystyle \frac {
\sqrt{\pi }\,2^{(2\,n - 2\,a - 2\,b - 2)}\,(-1)^{n}\,\mathrm{sin}
(\pi \,(a + b))\,\mathrm{sin}(2\,\pi \,b)}{\mathrm{sin}(\pi \,b)}
}  \\
\mbox{} - {\displaystyle \frac {1}{2}} \,\sqrt{\pi }\,2^{(2\,n - 
2\,b - 2\,a - 1)}\,(-1)^{n}\,\mathrm{sin}(\pi \,(2\,b + a)))
\Gamma (a + 2\,b - n + 1)\,(-1)^{n}\,\pi  \\
\Gamma (n - a) \left/ {\vrule height0.80em width0em depth0.80em}
 \right. \!  \! (\mathrm{sin}(\pi \,b)\,\mathrm{sin}(\pi \,a)\,
\mathrm{sin}(\pi \,(a + b + {\displaystyle \frac {1}{2}} ))\,
\Gamma ({\displaystyle \frac {3}{2}}  - n + a + b)\,\Gamma (n - 2
\,a) 
\Gamma (b)\,\Gamma (a)) }
}
\end{maplelatex}

\begin{maplelatex}
\mapleinline{inert}{2d}{216, ":   ",
(1/2*n-b-1/2)*sin(Pi*(n-b))/sin(Pi*(b+1))*GAMMA(1/2*a-1/2*n+1)/GAMMA(-
1/2*n-1/2*a+1)*GAMMA(n-b-1)/GAMMA(b)*Sum(GAMMA(b+1+k-n)*GAMMA(-1/2*a-1
/2*n+1+k)/GAMMA(1-b+k)/GAMMA(1+1/2*a-1/2*n+k),k = 1 ..
n-2)+(-1/2*Pi^(1/2)*2^(a+n-2*b-2)*sin(2*Pi*b)/sin(Pi*b)+1/2*Pi^(1/2)*2
^(n+a-1)*4^(-b)*sin(Pi*(1/2*a+1/2*n-2*b))/sin(Pi*(1/2*a+1/2*n-b)))/sin
(Pi*b)*GAMMA(1/2*a+1/2*n-b-1/2)/GAMMA(1/2*a+1/2*n-b)*GAMMA(1/2*n+1/2*a
)*GAMMA(1/2*a-1/2*n+1)*GAMMA(-n+2*b+2)/GAMMA(b-n+2)/GAMMA(b)/GAMMA(a);
}{%
\maplemultiline{
216\mbox{:~~~} ({\displaystyle \frac {n}{2}}  - b - 
{\displaystyle \frac {1}{2}} )\,\mathrm{sin}(\pi \,(n - b))\,
\Gamma ({\displaystyle \frac {a}{2}}  - {\displaystyle \frac {n}{
2}}  + 1)\,\Gamma (n - b - 1) \\
 \left(  \! {\displaystyle \sum _{k=1}^{n - 2}} \,{\displaystyle 
\frac {\Gamma (b + 1 + k - n)\,\Gamma ( - {\displaystyle \frac {a
}{2}}  - {\displaystyle \frac {n}{2}}  + 1 + k)}{\Gamma (1 - b + 
k)\,\Gamma (1 + {\displaystyle \frac {a}{2}}  - {\displaystyle 
\frac {n}{2}}  + k)}}  \!  \right)  \left/ {\vrule 
height0.80em width0em depth0.80em} \right. \!  \! (\mathrm{sin}(
\pi \,(b + 1)) 
\Gamma ( - {\displaystyle \frac {n}{2}}  - {\displaystyle \frac {
a}{2}}  + 1)\,\Gamma (b))\mbox{} +  \\
 \left(  \!  - {\displaystyle \frac {1}{2}} \,{\displaystyle 
\frac {\sqrt{\pi }\,2^{(a + n - 2\,b - 2)}\,\mathrm{sin}(2\,\pi 
\,b)}{\mathrm{sin}(\pi \,b)}}  + {\displaystyle \frac {1}{2}} \,
{\displaystyle \frac {\sqrt{\pi }\,2^{(n + a - 1)}\,4^{( - b)}\,
\mathrm{sin}(\pi \,({\displaystyle \frac {a}{2}}  + 
{\displaystyle \frac {n}{2}}  - 2\,b))}{\mathrm{sin}(\pi \,(
{\displaystyle \frac {a}{2}}  + {\displaystyle \frac {n}{2}}  - b
))}}  \!  \right)  \\
\Gamma ({\displaystyle \frac {a}{2}}  + {\displaystyle \frac {n}{
2}}  - b - {\displaystyle \frac {1}{2}} )\,\Gamma (
{\displaystyle \frac {n}{2}}  + {\displaystyle \frac {a}{2}} )\,
\Gamma ({\displaystyle \frac {a}{2}}  - {\displaystyle \frac {n}{
2}}  + 1)\,\Gamma ( - n + 2\,b + 2) \left/ {\vrule 
height0.80em width0em depth0.80em} \right. \!  \! (\mathrm{sin}(
\pi \,b) 
\Gamma ({\displaystyle \frac {a}{2}}  + {\displaystyle \frac {n}{
2}}  - b)\,\Gamma (b - n + 2)\,\Gamma (b)\,\Gamma (a)) }
}
\end{maplelatex}

\begin{maplelatex}
\mapleinline{inert}{2d}{217, ":   ",
-1/2*sin(Pi*(2*b+2*a))*sin(Pi*(a+b))*sin(Pi*(b+1))*Sum(GAMMA(1+b+a-n+k
)*GAMMA(-b+k)/GAMMA(-a-b+1+k)/GAMMA(2+b-n+k),k = 1 ..
n-2)/Pi^2/sin(Pi*(b+a+1))/GAMMA(a)*GAMMA(n-b-1-a)*GAMMA(2*a+b)/GAMMA(a
+b)*GAMMA(n-2*b-2*a)*GAMMA(-n+a+2*b+2)*GAMMA(b-n+2)+1/2*4^(-a)*(-sin(P
i*(3*b+2*a))+3*sin(Pi*b)-sin(Pi*(4*a+3*b))+3*sin(Pi*(2*a+b)))/Pi^(1/2)
/(3*sin(Pi*(a+b))-sin(Pi*(3*a+3*b)))*GAMMA(2*a+b)*GAMMA(-n+a+2*b+2)*GA
MMA(1/2-a)/GAMMA(2+a-n+b)/GAMMA(a+b)/GAMMA(-n+2*b+2)*GAMMA(b-n+2);}{%
\maplemultiline{
217\mbox{:~~~}  - {\displaystyle \frac {1}{2}} \mathrm{
sin}(\pi \,(2\,b + 2\,a))\,\mathrm{sin}(\pi \,(a + b))\,\mathrm{
sin}(\pi \,(b + 1)) \\
 \left(  \! {\displaystyle \sum _{k=1}^{n - 2}} \,{\displaystyle 
\frac {\Gamma (1 + b + a - n + k)\,\Gamma ( - b + k)}{\Gamma ( - 
a - b + 1 + k)\,\Gamma (2 + b - n + k)}}  \!  \right) \,\Gamma (n
 - b - 1 - a)\,\Gamma (2\,a + b) \\
\Gamma (n - 2\,b - 2\,a)\,\Gamma ( - n + a + 2\,b + 2)\,\Gamma (b
 - n + 2) \left/ {\vrule height0.43em width0em depth0.43em}
 \right. \!  \! (\pi ^{2}\,\mathrm{sin}(\pi \,(b + a + 1))\,
\Gamma (a) \,
\Gamma (a + b))\mbox{} \\ + {\displaystyle \frac {1}{2}} 4^{( - a)}
( - \mathrm{sin}(\pi \,(3\,b + 2\,a)) + 3\,\mathrm{sin}(\pi \,b)
 - \mathrm{sin}(\pi \,(4\,a + 3\,b)) + 3\,\mathrm{sin}(\pi \,(2\,
a + b))) \\
\Gamma (2\,a + b)\,\Gamma ( - n + a + 2\,b + 2)\,\Gamma (
{\displaystyle \frac {1}{2}}  - a)\,\Gamma (b - n + 2) \left/ 
{\vrule height0.41em width0em depth0.41em} \right. \!  \! (\sqrt{
\pi } \\
(3\,\mathrm{sin}(\pi \,(a + b)) - \mathrm{sin}(\pi \,(3\,a + 3\,b
)))\,\Gamma (2 + a - n + b)\,\Gamma (a + b)\,\Gamma ( - n + 2\,b
 + 2)) }
}
\end{maplelatex}

\begin{maplelatex}
\mapleinline{inert}{2d}{218, ":   ",
1/2*Sum(GAMMA(a-1+k)*GAMMA(-1/2*b-1/2*n+1+k)/GAMMA(3-n-a+k)/GAMMA(1+1/
2*b-1/2*n+k),k = 1 ..
n-2)*Pi*(2*a+n-3)/sin(Pi*a)/GAMMA(a)/GAMMA(n+a-2)/GAMMA(1/2*n+a-1/2*b-
1)*GAMMA(1/2*n+1/2*b+a-1)-2*2^(b-n-2*a)*(2*sin(1/2*Pi*(4*a-b+n))*sin(P
i*a)-sin(2*Pi*a)*sin(1/2*Pi*(n+2*a-b)))*GAMMA(1/2*n+1/2*b+a-1)*GAMMA(-
1/2*n-a+3/2+1/2*b)*GAMMA(n+2*a-2)*Pi^(1/2)/sin(1/2*Pi*(n-b))/sin(Pi*a)
^2/GAMMA(n+a-2)/GAMMA(b)/GAMMA(a);}{%
\maplemultiline{
218\mbox{:~~~}
{\displaystyle \frac {1}{2}} \,{\displaystyle \frac { \left(  \! 
{\displaystyle \sum _{k=1}^{n - 2}} \,{\displaystyle \frac {
\Gamma (a - 1 + k)\,\Gamma ( - {\displaystyle \frac {b}{2}}  - 
{\displaystyle \frac {n}{2}}  + 1 + k)}{\Gamma (3 - n - a + k)\,
\Gamma (1 + {\displaystyle \frac {b}{2}}  - {\displaystyle 
\frac {n}{2}}  + k)}}  \!  \right) \,\pi \,(2\,a + n - 3)\,\Gamma
 ({\displaystyle \frac {n}{2}}  + {\displaystyle \frac {b}{2}} 
 + a - 1)}{\mathrm{sin}(\pi \,a)\,\Gamma (a)\,\Gamma (n + a - 2)
\,\Gamma ({\displaystyle \frac {n}{2}}  + a - {\displaystyle 
\frac {b}{2}}  - 1)}}  \\
\mbox{} - 2\,2^{(b - n - 2\,a)}\,(2\,\mathrm{sin}({\displaystyle 
\frac {\pi \,(4\,a - b + n)}{2}} )\,\mathrm{sin}(\pi \,a) - 
\mathrm{sin}(2\,\pi \,a)\,\mathrm{sin}({\displaystyle \frac {\pi 
\,(n + 2\,a - b)}{2}} )) \\
\Gamma ({\displaystyle \frac {n}{2}}  + {\displaystyle \frac {b}{
2}}  + a - 1)\,\Gamma ( - {\displaystyle \frac {n}{2}}  - a + 
{\displaystyle \frac {3}{2}}  + {\displaystyle \frac {b}{2}} )\,
\Gamma (n + 2\,a - 2)\,\sqrt{\pi } \left/ {\vrule 
height0.80em width0em depth0.80em} \right. \!  \! (\mathrm{sin}(
{\displaystyle \frac {\pi \,(n - b)}{2}} )
\mathrm{sin}(\pi \,a)^{2}\,\Gamma (n + a - 2)\,\Gamma (b)\,\Gamma
 (a)) }
}
\end{maplelatex}

\begin{maplelatex}
\mapleinline{inert}{2d}{219, ":   ",
-1/2*(2*b-1)*Sum(GAMMA(3/2-b-n+1/2*a+k)*GAMMA(-1/2*a-n+3/2+k)/GAMMA(3/
2-1/2*a+b-n+k)/GAMMA(3/2+1/2*a-n+k),k = 1 ..
2*n-3)*sin(1/2*Pi*(1+a-2*b))*sin(1/2*Pi*(a-1))/sin(1/2*Pi*(2*n+a-2*b-1
))/Pi*GAMMA(5/2+1/2*a-n)*GAMMA(n-3/2+b-1/2*a)*GAMMA(n-3/2+1/2*a)/GAMMA
(n+1/2*a-b-1/2)+1/16*Pi^(1/2)*GAMMA(5/2+1/2*a-n)*GAMMA(1+a-2*b)*GAMMA(
b+1/2)*GAMMA(n-3/2+1/2*a)*4^b*(-sin(2*Pi*(a-b))+sin(2*Pi*(a-2*b))+sin(
2*Pi*b))/GAMMA(-b+5/2+1/2*a-n)/GAMMA(n+1/2*a-b-1/2)/GAMMA(b)/GAMMA(a)/
cos(Pi*(1/2*a-b))^2/sin(Pi*(a-b))/sin(Pi*b);}{%
\maplemultiline{
219\mbox{:~~~}  - {\displaystyle \frac {1}{2}} (2\,b - 1
)\, \left(  \! {\displaystyle \sum _{k=1}^{2\,n - 3}} \,
{\displaystyle \frac {\Gamma ({\displaystyle \frac {3}{2}}  - b
 - n + {\displaystyle \frac {a}{2}}  + k)\,\Gamma ( - 
{\displaystyle \frac {a}{2}}  - n + {\displaystyle \frac {3}{2}} 
 + k)}{\Gamma ({\displaystyle \frac {3}{2}}  - {\displaystyle 
\frac {a}{2}}  + b - n + k)\,\Gamma ({\displaystyle \frac {3}{2}
}  + {\displaystyle \frac {a}{2}}  - n + k)}}  \!  \right)  \\
\mathrm{sin}({\displaystyle \frac {\pi \,(1 + a - 2\,b)}{2}} )\,
\mathrm{sin}({\displaystyle \frac {\pi \,(a - 1)}{2}} )\,\Gamma (
{\displaystyle \frac {5}{2}}  + {\displaystyle \frac {a}{2}}  - n
)\,\Gamma (n - {\displaystyle \frac {3}{2}}  + b - 
{\displaystyle \frac {a}{2}} )\,\Gamma (n - {\displaystyle 
\frac {3}{2}}  + {\displaystyle \frac {a}{2}} ) \\
 \left/ {\vrule height0.80em width0em depth0.80em} \right. \! 
 \! (\mathrm{sin}({\displaystyle \frac {\pi \,(2\,n + a - 2\,b - 
1)}{2}} )\,\pi \,\Gamma (n + {\displaystyle \frac {a}{2}}  - b - 
{\displaystyle \frac {1}{2}} ))\mbox{} + {\displaystyle \frac {1
}{16}} \sqrt{\pi }\,\Gamma ({\displaystyle \frac {5}{2}}  + 
{\displaystyle \frac {a}{2}}  - n) \\
\Gamma (1 + a - 2\,b)\,\Gamma (b + {\displaystyle \frac {1}{2}} )
\,\Gamma (n - {\displaystyle \frac {3}{2}}  + {\displaystyle 
\frac {a}{2}} )\,4^{b} \\
( - \mathrm{sin}(2\,\pi \,(a - b)) + \mathrm{sin}(2\,\pi \,(a - 2
\,b)) + \mathrm{sin}(2\,\pi \,b)) \left/ {\vrule 
height0.80em width0em depth0.80em} \right. \!  \! (\Gamma ( - b
 + {\displaystyle \frac {5}{2}}  + {\displaystyle \frac {a}{2}} 
 - n) \\
\Gamma (n + {\displaystyle \frac {a}{2}}  - b - {\displaystyle 
\frac {1}{2}} )\,\Gamma (b)\,\Gamma (a)\,\mathrm{cos}(\pi \,(
{\displaystyle \frac {a}{2}}  - b))^{2}\,\mathrm{sin}(\pi \,(a - 
b))\,\mathrm{sin}(\pi \,b)) }
}
\end{maplelatex}

\begin{maplelatex}
\mapleinline{inert}{2d}{220, ":   ",
-1/2*(2*b-1)*Sum(GAMMA(3/2-b-n+1/2*a+k)*GAMMA(-1/2*a-n+3/2+k)/GAMMA(3/
2-1/2*a+b-n+k)/GAMMA(3/2+1/2*a-n+k),k = 1 ..
2*n-3)*sin(1/2*Pi*(1+a-2*b))*sin(1/2*Pi*(a-1))/sin(1/2*Pi*(a-3+2*n))/P
i*GAMMA(n-3/2+b-1/2*a)*GAMMA(b+1/2*a-n+3/2)+(sin(2*Pi*a-4*Pi*b)+sin(2*
Pi*b)-sin(-2*Pi*b+2*Pi*a))/(-8*cos(Pi*a-3*Pi*b)+8*cos(-Pi*b+2*Pi*a)-8*
cos(-3*Pi*b+2*Pi*a)+8*cos(Pi*a+Pi*b))*Pi^(1/2)*2^(2*b+2)*GAMMA(b+1/2)*
GAMMA(b+1/2*a-n+3/2)*GAMMA(1+a-2*b)/GAMMA(-b+5/2+1/2*a-n)/GAMMA(b)/GAM
MA(a);}{%
\maplemultiline{
220\mbox{:~~~}  - {\displaystyle \frac {1}{2}} (2\,b - 1
)\, \left(  \! {\displaystyle \sum _{k=1}^{2\,n - 3}} \,
{\displaystyle \frac {\Gamma ({\displaystyle \frac {3}{2}}  - b
 - n + {\displaystyle \frac {a}{2}}  + k)\,\Gamma ( - 
{\displaystyle \frac {a}{2}}  - n + {\displaystyle \frac {3}{2}} 
 + k)}{\Gamma ({\displaystyle \frac {3}{2}}  - {\displaystyle 
\frac {a}{2}}  + b - n + k)\,\Gamma ({\displaystyle \frac {3}{2}
}  + {\displaystyle \frac {a}{2}}  - n + k)}}  \!  \right)  \\
\mathrm{sin}({\displaystyle \frac {\pi \,(1 + a - 2\,b)}{2}} )\,
\mathrm{sin}({\displaystyle \frac {\pi \,(a - 1)}{2}} )\,\Gamma (
n - {\displaystyle \frac {3}{2}}  + b - {\displaystyle \frac {a}{
2}} )\,\Gamma (b + {\displaystyle \frac {a}{2}}  - n + 
{\displaystyle \frac {3}{2}} ) \left/ {\vrule 
height0.80em width0em depth0.80em} \right. \!  \! ( \\
\mathrm{sin}({\displaystyle \frac {\pi \,(a - 3 + 2\,n)}{2}} )\,
\pi )\mbox{} + (\mathrm{sin}(2\,\pi \,a - 4\,\pi \,b) + \mathrm{
sin}(2\,\pi \,b) - \mathrm{sin}( - 2\,\pi \,b + 2\,\pi \,a))\,
\sqrt{\pi } \\
2^{(2\,b + 2)}\,\Gamma (b + {\displaystyle \frac {1}{2}} )\,
\Gamma (b + {\displaystyle \frac {a}{2}}  - n + {\displaystyle 
\frac {3}{2}} )\,\Gamma (1 + a - 2\,b) \left/ {\vrule 
height0.80em width0em depth0.80em} \right. \!  \! (( - 8\,
\mathrm{cos}(\pi \,a - 3\,\pi \,b) \\
\mbox{} + 8\,\mathrm{cos}( - \pi \,b + 2\,\pi \,a) - 8\,\mathrm{
cos}( - 3\,\pi \,b + 2\,\pi \,a) + 8\,\mathrm{cos}(\pi \,a + \pi 
\,b)) 
\Gamma ( - b + {\displaystyle \frac {5}{2}}  + {\displaystyle 
\frac {a}{2}}  - n)\,\Gamma (b)\,\Gamma (a)) }
}
\end{maplelatex}

\end{maplegroup}

%% file: AppendixB221to240.tex
\begin{maplegroup}
\mapleresult
\begin{maplelatex}
\mapleinline{inert}{2d}{221, ":   ",
1/2*Sum(GAMMA(a-1+k)*GAMMA(b-1+k)/GAMMA(4-a-2*n+k)/GAMMA(4-b-2*n+k),k
= 1 ..
2*n-3)/GAMMA(2*n-4+a)/GAMMA(b-4+2*n)*GAMMA(-a+1)*GAMMA(-b+1)+256*(sin(
2*Pi*b+2*Pi*a)-sin(4*Pi*a)+sin(-2*Pi*b+2*Pi*a))/(-cos(Pi*a)+2*cos(Pi*a
-2*Pi*b)+cos(3*Pi*a)-cos(3*Pi*a-2*Pi*b)-cos(Pi*a+2*Pi*b))*Pi^(1/2)*2^(
-2*b-2*a)*2^(-4*n)*GAMMA(b+a-4+2*n)/GAMMA(2*n-4+a)/GAMMA(b-4+2*n)*GAMM
A(2*a-4+2*n)*GAMMA(9/2-2*n-b-a)/GAMMA(-2*b+5-2*n)*GAMMA(-b+1)/GAMMA(a)
;}{%
\maplemultiline{
221\mbox{:~~~} {\displaystyle \frac {1}{2}} \,
{\displaystyle \frac { \left(  \! {\displaystyle \sum _{k=1}^{2\,
n - 3}} \,{\displaystyle \frac {\Gamma (a - 1 + k)\,\Gamma (b - 1
 + k)}{\Gamma (4 - a - 2\,n + k)\,\Gamma (4 - b - 2\,n + k)}} 
 \!  \right) \,\Gamma ( - a + 1)\,\Gamma ( - b + 1)}{\Gamma (2\,n
 - 4 + a)\,\Gamma (b - 4 + 2\,n)}}  + 256 \\
(\mathrm{sin}(2\,\pi \,b + 2\,\pi \,a) - \mathrm{sin}(4\,\pi \,a)
 + \mathrm{sin}( - 2\,\pi \,b + 2\,\pi \,a))\,\sqrt{\pi }\,2^{(
 - 2\,b - 2\,a)}\,2^{( - 4\,n)} \\
\Gamma (b + a - 4 + 2\,n)\,\Gamma (2\,a - 4 + 2\,n)\,\Gamma (
{\displaystyle \frac {9}{2}}  - 2\,n - b - a)\,\Gamma ( - b + 1)/
(( \\
 - \mathrm{cos}(\pi \,a) + 2\,\mathrm{cos}(\pi \,a - 2\,\pi \,b)
 + \mathrm{cos}(3\,\pi \,a) - \mathrm{cos}(3\,\pi \,a - 2\,\pi \,
b) - \mathrm{cos}(\pi \,a + 2\,\pi \,b) ) \\
\Gamma (2\,n - 4 + a)\,\Gamma (b - 4 + 2\,n)\,\Gamma ( - 2\,b + 
5 - 2\,n)\,\Gamma (a)) }
}
\end{maplelatex}

\begin{maplelatex}
\mapleinline{inert}{2d}{222, ":   ",
-2*4^(-a)*Sum(GAMMA(-a+b+k)*GAMMA(3-b+k-2*n)/GAMMA(3+a-2*n-b+k)/GAMMA(
b+k),k = 1 ..
2*n-3)*sin(Pi*(a-b))*sin(Pi*b)/Pi^(3/2)*GAMMA(b+1)*GAMMA(3+2*a-2*n-b)*
GAMMA(2*n-3+b)/GAMMA(a-1/2)*GAMMA(a-b)+1/2*1/GAMMA(-2*b+3+2*a-2*n)*GAM
MA(3+2*a-2*n-b)/GAMMA(2*b-3+2*n)*GAMMA(2*n-3+b)*GAMMA(b+1)*GAMMA(a-b)/
GAMMA(a);}{%
\maplemultiline{
222\mbox{:~~~}  - 2\,4^{( - a)}\, \left(  \! 
{\displaystyle \sum _{k=1}^{2\,n - 3}} \,{\displaystyle \frac {
\Gamma ( - a + b + k)\,\Gamma (3 - b + k - 2\,n)}{\Gamma (3 + a
 - 2\,n - b + k)\,\Gamma (b + k)}}  \!  \right) \,\mathrm{sin}(
\pi \,(a - b))\,\mathrm{sin}(\pi \,b)\,\Gamma (b + 1) \\
\Gamma (3 + 2\,a - 2\,n - b)\,\Gamma (2\,n - 3 + b)\,\Gamma (a - 
b) \left/ {\vrule height0.80em width0em depth0.80em} \right. \! 
 \! (\pi ^{(3/2)}\,\Gamma (a - {\displaystyle \frac {1}{2}} ))
 \\
\mbox{} + {\displaystyle \frac {1}{2}} \,{\displaystyle \frac {
\Gamma (3 + 2\,a - 2\,n - b)\,\Gamma (2\,n - 3 + b)\,\Gamma (b + 
1)\,\Gamma (a - b)}{\Gamma ( - 2\,b + 3 + 2\,a - 2\,n)\,\Gamma (2
\,b - 3 + 2\,n)\,\Gamma (a)}}  }
}
\end{maplelatex}

\begin{maplelatex}
\mapleinline{inert}{2d}{223, ":   ",
-1/2*Sum(GAMMA(-1/2*a-1/2+b+k)*GAMMA(3-b+k-2*n)/GAMMA(7/2+1/2*a-2*n-b+
k)/GAMMA(b+k),k = 1 ..
2*n-3)/Pi/GAMMA(2*n-3-a+b)*GAMMA(1/2*a+b+1/2)*GAMMA(1/2*a-b+1/2)*GAMMA
(2*n-3+b)/GAMMA(a)*sin(Pi*b)*cos(1/2*Pi*a-Pi*b)/sin(Pi*a-Pi*b)+(sin(Pi
*(1/2*a-b+1/2))/sin(1/2*Pi*(1+a))-1/2*(sin(Pi*a)+2*sin(2*Pi*b)+2*sin(P
i*a-2*Pi*b)-sin(Pi*a-4*Pi*b)+sin(2*Pi*a-4*Pi*b))*sin(Pi*a-2*Pi*b)/sin(
Pi*a-Pi*b)/(sin(Pi*a-4*Pi*b)+sin(Pi*a)+sin(2*Pi*a-4*Pi*b)))*GAMMA(2*b-
3-a+2*n)/GAMMA(1/2*a+1/2)^2*GAMMA(1/2*a+b+1/2)*GAMMA(1/2*a-b+1/2)/GAMM
A(2*b-3+2*n)*GAMMA(2*n-3+b)/GAMMA(2*n-3-a+b);}{%
\maplemultiline{
223\mbox{:~~~}  - {\displaystyle \frac {1}{2}}  \left( 
 \! {\displaystyle \sum _{k=1}^{2\,n - 3}} \,{\displaystyle 
\frac {\Gamma ( - {\displaystyle \frac {a}{2}}  - {\displaystyle 
\frac {1}{2}}  + b + k)\,\Gamma (3 - b + k - 2\,n)}{\Gamma (
{\displaystyle \frac {7}{2}}  + {\displaystyle \frac {a}{2}}  - 2
\,n - b + k)\,\Gamma (b + k)}}  \!  \right) \,\Gamma (
{\displaystyle \frac {a}{2}}  + b + {\displaystyle \frac {1}{2}} 
)\,\Gamma ({\displaystyle \frac {a}{2}}  - b + {\displaystyle 
\frac {1}{2}} ) \\
\Gamma (2\,n - 3 + b)\,\mathrm{sin}(\pi \,b)\,\mathrm{cos}(
{\displaystyle \frac {1}{2}} \,\pi \,a - \pi \,b)/(\pi \,\Gamma (
2\,n - 3 - a + b)\,\Gamma (a)\,\mathrm{sin}(\pi \,a - \pi \,b))
\mbox{} +  \left( {\vrule height1.67em width0em depth1.67em}
 \right. \!  \! {\displaystyle \frac {\mathrm{sin}(\pi \,(
{\displaystyle \frac {a}{2}}  - b + {\displaystyle \frac {1}{2}} 
))}{\mathrm{sin}({\displaystyle \frac {\pi \,(1 + a)}{2}} )}} \\  - 
{\displaystyle \frac {1}{2}}
(\mathrm{sin}(\pi \,a) + 2\,\mathrm{sin}(2\,\pi \,b) + 2\,
\mathrm{sin}(\pi \,a - 2\,\pi \,b) - \mathrm{sin}(\pi \,a - 4\,
\pi \,b) + \mathrm{sin}(2\,\pi \,a - 4\,\pi \,b)) \\
\mathrm{sin}(\pi \,a - 2\,\pi \,b)/(\mathrm{sin}(\pi \,a - \pi \,
b)\,(\mathrm{sin}(\pi \,a - 4\,\pi \,b) + \mathrm{sin}(\pi \,a)
 + \mathrm{sin}(2\,\pi \,a - 4\,\pi \,b))) \! \! \left. {\vrule 
height1.67em width0em depth1.67em} \right)  \\
\Gamma (2\,b - 3 - a + 2\,n)\,\Gamma ({\displaystyle \frac {a}{2}
}  + b + {\displaystyle \frac {1}{2}} )\,\Gamma ({\displaystyle 
\frac {a}{2}}  - b + {\displaystyle \frac {1}{2}} )\,\Gamma (2\,n
 - 3 + b) \left/ {\vrule height0.87em width0em depth0.87em}
 \right. \!  \! (\Gamma ({\displaystyle \frac {a}{2}}  + 
{\displaystyle \frac {1}{2}} )^{2} 
\Gamma (2\,b - 3 + 2\,n)\,\Gamma (2\,n - 3 - a + b)) }
}
\end{maplelatex}

\begin{maplelatex}
\mapleinline{inert}{2d}{224, ":   ",
1/2*sin(Pi*(a-b))*sin(Pi*a)*Sum(GAMMA(-a+2+k-n)*GAMMA(b-1+k)/GAMMA(a+k
)/GAMMA(3-b-n+k),k = 1 ..
n-2)/sin(Pi*b)/Pi*GAMMA(1+a)*GAMMA(1+a-b)/GAMMA(b)/GAMMA(n+2*b-3)*GAMM
A(n+a-2)*GAMMA(2*b-a)-1/Pi^(1/2)/GAMMA(n+2*a-2)/GAMMA(b)/sin(Pi*b)^2/s
in(Pi*(a+b))*2^(-1+2*a-2*b)*GAMMA(a-b+1/2)*GAMMA(n+a-2)*GAMMA(2*b-a)*G
AMMA(1+a)*(sin(2*Pi*b)*sin(Pi*(a-b))*sin(Pi*(a+b))-2*sin(Pi*a)*sin(Pi*
b)*sin(2*Pi*b)*cos(Pi*(a-b))-2*sin(Pi*a)*sin(Pi*b)*cos(2*Pi*b)*sin(Pi*
(a-b)));}{%
\maplemultiline{
224\mbox{:~~~} {\displaystyle \frac {1}{2}} \mathrm{sin}
(\pi \,(a - b))\,\mathrm{sin}(\pi \,a)\, \left(  \! 
{\displaystyle \sum _{k=1}^{n - 2}} \,{\displaystyle \frac {
\Gamma ( - a + 2 + k - n)\,\Gamma (b - 1 + k)}{\Gamma (a + k)\,
\Gamma (3 - b - n + k)}}  \!  \right) \,\Gamma (1 + a) \\
\Gamma (1 + a - b)\,\Gamma (n + a - 2)\,\Gamma (2\,b - a)/(
\mathrm{sin}(\pi \,b)\,\pi \,\Gamma (b)\,\Gamma (n + 2\,b - 3))
\mbox{} -  \\
2^{( - 1 + 2\,a - 2\,b)}\,\Gamma (a - b + {\displaystyle \frac {1
}{2}} )\,\Gamma (n + a - 2)\,\Gamma (2\,b - a)\,\Gamma (1 + a)(
\mathrm{sin}(2\,\pi \,b)\,\mathrm{sin}(\pi \,(a - b))\,\mathrm{
sin}(\pi \,(a + b)) \\
\mbox{} - 2\,\mathrm{sin}(\pi \,a)\,\mathrm{sin}(\pi \,b)\,
\mathrm{sin}(2\,\pi \,b)\,\mathrm{cos}(\pi \,(a - b))
\mbox{} - 2\,\mathrm{sin}(\pi \,a)\,\mathrm{sin}(\pi \,b)\,
\mathrm{cos}(2\,\pi \,b)\,\mathrm{sin}(\pi \,(a - b))) \left/ 
{\vrule height0.44em width0em depth0.44em} \right. \! \\ \! (\sqrt{
\pi }\,\Gamma (n + 2\,a - 2)\,\Gamma (b) 
\mathrm{sin}(\pi \,b)^{2}\,\mathrm{sin}(\pi \,(a + b))) }
}
\end{maplelatex}

\begin{maplelatex}
\mapleinline{inert}{2d}{225, ":   ",
1/2*Sum(GAMMA(-a+2+k-n)*GAMMA(-b+a+k)/GAMMA(a+k)/GAMMA(-a+b+2+k-n),k =
1 ..
n-2)/Pi*GAMMA(1+a)/GAMMA(1+a-b)*GAMMA(n+a-2)/GAMMA(n+a-b-1)/GAMMA(2*a-
2*b+n-1)*GAMMA(-1-b+2*a+n)*GAMMA(b)*GAMMA(2-2*b)*sin(Pi*b)*sin(Pi*a)/s
in(Pi*a-Pi*b)+(-Pi*cos(-2*Pi*b+2*Pi*a)/sin(2*Pi*b)/sin(-Pi*b+2*Pi*a)*s
in(Pi*b)-Pi*sin(-2*Pi*b+2*Pi*a)*cos(Pi*b)/sin(2*Pi*b)/sin(-Pi*b+2*Pi*a
))*GAMMA(-1-b+2*a+n)/GAMMA(n+2*a-2)/GAMMA(-n-a+3)/GAMMA(1+a-b)*GAMMA(2
+b-a-n)/GAMMA(b)*GAMMA(1+a)-1/2*GAMMA(2+b-a-n)*Pi^(1/2)*2^(2*b-2)*GAMM
A(b-a)*GAMMA(b-1/2)*GAMMA(n+a-2)/GAMMA(-2*a+2*b-n+1)/GAMMA(b-1)/GAMMA(
-a)/sin(Pi*(a-1))*GAMMA(2-2*b)*GAMMA(-1-b+2*a+n)/GAMMA(2*a-2*b+n)/GAMM
A(-b+2)/GAMMA(n+2*a-2);}{%
\maplemultiline{
225\mbox{:~~~} {\displaystyle \frac {1}{2}}  \left(  \! 
{\displaystyle \sum _{k=1}^{n - 2}} \,{\displaystyle \frac {
\Gamma ( - a + 2 + k - n)\,\Gamma ( - b + a + k)}{\Gamma (a + k)
\,\Gamma ( - a + b + 2 + k - n)}}  \!  \right) \,\Gamma (1 + a)\,
\Gamma (n + a - 2) \\
\Gamma ( - 1 - b + 2\,a + n)\,\Gamma (b)\,\Gamma (2 - 2\,b)\,
\mathrm{sin}(\pi \,b)\,\mathrm{sin}(\pi \,a)/(\pi \,\Gamma (1 + a
 - b) \\
\Gamma (n + a - b - 1)\,\Gamma (2\,a - 2\,b + n - 1)\,\mathrm{sin
}(\pi \,a - \pi \,b))\mbox{} +  \\
( - {\displaystyle \frac {\pi \,\mathrm{cos}( - 2\,\pi \,b + 2\,
\pi \,a)\,\mathrm{sin}(\pi \,b)}{\mathrm{sin}(2\,\pi \,b)\,
\mathrm{sin}( - \pi \,b + 2\,\pi \,a)}}  - {\displaystyle \frac {
\pi \,\mathrm{sin}( - 2\,\pi \,b + 2\,\pi \,a)\,\mathrm{cos}(\pi 
\,b)}{\mathrm{sin}(2\,\pi \,b)\,\mathrm{sin}( - \pi \,b + 2\,\pi 
\,a)}} ) \\
\Gamma ( - 1 - b + 2\,a + n)\,\Gamma (2 + b - a - n)\,\Gamma (1
 + a)/(\Gamma (n + 2\,a - 2)\,\Gamma ( - n - a + 3) \\
\Gamma (1 + a - b)\,\Gamma (b))\mbox{} - {\displaystyle \frac {1
}{2}} \Gamma (2 + b - a - n)\,\sqrt{\pi }\,2^{(2\,b - 2)}\,\Gamma
 (b - a)\,\Gamma (b - {\displaystyle \frac {1}{2}} )\,\Gamma (n
 + a - 2) \\
\Gamma (2 - 2\,b)\,\Gamma ( - 1 - b + 2\,a + n)/(\Gamma ( - 2\,a
 + 2\,b - n + 1)\,\Gamma (b - 1)\,\Gamma ( - a)\,\mathrm{sin}(\pi
 \,(a - 1)) \\
\Gamma (2\,a - 2\,b + n)\,\Gamma ( - b + 2)\,\Gamma (n + 2\,a - 2
)) }
}
\end{maplelatex}

\begin{maplelatex}
\mapleinline{inert}{2d}{226, ":   ",
-1/2*Sum(GAMMA(a-1+k)*GAMMA(-1/2*b-1/2*n+1+k)/GAMMA(3-n-a+k)/GAMMA(1+1
/2*b-1/2*n+k),k = 1 ..
n-2)*sin(1/2*Pi*(b+n))*(2*a+n-3)*Pi/sin(1/2*Pi*(n-b))/sin(Pi*(n+a-2))/
GAMMA(a)/GAMMA(-1/2*b+1/2*n-1)*GAMMA(1/2*n+1/2*b-1)/GAMMA(n+a-2)-1/8*2
^(b-n-2*a+2)*(sin(3*Pi*a)-3*sin(Pi*a)+4*sin(Pi*(a-b))*sin(Pi*a)^2*(-1)
^n)*GAMMA(1/2*n+a-1/2*b-1)*GAMMA(-1/2*n-a+3/2+1/2*b)*GAMMA(n+2*a-2)*GA
MMA(1/2*b-1/2*n+2)*GAMMA(1/2*n+1/2*b-1)/Pi^(1/2)/sin(Pi*a)^2/sin(1/2*P
i*(n-b))/GAMMA(n+a-2)/GAMMA(b)/GAMMA(a);}{%
\maplemultiline{
226\mbox{:~~~}  - {\displaystyle \frac {1}{2}} \,
{\displaystyle \frac { \left(  \! {\displaystyle \sum _{k=1}^{n
 - 2}} \,{\displaystyle \frac {\Gamma (a - 1 + k)\,\Gamma ( - 
{\displaystyle \frac {b}{2}}  - {\displaystyle \frac {n}{2}}  + 1
 + k)}{\Gamma (3 - n - a + k)\,\Gamma (1 + {\displaystyle \frac {
b}{2}}  - {\displaystyle \frac {n}{2}}  + k)}}  \!  \right) \,
\mathrm{sin}({\displaystyle \frac {\pi \,(b + n)}{2}} )\,(2\,a + 
n - 3)\,\pi \,\Gamma ({\displaystyle \frac {n}{2}}  + 
{\displaystyle \frac {b}{2}}  - 1)}{\mathrm{sin}({\displaystyle 
\frac {\pi \,(n - b)}{2}} )\,\mathrm{sin}(\pi \,(n + a - 2))\,
\Gamma (a)\,\Gamma ( - {\displaystyle \frac {b}{2}}  + 
{\displaystyle \frac {n}{2}}  - 1)\,\Gamma (n + a - 2)}}  \\
\mbox{} - {\displaystyle \frac {1}{8}} 2^{(b - n - 2\,a + 2)}\,(
\mathrm{sin}(3\,\pi \,a) - 3\,\mathrm{sin}(\pi \,a) + 4\,\mathrm{
sin}(\pi \,(a - b))\,\mathrm{sin}(\pi \,a)^{2}\,(-1)^{n}) \\
\Gamma ({\displaystyle \frac {n}{2}}  + a - {\displaystyle 
\frac {b}{2}}  - 1)\,\Gamma ( - {\displaystyle \frac {n}{2}}  - a
 + {\displaystyle \frac {3}{2}}  + {\displaystyle \frac {b}{2}} )
\,\Gamma (n + 2\,a - 2)\,\Gamma ({\displaystyle \frac {b}{2}}  - 
{\displaystyle \frac {n}{2}}  + 2)\,\Gamma ({\displaystyle 
\frac {n}{2}}  + {\displaystyle \frac {b}{2}}  - 1) \\
 \left/ {\vrule height0.80em width0em depth0.80em} \right. \! 
 \! (\sqrt{\pi }\,\mathrm{sin}(\pi \,a)^{2}\,\mathrm{sin}(
{\displaystyle \frac {\pi \,(n - b)}{2}} )\,\Gamma (n + a - 2)\,
\Gamma (b)\,\Gamma (a)) }
}
\end{maplelatex}

\begin{maplelatex}
\mapleinline{inert}{2d}{227, ":   ",
(1/2*b-1/2)*sin(Pi*(a+b))/sin(Pi*a)*GAMMA(b+n+a-3)/GAMMA(n+a-3)*GAMMA(
-1+a+b)/GAMMA(a)*Sum(GAMMA(a-1+k)*GAMMA(-a-b+3+k-n)/GAMMA(3-n-a+k)/GAM
MA(-1+b+a+k),k = 1 ..
n-2)+1/8*Pi^(1/2)*4^b*(2*sin(Pi*(a+b))*sin(Pi*a)+sin(Pi*b)*sin(2*Pi*(a
+b)))/sin(Pi*b)/sin(2*Pi*(a+b))/sin(Pi*a)^2*GAMMA(b+n+a-3)/GAMMA(n+2*a
+2*b-4)*sin(2*Pi*a)*GAMMA(2*a+n-3)/GAMMA(n+a-3)*GAMMA(-1+a+b)/GAMMA(b-
1)*GAMMA(b-1/2)/GAMMA(a);}{%
\maplemultiline{
227\mbox{:~~~} ({\displaystyle \frac {b}{2}}  - 
{\displaystyle \frac {1}{2}} )\,\mathrm{sin}(\pi \,(a + b))\,
\Gamma (b + n + a - 3)\,\Gamma ( - 1 + a + b) \\
 \left(  \! {\displaystyle \sum _{k=1}^{n - 2}} \,{\displaystyle 
\frac {\Gamma (a - 1 + k)\,\Gamma ( - a - b + 3 + k - n)}{\Gamma 
(3 - n - a + k)\,\Gamma ( - 1 + b + a + k)}}  \!  \right) /(
\mathrm{sin}(\pi \,a)\,\Gamma (n + a - 3)\,\Gamma (a))\mbox{} + 
{\displaystyle \frac {1}{8}} \sqrt{\pi }\,4^{b} \\
(2\,\mathrm{sin}(\pi \,(a + b))\,\mathrm{sin}(\pi \,a) + \mathrm{
sin}(\pi \,b)\,\mathrm{sin}(2\,\pi \,(a + b)))\,\Gamma (b + n + a
 - 3)\,\mathrm{sin}(2\,\pi \,a) \\
\Gamma (2\,a + n - 3)\,\Gamma ( - 1 + a + b)\,\Gamma (b - 
{\displaystyle \frac {1}{2}} ) \left/ {\vrule 
height0.44em width0em depth0.44em} \right. \!  \! (\mathrm{sin}(
\pi \,b)\,\mathrm{sin}(2\,\pi \,(a + b))\,\mathrm{sin}(\pi \,a)^{
2} \\
\Gamma (n + 2\,a + 2\,b - 4)\,\Gamma (n + a - 3)\,\Gamma (b - 1)
\,\Gamma (a)) }
}
\end{maplelatex}

\begin{maplelatex}
\mapleinline{inert}{2d}{228, ":   ",
1/2*sin(1/2*Pi*(2*a+n-b-3))*sin(Pi*(n+a-1))*sin(1/2*Pi*(n-b-1))*Sum(GA
MMA(-a+2+k-n)*GAMMA(1/2*b+1/2-1/2*n+k)/GAMMA(a+k)/GAMMA(3/2-1/2*b-1/2*
n+k),k = 1 ..
n-2)/sin(1/2*Pi*(n+b-3))/Pi^2/GAMMA(b)*GAMMA(n+a-2)/GAMMA(1/2*b-1/2+1/
2*n)*GAMMA(a-1/2*b+1/2*n-1/2)*GAMMA(a+1/2*b+1/2*n-1/2)*GAMMA(b+3-a-n)*
GAMMA(1/2*n-1/2*b-1/2)-1/32*4^a*2^(n-b)*(-sin(Pi*(n+2*a+b))-sin(Pi*(n+
2*a-b))-2*sin(2*Pi*(a+n))+2*sin(Pi*(b-n))+sin(2*Pi*b))/Pi^(3/2)/(cos(P
i*a)+cos(Pi*(b+a+n)))*GAMMA(1/2*n+a-1/2*b-1)/GAMMA(n+2*a-2)*GAMMA(1/2*
n-1/2*b-1/2)*GAMMA(a+1/2*b+1/2*n-1/2)*GAMMA(b+3-a-n)/GAMMA(1/2*b-1/2+1
/2*n)*GAMMA(n+a-2);}{%
\maplemultiline{
228\mbox{:~~~} {\displaystyle \frac {1}{2}} \mathrm{sin}
({\displaystyle \frac {\pi \,(2\,a + n - b - 3)}{2}} )\,\mathrm{
sin}(\pi \,(n + a - 1))\,\mathrm{sin}({\displaystyle \frac {\pi 
\,(n - b - 1)}{2}} ) \\
 \left(  \! {\displaystyle \sum _{k=1}^{n - 2}} \,{\displaystyle 
\frac {\Gamma ( - a + 2 + k - n)\,\Gamma ({\displaystyle \frac {b
}{2}}  + {\displaystyle \frac {1}{2}}  - {\displaystyle \frac {n
}{2}}  + k)}{\Gamma (a + k)\,\Gamma ({\displaystyle \frac {3}{2}
}  - {\displaystyle \frac {b}{2}}  - {\displaystyle \frac {n}{2}
}  + k)}}  \!  \right) \,\Gamma (n + a - 2)\,\Gamma (a - 
{\displaystyle \frac {b}{2}}  + {\displaystyle \frac {n}{2}}  - 
{\displaystyle \frac {1}{2}} ) \\
\Gamma (a + {\displaystyle \frac {b}{2}}  + {\displaystyle 
\frac {n}{2}}  - {\displaystyle \frac {1}{2}} )\,\Gamma (b + 3 - 
a - n)\,\Gamma ({\displaystyle \frac {n}{2}}  - {\displaystyle 
\frac {b}{2}}  - {\displaystyle \frac {1}{2}} ) \left/ {\vrule 
height0.80em width0em depth0.80em} \right. \!  \! (\mathrm{sin}(
{\displaystyle \frac {\pi \,(n + b - 3)}{2}} )\,\pi ^{2}\,\Gamma 
(b) \\
\Gamma ({\displaystyle \frac {b}{2}}  - {\displaystyle \frac {1}{
2}}  + {\displaystyle \frac {n}{2}} ))\mbox{} - {\displaystyle 
\frac {1}{32}} 4^{a}\,2^{(n - b)}( - \mathrm{sin}(\pi \,(n + 2\,a
 + b)) - \mathrm{sin}(\pi \,(n + 2\,a - b)) \\
\mbox{} - 2\,\mathrm{sin}(2\,\pi \,(a + n)) + 2\,\mathrm{sin}(\pi
 \,(b - n)) + \mathrm{sin}(2\,\pi \,b))\Gamma ({\displaystyle 
\frac {n}{2}}  + a - {\displaystyle \frac {b}{2}}  - 1) \\
\Gamma ({\displaystyle \frac {n}{2}}  - {\displaystyle \frac {b}{
2}}  - {\displaystyle \frac {1}{2}} )\,\Gamma (a + 
{\displaystyle \frac {b}{2}}  + {\displaystyle \frac {n}{2}}  - 
{\displaystyle \frac {1}{2}} )\,\Gamma (b + 3 - a - n)\,\Gamma (n
 + a - 2) \left/ {\vrule height0.80em width0em depth0.80em}
 \right. \!  \! (\pi ^{(3/2)} \\
(\mathrm{cos}(\pi \,a) + \mathrm{cos}(\pi \,(b + a + n)))\,\Gamma
 (n + 2\,a - 2)\,\Gamma ({\displaystyle \frac {b}{2}}  - 
{\displaystyle \frac {1}{2}}  + {\displaystyle \frac {n}{2}} ))
 }
}
\end{maplelatex}

\begin{maplelatex}
\mapleinline{inert}{2d}{229, ":   ",
1/2*Sum(GAMMA(a-1+k)*GAMMA(-a+b+2+k-n)/GAMMA(3-n-a+k)/GAMMA(-b+a+k),k
= 1 ..
n-2)*b*Pi/sin(Pi*(b-a))/GAMMA(b-a+1)/GAMMA(a)/GAMMA(b+3-a-n)/GAMMA(-n+
2*b-2*a+3)*GAMMA(-n-a+2*b+3)*GAMMA(-n-2*a+4)+(b*(-1/2*cos(Pi*a)+1/2*co
s(3*Pi*a-2*Pi*b))/sin(Pi*(-2*a+b))/(b-a)*GAMMA(-2*b+1)*GAMMA(b)^2/sin(
Pi*(b-a))^2/GAMMA(b+3-a-n)*GAMMA(-n-a+2*b+3)*sin(Pi*b)+b*GAMMA(b)*Pi^(
3/2)*4^(-b)*sin(Pi*(2*b-2*a))/(b-a)/sin(Pi*a)/(2*b-1)/cos(Pi*b)/sin(Pi
*(b-a))^2/GAMMA(b+3-a-n)*GAMMA(-n-a+2*b+3)*sin(Pi*b)/GAMMA(b-1/2))/GAM
MA(b-a)/GAMMA(a);}{%
\maplemultiline{
229\mbox{:~~~} {\displaystyle \frac {1}{2}} \,
{\displaystyle \frac { \left(  \! {\displaystyle \sum _{k=1}^{n
 - 2}} \,{\displaystyle \frac {\Gamma (a - 1 + k)\,\Gamma ( - a
 + b + 2 + k - n)}{\Gamma (3 - n - a + k)\,\Gamma ( - b + a + k)}
}  \!  \right) \,b\,\pi \,\Gamma ( - n - a + 2\,b + 3)\,\Gamma (
 - n - 2\,a + 4)}{\mathrm{sin}(\pi \,(b - a))\,\Gamma (b - a + 1)
\,\Gamma (a)\,\Gamma (b + 3 - a - n)\,\Gamma ( - n + 2\,b - 2\,a
 + 3)}}  \\
\mbox{} +  \left( {\vrule height1.67em width0em depth1.67em}
 \right. \!  \! {\displaystyle \frac {b\,( - {\displaystyle 
\frac {1}{2}} \,\mathrm{cos}(\pi \,a) + {\displaystyle \frac {1}{
2}} \,\mathrm{cos}(3\,\pi \,a - 2\,\pi \,b))\,\Gamma ( - 2\,b + 1
)\,\Gamma (b)^{2}\,\Gamma ( - n - a + 2\,b + 3)\,\mathrm{sin}(\pi
 \,b)}{\mathrm{sin}(\pi \,( - 2\,a + b))\,(b - a)\,\mathrm{sin}(
\pi \,(b - a))^{2}\,\Gamma (b + 3 - a - n)}}  \\
\mbox{} + {\displaystyle \frac {b\,\Gamma (b)\,\pi ^{(3/2)}\,4^{(
 - b)}\,\mathrm{sin}(\pi \,(2\,b - 2\,a))\,\Gamma ( - n - a + 2\,
b + 3)\,\mathrm{sin}(\pi \,b)}{(b - a)\,\mathrm{sin}(\pi \,a)\,(2
\,b - 1)\,\mathrm{cos}(\pi \,b)\,\mathrm{sin}(\pi \,(b - a))^{2}
\,\Gamma (b + 3 - a - n)\,\Gamma (b - {\displaystyle \frac {1}{2}
} )}}  \! \! \left. {\vrule height1.67em width0em depth1.67em}
 \right)  \left/ {\vrule height0.37em width0em depth0.37em}
 \right. \!  \! (
\Gamma (b - a)\,\Gamma (a)) }
}
\end{maplelatex}

\begin{maplelatex}
\mapleinline{inert}{2d}{230, ":   ",
1/2*b*sin(1/2*Pi*(1+a-b))*sin(1/2*Pi*(n-b-1))*Sum(GAMMA(-1/2*a-1/2*n+1
+k)*GAMMA(1/2*b+1/2-1/2*n+k)/GAMMA(1+1/2*a-1/2*n+k)/GAMMA(3/2-1/2*b-1/
2*n+k),k = 1 ..
n-2)/sin(1/2*Pi*(n+b-3))/Pi/GAMMA(1/2*b-1/2+1/2*n)*GAMMA(1/2*a+1/2-1/2
*b)*GAMMA(1/2*n-1/2*b-1/2)*GAMMA(1/2*a+1/2+1/2*b)+1/16*1/Pi^(1/2)*GAMM
A(b+1)/GAMMA(1/2*b-1/2+1/2*n)*GAMMA(1/2*a+1/2+1/2*b)*GAMMA(1/2*n-1/2*b
-1/2)/GAMMA(a)*GAMMA(1/2*a-1/2*b)/cos(1/2*Pi*(b+n))/cos(1/2*Pi*(a+b))*
2^(a-b)*(sin(Pi*(a+b))+sin(Pi*(a-b))+2*sin(Pi*(n-b))+2*sin(Pi*(a+n))-s
in(2*Pi*b))/sin(1/2*Pi*(a+n));}{%
\maplemultiline{
230\mbox{:~~~} {\displaystyle \frac {1}{2}} b\,\mathrm{
sin}({\displaystyle \frac {\pi \,(1 + a - b)}{2}} )\,\mathrm{sin}
({\displaystyle \frac {\pi \,(n - b - 1)}{2}} ) \\
 \left(  \! {\displaystyle \sum _{k=1}^{n - 2}} \,{\displaystyle 
\frac {\Gamma ( - {\displaystyle \frac {a}{2}}  - {\displaystyle 
\frac {n}{2}}  + 1 + k)\,\Gamma ({\displaystyle \frac {b}{2}}  + 
{\displaystyle \frac {1}{2}}  - {\displaystyle \frac {n}{2}}  + k
)}{\Gamma (1 + {\displaystyle \frac {a}{2}}  - {\displaystyle 
\frac {n}{2}}  + k)\,\Gamma ({\displaystyle \frac {3}{2}}  - 
{\displaystyle \frac {b}{2}}  - {\displaystyle \frac {n}{2}}  + k
)}}  \!  \right) \,\Gamma ({\displaystyle \frac {a}{2}}  + 
{\displaystyle \frac {1}{2}}  - {\displaystyle \frac {b}{2}} )\,
\Gamma ({\displaystyle \frac {n}{2}}  - {\displaystyle \frac {b}{
2}}  - {\displaystyle \frac {1}{2}} ) \\
\Gamma ({\displaystyle \frac {a}{2}}  + {\displaystyle \frac {1}{
2}}  + {\displaystyle \frac {b}{2}} ) \left/ {\vrule 
height0.80em width0em depth0.80em} \right. \!  \! (\mathrm{sin}(
{\displaystyle \frac {\pi \,(n + b - 3)}{2}} )\,\pi \,\Gamma (
{\displaystyle \frac {b}{2}}  - {\displaystyle \frac {1}{2}}  + 
{\displaystyle \frac {n}{2}} ))\mbox{} + {\displaystyle \frac {1
}{16}} \Gamma (b + 1) \\
\Gamma ({\displaystyle \frac {a}{2}}  + {\displaystyle \frac {1}{
2}}  + {\displaystyle \frac {b}{2}} )\,\Gamma ({\displaystyle 
\frac {n}{2}}  - {\displaystyle \frac {b}{2}}  - {\displaystyle 
\frac {1}{2}} )\,\Gamma ({\displaystyle \frac {a}{2}}  - 
{\displaystyle \frac {b}{2}} )\,2^{(a - b)} \\
(\mathrm{sin}(\pi \,(a + b)) + \mathrm{sin}(\pi \,(a - b)) + 2\,
\mathrm{sin}(\pi \,(n - b)) + 2\,\mathrm{sin}(\pi \,(a + n)) - 
\mathrm{sin}(2\,\pi \,b)) \\
 \left/ {\vrule height0.80em width0em depth0.80em} \right. \! 
 \! (\sqrt{\pi }\,\Gamma ({\displaystyle \frac {b}{2}}  - 
{\displaystyle \frac {1}{2}}  + {\displaystyle \frac {n}{2}} )\,
\Gamma (a)\,\mathrm{cos}({\displaystyle \frac {\pi \,(b + n)}{2}
} )\,\mathrm{cos}({\displaystyle \frac {\pi \,(a + b)}{2}} )\,
\mathrm{sin}({\displaystyle \frac {\pi \,(a + n)}{2}} )) }
}
\end{maplelatex}

\begin{maplelatex}
\mapleinline{inert}{2d}{231, ":   ",
-1/2*Sum(GAMMA(-a-b+1+k)*GAMMA(b+1+k-n)/GAMMA(1-b+k)/GAMMA(1+b+a-n+k),
k = 1 ..
n-2)*(a-1)/Pi*GAMMA(-1+a+b)*GAMMA(a+2*b-n+1)*GAMMA(n-2*b)/GAMMA(b)*sin
(Pi*a+Pi*b)*sin(Pi*n-2*Pi*b)/sin(Pi*b)-1/8*Pi^(3/2)*GAMMA(a-1/2)*GAMMA
(-1+a+b)*4^a*(1/4*sin(2*Pi*a)-1/2*sin(2*Pi*b)+1/4*sin(2*Pi*(2*b+a))-1/
2*sin(2*Pi*(a+b))+1/4*sin(4*Pi*b))/GAMMA(a-1)/GAMMA(-a-2*b+n)/GAMMA(-n
+2*b+2*a)/GAMMA(b)/sin(Pi*(-a-2*b))^2/sin(Pi*a)/sin(Pi*(n-b))/sin(Pi*b
);}{%
\maplemultiline{
231\mbox{:~~~}  - {\displaystyle \frac {1}{2}}  \left( 
 \! {\displaystyle \sum _{k=1}^{n - 2}} \,{\displaystyle \frac {
\Gamma ( - a - b + 1 + k)\,\Gamma (b + 1 + k - n)}{\Gamma (1 - b
 + k)\,\Gamma (1 + b + a - n + k)}}  \!  \right) \,(a - 1)\,
\Gamma ( - 1 + a + b) \\
\Gamma (a + 2\,b - n + 1)\,\Gamma (n - 2\,b)\,\mathrm{sin}(\pi \,
a + \pi \,b)\,\mathrm{sin}(\pi \,n - 2\,\pi \,b)/(\pi \,\Gamma (b
)\,\mathrm{sin}(\pi \,b))\mbox{} - {\displaystyle \frac {1}{8}} 
\pi ^{(3/2)}\,\Gamma (a - {\displaystyle \frac {1}{2}} )\,\Gamma 
( - 1 + a + b)\,4^{a} \\ (
{\displaystyle \frac {1}{4}} \,\mathrm{sin}(2\,\pi \,a) - 
{\displaystyle \frac {1}{2}} \,\mathrm{sin}(2\,\pi \,b) + 
{\displaystyle \frac {1}{4}} \,\mathrm{sin}(2\,\pi \,(2\,b + a))
 - {\displaystyle \frac {1}{2}} \,\mathrm{sin}(2\,\pi \,(a + b))
 + {\displaystyle \frac {1}{4}} \,\mathrm{sin}(4\,\pi \,b)) \\
 \left/ {\vrule height0.44em width0em depth0.44em} \right. \! 
 \! (\Gamma (a - 1)\,\Gamma ( - a - 2\,b + n)\,\Gamma ( - n + 2\,
b + 2\,a)\,\Gamma (b)\,\mathrm{sin}(\pi \,( - a - 2\,b))^{2}\,
\mathrm{sin}(\pi \,a)
\mathrm{sin}(\pi \,(n - b))\,\mathrm{sin}(\pi \,b)) }
}
\end{maplelatex}

\begin{maplelatex}
\mapleinline{inert}{2d}{232, ":   ",
-1/2*Sum(GAMMA(-1/2*b-n+3/2+k)*GAMMA(1+1/2*a-n+k)/GAMMA(3/2+1/2*b-n+k)
/GAMMA(2-1/2*a-n+k),k = 1 ..
2*n-3)*a/Pi/GAMMA(1/2*a+1/2-1/2*b)*GAMMA(n-1-1/2*a)*GAMMA(2-1/2*a-n)*G
AMMA(1/2*a+1/2+1/2*b)*sin(1/2*Pi*a)*cos(1/2*Pi*b)/cos(Pi*n+1/2*Pi*b)+1
/16*2^(b-a)*Pi^(3/2)/GAMMA(1/2*a+1-1/2*b)/GAMMA(-1+1/2*a+n)/GAMMA(2-n+
1/2*a)/GAMMA(b)*GAMMA(1+a)*(sin(2*Pi*a)-sin(Pi*(a+b))-2*sin(Pi*a)+2*si
n(Pi*b)+sin(Pi*(a-b))+16*cos(1/2*Pi*(a+b))*sin(1/2*Pi*(a-b))*sin(1/2*P
i*(a-2*n))*sin(1/2*Pi*(2*n+a)))/sin(1/2*Pi*(2*n+a))/sin(1/2*Pi*(a-2*n)
)/sin(1/2*Pi*(a-b))/cos(1/2*Pi*(b+2*n))/cos(1/2*Pi*(a+b))*GAMMA(1/2*a+
1/2+1/2*b);}{%
\maplemultiline{
232\mbox{:~~~}  - {\displaystyle \frac {1}{2}}  \left( 
 \! {\displaystyle \sum _{k=1}^{2\,n - 3}} \,{\displaystyle 
\frac {\Gamma ( - {\displaystyle \frac {b}{2}}  - n + 
{\displaystyle \frac {3}{2}}  + k)\,\Gamma (1 + {\displaystyle 
\frac {a}{2}}  - n + k)}{\Gamma ({\displaystyle \frac {3}{2}}  + 
{\displaystyle \frac {b}{2}}  - n + k)\,\Gamma (2 - 
{\displaystyle \frac {a}{2}}  - n + k)}}  \!  \right) \,a\,\Gamma
 (n - 1 - {\displaystyle \frac {a}{2}} )\,\Gamma (2 - 
{\displaystyle \frac {a}{2}}  - n) \\
\Gamma ({\displaystyle \frac {a}{2}}  + {\displaystyle \frac {1}{
2}}  + {\displaystyle \frac {b}{2}} )\,\mathrm{sin}(
{\displaystyle \frac {\pi \,a}{2}} )\,\mathrm{cos}(
{\displaystyle \frac {\pi \,b}{2}} ) \left/ {\vrule 
height0.80em width0em depth0.80em} \right. \!  \! (\pi \,\Gamma (
{\displaystyle \frac {a}{2}}  + {\displaystyle \frac {1}{2}}  - 
{\displaystyle \frac {b}{2}} )\,\mathrm{cos}(\pi \,n + 
{\displaystyle \frac {1}{2}} \,\pi \,b))\mbox{} \\ + {\displaystyle 
\frac {1}{16}} \,
2^{(b - a)}\,\pi ^{(3/2)}\,\Gamma (1 + a)(\mathrm{sin}(2\,\pi \,a
) - \mathrm{sin}(\pi \,(a + b)) - 2\,\mathrm{sin}(\pi \,a) + 2\,
\mathrm{sin}(\pi \,b) 
\mbox{} + \mathrm{sin}(\pi \,(a - b)) \\
\mbox{} + 16\,\mathrm{cos}({\displaystyle \frac {\pi \,(a + b)}{2
}} )\,\mathrm{sin}({\displaystyle \frac {\pi \,(a - b)}{2}} )\,
\mathrm{sin}({\displaystyle \frac {\pi \,(a - 2\,n)}{2}} )\,
\mathrm{sin}({\displaystyle \frac {\pi \,(2\,n + a)}{2}} )) \\
\Gamma ({\displaystyle \frac {a}{2}}  + {\displaystyle \frac {1}{
2}}  + {\displaystyle \frac {b}{2}} ) \left/ {\vrule 
height0.80em width0em depth0.80em} \right. \!  \! (\Gamma (
{\displaystyle \frac {a}{2}}  + 1 - {\displaystyle \frac {b}{2}} 
)\,\Gamma ( - 1 + {\displaystyle \frac {a}{2}}  + n)\,\Gamma (2
 - n + {\displaystyle \frac {a}{2}} )\,\Gamma (b) \\
\mathrm{sin}({\displaystyle \frac {\pi \,(2\,n + a)}{2}} )\,
\mathrm{sin}({\displaystyle \frac {\pi \,(a - 2\,n)}{2}} )\,
\mathrm{sin}({\displaystyle \frac {\pi \,(a - b)}{2}} )\,\mathrm{
cos}({\displaystyle \frac {\pi \,(b + 2\,n)}{2}} )\,\mathrm{cos}(
{\displaystyle \frac {\pi \,(a + b)}{2}} )) }
}
\end{maplelatex}

\begin{maplelatex}
\mapleinline{inert}{2d}{233, ":   ",
1/2*Sum(GAMMA(1-1/2*a-n+1/2*b+k)*GAMMA(1+1/2*a-n+k)/GAMMA(2+1/2*a-n-1/
2*b+k)/GAMMA(2-1/2*a-n+k),k = 1 ..
2*n-3)*a/Pi*GAMMA(n-1+1/2*a-1/2*b)/GAMMA(1/2*b)^2*GAMMA(1/2*a+1/2*b+2-
n)*GAMMA(n-1-1/2*a)*GAMMA(2-1/2*a-n)*sin(1/2*Pi*a)*sin(-1/2*Pi*b+1/2*P
i*a)/sin(1/2*Pi*b)-1/16*Pi/GAMMA(b)/GAMMA(-1+1/2*a+n)/GAMMA(2-n+1/2*a)
*GAMMA(1+a)/sin(1/2*Pi*b)/sin(1/2*Pi*(2*a-b))/sin(1/2*Pi*a)^2/cos(1/2*
Pi*b)*(-sin(Pi*(2*a-b))+sin(Pi*b)+2*sin(Pi*(a-b))-sin(2*Pi*a)+2*sin(Pi
*a)-16*sin(1/2*Pi*(2*a-b))*sin(1/2*Pi*a)^2*cos(1/2*Pi*b))*GAMMA(n-1+1/
2*a-1/2*b)/GAMMA(1+a-b)*GAMMA(1/2*a+1/2*b+2-n);}{%
\maplemultiline{
233\mbox{:~~~} {\displaystyle \frac {1}{2}}  \left(  \! 
{\displaystyle \sum _{k=1}^{2\,n - 3}} \,{\displaystyle \frac {
\Gamma (1 - {\displaystyle \frac {a}{2}}  - n + {\displaystyle 
\frac {b}{2}}  + k)\,\Gamma (1 + {\displaystyle \frac {a}{2}}  - 
n + k)}{\Gamma (2 + {\displaystyle \frac {a}{2}}  - n - 
{\displaystyle \frac {b}{2}}  + k)\,\Gamma (2 - {\displaystyle 
\frac {a}{2}}  - n + k)}}  \!  \right) \,a\,\Gamma (n - 1 + 
{\displaystyle \frac {a}{2}}  - {\displaystyle \frac {b}{2}} )
 \\
\Gamma ({\displaystyle \frac {a}{2}}  + {\displaystyle \frac {b}{
2}}  + 2 - n)\,\Gamma (n - 1 - {\displaystyle \frac {a}{2}} )\,
\Gamma (2 - {\displaystyle \frac {a}{2}}  - n)\,\mathrm{sin}(
{\displaystyle \frac {\pi \,a}{2}} )\,\mathrm{sin}( - 
{\displaystyle \frac {1}{2}} \,\pi \,b + {\displaystyle \frac {1
}{2}} \,\pi \,a) \left/ {\vrule height0.87em width0em depth0.87em
} \right. \!  \! \\ (
\pi \,\Gamma ({\displaystyle \frac {b}{2}} )^{2}\,\mathrm{sin}(
{\displaystyle \frac {\pi \,b}{2}} ))\mbox{} - {\displaystyle 
\frac {1}{16}} \pi \,\Gamma (1 + a)( - \mathrm{sin}(\pi \,(2\,a
 - b)) + \mathrm{sin}(\pi \,b) + 2\,\mathrm{sin}(\pi \,(a - b))
 \\
\mbox{} - \mathrm{sin}(2\,\pi \,a) + 2\,\mathrm{sin}(\pi \,a) - 
16\,\mathrm{sin}({\displaystyle \frac {\pi \,(2\,a - b)}{2}} )\,
\mathrm{sin}({\displaystyle \frac {\pi \,a}{2}} )^{2}\,\mathrm{
cos}({\displaystyle \frac {\pi \,b}{2}} )) \\
\Gamma (n - 1 + {\displaystyle \frac {a}{2}}  - {\displaystyle 
\frac {b}{2}} )\,\Gamma ({\displaystyle \frac {a}{2}}  + 
{\displaystyle \frac {b}{2}}  + 2 - n) \left/ {\vrule 
height0.87em width0em depth0.87em} \right. \!  \! (\Gamma (b)\,
\Gamma ( - 1 + {\displaystyle \frac {a}{2}}  + n)\,\Gamma (2 - n
 + {\displaystyle \frac {a}{2}} ) \\
\mathrm{sin}({\displaystyle \frac {\pi \,b}{2}} )\,\mathrm{sin}(
{\displaystyle \frac {\pi \,(2\,a - b)}{2}} )\,\mathrm{sin}(
{\displaystyle \frac {\pi \,a}{2}} )^{2}\,\mathrm{cos}(
{\displaystyle \frac {\pi \,b}{2}} )\,\Gamma (1 + a - b)) }
}
\end{maplelatex}

\begin{maplelatex}
\mapleinline{inert}{2d}{234, ":   ",
1/2*b*(n-a-1)*sin(Pi*a)*sin(Pi*(-1+1/2*b-a))*Sum(GAMMA(1/2*b-a+k)*GAMM
A(a+2-2*n+k)/GAMMA(3+a-2*n-1/2*b+k)/GAMMA(1-a+k),k = 1 ..
2*n-3)/Pi/sin(1/2*Pi*(4*n-6-2*a+b))/GAMMA(-a+2*n-2+1/2*b)*GAMMA(-1/2*b
+a)*GAMMA(2*n-2-a)*GAMMA(-a+1)+1/8*Pi^(3/2)*GAMMA(1/2*b+1)*GAMMA(2*a-2
*n+3)*GAMMA(-1/2*b+a)*(1/2*cos(3*Pi*a)-1/2*cos(Pi*a)-cos(Pi*a+Pi*b)-1/
2*cos(Pi*a-2*Pi*b)+cos(3*Pi*a-Pi*b)+1/2*cos(3*Pi*a+Pi*b)-1/2*cos(5*Pi*
a-2*Pi*b)+cos(3*Pi*a-2*Pi*b)-1/2*cos(5*Pi*a-Pi*b))*2^(-b)/GAMMA(3+2*a-
2*n-b)/GAMMA(-a+2*n-2+1/2*b)/GAMMA(3-2*n+a)/GAMMA(1/2*b+1/2)/GAMMA(a)/
sin(Pi*(a-b))/sin(Pi*(-1/2*b+a))/sin(Pi*(-1/2*b+2*a))/sin(Pi*a)^2/cos(
1/2*Pi*b);}{%
\maplemultiline{
234\mbox{:~~~} {\displaystyle \frac {1}{2}} b\,(n - a - 
1)\,\mathrm{sin}(\pi \,a)\,\mathrm{sin}(\pi \,( - 1 + 
{\displaystyle \frac {b}{2}}  - a)) \\
 \left(  \! {\displaystyle \sum _{k=1}^{2\,n - 3}} \,
{\displaystyle \frac {\Gamma ({\displaystyle \frac {b}{2}}  - a
 + k)\,\Gamma (a + 2 - 2\,n + k)}{\Gamma (3 + a - 2\,n - 
{\displaystyle \frac {b}{2}}  + k)\,\Gamma (1 - a + k)}}  \! 
 \right) \,\Gamma ( - {\displaystyle \frac {b}{2}}  + a)\,\Gamma 
(2\,n - 2 - a)\,\Gamma ( - a + 1) \\
 \left/ {\vrule height0.80em width0em depth0.80em} \right. \! 
 \! (\pi \,\mathrm{sin}({\displaystyle \frac {\pi \,(4\,n - 6 - 2
\,a + b)}{2}} )\,\Gamma ( - a + 2\,n - 2 + {\displaystyle \frac {
b}{2}} ))\mbox{} + {\displaystyle \frac {1}{8}} \pi ^{(3/2)}\,
\Gamma ({\displaystyle \frac {b}{2}}  + 1) \\
\Gamma (2\,a - 2\,n + 3)\,\Gamma ( - {\displaystyle \frac {b}{2}
}  + a)({\displaystyle \frac {1}{2}} \,\mathrm{cos}(3\,\pi \,a)
 - {\displaystyle \frac {1}{2}} \,\mathrm{cos}(\pi \,a) - 
\mathrm{cos}(\pi \,a + \pi \,b) \\
\mbox{} - {\displaystyle \frac {1}{2}} \,\mathrm{cos}(\pi \,a - 2
\,\pi \,b) + \mathrm{cos}(3\,\pi \,a - \pi \,b) + {\displaystyle 
\frac {1}{2}} \,\mathrm{cos}(3\,\pi \,a + \pi \,b) - 
{\displaystyle \frac {1}{2}} \,\mathrm{cos}(5\,\pi \,a - 2\,\pi 
\,b) \\
\mbox{} + \mathrm{cos}(3\,\pi \,a - 2\,\pi \,b) - {\displaystyle 
\frac {1}{2}} \,\mathrm{cos}(5\,\pi \,a - \pi \,b))2^{( - b)}
 \left/ {\vrule height0.80em width0em depth0.80em} \right. \! 
 \! (\Gamma (3 + 2\,a - 2\,n - b) \\
\Gamma ( - a + 2\,n - 2 + {\displaystyle \frac {b}{2}} )\,\Gamma 
(3 - 2\,n + a)\,\Gamma ({\displaystyle \frac {b}{2}}  + 
{\displaystyle \frac {1}{2}} )\,\Gamma (a)\,\mathrm{sin}(\pi \,(a
 - b))\,\mathrm{sin}(\pi \,( - {\displaystyle \frac {b}{2}}  + a)
) \\
\mathrm{sin}(\pi \,( - {\displaystyle \frac {b}{2}}  + 2\,a))\,
\mathrm{sin}(\pi \,a)^{2}\,\mathrm{cos}({\displaystyle \frac {\pi
 \,b}{2}} )) }
}
\end{maplelatex}

\begin{maplelatex}
\mapleinline{inert}{2d}{235, ":   ",
1/2*sin(1/2*Pi*a)*sin(Pi*(n+b-2))*Sum(GAMMA(3-b+k-2*n)*GAMMA(1+1/2*a-n
+k)/GAMMA(b+k)/GAMMA(2-1/2*a-n+k),k = 1 ..
2*n-3)/Pi^2/GAMMA(a)/GAMMA(2-n+1/2*a-b)*GAMMA(2*n-3+b)*GAMMA(1/2*a+n+b
-1)*GAMMA(-b+4+a-2*n)*GAMMA(2-1/2*a-n)*GAMMA(n-1-1/2*a)+1/256*Pi^(1/2)
*GAMMA(2*n-3+b)*GAMMA(1/2*a+n+b-1)*GAMMA(-b+4+a-2*n)*(cos(1/2*Pi*(-4*b
+5*a))-cos(1/2*Pi*(3*a+4*b))-cos(1/2*Pi*(a+8*b))-2*cos(1/2*Pi*(-4*b+3*
a))-cos(1/2*Pi*(-8*b+3*a))+2*cos(1/2*Pi*(a-8*b))+2*cos(1/2*Pi*(a+4*b))
+cos(3/2*Pi*a)-cos(1/2*Pi*a))*2^(-a+2*n+2*b)/GAMMA(-b+5/2+1/2*a-n)/GAM
MA(2*b-3+2*n)/GAMMA(2-n+1/2*a)/GAMMA(-1+1/2*a+n)/cos(1/2*Pi*(a-2*b))/s
in(1/2*Pi*a)^2/sin(1/2*Pi*(2*b+a))/sin(1/2*Pi*(a-4*b));}{%
\maplemultiline{
235\mbox{:~~~} {\displaystyle \frac {1}{2}} \mathrm{sin}
({\displaystyle \frac {\pi \,a}{2}} )\,\mathrm{sin}(\pi \,(n + b
 - 2))\, \left(  \! {\displaystyle \sum _{k=1}^{2\,n - 3}} \,
{\displaystyle \frac {\Gamma (3 - b + k - 2\,n)\,\Gamma (1 + 
{\displaystyle \frac {a}{2}}  - n + k)}{\Gamma (b + k)\,\Gamma (2
 - {\displaystyle \frac {a}{2}}  - n + k)}}  \!  \right)  \\
\Gamma (2\,n - 3 + b)\,\Gamma ({\displaystyle \frac {a}{2}}  + n
 + b - 1)\,\Gamma ( - b + 4 + a - 2\,n)\,\Gamma (2 - 
{\displaystyle \frac {a}{2}}  - n)\,\Gamma (n - 1 - 
{\displaystyle \frac {a}{2}} ) \\
 \left/ {\vrule height0.80em width0em depth0.80em} \right. \! 
 \! (\pi ^{2}\,\Gamma (a)\,\Gamma (2 - n + {\displaystyle \frac {
a}{2}}  - b))\mbox{} + {\displaystyle \frac {1}{256}} \sqrt{\pi }
\,\Gamma (2\,n - 3 + b)\,\Gamma ({\displaystyle \frac {a}{2}}  + 
n + b - 1) \\
\Gamma ( - b + 4 + a - 2\,n)(\mathrm{cos}({\displaystyle \frac {
\pi \,( - 4\,b + 5\,a)}{2}} ) - \mathrm{cos}({\displaystyle 
\frac {\pi \,(3\,a + 4\,b)}{2}} ) - \mathrm{cos}({\displaystyle 
\frac {\pi \,(a + 8\,b)}{2}} ) \\
\mbox{} - 2\,\mathrm{cos}({\displaystyle \frac {\pi \,( - 4\,b + 
3\,a)}{2}} ) - \mathrm{cos}({\displaystyle \frac {\pi \,( - 8\,b
 + 3\,a)}{2}} ) + 2\,\mathrm{cos}({\displaystyle \frac {\pi \,(a
 - 8\,b)}{2}} ) \\
\mbox{} + 2\,\mathrm{cos}({\displaystyle \frac {\pi \,(a + 4\,b)
}{2}} ) + \mathrm{cos}({\displaystyle \frac {3\,\pi \,a}{2}} ) - 
\mathrm{cos}({\displaystyle \frac {\pi \,a}{2}} ))2^{( - a + 2\,n
 + 2\,b)} \left/ {\vrule height0.80em width0em depth0.80em}
 \right. \!  \! \\ (
\Gamma ( - b + {\displaystyle \frac {5}{2}}  + {\displaystyle 
\frac {a}{2}}  - n)\,\Gamma (2\,b - 3 + 2\,n)\,\Gamma (2 - n + 
{\displaystyle \frac {a}{2}} )\,\Gamma ( - 1 + {\displaystyle 
\frac {a}{2}}  + n)\,\mathrm{cos}({\displaystyle \frac {\pi \,(a
 - 2\,b)}{2}} ) \\
\mathrm{sin}({\displaystyle \frac {\pi \,a}{2}} )^{2}\,\mathrm{
sin}({\displaystyle \frac {\pi \,(2\,b + a)}{2}} )\,\mathrm{sin}(
{\displaystyle \frac {\pi \,(a - 4\,b)}{2}} )) }
}
\end{maplelatex}

\begin{maplelatex}
\mapleinline{inert}{2d}{236, ":   ",
-1/2*Sum(GAMMA(3/2-b-n+1/2*a+k)*GAMMA(-1/2*a-n+3/2+k)/GAMMA(3/2-1/2*a+
b-n+k)/GAMMA(3/2+1/2*a-n+k),k = 1 ..
2*n-3)/Pi^2/GAMMA(1+a-2*b)/GAMMA(n-3/2-1/2*a)*GAMMA(n-3/2+1/2*a)*GAMMA
(n-3/2+b-1/2*a)*GAMMA(-1/2*a+b+3/2-n)*GAMMA(1+a-b)*GAMMA(b)*GAMMA(2-2*
b)*sin(Pi*b)*cos(1/2*Pi*a-Pi*b)/((-1)^n)+(-1/32*4^b/Pi^(3/2)/GAMMA(n-3
/2-1/2*a)*GAMMA(n-3/2+1/2*a)*GAMMA(n-3/2+b-1/2*a)*GAMMA(-1/2*a+b+3/2-n
)*GAMMA(1+a-b)*GAMMA(b-1/2)*GAMMA(2-2*b)*(2*sin(Pi*a)+sin(2*Pi*b)+2*si
n(Pi*a-2*Pi*b)+sin(-2*Pi*b+2*Pi*a)+sin(2*Pi*a-4*Pi*b))/cos(1/2*Pi*(-a+
2*n))/sin(Pi*(a-b))-Pi^2*(sin(1/2*Pi*a)*sin(Pi*b)-cos(1/2*Pi*a)*cos(Pi
*b))/cos(1/2*Pi*a)/sin(Pi*b)/sin(1/2*Pi*(-a+2*n+1-2*b))/GAMMA(n-3/2-1/
2*a)/GAMMA(-b+5/2+1/2*a-n)/GAMMA(n+1/2*a-b-1/2)*GAMMA(n-3/2+1/2*a)*GAM
MA(1+a-b)/GAMMA(b))/GAMMA(a);}{%
\maplemultiline{
236\mbox{:~~~}  - {\displaystyle \frac {1}{2}}  \left( 
 \! {\displaystyle \sum _{k=1}^{2\,n - 3}} \,{\displaystyle 
\frac {\Gamma ({\displaystyle \frac {3}{2}}  - b - n + 
{\displaystyle \frac {a}{2}}  + k)\,\Gamma ( - {\displaystyle 
\frac {a}{2}}  - n + {\displaystyle \frac {3}{2}}  + k)}{\Gamma (
{\displaystyle \frac {3}{2}}  - {\displaystyle \frac {a}{2}}  + b
 - n + k)\,\Gamma ({\displaystyle \frac {3}{2}}  + 
{\displaystyle \frac {a}{2}}  - n + k)}}  \!  \right) \,\Gamma (n
 - {\displaystyle \frac {3}{2}}  + {\displaystyle \frac {a}{2}} )
 \\
\Gamma (n - {\displaystyle \frac {3}{2}}  + b - {\displaystyle 
\frac {a}{2}} )\,\Gamma ( - {\displaystyle \frac {a}{2}}  + b + 
{\displaystyle \frac {3}{2}}  - n)\,\Gamma (1 + a - b)\,\Gamma (b
)\,\Gamma (2 - 2\,b)\,\mathrm{sin}(\pi \,b) \\
\mathrm{cos}({\displaystyle \frac {1}{2}} \,\pi \,a - \pi \,b)
 \left/ {\vrule height0.80em width0em depth0.80em} \right. \! 
 \! (\pi ^{2}\,\Gamma (1 + a - 2\,b)\,\Gamma (n - {\displaystyle 
\frac {3}{2}}  - {\displaystyle \frac {a}{2}} )\,(-1)^{n})\mbox{}
 + ( - {\displaystyle \frac {1}{32}} 4^{b} \\
\Gamma (n - {\displaystyle \frac {3}{2}}  + {\displaystyle 
\frac {a}{2}} )\,\Gamma (n - {\displaystyle \frac {3}{2}}  + b - 
{\displaystyle \frac {a}{2}} )\,\Gamma ( - {\displaystyle \frac {
a}{2}}  + b + {\displaystyle \frac {3}{2}}  - n)\,\Gamma (1 + a
 - b)\,\Gamma (b - {\displaystyle \frac {1}{2}} )\,\Gamma (2 - 2
\,b) \\
(2\,\mathrm{sin}(\pi \,a) + \mathrm{sin}(2\,\pi \,b) + 2\,
\mathrm{sin}(\pi \,a - 2\,\pi \,b) + \mathrm{sin}( - 2\,\pi \,b
 + 2\,\pi \,a) \\
\mbox{} + \mathrm{sin}(2\,\pi \,a - 4\,\pi \,b)) \left/ {\vrule 
height0.80em width0em depth0.80em} \right. \!  \! (\pi ^{(3/2)}\,
\Gamma (n - {\displaystyle \frac {3}{2}}  - {\displaystyle 
\frac {a}{2}} )\,\mathrm{cos}({\displaystyle \frac {\pi \,( - a
 + 2\,n)}{2}} )\,\mathrm{sin}(\pi \,(a - b))) \\
\mbox{} - \pi ^{2}\,(\mathrm{sin}({\displaystyle \frac {\pi \,a}{
2}} )\,\mathrm{sin}(\pi \,b) - \mathrm{cos}({\displaystyle 
\frac {\pi \,a}{2}} )\,\mathrm{cos}(\pi \,b))\,\Gamma (n - 
{\displaystyle \frac {3}{2}}  + {\displaystyle \frac {a}{2}} )\,
\Gamma (1 + a - b) \left/ {\vrule 
height0.80em width0em depth0.80em} \right. \!  \! \\ (
\mathrm{cos}({\displaystyle \frac {\pi \,a}{2}} )\,\mathrm{sin}(
\pi \,b)\,\mathrm{sin}({\displaystyle \frac {\pi \,( - a + 2\,n
 + 1 - 2\,b)}{2}} )\,\Gamma (n - {\displaystyle \frac {3}{2}}  - 
{\displaystyle \frac {a}{2}} )\,\Gamma ( - b + {\displaystyle 
\frac {5}{2}}  + {\displaystyle \frac {a}{2}}  - n) 
\Gamma (n + {\displaystyle \frac {a}{2}}  - b - {\displaystyle 
\frac {1}{2}} )\,\Gamma (b)))/\Gamma (a) }
}
\end{maplelatex}

\begin{maplelatex}
\mapleinline{inert}{2d}{237, ":   ",
1/2*b*sin(Pi*(n-a+b-1))*(-1)^n*sin(Pi*a)*Sum(GAMMA(-a+b+k)*GAMMA(a+2-2
*n+k)/GAMMA(3+a-2*n-b+k)/GAMMA(1-a+k),k = 1 ..
2*n-3)/sin(Pi*(2*n-2-a))/Pi/GAMMA(2-2*n+2*a)*GAMMA(3+2*a-2*n-b)*GAMMA(
-a+1)*GAMMA(a-b)+1/8*4^(-b)/Pi^(3/2)*GAMMA(b+1)/GAMMA(-2*b+3+2*a-2*n)*
GAMMA(3+2*a-2*n-b)*GAMMA(-b+1/2)*GAMMA(-a+1)*GAMMA(a-b)*(sin(2*Pi*b)-s
in(4*Pi*a-2*Pi*b)+2*sin(-2*Pi*b+2*Pi*a)-sin(4*Pi*a)+2*sin(2*Pi*a))/sin
(Pi*a)/sin(Pi*(-2*a+b))-(-cot(Pi*a-Pi*b)+cot(Pi*b))*b/GAMMA(-2*b+3+2*a
-2*n)/GAMMA(2*b)*GAMMA(3+2*a-2*n-b)*GAMMA(b)^2*GAMMA(a-b)/GAMMA(a)*sin
(Pi*a-Pi*b)*sin(Pi*b)/sin(2*Pi*b-Pi*a+2*Pi*n);}{%
\maplemultiline{
237\mbox{:~~~} {\displaystyle \frac {1}{2}} b\,\mathrm{
sin}(\pi \,(n - a + b - 1))\,(-1)^{n}\,\mathrm{sin}(\pi \,a) \\
 \left(  \! {\displaystyle \sum _{k=1}^{2\,n - 3}} \,
{\displaystyle \frac {\Gamma ( - a + b + k)\,\Gamma (a + 2 - 2\,n
 + k)}{\Gamma (3 + a - 2\,n - b + k)\,\Gamma (1 - a + k)}}  \! 
 \right) \,\Gamma (3 + 2\,a - 2\,n - b)\,\Gamma ( - a + 1)\,
\Gamma (a - b) \\
/(\mathrm{sin}(\pi \,(2\,n - 2 - a))\,\pi \,\Gamma (2 - 2\,n + 2
\,a))\mbox{} + {\displaystyle \frac {1}{8}} 4^{( - b)}\,\Gamma (b
 + 1)\,\Gamma (3 + 2\,a - 2\,n - b) 
\Gamma ( - b + {\displaystyle \frac {1}{2}} )\,\Gamma ( - a + 1)
\,\Gamma (a - b) \\
(\mathrm{sin}(2\,\pi \,b) - \mathrm{sin}(4\,\pi \,a - 2\,\pi \,b)
 + 2\,\mathrm{sin}( - 2\,\pi \,b + 2\,\pi \,a) - \mathrm{sin}(4\,
\pi \,a) + 2\,\mathrm{sin}(2\,\pi \,a)) \\
 \left/ {\vrule height0.63em width0em depth0.63em} \right. \! 
 \! (\pi ^{(3/2)}\,\Gamma ( - 2\,b + 3 + 2\,a - 2\,n)\,\mathrm{
sin}(\pi \,a)\,\mathrm{sin}(\pi \,( - 2\,a + b)))\mbox{} \\ - 
{\displaystyle \frac {( - \mathrm{cot}(\pi \,a - \pi \,b) + 
\mathrm{cot}(\pi \,b))\,b\,\Gamma (3 + 2\,a - 2\,n - b)\,\Gamma (
b)^{2}\,\Gamma (a - b)\,\mathrm{sin}(\pi \,a - \pi \,b)\,\mathrm{
sin}(\pi \,b)}{\Gamma ( - 2\,b + 3 + 2\,a - 2\,n)\,\Gamma (2\,b)
\,\Gamma (a)\,\mathrm{sin}(2\,\pi \,b - \pi \,a + 2\,\pi \,n)}} 
 }
}
\end{maplelatex}

\begin{maplelatex}
\mapleinline{inert}{2d}{238, ":   ",
-1/2*Sum(GAMMA(b-1+k)*GAMMA(-1/2*a-1/2*n+1+k)/GAMMA(1+1/2*a-1/2*n+k)/G
AMMA(3-b-n+k),k = 1 ..
n-2)*sin(1/2*Pi*(a+n))/Pi*(-1/2*sin(Pi*(a+b))+1/2*sin(Pi*(a-b)))/((-1)
^n)/sin(Pi*b)/sin(1/2*Pi*(2*a+1))*GAMMA(1+a-b)/GAMMA(n+b-2)/GAMMA(3-2*
b-n+a)*GAMMA(1/2*a-1/2*n+1)*GAMMA(1/2*n+1/2*a)*GAMMA(-n-2*b+4)+(sin(1/
2*Pi*(1+2*b))*Pi/sin(Pi*(1/2*n+b+1/2*a))/GAMMA(-1/2*n-b+2+1/2*a)^2*GAM
MA(1/2*a-1/2*n+1)/GAMMA(-1/2*n-1/2*a+1)*GAMMA(-b+3-n)-1/2*Pi*GAMMA(1/2
*n+1/2*a)/GAMMA(n+b-2)/GAMMA(-1/2*n-b+2+1/2*a)^2*GAMMA(1/2*a-1/2*n+1)*
(1/4*cos(2*Pi*b)+1/4*cos(-2*Pi*b+2*Pi*a)-1/4*cos(-Pi*a+Pi*n)-1/4*cos(P
i*a+Pi*n))/((-1)^n)/sin(Pi*b)/sin(1/2*Pi*n+Pi*b+1/2*Pi*a)/sin(1/2*Pi*n
-1/2*Pi*a+Pi*b)/cos(Pi*a))*GAMMA(1+a-b)/GAMMA(a);}{%
\maplemultiline{
238\mbox{:~~~}  - {\displaystyle \frac {1}{2}}  \left( 
 \! {\displaystyle \sum _{k=1}^{n - 2}} \,{\displaystyle \frac {
\Gamma (b - 1 + k)\,\Gamma ( - {\displaystyle \frac {a}{2}}  - 
{\displaystyle \frac {n}{2}}  + 1 + k)}{\Gamma (1 + 
{\displaystyle \frac {a}{2}}  - {\displaystyle \frac {n}{2}}  + k
)\,\Gamma (3 - b - n + k)}}  \!  \right) \,\mathrm{sin}(
{\displaystyle \frac {\pi \,(a + n)}{2}} ) \\
( - {\displaystyle \frac {1}{2}} \,\mathrm{sin}(\pi \,(a + b)) + 
{\displaystyle \frac {1}{2}} \,\mathrm{sin}(\pi \,(a - b)))\,
\Gamma (1 + a - b)\,\Gamma ({\displaystyle \frac {a}{2}}  - 
{\displaystyle \frac {n}{2}}  + 1)\,\Gamma ({\displaystyle 
\frac {n}{2}}  + {\displaystyle \frac {a}{2}} ) \\
\Gamma ( - n - 2\,b + 4) \left/ {\vrule 
height0.80em width0em depth0.80em} \right. \!  \! (\pi \,(-1)^{n}
\,\mathrm{sin}(\pi \,b)\,\mathrm{sin}({\displaystyle \frac {\pi 
\,(2\,a + 1)}{2}} )\,\Gamma (n + b - 2) \\
\Gamma (3 - 2\,b - n + a))\mbox{} +  \left( {\vrule 
height1.74em width0em depth1.74em} \right. \!  \! {\displaystyle 
\frac {\mathrm{sin}({\displaystyle \frac {\pi \,(1 + 2\,b)}{2}} )
\,\pi \,\Gamma ({\displaystyle \frac {a}{2}}  - {\displaystyle 
\frac {n}{2}}  + 1)\,\Gamma ( - b + 3 - n)}{\mathrm{sin}(\pi \,(
{\displaystyle \frac {n}{2}}  + b + {\displaystyle \frac {a}{2}} 
))\,\Gamma ( - {\displaystyle \frac {n}{2}}  - b + 2 + 
{\displaystyle \frac {a}{2}} )^{2}\,\Gamma ( - {\displaystyle 
\frac {n}{2}}  - {\displaystyle \frac {a}{2}}  + 1)}} \\  - 
{\displaystyle \frac {1}{2}} \pi
\Gamma ({\displaystyle \frac {n}{2}}  + {\displaystyle \frac {a}{
2}} )\,\Gamma ({\displaystyle \frac {a}{2}}  - {\displaystyle 
\frac {n}{2}}  + 1)
({\displaystyle \frac {1}{4}} \,\mathrm{cos}(2\,\pi \,b) + 
{\displaystyle \frac {1}{4}} \,\mathrm{cos}( - 2\,\pi \,b + 2\,
\pi \,a) - {\displaystyle \frac {1}{4}} \,\mathrm{cos}( - \pi \,a
 + \pi \,n) - {\displaystyle \frac {1}{4}} \,\mathrm{cos}(\pi \,a
 + \pi \,n)) \\
 \left/ {\vrule height0.80em width0em depth0.80em} \right. \! 
 \! (\Gamma (n + b - 2)\,\Gamma ( - {\displaystyle \frac {n}{2}} 
 - b + 2 + {\displaystyle \frac {a}{2}} )^{2}\,(-1)^{n}\,\mathrm{
sin}(\pi \,b)\,\mathrm{sin}({\displaystyle \frac {\pi \,n}{2}} 
 + \pi \,b + {\displaystyle \frac {\pi \,a}{2}} ) 
\mathrm{sin}({\displaystyle \frac {\pi \,n}{2}}  - 
{\displaystyle \frac {\pi \,a}{2}}  + \pi \,b)\,\mathrm{cos}(\pi 
\,a)) \! \! \left. {\vrule height1.74em width0em depth1.74em}
 \right) \Gamma (1 + a - b)/\Gamma (a) }
}
\end{maplelatex}

\begin{maplelatex}
\mapleinline{inert}{2d}{239, ":   ",
1/2*(sin(2*Pi*a)*sin(Pi*(2*b+a))*sin(Pi*(a+b))+sin(Pi*a)*sin(2*Pi*(a+b
))*sin(Pi*(a-b)))/sin(Pi*a)*Sum(GAMMA(-a+k)*GAMMA(2+a-2*n+b+k)/GAMMA(3
-2*n+a+k)/GAMMA(-a-b+1+k),k = 1 ..
2*n-3)/sin(Pi*(2*b+a))/sin(Pi*(2*a+b))/GAMMA(a+2+2*b-2*n)*GAMMA(-b+1-a
)*2^(2*b-1)*Pi^(3/2)/GAMMA(-b+1)/sin(1/2*Pi*(1+2*b))/GAMMA(-b+1/2)/sin
(Pi*(a-b))/GAMMA(b-a)/GAMMA(2*n-2-a)+(2*sin(3*Pi*a+Pi*b)-sin(3*Pi*a-Pi
*b)+4*sin(Pi*a+Pi*b)-2*sin(Pi*a-Pi*b)-sin(3*Pi*b+3*Pi*a)+sin(5*Pi*a+5*
Pi*b)-2*sin(3*Pi*b+5*Pi*a)-2*sin(Pi*a+3*Pi*b)+sin(5*Pi*a+Pi*b))/(-cos(
-2*Pi*b+2*Pi*a)-cos(4*Pi*a-2*Pi*b)+1-cos(2*Pi*a)-cos(2*Pi*b+2*Pi*a)+2*
cos(4*Pi*b)-cos(2*Pi*a+4*Pi*b)-cos(4*Pi*a+6*Pi*b)+cos(6*Pi*b+6*Pi*a)+c
os(4*Pi*b+6*Pi*a)-cos(4*Pi*b+8*Pi*a)-cos(4*Pi*a+4*Pi*b)-cos(4*Pi*a)+co
s(6*Pi*a)+2*cos(2*Pi*b)+cos(2*Pi*b+6*Pi*a))/GAMMA(2*n-2-a)/sin(Pi*b)^2
/GAMMA(-b+1)^2/GAMMA(b-a)/GAMMA(a+2+2*b-2*n)/GAMMA(2*n-2*b-2*a-1)/GAMM
A(2*a-2*n+3)*Pi^3*GAMMA(-b+1-a);}{%
\maplemultiline{
239\mbox{:~~~} \\ {\displaystyle \frac {1}{2}}
(\mathrm{sin}(2\,\pi \,a)\,\mathrm{sin}(\pi \,(2\,b + a))\,
\mathrm{sin}(\pi \,(a + b)) + \mathrm{sin}(\pi \,a)\,\mathrm{sin}
(2\,\pi \,(a + b))\,\mathrm{sin}(\pi \,(a - b))) \\
 \left(  \! {\displaystyle \sum _{k=1}^{2\,n - 3}} \,
{\displaystyle \frac {\Gamma ( - a + k)\,\Gamma (2 + a - 2\,n + b
 + k)}{\Gamma (3 - 2\,n + a + k)\,\Gamma ( - a - b + 1 + k)}} 
 \!  \right) \,\Gamma ( - b + 1 - a)\,2^{(2\,b - 1)}\,\pi ^{(3/2)
} \left/ {\vrule height0.80em width0em depth0.80em} \right. \! 
 \! \\ (
\mathrm{sin}(\pi \,a)\,\mathrm{sin}(\pi \,(2\,b + a))\,\mathrm{
sin}(\pi \,(2\,a + b))\,\Gamma (a + 2 + 2\,b - 2\,n)\,\Gamma ( - 
b + 1) \\
\mathrm{sin}({\displaystyle \frac {\pi \,(1 + 2\,b)}{2}} )\,
\Gamma ( - b + {\displaystyle \frac {1}{2}} )\,\mathrm{sin}(\pi 
\,(a - b))\,\Gamma (b - a)\,\Gamma (2\,n - 2 - a))\mbox{} + \\ (
2\,\mathrm{sin}(3\,\pi \,a + \pi \,b) - \mathrm{sin}(3\,\pi \,a
 - \pi \,b) + 4\,\mathrm{sin}(\pi \,a + \pi \,b) - 2\,\mathrm{sin
}(\pi \,a - \pi \,b) \\
\mbox{} - \mathrm{sin}(3\,\pi \,b + 3\,\pi \,a) + \mathrm{sin}(5
\,\pi \,a + 5\,\pi \,b) - 2\,\mathrm{sin}(3\,\pi \,b + 5\,\pi \,a
) - 2\,\mathrm{sin}(\pi \,a + 3\,\pi \,b) \\
\mbox{} + \mathrm{sin}(5\,\pi \,a + \pi \,b))\pi ^{3}\,\Gamma (
 - b + 1 - a) \left/ {\vrule height0.43em width0em depth0.43em}
 \right. \!  \! (( - \mathrm{cos}( - 2\,\pi \,b + 2\,\pi \,a) \\
\mbox{} - \mathrm{cos}(4\,\pi \,a - 2\,\pi \,b) + 1 - \mathrm{cos
}(2\,\pi \,a) - \mathrm{cos}(2\,\pi \,b + 2\,\pi \,a) + 2\,
\mathrm{cos}(4\,\pi \,b) \\
\mbox{} - \mathrm{cos}(2\,\pi \,a + 4\,\pi \,b) - \mathrm{cos}(4
\,\pi \,a + 6\,\pi \,b) + \mathrm{cos}(6\,\pi \,b + 6\,\pi \,a)
 + \mathrm{cos}(4\,\pi \,b + 6\,\pi \,a) \\
\mbox{} - \mathrm{cos}(4\,\pi \,b + 8\,\pi \,a) - \mathrm{cos}(4
\,\pi \,a + 4\,\pi \,b) - \mathrm{cos}(4\,\pi \,a) + \mathrm{cos}
(6\,\pi \,a) + 2\,\mathrm{cos}(2\,\pi \,b) \\
\mbox{} + \mathrm{cos}(2\,\pi \,b + 6\,\pi \,a))\Gamma (2\,n - 2
 - a)\,\mathrm{sin}(\pi \,b)^{2}\,\Gamma ( - b + 1)^{2}\,\Gamma (
b - a) \\
\Gamma (a + 2 + 2\,b - 2\,n)\,\Gamma (2\,n - 2\,b - 2\,a - 1)\,
\Gamma (2\,a - 2\,n + 3)) }
}
\end{maplelatex}

\mapleinline{inert}{2d}{240, ":   ",
-1/4*Sum(GAMMA(1/2*b-a+k)*GAMMA(a+2-2*n+k)/GAMMA(3+a-2*n-1/2*b+k)/GAMM
A(1-a+k),k = 1 ..
2*n-3)*2^b*Pi^(1/2)*GAMMA(1+a-b)*GAMMA(1/2*b+1/2)/sin(Pi*(-1/2*b+a))/G
AMMA(-a+2*n-2+1/2*b)*sin(1/2*Pi*(-2*a+2*n+b))*sin(Pi*(n-a))*sin(Pi*b)/
GAMMA(a)/GAMMA(-2*n+2+a)/sin(Pi*a)^2/sin(1/2*Pi*b)+(sin(Pi*b)*sin(-1/2
*Pi*b+Pi*a)/sin(Pi*a)/sin(1/2*Pi*b)+sin(2*Pi*a)/sin(Pi*(2-2*n+2*a))/si
n(Pi*a)*sin(Pi*(b-a))-1/16*sin(Pi*(-b+3+2*a))/sin(Pi*(-1/2*b+a))*(sin(
4*Pi*a-2*Pi*b)+sin(-Pi*b+4*Pi*a)-2*sin(2*Pi*a)-2*sin(-Pi*b+2*Pi*a)+sin
(Pi*b))*sin(Pi*b)/sin(-Pi*b+2*Pi*a)/sin(-1/2*Pi*b+2*Pi*a)/sin(1/2*Pi*b
)/cos(1/2*Pi*b)/sin(Pi*a))/GAMMA(2*n-2*a-1)*GAMMA(1+a-b)*GAMMA(-2+b-2*
a+2*n)/GAMMA(-a+2*n-2+1/2*b)*GAMMA(1/2*b)/GAMMA(a)*GAMMA(2*n-1-a);}{%
\maplemultiline{
240\mbox{:~~~}  - {\displaystyle \frac {1}{4}}  \left( 
 \! {\displaystyle \sum _{k=1}^{2\,n - 3}} \,{\displaystyle 
\frac {\Gamma ({\displaystyle \frac {b}{2}}  - a + k)\,\Gamma (a
 + 2 - 2\,n + k)}{\Gamma (3 + a - 2\,n - {\displaystyle \frac {b
}{2}}  + k)\,\Gamma (1 - a + k)}}  \!  \right) \,2^{b}\,\sqrt{\pi
 }\,\Gamma (1 + a - b)\,\Gamma ({\displaystyle \frac {b}{2}}  + 
{\displaystyle \frac {1}{2}} ) \\
\mathrm{sin}({\displaystyle \frac {\pi \,( - 2\,a + 2\,n + b)}{2}
} )\,\mathrm{sin}(\pi \,(n - a))\,\mathrm{sin}(\pi \,b) \left/ 
{\vrule height0.80em width0em depth0.80em} \right. \!  \! (
\mathrm{sin}(\pi \,( - {\displaystyle \frac {b}{2}}  + a))
\Gamma ( - a + 2\,n - 2 + {\displaystyle \frac {b}{2}} )\,\Gamma 
(a)\,\Gamma ( - 2\,n + 2 + a)\,\mathrm{sin}(\pi \,a)^{2}\,
\mathrm{sin}({\displaystyle \frac {\pi \,b}{2}} ))\mbox{} + 
 \\  \left( {\vrule height1.67em width0em depth1.67em} \right. \! 
 \!
{\displaystyle \frac {\mathrm{sin}(\pi \,b)\,\mathrm{sin}( - 
{\displaystyle \frac {1}{2}} \,\pi \,b + \pi \,a)}{\mathrm{sin}(
\pi \,a)\,\mathrm{sin}({\displaystyle \frac {\pi \,b}{2}} )}}  + 
{\displaystyle \frac {\mathrm{sin}(2\,\pi \,a)\,\mathrm{sin}(\pi 
\,(b - a))}{\mathrm{sin}(\pi \,(2 - 2\,n + 2\,a))\,\mathrm{sin}(
\pi \,a)}}  - {\displaystyle \frac {1}{16}}
\mathrm{sin}(\pi \,( - b + 3 + 2\,a))(\mathrm{sin}(4\,\pi \,a - 2
\,\pi \,b) + \mathrm{sin}( - \pi \,b + 4\,\pi \,a) \\ 
- 2\,\mathrm{sin}(2\,\pi \,a) 
\mbox{} - 2\,\mathrm{sin}( - \pi \,b + 2\,\pi \,a) + \mathrm{sin}
(\pi \,b))\mathrm{sin}(\pi \,b)/(\mathrm{sin}(\pi \,( - 
{\displaystyle \frac {b}{2}}  + a))\,\mathrm{sin}( - \pi \,b + 2
\,\pi \,a) \\
\mathrm{sin}( - {\displaystyle \frac {1}{2}} \,\pi \,b + 2\,\pi 
\,a)\,\mathrm{sin}({\displaystyle \frac {\pi \,b}{2}} )\,\mathrm{
cos}({\displaystyle \frac {\pi \,b}{2}} )\,\mathrm{sin}(\pi \,a))
 \! \! \left. {\vrule height1.67em width0em depth1.67em} \right) 
\Gamma (1 + a - b)\,\Gamma ( - 2 + b - 2\,a + 2\,n) \\
\Gamma ({\displaystyle \frac {b}{2}} )\,\Gamma (2\,n - 1 - a)
 \left/ {\vrule height0.80em width0em depth0.80em} \right. \! 
 \! (\Gamma (2\,n - 2\,a - 1)\,\Gamma ( - a + 2\,n - 2 + 
{\displaystyle \frac {b}{2}} )\,\Gamma (a)) }
}

\begin{mapleinput}
\end{mapleinput}

\end{maplegroup}

%% file: AppendixB241to260.tex
\begin{maplegroup}
\mapleresult
\begin{maplelatex}
\mapleinline{inert}{2d}{241, ":   ",
1/2*Sum(GAMMA(b-1+k)*GAMMA(-a-b+3+k-n)/GAMMA(3-b-n+k)/GAMMA(-1+b+a+k),
k = 1 ..
n-2)*(2*a-1)/Pi*GAMMA(1+a-b)*GAMMA(-1+a+b)*sin(Pi*a+Pi*b)-1/2*1/Pi^(1/
2)*cos(Pi*b)*4^a/GAMMA(a)*GAMMA(a+1/2)*GAMMA(-1+a+b)*GAMMA(n+2*b-3)/GA
MMA(n+2*a+2*b-4)*GAMMA(1+a-b);}{%
\maplemultiline{
241\mbox{:~~~} {\displaystyle \frac {1}{2}}  \left(  \! 
{\displaystyle \sum _{k=1}^{n - 2}} \,{\displaystyle \frac {
\Gamma (b - 1 + k)\,\Gamma ( - a - b + 3 + k - n)}{\Gamma (3 - b
 - n + k)\,\Gamma ( - 1 + b + a + k)}}  \!  \right) \,(2\,a - 1)
\,\Gamma (1 + a - b)
\Gamma ( - 1 + a + b)\,\mathrm{sin}(\pi \,a + \pi \,b)/\pi  \\
\mbox{} - {\displaystyle \frac {1}{2}} \,{\displaystyle \frac {
\mathrm{cos}(\pi \,b)\,4^{a}\,\Gamma (a + {\displaystyle \frac {1
}{2}} )\,\Gamma ( - 1 + a + b)\,\Gamma (n + 2\,b - 3)\,\Gamma (1
 + a - b)}{\sqrt{\pi }\,\Gamma (a)\,\Gamma (n + 2\,a + 2\,b - 4)}
}  }
}
\end{maplelatex}

\begin{maplelatex}
\mapleinline{inert}{2d}{242, ":   ",
-1/2*(b-1)*Sum(GAMMA(a-1+k)*GAMMA(1/2*b+2-a+k-n)/GAMMA(3-n-a+k)/GAMMA(
a-1/2*b+k),k = 1 ..
n-2)/GAMMA(a-1)*GAMMA(a-2+n-1/2*b)/GAMMA(n+a-2)*GAMMA(-1/2*b+a)*sin(-1
/2*Pi*b+Pi*n+Pi*a)/sin(Pi*a+Pi*n)-(4*sin(1/2*Pi*(2*a+b))*sin(Pi*a)*cos
(1/2*Pi*b)-sin(2*Pi*a)*sin(Pi*b))*Pi^(1/2)*2^(-b)*GAMMA(1/2*b)*GAMMA(-
1/2*b+a)/GAMMA(1/2*b-1/2)/GAMMA(n+a-2)/cos(1/2*Pi*b)/sin(Pi*b)/sin(Pi*
a)^2/GAMMA(n-b-2+2*a)*sin(1/2*Pi*b)/GAMMA(a-1)*GAMMA(a-2+n-1/2*b)*GAMM
A(2*a+n-3)-Pi*(cot(Pi*a)-cot(Pi*b))*GAMMA(1/2*b)^2/GAMMA(n+a-2)/GAMMA(
-a+1+1/2*b)/sin(Pi*a-Pi*b)/GAMMA(n-b-2+2*a)*sin(1/2*Pi*b)/GAMMA(a-1)*G
AMMA(a-2+n-1/2*b)*GAMMA(2*a+n-3)/GAMMA(b-1);}{%
\maplemultiline{
242\mbox{:~~~}  - {\displaystyle \frac {1}{2}} (b - 1)\,
 \left(  \! {\displaystyle \sum _{k=1}^{n - 2}} \,{\displaystyle 
\frac {\Gamma (a - 1 + k)\,\Gamma ({\displaystyle \frac {b}{2}} 
 + 2 - a + k - n)}{\Gamma (3 - n - a + k)\,\Gamma (a - 
{\displaystyle \frac {b}{2}}  + k)}}  \!  \right) \,\Gamma (a - 2
 + n - {\displaystyle \frac {b}{2}} ) \\
\Gamma ( - {\displaystyle \frac {b}{2}}  + a)\,\mathrm{sin}( - 
{\displaystyle \frac {\pi \,b}{2}}  + \pi \,n + \pi \,a)/(\Gamma 
(a - 1)\,\Gamma (n + a - 2)\,\mathrm{sin}(\pi \,a + \pi \,n))
\mbox{} -  \\
(4\,\mathrm{sin}({\displaystyle \frac {\pi \,(2\,a + b)}{2}} )\,
\mathrm{sin}(\pi \,a)\,\mathrm{cos}({\displaystyle \frac {\pi \,b
}{2}} ) - \mathrm{sin}(2\,\pi \,a)\,\mathrm{sin}(\pi \,b))\,
\sqrt{\pi }\,2^{( - b)}\,\Gamma ({\displaystyle \frac {b}{2}} )
 \\
\Gamma ( - {\displaystyle \frac {b}{2}}  + a)\,\mathrm{sin}(
{\displaystyle \frac {\pi \,b}{2}} )\,\Gamma (a - 2 + n - 
{\displaystyle \frac {b}{2}} )\,\Gamma (2\,a + n - 3) \left/ 
{\vrule height0.80em width0em depth0.80em} \right. \!  \! (\Gamma
 ({\displaystyle \frac {b}{2}}  - {\displaystyle \frac {1}{2}} )
\,\Gamma (n + a - 2) \\
\mathrm{cos}({\displaystyle \frac {\pi \,b}{2}} )\,\mathrm{sin}(
\pi \,b)\,\mathrm{sin}(\pi \,a)^{2}\,\Gamma (n - b - 2 + 2\,a)\,
\Gamma (a - 1)) \\
\mbox{} - {\displaystyle \frac {\pi \,(\mathrm{cot}(\pi \,a) - 
\mathrm{cot}(\pi \,b))\,\Gamma ({\displaystyle \frac {b}{2}} )^{2
}\,\mathrm{sin}({\displaystyle \frac {\pi \,b}{2}} )\,\Gamma (a
 - 2 + n - {\displaystyle \frac {b}{2}} )\,\Gamma (2\,a + n - 3)
}{\Gamma (n + a - 2)\,\Gamma ( - a + 1 + {\displaystyle \frac {b
}{2}} )\,\mathrm{sin}(\pi \,a - \pi \,b)\,\Gamma (n - b - 2 + 2\,
a)\,\Gamma (a - 1)\,\Gamma (b - 1)}}  }
}
\end{maplelatex}

\begin{maplelatex}
\mapleinline{inert}{2d}{243, ":   ",
-1/2*Sum(GAMMA(a-1+k)*GAMMA(b-1+k)/GAMMA(4-a-2*n+k)/GAMMA(4-b-2*n+k),k
= 1 ..
2*n-3)*Pi^2/GAMMA(a-1)/GAMMA(b-1)/GAMMA(a-3+2*n)/GAMMA(2*n-3+b)/sin(Pi
*b)/sin(Pi*a)+((cot(Pi*a)-cot(Pi*b))/sin(Pi*(a-b))*2^(8-2*a-4*n-2*b)*G
AMMA(9/2-2*n-b-a)/GAMMA(5-a-b-2*n)/sin(Pi*a)/GAMMA(a-3+2*n)/GAMMA(-2*a
+5-2*n)/GAMMA(-2*b+5-2*n)*GAMMA(-b+4-2*n)/GAMMA(a-1)*Pi^(5/2)-32*GAMMA
(b+a-4+2*n)*GAMMA(2*a-4+2*n)*GAMMA(9/2-2*n-b-a)*2^(-2*b-2*a)*16^(-n)*(
sin(2*Pi*b+2*Pi*a)-sin(4*Pi*a)+sin(-2*Pi*b+2*Pi*a))/GAMMA(a-3+2*n)/GAM
MA(-2*b+5-2*n)/GAMMA(2*n-3+b)/sin(Pi*a)^2/sin(Pi*b)^2/sin(Pi*(a-b))/GA
MMA(a-1)*Pi^(3/2))/GAMMA(b-1);}{%
\maplemultiline{
243\mbox{:~~~}  - {\displaystyle \frac {1}{2}} \,
{\displaystyle \frac { \left(  \! {\displaystyle \sum _{k=1}^{2\,
n - 3}} \,{\displaystyle \frac {\Gamma (a - 1 + k)\,\Gamma (b - 1
 + k)}{\Gamma (4 - a - 2\,n + k)\,\Gamma (4 - b - 2\,n + k)}} 
 \!  \right) \,\pi ^{2}}{\Gamma (a - 1)\,\Gamma (b - 1)\,\Gamma (
a - 3 + 2\,n)\,\Gamma (2\,n - 3 + b)\,\mathrm{sin}(\pi \,b)\,
\mathrm{sin}(\pi \,a)}}  + \\ ( 
(\mathrm{cot}(\pi \,a) - \mathrm{cot}(\pi \,b))\,2^{(8 - 2\,a - 4
\,n - 2\,b)}\,\Gamma ({\displaystyle \frac {9}{2}}  - 2\,n - b - 
a)\,\Gamma ( - b + 4 - 2\,n)\,\pi ^{(5/2)} \left/ {\vrule 
height0.37em width0em depth0.37em} \right. \!  \!  \\
(\mathrm{sin}(\pi \,(a - b))\,\Gamma (5 - a - b - 2\,n)\,\mathrm{
sin}(\pi \,a)\,\Gamma (a - 3 + 2\,n)\,\Gamma ( - 2\,a + 5 - 2\,n)
 \\
\Gamma ( - 2\,b + 5 - 2\,n)\,\Gamma (a - 1))\mbox{} - 32\,\Gamma 
(b + a - 4 + 2\,n)\,\Gamma (2\,a - 4 + 2\,n) \\
\Gamma ({\displaystyle \frac {9}{2}}  - 2\,n - b - a)\,2^{( - 2\,
b - 2\,a)}\,16^{( - n)} \\
(\mathrm{sin}(2\,\pi \,b + 2\,\pi \,a) - \mathrm{sin}(4\,\pi \,a)
 + \mathrm{sin}( - 2\,\pi \,b + 2\,\pi \,a))\,\pi ^{(3/2)}
 \left/ {\vrule height0.44em width0em depth0.44em} \right. \! 
 \! (\Gamma (a - 3 + 2\,n) \\
\Gamma ( - 2\,b + 5 - 2\,n)\,\Gamma (2\,n - 3 + b)\,\mathrm{sin}(
\pi \,a)^{2}\,\mathrm{sin}(\pi \,b)^{2}\,\mathrm{sin}(\pi \,(a - 
b))\,\Gamma (a - 1))) \left/ {\vrule 
height0.37em width0em depth0.37em} \right. \!  \!
\Gamma (b - 1) }
}
\end{maplelatex}

\begin{maplelatex}
\mapleinline{inert}{2d}{244, ":   ",
Sum(GAMMA(a+k+1-n-c)*GAMMA(-b-n+1+a+k)/GAMMA(k+1)/GAMMA(a+1+k-n),k = 0
..
n)*GAMMA(a)*GAMMA(b)*GAMMA(a+2-b-c-n)*GAMMA(2+n)*(n-a)*(-b*c+c+a*c+a*b
-1-2*a+b+a*n-a^2)*Pi/GAMMA(-c+1)/GAMMA(a-n+1)/GAMMA(b-a+n)/GAMMA(-b-n+
1+a)/GAMMA(2+a-b)/GAMMA(2+a-c)/sin(Pi*a)/cos(Pi*a)-Sum(GAMMA(b-a-1+k)*
GAMMA(a+k+1-n-c)/GAMMA(k+1)/GAMMA(-c-n+2+k),k = 0 ..
n)*GAMMA(b)*GAMMA(a+2-b-c-n)*(b+c-2+n-a)*GAMMA(2+n)*cos(Pi*b)*(n-a)*(-
b*c+c+a*c+a*b-1-2*a+b+a*n-a^2)*Pi/GAMMA(n-a)/GAMMA(a-n+1)^2/GAMMA(b-a+
n)/GAMMA(2+a-c)/sin(Pi*a)/(b-1)/cos(Pi*a)+GAMMA(a)*GAMMA(b)*GAMMA(a+2-
b-c-n)*GAMMA(2+n)*(n-a)*Pi/GAMMA(n)/GAMMA(-c+1)/GAMMA(a-n+1)/GAMMA(b-a
+n)/GAMMA(1+a)/GAMMA(-b-n+1+a)/n/sin(Pi*a)/cos(Pi*a)+GAMMA(1+a)*GAMMA(
1+a-c-n)*GAMMA(2+n)*Pi^2*(n-a)*(-b*c+c+a*c+a*b-1-2*a+b+a*n-a^2)/GAMMA(
n-a)/GAMMA(a-n+1)^2/GAMMA(b-a+n)/GAMMA(2+a-b)/GAMMA(2+a-c)/sin(Pi*(a-b
))/sin(Pi*a)/a+GAMMA(b)*GAMMA(a+2-b-c-n)*(b+c-2+n-a)*GAMMA(2+n)*cos(Pi
*b)*(n-a)*Pi/GAMMA(n-a)/GAMMA(n)/GAMMA(a-n+1)^2/GAMMA(2-c)/sin(Pi*a)/(
b-1)/n/cos(Pi*a);}{%
\maplemultiline{
244\mbox{:~~~}  \left(  \! {\displaystyle \sum _{k=0}^{n
}} \,{\displaystyle \frac {\Gamma (a + k + 1 - n - c)\,\Gamma (
 - b - n + 1 + a + k)}{\Gamma (k + 1)\,\Gamma (a + 1 + k - n)}} 
 \!  \right) \,\Gamma (a)\,\Gamma (b)\,\mathrm{\%1}\,\Gamma (2 + 
n) \\
(n - a)\,\mathrm{\%2}\,\pi /(\Gamma ( - c + 1)\,\Gamma (a - n + 1
)\,\Gamma (b - a + n)\,\Gamma ( - b - n + 1 + a)\,\Gamma (2 + a
 - b) \\
\Gamma (2 + a - c)\,\mathrm{sin}(\pi \,a)\,\mathrm{cos}(\pi \,a))
\mbox{} -  \left(  \! {\displaystyle \sum _{k=0}^{n}} \,
{\displaystyle \frac {\Gamma (b - a - 1 + k)\,\Gamma (a + k + 1
 - n - c)}{\Gamma (k + 1)\,\Gamma ( - c - n + 2 + k)}}  \! 
 \right) \,\Gamma (b) \\
\mathrm{\%1}\,(b + c - 2 + n - a)\,\Gamma (2 + n)\,\mathrm{cos}(
\pi \,b)\,(n - a)\,\mathrm{\%2}\,\pi  \left/ {\vrule 
height0.44em width0em depth0.44em} \right. \!  \! (\Gamma (n - a)
\,\Gamma (a - n + 1)^{2} \\
\Gamma (b - a + n)\,\Gamma (2 + a - c)\,\mathrm{sin}(\pi \,a)\,(b
 - 1)\,\mathrm{cos}(\pi \,a))\mbox{} \\ + {\displaystyle \frac {
\Gamma (a)\,\Gamma (b)\,\mathrm{\%1}\,\Gamma (2 + n)\,(n - a)\,
\pi }{\Gamma (n)\,\Gamma ( - c + 1)\,\Gamma (a - n + 1)\,\Gamma (
b - a + n)\,\Gamma (1 + a)\,\Gamma ( - b - n + 1 + a)\,n\,
\mathrm{sin}(\pi \,a)\,\mathrm{cos}(\pi \,a)}}  \\
\mbox{} + {\displaystyle \frac {\Gamma (1 + a)\,\Gamma (1 + a - c
 - n)\,\Gamma (2 + n)\,\pi ^{2}\,(n - a)\,\mathrm{\%2}}{\Gamma (n
 - a)\,\Gamma (a - n + 1)^{2}\,\Gamma (b - a + n)\,\Gamma (2 + a
 - b)\,\Gamma (2 + a - c)\,\mathrm{sin}(\pi \,(a - b))\,\mathrm{
sin}(\pi \,a)\,a}}  \\
\mbox{} + {\displaystyle \frac {\Gamma (b)\,\mathrm{\%1}\,(b + c
 - 2 + n - a)\,\Gamma (2 + n)\,\mathrm{cos}(\pi \,b)\,(n - a)\,
\pi }{\Gamma (n - a)\,\Gamma (n)\,\Gamma (a - n + 1)^{2}\,\Gamma 
(2 - c)\,\mathrm{sin}(\pi \,a)\,(b - 1)\,n\,\mathrm{cos}(\pi \,a)
}}  \\
\mathrm{\%1} := \Gamma (a + 2 - b - c - n) \\
\mathrm{\%2} :=  - b\,c + c + a\,c + a\,b - 1 - 2\,a + b + a\,n
 - a^{2} }
}
\end{maplelatex}

\begin{maplelatex}
\mapleinline{inert}{2d}{245, ":   ",
(a*n+4-2*c-n-2*b+b*c)*(-2+a)*Sum(GAMMA(-2+b+k)*GAMMA(c-2+k)/GAMMA(k+1)
/GAMMA(c-2-a+k+b),k = 0 ..
n)/GAMMA(c-1+n)*GAMMA(2+n)*GAMMA(b-3+c+n-a)/GAMMA(b-1+n)-(-2+a)*(n+1)/
(b-3+c+n-a)+sin(Pi*a-Pi*c)/sin(Pi*a)*GAMMA(b-2)*GAMMA(2+n)*GAMMA(b-3+c
+n-a)/GAMMA(-2+a)*GAMMA(c-2)/GAMMA(b-1+n)/GAMMA(c-1+n)*GAMMA(1+a-b)*GA
MMA(1+a-c)*(a*n+4-2*c-n-2*b+b*c)*sin(Pi*a-Pi*b)/Pi;}{%
\maplemultiline{
245\mbox{:~~~} (a\,n + 4 - 2\,c - n - 2\,b + b\,c)\,( - 
2 + a)\, \left(  \! {\displaystyle \sum _{k=0}^{n}} \,
{\displaystyle \frac {\Gamma ( - 2 + b + k)\,\Gamma (c - 2 + k)}{
\Gamma (k + 1)\,\Gamma (c - 2 - a + k + b)}}  \!  \right)  \\
\Gamma (2 + n)\,\Gamma (b - 3 + c + n - a)/(\Gamma (c - 1 + n)\,
\Gamma (b - 1 + n))\mbox{} - {\displaystyle \frac {( - 2 + a)\,(n
 + 1)}{b - 3 + c + n - a}}  +  \\
\mathrm{sin}(\pi \,a - \pi \,c)\,\Gamma (b - 2)\,\Gamma (2 + n)\,
\Gamma (b - 3 + c + n - a)\,\Gamma (c - 2)\,\Gamma (1 + a - b)
 \\
\Gamma (1 + a - c)\,(a\,n + 4 - 2\,c - n - 2\,b + b\,c)\,\mathrm{
sin}(\pi \,a - \pi \,b)/(\mathrm{sin}(\pi \,a)\,\Gamma ( - 2 + a)
 \\
\Gamma (b - 1 + n)\,\Gamma (c - 1 + n)\,\pi ) }
}
\end{maplelatex}

\begin{maplelatex}
\mapleinline{inert}{2d}{246, ":   ",
-Sum(GAMMA(a+k+1-n-c)*GAMMA(-b-n+1+a+k)/GAMMA(k+1)/GAMMA(a+1+k-n),k =
0 ..
n)*GAMMA(a)*GAMMA(b)*GAMMA(c-a+n)*GAMMA(3-c)*GAMMA(3-b+a-c)*(n-a)*(-b*
c+c+a*c+a*b-1-2*a+b+a*n-a^2)*Pi/GAMMA(n+1-a)/GAMMA(-c+1)/GAMMA(a-n+1)/
GAMMA(b-a+n)/GAMMA(-b-n+1+a)/GAMMA(2+a-b)/GAMMA(2+a-c)/(b+c-2+n-a)/sin
(Pi*a)/cos(Pi*a)+Sum(GAMMA(b-a-1+k)*GAMMA(a+k+1-n-c)/GAMMA(k+1)/GAMMA(
-c-n+2+k),k = 0 ..
n)*GAMMA(a-b)*GAMMA(b)*GAMMA(c-a+n)*GAMMA(3-c)*GAMMA(3-b+a-c)*cos(Pi*b
)*(1+a-b)*(a-b)*(n-a)*(-b*c+c+a*c+a*b-1-2*a+b+a*n-a^2)*Pi/GAMMA(n-a)/G
AMMA(n+1-a)/GAMMA(a-n+1)^2/GAMMA(b-a+n)/GAMMA(2+a-b)/GAMMA(2+a-c)/sin(
Pi*a)/(b-1)/cos(Pi*a)+sin(Pi*a)*(-b*c+c+a*c+a*b-1-2*a+b+a*n-a^2)*Pi/si
n(Pi*a-Pi*b)/sin(Pi*c-Pi*a+Pi*n)/GAMMA(2+a-b)/GAMMA(a-c+3-b-n)*GAMMA(3
-c)/GAMMA(2+a-c)/GAMMA(b-a+n)*GAMMA(3-b+a-c)*GAMMA(a)+(-(c-1)*(c-2)*(b
-1)/(b+c-2+n-a)/Pi/a*sin(Pi*a-Pi*b)+cos(Pi*b)*(c-2)/Pi*sin(Pi*a))*GAMM
A(c-a+n)/GAMMA(n+1)*GAMMA(b-1)/cos(Pi*a)*GAMMA(3-b+a-c);}{%
\maplemultiline{
246\mbox{:~~~}  -  \left(  \! {\displaystyle \sum _{k=0}
^{n}} \,{\displaystyle \frac {\Gamma (a + k + 1 - n - c)\,\Gamma 
( - b - n + 1 + a + k)}{\Gamma (k + 1)\,\Gamma (a + 1 + k - n)}} 
 \!  \right) \,\Gamma (a)\,\Gamma (b)\,\Gamma (c - a + n) \\
\Gamma (3 - c)\,\Gamma (3 - b + a - c)\,(n - a)\,\mathrm{\%1}\,
\pi /(\Gamma (n + 1 - a)\,\Gamma ( - c + 1)\,\Gamma (a - n + 1)
 \\
\Gamma (b - a + n)\,\Gamma ( - b - n + 1 + a)\,\Gamma (2 + a - b)
\,\Gamma (2 + a - c)\,(b + c - 2 + n - a) \\
\mathrm{sin}(\pi \,a)\,\mathrm{cos}(\pi \,a))\mbox{} +  \left( 
 \! {\displaystyle \sum _{k=0}^{n}} \,{\displaystyle \frac {
\Gamma (b - a - 1 + k)\,\Gamma (a + k + 1 - n - c)}{\Gamma (k + 1
)\,\Gamma ( - c - n + 2 + k)}}  \!  \right) \,\Gamma (a - b)\,
\Gamma (b) \\
\Gamma (c - a + n)\,\Gamma (3 - c)\,\Gamma (3 - b + a - c)\,
\mathrm{cos}(\pi \,b)\,(1 + a - b)\,(a - b)\,(n - a)\,\mathrm{\%1
}\,\pi  \left/ {\vrule height0.44em width0em depth0.44em}
 \right. \!  \!  \\
(\Gamma (n - a)\,\Gamma (n + 1 - a)\,\Gamma (a - n + 1)^{2}\,
\Gamma (b - a + n)\,\Gamma (2 + a - b)\,\Gamma (2 + a - c)\,
\mathrm{sin}(\pi \,a) \\
(b - 1)\,\mathrm{cos}(\pi \,a))\mbox{} + \mathrm{sin}(\pi \,a)\,
\mathrm{\%1}\,\pi \,\Gamma (3 - c)\,\Gamma (3 - b + a - c)\,
\Gamma (a)/(\mathrm{sin}(\pi \,a - \pi \,b) \\
\mathrm{sin}(\pi \,c - \pi \,a + \pi \,n)\,\Gamma (2 + a - b)\,
\Gamma (a - c + 3 - b - n)\,\Gamma (2 + a - c)\,\Gamma (b - a + n
))\mbox{} +  \\
( - {\displaystyle \frac {(c - 1)\,(c - 2)\,(b - 1)\,\mathrm{sin}
(\pi \,a - \pi \,b)}{(b + c - 2 + n - a)\,\pi \,a}}  + 
{\displaystyle \frac {\mathrm{cos}(\pi \,b)\,(c - 2)\,\mathrm{sin
}(\pi \,a)}{\pi }} )\,\Gamma (c - a + n) \\
\Gamma (b - 1)\,\Gamma (3 - b + a - c)/(\Gamma (n + 1)\,\mathrm{
cos}(\pi \,a)) \\
\mathrm{\%1} :=  - b\,c + c + a\,c + a\,b - 1 - 2\,a + b + a\,n
 - a^{2} }
}
\end{maplelatex}

\begin{maplelatex}
\mapleinline{inert}{2d}{247, ":   ",
-(-n*c+b*n+n+a*n+a+a*b-a*c-2*c+1-b*c+b+c^2)*cos(Pi*(b-1))*GAMMA(c)*GAM
MA(-b+1)*Sum(GAMMA(1+a-c+k)*GAMMA(c-1-n-b+k)/GAMMA(k+1)/GAMMA(a+1+k-n)
,k = 0 ..
n)/GAMMA(c-b)/Pi*sin(Pi*(b-c+a))*GAMMA(n+2-c+b+a)/GAMMA(2-c+a+n)/b/cos
(Pi*(b-c+a))+(-n*c+b*n+n+a*n+a+a*b-a*c-2*c+1-b*c+b+c^2)*GAMMA(c-1-b-a)
*GAMMA(c-1-n)*GAMMA(c)*GAMMA(-b+1)*Sum(GAMMA(c-1-n-b+k)*GAMMA(-1-n+c-a
+k)/GAMMA(k+1)/GAMMA(-b+c-a-n+k),k = 0 ..
n)/GAMMA(a)/GAMMA(c-a)/GAMMA(c-b)/GAMMA(c-n)/Pi*GAMMA(n+2-c+b+a)*sin(P
i*(c-a))/cos(Pi*(b-c+a))-cos(Pi*(b-1))*GAMMA(c)/GAMMA(1+a)/GAMMA(n)/Pi
*sin(Pi*(b-c+a))/b/n/cos(Pi*(b-c+a))*GAMMA(n+2-c+b+a)*GAMMA(-b+1)-(-n*
c+b*n+n+a*n+a+a*b-a*c-2*c+1-b*c+b+c^2)*GAMMA(c-b-1-n)*GAMMA(c-1-b-a)*G
AMMA(c)/(n+1-c)/sin(Pi*(-c+a))/GAMMA(2-c+a+n)/GAMMA(c-1-n)*sin(Pi*(b-c
+a))/GAMMA(c-b)/GAMMA(c-a)*GAMMA(n+2-c+b+a)+GAMMA(c-1-n)*GAMMA(c)*GAMM
A(-b+1)/GAMMA(n)/GAMMA(c-n)/(-c+b+1+a)/Pi*sin(Pi*(c-a))/n/cos(Pi*(b-c+
a))*GAMMA(n+2-c+b+a)/GAMMA(a);}{%
\maplemultiline{
247\mbox{:~~~}  - \mathrm{\%2}\,\mathrm{cos}(\pi \,(b - 
1))\,\Gamma (c)\,\Gamma ( - b + 1) \\
 \left(  \! {\displaystyle \sum _{k=0}^{n}} \,{\displaystyle 
\frac {\Gamma (1 + a - c + k)\,\Gamma (c - 1 - n - b + k)}{\Gamma
 (k + 1)\,\Gamma (a + 1 + k - n)}}  \!  \right) \,\mathrm{sin}(
\pi \,(b - c + a))\,\mathrm{\%1}/(\Gamma (c - b)\,\pi  \\
\Gamma (2 - c + a + n)\,b\,\mathrm{cos}(\pi \,(b - c + a)))
\mbox{} + \mathrm{\%2}\,\Gamma (c - 1 - b - a)\,\Gamma (c - 1 - n
)\,\Gamma (c) \\
\Gamma ( - b + 1)\, \left(  \! {\displaystyle \sum _{k=0}^{n}} \,
{\displaystyle \frac {\Gamma (c - 1 - n - b + k)\,\Gamma ( - 1 - 
n + c - a + k)}{\Gamma (k + 1)\,\Gamma ( - b + c - a - n + k)}} 
 \!  \right) \,\mathrm{\%1}\,\mathrm{sin}(\pi \,(c - a))/( \\
\Gamma (a)\,\Gamma (c - a)\,\Gamma (c - b)\,\Gamma (c - n)\,\pi 
\,\mathrm{cos}(\pi \,(b - c + a))) \\
\mbox{} - {\displaystyle \frac {\mathrm{cos}(\pi \,(b - 1))\,
\Gamma (c)\,\mathrm{sin}(\pi \,(b - c + a))\,\mathrm{\%1}\,\Gamma
 ( - b + 1)}{\Gamma (1 + a)\,\Gamma (n)\,\pi \,b\,n\,\mathrm{cos}
(\pi \,(b - c + a))}}  \\
\mbox{} - {\displaystyle \frac {\mathrm{\%2}\,\Gamma (c - b - 1
 - n)\,\Gamma (c - 1 - b - a)\,\Gamma (c)\,\mathrm{sin}(\pi \,(b
 - c + a))\,\mathrm{\%1}}{(n + 1 - c)\,\mathrm{sin}(\pi \,( - c
 + a))\,\Gamma (2 - c + a + n)\,\Gamma (c - 1 - n)\,\Gamma (c - b
)\,\Gamma (c - a)}}  \\
\mbox{} + {\displaystyle \frac {\Gamma (c - 1 - n)\,\Gamma (c)\,
\Gamma ( - b + 1)\,\mathrm{sin}(\pi \,(c - a))\,\mathrm{\%1}}{
\Gamma (n)\,\Gamma (c - n)\,( - c + b + 1 + a)\,\pi \,n\,\mathrm{
cos}(\pi \,(b - c + a))\,\Gamma (a)}}  \\
\mathrm{\%1} := \Gamma (n + 2 - c + b + a) \\
\mathrm{\%2} :=  - n\,c + b\,n + n + a\,n + a + a\,b - a\,c - 2\,
c + 1 - b\,c + b + c^{2} }
}
\end{maplelatex}

\begin{maplelatex}
\mapleinline{inert}{2d}{248, ":   ",
(n*c-2*n-a+2*c-1+b*c+a*c-b-a*b-c^2)*(-c+b+a+n+1)*Sum(GAMMA(-1-n+c-a+k)
*GAMMA(c-1-n-b+k)/GAMMA(k+1)/GAMMA(c-n-1+k),k = 0 ..
n)*GAMMA(n)/GAMMA(c-b)/GAMMA(c-a)*GAMMA(c)*cos(Pi*(b-c+a))*sin(Pi*b)/c
os(Pi*b)/sin(Pi*(b-c+a))+(c-1)*(c-2)*sin(Pi*(-c+a))*Sum(GAMMA(c-1-n-b+
k)*GAMMA(1+a-c+k)/GAMMA(k+1)/GAMMA(-b+2+k-n),k = 0 ..
n)*(-c+b+a+n+1)*(-c+b+a+n)*(n*c-2*n-a+2*c-1+b*c+a*c-b-a*b-c^2)/cos(Pi*
b)/(a-1)/sin(Pi*(b-c+a))/GAMMA(2-c+a+n)*GAMMA(n)/GAMMA(c-b)*GAMMA(-b+1
)-(n*c-2*n-a+2*c-1+b*c+a*c-b-a*b-c^2)*GAMMA(-b+1)/GAMMA(a)*GAMMA(c-b-1
-n)/GAMMA(2-c+a+n)*GAMMA(n+2-c+b+a)*GAMMA(n)/GAMMA(c-b)/GAMMA(c-a)*GAM
MA(c)*sin(Pi*b)/((-1)^n)/sin(Pi*(-c+a))+(-c+b+a+n+1)*(c-1)/n*cos(Pi*(b
-c+a))*sin(Pi*b)/cos(Pi*b)/sin(Pi*(b-c+a))+(c-1)*(c-2)*sin(Pi*(-c+a))*
(-c+b+a+n+1)*(-c+b+a+n)/n/cos(Pi*b)/(b-1)/(a-1)/sin(Pi*(b-c+a));}{%
\maplemultiline{
248\mbox{:~~~} \mathrm{\%1}\,( - c + b + a + n + 1)\,
 \left(  \! {\displaystyle \sum _{k=0}^{n}} \,{\displaystyle 
\frac {\Gamma ( - 1 - n + c - a + k)\,\Gamma (c - 1 - n - b + k)
}{\Gamma (k + 1)\,\Gamma (c - n - 1 + k)}}  \!  \right) \,\Gamma 
(n) \\
\Gamma (c)\,\mathrm{cos}(\pi \,(b - c + a))\,\mathrm{sin}(\pi \,b
)/(\Gamma (c - b)\,\Gamma (c - a)\,\mathrm{cos}(\pi \,b)\,
\mathrm{sin}(\pi \,(b - c + a)))\mbox{} +  \\
(c - 1)\,(c - 2)\,\mathrm{sin}(\pi \,( - c + a))\, \left(  \! 
{\displaystyle \sum _{k=0}^{n}} \,{\displaystyle \frac {\Gamma (c
 - 1 - n - b + k)\,\Gamma (1 + a - c + k)}{\Gamma (k + 1)\,\Gamma
 ( - b + 2 + k - n)}}  \!  \right)  \\
( - c + b + a + n + 1)\,( - c + b + a + n)\,\mathrm{\%1}\,\Gamma 
(n)\,\Gamma ( - b + 1)/(\mathrm{cos}(\pi \,b)\,(a - 1) \\
\mathrm{sin}(\pi \,(b - c + a))\,\Gamma (2 - c + a + n)\,\Gamma (
c - b)) \\
\mbox{} - {\displaystyle \frac {\mathrm{\%1}\,\Gamma ( - b + 1)\,
\Gamma (c - b - 1 - n)\,\Gamma (n + 2 - c + b + a)\,\Gamma (n)\,
\Gamma (c)\,\mathrm{sin}(\pi \,b)}{\Gamma (a)\,\Gamma (2 - c + a
 + n)\,\Gamma (c - b)\,\Gamma (c - a)\,(-1)^{n}\,\mathrm{sin}(\pi
 \,( - c + a))}}  \\
\mbox{} + {\displaystyle \frac {( - c + b + a + n + 1)\,(c - 1)\,
\mathrm{cos}(\pi \,(b - c + a))\,\mathrm{sin}(\pi \,b)}{n\,
\mathrm{cos}(\pi \,b)\,\mathrm{sin}(\pi \,(b - c + a))}}  \\
\mbox{} + {\displaystyle \frac {(c - 1)\,(c - 2)\,\mathrm{sin}(
\pi \,( - c + a))\,( - c + b + a + n + 1)\,( - c + b + a + n)}{n
\,\mathrm{cos}(\pi \,b)\,(b - 1)\,(a - 1)\,\mathrm{sin}(\pi \,(b
 - c + a))}}  \\
\mathrm{\%1} := n\,c - 2\,n - a + 2\,c - 1 + b\,c + a\,c - b - a
\,b - c^{2} }
}
\end{maplelatex}

\begin{maplelatex}
\mapleinline{inert}{2d}{253, ":   ",
-2/3*GAMMA(a-1/2*n)*GAMMA(a-1/3-1/2*n)*GAMMA(2*a-1/3)*GAMMA(2/3)*GAMMA
(-a-1/2*n+5/6)*V7(a-1/3-1/2*n,n+1)*GAMMA(1/6+a+1/2*n)/GAMMA(a+1/2*n-5/
6)/GAMMA(a-1/2*n+5/6)/GAMMA(-n+1/2)/GAMMA(1/2+a-1/2*n)/GAMMA(1/2+n)/GA
MMA(3*a-1/2*n)/Pi^(1/2)*GAMMA(2*a)+2/3*GAMMA(a-1/2*n)*GAMMA(a-1/3-1/2*
n)*GAMMA(2*a-1/3)*GAMMA(2/3)*V8(a,n)/GAMMA(a+1/2*n-5/6)/GAMMA(-a+5/6+1
/2*n)/GAMMA(2*a+5/6-n)/GAMMA(2*a+1/6-n);}{%
\maplemultiline{
253\mbox{:~~~}  - {\displaystyle \frac {2}{3}} \Gamma (a
 - {\displaystyle \frac {n}{2}} )\,\Gamma (a - {\displaystyle 
\frac {1}{3}}  - {\displaystyle \frac {n}{2}} )\,\Gamma (2\,a - 
{\displaystyle \frac {1}{3}} )\,\Gamma ({\displaystyle \frac {2}{
3}} )\,\Gamma ( - a - {\displaystyle \frac {n}{2}}  + 
{\displaystyle \frac {5}{6}} )\,\mathrm{V7}(a - {\displaystyle 
\frac {1}{3}}  - {\displaystyle \frac {n}{2}} , \,n + 1) \\
\Gamma ({\displaystyle \frac {1}{6}}  + a + {\displaystyle 
\frac {n}{2}} )\,\Gamma (2\,a) \left/ {\vrule 
height0.80em width0em depth0.80em} \right. \!  \! (\Gamma (a + 
{\displaystyle \frac {n}{2}}  - {\displaystyle \frac {5}{6}} )\,
\Gamma (a - {\displaystyle \frac {n}{2}}  + {\displaystyle 
\frac {5}{6}} )\,\Gamma ( - n + {\displaystyle \frac {1}{2}} )\,
\Gamma ({\displaystyle \frac {1}{2}}  + a - {\displaystyle 
\frac {n}{2}} ) \\
\Gamma ({\displaystyle \frac {1}{2}}  + n)\,\Gamma (3\,a - 
{\displaystyle \frac {n}{2}} )\,\sqrt{\pi })
\mbox{} + {\displaystyle \frac {2}{3}} \,{\displaystyle \frac {
\Gamma (a - {\displaystyle \frac {n}{2}} )\,\Gamma (a - 
{\displaystyle \frac {1}{3}}  - {\displaystyle \frac {n}{2}} )\,
\Gamma (2\,a - {\displaystyle \frac {1}{3}} )\,\Gamma (
{\displaystyle \frac {2}{3}} )\,\mathrm{V8}(a, \,n)}{\Gamma (a + 
{\displaystyle \frac {n}{2}}  - {\displaystyle \frac {5}{6}} )\,
\Gamma ( - a + {\displaystyle \frac {5}{6}}  + {\displaystyle 
\frac {n}{2}} )\,\Gamma (2\,a + {\displaystyle \frac {5}{6}}  - n
)\,\Gamma (2\,a + {\displaystyle \frac {1}{6}}  - n)}}  }
}
\end{maplelatex}

\begin{maplelatex}
\mapleinline{inert}{2d}{254, ":   ",
1/2*GAMMA(2*a-1/3)*GAMMA(a+1/3-1/2*n)*GAMMA(a-1/3-1/2*n)*GAMMA(a-1/2*n
)*GAMMA(4/3+a+1/2*n)*V7(a-1/3-1/2*n,n+1)/GAMMA(1/2+n)/GAMMA(a-1/2*n+5/
6)/GAMMA(a+1/6-1/2*n)/GAMMA(1/2+a-1/2*n)/GAMMA(3*a-1/2*n)-1/2*GAMMA(2*
a-1/3)*GAMMA(a+1/3-1/2*n)*GAMMA(a-1/3-1/2*n)*GAMMA(a-1/2*n)*GAMMA(4/3+
a+1/2*n)*V8(a,n)/GAMMA(1/2+n)/Pi^(1/2)/GAMMA(2*a+5/6-n)/GAMMA(2*a+1/6-
n)/GAMMA(2*a);}{%
\maplemultiline{
254\mbox{:~~~} {\displaystyle \frac {1}{2}} 
{\displaystyle \frac {\Gamma (2\,a - {\displaystyle \frac {1}{3}
} )\,\Gamma (a + {\displaystyle \frac {1}{3}}  - {\displaystyle 
\frac {n}{2}} )\,\Gamma (a - {\displaystyle \frac {1}{3}}  - 
{\displaystyle \frac {n}{2}} )\,\Gamma (a - {\displaystyle 
\frac {n}{2}} )\,\Gamma ({\displaystyle \frac {4}{3}}  + a + 
{\displaystyle \frac {n}{2}} )\,\mathrm{V7}(a - {\displaystyle 
\frac {1}{3}}  - {\displaystyle \frac {n}{2}} , \,n + 1)}{\Gamma 
({\displaystyle \frac {1}{2}}  + n)\,\Gamma (a - {\displaystyle 
\frac {n}{2}}  + {\displaystyle \frac {5}{6}} )\,\Gamma (a + 
{\displaystyle \frac {1}{6}}  - {\displaystyle \frac {n}{2}} )\,
\Gamma ({\displaystyle \frac {1}{2}}  + a - {\displaystyle 
\frac {n}{2}} )\,\Gamma (3\,a - {\displaystyle \frac {n}{2}} )}} 
 \\
\mbox{} - {\displaystyle \frac {1}{2}} \,{\displaystyle \frac {
\Gamma (2\,a - {\displaystyle \frac {1}{3}} )\,\Gamma (a + 
{\displaystyle \frac {1}{3}}  - {\displaystyle \frac {n}{2}} )\,
\Gamma (a - {\displaystyle \frac {1}{3}}  - {\displaystyle 
\frac {n}{2}} )\,\Gamma (a - {\displaystyle \frac {n}{2}} )\,
\Gamma ({\displaystyle \frac {4}{3}}  + a + {\displaystyle 
\frac {n}{2}} )\,\mathrm{V8}(a, \,n)}{\Gamma ({\displaystyle 
\frac {1}{2}}  + n)\,\sqrt{\pi }\,\Gamma (2\,a + {\displaystyle 
\frac {5}{6}}  - n)\,\Gamma (2\,a + {\displaystyle \frac {1}{6}} 
 - n)\,\Gamma (2\,a)}}  }
}
\end{maplelatex}

\begin{maplelatex}
\mapleinline{inert}{2d}{255, ":   ",
1/2*GAMMA(2*a+1/3)*GAMMA(a-1/3-1/2*n)*GAMMA(a+1/3-1/2*n)*GAMMA(a-1/2*n
)*GAMMA(2/3+a+1/2*n)*V7(a-1/3-1/2*n,n+1)/GAMMA(1/2+n)/GAMMA(a+1/6-1/2*
n)/GAMMA(a-1/2*n+5/6)/GAMMA(1/2+a-1/2*n)/GAMMA(3*a-1/2*n)-1/2*GAMMA(2*
a+1/3)*GAMMA(a-1/3-1/2*n)*GAMMA(a+1/3-1/2*n)*GAMMA(a-1/2*n)*GAMMA(2/3+
a+1/2*n)*V8(a,n)/GAMMA(1/2+n)/Pi^(1/2)/GAMMA(2*a+1/6-n)/GAMMA(2*a+5/6-
n)/GAMMA(2*a);}{%
\maplemultiline{
255\mbox{:~~~} {\displaystyle \frac {1}{2}}  
{\displaystyle \frac {\Gamma (2\,a + {\displaystyle \frac {1}{3}
} )\,\Gamma (a - {\displaystyle \frac {1}{3}}  - {\displaystyle 
\frac {n}{2}} )\,\Gamma (a + {\displaystyle \frac {1}{3}}  - 
{\displaystyle \frac {n}{2}} )\,\Gamma (a - {\displaystyle 
\frac {n}{2}} )\,\Gamma ({\displaystyle \frac {2}{3}}  + a + 
{\displaystyle \frac {n}{2}} )\,\mathrm{V7}(a - {\displaystyle 
\frac {1}{3}}  - {\displaystyle \frac {n}{2}} , \,n + 1)}{\Gamma 
({\displaystyle \frac {1}{2}}  + n)\,\Gamma (a + {\displaystyle 
\frac {1}{6}}  - {\displaystyle \frac {n}{2}} )\,\Gamma (a - 
{\displaystyle \frac {n}{2}}  + {\displaystyle \frac {5}{6}} )\,
\Gamma ({\displaystyle \frac {1}{2}}  + a - {\displaystyle 
\frac {n}{2}} )\,\Gamma (3\,a - {\displaystyle \frac {n}{2}} )}} 
 \\
\mbox{} - {\displaystyle \frac {1}{2}} \,{\displaystyle \frac {
\Gamma (2\,a + {\displaystyle \frac {1}{3}} )\,\Gamma (a - 
{\displaystyle \frac {1}{3}}  - {\displaystyle \frac {n}{2}} )\,
\Gamma (a + {\displaystyle \frac {1}{3}}  - {\displaystyle 
\frac {n}{2}} )\,\Gamma (a - {\displaystyle \frac {n}{2}} )\,
\Gamma ({\displaystyle \frac {2}{3}}  + a + {\displaystyle 
\frac {n}{2}} )\,\mathrm{V8}(a, \,n)}{\Gamma ({\displaystyle 
\frac {1}{2}}  + n)\,\sqrt{\pi }\,\Gamma (2\,a + {\displaystyle 
\frac {1}{6}}  - n)\,\Gamma (2\,a + {\displaystyle \frac {5}{6}} 
 - n)\,\Gamma (2\,a)}}  }
}
\end{maplelatex}

\begin{maplelatex}
\mapleinline{inert}{2d}{256, ":   ",
-2/9*GAMMA(2*a+1/3)*GAMMA(a-1/3-1/2*n)*GAMMA(2*a-1/3)*GAMMA(a-1/2*n)*P
i^(1/2)*3^(1/2)/GAMMA(2/3)*GAMMA(-a-1/2*n+1/2)*V7(a-1/3-1/2*n,n+1)*GAM
MA(a+1/2+1/2*n)/GAMMA(a-1/2+1/2*n)/GAMMA(a-1/2*n+5/6)/GAMMA(a+1/6-1/2*
n)/GAMMA(-n+1/2)/GAMMA(1/2+n)/GAMMA(3*a-1/2*n)+2/9*GAMMA(2*a+1/3)*GAMM
A(a-1/3-1/2*n)*GAMMA(2*a-1/3)*GAMMA(a-1/2*n)*Pi*3^(1/2)/GAMMA(2/3)*V8(
a,n)/GAMMA(a-1/2+1/2*n)/GAMMA(2*a+5/6-n)/GAMMA(1/2-a+1/2*n)/GAMMA(2*a+
1/6-n)/GAMMA(2*a);}{%
\maplemultiline{
256\mbox{:~~~}  - {\displaystyle \frac {2}{9}} \Gamma (2
\,a + {\displaystyle \frac {1}{3}} )\,\Gamma (a - {\displaystyle 
\frac {1}{3}}  - {\displaystyle \frac {n}{2}} )\,\Gamma (2\,a - 
{\displaystyle \frac {1}{3}} )\,\Gamma (a - {\displaystyle 
\frac {n}{2}} )\,\sqrt{\pi }\,\sqrt{3}\,\Gamma ( - a - 
{\displaystyle \frac {n}{2}}  + {\displaystyle \frac {1}{2}} )
 \\
\mathrm{V7}(a - {\displaystyle \frac {1}{3}}  - {\displaystyle 
\frac {n}{2}} , \,n + 1)\,\Gamma (a + {\displaystyle \frac {1}{2}
}  + {\displaystyle \frac {n}{2}} ) \left/ {\vrule 
height0.80em width0em depth0.80em} \right. \!  \! (\Gamma (
{\displaystyle \frac {2}{3}} )\,\Gamma (a - {\displaystyle 
\frac {1}{2}}  + {\displaystyle \frac {n}{2}} )\,\Gamma (a - 
{\displaystyle \frac {n}{2}}  + {\displaystyle \frac {5}{6}} )
 \\
\Gamma (a + {\displaystyle \frac {1}{6}}  - {\displaystyle 
\frac {n}{2}} )\,\Gamma ( - n + {\displaystyle \frac {1}{2}} )\,
\Gamma ({\displaystyle \frac {1}{2}}  + n)\,\Gamma (3\,a - 
{\displaystyle \frac {n}{2}} )) \\
\mbox{} + {\displaystyle \frac {2}{9}} \,{\displaystyle \frac {
\Gamma (2\,a + {\displaystyle \frac {1}{3}} )\,\Gamma (a - 
{\displaystyle \frac {1}{3}}  - {\displaystyle \frac {n}{2}} )\,
\Gamma (2\,a - {\displaystyle \frac {1}{3}} )\,\Gamma (a - 
{\displaystyle \frac {n}{2}} )\,\pi \,\sqrt{3}\,\mathrm{V8}(a, \,
n)}{\Gamma ({\displaystyle \frac {2}{3}} )\,\Gamma (a - 
{\displaystyle \frac {1}{2}}  + {\displaystyle \frac {n}{2}} )\,
\Gamma (2\,a + {\displaystyle \frac {5}{6}}  - n)\,\Gamma (
{\displaystyle \frac {1}{2}}  - a + {\displaystyle \frac {n}{2}} 
)\,\Gamma (2\,a + {\displaystyle \frac {1}{6}}  - n)\,\Gamma (2\,
a)}}  }
}
\end{maplelatex}

\begin{maplelatex}
\mapleinline{inert}{2d}{257, ":   ",
-GAMMA(2*a-1/3)*GAMMA(a+1/3-1/2*n)*GAMMA(2*a+1/3)*GAMMA(a-1/2*n)*GAMMA
(2/3)*GAMMA(-a-1/2*n+1/2)*V7(a-1/3-1/2*n,n+1)*GAMMA(a+1/2+1/2*n)/GAMMA
(a-1/2+1/2*n)/GAMMA(a+1/6-1/2*n)/GAMMA(a-1/2*n+5/6)/GAMMA(-n+1/2)/GAMM
A(1/2+n)/GAMMA(3*a-1/2*n)/Pi^(1/2)+GAMMA(2*a-1/3)*GAMMA(a+1/3-1/2*n)*G
AMMA(2*a+1/3)*GAMMA(a-1/2*n)*GAMMA(2/3)*V8(a,n)/GAMMA(a-1/2+1/2*n)/GAM
MA(2*a+1/6-n)/GAMMA(1/2-a+1/2*n)/GAMMA(2*a+5/6-n)/GAMMA(2*a);}{%
\maplemultiline{
257\mbox{:~~~}  - \Gamma (2\,a - {\displaystyle \frac {1
}{3}} )\,\Gamma (a + {\displaystyle \frac {1}{3}}  - 
{\displaystyle \frac {n}{2}} )\,\Gamma (2\,a + {\displaystyle 
\frac {1}{3}} )\,\Gamma (a - {\displaystyle \frac {n}{2}} )\,
\Gamma ({\displaystyle \frac {2}{3}} )\,\Gamma ( - a - 
{\displaystyle \frac {n}{2}}  + {\displaystyle \frac {1}{2}} )
 \\
\mathrm{V7}(a - {\displaystyle \frac {1}{3}}  - {\displaystyle 
\frac {n}{2}} , \,n + 1)\,\Gamma (a + {\displaystyle \frac {1}{2}
}  + {\displaystyle \frac {n}{2}} ) \left/ {\vrule 
height0.80em width0em depth0.80em} \right. \!  \! (\Gamma (a - 
{\displaystyle \frac {1}{2}}  + {\displaystyle \frac {n}{2}} )\,
\Gamma (a + {\displaystyle \frac {1}{6}}  - {\displaystyle 
\frac {n}{2}} )\,\Gamma (a - {\displaystyle \frac {n}{2}}  + 
{\displaystyle \frac {5}{6}} ) \\
\Gamma ( - n + {\displaystyle \frac {1}{2}} )\,\Gamma (
{\displaystyle \frac {1}{2}}  + n)\,\Gamma (3\,a - 
{\displaystyle \frac {n}{2}} )\,\sqrt{\pi })
\mbox{} + {\displaystyle \frac {\Gamma (2\,a - {\displaystyle 
\frac {1}{3}} )\,\Gamma (a + {\displaystyle \frac {1}{3}}  - 
{\displaystyle \frac {n}{2}} )\,\Gamma (2\,a + {\displaystyle 
\frac {1}{3}} )\,\Gamma (a - {\displaystyle \frac {n}{2}} )\,
\Gamma ({\displaystyle \frac {2}{3}} )\,\mathrm{V8}(a, \,n)}{
\Gamma (a - {\displaystyle \frac {1}{2}}  + {\displaystyle 
\frac {n}{2}} )\,\Gamma (2\,a + {\displaystyle \frac {1}{6}}  - n
)\,\Gamma ({\displaystyle \frac {1}{2}}  - a + {\displaystyle 
\frac {n}{2}} )\,\Gamma (2\,a + {\displaystyle \frac {5}{6}}  - n
)\,\Gamma (2\,a)}}  }
}
\end{maplelatex}

\begin{maplelatex}
\mapleinline{inert}{2d}{258, ":   ",
-GAMMA(2*a+1/3)*GAMMA(2*a-1/3)*GAMMA(a-1/2*n)*GAMMA(1-1/2*n-a)*GAMMA(-
a-1/2*n+1/2)*V7(a-1/3-1/2*n,n+1)*GAMMA(a+1/2+1/2*n)/GAMMA(a-1/2+1/2*n)
/GAMMA(-n+1/2)/GAMMA(a-1/2*n+5/6)/GAMMA(a+1/6-1/2*n)/Pi^(1/2)/GAMMA(1/
2+n)/GAMMA(3*a-1/2*n)*GAMMA(2*a)+GAMMA(2*a+1/3)*GAMMA(2*a-1/3)*GAMMA(a
-1/2*n)*GAMMA(1-1/2*n-a)*V8(a,n)/GAMMA(a-1/2+1/2*n)/GAMMA(2*a+5/6-n)/G
AMMA(2*a+1/6-n)/GAMMA(1/2-a+1/2*n);}{%
\maplemultiline{
258\mbox{:~~~}  - \Gamma (2\,a + {\displaystyle \frac {1
}{3}} )\,\Gamma (2\,a - {\displaystyle \frac {1}{3}} )\,\Gamma (a
 - {\displaystyle \frac {n}{2}} )\,\Gamma (1 - {\displaystyle 
\frac {n}{2}}  - a)\,\Gamma ( - a - {\displaystyle \frac {n}{2}} 
 + {\displaystyle \frac {1}{2}} ) \\
\mathrm{V7}(a - {\displaystyle \frac {1}{3}}  - {\displaystyle 
\frac {n}{2}} , \,n + 1)\,\Gamma (a + {\displaystyle \frac {1}{2}
}  + {\displaystyle \frac {n}{2}} )\,\Gamma (2\,a) \left/ 
{\vrule height0.80em width0em depth0.80em} \right. \!  \! (\Gamma
 (a - {\displaystyle \frac {1}{2}}  + {\displaystyle \frac {n}{2}
} )\,\Gamma ( - n + {\displaystyle \frac {1}{2}} ) \\
\Gamma (a - {\displaystyle \frac {n}{2}}  + {\displaystyle 
\frac {5}{6}} )\,\Gamma (a + {\displaystyle \frac {1}{6}}  - 
{\displaystyle \frac {n}{2}} )\,\sqrt{\pi }\,\Gamma (
{\displaystyle \frac {1}{2}}  + n)\,\Gamma (3\,a - 
{\displaystyle \frac {n}{2}} ))
\mbox{} + {\displaystyle \frac {\Gamma (2\,a + {\displaystyle 
\frac {1}{3}} )\,\Gamma (2\,a - {\displaystyle \frac {1}{3}} )\,
\Gamma (a - {\displaystyle \frac {n}{2}} )\,\Gamma (1 - 
{\displaystyle \frac {n}{2}}  - a)\,\mathrm{V8}(a, \,n)}{\Gamma (
a - {\displaystyle \frac {1}{2}}  + {\displaystyle \frac {n}{2}} 
)\,\Gamma (2\,a + {\displaystyle \frac {5}{6}}  - n)\,\Gamma (2\,
a + {\displaystyle \frac {1}{6}}  - n)\,\Gamma ({\displaystyle 
\frac {1}{2}}  - a + {\displaystyle \frac {n}{2}} )}}  }
}
\end{maplelatex}

\begin{maplelatex}
\mapleinline{inert}{2d}{259, ":   ",
GAMMA(n+1)/GAMMA(a)/GAMMA(c)*GAMMA(b)*GAMMA(e)/GAMMA(e-a)/GAMMA(e-c)*S
um(1/GAMMA(1+j)/GAMMA(n-j+1)*GAMMA(a+j)*GAMMA(c+j)/GAMMA(b+j)*GAMMA(e-
j-a-c),j = 0 .. n);}{%
\[
259\mbox{:~~~} \,{\displaystyle \frac {\Gamma (n + 1)\,
\Gamma (b)\,\Gamma (e)\, \left(  \! {\displaystyle \sum _{j=0}^{n
}} \,{\displaystyle \frac {\Gamma (a + j)\,\Gamma (c + j)\,\Gamma
 (e - j - a - c)}{\Gamma (1 + j)\,\Gamma (n - j + 1)\,\Gamma (b
 + j)}}  \!  \right) }{\Gamma (a)\,\Gamma (c)\,\Gamma (e - a)\,
\Gamma (e - c)}} 
\]
}
\end{maplelatex}

\begin{maplelatex}
\mapleinline{inert}{2d}{260, ":   ",
Sum(GAMMA(a+k)/GAMMA(a)*GAMMA(b+k)/GAMMA(b)/GAMMA(c+k)*GAMMA(c)/GAMMA(
e+k)*GAMMA(e)/GAMMA(k+1)*(-1)^k/GAMMA(n-k+1)*GAMMA(n+1),k = 0 ..
n);}{%
\[
260\mbox{:~~~} \,{\displaystyle \sum _{k=0}^{n}} \,
{\displaystyle \frac {\Gamma (a + k)\,\Gamma (b + k)\,\Gamma (c)
\,\Gamma (e)\,(-1)^{k}\,\Gamma (n + 1)}{\Gamma (a)\,\Gamma (b)\,
\Gamma (c + k)\,\Gamma (e + k)\,\Gamma (k + 1)\,\Gamma (n - k + 1
)}} 
\]
}
\end{maplelatex}

\end{maplegroup}

%% file: AppendixB261to290.tex
\begin{maplegroup}
\mapleresult
\begin{maplelatex}
\mapleinline{inert}{2d}{267, ":        ",
-GAMMA(2*a)*GAMMA(a-1/3-1/2*n)*GAMMA(2*a-1/3)*GAMMA(-a-1/2*n+5/6)*V7(a
-1/3-1/2*n,n+1)*GAMMA(2*a+5/6-n)/GAMMA(a+1/2*n-5/6)/GAMMA(a-1/2*n+5/6)
/GAMMA(-n+1/2)/GAMMA(1/2+a-1/2*n)/GAMMA(1/2+n)/GAMMA(3*a-1/2*n)/Pi^(1/
2)*GAMMA(4/3-a+1/2*n)*GAMMA(4/3-a-1/2*n)/GAMMA(-a-1/2*n+11/6)*GAMMA(1/
6+a+1/2*n)+GAMMA(a-1/3-1/2*n)*GAMMA(2*a-1/3)*V8(a,n)/GAMMA(a+1/2*n-5/6
)/GAMMA(-a+5/6+1/2*n)/GAMMA(2*a+1/6-n)*GAMMA(4/3-a+1/2*n)*GAMMA(4/3-a-
1/2*n)/GAMMA(-a-1/2*n+11/6);}{%
\maplemultiline{
267\mbox{:     ~~~}    - \Gamma (2\,a)\,\Gamma (a - 
{\displaystyle \frac {1}{3}}  - {\displaystyle \frac {n}{2}} )\,
\Gamma (2\,a - {\displaystyle \frac {1}{3}} )\,\Gamma ( - a - 
{\displaystyle \frac {n}{2}}  + {\displaystyle \frac {5}{6}} )\,
\mathrm{V7}(a - {\displaystyle \frac {1}{3}}  - {\displaystyle 
\frac {n}{2}} , \,n + 1) \\
\Gamma (2\,a + {\displaystyle \frac {5}{6}}  - n)\,\Gamma (
{\displaystyle \frac {4}{3}}  - a + {\displaystyle \frac {n}{2}} 
)\,\Gamma ({\displaystyle \frac {4}{3}}  - a - {\displaystyle 
\frac {n}{2}} )\,\Gamma ({\displaystyle \frac {1}{6}}  + a + 
{\displaystyle \frac {n}{2}} ) \left/ {\vrule 
height0.80em width0em depth0.80em} \right. \!  \! (\Gamma (a + 
{\displaystyle \frac {n}{2}}  - {\displaystyle \frac {5}{6}} )
 \\
\Gamma (a - {\displaystyle \frac {n}{2}}  + {\displaystyle 
\frac {5}{6}} )\,\Gamma ( - n + {\displaystyle \frac {1}{2}} )\,
\Gamma ({\displaystyle \frac {1}{2}}  + a - {\displaystyle 
\frac {n}{2}} )\,\Gamma ({\displaystyle \frac {1}{2}}  + n)\,
\Gamma (3\,a - {\displaystyle \frac {n}{2}} )\,\sqrt{\pi }\,
\Gamma ( - a - {\displaystyle \frac {n}{2}}  + {\displaystyle 
\frac {11}{6}} )) \\
\mbox{} + {\displaystyle \frac {\Gamma (a - {\displaystyle 
\frac {1}{3}}  - {\displaystyle \frac {n}{2}} )\,\Gamma (2\,a - 
{\displaystyle \frac {1}{3}} )\,\mathrm{V8}(a, \,n)\,\Gamma (
{\displaystyle \frac {4}{3}}  - a + {\displaystyle \frac {n}{2}} 
)\,\Gamma ({\displaystyle \frac {4}{3}}  - a - {\displaystyle 
\frac {n}{2}} )}{\Gamma (a + {\displaystyle \frac {n}{2}}  - 
{\displaystyle \frac {5}{6}} )\,\Gamma ( - a + {\displaystyle 
\frac {5}{6}}  + {\displaystyle \frac {n}{2}} )\,\Gamma (2\,a + 
{\displaystyle \frac {1}{6}}  - n)\,\Gamma ( - a - 
{\displaystyle \frac {n}{2}}  + {\displaystyle \frac {11}{6}} )}
}  }
}
\end{maplelatex}

\begin{maplelatex}
\mapleinline{inert}{2d}{268, ":        ",
-2/9*GAMMA(2*a)*GAMMA(a-1/3-1/2*n)*GAMMA(2*a-1/3)/GAMMA(2/3)*GAMMA(-a-
1/2*n+5/6)*V7(a-1/3-1/2*n,n+1)/GAMMA(a+1/2*n-5/6)/GAMMA(-n+1/2)/GAMMA(
1/2+a-1/2*n)/GAMMA(1/2+n)/GAMMA(3*a-1/2*n)*Pi^(1/2)*GAMMA(4/3-a-1/2*n)
*3^(1/2)/GAMMA(-a-1/2*n+11/6)*GAMMA(1/6+a+1/2*n)+2/9*GAMMA(a-1/3-1/2*n
)*GAMMA(2*a-1/3)/GAMMA(2/3)*V8(a,n)/GAMMA(a+1/2*n-5/6)/GAMMA(-a+5/6+1/
2*n)/GAMMA(2*a+5/6-n)/GAMMA(2*a+1/6-n)*GAMMA(a-1/2*n+5/6)*GAMMA(4/3-a-
1/2*n)*Pi*3^(1/2)/GAMMA(-a-1/2*n+11/6);}{%
\maplemultiline{
268\mbox{:     ~~~}    - {\displaystyle \frac {2}{9}} \Gamma (2
\,a)\,\Gamma (a - {\displaystyle \frac {1}{3}}  - {\displaystyle 
\frac {n}{2}} )\,\Gamma (2\,a - {\displaystyle \frac {1}{3}} )\,
\Gamma ( - a - {\displaystyle \frac {n}{2}}  + {\displaystyle 
\frac {5}{6}} )\,\mathrm{V7}(a - {\displaystyle \frac {1}{3}}  - 
{\displaystyle \frac {n}{2}} , \,n + 1)\,\sqrt{\pi } \\
\Gamma ({\displaystyle \frac {4}{3}}  - a - {\displaystyle 
\frac {n}{2}} )\,\sqrt{3}\,\Gamma ({\displaystyle \frac {1}{6}} 
 + a + {\displaystyle \frac {n}{2}} ) \left/ {\vrule 
height0.80em width0em depth0.80em} \right. \!  \! (\Gamma (
{\displaystyle \frac {2}{3}} )\,\Gamma (a + {\displaystyle 
\frac {n}{2}}  - {\displaystyle \frac {5}{6}} )\,\Gamma ( - n + 
{\displaystyle \frac {1}{2}} )\,\Gamma ({\displaystyle \frac {1}{
2}}  + a - {\displaystyle \frac {n}{2}} ) \\
\Gamma ({\displaystyle \frac {1}{2}}  + n)\,\Gamma (3\,a - 
{\displaystyle \frac {n}{2}} )\,\Gamma ( - a - {\displaystyle 
\frac {n}{2}}  + {\displaystyle \frac {11}{6}} ))\mbox{} +  \\
{\displaystyle \frac {2}{9}} \,{\displaystyle \frac {\Gamma (a - 
{\displaystyle \frac {1}{3}}  - {\displaystyle \frac {n}{2}} )\,
\Gamma (2\,a - {\displaystyle \frac {1}{3}} )\,\mathrm{V8}(a, \,n
)\,\Gamma (a - {\displaystyle \frac {n}{2}}  + {\displaystyle 
\frac {5}{6}} )\,\Gamma ({\displaystyle \frac {4}{3}}  - a - 
{\displaystyle \frac {n}{2}} )\,\pi \,\sqrt{3}}{\Gamma (
{\displaystyle \frac {2}{3}} )\,\Gamma (a + {\displaystyle 
\frac {n}{2}}  - {\displaystyle \frac {5}{6}} )\,\Gamma ( - a + 
{\displaystyle \frac {5}{6}}  + {\displaystyle \frac {n}{2}} )\,
\Gamma (2\,a + {\displaystyle \frac {5}{6}}  - n)\,\Gamma (2\,a
 + {\displaystyle \frac {1}{6}}  - n)\,\Gamma ( - a - 
{\displaystyle \frac {n}{2}}  + {\displaystyle \frac {11}{6}} )}
}  }
}
\end{maplelatex}

\begin{maplelatex}
\mapleinline{inert}{2d}{269, ":        ",
-GAMMA(a-1/3-1/2*n)*GAMMA(2*a-1/3)*GAMMA(-a-1/2*n+5/6)*V7(a-1/3-1/2*n,
n+1)/GAMMA(a+1/2*n-5/6)/GAMMA(a-1/2*n+5/6)/GAMMA(-n+1/2)/GAMMA(1/2+a-1
/2*n)/GAMMA(1/2+n)/GAMMA(3*a-1/2*n)/Pi^(1/2)*GAMMA(4/3-a-1/2*n)*GAMMA(
1/6+a+1/2*n)*GAMMA(2*a+1/3)*GAMMA(2*a)+GAMMA(a-1/3-1/2*n)*GAMMA(2*a-1/
3)*V8(a,n)/GAMMA(a+1/2*n-5/6)/GAMMA(-a+5/6+1/2*n)/GAMMA(2*a+5/6-n)/GAM
MA(2*a+1/6-n)*GAMMA(4/3-a-1/2*n)*GAMMA(2*a+1/3);}{%
\maplemultiline{
269\mbox{:     ~~~}    - \Gamma (a - {\displaystyle \frac {1}{3
}}  - {\displaystyle \frac {n}{2}} )\,\Gamma (2\,a - 
{\displaystyle \frac {1}{3}} )\,\Gamma ( - a - {\displaystyle 
\frac {n}{2}}  + {\displaystyle \frac {5}{6}} )\,\mathrm{V7}(a - 
{\displaystyle \frac {1}{3}}  - {\displaystyle \frac {n}{2}} , \,
n + 1)\,\Gamma ({\displaystyle \frac {4}{3}}  - a - 
{\displaystyle \frac {n}{2}} ) \\
\Gamma ({\displaystyle \frac {1}{6}}  + a + {\displaystyle 
\frac {n}{2}} )\,\Gamma (2\,a + {\displaystyle \frac {1}{3}} )\,
\Gamma (2\,a) \left/ {\vrule height0.80em width0em depth0.80em}
 \right. \!  \! (\Gamma (a + {\displaystyle \frac {n}{2}}  - 
{\displaystyle \frac {5}{6}} )\,\Gamma (a - {\displaystyle 
\frac {n}{2}}  + {\displaystyle \frac {5}{6}} )\,\Gamma ( - n + 
{\displaystyle \frac {1}{2}} ) \\
\Gamma ({\displaystyle \frac {1}{2}}  + a - {\displaystyle 
\frac {n}{2}} )\,\Gamma ({\displaystyle \frac {1}{2}}  + n)\,
\Gamma (3\,a - {\displaystyle \frac {n}{2}} )\,\sqrt{\pi }) \\
\mbox{} + {\displaystyle \frac {\Gamma (a - {\displaystyle 
\frac {1}{3}}  - {\displaystyle \frac {n}{2}} )\,\Gamma (2\,a - 
{\displaystyle \frac {1}{3}} )\,\mathrm{V8}(a, \,n)\,\Gamma (
{\displaystyle \frac {4}{3}}  - a - {\displaystyle \frac {n}{2}} 
)\,\Gamma (2\,a + {\displaystyle \frac {1}{3}} )}{\Gamma (a + 
{\displaystyle \frac {n}{2}}  - {\displaystyle \frac {5}{6}} )\,
\Gamma ( - a + {\displaystyle \frac {5}{6}}  + {\displaystyle 
\frac {n}{2}} )\,\Gamma (2\,a + {\displaystyle \frac {5}{6}}  - n
)\,\Gamma (2\,a + {\displaystyle \frac {1}{6}}  - n)}}  }
}
\end{maplelatex}

\begin{maplelatex}
\mapleinline{inert}{2d}{270, ":        ",
-2/3*GAMMA(2*a)*GAMMA(a-1/3-1/2*n)*GAMMA(2*a-1/3)*GAMMA(2/3)*GAMMA(-a-
1/2*n+5/6)*V7(a-1/3-1/2*n,n+1)/GAMMA(a+1/2*n-5/6)/GAMMA(a-1/2*n+5/6)/G
AMMA(-n+1/2)/GAMMA(1/2+n)/GAMMA(3*a-1/2*n)/Pi^(1/2)*GAMMA(4/3-a-1/2*n)
/GAMMA(-a-1/2*n+11/6)*GAMMA(1/6+a+1/2*n)+2/3*GAMMA(a-1/3-1/2*n)*GAMMA(
2*a-1/3)*GAMMA(2/3)*V8(a,n)/GAMMA(a+1/2*n-5/6)/GAMMA(-a+5/6+1/2*n)/GAM
MA(2*a+5/6-n)/GAMMA(2*a+1/6-n)*GAMMA(1/2+a-1/2*n)*GAMMA(4/3-a-1/2*n)/G
AMMA(-a-1/2*n+11/6);}{%
\maplemultiline{
270\mbox{:     ~~~}    - {\displaystyle \frac {2}{3}} \Gamma (2
\,a)\,\Gamma (a - {\displaystyle \frac {1}{3}}  - {\displaystyle 
\frac {n}{2}} )\,\Gamma (2\,a - {\displaystyle \frac {1}{3}} )\,
\Gamma ({\displaystyle \frac {2}{3}} )\,\Gamma ( - a - 
{\displaystyle \frac {n}{2}}  + {\displaystyle \frac {5}{6}} )\,
\mathrm{V7}(a - {\displaystyle \frac {1}{3}}  - {\displaystyle 
\frac {n}{2}} , \,n + 1) \\
\Gamma ({\displaystyle \frac {4}{3}}  - a - {\displaystyle 
\frac {n}{2}} )\,\Gamma ({\displaystyle \frac {1}{6}}  + a + 
{\displaystyle \frac {n}{2}} ) \left/ {\vrule 
height0.80em width0em depth0.80em} \right. \!  \! (\Gamma (a + 
{\displaystyle \frac {n}{2}}  - {\displaystyle \frac {5}{6}} )\,
\Gamma (a - {\displaystyle \frac {n}{2}}  + {\displaystyle 
\frac {5}{6}} )\,\Gamma ( - n + {\displaystyle \frac {1}{2}} )\,
\Gamma ({\displaystyle \frac {1}{2}}  + n) \\
\Gamma (3\,a - {\displaystyle \frac {n}{2}} )\,\sqrt{\pi }\,
\Gamma ( - a - {\displaystyle \frac {n}{2}}  + {\displaystyle 
\frac {11}{6}} )) \\
\mbox{} + {\displaystyle \frac {2}{3}} \,{\displaystyle \frac {
\Gamma (a - {\displaystyle \frac {1}{3}}  - {\displaystyle 
\frac {n}{2}} )\,\Gamma (2\,a - {\displaystyle \frac {1}{3}} )\,
\Gamma ({\displaystyle \frac {2}{3}} )\,\mathrm{V8}(a, \,n)\,
\Gamma ({\displaystyle \frac {1}{2}}  + a - {\displaystyle 
\frac {n}{2}} )\,\Gamma ({\displaystyle \frac {4}{3}}  - a - 
{\displaystyle \frac {n}{2}} )}{\Gamma (a + {\displaystyle 
\frac {n}{2}}  - {\displaystyle \frac {5}{6}} )\,\Gamma ( - a + 
{\displaystyle \frac {5}{6}}  + {\displaystyle \frac {n}{2}} )\,
\Gamma (2\,a + {\displaystyle \frac {5}{6}}  - n)\,\Gamma (2\,a
 + {\displaystyle \frac {1}{6}}  - n)\,\Gamma ( - a - 
{\displaystyle \frac {n}{2}}  + {\displaystyle \frac {11}{6}} )}
}  }
}
\end{maplelatex}

\begin{maplelatex}
\mapleinline{inert}{2d}{271, ":        ",
-2/9*GAMMA(2*a)*GAMMA(a-1/3-1/2*n)*GAMMA(2*a-1/3)/GAMMA(2/3)*GAMMA(-a-
1/2*n+5/6)*V7(a-1/3-1/2*n,n+1)*GAMMA(2*a+5/6-n)/GAMMA(a+1/2*n-5/6)/GAM
MA(a-1/2*n+5/6)/GAMMA(-n+1/2)/GAMMA(1/2+a-1/2*n)/GAMMA(1/2+n)/GAMMA(3*
a-1/2*n)*Pi^(1/2)*GAMMA(4/3-2*a)*3^(1/2)/GAMMA(-a-1/2*n+11/6)*GAMMA(1/
6+a+1/2*n)+2/9*GAMMA(a-1/3-1/2*n)*GAMMA(2*a-1/3)/GAMMA(2/3)*V8(a,n)/GA
MMA(a+1/2*n-5/6)/GAMMA(-a+5/6+1/2*n)/GAMMA(2*a+1/6-n)*GAMMA(4/3-2*a)*P
i*3^(1/2)/GAMMA(-a-1/2*n+11/6);}{%
\maplemultiline{
271\mbox{:     ~~~}    - {\displaystyle \frac {2}{9}} \Gamma (2
\,a)\,\Gamma (a - {\displaystyle \frac {1}{3}}  - {\displaystyle 
\frac {n}{2}} )\,\Gamma (2\,a - {\displaystyle \frac {1}{3}} )\,
\Gamma ( - a - {\displaystyle \frac {n}{2}}  + {\displaystyle 
\frac {5}{6}} )\,\mathrm{V7}(a - {\displaystyle \frac {1}{3}}  - 
{\displaystyle \frac {n}{2}} , \,n + 1) \\
\Gamma (2\,a + {\displaystyle \frac {5}{6}}  - n)\,\sqrt{\pi }\,
\Gamma ({\displaystyle \frac {4}{3}}  - 2\,a)\,\sqrt{3}\,\Gamma (
{\displaystyle \frac {1}{6}}  + a + {\displaystyle \frac {n}{2}} 
) \left/ {\vrule height0.80em width0em depth0.80em} \right. \! 
 \! (\Gamma ({\displaystyle \frac {2}{3}} )\,\Gamma (a + 
{\displaystyle \frac {n}{2}}  - {\displaystyle \frac {5}{6}} )
 \\
\Gamma (a - {\displaystyle \frac {n}{2}}  + {\displaystyle 
\frac {5}{6}} )\,\Gamma ( - n + {\displaystyle \frac {1}{2}} )\,
\Gamma ({\displaystyle \frac {1}{2}}  + a - {\displaystyle 
\frac {n}{2}} )\,\Gamma ({\displaystyle \frac {1}{2}}  + n)\,
\Gamma (3\,a - {\displaystyle \frac {n}{2}} )\,\Gamma ( - a - 
{\displaystyle \frac {n}{2}}  + {\displaystyle \frac {11}{6}} ))
 \\
\mbox{} + {\displaystyle \frac {2}{9}} \,{\displaystyle \frac {
\Gamma (a - {\displaystyle \frac {1}{3}}  - {\displaystyle 
\frac {n}{2}} )\,\Gamma (2\,a - {\displaystyle \frac {1}{3}} )\,
\mathrm{V8}(a, \,n)\,\Gamma ({\displaystyle \frac {4}{3}}  - 2\,a
)\,\pi \,\sqrt{3}}{\Gamma ({\displaystyle \frac {2}{3}} )\,\Gamma
 (a + {\displaystyle \frac {n}{2}}  - {\displaystyle \frac {5}{6}
} )\,\Gamma ( - a + {\displaystyle \frac {5}{6}}  + 
{\displaystyle \frac {n}{2}} )\,\Gamma (2\,a + {\displaystyle 
\frac {1}{6}}  - n)\,\Gamma ( - a - {\displaystyle \frac {n}{2}} 
 + {\displaystyle \frac {11}{6}} )}}  }
}
\end{maplelatex}

\begin{maplelatex}
\mapleinline{inert}{2d}{272, ":        ",
-GAMMA(2*a)*GAMMA(a-1/3-1/2*n)*GAMMA(2*a-1/3)*GAMMA(-a-1/2*n+5/6)*V7(a
-1/3-1/2*n,n+1)*GAMMA(2*a+5/6-n)/GAMMA(a+1/2*n-5/6)/GAMMA(a-1/2*n+5/6)
/GAMMA(-n+1/2)/GAMMA(1/2+a-1/2*n)/GAMMA(1/2+n)/GAMMA(3*a-1/2*n)/Pi^(1/
2)*GAMMA(-a+5/6+1/2*n)*GAMMA(1/6+a+1/2*n)+GAMMA(a-1/3-1/2*n)*GAMMA(2*a
-1/3)*V8(a,n)/GAMMA(a+1/2*n-5/6)/GAMMA(2*a+1/6-n);}{%
\maplemultiline{
272\mbox{:     ~~~}    - \Gamma (2\,a)\,\Gamma (a - 
{\displaystyle \frac {1}{3}}  - {\displaystyle \frac {n}{2}} )\,
\Gamma (2\,a - {\displaystyle \frac {1}{3}} )\,\Gamma ( - a - 
{\displaystyle \frac {n}{2}}  + {\displaystyle \frac {5}{6}} )\,
\mathrm{V7}(a - {\displaystyle \frac {1}{3}}  - {\displaystyle 
\frac {n}{2}} , \,n + 1) \\
\Gamma (2\,a + {\displaystyle \frac {5}{6}}  - n)\,\Gamma ( - a
 + {\displaystyle \frac {5}{6}}  + {\displaystyle \frac {n}{2}} )
\,\Gamma ({\displaystyle \frac {1}{6}}  + a + {\displaystyle 
\frac {n}{2}} ) \left/ {\vrule height0.80em width0em depth0.80em}
 \right. \!  \! (\Gamma (a + {\displaystyle \frac {n}{2}}  - 
{\displaystyle \frac {5}{6}} )\,\Gamma (a - {\displaystyle 
\frac {n}{2}}  + {\displaystyle \frac {5}{6}} ) \\
\Gamma ( - n + {\displaystyle \frac {1}{2}} )\,\Gamma (
{\displaystyle \frac {1}{2}}  + a - {\displaystyle \frac {n}{2}} 
)\,\Gamma ({\displaystyle \frac {1}{2}}  + n)\,\Gamma (3\,a - 
{\displaystyle \frac {n}{2}} )\,\sqrt{\pi }) 
\mbox{} + {\displaystyle \frac {\Gamma (a - {\displaystyle 
\frac {1}{3}}  - {\displaystyle \frac {n}{2}} )\,\Gamma (2\,a - 
{\displaystyle \frac {1}{3}} )\,\mathrm{V8}(a, \,n)}{\Gamma (a + 
{\displaystyle \frac {n}{2}}  - {\displaystyle \frac {5}{6}} )\,
\Gamma (2\,a + {\displaystyle \frac {1}{6}}  - n)}}  }
}
\end{maplelatex}

\begin{maplelatex}
\mapleinline{inert}{2d}{273, ":        ",
-2/3*GAMMA(2*a)*GAMMA(a-1/3-1/2*n)*GAMMA(2*a-1/3)*GAMMA(2/3)*GAMMA(-a-
1/2*n+5/6)*V7(a-1/3-1/2*n,n+1)*GAMMA(2*a+5/6-n)/GAMMA(a+1/2*n-5/6)/GAM
MA(a-1/2*n+5/6)/GAMMA(-n+1/2)/GAMMA(1/2+a-1/2*n)/GAMMA(1/2+n)/GAMMA(3*
a-1/2*n)/Pi^(1/2)*GAMMA(1-2*a)/GAMMA(-a-1/2*n+11/6)*GAMMA(1/6+a+1/2*n)
+2/3*GAMMA(a-1/3-1/2*n)*GAMMA(2*a-1/3)*GAMMA(2/3)*V8(a,n)/GAMMA(a+1/2*
n-5/6)/GAMMA(-a+5/6+1/2*n)/GAMMA(2*a+1/6-n)*GAMMA(1-2*a)/GAMMA(-a-1/2*
n+11/6);}{%
\maplemultiline{
273\mbox{:     ~~~}    - {\displaystyle \frac {2}{3}} \Gamma (2
\,a)\,\Gamma (a - {\displaystyle \frac {1}{3}}  - {\displaystyle 
\frac {n}{2}} )\,\Gamma (2\,a - {\displaystyle \frac {1}{3}} )\,
\Gamma ({\displaystyle \frac {2}{3}} )\,\Gamma ( - a - 
{\displaystyle \frac {n}{2}}  + {\displaystyle \frac {5}{6}} )\,
\mathrm{V7}(a - {\displaystyle \frac {1}{3}}  - {\displaystyle 
\frac {n}{2}} , \,n + 1) \\
\Gamma (2\,a + {\displaystyle \frac {5}{6}}  - n)\,\Gamma (1 - 2
\,a)\,\Gamma ({\displaystyle \frac {1}{6}}  + a + {\displaystyle 
\frac {n}{2}} ) \left/ {\vrule height0.80em width0em depth0.80em}
 \right. \!  \! (\Gamma (a + {\displaystyle \frac {n}{2}}  - 
{\displaystyle \frac {5}{6}} )\,\Gamma (a - {\displaystyle 
\frac {n}{2}}  + {\displaystyle \frac {5}{6}} ) \\
\Gamma ( - n + {\displaystyle \frac {1}{2}} )\,\Gamma (
{\displaystyle \frac {1}{2}}  + a - {\displaystyle \frac {n}{2}} 
)\,\Gamma ({\displaystyle \frac {1}{2}}  + n)\,\Gamma (3\,a - 
{\displaystyle \frac {n}{2}} )\,\sqrt{\pi }\,\Gamma ( - a - 
{\displaystyle \frac {n}{2}}  + {\displaystyle \frac {11}{6}} ))
 \\
\mbox{} + {\displaystyle \frac {2}{3}} \,{\displaystyle \frac {
\Gamma (a - {\displaystyle \frac {1}{3}}  - {\displaystyle 
\frac {n}{2}} )\,\Gamma (2\,a - {\displaystyle \frac {1}{3}} )\,
\Gamma ({\displaystyle \frac {2}{3}} )\,\mathrm{V8}(a, \,n)\,
\Gamma (1 - 2\,a)}{\Gamma (a + {\displaystyle \frac {n}{2}}  - 
{\displaystyle \frac {5}{6}} )\,\Gamma ( - a + {\displaystyle 
\frac {5}{6}}  + {\displaystyle \frac {n}{2}} )\,\Gamma (2\,a + 
{\displaystyle \frac {1}{6}}  - n)\,\Gamma ( - a - 
{\displaystyle \frac {n}{2}}  + {\displaystyle \frac {11}{6}} )}
}  }
}
\end{maplelatex}

\begin{maplelatex}
\mapleinline{inert}{2d}{274, ":        ",
-2/9*GAMMA(2*a)*GAMMA(a-1/3-1/2*n)*GAMMA(2*a-1/3)/GAMMA(2/3)*GAMMA(-a-
1/2*n+5/6)*V7(a-1/3-1/2*n,n+1)/GAMMA(a+1/2*n-5/6)/GAMMA(a-1/2*n+5/6)/G
AMMA(-n+1/2)/GAMMA(1/2+a-1/2*n)/GAMMA(1/2+n)/GAMMA(3*a-1/2*n)*Pi^(1/2)
*GAMMA(a+1/3-1/2*n)*3^(1/2)*GAMMA(1/6+a+1/2*n)+2/9*GAMMA(a-1/3-1/2*n)*
GAMMA(2*a-1/3)/GAMMA(2/3)*V8(a,n)/GAMMA(a+1/2*n-5/6)/GAMMA(-a+5/6+1/2*
n)/GAMMA(2*a+5/6-n)/GAMMA(2*a+1/6-n)*GAMMA(a+1/3-1/2*n)*Pi*3^(1/2);}{%
\maplemultiline{
274\mbox{:     ~~~}    - {\displaystyle \frac {2}{9}} \Gamma (2
\,a)\,\Gamma (a - {\displaystyle \frac {1}{3}}  - {\displaystyle 
\frac {n}{2}} )\,\Gamma (2\,a - {\displaystyle \frac {1}{3}} )\,
\Gamma ( - a - {\displaystyle \frac {n}{2}}  + {\displaystyle 
\frac {5}{6}} )\,\mathrm{V7}(a - {\displaystyle \frac {1}{3}}  - 
{\displaystyle \frac {n}{2}} , \,n + 1)\,\sqrt{\pi } \\
\Gamma (a + {\displaystyle \frac {1}{3}}  - {\displaystyle 
\frac {n}{2}} )\,\sqrt{3}\,\Gamma ({\displaystyle \frac {1}{6}} 
 + a + {\displaystyle \frac {n}{2}} ) \left/ {\vrule 
height0.80em width0em depth0.80em} \right. \!  \! (\Gamma (
{\displaystyle \frac {2}{3}} )\,\Gamma (a + {\displaystyle 
\frac {n}{2}}  - {\displaystyle \frac {5}{6}} )\,\Gamma (a - 
{\displaystyle \frac {n}{2}}  + {\displaystyle \frac {5}{6}} )\,
\Gamma ( - n + {\displaystyle \frac {1}{2}} ) \\
\Gamma ({\displaystyle \frac {1}{2}}  + a - {\displaystyle 
\frac {n}{2}} )\,\Gamma ({\displaystyle \frac {1}{2}}  + n)\,
\Gamma (3\,a - {\displaystyle \frac {n}{2}} )) \\
\mbox{} + {\displaystyle \frac {2}{9}} \,{\displaystyle \frac {
\Gamma (a - {\displaystyle \frac {1}{3}}  - {\displaystyle 
\frac {n}{2}} )\,\Gamma (2\,a - {\displaystyle \frac {1}{3}} )\,
\mathrm{V8}(a, \,n)\,\Gamma (a + {\displaystyle \frac {1}{3}}  - 
{\displaystyle \frac {n}{2}} )\,\pi \,\sqrt{3}}{\Gamma (
{\displaystyle \frac {2}{3}} )\,\Gamma (a + {\displaystyle 
\frac {n}{2}}  - {\displaystyle \frac {5}{6}} )\,\Gamma ( - a + 
{\displaystyle \frac {5}{6}}  + {\displaystyle \frac {n}{2}} )\,
\Gamma (2\,a + {\displaystyle \frac {5}{6}}  - n)\,\Gamma (2\,a
 + {\displaystyle \frac {1}{6}}  - n)}}  }
}
\end{maplelatex}

\begin{maplelatex}
\mapleinline{inert}{2d}{275, ":        ",
-4/27*GAMMA(2*a)*GAMMA(a-1/3-1/2*n)*GAMMA(2*a-1/3)*GAMMA(-a-1/2*n+5/6)
*V7(a-1/3-1/2*n,n+1)/GAMMA(a+1/2*n-5/6)/GAMMA(a-1/2*n+5/6)/GAMMA(1/2+a
-1/2*n)/GAMMA(1/2+n)/GAMMA(3*a-1/2*n)*Pi^(1/2)*3^(1/2)/GAMMA(-a-1/2*n+
11/6)*GAMMA(1/6+a+1/2*n)+4/27*GAMMA(a-1/3-1/2*n)*GAMMA(2*a-1/3)*V8(a,n
)/GAMMA(a+1/2*n-5/6)/GAMMA(-a+5/6+1/2*n)/GAMMA(2*a+5/6-n)/GAMMA(2*a+1/
6-n)*GAMMA(-n+1/2)*Pi*3^(1/2)/GAMMA(-a-1/2*n+11/6);}{%
\maplemultiline{
275\mbox{:     ~~~}    - {\displaystyle \frac {4}{27}} \,
{\displaystyle \frac {\Gamma (2\,a)\,\Gamma (a - {\displaystyle 
\frac {1}{3}}  - {\displaystyle \frac {n}{2}} )\,\Gamma (2\,a - 
{\displaystyle \frac {1}{3}} )\,\Gamma ( - a - {\displaystyle 
\frac {n}{2}}  + {\displaystyle \frac {5}{6}} )\,\mathrm{V7}(a - 
{\displaystyle \frac {1}{3}}  - {\displaystyle \frac {n}{2}} , \,
n + 1)\,\sqrt{\pi }\,\sqrt{3}\,\Gamma ({\displaystyle \frac {1}{6
}}  + a + {\displaystyle \frac {n}{2}} )}{\Gamma (a + 
{\displaystyle \frac {n}{2}}  - {\displaystyle \frac {5}{6}} )\,
\Gamma (a - {\displaystyle \frac {n}{2}}  + {\displaystyle 
\frac {5}{6}} )\,\Gamma ({\displaystyle \frac {1}{2}}  + a - 
{\displaystyle \frac {n}{2}} )\,\Gamma ({\displaystyle \frac {1}{
2}}  + n)\,\Gamma (3\,a - {\displaystyle \frac {n}{2}} )\,\Gamma 
( - a - {\displaystyle \frac {n}{2}}  + {\displaystyle \frac {11
}{6}} )}}  \\
\mbox{} + {\displaystyle \frac {4}{27}} \,{\displaystyle \frac {
\Gamma (a - {\displaystyle \frac {1}{3}}  - {\displaystyle 
\frac {n}{2}} )\,\Gamma (2\,a - {\displaystyle \frac {1}{3}} )\,
\mathrm{V8}(a, \,n)\,\Gamma ( - n + {\displaystyle \frac {1}{2}} 
)\,\pi \,\sqrt{3}}{\Gamma (a + {\displaystyle \frac {n}{2}}  - 
{\displaystyle \frac {5}{6}} )\,\Gamma ( - a + {\displaystyle 
\frac {5}{6}}  + {\displaystyle \frac {n}{2}} )\,\Gamma (2\,a + 
{\displaystyle \frac {5}{6}}  - n)\,\Gamma (2\,a + 
{\displaystyle \frac {1}{6}}  - n)\,\Gamma ( - a - 
{\displaystyle \frac {n}{2}}  + {\displaystyle \frac {11}{6}} )}
}  }
}
\end{maplelatex}

\begin{maplelatex}
\mapleinline{inert}{2d}{276, ":        ",
-GAMMA(a+1/3-1/2*n)*GAMMA(a-1/3-1/2*n)*GAMMA(a-1/2*n)*V8(a,n)/Pi/GAMMA
(2*a+5/6-n)/GAMMA(1/2+n)/GAMMA(2*a)/GAMMA(2*a+1/6-n)*GAMMA(4/3-a+1/2*n
)*GAMMA(2/3+a+1/2*n)*GAMMA(3*a+1/2-1/2*n)+GAMMA(a+1/3-1/2*n)*GAMMA(a-1
/3-1/2*n)*GAMMA(a-1/2*n)/Pi^(1/2)*V7(a-1/3-1/2*n,n+1)/GAMMA(a-1/2*n+5/
6)/GAMMA(a+1/6-1/2*n)/GAMMA(1/2+a-1/2*n)/GAMMA(1/2+n)*GAMMA(4/3-a+1/2*
n)*GAMMA(2/3+a+1/2*n)*GAMMA(3*a+1/2-1/2*n)/GAMMA(3*a-1/2*n);}{%
\maplemultiline{
276\mbox{:     ~~~}    - {\displaystyle \frac {\Gamma (a + 
{\displaystyle \frac {1}{3}}  - {\displaystyle \frac {n}{2}} )\,
\Gamma (a - {\displaystyle \frac {1}{3}}  - {\displaystyle 
\frac {n}{2}} )\,\Gamma (a - {\displaystyle \frac {n}{2}} )\,
\mathrm{V8}(a, \,n)\,\Gamma ({\displaystyle \frac {4}{3}}  - a + 
{\displaystyle \frac {n}{2}} )\,\Gamma ({\displaystyle \frac {2}{
3}}  + a + {\displaystyle \frac {n}{2}} )\,\Gamma (3\,a + 
{\displaystyle \frac {1}{2}}  - {\displaystyle \frac {n}{2}} )}{
\pi \,\Gamma (2\,a + {\displaystyle \frac {5}{6}}  - n)\,\Gamma (
{\displaystyle \frac {1}{2}}  + n)\,\Gamma (2\,a)\,\Gamma (2\,a
 + {\displaystyle \frac {1}{6}}  - n)}}  \\
\mbox{} + \Gamma (a + {\displaystyle \frac {1}{3}}  - 
{\displaystyle \frac {n}{2}} )\,\Gamma (a - {\displaystyle 
\frac {1}{3}}  - {\displaystyle \frac {n}{2}} )\,\Gamma (a - 
{\displaystyle \frac {n}{2}} )\,\mathrm{V7}(a - {\displaystyle 
\frac {1}{3}}  - {\displaystyle \frac {n}{2}} , \,n + 1)\,\Gamma 
({\displaystyle \frac {4}{3}}  - a + {\displaystyle \frac {n}{2}
} ) \\
\Gamma ({\displaystyle \frac {2}{3}}  + a + {\displaystyle 
\frac {n}{2}} )\,\Gamma (3\,a + {\displaystyle \frac {1}{2}}  - 
{\displaystyle \frac {n}{2}} ) \left/ {\vrule 
height0.80em width0em depth0.80em} \right. \!  \! (\sqrt{\pi }\,
\Gamma (a - {\displaystyle \frac {n}{2}}  + {\displaystyle 
\frac {5}{6}} )\,\Gamma (a + {\displaystyle \frac {1}{6}}  - 
{\displaystyle \frac {n}{2}} )\,\Gamma ({\displaystyle \frac {1}{
2}}  + a - {\displaystyle \frac {n}{2}} )
\Gamma ({\displaystyle \frac {1}{2}}  + n)\,\Gamma (3\,a - 
{\displaystyle \frac {n}{2}} )) }
}
\end{maplelatex}

\begin{maplelatex}
\mapleinline{inert}{2d}{277, ":        ",
-GAMMA(a+1/3-1/2*n)*GAMMA(a-1/3-1/2*n)*GAMMA(a-1/2*n)*V8(a,n)/Pi/GAMMA
(2*a+5/6-n)/GAMMA(1/2+n)/GAMMA(2*a)/GAMMA(2*a+1/6-n)*GAMMA(a-1/2*n+5/6
)*GAMMA(2/3+a+1/2*n)*GAMMA(1+a+1/2*n)+GAMMA(a+1/3-1/2*n)*GAMMA(a-1/3-1
/2*n)*GAMMA(a-1/2*n)/Pi^(1/2)*V7(a-1/3-1/2*n,n+1)/GAMMA(a+1/6-1/2*n)/G
AMMA(1/2+a-1/2*n)/GAMMA(1/2+n)*GAMMA(2/3+a+1/2*n)*GAMMA(1+a+1/2*n)/GAM
MA(3*a-1/2*n);}{%
\maplemultiline{
277\mbox{:     ~~~}    - {\displaystyle \frac {\Gamma (a + 
{\displaystyle \frac {1}{3}}  - {\displaystyle \frac {n}{2}} )\,
\Gamma (a - {\displaystyle \frac {1}{3}}  - {\displaystyle 
\frac {n}{2}} )\,\Gamma (a - {\displaystyle \frac {n}{2}} )\,
\mathrm{V8}(a, \,n)\,\Gamma (a - {\displaystyle \frac {n}{2}}  + 
{\displaystyle \frac {5}{6}} )\,\Gamma ({\displaystyle \frac {2}{
3}}  + a + {\displaystyle \frac {n}{2}} )\,\Gamma (1 + a + 
{\displaystyle \frac {n}{2}} )}{\pi \,\Gamma (2\,a + 
{\displaystyle \frac {5}{6}}  - n)\,\Gamma ({\displaystyle 
\frac {1}{2}}  + n)\,\Gamma (2\,a)\,\Gamma (2\,a + 
{\displaystyle \frac {1}{6}}  - n)}}  \\
\mbox{} + {\displaystyle \frac {\Gamma (a + {\displaystyle 
\frac {1}{3}}  - {\displaystyle \frac {n}{2}} )\,\Gamma (a - 
{\displaystyle \frac {1}{3}}  - {\displaystyle \frac {n}{2}} )\,
\Gamma (a - {\displaystyle \frac {n}{2}} )\,\mathrm{V7}(a - 
{\displaystyle \frac {1}{3}}  - {\displaystyle \frac {n}{2}} , \,
n + 1)\,\Gamma ({\displaystyle \frac {2}{3}}  + a + 
{\displaystyle \frac {n}{2}} )\,\Gamma (1 + a + {\displaystyle 
\frac {n}{2}} )}{\sqrt{\pi }\,\Gamma (a + {\displaystyle \frac {1
}{6}}  - {\displaystyle \frac {n}{2}} )\,\Gamma ({\displaystyle 
\frac {1}{2}}  + a - {\displaystyle \frac {n}{2}} )\,\Gamma (
{\displaystyle \frac {1}{2}}  + n)\,\Gamma (3\,a - 
{\displaystyle \frac {n}{2}} )}}  }
}
\end{maplelatex}

\begin{maplelatex}
\mapleinline{inert}{2d}{278, ":        ",
-GAMMA(a+1/3-1/2*n)*GAMMA(a-1/3-1/2*n)*GAMMA(a-1/2*n)*GAMMA(4/3+a+1/2*
n)*V8(a,n)/Pi/GAMMA(2*a+5/6-n)/GAMMA(1/2+n)/GAMMA(2*a)/GAMMA(2*a+1/6-n
)*GAMMA(1/2+a-1/2*n)*GAMMA(2/3+a+1/2*n)+GAMMA(a+1/3-1/2*n)*GAMMA(a-1/3
-1/2*n)*GAMMA(a-1/2*n)*GAMMA(4/3+a+1/2*n)/Pi^(1/2)*V7(a-1/3-1/2*n,n+1)
/GAMMA(a-1/2*n+5/6)/GAMMA(a+1/6-1/2*n)/GAMMA(1/2+n)*GAMMA(2/3+a+1/2*n)
/GAMMA(3*a-1/2*n);}{%
\maplemultiline{
278\mbox{:     ~~~}    - {\displaystyle \frac {\Gamma (a + 
{\displaystyle \frac {1}{3}}  - {\displaystyle \frac {n}{2}} )\,
\Gamma (a - {\displaystyle \frac {1}{3}}  - {\displaystyle 
\frac {n}{2}} )\,\Gamma (a - {\displaystyle \frac {n}{2}} )\,
\Gamma ({\displaystyle \frac {4}{3}}  + a + {\displaystyle 
\frac {n}{2}} )\,\mathrm{V8}(a, \,n)\,\Gamma ({\displaystyle 
\frac {1}{2}}  + a - {\displaystyle \frac {n}{2}} )\,\Gamma (
{\displaystyle \frac {2}{3}}  + a + {\displaystyle \frac {n}{2}} 
)}{\pi \,\Gamma (2\,a + {\displaystyle \frac {5}{6}}  - n)\,
\Gamma ({\displaystyle \frac {1}{2}}  + n)\,\Gamma (2\,a)\,\Gamma
 (2\,a + {\displaystyle \frac {1}{6}}  - n)}}  \\
\mbox{} + {\displaystyle \frac {\Gamma (a + {\displaystyle 
\frac {1}{3}}  - {\displaystyle \frac {n}{2}} )\,\Gamma (a - 
{\displaystyle \frac {1}{3}}  - {\displaystyle \frac {n}{2}} )\,
\Gamma (a - {\displaystyle \frac {n}{2}} )\,\Gamma (
{\displaystyle \frac {4}{3}}  + a + {\displaystyle \frac {n}{2}} 
)\,\mathrm{V7}(a - {\displaystyle \frac {1}{3}}  - 
{\displaystyle \frac {n}{2}} , \,n + 1)\,\Gamma ({\displaystyle 
\frac {2}{3}}  + a + {\displaystyle \frac {n}{2}} )}{\sqrt{\pi }
\,\Gamma (a - {\displaystyle \frac {n}{2}}  + {\displaystyle 
\frac {5}{6}} )\,\Gamma (a + {\displaystyle \frac {1}{6}}  - 
{\displaystyle \frac {n}{2}} )\,\Gamma ({\displaystyle \frac {1}{
2}}  + n)\,\Gamma (3\,a - {\displaystyle \frac {n}{2}} )}}  }
}
\end{maplelatex}

\begin{maplelatex}
\mapleinline{inert}{2d}{279, ":        ",
-GAMMA(a+1/3-1/2*n)*GAMMA(a-1/3-1/2*n)*GAMMA(a-1/2*n)*V8(a,n)/Pi/GAMMA
(2*a+5/6-n)/GAMMA(1/2+n)/GAMMA(2*a)/GAMMA(2*a+1/6-n)*GAMMA(1-a+1/2*n)*
GAMMA(3*a+1/2-1/2*n)*GAMMA(1+a+1/2*n)+GAMMA(a+1/3-1/2*n)*GAMMA(a-1/3-1
/2*n)*GAMMA(a-1/2*n)/Pi^(1/2)*V7(a-1/3-1/2*n,n+1)/GAMMA(a-1/2*n+5/6)/G
AMMA(a+1/6-1/2*n)/GAMMA(1/2+a-1/2*n)/GAMMA(1/2+n)*GAMMA(1-a+1/2*n)*GAM
MA(3*a+1/2-1/2*n)*GAMMA(1+a+1/2*n)/GAMMA(3*a-1/2*n);}{%
\maplemultiline{
279\mbox{:     ~~~}    - {\displaystyle \frac {\Gamma (a + 
{\displaystyle \frac {1}{3}}  - {\displaystyle \frac {n}{2}} )\,
\Gamma (a - {\displaystyle \frac {1}{3}}  - {\displaystyle 
\frac {n}{2}} )\,\Gamma (a - {\displaystyle \frac {n}{2}} )\,
\mathrm{V8}(a, \,n)\,\Gamma (1 - a + {\displaystyle \frac {n}{2}
} )\,\Gamma (3\,a + {\displaystyle \frac {1}{2}}  - 
{\displaystyle \frac {n}{2}} )\,\Gamma (1 + a + {\displaystyle 
\frac {n}{2}} )}{\pi \,\Gamma (2\,a + {\displaystyle \frac {5}{6}
}  - n)\,\Gamma ({\displaystyle \frac {1}{2}}  + n)\,\Gamma (2\,a
)\,\Gamma (2\,a + {\displaystyle \frac {1}{6}}  - n)}}  \\
\mbox{} + \Gamma (a + {\displaystyle \frac {1}{3}}  - 
{\displaystyle \frac {n}{2}} )\,\Gamma (a - {\displaystyle 
\frac {1}{3}}  - {\displaystyle \frac {n}{2}} )\,\Gamma (a - 
{\displaystyle \frac {n}{2}} )\,\mathrm{V7}(a - {\displaystyle 
\frac {1}{3}}  - {\displaystyle \frac {n}{2}} , \,n + 1)\,\Gamma 
(1 - a + {\displaystyle \frac {n}{2}} ) \\
\Gamma (3\,a + {\displaystyle \frac {1}{2}}  - {\displaystyle 
\frac {n}{2}} )\,\Gamma (1 + a + {\displaystyle \frac {n}{2}} )
 \left/ {\vrule height0.80em width0em depth0.80em} \right. \! 
 \! (\sqrt{\pi }\,\Gamma (a - {\displaystyle \frac {n}{2}}  + 
{\displaystyle \frac {5}{6}} )\,\Gamma (a + {\displaystyle 
\frac {1}{6}}  - {\displaystyle \frac {n}{2}} )\,\Gamma (
{\displaystyle \frac {1}{2}}  + a - {\displaystyle \frac {n}{2}} 
)
\Gamma ({\displaystyle \frac {1}{2}}  + n)\,\Gamma (3\,a - 
{\displaystyle \frac {n}{2}} )) }
}
\end{maplelatex}

\begin{maplelatex}
\mapleinline{inert}{2d}{280, ":        ",
-GAMMA(a+1/3-1/2*n)*GAMMA(a-1/3-1/2*n)*GAMMA(a-1/2*n)*GAMMA(4/3+a+1/2*
n)*V8(a,n)/Pi/GAMMA(2*a+5/6-n)/GAMMA(1/2+n)/GAMMA(2*a)/GAMMA(2*a+1/6-n
)*GAMMA(2/3-a+1/2*n)*GAMMA(3*a+1/2-1/2*n)+GAMMA(a+1/3-1/2*n)*GAMMA(a-1
/3-1/2*n)*GAMMA(a-1/2*n)*GAMMA(4/3+a+1/2*n)/Pi^(1/2)*V7(a-1/3-1/2*n,n+
1)/GAMMA(a-1/2*n+5/6)/GAMMA(a+1/6-1/2*n)/GAMMA(1/2+a-1/2*n)/GAMMA(1/2+
n)*GAMMA(2/3-a+1/2*n)*GAMMA(3*a+1/2-1/2*n)/GAMMA(3*a-1/2*n);}{%
\maplemultiline{
280\mbox{:     ~~~}    - {\displaystyle \frac {\Gamma (a + 
{\displaystyle \frac {1}{3}}  - {\displaystyle \frac {n}{2}} )\,
\Gamma (a - {\displaystyle \frac {1}{3}}  - {\displaystyle 
\frac {n}{2}} )\,\Gamma (a - {\displaystyle \frac {n}{2}} )\,
\Gamma ({\displaystyle \frac {4}{3}}  + a + {\displaystyle 
\frac {n}{2}} )\,\mathrm{V8}(a, \,n)\,\Gamma ({\displaystyle 
\frac {2}{3}}  - a + {\displaystyle \frac {n}{2}} )\,\Gamma (3\,a
 + {\displaystyle \frac {1}{2}}  - {\displaystyle \frac {n}{2}} )
}{\pi \,\Gamma (2\,a + {\displaystyle \frac {5}{6}}  - n)\,\Gamma
 ({\displaystyle \frac {1}{2}}  + n)\,\Gamma (2\,a)\,\Gamma (2\,a
 + {\displaystyle \frac {1}{6}}  - n)}}  \\
\mbox{} + \Gamma (a + {\displaystyle \frac {1}{3}}  - 
{\displaystyle \frac {n}{2}} )\,\Gamma (a - {\displaystyle 
\frac {1}{3}}  - {\displaystyle \frac {n}{2}} )\,\Gamma (a - 
{\displaystyle \frac {n}{2}} )\,\Gamma ({\displaystyle \frac {4}{
3}}  + a + {\displaystyle \frac {n}{2}} )\,\mathrm{V7}(a - 
{\displaystyle \frac {1}{3}}  - {\displaystyle \frac {n}{2}} , \,
n + 1) \\
\Gamma ({\displaystyle \frac {2}{3}}  - a + {\displaystyle 
\frac {n}{2}} )\,\Gamma (3\,a + {\displaystyle \frac {1}{2}}  - 
{\displaystyle \frac {n}{2}} ) \left/ {\vrule 
height0.80em width0em depth0.80em} \right. \!  \! (\sqrt{\pi }\,
\Gamma (a - {\displaystyle \frac {n}{2}}  + {\displaystyle 
\frac {5}{6}} )\,\Gamma (a + {\displaystyle \frac {1}{6}}  - 
{\displaystyle \frac {n}{2}} )\,\Gamma ({\displaystyle \frac {1}{
2}}  + a - {\displaystyle \frac {n}{2}} ) 
\Gamma ({\displaystyle \frac {1}{2}}  + n)\,\Gamma (3\,a - 
{\displaystyle \frac {n}{2}} )) }
}
\end{maplelatex}

\begin{maplelatex}
\mapleinline{inert}{2d}{281, ":        ",
-1/2*GAMMA(a+1/3-1/2*n)*GAMMA(a-1/3-1/2*n)*GAMMA(a-1/2*n)*V8(a,n)/Pi^(
1/2)/GAMMA(2*a+5/6-n)/GAMMA(2*a)/GAMMA(2*a+1/6-n)*GAMMA(3*a+1/2-1/2*n)
+1/2*GAMMA(a+1/3-1/2*n)*GAMMA(a-1/3-1/2*n)*GAMMA(a-1/2*n)*V7(a-1/3-1/2
*n,n+1)/GAMMA(a-1/2*n+5/6)/GAMMA(a+1/6-1/2*n)/GAMMA(1/2+a-1/2*n)*GAMMA
(3*a+1/2-1/2*n)/GAMMA(3*a-1/2*n);}{%
\maplemultiline{
281\mbox{:     ~~~}    - {\displaystyle \frac {1}{2}} \,
{\displaystyle \frac {\Gamma (a + {\displaystyle \frac {1}{3}} 
 - {\displaystyle \frac {n}{2}} )\,\Gamma (a - {\displaystyle 
\frac {1}{3}}  - {\displaystyle \frac {n}{2}} )\,\Gamma (a - 
{\displaystyle \frac {n}{2}} )\,\mathrm{V8}(a, \,n)\,\Gamma (3\,a
 + {\displaystyle \frac {1}{2}}  - {\displaystyle \frac {n}{2}} )
}{\sqrt{\pi }\,\Gamma (2\,a + {\displaystyle \frac {5}{6}}  - n)
\,\Gamma (2\,a)\,\Gamma (2\,a + {\displaystyle \frac {1}{6}}  - n
)}}  \\
\mbox{} + {\displaystyle \frac {1}{2}} \,{\displaystyle \frac {
\Gamma (a + {\displaystyle \frac {1}{3}}  - {\displaystyle 
\frac {n}{2}} )\,\Gamma (a - {\displaystyle \frac {1}{3}}  - 
{\displaystyle \frac {n}{2}} )\,\Gamma (a - {\displaystyle 
\frac {n}{2}} )\,\mathrm{V7}(a - {\displaystyle \frac {1}{3}}  - 
{\displaystyle \frac {n}{2}} , \,n + 1)\,\Gamma (3\,a + 
{\displaystyle \frac {1}{2}}  - {\displaystyle \frac {n}{2}} )}{
\Gamma (a - {\displaystyle \frac {n}{2}}  + {\displaystyle 
\frac {5}{6}} )\,\Gamma (a + {\displaystyle \frac {1}{6}}  - 
{\displaystyle \frac {n}{2}} )\,\Gamma ({\displaystyle \frac {1}{
2}}  + a - {\displaystyle \frac {n}{2}} )\,\Gamma (3\,a - 
{\displaystyle \frac {n}{2}} )}}  }
}
\end{maplelatex}

\begin{maplelatex}
\mapleinline{inert}{2d}{282, ":        ",
-GAMMA(a+1/3-1/2*n)*GAMMA(a-1/3-1/2*n)*GAMMA(a-1/2*n)*GAMMA(4/3+a+1/2*
n)*V8(a,n)/Pi/GAMMA(2*a+5/6-n)/GAMMA(1/2+n)/GAMMA(2*a)/GAMMA(2*a+1/6-n
)*GAMMA(a+1/6-1/2*n)*GAMMA(1+a+1/2*n)+GAMMA(a+1/3-1/2*n)*GAMMA(a-1/3-1
/2*n)*GAMMA(a-1/2*n)*GAMMA(4/3+a+1/2*n)/Pi^(1/2)*V7(a-1/3-1/2*n,n+1)/G
AMMA(a-1/2*n+5/6)/GAMMA(1/2+a-1/2*n)/GAMMA(1/2+n)*GAMMA(1+a+1/2*n)/GAM
MA(3*a-1/2*n);}{%
\maplemultiline{
282\mbox{:     ~~~}    - {\displaystyle \frac {\Gamma (a + 
{\displaystyle \frac {1}{3}}  - {\displaystyle \frac {n}{2}} )\,
\Gamma (a - {\displaystyle \frac {1}{3}}  - {\displaystyle 
\frac {n}{2}} )\,\Gamma (a - {\displaystyle \frac {n}{2}} )\,
\Gamma ({\displaystyle \frac {4}{3}}  + a + {\displaystyle 
\frac {n}{2}} )\,\mathrm{V8}(a, \,n)\,\Gamma (a + {\displaystyle 
\frac {1}{6}}  - {\displaystyle \frac {n}{2}} )\,\Gamma (1 + a + 
{\displaystyle \frac {n}{2}} )}{\pi \,\Gamma (2\,a + 
{\displaystyle \frac {5}{6}}  - n)\,\Gamma ({\displaystyle 
\frac {1}{2}}  + n)\,\Gamma (2\,a)\,\Gamma (2\,a + 
{\displaystyle \frac {1}{6}}  - n)}}  \\
\mbox{} + {\displaystyle \frac {\Gamma (a + {\displaystyle 
\frac {1}{3}}  - {\displaystyle \frac {n}{2}} )\,\Gamma (a - 
{\displaystyle \frac {1}{3}}  - {\displaystyle \frac {n}{2}} )\,
\Gamma (a - {\displaystyle \frac {n}{2}} )\,\Gamma (
{\displaystyle \frac {4}{3}}  + a + {\displaystyle \frac {n}{2}} 
)\,\mathrm{V7}(a - {\displaystyle \frac {1}{3}}  - 
{\displaystyle \frac {n}{2}} , \,n + 1)\,\Gamma (1 + a + 
{\displaystyle \frac {n}{2}} )}{\sqrt{\pi }\,\Gamma (a - 
{\displaystyle \frac {n}{2}}  + {\displaystyle \frac {5}{6}} )\,
\Gamma ({\displaystyle \frac {1}{2}}  + a - {\displaystyle 
\frac {n}{2}} )\,\Gamma ({\displaystyle \frac {1}{2}}  + n)\,
\Gamma (3\,a - {\displaystyle \frac {n}{2}} )}}  }
}
\end{maplelatex}

\begin{maplelatex}
\mapleinline{inert}{2d}{283, ":        ",
-1/2*GAMMA(a+1/3-1/2*n)*GAMMA(a-1/3-1/2*n)*GAMMA(a-1/2*n)*V8(a,n)/Pi^(
1/2)/GAMMA(2*a+5/6-n)/GAMMA(1/2+n)/GAMMA(2*a+1/6-n)*GAMMA(1+a+1/2*n)+1
/2*GAMMA(a+1/3-1/2*n)*GAMMA(a-1/3-1/2*n)*GAMMA(a-1/2*n)*V7(a-1/3-1/2*n
,n+1)*GAMMA(2*a)/GAMMA(a-1/2*n+5/6)/GAMMA(a+1/6-1/2*n)/GAMMA(1/2+a-1/2
*n)/GAMMA(1/2+n)*GAMMA(1+a+1/2*n)/GAMMA(3*a-1/2*n);}{%
\maplemultiline{
283\mbox{:     ~~~}    - {\displaystyle \frac {1}{2}} \,
{\displaystyle \frac {\Gamma (a + {\displaystyle \frac {1}{3}} 
 - {\displaystyle \frac {n}{2}} )\,\Gamma (a - {\displaystyle 
\frac {1}{3}}  - {\displaystyle \frac {n}{2}} )\,\Gamma (a - 
{\displaystyle \frac {n}{2}} )\,\mathrm{V8}(a, \,n)\,\Gamma (1 + 
a + {\displaystyle \frac {n}{2}} )}{\sqrt{\pi }\,\Gamma (2\,a + 
{\displaystyle \frac {5}{6}}  - n)\,\Gamma ({\displaystyle 
\frac {1}{2}}  + n)\,\Gamma (2\,a + {\displaystyle \frac {1}{6}} 
 - n)}}  +  \\
{\displaystyle \frac {1}{2}} \,{\displaystyle \frac {\Gamma (a + 
{\displaystyle \frac {1}{3}}  - {\displaystyle \frac {n}{2}} )\,
\Gamma (a - {\displaystyle \frac {1}{3}}  - {\displaystyle 
\frac {n}{2}} )\,\Gamma (a - {\displaystyle \frac {n}{2}} )\,
\mathrm{V7}(a - {\displaystyle \frac {1}{3}}  - {\displaystyle 
\frac {n}{2}} , \,n + 1)\,\Gamma (2\,a)\,\Gamma (1 + a + 
{\displaystyle \frac {n}{2}} )}{\Gamma (a - {\displaystyle 
\frac {n}{2}}  + {\displaystyle \frac {5}{6}} )\,\Gamma (a + 
{\displaystyle \frac {1}{6}}  - {\displaystyle \frac {n}{2}} )\,
\Gamma ({\displaystyle \frac {1}{2}}  + a - {\displaystyle 
\frac {n}{2}} )\,\Gamma ({\displaystyle \frac {1}{2}}  + n)\,
\Gamma (3\,a - {\displaystyle \frac {n}{2}} )}}  }
}
\end{maplelatex}

\begin{maplelatex}
\mapleinline{inert}{2d}{284, ":        ",
GAMMA(2*a+1/3)*GAMMA(2*a-1/3)*GAMMA(a-1/2*n)*GAMMA(2/3)*V8(a,n)/GAMMA(
a-1/2+1/2*n)/GAMMA(2*a+5/6-n)/GAMMA(1/2-a+1/2*n)/GAMMA(2*a)/GAMMA(2*a+
1/6-n)*GAMMA(a-1/2*n+5/6)*GAMMA(1-1/2*n-a)/GAMMA(3/2-a-1/2*n)-GAMMA(2*
a+1/3)*GAMMA(2*a-1/3)*GAMMA(a-1/2*n)/Pi^(1/2)*GAMMA(2/3)*GAMMA(-a-1/2*
n+1/2)*V7(a-1/3-1/2*n,n+1)*GAMMA(a+1/2+1/2*n)/GAMMA(a-1/2+1/2*n)/GAMMA
(a+1/6-1/2*n)/GAMMA(-n+1/2)/GAMMA(1/2+n)*GAMMA(1-1/2*n-a)/GAMMA(3/2-a-
1/2*n)/GAMMA(3*a-1/2*n);}{%
\maplemultiline{
284\mbox{:     ~~~}   {\displaystyle \frac {\Gamma (2\,a + 
{\displaystyle \frac {1}{3}} )\,\Gamma (2\,a - {\displaystyle 
\frac {1}{3}} )\,\Gamma (a - {\displaystyle \frac {n}{2}} )\,
\Gamma ({\displaystyle \frac {2}{3}} )\,\mathrm{V8}(a, \,n)\,
\Gamma (a - {\displaystyle \frac {n}{2}}  + {\displaystyle 
\frac {5}{6}} )\,\Gamma (1 - {\displaystyle \frac {n}{2}}  - a)}{
\Gamma (a - {\displaystyle \frac {1}{2}}  + {\displaystyle 
\frac {n}{2}} )\,\Gamma (2\,a + {\displaystyle \frac {5}{6}}  - n
)\,\Gamma ({\displaystyle \frac {1}{2}}  - a + {\displaystyle 
\frac {n}{2}} )\,\Gamma (2\,a)\,\Gamma (2\,a + {\displaystyle 
\frac {1}{6}}  - n)\,\Gamma ({\displaystyle \frac {3}{2}}  - a - 
{\displaystyle \frac {n}{2}} )}}  -  \\
\Gamma (2\,a + {\displaystyle \frac {1}{3}} )\,\Gamma (2\,a - 
{\displaystyle \frac {1}{3}} )\,\Gamma (a - {\displaystyle 
\frac {n}{2}} )\,\Gamma ({\displaystyle \frac {2}{3}} )\,\Gamma (
 - a - {\displaystyle \frac {n}{2}}  + {\displaystyle \frac {1}{2
}} )\,\mathrm{V7}(a - {\displaystyle \frac {1}{3}}  - 
{\displaystyle \frac {n}{2}} , \,n + 1) \\
\Gamma (a + {\displaystyle \frac {1}{2}}  + {\displaystyle 
\frac {n}{2}} )\,\Gamma (1 - {\displaystyle \frac {n}{2}}  - a)
 \left/ {\vrule height0.80em width0em depth0.80em} \right. \! 
 \! (\sqrt{\pi }\,\Gamma (a - {\displaystyle \frac {1}{2}}  + 
{\displaystyle \frac {n}{2}} )\,\Gamma (a + {\displaystyle 
\frac {1}{6}}  - {\displaystyle \frac {n}{2}} )\,\Gamma ( - n + 
{\displaystyle \frac {1}{2}} ) 
\Gamma ({\displaystyle \frac {1}{2}}  + n)\,\Gamma (
{\displaystyle \frac {3}{2}}  - a - {\displaystyle \frac {n}{2}} 
)\,\Gamma (3\,a - {\displaystyle \frac {n}{2}} )) }
}
\end{maplelatex}

\begin{maplelatex}
\mapleinline{inert}{2d}{285, ":        ",
GAMMA(2*a+1/3)*GAMMA(2*a-1/3)*GAMMA(a-1/2*n)*V8(a,n)/GAMMA(a-1/2+1/2*n
)/GAMMA(2*a+5/6-n)/GAMMA(1/2-a+1/2*n)/GAMMA(2*a)/GAMMA(2*a+1/6-n)*GAMM
A(1-a+1/2*n)*GAMMA(1-1/2*n-a)*GAMMA(1/2+2*a-n)/GAMMA(3/2-a-1/2*n)-GAMM
A(2*a+1/3)*GAMMA(2*a-1/3)*GAMMA(a-1/2*n)/Pi^(1/2)*GAMMA(-a-1/2*n+1/2)*
V7(a-1/3-1/2*n,n+1)*GAMMA(a+1/2+1/2*n)/GAMMA(a-1/2+1/2*n)/GAMMA(a-1/2*
n+5/6)/GAMMA(a+1/6-1/2*n)/GAMMA(-n+1/2)/GAMMA(1/2+n)*GAMMA(1-a+1/2*n)*
GAMMA(1-1/2*n-a)*GAMMA(1/2+2*a-n)/GAMMA(3/2-a-1/2*n)/GAMMA(3*a-1/2*n);
}{%
\maplemultiline{
285\mbox{:     ~~~}    \\
{\displaystyle \frac {\Gamma (2\,a + {\displaystyle \frac {1}{3}
} )\,\Gamma (2\,a - {\displaystyle \frac {1}{3}} )\,\Gamma (a - 
{\displaystyle \frac {n}{2}} )\,\mathrm{V8}(a, \,n)\,\Gamma (1 - 
a + {\displaystyle \frac {n}{2}} )\,\Gamma (1 - {\displaystyle 
\frac {n}{2}}  - a)\,\Gamma ({\displaystyle \frac {1}{2}}  + 2\,a
 - n)}{\Gamma (a - {\displaystyle \frac {1}{2}}  + 
{\displaystyle \frac {n}{2}} )\,\Gamma (2\,a + {\displaystyle 
\frac {5}{6}}  - n)\,\Gamma ({\displaystyle \frac {1}{2}}  - a + 
{\displaystyle \frac {n}{2}} )\,\Gamma (2\,a)\,\Gamma (2\,a + 
{\displaystyle \frac {1}{6}}  - n)\,\Gamma ({\displaystyle 
\frac {3}{2}}  - a - {\displaystyle \frac {n}{2}} )}}  \\
\mbox{} - \Gamma (2\,a + {\displaystyle \frac {1}{3}} )\,\Gamma (
2\,a - {\displaystyle \frac {1}{3}} )\,\Gamma (a - 
{\displaystyle \frac {n}{2}} )\,\Gamma ( - a - {\displaystyle 
\frac {n}{2}}  + {\displaystyle \frac {1}{2}} )\,\mathrm{V7}(a - 
{\displaystyle \frac {1}{3}}  - {\displaystyle \frac {n}{2}} , \,
n + 1) \\
\Gamma (a + {\displaystyle \frac {1}{2}}  + {\displaystyle 
\frac {n}{2}} )\,\Gamma (1 - a + {\displaystyle \frac {n}{2}} )\,
\Gamma (1 - {\displaystyle \frac {n}{2}}  - a)\,\Gamma (
{\displaystyle \frac {1}{2}}  + 2\,a - n) \left/ {\vrule 
height0.80em width0em depth0.80em} \right. \!  \! (\sqrt{\pi }\,
\Gamma (a - {\displaystyle \frac {1}{2}}  + {\displaystyle 
\frac {n}{2}} ) \\
\Gamma (a - {\displaystyle \frac {n}{2}}  + {\displaystyle 
\frac {5}{6}} )\,\Gamma (a + {\displaystyle \frac {1}{6}}  - 
{\displaystyle \frac {n}{2}} )\,\Gamma ( - n + {\displaystyle 
\frac {1}{2}} )\,\Gamma ({\displaystyle \frac {1}{2}}  + n)\,
\Gamma ({\displaystyle \frac {3}{2}}  - a - {\displaystyle 
\frac {n}{2}} )\,\Gamma (3\,a - {\displaystyle \frac {n}{2}} ))
 }
}
\end{maplelatex}

\begin{maplelatex}
\mapleinline{inert}{2d}{286, ":        ",
2/9*GAMMA(2*a+1/3)*GAMMA(2*a-1/3)*GAMMA(a-1/2*n)*Pi*3^(1/2)/GAMMA(2/3)
*V8(a,n)/GAMMA(a-1/2+1/2*n)/GAMMA(2*a+5/6-n)/GAMMA(1/2-a+1/2*n)/GAMMA(
2*a)/GAMMA(2*a+1/6-n)*GAMMA(a+1/6-1/2*n)*GAMMA(1-1/2*n-a)/GAMMA(3/2-a-
1/2*n)-2/9*GAMMA(2*a+1/3)*GAMMA(2*a-1/3)*GAMMA(a-1/2*n)*Pi^(1/2)*3^(1/
2)/GAMMA(2/3)*GAMMA(-a-1/2*n+1/2)*V7(a-1/3-1/2*n,n+1)*GAMMA(a+1/2+1/2*
n)/GAMMA(a-1/2+1/2*n)/GAMMA(a-1/2*n+5/6)/GAMMA(-n+1/2)/GAMMA(1/2+n)*GA
MMA(1-1/2*n-a)/GAMMA(3/2-a-1/2*n)/GAMMA(3*a-1/2*n);}{%
\maplemultiline{
286\mbox{:     ~~~}   {\displaystyle \frac {2}{9}} 
{\displaystyle \frac {\Gamma (2\,a + {\displaystyle \frac {1}{3}
} )\,\Gamma (2\,a - {\displaystyle \frac {1}{3}} )\,\Gamma (a - 
{\displaystyle \frac {n}{2}} )\,\pi \,\sqrt{3}\,\mathrm{V8}(a, \,
n)\,\Gamma (a + {\displaystyle \frac {1}{6}}  - {\displaystyle 
\frac {n}{2}} )\,\Gamma (1 - {\displaystyle \frac {n}{2}}  - a)}{
\Gamma ({\displaystyle \frac {2}{3}} )\,\Gamma (a - 
{\displaystyle \frac {1}{2}}  + {\displaystyle \frac {n}{2}} )\,
\Gamma (2\,a + {\displaystyle \frac {5}{6}}  - n)\,\Gamma (
{\displaystyle \frac {1}{2}}  - a + {\displaystyle \frac {n}{2}} 
)\,\Gamma (2\,a)\,\Gamma (2\,a + {\displaystyle \frac {1}{6}}  - 
n)\,\Gamma ({\displaystyle \frac {3}{2}}  - a - {\displaystyle 
\frac {n}{2}} )}}  \\
\mbox{} - {\displaystyle \frac {2}{9}} \Gamma (2\,a + 
{\displaystyle \frac {1}{3}} )\,\Gamma (2\,a - {\displaystyle 
\frac {1}{3}} )\,\Gamma (a - {\displaystyle \frac {n}{2}} )\,
\sqrt{\pi }\,\sqrt{3}\,\Gamma ( - a - {\displaystyle \frac {n}{2}
}  + {\displaystyle \frac {1}{2}} )\,\mathrm{V7}(a - 
{\displaystyle \frac {1}{3}}  - {\displaystyle \frac {n}{2}} , \,
n + 1) \\
\Gamma (a + {\displaystyle \frac {1}{2}}  + {\displaystyle 
\frac {n}{2}} )\,\Gamma (1 - {\displaystyle \frac {n}{2}}  - a)
 \left/ {\vrule height0.80em width0em depth0.80em} \right. \! 
 \! (\Gamma ({\displaystyle \frac {2}{3}} )\,\Gamma (a - 
{\displaystyle \frac {1}{2}}  + {\displaystyle \frac {n}{2}} )\,
\Gamma (a - {\displaystyle \frac {n}{2}}  + {\displaystyle 
\frac {5}{6}} )\,\Gamma ( - n + {\displaystyle \frac {1}{2}} )
\Gamma ({\displaystyle \frac {1}{2}}  + n)\,\Gamma (
{\displaystyle \frac {3}{2}}  - a - {\displaystyle \frac {n}{2}} 
)\,\Gamma (3\,a - {\displaystyle \frac {n}{2}} )) }
}
\end{maplelatex}

\begin{maplelatex}
\mapleinline{inert}{2d}{287, ":        ",
GAMMA(2*a+1/3)*GAMMA(2*a-1/3)*GAMMA(a-1/2*n)*GAMMA(2/3)*V8(a,n)/GAMMA(
a-1/2+1/2*n)/GAMMA(2*a+5/6-n)/GAMMA(1/2-a+1/2*n)/GAMMA(2*a)/GAMMA(2*a+
1/6-n)*GAMMA(4/3-2*a)*GAMMA(1/2+2*a-n)/GAMMA(3/2-a-1/2*n)-GAMMA(2*a+1/
3)*GAMMA(2*a-1/3)*GAMMA(a-1/2*n)/Pi^(1/2)*GAMMA(2/3)*GAMMA(-a-1/2*n+1/
2)*V7(a-1/3-1/2*n,n+1)*GAMMA(a+1/2+1/2*n)/GAMMA(a-1/2+1/2*n)/GAMMA(a-1
/2*n+5/6)/GAMMA(a+1/6-1/2*n)/GAMMA(-n+1/2)/GAMMA(1/2+n)*GAMMA(4/3-2*a)
*GAMMA(1/2+2*a-n)/GAMMA(3/2-a-1/2*n)/GAMMA(3*a-1/2*n);}{%
\maplemultiline{
287\mbox{:     ~~~}   {\displaystyle \frac {\Gamma (2\,a + 
{\displaystyle \frac {1}{3}} )\,\Gamma (2\,a - {\displaystyle 
\frac {1}{3}} )\,\Gamma (a - {\displaystyle \frac {n}{2}} )\,
\Gamma ({\displaystyle \frac {2}{3}} )\,\mathrm{V8}(a, \,n)\,
\Gamma ({\displaystyle \frac {4}{3}}  - 2\,a)\,\Gamma (
{\displaystyle \frac {1}{2}}  + 2\,a - n)}{\Gamma (a - 
{\displaystyle \frac {1}{2}}  + {\displaystyle \frac {n}{2}} )\,
\Gamma (2\,a + {\displaystyle \frac {5}{6}}  - n)\,\Gamma (
{\displaystyle \frac {1}{2}}  - a + {\displaystyle \frac {n}{2}} 
)\,\Gamma (2\,a)\,\Gamma (2\,a + {\displaystyle \frac {1}{6}}  - 
n)\,\Gamma ({\displaystyle \frac {3}{2}}  - a - {\displaystyle 
\frac {n}{2}} )}}  -  \\
\Gamma (2\,a + {\displaystyle \frac {1}{3}} )\,\Gamma (2\,a - 
{\displaystyle \frac {1}{3}} )\,\Gamma (a - {\displaystyle 
\frac {n}{2}} )\,\Gamma ({\displaystyle \frac {2}{3}} )\,\Gamma (
 - a - {\displaystyle \frac {n}{2}}  + {\displaystyle \frac {1}{2
}} )\,\mathrm{V7}(a - {\displaystyle \frac {1}{3}}  - 
{\displaystyle \frac {n}{2}} , \,n + 1) \\
\Gamma (a + {\displaystyle \frac {1}{2}}  + {\displaystyle 
\frac {n}{2}} )\,\Gamma ({\displaystyle \frac {4}{3}}  - 2\,a)\,
\Gamma ({\displaystyle \frac {1}{2}}  + 2\,a - n) \left/ {\vrule 
height0.80em width0em depth0.80em} \right. \!  \! (\sqrt{\pi }\,
\Gamma (a - {\displaystyle \frac {1}{2}}  + {\displaystyle 
\frac {n}{2}} )\,\Gamma (a - {\displaystyle \frac {n}{2}}  + 
{\displaystyle \frac {5}{6}} ) \\
\Gamma (a + {\displaystyle \frac {1}{6}}  - {\displaystyle 
\frac {n}{2}} )\,\Gamma ( - n + {\displaystyle \frac {1}{2}} )\,
\Gamma ({\displaystyle \frac {1}{2}}  + n)\,\Gamma (
{\displaystyle \frac {3}{2}}  - a - {\displaystyle \frac {n}{2}} 
)\,\Gamma (3\,a - {\displaystyle \frac {n}{2}} )) }
}
\end{maplelatex}

\begin{maplelatex}
\mapleinline{inert}{2d}{288, ":        ",
2/9*GAMMA(2*a+1/3)*GAMMA(2*a-1/3)*GAMMA(a-1/2*n)*Pi*3^(1/2)*V8(a,n)/GA
MMA(a-1/2+1/2*n)/GAMMA(2*a+5/6-n)/GAMMA(1/2-a+1/2*n)/GAMMA(2*a)/GAMMA(
2*a+1/6-n)*GAMMA(-n+1/2)/GAMMA(3/2-a-1/2*n)-2/9*GAMMA(2*a+1/3)*GAMMA(2
*a-1/3)*GAMMA(a-1/2*n)*Pi^(1/2)*3^(1/2)*GAMMA(-a-1/2*n+1/2)*V7(a-1/3-1
/2*n,n+1)*GAMMA(a+1/2+1/2*n)/GAMMA(a-1/2+1/2*n)/GAMMA(a-1/2*n+5/6)/GAM
MA(a+1/6-1/2*n)/GAMMA(1/2+n)/GAMMA(3/2-a-1/2*n)/GAMMA(3*a-1/2*n);}{%
\maplemultiline{
288\mbox{:     ~~~}   {\displaystyle \frac {2}{9}} \,
{\displaystyle \frac {\Gamma (2\,a + {\displaystyle \frac {1}{3}
} )\,\Gamma (2\,a - {\displaystyle \frac {1}{3}} )\,\Gamma (a - 
{\displaystyle \frac {n}{2}} )\,\pi \,\sqrt{3}\,\mathrm{V8}(a, \,
n)\,\Gamma ( - n + {\displaystyle \frac {1}{2}} )}{\Gamma (a - 
{\displaystyle \frac {1}{2}}  + {\displaystyle \frac {n}{2}} )\,
\Gamma (2\,a + {\displaystyle \frac {5}{6}}  - n)\,\Gamma (
{\displaystyle \frac {1}{2}}  - a + {\displaystyle \frac {n}{2}} 
)\,\Gamma (2\,a)\,\Gamma (2\,a + {\displaystyle \frac {1}{6}}  - 
n)\,\Gamma ({\displaystyle \frac {3}{2}}  - a - {\displaystyle 
\frac {n}{2}} )}}  \\
\mbox{} - {\displaystyle \frac {2}{9}} \,{\displaystyle \frac {
\Gamma (2\,a + {\displaystyle \frac {1}{3}} )\,\Gamma (2\,a - 
{\displaystyle \frac {1}{3}} )\,\Gamma (a - {\displaystyle 
\frac {n}{2}} )\,\sqrt{\pi }\,\sqrt{3}\,\Gamma ( - a - 
{\displaystyle \frac {n}{2}}  + {\displaystyle \frac {1}{2}} )\,
\mathrm{V7}(a - {\displaystyle \frac {1}{3}}  - {\displaystyle 
\frac {n}{2}} , \,n + 1)\,\Gamma (a + {\displaystyle \frac {1}{2}
}  + {\displaystyle \frac {n}{2}} )}{\Gamma (a - {\displaystyle 
\frac {1}{2}}  + {\displaystyle \frac {n}{2}} )\,\Gamma (a - 
{\displaystyle \frac {n}{2}}  + {\displaystyle \frac {5}{6}} )\,
\Gamma (a + {\displaystyle \frac {1}{6}}  - {\displaystyle 
\frac {n}{2}} )\,\Gamma ({\displaystyle \frac {1}{2}}  + n)\,
\Gamma ({\displaystyle \frac {3}{2}}  - a - {\displaystyle 
\frac {n}{2}} )\,\Gamma (3\,a - {\displaystyle \frac {n}{2}} )}} 
 }
}
\end{maplelatex}

\begin{maplelatex}
\mapleinline{inert}{2d}{289, ":        ",
GAMMA(2*a+1/3)*GAMMA(2*a-1/3)*GAMMA(a-1/2*n)*V8(a,n)/GAMMA(a-1/2+1/2*n
)/GAMMA(2*a+5/6-n)/GAMMA(2*a)/GAMMA(2*a+1/6-n)*GAMMA(1/2+2*a-n)-GAMMA(
2*a+1/3)*GAMMA(2*a-1/3)*GAMMA(a-1/2*n)/Pi^(1/2)*GAMMA(-a-1/2*n+1/2)*V7
(a-1/3-1/2*n,n+1)*GAMMA(a+1/2+1/2*n)/GAMMA(a-1/2+1/2*n)/GAMMA(a-1/2*n+
5/6)/GAMMA(a+1/6-1/2*n)/GAMMA(-n+1/2)/GAMMA(1/2+n)*GAMMA(1/2-a+1/2*n)*
GAMMA(1/2+2*a-n)/GAMMA(3*a-1/2*n);}{%
\maplemultiline{
289\mbox{:     ~~~}   {\displaystyle \frac {\Gamma (2\,a + 
{\displaystyle \frac {1}{3}} )\,\Gamma (2\,a - {\displaystyle 
\frac {1}{3}} )\,\Gamma (a - {\displaystyle \frac {n}{2}} )\,
\mathrm{V8}(a, \,n)\,\Gamma ({\displaystyle \frac {1}{2}}  + 2\,a
 - n)}{\Gamma (a - {\displaystyle \frac {1}{2}}  + 
{\displaystyle \frac {n}{2}} )\,\Gamma (2\,a + {\displaystyle 
\frac {5}{6}}  - n)\,\Gamma (2\,a)\,\Gamma (2\,a + 
{\displaystyle \frac {1}{6}}  - n)}}  - \Gamma (2\,a + 
{\displaystyle \frac {1}{3}} )\,\Gamma (2\,a - {\displaystyle 
\frac {1}{3}} ) \\
\Gamma (a - {\displaystyle \frac {n}{2}} )\,\Gamma ( - a - 
{\displaystyle \frac {n}{2}}  + {\displaystyle \frac {1}{2}} )\,
\mathrm{V7}(a - {\displaystyle \frac {1}{3}}  - {\displaystyle 
\frac {n}{2}} , \,n + 1)\,\Gamma (a + {\displaystyle \frac {1}{2}
}  + {\displaystyle \frac {n}{2}} )\,\Gamma ({\displaystyle 
\frac {1}{2}}  - a + {\displaystyle \frac {n}{2}} ) \\
\Gamma ({\displaystyle \frac {1}{2}}  + 2\,a - n) \left/ {\vrule 
height0.80em width0em depth0.80em} \right. \!  \! (\sqrt{\pi }\,
\Gamma (a - {\displaystyle \frac {1}{2}}  + {\displaystyle 
\frac {n}{2}} )\,\Gamma (a - {\displaystyle \frac {n}{2}}  + 
{\displaystyle \frac {5}{6}} )\,\Gamma (a + {\displaystyle 
\frac {1}{6}}  - {\displaystyle \frac {n}{2}} )\,\Gamma ( - n + 
{\displaystyle \frac {1}{2}} )
\Gamma ({\displaystyle \frac {1}{2}}  + n)\,\Gamma (3\,a - 
{\displaystyle \frac {n}{2}} )) }
}
\end{maplelatex}

\mapleinline{inert}{2d}{290, ":        ",
2/9*GAMMA(2*a+1/3)*GAMMA(2*a-1/3)*GAMMA(a-1/2*n)*Pi*3^(1/2)/GAMMA(2/3)
*V8(a,n)/GAMMA(a-1/2+1/2*n)/GAMMA(2*a+5/6-n)/GAMMA(1/2-a+1/2*n)/GAMMA(
2*a)/GAMMA(2*a+1/6-n)*GAMMA(2/3-2*a)*GAMMA(1/2+2*a-n)/GAMMA(3/2-a-1/2*
n)-2/9*GAMMA(2*a+1/3)*GAMMA(2*a-1/3)*GAMMA(a-1/2*n)*Pi^(1/2)*3^(1/2)/G
AMMA(2/3)*GAMMA(-a-1/2*n+1/2)*V7(a-1/3-1/2*n,n+1)*GAMMA(a+1/2+1/2*n)/G
AMMA(a-1/2+1/2*n)/GAMMA(a-1/2*n+5/6)/GAMMA(a+1/6-1/2*n)/GAMMA(-n+1/2)/
GAMMA(1/2+n)*GAMMA(2/3-2*a)*GAMMA(1/2+2*a-n)/GAMMA(3/2-a-1/2*n)/GAMMA(
3*a-1/2*n);}{%
\maplemultiline{
290\mbox{:~~~}   {\displaystyle \frac {2}{9}} 
{\displaystyle \frac {\Gamma (2\,a + {\displaystyle \frac {1}{3}
} )\,\Gamma (2\,a - {\displaystyle \frac {1}{3}} )\,\Gamma (a - 
{\displaystyle \frac {n}{2}} )\,\pi \,\sqrt{3}\,\mathrm{V8}(a, \,
n)\,\Gamma ({\displaystyle \frac {2}{3}}  - 2\,a)\,\Gamma (
{\displaystyle \frac {1}{2}}  + 2\,a - n)}{\Gamma (
{\displaystyle \frac {2}{3}} )\,\Gamma (a - {\displaystyle 
\frac {1}{2}}  + {\displaystyle \frac {n}{2}} )\,\Gamma (2\,a + 
{\displaystyle \frac {5}{6}}  - n)\,\Gamma ({\displaystyle 
\frac {1}{2}}  - a + {\displaystyle \frac {n}{2}} )\,\Gamma (2\,a
)\,\Gamma (2\,a + {\displaystyle \frac {1}{6}}  - n)\,\Gamma (
{\displaystyle \frac {3}{2}}  - a - {\displaystyle \frac {n}{2}} 
)}}  \\
\mbox{} - {\displaystyle \frac {2}{9}} \Gamma (2\,a + 
{\displaystyle \frac {1}{3}} )\,\Gamma (2\,a - {\displaystyle 
\frac {1}{3}} )\,\Gamma (a - {\displaystyle \frac {n}{2}} )\,
\sqrt{\pi }\,\sqrt{3}\,\Gamma ( - a - {\displaystyle \frac {n}{2}
}  + {\displaystyle \frac {1}{2}} )\,\mathrm{V7}(a - 
{\displaystyle \frac {1}{3}}  - {\displaystyle \frac {n}{2}} , \,
n + 1) \\
\Gamma (a + {\displaystyle \frac {1}{2}}  + {\displaystyle 
\frac {n}{2}} )\,\Gamma ({\displaystyle \frac {2}{3}}  - 2\,a)\,
\Gamma ({\displaystyle \frac {1}{2}}  + 2\,a - n) \left/ {\vrule 
height0.80em width0em depth0.80em} \right. \!  \! (\Gamma (
{\displaystyle \frac {2}{3}} )\,\Gamma (a - {\displaystyle 
\frac {1}{2}}  + {\displaystyle \frac {n}{2}} )\,\Gamma (a - 
{\displaystyle \frac {n}{2}}  + {\displaystyle \frac {5}{6}} )
 \\
\Gamma (a + {\displaystyle \frac {1}{6}}  - {\displaystyle 
\frac {n}{2}} )\,\Gamma ( - n + {\displaystyle \frac {1}{2}} )\,
\Gamma ({\displaystyle \frac {1}{2}}  + n)\,\Gamma (
{\displaystyle \frac {3}{2}}  - a - {\displaystyle \frac {n}{2}} 
)\,\Gamma (3\,a - {\displaystyle \frac {n}{2}} )) }
}

\begin{mapleinput}
\end{mapleinput}

\end{maplegroup}

%% file: AppendixB291to320.tex
\begin{maplegroup}
\mapleresult
\begin{maplelatex}
\mapleinline{inert}{2d}{293, ":   ",
V7(a-1/3-1/2*n,n+1)*2^(6*a-n-2)/((-1)^(-n))/Pi^2/GAMMA(6*a-n-1)*GAMMA(
1-1/2*n-a)*GAMMA(a-1/2*n)*GAMMA(2*a+1/3)*GAMMA(2*a-1/3)*GAMMA(2*a)-1/4
*(2*sin(1/6*Pi*(-12*a+6*n+1))+1)*V8(a,n)*GAMMA(2*a-n)*GAMMA(1-1/2*n-a)
*GAMMA(a-1/2*n+5/6)*GAMMA(2*a+1/3)*GAMMA(2*a-1/3)/cos(1/2*Pi*(-6*a+n))
/Pi^(1/2)/(2^(2*a-n-1))/GAMMA(3*a-1/2*n-1/2)/GAMMA(2*a+1/6-n)/GAMMA(2*
a+5/6-n)/GAMMA(-a+5/6+1/2*n);}{%
\maplemultiline{
293\mbox{:~~~} {\displaystyle \frac {\mathrm{V7}(a - 
{\displaystyle \frac {1}{3}}  - {\displaystyle \frac {n}{2}} , \,
n + 1)\,2^{(6\,a - n - 2)}\,\Gamma (1 - {\displaystyle \frac {n}{
2}}  - a)\,\Gamma (a - {\displaystyle \frac {n}{2}} )\,\Gamma (2
\,a + {\displaystyle \frac {1}{3}} )\,\Gamma (2\,a - 
{\displaystyle \frac {1}{3}} )\,\Gamma (2\,a)}{(-1)^{( - n)}\,\pi
 ^{2}\,\Gamma (6\,a - n - 1)}}  \\
\mbox{} - {\displaystyle \frac {1}{4}} (2\,\mathrm{sin}(
{\displaystyle \frac {\pi \,( - 12\,a + 6\,n + 1)}{6}} ) + 1)\,
\mathrm{V8}(a, \,n)\,\Gamma (2\,a - n)\,\Gamma (1 - 
{\displaystyle \frac {n}{2}}  - a)\,\Gamma (a - {\displaystyle 
\frac {n}{2}}  + {\displaystyle \frac {5}{6}} ) \\
\Gamma (2\,a + {\displaystyle \frac {1}{3}} )\,\Gamma (2\,a - 
{\displaystyle \frac {1}{3}} ) \left/ {\vrule 
height0.80em width0em depth0.80em} \right. \!  \! (\mathrm{cos}(
{\displaystyle \frac {\pi \,( - 6\,a + n)}{2}} )\,\sqrt{\pi }\,2
^{(2\,a - n - 1)}\,\Gamma (3\,a - {\displaystyle \frac {n}{2}} 
 - {\displaystyle \frac {1}{2}} ) \\
\Gamma (2\,a + {\displaystyle \frac {1}{6}}  - n)\,\Gamma (2\,a
 + {\displaystyle \frac {5}{6}}  - n)\,\Gamma ( - a + 
{\displaystyle \frac {5}{6}}  + {\displaystyle \frac {n}{2}} ))
 }
}
\end{maplelatex}

\begin{maplelatex}
\mapleinline{inert}{2d}{294, ":   ",
GAMMA(2*a-1/3)*GAMMA(2*a+1/3)*GAMMA(2/3)*GAMMA(a+1/3-1/2*n)*GAMMA(a-1/
2*n)/Pi^(1/2)/GAMMA(1/2+n)/GAMMA(2*a-n-1/6)/GAMMA(a+1/6-1/2*n)/GAMMA(3
*a-1/2*n)*V7(a-1/3-1/2*n,n+1)+1/24*V8(a,n)*GAMMA(2*a-n)*GAMMA(2*a+2/3-
n)*GAMMA(2*a+1/3)*GAMMA(2*a-1/3)*GAMMA(2/3)*(2*sin(1/6*Pi*(-12*a+6*n+1
))+1)*(-12*a+6*n+1)/sin(1/6*Pi*(-12*a+6*n+1))/(2^(2*a-n-1))/(2^(-1/3+2
*a-n))/GAMMA(2*a+1/6-n)/GAMMA(2*a+5/6-n)^2/GAMMA(1/2+n)/GAMMA(2*a);}{%
\maplemultiline{
294\mbox{:~~~} {\displaystyle \frac {\Gamma (2\,a - 
{\displaystyle \frac {1}{3}} )\,\Gamma (2\,a + {\displaystyle 
\frac {1}{3}} )\,\Gamma ({\displaystyle \frac {2}{3}} )\,\Gamma (
a + {\displaystyle \frac {1}{3}}  - {\displaystyle \frac {n}{2}} 
)\,\Gamma (a - {\displaystyle \frac {n}{2}} )\,\mathrm{V7}(a - 
{\displaystyle \frac {1}{3}}  - {\displaystyle \frac {n}{2}} , \,
n + 1)}{\sqrt{\pi }\,\Gamma ({\displaystyle \frac {1}{2}}  + n)\,
\Gamma (2\,a - n - {\displaystyle \frac {1}{6}} )\,\Gamma (a + 
{\displaystyle \frac {1}{6}}  - {\displaystyle \frac {n}{2}} )\,
\Gamma (3\,a - {\displaystyle \frac {n}{2}} )}}  + 
{\displaystyle \frac {1}{24}}  \\
\mathrm{V8}(a, \,n)\,\Gamma (2\,a - n)\,\Gamma (2\,a + 
{\displaystyle \frac {2}{3}}  - n)\,\Gamma (2\,a + 
{\displaystyle \frac {1}{3}} )\,\Gamma (2\,a - {\displaystyle 
\frac {1}{3}} )\,\Gamma ({\displaystyle \frac {2}{3}} ) \\
(2\,\mathrm{sin}({\displaystyle \frac {\pi \,( - 12\,a + 6\,n + 1
)}{6}} ) + 1)\,( - 12\,a + 6\,n + 1) \left/ {\vrule 
height0.80em width0em depth0.80em} \right. \!  \! (\mathrm{sin}(
{\displaystyle \frac {\pi \,( - 12\,a + 6\,n + 1)}{6}} ) \\
2^{(2\,a - n - 1)}\,2^{( - 1/3 + 2\,a - n)}\,\Gamma (2\,a + 
{\displaystyle \frac {1}{6}}  - n)\,\Gamma (2\,a + 
{\displaystyle \frac {5}{6}}  - n)^{2}\,\Gamma ({\displaystyle 
\frac {1}{2}}  + n)\,\Gamma (2\,a)) }
}
\end{maplelatex}

\begin{maplelatex}
\mapleinline{inert}{2d}{295, ":   ",
2/3*V7(a-1/3-1/2*n,n+1)*3^(1/2)*Pi^(1/2)/GAMMA(3*a-1/2*n)/GAMMA(1/2+n)
/GAMMA(a+1/6-1/2*n)/GAMMA(2*a-n-1/6)*GAMMA(a-1/2*n)*GAMMA(a+1/3-1/2*n)
*GAMMA(2*a+1/3)/GAMMA(2/3)*GAMMA(2*a)-1/18*3^(1/2)*Pi^(3/2)*(-1)^(-n)*
V8(a,n)*GAMMA(2*a+2/3-n)*GAMMA(a-1/2*n)*GAMMA(2*a+1/3)*(-sin(1/6*Pi*(9
*n-6*a+1))+sin(1/6*Pi*(6*a-1+3*n)))*(-12*a+6*n+1)/sin(1/6*Pi*(-12*a+6*
n+1))/(2^(-1/3+2*a-n))/GAMMA(2*a+1/6-n)/GAMMA(2*a+5/6-n)^2/GAMMA(1/2-a
+1/2*n)/GAMMA(1/2+n)/GAMMA(2/3);}{%
\maplemultiline{
295\mbox{:~~~} {\displaystyle \frac {2}{3}} \,
{\displaystyle \frac {\mathrm{V7}(a - {\displaystyle \frac {1}{3}
}  - {\displaystyle \frac {n}{2}} , \,n + 1)\,\sqrt{3}\,\sqrt{\pi
 }\,\Gamma (a - {\displaystyle \frac {n}{2}} )\,\Gamma (a + 
{\displaystyle \frac {1}{3}}  - {\displaystyle \frac {n}{2}} )\,
\Gamma (2\,a + {\displaystyle \frac {1}{3}} )\,\Gamma (2\,a)}{
\Gamma (3\,a - {\displaystyle \frac {n}{2}} )\,\Gamma (
{\displaystyle \frac {1}{2}}  + n)\,\Gamma (a + {\displaystyle 
\frac {1}{6}}  - {\displaystyle \frac {n}{2}} )\,\Gamma (2\,a - n
 - {\displaystyle \frac {1}{6}} )\,\Gamma ({\displaystyle \frac {
2}{3}} )}}  - {\displaystyle \frac {1}{18}}  \\
\sqrt{3}\,\pi ^{(3/2)}\,(-1)^{( - n)}\,\mathrm{V8}(a, \,n)\,
\Gamma (2\,a + {\displaystyle \frac {2}{3}}  - n)\,\Gamma (a - 
{\displaystyle \frac {n}{2}} )\,\Gamma (2\,a + {\displaystyle 
\frac {1}{3}} ) \\
( - \mathrm{sin}({\displaystyle \frac {\pi \,(9\,n - 6\,a + 1)}{6
}} ) + \mathrm{sin}({\displaystyle \frac {\pi \,(6\,a - 1 + 3\,n)
}{6}} ))\,( - 12\,a + 6\,n + 1) \left/ {\vrule 
height0.84em width0em depth0.84em} \right. \!  \! ( \\
\mathrm{sin}({\displaystyle \frac {\pi \,( - 12\,a + 6\,n + 1)}{6
}} )\,2^{( - 1/3 + 2\,a - n)}\,\Gamma (2\,a + {\displaystyle 
\frac {1}{6}}  - n)\,\Gamma (2\,a + {\displaystyle \frac {5}{6}} 
 - n)^{2}\,\Gamma ({\displaystyle \frac {1}{2}}  - a + 
{\displaystyle \frac {n}{2}} ) \\
\Gamma ({\displaystyle \frac {1}{2}}  + n)\,\Gamma (
{\displaystyle \frac {2}{3}} )) }
}
\end{maplelatex}

\begin{maplelatex}
\mapleinline{inert}{2d}{298, ":   ",
-1/2*1/((-1)^(-n))/Pi^(3/2)/GAMMA(3*a-1/2*n)^2/GAMMA(1/2+a-1/2*n)/GAMM
A(a+1/6-1/2*n)*V7(a-1/3-1/2*n,n+1)*(-6*a+n+2)*GAMMA(a-1/3-1/2*n)*GAMMA
(7/6-a-1/2*n)*GAMMA(a+1/3-1/2*n)*GAMMA(a-1/2*n)*GAMMA(2*a+1/3)*GAMMA(2
*a-1/3)*GAMMA(2*a)-1/2*V8(a,n)*(sin(3/2*Pi*(n-2*a))-2*sin(1/2*Pi*(6*a+
n)))*GAMMA(2*a+2/3-n)*GAMMA(a-1/3-1/2*n)*GAMMA(7/6-a-1/2*n)*GAMMA(a-1/
2*n)*GAMMA(2*a+1/3)*GAMMA(2*a-1/3)/sin(1/2*Pi*(-6*a+n))/Pi^(3/2)/(2^(2
*a+2/3-n))/GAMMA(3*a-1/2*n-1)/GAMMA(2*a+5/6-n)/GAMMA(2*a+1/6-n);}{%
\maplemultiline{
298\mbox{:~~~}  - {\displaystyle \frac {1}{2}} \mathrm{
V7}(a - {\displaystyle \frac {1}{3}}  - {\displaystyle \frac {n}{
2}} , \,n + 1)\,( - 6\,a + n + 2)\,\Gamma (a - {\displaystyle 
\frac {1}{3}}  - {\displaystyle \frac {n}{2}} )\,\Gamma (
{\displaystyle \frac {7}{6}}  - a - {\displaystyle \frac {n}{2}} 
)\,\Gamma (a + {\displaystyle \frac {1}{3}}  - {\displaystyle 
\frac {n}{2}} ) \\
\Gamma (a - {\displaystyle \frac {n}{2}} )\,\Gamma (2\,a + 
{\displaystyle \frac {1}{3}} )\,\Gamma (2\,a - {\displaystyle 
\frac {1}{3}} )\,\Gamma (2\,a) \left/ {\vrule 
height0.80em width0em depth0.80em} \right. \!  \! ((-1)^{( - n)}
\,\pi ^{(3/2)}\,\Gamma (3\,a - {\displaystyle \frac {n}{2}} )^{2}
 \\
\Gamma ({\displaystyle \frac {1}{2}}  + a - {\displaystyle 
\frac {n}{2}} )\,\Gamma (a + {\displaystyle \frac {1}{6}}  - 
{\displaystyle \frac {n}{2}} ))\mbox{} - {\displaystyle \frac {1
}{2}} \mathrm{V8}(a, \,n)\,(\mathrm{sin}({\displaystyle \frac {3
\,\pi \,(n - 2\,a)}{2}} ) - 2\,\mathrm{sin}({\displaystyle 
\frac {\pi \,(6\,a + n)}{2}} )) \\
\Gamma (2\,a + {\displaystyle \frac {2}{3}}  - n)\,\Gamma (a - 
{\displaystyle \frac {1}{3}}  - {\displaystyle \frac {n}{2}} )\,
\Gamma ({\displaystyle \frac {7}{6}}  - a - {\displaystyle 
\frac {n}{2}} )\,\Gamma (a - {\displaystyle \frac {n}{2}} )\,
\Gamma (2\,a + {\displaystyle \frac {1}{3}} )\,\Gamma (2\,a - 
{\displaystyle \frac {1}{3}} ) \left/ {\vrule 
height0.80em width0em depth0.80em} \right. \!  \! ( \\
\mathrm{sin}({\displaystyle \frac {\pi \,( - 6\,a + n)}{2}} )\,
\pi ^{(3/2)}\,2^{(2\,a + 2/3 - n)}\,\Gamma (3\,a - 
{\displaystyle \frac {n}{2}}  - 1)\,\Gamma (2\,a + 
{\displaystyle \frac {5}{6}}  - n)\,\Gamma (2\,a + 
{\displaystyle \frac {1}{6}}  - n)) }
}
\end{maplelatex}

\begin{maplelatex}
\mapleinline{inert}{2d}{299, ":   ",
V7(a-1/3-1/2*n,n+1)/Pi^(1/2)/GAMMA(3*a-1/2*n)/GAMMA(1/2+n)/GAMMA(1/2+a
-1/2*n)/GAMMA(2*a-1/2-n)*GAMMA(a-1/3-1/2*n)*GAMMA(a+1/3-1/2*n)*GAMMA(2
*a+1/3)*GAMMA(2/3)*GAMMA(2*a)+1/2*V8(a,n)*(-sin(1/6*Pi*(9*n-6*a+1))+si
n(1/6*Pi*(6*a-1+3*n)))*(-1)^(-n)*Pi^(1/2)*GAMMA(2*a+2/3-n)*GAMMA(a-1/3
-1/2*n)*GAMMA(2*a+1/3)*GAMMA(2/3)/cos(Pi*(n-2*a))/(2^(-1/3+2*a-n))/GAM
MA(2*a-1/2-n)/GAMMA(2*a+5/6-n)/GAMMA(-a+5/6+1/2*n)/GAMMA(2*a+1/6-n)/GA
MMA(1/2+n);}{%
\maplemultiline{
299\mbox{:~~~} {\displaystyle \frac {\mathrm{V7}(a - 
{\displaystyle \frac {1}{3}}  - {\displaystyle \frac {n}{2}} , \,
n + 1)\,\Gamma (a - {\displaystyle \frac {1}{3}}  - 
{\displaystyle \frac {n}{2}} )\,\Gamma (a + {\displaystyle 
\frac {1}{3}}  - {\displaystyle \frac {n}{2}} )\,\Gamma (2\,a + 
{\displaystyle \frac {1}{3}} )\,\Gamma ({\displaystyle \frac {2}{
3}} )\,\Gamma (2\,a)}{\sqrt{\pi }\,\Gamma (3\,a - {\displaystyle 
\frac {n}{2}} )\,\Gamma ({\displaystyle \frac {1}{2}}  + n)\,
\Gamma ({\displaystyle \frac {1}{2}}  + a - {\displaystyle 
\frac {n}{2}} )\,\Gamma (2\,a - {\displaystyle \frac {1}{2}}  - n
)}}  + {\displaystyle \frac {1}{2}}  \\
\mathrm{V8}(a, \,n)\,( - \mathrm{sin}({\displaystyle \frac {\pi 
\,(9\,n - 6\,a + 1)}{6}} ) + \mathrm{sin}({\displaystyle \frac {
\pi \,(6\,a - 1 + 3\,n)}{6}} ))\,(-1)^{( - n)}\,\sqrt{\pi } \\
\Gamma (2\,a + {\displaystyle \frac {2}{3}}  - n)\,\Gamma (a - 
{\displaystyle \frac {1}{3}}  - {\displaystyle \frac {n}{2}} )\,
\Gamma (2\,a + {\displaystyle \frac {1}{3}} )\,\Gamma (
{\displaystyle \frac {2}{3}} ) \left/ {\vrule 
height0.80em width0em depth0.80em} \right. \!  \! (\mathrm{cos}(
\pi \,(n - 2\,a))\,2^{( - 1/3 + 2\,a - n)} \\
\Gamma (2\,a - {\displaystyle \frac {1}{2}}  - n)\,\Gamma (2\,a
 + {\displaystyle \frac {5}{6}}  - n)\,\Gamma ( - a + 
{\displaystyle \frac {5}{6}}  + {\displaystyle \frac {n}{2}} )\,
\Gamma (2\,a + {\displaystyle \frac {1}{6}}  - n)\,\Gamma (
{\displaystyle \frac {1}{2}}  + n)) }
}
\end{maplelatex}

\begin{maplelatex}
\mapleinline{inert}{2d}{300, ":   ",
V7(a-1/3-1/2*n,n+1)/Pi^(1/2)/GAMMA(3*a-1/2*n)/GAMMA(1/2+n)/GAMMA(2*a-1
/2-n)*GAMMA(a-1/3-1/2*n)*GAMMA(a+1/2+1/2*n)*GAMMA(a+1/3-1/2*n)*GAMMA(2
*a)-1/2*V8(a,n)*(sin(1/6*Pi*(9*n-6*a+1))-sin(1/6*Pi*(6*a-1+3*n)))*(-1)
^(-n)*Pi^(1/2)*GAMMA(2*a+2/3-n)*GAMMA(1/2+a-1/2*n)*GAMMA(a-1/3-1/2*n)*
GAMMA(a+1/2+1/2*n)/cos(Pi*(n-2*a))/(2^(-1/3+2*a-n))/GAMMA(2*a-1/2-n)/G
AMMA(2*a+1/6-n)/GAMMA(2*a+5/6-n)/GAMMA(-a+5/6+1/2*n)/GAMMA(1/2+n);}{%
\maplemultiline{
300\mbox{:~~~} {\displaystyle \frac {\mathrm{V7}(a - 
{\displaystyle \frac {1}{3}}  - {\displaystyle \frac {n}{2}} , \,
n + 1)\,\Gamma (a - {\displaystyle \frac {1}{3}}  - 
{\displaystyle \frac {n}{2}} )\,\Gamma (a + {\displaystyle 
\frac {1}{2}}  + {\displaystyle \frac {n}{2}} )\,\Gamma (a + 
{\displaystyle \frac {1}{3}}  - {\displaystyle \frac {n}{2}} )\,
\Gamma (2\,a)}{\sqrt{\pi }\,\Gamma (3\,a - {\displaystyle \frac {
n}{2}} )\,\Gamma ({\displaystyle \frac {1}{2}}  + n)\,\Gamma (2\,
a - {\displaystyle \frac {1}{2}}  - n)}}  - {\displaystyle 
\frac {1}{2}}  \\
\mathrm{V8}(a, \,n)\,(\mathrm{sin}({\displaystyle \frac {\pi \,(9
\,n - 6\,a + 1)}{6}} ) - \mathrm{sin}({\displaystyle \frac {\pi 
\,(6\,a - 1 + 3\,n)}{6}} ))\,(-1)^{( - n)}\,\sqrt{\pi } \\
\Gamma (2\,a + {\displaystyle \frac {2}{3}}  - n)\,\Gamma (
{\displaystyle \frac {1}{2}}  + a - {\displaystyle \frac {n}{2}} 
)\,\Gamma (a - {\displaystyle \frac {1}{3}}  - {\displaystyle 
\frac {n}{2}} )\,\Gamma (a + {\displaystyle \frac {1}{2}}  + 
{\displaystyle \frac {n}{2}} ) \left/ {\vrule 
height0.80em width0em depth0.80em} \right. \!  \! (\mathrm{cos}(
\pi \,(n - 2\,a)) \\
2^{( - 1/3 + 2\,a - n)}\,\Gamma (2\,a - {\displaystyle \frac {1}{
2}}  - n)\,\Gamma (2\,a + {\displaystyle \frac {1}{6}}  - n)\,
\Gamma (2\,a + {\displaystyle \frac {5}{6}}  - n)\,\Gamma ( - a
 + {\displaystyle \frac {5}{6}}  + {\displaystyle \frac {n}{2}} )
\,\Gamma ({\displaystyle \frac {1}{2}}  + n)) }
}
\end{maplelatex}

\begin{maplelatex}
\mapleinline{inert}{2d}{301, ":   ",
V7(a-1/3-1/2*n,n+1)*2^(6*a-n-2)/((-1)^(-n))/Pi^2/GAMMA(6*a-n-1)*GAMMA(
2/3-a-1/2*n)*GAMMA(a+1/3-1/2*n)*GAMMA(2*a+1/3)*GAMMA(2*a-1/3)*GAMMA(2*
a)-1/4*V8(a,n)*(2*sin(1/6*Pi*(-12*a+6*n+1))+1)*GAMMA(2*a+2/3-n)*GAMMA(
1/2+a-1/2*n)*GAMMA(2/3-a-1/2*n)*GAMMA(2*a+1/3)*GAMMA(2*a-1/3)/cos(1/2*
Pi*(-6*a+n))/Pi^(1/2)/(2^(-1/3+2*a-n))/GAMMA(3*a-1/2*n-1/2)/GAMMA(2*a+
5/6-n)/GAMMA(-a+5/6+1/2*n)/GAMMA(2*a+1/6-n);}{%
\maplemultiline{
301\mbox{:~~~} {\displaystyle \frac {\mathrm{V7}(a - 
{\displaystyle \frac {1}{3}}  - {\displaystyle \frac {n}{2}} , \,
n + 1)\,2^{(6\,a - n - 2)}\,\Gamma ({\displaystyle \frac {2}{3}} 
 - a - {\displaystyle \frac {n}{2}} )\,\Gamma (a + 
{\displaystyle \frac {1}{3}}  - {\displaystyle \frac {n}{2}} )\,
\Gamma (2\,a + {\displaystyle \frac {1}{3}} )\,\Gamma (2\,a - 
{\displaystyle \frac {1}{3}} )\,\Gamma (2\,a)}{(-1)^{( - n)}\,\pi
 ^{2}\,\Gamma (6\,a - n - 1)}}  \\
\mbox{} - {\displaystyle \frac {1}{4}} \mathrm{V8}(a, \,n)\,(2\,
\mathrm{sin}({\displaystyle \frac {\pi \,( - 12\,a + 6\,n + 1)}{6
}} ) + 1)\,\Gamma (2\,a + {\displaystyle \frac {2}{3}}  - n)\,
\Gamma ({\displaystyle \frac {1}{2}}  + a - {\displaystyle 
\frac {n}{2}} ) \\
\Gamma ({\displaystyle \frac {2}{3}}  - a - {\displaystyle 
\frac {n}{2}} )\,\Gamma (2\,a + {\displaystyle \frac {1}{3}} )\,
\Gamma (2\,a - {\displaystyle \frac {1}{3}} ) \left/ {\vrule 
height0.80em width0em depth0.80em} \right. \!  \! (\mathrm{cos}(
{\displaystyle \frac {\pi \,( - 6\,a + n)}{2}} )\,\sqrt{\pi }\,2
^{( - 1/3 + 2\,a - n)} \\
\Gamma (3\,a - {\displaystyle \frac {n}{2}}  - {\displaystyle 
\frac {1}{2}} )\,\Gamma (2\,a + {\displaystyle \frac {5}{6}}  - n
)\,\Gamma ( - a + {\displaystyle \frac {5}{6}}  + {\displaystyle 
\frac {n}{2}} )\,\Gamma (2\,a + {\displaystyle \frac {1}{6}}  - n
)) }
}
\end{maplelatex}

\begin{maplelatex}
\mapleinline{inert}{2d}{302, ":   ",
1/2*V7(a-1/3-1/2*n,n+1)/Pi^(1/2)*2^(6*a-n-1)/GAMMA(6*a-n-1)*GAMMA(2*a+
1/3)*GAMMA(2*a-1/3)/GAMMA(1/2+n)*GAMMA(2*a)-1/4*V8(a,n)*(2*sin(1/6*Pi*
(-12*a+6*n+1))+1)*(-1)^(-n)*Pi^(1/2)*GAMMA(a-1/2*n+5/6)*GAMMA(1/2+a-1/
2*n)*GAMMA(2*a+1/3)*GAMMA(2*a-1/3)/cos(1/2*Pi*(-6*a+n))/GAMMA(3*a-1/2*
n-1/2)/GAMMA(2*a+1/6-n)/GAMMA(2*a+5/6-n)/GAMMA(-a+5/6+1/2*n)/GAMMA(1/2
+n);}{%
\maplemultiline{
302\mbox{:~~~} {\displaystyle \frac {1}{2}} \,
{\displaystyle \frac {\mathrm{V7}(a - {\displaystyle \frac {1}{3}
}  - {\displaystyle \frac {n}{2}} , \,n + 1)\,2^{(6\,a - n - 1)}
\,\Gamma (2\,a + {\displaystyle \frac {1}{3}} )\,\Gamma (2\,a - 
{\displaystyle \frac {1}{3}} )\,\Gamma (2\,a)}{\sqrt{\pi }\,
\Gamma (6\,a - n - 1)\,\Gamma ({\displaystyle \frac {1}{2}}  + n)
}}  - {\displaystyle \frac {1}{4}} \mathrm{V8}(a, \,n) \\
(2\,\mathrm{sin}({\displaystyle \frac {\pi \,( - 12\,a + 6\,n + 1
)}{6}} ) + 1)\,(-1)^{( - n)}\,\sqrt{\pi }\,\Gamma (a - 
{\displaystyle \frac {n}{2}}  + {\displaystyle \frac {5}{6}} )\,
\Gamma ({\displaystyle \frac {1}{2}}  + a - {\displaystyle 
\frac {n}{2}} )\,\Gamma (2\,a + {\displaystyle \frac {1}{3}} )
 \\
\Gamma (2\,a - {\displaystyle \frac {1}{3}} ) \left/ {\vrule 
height0.80em width0em depth0.80em} \right. \!  \! (\mathrm{cos}(
{\displaystyle \frac {\pi \,( - 6\,a + n)}{2}} )\,\Gamma (3\,a - 
{\displaystyle \frac {n}{2}}  - {\displaystyle \frac {1}{2}} )\,
\Gamma (2\,a + {\displaystyle \frac {1}{6}}  - n)\,\Gamma (2\,a
 + {\displaystyle \frac {5}{6}}  - n) \\
\Gamma ( - a + {\displaystyle \frac {5}{6}}  + {\displaystyle 
\frac {n}{2}} )\,\Gamma ({\displaystyle \frac {1}{2}}  + n)) }
}
\end{maplelatex}

\begin{maplelatex}
\mapleinline{inert}{2d}{303, ":   ",
GAMMA(2*a-1/3)*GAMMA(a-1/3-1/2*n)*GAMMA(a-1/2*n)/Pi^(1/2)*GAMMA(5/6+a+
1/2*n)/GAMMA(1/2+n)/GAMMA(3*a-1/2*n)/GAMMA(2*a-5/6-n)*V7(a-1/3-1/2*n,n
+1)+1/24*V8(a,n)*GAMMA(2*a-n)*GAMMA(2*a-2/3-n)*GAMMA(a-1/2*n+5/6)*GAMM
A(5/6+a+1/2*n)*GAMMA(2*a-1/3)*(2*sin(1/6*Pi*(-12*a-1+6*n))-1)*(-12*a+5
+6*n)/sin(1/6*Pi*(-12*a-1+6*n))/(2^(2*a-n-1))/(2^(2*a-5/3-n))/GAMMA(2*
a+1/6-n)^2/GAMMA(2*a+5/6-n)/GAMMA(1/2+n)/GAMMA(2*a);}{%
\maplemultiline{
303\mbox{:~~~} {\displaystyle \frac {\Gamma (2\,a - 
{\displaystyle \frac {1}{3}} )\,\Gamma (a - {\displaystyle 
\frac {1}{3}}  - {\displaystyle \frac {n}{2}} )\,\Gamma (a - 
{\displaystyle \frac {n}{2}} )\,\Gamma ({\displaystyle \frac {5}{
6}}  + a + {\displaystyle \frac {n}{2}} )\,\mathrm{V7}(a - 
{\displaystyle \frac {1}{3}}  - {\displaystyle \frac {n}{2}} , \,
n + 1)}{\sqrt{\pi }\,\Gamma ({\displaystyle \frac {1}{2}}  + n)\,
\Gamma (3\,a - {\displaystyle \frac {n}{2}} )\,\Gamma (2\,a - 
{\displaystyle \frac {5}{6}}  - n)}}  + {\displaystyle \frac {1}{
24}}  \\
\mathrm{V8}(a, \,n)\,\Gamma (2\,a - n)\,\Gamma (2\,a - 
{\displaystyle \frac {2}{3}}  - n)\,\Gamma (a - {\displaystyle 
\frac {n}{2}}  + {\displaystyle \frac {5}{6}} )\,\Gamma (
{\displaystyle \frac {5}{6}}  + a + {\displaystyle \frac {n}{2}} 
)\,\Gamma (2\,a - {\displaystyle \frac {1}{3}} ) \\
(2\,\mathrm{sin}({\displaystyle \frac {\pi \,( - 12\,a - 1 + 6\,n
)}{6}} ) - 1)\,( - 12\,a + 5 + 6\,n) \left/ {\vrule 
height0.80em width0em depth0.80em} \right. \!  \! (\mathrm{sin}(
{\displaystyle \frac {\pi \,( - 12\,a - 1 + 6\,n)}{6}} ) \\
2^{(2\,a - n - 1)}\,2^{(2\,a - 5/3 - n)}\,\Gamma (2\,a + 
{\displaystyle \frac {1}{6}}  - n)^{2}\,\Gamma (2\,a + 
{\displaystyle \frac {5}{6}}  - n)\,\Gamma ({\displaystyle 
\frac {1}{2}}  + n)\,\Gamma (2\,a)) }
}
\end{maplelatex}

\begin{maplelatex}
\mapleinline{inert}{2d}{304, ":   ",
GAMMA(2*a+1/3)*GAMMA(1/6+a+1/2*n)*GAMMA(a+1/3-1/2*n)*GAMMA(a-1/2*n)/Pi
^(1/2)/GAMMA(1/2+n)/GAMMA(3*a-1/2*n)/GAMMA(2*a-n-1/6)*V7(a-1/3-1/2*n,n
+1)+1/24*V8(a,n)*GAMMA(2*a+2/3-n)*GAMMA(2*a-n)*GAMMA(2*a+1/3)*GAMMA(a+
1/6-1/2*n)*GAMMA(1/6+a+1/2*n)*(2*sin(1/6*Pi*(-12*a+6*n+1))+1)*(-12*a+6
*n+1)/sin(1/6*Pi*(-12*a+6*n+1))/(2^(-1/3+2*a-n))/(2^(2*a-n-1))/GAMMA(2
*a+1/6-n)/GAMMA(2*a+5/6-n)^2/GAMMA(1/2+n)/GAMMA(2*a);}{%
\maplemultiline{
304\mbox{:~~~} {\displaystyle \frac {\Gamma (2\,a + 
{\displaystyle \frac {1}{3}} )\,\Gamma ({\displaystyle \frac {1}{
6}}  + a + {\displaystyle \frac {n}{2}} )\,\Gamma (a + 
{\displaystyle \frac {1}{3}}  - {\displaystyle \frac {n}{2}} )\,
\Gamma (a - {\displaystyle \frac {n}{2}} )\,\mathrm{V7}(a - 
{\displaystyle \frac {1}{3}}  - {\displaystyle \frac {n}{2}} , \,
n + 1)}{\sqrt{\pi }\,\Gamma ({\displaystyle \frac {1}{2}}  + n)\,
\Gamma (3\,a - {\displaystyle \frac {n}{2}} )\,\Gamma (2\,a - n
 - {\displaystyle \frac {1}{6}} )}}  + {\displaystyle \frac {1}{
24}}  \\
\mathrm{V8}(a, \,n)\,\Gamma (2\,a + {\displaystyle \frac {2}{3}} 
 - n)\,\Gamma (2\,a - n)\,\Gamma (2\,a + {\displaystyle \frac {1
}{3}} )\,\Gamma (a + {\displaystyle \frac {1}{6}}  - 
{\displaystyle \frac {n}{2}} )\,\Gamma ({\displaystyle \frac {1}{
6}}  + a + {\displaystyle \frac {n}{2}} ) \\
(2\,\mathrm{sin}({\displaystyle \frac {\pi \,( - 12\,a + 6\,n + 1
)}{6}} ) + 1)\,( - 12\,a + 6\,n + 1) \left/ {\vrule 
height0.80em width0em depth0.80em} \right. \!  \! (\mathrm{sin}(
{\displaystyle \frac {\pi \,( - 12\,a + 6\,n + 1)}{6}} ) \\
2^{( - 1/3 + 2\,a - n)}\,2^{(2\,a - n - 1)}\,\Gamma (2\,a + 
{\displaystyle \frac {1}{6}}  - n)\,\Gamma (2\,a + 
{\displaystyle \frac {5}{6}}  - n)^{2}\,\Gamma ({\displaystyle 
\frac {1}{2}}  + n)\,\Gamma (2\,a)) }
}
\end{maplelatex}

\begin{maplelatex}
\mapleinline{inert}{2d}{305, ":   ",
-cos(1/2*Pi*(-6*a+n))*V7(a-1/3-1/2*n,n+1)/((-1)^(-n))/Pi^(5/2)*GAMMA(2
/3-a-1/2*n)*GAMMA(4/3-a-1/2*n)*GAMMA(1/2-a+1/2*n)/GAMMA(3*a-1/2*n)*GAM
MA(2*a+1/3)*GAMMA(2*a-1/3)*GAMMA(2*a)+1/4*(2*sin(1/6*Pi*(-12*a+6*n+1))
+1)*V8(a,n)*GAMMA(a-1/2*n+5/6)*GAMMA(2/3-a-1/2*n)*GAMMA(4/3-a-1/2*n)*G
AMMA(2*a+1/3)*GAMMA(2*a-1/3)/cos(1/2*Pi*(n-2*a))/Pi/GAMMA(2*a+5/6-n)/G
AMMA(2*a+1/6-n)/GAMMA(-a+5/6+1/2*n);}{%
\maplemultiline{
305\mbox{:~~~}  - \mathrm{cos}({\displaystyle \frac {\pi
 \,( - 6\,a + n)}{2}} )\,\mathrm{V7}(a - {\displaystyle \frac {1
}{3}}  - {\displaystyle \frac {n}{2}} , \,n + 1)\,\Gamma (
{\displaystyle \frac {2}{3}}  - a - {\displaystyle \frac {n}{2}} 
)\,\Gamma ({\displaystyle \frac {4}{3}}  - a - {\displaystyle 
\frac {n}{2}} ) \\
\Gamma ({\displaystyle \frac {1}{2}}  - a + {\displaystyle 
\frac {n}{2}} )\,\Gamma (2\,a + {\displaystyle \frac {1}{3}} )\,
\Gamma (2\,a - {\displaystyle \frac {1}{3}} )\,\Gamma (2\,a)
 \left/ {\vrule height0.80em width0em depth0.80em} \right. \! 
 \! ((-1)^{( - n)}\,\pi ^{(5/2)}\,\Gamma (3\,a - {\displaystyle 
\frac {n}{2}} ))\mbox{} + {\displaystyle \frac {1}{4}}  \\
(2\,\mathrm{sin}({\displaystyle \frac {\pi \,( - 12\,a + 6\,n + 1
)}{6}} ) + 1)\,\mathrm{V8}(a, \,n)\,\Gamma (a - {\displaystyle 
\frac {n}{2}}  + {\displaystyle \frac {5}{6}} )\,\Gamma (
{\displaystyle \frac {2}{3}}  - a - {\displaystyle \frac {n}{2}} 
)\,\Gamma ({\displaystyle \frac {4}{3}}  - a - {\displaystyle 
\frac {n}{2}} ) \\
\Gamma (2\,a + {\displaystyle \frac {1}{3}} )\,\Gamma (2\,a - 
{\displaystyle \frac {1}{3}} ) \left/ {\vrule 
height0.80em width0em depth0.80em} \right. \!  \! (\mathrm{cos}(
{\displaystyle \frac {\pi \,(n - 2\,a)}{2}} )\,\pi \,\Gamma (2\,a
 + {\displaystyle \frac {5}{6}}  - n)\,\Gamma (2\,a + 
{\displaystyle \frac {1}{6}}  - n) \\
\Gamma ( - a + {\displaystyle \frac {5}{6}}  + {\displaystyle 
\frac {n}{2}} )) }
}
\end{maplelatex}

\begin{maplelatex}
\mapleinline{inert}{2d}{306, ":   ",
-cos(1/2*Pi*(-6*a+n))*V7(a-1/3-1/2*n,n+1)/Pi/sin(1/3*Pi*(6*a+1))/((-1)
^(-n))/GAMMA(3*a-1/2*n)*GAMMA(4/3-a-1/2*n)*GAMMA(2*a-1/3)*GAMMA(2*a)+1
/4*(2*sin(1/6*Pi*(-12*a+6*n+1))+1)*V8(a,n)*GAMMA(4/3-a-1/2*n)*GAMMA(2*
a-1/3)*GAMMA(1/2+a-1/2*n)*GAMMA(a-1/2*n+5/6)/sin(1/3*Pi*(6*a+1))/Pi^(1
/2)/GAMMA(2*a+5/6-n)/GAMMA(2*a+1/6-n)/GAMMA(-a+5/6+1/2*n);}{%
\maplemultiline{
306\mbox{:~~~}  - {\displaystyle \frac {\mathrm{cos}(
{\displaystyle \frac {\pi \,( - 6\,a + n)}{2}} )\,\mathrm{V7}(a
 - {\displaystyle \frac {1}{3}}  - {\displaystyle \frac {n}{2}} 
, \,n + 1)\,\Gamma ({\displaystyle \frac {4}{3}}  - a - 
{\displaystyle \frac {n}{2}} )\,\Gamma (2\,a - {\displaystyle 
\frac {1}{3}} )\,\Gamma (2\,a)}{\pi \,\mathrm{sin}(
{\displaystyle \frac {\pi \,(6\,a + 1)}{3}} )\,(-1)^{( - n)}\,
\Gamma (3\,a - {\displaystyle \frac {n}{2}} )}}  + 
{\displaystyle \frac {1}{4}}  \\
(2\,\mathrm{sin}({\displaystyle \frac {\pi \,( - 12\,a + 6\,n + 1
)}{6}} ) + 1)\,\mathrm{V8}(a, \,n)\,\Gamma ({\displaystyle 
\frac {4}{3}}  - a - {\displaystyle \frac {n}{2}} )\,\Gamma (2\,a
 - {\displaystyle \frac {1}{3}} )\,\Gamma ({\displaystyle \frac {
1}{2}}  + a - {\displaystyle \frac {n}{2}} ) \\
\Gamma (a - {\displaystyle \frac {n}{2}}  + {\displaystyle 
\frac {5}{6}} ) \left/ {\vrule height0.80em width0em depth0.80em}
 \right. \!  \! (\mathrm{sin}({\displaystyle \frac {\pi \,(6\,a
 + 1)}{3}} )\,\sqrt{\pi }\,\Gamma (2\,a + {\displaystyle \frac {5
}{6}}  - n)\,\Gamma (2\,a + {\displaystyle \frac {1}{6}}  - n)
 \\
\Gamma ( - a + {\displaystyle \frac {5}{6}}  + {\displaystyle 
\frac {n}{2}} )) }
}
\end{maplelatex}

\begin{maplelatex}
\mapleinline{inert}{2d}{307, ":   ",
V7(a-1/3-1/2*n,n+1)*2^(6*a-n-2)/((-1)^(-n))/Pi^2/GAMMA(6*a-n-1)*GAMMA(
2*a+1/3)*GAMMA(2*a-1/3)*GAMMA(4/3-a-1/2*n)*GAMMA(a-1/3-1/2*n)*GAMMA(2*
a)-1/4*V8(a,n)*(2*sin(1/6*Pi*(-12*a+6*n+1))+1)*GAMMA(1/2+a-1/2*n)*GAMM
A(a-1/2*n+5/6)*GAMMA(4/3-a-1/2*n)*GAMMA(a-1/3-1/2*n)*GAMMA(2*a+1/3)*GA
MMA(2*a-1/3)/cos(1/2*Pi*(-6*a+n))/Pi/GAMMA(3*a-1/2*n-1/2)/GAMMA(2*a+5/
6-n)/GAMMA(2*a+1/6-n)/GAMMA(-a+5/6+1/2*n);}{%
\maplemultiline{
307\mbox{:~~~} {\displaystyle \frac {\mathrm{V7}(a - 
{\displaystyle \frac {1}{3}}  - {\displaystyle \frac {n}{2}} , \,
n + 1)\,2^{(6\,a - n - 2)}\,\Gamma (2\,a + {\displaystyle \frac {
1}{3}} )\,\Gamma (2\,a - {\displaystyle \frac {1}{3}} )\,\Gamma (
{\displaystyle \frac {4}{3}}  - a - {\displaystyle \frac {n}{2}} 
)\,\Gamma (a - {\displaystyle \frac {1}{3}}  - {\displaystyle 
\frac {n}{2}} )\,\Gamma (2\,a)}{(-1)^{( - n)}\,\pi ^{2}\,\Gamma (
6\,a - n - 1)}}  \\
\mbox{} - {\displaystyle \frac {1}{4}} \mathrm{V8}(a, \,n)\,(2\,
\mathrm{sin}({\displaystyle \frac {\pi \,( - 12\,a + 6\,n + 1)}{6
}} ) + 1)\,\Gamma ({\displaystyle \frac {1}{2}}  + a - 
{\displaystyle \frac {n}{2}} )\,\Gamma (a - {\displaystyle 
\frac {n}{2}}  + {\displaystyle \frac {5}{6}} ) \\
\Gamma ({\displaystyle \frac {4}{3}}  - a - {\displaystyle 
\frac {n}{2}} )\,\Gamma (a - {\displaystyle \frac {1}{3}}  - 
{\displaystyle \frac {n}{2}} )\,\Gamma (2\,a + {\displaystyle 
\frac {1}{3}} )\,\Gamma (2\,a - {\displaystyle \frac {1}{3}} )
 \left/ {\vrule height0.80em width0em depth0.80em} \right. \! 
 \! (\mathrm{cos}({\displaystyle \frac {\pi \,( - 6\,a + n)}{2}} 
)\,\pi  \\
\Gamma (3\,a - {\displaystyle \frac {n}{2}}  - {\displaystyle 
\frac {1}{2}} )\,\Gamma (2\,a + {\displaystyle \frac {5}{6}}  - n
)\,\Gamma (2\,a + {\displaystyle \frac {1}{6}}  - n)\,\Gamma ( - 
a + {\displaystyle \frac {5}{6}}  + {\displaystyle \frac {n}{2}} 
)) }
}
\end{maplelatex}

\begin{maplelatex}
\mapleinline{inert}{2d}{308, ":   ",
-cos(1/2*Pi*(-6*a+n))*V7(a-1/3-1/2*n,n+1)/((-1)^(-n))/Pi^(5/2)*GAMMA(4
/3-a-1/2*n)*GAMMA(1/6-a+1/2*n)/GAMMA(3*a-1/2*n)*GAMMA(1-1/2*n-a)*GAMMA
(2*a+1/3)*GAMMA(2*a-1/3)*GAMMA(2*a)+1/4*(2*sin(1/6*Pi*(-12*a+6*n+1))+1
)*V8(a,n)*GAMMA(1/2+a-1/2*n)*GAMMA(4/3-a-1/2*n)*GAMMA(1-1/2*n-a)*GAMMA
(2*a+1/3)*GAMMA(2*a-1/3)/sin(1/6*Pi*(1-6*a+3*n))/Pi/GAMMA(2*a+5/6-n)/G
AMMA(2*a+1/6-n)/GAMMA(-a+5/6+1/2*n);}{%
\maplemultiline{
308\mbox{:~~~}  - \mathrm{cos}({\displaystyle \frac {\pi
 \,( - 6\,a + n)}{2}} )\,\mathrm{V7}(a - {\displaystyle \frac {1
}{3}}  - {\displaystyle \frac {n}{2}} , \,n + 1)\,\Gamma (
{\displaystyle \frac {4}{3}}  - a - {\displaystyle \frac {n}{2}} 
)\,\Gamma ({\displaystyle \frac {1}{6}}  - a + {\displaystyle 
\frac {n}{2}} ) \\
\Gamma (1 - {\displaystyle \frac {n}{2}}  - a)\,\Gamma (2\,a + 
{\displaystyle \frac {1}{3}} )\,\Gamma (2\,a - {\displaystyle 
\frac {1}{3}} )\,\Gamma (2\,a) \left/ {\vrule 
height0.80em width0em depth0.80em} \right. \!  \! ((-1)^{( - n)}
\,\pi ^{(5/2)}\,\Gamma (3\,a - {\displaystyle \frac {n}{2}} ))
\mbox{} + {\displaystyle \frac {1}{4}}  \\
(2\,\mathrm{sin}({\displaystyle \frac {\pi \,( - 12\,a + 6\,n + 1
)}{6}} ) + 1)\,\mathrm{V8}(a, \,n)\,\Gamma ({\displaystyle 
\frac {1}{2}}  + a - {\displaystyle \frac {n}{2}} )\,\Gamma (
{\displaystyle \frac {4}{3}}  - a - {\displaystyle \frac {n}{2}} 
)\,\Gamma (1 - {\displaystyle \frac {n}{2}}  - a) \\
\Gamma (2\,a + {\displaystyle \frac {1}{3}} )\,\Gamma (2\,a - 
{\displaystyle \frac {1}{3}} ) \left/ {\vrule 
height0.80em width0em depth0.80em} \right. \!  \! (\mathrm{sin}(
{\displaystyle \frac {\pi \,(1 - 6\,a + 3\,n)}{6}} )\,\pi \,
\Gamma (2\,a + {\displaystyle \frac {5}{6}}  - n)\,\Gamma (2\,a
 + {\displaystyle \frac {1}{6}}  - n) \\
\Gamma ( - a + {\displaystyle \frac {5}{6}}  + {\displaystyle 
\frac {n}{2}} )) }
}
\end{maplelatex}

\begin{maplelatex}
\mapleinline{inert}{2d}{309, ":   ",
-cos(1/2*Pi*(-6*a+n))*V7(a-1/3-1/2*n,n+1)/Pi/sin(1/3*Pi*(6*a-1))/((-1)
^(-n))*GAMMA(2/3-a-1/2*n)/GAMMA(3*a-1/2*n)*GAMMA(2*a+1/3)*GAMMA(2*a)+1
/4*(2*sin(1/6*Pi*(-12*a+6*n+1))+1)*V8(a,n)*GAMMA(a-1/2*n+5/6)*GAMMA(2/
3-a-1/2*n)*GAMMA(1/2+a-1/2*n)*GAMMA(2*a+1/3)/sin(1/3*Pi*(6*a-1))/Pi^(1
/2)/GAMMA(2*a+5/6-n)/GAMMA(-a+5/6+1/2*n)/GAMMA(2*a+1/6-n);}{%
\maplemultiline{
309\mbox{:~~~}  - {\displaystyle \frac {\mathrm{cos}(
{\displaystyle \frac {\pi \,( - 6\,a + n)}{2}} )\,\mathrm{V7}(a
 - {\displaystyle \frac {1}{3}}  - {\displaystyle \frac {n}{2}} 
, \,n + 1)\,\Gamma ({\displaystyle \frac {2}{3}}  - a - 
{\displaystyle \frac {n}{2}} )\,\Gamma (2\,a + {\displaystyle 
\frac {1}{3}} )\,\Gamma (2\,a)}{\pi \,\mathrm{sin}(
{\displaystyle \frac {\pi \,(6\,a - 1)}{3}} )\,(-1)^{( - n)}\,
\Gamma (3\,a - {\displaystyle \frac {n}{2}} )}}  + 
{\displaystyle \frac {1}{4}}  \\
(2\,\mathrm{sin}({\displaystyle \frac {\pi \,( - 12\,a + 6\,n + 1
)}{6}} ) + 1)\,\mathrm{V8}(a, \,n)\,\Gamma (a - {\displaystyle 
\frac {n}{2}}  + {\displaystyle \frac {5}{6}} )\,\Gamma (
{\displaystyle \frac {2}{3}}  - a - {\displaystyle \frac {n}{2}} 
)\,\Gamma ({\displaystyle \frac {1}{2}}  + a - {\displaystyle 
\frac {n}{2}} ) \\
\Gamma (2\,a + {\displaystyle \frac {1}{3}} ) \left/ {\vrule 
height0.80em width0em depth0.80em} \right. \!  \! (\mathrm{sin}(
{\displaystyle \frac {\pi \,(6\,a - 1)}{3}} )\,\sqrt{\pi }\,
\Gamma (2\,a + {\displaystyle \frac {5}{6}}  - n)\,\Gamma ( - a
 + {\displaystyle \frac {5}{6}}  + {\displaystyle \frac {n}{2}} )
\,\Gamma (2\,a + {\displaystyle \frac {1}{6}}  - n)
) }
}
\end{maplelatex}

\begin{maplelatex}
\mapleinline{inert}{2d}{310, ":   ",
-cos(1/2*Pi*(-6*a+n))*V7(a-1/3-1/2*n,n+1)/((-1)^(-n))/Pi^(5/2)*GAMMA(2
/3-a-1/2*n)*GAMMA(-a+5/6+1/2*n)/GAMMA(3*a-1/2*n)*GAMMA(1-1/2*n-a)*GAMM
A(2*a+1/3)*GAMMA(2*a-1/3)*GAMMA(2*a)+1/4*(2*sin(1/6*Pi*(-12*a+6*n+1))+
1)*V8(a,n)*GAMMA(2/3-a-1/2*n)*GAMMA(a-1/2*n+5/6)*GAMMA(1/2+a-1/2*n)*GA
MMA(1-1/2*n-a)*GAMMA(2*a+1/3)*GAMMA(2*a-1/3)/Pi^2/GAMMA(2*a+1/6-n)/GAM
MA(2*a+5/6-n);}{%
\maplemultiline{
310\mbox{:~~~}  - \mathrm{cos}({\displaystyle \frac {\pi
 \,( - 6\,a + n)}{2}} )\,\mathrm{V7}(a - {\displaystyle \frac {1
}{3}}  - {\displaystyle \frac {n}{2}} , \,n + 1)\,\Gamma (
{\displaystyle \frac {2}{3}}  - a - {\displaystyle \frac {n}{2}} 
)\,\Gamma ( - a + {\displaystyle \frac {5}{6}}  + {\displaystyle 
\frac {n}{2}} ) \\
\Gamma (1 - {\displaystyle \frac {n}{2}}  - a)\,\Gamma (2\,a + 
{\displaystyle \frac {1}{3}} )\,\Gamma (2\,a - {\displaystyle 
\frac {1}{3}} )\,\Gamma (2\,a) \left/ {\vrule 
height0.80em width0em depth0.80em} \right. \!  \! ((-1)^{( - n)}
\,\pi ^{(5/2)}\,\Gamma (3\,a - {\displaystyle \frac {n}{2}} ))
\mbox{} + {\displaystyle \frac {1}{4}}  \\
(2\,\mathrm{sin}({\displaystyle \frac {\pi \,( - 12\,a + 6\,n + 1
)}{6}} ) + 1)\,\mathrm{V8}(a, \,n)\,\Gamma ({\displaystyle 
\frac {2}{3}}  - a - {\displaystyle \frac {n}{2}} )\,\Gamma (a - 
{\displaystyle \frac {n}{2}}  + {\displaystyle \frac {5}{6}} )\,
\Gamma ({\displaystyle \frac {1}{2}}  + a - {\displaystyle 
\frac {n}{2}} ) \\
\Gamma (1 - {\displaystyle \frac {n}{2}}  - a)\,\Gamma (2\,a + 
{\displaystyle \frac {1}{3}} )\,\Gamma (2\,a - {\displaystyle 
\frac {1}{3}} ) \left/ {\vrule height0.80em width0em depth0.80em}
 \right. \!  \! (\pi ^{2}\,\Gamma (2\,a + {\displaystyle \frac {1
}{6}}  - n)\,\Gamma (2\,a + {\displaystyle \frac {5}{6}}  - n))
 }
}
\end{maplelatex}

\begin{maplelatex}
\mapleinline{inert}{2d}{311, ":   ",
-cos(1/2*Pi*(-6*a+n))*V7(a-1/3-1/2*n,n+1)/Pi/sin(2*Pi*a)/((-1)^(-n))/G
AMMA(3*a-1/2*n)*GAMMA(1-1/2*n-a)*GAMMA(2*a+1/3)*GAMMA(2*a-1/3)+1/4*(2*
sin(1/6*Pi*(-12*a+6*n+1))+1)*V8(a,n)*GAMMA(1/2+a-1/2*n)*GAMMA(a-1/2*n+
5/6)*GAMMA(1-1/2*n-a)*GAMMA(2*a+1/3)*GAMMA(2*a-1/3)/sin(2*Pi*a)/Pi^(1/
2)/GAMMA(2*a+5/6-n)/GAMMA(2*a+1/6-n)/GAMMA(-a+5/6+1/2*n)/GAMMA(2*a);}{
\maplemultiline{
311\mbox{:~~~}  - {\displaystyle \frac {\mathrm{cos}(
{\displaystyle \frac {\pi \,( - 6\,a + n)}{2}} )\,\mathrm{V7}(a
 - {\displaystyle \frac {1}{3}}  - {\displaystyle \frac {n}{2}} 
, \,n + 1)\,\Gamma (1 - {\displaystyle \frac {n}{2}}  - a)\,
\Gamma (2\,a + {\displaystyle \frac {1}{3}} )\,\Gamma (2\,a - 
{\displaystyle \frac {1}{3}} )}{\pi \,\mathrm{sin}(2\,\pi \,a)\,(
-1)^{( - n)}\,\Gamma (3\,a - {\displaystyle \frac {n}{2}} )}}  + 
 \\
{\displaystyle \frac {1}{4}} (2\,\mathrm{sin}({\displaystyle 
\frac {\pi \,( - 12\,a + 6\,n + 1)}{6}} ) + 1)\,\mathrm{V8}(a, \,
n)\,\Gamma ({\displaystyle \frac {1}{2}}  + a - {\displaystyle 
\frac {n}{2}} )\,\Gamma (a - {\displaystyle \frac {n}{2}}  + 
{\displaystyle \frac {5}{6}} ) \\
\Gamma (1 - {\displaystyle \frac {n}{2}}  - a)\,\Gamma (2\,a + 
{\displaystyle \frac {1}{3}} )\,\Gamma (2\,a - {\displaystyle 
\frac {1}{3}} ) \left/ {\vrule height0.80em width0em depth0.80em}
 \right. \!  \! (\mathrm{sin}(2\,\pi \,a)\,\sqrt{\pi }\,\Gamma (2
\,a + {\displaystyle \frac {5}{6}}  - n) \\
\Gamma (2\,a + {\displaystyle \frac {1}{6}}  - n)\,\Gamma ( - a
 + {\displaystyle \frac {5}{6}}  + {\displaystyle \frac {n}{2}} )
\,\Gamma (2\,a)) }
}
\end{maplelatex}

\begin{maplelatex}
\mapleinline{inert}{2d}{312, ":   ",
-2/3*sin(1/6*Pi*(-12*a+6*n+1))*3^(1/2)*V7(a-1/3-1/2*n,n+1)/Pi^(1/2)*GA
MMA(2/3+a+1/2*n)/GAMMA(a+1/6-1/2*n)/GAMMA(3*a-1/2*n)*GAMMA(1/2-a+1/2*n
)/GAMMA(1/2+n)*GAMMA(a-1/2*n)*GAMMA(a+1/3-1/2*n)*GAMMA(2*a+1/3)/GAMMA(
2/3)+1/6*V8(a,n)*(2*sin(1/6*Pi*(-12*a+6*n+1))+1)*3^(1/2)*Pi^(1/2)*GAMM
A(2*a+2/3-n)*GAMMA(2/3+a+1/2*n)*GAMMA(a-1/2*n)*GAMMA(2*a+1/3)/cos(1/2*
Pi*(n-2*a))/(2^(-1/3+2*a-n))/GAMMA(2*a+1/6-n)/GAMMA(2*a+5/6-n)/GAMMA(1
/2+n)/GAMMA(2/3)/GAMMA(2*a);}{%
\maplemultiline{
312\mbox{:~~~}  - {\displaystyle \frac {2}{3}} \mathrm{
sin}({\displaystyle \frac {\pi \,( - 12\,a + 6\,n + 1)}{6}} )\,
\sqrt{3}\,\mathrm{V7}(a - {\displaystyle \frac {1}{3}}  - 
{\displaystyle \frac {n}{2}} , \,n + 1)\,\Gamma ({\displaystyle 
\frac {2}{3}}  + a + {\displaystyle \frac {n}{2}} )\,\Gamma (
{\displaystyle \frac {1}{2}}  - a + {\displaystyle \frac {n}{2}} 
) \\
\Gamma (a - {\displaystyle \frac {n}{2}} )\,\Gamma (a + 
{\displaystyle \frac {1}{3}}  - {\displaystyle \frac {n}{2}} )\,
\Gamma (2\,a + {\displaystyle \frac {1}{3}} ) \left/ {\vrule 
height0.80em width0em depth0.80em} \right. \!  \! (\sqrt{\pi }\,
\Gamma (a + {\displaystyle \frac {1}{6}}  - {\displaystyle 
\frac {n}{2}} )\,\Gamma (3\,a - {\displaystyle \frac {n}{2}} )\,
\Gamma ({\displaystyle \frac {1}{2}}  + n) \\
\Gamma ({\displaystyle \frac {2}{3}} ))\mbox{} + {\displaystyle 
\frac {1}{6}} \mathrm{V8}(a, \,n)\,(2\,\mathrm{sin}(
{\displaystyle \frac {\pi \,( - 12\,a + 6\,n + 1)}{6}} ) + 1)\,
\sqrt{3}\,\sqrt{\pi }\,\Gamma (2\,a + {\displaystyle \frac {2}{3}
}  - n) \\
\Gamma ({\displaystyle \frac {2}{3}}  + a + {\displaystyle 
\frac {n}{2}} )\,\Gamma (a - {\displaystyle \frac {n}{2}} )\,
\Gamma (2\,a + {\displaystyle \frac {1}{3}} ) \left/ {\vrule 
height0.80em width0em depth0.80em} \right. \!  \! (\mathrm{cos}(
{\displaystyle \frac {\pi \,(n - 2\,a)}{2}} )\,2^{( - 1/3 + 2\,a
 - n)} \\
\Gamma (2\,a + {\displaystyle \frac {1}{6}}  - n)\,\Gamma (2\,a
 + {\displaystyle \frac {5}{6}}  - n)\,\Gamma ({\displaystyle 
\frac {1}{2}}  + n)\,\Gamma ({\displaystyle \frac {2}{3}} )\,
\Gamma (2\,a)) }
}
\end{maplelatex}

\begin{maplelatex}
\mapleinline{inert}{2d}{313, ":   ",
2/3*Pi^(1/2)*3^(1/2)*sin(1/6*Pi*(-12*a+6*n+1))*V7(a-1/3-1/2*n,n+1)/sin
(1/2*Pi*(n-2*a))*GAMMA(1/6+a+1/2*n)/GAMMA(3*a-1/2*n)/GAMMA(1/2+n)/GAMM
A(a+1/6-1/2*n)*GAMMA(a+1/3-1/2*n)*GAMMA(2*a+1/3)/GAMMA(2/3)-1/6*3^(1/2
)*(2*sin(1/6*Pi*(-12*a+6*n+1))+1)*V8(a,n)*Pi^(1/2)*GAMMA(1/6+a+1/2*n)*
GAMMA(1/2+a-1/2*n)*GAMMA(2*a+2/3-n)*GAMMA(2*a+1/3)/sin(1/2*Pi*(n-2*a))
/(2^(-1/3+2*a-n))/GAMMA(2*a+1/6-n)/GAMMA(2*a+5/6-n)/GAMMA(1/2+n)/GAMMA
(2/3)/GAMMA(2*a);}{%
\maplemultiline{
313\mbox{:~~~} {\displaystyle \frac {2}{3}} \sqrt{\pi }
\,\sqrt{3}\,\mathrm{sin}({\displaystyle \frac {\pi \,( - 12\,a + 
6\,n + 1)}{6}} )\,\mathrm{V7}(a - {\displaystyle \frac {1}{3}} 
 - {\displaystyle \frac {n}{2}} , \,n + 1)\,\Gamma (
{\displaystyle \frac {1}{6}}  + a + {\displaystyle \frac {n}{2}} 
) \\
\Gamma (a + {\displaystyle \frac {1}{3}}  - {\displaystyle 
\frac {n}{2}} )\,\Gamma (2\,a + {\displaystyle \frac {1}{3}} )
 \left/ {\vrule height0.80em width0em depth0.80em} \right. \! 
 \! (\mathrm{sin}({\displaystyle \frac {\pi \,(n - 2\,a)}{2}} )\,
\Gamma (3\,a - {\displaystyle \frac {n}{2}} )\,\Gamma (
{\displaystyle \frac {1}{2}}  + n)\,\Gamma (a + {\displaystyle 
\frac {1}{6}}  - {\displaystyle \frac {n}{2}} ) \\
\Gamma ({\displaystyle \frac {2}{3}} ))\mbox{} - {\displaystyle 
\frac {1}{6}} \sqrt{3}\,(2\,\mathrm{sin}({\displaystyle \frac {
\pi \,( - 12\,a + 6\,n + 1)}{6}} ) + 1)\,\mathrm{V8}(a, \,n)\,
\sqrt{\pi }\,\Gamma ({\displaystyle \frac {1}{6}}  + a + 
{\displaystyle \frac {n}{2}} ) \\
\Gamma ({\displaystyle \frac {1}{2}}  + a - {\displaystyle 
\frac {n}{2}} )\,\Gamma (2\,a + {\displaystyle \frac {2}{3}}  - n
)\,\Gamma (2\,a + {\displaystyle \frac {1}{3}} ) \left/ {\vrule 
height0.80em width0em depth0.80em} \right. \!  \! (\mathrm{sin}(
{\displaystyle \frac {\pi \,(n - 2\,a)}{2}} )\,2^{( - 1/3 + 2\,a
 - n)} \\
\Gamma (2\,a + {\displaystyle \frac {1}{6}}  - n)\,\Gamma (2\,a
 + {\displaystyle \frac {5}{6}}  - n)\,\Gamma ({\displaystyle 
\frac {1}{2}}  + n)\,\Gamma ({\displaystyle \frac {2}{3}} )\,
\Gamma (2\,a)) }
}
\end{maplelatex}

\begin{maplelatex}
\mapleinline{inert}{2d}{314, ":   ",
-2/3*sin(1/6*Pi*(-12*a+6*n+1))*3^(1/2)*V7(a-1/3-1/2*n,n+1)/Pi^(1/2)/GA
MMA(3*a-1/2*n)*GAMMA(a+1/3-1/2*n)/GAMMA(a+1/6-1/2*n)*GAMMA(a-1/2*n)*GA
MMA(2*a+1/3)+1/3*3^(1/2)*(2*sin(1/6*Pi*(-12*a+6*n+1))+1)*V8(a,n)*GAMMA
(2*a+2/3-n)*GAMMA(2*a-n)*GAMMA(2*a+1/3)/(2^(-1/3+2*a-n))/(2^(2*a-n))/G
AMMA(2*a+5/6-n)/GAMMA(2*a+1/6-n)/GAMMA(2*a);}{%
\maplemultiline{
314\mbox{:~~~}  - {\displaystyle \frac {2}{3}} \,
{\displaystyle \frac {\mathrm{sin}({\displaystyle \frac {\pi \,(
 - 12\,a + 6\,n + 1)}{6}} )\,\sqrt{3}\,\mathrm{V7}(a - 
{\displaystyle \frac {1}{3}}  - {\displaystyle \frac {n}{2}} , \,
n + 1)\,\Gamma (a + {\displaystyle \frac {1}{3}}  - 
{\displaystyle \frac {n}{2}} )\,\Gamma (a - {\displaystyle 
\frac {n}{2}} )\,\Gamma (2\,a + {\displaystyle \frac {1}{3}} )}{
\sqrt{\pi }\,\Gamma (3\,a - {\displaystyle \frac {n}{2}} )\,
\Gamma (a + {\displaystyle \frac {1}{6}}  - {\displaystyle 
\frac {n}{2}} )}}  \\
\mbox{} + {\displaystyle \frac {1}{3}}
{\displaystyle \frac {\sqrt{3}\,(2\,\mathrm{sin}({\displaystyle 
\frac {\pi \,( - 12\,a + 6\,n + 1)}{6}} ) + 1)\,\mathrm{V8}(a, \,
n)\,\Gamma (2\,a + {\displaystyle \frac {2}{3}}  - n)\,\Gamma (2
\,a - n)\,\Gamma (2\,a + {\displaystyle \frac {1}{3}} )}{2^{( - 1
/3 + 2\,a - n)}\,2^{(2\,a - n)}\,\Gamma (2\,a + {\displaystyle 
\frac {5}{6}}  - n)\,\Gamma (2\,a + {\displaystyle \frac {1}{6}} 
 - n)\,\Gamma (2\,a)}}  }
}
\end{maplelatex}

\begin{maplelatex}
\mapleinline{inert}{2d}{315, ":   ",
-2/sin(1/3*Pi*(6*a+1))/(2^(1/3+2*a+n))/GAMMA(3*a-1/2*n)/GAMMA(1/2+n)/G
AMMA(a+1/6-1/2*n)*sin(1/6*Pi*(-12*a+6*n+1))*V7(a-1/3-1/2*n,n+1)*GAMMA(
1/3+2*a+n)*GAMMA(a-1/2*n)*GAMMA(a+1/3-1/2*n)+1/2*(2*sin(1/6*Pi*(-12*a+
6*n+1))+1)*V8(a,n)*Pi^(1/2)*GAMMA(2*a+2/3-n)*GAMMA(2*a-n)*GAMMA(1/3+2*
a+n)/sin(1/3*Pi*(6*a+1))/(2^(2*a-n-1))/(2^(-1/3+2*a-n))/(2^(1/3+2*a+n)
)/GAMMA(2*a+1/6-n)/GAMMA(1/2+n)/GAMMA(2*a+5/6-n)/GAMMA(2*a);}{%
\maplemultiline{
315\mbox{:~~~}  - {\displaystyle \frac {2\,\mathrm{sin}(
{\displaystyle \frac {\pi \,( - 12\,a + 6\,n + 1)}{6}} )\,
\mathrm{V7}(a - {\displaystyle \frac {1}{3}}  - {\displaystyle 
\frac {n}{2}} , \,n + 1)\,\Gamma ({\displaystyle \frac {1}{3}} 
 + 2\,a + n)\,\Gamma (a - {\displaystyle \frac {n}{2}} )\,\Gamma 
(a + {\displaystyle \frac {1}{3}}  - {\displaystyle \frac {n}{2}
} )}{\mathrm{sin}({\displaystyle \frac {\pi \,(6\,a + 1)}{3}} )\,
2^{(1/3 + 2\,a + n)}\,\Gamma (3\,a - {\displaystyle \frac {n}{2}
} )\,\Gamma ({\displaystyle \frac {1}{2}}  + n)\,\Gamma (a + 
{\displaystyle \frac {1}{6}}  - {\displaystyle \frac {n}{2}} )}} 
 \\
\mbox{} + {\displaystyle \frac {1}{2}} (2\,\mathrm{sin}(
{\displaystyle \frac {\pi \,( - 12\,a + 6\,n + 1)}{6}} ) + 1)\,
\mathrm{V8}(a, \,n)\,\sqrt{\pi }\,\Gamma (2\,a + {\displaystyle 
\frac {2}{3}}  - n)\,\Gamma (2\,a - n) \\
\Gamma ({\displaystyle \frac {1}{3}}  + 2\,a + n) \left/ {\vrule 
height0.80em width0em depth0.80em} \right. \!  \! (\mathrm{sin}(
{\displaystyle \frac {\pi \,(6\,a + 1)}{3}} )\,2^{(2\,a - n - 1)}
\,2^{( - 1/3 + 2\,a - n)}\,2^{(1/3 + 2\,a + n)} \\
\Gamma (2\,a + {\displaystyle \frac {1}{6}}  - n)\,\Gamma (
{\displaystyle \frac {1}{2}}  + n)\,\Gamma (2\,a + 
{\displaystyle \frac {5}{6}}  - n)\,\Gamma (2\,a)) }
}
\end{maplelatex}

\begin{maplelatex}
\mapleinline{inert}{2d}{316, ":   ",
GAMMA(2*a+1/3)*GAMMA(a+1/3-1/2*n)*GAMMA(a-1/2*n)/Pi^(1/2)*V7(a-1/3-1/2
*n,n+1)/GAMMA(2*a-n-1/6)/GAMMA(a+1/6-1/2*n)/GAMMA(1/2+n)*GAMMA(a-1/3-1
/2*n)*GAMMA(2/3+a+1/2*n)/GAMMA(3*a-1/2*n)+1/24*V8(a,n)*GAMMA(2*a+2/3-n
)*GAMMA(2/3+a+1/2*n)*GAMMA(a-1/3-1/2*n)*GAMMA(2*a-n)*GAMMA(2*a+1/3)*(2
*sin(1/6*Pi*(-12*a+6*n+1))+1)*(-12*a+6*n+1)/sin(1/6*Pi*(-12*a+6*n+1))/
(2^(2*a-n-1))/(2^(-1/3+2*a-n))/GAMMA(2*a+1/6-n)/GAMMA(2*a+5/6-n)^2/GAM
MA(1/2+n)/GAMMA(2*a);}{%
\maplemultiline{
316\mbox{:~~~}  \\
{\displaystyle \frac {\Gamma (2\,a + {\displaystyle \frac {1}{3}
} )\,\Gamma (a + {\displaystyle \frac {1}{3}}  - {\displaystyle 
\frac {n}{2}} )\,\Gamma (a - {\displaystyle \frac {n}{2}} )\,
\mathrm{V7}(a - {\displaystyle \frac {1}{3}}  - {\displaystyle 
\frac {n}{2}} , \,n + 1)\,\Gamma (a - {\displaystyle \frac {1}{3}
}  - {\displaystyle \frac {n}{2}} )\,\Gamma ({\displaystyle 
\frac {2}{3}}  + a + {\displaystyle \frac {n}{2}} )}{\sqrt{\pi }
\,\Gamma (2\,a - n - {\displaystyle \frac {1}{6}} )\,\Gamma (a + 
{\displaystyle \frac {1}{6}}  - {\displaystyle \frac {n}{2}} )\,
\Gamma ({\displaystyle \frac {1}{2}}  + n)\,\Gamma (3\,a - 
{\displaystyle \frac {n}{2}} )}}  \\
\mbox{} + {\displaystyle \frac {1}{24}} \mathrm{V8}(a, \,n)\,
\Gamma (2\,a + {\displaystyle \frac {2}{3}}  - n)\,\Gamma (
{\displaystyle \frac {2}{3}}  + a + {\displaystyle \frac {n}{2}} 
)\,\Gamma (a - {\displaystyle \frac {1}{3}}  - {\displaystyle 
\frac {n}{2}} )\,\Gamma (2\,a - n)\,\Gamma (2\,a + 
{\displaystyle \frac {1}{3}} ) \\
(2\,\mathrm{sin}({\displaystyle \frac {\pi \,( - 12\,a + 6\,n + 1
)}{6}} ) + 1)\,( - 12\,a + 6\,n + 1) \left/ {\vrule 
height0.80em width0em depth0.80em} \right. \!  \! (\mathrm{sin}(
{\displaystyle \frac {\pi \,( - 12\,a + 6\,n + 1)}{6}} ) \\
2^{(2\,a - n - 1)}\,2^{( - 1/3 + 2\,a - n)}\,\Gamma (2\,a + 
{\displaystyle \frac {1}{6}}  - n)\,\Gamma (2\,a + 
{\displaystyle \frac {5}{6}}  - n)^{2}\,\Gamma ({\displaystyle 
\frac {1}{2}}  + n)\,\Gamma (2\,a)) }
}
\end{maplelatex}

\begin{maplelatex}
\mapleinline{inert}{2d}{317, ":   ",
-sin(1/6*Pi*(-12*a+6*n+1))*V7(a-1/3-1/2*n,n+1)/Pi^(3/2)/GAMMA(3*a-1/2*
n)*GAMMA(2/3+a+1/2*n)*GAMMA(1/6-a+1/2*n)/GAMMA(1/2+n)/GAMMA(a+1/6-1/2*
n)*GAMMA(a-1/2*n)*GAMMA(a+1/3-1/2*n)*GAMMA(2*a+1/3)*GAMMA(2/3)+1/4*(2*
sin(1/6*Pi*(-12*a+6*n+1))+1)*V8(a,n)*GAMMA(2*a-n)*GAMMA(2/3+a+1/2*n)*G
AMMA(a+1/3-1/2*n)*GAMMA(2*a+1/3)*GAMMA(2/3)/sin(1/6*Pi*(1-6*a+3*n))/Pi
^(1/2)/(2^(2*a-n-1))/GAMMA(2*a+1/6-n)/GAMMA(2*a+5/6-n)/GAMMA(1/2+n)/GA
MMA(2*a);}{%
\maplemultiline{
317\mbox{:~~~}  - \mathrm{sin}({\displaystyle \frac {\pi
 \,( - 12\,a + 6\,n + 1)}{6}} )\,\mathrm{V7}(a - {\displaystyle 
\frac {1}{3}}  - {\displaystyle \frac {n}{2}} , \,n + 1)\,\Gamma 
({\displaystyle \frac {2}{3}}  + a + {\displaystyle \frac {n}{2}
} )\,\Gamma ({\displaystyle \frac {1}{6}}  - a + {\displaystyle 
\frac {n}{2}} ) \\
\Gamma (a - {\displaystyle \frac {n}{2}} )\,\Gamma (a + 
{\displaystyle \frac {1}{3}}  - {\displaystyle \frac {n}{2}} )\,
\Gamma (2\,a + {\displaystyle \frac {1}{3}} )\,\Gamma (
{\displaystyle \frac {2}{3}} ) \left/ {\vrule 
height0.80em width0em depth0.80em} \right. \!  \! (\pi ^{(3/2)}\,
\Gamma (3\,a - {\displaystyle \frac {n}{2}} )\,\Gamma (
{\displaystyle \frac {1}{2}}  + n) \\
\Gamma (a + {\displaystyle \frac {1}{6}}  - {\displaystyle 
\frac {n}{2}} ))\mbox{} + {\displaystyle \frac {1}{4}} (2\,
\mathrm{sin}({\displaystyle \frac {\pi \,( - 12\,a + 6\,n + 1)}{6
}} ) + 1)\,\mathrm{V8}(a, \,n)\,\Gamma (2\,a - n) \\
\Gamma ({\displaystyle \frac {2}{3}}  + a + {\displaystyle 
\frac {n}{2}} )\,\Gamma (a + {\displaystyle \frac {1}{3}}  - 
{\displaystyle \frac {n}{2}} )\,\Gamma (2\,a + {\displaystyle 
\frac {1}{3}} )\,\Gamma ({\displaystyle \frac {2}{3}} ) \left/ 
{\vrule height0.80em width0em depth0.80em} \right. \!  \! (
\mathrm{sin}({\displaystyle \frac {\pi \,(1 - 6\,a + 3\,n)}{6}} )
\,\sqrt{\pi } \\
2^{(2\,a - n - 1)}\,\Gamma (2\,a + {\displaystyle \frac {1}{6}} 
 - n)\,\Gamma (2\,a + {\displaystyle \frac {5}{6}}  - n)\,\Gamma 
({\displaystyle \frac {1}{2}}  + n)\,\Gamma (2\,a)) }
}
\end{maplelatex}

\begin{maplelatex}
\mapleinline{inert}{2d}{318, ":   ",
sin(1/6*Pi*(-12*a+6*n+1))*V7(a-1/3-1/2*n,n+1)/sin(1/6*Pi*(-6*a-2+3*n))
/Pi^(1/2)/GAMMA(3*a-1/2*n)*GAMMA(1/6+a+1/2*n)/GAMMA(1/2+n)/GAMMA(a+1/6
-1/2*n)*GAMMA(a-1/2*n)*GAMMA(2*a+1/3)*GAMMA(2/3)-1/4*(2*sin(1/6*Pi*(-1
2*a+6*n+1))+1)*V8(a,n)*GAMMA(2*a-n)*GAMMA(1/6+a+1/2*n)*GAMMA(a-1/2*n+5
/6)*GAMMA(2*a+1/3)*GAMMA(2/3)/sin(1/6*Pi*(-6*a-2+3*n))/Pi^(1/2)/(2^(2*
a-n-1))/GAMMA(2*a+5/6-n)/GAMMA(1/2+n)/GAMMA(2*a+1/6-n)/GAMMA(2*a);}{%
\maplemultiline{
318\mbox{:~~~} {\displaystyle \frac {\mathrm{sin}(
{\displaystyle \frac {\pi \,( - 12\,a + 6\,n + 1)}{6}} )\,
\mathrm{V7}(a - {\displaystyle \frac {1}{3}}  - {\displaystyle 
\frac {n}{2}} , \,n + 1)\,\Gamma ({\displaystyle \frac {1}{6}} 
 + a + {\displaystyle \frac {n}{2}} )\,\Gamma (a - 
{\displaystyle \frac {n}{2}} )\,\Gamma (2\,a + {\displaystyle 
\frac {1}{3}} )\,\Gamma ({\displaystyle \frac {2}{3}} )}{\mathrm{
sin}({\displaystyle \frac {\pi \,( - 6\,a - 2 + 3\,n)}{6}} )\,
\sqrt{\pi }\,\Gamma (3\,a - {\displaystyle \frac {n}{2}} )\,
\Gamma ({\displaystyle \frac {1}{2}}  + n)\,\Gamma (a + 
{\displaystyle \frac {1}{6}}  - {\displaystyle \frac {n}{2}} )}} 
 \\
\mbox{} - {\displaystyle \frac {1}{4}} (2\,\mathrm{sin}(
{\displaystyle \frac {\pi \,( - 12\,a + 6\,n + 1)}{6}} ) + 1)\,
\mathrm{V8}(a, \,n)\,\Gamma (2\,a - n)\,\Gamma ({\displaystyle 
\frac {1}{6}}  + a + {\displaystyle \frac {n}{2}} ) \\
\Gamma (a - {\displaystyle \frac {n}{2}}  + {\displaystyle 
\frac {5}{6}} )\,\Gamma (2\,a + {\displaystyle \frac {1}{3}} )\,
\Gamma ({\displaystyle \frac {2}{3}} ) \left/ {\vrule 
height0.80em width0em depth0.80em} \right. \!  \! (\mathrm{sin}(
{\displaystyle \frac {\pi \,( - 6\,a - 2 + 3\,n)}{6}} )\,\sqrt{
\pi }\,2^{(2\,a - n - 1)} \\
\Gamma (2\,a + {\displaystyle \frac {5}{6}}  - n)\,\Gamma (
{\displaystyle \frac {1}{2}}  + n)\,\Gamma (2\,a + 
{\displaystyle \frac {1}{6}}  - n)\,\Gamma (2\,a)) }
}
\end{maplelatex}

\begin{maplelatex}
\mapleinline{inert}{2d}{319, ":   ",
-(-1+cos(Pi*(-6*a+n)))*V7(a-1/3-1/2*n,n+1)*GAMMA(1-3*a+1/2*n)*GAMMA(3/
2-a-1/2*n)*GAMMA(a-1/3-1/2*n)*GAMMA(a-1/2*n)*GAMMA(2*a+1/3)*GAMMA(2*a-
1/3)*GAMMA(2*a)/((-1)^(-n))/Pi^(3/2)/GAMMA(a+1/2*n-5/6)/GAMMA(1/2+a-1/
2*n)/GAMMA(a-1/2*n+5/6)/GAMMA(a+1/6-1/2*n)/(sin(1/3*Pi*(3*n+1))+sin(2*
Pi*a))+1/2*(sin(3/2*Pi*(n-2*a))-2*sin(1/2*Pi*(6*a+n)))*V8(a,n)*GAMMA(3
/2-a-1/2*n)*GAMMA(a-1/3-1/2*n)*GAMMA(a-1/2*n)*GAMMA(2*a+1/3)*GAMMA(2*a
-1/3)/Pi/GAMMA(2*a+5/6-n)/GAMMA(2*a+1/6-n)/GAMMA(a+1/2*n-5/6)/(sin(1/3
*Pi*(3*n+1))+sin(2*Pi*a));}{%
\maplemultiline{
319\mbox{:~~~}  - ( - 1 + \mathrm{cos}(\pi \,( - 6\,a + 
n)))\,\mathrm{V7}(a - {\displaystyle \frac {1}{3}}  - 
{\displaystyle \frac {n}{2}} , \,n + 1)\,\Gamma (1 - 3\,a + 
{\displaystyle \frac {n}{2}} )\,\Gamma ({\displaystyle \frac {3}{
2}}  - a - {\displaystyle \frac {n}{2}} ) \\
\Gamma (a - {\displaystyle \frac {1}{3}}  - {\displaystyle 
\frac {n}{2}} )\,\Gamma (a - {\displaystyle \frac {n}{2}} )\,
\Gamma (2\,a + {\displaystyle \frac {1}{3}} )\,\Gamma (2\,a - 
{\displaystyle \frac {1}{3}} )\,\Gamma (2\,a) \left/ {\vrule 
height0.80em width0em depth0.80em} \right. \!  \! ((-1)^{( - n)}
\,\pi ^{(3/2)} \\
\Gamma (a + {\displaystyle \frac {n}{2}}  - {\displaystyle 
\frac {5}{6}} )\,\Gamma ({\displaystyle \frac {1}{2}}  + a - 
{\displaystyle \frac {n}{2}} )\,\Gamma (a - {\displaystyle 
\frac {n}{2}}  + {\displaystyle \frac {5}{6}} )\,\Gamma (a + 
{\displaystyle \frac {1}{6}}  - {\displaystyle \frac {n}{2}} )
 \\
(\mathrm{sin}({\displaystyle \frac {\pi \,(3\,n + 1)}{3}} ) + 
\mathrm{sin}(2\,\pi \,a)))\mbox{} + {\displaystyle \frac {1}{2}} 
(\mathrm{sin}({\displaystyle \frac {3\,\pi \,(n - 2\,a)}{2}} ) - 
2\,\mathrm{sin}({\displaystyle \frac {\pi \,(6\,a + n)}{2}} ))
 \\
\mathrm{V8}(a, \,n)\,\Gamma ({\displaystyle \frac {3}{2}}  - a - 
{\displaystyle \frac {n}{2}} )\,\Gamma (a - {\displaystyle 
\frac {1}{3}}  - {\displaystyle \frac {n}{2}} )\,\Gamma (a - 
{\displaystyle \frac {n}{2}} )\,\Gamma (2\,a + {\displaystyle 
\frac {1}{3}} )\,\Gamma (2\,a - {\displaystyle \frac {1}{3}} )
 \left/ {\vrule height0.80em width0em depth0.80em} \right. \! 
 \! (\pi  \\
\Gamma (2\,a + {\displaystyle \frac {5}{6}}  - n)\,\Gamma (2\,a
 + {\displaystyle \frac {1}{6}}  - n)\,\Gamma (a + 
{\displaystyle \frac {n}{2}}  - {\displaystyle \frac {5}{6}} )\,(
\mathrm{sin}({\displaystyle \frac {\pi \,(3\,n + 1)}{3}} ) + 
\mathrm{sin}(2\,\pi \,a))) }
}
\end{maplelatex}

\mapleinline{inert}{2d}{320, ":   ",
-sin(1/2*Pi*(-6*a+n))*V7(a-1/3-1/2*n,n+1)*GAMMA(a-1/3-1/2*n)*GAMMA(a+1
/3-1/2*n)*GAMMA(a-1/2*n)*GAMMA(2*a-1/3)*GAMMA(2*a)/((-1)^(-n))/GAMMA(a
+1/2*n-5/6)/GAMMA(3*a-1/2*n)/GAMMA(1/2+a-1/2*n)/GAMMA(a-1/2*n+5/6)/GAM
MA(a+1/6-1/2*n)/(sin(1/6*Pi*(-6*a+2+3*n))+sin(1/2*Pi*(6*a+n)))+1/4*(si
n(3/2*Pi*(n-2*a))-2*sin(1/2*Pi*(6*a+n)))*V8(a,n)*GAMMA(a-1/3-1/2*n)*GA
MMA(a+1/3-1/2*n)*GAMMA(a-1/2*n)*GAMMA(2*a-1/3)/Pi^(1/2)/GAMMA(2*a+5/6-
n)/GAMMA(2*a+1/6-n)/GAMMA(a+1/2*n-5/6)/(sin(1/6*Pi*(-6*a+2+3*n))+sin(1
/2*Pi*(6*a+n)));}{%
\maplemultiline{
320\mbox{:~~~}  - \mathrm{sin}({\displaystyle \frac {\pi
 \,( - 6\,a + n)}{2}} )\,\mathrm{V7}(a - {\displaystyle \frac {1
}{3}}  - {\displaystyle \frac {n}{2}} , \,n + 1)\,\Gamma (a - 
{\displaystyle \frac {1}{3}}  - {\displaystyle \frac {n}{2}} )\,
\Gamma (a + {\displaystyle \frac {1}{3}}  - {\displaystyle 
\frac {n}{2}} )\,\Gamma (a - {\displaystyle \frac {n}{2}} ) \\
\Gamma (2\,a - {\displaystyle \frac {1}{3}} )\,\Gamma (2\,a)
 \left/ {\vrule height0.80em width0em depth0.80em} \right. \! 
 \! ((-1)^{( - n)}\,\Gamma (a + {\displaystyle \frac {n}{2}}  - 
{\displaystyle \frac {5}{6}} )\,\Gamma (3\,a - {\displaystyle 
\frac {n}{2}} )\,\Gamma ({\displaystyle \frac {1}{2}}  + a - 
{\displaystyle \frac {n}{2}} ) \\
\Gamma (a - {\displaystyle \frac {n}{2}}  + {\displaystyle 
\frac {5}{6}} )\,\Gamma (a + {\displaystyle \frac {1}{6}}  - 
{\displaystyle \frac {n}{2}} )\,(\mathrm{sin}({\displaystyle 
\frac {\pi \,( - 6\,a + 2 + 3\,n)}{6}} ) + \mathrm{sin}(
{\displaystyle \frac {\pi \,(6\,a + n)}{2}} )))\mbox{} + 
{\displaystyle \frac {1}{4}}  \\
(\mathrm{sin}({\displaystyle \frac {3\,\pi \,(n - 2\,a)}{2}} ) - 
2\,\mathrm{sin}({\displaystyle \frac {\pi \,(6\,a + n)}{2}} ))\,
\mathrm{V8}(a, \,n)\,\Gamma (a - {\displaystyle \frac {1}{3}}  - 
{\displaystyle \frac {n}{2}} )\,\Gamma (a + {\displaystyle 
\frac {1}{3}}  - {\displaystyle \frac {n}{2}} ) \\
\Gamma (a - {\displaystyle \frac {n}{2}} )\,\Gamma (2\,a - 
{\displaystyle \frac {1}{3}} ) \left/ {\vrule 
height0.80em width0em depth0.80em} \right. \!  \! (\sqrt{\pi }\,
\Gamma (2\,a + {\displaystyle \frac {5}{6}}  - n)\,\Gamma (2\,a
 + {\displaystyle \frac {1}{6}}  - n)\,\Gamma (a + 
{\displaystyle \frac {n}{2}}  - {\displaystyle \frac {5}{6}} )
 \\
(\mathrm{sin}({\displaystyle \frac {\pi \,( - 6\,a + 2 + 3\,n)}{6
}} ) + \mathrm{sin}({\displaystyle \frac {\pi \,(6\,a + n)}{2}} )
)) }
}

\begin{mapleinput}
\end{mapleinput}

\end{maplegroup}

%% file: AppendixB321to350.tex
\begin{maplegroup}
\mapleresult
\begin{maplelatex}
\mapleinline{inert}{2d}{321, ":   ",
-2*sin(1/2*Pi*(-6*a+n))*V7(a-1/3-1/2*n,n+1)*GAMMA(a-1/3-1/2*n)*GAMMA(7
/6-a-1/2*n)*GAMMA(a+1/3-1/2*n)*GAMMA(2*a+1/3)*GAMMA(2*a-1/3)*GAMMA(2*a
)/((-1)^(-n))/Pi^(1/2)/GAMMA(a+1/2*n-5/6)/GAMMA(3*a-1/2*n)/GAMMA(1/2+a
-1/2*n)/GAMMA(a-1/2*n+5/6)/GAMMA(a+1/6-1/2*n)/(sin(1/3*Pi*(6*a-1))-sin
(1/3*Pi*(3*n-1)))-1/2*(2*sin(1/2*Pi*(6*a+n))-sin(3/2*Pi*(n-2*a)))*V8(a
,n)*GAMMA(a-1/3-1/2*n)*GAMMA(7/6-a-1/2*n)*GAMMA(a+1/3-1/2*n)*GAMMA(2*a
+1/3)*GAMMA(2*a-1/3)/Pi/GAMMA(2*a+5/6-n)/GAMMA(2*a+1/6-n)/GAMMA(a+1/2*
n-5/6)/(sin(1/3*Pi*(6*a-1))-sin(1/3*Pi*(3*n-1)));}{%
\maplemultiline{
321\mbox{:~~~}  - 2\,\mathrm{sin}({\displaystyle \frac {
\pi \,( - 6\,a + n)}{2}} )\,\mathrm{V7}(a - {\displaystyle 
\frac {1}{3}}  - {\displaystyle \frac {n}{2}} , \,n + 1)\,\Gamma 
(a - {\displaystyle \frac {1}{3}}  - {\displaystyle \frac {n}{2}
} )\,\Gamma ({\displaystyle \frac {7}{6}}  - a - {\displaystyle 
\frac {n}{2}} ) \\
\Gamma (a + {\displaystyle \frac {1}{3}}  - {\displaystyle 
\frac {n}{2}} )\,\Gamma (2\,a + {\displaystyle \frac {1}{3}} )\,
\Gamma (2\,a - {\displaystyle \frac {1}{3}} )\,\Gamma (2\,a)
 \left/ {\vrule height0.80em width0em depth0.80em} \right. \! 
 \! ((-1)^{( - n)}\,\sqrt{\pi }\,\Gamma (a + {\displaystyle 
\frac {n}{2}}  - {\displaystyle \frac {5}{6}} ) \\
\Gamma (3\,a - {\displaystyle \frac {n}{2}} )\,\Gamma (
{\displaystyle \frac {1}{2}}  + a - {\displaystyle \frac {n}{2}} 
)\,\Gamma (a - {\displaystyle \frac {n}{2}}  + {\displaystyle 
\frac {5}{6}} )\,\Gamma (a + {\displaystyle \frac {1}{6}}  - 
{\displaystyle \frac {n}{2}} ) \\
(\mathrm{sin}({\displaystyle \frac {\pi \,(6\,a - 1)}{3}} ) - 
\mathrm{sin}({\displaystyle \frac {\pi \,(3\,n - 1)}{3}} )))
\mbox{} - {\displaystyle \frac {1}{2}} (2\,\mathrm{sin}(
{\displaystyle \frac {\pi \,(6\,a + n)}{2}} ) - \mathrm{sin}(
{\displaystyle \frac {3\,\pi \,(n - 2\,a)}{2}} )) \\
\mathrm{V8}(a, \,n)\,\Gamma (a - {\displaystyle \frac {1}{3}}  - 
{\displaystyle \frac {n}{2}} )\,\Gamma ({\displaystyle \frac {7}{
6}}  - a - {\displaystyle \frac {n}{2}} )\,\Gamma (a + 
{\displaystyle \frac {1}{3}}  - {\displaystyle \frac {n}{2}} )\,
\Gamma (2\,a + {\displaystyle \frac {1}{3}} )\,\Gamma (2\,a - 
{\displaystyle \frac {1}{3}} ) \left/ {\vrule 
height0.80em width0em depth0.80em} \right. \!  \! (\pi  \\
\Gamma (2\,a + {\displaystyle \frac {5}{6}}  - n)\,\Gamma (2\,a
 + {\displaystyle \frac {1}{6}}  - n)\,\Gamma (a + 
{\displaystyle \frac {n}{2}}  - {\displaystyle \frac {5}{6}} )\,(
\mathrm{sin}({\displaystyle \frac {\pi \,(6\,a - 1)}{3}} ) - 
\mathrm{sin}({\displaystyle \frac {\pi \,(3\,n - 1)}{3}} ))) }
}
\end{maplelatex}

\begin{maplelatex}
\mapleinline{inert}{2d}{322, ":   ",
1/2*1/((-1)^(-n))/sin(1/6*Pi*(1+6*a+3*n))/Pi^(1/2)/GAMMA(a+1/2*n-5/6)/
GAMMA(1/2+a-1/2*n)/GAMMA(3*a-1/2*n)^2/GAMMA(a-1/2*n+5/6)*(-6*a+n+2)*V7
(a-1/3-1/2*n,n+1)*GAMMA(a-1/3-1/2*n)*GAMMA(a+1/3-1/2*n)*GAMMA(a-1/2*n)
*GAMMA(2*a+1/3)*GAMMA(2*a-1/3)*GAMMA(2*a)+1/2*V8(a,n)*(2*sin(1/2*Pi*(6
*a+n))-sin(3/2*Pi*(n-2*a)))*GAMMA(2*a-2/3-n)*GAMMA(a+1/3-1/2*n)*GAMMA(
a-1/2*n)*GAMMA(2*a+1/3)*GAMMA(2*a-1/3)/Pi^(1/2)/(2^(2*a-5/3-n))/GAMMA(
3*a-1/2*n-1)/GAMMA(2*a+5/6-n)/GAMMA(2*a+1/6-n)/GAMMA(a+1/2*n-5/6)/(sin
(1/3*Pi*(12*a-1))-sin(1/3*Pi*(-6*a-1+3*n)));}{%
\maplemultiline{
322\mbox{:~~~} {\displaystyle \frac {1}{2}} ( - 6\,a + n
 + 2)\,\mathrm{V7}(a - {\displaystyle \frac {1}{3}}  - 
{\displaystyle \frac {n}{2}} , \,n + 1)\,\Gamma (a - 
{\displaystyle \frac {1}{3}}  - {\displaystyle \frac {n}{2}} )\,
\Gamma (a + {\displaystyle \frac {1}{3}}  - {\displaystyle 
\frac {n}{2}} )\,\Gamma (a - {\displaystyle \frac {n}{2}} ) \\
\Gamma (2\,a + {\displaystyle \frac {1}{3}} )\,\Gamma (2\,a - 
{\displaystyle \frac {1}{3}} )\,\Gamma (2\,a) \left/ {\vrule 
height0.80em width0em depth0.80em} \right. \!  \! ((-1)^{( - n)}
\,\mathrm{sin}({\displaystyle \frac {\pi \,(1 + 6\,a + 3\,n)}{6}
} )\,\sqrt{\pi }\,\Gamma (a + {\displaystyle \frac {n}{2}}  - 
{\displaystyle \frac {5}{6}} ) \\
\Gamma ({\displaystyle \frac {1}{2}}  + a - {\displaystyle 
\frac {n}{2}} )\,\Gamma (3\,a - {\displaystyle \frac {n}{2}} )^{2
}\,\Gamma (a - {\displaystyle \frac {n}{2}}  + {\displaystyle 
\frac {5}{6}} ))\mbox{} + {\displaystyle \frac {1}{2}} \mathrm{V8
}(a, \,n) \\
(2\,\mathrm{sin}({\displaystyle \frac {\pi \,(6\,a + n)}{2}} ) - 
\mathrm{sin}({\displaystyle \frac {3\,\pi \,(n - 2\,a)}{2}} ))\,
\Gamma (2\,a - {\displaystyle \frac {2}{3}}  - n)\,\Gamma (a + 
{\displaystyle \frac {1}{3}}  - {\displaystyle \frac {n}{2}} )\,
\Gamma (a - {\displaystyle \frac {n}{2}} ) \\
\Gamma (2\,a + {\displaystyle \frac {1}{3}} )\,\Gamma (2\,a - 
{\displaystyle \frac {1}{3}} ) \left/ {\vrule 
height0.80em width0em depth0.80em} \right. \!  \! (\sqrt{\pi }\,2
^{(2\,a - 5/3 - n)}\,\Gamma (3\,a - {\displaystyle \frac {n}{2}} 
 - 1)\,\Gamma (2\,a + {\displaystyle \frac {5}{6}}  - n) \\
\Gamma (2\,a + {\displaystyle \frac {1}{6}}  - n)\,\Gamma (a + 
{\displaystyle \frac {n}{2}}  - {\displaystyle \frac {5}{6}} )\,(
\mathrm{sin}({\displaystyle \frac {\pi \,(12\,a - 1)}{3}} ) - 
\mathrm{sin}({\displaystyle \frac {\pi \,( - 6\,a - 1 + 3\,n)}{3}
} ))) }
}
\end{maplelatex}

\begin{maplelatex}
\mapleinline{inert}{2d}{323, ":   ",
1/2*1/Pi/sin(2*Pi*a)/((-1)^(-n))/GAMMA(3*a-1/2*n)/GAMMA(1/2+a-1/2*n)/G
AMMA(a-1/2*n+5/6)/GAMMA(a+1/6-1/2*n)*sin(1/2*Pi*(-6*a+n))*V7(a-1/3-1/2
*n,n+1)*GAMMA(3/2-a-1/2*n)*GAMMA(a-1/3-1/2*n)*GAMMA(a+1/3-1/2*n)*GAMMA
(a-1/2*n)*GAMMA(2*a+1/3)*GAMMA(2*a-1/3)+1/8*(2*sin(1/2*Pi*(6*a+n))-sin
(1/2*Pi*(n-2*a))+sin(1/6*Pi*(-6*a-2+3*n))-sin(3/2*Pi*(n-2*a))+sin(1/6*
Pi*(-6*a+2+3*n)))*V8(a,n)*GAMMA(3/2-a-1/2*n)*GAMMA(a-1/3-1/2*n)*GAMMA(
a+1/3-1/2*n)*GAMMA(a-1/2*n)*GAMMA(2*a+1/3)*GAMMA(2*a-1/3)/sin(2*Pi*a)/
Pi^(3/2)/GAMMA(2*a+5/6-n)/GAMMA(2*a+1/6-n)/GAMMA(2*a);}{%
\maplemultiline{
323\mbox{:~~~} {\displaystyle \frac {1}{2}} \mathrm{sin}
({\displaystyle \frac {\pi \,( - 6\,a + n)}{2}} )\,\mathrm{V7}(a
 - {\displaystyle \frac {1}{3}}  - {\displaystyle \frac {n}{2}} 
, \,n + 1)\,\Gamma ({\displaystyle \frac {3}{2}}  - a - 
{\displaystyle \frac {n}{2}} )\,\Gamma (a - {\displaystyle 
\frac {1}{3}}  - {\displaystyle \frac {n}{2}} ) \\
\Gamma (a + {\displaystyle \frac {1}{3}}  - {\displaystyle 
\frac {n}{2}} )\,\Gamma (a - {\displaystyle \frac {n}{2}} )\,
\Gamma (2\,a + {\displaystyle \frac {1}{3}} )\,\Gamma (2\,a - 
{\displaystyle \frac {1}{3}} ) \left/ {\vrule 
height0.80em width0em depth0.80em} \right. \!  \! (\pi \,\mathrm{
sin}(2\,\pi \,a)\,(-1)^{( - n)} \\
\Gamma (3\,a - {\displaystyle \frac {n}{2}} )\,\Gamma (
{\displaystyle \frac {1}{2}}  + a - {\displaystyle \frac {n}{2}} 
)\,\Gamma (a - {\displaystyle \frac {n}{2}}  + {\displaystyle 
\frac {5}{6}} )\,\Gamma (a + {\displaystyle \frac {1}{6}}  - 
{\displaystyle \frac {n}{2}} ))\mbox{} + {\displaystyle \frac {1
}{8}} (2\,\mathrm{sin}({\displaystyle \frac {\pi \,(6\,a + n)}{2}
} ) \\
\mbox{} - \mathrm{sin}({\displaystyle \frac {\pi \,(n - 2\,a)}{2}
} ) + \mathrm{sin}({\displaystyle \frac {\pi \,( - 6\,a - 2 + 3\,
n)}{6}} ) - \mathrm{sin}({\displaystyle \frac {3\,\pi \,(n - 2\,a
)}{2}} ) \\
\mbox{} + \mathrm{sin}({\displaystyle \frac {\pi \,( - 6\,a + 2
 + 3\,n)}{6}} ))\mathrm{V8}(a, \,n)\,\Gamma ({\displaystyle 
\frac {3}{2}}  - a - {\displaystyle \frac {n}{2}} )\,\Gamma (a - 
{\displaystyle \frac {1}{3}}  - {\displaystyle \frac {n}{2}} )\,
\Gamma (a + {\displaystyle \frac {1}{3}}  - {\displaystyle 
\frac {n}{2}} ) \\
\Gamma (a - {\displaystyle \frac {n}{2}} )\,\Gamma (2\,a + 
{\displaystyle \frac {1}{3}} )\,\Gamma (2\,a - {\displaystyle 
\frac {1}{3}} ) \left/ {\vrule height0.80em width0em depth0.80em}
 \right. \!  \! (\mathrm{sin}(2\,\pi \,a)\,\pi ^{(3/2)}\,\Gamma (
2\,a + {\displaystyle \frac {5}{6}}  - n) \\
\Gamma (2\,a + {\displaystyle \frac {1}{6}}  - n)\,\Gamma (2\,a))
 }
}
\end{maplelatex}

\begin{maplelatex}
\mapleinline{inert}{2d}{324, ":   ",
-1/sin(1/6*Pi*(-6*a+2+3*n))/((-1)^(-n))/Pi^(3/2)/GAMMA(3*a-1/2*n)/GAMM
A(1/2+a-1/2*n)/GAMMA(a-1/2*n+5/6)/GAMMA(a+1/6-1/2*n)*sin(1/2*Pi*(-6*a+
n))*V7(a-1/3-1/2*n,n+1)*GAMMA(3/2-a-1/2*n)*GAMMA(7/6-a-1/2*n)*GAMMA(a+
1/3-1/2*n)*GAMMA(a-1/2*n)*GAMMA(2*a+1/3)*GAMMA(2*a-1/3)*GAMMA(2*a)+1/4
*(-2*sin(1/2*Pi*(6*a+n))+sin(3/2*Pi*(n-2*a)))*V8(a,n)*GAMMA(3/2-a-1/2*
n)*GAMMA(7/6-a-1/2*n)*GAMMA(a+1/3-1/2*n)*GAMMA(a-1/2*n)*GAMMA(2*a+1/3)
*GAMMA(2*a-1/3)/sin(1/6*Pi*(-6*a+2+3*n))/Pi^2/GAMMA(2*a+5/6-n)/GAMMA(2
*a+1/6-n);}{%
\maplemultiline{
324\mbox{:~~~}  - \mathrm{sin}({\displaystyle \frac {\pi
 \,( - 6\,a + n)}{2}} )\,\mathrm{V7}(a - {\displaystyle \frac {1
}{3}}  - {\displaystyle \frac {n}{2}} , \,n + 1)\,\Gamma (
{\displaystyle \frac {3}{2}}  - a - {\displaystyle \frac {n}{2}} 
)\,\Gamma ({\displaystyle \frac {7}{6}}  - a - {\displaystyle 
\frac {n}{2}} ) \\
\Gamma (a + {\displaystyle \frac {1}{3}}  - {\displaystyle 
\frac {n}{2}} )\,\Gamma (a - {\displaystyle \frac {n}{2}} )\,
\Gamma (2\,a + {\displaystyle \frac {1}{3}} )\,\Gamma (2\,a - 
{\displaystyle \frac {1}{3}} )\,\Gamma (2\,a) \left/ {\vrule 
height0.80em width0em depth0.80em} \right. \!  \! (\mathrm{sin}(
{\displaystyle \frac {\pi \,( - 6\,a + 2 + 3\,n)}{6}} ) \\
(-1)^{( - n)}\,\pi ^{(3/2)}\,\Gamma (3\,a - {\displaystyle 
\frac {n}{2}} )\,\Gamma ({\displaystyle \frac {1}{2}}  + a - 
{\displaystyle \frac {n}{2}} )\,\Gamma (a - {\displaystyle 
\frac {n}{2}}  + {\displaystyle \frac {5}{6}} )\,\Gamma (a + 
{\displaystyle \frac {1}{6}}  - {\displaystyle \frac {n}{2}} ))
\mbox{} + {\displaystyle \frac {1}{4}}  \\
( - 2\,\mathrm{sin}({\displaystyle \frac {\pi \,(6\,a + n)}{2}} )
 + \mathrm{sin}({\displaystyle \frac {3\,\pi \,(n - 2\,a)}{2}} ))
\,\mathrm{V8}(a, \,n)\,\Gamma ({\displaystyle \frac {3}{2}}  - a
 - {\displaystyle \frac {n}{2}} )\,\Gamma ({\displaystyle \frac {
7}{6}}  - a - {\displaystyle \frac {n}{2}} ) \\
\Gamma (a + {\displaystyle \frac {1}{3}}  - {\displaystyle 
\frac {n}{2}} )\,\Gamma (a - {\displaystyle \frac {n}{2}} )\,
\Gamma (2\,a + {\displaystyle \frac {1}{3}} )\,\Gamma (2\,a - 
{\displaystyle \frac {1}{3}} ) \left/ {\vrule 
height0.80em width0em depth0.80em} \right. \!  \! (\mathrm{sin}(
{\displaystyle \frac {\pi \,( - 6\,a + 2 + 3\,n)}{6}} )\,\pi ^{2}
 \\
\Gamma (2\,a + {\displaystyle \frac {5}{6}}  - n)\,\Gamma (2\,a
 + {\displaystyle \frac {1}{6}}  - n)) }
}
\end{maplelatex}

\begin{maplelatex}
\mapleinline{inert}{2d}{325, ":   ",
-1/2*1/((-1)^(-n))/Pi^(3/2)/GAMMA(3*a-1/2*n)^2/GAMMA(a-1/2*n+5/6)/GAMM
A(a+1/6-1/2*n)*(-6*a+n+2)*V7(a-1/3-1/2*n,n+1)*GAMMA(3/2-a-1/2*n)*GAMMA
(a-1/2*n)*GAMMA(a-1/3-1/2*n)*GAMMA(a+1/3-1/2*n)*GAMMA(2*a+1/3)*GAMMA(2
*a-1/3)*GAMMA(2*a)-1/4*V8(a,n)*(-2*sin(1/2*Pi*(6*a+n))+sin(3/2*Pi*(n-2
*a)))*GAMMA(2*a-n)*GAMMA(3/2-a-1/2*n)*GAMMA(a-1/3-1/2*n)*GAMMA(a+1/3-1
/2*n)*GAMMA(2*a+1/3)*GAMMA(2*a-1/3)/sin(1/2*Pi*(-6*a+n))/Pi^(3/2)/(2^(
2*a-n-1))/GAMMA(3*a-1/2*n-1)/GAMMA(2*a+5/6-n)/GAMMA(2*a+1/6-n);}{%
\maplemultiline{
325\mbox{:~~~}  - {\displaystyle \frac {1}{2}} ( - 6\,a
 + n + 2)\,\mathrm{V7}(a - {\displaystyle \frac {1}{3}}  - 
{\displaystyle \frac {n}{2}} , \,n + 1)\,\Gamma ({\displaystyle 
\frac {3}{2}}  - a - {\displaystyle \frac {n}{2}} )\,\Gamma (a - 
{\displaystyle \frac {n}{2}} )\,\Gamma (a - {\displaystyle 
\frac {1}{3}}  - {\displaystyle \frac {n}{2}} ) \\
\Gamma (a + {\displaystyle \frac {1}{3}}  - {\displaystyle 
\frac {n}{2}} )\,\Gamma (2\,a + {\displaystyle \frac {1}{3}} )\,
\Gamma (2\,a - {\displaystyle \frac {1}{3}} )\,\Gamma (2\,a)
 \left/ {\vrule height0.80em width0em depth0.80em} \right. \! 
 \! ((-1)^{( - n)}\,\pi ^{(3/2)}\,\Gamma (3\,a - {\displaystyle 
\frac {n}{2}} )^{2} \\
\Gamma (a - {\displaystyle \frac {n}{2}}  + {\displaystyle 
\frac {5}{6}} )\,\Gamma (a + {\displaystyle \frac {1}{6}}  - 
{\displaystyle \frac {n}{2}} ))\mbox{} - {\displaystyle \frac {1
}{4}} \mathrm{V8}(a, \,n)\,( - 2\,\mathrm{sin}({\displaystyle 
\frac {\pi \,(6\,a + n)}{2}} ) + \mathrm{sin}({\displaystyle 
\frac {3\,\pi \,(n - 2\,a)}{2}} )) \\
\Gamma (2\,a - n)\,\Gamma ({\displaystyle \frac {3}{2}}  - a - 
{\displaystyle \frac {n}{2}} )\,\Gamma (a - {\displaystyle 
\frac {1}{3}}  - {\displaystyle \frac {n}{2}} )\,\Gamma (a + 
{\displaystyle \frac {1}{3}}  - {\displaystyle \frac {n}{2}} )\,
\Gamma (2\,a + {\displaystyle \frac {1}{3}} )\,\Gamma (2\,a - 
{\displaystyle \frac {1}{3}} ) \left/ {\vrule 
height0.80em width0em depth0.80em} \right. \!  \! ( \\
\mathrm{sin}({\displaystyle \frac {\pi \,( - 6\,a + n)}{2}} )\,
\pi ^{(3/2)}\,2^{(2\,a - n - 1)}\,\Gamma (3\,a - {\displaystyle 
\frac {n}{2}}  - 1)\,\Gamma (2\,a + {\displaystyle \frac {5}{6}} 
 - n)\,\Gamma (2\,a + {\displaystyle \frac {1}{6}}  - n)) }
}
\end{maplelatex}

\begin{maplelatex}
\mapleinline{inert}{2d}{326, ":   ",
1/2*1/Pi/sin(1/3*Pi*(6*a-1))/((-1)^(-n))/GAMMA(3*a-1/2*n)/GAMMA(1/2+a-
1/2*n)/GAMMA(a-1/2*n+5/6)/GAMMA(a+1/6-1/2*n)*sin(1/2*Pi*(-6*a+n))*V7(a
-1/3-1/2*n,n+1)*GAMMA(a-1/3-1/2*n)*GAMMA(7/6-a-1/2*n)*GAMMA(a+1/3-1/2*
n)*GAMMA(a-1/2*n)*GAMMA(2*a+1/3)*GAMMA(2*a)-1/8*(sin(1/2*Pi*(n-2*a))-s
in(1/6*Pi*(-6*a-2+3*n))-2*sin(1/2*Pi*(6*a+n))+sin(3/2*Pi*(n-2*a))-sin(
1/6*Pi*(-6*a+2+3*n)))*V8(a,n)*GAMMA(a-1/3-1/2*n)*GAMMA(7/6-a-1/2*n)*GA
MMA(a+1/3-1/2*n)*GAMMA(a-1/2*n)*GAMMA(2*a+1/3)/sin(1/3*Pi*(6*a-1))/Pi^
(3/2)/GAMMA(2*a+1/6-n)/GAMMA(2*a+5/6-n);}{%
\maplemultiline{
326\mbox{:~~~} {\displaystyle \frac {1}{2}} \mathrm{sin}
({\displaystyle \frac {\pi \,( - 6\,a + n)}{2}} )\,\mathrm{V7}(a
 - {\displaystyle \frac {1}{3}}  - {\displaystyle \frac {n}{2}} 
, \,n + 1)\,\Gamma (a - {\displaystyle \frac {1}{3}}  - 
{\displaystyle \frac {n}{2}} )\,\Gamma ({\displaystyle \frac {7}{
6}}  - a - {\displaystyle \frac {n}{2}} ) \\
\Gamma (a + {\displaystyle \frac {1}{3}}  - {\displaystyle 
\frac {n}{2}} )\,\Gamma (a - {\displaystyle \frac {n}{2}} )\,
\Gamma (2\,a + {\displaystyle \frac {1}{3}} )\,\Gamma (2\,a)
 \left/ {\vrule height0.80em width0em depth0.80em} \right. \! 
 \! (\pi \,\mathrm{sin}({\displaystyle \frac {\pi \,(6\,a - 1)}{3
}} )\,(-1)^{( - n)} \\
\Gamma (3\,a - {\displaystyle \frac {n}{2}} )\,\Gamma (
{\displaystyle \frac {1}{2}}  + a - {\displaystyle \frac {n}{2}} 
)\,\Gamma (a - {\displaystyle \frac {n}{2}}  + {\displaystyle 
\frac {5}{6}} )\,\Gamma (a + {\displaystyle \frac {1}{6}}  - 
{\displaystyle \frac {n}{2}} ))\mbox{} - {\displaystyle \frac {1
}{8}} (\mathrm{sin}({\displaystyle \frac {\pi \,(n - 2\,a)}{2}} )
 \\
\mbox{} - \mathrm{sin}({\displaystyle \frac {\pi \,( - 6\,a - 2
 + 3\,n)}{6}} ) - 2\,\mathrm{sin}({\displaystyle \frac {\pi \,(6
\,a + n)}{2}} ) + \mathrm{sin}({\displaystyle \frac {3\,\pi \,(n
 - 2\,a)}{2}} ) \\
\mbox{} - \mathrm{sin}({\displaystyle \frac {\pi \,( - 6\,a + 2
 + 3\,n)}{6}} ))\mathrm{V8}(a, \,n)\,\Gamma (a - {\displaystyle 
\frac {1}{3}}  - {\displaystyle \frac {n}{2}} )\,\Gamma (
{\displaystyle \frac {7}{6}}  - a - {\displaystyle \frac {n}{2}} 
)\,\Gamma (a + {\displaystyle \frac {1}{3}}  - {\displaystyle 
\frac {n}{2}} ) \\
\Gamma (a - {\displaystyle \frac {n}{2}} )\,\Gamma (2\,a + 
{\displaystyle \frac {1}{3}} ) \left/ {\vrule 
height0.80em width0em depth0.80em} \right. \!  \! (\mathrm{sin}(
{\displaystyle \frac {\pi \,(6\,a - 1)}{3}} )\,\pi ^{(3/2)}\,
\Gamma (2\,a + {\displaystyle \frac {1}{6}}  - n)\,\Gamma (2\,a
 + {\displaystyle \frac {5}{6}}  - n)) }
}
\end{maplelatex}

\begin{maplelatex}
\mapleinline{inert}{2d}{327, ":   ",
(3/2*a-1/4*n-1/2)/GAMMA(a-1/2*n+5/6)/GAMMA(1/2+a-1/2*n)/GAMMA(a+1/6-1/
2*n)*GAMMA(a-1/2*n)*GAMMA(2*a+1/3)*GAMMA(2*a-1/3)*GAMMA(a-1/3-1/2*n)*G
AMMA(a+1/3-1/2*n)/GAMMA(3*a-1/2*n)^2/GAMMA(1/2+n)*GAMMA(2*a)*V7(a-1/3-
1/2*n,n+1)+1/8*V8(a,n)*(2*sin(1/2*Pi*(6*a+n))-sin(3/2*Pi*(n-2*a)))*(-1
)^(-n)*GAMMA(a-1/3-1/2*n)*GAMMA(a+1/3-1/2*n)*GAMMA(a-1/2*n)*GAMMA(2*a+
1/3)*GAMMA(2*a-1/3)/sin(1/2*Pi*(-6*a+n))/Pi^(1/2)/GAMMA(3*a-1/2*n-1)/G
AMMA(2*a+5/6-n)/GAMMA(2*a+1/6-n)/GAMMA(1/2+n);}{%
\maplemultiline{
327\mbox{:~~~} ({\displaystyle \frac {3\,a}{2}}  - 
{\displaystyle \frac {n}{4}}  - {\displaystyle \frac {1}{2}} )\,
\Gamma (a - {\displaystyle \frac {n}{2}} )\,\Gamma (2\,a + 
{\displaystyle \frac {1}{3}} )\,\Gamma (2\,a - {\displaystyle 
\frac {1}{3}} )\,\Gamma (a - {\displaystyle \frac {1}{3}}  - 
{\displaystyle \frac {n}{2}} )\,\Gamma (a + {\displaystyle 
\frac {1}{3}}  - {\displaystyle \frac {n}{2}} )\,\Gamma (2\,a)
 \\
\mathrm{V7}(a - {\displaystyle \frac {1}{3}}  - {\displaystyle 
\frac {n}{2}} , \,n + 1) \left/ {\vrule 
height0.80em width0em depth0.80em} \right. \!  \! (\Gamma (a - 
{\displaystyle \frac {n}{2}}  + {\displaystyle \frac {5}{6}} )\,
\Gamma ({\displaystyle \frac {1}{2}}  + a - {\displaystyle 
\frac {n}{2}} )\,\Gamma (a + {\displaystyle \frac {1}{6}}  - 
{\displaystyle \frac {n}{2}} )\,\Gamma (3\,a - {\displaystyle 
\frac {n}{2}} )^{2} \\
\Gamma ({\displaystyle \frac {1}{2}}  + n))\mbox{} + 
{\displaystyle \frac {1}{8}} \mathrm{V8}(a, \,n)\,(2\,\mathrm{sin
}({\displaystyle \frac {\pi \,(6\,a + n)}{2}} ) - \mathrm{sin}(
{\displaystyle \frac {3\,\pi \,(n - 2\,a)}{2}} ))\,(-1)^{( - n)}
 \\
\Gamma (a - {\displaystyle \frac {1}{3}}  - {\displaystyle 
\frac {n}{2}} )\,\Gamma (a + {\displaystyle \frac {1}{3}}  - 
{\displaystyle \frac {n}{2}} )\,\Gamma (a - {\displaystyle 
\frac {n}{2}} )\,\Gamma (2\,a + {\displaystyle \frac {1}{3}} )\,
\Gamma (2\,a - {\displaystyle \frac {1}{3}} ) \left/ {\vrule 
height0.80em width0em depth0.80em} \right. \!  \! (\mathrm{sin}(
{\displaystyle \frac {\pi \,( - 6\,a + n)}{2}} ) \\
\sqrt{\pi }\,\Gamma (3\,a - {\displaystyle \frac {n}{2}}  - 1)\,
\Gamma (2\,a + {\displaystyle \frac {5}{6}}  - n)\,\Gamma (2\,a
 + {\displaystyle \frac {1}{6}}  - n)\,\Gamma ({\displaystyle 
\frac {1}{2}}  + n)) }
}
\end{maplelatex}

\begin{maplelatex}
\mapleinline{inert}{2d}{328, ":   ",
-2/9*cos(Pi*(n-2*a))*3^(1/2)*V7(a-1/3-1/2*n,n+1)/Pi^(1/2)/GAMMA(3*a-1/
2*n)*GAMMA(1+a+1/2*n)*GAMMA(1/6-a+1/2*n)/GAMMA(1/2+n)/GAMMA(1/2+a-1/2*
n)*GAMMA(a-1/3-1/2*n)*GAMMA(a+1/3-1/2*n)/GAMMA(2/3)*GAMMA(2*a)-1/9*V8(
a,n)*(sin(1/6*Pi*(6*a-1+3*n))-sin(1/6*Pi*(9*n-6*a+1)))*3^(1/2)*Pi*(-1)
^(-n)*GAMMA(1+a+1/2*n)*GAMMA(a-1/3-1/2*n)*GAMMA(a+1/3-1/2*n)/sin(1/6*P
i*(1-6*a+3*n))/GAMMA(1/2+n)/GAMMA(2*a+5/6-n)/GAMMA(2*a+1/6-n)/GAMMA(-a
+5/6+1/2*n)/GAMMA(2/3);}{%
\maplemultiline{
328\mbox{:~~~}  - {\displaystyle \frac {2}{9}} \mathrm{
cos}(\pi \,(n - 2\,a))\,\sqrt{3}\,\mathrm{V7}(a - {\displaystyle 
\frac {1}{3}}  - {\displaystyle \frac {n}{2}} , \,n + 1)\,\Gamma 
(1 + a + {\displaystyle \frac {n}{2}} )\,\Gamma ({\displaystyle 
\frac {1}{6}}  - a + {\displaystyle \frac {n}{2}} ) \\
\Gamma (a - {\displaystyle \frac {1}{3}}  - {\displaystyle 
\frac {n}{2}} )\,\Gamma (a + {\displaystyle \frac {1}{3}}  - 
{\displaystyle \frac {n}{2}} )\,\Gamma (2\,a) \left/ {\vrule 
height0.80em width0em depth0.80em} \right. \!  \! (\sqrt{\pi }\,
\Gamma (3\,a - {\displaystyle \frac {n}{2}} )\,\Gamma (
{\displaystyle \frac {1}{2}}  + n)\,\Gamma ({\displaystyle 
\frac {1}{2}}  + a - {\displaystyle \frac {n}{2}} ) \\
\Gamma ({\displaystyle \frac {2}{3}} ))\mbox{} - {\displaystyle 
\frac {1}{9}} \mathrm{V8}(a, \,n)\,(\mathrm{sin}({\displaystyle 
\frac {\pi \,(6\,a - 1 + 3\,n)}{6}} ) - \mathrm{sin}(
{\displaystyle \frac {\pi \,(9\,n - 6\,a + 1)}{6}} ))\,\sqrt{3}\,
\pi \,(-1)^{( - n)} \\
\Gamma (1 + a + {\displaystyle \frac {n}{2}} )\,\Gamma (a - 
{\displaystyle \frac {1}{3}}  - {\displaystyle \frac {n}{2}} )\,
\Gamma (a + {\displaystyle \frac {1}{3}}  - {\displaystyle 
\frac {n}{2}} ) \left/ {\vrule height0.80em width0em depth0.80em}
 \right. \!  \! (\mathrm{sin}({\displaystyle \frac {\pi \,(1 - 6
\,a + 3\,n)}{6}} )\,\Gamma ({\displaystyle \frac {1}{2}}  + n)
 \\
\Gamma (2\,a + {\displaystyle \frac {5}{6}}  - n)\,\Gamma (2\,a
 + {\displaystyle \frac {1}{6}}  - n)\,\Gamma ( - a + 
{\displaystyle \frac {5}{6}}  + {\displaystyle \frac {n}{2}} )\,
\Gamma ({\displaystyle \frac {2}{3}} )) }
}
\end{maplelatex}

\begin{maplelatex}
\mapleinline{inert}{2d}{329, ":   ",
2/9*cos(Pi*(n-2*a))*3^(1/2)*Pi^(1/2)*V7(a-1/3-1/2*n,n+1)/sin(1/6*Pi*(-
6*a-2+3*n))*GAMMA(a+1/2+1/2*n)/GAMMA(3*a-1/2*n)/GAMMA(1/2+n)/GAMMA(1/2
+a-1/2*n)*GAMMA(a-1/3-1/2*n)/GAMMA(2/3)*GAMMA(2*a)+1/9*V8(a,n)*(sin(1/
6*Pi*(6*a-1+3*n))-sin(1/6*Pi*(9*n-6*a+1)))*3^(1/2)*(-1)^(-n)*Pi*GAMMA(
a-1/2*n+5/6)*GAMMA(a+1/2+1/2*n)*GAMMA(a-1/3-1/2*n)/sin(1/6*Pi*(-6*a-2+
3*n))/GAMMA(2/3)/GAMMA(2*a+5/6-n)/GAMMA(2*a+1/6-n)/GAMMA(-a+5/6+1/2*n)
/GAMMA(1/2+n);}{%
\maplemultiline{
329\mbox{:~~~} 
{\displaystyle \frac {2}{9}} \,{\displaystyle \frac {\mathrm{cos}
(\pi \,(n - 2\,a))\,\sqrt{3}\,\sqrt{\pi }\,\mathrm{V7}(a - 
{\displaystyle \frac {1}{3}}  - {\displaystyle \frac {n}{2}} , \,
n + 1)\,\Gamma (a + {\displaystyle \frac {1}{2}}  + 
{\displaystyle \frac {n}{2}} )\,\Gamma (a - {\displaystyle 
\frac {1}{3}}  - {\displaystyle \frac {n}{2}} )\,\Gamma (2\,a)}{
\mathrm{sin}({\displaystyle \frac {\pi \,( - 6\,a - 2 + 3\,n)}{6}
} )\,\Gamma (3\,a - {\displaystyle \frac {n}{2}} )\,\Gamma (
{\displaystyle \frac {1}{2}}  + n)\,\Gamma ({\displaystyle 
\frac {1}{2}}  + a - {\displaystyle \frac {n}{2}} )\,\Gamma (
{\displaystyle \frac {2}{3}} )}}  \\
\mbox{} + {\displaystyle \frac {1}{9}} \mathrm{V8}(a, \,n)\,(
\mathrm{sin}({\displaystyle \frac {\pi \,(6\,a - 1 + 3\,n)}{6}} )
 - \mathrm{sin}({\displaystyle \frac {\pi \,(9\,n - 6\,a + 1)}{6}
} ))\,\sqrt{3}\,(-1)^{( - n)}\,\pi  \\
\Gamma (a - {\displaystyle \frac {n}{2}}  + {\displaystyle 
\frac {5}{6}} )\,\Gamma (a + {\displaystyle \frac {1}{2}}  + 
{\displaystyle \frac {n}{2}} )\,\Gamma (a - {\displaystyle 
\frac {1}{3}}  - {\displaystyle \frac {n}{2}} ) \left/ {\vrule 
height0.80em width0em depth0.80em} \right. \!  \! (\mathrm{sin}(
{\displaystyle \frac {\pi \,( - 6\,a - 2 + 3\,n)}{6}} )\,\Gamma (
{\displaystyle \frac {2}{3}} ) \\
\Gamma (2\,a + {\displaystyle \frac {5}{6}}  - n)\,\Gamma (2\,a
 + {\displaystyle \frac {1}{6}}  - n)\,\Gamma ( - a + 
{\displaystyle \frac {5}{6}}  + {\displaystyle \frac {n}{2}} )\,
\Gamma ({\displaystyle \frac {1}{2}}  + n)) }
}
\end{maplelatex}

\begin{maplelatex}
\mapleinline{inert}{2d}{330, ":   ",
2/9*3^(1/2)*Pi^(1/2)*V7(a-1/3-1/2*n,n+1)/GAMMA(3*a-1/2*n)/GAMMA(1/2+n)
/GAMMA(1/2+a-1/2*n)/GAMMA(2*a-1/2-n)*GAMMA(a-1/3-1/2*n)*GAMMA(a+1/3-1/
2*n)*GAMMA(2*a-1/3)/GAMMA(2/3)*GAMMA(2*a)+1/9*V8(a,n)*(sin(1/6*Pi*(6*a
-1+3*n))-sin(1/6*Pi*(9*n-6*a+1)))*3^(1/2)*(-1)^(-n)*Pi^(3/2)*GAMMA(2*a
+2/3-n)*GAMMA(a-1/3-1/2*n)*GAMMA(2*a-1/3)/cos(Pi*(n-2*a))/(2^(-1/3+2*a
-n))/GAMMA(2*a-1/2-n)/GAMMA(2*a+5/6-n)/GAMMA(2*a+1/6-n)/GAMMA(-a+5/6+1
/2*n)/GAMMA(1/2+n)/GAMMA(2/3);}{%
\maplemultiline{
330\mbox{:~~~} {\displaystyle \frac {2}{9}} \,
{\displaystyle \frac {\sqrt{3}\,\sqrt{\pi }\,\mathrm{V7}(a - 
{\displaystyle \frac {1}{3}}  - {\displaystyle \frac {n}{2}} , \,
n + 1)\,\Gamma (a - {\displaystyle \frac {1}{3}}  - 
{\displaystyle \frac {n}{2}} )\,\Gamma (a + {\displaystyle 
\frac {1}{3}}  - {\displaystyle \frac {n}{2}} )\,\Gamma (2\,a - 
{\displaystyle \frac {1}{3}} )\,\Gamma (2\,a)}{\Gamma (3\,a - 
{\displaystyle \frac {n}{2}} )\,\Gamma ({\displaystyle \frac {1}{
2}}  + n)\,\Gamma ({\displaystyle \frac {1}{2}}  + a - 
{\displaystyle \frac {n}{2}} )\,\Gamma (2\,a - {\displaystyle 
\frac {1}{2}}  - n)\,\Gamma ({\displaystyle \frac {2}{3}} )}}  + 
{\displaystyle \frac {1}{9}}  \\
\mathrm{V8}(a, \,n)\,(\mathrm{sin}({\displaystyle \frac {\pi \,(6
\,a - 1 + 3\,n)}{6}} ) - \mathrm{sin}({\displaystyle \frac {\pi 
\,(9\,n - 6\,a + 1)}{6}} ))\,\sqrt{3}\,(-1)^{( - n)}\,\pi ^{(3/2)
} \\
\Gamma (2\,a + {\displaystyle \frac {2}{3}}  - n)\,\Gamma (a - 
{\displaystyle \frac {1}{3}}  - {\displaystyle \frac {n}{2}} )\,
\Gamma (2\,a - {\displaystyle \frac {1}{3}} ) \left/ {\vrule 
height0.80em width0em depth0.80em} \right. \!  \! (\mathrm{cos}(
\pi \,(n - 2\,a))\,2^{( - 1/3 + 2\,a - n)} \\
\Gamma (2\,a - {\displaystyle \frac {1}{2}}  - n)\,\Gamma (2\,a
 + {\displaystyle \frac {5}{6}}  - n)\,\Gamma (2\,a + 
{\displaystyle \frac {1}{6}}  - n)\,\Gamma ( - a + 
{\displaystyle \frac {5}{6}}  + {\displaystyle \frac {n}{2}} )\,
\Gamma ({\displaystyle \frac {1}{2}}  + n)\,\Gamma (
{\displaystyle \frac {2}{3}} )) }
}
\end{maplelatex}

\begin{maplelatex}
\mapleinline{inert}{2d}{331, ":   ",
-2/9*cos(Pi*(n-2*a))*3^(1/2)*V7(a-1/3-1/2*n,n+1)/Pi^(1/2)/GAMMA(3*a-1/
2*n)/GAMMA(1/2+a-1/2*n)*GAMMA(a-1/3-1/2*n)*GAMMA(a+1/3-1/2*n)*GAMMA(2*
a)-1/9*(-1)^(-n)*3^(1/2)*(sin(1/6*Pi*(6*a-1+3*n))-sin(1/6*Pi*(9*n-6*a+
1)))*V8(a,n)*Pi^(1/2)*GAMMA(2*a+2/3-n)*GAMMA(a-1/3-1/2*n)/(2^(-1/3+2*a
-n))/GAMMA(2*a+5/6-n)/GAMMA(2*a+1/6-n)/GAMMA(-a+5/6+1/2*n);}{%
\maplemultiline{
331\mbox{:~~~}  - {\displaystyle \frac {2}{9}} \,
{\displaystyle \frac {\mathrm{cos}(\pi \,(n - 2\,a))\,\sqrt{3}\,
\mathrm{V7}(a - {\displaystyle \frac {1}{3}}  - {\displaystyle 
\frac {n}{2}} , \,n + 1)\,\Gamma (a - {\displaystyle \frac {1}{3}
}  - {\displaystyle \frac {n}{2}} )\,\Gamma (a + {\displaystyle 
\frac {1}{3}}  - {\displaystyle \frac {n}{2}} )\,\Gamma (2\,a)}{
\sqrt{\pi }\,\Gamma (3\,a - {\displaystyle \frac {n}{2}} )\,
\Gamma ({\displaystyle \frac {1}{2}}  + a - {\displaystyle 
\frac {n}{2}} )}}  \\
\mbox{} - {\displaystyle \frac {1}{9}} (-1)^{( - n)}\,\sqrt{3}\,(
\mathrm{sin}({\displaystyle \frac {\pi \,(6\,a - 1 + 3\,n)}{6}} )
 - \mathrm{sin}({\displaystyle \frac {\pi \,(9\,n - 6\,a + 1)}{6}
} ))\,\mathrm{V8}(a, \,n)\,\sqrt{\pi } \\
\Gamma (2\,a + {\displaystyle \frac {2}{3}}  - n)\,\Gamma (a - 
{\displaystyle \frac {1}{3}}  - {\displaystyle \frac {n}{2}} )
 \left/ {\vrule height0.80em width0em depth0.80em} \right. \! 
 \! (2^{( - 1/3 + 2\,a - n)}\,\Gamma (2\,a + {\displaystyle 
\frac {5}{6}}  - n)\,\Gamma (2\,a + {\displaystyle \frac {1}{6}} 
 - n) 
\Gamma ( - a + {\displaystyle \frac {5}{6}}  + {\displaystyle 
\frac {n}{2}} )) }
}
\end{maplelatex}

\begin{maplelatex}
\mapleinline{inert}{2d}{332, ":   ",
-2*cos(Pi*(n-2*a))*V7(a-1/3-1/2*n,n+1)/sin(2*Pi*a)/(2^(2*a+1+n))*GAMMA
(2*a+1+n)/GAMMA(3*a-1/2*n)/GAMMA(1/2+n)/GAMMA(1/2+a-1/2*n)*GAMMA(a-1/3
-1/2*n)*GAMMA(a+1/3-1/2*n)-(-1)^(-n)*(sin(1/6*Pi*(6*a-1+3*n))-sin(1/6*
Pi*(9*n-6*a+1)))*V8(a,n)*Pi*GAMMA(2*a+2/3-n)*GAMMA(a-1/3-1/2*n)*GAMMA(
2*a+1+n)/sin(2*Pi*a)/(2^(-1/3+2*a-n))/(2^(2*a+1+n))/GAMMA(2*a+5/6-n)/G
AMMA(2*a+1/6-n)/GAMMA(-a+5/6+1/2*n)/GAMMA(1/2+n)/GAMMA(2*a);}{%
\maplemultiline{
332\mbox{:~~~} 
 - {\displaystyle \frac {2\,\mathrm{cos}(\pi \,(n - 2\,a))\,
\mathrm{V7}(a - {\displaystyle \frac {1}{3}}  - {\displaystyle 
\frac {n}{2}} , \,n + 1)\,\Gamma (2\,a + 1 + n)\,\Gamma (a - 
{\displaystyle \frac {1}{3}}  - {\displaystyle \frac {n}{2}} )\,
\Gamma (a + {\displaystyle \frac {1}{3}}  - {\displaystyle 
\frac {n}{2}} )}{\mathrm{sin}(2\,\pi \,a)\,2^{(2\,a + 1 + n)}\,
\Gamma (3\,a - {\displaystyle \frac {n}{2}} )\,\Gamma (
{\displaystyle \frac {1}{2}}  + n)\,\Gamma ({\displaystyle 
\frac {1}{2}}  + a - {\displaystyle \frac {n}{2}} )}}  \\
\mbox{} - (-1)^{( - n)}\,(\mathrm{sin}({\displaystyle \frac {\pi 
\,(6\,a - 1 + 3\,n)}{6}} ) - \mathrm{sin}({\displaystyle \frac {
\pi \,(9\,n - 6\,a + 1)}{6}} ))\,\mathrm{V8}(a, \,n)\,\pi  \\
\Gamma (2\,a + {\displaystyle \frac {2}{3}}  - n)\,\Gamma (a - 
{\displaystyle \frac {1}{3}}  - {\displaystyle \frac {n}{2}} )\,
\Gamma (2\,a + 1 + n) \left/ {\vrule 
height0.80em width0em depth0.80em} \right. \!  \! (\mathrm{sin}(2
\,\pi \,a)\,2^{( - 1/3 + 2\,a - n)} \\
2^{(2\,a + 1 + n)}\,\Gamma (2\,a + {\displaystyle \frac {5}{6}} 
 - n)\,\Gamma (2\,a + {\displaystyle \frac {1}{6}}  - n)\,\Gamma 
( - a + {\displaystyle \frac {5}{6}}  + {\displaystyle \frac {n}{
2}} )\,\Gamma ({\displaystyle \frac {1}{2}}  + n)\,\Gamma (2\,a))
 }
}
\end{maplelatex}

\begin{maplelatex}
\mapleinline{inert}{2d}{333, ":   ",
V7(a-1/3-1/2*n,n+1)/Pi^(1/2)*GAMMA(1+a+1/2*n)*GAMMA(a-1/2*n)/GAMMA(3*a
-1/2*n)/GAMMA(1/2+n)/GAMMA(1/2+a-1/2*n)/GAMMA(2*a-1/2-n)*GAMMA(a-1/3-1
/2*n)*GAMMA(a+1/3-1/2*n)*GAMMA(2*a)+1/2*V8(a,n)*(sin(1/6*Pi*(6*a-1+3*n
))-sin(1/6*Pi*(9*n-6*a+1)))*(-1)^(-n)*Pi^(1/2)*GAMMA(2*a+2/3-n)*GAMMA(
1+a+1/2*n)*GAMMA(a-1/2*n)*GAMMA(a-1/3-1/2*n)/cos(Pi*(n-2*a))/(2^(-1/3+
2*a-n))/GAMMA(2*a-1/2-n)/GAMMA(2*a+1/6-n)/GAMMA(2*a+5/6-n)/GAMMA(-a+5/
6+1/2*n)/GAMMA(1/2+n);}{%
\maplemultiline{
333\mbox{:~~~} {\displaystyle \frac {\mathrm{V7}(a - 
{\displaystyle \frac {1}{3}}  - {\displaystyle \frac {n}{2}} , \,
n + 1)\,\Gamma (1 + a + {\displaystyle \frac {n}{2}} )\,\Gamma (a
 - {\displaystyle \frac {n}{2}} )\,\Gamma (a - {\displaystyle 
\frac {1}{3}}  - {\displaystyle \frac {n}{2}} )\,\Gamma (a + 
{\displaystyle \frac {1}{3}}  - {\displaystyle \frac {n}{2}} )\,
\Gamma (2\,a)}{\sqrt{\pi }\,\Gamma (3\,a - {\displaystyle \frac {
n}{2}} )\,\Gamma ({\displaystyle \frac {1}{2}}  + n)\,\Gamma (
{\displaystyle \frac {1}{2}}  + a - {\displaystyle \frac {n}{2}} 
)\,\Gamma (2\,a - {\displaystyle \frac {1}{2}}  - n)}}  \\
\mbox{} + {\displaystyle \frac {1}{2}} \mathrm{V8}(a, \,n)\,(
\mathrm{sin}({\displaystyle \frac {\pi \,(6\,a - 1 + 3\,n)}{6}} )
 - \mathrm{sin}({\displaystyle \frac {\pi \,(9\,n - 6\,a + 1)}{6}
} ))\,(-1)^{( - n)}\,\sqrt{\pi } \\
\Gamma (2\,a + {\displaystyle \frac {2}{3}}  - n)\,\Gamma (1 + a
 + {\displaystyle \frac {n}{2}} )\,\Gamma (a - {\displaystyle 
\frac {n}{2}} )\,\Gamma (a - {\displaystyle \frac {1}{3}}  - 
{\displaystyle \frac {n}{2}} ) \left/ {\vrule 
height0.80em width0em depth0.80em} \right. \!  \! (\mathrm{cos}(
\pi \,(n - 2\,a)) \\
2^{( - 1/3 + 2\,a - n)}\,\Gamma (2\,a - {\displaystyle \frac {1}{
2}}  - n)\,\Gamma (2\,a + {\displaystyle \frac {1}{6}}  - n)\,
\Gamma (2\,a + {\displaystyle \frac {5}{6}}  - n)\,\Gamma ( - a
 + {\displaystyle \frac {5}{6}}  + {\displaystyle \frac {n}{2}} )
\,\Gamma ({\displaystyle \frac {1}{2}}  + n)) }
}
\end{maplelatex}

\begin{maplelatex}
\mapleinline{inert}{2d}{334, ":   ",
-cos(Pi*(n-2*a))*V7(a-1/3-1/2*n,n+1)/Pi^(3/2)*GAMMA(1+a+1/2*n)/GAMMA(3
*a-1/2*n)*GAMMA(-a+5/6+1/2*n)*GAMMA(a-1/3-1/2*n)/GAMMA(1/2+n)/GAMMA(1/
2+a-1/2*n)*GAMMA(a+1/3-1/2*n)*GAMMA(2/3)*GAMMA(2*a)-1/2*(-1)^(-n)*(sin
(1/6*Pi*(6*a-1+3*n))-sin(1/6*Pi*(9*n-6*a+1)))*V8(a,n)*GAMMA(2*a+2/3-n)
*GAMMA(1+a+1/2*n)*GAMMA(a-1/3-1/2*n)*GAMMA(2/3)/Pi^(1/2)/(2^(-1/3+2*a-
n))/GAMMA(2*a+5/6-n)/GAMMA(2*a+1/6-n)/GAMMA(1/2+n);}{%
\maplemultiline{
334\mbox{:~~~}  - \mathrm{cos}(\pi \,(n - 2\,a))\,
\mathrm{V7}(a - {\displaystyle \frac {1}{3}}  - {\displaystyle 
\frac {n}{2}} , \,n + 1)\,\Gamma (1 + a + {\displaystyle \frac {n
}{2}} )\,\Gamma ( - a + {\displaystyle \frac {5}{6}}  + 
{\displaystyle \frac {n}{2}} ) \\
\Gamma (a - {\displaystyle \frac {1}{3}}  - {\displaystyle 
\frac {n}{2}} )\,\Gamma (a + {\displaystyle \frac {1}{3}}  - 
{\displaystyle \frac {n}{2}} )\,\Gamma ({\displaystyle \frac {2}{
3}} )\,\Gamma (2\,a) \left/ {\vrule 
height0.80em width0em depth0.80em} \right. \!  \! (\pi ^{(3/2)}\,
\Gamma (3\,a - {\displaystyle \frac {n}{2}} )\,\Gamma (
{\displaystyle \frac {1}{2}}  + n) \\
\Gamma ({\displaystyle \frac {1}{2}}  + a - {\displaystyle 
\frac {n}{2}} ))\mbox{} - {\displaystyle \frac {1}{2}} (-1)^{( - 
n)}\,(\mathrm{sin}({\displaystyle \frac {\pi \,(6\,a - 1 + 3\,n)
}{6}} ) - \mathrm{sin}({\displaystyle \frac {\pi \,(9\,n - 6\,a
 + 1)}{6}} )) \\
\mathrm{V8}(a, \,n)\,\Gamma (2\,a + {\displaystyle \frac {2}{3}} 
 - n)\,\Gamma (1 + a + {\displaystyle \frac {n}{2}} )\,\Gamma (a
 - {\displaystyle \frac {1}{3}}  - {\displaystyle \frac {n}{2}} )
\,\Gamma ({\displaystyle \frac {2}{3}} ) \left/ {\vrule 
height0.80em width0em depth0.80em} \right. \!  \! (\sqrt{\pi }\,2
^{( - 1/3 + 2\,a - n)} \\
\Gamma (2\,a + {\displaystyle \frac {5}{6}}  - n)\,\Gamma (2\,a
 + {\displaystyle \frac {1}{6}}  - n)\,\Gamma ({\displaystyle 
\frac {1}{2}}  + n)) }
}
\end{maplelatex}

\begin{maplelatex}
\mapleinline{inert}{2d}{335, ":   ",
cos(Pi*(n-2*a))*V7(a-1/3-1/2*n,n+1)/sin(1/6*Pi*(-6*a+2+3*n))/Pi^(1/2)/
GAMMA(3*a-1/2*n)*GAMMA(a+1/2+1/2*n)/GAMMA(1/2+n)*GAMMA(a+1/3-1/2*n)/GA
MMA(1/2+a-1/2*n)*GAMMA(2/3)*GAMMA(2*a)-1/2*(-1)^(-n)*(-sin(1/6*Pi*(6*a
-1+3*n))+sin(1/6*Pi*(9*n-6*a+1)))*V8(a,n)*Pi^(1/2)*GAMMA(2*a+2/3-n)*GA
MMA(a+1/2+1/2*n)*GAMMA(2/3)/sin(1/6*Pi*(-6*a+2+3*n))/(2^(-1/3+2*a-n))/
GAMMA(1/2+n)/GAMMA(2*a+5/6-n)/GAMMA(2*a+1/6-n)/GAMMA(-a+5/6+1/2*n);}{%
\maplemultiline{
335\mbox{:~~~} {\displaystyle \frac {\mathrm{cos}(\pi \,
(n - 2\,a))\,\mathrm{V7}(a - {\displaystyle \frac {1}{3}}  - 
{\displaystyle \frac {n}{2}} , \,n + 1)\,\Gamma (a + 
{\displaystyle \frac {1}{2}}  + {\displaystyle \frac {n}{2}} )\,
\Gamma (a + {\displaystyle \frac {1}{3}}  - {\displaystyle 
\frac {n}{2}} )\,\Gamma ({\displaystyle \frac {2}{3}} )\,\Gamma (
2\,a)}{\mathrm{sin}({\displaystyle \frac {\pi \,( - 6\,a + 2 + 3
\,n)}{6}} )\,\sqrt{\pi }\,\Gamma (3\,a - {\displaystyle \frac {n
}{2}} )\,\Gamma ({\displaystyle \frac {1}{2}}  + n)\,\Gamma (
{\displaystyle \frac {1}{2}}  + a - {\displaystyle \frac {n}{2}} 
)}}  -  \\
{\displaystyle \frac {1}{2}} (-1)^{( - n)}\,( - \mathrm{sin}(
{\displaystyle \frac {\pi \,(6\,a - 1 + 3\,n)}{6}} ) + \mathrm{
sin}({\displaystyle \frac {\pi \,(9\,n - 6\,a + 1)}{6}} ))\,
\mathrm{V8}(a, \,n)\,\sqrt{\pi } \\
\Gamma (2\,a + {\displaystyle \frac {2}{3}}  - n)\,\Gamma (a + 
{\displaystyle \frac {1}{2}}  + {\displaystyle \frac {n}{2}} )\,
\Gamma ({\displaystyle \frac {2}{3}} ) \left/ {\vrule 
height0.80em width0em depth0.80em} \right. \!  \! (\mathrm{sin}(
{\displaystyle \frac {\pi \,( - 6\,a + 2 + 3\,n)}{6}} )\,2^{( - 1
/3 + 2\,a - n)} \\
\Gamma ({\displaystyle \frac {1}{2}}  + n)\,\Gamma (2\,a + 
{\displaystyle \frac {5}{6}}  - n)\,\Gamma (2\,a + 
{\displaystyle \frac {1}{6}}  - n)\,\Gamma ( - a + 
{\displaystyle \frac {5}{6}}  + {\displaystyle \frac {n}{2}} ))
 }
}
\end{maplelatex}

\begin{maplelatex}
\mapleinline{inert}{2d}{336, ":   ",
2/3*sin(1/6*Pi*(-12*a-1+6*n))*V7(a-1/3-1/2*n,n+1)/Pi^(3/2)*GAMMA(a-1/2
*n)*GAMMA(a-1/3-1/2*n)*GAMMA(4/3+a+1/2*n)/GAMMA(1/2+n)/GAMMA(a-1/2*n+5
/6)/GAMMA(3*a-1/2*n)*GAMMA(1/2-a+1/2*n)*GAMMA(2*a-1/3)*GAMMA(2/3)-1/3*
(2*sin(1/6*Pi*(-12*a-1+6*n))-1)*V8(a,n)*GAMMA(2*a-2/3-n)*GAMMA(a-1/2*n
)*GAMMA(4/3+a+1/2*n)*GAMMA(2*a-1/3)*GAMMA(2/3)/cos(1/2*Pi*(n-2*a))/Pi^
(1/2)/(2^(2*a-2/3-n))/GAMMA(2*a+1/6-n)/GAMMA(2*a+5/6-n)/GAMMA(1/2+n)/G
AMMA(2*a);}{%
\maplemultiline{
336\mbox{:~~~} {\displaystyle \frac {2}{3}} \mathrm{sin}
({\displaystyle \frac {\pi \,( - 12\,a - 1 + 6\,n)}{6}} )\,
\mathrm{V7}(a - {\displaystyle \frac {1}{3}}  - {\displaystyle 
\frac {n}{2}} , \,n + 1)\,\Gamma (a - {\displaystyle \frac {n}{2}
} )\,\Gamma (a - {\displaystyle \frac {1}{3}}  - {\displaystyle 
\frac {n}{2}} ) \\
\Gamma ({\displaystyle \frac {4}{3}}  + a + {\displaystyle 
\frac {n}{2}} )\,\Gamma ({\displaystyle \frac {1}{2}}  - a + 
{\displaystyle \frac {n}{2}} )\,\Gamma (2\,a - {\displaystyle 
\frac {1}{3}} )\,\Gamma ({\displaystyle \frac {2}{3}} ) \left/ 
{\vrule height0.80em width0em depth0.80em} \right. \!  \! (\pi ^{
(3/2)}\,\Gamma ({\displaystyle \frac {1}{2}}  + n)\,\Gamma (a - 
{\displaystyle \frac {n}{2}}  + {\displaystyle \frac {5}{6}} )
 \\
\Gamma (3\,a - {\displaystyle \frac {n}{2}} ))\mbox{} - 
{\displaystyle \frac {1}{3}} (2\,\mathrm{sin}({\displaystyle 
\frac {\pi \,( - 12\,a - 1 + 6\,n)}{6}} ) - 1)\,\mathrm{V8}(a, \,
n)\,\Gamma (2\,a - {\displaystyle \frac {2}{3}}  - n)\,\Gamma (a
 - {\displaystyle \frac {n}{2}} ) \\
\Gamma ({\displaystyle \frac {4}{3}}  + a + {\displaystyle 
\frac {n}{2}} )\,\Gamma (2\,a - {\displaystyle \frac {1}{3}} )\,
\Gamma ({\displaystyle \frac {2}{3}} ) \left/ {\vrule 
height0.80em width0em depth0.80em} \right. \!  \! (\mathrm{cos}(
{\displaystyle \frac {\pi \,(n - 2\,a)}{2}} )\,\sqrt{\pi }\,2^{(2
\,a - 2/3 - n)} \\
\Gamma (2\,a + {\displaystyle \frac {1}{6}}  - n)\,\Gamma (2\,a
 + {\displaystyle \frac {5}{6}}  - n)\,\Gamma ({\displaystyle 
\frac {1}{2}}  + n)\,\Gamma (2\,a)) }
}
\end{maplelatex}

\begin{maplelatex}
\mapleinline{inert}{2d}{337, ":   ",
2/9*sin(1/6*Pi*(-12*a-1+6*n))*3^(1/2)*V7(a-1/3-1/2*n,n+1)/Pi^(1/2)/GAM
MA(3*a-1/2*n)*GAMMA(-a+5/6+1/2*n)*GAMMA(a-1/2*n)*GAMMA(a-1/3-1/2*n)*GA
MMA(4/3+a+1/2*n)/GAMMA(1/2+n)/GAMMA(a-1/2*n+5/6)*GAMMA(2*a-1/3)/GAMMA(
2/3)+1/9*3^(1/2)*(2*sin(1/6*Pi*(-12*a-1+6*n))-1)*V8(a,n)*Pi^(1/2)*GAMM
A(2*a-n)*GAMMA(4/3+a+1/2*n)*GAMMA(a-1/3-1/2*n)*GAMMA(2*a-1/3)/sin(1/6*
Pi*(-6*a-1+3*n))/(2^(2*a-n))/GAMMA(2*a+1/6-n)/GAMMA(2*a+5/6-n)/GAMMA(1
/2+n)/GAMMA(2/3)/GAMMA(2*a);}{%
\maplemultiline{
337\mbox{:~~~} {\displaystyle \frac {2}{9}} \mathrm{sin}
({\displaystyle \frac {\pi \,( - 12\,a - 1 + 6\,n)}{6}} )\,\sqrt{
3}\,\mathrm{V7}(a - {\displaystyle \frac {1}{3}}  - 
{\displaystyle \frac {n}{2}} , \,n + 1)\,\Gamma ( - a + 
{\displaystyle \frac {5}{6}}  + {\displaystyle \frac {n}{2}} )\,
\Gamma (a - {\displaystyle \frac {n}{2}} ) \\
\Gamma (a - {\displaystyle \frac {1}{3}}  - {\displaystyle 
\frac {n}{2}} )\,\Gamma ({\displaystyle \frac {4}{3}}  + a + 
{\displaystyle \frac {n}{2}} )\,\Gamma (2\,a - {\displaystyle 
\frac {1}{3}} ) \left/ {\vrule height0.80em width0em depth0.80em}
 \right. \!  \! (\sqrt{\pi }\,\Gamma (3\,a - {\displaystyle 
\frac {n}{2}} )\,\Gamma ({\displaystyle \frac {1}{2}}  + n)\,
\Gamma (a - {\displaystyle \frac {n}{2}}  + {\displaystyle 
\frac {5}{6}} ) \\
\Gamma ({\displaystyle \frac {2}{3}} ))\mbox{} + {\displaystyle 
\frac {1}{9}} \sqrt{3}\,(2\,\mathrm{sin}({\displaystyle \frac {
\pi \,( - 12\,a - 1 + 6\,n)}{6}} ) - 1)\,\mathrm{V8}(a, \,n)\,
\sqrt{\pi }\,\Gamma (2\,a - n)\,\Gamma ({\displaystyle \frac {4}{
3}}  + a + {\displaystyle \frac {n}{2}} ) \\
\Gamma (a - {\displaystyle \frac {1}{3}}  - {\displaystyle 
\frac {n}{2}} )\,\Gamma (2\,a - {\displaystyle \frac {1}{3}} )
 \left/ {\vrule height0.80em width0em depth0.80em} \right. \! 
 \! (\mathrm{sin}({\displaystyle \frac {\pi \,( - 6\,a - 1 + 3\,n
)}{6}} )\,2^{(2\,a - n)}\,\Gamma (2\,a + {\displaystyle \frac {1
}{6}}  - n) \\
\Gamma (2\,a + {\displaystyle \frac {5}{6}}  - n)\,\Gamma (
{\displaystyle \frac {1}{2}}  + n)\,\Gamma ({\displaystyle 
\frac {2}{3}} )\,\Gamma (2\,a)) }
}
\end{maplelatex}

\begin{maplelatex}
\mapleinline{inert}{2d}{338, ":   ",
2*sin(1/6*Pi*(-12*a-1+6*n))*V7(a-1/3-1/2*n,n+1)/sin(1/3*Pi*(6*a-1))/(2
^(5/3+2*a+n))*GAMMA(5/3+2*a+n)/GAMMA(3*a-1/2*n)/GAMMA(1/2+n)/GAMMA(a-1
/2*n+5/6)*GAMMA(a-1/2*n)*GAMMA(a-1/3-1/2*n)-1/2*(2*sin(1/6*Pi*(-12*a-1
+6*n))-1)*V8(a,n)*Pi^(1/2)*GAMMA(2*a-2/3-n)*GAMMA(2*a-n)*GAMMA(5/3+2*a
+n)/sin(1/3*Pi*(6*a-1))/(2^(2*a-n-1))/(2^(2*a-5/3-n))/(2^(5/3+2*a+n))/
GAMMA(2*a+1/6-n)/GAMMA(2*a+5/6-n)/GAMMA(1/2+n)/GAMMA(2*a);}{%
\maplemultiline{
338\mbox{:~~~} {\displaystyle \frac {2\,\mathrm{sin}(
{\displaystyle \frac {\pi \,( - 12\,a - 1 + 6\,n)}{6}} )\,
\mathrm{V7}(a - {\displaystyle \frac {1}{3}}  - {\displaystyle 
\frac {n}{2}} , \,n + 1)\,\Gamma ({\displaystyle \frac {5}{3}} 
 + 2\,a + n)\,\Gamma (a - {\displaystyle \frac {n}{2}} )\,\Gamma 
(a - {\displaystyle \frac {1}{3}}  - {\displaystyle \frac {n}{2}
} )}{\mathrm{sin}({\displaystyle \frac {\pi \,(6\,a - 1)}{3}} )\,
2^{(5/3 + 2\,a + n)}\,\Gamma (3\,a - {\displaystyle \frac {n}{2}
} )\,\Gamma ({\displaystyle \frac {1}{2}}  + n)\,\Gamma (a - 
{\displaystyle \frac {n}{2}}  + {\displaystyle \frac {5}{6}} )}} 
 \\
\mbox{} - {\displaystyle \frac {1}{2}} (2\,\mathrm{sin}(
{\displaystyle \frac {\pi \,( - 12\,a - 1 + 6\,n)}{6}} ) - 1)\,
\mathrm{V8}(a, \,n)\,\sqrt{\pi }\,\Gamma (2\,a - {\displaystyle 
\frac {2}{3}}  - n)\,\Gamma (2\,a - n) \\
\Gamma ({\displaystyle \frac {5}{3}}  + 2\,a + n) \left/ {\vrule 
height0.80em width0em depth0.80em} \right. \!  \! (\mathrm{sin}(
{\displaystyle \frac {\pi \,(6\,a - 1)}{3}} )\,2^{(2\,a - n - 1)}
\,2^{(2\,a - 5/3 - n)}\,2^{(5/3 + 2\,a + n)} \\
\Gamma (2\,a + {\displaystyle \frac {1}{6}}  - n)\,\Gamma (2\,a
 + {\displaystyle \frac {5}{6}}  - n)\,\Gamma ({\displaystyle 
\frac {1}{2}}  + n)\,\Gamma (2\,a)) }
}
\end{maplelatex}

\begin{maplelatex}
\mapleinline{inert}{2d}{339, ":   ",
GAMMA(2*a-1/3)*GAMMA(a-1/3-1/2*n)*GAMMA(a-1/2*n)*GAMMA(4/3+a+1/2*n)/Pi
^(1/2)*V7(a-1/3-1/2*n,n+1)/GAMMA(2*a-5/6-n)/GAMMA(a-1/2*n+5/6)/GAMMA(1
/2+n)*GAMMA(a+1/3-1/2*n)/GAMMA(3*a-1/2*n)+1/24*V8(a,n)*GAMMA(a+1/3-1/2
*n)*GAMMA(4/3+a+1/2*n)*GAMMA(2*a-2/3-n)*GAMMA(2*a-n)*GAMMA(2*a-1/3)*(2
*sin(1/6*Pi*(-12*a-1+6*n))-1)*(-12*a+5+6*n)/sin(1/6*Pi*(-12*a-1+6*n))/
(2^(2*a-n-1))/(2^(2*a-5/3-n))/GAMMA(2*a+1/6-n)^2/GAMMA(2*a+5/6-n)/GAMM
A(1/2+n)/GAMMA(2*a);}{%
\maplemultiline{
339\mbox{:~~~} 
{\displaystyle \frac {\Gamma (2\,a - {\displaystyle \frac {1}{3}
} )\,\Gamma (a - {\displaystyle \frac {1}{3}}  - {\displaystyle 
\frac {n}{2}} )\,\Gamma (a - {\displaystyle \frac {n}{2}} )\,
\Gamma ({\displaystyle \frac {4}{3}}  + a + {\displaystyle 
\frac {n}{2}} )\,\mathrm{V7}(a - {\displaystyle \frac {1}{3}}  - 
{\displaystyle \frac {n}{2}} , \,n + 1)\,\Gamma (a + 
{\displaystyle \frac {1}{3}}  - {\displaystyle \frac {n}{2}} )}{
\sqrt{\pi }\,\Gamma (2\,a - {\displaystyle \frac {5}{6}}  - n)\,
\Gamma (a - {\displaystyle \frac {n}{2}}  + {\displaystyle 
\frac {5}{6}} )\,\Gamma ({\displaystyle \frac {1}{2}}  + n)\,
\Gamma (3\,a - {\displaystyle \frac {n}{2}} )}}  \\
\mbox{} + {\displaystyle \frac {1}{24}} \mathrm{V8}(a, \,n)\,
\Gamma (a + {\displaystyle \frac {1}{3}}  - {\displaystyle 
\frac {n}{2}} )\,\Gamma ({\displaystyle \frac {4}{3}}  + a + 
{\displaystyle \frac {n}{2}} )\,\Gamma (2\,a - {\displaystyle 
\frac {2}{3}}  - n)\,\Gamma (2\,a - n)\,\Gamma (2\,a - 
{\displaystyle \frac {1}{3}} ) \\
(2\,\mathrm{sin}({\displaystyle \frac {\pi \,( - 12\,a - 1 + 6\,n
)}{6}} ) - 1)\,( - 12\,a + 5 + 6\,n) \left/ {\vrule 
height0.80em width0em depth0.80em} \right. \!  \! (\mathrm{sin}(
{\displaystyle \frac {\pi \,( - 12\,a - 1 + 6\,n)}{6}} ) \\
2^{(2\,a - n - 1)}\,2^{(2\,a - 5/3 - n)}\,\Gamma (2\,a + 
{\displaystyle \frac {1}{6}}  - n)^{2}\,\Gamma (2\,a + 
{\displaystyle \frac {5}{6}}  - n)\,\Gamma ({\displaystyle 
\frac {1}{2}}  + n)\,\Gamma (2\,a)) }
}
\end{maplelatex}

\begin{maplelatex}
\mapleinline{inert}{2d}{340, ":   ",
4/27*sin(1/6*Pi*(-12*a-1+6*n))*3^(1/2)*V7(a-1/3-1/2*n,n+1)/Pi^(1/2)/GA
MMA(3*a-1/2*n)/GAMMA(a-1/2*n+5/6)*GAMMA(a-1/2*n)*GAMMA(a-1/3-1/2*n)*GA
MMA(2*a-1/3)-1/27*3^(1/2)*(2*sin(1/6*Pi*(-12*a-1+6*n))-1)*V8(a,n)*GAMM
A(2*a-2/3-n)*GAMMA(2*a-n)*GAMMA(2*a-1/3)/(2^(2*a-n-1))/(2^(2*a-5/3-n))
/GAMMA(2*a+1/6-n)/GAMMA(2*a+5/6-n)/GAMMA(2*a);}{%
\maplemultiline{
340\mbox{:~~~} {\displaystyle \frac {4}{27}} \,
{\displaystyle \frac {\mathrm{sin}({\displaystyle \frac {\pi \,(
 - 12\,a - 1 + 6\,n)}{6}} )\,\sqrt{3}\,\mathrm{V7}(a - 
{\displaystyle \frac {1}{3}}  - {\displaystyle \frac {n}{2}} , \,
n + 1)\,\Gamma (a - {\displaystyle \frac {n}{2}} )\,\Gamma (a - 
{\displaystyle \frac {1}{3}}  - {\displaystyle \frac {n}{2}} )\,
\Gamma (2\,a - {\displaystyle \frac {1}{3}} )}{\sqrt{\pi }\,
\Gamma (3\,a - {\displaystyle \frac {n}{2}} )\,\Gamma (a - 
{\displaystyle \frac {n}{2}}  + {\displaystyle \frac {5}{6}} )}} 
 \\
\mbox{} - {\displaystyle \frac {1}{27}}
{\displaystyle \frac {\sqrt{3}\,(2\,\mathrm{sin}({\displaystyle 
\frac {\pi \,( - 12\,a - 1 + 6\,n)}{6}} ) - 1)\,\mathrm{V8}(a, \,
n)\,\Gamma (2\,a - {\displaystyle \frac {2}{3}}  - n)\,\Gamma (2
\,a - n)\,\Gamma (2\,a - {\displaystyle \frac {1}{3}} )}{2^{(2\,a
 - n - 1)}\,2^{(2\,a - 5/3 - n)}\,\Gamma (2\,a + {\displaystyle 
\frac {1}{6}}  - n)\,\Gamma (2\,a + {\displaystyle \frac {5}{6}} 
 - n)\,\Gamma (2\,a)}}  }
}
\end{maplelatex}

\begin{maplelatex}
\mapleinline{inert}{2d}{341, ":   ",
-2/3*sin(1/6*Pi*(-12*a-1+6*n))*V7(a-1/3-1/2*n,n+1)/sin(1/2*Pi*(n-2*a))
/Pi^(1/2)/GAMMA(3*a-1/2*n)/GAMMA(1/2+n)/GAMMA(a-1/2*n+5/6)*GAMMA(5/6+a
+1/2*n)*GAMMA(a-1/3-1/2*n)*GAMMA(2*a-1/3)*GAMMA(2/3)+1/6*(2*sin(1/6*Pi
*(-12*a-1+6*n))-1)*V8(a,n)*GAMMA(2*a-2/3-n)*GAMMA(1/2+a-1/2*n)*GAMMA(5
/6+a+1/2*n)*GAMMA(2*a-1/3)*GAMMA(2/3)/Pi^(1/2)/sin(1/2*Pi*(n-2*a))/(2^
(2*a-5/3-n))/GAMMA(2*a+1/6-n)/GAMMA(2*a+5/6-n)/GAMMA(1/2+n)/GAMMA(2*a)
;}{%
\maplemultiline{
341\mbox{:~~~}  - {\displaystyle \frac {2}{3}} \,
{\displaystyle \frac {\mathrm{sin}({\displaystyle \frac {\pi \,(
 - 12\,a - 1 + 6\,n)}{6}} )\,\mathrm{V7}(a - {\displaystyle 
\frac {1}{3}}  - {\displaystyle \frac {n}{2}} , \,n + 1)\,\Gamma 
({\displaystyle \frac {5}{6}}  + a + {\displaystyle \frac {n}{2}
} )\,\Gamma (a - {\displaystyle \frac {1}{3}}  - {\displaystyle 
\frac {n}{2}} )\,\Gamma (2\,a - {\displaystyle \frac {1}{3}} )\,
\Gamma ({\displaystyle \frac {2}{3}} )}{\mathrm{sin}(
{\displaystyle \frac {\pi \,(n - 2\,a)}{2}} )\,\sqrt{\pi }\,
\Gamma (3\,a - {\displaystyle \frac {n}{2}} )\,\Gamma (
{\displaystyle \frac {1}{2}}  + n)\,\Gamma (a - {\displaystyle 
\frac {n}{2}}  + {\displaystyle \frac {5}{6}} )}}  \\
\mbox{} + {\displaystyle \frac {1}{6}} (2\,\mathrm{sin}(
{\displaystyle \frac {\pi \,( - 12\,a - 1 + 6\,n)}{6}} ) - 1)\,
\mathrm{V8}(a, \,n)\,\Gamma (2\,a - {\displaystyle \frac {2}{3}} 
 - n)\,\Gamma ({\displaystyle \frac {1}{2}}  + a - 
{\displaystyle \frac {n}{2}} ) \\
\Gamma ({\displaystyle \frac {5}{6}}  + a + {\displaystyle 
\frac {n}{2}} )\,\Gamma (2\,a - {\displaystyle \frac {1}{3}} )\,
\Gamma ({\displaystyle \frac {2}{3}} ) \left/ {\vrule 
height0.80em width0em depth0.80em} \right. \!  \! (\sqrt{\pi }\,
\mathrm{sin}({\displaystyle \frac {\pi \,(n - 2\,a)}{2}} )\,2^{(2
\,a - 5/3 - n)} \\
\Gamma (2\,a + {\displaystyle \frac {1}{6}}  - n)\,\Gamma (2\,a
 + {\displaystyle \frac {5}{6}}  - n)\,\Gamma ({\displaystyle 
\frac {1}{2}}  + n)\,\Gamma (2\,a)) }
}
\end{maplelatex}

\begin{maplelatex}
\mapleinline{inert}{2d}{342, ":   ",
2/3*GAMMA(2*a-1/3)*GAMMA(a-1/3-1/2*n)*GAMMA(a-1/2*n)/Pi^(1/2)*V7(a-1/3
-1/2*n,n+1)*GAMMA(2*a)/GAMMA(2*a-5/6-n)/GAMMA(a-1/2*n+5/6)/GAMMA(1/2+n
)*GAMMA(2/3)/GAMMA(3*a-1/2*n)+1/36*V8(a,n)*GAMMA(2*a-n)*GAMMA(2*a-2/3-
n)*GAMMA(2*a-1/3)*GAMMA(2/3)*(2*sin(1/6*Pi*(-12*a-1+6*n))-1)*(-12*a+5+
6*n)/sin(1/6*Pi*(-12*a-1+6*n))/(2^(2*a-n-1))/(2^(2*a-5/3-n))/GAMMA(2*a
+1/6-n)^2/GAMMA(2*a+5/6-n)/GAMMA(1/2+n);}{%
\maplemultiline{
342\mbox{:~~~} {\displaystyle \frac {2}{3}} \,
{\displaystyle \frac {\Gamma (2\,a - {\displaystyle \frac {1}{3}
} )\,\Gamma (a - {\displaystyle \frac {1}{3}}  - {\displaystyle 
\frac {n}{2}} )\,\Gamma (a - {\displaystyle \frac {n}{2}} )\,
\mathrm{V7}(a - {\displaystyle \frac {1}{3}}  - {\displaystyle 
\frac {n}{2}} , \,n + 1)\,\Gamma (2\,a)\,\Gamma ({\displaystyle 
\frac {2}{3}} )}{\sqrt{\pi }\,\Gamma (2\,a - {\displaystyle 
\frac {5}{6}}  - n)\,\Gamma (a - {\displaystyle \frac {n}{2}}  + 
{\displaystyle \frac {5}{6}} )\,\Gamma ({\displaystyle \frac {1}{
2}}  + n)\,\Gamma (3\,a - {\displaystyle \frac {n}{2}} )}}  + 
{\displaystyle \frac {1}{36}}  \\
\mathrm{V8}(a, \,n)\,\Gamma (2\,a - n)\,\Gamma (2\,a - 
{\displaystyle \frac {2}{3}}  - n)\,\Gamma (2\,a - 
{\displaystyle \frac {1}{3}} )\,\Gamma ({\displaystyle \frac {2}{
3}} )\,(2\,\mathrm{sin}({\displaystyle \frac {\pi \,( - 12\,a - 1
 + 6\,n)}{6}} ) - 1) \\
( - 12\,a + 5 + 6\,n) \left/ {\vrule 
height0.80em width0em depth0.80em} \right. \!  \! (\mathrm{sin}(
{\displaystyle \frac {\pi \,( - 12\,a - 1 + 6\,n)}{6}} )\,2^{(2\,
a - n - 1)}\,2^{(2\,a - 5/3 - n)} \\
\Gamma (2\,a + {\displaystyle \frac {1}{6}}  - n)^{2}\,\Gamma (2
\,a + {\displaystyle \frac {5}{6}}  - n)\,\Gamma ({\displaystyle 
\frac {1}{2}}  + n)) }
}
\end{maplelatex}

\begin{maplelatex}
\mapleinline{inert}{2d}{343, ":   ",
-2/9*sin(1/6*Pi*(-12*a-1+6*n))*3^(1/2)*V7(a-1/3-1/2*n,n+1)*Pi^(1/2)/si
n(1/6*Pi*(-6*a+2+3*n))/GAMMA(3*a-1/2*n)/GAMMA(1/2+n)/GAMMA(a-1/2*n+5/6
)*GAMMA(5/6+a+1/2*n)*GAMMA(a-1/2*n)*GAMMA(2*a-1/3)/GAMMA(2/3)+1/18*V8(
a,n)*(2*sin(1/6*Pi*(-12*a-1+6*n))-1)*3^(1/2)*Pi^(1/2)*GAMMA(2*a-n)*GAM
MA(a+1/6-1/2*n)*GAMMA(5/6+a+1/2*n)*GAMMA(2*a-1/3)/sin(1/6*Pi*(-6*a+2+3
*n))/(2^(2*a-n-1))/GAMMA(1/2+n)/GAMMA(2*a+1/6-n)/GAMMA(2*a+5/6-n)/GAMM
A(2/3)/GAMMA(2*a);}{%
\maplemultiline{
343\mbox{:~~~}  - {\displaystyle \frac {2}{9}} \,
{\displaystyle \frac {\mathrm{sin}({\displaystyle \frac {\pi \,(
 - 12\,a - 1 + 6\,n)}{6}} )\,\sqrt{3}\,\mathrm{V7}(a - 
{\displaystyle \frac {1}{3}}  - {\displaystyle \frac {n}{2}} , \,
n + 1)\,\sqrt{\pi }\,\Gamma ({\displaystyle \frac {5}{6}}  + a + 
{\displaystyle \frac {n}{2}} )\,\Gamma (a - {\displaystyle 
\frac {n}{2}} )\,\Gamma (2\,a - {\displaystyle \frac {1}{3}} )}{
\mathrm{sin}({\displaystyle \frac {\pi \,( - 6\,a + 2 + 3\,n)}{6}
} )\,\Gamma (3\,a - {\displaystyle \frac {n}{2}} )\,\Gamma (
{\displaystyle \frac {1}{2}}  + n)\,\Gamma (a - {\displaystyle 
\frac {n}{2}}  + {\displaystyle \frac {5}{6}} )\,\Gamma (
{\displaystyle \frac {2}{3}} )}}  \\
\mbox{} + {\displaystyle \frac {1}{18}} \mathrm{V8}(a, \,n)\,(2\,
\mathrm{sin}({\displaystyle \frac {\pi \,( - 12\,a - 1 + 6\,n)}{6
}} ) - 1)\,\sqrt{3}\,\sqrt{\pi }\,\Gamma (2\,a - n)\,\Gamma (a + 
{\displaystyle \frac {1}{6}}  - {\displaystyle \frac {n}{2}} )
 \\
\Gamma ({\displaystyle \frac {5}{6}}  + a + {\displaystyle 
\frac {n}{2}} )\,\Gamma (2\,a - {\displaystyle \frac {1}{3}} )
 \left/ {\vrule height0.80em width0em depth0.80em} \right. \! 
 \! (\mathrm{sin}({\displaystyle \frac {\pi \,( - 6\,a + 2 + 3\,n
)}{6}} )\,2^{(2\,a - n - 1)}\,\Gamma ({\displaystyle \frac {1}{2}
}  + n) \\
\Gamma (2\,a + {\displaystyle \frac {1}{6}}  - n)\,\Gamma (2\,a
 + {\displaystyle \frac {5}{6}}  - n)\,\Gamma ({\displaystyle 
\frac {2}{3}} )\,\Gamma (2\,a)) }
}
\end{maplelatex}

\begin{maplelatex}
\mapleinline{inert}{2d}{344, ":   ",
2/9*GAMMA(2*a-1/3)*GAMMA(a-1/3-1/2*n)*GAMMA(a-1/2*n)*Pi^(1/2)/GAMMA(2*
a-5/6-n)/GAMMA(a-1/2*n+5/6)/GAMMA(1/2+n)*3^(1/2)/GAMMA(2/3)/GAMMA(3*a-
1/2*n)*GAMMA(2*a+1/3)*V7(a-1/3-1/2*n,n+1)+1/108*V8(a,n)*3^(1/2)*Pi*GAM
MA(2*a-2/3-n)*GAMMA(2*a-n)*GAMMA(2*a+1/3)*GAMMA(2*a-1/3)*(2*sin(1/6*Pi
*(-12*a-1+6*n))-1)*(-12*a+5+6*n)/sin(1/6*Pi*(-12*a-1+6*n))/(2^(2*a-n-1
))/(2^(2*a-5/3-n))/GAMMA(2*a+1/6-n)^2/GAMMA(2*a+5/6-n)/GAMMA(1/2+n)/GA
MMA(2/3)/GAMMA(2*a);}{%
\maplemultiline{
344\mbox{:~~~} {\displaystyle \frac {2}{9}} \,
{\displaystyle \frac {\Gamma (2\,a - {\displaystyle \frac {1}{3}
} )\,\Gamma (a - {\displaystyle \frac {1}{3}}  - {\displaystyle 
\frac {n}{2}} )\,\Gamma (a - {\displaystyle \frac {n}{2}} )\,
\sqrt{\pi }\,\sqrt{3}\,\Gamma (2\,a + {\displaystyle \frac {1}{3}
} )\,\mathrm{V7}(a - {\displaystyle \frac {1}{3}}  - 
{\displaystyle \frac {n}{2}} , \,n + 1)}{\Gamma (2\,a - 
{\displaystyle \frac {5}{6}}  - n)\,\Gamma (a - {\displaystyle 
\frac {n}{2}}  + {\displaystyle \frac {5}{6}} )\,\Gamma (
{\displaystyle \frac {1}{2}}  + n)\,\Gamma ({\displaystyle 
\frac {2}{3}} )\,\Gamma (3\,a - {\displaystyle \frac {n}{2}} )}} 
 +  \\
{\displaystyle \frac {1}{108}} \mathrm{V8}(a, \,n)\,\sqrt{3}\,\pi
 \,\Gamma (2\,a - {\displaystyle \frac {2}{3}}  - n)\,\Gamma (2\,
a - n)\,\Gamma (2\,a + {\displaystyle \frac {1}{3}} )\,\Gamma (2
\,a - {\displaystyle \frac {1}{3}} ) \\
(2\,\mathrm{sin}({\displaystyle \frac {\pi \,( - 12\,a - 1 + 6\,n
)}{6}} ) - 1)\,( - 12\,a + 5 + 6\,n) \left/ {\vrule 
height0.80em width0em depth0.80em} \right. \!  \! (\mathrm{sin}(
{\displaystyle \frac {\pi \,( - 12\,a - 1 + 6\,n)}{6}} ) \\
2^{(2\,a - n - 1)}\,2^{(2\,a - 5/3 - n)}\,\Gamma (2\,a + 
{\displaystyle \frac {1}{6}}  - n)^{2}\,\Gamma (2\,a + 
{\displaystyle \frac {5}{6}}  - n)\,\Gamma ({\displaystyle 
\frac {1}{2}}  + n)\,\Gamma ({\displaystyle \frac {2}{3}} )\,
\Gamma (2\,a)) }
}
\end{maplelatex}

\begin{maplelatex}
\mapleinline{inert}{2d}{345, ":   ",
1/3*(3^(1/2)*sin(1/6*Pi*(-18*a+2+3*n))+3^(1/2)*sin(1/2*Pi*(2*a+n))-3*c
os(1/2*Pi*(-6*a+n)))/Pi^(3/2)*V7(a-1/3-1/2*n,n+1)/sin(1/3*Pi*(6*a-1))/
((-1)^(-n))/GAMMA(3*a-1/2*n)*GAMMA(-a+5/6+1/2*n)*GAMMA(a+1/3-1/2*n)/GA
MMA(1/2+a-1/2*n)*GAMMA(7/6-a-1/2*n)*GAMMA(a-1/2*n)*GAMMA(2*a+1/3)/GAMM
A(2/3)*GAMMA(2*a)-1/3*V8(a,n)*(3^(1/2)*sin(Pi*(n-2*a))+3^(1/2)*sin(1/3
*Pi*(12*a+1))+3^(1/2)*sin(1/3*Pi*(6*a-1+3*n))-3*sin(1/6*Pi*(-12*a+6*n+
1)))*GAMMA(2*a+2/3-n)*GAMMA(7/6-a-1/2*n)*GAMMA(a-1/2*n)*GAMMA(2*a+1/3)
/Pi^(1/2)/sin(1/3*Pi*(6*a-1))/(2^(2*a+2/3-n))/GAMMA(2*a+1/6-n)/GAMMA(2
*a+5/6-n)/GAMMA(2/3);}{%
\maplemultiline{
345\mbox{:~~~} {\displaystyle \frac {1}{3}} (\sqrt{3}\,
\mathrm{sin}({\displaystyle \frac {\pi \,( - 18\,a + 2 + 3\,n)}{6
}} ) + \sqrt{3}\,\mathrm{sin}({\displaystyle \frac {\pi \,(2\,a
 + n)}{2}} ) - 3\,\mathrm{cos}({\displaystyle \frac {\pi \,( - 6
\,a + n)}{2}} )) \\
\mathrm{V7}(a - {\displaystyle \frac {1}{3}}  - {\displaystyle 
\frac {n}{2}} , \,n + 1)\,\Gamma ( - a + {\displaystyle \frac {5
}{6}}  + {\displaystyle \frac {n}{2}} )\,\Gamma (a + 
{\displaystyle \frac {1}{3}}  - {\displaystyle \frac {n}{2}} )\,
\Gamma ({\displaystyle \frac {7}{6}}  - a - {\displaystyle 
\frac {n}{2}} )\,\Gamma (a - {\displaystyle \frac {n}{2}} ) \\
\Gamma (2\,a + {\displaystyle \frac {1}{3}} )\,\Gamma (2\,a)
 \left/ {\vrule height0.80em width0em depth0.80em} \right. \! 
 \! (\pi ^{(3/2)}\,\mathrm{sin}({\displaystyle \frac {\pi \,(6\,a
 - 1)}{3}} )\,(-1)^{( - n)}\,\Gamma (3\,a - {\displaystyle 
\frac {n}{2}} )\,\Gamma ({\displaystyle \frac {1}{2}}  + a - 
{\displaystyle \frac {n}{2}} ) \\
\Gamma ({\displaystyle \frac {2}{3}} ))\mbox{} - {\displaystyle 
\frac {1}{3}} \mathrm{V8}(a, \,n)(\sqrt{3}\,\mathrm{sin}(\pi \,(n
 - 2\,a)) + \sqrt{3}\,\mathrm{sin}({\displaystyle \frac {\pi \,(
12\,a + 1)}{3}} ) \\
\mbox{} + \sqrt{3}\,\mathrm{sin}({\displaystyle \frac {\pi \,(6\,
a - 1 + 3\,n)}{3}} ) - 3\,\mathrm{sin}({\displaystyle \frac {\pi 
\,( - 12\,a + 6\,n + 1)}{6}} ))\Gamma (2\,a + {\displaystyle 
\frac {2}{3}}  - n) \\
\Gamma ({\displaystyle \frac {7}{6}}  - a - {\displaystyle 
\frac {n}{2}} )\,\Gamma (a - {\displaystyle \frac {n}{2}} )\,
\Gamma (2\,a + {\displaystyle \frac {1}{3}} ) \left/ {\vrule 
height0.80em width0em depth0.80em} \right. \!  \! (\sqrt{\pi }\,
\mathrm{sin}({\displaystyle \frac {\pi \,(6\,a - 1)}{3}} )\,2^{(2
\,a + 2/3 - n)} \\
\Gamma (2\,a + {\displaystyle \frac {1}{6}}  - n)\,\Gamma (2\,a
 + {\displaystyle \frac {5}{6}}  - n)\,\Gamma ({\displaystyle 
\frac {2}{3}} )) }
}
\end{maplelatex}

\begin{maplelatex}
\mapleinline{inert}{2d}{346, ":   ",
V7(a-1/3-1/2*n,n+1)*(sin(1/2*Pi*(2*a+n))-sin(1/2*Pi*(-2*a+3*n)))*GAMMA
(a+1/3-1/2*n)*GAMMA(7/6-a-1/2*n)*GAMMA(2*a+1/3)*GAMMA(2*a-1/3)*GAMMA(2
*a)/Pi^(1/2)/GAMMA(1/3+a+1/2*n)/GAMMA(3*a-1/2*n)/GAMMA(a+1/6-1/2*n)/GA
MMA(1/2+a-1/2*n)/(sin(1/6*Pi*(12*a-1))-sin(1/6*Pi*(-1+6*n)))-V8(a,n)*(
sin(1/6*Pi*(-18*a+2+3*n))+sin(1/2*Pi*(2*a+n))-sin(1/2*Pi*(-2*a+3*n))+s
in(1/6*Pi*(6*a+3*n-2)))*GAMMA(2*a+2/3-n)*GAMMA(7/6-a-1/2*n)*GAMMA(2*a+
1/3)*GAMMA(2*a-1/3)/Pi^(1/2)/(2^(2*a+2/3-n))/GAMMA(1/3+a+1/2*n)/GAMMA(
2*a+1/6-n)/GAMMA(2*a+5/6-n)/(sin(1/6*Pi*(12*a-1))-sin(1/6*Pi*(-1+6*n))
);}{%
\maplemultiline{
346\mbox{:~~~} \mathrm{V7}(a - {\displaystyle \frac {1}{
3}}  - {\displaystyle \frac {n}{2}} , \,n + 1)\,(\mathrm{sin}(
{\displaystyle \frac {\pi \,(2\,a + n)}{2}} ) - \mathrm{sin}(
{\displaystyle \frac {\pi \,( - 2\,a + 3\,n)}{2}} ))\,\Gamma (a
 + {\displaystyle \frac {1}{3}}  - {\displaystyle \frac {n}{2}} )
 \\
\Gamma ({\displaystyle \frac {7}{6}}  - a - {\displaystyle 
\frac {n}{2}} )\,\Gamma (2\,a + {\displaystyle \frac {1}{3}} )\,
\Gamma (2\,a - {\displaystyle \frac {1}{3}} )\,\Gamma (2\,a)
 \left/ {\vrule height0.80em width0em depth0.80em} \right. \! 
 \! (\sqrt{\pi }\,\Gamma ({\displaystyle \frac {1}{3}}  + a + 
{\displaystyle \frac {n}{2}} )\,\Gamma (3\,a - {\displaystyle 
\frac {n}{2}} ) \\
\Gamma (a + {\displaystyle \frac {1}{6}}  - {\displaystyle 
\frac {n}{2}} )\,\Gamma ({\displaystyle \frac {1}{2}}  + a - 
{\displaystyle \frac {n}{2}} )\,(\mathrm{sin}({\displaystyle 
\frac {\pi \,(12\,a - 1)}{6}} ) - \mathrm{sin}({\displaystyle 
\frac {\pi \,( - 1 + 6\,n)}{6}} )))\mbox{} - \mathrm{V8}(a, \,n)(
 \\
\mathrm{sin}({\displaystyle \frac {\pi \,( - 18\,a + 2 + 3\,n)}{6
}} ) + \mathrm{sin}({\displaystyle \frac {\pi \,(2\,a + n)}{2}} )
 - \mathrm{sin}({\displaystyle \frac {\pi \,( - 2\,a + 3\,n)}{2}
} ) \\
\mbox{} + \mathrm{sin}({\displaystyle \frac {\pi \,(6\,a + 3\,n
 - 2)}{6}} ))\Gamma (2\,a + {\displaystyle \frac {2}{3}}  - n)\,
\Gamma ({\displaystyle \frac {7}{6}}  - a - {\displaystyle 
\frac {n}{2}} )\,\Gamma (2\,a + {\displaystyle \frac {1}{3}} )\,
\Gamma (2\,a - {\displaystyle \frac {1}{3}} ) \left/ {\vrule 
height0.80em width0em depth0.80em} \right. \!  \!  \\
(\sqrt{\pi }\,2^{(2\,a + 2/3 - n)}\,\Gamma ({\displaystyle 
\frac {1}{3}}  + a + {\displaystyle \frac {n}{2}} )\,\Gamma (2\,a
 + {\displaystyle \frac {1}{6}}  - n)\,\Gamma (2\,a + 
{\displaystyle \frac {5}{6}}  - n) \\
(\mathrm{sin}({\displaystyle \frac {\pi \,(12\,a - 1)}{6}} ) - 
\mathrm{sin}({\displaystyle \frac {\pi \,( - 1 + 6\,n)}{6}} )))
 }
}
\end{maplelatex}

\begin{maplelatex}
\mapleinline{inert}{2d}{347, ":   ",
1/8*V7(a-1/3-1/2*n,n+1)*GAMMA(2*a+1/3)*GAMMA(2*a-1/3)*GAMMA(3/2-2*a+n)
*GAMMA(a-1/3-1/2*n)*GAMMA(a+1/3-1/2*n)*GAMMA(a-1/2*n)*GAMMA(1/2-a+1/2*
n)*GAMMA(2*a)*(2*sin(Pi*(-6*a+n))+(-1)^(-n)*sin(6*Pi*a))*(-6*a+n+2)/si
n(1/2*Pi*(-6*a+n))/((-1)^(-n))/Pi^(5/2)/GAMMA(3*a-1/2*n)^2/GAMMA(1/2+n
)+1/4*(1+2*sin(1/6*Pi*(-12*a+6*n+1)))*V8(a,n)*GAMMA(2*a+2/3-n)*GAMMA(3
/2-2*a+n)*GAMMA(a-1/3-1/2*n)*GAMMA(a-1/2*n)*GAMMA(2*a+1/3)*GAMMA(2*a-1
/3)/cos(1/2*Pi*(n-2*a))/Pi^(1/2)/(2^(-1/3+2*a-n))/GAMMA(2*a+5/6-n)/GAM
MA(2*a+1/6-n)/GAMMA(-a+5/6+1/2*n)/GAMMA(1/2+n)/GAMMA(3*a-1/2*n-1);}{%
\maplemultiline{
347\mbox{:~~~} {\displaystyle \frac {1}{8}} \mathrm{V7}(
a - {\displaystyle \frac {1}{3}}  - {\displaystyle \frac {n}{2}} 
, \,n + 1)\,\Gamma (2\,a + {\displaystyle \frac {1}{3}} )\,\Gamma
 (2\,a - {\displaystyle \frac {1}{3}} )\,\Gamma ({\displaystyle 
\frac {3}{2}}  - 2\,a + n)\,\Gamma (a - {\displaystyle \frac {1}{
3}}  - {\displaystyle \frac {n}{2}} ) \\
\Gamma (a + {\displaystyle \frac {1}{3}}  - {\displaystyle 
\frac {n}{2}} )\,\Gamma (a - {\displaystyle \frac {n}{2}} )\,
\Gamma ({\displaystyle \frac {1}{2}}  - a + {\displaystyle 
\frac {n}{2}} )\,\Gamma (2\,a) \\
(2\,\mathrm{sin}(\pi \,( - 6\,a + n)) + (-1)^{( - n)}\,\mathrm{
sin}(6\,\pi \,a))\,( - 6\,a + n + 2) \left/ {\vrule 
height0.80em width0em depth0.80em} \right. \!  \! (\mathrm{sin}(
{\displaystyle \frac {\pi \,( - 6\,a + n)}{2}} ) \\
(-1)^{( - n)}\,\pi ^{(5/2)}\,\Gamma (3\,a - {\displaystyle 
\frac {n}{2}} )^{2}\,\Gamma ({\displaystyle \frac {1}{2}}  + n))
\mbox{} + {\displaystyle \frac {1}{4}} (1 + 2\,\mathrm{sin}(
{\displaystyle \frac {\pi \,( - 12\,a + 6\,n + 1)}{6}} ))\,
\mathrm{V8}(a, \,n) \\
\Gamma (2\,a + {\displaystyle \frac {2}{3}}  - n)\,\Gamma (
{\displaystyle \frac {3}{2}}  - 2\,a + n)\,\Gamma (a - 
{\displaystyle \frac {1}{3}}  - {\displaystyle \frac {n}{2}} )\,
\Gamma (a - {\displaystyle \frac {n}{2}} )\,\Gamma (2\,a + 
{\displaystyle \frac {1}{3}} )\,\Gamma (2\,a - {\displaystyle 
\frac {1}{3}} ) \left/ {\vrule height0.80em width0em depth0.80em}
 \right. \!  \! ( \\
\mathrm{cos}({\displaystyle \frac {\pi \,(n - 2\,a)}{2}} )\,
\sqrt{\pi }\,2^{( - 1/3 + 2\,a - n)}\,\Gamma (2\,a + 
{\displaystyle \frac {5}{6}}  - n)\,\Gamma (2\,a + 
{\displaystyle \frac {1}{6}}  - n)\,\Gamma ( - a + 
{\displaystyle \frac {5}{6}}  + {\displaystyle \frac {n}{2}} )
 \\
\Gamma ({\displaystyle \frac {1}{2}}  + n)\,\Gamma (3\,a - 
{\displaystyle \frac {n}{2}}  - 1)) }
}
\end{maplelatex}

\begin{maplelatex}
\mapleinline{inert}{2d}{348, ":   ",
-1/4*(2*sin(1/6*Pi*(-12*a-1+6*n))-1)*V7(a-1/3-1/2*n,n+1)*GAMMA(-a+5/6+
1/2*n)*GAMMA(1/2-a+1/2*n)*GAMMA(a+1/3-1/2*n)*GAMMA(a-1/3-1/2*n)*GAMMA(
2*a+1/3)*GAMMA(2/3)*GAMMA(2*a)/sin(1/6*Pi*(6*a-1+3*n))/((-1)^(-n))/Pi^
(5/2)/GAMMA(a-1/6+1/2*n)/GAMMA(3*a-1/2*n)+1/4*(1+2*sin(1/6*Pi*(-12*a+6
*n+1)))*V8(a,n)*GAMMA(2*a+2/3-n)*GAMMA(a-1/3-1/2*n)*GAMMA(2*a+1/3)*GAM
MA(2/3)/cos(1/2*Pi*(n-2*a))/Pi^(1/2)/(2^(-1/3+2*a-n))/GAMMA(2*a+5/6-n)
/GAMMA(2*a+1/6-n)/GAMMA(a-1/6+1/2*n);}{%
\maplemultiline{
348\mbox{:~~~}  - {\displaystyle \frac {1}{4}} (2\,
\mathrm{sin}({\displaystyle \frac {\pi \,( - 12\,a - 1 + 6\,n)}{6
}} ) - 1)\,\mathrm{V7}(a - {\displaystyle \frac {1}{3}}  - 
{\displaystyle \frac {n}{2}} , \,n + 1)\,\Gamma ( - a + 
{\displaystyle \frac {5}{6}}  + {\displaystyle \frac {n}{2}} )
 \\
\Gamma ({\displaystyle \frac {1}{2}}  - a + {\displaystyle 
\frac {n}{2}} )\,\Gamma (a + {\displaystyle \frac {1}{3}}  - 
{\displaystyle \frac {n}{2}} )\,\Gamma (a - {\displaystyle 
\frac {1}{3}}  - {\displaystyle \frac {n}{2}} )\,\Gamma (2\,a + 
{\displaystyle \frac {1}{3}} )\,\Gamma ({\displaystyle \frac {2}{
3}} )\,\Gamma (2\,a) \left/ {\vrule 
height0.80em width0em depth0.80em} \right. \!  \! \\ (
\mathrm{sin}({\displaystyle \frac {\pi \,(6\,a - 1 + 3\,n)}{6}} )
\,(-1)^{( - n)}\,\pi ^{(5/2)}\,\Gamma (a - {\displaystyle \frac {
1}{6}}  + {\displaystyle \frac {n}{2}} )\,\Gamma (3\,a - 
{\displaystyle \frac {n}{2}} ))\mbox{} \\ + {\displaystyle \frac {1
}{4}} \,{\displaystyle \frac {(1 + 2\,\mathrm{sin}(
{\displaystyle \frac {\pi \,( - 12\,a + 6\,n + 1)}{6}} ))\,
\mathrm{V8}(a, \,n)\,\Gamma (2\,a + {\displaystyle \frac {2}{3}} 
 - n)\,\Gamma (a - {\displaystyle \frac {1}{3}}  - 
{\displaystyle \frac {n}{2}} )\,\Gamma (2\,a + {\displaystyle 
\frac {1}{3}} )\,\Gamma ({\displaystyle \frac {2}{3}} )}{\mathrm{
cos}({\displaystyle \frac {\pi \,(n - 2\,a)}{2}} )\,\sqrt{\pi }\,
2^{( - 1/3 + 2\,a - n)}\,\Gamma (2\,a + {\displaystyle \frac {5}{
6}}  - n)\,\Gamma (2\,a + {\displaystyle \frac {1}{6}}  - n)\,
\Gamma (a - {\displaystyle \frac {1}{6}}  + {\displaystyle 
\frac {n}{2}} )}}  }
}
\end{maplelatex}

\begin{maplelatex}
\mapleinline{inert}{2d}{349, ":   ",
-1/4*(2*sin(1/6*Pi*(-12*a-1+6*n))-1)*V7(a-1/3-1/2*n,n+1)*GAMMA(2*a+1/6
-n)*GAMMA(-a+5/6+1/2*n)*GAMMA(1/2-a+1/2*n)*GAMMA(a+1/3-1/2*n)*GAMMA(1/
6-a+1/2*n)*GAMMA(2*a+1/3)*GAMMA(2*a)/sin(1/6*Pi*(6*a-1+3*n))/((-1)^(-n
))/Pi^(5/2)/GAMMA(a-1/6+1/2*n)/GAMMA(3*a-1/2*n)+V8(a,n)*GAMMA(a+1/3-1/
2*n)*GAMMA(2*a+1/3)/GAMMA(2*a+5/6-n)/GAMMA(a-1/6+1/2*n);}{%
\maplemultiline{
349\mbox{:~~~}  - {\displaystyle \frac {1}{4}} (2\,
\mathrm{sin}({\displaystyle \frac {\pi \,( - 12\,a - 1 + 6\,n)}{6
}} ) - 1)\,\mathrm{V7}(a - {\displaystyle \frac {1}{3}}  - 
{\displaystyle \frac {n}{2}} , \,n + 1)\,\Gamma (2\,a + 
{\displaystyle \frac {1}{6}}  - n) \\
\Gamma ( - a + {\displaystyle \frac {5}{6}}  + {\displaystyle 
\frac {n}{2}} )\,\Gamma ({\displaystyle \frac {1}{2}}  - a + 
{\displaystyle \frac {n}{2}} )\,\Gamma (a + {\displaystyle 
\frac {1}{3}}  - {\displaystyle \frac {n}{2}} )\,\Gamma (
{\displaystyle \frac {1}{6}}  - a + {\displaystyle \frac {n}{2}} 
)\,\Gamma (2\,a + {\displaystyle \frac {1}{3}} )\,\Gamma (2\,a)
 \left/ {\vrule height0.80em width0em depth0.80em} \right. \! 
 \! \\ (
\mathrm{sin}({\displaystyle \frac {\pi \,(6\,a - 1 + 3\,n)}{6}} )
\,(-1)^{( - n)}\,\pi ^{(5/2)}\,\Gamma (a - {\displaystyle \frac {
1}{6}}  + {\displaystyle \frac {n}{2}} )\,\Gamma (3\,a - 
{\displaystyle \frac {n}{2}} )) 
\mbox{} + {\displaystyle \frac {\mathrm{V8}(a, \,n)\,\Gamma (a + 
{\displaystyle \frac {1}{3}}  - {\displaystyle \frac {n}{2}} )\,
\Gamma (2\,a + {\displaystyle \frac {1}{3}} )}{\Gamma (2\,a + 
{\displaystyle \frac {5}{6}}  - n)\,\Gamma (a - {\displaystyle 
\frac {1}{6}}  + {\displaystyle \frac {n}{2}} )}}  }
}
\end{maplelatex}

\mapleinline{inert}{2d}{350, ":   ",
1/2*(-sin(1/3*Pi*(6*a+1))-sin(2*Pi*(-3*a+n))-sin(1/3*Pi*(6*a-1))+sin(2
*Pi*a))*V7(a-1/3-1/2*n,n+1)*GAMMA(a-1/3-1/2*n)*GAMMA(3/2-2*a+n)*GAMMA(
a+1/3-1/2*n)*GAMMA(a-1/2*n)*GAMMA(2*a+1/3)*GAMMA(2*a)/Pi^(3/2)/GAMMA(2
*a-n-1/6)/GAMMA(3*a-1/2*n)/GAMMA(1/2+n)/(sin(1/3*Pi*(1+3*n))+sin(Pi*(-
4*a+n)))-1/12*V8(a,n)*(-1)^(-n)*GAMMA(2*a-n)*GAMMA(2*a+2/3-n)*GAMMA(a-
1/3-1/2*n)*GAMMA(3/2-2*a+n)*GAMMA(2*a+1/3)*(sin(1/6*Pi*(15*n-30*a-2))+
sin(1/6*Pi*(-6*a+2+3*n))+sin(1/6*Pi*(-18*a-2+9*n))-sin(1/6*Pi*(18*a+2+
3*n))+sin(1/6*Pi*(2-18*a+9*n))-sin(1/6*Pi*(-2+6*a+9*n))-sin(1/6*Pi*(-2
+3*n+18*a))+sin(1/6*Pi*(-30*a-2+3*n)))*(-12*a+6*n+1)/(2^(2*a+2/3-n))/(
2^(2*a-n-1))/GAMMA(2*a+5/6-n)^2/GAMMA(1/2+n)/GAMMA(2*a+1/6-n)/GAMMA(-a
+5/6+1/2*n)/(sin(1/3*Pi*(1+3*n))+sin(Pi*(-4*a+n)));}{%
\maplemultiline{
350\mbox{:~~~} {\displaystyle \frac {1}{2}} ( - \mathrm{
sin}({\displaystyle \frac {\pi \,(6\,a + 1)}{3}} ) - \mathrm{sin}
(2\,\pi \,( - 3\,a + n)) - \mathrm{sin}({\displaystyle \frac {\pi
 \,(6\,a - 1)}{3}} ) + \mathrm{sin}(2\,\pi \,a)) \\
\mathrm{V7}(a - {\displaystyle \frac {1}{3}}  - {\displaystyle 
\frac {n}{2}} , \,n + 1)\,\Gamma (a - {\displaystyle \frac {1}{3}
}  - {\displaystyle \frac {n}{2}} )\,\Gamma ({\displaystyle 
\frac {3}{2}}  - 2\,a + n)\,\Gamma (a + {\displaystyle \frac {1}{
3}}  - {\displaystyle \frac {n}{2}} )\,\Gamma (a - 
{\displaystyle \frac {n}{2}} ) \\
\Gamma (2\,a + {\displaystyle \frac {1}{3}} )\,\Gamma (2\,a)
 \left/ {\vrule height0.80em width0em depth0.80em} \right. \! 
 \! (\pi ^{(3/2)}\,\Gamma (2\,a - n - {\displaystyle \frac {1}{6}
} )\,\Gamma (3\,a - {\displaystyle \frac {n}{2}} )\,\Gamma (
{\displaystyle \frac {1}{2}}  + n) \\
(\mathrm{sin}({\displaystyle \frac {\pi \,(1 + 3\,n)}{3}} ) + 
\mathrm{sin}(\pi \,( - 4\,a + n))))\mbox{} - {\displaystyle 
\frac {1}{12}} \mathrm{V8}(a, \,n)\,(-1)^{( - n)}\,\Gamma (2\,a
 - n) \\
\Gamma (2\,a + {\displaystyle \frac {2}{3}}  - n)\,\Gamma (a - 
{\displaystyle \frac {1}{3}}  - {\displaystyle \frac {n}{2}} )\,
\Gamma ({\displaystyle \frac {3}{2}}  - 2\,a + n)\,\Gamma (2\,a
 + {\displaystyle \frac {1}{3}} )(\mathrm{sin}({\displaystyle 
\frac {\pi \,(15\,n - 30\,a - 2)}{6}} ) \\
\mbox{} + \mathrm{sin}({\displaystyle \frac {\pi \,( - 6\,a + 2
 + 3\,n)}{6}} ) + \mathrm{sin}({\displaystyle \frac {\pi \,( - 18
\,a - 2 + 9\,n)}{6}} ) - \mathrm{sin}({\displaystyle \frac {\pi 
\,(18\,a + 2 + 3\,n)}{6}} ) \\
\mbox{} + \mathrm{sin}({\displaystyle \frac {\pi \,(2 - 18\,a + 9
\,n)}{6}} ) - \mathrm{sin}({\displaystyle \frac {\pi \,( - 2 + 6
\,a + 9\,n)}{6}} ) - \mathrm{sin}({\displaystyle \frac {\pi \,(
 - 2 + 3\,n + 18\,a)}{6}} ) \\
\mbox{} + \mathrm{sin}({\displaystyle \frac {\pi \,( - 30\,a - 2
 + 3\,n)}{6}} ))( - 12\,a + 6\,n + 1) \left/ {\vrule 
height0.80em width0em depth0.80em} \right. \!  \! (2^{(2\,a + 2/3
 - n)}\,2^{(2\,a - n - 1)} \\
\Gamma (2\,a + {\displaystyle \frac {5}{6}}  - n)^{2}\,\Gamma (
{\displaystyle \frac {1}{2}}  + n)\,\Gamma (2\,a + 
{\displaystyle \frac {1}{6}}  - n)\,\Gamma ( - a + 
{\displaystyle \frac {5}{6}}  + {\displaystyle \frac {n}{2}} )
(\mathrm{sin}({\displaystyle \frac {\pi \,(1 + 3\,n)}{3}} ) + 
\mathrm{sin}(\pi \,( - 4\,a + n)))) }
}

\begin{mapleinput}
\end{mapleinput}

\end{maplegroup}

%% file: AppendixB351to380.tex
\begin{maplegroup}
\mapleinline{inert}{2d}{351, ":   ",
4/3*V7(a-1/3-1/2*n,n+1)*(-3^(1/2)*sin(1/6*Pi*(-18*a+3*n+2))-3^(1/2)*si
n(1/2*Pi*(2*a+n))+3*cos(1/2*Pi*(-6*a+n)))*GAMMA(2*a+1/6-n)*GAMMA(-a+5/
6+1/2*n)*GAMMA(a+1/3-1/2*n)*GAMMA(2*a+1/3)/Pi^(1/2)/((-1)^(-n))/GAMMA(
3*a-1/2*n)/GAMMA(1/2+a-1/2*n)/GAMMA(2/3)/(2*sin(1/6*Pi*(24*a+1))-1)+4/
3*V8(a,n)*(3^(1/2)*sin(Pi*(n-2*a))+3^(1/2)*sin(1/3*Pi*(12*a+1))+3^(1/2
)*sin(1/3*Pi*(6*a-1+3*n))-3*sin(1/6*Pi*(-12*a+6*n+1)))*Pi^(1/2)*GAMMA(
2*a+2/3-n)*GAMMA(2*a+1/3)/(2^(2*a+2/3-n))/GAMMA(2*a+5/6-n)/GAMMA(2/3)/
GAMMA(2*a)/(2*sin(1/6*Pi*(24*a+1))-1);}{%
\maplemultiline{
351\mbox{:~~~} {\displaystyle \frac {4}{3}} \mathrm{V7}(
a - {\displaystyle \frac {1}{3}}  - {\displaystyle \frac {n}{2}} 
, \,n + 1) \\
( - \sqrt{3}\,\mathrm{sin}({\displaystyle \frac {\pi \,( - 18\,a
 + 3\,n + 2)}{6}} ) - \sqrt{3}\,\mathrm{sin}({\displaystyle 
\frac {\pi \,(2\,a + n)}{2}} ) + 3\,\mathrm{cos}({\displaystyle 
\frac {\pi \,( - 6\,a + n)}{2}} )) \\
\Gamma (2\,a + {\displaystyle \frac {1}{6}}  - n)\,\Gamma ( - a
 + {\displaystyle \frac {5}{6}}  + {\displaystyle \frac {n}{2}} )
\,\Gamma (a + {\displaystyle \frac {1}{3}}  - {\displaystyle 
\frac {n}{2}} )\,\Gamma (2\,a + {\displaystyle \frac {1}{3}} )
 \left/ {\vrule height0.80em width0em depth0.80em} \right. \! 
 \! (\sqrt{\pi }\,(-1)^{( - n)} \\
\Gamma (3\,a - {\displaystyle \frac {n}{2}} )\,\Gamma (
{\displaystyle \frac {1}{2}}  + a - {\displaystyle \frac {n}{2}} 
)\,\Gamma ({\displaystyle \frac {2}{3}} )\,(2\,\mathrm{sin}(
{\displaystyle \frac {\pi \,(24\,a + 1)}{6}} ) - 1))\mbox{} + 
{\displaystyle \frac {4}{3}} \mathrm{V8}(a, \,n)( \\
\sqrt{3}\,\mathrm{sin}(\pi \,(n - 2\,a)) + \sqrt{3}\,\mathrm{sin}
({\displaystyle \frac {\pi \,(12\,a + 1)}{3}} ) + \sqrt{3}\,
\mathrm{sin}({\displaystyle \frac {\pi \,(6\,a - 1 + 3\,n)}{3}} )
 \\
\mbox{} - 3\,\mathrm{sin}({\displaystyle \frac {\pi \,( - 12\,a
 + 6\,n + 1)}{6}} ))\sqrt{\pi }\,\Gamma (2\,a + {\displaystyle 
\frac {2}{3}}  - n)\,\Gamma (2\,a + {\displaystyle \frac {1}{3}} 
) \left/ {\vrule height0.80em width0em depth0.80em} \right. \! 
 \! (2^{(2\,a + 2/3 - n)} \\
\Gamma (2\,a + {\displaystyle \frac {5}{6}}  - n)\,\Gamma (
{\displaystyle \frac {2}{3}} )\,\Gamma (2\,a)\,(2\,\mathrm{sin}(
{\displaystyle \frac {\pi \,(24\,a + 1)}{6}} ) - 1)) }
}

\mapleresult
\begin{maplelatex}
\mapleinline{inert}{2d}{352, ":   ",
-1/2*(sin(1/3*Pi*(6*a+1))+sin(2*Pi*(-3*a+n))+sin(1/3*Pi*(6*a-1))-sin(2
*Pi*a))*V7(a-1/3-1/2*n,n+1)*2^(6*a-n-2)*GAMMA(a-1/3-1/2*n)*GAMMA(a+1/3
-1/2*n)*GAMMA(3/2-2*a+n)*GAMMA(2*a+1/3)*GAMMA(2*a-1/3)*GAMMA(2*a)/Pi^2
/GAMMA(6*a-n-1)/GAMMA(1/2+n)/(sin(Pi*(-4*a+n))+sin(2*Pi*a))+1/4*(2*sin
(1/3*Pi*(6*a-1))-sin(Pi*n)-2*sin(2*Pi*a)+sin(2*Pi*(-3*a+n))+2*sin(1/3*
Pi*(6*a+1)))*V8(a,n)*GAMMA(2*a-2/3-n)*GAMMA(2*a+2/3-n)*GAMMA(1/2+a-1/2
*n)*GAMMA(3/2-2*a+n)*GAMMA(2*a+1/3)*GAMMA(2*a-1/3)/Pi/(2^(2*a-5/3-n))/
(2^(-1/3+2*a-n))/GAMMA(2*a+1/6-n)/GAMMA(2*a+5/6-n)/GAMMA(3*a-1/2*n-1/2
)/GAMMA(1/2+n)/(sin(Pi*(-4*a+n))+sin(2*Pi*a));}{%
\maplemultiline{
352\mbox{:~~~}  - {\displaystyle \frac {1}{2}} (\mathrm{
sin}({\displaystyle \frac {\pi \,(6\,a + 1)}{3}} ) + \mathrm{sin}
(2\,\pi \,( - 3\,a + n)) + \mathrm{sin}({\displaystyle \frac {\pi
 \,(6\,a - 1)}{3}} ) - \mathrm{sin}(2\,\pi \,a)) \\
\mathrm{V7}(a - {\displaystyle \frac {1}{3}}  - {\displaystyle 
\frac {n}{2}} , \,n + 1)\,2^{(6\,a - n - 2)}\,\Gamma (a - 
{\displaystyle \frac {1}{3}}  - {\displaystyle \frac {n}{2}} )\,
\Gamma (a + {\displaystyle \frac {1}{3}}  - {\displaystyle 
\frac {n}{2}} )\,\Gamma ({\displaystyle \frac {3}{2}}  - 2\,a + n
) \\
\Gamma (2\,a + {\displaystyle \frac {1}{3}} )\,\Gamma (2\,a - 
{\displaystyle \frac {1}{3}} )\,\Gamma (2\,a) \left/ {\vrule 
height0.80em width0em depth0.80em} \right. \!  \! (\pi ^{2}\,
\Gamma (6\,a - n - 1)\,\Gamma ({\displaystyle \frac {1}{2}}  + n)
 \\
(\mathrm{sin}(\pi \,( - 4\,a + n)) + \mathrm{sin}(2\,\pi \,a)))
\mbox{} + {\displaystyle \frac {1}{4}} (2\,\mathrm{sin}(
{\displaystyle \frac {\pi \,(6\,a - 1)}{3}} ) - \mathrm{sin}(\pi 
\,n) - 2\,\mathrm{sin}(2\,\pi \,a) \\
\mbox{} + \mathrm{sin}(2\,\pi \,( - 3\,a + n)) + 2\,\mathrm{sin}(
{\displaystyle \frac {\pi \,(6\,a + 1)}{3}} ))\mathrm{V8}(a, \,n)
\,\Gamma (2\,a - {\displaystyle \frac {2}{3}}  - n)\,\Gamma (2\,a
 + {\displaystyle \frac {2}{3}}  - n) \\
\Gamma ({\displaystyle \frac {1}{2}}  + a - {\displaystyle 
\frac {n}{2}} )\,\Gamma ({\displaystyle \frac {3}{2}}  - 2\,a + n
)\,\Gamma (2\,a + {\displaystyle \frac {1}{3}} )\,\Gamma (2\,a - 
{\displaystyle \frac {1}{3}} ) \left/ {\vrule 
height0.80em width0em depth0.80em} \right. \!  \! (\pi \,2^{(2\,a
 - 5/3 - n)} \\
2^{( - 1/3 + 2\,a - n)}\,\Gamma (2\,a + {\displaystyle \frac {1}{
6}}  - n)\,\Gamma (2\,a + {\displaystyle \frac {5}{6}}  - n)\,
\Gamma (3\,a - {\displaystyle \frac {n}{2}}  - {\displaystyle 
\frac {1}{2}} )\,\Gamma ({\displaystyle \frac {1}{2}}  + n)
(\mathrm{sin}(\pi \,( - 4\,a + n)) + \mathrm{sin}(2\,\pi \,a)))
 }
}
\end{maplelatex}

\begin{maplelatex}
\mapleinline{inert}{2d}{353, ":   ",
-2*Pi^(1/2)*(sin(1/2*Pi*(2*a+n))-sin(1/2*Pi*(-2*a+3*n)))*V7(a-1/3-1/2*
n,n+1)*GAMMA(2*a+1/6-n)*GAMMA(2*a+1/3)*GAMMA(2*a)/GAMMA(1/3+a+1/2*n)/G
AMMA(3*a-1/2*n)/GAMMA(1/2+a-1/2*n)/GAMMA(a+1/6-1/2*n)/(sin(1/6*Pi*(-18
*a+3*n+2))+sin(1/2*Pi*(-2*a+3*n))-sin(1/6*Pi*(2+6*a+3*n)))+(sin(1/6*Pi
*(-18*a+3*n+2))-sin(1/2*Pi*(-2*a+3*n))-sin(1/6*Pi*(2+6*a+3*n))+2*sin(1
/2*Pi*(2*a+n)))*V8(a,n)*GAMMA(a-1/2*n+5/6)*GAMMA(2*a+1/3)/GAMMA(2*a+5/
6-n)/GAMMA(1/3+a+1/2*n)/(sin(1/6*Pi*(-18*a+3*n+2))+sin(1/2*Pi*(-2*a+3*
n))-sin(1/6*Pi*(2+6*a+3*n)));}{%
\maplemultiline{
353\mbox{:~~~}  - 2\,\sqrt{\pi }\,(\mathrm{sin}(
{\displaystyle \frac {\pi \,(2\,a + n)}{2}} ) - \mathrm{\%1})\,
\mathrm{V7}(a - {\displaystyle \frac {1}{3}}  - {\displaystyle 
\frac {n}{2}} , \,n + 1)\,\Gamma (2\,a + {\displaystyle \frac {1
}{6}}  - n)\,\Gamma (2\,a + {\displaystyle \frac {1}{3}} ) \\
\Gamma (2\,a) \left/ {\vrule height0.80em width0em depth0.80em}
 \right. \!  \! (\Gamma ({\displaystyle \frac {1}{3}}  + a + 
{\displaystyle \frac {n}{2}} )\,\Gamma (3\,a - {\displaystyle 
\frac {n}{2}} )\,\Gamma ({\displaystyle \frac {1}{2}}  + a - 
{\displaystyle \frac {n}{2}} )\,\Gamma (a + {\displaystyle 
\frac {1}{6}}  - {\displaystyle \frac {n}{2}} ) \\
(\mathrm{sin}({\displaystyle \frac {\pi \,( - 18\,a + 3\,n + 2)}{
6}} ) + \mathrm{\%1} - \mathrm{sin}({\displaystyle \frac {\pi \,(
2 + 6\,a + 3\,n)}{6}} )))\mbox{} +  \\
(\mathrm{sin}({\displaystyle \frac {\pi \,( - 18\,a + 3\,n + 2)}{
6}} ) - \mathrm{\%1} - \mathrm{sin}({\displaystyle \frac {\pi \,(
2 + 6\,a + 3\,n)}{6}} ) + 2\,\mathrm{sin}({\displaystyle \frac {
\pi \,(2\,a + n)}{2}} )) \\
\mathrm{V8}(a, \,n)\,\Gamma (a - {\displaystyle \frac {n}{2}}  + 
{\displaystyle \frac {5}{6}} )\,\Gamma (2\,a + {\displaystyle 
\frac {1}{3}} ) \left/ {\vrule height0.80em width0em depth0.80em}
 \right. \!  \! (\Gamma (2\,a + {\displaystyle \frac {5}{6}}  - n
)\,\Gamma ({\displaystyle \frac {1}{3}}  + a + {\displaystyle 
\frac {n}{2}} ) \\
(\mathrm{sin}({\displaystyle \frac {\pi \,( - 18\,a + 3\,n + 2)}{
6}} ) + \mathrm{\%1} - \mathrm{sin}({\displaystyle \frac {\pi \,(
2 + 6\,a + 3\,n)}{6}} ))) \\
\mathrm{\%1} := \mathrm{sin}({\displaystyle \frac {\pi \,( - 2\,a
 + 3\,n)}{2}} ) }
}
\end{maplelatex}

\begin{maplelatex}
\mapleinline{inert}{2d}{354, ":   ",
2*V7(a-1/3-1/2*n,n+1)*GAMMA(2*a+1/6-n)*GAMMA(a+1/3-1/2*n)*GAMMA(2/3)*G
AMMA(2*a)*(sin(1/6*Pi*(1+6*a+3*n))-sin(1/6*Pi*(-18*a+1+3*n))+3^(1/2)*s
in(1/2*Pi*(-6*a+n)))/Pi^(1/2)/((-1)^(-n))/GAMMA(1/2+a-1/2*n)/GAMMA(a+1
/6-1/2*n)/GAMMA(3*a-1/2*n)/(sin(1/6*Pi*(-6*a+2+3*n))+sin(1/2*Pi*(6*a+n
))+sin(1/6*Pi*(-6*a-2+3*n))+sin(1/2*Pi*(-10*a+n)))-1/2*V8(a,n)*GAMMA(2
*a+2/3-n)*GAMMA(2/3)*(sin(1/6*Pi*(18*a+3*n+1))+sin(1/6*Pi*(9*n+6*a+1))
-sin(1/6*Pi*(-6*a-1+3*n))-sin(1/6*Pi*(-30*a+3*n+1))-sin(1/6*Pi*(-18*a+
1+9*n))+sin(1/6*Pi*(18*a-1+3*n))-2*3^(1/2)*sin(1/2*Pi*(6*a+n))-3^(1/2)
*sin(1/6*Pi*(-6*a+2+3*n))+3^(1/2)*sin(3/2*Pi*(n-2*a))+3^(1/2)*sin(1/2*
Pi*(n-2*a))-3^(1/2)*sin(1/6*Pi*(-6*a-2+3*n)))/Pi^(1/2)/(2^(-1/3+2*a-n)
)/GAMMA(2*a+5/6-n)/(sin(1/6*Pi*(-6*a+2+3*n))+sin(1/2*Pi*(6*a+n))+sin(1
/6*Pi*(-6*a-2+3*n))+sin(1/2*Pi*(-10*a+n)));}{%
\maplemultiline{
354\mbox{:~~~} 2\,\mathrm{V7}(a - {\displaystyle \frac {
1}{3}}  - {\displaystyle \frac {n}{2}} , \,n + 1)\,\Gamma (2\,a
 + {\displaystyle \frac {1}{6}}  - n)\,\Gamma (a + 
{\displaystyle \frac {1}{3}}  - {\displaystyle \frac {n}{2}} )\,
\Gamma ({\displaystyle \frac {2}{3}} )\,\Gamma (2\,a) \\
(\mathrm{sin}({\displaystyle \frac {\pi \,(1 + 6\,a + 3\,n)}{6}} 
) - \mathrm{sin}({\displaystyle \frac {\pi \,( - 18\,a + 1 + 3\,n
)}{6}} ) + \sqrt{3}\,\mathrm{sin}({\displaystyle \frac {\pi \,(
 - 6\,a + n)}{2}} )) \left/ {\vrule 
height0.80em width0em depth0.80em} \right. \!  \! \\ (
\sqrt{\pi }\,(-1)^{( - n)}\,\Gamma ({\displaystyle \frac {1}{2}} 
 + a - {\displaystyle \frac {n}{2}} )\,\Gamma (a + 
{\displaystyle \frac {1}{6}}  - {\displaystyle \frac {n}{2}} )\,
\Gamma (3\,a - {\displaystyle \frac {n}{2}} )(\mathrm{sin}(
{\displaystyle \frac {\pi \,( - 6\,a + 2 + 3\,n)}{6}} ) \\
\mbox{} + \mathrm{sin}({\displaystyle \frac {\pi \,(6\,a + n)}{2}
} ) + \mathrm{sin}({\displaystyle \frac {\pi \,( - 6\,a - 2 + 3\,
n)}{6}} ) + \mathrm{sin}({\displaystyle \frac {\pi \,( - 10\,a + 
n)}{2}} )))\mbox{} - {\displaystyle \frac {1}{2}} \mathrm{V8}(a, 
\,n) \\
\Gamma (2\,a + {\displaystyle \frac {2}{3}}  - n)\,\Gamma (
{\displaystyle \frac {2}{3}} )(\mathrm{sin}({\displaystyle 
\frac {\pi \,(18\,a + 3\,n + 1)}{6}} ) + \mathrm{sin}(
{\displaystyle \frac {\pi \,(9\,n + 6\,a + 1)}{6}} ) \\
\mbox{} - \mathrm{sin}({\displaystyle \frac {\pi \,( - 6\,a - 1
 + 3\,n)}{6}} ) - \mathrm{sin}({\displaystyle \frac {\pi \,( - 30
\,a + 3\,n + 1)}{6}} ) - \mathrm{sin}({\displaystyle \frac {\pi 
\,( - 18\,a + 1 + 9\,n)}{6}} ) \\
\mbox{} + \mathrm{sin}({\displaystyle \frac {\pi \,(18\,a - 1 + 3
\,n)}{6}} ) - 2\,\sqrt{3}\,\mathrm{sin}({\displaystyle \frac {\pi
 \,(6\,a + n)}{2}} ) - \sqrt{3}\,\mathrm{sin}({\displaystyle 
\frac {\pi \,( - 6\,a + 2 + 3\,n)}{6}} ) \\
\mbox{} + \sqrt{3}\,\mathrm{sin}({\displaystyle \frac {3\,\pi \,(
n - 2\,a)}{2}} ) + \sqrt{3}\,\mathrm{sin}({\displaystyle \frac {
\pi \,(n - 2\,a)}{2}} ) - \sqrt{3}\,\mathrm{sin}({\displaystyle 
\frac {\pi \,( - 6\,a - 2 + 3\,n)}{6}} )) \\
 \left/ {\vrule height0.80em width0em depth0.80em} \right. \! 
 \! (\sqrt{\pi }\,2^{( - 1/3 + 2\,a - n)}\,\Gamma (2\,a + 
{\displaystyle \frac {5}{6}}  - n)(\mathrm{sin}({\displaystyle 
\frac {\pi \,( - 6\,a + 2 + 3\,n)}{6}} ) \\
\mbox{} + \mathrm{sin}({\displaystyle \frac {\pi \,(6\,a + n)}{2}
} ) + \mathrm{sin}({\displaystyle \frac {\pi \,( - 6\,a - 2 + 3\,
n)}{6}} ) + \mathrm{sin}({\displaystyle \frac {\pi \,( - 10\,a + 
n)}{2}} ))) }
}
\end{maplelatex}

\begin{maplelatex}
\mapleinline{inert}{2d}{355, ":   ",
1/8*V7(a-1/3-1/2*n,n+1)*GAMMA(a-1/3-1/2*n)*GAMMA(a+1/3-1/2*n)*GAMMA(11
/6+n-2*a)*GAMMA(a-1/2*n)*GAMMA(1/6-a+1/2*n)*GAMMA(2*a+1/3)*GAMMA(2*a-1
/3)*GAMMA(2*a)*(sin(1/3*Pi*(6*a+1))+sin(2*Pi*(-3*a+n))+sin(1/3*Pi*(6*a
-1))-sin(2*Pi*a))*(-6*a+n+2)/sin(1/2*Pi*(-6*a+n))/Pi^(5/2)/GAMMA(3*a-1
/2*n)^2/GAMMA(1/2+n)-6*Pi*(-1)^(-n)*(sin(1/6*Pi*(9*n-6*a+1))-sin(1/6*P
i*(6*a-1+3*n)))*V8(a,n)*GAMMA(a-1/3-1/2*n)*GAMMA(a+1/3-1/2*n)*GAMMA(a-
1/2*n)*GAMMA(2*a+1/3)*GAMMA(2*a-1/3)/GAMMA(2*a-5/6-n)^2/GAMMA(2*a+5/6-
n)/GAMMA(-a+5/6+1/2*n)/GAMMA(1/2+n)/GAMMA(1/2-a+1/2*n)/GAMMA(3*a-1/2*n
-1)/(-cos(3/2*Pi*(n-2*a))+sin(1/6*Pi*(1-6*a+3*n)))/(-12*a+5+6*n);}{%
\maplemultiline{
355\mbox{:~~~} {\displaystyle \frac {1}{8}} \mathrm{V7}(
a - {\displaystyle \frac {1}{3}}  - {\displaystyle \frac {n}{2}} 
, \,n + 1)\,\Gamma (a - {\displaystyle \frac {1}{3}}  - 
{\displaystyle \frac {n}{2}} )\,\Gamma (a + {\displaystyle 
\frac {1}{3}}  - {\displaystyle \frac {n}{2}} )\,\Gamma (
{\displaystyle \frac {11}{6}}  + n - 2\,a)\,\Gamma (a - 
{\displaystyle \frac {n}{2}} ) \\
\Gamma ({\displaystyle \frac {1}{6}}  - a + {\displaystyle 
\frac {n}{2}} )\,\Gamma (2\,a + {\displaystyle \frac {1}{3}} )\,
\Gamma (2\,a - {\displaystyle \frac {1}{3}} )\,\Gamma (2\,a) \\
(\mathrm{sin}({\displaystyle \frac {\pi \,(6\,a + 1)}{3}} ) + 
\mathrm{sin}(2\,\pi \,( - 3\,a + n)) + \mathrm{sin}(
{\displaystyle \frac {\pi \,(6\,a - 1)}{3}} ) - \mathrm{sin}(2\,
\pi \,a)) \\
( - 6\,a + n + 2) \left/ {\vrule 
height0.80em width0em depth0.80em} \right. \!  \! (\mathrm{sin}(
{\displaystyle \frac {\pi \,( - 6\,a + n)}{2}} )\,\pi ^{(5/2)}\,
\Gamma (3\,a - {\displaystyle \frac {n}{2}} )^{2}\,\Gamma (
{\displaystyle \frac {1}{2}}  + n))\mbox{} - 6\,\pi \,(-1)^{( - n
)} \\
(\mathrm{sin}({\displaystyle \frac {\pi \,(9\,n - 6\,a + 1)}{6}} 
) - \mathrm{sin}({\displaystyle \frac {\pi \,(6\,a - 1 + 3\,n)}{6
}} ))\,\mathrm{V8}(a, \,n)\,\Gamma (a - {\displaystyle \frac {1}{
3}}  - {\displaystyle \frac {n}{2}} ) \\
\Gamma (a + {\displaystyle \frac {1}{3}}  - {\displaystyle 
\frac {n}{2}} )\,\Gamma (a - {\displaystyle \frac {n}{2}} )\,
\Gamma (2\,a + {\displaystyle \frac {1}{3}} )\,\Gamma (2\,a - 
{\displaystyle \frac {1}{3}} ) \left/ {\vrule 
height0.80em width0em depth0.80em} \right. \!  \! (\Gamma (2\,a
 - {\displaystyle \frac {5}{6}}  - n)^{2}\,\Gamma (2\,a + 
{\displaystyle \frac {5}{6}}  - n) \\
\Gamma ( - a + {\displaystyle \frac {5}{6}}  + {\displaystyle 
\frac {n}{2}} )\,\Gamma ({\displaystyle \frac {1}{2}}  + n)\,
\Gamma ({\displaystyle \frac {1}{2}}  - a + {\displaystyle 
\frac {n}{2}} )\,\Gamma (3\,a - {\displaystyle \frac {n}{2}}  - 1
) \\
( - \mathrm{cos}({\displaystyle \frac {3\,\pi \,(n - 2\,a)}{2}} )
 + \mathrm{sin}({\displaystyle \frac {\pi \,(1 - 6\,a + 3\,n)}{6}
} ))\,( - 12\,a + 5 + 6\,n)) }
}
\end{maplelatex}

\begin{maplelatex}
\mapleinline{inert}{2d}{356, ":   ",
-1/4*2^(6*a-n-1)*(sin(1/3*Pi*(6*a+1))+sin(2*Pi*(-3*a+n))+sin(1/3*Pi*(6
*a-1))-sin(2*Pi*a))*V7(a-1/3-1/2*n,n+1)*GAMMA(a-1/2*n)*GAMMA(a-1/3-1/2
*n)*GAMMA(11/6+n-2*a)*GAMMA(2*a+1/3)*GAMMA(2*a-1/3)*GAMMA(2*a)/Pi^2/GA
MMA(1/2+n)/GAMMA(6*a-n-1)/(sin(1/3*Pi*(6*a-1))+sin(1/3*Pi*(-12*a-1+3*n
)))+3*(-(-1)^(-n)*sin(1/6*Pi*(-6*a+2+3*n))-(-1)^(-n)*sin(1/6*Pi*(-2+3*
n+18*a))-(-1)^(-n)*sin(1/6*Pi*(-18*a+9*n-2))+(-1)^(-n)*sin(1/6*Pi*(9*n
+6*a-2))+(-1)^(-n)*sin(1/6*Pi*(15*n-30*a-2))+(-1)^(-n)*sin(1/6*Pi*(-30
*a+3*n-2))+(-1)^(-n)*sin(1/6*Pi*(-18*a+2+9*n))+(-1)^(-n)*sin(1/6*Pi*(1
8*a+3*n+2))+4*sin(1/6*Pi*(-18*a-2+3*n))-4*sin(1/6*Pi*(6*a+3*n-2)))*V8(
a,n)*Pi^(1/2)*GAMMA(2*a-n)*GAMMA(a-1/2*n+5/6)*GAMMA(a-1/3-1/2*n)*GAMMA
(2*a+1/3)*GAMMA(2*a-1/3)/(2^(2*a-n-1))/GAMMA(2*a-5/6-n)^2/GAMMA(-a+5/6
+1/2*n)/GAMMA(2*a+5/6-n)/GAMMA(3*a-1/2*n-1/2)/GAMMA(1/2+n)/(sin(1/3*Pi
*(6*a-1))+sin(Pi*n)+sin(1/3*Pi*(-12*a-1+3*n))+sin(2*Pi*(-3*a+n)))/(-12
*a+5+6*n);}{%
\maplemultiline{
356\mbox{:~~~}  - {\displaystyle \frac {1}{4}} 2^{(6\,a
 - n - 1)}\,(\mathrm{sin}({\displaystyle \frac {\pi \,(6\,a + 1)
}{3}} ) + \mathrm{sin}(2\,\pi \,( - 3\,a + n)) + \mathrm{sin}(
{\displaystyle \frac {\pi \,(6\,a - 1)}{3}} ) - \mathrm{sin}(2\,
\pi \,a)) \\
\mathrm{V7}(a - {\displaystyle \frac {1}{3}}  - {\displaystyle 
\frac {n}{2}} , \,n + 1)\,\Gamma (a - {\displaystyle \frac {n}{2}
} )\,\Gamma (a - {\displaystyle \frac {1}{3}}  - {\displaystyle 
\frac {n}{2}} )\,\Gamma ({\displaystyle \frac {11}{6}}  + n - 2\,
a)\,\Gamma (2\,a + {\displaystyle \frac {1}{3}} )\,\Gamma (2\,a
 - {\displaystyle \frac {1}{3}} )\,\Gamma (2\,a) \\
 \left/ {\vrule height0.80em width0em depth0.80em} \right. \! 
 \! (\pi ^{2}\,\Gamma ({\displaystyle \frac {1}{2}}  + n)\,\Gamma
 (6\,a - n - 1)\,(\mathrm{sin}({\displaystyle \frac {\pi \,(6\,a
 - 1)}{3}} ) + \mathrm{sin}({\displaystyle \frac {\pi \,( - 12\,a
 - 1 + 3\,n)}{3}} )))\mbox{} + 3( \\
 - (-1)^{( - n)}\,\mathrm{sin}({\displaystyle \frac {\pi \,( - 6
\,a + 2 + 3\,n)}{6}} ) - (-1)^{( - n)}\,\mathrm{sin}(
{\displaystyle \frac {\pi \,( - 2 + 3\,n + 18\,a)}{6}} ) \\
\mbox{} - (-1)^{( - n)}\,\mathrm{sin}({\displaystyle \frac {\pi 
\,( - 18\,a + 9\,n - 2)}{6}} ) + (-1)^{( - n)}\,\mathrm{sin}(
{\displaystyle \frac {\pi \,(9\,n + 6\,a - 2)}{6}} ) \\
\mbox{} + (-1)^{( - n)}\,\mathrm{sin}({\displaystyle \frac {\pi 
\,(15\,n - 30\,a - 2)}{6}} ) + (-1)^{( - n)}\,\mathrm{sin}(
{\displaystyle \frac {\pi \,( - 30\,a + 3\,n - 2)}{6}} ) \\
\mbox{} + (-1)^{( - n)}\,\mathrm{sin}({\displaystyle \frac {\pi 
\,( - 18\,a + 2 + 9\,n)}{6}} ) + (-1)^{( - n)}\,\mathrm{sin}(
{\displaystyle \frac {\pi \,(18\,a + 3\,n + 2)}{6}} ) \\
\mbox{} + 4\,\mathrm{sin}({\displaystyle \frac {\pi \,( - 18\,a
 - 2 + 3\,n)}{6}} ) - 4\,\mathrm{sin}({\displaystyle \frac {\pi 
\,(6\,a + 3\,n - 2)}{6}} ))\mathrm{V8}(a, \,n)\,\sqrt{\pi }\,
\Gamma (2\,a - n) \\
\Gamma (a - {\displaystyle \frac {n}{2}}  + {\displaystyle 
\frac {5}{6}} )\,\Gamma (a - {\displaystyle \frac {1}{3}}  - 
{\displaystyle \frac {n}{2}} )\,\Gamma (2\,a + {\displaystyle 
\frac {1}{3}} )\,\Gamma (2\,a - {\displaystyle \frac {1}{3}} )
 \left/ {\vrule height0.80em width0em depth0.80em} \right. \! 
 \! (2^{(2\,a - n - 1)}\,\Gamma (2\,a - {\displaystyle \frac {5}{
6}}  - n)^{2} \\
\Gamma ( - a + {\displaystyle \frac {5}{6}}  + {\displaystyle 
\frac {n}{2}} )\,\Gamma (2\,a + {\displaystyle \frac {5}{6}}  - n
)\,\Gamma (3\,a - {\displaystyle \frac {n}{2}}  - {\displaystyle 
\frac {1}{2}} )\,\Gamma ({\displaystyle \frac {1}{2}}  + n) \\
(\mathrm{sin}({\displaystyle \frac {\pi \,(6\,a - 1)}{3}} ) + 
\mathrm{sin}(\pi \,n) + \mathrm{sin}({\displaystyle \frac {\pi \,
( - 12\,a - 1 + 3\,n)}{3}} ) + \mathrm{sin}(2\,\pi \,( - 3\,a + n
))) ( - 12\,a + 5 + 6\,n)) }
}
\end{maplelatex}

\begin{maplelatex}
\mapleinline{inert}{2d}{357, ":   ",
1/2*(2*sin(1/3*Pi*(1-6*a+3*n))-(-1)^(-n)*sin(1/3*Pi*(6*a+1))-(-1)^(-n)
*sin(2*Pi*(-3*a+n))-(-1)^(-n)*sin(1/3*Pi*(-6*a+6*n+1))+(-1)^(-n)*sin(2
*Pi*a))*V7(a-1/3-1/2*n,n+1)*GAMMA(a+1/3-1/2*n)*GAMMA(a-1/3-1/2*n)*GAMM
A(11/6+n-2*a)*GAMMA(a-1/2*n)*GAMMA(2*a+1/3)*GAMMA(2*a-1/3)/((-1)^(-n))
/Pi^(3/2)/GAMMA(3*a-1/2*n)/GAMMA(1/2+n)/GAMMA(2*a-n-1/6)/(sin(1/3*Pi*(
-12*a-1+3*n))-sin(1/3*Pi*(3*n-1)))+Pi*V8(a,n)*GAMMA(a-1/3-1/2*n)*GAMMA
(2*a-n)*GAMMA(2*a+1/3)*GAMMA(2*a-1/3)*GAMMA(2*a+2/3-n)*((-1)^(-n)*sin(
1/6*Pi*(15*n-30*a-2))+(-1)^(-n)*sin(1/6*Pi*(-30*a+3*n-2))+(-1)^(-n)*si
n(1/6*Pi*(-18*a+2+9*n))+(-1)^(-n)*sin(1/6*Pi*(18*a+3*n+2))+2*(-1)^(-n)
*sin(1/2*Pi*(n-2*a))-(-1)^(-n)*sin(1/6*Pi*(-6*a-2+3*n))+(-1)^(-n)*sin(
1/2*Pi*(5*n-2*a))-2*sin(1/2*Pi*(-2*a+3*n))-(-1)^(-n)*sin(1/2*Pi*(6*a+n
))+2*sin(1/6*Pi*(-18*a-2+3*n))-(-1)^(-n)*sin(1/6*Pi*(9*n+6*a-2)))*(-12
*a+5+6*n)/(2^(2*a-n-1))/(2^(2*a+2/3-n))/GAMMA(1/2+n)/GAMMA(2*a-n-1/6)^
2/GAMMA(2*a+1/6-n)^2/GAMMA(-a+5/6+1/2*n)/GAMMA(2*a)/(-sin(1/3*Pi*(6*a-
1))-sin(2*Pi*(-3*a+n))-sin(1/3*Pi*(6*a+1))+sin(2*Pi*(n-a)))/(-12*a+6*n
+1);}{%
\maplemultiline{
357\mbox{:~~~} {\displaystyle \frac {1}{2}} (2\,\mathrm{
sin}({\displaystyle \frac {\pi \,(1 - 6\,a + 3\,n)}{3}} ) - (-1)
^{( - n)}\,\mathrm{sin}({\displaystyle \frac {\pi \,(6\,a + 1)}{3
}} ) \\
\mbox{} - (-1)^{( - n)}\,\mathrm{sin}(2\,\pi \,( - 3\,a + n)) - (
-1)^{( - n)}\,\mathrm{sin}({\displaystyle \frac {\pi \,( - 6\,a
 + 6\,n + 1)}{3}} ) + (-1)^{( - n)}\,\mathrm{sin}(2\,\pi \,a))
 \\
\mathrm{V7}(a - {\displaystyle \frac {1}{3}}  - {\displaystyle 
\frac {n}{2}} , \,n + 1)\,\Gamma (a + {\displaystyle \frac {1}{3}
}  - {\displaystyle \frac {n}{2}} )\,\Gamma (a - {\displaystyle 
\frac {1}{3}}  - {\displaystyle \frac {n}{2}} )\,\Gamma (
{\displaystyle \frac {11}{6}}  + n - 2\,a)\,\Gamma (a - 
{\displaystyle \frac {n}{2}} ) \\
\Gamma (2\,a + {\displaystyle \frac {1}{3}} )\,\Gamma (2\,a - 
{\displaystyle \frac {1}{3}} ) \left/ {\vrule 
height0.80em width0em depth0.80em} \right. \!  \! ((-1)^{( - n)}
\,\pi ^{(3/2)}\,\Gamma (3\,a - {\displaystyle \frac {n}{2}} )\,
\Gamma ({\displaystyle \frac {1}{2}}  + n)\,\Gamma (2\,a - n - 
{\displaystyle \frac {1}{6}} ) \\
(\mathrm{sin}({\displaystyle \frac {\pi \,( - 12\,a - 1 + 3\,n)}{
3}} ) - \mathrm{sin}({\displaystyle \frac {\pi \,(3\,n - 1)}{3}} 
)))\mbox{} + \pi \,\mathrm{V8}(a, \,n)\,\Gamma (a - 
{\displaystyle \frac {1}{3}}  - {\displaystyle \frac {n}{2}} )\,
\Gamma (2\,a - n) \\
\Gamma (2\,a + {\displaystyle \frac {1}{3}} )\,\Gamma (2\,a - 
{\displaystyle \frac {1}{3}} )\,\Gamma (2\,a + {\displaystyle 
\frac {2}{3}}  - n)((-1)^{( - n)}\,\mathrm{sin}({\displaystyle 
\frac {\pi \,(15\,n - 30\,a - 2)}{6}} ) \\
\mbox{} + (-1)^{( - n)}\,\mathrm{sin}({\displaystyle \frac {\pi 
\,( - 30\,a + 3\,n - 2)}{6}} ) + (-1)^{( - n)}\,\mathrm{sin}(
{\displaystyle \frac {\pi \,( - 18\,a + 2 + 9\,n)}{6}} ) \\
\mbox{} + (-1)^{( - n)}\,\mathrm{sin}({\displaystyle \frac {\pi 
\,(18\,a + 3\,n + 2)}{6}} ) + 2\,(-1)^{( - n)}\,\mathrm{sin}(
{\displaystyle \frac {\pi \,(n - 2\,a)}{2}} ) \\
\mbox{} - (-1)^{( - n)}\,\mathrm{sin}({\displaystyle \frac {\pi 
\,( - 6\,a - 2 + 3\,n)}{6}} ) + (-1)^{( - n)}\,\mathrm{sin}(
{\displaystyle \frac {\pi \,(5\,n - 2\,a)}{2}} ) \\
\mbox{} - 2\,\mathrm{sin}({\displaystyle \frac {\pi \,( - 2\,a + 
3\,n)}{2}} ) - (-1)^{( - n)}\,\mathrm{sin}({\displaystyle \frac {
\pi \,(6\,a + n)}{2}} ) + 2\,\mathrm{sin}({\displaystyle \frac {
\pi \,( - 18\,a - 2 + 3\,n)}{6}} ) \\
\mbox{} - (-1)^{( - n)}\,\mathrm{sin}({\displaystyle \frac {\pi 
\,(9\,n + 6\,a - 2)}{6}} ))( - 12\,a + 5 + 6\,n) \left/ {\vrule 
height0.80em width0em depth0.80em} \right. \!  \! (2^{(2\,a - n
 - 1)} \\
2^{(2\,a + 2/3 - n)}\,\Gamma ({\displaystyle \frac {1}{2}}  + n)
\,\Gamma (2\,a - n - {\displaystyle \frac {1}{6}} )^{2}\,\Gamma (
2\,a + {\displaystyle \frac {1}{6}}  - n)^{2}\,\Gamma ( - a + 
{\displaystyle \frac {5}{6}}  + {\displaystyle \frac {n}{2}} )\,
\Gamma (2\,a) \\
( - \mathrm{sin}({\displaystyle \frac {\pi \,(6\,a - 1)}{3}} ) - 
\mathrm{sin}(2\,\pi \,( - 3\,a + n)) - \mathrm{sin}(
{\displaystyle \frac {\pi \,(6\,a + 1)}{3}} ) + \mathrm{sin}(2\,
\pi \,(n - a))) \\
( - 12\,a + 6\,n + 1)) }
}
\end{maplelatex}

\begin{maplelatex}
\mapleinline{inert}{2d}{358, ":   ",
-(-sin(1/6*Pi*(-18*a+3*n+2))-sin(1/2*Pi*(2*a+n))+3^(1/2)*cos(1/2*Pi*(-
6*a+n)))*V7(a-1/3-1/2*n,n+1)*GAMMA(2/3-a-1/2*n)*GAMMA(a+1/3-1/2*n)*GAM
MA(2*a+1/3)*GAMMA(2/3)*GAMMA(2*a)/Pi^(3/2)/((-1)^(-n))/GAMMA(1/2+a-1/2
*n)/GAMMA(3*a-1/2*n)/(sin(1/6*Pi*(-18*a-2+3*n))+sin(1/2*Pi*(2*a+n)))+1
/2*V8(a,n)*(-sin(Pi*(n-2*a))-sin(1/3*Pi*(12*a+1))-sin(1/3*Pi*(6*a-1+3*
n))+3^(1/2)*sin(1/6*Pi*(-12*a+6*n+1)))*GAMMA(2*a+2/3-n)*GAMMA(2/3-a-1/
2*n)*GAMMA(a+1/6-1/2*n)*GAMMA(2*a+1/3)*GAMMA(2/3)/Pi^(3/2)/sin(1/3*Pi*
(6*a-1))/(2^(2*a+2/3-n))/GAMMA(2*a+1/6-n)/GAMMA(2*a+5/6-n);}{%
\maplemultiline{
358\mbox{:~~~}  - ( - \mathrm{sin}({\displaystyle 
\frac {\pi \,( - 18\,a + 3\,n + 2)}{6}} ) - \mathrm{sin}(
{\displaystyle \frac {\pi \,(2\,a + n)}{2}} ) + \sqrt{3}\,
\mathrm{cos}({\displaystyle \frac {\pi \,( - 6\,a + n)}{2}} ))
 \\
\mathrm{V7}(a - {\displaystyle \frac {1}{3}}  - {\displaystyle 
\frac {n}{2}} , \,n + 1)\,\Gamma ({\displaystyle \frac {2}{3}} 
 - a - {\displaystyle \frac {n}{2}} )\,\Gamma (a + 
{\displaystyle \frac {1}{3}}  - {\displaystyle \frac {n}{2}} )\,
\Gamma (2\,a + {\displaystyle \frac {1}{3}} )\,\Gamma (
{\displaystyle \frac {2}{3}} )\,\Gamma (2\,a) \left/ {\vrule 
height0.80em width0em depth0.80em} \right. \!  \! ( \\
\pi ^{(3/2)}\,(-1)^{( - n)}\,\Gamma ({\displaystyle \frac {1}{2}
}  + a - {\displaystyle \frac {n}{2}} )\,\Gamma (3\,a - 
{\displaystyle \frac {n}{2}} ) \\
(\mathrm{sin}({\displaystyle \frac {\pi \,( - 18\,a - 2 + 3\,n)}{
6}} ) + \mathrm{sin}({\displaystyle \frac {\pi \,(2\,a + n)}{2}} 
)))\mbox{} + {\displaystyle \frac {1}{2}} \mathrm{V8}(a, \,n)( - 
\mathrm{sin}(\pi \,(n - 2\,a)) \\
\mbox{} - \mathrm{sin}({\displaystyle \frac {\pi \,(12\,a + 1)}{3
}} ) - \mathrm{sin}({\displaystyle \frac {\pi \,(6\,a - 1 + 3\,n)
}{3}} ) + \sqrt{3}\,\mathrm{sin}({\displaystyle \frac {\pi \,( - 
12\,a + 6\,n + 1)}{6}} )) \\
\Gamma (2\,a + {\displaystyle \frac {2}{3}}  - n)\,\Gamma (
{\displaystyle \frac {2}{3}}  - a - {\displaystyle \frac {n}{2}} 
)\,\Gamma (a + {\displaystyle \frac {1}{6}}  - {\displaystyle 
\frac {n}{2}} )\,\Gamma (2\,a + {\displaystyle \frac {1}{3}} )\,
\Gamma ({\displaystyle \frac {2}{3}} ) \left/ {\vrule 
height0.80em width0em depth0.80em} \right. \!  \! (\pi ^{(3/2)}
 \\
\mathrm{sin}({\displaystyle \frac {\pi \,(6\,a - 1)}{3}} )\,2^{(2
\,a + 2/3 - n)}\,\Gamma (2\,a + {\displaystyle \frac {1}{6}}  - n
)\,\Gamma (2\,a + {\displaystyle \frac {5}{6}}  - n)) }
}
\end{maplelatex}

\begin{maplelatex}
\mapleinline{inert}{2d}{359, ":   ",
1/3*V7(a-1/3-1/2*n,n+1)*(3^(1/2)*sin(1/6*Pi*(-18*a+3*n+2))+3^(1/2)*sin
(1/2*Pi*(2*a+n))-3*cos(1/2*Pi*(-6*a+n)))*GAMMA(-a+5/6+1/2*n)*GAMMA(a+1
/3-1/2*n)*GAMMA(2/3-a-1/2*n)*GAMMA(2*a+1/3)*GAMMA(2*a)/Pi^(3/2)/sin(1/
3*Pi*(6*a-1))/((-1)^(-n))/GAMMA(3*a-1/2*n)/GAMMA(2/3)-1/3*V8(a,n)*(3^(
1/2)*sin(Pi*(n-2*a))+3^(1/2)*sin(1/3*Pi*(12*a+1))+3^(1/2)*sin(1/3*Pi*(
6*a-1+3*n))-3*sin(1/6*Pi*(-12*a+6*n+1)))*GAMMA(2*a+2/3-n)*GAMMA(1/2+a-
1/2*n)*GAMMA(2/3-a-1/2*n)*GAMMA(2*a+1/3)/Pi^(1/2)/sin(1/3*Pi*(6*a-1))/
(2^(2*a+2/3-n))/GAMMA(2*a+1/6-n)/GAMMA(2*a+5/6-n)/GAMMA(2/3);}{%
\maplemultiline{
359\mbox{:~~~} {\displaystyle \frac {1}{3}} \mathrm{V7}(
a - {\displaystyle \frac {1}{3}}  - {\displaystyle \frac {n}{2}} 
, \,n + 1) \\
(\sqrt{3}\,\mathrm{sin}({\displaystyle \frac {\pi \,( - 18\,a + 3
\,n + 2)}{6}} ) + \sqrt{3}\,\mathrm{sin}({\displaystyle \frac {
\pi \,(2\,a + n)}{2}} ) - 3\,\mathrm{cos}({\displaystyle \frac {
\pi \,( - 6\,a + n)}{2}} )) \\
\Gamma ( - a + {\displaystyle \frac {5}{6}}  + {\displaystyle 
\frac {n}{2}} )\,\Gamma (a + {\displaystyle \frac {1}{3}}  - 
{\displaystyle \frac {n}{2}} )\,\Gamma ({\displaystyle \frac {2}{
3}}  - a - {\displaystyle \frac {n}{2}} )\,\Gamma (2\,a + 
{\displaystyle \frac {1}{3}} )\,\Gamma (2\,a) \left/ {\vrule 
height0.80em width0em depth0.80em} \right. \!  \! (\pi ^{(3/2)}
 \\
\mathrm{sin}({\displaystyle \frac {\pi \,(6\,a - 1)}{3}} )\,(-1)
^{( - n)}\,\Gamma (3\,a - {\displaystyle \frac {n}{2}} )\,\Gamma 
({\displaystyle \frac {2}{3}} ))\mbox{} - {\displaystyle \frac {1
}{3}} \mathrm{V8}(a, \,n)(\sqrt{3}\,\mathrm{sin}(\pi \,(n - 2\,a)
) \\
\mbox{} + \sqrt{3}\,\mathrm{sin}({\displaystyle \frac {\pi \,(12
\,a + 1)}{3}} ) + \sqrt{3}\,\mathrm{sin}({\displaystyle \frac {
\pi \,(6\,a - 1 + 3\,n)}{3}} ) - 3\,\mathrm{sin}({\displaystyle 
\frac {\pi \,( - 12\,a + 6\,n + 1)}{6}} )) \\
\Gamma (2\,a + {\displaystyle \frac {2}{3}}  - n)\,\Gamma (
{\displaystyle \frac {1}{2}}  + a - {\displaystyle \frac {n}{2}} 
)\,\Gamma ({\displaystyle \frac {2}{3}}  - a - {\displaystyle 
\frac {n}{2}} )\,\Gamma (2\,a + {\displaystyle \frac {1}{3}} )
 \left/ {\vrule height0.80em width0em depth0.80em} \right. \! 
 \! (\sqrt{\pi }\,\mathrm{sin}({\displaystyle \frac {\pi \,(6\,a
 - 1)}{3}} ) \\
2^{(2\,a + 2/3 - n)}\,\Gamma (2\,a + {\displaystyle \frac {1}{6}
}  - n)\,\Gamma (2\,a + {\displaystyle \frac {5}{6}}  - n)\,
\Gamma ({\displaystyle \frac {2}{3}} )) }
}
\end{maplelatex}

\begin{maplelatex}
\mapleinline{inert}{2d}{360, ":   ",
1/3*V7(a-1/3-1/2*n,n+1)/Pi^(1/2)*(3^(1/2)*sin(1/6*Pi*(-18*a+3*n+2))+3^
(1/2)*sin(1/2*Pi*(2*a+n))-3*cos(1/2*Pi*(-6*a+n)))/sin(1/3*Pi*(6*a-1))/
GAMMA(1/2+n)/GAMMA(3*a-1/2*n)*GAMMA(-a+5/6+1/2*n)*GAMMA(a+1/3-1/2*n)/G
AMMA(1/2+a-1/2*n)*GAMMA(2*a+1/3)*GAMMA(2*a)-1/6*(-1)^(-n)*(3^(1/2)*sin
(Pi*(n-2*a))+3^(1/2)*sin(1/3*Pi*(12*a+1))+3^(1/2)*sin(1/3*Pi*(6*a-1+3*
n))-3*sin(1/6*Pi*(-12*a+6*n+1)))*V8(a,n)*Pi^(1/2)*GAMMA(2*a+2/3-n)*GAM
MA(2*a+1/3)/sin(1/3*Pi*(6*a-1))/(2^(-1/3+2*a-n))/GAMMA(1/2+n)/GAMMA(2*
a+5/6-n)/GAMMA(2*a+1/6-n);}{%
\maplemultiline{
360\mbox{:~~~} {\displaystyle \frac {1}{3}} \mathrm{V7}(
a - {\displaystyle \frac {1}{3}}  - {\displaystyle \frac {n}{2}} 
, \,n + 1) \\
(\sqrt{3}\,\mathrm{sin}({\displaystyle \frac {\pi \,( - 18\,a + 3
\,n + 2)}{6}} ) + \sqrt{3}\,\mathrm{sin}({\displaystyle \frac {
\pi \,(2\,a + n)}{2}} ) - 3\,\mathrm{cos}({\displaystyle \frac {
\pi \,( - 6\,a + n)}{2}} )) \\
\Gamma ( - a + {\displaystyle \frac {5}{6}}  + {\displaystyle 
\frac {n}{2}} )\,\Gamma (a + {\displaystyle \frac {1}{3}}  - 
{\displaystyle \frac {n}{2}} )\,\Gamma (2\,a + {\displaystyle 
\frac {1}{3}} )\,\Gamma (2\,a) \left/ {\vrule 
height0.80em width0em depth0.80em} \right. \!  \! (\sqrt{\pi }\,
\mathrm{sin}({\displaystyle \frac {\pi \,(6\,a - 1)}{3}} )\,
\Gamma ({\displaystyle \frac {1}{2}}  + n) \\
\Gamma (3\,a - {\displaystyle \frac {n}{2}} )\,\Gamma (
{\displaystyle \frac {1}{2}}  + a - {\displaystyle \frac {n}{2}} 
))\mbox{} - {\displaystyle \frac {1}{6}} (-1)^{( - n)}(\sqrt{3}\,
\mathrm{sin}(\pi \,(n - 2\,a)) + \sqrt{3}\,\mathrm{sin}(
{\displaystyle \frac {\pi \,(12\,a + 1)}{3}} ) \\
\mbox{} + \sqrt{3}\,\mathrm{sin}({\displaystyle \frac {\pi \,(6\,
a - 1 + 3\,n)}{3}} ) - 3\,\mathrm{sin}({\displaystyle \frac {\pi 
\,( - 12\,a + 6\,n + 1)}{6}} ))\mathrm{V8}(a, \,n)\,\sqrt{\pi }
 \\
\Gamma (2\,a + {\displaystyle \frac {2}{3}}  - n)\,\Gamma (2\,a
 + {\displaystyle \frac {1}{3}} ) \left/ {\vrule 
height0.80em width0em depth0.80em} \right. \!  \! (\mathrm{sin}(
{\displaystyle \frac {\pi \,(6\,a - 1)}{3}} )\,2^{( - 1/3 + 2\,a
 - n)}\,\Gamma ({\displaystyle \frac {1}{2}}  + n) \\
\Gamma (2\,a + {\displaystyle \frac {5}{6}}  - n)\,\Gamma (2\,a
 + {\displaystyle \frac {1}{6}}  - n)) }
}
\end{maplelatex}

\begin{maplelatex}
\mapleinline{inert}{2d}{361, ":   ",
-1/2*(2*sin(Pi*(-6*a+n))+(-1)^(-n)*sin(6*Pi*a))*V7(a-1/3-1/2*n,n+1)*2^
(6*a-n-2)*GAMMA(a+1/3-1/2*n)*GAMMA(a-1/2*n)*GAMMA(7/6-2*a+n)*GAMMA(2*a
+1/3)*GAMMA(2*a-1/3)*GAMMA(2*a)/((-1)^(-n))/Pi^2/GAMMA(6*a-n-1)/GAMMA(
1/2+n)/(sin(1/3*Pi*(-12*a+1+3*n))+sin(1/3*Pi*(6*a+1)))+1/6*V8(a,n)*Pi*
GAMMA(2*a-n)*GAMMA(2*a+2/3-n)*GAMMA(2*a+1/3)*GAMMA(2*a-1/3)*(-sin(1/6*
Pi*(6*a+3*n-2))+sin(1/6*Pi*(-30*a+9*n-2))+sin(1/2*Pi*(-6*a+n)))*(-12*a
+6*n+1)/(2^(2*a-n-1))/(2^(-1/3+2*a-n))/GAMMA(2*a+1/6-n)/GAMMA(2*a+5/6-
n)^2/GAMMA(-a+5/6+1/2*n)/GAMMA(3*a-1/2*n-1/2)/GAMMA(1/2+n)/(sin(1/3*Pi
*(6*a+1))+sin(2*Pi*(-3*a+n))+sin(1/3*Pi*(-12*a+1+3*n))+sin(Pi*n));}{%
\maplemultiline{
361\mbox{:~~~}  - {\displaystyle \frac {1}{2}} (2\,
\mathrm{sin}(\pi \,( - 6\,a + n)) + (-1)^{( - n)}\,\mathrm{sin}(6
\,\pi \,a))\,\mathrm{V7}(a - {\displaystyle \frac {1}{3}}  - 
{\displaystyle \frac {n}{2}} , \,n + 1)\,2^{(6\,a - n - 2)} \\
\Gamma (a + {\displaystyle \frac {1}{3}}  - {\displaystyle 
\frac {n}{2}} )\,\Gamma (a - {\displaystyle \frac {n}{2}} )\,
\Gamma ({\displaystyle \frac {7}{6}}  - 2\,a + n)\,\Gamma (2\,a
 + {\displaystyle \frac {1}{3}} )\,\Gamma (2\,a - {\displaystyle 
\frac {1}{3}} )\,\Gamma (2\,a) \left/ {\vrule 
height0.80em width0em depth0.80em} \right. \!  \! ((-1)^{( - n)}
 \\
\pi ^{2}\,\Gamma (6\,a - n - 1)\,\Gamma ({\displaystyle \frac {1
}{2}}  + n)\,(\mathrm{sin}({\displaystyle \frac {\pi \,( - 12\,a
 + 1 + 3\,n)}{3}} ) + \mathrm{sin}({\displaystyle \frac {\pi \,(6
\,a + 1)}{3}} )))\mbox{} + {\displaystyle \frac {1}{6}}  \\
\mathrm{V8}(a, \,n)\,\pi \,\Gamma (2\,a - n)\,\Gamma (2\,a + 
{\displaystyle \frac {2}{3}}  - n)\,\Gamma (2\,a + 
{\displaystyle \frac {1}{3}} )\,\Gamma (2\,a - {\displaystyle 
\frac {1}{3}} ) \\
( - \mathrm{sin}({\displaystyle \frac {\pi \,(6\,a + 3\,n - 2)}{6
}} ) + \mathrm{sin}({\displaystyle \frac {\pi \,( - 30\,a + 9\,n
 - 2)}{6}} ) + \mathrm{sin}({\displaystyle \frac {\pi \,( - 6\,a
 + n)}{2}} )) \\
( - 12\,a + 6\,n + 1) \left/ {\vrule 
height0.80em width0em depth0.80em} \right. \!  \! (2^{(2\,a - n
 - 1)}\,2^{( - 1/3 + 2\,a - n)}\,\Gamma (2\,a + {\displaystyle 
\frac {1}{6}}  - n)\,\Gamma (2\,a + {\displaystyle \frac {5}{6}} 
 - n)^{2} \\
\Gamma ( - a + {\displaystyle \frac {5}{6}}  + {\displaystyle 
\frac {n}{2}} )\,\Gamma (3\,a - {\displaystyle \frac {n}{2}}  - 
{\displaystyle \frac {1}{2}} )\,\Gamma ({\displaystyle \frac {1}{
2}}  + n) \\
(\mathrm{sin}({\displaystyle \frac {\pi \,(6\,a + 1)}{3}} ) + 
\mathrm{sin}(2\,\pi \,( - 3\,a + n)) + \mathrm{sin}(
{\displaystyle \frac {\pi \,( - 12\,a + 1 + 3\,n)}{3}} ) + 
\mathrm{sin}(\pi \,n))) }
}
\end{maplelatex}

\begin{maplelatex}
\mapleinline{inert}{2d}{362, ":   ",
-1/4*V7(a-1/3-1/2*n,n+1)*GAMMA(a-1/3-1/2*n)*GAMMA(a+1/3-1/2*n)*GAMMA(a
-1/2*n)*GAMMA(2*a+1/3)*GAMMA(2*a-1/3)*GAMMA(2*a)*(2*sin(Pi*(-6*a+n))+(
-1)^(-n)*sin(6*Pi*a))*(-6*a+n+2)/((-1)^(-n))/Pi^(3/2)/GAMMA(3*a-1/2*n)
^2/GAMMA(3*a-1/2*n-1/2)/sin(Pi*(-6*a+n))-1/4*2^(6*a-3-n)*(2*sin(1/6*Pi
*(-12*a+6*n+1))+1)*V8(a,n)*GAMMA(2*a-n)*GAMMA(2*a+2/3-n)*GAMMA(a-1/3-1
/2*n)*GAMMA(2*a+1/3)*GAMMA(2*a-1/3)/(2^(-1/3+2*a-n))/(2^(2*a-n-1))/Pi^
(1/2)/cos(1/2*Pi*(-6*a+n))/GAMMA(6*a-n-2)/GAMMA(2*a+1/6-n)/GAMMA(2*a+5
/6-n)/GAMMA(-a+5/6+1/2*n);}{%
\maplemultiline{
362\mbox{:~~~}  - {\displaystyle \frac {1}{4}} \mathrm{
V7}(a - {\displaystyle \frac {1}{3}}  - {\displaystyle \frac {n}{
2}} , \,n + 1)\,\Gamma (a - {\displaystyle \frac {1}{3}}  - 
{\displaystyle \frac {n}{2}} )\,\Gamma (a + {\displaystyle 
\frac {1}{3}}  - {\displaystyle \frac {n}{2}} )\,\Gamma (a - 
{\displaystyle \frac {n}{2}} )\,\Gamma (2\,a + {\displaystyle 
\frac {1}{3}} ) \\
\Gamma (2\,a - {\displaystyle \frac {1}{3}} )\,\Gamma (2\,a)\,(2
\,\mathrm{sin}(\pi \,( - 6\,a + n)) + (-1)^{( - n)}\,\mathrm{sin}
(6\,\pi \,a))\,( - 6\,a + n + 2) \left/ {\vrule 
height0.80em width0em depth0.80em} \right. \!  \! ( \\
(-1)^{( - n)}\,\pi ^{(3/2)}\,\Gamma (3\,a - {\displaystyle 
\frac {n}{2}} )^{2}\,\Gamma (3\,a - {\displaystyle \frac {n}{2}} 
 - {\displaystyle \frac {1}{2}} )\,\mathrm{sin}(\pi \,( - 6\,a + 
n)))\mbox{} - {\displaystyle \frac {1}{4}} 2^{(6\,a - 3 - n)} \\
(2\,\mathrm{sin}({\displaystyle \frac {\pi \,( - 12\,a + 6\,n + 1
)}{6}} ) + 1)\,\mathrm{V8}(a, \,n)\,\Gamma (2\,a - n)\,\Gamma (2
\,a + {\displaystyle \frac {2}{3}}  - n)\,\Gamma (a - 
{\displaystyle \frac {1}{3}}  - {\displaystyle \frac {n}{2}} )
 \\
\Gamma (2\,a + {\displaystyle \frac {1}{3}} )\,\Gamma (2\,a - 
{\displaystyle \frac {1}{3}} ) \left/ {\vrule 
height0.80em width0em depth0.80em} \right. \!  \! (2^{( - 1/3 + 2
\,a - n)}\,2^{(2\,a - n - 1)}\,\sqrt{\pi }\,\mathrm{cos}(
{\displaystyle \frac {\pi \,( - 6\,a + n)}{2}} ) \\
\Gamma (6\,a - n - 2)\,\Gamma (2\,a + {\displaystyle \frac {1}{6}
}  - n)\,\Gamma (2\,a + {\displaystyle \frac {5}{6}}  - n)\,
\Gamma ( - a + {\displaystyle \frac {5}{6}}  + {\displaystyle 
\frac {n}{2}} )) }
}
\end{maplelatex}

\begin{maplelatex}
\mapleinline{inert}{2d}{363, ":   ",
1/8*V7(a-1/3-1/2*n,n+1)*GAMMA(7/6-2*a+n)*GAMMA(-a+5/6+1/2*n)*GAMMA(a-1
/3-1/2*n)*GAMMA(a+1/3-1/2*n)*GAMMA(a-1/2*n)*GAMMA(2*a+1/3)*GAMMA(2*a-1
/3)*GAMMA(2*a)*(2*sin(Pi*(-6*a+n))+(-1)^(-n)*sin(6*Pi*a))*(-6*a+n+2)/s
in(1/2*Pi*(-6*a+n))/((-1)^(-n))/Pi^(5/2)/GAMMA(3*a-1/2*n)^2/GAMMA(1/2+
n)+1/24*V8(a,n)*GAMMA(2*a-n)*GAMMA(2*a+2/3-n)*GAMMA(a-1/3-1/2*n)*GAMMA
(2*a+1/3)*GAMMA(2*a-1/3)*(2*sin(1/6*Pi*(-12*a+6*n+1))+1)*(-12*a+6*n+1)
/sin(1/6*Pi*(-12*a+6*n+1))/(2^(2*a-n-1))/(2^(-1/3+2*a-n))/GAMMA(2*a+1/
6-n)/GAMMA(2*a+5/6-n)^2/GAMMA(3*a-1/2*n-1)/GAMMA(1/2+n);}{%
\maplemultiline{
363\mbox{:~~~} {\displaystyle \frac {1}{8}} \mathrm{V7}(
a - {\displaystyle \frac {1}{3}}  - {\displaystyle \frac {n}{2}} 
, \,n + 1)\,\Gamma ({\displaystyle \frac {7}{6}}  - 2\,a + n)\,
\Gamma ( - a + {\displaystyle \frac {5}{6}}  + {\displaystyle 
\frac {n}{2}} )\,\Gamma (a - {\displaystyle \frac {1}{3}}  - 
{\displaystyle \frac {n}{2}} )\,\Gamma (a + {\displaystyle 
\frac {1}{3}}  - {\displaystyle \frac {n}{2}} ) \\
\Gamma (a - {\displaystyle \frac {n}{2}} )\,\Gamma (2\,a + 
{\displaystyle \frac {1}{3}} )\,\Gamma (2\,a - {\displaystyle 
\frac {1}{3}} )\,\Gamma (2\,a)\,(2\,\mathrm{sin}(\pi \,( - 6\,a
 + n)) + (-1)^{( - n)}\,\mathrm{sin}(6\,\pi \,a)) \\
( - 6\,a + n + 2) \left/ {\vrule 
height0.80em width0em depth0.80em} \right. \!  \! (\mathrm{sin}(
{\displaystyle \frac {\pi \,( - 6\,a + n)}{2}} )\,(-1)^{( - n)}\,
\pi ^{(5/2)}\,\Gamma (3\,a - {\displaystyle \frac {n}{2}} )^{2}\,
\Gamma ({\displaystyle \frac {1}{2}}  + n))\mbox{} + 
{\displaystyle \frac {1}{24}}  \\
\mathrm{V8}(a, \,n)\,\Gamma (2\,a - n)\,\Gamma (2\,a + 
{\displaystyle \frac {2}{3}}  - n)\,\Gamma (a - {\displaystyle 
\frac {1}{3}}  - {\displaystyle \frac {n}{2}} )\,\Gamma (2\,a + 
{\displaystyle \frac {1}{3}} )\,\Gamma (2\,a - {\displaystyle 
\frac {1}{3}} ) \\
(2\,\mathrm{sin}({\displaystyle \frac {\pi \,( - 12\,a + 6\,n + 1
)}{6}} ) + 1)\,( - 12\,a + 6\,n + 1) \left/ {\vrule 
height0.80em width0em depth0.80em} \right. \!  \! (\mathrm{sin}(
{\displaystyle \frac {\pi \,( - 12\,a + 6\,n + 1)}{6}} ) \\
2^{(2\,a - n - 1)}\,2^{( - 1/3 + 2\,a - n)}\,\Gamma (2\,a + 
{\displaystyle \frac {1}{6}}  - n)\,\Gamma (2\,a + 
{\displaystyle \frac {5}{6}}  - n)^{2}\,\Gamma (3\,a - 
{\displaystyle \frac {n}{2}}  - 1)\,\Gamma ({\displaystyle 
\frac {1}{2}}  + n)) }
}
\end{maplelatex}

\begin{maplelatex}
\mapleinline{inert}{2d}{364, ":   ",
-1/4*(2*sin(Pi*(-6*a+n))+(-1)^(-n)*sin(6*Pi*a))*V7(a-1/3-1/2*n,n+1)*GA
MMA(11/6+n-2*a)*GAMMA(3/2-2*a+n)*GAMMA(a-1/3-1/2*n)*GAMMA(a+1/3-1/2*n)
*GAMMA(a-1/2*n)*GAMMA(2*a-1/3)*GAMMA(2*a)/sin(1/3*Pi*(6*a+1))/((-1)^(-
n))/Pi^(5/2)/GAMMA(3*a-1/2*n)/GAMMA(1/2+n)-1/12*V8(a,n)*GAMMA(2*a+2/3-
n)*GAMMA(a-1/3-1/2*n)*GAMMA(3/2-2*a+n)*GAMMA(2*a-n)*GAMMA(2*a-1/3)*(-s
in(1/6*Pi*(6*a+3*n-2))+sin(1/6*Pi*(-30*a+9*n-2))+sin(1/2*Pi*(-6*a+n)))
*(-12*a+5+6*n)/(2^(-1/3+2*a-n))/(2^(2*a-n-1))/GAMMA(2*a+1/6-n)^2/GAMMA
(-a+5/6+1/2*n)/GAMMA(2*a+5/6-n)/GAMMA(1/2+n)/(sin(1/3*Pi*(3*n-1))+sin(
Pi*(-4*a+n)));}{%
\maplemultiline{
364\mbox{:~~~}  - {\displaystyle \frac {1}{4}} (2\,
\mathrm{sin}(\pi \,( - 6\,a + n)) + (-1)^{( - n)}\,\mathrm{sin}(6
\,\pi \,a))\,\mathrm{V7}(a - {\displaystyle \frac {1}{3}}  - 
{\displaystyle \frac {n}{2}} , \,n + 1) \\
\Gamma ({\displaystyle \frac {11}{6}}  + n - 2\,a)\,\Gamma (
{\displaystyle \frac {3}{2}}  - 2\,a + n)\,\Gamma (a - 
{\displaystyle \frac {1}{3}}  - {\displaystyle \frac {n}{2}} )\,
\Gamma (a + {\displaystyle \frac {1}{3}}  - {\displaystyle 
\frac {n}{2}} )\,\Gamma (a - {\displaystyle \frac {n}{2}} )\,
\Gamma (2\,a - {\displaystyle \frac {1}{3}} ) \\
\Gamma (2\,a) \left/ {\vrule height0.80em width0em depth0.80em}
 \right. \!  \! (\mathrm{sin}({\displaystyle \frac {\pi \,(6\,a
 + 1)}{3}} )\,(-1)^{( - n)}\,\pi ^{(5/2)}\,\Gamma (3\,a - 
{\displaystyle \frac {n}{2}} )\,\Gamma ({\displaystyle \frac {1}{
2}}  + n))\mbox{} - {\displaystyle \frac {1}{12}} \mathrm{V8}(a, 
\,n) \\
\Gamma (2\,a + {\displaystyle \frac {2}{3}}  - n)\,\Gamma (a - 
{\displaystyle \frac {1}{3}}  - {\displaystyle \frac {n}{2}} )\,
\Gamma ({\displaystyle \frac {3}{2}}  - 2\,a + n)\,\Gamma (2\,a
 - n)\,\Gamma (2\,a - {\displaystyle \frac {1}{3}} ) \\
( - \mathrm{sin}({\displaystyle \frac {\pi \,(6\,a + 3\,n - 2)}{6
}} ) + \mathrm{sin}({\displaystyle \frac {\pi \,( - 30\,a + 9\,n
 - 2)}{6}} ) + \mathrm{sin}({\displaystyle \frac {\pi \,( - 6\,a
 + n)}{2}} )) \\
( - 12\,a + 5 + 6\,n) \left/ {\vrule 
height0.80em width0em depth0.80em} \right. \!  \! (2^{( - 1/3 + 2
\,a - n)}\,2^{(2\,a - n - 1)}\,\Gamma (2\,a + {\displaystyle 
\frac {1}{6}}  - n)^{2}\,\Gamma ( - a + {\displaystyle \frac {5}{
6}}  + {\displaystyle \frac {n}{2}} ) \\
\Gamma (2\,a + {\displaystyle \frac {5}{6}}  - n)\,\Gamma (
{\displaystyle \frac {1}{2}}  + n)\,(\mathrm{sin}({\displaystyle 
\frac {\pi \,(3\,n - 1)}{3}} ) + \mathrm{sin}(\pi \,( - 4\,a + n)
))) }
}
\end{maplelatex}

\begin{maplelatex}
\mapleinline{inert}{2d}{366, ":   ",
GAMMA(1/2*c+1/2*b+n+m-1/2-a)*GAMMA(1+a-b)*GAMMA(1+a-c)*WPMR(2*m+2*n-2,
2*n-1,c,b,a)*GAMMA(2*n-2*a)/GAMMA(-a+2*n)/GAMMA(1/2*c+1/2*b+n+m-1/2)/G
AMMA(-c+1)/GAMMA(-b+1)+GAMMA(1/2*c+1/2*b+n+m-1/2-a)*GAMMA(1+a-b)*GAMMA
(1+a-c)*GAMMA(2*a-2*n)*WPMR(2*m+2*n-2,2*n-1,c-2*a+2*n,-2*a+2*n+b,-a+2*
n)/GAMMA(-b+2*a-2*n+1)/GAMMA(2*a-2*n+1-c)/GAMMA(-1/2+1/2*c+1/2*b+3*n+m
-2*a)/GAMMA(a);}{%
\maplemultiline{
366\mbox{:~~~} \Gamma ({\displaystyle \frac {c}{2}}  + 
{\displaystyle \frac {b}{2}}  + n + m - {\displaystyle \frac {1}{
2}}  - a)\,\Gamma (1 + a - b)\,\Gamma (1 + a - c) \\
\mathrm{WPMR}(2\,m + 2\,n - 2, \,2\,n - 1, \,c, \,b, \,a)\,\Gamma
 (2\,n - 2\,a) \left/ {\vrule height0.80em width0em depth0.80em}
 \right. \!  \! (\Gamma ( - a + 2\,n) \\
\Gamma ({\displaystyle \frac {c}{2}}  + {\displaystyle \frac {b}{
2}}  + n + m - {\displaystyle \frac {1}{2}} )\,\Gamma ( - c + 1)
\,\Gamma ( - b + 1))\mbox{} + \Gamma ({\displaystyle \frac {c}{2}
}  + {\displaystyle \frac {b}{2}}  + n + m - {\displaystyle 
\frac {1}{2}}  - a) \\
\Gamma (1 + a - b)\,\Gamma (1 + a - c)\,\Gamma (2\,a - 2\,n) \\
\mathrm{WPMR}(2\,m + 2\,n - 2, \,2\,n - 1, \,c - 2\,a + 2\,n, \,
 - 2\,a + 2\,n + b, \, - a + 2\,n) \left/ {\vrule 
height0.80em width0em depth0.80em} \right. \!  \! ( \\
\Gamma ( - b + 2\,a - 2\,n + 1)\,\Gamma (2\,a - 2\,n + 1 - c)\,
\Gamma ( - {\displaystyle \frac {1}{2}}  + {\displaystyle \frac {
c}{2}}  + {\displaystyle \frac {b}{2}}  + 3\,n + m - 2\,a)\,
\Gamma (a)) }
}
\end{maplelatex}

\begin{maplelatex}
\mapleinline{inert}{2d}{367, ":   ",
GAMMA(1/2*c+1/2*b+n+m-a)*GAMMA(1+a-b)*GAMMA(1+a-c)*WPMR(2*m+2*n-1,2*n-
1,c,b,a)*GAMMA(2*n-2*a)/GAMMA(-a+2*n)/GAMMA(1/2*c+1/2*b+n+m)/GAMMA(-c+
1)/GAMMA(-b+1)+GAMMA(1/2*c+1/2*b+n+m-a)*GAMMA(1+a-b)*GAMMA(1+a-c)*GAMM
A(2*a-2*n)*WPMR(2*m+2*n-1,2*n-1,c-2*a+2*n,-2*a+2*n+b,-a+2*n)/GAMMA(-b+
2*a-2*n+1)/GAMMA(2*a-2*n+1-c)/GAMMA(1/2*c+1/2*b+3*n+m-2*a)/GAMMA(a);}{
\maplemultiline{
367\mbox{:~~~} \Gamma ({\displaystyle \frac {c}{2}}  + 
{\displaystyle \frac {b}{2}}  + n + m - a)\,\Gamma (1 + a - b)\,
\Gamma (1 + a - c) \\
\mathrm{WPMR}(2\,m + 2\,n - 1, \,2\,n - 1, \,c, \,b, \,a)\,\Gamma
 (2\,n - 2\,a) \left/ {\vrule height0.80em width0em depth0.80em}
 \right. \!  \! (\Gamma ( - a + 2\,n) \\
\Gamma ({\displaystyle \frac {c}{2}}  + {\displaystyle \frac {b}{
2}}  + n + m)\,\Gamma ( - c + 1)\,\Gamma ( - b + 1))\mbox{} + 
\Gamma ({\displaystyle \frac {c}{2}}  + {\displaystyle \frac {b}{
2}}  + n + m - a)\,\Gamma (1 + a - b) \\
\Gamma (1 + a - c)\,\Gamma (2\,a - 2\,n) \\
\mathrm{WPMR}(2\,m + 2\,n - 1, \,2\,n - 1, \,c - 2\,a + 2\,n, \,
 - 2\,a + 2\,n + b, \, - a + 2\,n) \left/ {\vrule 
height0.80em width0em depth0.80em} \right. \!  \! ( \\
\Gamma ( - b + 2\,a - 2\,n + 1)\,\Gamma (2\,a - 2\,n + 1 - c)\,
\Gamma ({\displaystyle \frac {c}{2}}  + {\displaystyle \frac {b}{
2}}  + 3\,n + m - 2\,a)\,\Gamma (a)) }
}
\end{maplelatex}

\begin{maplelatex}
\mapleinline{inert}{2d}{368, ":   ",
GAMMA(1/2*c+1/2*b+n+m-a)*GAMMA(1+a-b)*GAMMA(1+a-c)*WPMR(2*m+2*n-1,2*n,
c,b,a)*GAMMA(1-2*a+2*n)/GAMMA(1-a+2*n)/GAMMA(1/2*c+1/2*b+n+m)/GAMMA(-c
+1)/GAMMA(-b+1)+GAMMA(1/2*c+1/2*b+n+m-a)*GAMMA(1+a-b)*GAMMA(1+a-c)*GAM
MA(2*a-2*n-1)*WPMR(2*m+2*n-1,2*n,1+b-2*a+2*n,1-2*a+c+2*n,1-a+2*n)/GAMM
A(-b+2*a-2*n)/GAMMA(2*a-2*n-c)/GAMMA(1+1/2*c+1/2*b+3*n+m-2*a)/GAMMA(a)
;}{%
\maplemultiline{
368\mbox{:~~~} \Gamma ({\displaystyle \frac {c}{2}}  + 
{\displaystyle \frac {b}{2}}  + n + m - a)\,\Gamma (1 + a - b)\,
\Gamma (1 + a - c) \\
\mathrm{WPMR}(2\,m + 2\,n - 1, \,2\,n, \,c, \,b, \,a)\,\Gamma (1
 - 2\,a + 2\,n) \left/ {\vrule height0.80em width0em depth0.80em}
 \right. \!  \! (\Gamma (1 - a + 2\,n) \\
\Gamma ({\displaystyle \frac {c}{2}}  + {\displaystyle \frac {b}{
2}}  + n + m)\,\Gamma ( - c + 1)\,\Gamma ( - b + 1))\mbox{} + 
\Gamma ({\displaystyle \frac {c}{2}}  + {\displaystyle \frac {b}{
2}}  + n + m - a)\,\Gamma (1 + a - b) \\
\Gamma (1 + a - c)\,\Gamma (2\,a - 2\,n - 1) \\
\mathrm{WPMR}(2\,m + 2\,n - 1, \,2\,n, \,1 + b - 2\,a + 2\,n, \,1
 - 2\,a + c + 2\,n, \,1 - a + 2\,n) \left/ {\vrule 
height0.80em width0em depth0.80em} \right. \!  \!  \\
(\Gamma ( - b + 2\,a - 2\,n)\,\Gamma (2\,a - 2\,n - c)\,\Gamma (1
 + {\displaystyle \frac {c}{2}}  + {\displaystyle \frac {b}{2}} 
 + 3\,n + m - 2\,a)\,\Gamma (a)) }
}
\end{maplelatex}

\begin{maplelatex}
\mapleinline{inert}{2d}{375, ":   ",
(GAMMA(-a+2*m)*WPMR(2*m+2*n-2,2*n-1,c+1-2*a-2*n+b,-b+c+2*m-2*a,c-a)*GA
MMA(2*n+2*a-2*c)/GAMMA(a-c+2*n)/GAMMA(c-2*a+2*m)/GAMMA(-c+2*a+2*n-b)/G
AMMA(b-c+2*a-2*m+1)*GAMMA(1+a-2*m)*GAMMA(-c+2*a+2*n)+GAMMA(-a+2*m)*GAM
MA(-2*a+2*c-2*n)*WPMR(2*m+2*n-2,2*n-1,1+b-c,-b-c+2*n+2*m,a-c+2*n)/GAMM
A(b+c-2*m-2*n+1)/GAMMA(c-b)/GAMMA(-c+2*n+2*m)*GAMMA(1+a-2*m)*GAMMA(-c+
2*a+2*n)/GAMMA(c-a))*GAMMA(c)/GAMMA(a);}{%
\maplemultiline{
375\mbox{:~~~} (\Gamma ( - a + 2\,m) \\
\mathrm{WPMR}(2\,m + 2\,n - 2, \,2\,n - 1, \,c + 1 - 2\,a - 2\,n
 + b, \, - b + c + 2\,m - 2\,a, \,c - a) \\
\Gamma (2\,n + 2\,a - 2\,c)\,\Gamma (1 + a - 2\,m)\,\Gamma ( - c
 + 2\,a + 2\,n)/(\Gamma (a - c + 2\,n) \\
\Gamma (c - 2\,a + 2\,m)\,\Gamma ( - c + 2\,a + 2\,n - b)\,\Gamma
 (b - c + 2\,a - 2\,m + 1))\mbox{} + \Gamma ( - a + 2\,m) \\
\Gamma ( - 2\,a + 2\,c - 2\,n) \\
\mathrm{WPMR}(2\,m + 2\,n - 2, \,2\,n - 1, \,1 + b - c, \, - b - 
c + 2\,n + 2\,m, \,a - c + 2\,n) \\
\Gamma (1 + a - 2\,m)\,\Gamma ( - c + 2\,a + 2\,n)/(\Gamma (b + c
 - 2\,m - 2\,n + 1)\,\Gamma (c - b) \\
\Gamma ( - c + 2\,n + 2\,m)\,\Gamma (c - a)))\Gamma (c)/\Gamma (a
) }
}
\end{maplelatex}

\begin{maplelatex}
\mapleinline{inert}{2d}{376, ":   ",
GAMMA(a+1-b+c-2*n)*GAMMA(-c+b+2*n)*WPMR(2*m+2*n-2,2*n-1,2*a+2-b+2*c-4*
n-2*m,1+2*c-b-2*n,c)/GAMMA(-a-1+b-c+2*n+2*m)*GAMMA(2*n-2*c)/GAMMA(-c+2
*n)/GAMMA(a+1-b+2*c-2*n)/GAMMA(b-2*c+2*n)*GAMMA(-a-1+b+2*n+2*m)+GAMMA(
a+1-b+c-2*n)*GAMMA(-c+b+2*n)*GAMMA(2*c-2*n)*WPMR(2*m+2*n-2,2*n-1,2*a+2
-b-2*n-2*m,-b+1,-c+2*n)/GAMMA(-a-1+b-c+2*n+2*m)/GAMMA(c)/GAMMA(2*n-1-2
*a+b+2*m)/GAMMA(1+a-b)*GAMMA(-2*a-1+b-2*c+4*n+2*m)*GAMMA(-a-1+b+2*n+2*
m)/GAMMA(b);}{%
\maplemultiline{
376\mbox{:~~~} \Gamma (a + 1 - b + c - 2\,n)\,\Gamma (
 - c + b + 2\,n) \\
\mathrm{WPMR}(2\,m + 2\,n - 2, \,2\,n - 1, \,2\,a + 2 - b + 2\,c
 - 4\,n - 2\,m, \,1 + 2\,c - b - 2\,n, \,c) \\
\Gamma (2\,n - 2\,c)\,\Gamma ( - a - 1 + b + 2\,n + 2\,m)/(\Gamma
 ( - a - 1 + b - c + 2\,n + 2\,m)\,\Gamma ( - c + 2\,n) \\
\Gamma (a + 1 - b + 2\,c - 2\,n)\,\Gamma (b - 2\,c + 2\,n))
\mbox{} + \Gamma (a + 1 - b + c - 2\,n)\,\Gamma ( - c + b + 2\,n)
 \\
\Gamma (2\,c - 2\,n) 
\mathrm{WPMR}(2\,m + 2\,n - 2, \,2\,n - 1, \,2\,a + 2 - b - 2\,n
 - 2\,m, \, - b + 1, \, - c + 2\,n) \\
\Gamma ( - 2\,a - 1 + b - 2\,c + 4\,n + 2\,m)\,\Gamma ( - a - 1
 + b + 2\,n + 2\,m)/( \\
\Gamma ( - a - 1 + b - c + 2\,n + 2\,m)\,\Gamma (c)\,\Gamma (2\,n
 - 1 - 2\,a + b + 2\,m)\,\Gamma (1 + a - b)\,\Gamma (b)) }
}
\end{maplelatex}

\begin{maplelatex}
\mapleinline{inert}{2d}{377, ":   ",
GAMMA(-c+1)*GAMMA(1/2*a+1/2*b+n)*WPMR(2*m+2*n-2,2*n-1,-2*c+2+a-2*m,-b+
1,n+1/2*a-1/2*b)*GAMMA(c-1+2*m-1/2*b-1/2*a+n)/GAMMA(2*m+c-1)*GAMMA(b-a
)/GAMMA(-1/2*a+1/2*b+n)/GAMMA(-c+1+1/2*a-1/2*b+n)/GAMMA(2*c-1-a+2*m)/G
AMMA(b)*GAMMA(2*c+2*m-1)+GAMMA(-c+1)*GAMMA(1/2*a+1/2*b+n)*GAMMA(a-b)*W
PMR(2*m+2*n-2,2*n-1,-2*c+2-2*m+b,-a+1,-1/2*a+1/2*b+n)*GAMMA(c-1+2*m-1/
2*b-1/2*a+n)/GAMMA(2*m+c-1)/GAMMA(n+1/2*a-1/2*b)/GAMMA(-b-1+2*c+2*m)/G
AMMA(1-c-1/2*a+1/2*b+n)*GAMMA(2*c+2*m-1)/GAMMA(a);}{%
\maplemultiline{
377\mbox{:~~~} \Gamma ( - c + 1)\,\Gamma (
{\displaystyle \frac {a}{2}}  + {\displaystyle \frac {b}{2}}  + n
) \,
\mathrm{WPMR}(2\,m + 2\,n - 2, \,2\,n - 1, \, - 2\,c + 2 + a - 2
\,m, \, - b + 1, \,n + {\displaystyle \frac {a}{2}}  - 
{\displaystyle \frac {b}{2}} ) \\
\Gamma (c - 1 + 2\,m - {\displaystyle \frac {b}{2}}  - 
{\displaystyle \frac {a}{2}}  + n)\,\Gamma (b - a)\,\Gamma (2\,c
 + 2\,m - 1) \left/ {\vrule height0.80em width0em depth0.80em}
 \right. \!  \! (\Gamma (2\,m + c - 1) \\
\Gamma ( - {\displaystyle \frac {a}{2}}  + {\displaystyle \frac {
b}{2}}  + n)\,\Gamma ( - c + 1 + {\displaystyle \frac {a}{2}}  - 
{\displaystyle \frac {b}{2}}  + n)\,\Gamma (2\,c - 1 - a + 2\,m)
\,\Gamma (b))\mbox{} + \Gamma ( - c + 1) \\
\Gamma ({\displaystyle \frac {a}{2}}  + {\displaystyle \frac {b}{
2}}  + n)\,\Gamma (a - b) \,
\mathrm{WPMR}(2\,m + 2\,n - 2, \,2\,n - 1, \, - 2\,c + 2 - 2\,m
 + b, \, - a + 1, \, - {\displaystyle \frac {a}{2}}  + 
{\displaystyle \frac {b}{2}}  + n) \\
\Gamma (c - 1 + 2\,m - {\displaystyle \frac {b}{2}}  - 
{\displaystyle \frac {a}{2}}  + n)\,\Gamma (2\,c + 2\,m - 1)
 \left/ {\vrule height0.80em width0em depth0.80em} \right. \! 
 \! (\Gamma (2\,m + c - 1)\,\Gamma (n + {\displaystyle \frac {a}{
2}}  - {\displaystyle \frac {b}{2}} ) \\
\Gamma ( - b - 1 + 2\,c + 2\,m)\,\Gamma (1 - c - {\displaystyle 
\frac {a}{2}}  + {\displaystyle \frac {b}{2}}  + n)\,\Gamma (a))
 }
}
\end{maplelatex}

\begin{maplelatex}
\mapleinline{inert}{2d}{378, ":   ",
(GAMMA(1-a+2*m)*WPMR(2*m+2*n-1,2*n-1,c+1-2*a-2*n+b,-b+c+2*m-2*a+1,c-a)
*GAMMA(2*n+2*a-2*c)/GAMMA(a-c+2*n)/GAMMA(c-2*a+1+2*m)/GAMMA(-c+2*a+2*n
-b)/GAMMA(b-c+2*a-2*m)*GAMMA(a-2*m)*GAMMA(-c+2*a+2*n)+GAMMA(1-a+2*m)*G
AMMA(-2*a+2*c-2*n)*WPMR(2*m+2*n-1,2*n-1,1+b-c,1-b-c+2*m+2*n,a-c+2*n)/G
AMMA(b+c-2*m-2*n)/GAMMA(c-b)/GAMMA(-c+1+2*n+2*m)*GAMMA(a-2*m)*GAMMA(-c
+2*a+2*n)/GAMMA(c-a))*GAMMA(c)/GAMMA(a);}{%
\maplemultiline{
378\mbox{:~~~} (\Gamma (1 - a + 2\,m) \,
\mathrm{WPMR}(2\,m + 2\,n - 1, \,2\,n - 1, \,c + 1 - 2\,a - 2\,n
 + b, \, - b + c + 2\,m - 2\,a + 1, \,c - a) \\
\Gamma (2\,n + 2\,a - 2\,c)\,\Gamma (a - 2\,m)\,\Gamma ( - c + 2
\,a + 2\,n)/(\Gamma (a - c + 2\,n) \\
\Gamma (c - 2\,a + 1 + 2\,m)\,\Gamma ( - c + 2\,a + 2\,n - b)\,
\Gamma (b - c + 2\,a - 2\,m))\mbox{} + \Gamma (1 - a + 2\,m) \\
\Gamma ( - 2\,a + 2\,c - 2\,n) 
\mathrm{WPMR}(2\,m + 2\,n - 1, \,2\,n - 1, \,1 + b - c, \,1 - b
 - c + 2\,m + 2\,n, \,a - c + 2\,n) \\
\Gamma (a - 2\,m)\,\Gamma ( - c + 2\,a + 2\,n)/(\Gamma (b + c - 2
\,m - 2\,n)\,\Gamma (c - b)\,\Gamma ( - c + 1 + 2\,n + 2\,m)
\Gamma (c - a)))\Gamma (c)/\Gamma (a) }
}
\end{maplelatex}

\begin{maplelatex}
\mapleinline{inert}{2d}{379, ":   ",
GAMMA(a+1-b+c-2*n)*GAMMA(-c+b+2*n)*WPMR(2*m+2*n-1,2*n-1,2*a-b+2*c-4*n+
1-2*m,1+2*c-b-2*n,c)/GAMMA(-a+b-c+2*n+2*m)*GAMMA(2*n-2*c)/GAMMA(-c+2*n
)/GAMMA(a+1-b+2*c-2*n)/GAMMA(b-2*c+2*n)*GAMMA(-a+b+2*n+2*m)+GAMMA(a+1-
b+c-2*n)*GAMMA(-c+b+2*n)*GAMMA(2*c-2*n)*WPMR(2*m+2*n-1,2*n-1,2*a-b-2*n
+1-2*m,-b+1,-c+2*n)/GAMMA(-a+b-c+2*n+2*m)/GAMMA(c)/GAMMA(2*n-2*a+b+2*m
)/GAMMA(1+a-b)*GAMMA(-a+b+2*n+2*m)*GAMMA(-2*a+b-2*c+4*n+2*m)/GAMMA(b);
}{%
\maplemultiline{
379\mbox{:~~~} \Gamma (a + 1 - b + c - 2\,n)\,\Gamma (
 - c + b + 2\,n) \\
\mathrm{WPMR}(2\,m + 2\,n - 1, \,2\,n - 1, \,2\,a - b + 2\,c - 4
\,n + 1 - 2\,m, \,1 + 2\,c - b - 2\,n, \,c) \\
\Gamma (2\,n - 2\,c)\,\Gamma ( - a + b + 2\,n + 2\,m)/(\Gamma (
 - a + b - c + 2\,n + 2\,m)\,\Gamma ( - c + 2\,n) \\
\Gamma (a + 1 - b + 2\,c - 2\,n)\,\Gamma (b - 2\,c + 2\,n))
\mbox{} + \Gamma (a + 1 - b + c - 2\,n)\,\Gamma ( - c + b + 2\,n)
 \\
\Gamma (2\,c - 2\,n) \,
\mathrm{WPMR}(2\,m + 2\,n - 1, \,2\,n - 1, \,2\,a - b - 2\,n + 1
 - 2\,m, \, - b + 1, \, - c + 2\,n) \\
\Gamma ( - a + b + 2\,n + 2\,m)\,\Gamma ( - 2\,a + b - 2\,c + 4\,
n + 2\,m)/(\Gamma ( - a + b - c + 2\,n + 2\,m) \\
\Gamma (c)\,\Gamma (2\,n - 2\,a + b + 2\,m)\,\Gamma (1 + a - b)\,
\Gamma (b)) }
}
\end{maplelatex}

\begin{maplelatex}
\mapleinline{inert}{2d}{380, ":   ",
GAMMA(-c+1)*GAMMA(1/2*a+1/2*b+n)*WPMR(2*m+2*n-1,2*n-1,-2*c+1+a-2*m,-b+
1,n+1/2*a-1/2*b)*GAMMA(-1/2*a-1/2*b+n+c+2*m)/GAMMA(2*m+c)*GAMMA(b-a)/G
AMMA(-1/2*a+1/2*b+n)/GAMMA(-c+1+1/2*a-1/2*b+n)/GAMMA(2*c-a+2*m)/GAMMA(
b)*GAMMA(2*c+2*m)+GAMMA(-c+1)*GAMMA(1/2*a+1/2*b+n)*GAMMA(a-b)*WPMR(2*m
+2*n-1,2*n-1,-2*c+1-2*m+b,-a+1,-1/2*a+1/2*b+n)*GAMMA(-1/2*a-1/2*b+n+c+
2*m)/GAMMA(2*m+c)/GAMMA(n+1/2*a-1/2*b)/GAMMA(-b+2*c+2*m)/GAMMA(1-c-1/2
*a+1/2*b+n)*GAMMA(2*c+2*m)/GAMMA(a);}{%
\maplemultiline{
380\mbox{:~~~} \Gamma ( - c + 1)\,\Gamma (
{\displaystyle \frac {a}{2}}  + {\displaystyle \frac {b}{2}}  + n
) \,
\mathrm{WPMR}(2\,m + 2\,n - 1, \,2\,n - 1, \, - 2\,c + 1 + a - 2
\,m, \, - b + 1, \,n + {\displaystyle \frac {a}{2}}  - 
{\displaystyle \frac {b}{2}} ) \\
\Gamma ( - {\displaystyle \frac {a}{2}}  - {\displaystyle \frac {
b}{2}}  + n + c + 2\,m)\,\Gamma (b - a)\,\Gamma (2\,c + 2\,m)
 \left/ {\vrule height0.80em width0em depth0.80em} \right. \! 
 \! (\Gamma (2\,m + c)\,\Gamma ( - {\displaystyle \frac {a}{2}} 
 + {\displaystyle \frac {b}{2}}  + n) \\
\Gamma ( - c + 1 + {\displaystyle \frac {a}{2}}  - 
{\displaystyle \frac {b}{2}}  + n)\,\Gamma (2\,c - a + 2\,m)\,
\Gamma (b))\mbox{} + \Gamma ( - c + 1)\,\Gamma ({\displaystyle 
\frac {a}{2}}  + {\displaystyle \frac {b}{2}}  + n)\,\Gamma (a - 
b) \\
\mathrm{WPMR}(2\,m + 2\,n - 1, \,2\,n - 1, \, - 2\,c + 1 - 2\,m
 + b, \, - a + 1, \, - {\displaystyle \frac {a}{2}}  + 
{\displaystyle \frac {b}{2}}  + n) \\
\Gamma ( - {\displaystyle \frac {a}{2}}  - {\displaystyle \frac {
b}{2}}  + n + c + 2\,m)\,\Gamma (2\,c + 2\,m) \left/ {\vrule 
height0.80em width0em depth0.80em} \right. \!  \! (\Gamma (2\,m
 + c)\,\Gamma (n + {\displaystyle \frac {a}{2}}  - 
{\displaystyle \frac {b}{2}} ) \\
\Gamma ( - b + 2\,c + 2\,m)\,\Gamma (1 - c - {\displaystyle 
\frac {a}{2}}  + {\displaystyle \frac {b}{2}}  + n)\,\Gamma (a))
 }
}
\end{maplelatex}

\end{maplegroup}

%% file: AppendixB381to410.tex
\begin{maplegroup}
\mapleresult
\begin{maplelatex}
\mapleinline{inert}{2d}{381, ":   ",
GAMMA(-a+2*m)*GAMMA(-2*a+2*c-2*n-1)/GAMMA(a)/GAMMA(c-a)/GAMMA(b+c-2*m-
2*n)/GAMMA(c-b)/GAMMA(-c+1+2*n+2*m)*GAMMA(1+a-2*m)*GAMMA(-c+2*a+2*n+1)
*GAMMA(c)*WPMR(2*m+2*n-1,2*n,1-b-c+2*m+2*n,1+b-c,1+a-c+2*n)+GAMMA(-a+2
*m)*WPMR(2*m+2*n-1,2*n,c-2*a-2*n+b,-b+c+2*m-2*a,c-a)/GAMMA(a)*GAMMA(1+
2*a-2*c+2*n)/GAMMA(1+a-c+2*n)/GAMMA(c-2*a+2*m)/GAMMA(-c+2*a+2*n-b+1)/G
AMMA(b-c+2*a-2*m+1)*GAMMA(1+a-2*m)*GAMMA(-c+2*a+2*n+1)*GAMMA(c);}{%
\maplemultiline{
381\mbox{:~~~} \Gamma ( - a + 2\,m)\,\Gamma ( - 2\,a + 2
\,c - 2\,n - 1)\,\Gamma (1 + a - 2\,m)\,\Gamma ( - c + 2\,a + 2\,
n + 1)\,\Gamma (c) \\
\mathrm{WPMR}(2\,m + 2\,n - 1, \,2\,n, \,1 - b - c + 2\,m + 2\,n
, \,1 + b - c, \,1 + a - c + 2\,n)/(\Gamma (a) \\
\Gamma (c - a)\,\Gamma (b + c - 2\,m - 2\,n)\,\Gamma (c - b)\,
\Gamma ( - c + 1 + 2\,n + 2\,m))\mbox{} + \Gamma ( - a + 2\,m)
 \\
\mathrm{WPMR}(2\,m + 2\,n - 1, \,2\,n, \,c - 2\,a - 2\,n + b, \,
 - b + c + 2\,m - 2\,a, \,c - a) \\
\Gamma (1 + 2\,a - 2\,c + 2\,n)\,\Gamma (1 + a - 2\,m)\,\Gamma (
 - c + 2\,a + 2\,n + 1)\,\Gamma (c)/(\Gamma (a) \\
\Gamma (1 + a - c + 2\,n)\,\Gamma (c - 2\,a + 2\,m)\,\Gamma ( - c
 + 2\,a + 2\,n - b + 1) \\
\Gamma (b - c + 2\,a - 2\,m + 1)) }
}
\end{maplelatex}

\begin{maplelatex}
\mapleinline{inert}{2d}{382, ":   ",
GAMMA(a-b+c-2*n)*GAMMA(1-c+b+2*n)*GAMMA(2*c-2*n-1)/GAMMA(-a+b-c+2*n+2*
m)/GAMMA(c)/GAMMA(b)/GAMMA(2*n-2*a+b+2*m)/GAMMA(1+a-b)*GAMMA(-a+b+2*n+
2*m)*GAMMA(-2*a+b-2*c+4*n+2*m+1)*WPMR(2*m+2*n-1,2*n,-b+1,2*a-b-2*n+1-2
*m,1-c+2*n)+GAMMA(a-b+c-2*n)*GAMMA(1-c+b+2*n)*WPMR(2*m+2*n-1,2*n,2*a-b
+2*c-4*n-2*m,2*c-b-2*n,c)/GAMMA(-a+b-c+2*n+2*m)*GAMMA(1-2*c+2*n)/GAMMA
(1-c+2*n)/GAMMA(a-b+2*c-2*n)/GAMMA(b-2*c+2*n+1)*GAMMA(-a+b+2*n+2*m);}{
\maplemultiline{
382\mbox{:~~~} \Gamma (a - b + c - 2\,n)\,\Gamma (1 - c
 + b + 2\,n)\,\Gamma (2\,c - 2\,n - 1)\,\Gamma ( - a + b + 2\,n
 + 2\,m) \\
\Gamma ( - 2\,a + b - 2\,c + 4\,n + 2\,m + 1) \\
\mathrm{WPMR}(2\,m + 2\,n - 1, \,2\,n, \, - b + 1, \,2\,a - b - 2
\,n + 1 - 2\,m, \,1 - c + 2\,n)/( \\
\Gamma ( - a + b - c + 2\,n + 2\,m)\,\Gamma (c)\,\Gamma (b)\,
\Gamma (2\,n - 2\,a + b + 2\,m)\,\Gamma (1 + a - b))\mbox{} + 
 \\
\Gamma (a - b + c - 2\,n)\,\Gamma (1 - c + b + 2\,n) \\
\mathrm{WPMR}(2\,m + 2\,n - 1, \,2\,n, \,2\,a - b + 2\,c - 4\,n
 - 2\,m, \,2\,c - b - 2\,n, \,c) \\
\Gamma (1 - 2\,c + 2\,n)\,\Gamma ( - a + b + 2\,n + 2\,m)/(\Gamma
 ( - a + b - c + 2\,n + 2\,m)\,\Gamma (1 - c + 2\,n) \\
\Gamma (a - b + 2\,c - 2\,n)\,\Gamma (b - 2\,c + 2\,n + 1)) }
}
\end{maplelatex}

\begin{maplelatex}
\mapleinline{inert}{2d}{383, ":   ",
GAMMA(-c+1)*GAMMA(1/2+1/2*a+n+1/2*b)*GAMMA(a-b)/GAMMA(2*m+c-1)*GAMMA(c
-1/2-1/2*a+2*m-1/2*b+n)/GAMMA(n+1/2*a+1/2-1/2*b)/GAMMA(a)/GAMMA(-b-1+2
*c+2*m)/GAMMA(3/2-c-1/2*a+1/2*b+n)*GAMMA(2*c+2*m-1)*WPMR(2*m+2*n-1,2*n
,-a+1,-2*c+2-2*m+b,1/2-1/2*a+n+1/2*b)+GAMMA(-c+1)*GAMMA(1/2+1/2*a+n+1/
2*b)*WPMR(2*m+2*n-1,2*n,-2*c+2+a-2*m,-b+1,n+1/2*a+1/2-1/2*b)*GAMMA(c-1
/2-1/2*a+2*m-1/2*b+n)/GAMMA(2*m+c-1)*GAMMA(b-a)/GAMMA(1/2-1/2*a+n+1/2*
b)/GAMMA(-c+3/2+1/2*a-1/2*b+n)/GAMMA(2*c-1-a+2*m)/GAMMA(b)*GAMMA(2*c+2
*m-1);}{%
\maplemultiline{
383\mbox{:~~~} \Gamma ( - c + 1)\,\Gamma (
{\displaystyle \frac {1}{2}}  + {\displaystyle \frac {a}{2}}  + n
 + {\displaystyle \frac {b}{2}} )\,\Gamma (a - b)\,\Gamma (c - 
{\displaystyle \frac {1}{2}}  - {\displaystyle \frac {a}{2}}  + 2
\,m - {\displaystyle \frac {b}{2}}  + n) \\
\Gamma (2\,c + 2\,m - 1) \\
\mathrm{WPMR}(2\,m + 2\,n - 1, \,2\,n, \, - a + 1, \, - 2\,c + 2
 - 2\,m + b, \,{\displaystyle \frac {1}{2}}  - {\displaystyle 
\frac {a}{2}}  + n + {\displaystyle \frac {b}{2}} ) \left/ 
{\vrule height0.80em width0em depth0.80em} \right. \!  \! ( \\
\Gamma (2\,m + c - 1)\,\Gamma (n + {\displaystyle \frac {a}{2}} 
 + {\displaystyle \frac {1}{2}}  - {\displaystyle \frac {b}{2}} )
\,\Gamma (a)\,\Gamma ( - b - 1 + 2\,c + 2\,m)\,\Gamma (
{\displaystyle \frac {3}{2}}  - c - {\displaystyle \frac {a}{2}} 
 + {\displaystyle \frac {b}{2}}  + n)) \\
\mbox{} + \Gamma ( - c + 1)\,\Gamma ({\displaystyle \frac {1}{2}
}  + {\displaystyle \frac {a}{2}}  + n + {\displaystyle \frac {b
}{2}} ) \\
\mathrm{WPMR}(2\,m + 2\,n - 1, \,2\,n, \, - 2\,c + 2 + a - 2\,m, 
\, - b + 1, \,n + {\displaystyle \frac {a}{2}}  + {\displaystyle 
\frac {1}{2}}  - {\displaystyle \frac {b}{2}} ) \\
\Gamma (c - {\displaystyle \frac {1}{2}}  - {\displaystyle 
\frac {a}{2}}  + 2\,m - {\displaystyle \frac {b}{2}}  + n)\,
\Gamma (b - a)\,\Gamma (2\,c + 2\,m - 1) \left/ {\vrule 
height0.80em width0em depth0.80em} \right. \!  \! (\Gamma (2\,m
 + c - 1) \\
\Gamma ({\displaystyle \frac {1}{2}}  - {\displaystyle \frac {a}{
2}}  + n + {\displaystyle \frac {b}{2}} )\,\Gamma ( - c + 
{\displaystyle \frac {3}{2}}  + {\displaystyle \frac {a}{2}}  - 
{\displaystyle \frac {b}{2}}  + n)\,\Gamma (2\,c - 1 - a + 2\,m)
\,\Gamma (b)) }
}
\end{maplelatex}

\begin{maplelatex}
\mapleinline{inert}{2d}{384, ":   ",
GAMMA(-1/2-1/2*b+1/2*c+n+m)*GAMMA(1+b-c)*Pi/sin(Pi*(a-b))/GAMMA(a-b)*G
AMMA(2*a-2*n)/GAMMA(b)/GAMMA(a-2*n+1)/GAMMA(2*a-2*n+1-c)/GAMMA(-1/2+1/
2*c+1/2*b+3*n+m-2*a)*WPMR(2*m+2*n-2,2*n-1,c-2*a+2*n,-2*a+2*n+b,-a+2*n)
+(cos(Pi*(-3*a+b+2*n))-cos(Pi*(2*n-a+b)))*WPMR(2*m+2*n-2,2*n-1,c,b,a)*
GAMMA(-b+2*a-2*n+1)*GAMMA(-1/2-1/2*b+1/2*c+n+m)*GAMMA(1+b-c)*GAMMA(a)/
GAMMA(1/2*c+1/2*b+n+m-1/2)/GAMMA(2*a-2*n+1)/GAMMA(-c+1)/GAMMA(a-b)/(co
s(Pi*(-3*a+b+2*n))-cos(Pi*(-a-b+2*n)));}{%
\maplemultiline{
384\mbox{:~~~} \Gamma ( - {\displaystyle \frac {1}{2}} 
 - {\displaystyle \frac {b}{2}}  + {\displaystyle \frac {c}{2}} 
 + n + m)\,\Gamma (1 + b - c)\,\pi \,\Gamma (2\,a - 2\,n) \\
\mathrm{WPMR}(2\,m + 2\,n - 2, \,2\,n - 1, \,c - 2\,a + 2\,n, \,
 - 2\,a + 2\,n + b, \, - a + 2\,n) \left/ {\vrule 
height0.80em width0em depth0.80em} \right. \!  \! ( \\
\mathrm{sin}(\pi \,(a - b))\,\Gamma (a - b)\,\Gamma (b)\,\Gamma (
a - 2\,n + 1)\,\Gamma (2\,a - 2\,n + 1 - c) \\
\Gamma ( - {\displaystyle \frac {1}{2}}  + {\displaystyle \frac {
c}{2}}  + {\displaystyle \frac {b}{2}}  + 3\,n + m - 2\,a))
\mbox{} + (\mathrm{cos}(\pi \,( - 3\,a + b + 2\,n)) - \mathrm{cos
}(\pi \,(2\,n - a + b))) \\
\mathrm{WPMR}(2\,m + 2\,n - 2, \,2\,n - 1, \,c, \,b, \,a)\,\Gamma
 ( - b + 2\,a - 2\,n + 1) \\
\Gamma ( - {\displaystyle \frac {1}{2}}  - {\displaystyle \frac {
b}{2}}  + {\displaystyle \frac {c}{2}}  + n + m)\,\Gamma (1 + b
 - c)\,\Gamma (a) \left/ {\vrule 
height0.80em width0em depth0.80em} \right. \!  \! (\Gamma (
{\displaystyle \frac {c}{2}}  + {\displaystyle \frac {b}{2}}  + n
 + m - {\displaystyle \frac {1}{2}} ) \\
\Gamma (2\,a - 2\,n + 1)\,\Gamma ( - c + 1)\,\Gamma (a - b) \\
(\mathrm{cos}(\pi \,( - 3\,a + b + 2\,n)) - \mathrm{cos}(\pi \,(
 - a - b + 2\,n)))) }
}
\end{maplelatex}

\begin{maplelatex}
\mapleinline{inert}{2d}{385, ":   ",
GAMMA(-1/2*b+n+1/2*c+m)*GAMMA(1+b-c)*Pi/sin(Pi*(a-b))/GAMMA(a-b)*GAMMA
(2*a-2*n)/GAMMA(b)/GAMMA(a-2*n+1)/GAMMA(2*a-2*n+1-c)/GAMMA(1/2*c+1/2*b
+3*n+m-2*a)*WPMR(2*m+2*n-1,2*n-1,c-2*a+2*n,-2*a+2*n+b,-a+2*n)+(cos(Pi*
(-3*a+b+2*n))-cos(Pi*(2*n-a+b)))*WPMR(2*m+2*n-1,2*n-1,c,b,a)*GAMMA(-b+
2*a-2*n+1)*GAMMA(1+b-c)*GAMMA(-1/2*b+n+1/2*c+m)*GAMMA(a)/GAMMA(1/2*c+1
/2*b+n+m)/GAMMA(2*a-2*n+1)/GAMMA(-c+1)/GAMMA(a-b)/(cos(Pi*(-3*a+b+2*n)
)-cos(Pi*(-a-b+2*n)));}{%
\maplemultiline{
385\mbox{:~~~} \Gamma ( - {\displaystyle \frac {b}{2}} 
 + n + {\displaystyle \frac {c}{2}}  + m)\,\Gamma (1 + b - c)\,
\pi \,\Gamma (2\,a - 2\,n) \\
\mathrm{WPMR}(2\,m + 2\,n - 1, \,2\,n - 1, \,c - 2\,a + 2\,n, \,
 - 2\,a + 2\,n + b, \, - a + 2\,n) \left/ {\vrule 
height0.80em width0em depth0.80em} \right. \!  \! ( \\
\mathrm{sin}(\pi \,(a - b))\,\Gamma (a - b)\,\Gamma (b)\,\Gamma (
a - 2\,n + 1)\,\Gamma (2\,a - 2\,n + 1 - c) \\
\Gamma ({\displaystyle \frac {c}{2}}  + {\displaystyle \frac {b}{
2}}  + 3\,n + m - 2\,a))\mbox{} + (\mathrm{cos}(\pi \,( - 3\,a + 
b + 2\,n)) - \mathrm{cos}(\pi \,(2\,n - a + b))) \\
\mathrm{WPMR}(2\,m + 2\,n - 1, \,2\,n - 1, \,c, \,b, \,a)\,\Gamma
 ( - b + 2\,a - 2\,n + 1)\,\Gamma (1 + b - c) \\
\Gamma ( - {\displaystyle \frac {b}{2}}  + n + {\displaystyle 
\frac {c}{2}}  + m)\,\Gamma (a) \left/ {\vrule 
height0.80em width0em depth0.80em} \right. \!  \! (\Gamma (
{\displaystyle \frac {c}{2}}  + {\displaystyle \frac {b}{2}}  + n
 + m)\,\Gamma (2\,a - 2\,n + 1)\,\Gamma ( - c + 1) \\
\Gamma (a - b)\,(\mathrm{cos}(\pi \,( - 3\,a + b + 2\,n)) - 
\mathrm{cos}(\pi \,( - a - b + 2\,n)))) }
}
\end{maplelatex}

\begin{maplelatex}
\mapleinline{inert}{2d}{386, ":   ",
GAMMA(-1/2*b+n+1/2*c+m)*GAMMA(1+b-c)*Pi/sin(Pi*(a-b))/GAMMA(a-b)*GAMMA
(2*a-2*n-1)/GAMMA(b)/GAMMA(a-2*n)/GAMMA(2*a-2*n-c)/GAMMA(1+1/2*c+1/2*b
+3*n+m-2*a)*WPMR(2*m+2*n-1,2*n,1+b-2*a+2*n,1-2*a+c+2*n,1-a+2*n)+(cos(P
i*(-3*a+b+2*n))-cos(Pi*(2*n-a+b)))*WPMR(2*m+2*n-1,2*n,c,b,a)*GAMMA(-b+
2*a-2*n)*GAMMA(-1/2*b+n+1/2*c+m)*GAMMA(1+b-c)*GAMMA(a)/GAMMA(2*a-2*n)/
GAMMA(1/2*c+1/2*b+n+m)/GAMMA(-c+1)/GAMMA(a-b)/(cos(Pi*(-3*a+b+2*n))-co
s(Pi*(-a-b+2*n)));}{%
\maplemultiline{
386\mbox{:~~~} \Gamma ( - {\displaystyle \frac {b}{2}} 
 + n + {\displaystyle \frac {c}{2}}  + m)\,\Gamma (1 + b - c)\,
\pi \,\Gamma (2\,a - 2\,n - 1) \\
\mathrm{WPMR}(2\,m + 2\,n - 1, \,2\,n, \,1 + b - 2\,a + 2\,n, \,1
 - 2\,a + c + 2\,n, \,1 - a + 2\,n) \left/ {\vrule 
height0.80em width0em depth0.80em} \right. \!  \!  \\
(\mathrm{sin}(\pi \,(a - b))\,\Gamma (a - b)\,\Gamma (b)\,\Gamma 
(a - 2\,n)\,\Gamma (2\,a - 2\,n - c) \\
\Gamma (1 + {\displaystyle \frac {c}{2}}  + {\displaystyle 
\frac {b}{2}}  + 3\,n + m - 2\,a))\mbox{} + (\mathrm{cos}(\pi \,(
 - 3\,a + b + 2\,n)) - \mathrm{cos}(\pi \,(2\,n - a + b))) \\
\mathrm{WPMR}(2\,m + 2\,n - 1, \,2\,n, \,c, \,b, \,a)\,\Gamma (
 - b + 2\,a - 2\,n)\,\Gamma ( - {\displaystyle \frac {b}{2}}  + n
 + {\displaystyle \frac {c}{2}}  + m) \\
\Gamma (1 + b - c)\,\Gamma (a) \left/ {\vrule 
height0.80em width0em depth0.80em} \right. \!  \! (\Gamma (2\,a
 - 2\,n)\,\Gamma ({\displaystyle \frac {c}{2}}  + {\displaystyle 
\frac {b}{2}}  + n + m)\,\Gamma ( - c + 1)\,\Gamma (a - b) \\
(\mathrm{cos}(\pi \,( - 3\,a + b + 2\,n)) - \mathrm{cos}(\pi \,(
 - a - b + 2\,n)))) }
}
\end{maplelatex}

\begin{maplelatex}
\mapleinline{inert}{2d}{387, ":   ",
-1/2*WPMR(2*m+2*n-2,2*n-1,2*a+2-b-2*n-2*m,-b+1,-c+2*n)*GAMMA(b-2*c+2*n
)*GAMMA(-2*a-1+b-2*c+4*n+2*m)*GAMMA(-c+2*n)*GAMMA(2*c-2*n)*GAMMA(-a+2*
m+c)*GAMMA(a)*(cos(Pi*(b-c))-cos(Pi*(c+b)))/sin(Pi*(-c+b+2*n))/Pi/GAMM
A(-c+b+2*n-1)/GAMMA(2*n-1-2*a+b+2*m)/GAMMA(-a-1+b-c+2*n+2*m)/GAMMA(1+a
-b)-WPMR(2*m+2*n-2,2*n-1,2*a+2-b+2*c-4*n-2*m,1+2*c-b-2*n,c)*Pi/sin(Pi*
(-c+b+2*n))/GAMMA(-c+b+2*n-1)/GAMMA(a+1-b+2*c-2*n)*GAMMA(2*n-2*c)/GAMM
A(-a-1+b-c+2*n+2*m)*GAMMA(-a+2*m+c)/GAMMA(-c+1)/GAMMA(-b+1)*GAMMA(a)-W
PMR(2*m+2*n-2,2*m-1,b,2*n-1-2*a+b+2*m,-a-1+b-c+2*n+2*m)*sin(Pi*(2*n-2*
a+b+2*m))/Pi/GAMMA(-a-1+b+2*n+2*m)/GAMMA(a+1-b+c-2*n)*GAMMA(2*a+2-2*b+
2*c-4*n-2*m)/GAMMA(-c+b+2*n-1)*GAMMA(b-2*c+2*n)*GAMMA(-2*a-1+b-2*c+4*n
+2*m)*GAMMA(-a+2*m+c)*GAMMA(c)/GAMMA(-b+1)*GAMMA(a)+GAMMA(a)*GAMMA(-a+
2*m+c)*GAMMA(c)*WPMR(2*m+2*n-2,2*m-1,2*a+2-b+2*c-4*n-2*m,1+2*c-b-2*n,a
+1-b+c-2*n)*GAMMA(2*m-2*a-2+2*b-2*c+4*n)/GAMMA(-c+b+2*n-1)/GAMMA(2*n-1
-2*a+b+2*m)/GAMMA(-a-1+b-c+2*n+2*m)/GAMMA(a+1-b+2*c-2*n);}{%
\maplemultiline{
387\mbox{:~~~}  - {\displaystyle \frac {1}{2}} \mathrm{
WPMR}(2\,m + 2\,n - 2, \,2\,n - 1, \,2\,a + 2 - b - 2\,n - 2\,m, 
\, - b + 1, \, - c + 2\,n) \\
\Gamma (b - 2\,c + 2\,n)\,\Gamma ( - 2\,a - 1 + b - 2\,c + 4\,n
 + 2\,m)\,\Gamma ( - c + 2\,n)\,\Gamma (2\,c - 2\,n) \\
\Gamma ( - a + 2\,m + c)\,\Gamma (a)\,(\mathrm{cos}(\pi \,(b - c)
) - \mathrm{cos}(\pi \,(c + b)))/(\mathrm{sin}(\pi \,( - c + b + 
2\,n))\,\pi \,\mathrm{\%2} \\
\Gamma (2\,n - 1 - 2\,a + b + 2\,m)\,\Gamma (\mathrm{\%1})\,
\Gamma (1 + a - b))\mbox{} -  \\
\mathrm{WPMR}(2\,m + 2\,n - 2, \,2\,n - 1, \,2\,a + 2 - b + 2\,c
 - 4\,n - 2\,m, \,1 + 2\,c - b - 2\,n, \,c) \\
\pi \,\Gamma (2\,n - 2\,c)\,\Gamma ( - a + 2\,m + c)\,\Gamma (a)/
(\mathrm{sin}(\pi \,( - c + b + 2\,n))\,\mathrm{\%2} \\
\Gamma (a + 1 - b + 2\,c - 2\,n)\,\Gamma (\mathrm{\%1})\,\Gamma (
 - c + 1)\,\Gamma ( - b + 1))\mbox{} -  \\
\mathrm{WPMR}(2\,m + 2\,n - 2, \,2\,m - 1, \,b, \,2\,n - 1 - 2\,a
 + b + 2\,m, \,\mathrm{\%1}) \\
\mathrm{sin}(\pi \,(2\,n - 2\,a + b + 2\,m))\,\Gamma (2\,a + 2 - 
2\,b + 2\,c - 4\,n - 2\,m)\,\Gamma (b - 2\,c + 2\,n) \\
\Gamma ( - 2\,a - 1 + b - 2\,c + 4\,n + 2\,m)\,\Gamma ( - a + 2\,
m + c)\,\Gamma (c)\,\Gamma (a)/(\pi  \\
\Gamma ( - a - 1 + b + 2\,n + 2\,m)\,\Gamma (a + 1 - b + c - 2\,n
)\,\mathrm{\%2}\,\Gamma ( - b + 1))\mbox{} + \Gamma (a) \\
\Gamma ( - a + 2\,m + c)\,\Gamma (c)\mathrm{WPMR}(2\,m + 2\,n - 2
, \,2\,m - 1, \,2\,a + 2 - b + 2\,c - 4\,n - 2\,m,  \\
1 + 2\,c - b - 2\,n, \,a + 1 - b + c - 2\,n)\Gamma (2\,m - 2\,a
 - 2 + 2\,b - 2\,c + 4\,n)/(\mathrm{\%2} \\
\Gamma (2\,n - 1 - 2\,a + b + 2\,m)\,\Gamma (\mathrm{\%1})\,
\Gamma (a + 1 - b + 2\,c - 2\,n)) \\
\mathrm{\%1} :=  - a - 1 + b - c + 2\,n + 2\,m \\
\mathrm{\%2} := \Gamma ( - c + b + 2\,n - 1) }
}
\end{maplelatex}

\begin{maplelatex}
\mapleinline{inert}{2d}{388, ":   ",
-WPMR(2*m+2*n-2,2*n-1,-2*c+2-2*m+b,-a+1,-1/2*a+1/2*b+n)*Pi/sin(1/2*Pi*
(a+b+2*n))/GAMMA(n-1+1/2*a+1/2*b)/GAMMA(-b-1+2*c+2*m)/GAMMA(2*m+c-1)*G
AMMA(a-b)/GAMMA(-b+1)*GAMMA(1/2*a+1/2*b+n-c)*GAMMA(2*c+2*m-1/2*a-1/2*b
-n)/GAMMA(1/2*a-1/2*b-n+1)*GAMMA(c-1+2*m-1/2*b-1/2*a+n)/GAMMA(1-c-1/2*
a+1/2*b+n)-1/2*WPMR(2*m+2*n-2,2*n-1,-2*c+2+a-2*m,-b+1,n+1/2*a-1/2*b)*G
AMMA(n+1/2*a-1/2*b)*GAMMA(2*c+2*m-1/2*a-1/2*b-n)*GAMMA(c-1+2*m-1/2*b-1
/2*a+n)*GAMMA(1/2*a+1/2*b+n-c)*GAMMA(b-a)*GAMMA(a)*(cos(1/2*Pi*(-a-b+2
*n))-cos(1/2*Pi*(-a+3*b+2*n)))/Pi/sin(1/2*Pi*(a+b+2*n))/GAMMA(n-1+1/2*
a+1/2*b)/GAMMA(2*c-1-a+2*m)/GAMMA(-c+1+1/2*a-1/2*b+n)/GAMMA(2*m+c-1)+G
AMMA(1/2*a+1/2*b+n-c)*GAMMA(2*c+2*m-1/2*a-1/2*b-n)*GAMMA(-1/2*a+1/2*b+
n)*GAMMA(2*c-2+2*m)/GAMMA(n-1+1/2*a+1/2*b)/GAMMA(2*c-1-a+2*m)/GAMMA(2*
m+c-1)/GAMMA(-b-1+2*c+2*m)/GAMMA(1-c-1/2*a+1/2*b+n)*GAMMA(c-1+2*m-1/2*
b-1/2*a+n)*WPMR(2*m+2*n-2,2*m-1,-2*c+2-2*m+b,-a+1,-c+1)+sin(Pi*(-2*c+a
-2*m))*WPMR(2*m+2*n-2,2*m-1,b,2*c-1-a+2*m,2*m+c-1)/Pi*GAMMA(-2*c+2-2*m
)/GAMMA(-1/2*a+1/2*b+n-1+c+2*m)*GAMMA(-1/2*a+1/2*b+n)*GAMMA(c-1+2*m-1/
2*b-1/2*a+n)/GAMMA(n-1+1/2*a+1/2*b)*GAMMA(2*c+2*m-1/2*a-1/2*b-n)*GAMMA
(1/2*a+1/2*b+n-c)/GAMMA(-c+1)/GAMMA(-b+1)*GAMMA(a);}{%
\maplemultiline{
388\mbox{:~~~}  - \mathrm{WPMR}(2\,m + 2\,n - 2, \,2\,n
 - 1, \, - 2\,c + 2 - 2\,m + b, \, - a + 1, \, - {\displaystyle 
\frac {a}{2}}  + {\displaystyle \frac {b}{2}}  + n)\,\pi  \\
\Gamma (a - b)\,\mathrm{\%1}\,\mathrm{\%2}\,\mathrm{\%4} \left/ 
{\vrule height0.80em width0em depth0.80em} \right. \!  \! (
\mathrm{sin}({\displaystyle \frac {\pi \,(a + b + 2\,n)}{2}} )\,
\mathrm{\%3}\,\Gamma ( - b - 1 + 2\,c + 2\,m) \\
\Gamma (2\,m + c - 1)\,\Gamma ( - b + 1)\,\Gamma ({\displaystyle 
\frac {a}{2}}  - {\displaystyle \frac {b}{2}}  - n + 1)\,\Gamma (
1 - c - {\displaystyle \frac {a}{2}}  + {\displaystyle \frac {b}{
2}}  + n))\mbox{} - {\displaystyle \frac {1}{2}}  \\
\mathrm{WPMR}(2\,m + 2\,n - 2, \,2\,n - 1, \, - 2\,c + 2 + a - 2
\,m, \, - b + 1, \,n + {\displaystyle \frac {a}{2}}  - 
{\displaystyle \frac {b}{2}} ) \\
\Gamma (n + {\displaystyle \frac {a}{2}}  - {\displaystyle 
\frac {b}{2}} )\,\mathrm{\%2}\,\mathrm{\%4}\,\mathrm{\%1}\,\Gamma
 (b - a)\,\Gamma (a) \,
(\mathrm{cos}({\displaystyle \frac {\pi \,( - a - b + 2\,n)}{2}} 
) - \mathrm{cos}({\displaystyle \frac {\pi \,( - a + 3\,b + 2\,n)
}{2}} )) \\ \left/ {\vrule height0.80em width0em depth0.80em}
 \right. \!  \! (\pi \,\mathrm{sin}({\displaystyle \frac {\pi \,(
a + b + 2\,n)}{2}} ) 
\mathrm{\%3}\,\Gamma (2\,c - 1 - a + 2\,m)\,\Gamma ( - c + 1 + 
{\displaystyle \frac {a}{2}}  - {\displaystyle \frac {b}{2}}  + n
)\,\Gamma (2\,m + c - 1))\mbox{} \\ + \mathrm{\%1}\,\mathrm{\%2} 
\Gamma ( - {\displaystyle \frac {a}{2}}  + {\displaystyle \frac {
b}{2}}  + n)\,\Gamma (2\,c - 2 + 2\,m)\,\mathrm{\%4} 
\mathrm{WPMR}(2\,m + 2\,n - 2, \,2\,m - 1, \, - 2\,c + 2 - 2\,m
 + b, \, - a + 1, \, - c + 1) \left/ {\vrule 
height0.80em width0em depth0.80em} \right. \!  \! \\ (\mathrm{\%3}
\,
\Gamma (2\,c - 1 - a + 2\,m)\,\Gamma (2\,m + c - 1)\,\Gamma ( - b
 - 1 + 2\,c + 2\,m)\,\Gamma (1 - c - {\displaystyle \frac {a}{2}
}  + {\displaystyle \frac {b}{2}}  + n)) \\
\mbox{} + \mathrm{sin}(\pi \,( - 2\,c + a - 2\,m)) \,
\mathrm{WPMR}(2\,m + 2\,n - 2, \,2\,m - 1, \,b, \,2\,c - 1 - a + 
2\,m, \,2\,m + c - 1) \\
\Gamma ( - 2\,c + 2 - 2\,m)\,\Gamma ( - {\displaystyle \frac {a}{
2}}  + {\displaystyle \frac {b}{2}}  + n)\,\mathrm{\%4}\,\mathrm{
\%2}\,\mathrm{\%1}\,\Gamma (a) \left/ {\vrule 
height0.80em width0em depth0.80em} \right. \!  \! (\pi  
\Gamma ( - {\displaystyle \frac {a}{2}}  + {\displaystyle \frac {
b}{2}}  + n - 1 + c + 2\,m)\,\mathrm{\%3}\,\Gamma ( - c + 1)\,
\Gamma ( - b + 1)) \\
\mathrm{\%1} := \Gamma ({\displaystyle \frac {a}{2}}  + 
{\displaystyle \frac {b}{2}}  + n - c) \\
\mathrm{\%2} := \Gamma (2\,c + 2\,m - {\displaystyle \frac {a}{2}
}  - {\displaystyle \frac {b}{2}}  - n) \\
\mathrm{\%3} := \Gamma (n - 1 + {\displaystyle \frac {a}{2}}  + 
{\displaystyle \frac {b}{2}} ) \\
\mathrm{\%4} := \Gamma (c - 1 + 2\,m - {\displaystyle \frac {b}{2
}}  - {\displaystyle \frac {a}{2}}  + n) }
}
\end{maplelatex}

\begin{maplelatex}
\mapleinline{inert}{2d}{389, ":   ",
-WPMR(2*m+2*n-1,2*n-1,2*a-b+2*c-4*n+1-2*m,1+2*c-b-2*n,c)*Pi/sin(Pi*(-c
+b+2*n))/GAMMA(-c+b+2*n-1)/GAMMA(a+1-b+2*c-2*n)*GAMMA(2*n-2*c)/GAMMA(-
a+b-c+2*n+2*m)*GAMMA(1-a+2*m+c)/GAMMA(-c+1)/GAMMA(-b+1)*GAMMA(a)-1/2*W
PMR(2*m+2*n-1,2*n-1,2*a-b-2*n+1-2*m,-b+1,-c+2*n)*GAMMA(1-a+2*m+c)*GAMM
A(b-2*c+2*n)*GAMMA(-2*a+b-2*c+4*n+2*m)*GAMMA(-c+2*n)*GAMMA(2*c-2*n)*GA
MMA(a)*(cos(Pi*(b-c))-cos(Pi*(c+b)))/sin(Pi*(-c+b+2*n))/Pi/GAMMA(-c+b+
2*n-1)/GAMMA(2*n-2*a+b+2*m)/GAMMA(-a+b-c+2*n+2*m)/GAMMA(1+a-b)+WPMR(2*
m+2*n-1,2*m,2*n-2*a+b+2*m,b,-a+b-c+2*n+2*m)*sin(Pi*(2*n-2*a+b+2*m))/Pi
*GAMMA(1-a+2*m+c)/GAMMA(-a+b+2*n+2*m)/GAMMA(a+1-b+c-2*n)*GAMMA(2*a+1-2
*b+2*c-4*n-2*m)*GAMMA(b-2*c+2*n)*GAMMA(-2*a+b-2*c+4*n+2*m)/GAMMA(-c+b+
2*n-1)*GAMMA(c)/GAMMA(-b+1)*GAMMA(a)+GAMMA(a)*GAMMA(1-a+2*m+c)*GAMMA(c
)*WPMR(2*m+2*n-1,2*m,2*a-b+2*c-4*n+1-2*m,1+2*c-b-2*n,a+1-b+c-2*n)*GAMM
A(-1-2*a+2*b-2*c+4*n+2*m)/GAMMA(-c+b+2*n-1)/GAMMA(2*n-2*a+b+2*m)/GAMMA
(-a+b-c+2*n+2*m)/GAMMA(a+1-b+2*c-2*n);}{%
\maplemultiline{
389\mbox{:~~~}  -  \\
\mathrm{WPMR}(2\,m + 2\,n - 1, \,2\,n - 1, \,2\,a - b + 2\,c - 4
\,n + 1 - 2\,m, \,1 + 2\,c - b - 2\,n, \,c) \\
\pi \,\Gamma (2\,n - 2\,c)\,\mathrm{\%4}\,\Gamma (a)/(\mathrm{sin
}(\pi \,( - c + b + 2\,n))\,\mathrm{\%3}\,\Gamma (a + 1 - b + 2\,
c - 2\,n)\,\Gamma (\mathrm{\%1}) \\
\Gamma ( - c + 1)\,\Gamma ( - b + 1))\mbox{} - {\displaystyle 
\frac {1}{2}}  \\
\mathrm{WPMR}(2\,m + 2\,n - 1, \,2\,n - 1, \,2\,a - b - 2\,n + 1
 - 2\,m, \, - b + 1, \, - c + 2\,n)\,\mathrm{\%4} \\
\Gamma (b - 2\,c + 2\,n)\,\Gamma ( - 2\,a + b - 2\,c + 4\,n + 2\,
m)\,\Gamma ( - c + 2\,n)\,\Gamma (2\,c - 2\,n)\,\Gamma (a) \\
(\mathrm{cos}(\pi \,(b - c)) - \mathrm{cos}(\pi \,(c + b)))/(
\mathrm{sin}(\pi \,( - c + b + 2\,n))\,\pi \,\mathrm{\%3}\,\Gamma
 (\mathrm{\%2})\,\Gamma (\mathrm{\%1}) \\
\Gamma (1 + a - b))\mbox{} + \mathrm{WPMR}(2\,m + 2\,n - 1, \,2\,
m, \,\mathrm{\%2}, \,b, \,\mathrm{\%1})\,\mathrm{sin}(\pi \,
\mathrm{\%2})\,\mathrm{\%4} \\
\Gamma (2\,a + 1 - 2\,b + 2\,c - 4\,n - 2\,m)\,\Gamma (b - 2\,c
 + 2\,n)\,\Gamma ( - 2\,a + b - 2\,c + 4\,n + 2\,m) \\
\Gamma (c)\,\Gamma (a)/(\pi \,\Gamma ( - a + b + 2\,n + 2\,m)\,
\Gamma (a + 1 - b + c - 2\,n)\,\mathrm{\%3}\,\Gamma ( - b + 1))
\mbox{} + \Gamma (a) \\
\mathrm{\%4}\,\Gamma (c)\mathrm{WPMR}(2\,m + 2\,n - 1, \,2\,m, \,
2\,a - b + 2\,c - 4\,n + 1 - 2\,m,  \\
1 + 2\,c - b - 2\,n, \,a + 1 - b + c - 2\,n)\Gamma ( - 1 - 2\,a
 + 2\,b - 2\,c + 4\,n + 2\,m)/(\mathrm{\%3} \\
\Gamma (\mathrm{\%2})\,\Gamma (\mathrm{\%1})\,\Gamma (a + 1 - b
 + 2\,c - 2\,n)) \\
\mathrm{\%1} :=  - a + b - c + 2\,n + 2\,m \\
\mathrm{\%2} := 2\,n - 2\,a + b + 2\,m \\
\mathrm{\%3} := \Gamma ( - c + b + 2\,n - 1) \\
\mathrm{\%4} := \Gamma (1 - a + 2\,m + c) }
}
\end{maplelatex}

\begin{maplelatex}
\mapleinline{inert}{2d}{390, ":   ",
-1/2*WPMR(2*m+2*n-1,2*n-1,-2*c+1+a-2*m,-b+1,n+1/2*a-1/2*b)*GAMMA(1/2*a
+1/2*b+n-c)*GAMMA(-1/2*a-1/2*b+n+c+2*m)*GAMMA(n+1/2*a-1/2*b)*GAMMA(1+2
*c+2*m-1/2*a-1/2*b-n)*GAMMA(b-a)*GAMMA(a)*(cos(1/2*Pi*(-a-b+2*n))-cos(
1/2*Pi*(-a+3*b+2*n)))/Pi/sin(1/2*Pi*(a+b+2*n))/GAMMA(n-1+1/2*a+1/2*b)/
GAMMA(2*c-a+2*m)/GAMMA(-c+1+1/2*a-1/2*b+n)/GAMMA(2*m+c)-WPMR(2*m+2*n-1
,2*n-1,-2*c+1-2*m+b,-a+1,-1/2*a+1/2*b+n)*Pi/sin(1/2*Pi*(a+b+2*n))/GAMM
A(n-1+1/2*a+1/2*b)/GAMMA(-b+2*c+2*m)/GAMMA(1-c-1/2*a+1/2*b+n)*GAMMA(-1
/2*a-1/2*b+n+c+2*m)/GAMMA(1/2*a-1/2*b-n+1)*GAMMA(1+2*c+2*m-1/2*a-1/2*b
-n)*GAMMA(1/2*a+1/2*b+n-c)/GAMMA(2*m+c)*GAMMA(a-b)/GAMMA(-b+1)+GAMMA(1
/2*a+1/2*b+n-c)*GAMMA(1+2*c+2*m-1/2*a-1/2*b-n)*GAMMA(-1/2*a+1/2*b+n)*G
AMMA(2*c+2*m-1)/GAMMA(n-1+1/2*a+1/2*b)/GAMMA(2*c-a+2*m)/GAMMA(2*m+c)/G
AMMA(-b+2*c+2*m)/GAMMA(1-c-1/2*a+1/2*b+n)*GAMMA(-1/2*a-1/2*b+n+c+2*m)*
WPMR(2*m+2*n-1,2*m,-a+1,-2*c+1-2*m+b,-c+1)-sin(Pi*(-2*c+a-2*m))*WPMR(2
*m+2*n-1,2*m,b,2*c-a+2*m,2*m+c)/Pi/GAMMA(-1/2*a+1/2*b+n+c+2*m)*GAMMA(-
1/2*a-1/2*b+n+c+2*m)*GAMMA(-1/2*a+1/2*b+n)/GAMMA(n-1+1/2*a+1/2*b)*GAMM
A(1+2*c+2*m-1/2*a-1/2*b-n)*GAMMA(1/2*a+1/2*b+n-c)*GAMMA(-2*c+1-2*m)/GA
MMA(-c+1)/GAMMA(-b+1)*GAMMA(a);}{%
\maplemultiline{
390\mbox{:~~~}  - {\displaystyle \frac {1}{2}} \mathrm{
WPMR}(2\,m + 2\,n - 1, \,2\,n - 1, \, - 2\,c + 1 + a - 2\,m, \,
 - b + 1, \,n + {\displaystyle \frac {a}{2}}  - {\displaystyle 
\frac {b}{2}} )\,\mathrm{\%1} \\
\mathrm{\%4}\,\Gamma (n + {\displaystyle \frac {a}{2}}  - 
{\displaystyle \frac {b}{2}} )\,\mathrm{\%2}\,\Gamma (b - a)\,
\Gamma (a) \\
(\mathrm{cos}({\displaystyle \frac {\pi \,( - a - b + 2\,n)}{2}} 
) - \mathrm{cos}({\displaystyle \frac {\pi \,( - a + 3\,b + 2\,n)
}{2}} )) \left/ {\vrule height0.80em width0em depth0.80em}
 \right. \!  \! (\pi \,\mathrm{sin}({\displaystyle \frac {\pi \,(
a + b + 2\,n)}{2}} ) \\
\mathrm{\%3}\,\Gamma (2\,c - a + 2\,m)\,\Gamma ( - c + 1 + 
{\displaystyle \frac {a}{2}}  - {\displaystyle \frac {b}{2}}  + n
)\,\Gamma (2\,m + c))\mbox{} -  \\
\mathrm{WPMR}(2\,m + 2\,n - 1, \,2\,n - 1, \, - 2\,c + 1 - 2\,m
 + b, \, - a + 1, \, - {\displaystyle \frac {a}{2}}  + 
{\displaystyle \frac {b}{2}}  + n)\,\pi \,\mathrm{\%4}\,\mathrm{
\%2} \\
\mathrm{\%1}\,\Gamma (a - b) \left/ {\vrule 
height0.80em width0em depth0.80em} \right. \!  \! (\mathrm{sin}(
{\displaystyle \frac {\pi \,(a + b + 2\,n)}{2}} )\,\mathrm{\%3}\,
\Gamma ( - b + 2\,c + 2\,m) \\
\Gamma (1 - c - {\displaystyle \frac {a}{2}}  + {\displaystyle 
\frac {b}{2}}  + n)\,\Gamma ({\displaystyle \frac {a}{2}}  - 
{\displaystyle \frac {b}{2}}  - n + 1)\,\Gamma (2\,m + c)\,\Gamma
 ( - b + 1))\mbox{} + \mathrm{\%1}\,\mathrm{\%2} \\
\Gamma ( - {\displaystyle \frac {a}{2}}  + {\displaystyle \frac {
b}{2}}  + n)\,\Gamma (2\,c + 2\,m - 1)\,\mathrm{\%4} \\
\mathrm{WPMR}(2\,m + 2\,n - 1, \,2\,m, \, - a + 1, \, - 2\,c + 1
 - 2\,m + b, \, - c + 1) \left/ {\vrule 
height0.80em width0em depth0.80em} \right. \!  \! (\mathrm{\%3}
 \\
\Gamma (2\,c - a + 2\,m)\,\Gamma (2\,m + c)\,\Gamma ( - b + 2\,c
 + 2\,m)\,\Gamma (1 - c - {\displaystyle \frac {a}{2}}  + 
{\displaystyle \frac {b}{2}}  + n))\mbox{} -  \\
\mathrm{sin}(\pi \,( - 2\,c + a - 2\,m))\,\mathrm{WPMR}(2\,m + 2
\,n - 1, \,2\,m, \,b, \,2\,c - a + 2\,m, \,2\,m + c)\,\mathrm{\%4
} \\
\Gamma ( - {\displaystyle \frac {a}{2}}  + {\displaystyle \frac {
b}{2}}  + n)\,\mathrm{\%2}\,\mathrm{\%1}\,\Gamma ( - 2\,c + 1 - 2
\,m)\,\Gamma (a) \left/ {\vrule height0.80em width0em depth0.80em
} \right. \!  \! (\pi \,\Gamma ( - {\displaystyle \frac {a}{2}} 
 + {\displaystyle \frac {b}{2}}  + n + c + 2\,m) \\
\mathrm{\%3}\,\Gamma ( - c + 1)\,\Gamma ( - b + 1)) \\
\mathrm{\%1} := \Gamma ({\displaystyle \frac {a}{2}}  + 
{\displaystyle \frac {b}{2}}  + n - c) \\
\mathrm{\%2} := \Gamma (1 + 2\,c + 2\,m - {\displaystyle \frac {a
}{2}}  - {\displaystyle \frac {b}{2}}  - n) \\
\mathrm{\%3} := \Gamma (n - 1 + {\displaystyle \frac {a}{2}}  + 
{\displaystyle \frac {b}{2}} ) \\
\mathrm{\%4} := \Gamma ( - {\displaystyle \frac {a}{2}}  - 
{\displaystyle \frac {b}{2}}  + n + c + 2\,m) }
}
\end{maplelatex}

\begin{maplelatex}
\mapleinline{inert}{2d}{391, ":   ",
1/2*GAMMA(a)*GAMMA(-a+2*m+c)*GAMMA(2*c-2*n-1)*GAMMA(1-c+2*n)*GAMMA(-2*
a+b-2*c+4*n+2*m+1)*GAMMA(b-2*c+2*n+1)*WPMR(2*m+2*n-1,2*n,-b+1,2*a-b-2*
n+1-2*m,1-c+2*n)*(cos(Pi*(b-c))-cos(Pi*(c+b)))/Pi/sin(Pi*(-c+b+2*n))/G
AMMA(-c+b+2*n)/GAMMA(-a+b-c+2*n+2*m)/GAMMA(2*n-2*a+b+2*m)/GAMMA(1+a-b)
+GAMMA(a)*GAMMA(-a+2*m+c)*Pi/sin(Pi*(-c+b+2*n))/GAMMA(-c+b+2*n)/GAMMA(
-a+b-c+2*n+2*m)*WPMR(2*m+2*n-1,2*n,2*a-b+2*c-4*n-2*m,2*c-b-2*n,c)*GAMM
A(1-2*c+2*n)/GAMMA(-c+1)/GAMMA(-b+1)/GAMMA(a-b+2*c-2*n)+WPMR(2*m+2*n-1
,2*m-1,b,2*n-2*a+b+2*m,-a+b-c+2*n+2*m)*sin(Pi*(2*n-2*a+b+2*m))/Pi/GAMM
A(a-b+c-2*n)/GAMMA(-a+b+2*n+2*m)*GAMMA(2*a-2*b+2*c-4*n-2*m)/GAMMA(-c+b
+2*n)*GAMMA(b-2*c+2*n+1)*GAMMA(-2*a+b-2*c+4*n+2*m+1)*GAMMA(-a+2*m+c)*G
AMMA(c)/GAMMA(-b+1)*GAMMA(a)+GAMMA(a)*GAMMA(-a+2*m+c)*GAMMA(c)*WPMR(2*
m+2*n-1,2*m-1,2*a-b+2*c-4*n-2*m,2*c-b-2*n,a-b+c-2*n)*GAMMA(2*m-2*a+2*b
-2*c+4*n)/GAMMA(-c+b+2*n)/GAMMA(2*n-2*a+b+2*m)/GAMMA(-a+b-c+2*n+2*m)/G
AMMA(a-b+2*c-2*n);}{%
\maplemultiline{
391\mbox{:~~~} {\displaystyle \frac {1}{2}} \Gamma (a)\,
\Gamma ( - a + 2\,m + c)\,\Gamma (2\,c - 2\,n - 1)\,\Gamma (1 - c
 + 2\,n) \\
\Gamma ( - 2\,a + b - 2\,c + 4\,n + 2\,m + 1)\,\Gamma (b - 2\,c
 + 2\,n + 1) \\
\mathrm{WPMR}(2\,m + 2\,n - 1, \,2\,n, \, - b + 1, \,2\,a - b - 2
\,n + 1 - 2\,m, \,1 - c + 2\,n) \\
(\mathrm{cos}(\pi \,(b - c)) - \mathrm{cos}(\pi \,(c + b)))/(\pi 
\,\mathrm{sin}(\pi \,( - c + b + 2\,n))\,\Gamma ( - c + b + 2\,n)
\,\Gamma (\mathrm{\%1}) \\
\Gamma (\mathrm{\%2})\,\Gamma (1 + a - b))\mbox{} + \Gamma (a)\,
\Gamma ( - a + 2\,m + c)\,\pi  \\
\mathrm{WPMR}(2\,m + 2\,n - 1, \,2\,n, \,2\,a - b + 2\,c - 4\,n
 - 2\,m, \,2\,c - b - 2\,n, \,c) \\
\Gamma (1 - 2\,c + 2\,n)/(\mathrm{sin}(\pi \,( - c + b + 2\,n))\,
\Gamma ( - c + b + 2\,n)\,\Gamma (\mathrm{\%1})\,\Gamma ( - c + 1
)\,\Gamma ( - b + 1) \\
\Gamma (a - b + 2\,c - 2\,n))\mbox{} + \mathrm{WPMR}(2\,m + 2\,n
 - 1, \,2\,m - 1, \,b, \,\mathrm{\%2}, \,\mathrm{\%1})\,\mathrm{
sin}(\pi \,\mathrm{\%2}) \\
\Gamma (2\,a - 2\,b + 2\,c - 4\,n - 2\,m)\,\Gamma (b - 2\,c + 2\,
n + 1) \\
\Gamma ( - 2\,a + b - 2\,c + 4\,n + 2\,m + 1)\,\Gamma ( - a + 2\,
m + c)\,\Gamma (c)\,\Gamma (a)/(\pi \,\Gamma (a - b + c - 2\,n)
 \\
\Gamma ( - a + b + 2\,n + 2\,m)\,\Gamma ( - c + b + 2\,n)\,\Gamma
 ( - b + 1))\mbox{} + \Gamma (a)\,\Gamma ( - a + 2\,m + c)\,
\Gamma (c) \\
\mathrm{WPMR}(2\,m + 2\,n - 1, \,2\,m - 1, \,2\,a - b + 2\,c - 4
\,n - 2\,m, \,2\,c - b - 2\,n,  \\
a - b + c - 2\,n)\Gamma (2\,m - 2\,a + 2\,b - 2\,c + 4\,n)/(
\Gamma ( - c + b + 2\,n)\,\Gamma (\mathrm{\%2})\,\Gamma (\mathrm{
\%1}) \\
\Gamma (a - b + 2\,c - 2\,n)) \\
\mathrm{\%1} :=  - a + b - c + 2\,n + 2\,m \\
\mathrm{\%2} := 2\,n - 2\,a + b + 2\,m }
}
\end{maplelatex}

\begin{maplelatex}
\mapleinline{inert}{2d}{392, ":   ",
-Pi/cos(1/2*Pi*(a+b+2*n))*WPMR(2*m+2*n-1,2*n,-a+1,-2*c+2-2*m+b,1/2-1/2
*a+n+1/2*b)/GAMMA(-b-1+2*c+2*m)/GAMMA(1/2+1/2*a-n-1/2*b)*GAMMA(c-1/2-1
/2*a+2*m-1/2*b+n)/GAMMA(3/2-c-1/2*a+1/2*b+n)/GAMMA(-1/2+1/2*a+n+1/2*b)
*GAMMA(2*c+2*m-1/2-1/2*a-n-1/2*b)*GAMMA(1/2+1/2*b+n-c+1/2*a)/GAMMA(2*m
+c-1)*GAMMA(a-b)/GAMMA(-b+1)+1/2*WPMR(2*m+2*n-1,2*n,-2*c+2+a-2*m,-b+1,
n+1/2*a+1/2-1/2*b)*GAMMA(1/2+1/2*b+n-c+1/2*a)*GAMMA(2*c+2*m-1/2-1/2*a-
n-1/2*b)*GAMMA(c-1/2-1/2*a+2*m-1/2*b+n)*GAMMA(n+1/2*a+1/2-1/2*b)*GAMMA
(b-a)*GAMMA(a)*(-sin(1/2*Pi*(-a+3*b+2*n))+sin(1/2*Pi*(-a-b+2*n)))/cos(
1/2*Pi*(a+b+2*n))/Pi/GAMMA(-1/2+1/2*a+n+1/2*b)/GAMMA(-c+3/2+1/2*a-1/2*
b+n)/GAMMA(2*c-1-a+2*m)/GAMMA(2*m+c-1)+GAMMA(1/2+1/2*b+n-c+1/2*a)*GAMM
A(2*c+2*m-1/2-1/2*a-n-1/2*b)*WPPR(2*m-1,2*n,a,-2*c+2+a-2*m,1/2+1/2*a-n
-1/2*b)*GAMMA(a)*GAMMA(n+1/2*a+1/2-1/2*b)/GAMMA(-1/2+1/2*a+n+1/2*b)/GA
MMA(2*c-1-a+2*m)/GAMMA(1+a-b)/GAMMA(1+a-c);}{%
\maplemultiline{
392\mbox{:~~~}  - \pi \,\mathrm{WPMR}(2\,m + 2\,n - 1, 
\,2\,n, \, - a + 1, \, - 2\,c + 2 - 2\,m + b, \,{\displaystyle 
\frac {1}{2}}  - {\displaystyle \frac {a}{2}}  + n + 
{\displaystyle \frac {b}{2}} ) \\
\Gamma (c - {\displaystyle \frac {1}{2}}  - {\displaystyle 
\frac {a}{2}}  + 2\,m - {\displaystyle \frac {b}{2}}  + n)\,
\Gamma (2\,c + 2\,m - {\displaystyle \frac {1}{2}}  - 
{\displaystyle \frac {a}{2}}  - n - {\displaystyle \frac {b}{2}} 
)\,\Gamma ({\displaystyle \frac {1}{2}}  + {\displaystyle \frac {
b}{2}}  + n - c + {\displaystyle \frac {a}{2}} ) \\
\Gamma (a - b) \left/ {\vrule height0.80em width0em depth0.80em}
 \right. \!  \! (\mathrm{cos}({\displaystyle \frac {\pi \,(a + b
 + 2\,n)}{2}} )\,\Gamma ( - b - 1 + 2\,c + 2\,m)\,\Gamma (
{\displaystyle \frac {1}{2}}  + {\displaystyle \frac {a}{2}}  - n
 - {\displaystyle \frac {b}{2}} ) \\
\Gamma ({\displaystyle \frac {3}{2}}  - c - {\displaystyle 
\frac {a}{2}}  + {\displaystyle \frac {b}{2}}  + n)\,\Gamma ( - 
{\displaystyle \frac {1}{2}}  + {\displaystyle \frac {a}{2}}  + n
 + {\displaystyle \frac {b}{2}} )\,\Gamma (2\,m + c - 1)\,\Gamma 
( - b + 1))\mbox{} + {\displaystyle \frac {1}{2}}  \\
\mathrm{WPMR}(2\,m + 2\,n - 1, \,2\,n, \, - 2\,c + 2 + a - 2\,m, 
\, - b + 1, \,n + {\displaystyle \frac {a}{2}}  + {\displaystyle 
\frac {1}{2}}  - {\displaystyle \frac {b}{2}} ) \\
\Gamma ({\displaystyle \frac {1}{2}}  + {\displaystyle \frac {b}{
2}}  + n - c + {\displaystyle \frac {a}{2}} )\,\Gamma (2\,c + 2\,
m - {\displaystyle \frac {1}{2}}  - {\displaystyle \frac {a}{2}} 
 - n - {\displaystyle \frac {b}{2}} )\,\Gamma (c - 
{\displaystyle \frac {1}{2}}  - {\displaystyle \frac {a}{2}}  + 2
\,m - {\displaystyle \frac {b}{2}}  + n) \\
\Gamma (n + {\displaystyle \frac {a}{2}}  + {\displaystyle 
\frac {1}{2}}  - {\displaystyle \frac {b}{2}} )\,\Gamma (b - a)\,
\Gamma (a)\,( - \mathrm{sin}({\displaystyle \frac {\pi \,( - a + 
3\,b + 2\,n)}{2}} ) + \mathrm{sin}({\displaystyle \frac {\pi \,(
 - a - b + 2\,n)}{2}} )) \\
 \left/ {\vrule height0.80em width0em depth0.80em} \right. \! 
 \! (\mathrm{cos}({\displaystyle \frac {\pi \,(a + b + 2\,n)}{2}
} )\,\pi \,\Gamma ( - {\displaystyle \frac {1}{2}}  + 
{\displaystyle \frac {a}{2}}  + n + {\displaystyle \frac {b}{2}} 
)\,\Gamma ( - c + {\displaystyle \frac {3}{2}}  + {\displaystyle 
\frac {a}{2}}  - {\displaystyle \frac {b}{2}}  + n) \\
\Gamma (2\,c - 1 - a + 2\,m)\,\Gamma (2\,m + c - 1))\mbox{} + 
\Gamma ({\displaystyle \frac {1}{2}}  + {\displaystyle \frac {b}{
2}}  + n - c + {\displaystyle \frac {a}{2}} ) \\
\Gamma (2\,c + 2\,m - {\displaystyle \frac {1}{2}}  - 
{\displaystyle \frac {a}{2}}  - n - {\displaystyle \frac {b}{2}} 
) \\
\mathrm{WPPR}(2\,m - 1, \,2\,n, \,a, \, - 2\,c + 2 + a - 2\,m, \,
{\displaystyle \frac {1}{2}}  + {\displaystyle \frac {a}{2}}  - n
 - {\displaystyle \frac {b}{2}} )\,\Gamma (a)\,\Gamma (n + 
{\displaystyle \frac {a}{2}}  + {\displaystyle \frac {1}{2}}  - 
{\displaystyle \frac {b}{2}} ) \\
 \left/ {\vrule height0.80em width0em depth0.80em} \right. \! 
 \! (\Gamma ( - {\displaystyle \frac {1}{2}}  + {\displaystyle 
\frac {a}{2}}  + n + {\displaystyle \frac {b}{2}} )\,\Gamma (2\,c
 - 1 - a + 2\,m)\,\Gamma (1 + a - b)\,\Gamma (1 + a - c)) }
}
\end{maplelatex}

\begin{maplelatex}
\mapleinline{inert}{2d}{393, ":   ",
-Pi*WPMR(2*m+2*n-2,2*n-1,-c+1-a+2*m,-2*n+2-c+2*b-a,1+b-c)/sin(Pi*(2*n-
b+2*m))/GAMMA(1-c+b-a+2*m)/GAMMA(a-2*m+c)*GAMMA(2-c+2*b-2*n)*GAMMA(2*c
-2-2*b+2*n)/GAMMA(c-b)/GAMMA(c-a)*GAMMA(c)/GAMMA(a)+1/2*(cos(Pi*(2*c-3
*b+4*n+a))-cos(Pi*(a-b)))*WPMR(2*m+2*n-2,2*n-1,c-1-2*b+2*n-a+2*m,c-a,c
-1-b+2*n)*GAMMA(c+a-1-2*b+2*n)*GAMMA(c-1-b+2*n)*GAMMA(2-c+2*b-2*n)*GAM
MA(-2*n-2*c+2+2*b)*GAMMA(c)/Pi/sin(Pi*(2*n-b+2*m))/GAMMA(-c+2+2*b-2*n+
a-2*m)/GAMMA(c-1-b+2*n-a+2*m)/GAMMA(a);}{%
\maplemultiline{
393\mbox{:~~~}  - \pi  \\
\mathrm{WPMR}(2\,m + 2\,n - 2, \,2\,n - 1, \, - c + 1 - a + 2\,m
, \, - 2\,n + 2 - c + 2\,b - a, \,1 + b - c) \\
\Gamma (2 - c + 2\,b - 2\,n)\,\Gamma (2\,c - 2 - 2\,b + 2\,n)\,
\Gamma (c)/(\mathrm{sin}(\pi \,(2\,n - b + 2\,m)) \\
\Gamma (1 - c + b - a + 2\,m)\,\Gamma (a - 2\,m + c)\,\Gamma (c
 - b)\,\Gamma (c - a)\,\Gamma (a))\mbox{} + {\displaystyle 
\frac {1}{2}}  \\
(\mathrm{cos}(\pi \,(2\,c - 3\,b + 4\,n + a)) - \mathrm{cos}(\pi 
\,(a - b))) \\
\mathrm{WPMR}(2\,m + 2\,n - 2, \,2\,n - 1, \,c - 1 - 2\,b + 2\,n
 - a + 2\,m, \,c - a, \,c - 1 - b + 2\,n) \\
\Gamma (c + a - 1 - 2\,b + 2\,n)\,\Gamma (c - 1 - b + 2\,n)\,
\Gamma (2 - c + 2\,b - 2\,n) \\
\Gamma ( - 2\,n - 2\,c + 2 + 2\,b)\,\Gamma (c)/(\pi \,\mathrm{sin
}(\pi \,(2\,n - b + 2\,m)) \\
\Gamma ( - c + 2 + 2\,b - 2\,n + a - 2\,m)\,\Gamma (c - 1 - b + 2
\,n - a + 2\,m)\,\Gamma (a)) }
}
\end{maplelatex}

\begin{maplelatex}
\mapleinline{inert}{2d}{394, ":   ",
WPMR(2*m+2*n-2,2*n-1,-b+1,2*a+2-b-2*n-2*m,-1/2*b+1+a-1/2*c-m)*GAMMA(-1
/2*b-2*n+2+a-1/2*c-m)*GAMMA(-a-1+2*n+2*m)*GAMMA(b+2*n-2-2*a+c+2*m)*GAM
MA(2+2*a-2*n-2*m)/GAMMA(1/2*b-a+1/2*c+m)/GAMMA(-1/2*c+1/2*b+m)*GAMMA(1
/2*c+m+1/2*b)/GAMMA(1+a-b)/GAMMA(c)/GAMMA(b)+1/2*(cos(1/2*Pi*(c+3*b+8*
n-6*a+6*m))-cos(1/2*Pi*(-b+2*a+c-2*m)))*WPMR(2*m+2*n-2,2*n-1,-2*a+c+2*
n+2*m-1,c,1/2*b+2*n-1-a+1/2*c+m)*GAMMA(-2*n-b+2+2*a-c-2*m)*GAMMA(2*n-1
-2*a+b+2*m)*GAMMA(1/2*c+m+1/2*b)*GAMMA(2+2*a-2*n-2*m)*GAMMA(-a-1+2*n+2
*m)/Pi/sin(1/2*Pi*(b+4*n-2*a+c+2*m))/GAMMA(-a+c+2*n+2*m-1)/GAMMA(2*a-c
-2*n-2*m+2)/GAMMA(-1/2*c+1/2*b+m);}{%
\maplemultiline{
394\mbox{:~~~} \mathrm{WPMR}(2\,m + 2\,n - 2, \,2\,n - 1
, \, - b + 1, \,2\,a + 2 - b - 2\,n - 2\,m,  \\
 - {\displaystyle \frac {b}{2}}  + 1 + a - {\displaystyle \frac {
c}{2}}  - m)\Gamma ( - {\displaystyle \frac {b}{2}}  - 2\,n + 2
 + a - {\displaystyle \frac {c}{2}}  - m)\,\Gamma ( - a - 1 + 2\,
n + 2\,m) \\
\Gamma (b + 2\,n - 2 - 2\,a + c + 2\,m)\,\Gamma (2 + 2\,a - 2\,n
 - 2\,m)\,\Gamma ({\displaystyle \frac {c}{2}}  + m + 
{\displaystyle \frac {b}{2}} ) \left/ {\vrule 
height0.80em width0em depth0.80em} \right. \!  \! ( \\
\Gamma ({\displaystyle \frac {b}{2}}  - a + {\displaystyle 
\frac {c}{2}}  + m)\,\Gamma ( - {\displaystyle \frac {c}{2}}  + 
{\displaystyle \frac {b}{2}}  + m)\,\Gamma (1 + a - b)\,\Gamma (c
)\,\Gamma (b))\mbox{} + {\displaystyle \frac {1}{2}}  \\
(\mathrm{cos}({\displaystyle \frac {\pi \,(c + 3\,b + 8\,n - 6\,a
 + 6\,m)}{2}} ) - \mathrm{cos}({\displaystyle \frac {\pi \,( - b
 + 2\,a + c - 2\,m)}{2}} ))\mathrm{WPMR}( \\
2\,m + 2\,n - 2, \,2\,n - 1, \, - 2\,a + c + 2\,n + 2\,m - 1, \,c
, \,{\displaystyle \frac {b}{2}}  + 2\,n - 1 - a + 
{\displaystyle \frac {c}{2}}  + m) \\
\Gamma ( - 2\,n - b + 2 + 2\,a - c - 2\,m)\,\Gamma (2\,n - 1 - 2
\,a + b + 2\,m)\,\Gamma ({\displaystyle \frac {c}{2}}  + m + 
{\displaystyle \frac {b}{2}} ) \\
\Gamma (2 + 2\,a - 2\,n - 2\,m)\,\Gamma ( - a - 1 + 2\,n + 2\,m)
 \left/ {\vrule height0.80em width0em depth0.80em} \right. \! 
 \! (\pi  \\
\mathrm{sin}({\displaystyle \frac {\pi \,(b + 4\,n - 2\,a + c + 2
\,m)}{2}} )\,\Gamma ( - a + c + 2\,n + 2\,m - 1) \\
\Gamma (2\,a - c - 2\,n - 2\,m + 2)\,\Gamma ( - {\displaystyle 
\frac {c}{2}}  + {\displaystyle \frac {b}{2}}  + m)) }
}
\end{maplelatex}

\begin{maplelatex}
\mapleinline{inert}{2d}{395, ":   ",
Pi*WPMR(2*m+2*n-1,2*n-1,-a+2*m-c+2,-2*n+2-c+2*b-a,1+b-c)/sin(Pi*(2*n-b
+2*m))*GAMMA(2-c+2*b-2*n)*GAMMA(2*c-2-2*b+2*n)/GAMMA(-a+b+2*m-c+2)/GAM
MA(a-2*m-1+c)/GAMMA(c-b)/GAMMA(c-a)*GAMMA(c)/GAMMA(a)+1/2*(-cos(Pi*(2*
c-3*b+4*n+a))+cos(Pi*(a-b)))*WPMR(2*m+2*n-1,2*n-1,-a+2*n-2*b+2*m+c,c-a
,c-1-b+2*n)*GAMMA(c+a-1-2*b+2*n)*GAMMA(-2*n-2*c+2+2*b)*GAMMA(c-1-b+2*n
)*GAMMA(2-c+2*b-2*n)*GAMMA(c)/Pi/sin(Pi*(2*n-b+2*m))/GAMMA(a-2*n+2*b-2
*m-c+1)/GAMMA(-a+2*n-b+2*m+c)/GAMMA(a);}{%
\maplemultiline{
395\mbox{:~~~} \pi  \,
\mathrm{WPMR}(2\,m + 2\,n - 1, \,2\,n - 1, \, - a + 2\,m - c + 2
, \, - 2\,n + 2 - c + 2\,b - a, \,1 + b - c) \\
\Gamma (2 - c + 2\,b - 2\,n)\,\Gamma (2\,c - 2 - 2\,b + 2\,n)\,
\Gamma (c)/(\mathrm{sin}(\pi \,(2\,n - b + 2\,m)) \\
\Gamma ( - a + b + 2\,m - c + 2)\,\Gamma (a - 2\,m - 1 + c)\,
\Gamma (c - b)\,\Gamma (c - a)\,\Gamma (a))\mbox{} + 
{\displaystyle \frac {1}{2}}  \\
( - \mathrm{cos}(\pi \,(2\,c - 3\,b + 4\,n + a)) + \mathrm{cos}(
\pi \,(a - b))) \\
\mathrm{WPMR}(2\,m + 2\,n - 1, \,2\,n - 1, \, - a + 2\,n - 2\,b
 + 2\,m + c, \,c - a, \,c - 1 - b + 2\,n) \\
\Gamma (c + a - 1 - 2\,b + 2\,n)\,\Gamma ( - 2\,n - 2\,c + 2 + 2
\,b)\,\Gamma (c - 1 - b + 2\,n) \\
\Gamma (2 - c + 2\,b - 2\,n)\,\Gamma (c)/(\pi \,\mathrm{sin}(\pi 
\,(2\,n - b + 2\,m))\,\Gamma (a - 2\,n + 2\,b - 2\,m - c + 1) \\
\Gamma ( - a + 2\,n - b + 2\,m + c)\,\Gamma (a)) }
}
\end{maplelatex}

\begin{maplelatex}
\mapleinline{inert}{2d}{396, ":   ",
WPMR(2*m+2*n-1,2*n-1,-b+1,2*a-b-2*n+1-2*m,-1/2*b+1/2+a-1/2*c-m)/GAMMA(
1/2*b+1/2-a+1/2*c+m)*GAMMA(b+2*n-1-2*a+c+2*m)*GAMMA(-1/2*b-2*n+3/2+a-1
/2*c-m)*GAMMA(1+2*a-2*n-2*m)*GAMMA(-a+2*n+2*m)*GAMMA(1/2*b+1/2+1/2*c+m
)/GAMMA(-1/2*c+1/2*b+1/2+m)/GAMMA(1+a-b)/GAMMA(c)/GAMMA(b)+1/2*(sin(1/
2*Pi*(c+3*b+8*n-6*a+6*m))-sin(1/2*Pi*(-b+2*a+c-2*m)))*WPMR(2*m+2*n-1,2
*n-1,-2*a+c+2*n+2*m,c,1/2*b+2*n-1/2-a+1/2*c+m)*GAMMA(1+2*a-2*n-2*m)*GA
MMA(2*n-2*a+b+2*m)*GAMMA(-2*n-b+1+2*a-c-2*m)*GAMMA(-a+2*n+2*m)*GAMMA(1
/2*b+1/2+1/2*c+m)/cos(1/2*Pi*(b+4*n-2*a+c+2*m))/Pi/GAMMA(2*a-c-2*n-2*m
+1)/GAMMA(-a+c+2*n+2*m)/GAMMA(-1/2*c+1/2*b+1/2+m);}{%
\maplemultiline{
396\mbox{:~~~} \mathrm{WPMR}(2\,m + 2\,n - 1, \,2\,n - 1
, \, - b + 1, \,2\,a - b - 2\,n + 1 - 2\,m,  \\
 - {\displaystyle \frac {b}{2}}  + {\displaystyle \frac {1}{2}} 
 + a - {\displaystyle \frac {c}{2}}  - m)\Gamma (b + 2\,n - 1 - 2
\,a + c + 2\,m)\,\Gamma ( - {\displaystyle \frac {b}{2}}  - 2\,n
 + {\displaystyle \frac {3}{2}}  + a - {\displaystyle \frac {c}{2
}}  - m) \\
\Gamma (1 + 2\,a - 2\,n - 2\,m)\,\Gamma ( - a + 2\,n + 2\,m)\,
\Gamma ({\displaystyle \frac {b}{2}}  + {\displaystyle \frac {1}{
2}}  + {\displaystyle \frac {c}{2}}  + m) \left/ {\vrule 
height0.80em width0em depth0.80em} \right. \!  \! ( \\
\Gamma ({\displaystyle \frac {b}{2}}  + {\displaystyle \frac {1}{
2}}  - a + {\displaystyle \frac {c}{2}}  + m)\,\Gamma ( - 
{\displaystyle \frac {c}{2}}  + {\displaystyle \frac {b}{2}}  + 
{\displaystyle \frac {1}{2}}  + m)\,\Gamma (1 + a - b)\,\Gamma (c
)\,\Gamma (b))\mbox{} + {\displaystyle \frac {1}{2}}  \\
(\mathrm{sin}({\displaystyle \frac {\pi \,(c + 3\,b + 8\,n - 6\,a
 + 6\,m)}{2}} ) - \mathrm{sin}({\displaystyle \frac {\pi \,( - b
 + 2\,a + c - 2\,m)}{2}} )) \\
\mathrm{WPMR}(2\,m + 2\,n - 1, \,2\,n - 1, \, - 2\,a + c + 2\,n
 + 2\,m, \,c, \,{\displaystyle \frac {b}{2}}  + 2\,n - 
{\displaystyle \frac {1}{2}}  - a + {\displaystyle \frac {c}{2}} 
 + m) \\
\Gamma (1 + 2\,a - 2\,n - 2\,m)\,\Gamma (2\,n - 2\,a + b + 2\,m)
\,\Gamma ( - 2\,n - b + 1 + 2\,a - c - 2\,m) \\
\Gamma ( - a + 2\,n + 2\,m)\,\Gamma ({\displaystyle \frac {b}{2}
}  + {\displaystyle \frac {1}{2}}  + {\displaystyle \frac {c}{2}
}  + m) \left/ {\vrule height0.80em width0em depth0.80em}
 \right. \!  \! (\mathrm{cos}({\displaystyle \frac {\pi \,(b + 4
\,n - 2\,a + c + 2\,m)}{2}} )\,\pi  \\
\Gamma (2\,a - c - 2\,n - 2\,m + 1)\,\Gamma ( - a + c + 2\,n + 2
\,m)\,\Gamma ( - {\displaystyle \frac {c}{2}}  + {\displaystyle 
\frac {b}{2}}  + {\displaystyle \frac {1}{2}}  + m)) }
}
\end{maplelatex}

\begin{maplelatex}
\mapleinline{inert}{2d}{397, ":   ",
Pi*WPMR(2*m+2*n-1,2*n,1-a-c+2*b-2*n,-c+1-a+2*m,1+b-c)/sin(Pi*(2*n-b+2*
m))*GAMMA(1-c+2*b-2*n)*GAMMA(2*c-1-2*b+2*n)/GAMMA(1-c+b-a+2*m)/GAMMA(a
-2*m+c)/GAMMA(c-b)/GAMMA(c-a)*GAMMA(c)/GAMMA(a)-1/2*(cos(Pi*(2*c-3*b+4
*n+a))-cos(Pi*(a-b)))*WPMR(2*m+2*n-1,2*n,-a+2*n-2*b+2*m+c,c-a,c-b+2*n)
*GAMMA(c)*GAMMA(1-2*c+2*b-2*n)*GAMMA(a+c-2*b+2*n)*GAMMA(c-b+2*n)*GAMMA
(1-c+2*b-2*n)/Pi/sin(Pi*(2*n-b+2*m))/GAMMA(a)/GAMMA(-a+2*n-b+2*m+c)/GA
MMA(a-2*n+2*b-2*m-c+1);}{%
\maplemultiline{
397\mbox{:~~~} \pi \,\mathrm{WPMR}(2\,m + 2\,n - 1, \,2
\,n, \,1 - a - c + 2\,b - 2\,n, \, - c + 1 - a + 2\,m, \,1 + b - 
c) \\
\Gamma (1 - c + 2\,b - 2\,n)\,\Gamma (2\,c - 1 - 2\,b + 2\,n)\,
\Gamma (c)/(\mathrm{sin}(\pi \,(2\,n - b + 2\,m)) \\
\Gamma (1 - c + b - a + 2\,m)\,\Gamma (a - 2\,m + c)\,\Gamma (c
 - b)\,\Gamma (c - a)\,\Gamma (a))\mbox{} - {\displaystyle 
\frac {1}{2}}  \\
(\mathrm{cos}(\pi \,(2\,c - 3\,b + 4\,n + a)) - \mathrm{cos}(\pi 
\,(a - b))) \\
\mathrm{WPMR}(2\,m + 2\,n - 1, \,2\,n, \, - a + 2\,n - 2\,b + 2\,
m + c, \,c - a, \,c - b + 2\,n)\,\Gamma (c) \\
\Gamma (1 - 2\,c + 2\,b - 2\,n)\,\Gamma (a + c - 2\,b + 2\,n)\,
\Gamma (c - b + 2\,n)\,\Gamma (1 - c + 2\,b - 2\,n)/(\pi  \\
\mathrm{sin}(\pi \,(2\,n - b + 2\,m))\,\Gamma (a)\,\Gamma ( - a
 + 2\,n - b + 2\,m + c) \\
\Gamma (a - 2\,n + 2\,b - 2\,m - c + 1)) }
}
\end{maplelatex}

\begin{maplelatex}
\mapleinline{inert}{2d}{398, ":   ",
WPMR(2*m+2*n-1,2*n,2*a-b-2*n+1-2*m,-b+1,-1/2*b+1+a-1/2*c-m)*GAMMA(-1/2
*b-2*n+a-1/2*c-m+1)*GAMMA(b+2*n-1-2*a+c+2*m)*GAMMA(1+2*a-2*n-2*m)*GAMM
A(-a+2*n+2*m)/GAMMA(1/2*b-a+1/2*c+m)/GAMMA(-1/2*c+1/2*b+m)*GAMMA(1/2*c
+m+1/2*b)/GAMMA(1+a-b)/GAMMA(c)/GAMMA(b)-1/2*(cos(1/2*Pi*(c+3*b+8*n-6*
a+6*m))-cos(1/2*Pi*(-b+2*a+c-2*m)))*WPMR(2*m+2*n-1,2*n,-2*a+c+2*n+2*m,
c,1/2*b+2*n-a+1/2*c+m)*GAMMA(-2*n-b+1+2*a-c-2*m)*GAMMA(2*n-2*a+b+2*m)*
GAMMA(1+2*a-2*n-2*m)*GAMMA(-a+2*n+2*m)*GAMMA(1/2*c+m+1/2*b)/Pi/sin(1/2
*Pi*(b+4*n-2*a+c+2*m))/GAMMA(-a+c+2*n+2*m)/GAMMA(2*a-c-2*n-2*m+1)/GAMM
A(-1/2*c+1/2*b+m);}{%
\maplemultiline{
398\mbox{:~~~}
\mathrm{WPMR}(2\,m + 2\,n - 1, \,2\,n, \,2\,a - b - 2\,n + 1 - 2
\,m, \, - b + 1, \, - {\displaystyle \frac {b}{2}}  + 1 + a - 
{\displaystyle \frac {c}{2}}  - m) \\
\Gamma ( - {\displaystyle \frac {b}{2}}  - 2\,n + a - 
{\displaystyle \frac {c}{2}}  - m + 1)\,\Gamma (b + 2\,n - 1 - 2
\,a + c + 2\,m)\,\Gamma (1 + 2\,a - 2\,n - 2\,m) \\
\Gamma ( - a + 2\,n + 2\,m)\,\Gamma ({\displaystyle \frac {c}{2}
}  + m + {\displaystyle \frac {b}{2}} ) \left/ {\vrule 
height0.80em width0em depth0.80em} \right. \!  \! (\Gamma (
{\displaystyle \frac {b}{2}}  - a + {\displaystyle \frac {c}{2}} 
 + m)\,\Gamma ( - {\displaystyle \frac {c}{2}}  + {\displaystyle 
\frac {b}{2}}  + m) \\
\Gamma (1 + a - b)\,\Gamma (c)\,\Gamma (b))\mbox{} - 
{\displaystyle \frac {1}{2}}  \\
(\mathrm{cos}({\displaystyle \frac {\pi \,(c + 3\,b + 8\,n - 6\,a
 + 6\,m)}{2}} ) - \mathrm{cos}({\displaystyle \frac {\pi \,( - b
 + 2\,a + c - 2\,m)}{2}} )) \\
\mathrm{WPMR}(2\,m + 2\,n - 1, \,2\,n, \, - 2\,a + c + 2\,n + 2\,
m, \,c, \,{\displaystyle \frac {b}{2}}  + 2\,n - a + 
{\displaystyle \frac {c}{2}}  + m) \\
\Gamma ( - 2\,n - b + 1 + 2\,a - c - 2\,m)\,\Gamma (2\,n - 2\,a
 + b + 2\,m)\,\Gamma (1 + 2\,a - 2\,n - 2\,m) \\
\Gamma ( - a + 2\,n + 2\,m)\,\Gamma ({\displaystyle \frac {c}{2}
}  + m + {\displaystyle \frac {b}{2}} ) \left/ {\vrule 
height0.80em width0em depth0.80em} \right. \!  \! (\pi \,\mathrm{
sin}({\displaystyle \frac {\pi \,(b + 4\,n - 2\,a + c + 2\,m)}{2}
} ) \\
\Gamma ( - a + c + 2\,n + 2\,m)\,\Gamma (2\,a - c - 2\,n - 2\,m
 + 1)\,\Gamma ( - {\displaystyle \frac {c}{2}}  + {\displaystyle 
\frac {b}{2}}  + m)) }
}
\end{maplelatex}

\begin{maplelatex}
\mapleinline{inert}{2d}{399, ":   ",
WPMR(2*m+2*n-2,2*n-1,2*a+2-b-2*n-2*m,-b+1,-c+2*n)*GAMMA(a+1-b+c-2*n)*G
AMMA(-2*a-1+b-2*c+4*n+2*m)*GAMMA(2*c-2*n)*GAMMA(-a-1+2*n+2*m)*GAMMA(1+
b-c)*GAMMA(a)*(sin(Pi*(-a+2*n+2*b))-sin(Pi*(-2*c-a+2*b+2*n))-sin(Pi*(2
*n+a))+sin(Pi*(-2*c+a+2*n)))/GAMMA(a-b)^2/GAMMA(-a-1+b-c+2*n+2*m)/GAMM
A(c-2*n+1)/GAMMA(2*n-1-2*a+b+2*m)/GAMMA(b)/(sin(Pi*a)+sin(Pi*(a-2*b))-
sin(Pi*(a-2*c+4*n))+sin(Pi*(-a+2*b-2*c+4*n)))/(a-b)+1/2*(cos(Pi*(-a+2*
b-c+2*n))-cos(Pi*(a-c+2*n)))*WPMR(2*m+2*n-2,2*n-1,2*a+2-b+2*c-4*n-2*m,
1+2*c-b-2*n,c)*GAMMA(a+1-b+c-2*n)*GAMMA(2*n-2*c)*GAMMA(-a-1+2*n+2*m)*G
AMMA(1+b-c)*GAMMA(b-a+1)*GAMMA(a)/Pi/sin(Pi*(-c+b+2*n))/GAMMA(a+1-b+2*
c-2*n)/GAMMA(b-2*c+2*n)/GAMMA(-a-1+b-c+2*n+2*m)/GAMMA(-c+1)+1/4*WPMR(2
*m+2*n-2,2*m-1,b,2*n-1-2*a+b+2*m,-a-1+b-c+2*n+2*m)*GAMMA(2*a+2-2*b+2*c
-4*n-2*m)*GAMMA(-2*a-1+b-2*c+4*n+2*m)*GAMMA(-a-1+2*n+2*m)*GAMMA(1+b-c)
*GAMMA(c)*GAMMA(b-a+1)*GAMMA(a)*(-sin(Pi*(b+a-c-2*m))+sin(Pi*(-b+a-c-2
*m))+sin(Pi*(3*b-3*a-c+4*n+2*m))-sin(Pi*(b-3*a-c+4*n+2*m)))/Pi^3/GAMMA
(-a-1+b+2*n+2*m)+GAMMA(-a-1+2*n+2*m)*GAMMA(1+b-c)*GAMMA(b-a+1)*GAMMA(a
)*GAMMA(c)/GAMMA(b)*WPMR(2*m+2*n-2,2*m-1,2*a+2-b+2*c-4*n-2*m,1+2*c-b-2
*n,a+1-b+c-2*n)*GAMMA(2*m-2*a-2+2*b-2*c+4*n)/GAMMA(-a+b-c+2*n)/GAMMA(b
-2*c+2*n)/GAMMA(2*n-1-2*a+b+2*m)/GAMMA(a+1-b+2*c-2*n)/GAMMA(-a-1+b-c+2
*n+2*m);}{%
\maplemultiline{
399\mbox{:~~~} \mathrm{WPMR}(2\,m + 2\,n - 2, \,2\,n - 1
, \,2\,a + 2 - b - 2\,n - 2\,m, \, - b + 1, \, - c + 2\,n) \\
\Gamma (a + 1 - b + c - 2\,n)\,\Gamma ( - 2\,a - 1 + b - 2\,c + 4
\,n + 2\,m)\,\Gamma (2\,c - 2\,n)\,\mathrm{\%2} \\
\Gamma (1 + b - c)\,\Gamma (a)(\mathrm{sin}(\pi \,( - a + 2\,n + 
2\,b)) - \mathrm{sin}(\pi \,( - 2\,c - a + 2\,b + 2\,n)) \\
\mbox{} - \mathrm{sin}(\pi \,(2\,n + a)) + \mathrm{sin}(\pi \,(
 - 2\,c + a + 2\,n)))/(\Gamma (a - b)^{2}\,\Gamma (\mathrm{\%1})
\,\Gamma (c - 2\,n + 1) \\
\Gamma (2\,n - 1 - 2\,a + b + 2\,m)\,\Gamma (b) \\
(\mathrm{sin}(\pi \,a) + \mathrm{sin}(\pi \,(a - 2\,b)) - 
\mathrm{sin}(\pi \,(a - 2\,c + 4\,n)) + \mathrm{sin}(\pi \,( - a
 + 2\,b - 2\,c + 4\,n))) \\
(a - b))\mbox{} + {\displaystyle \frac {1}{2}} (\mathrm{cos}(\pi 
\,( - a + 2\,b - c + 2\,n)) - \mathrm{cos}(\pi \,(a - c + 2\,n)))
 \\
\mathrm{WPMR}(2\,m + 2\,n - 2, \,2\,n - 1, \,2\,a + 2 - b + 2\,c
 - 4\,n - 2\,m, \,1 + 2\,c - b - 2\,n, \,c) \\
\Gamma (a + 1 - b + c - 2\,n)\,\Gamma (2\,n - 2\,c)\,\mathrm{\%2}
\,\Gamma (1 + b - c)\,\Gamma (b - a + 1)\,\Gamma (a)/(\pi  \\
\mathrm{sin}(\pi \,( - c + b + 2\,n))\,\Gamma (a + 1 - b + 2\,c
 - 2\,n)\,\Gamma (b - 2\,c + 2\,n)\,\Gamma (\mathrm{\%1})\,\Gamma
 ( - c + 1))\mbox{} +  \\
{\displaystyle \frac {1}{4}} \mathrm{WPMR}(2\,m + 2\,n - 2, \,2\,
m - 1, \,b, \,2\,n - 1 - 2\,a + b + 2\,m, \,\mathrm{\%1}) \\
\Gamma (2\,a + 2 - 2\,b + 2\,c - 4\,n - 2\,m)\,\Gamma ( - 2\,a - 
1 + b - 2\,c + 4\,n + 2\,m)\,\mathrm{\%2} \\
\Gamma (1 + b - c)\,\Gamma (c)\,\Gamma (b - a + 1)\,\Gamma (a)(
 - \mathrm{sin}(\pi \,(b + a - c - 2\,m)) \\
\mbox{} + \mathrm{sin}(\pi \,( - b + a - c - 2\,m)) + \mathrm{sin
}(\pi \,(3\,b - 3\,a - c + 4\,n + 2\,m)) \\
\mbox{} - \mathrm{sin}(\pi \,(b - 3\,a - c + 4\,n + 2\,m)))/(\pi 
^{3}\,\Gamma ( - a - 1 + b + 2\,n + 2\,m))\mbox{} + \mathrm{\%2}
 \\
\Gamma (1 + b - c)\,\Gamma (b - a + 1)\,\Gamma (a)\,\Gamma (c)
\mathrm{WPMR}(2\,m + 2\,n - 2, \,2\,m - 1,  \\
2\,a + 2 - b + 2\,c - 4\,n - 2\,m, \,1 + 2\,c - b - 2\,n, \,a + 1
 - b + c - 2\,n) \\
\Gamma (2\,m - 2\,a - 2 + 2\,b - 2\,c + 4\,n)/(\Gamma (b)\,\Gamma
 ( - a + b - c + 2\,n)\,\Gamma (b - 2\,c + 2\,n) \\
\Gamma (2\,n - 1 - 2\,a + b + 2\,m)\,\Gamma (a + 1 - b + 2\,c - 2
\,n)\,\Gamma (\mathrm{\%1})) \\
\mathrm{\%1} :=  - a - 1 + b - c + 2\,n + 2\,m \\
\mathrm{\%2} := \Gamma ( - a - 1 + 2\,n + 2\,m) }
}
\end{maplelatex}

\begin{maplelatex}
\mapleinline{inert}{2d}{400, ":   ",
1/2*(cos(Pi*(-a+2*b-c+2*n))-cos(Pi*(a-c+2*n)))*WPMR(2*m+2*n-1,2*n-1,2*
a-b+2*c-4*n+1-2*m,1+2*c-b-2*n,c)*GAMMA(a+1-b+c-2*n)*GAMMA(2*n-2*c)*GAM
MA(-a+2*n+2*m)*GAMMA(1+b-c)*GAMMA(b-a+1)*GAMMA(a)/Pi/sin(Pi*(-c+b+2*n)
)/GAMMA(a+1-b+2*c-2*n)/GAMMA(b-2*c+2*n)/GAMMA(-a+b-c+2*n+2*m)/GAMMA(-c
+1)+WPMR(2*m+2*n-1,2*n-1,2*a-b-2*n+1-2*m,-b+1,-c+2*n)*GAMMA(-2*a+b-2*c
+4*n+2*m)*GAMMA(2*c-2*n)*GAMMA(a+1-b+c-2*n)*GAMMA(-a+2*n+2*m)*GAMMA(1+
b-c)*GAMMA(a)*(sin(Pi*(-a+2*n+2*b))-sin(Pi*(-2*c-a+2*b+2*n))-sin(Pi*(2
*n+a))+sin(Pi*(-2*c+a+2*n)))/GAMMA(a-b)^2/GAMMA(2*n-2*a+b+2*m)/GAMMA(c
-2*n+1)/GAMMA(-a+b-c+2*n+2*m)/GAMMA(b)/(sin(Pi*a)+sin(Pi*(a-2*b))-sin(
Pi*(a-2*c+4*n))+sin(Pi*(-a+2*b-2*c+4*n)))/(a-b)-1/4*WPMR(2*m+2*n-1,2*m
,2*n-2*a+b+2*m,b,-a+b-c+2*n+2*m)*GAMMA(2*a+1-2*b+2*c-4*n-2*m)*GAMMA(-2
*a+b-2*c+4*n+2*m)*GAMMA(-a+2*n+2*m)*GAMMA(1+b-c)*GAMMA(c)*GAMMA(b-a+1)
*GAMMA(a)*(-sin(Pi*(b+a-c-2*m))+sin(Pi*(-b+a-c-2*m))+sin(Pi*(3*b-3*a-c
+4*n+2*m))-sin(Pi*(b-3*a-c+4*n+2*m)))/Pi^3/GAMMA(-a+b+2*n+2*m)+GAMMA(-
a+2*n+2*m)*GAMMA(1+b-c)*GAMMA(b-a+1)*GAMMA(a)*GAMMA(c)/GAMMA(b)*WPMR(2
*m+2*n-1,2*m,2*a-b+2*c-4*n+1-2*m,1+2*c-b-2*n,a+1-b+c-2*n)*GAMMA(-1-2*a
+2*b-2*c+4*n+2*m)/GAMMA(-a+b-c+2*n)/GAMMA(b-2*c+2*n)/GAMMA(a+1-b+2*c-2
*n)/GAMMA(2*n-2*a+b+2*m)/GAMMA(-a+b-c+2*n+2*m);}{%
\maplemultiline{
400\mbox{:~~~} {\displaystyle \frac {1}{2}} (\mathrm{cos
}(\pi \,( - a + 2\,b - c + 2\,n)) - \mathrm{cos}(\pi \,(a - c + 2
\,n))) \\
\mathrm{WPMR}(2\,m + 2\,n - 1, \,2\,n - 1, \,2\,a - b + 2\,c - 4
\,n + 1 - 2\,m, \,1 + 2\,c - b - 2\,n, \,c) \\
\Gamma (a + 1 - b + c - 2\,n)\,\Gamma (2\,n - 2\,c)\,\mathrm{\%2}
\,\Gamma (1 + b - c)\,\Gamma (b - a + 1)\,\Gamma (a)/(\pi  \\
\mathrm{sin}(\pi \,( - c + b + 2\,n))\,\Gamma (a + 1 - b + 2\,c
 - 2\,n)\,\Gamma (b - 2\,c + 2\,n)\,\Gamma (\mathrm{\%1})\,\Gamma
 ( - c + 1))\mbox{} +  \\
\mathrm{WPMR}(2\,m + 2\,n - 1, \,2\,n - 1, \,2\,a - b - 2\,n + 1
 - 2\,m, \, - b + 1, \, - c + 2\,n) \\
\Gamma ( - 2\,a + b - 2\,c + 4\,n + 2\,m)\,\Gamma (2\,c - 2\,n)\,
\Gamma (a + 1 - b + c - 2\,n)\,\mathrm{\%2}\,\Gamma (1 + b - c)
 \\
\Gamma (a)(\mathrm{sin}(\pi \,( - a + 2\,n + 2\,b)) - \mathrm{sin
}(\pi \,( - 2\,c - a + 2\,b + 2\,n)) - \mathrm{sin}(\pi \,(2\,n
 + a)) \\
\mbox{} + \mathrm{sin}(\pi \,( - 2\,c + a + 2\,n)))/(\Gamma (a - 
b)^{2}\,\Gamma (2\,n - 2\,a + b + 2\,m)\,\Gamma (c - 2\,n + 1)\,
\Gamma (\mathrm{\%1}) \\
\Gamma (b) \,
(\mathrm{sin}(\pi \,a) + \mathrm{sin}(\pi \,(a - 2\,b)) - 
\mathrm{sin}(\pi \,(a - 2\,c + 4\,n)) + \mathrm{sin}(\pi \,( - a
 + 2\,b - 2\,c + 4\,n))) \\
(a - b))\mbox{} - {\displaystyle \frac {1}{4}} \mathrm{WPMR}(2\,m
 + 2\,n - 1, \,2\,m, \,2\,n - 2\,a + b + 2\,m, \,b, \,\mathrm{\%1
}) \\
\Gamma (2\,a + 1 - 2\,b + 2\,c - 4\,n - 2\,m)\,\Gamma ( - 2\,a + 
b - 2\,c + 4\,n + 2\,m)\,\mathrm{\%2}\,\Gamma (1 + b - c) \\
\Gamma (c)\,\Gamma (b - a + 1)\,\Gamma (a)( - \mathrm{sin}(\pi \,
(b + a - c - 2\,m)) + \mathrm{sin}(\pi \,( - b + a - c - 2\,m))
 \\
\mbox{} + \mathrm{sin}(\pi \,(3\,b - 3\,a - c + 4\,n + 2\,m)) - 
\mathrm{sin}(\pi \,(b - 3\,a - c + 4\,n + 2\,m)))/(\pi ^{3} \\
\Gamma ( - a + b + 2\,n + 2\,m))\mbox{} + \mathrm{\%2}\,\Gamma (1
 + b - c)\,\Gamma (b - a + 1)\,\Gamma (a)\,\Gamma (c)\mathrm{WPMR
}( \\
2\,m + 2\,n - 1, \,2\,m, \,2\,a - b + 2\,c - 4\,n + 1 - 2\,m, \,1
 + 2\,c - b - 2\,n,  \\
a + 1 - b + c - 2\,n)\Gamma ( - 1 - 2\,a + 2\,b - 2\,c + 4\,n + 2
\,m)/(\Gamma (b)\,\Gamma ( - a + b - c + 2\,n) \\
\Gamma (b - 2\,c + 2\,n)\,\Gamma (a + 1 - b + 2\,c - 2\,n)\,
\Gamma (2\,n - 2\,a + b + 2\,m)\,\Gamma (\mathrm{\%1})) \\
\mathrm{\%1} :=  - a + b - c + 2\,n + 2\,m \\
\mathrm{\%2} := \Gamma ( - a + 2\,n + 2\,m) }
}
\end{maplelatex}

\begin{maplelatex}
\mapleinline{inert}{2d}{401, ":   ",
WPMR(2*m+2*n-1,2*n,-b+1,2*a-b-2*n+1-2*m,1-c+2*n)*GAMMA(a-b+c-2*n)*GAMM
A(-2*a+b-2*c+4*n+2*m+1)*GAMMA(2*c-2*n-1)*GAMMA(-a+2*n+2*m)*GAMMA(1+b-c
)*GAMMA(a)*(-sin(Pi*(-a+2*n+2*b))+sin(Pi*(-2*c-a+2*b+2*n))+sin(Pi*(2*n
+a))-sin(Pi*(-2*c+a+2*n)))/GAMMA(a-b)^2/GAMMA(c-2*n)/GAMMA(2*n-2*a+b+2
*m)/GAMMA(-a+b-c+2*n+2*m)/GAMMA(b)/(sin(Pi*a)+sin(Pi*(a-2*b))-sin(Pi*(
a-2*c+4*n))+sin(Pi*(-a+2*b-2*c+4*n)))/(a-b)+1/2*(cos(Pi*(-a+2*b-c+2*n)
)-cos(Pi*(a-c+2*n)))*WPMR(2*m+2*n-1,2*n,2*a-b+2*c-4*n-2*m,2*c-b-2*n,c)
*GAMMA(1-2*c+2*n)*GAMMA(a-b+c-2*n)*GAMMA(-a+2*n+2*m)*GAMMA(1+b-c)*GAMM
A(b-a+1)*GAMMA(a)/Pi/sin(Pi*(-c+b+2*n))/GAMMA(a-b+2*c-2*n)/GAMMA(b-2*c
+2*n+1)/GAMMA(-a+b-c+2*n+2*m)/GAMMA(-c+1)-1/4*WPMR(2*m+2*n-1,2*m-1,b,2
*n-2*a+b+2*m,-a+b-c+2*n+2*m)*GAMMA(2*a-2*b+2*c-4*n-2*m)*GAMMA(-2*a+b-2
*c+4*n+2*m+1)*GAMMA(-a+2*n+2*m)*GAMMA(1+b-c)*GAMMA(c)*GAMMA(b-a+1)*GAM
MA(a)*(sin(Pi*(b+a-c-2*m))-sin(Pi*(-b+a-c-2*m))-sin(Pi*(3*b-3*a-c+4*n+
2*m))+sin(Pi*(b-3*a-c+4*n+2*m)))/Pi^3/GAMMA(-a+b+2*n+2*m)+GAMMA(-a+2*n
+2*m)*GAMMA(1+b-c)*GAMMA(b-a+1)*GAMMA(a)*GAMMA(c)/GAMMA(b)*WPMR(2*m+2*
n-1,2*m-1,2*a-b+2*c-4*n-2*m,2*c-b-2*n,a-b+c-2*n)*GAMMA(2*m-2*a+2*b-2*c
+4*n)/GAMMA(-a+b-c+2*n+1)/GAMMA(b-2*c+2*n+1)/GAMMA(2*n-2*a+b+2*m)/GAMM
A(a-b+2*c-2*n)/GAMMA(-a+b-c+2*n+2*m);}{%
\maplemultiline{
401\mbox{:~~~} \mathrm{WPMR}(2\,m + 2\,n - 1, \,2\,n, \,
 - b + 1, \,2\,a - b - 2\,n + 1 - 2\,m, \,1 - c + 2\,n) \\
\Gamma (a - b + c - 2\,n)\,\Gamma ( - 2\,a + b - 2\,c + 4\,n + 2
\,m + 1)\,\Gamma (2\,c - 2\,n - 1)\,\mathrm{\%2} \\
\Gamma (1 + b - c)\,\Gamma (a)( - \mathrm{sin}(\pi \,( - a + 2\,n
 + 2\,b)) + \mathrm{sin}(\pi \,( - 2\,c - a + 2\,b + 2\,n)) \\
\mbox{} + \mathrm{sin}(\pi \,(2\,n + a)) - \mathrm{sin}(\pi \,(
 - 2\,c + a + 2\,n)))/(\Gamma (a - b)^{2}\,\Gamma (c - 2\,n) \\
\Gamma (2\,n - 2\,a + b + 2\,m)\,\Gamma (\mathrm{\%1})\,\Gamma (b
) \\
(\mathrm{sin}(\pi \,a) + \mathrm{sin}(\pi \,(a - 2\,b)) - 
\mathrm{sin}(\pi \,(a - 2\,c + 4\,n)) + \mathrm{sin}(\pi \,( - a
 + 2\,b - 2\,c + 4\,n))) \\
(a - b))\mbox{} + {\displaystyle \frac {1}{2}} (\mathrm{cos}(\pi 
\,( - a + 2\,b - c + 2\,n)) - \mathrm{cos}(\pi \,(a - c + 2\,n)))
 \\
\mathrm{WPMR}(2\,m + 2\,n - 1, \,2\,n, \,2\,a - b + 2\,c - 4\,n
 - 2\,m, \,2\,c - b - 2\,n, \,c) \\
\Gamma (1 - 2\,c + 2\,n)\,\Gamma (a - b + c - 2\,n)\,\mathrm{\%2}
\,\Gamma (1 + b - c)\,\Gamma (b - a + 1)\,\Gamma (a)/(\pi  \\
\mathrm{sin}(\pi \,( - c + b + 2\,n))\,\Gamma (a - b + 2\,c - 2\,
n)\,\Gamma (b - 2\,c + 2\,n + 1)\,\Gamma (\mathrm{\%1})\,\Gamma (
 - c + 1))\mbox{} -  \\
{\displaystyle \frac {1}{4}} \mathrm{WPMR}(2\,m + 2\,n - 1, \,2\,
m - 1, \,b, \,2\,n - 2\,a + b + 2\,m, \,\mathrm{\%1}) \\
\Gamma (2\,a - 2\,b + 2\,c - 4\,n - 2\,m)\,\Gamma ( - 2\,a + b - 
2\,c + 4\,n + 2\,m + 1)\,\mathrm{\%2}\,\Gamma (1 + b - c) \\
\Gamma (c)\,\Gamma (b - a + 1)\,\Gamma (a)(\mathrm{sin}(\pi \,(b
 + a - c - 2\,m)) - \mathrm{sin}(\pi \,( - b + a - c - 2\,m)) \\
\mbox{} - \mathrm{sin}(\pi \,(3\,b - 3\,a - c + 4\,n + 2\,m)) + 
\mathrm{sin}(\pi \,(b - 3\,a - c + 4\,n + 2\,m)))/(\pi ^{3} \\
\Gamma ( - a + b + 2\,n + 2\,m))\mbox{} + \mathrm{\%2}\,\Gamma (1
 + b - c)\,\Gamma (b - a + 1)\,\Gamma (a)\,\Gamma (c)\mathrm{WPMR
}( \\
2\,m + 2\,n - 1, \,2\,m - 1, \,2\,a - b + 2\,c - 4\,n - 2\,m, \,2
\,c - b - 2\,n, \,a - b + c - 2\,n) \\
\Gamma (2\,m - 2\,a + 2\,b - 2\,c + 4\,n)/(\Gamma (b)\,\Gamma (
 - a + b - c + 2\,n + 1)\,\Gamma (b - 2\,c + 2\,n + 1) \\
\Gamma (2\,n - 2\,a + b + 2\,m)\,\Gamma (a - b + 2\,c - 2\,n)\,
\Gamma (\mathrm{\%1})) \\
\mathrm{\%1} :=  - a + b - c + 2\,n + 2\,m \\
\mathrm{\%2} := \Gamma ( - a + 2\,n + 2\,m) }
}
\end{maplelatex}

\begin{maplelatex}
\mapleinline{inert}{2d}{402, ":   ",
2*WPMR(2*m+2*n-1,2*n,1-a-c+2*b-2*n,-c+1-a+2*m,1+b-c)*Pi^2*GAMMA(2*c-1-
2*b+2*n)*GAMMA(1+a-b)*GAMMA(1+2*a-2*m)/GAMMA(a+c-2*b+2*n)/GAMMA(1-c+b-
a+2*m)/GAMMA(a-2*m+c)/GAMMA(1+a-2*m)/GAMMA(c-b)/GAMMA(-b+1)/GAMMA(a)/(
cos(Pi*(-c+b+2*m))-cos(Pi*(c-3*b+4*n+2*m)))+2*Pi*(cos(Pi*(2*c-3*b+4*n+
a))-cos(Pi*(a-b)))*WPMR(2*m+2*n-1,2*n,-a+2*n-2*b+2*m+c,c-a,c-b+2*n)*GA
MMA(1-2*c+2*b-2*n)*GAMMA(c-b+2*n)*GAMMA(1+a-b)*GAMMA(1+2*a-2*m)/GAMMA(
-a+2*n-b+2*m+c)/GAMMA(a-2*n+2*b-2*m-c+1)/GAMMA(1+a-2*m)/GAMMA(1+a-c)/G
AMMA(-b+1)/GAMMA(a)/(sin(Pi*(-3*b+4*n+2*m-a+2*c))+sin(Pi*(b+2*m+a-2*c)
)-sin(Pi*(-3*b+4*n+2*m+a))+sin(Pi*(-b+a-2*m)))+GAMMA(1+2*a-2*m)*GAMMA(
1+a-b)*GAMMA(-c+2*b-2*n)*WMMR(2*m-1,2*n,a+c-2*b+2*n,-a+2*n-2*b+2*m+c,c
-b+2*n)*GAMMA(c-a)/GAMMA(a)/GAMMA(a-2*n+2*b-2*m-c+1)/GAMMA(-c+b-2*n+1)
/GAMMA(2*c-2*b+2*n);}{%
\maplemultiline{
402\mbox{:~~~} 2\,\mathrm{WPMR}(2\,m + 2\,n - 1, \,2\,n
, \,1 - a - c + 2\,b - 2\,n, \, - c + 1 - a + 2\,m, \,1 + b - c)
 \\
\pi ^{2}\,\Gamma (2\,c - 1 - 2\,b + 2\,n)\,\Gamma (1 + a - b)\,
\Gamma (1 + 2\,a - 2\,m)/(\Gamma (a + c - 2\,b + 2\,n) \\
\Gamma (1 - c + b - a + 2\,m)\,\Gamma (a - 2\,m + c)\,\Gamma (1
 + a - 2\,m)\,\Gamma (c - b)\,\Gamma ( - b + 1)\,\Gamma (a) \\
(\mathrm{cos}(\pi \,( - c + b + 2\,m)) - \mathrm{cos}(\pi \,(c - 
3\,b + 4\,n + 2\,m))))\mbox{} + 2\,\pi  \\
(\mathrm{cos}(\pi \,(2\,c - 3\,b + 4\,n + a)) - \mathrm{cos}(\pi 
\,(a - b))) \\
\mathrm{WPMR}(2\,m + 2\,n - 1, \,2\,n, \, - a + 2\,n - 2\,b + 2\,
m + c, \,c - a, \,c - b + 2\,n) \\
\Gamma (1 - 2\,c + 2\,b - 2\,n)\,\Gamma (c - b + 2\,n)\,\Gamma (1
 + a - b)\,\Gamma (1 + 2\,a - 2\,m)/( \\
\Gamma ( - a + 2\,n - b + 2\,m + c)\,\Gamma (a - 2\,n + 2\,b - 2
\,m - c + 1)\,\Gamma (1 + a - 2\,m) \\
\Gamma (1 + a - c)\,\Gamma ( - b + 1)\,\Gamma (a)(\mathrm{sin}(
\pi \,( - 3\,b + 4\,n + 2\,m - a + 2\,c)) \\
\mbox{} + \mathrm{sin}(\pi \,(b + 2\,m + a - 2\,c)) - \mathrm{sin
}(\pi \,( - 3\,b + 4\,n + 2\,m + a)) \\
\mbox{} + \mathrm{sin}(\pi \,( - b + a - 2\,m))))\mbox{} + \Gamma
 (1 + 2\,a - 2\,m)\,\Gamma (1 + a - b)\,\Gamma ( - c + 2\,b - 2\,
n) \\
\mathrm{WMMR}(2\,m - 1, \,2\,n, \,a + c - 2\,b + 2\,n, \, - a + 2
\,n - 2\,b + 2\,m + c, \,c - b + 2\,n) \\
\Gamma (c - a)/(\Gamma (a)\,\Gamma (a - 2\,n + 2\,b - 2\,m - c + 
1)\,\Gamma ( - c + b - 2\,n + 1) \\
\Gamma (2\,c - 2\,b + 2\,n)) }
}
\end{maplelatex}

\begin{maplelatex}
\mapleinline{inert}{2d}{403, ":   ",
4*Pi^2*(cos(Pi*(c-b+2*n))-cos(Pi*(-3*c+a+4*m)))*WPMR(2*m+2*n-2,2*n-1,-
a+1,2*m-2*c+b,-1/2*a+1/2*b+n)*GAMMA(1+2*c-2*m)*GAMMA(1/2*a+1/2*b-n+1)*
GAMMA(a-b)/GAMMA(-c+1/2*b+1/2*a-n+1)/GAMMA(1/2*a-1/2*b-n+1)/GAMMA(-c+1
/2*b-1/2*a+n+2*m)/GAMMA(c+1-2*m)/GAMMA(1+2*c-2*m-b)/GAMMA(c)/GAMMA(a)/
(cos(1/2*Pi*(-3*b-4*m-3*a+6*c+2*n))-cos(1/2*Pi*(-3*b-12*m-3*a+10*c+2*n
))-cos(1/2*Pi*(b-4*m+a+2*c+2*n))-cos(1/2*Pi*(6*n+4*m+2*c-b-a))+cos(1/2
*Pi*(2*n+12*m+a+b-6*c))+cos(1/2*Pi*(6*n-4*m-a-b+6*c)))+WPMR(2*m+2*n-2,
2*n-1,-b+1,-2*c+a+2*m,n+1/2*a-1/2*b)*GAMMA(1/2*a+1/2*b-n+1)*GAMMA(n+1/
2*a-1/2*b)*GAMMA(1+2*c-2*m)*GAMMA(-2*c+a+2*m)*GAMMA(b-a)*(-cos(1/2*Pi*
(2*n-3*b-2*c+3*a+4*m))-cos(1/2*Pi*(6*n-b+6*c-3*a-4*m))+cos(1/2*Pi*(2*n
-3*b+6*c-a-4*m))+cos(1/2*Pi*(6*n-b-2*c+a+4*m))+cos(1/2*Pi*(10*c-5*a-12
*m+b+2*n))+cos(1/2*Pi*(-2*c+b-a+2*n+4*m))-cos(1/2*Pi*(2*c-a-4*m+b+2*n)
)-cos(1/2*Pi*(-10*c+3*a+12*m+b+2*n)))/GAMMA(-c-1/2*b+1/2*a+n+2*m)/GAMM
A(-c+1/2*b+1/2*a-n+1)/GAMMA(c+1-2*m)/GAMMA(c)/GAMMA(b)/(-cos(1/2*Pi*(-
3*b-4*m-3*a+6*c+2*n))+cos(1/2*Pi*(-3*b-12*m-3*a+10*c+2*n))+cos(1/2*Pi*
(b-4*m+a+2*c+2*n))+cos(1/2*Pi*(6*n+4*m+2*c-b-a))-cos(1/2*Pi*(2*n+12*m+
a+b-6*c))-cos(1/2*Pi*(6*n-4*m-a-b+6*c)))-2*WPMR(2*m+2*n-2,2*m-1,-a+1,2
*m-2*c+b,-c+2*m)*(c-m)*Pi/sin(1/2*Pi*(-2*c+b+a+2*n+4*m))*GAMMA(-1/2*a+
1/2*b+n)*GAMMA(c-1/2*b-1/2*a+n)*GAMMA(1/2*a+1/2*b-n+1)/GAMMA(-c+1/2*b-
1/2*a+n+2*m)*GAMMA(-2*m+2*c)^2/GAMMA(c+1-2*m)/GAMMA(1+2*c-a-2*m)/GAMMA
(1+2*c-2*m-b)/GAMMA(c)/GAMMA(b)/GAMMA(a)-1/2*WPMR(2*m+2*n-2,2*m-1,1+2*
c-a-2*m,b,c)*GAMMA(-1/2*a+1/2*b+n)*GAMMA(c-1/2*b-1/2*a+n)*GAMMA(1/2*a+
1/2*b-n+1)*(sin(Pi*(-3*c+a+4*m+b))+sin(Pi*(-c+a-b))-sin(Pi*(b-c+a))-si
n(Pi*(-3*c+a+4*m-b)))/Pi/GAMMA(c+1/2*b-1/2*a+n)/(cos(1/2*Pi*(-6*c+b+a+
2*n+8*m))-cos(1/2*Pi*(2*c+b+a+2*n)));}{%
\maplemultiline{
403\mbox{:~~~} 4\,\pi ^{2}\,(\mathrm{cos}(\pi \,(c - b
 + 2\,n)) - \mathrm{cos}(\pi \,( - 3\,c + a + 4\,m))) \\
\mathrm{WPMR}(2\,m + 2\,n - 2, \,2\,n - 1, \, - a + 1, \,2\,m - 2
\,c + b, \, - {\displaystyle \frac {a}{2}}  + {\displaystyle 
\frac {b}{2}}  + n) \\
\Gamma (1 + 2\,c - 2\,m)\,\mathrm{\%1}\,\Gamma (a - b) \left/ 
{\vrule height0.80em width0em depth0.80em} \right. \!  \! (\Gamma
 ( - c + {\displaystyle \frac {b}{2}}  + {\displaystyle \frac {a
}{2}}  - n + 1)\,\Gamma ({\displaystyle \frac {a}{2}}  - 
{\displaystyle \frac {b}{2}}  - n + 1) \\
\Gamma ( - c + {\displaystyle \frac {b}{2}}  - {\displaystyle 
\frac {a}{2}}  + n + 2\,m)\,\Gamma (c + 1 - 2\,m)\,\Gamma (1 + 2
\,c - 2\,m - b)\,\Gamma (c)\,\Gamma (a)( \\
\mathrm{cos}({\displaystyle \frac {\pi \,( - 3\,b - 4\,m - 3\,a
 + 6\,c + 2\,n)}{2}} ) - \mathrm{cos}({\displaystyle \frac {\pi 
\,( - 3\,b - 12\,m - 3\,a + 10\,c + 2\,n)}{2}} ) \\
\mbox{} - \mathrm{cos}({\displaystyle \frac {\pi \,(b - 4\,m + a
 + 2\,c + 2\,n)}{2}} ) - \mathrm{cos}({\displaystyle \frac {\pi 
\,(6\,n + 4\,m + 2\,c - b - a)}{2}} ) \\
\mbox{} + \mathrm{cos}({\displaystyle \frac {\pi \,(2\,n + 12\,m
 + a + b - 6\,c)}{2}} ) + \mathrm{cos}({\displaystyle \frac {\pi 
\,(6\,n - 4\,m - a - b + 6\,c)}{2}} )))\mbox{} +  \\
\mathrm{WPMR}(2\,m + 2\,n - 2, \,2\,n - 1, \, - b + 1, \, - 2\,c
 + a + 2\,m, \,n + {\displaystyle \frac {a}{2}}  - 
{\displaystyle \frac {b}{2}} )\,\mathrm{\%1} \\
\Gamma (n + {\displaystyle \frac {a}{2}}  - {\displaystyle 
\frac {b}{2}} )\,\Gamma (1 + 2\,c - 2\,m)\,\Gamma ( - 2\,c + a + 
2\,m)\,\Gamma (b - a)( \\
 - \mathrm{cos}({\displaystyle \frac {\pi \,(2\,n - 3\,b - 2\,c
 + 3\,a + 4\,m)}{2}} ) - \mathrm{cos}({\displaystyle \frac {\pi 
\,(6\,n - b + 6\,c - 3\,a - 4\,m)}{2}} ) \\
\mbox{} + \mathrm{cos}({\displaystyle \frac {\pi \,(2\,n - 3\,b
 + 6\,c - a - 4\,m)}{2}} ) + \mathrm{cos}({\displaystyle \frac {
\pi \,(6\,n - b - 2\,c + a + 4\,m)}{2}} ) \\
\mbox{} + \mathrm{cos}({\displaystyle \frac {\pi \,(10\,c - 5\,a
 - 12\,m + b + 2\,n)}{2}} ) + \mathrm{cos}({\displaystyle \frac {
\pi \,( - 2\,c + b - a + 2\,n + 4\,m)}{2}} ) \\
\mbox{} - \mathrm{cos}({\displaystyle \frac {\pi \,(2\,c - a - 4
\,m + b + 2\,n)}{2}} ) - \mathrm{cos}({\displaystyle \frac {\pi 
\,( - 10\,c + 3\,a + 12\,m + b + 2\,n)}{2}} )) \\
 \left/ {\vrule height0.80em width0em depth0.80em} \right. \! 
 \! (\Gamma ( - c - {\displaystyle \frac {b}{2}}  + 
{\displaystyle \frac {a}{2}}  + n + 2\,m)\,\Gamma ( - c + 
{\displaystyle \frac {b}{2}}  + {\displaystyle \frac {a}{2}}  - n
 + 1)\,\Gamma (c + 1 - 2\,m)\,\Gamma (c)\,\Gamma (b)( \\
 - \mathrm{cos}({\displaystyle \frac {\pi \,( - 3\,b - 4\,m - 3\,
a + 6\,c + 2\,n)}{2}} ) + \mathrm{cos}({\displaystyle \frac {\pi 
\,( - 3\,b - 12\,m - 3\,a + 10\,c + 2\,n)}{2}} ) \\
\mbox{} + \mathrm{cos}({\displaystyle \frac {\pi \,(b - 4\,m + a
 + 2\,c + 2\,n)}{2}} ) + \mathrm{cos}({\displaystyle \frac {\pi 
\,(6\,n + 4\,m + 2\,c - b - a)}{2}} ) \\
\mbox{} - \mathrm{cos}({\displaystyle \frac {\pi \,(2\,n + 12\,m
 + a + b - 6\,c)}{2}} ) - \mathrm{cos}({\displaystyle \frac {\pi 
\,(6\,n - 4\,m - a - b + 6\,c)}{2}} )))\mbox{} - 2 \\
\mathrm{WPMR}(2\,m + 2\,n - 2, \,2\,m - 1, \, - a + 1, \,2\,m - 2
\,c + b, \, - c + 2\,m)\,(c - m)\,\pi  \\
\Gamma ( - {\displaystyle \frac {a}{2}}  + {\displaystyle \frac {
b}{2}}  + n)\,\Gamma (c - {\displaystyle \frac {b}{2}}  - 
{\displaystyle \frac {a}{2}}  + n)\,\mathrm{\%1}\,\Gamma ( - 2\,m
 + 2\,c)^{2} \left/ {\vrule height0.80em width0em depth0.80em}
 \right. \!  \! ( \\
\mathrm{sin}({\displaystyle \frac {\pi \,( - 2\,c + b + a + 2\,n
 + 4\,m)}{2}} )\,\Gamma ( - c + {\displaystyle \frac {b}{2}}  - 
{\displaystyle \frac {a}{2}}  + n + 2\,m)\,\Gamma (c + 1 - 2\,m)
 \\
\Gamma (1 + 2\,c - a - 2\,m)\,\Gamma (1 + 2\,c - 2\,m - b)\,
\Gamma (c)\,\Gamma (b)\,\Gamma (a))\mbox{} - {\displaystyle 
\frac {1}{2}}  \\
\mathrm{WPMR}(2\,m + 2\,n - 2, \,2\,m - 1, \,1 + 2\,c - a - 2\,m
, \,b, \,c)\,\Gamma ( - {\displaystyle \frac {a}{2}}  + 
{\displaystyle \frac {b}{2}}  + n) \\
\Gamma (c - {\displaystyle \frac {b}{2}}  - {\displaystyle 
\frac {a}{2}}  + n)\,\mathrm{\%1}(\mathrm{sin}(\pi \,( - 3\,c + a
 + 4\,m + b)) + \mathrm{sin}(\pi \,( - c + a - b)) \\
\mbox{} - \mathrm{sin}(\pi \,(b - c + a)) - \mathrm{sin}(\pi \,(
 - 3\,c + a + 4\,m - b))) \left/ {\vrule 
height0.80em width0em depth0.80em} \right. \!  \! (\pi \,\Gamma (
c + {\displaystyle \frac {b}{2}}  - {\displaystyle \frac {a}{2}} 
 + n) \\
(\mathrm{cos}({\displaystyle \frac {\pi \,( - 6\,c + b + a + 2\,n
 + 8\,m)}{2}} ) - \mathrm{cos}({\displaystyle \frac {\pi \,(2\,c
 + b + a + 2\,n)}{2}} ))) \\
\mathrm{\%1} := \Gamma ({\displaystyle \frac {a}{2}}  + 
{\displaystyle \frac {b}{2}}  - n + 1) }
}
\end{maplelatex}

\begin{maplelatex}
\mapleinline{inert}{2d}{404, ":   ",
-WPMR(2*m+2*n-2,2*m-1,-a+1,-c+1,-b+1)*cos(1/2*Pi*(-a-c+2*n+2*m))*GAMMA
(1/2*a+1/2*c+3/2-n-m)^2*GAMMA(-1/2*a+b-1/2*c+1/2-n+m)*GAMMA(-1/2*c+1/2
*a-1/2+n+m)*GAMMA(-1/2*a+b-1/2*c-1/2+n+m)*GAMMA(1/2*c-1/2*a-1/2+n+m)*G
AMMA(2*b-1+2*m)^2/(-a-c-1+2*n+2*m)/Pi/(-1+b+m)/GAMMA(-1+2*m+b)/GAMMA(c
)/GAMMA(b)/GAMMA(a)/GAMMA(-a+2*b-1+2*m)/GAMMA(-c+2*b-1+2*m)+1/4*WPMR(2
*m+2*n-2,2*m-1,-a+2*b-1+2*m,-c+2*b-1+2*m,-1+2*m+b)*GAMMA(1/2*a+1/2*c+3
/2-n-m)*GAMMA(-1/2*a+b-1/2*c+1/2-n+m)*GAMMA(-1/2*a+b-1/2*c-1/2+n+m)*GA
MMA(1/2*c-1/2*a-1/2+n+m)*GAMMA(-1/2*c+1/2*a-1/2+n+m)*(sin(Pi*(-c+a-b))
+sin(Pi*(a-3*b-4*m+c))-sin(Pi*(b-c+a))-sin(Pi*(a-5*b-4*m+c)))/sin(2*Pi
*(b+m))/Pi^2/GAMMA(-1/2*c-1/2*a-3/2+n+3*m+2*b)-4*WPMR(2*m+2*n-2,2*n-1,
-a+1,-c+1,-1/2*a+b-1/2*c-1/2+n+m)*Pi*GAMMA(1/2*a+1/2*c+3/2-n-m)^2*GAMM
A(-1/2*c+1/2*a-1/2+n+m)*GAMMA(2*b-1+2*m)*GAMMA(a-2*b+c+1-2*m)*(cos(1/2
*Pi*(-3*a-c+2*n+6*m+6*b))+cos(1/2*Pi*(a-c+2*n-2*m-6*b))+cos(1/2*Pi*(a+
3*c+2*n-2*m-2*b))+cos(1/2*Pi*(-a+c+6*n+2*m-2*b)))/GAMMA(1/2*a-b+1/2*c+
3/2-n-m)/GAMMA(-1/2*c+1/2*a+3/2-n-m)/GAMMA(1/2*a-b+1/2*c+1/2+n-m)/GAMM
A(-1+2*m+b)/GAMMA(c)/GAMMA(b)/GAMMA(a)/(-sin(1/2*Pi*(2*b-a+3*c+6*n+2*m
))-sin(1/2*Pi*(6*b-3*a+c+2*n+6*m))+sin(1/2*Pi*(-6*b+a+c+2*n-2*m))-sin(
1/2*Pi*(2*b+a+c+2*n-2*m))+sin(1/2*Pi*(c-3*a+2*n+6*m+2*b))+sin(1/2*Pi*(
3*c-a+6*n+2*m-2*b)))/(-a-c-1+2*n+2*m)-WPMR(2*m+2*n-2,2*n-1,c+2-2*b-2*m
,a-2*b+2-2*m,1/2*a-b+1/2*c+1/2+n-m)*GAMMA(1/2*a+1/2*c+3/2-n-m)*GAMMA(-
1/2*c+1/2*a-1/2+n+m)*GAMMA(2*b-1+2*m)*GAMMA(a-2*b+2-2*m)*GAMMA(-a+2*b-
c-1+2*m)*(-sin(1/2*Pi*(3*a-8*b-6*m-c+2*n))-sin(1/2*Pi*(-a+4*b+2*m-c+2*
n))+sin(1/2*Pi*(-a-c+2*n+2*m))+sin(1/2*Pi*(-5*a+12*b+10*m-c+2*n))-sin(
1/2*Pi*(6*n-4*b+c+a-2*m))+sin(1/2*Pi*(2*n+3*c-a+2*m))+sin(1/2*Pi*(6*n+
4*b+c-3*a+6*m))-sin(1/2*Pi*(2*n-8*b+3*c+3*a-6*m)))/GAMMA(3/2+1/2*c+1/2
*a+n-m-2*b)/GAMMA(-1/2*c+1/2*a+3/2-n-m)/GAMMA(-1+2*m+b)/GAMMA(b)/GAMMA
(-c+2*b-1+2*m)/(sin(1/2*Pi*(2*b-a+3*c+6*n+2*m))+sin(1/2*Pi*(6*b-3*a+c+
2*n+6*m))-sin(1/2*Pi*(-6*b+a+c+2*n-2*m))+sin(1/2*Pi*(2*b+a+c+2*n-2*m))
-sin(1/2*Pi*(c-3*a+2*n+6*m+2*b))-sin(1/2*Pi*(3*c-a+6*n+2*m-2*b)));}{%
\maplemultiline{
404\mbox{:~~~}  - \mathrm{WPMR}(2\,m + 2\,n - 2, \,2\,m
 - 1, \, - a + 1, \, - c + 1, \, - b + 1) \\
\mathrm{cos}({\displaystyle \frac {\pi \,( - a - c + 2\,n + 2\,m)
}{2}} )\,\mathrm{\%2}^{2}\,\Gamma ( - {\displaystyle \frac {a}{2}
}  + b - {\displaystyle \frac {c}{2}}  + {\displaystyle \frac {1
}{2}}  - n + m)\,\mathrm{\%1} \\
\Gamma ( - {\displaystyle \frac {a}{2}}  + b - {\displaystyle 
\frac {c}{2}}  - {\displaystyle \frac {1}{2}}  + n + m)\,\Gamma (
{\displaystyle \frac {c}{2}}  - {\displaystyle \frac {a}{2}}  - 
{\displaystyle \frac {1}{2}}  + n + m)\,\Gamma (2\,b - 1 + 2\,m)
^{2}/( \\
( - a - c - 1 + 2\,n + 2\,m)\,\pi \,( - 1 + b + m)\,\Gamma ( - 1
 + 2\,m + b)\,\Gamma (c)\,\Gamma (b)\,\Gamma (a) \\
\Gamma ( - a + 2\,b - 1 + 2\,m)\,\Gamma ( - c + 2\,b - 1 + 2\,m))
\mbox{} + {\displaystyle \frac {1}{4}} \mathrm{WPMR}(2\,m + 2\,n
 - 2, \,2\,m - 1,  \\
 - a + 2\,b - 1 + 2\,m, \, - c + 2\,b - 1 + 2\,m, \, - 1 + 2\,m
 + b)\mathrm{\%2} \\
\Gamma ( - {\displaystyle \frac {a}{2}}  + b - {\displaystyle 
\frac {c}{2}}  + {\displaystyle \frac {1}{2}}  - n + m)\,\Gamma (
 - {\displaystyle \frac {a}{2}}  + b - {\displaystyle \frac {c}{2
}}  - {\displaystyle \frac {1}{2}}  + n + m)\,\Gamma (
{\displaystyle \frac {c}{2}}  - {\displaystyle \frac {a}{2}}  - 
{\displaystyle \frac {1}{2}}  + n + m)\,\mathrm{\%1}( \\
\mathrm{sin}(\pi \,( - c + a - b)) + \mathrm{sin}(\pi \,(a - 3\,b
 - 4\,m + c)) - \mathrm{sin}(\pi \,(b - c + a)) \\
\mbox{} - \mathrm{sin}(\pi \,(a - 5\,b - 4\,m + c))) \left/ 
{\vrule height0.80em width0em depth0.80em} \right. \!  \! (
\mathrm{sin}(2\,\pi \,(b + m))\,\pi ^{2}
\Gamma ( - {\displaystyle \frac {c}{2}}  - {\displaystyle \frac {
a}{2}}  - {\displaystyle \frac {3}{2}}  + n + 3\,m + 2\,b))
\mbox{} \\ - 4 \,
\mathrm{WPMR}(2\,m + 2\,n - 2, \,2\,n - 1, \, - a + 1, \, - c + 1
, \, - {\displaystyle \frac {a}{2}}  + b - {\displaystyle \frac {
c}{2}}  - {\displaystyle \frac {1}{2}}  + n + m)\,\pi \,\mathrm{
\%2}^{2}\,\mathrm{\%1} \\
\Gamma (2\,b - 1 + 2\,m)\,\Gamma (a - 2\,b + c + 1 - 2\,m)(
\mathrm{cos}({\displaystyle \frac {\pi \,( - 3\,a - c + 2\,n + 6
\,m + 6\,b)}{2}} ) \\
\mbox{} + \mathrm{cos}({\displaystyle \frac {\pi \,(a - c + 2\,n
 - 2\,m - 6\,b)}{2}} ) + \mathrm{cos}({\displaystyle \frac {\pi 
\,(a + 3\,c + 2\,n - 2\,m - 2\,b)}{2}} ) \\
\mbox{} + \mathrm{cos}({\displaystyle \frac {\pi \,( - a + c + 6
\,n + 2\,m - 2\,b)}{2}} )) \left/ {\vrule 
height0.80em width0em depth0.80em} \right. \!  \! (\Gamma (
{\displaystyle \frac {a}{2}}  - b + {\displaystyle \frac {c}{2}} 
 + {\displaystyle \frac {3}{2}}  - n - m) \\
\Gamma ( - {\displaystyle \frac {c}{2}}  + {\displaystyle \frac {
a}{2}}  + {\displaystyle \frac {3}{2}}  - n - m)\,\Gamma (
{\displaystyle \frac {a}{2}}  - b + {\displaystyle \frac {c}{2}} 
 + {\displaystyle \frac {1}{2}}  + n - m)\,\Gamma ( - 1 + 2\,m + 
b)\,\Gamma (c)\,\Gamma (b)\,\Gamma (a)( \\
 - \mathrm{sin}({\displaystyle \frac {\pi \,(2\,b - a + 3\,c + 6
\,n + 2\,m)}{2}} ) - \mathrm{sin}({\displaystyle \frac {\pi \,(6
\,b - 3\,a + c + 2\,n + 6\,m)}{2}} ) \\
\mbox{} + \mathrm{sin}({\displaystyle \frac {\pi \,( - 6\,b + a
 + c + 2\,n - 2\,m)}{2}} ) - \mathrm{sin}({\displaystyle \frac {
\pi \,(2\,b + a + c + 2\,n - 2\,m)}{2}} ) \\
\mbox{} + \mathrm{sin}({\displaystyle \frac {\pi \,(c - 3\,a + 2
\,n + 6\,m + 2\,b)}{2}} ) + \mathrm{sin}({\displaystyle \frac {
\pi \,(3\,c - a + 6\,n + 2\,m - 2\,b)}{2}} )) \\
( - a - c - 1 + 2\,n + 2\,m))\mbox{} - \mathrm{WPMR}(2\,m + 2\,n
 - 2, \,2\,n - 1, \,c + 2 - 2\,b - 2\,m,  \\
a - 2\,b + 2 - 2\,m, \,{\displaystyle \frac {a}{2}}  - b + 
{\displaystyle \frac {c}{2}}  + {\displaystyle \frac {1}{2}}  + n
 - m)\mathrm{\%2}\,\mathrm{\%1}\,\Gamma (2\,b - 1 + 2\,m) \\
\Gamma (a - 2\,b + 2 - 2\,m)\,\Gamma ( - a + 2\,b - c - 1 + 2\,m)
( \\
 - \mathrm{sin}({\displaystyle \frac {\pi \,(3\,a - 8\,b - 6\,m
 - c + 2\,n)}{2}} ) - \mathrm{sin}({\displaystyle \frac {\pi \,(
 - a + 4\,b + 2\,m - c + 2\,n)}{2}} ) \\
\mbox{} + \mathrm{sin}({\displaystyle \frac {\pi \,( - a - c + 2
\,n + 2\,m)}{2}} ) + \mathrm{sin}({\displaystyle \frac {\pi \,(
 - 5\,a + 12\,b + 10\,m - c + 2\,n)}{2}} ) \\
\mbox{} - \mathrm{sin}({\displaystyle \frac {\pi \,(6\,n - 4\,b
 + c + a - 2\,m)}{2}} ) + \mathrm{sin}({\displaystyle \frac {\pi 
\,(2\,n + 3\,c - a + 2\,m)}{2}} ) \\
\mbox{} + \mathrm{sin}({\displaystyle \frac {\pi \,(6\,n + 4\,b
 + c - 3\,a + 6\,m)}{2}} ) - \mathrm{sin}({\displaystyle \frac {
\pi \,(2\,n - 8\,b + 3\,c + 3\,a - 6\,m)}{2}} )) \left/ {\vrule 
height0.80em width0em depth0.80em} \right. \!  \!  \\
(\Gamma ({\displaystyle \frac {3}{2}}  + {\displaystyle \frac {c
}{2}}  + {\displaystyle \frac {a}{2}}  + n - m - 2\,b)\,\Gamma (
 - {\displaystyle \frac {c}{2}}  + {\displaystyle \frac {a}{2}} 
 + {\displaystyle \frac {3}{2}}  - n - m)\,\Gamma ( - 1 + 2\,m + 
b)\,\Gamma (b) \\
\Gamma ( - c + 2\,b - 1 + 2\,m)(\mathrm{sin}({\displaystyle 
\frac {\pi \,(2\,b - a + 3\,c + 6\,n + 2\,m)}{2}} ) \\
\mbox{} + \mathrm{sin}({\displaystyle \frac {\pi \,(6\,b - 3\,a
 + c + 2\,n + 6\,m)}{2}} ) - \mathrm{sin}({\displaystyle \frac {
\pi \,( - 6\,b + a + c + 2\,n - 2\,m)}{2}} ) \\
\mbox{} + \mathrm{sin}({\displaystyle \frac {\pi \,(2\,b + a + c
 + 2\,n - 2\,m)}{2}} ) - \mathrm{sin}({\displaystyle \frac {\pi 
\,(c - 3\,a + 2\,n + 6\,m + 2\,b)}{2}} ) \\
\mbox{} - \mathrm{sin}({\displaystyle \frac {\pi \,(3\,c - a + 6
\,n + 2\,m - 2\,b)}{2}} ))) \\
\mathrm{\%1} := \Gamma ( - {\displaystyle \frac {c}{2}}  + 
{\displaystyle \frac {a}{2}}  - {\displaystyle \frac {1}{2}}  + n
 + m) \\
\mathrm{\%2} := \Gamma ({\displaystyle \frac {a}{2}}  + 
{\displaystyle \frac {c}{2}}  + {\displaystyle \frac {3}{2}}  - n
 - m) }
}
\end{maplelatex}

\begin{maplelatex}
\mapleinline{inert}{2d}{405, ":   ",
-(cos(Pi*(b+2*n-3*a+c+4*m))-cos(Pi*(b+2*n+a-c)))*WPMR(2*m+2*n-2,2*n-1,
-b+3-2*n-c,-b+c+2*m-2*a,-b+1)*GAMMA(2*b+2*n-2)*GAMMA(3-c+2*a-2*n-2*m)*
GAMMA(-a-2+c+2*n+2*m)*GAMMA(c)*GAMMA(1+a-c)/GAMMA(b-c+2*a-2*m+1)/GAMMA
(-b+1-a+2*m)/GAMMA(c+b-2+2*n)/GAMMA(b)/GAMMA(a)/(cos(Pi*(-c+a))-cos(Pi
*(3*a-4*m-c)))+1/2*WPMR(2*m+2*n-2,2*n-1,1+b-c,b+2*n-2-2*a+c+2*m,b+2*n-
1)*GAMMA(-2*b+2-2*n)*GAMMA(b+2*n-1)*GAMMA(3-c+2*a-2*n-2*m)*GAMMA(-a-2+
c+2*n+2*m)*GAMMA(c)*GAMMA(1+a-c)*(-sin(Pi*(2*b+2*n-3*a+c+4*m))+sin(Pi*
(2*n-3*a+c+4*m))+sin(Pi*(2*b+2*n+a-c))-sin(Pi*(a-c+2*n)))/Pi/GAMMA(-a-
1+b+2*n+2*m)/GAMMA(-c-b+3+2*a-2*n-2*m)/GAMMA(c-b)/GAMMA(a)/(cos(Pi*(-c
+a))-cos(Pi*(3*a-4*m-c)))+GAMMA(-a-2+c+2*n+2*m)*GAMMA(3-c+2*a-2*n-2*m)
*GAMMA(1+a-c)*GAMMA(-b+1)/GAMMA(-c-b+3+2*a-2*n-2*m)*GAMMA(-2*m+2*a)/GA
MMA(b-c+2*a-2*m+1)/GAMMA(c-b)/GAMMA(-b+1-a+2*m)/GAMMA(c+b-2+2*n)/GAMMA
(a)*GAMMA(c)*WPMR(2*m+2*n-2,2*m-1,-b+3-2*n-c,-b+c+2*m-2*a,-a+2*m)-1/2*
WPMR(2*m+2*n-2,2*m-1,-c-b+3+2*a-2*n-2*m,c-b,a)*GAMMA(-2*a+2*m)*GAMMA(3
-c+2*a-2*n-2*m)*GAMMA(-a-2+c+2*n+2*m)*GAMMA(c)*GAMMA(1+a-c)*GAMMA(-b+1
)*(cos(2*Pi*(n-a+c+m))-cos(2*Pi*(b+n-a+m)))/Pi^2/GAMMA(-a+2*m)/GAMMA(1
+a-b);}{%
\maplemultiline{
405\mbox{:~~~}  - (\mathrm{cos}(\pi \,(b + 2\,n - 3\,a
 + c + 4\,m)) - \mathrm{cos}(\pi \,(b + 2\,n + a - c))) \\
\mathrm{WPMR}(2\,m + 2\,n - 2, \,2\,n - 1, \, - b + 3 - 2\,n - c
, \, - b + c + 2\,m - 2\,a, \, - b + 1) \\
\Gamma (2\,b + 2\,n - 2)\,\mathrm{\%2}\,\mathrm{\%1}\,\Gamma (c)
\,\Gamma (1 + a - c)/(\Gamma (b - c + 2\,a - 2\,m + 1) \\
\Gamma ( - b + 1 - a + 2\,m)\,\Gamma (c + b - 2 + 2\,n)\,\Gamma (
b)\,\Gamma (a) \\
(\mathrm{cos}(\pi \,( - c + a)) - \mathrm{cos}(\pi \,(3\,a - 4\,m
 - c))))\mbox{} + {\displaystyle \frac {1}{2}}  \\
\mathrm{WPMR}(2\,m + 2\,n - 2, \,2\,n - 1, \,1 + b - c, \,b + 2\,
n - 2 - 2\,a + c + 2\,m, \,b + 2\,n - 1) \\
\Gamma ( - 2\,b + 2 - 2\,n)\,\Gamma (b + 2\,n - 1)\,\mathrm{\%2}
\,\mathrm{\%1}\,\Gamma (c)\,\Gamma (1 + a - c)( \\
 - \mathrm{sin}(\pi \,(2\,b + 2\,n - 3\,a + c + 4\,m)) + \mathrm{
sin}(\pi \,(2\,n - 3\,a + c + 4\,m)) \\
\mbox{} + \mathrm{sin}(\pi \,(2\,b + 2\,n + a - c)) - \mathrm{sin
}(\pi \,(a - c + 2\,n)))/(\pi \,\Gamma ( - a - 1 + b + 2\,n + 2\,
m) \\
\Gamma ( - c - b + 3 + 2\,a - 2\,n - 2\,m)\,\Gamma (c - b)\,
\Gamma (a) \\
(\mathrm{cos}(\pi \,( - c + a)) - \mathrm{cos}(\pi \,(3\,a - 4\,m
 - c))))\mbox{} + \mathrm{\%1}\,\mathrm{\%2}\,\Gamma (1 + a - c)
\,\Gamma ( - b + 1) \\
\Gamma ( - 2\,m + 2\,a)\,\Gamma (c) \\
\mathrm{WPMR}(2\,m + 2\,n - 2, \,2\,m - 1, \, - b + 3 - 2\,n - c
, \, - b + c + 2\,m - 2\,a, \, - a + 2\,m)/( \\
\Gamma ( - c - b + 3 + 2\,a - 2\,n - 2\,m)\,\Gamma (b - c + 2\,a
 - 2\,m + 1)\,\Gamma (c - b) \\
\Gamma ( - b + 1 - a + 2\,m)\,\Gamma (c + b - 2 + 2\,n)\,\Gamma (
a))\mbox{} - {\displaystyle \frac {1}{2}}  \\
\mathrm{WPMR}(2\,m + 2\,n - 2, \,2\,m - 1, \, - c - b + 3 + 2\,a
 - 2\,n - 2\,m, \,c - b, \,a) \\
\Gamma ( - 2\,a + 2\,m)\,\mathrm{\%2}\,\mathrm{\%1}\,\Gamma (c)\,
\Gamma (1 + a - c)\,\Gamma ( - b + 1) \\
(\mathrm{cos}(2\,\pi \,(n - a + c + m)) - \mathrm{cos}(2\,\pi \,(
b + n - a + m)))/(\pi ^{2}\,\Gamma ( - a + 2\,m)\,\Gamma (1 + a
 - b)) \\
\mathrm{\%1} := \Gamma ( - a - 2 + c + 2\,n + 2\,m) \\
\mathrm{\%2} := \Gamma (3 - c + 2\,a - 2\,n - 2\,m) }
}
\end{maplelatex}

\begin{maplelatex}
\mapleinline{inert}{2d}{406, ":   ",
-WPMR(2*m+2*n-2,2*n-1,-b+4-c-2*n-2*m,1+b-c,1+a-c)*GAMMA(2*n-2*a-2+2*c)
*GAMMA(-c+2*a+2-2*n)*GAMMA(b-2+2*n+2*m)*GAMMA(c)^2*(-sin(Pi*(2*c-a+2*n
-b))+sin(Pi*(-a-b+2*n))+sin(Pi*(2*c-3*a+2*n+b))-sin(Pi*(-3*a+b+2*n)))/
GAMMA(b-3+c+2*n+2*m)/GAMMA(c-b)/GAMMA(c-a)/GAMMA(b)/GAMMA(a)/(-sin(2*P
i*(a-b))-sin(2*Pi*b)+sin(2*Pi*a))/(c-1)-GAMMA(2*a-2*n+1-c)^2*WPMR(2*m+
2*n-2,2*n-1,-b+2+c-2*a-2*m,c-2*a-1+2*n+b,-a-1+c+2*n)*GAMMA(2*a+2-2*c-2
*n)*GAMMA(b-2+2*n+2*m)*GAMMA(c)*(cos(Pi*(c-2*a+4*n-b))+2*cos(Pi*(b-c))
-cos(Pi*(3*c-4*a+4*n-b))-2*cos(Pi*(-c+2*a-b))+cos(Pi*(3*c-6*a+4*n+b))-
cos(Pi*(c-4*a+4*n+b)))*(-2*a+2*n-1+c)/GAMMA(-1+b-c+2*a+2*m)/GAMMA(-c+2
*a+2-2*n-b)/GAMMA(-c+2+a-2*n)/GAMMA(b)/GAMMA(a)/(-cos(Pi*(c-a+2*n-2*b)
)+cos(Pi*(c-3*a+2*n))+cos(Pi*(a+c-2*b+2*n))+cos(Pi*(2*b+2*n-a+c))-cos(
Pi*(2*b+2*n-3*a+c))-cos(Pi*(2*n+a+c)))-sin(Pi*c)*WPMR(2*m+2*n-2,2*m-1,
-b+4-c-2*n-2*m,1+b-c,-a+1)/Pi/(c-1)/GAMMA(-1+b-c+2*a+2*m)*GAMMA(2*m-2+
2*a)/GAMMA(-c+2*a+2-2*n-b)/GAMMA(b-3+c+2*n+2*m)*GAMMA(-c+2*a+2-2*n)*GA
MMA(b-2+2*n+2*m)/GAMMA(c-b)*GAMMA(c)^2*GAMMA(1+a-c)*GAMMA(-b+1)/GAMMA(
a)-1/4*WPMR(2*m+2*n-2,2*m-1,-c+2*a+2-2*n-b,-1+b-c+2*a+2*m,-1+a+2*m)*GA
MMA(2-2*a-2*m)*GAMMA(-c+2*a+2-2*n)*GAMMA(b-2+2*n+2*m)*GAMMA(-1+a+2*m)*
GAMMA(c)*GAMMA(1+a-c)*GAMMA(-b+1)*(-sin(Pi*(-3*a+2*c-2*m+2*n))+sin(Pi*
(-5*a+2*c-2*m+2*n))+sin(Pi*(a+2*b+2*m+2*n))-sin(Pi*(-a+2*b+2*m+2*n)))/
Pi^3/GAMMA(-c+2*a+2*m);}{%
\maplemultiline{
406\mbox{:~~~}  - \mathrm{WPMR}(2\,m + 2\,n - 2, \,2\,n
 - 1, \, - b + 4 - c - 2\,n - 2\,m, \,1 + b - c, \,1 + a - c) \\
\Gamma (2\,n - 2\,a - 2 + 2\,c)\,\Gamma ( - c + 2\,a + 2 - 2\,n)
\,\mathrm{\%1}\,\Gamma (c)^{2}( - \mathrm{sin}(\pi \,(2\,c - a + 
2\,n - b)) \\
\mbox{} + \mathrm{sin}(\pi \,( - a - b + 2\,n)) + \mathrm{sin}(
\pi \,(2\,c - 3\,a + 2\,n + b)) - \mathrm{sin}(\pi \,( - 3\,a + b
 + 2\,n)))/( \\
\Gamma (b - 3 + c + 2\,n + 2\,m)\,\Gamma (c - b)\,\Gamma (c - a)
\,\Gamma (b)\,\Gamma (a) \\
( - \mathrm{sin}(2\,\pi \,(a - b)) - \mathrm{sin}(2\,\pi \,b) + 
\mathrm{sin}(2\,\pi \,a))\,(c - 1))\mbox{} - \Gamma (2\,a - 2\,n
 + 1 - c)^{2}\mathrm{WPMR}( \\
2\,m + 2\,n - 2, \,2\,n - 1, \, - b + 2 + c - 2\,a - 2\,m, \,c - 
2\,a - 1 + 2\,n + b,  \\
 - a - 1 + c + 2\,n)\Gamma (2\,a + 2 - 2\,c - 2\,n)\,\mathrm{\%1}
\,\Gamma (c)(\mathrm{cos}(\pi \,(c - 2\,a + 4\,n - b)) \\
\mbox{} + 2\,\mathrm{cos}(\pi \,(b - c)) - \mathrm{cos}(\pi \,(3
\,c - 4\,a + 4\,n - b)) - 2\,\mathrm{cos}(\pi \,( - c + 2\,a - b)
) \\
\mbox{} + \mathrm{cos}(\pi \,(3\,c - 6\,a + 4\,n + b)) - \mathrm{
cos}(\pi \,(c - 4\,a + 4\,n + b)))( - 2\,a + 2\,n - 1 + c)/( \\
\Gamma ( - 1 + b - c + 2\,a + 2\,m)\,\Gamma ( - c + 2\,a + 2 - 2
\,n - b)\,\Gamma ( - c + 2 + a - 2\,n)\,\Gamma (b)\,\Gamma (a)(
 \\
 - \mathrm{cos}(\pi \,(c - a + 2\,n - 2\,b)) + \mathrm{cos}(\pi 
\,(c - 3\,a + 2\,n)) + \mathrm{cos}(\pi \,(a + c - 2\,b + 2\,n))
 \\
\mbox{} + \mathrm{cos}(\pi \,(2\,b + 2\,n - a + c)) - \mathrm{cos
}(\pi \,(2\,b + 2\,n - 3\,a + c)) - \mathrm{cos}(\pi \,(2\,n + a
 + c)))) \\
\mbox{} - \mathrm{sin}(\pi \,c) \,
\mathrm{WPMR}(2\,m + 2\,n - 2, \,2\,m - 1, \, - b + 4 - c - 2\,n
 - 2\,m, \,1 + b - c, \, - a + 1) \\
\Gamma (2\,m - 2 + 2\,a)\,\Gamma ( - c + 2\,a + 2 - 2\,n)\,
\mathrm{\%1}\,\Gamma (c)^{2}\,\Gamma (1 + a - c)\,\Gamma ( - b + 
1)/(\pi \,(c - 1) \\
\Gamma ( - 1 + b - c + 2\,a + 2\,m)\,\Gamma ( - c + 2\,a + 2 - 2
\,n - b)\,\Gamma (b - 3 + c + 2\,n + 2\,m) \\
\Gamma (c - b)\,\Gamma (a))\mbox{} - {\displaystyle \frac {1}{4}
} \mathrm{WPMR}(2\,m + 2\,n - 2, \,2\,m - 1, \, - c + 2\,a + 2 - 
2\,n - b,  \\
 - 1 + b - c + 2\,a + 2\,m, \, - 1 + a + 2\,m)\Gamma (2 - 2\,a - 
2\,m)\,\Gamma ( - c + 2\,a + 2 - 2\,n)\,\mathrm{\%1} \\
\Gamma ( - 1 + a + 2\,m)\,\Gamma (c)\,\Gamma (1 + a - c)\,\Gamma 
( - b + 1)( - \mathrm{sin}(\pi \,( - 3\,a + 2\,c - 2\,m + 2\,n))
 \\
\mbox{} + \mathrm{sin}(\pi \,( - 5\,a + 2\,c - 2\,m + 2\,n)) + 
\mathrm{sin}(\pi \,(a + 2\,b + 2\,m + 2\,n)) \\
\mbox{} - \mathrm{sin}(\pi \,( - a + 2\,b + 2\,m + 2\,n)))/(\pi 
^{3}\,\Gamma ( - c + 2\,a + 2\,m)) \\
\mathrm{\%1} := \Gamma (b - 2 + 2\,n + 2\,m) }
}
\end{maplelatex}

\begin{maplelatex}
\mapleinline{inert}{2d}{407, ":   ",
4*Pi^2*(cos(Pi*(-3*c+a+4*m))-cos(Pi*(c-b+2*n)))*WPMR(2*m+2*n-1,2*n-1,-
a+1,1+b-2*c+2*m,-1/2*a+1/2*b+n)*GAMMA(1/2*a+1/2*b-n+1)*GAMMA(-2*m+2*c)
*GAMMA(a-b)/GAMMA(-c+1/2*b+1/2*a-n+1)/GAMMA(-c+1/2*b-1/2*a+n+1+2*m)/GA
MMA(1/2*a-1/2*b-n+1)/GAMMA(2*c-2*m-b)/GAMMA(c-2*m)/GAMMA(c)/GAMMA(a)/(
cos(1/2*Pi*(-3*b-4*m-3*a+6*c+2*n))-cos(1/2*Pi*(-3*b-12*m-3*a+10*c+2*n)
)-cos(1/2*Pi*(b-4*m+a+2*c+2*n))-cos(1/2*Pi*(6*n+4*m+2*c-b-a))+cos(1/2*
Pi*(2*n+12*m+a+b-6*c))+cos(1/2*Pi*(6*n-4*m-a-b+6*c)))-WPMR(2*m+2*n-1,2
*n-1,-b+1,-2*c+a+2*m+1,n+1/2*a-1/2*b)*GAMMA(-2*c+a+2*m+1)*GAMMA(-2*m+2
*c)*GAMMA(n+1/2*a-1/2*b)*GAMMA(1/2*a+1/2*b-n+1)*GAMMA(b-a)*(-cos(1/2*P
i*(2*n-3*b-2*c+3*a+4*m))-cos(1/2*Pi*(6*n-b+6*c-3*a-4*m))+cos(1/2*Pi*(2
*n-3*b+6*c-a-4*m))+cos(1/2*Pi*(6*n-b-2*c+a+4*m))+cos(1/2*Pi*(10*c-5*a-
12*m+b+2*n))+cos(1/2*Pi*(-2*c+b-a+2*n+4*m))-cos(1/2*Pi*(2*c-a-4*m+b+2*
n))-cos(1/2*Pi*(-10*c+3*a+12*m+b+2*n)))/GAMMA(c-2*m)/GAMMA(1-c-1/2*b+1
/2*a+n+2*m)/GAMMA(-c+1/2*b+1/2*a-n+1)/GAMMA(c)/GAMMA(b)/(cos(1/2*Pi*(-
3*b-4*m-3*a+6*c+2*n))-cos(1/2*Pi*(-3*b-12*m-3*a+10*c+2*n))-cos(1/2*Pi*
(b-4*m+a+2*c+2*n))-cos(1/2*Pi*(6*n+4*m+2*c-b-a))+cos(1/2*Pi*(2*n+12*m+
a+b-6*c))+cos(1/2*Pi*(6*n-4*m-a-b+6*c)))-1/2*WPMR(2*m+2*n-1,2*m,b,2*c-
a-2*m,c)*GAMMA(-1/2*a+1/2*b+n)*GAMMA(c-1/2*b-1/2*a+n)*GAMMA(1/2*a+1/2*
b-n+1)*(sin(Pi*(-3*c+a+4*m+b))+sin(Pi*(-c+a-b))-sin(Pi*(b-c+a))-sin(Pi
*(-3*c+a+4*m-b)))/Pi/GAMMA(c+1/2*b-1/2*a+n)/(cos(1/2*Pi*(-6*c+b+a+2*n+
8*m))-cos(1/2*Pi*(2*c+b+a+2*n)))+Pi*WPMR(2*m+2*n-1,2*m,-a+1,1+b-2*c+2*
m,-c+1+2*m)/sin(1/2*Pi*(-2*c+b+a+2*n+4*m))/(-1+2*c-2*m)/GAMMA(2*c-2*m-
b)/GAMMA(c-2*m)/GAMMA(2*c-a-2*m)*GAMMA(-2*m+2*c)^2*GAMMA(-1/2*a+1/2*b+
n)*GAMMA(c-1/2*b-1/2*a+n)*GAMMA(1/2*a+1/2*b-n+1)/GAMMA(-c+1/2*b-1/2*a+
n+1+2*m)/GAMMA(c)/GAMMA(b)/GAMMA(a);}{%
\maplemultiline{
407\mbox{:~~~} 4\,\pi ^{2}\,(\mathrm{cos}(\pi \,( - 3\,c
 + a + 4\,m)) - \mathrm{cos}(\pi \,(c - b + 2\,n))) \\
\mathrm{WPMR}(2\,m + 2\,n - 1, \,2\,n - 1, \, - a + 1, \,1 + b - 
2\,c + 2\,m, \, - {\displaystyle \frac {a}{2}}  + {\displaystyle 
\frac {b}{2}}  + n)\,\mathrm{\%1} \\
\Gamma ( - 2\,m + 2\,c)\,\Gamma (a - b) \left/ {\vrule 
height0.80em width0em depth0.80em} \right. \!  \! (\Gamma ( - c
 + {\displaystyle \frac {b}{2}}  + {\displaystyle \frac {a}{2}} 
 - n + 1)\,\Gamma ( - c + {\displaystyle \frac {b}{2}}  - 
{\displaystyle \frac {a}{2}}  + n + 1 + 2\,m) \\
\Gamma ({\displaystyle \frac {a}{2}}  - {\displaystyle \frac {b}{
2}}  - n + 1)\,\Gamma (2\,c - 2\,m - b)\,\Gamma (c - 2\,m)\,
\Gamma (c)\,\Gamma (a)( \\
\mathrm{cos}({\displaystyle \frac {\pi \,( - 3\,b - 4\,m - 3\,a
 + 6\,c + 2\,n)}{2}} ) - \mathrm{cos}({\displaystyle \frac {\pi 
\,( - 3\,b - 12\,m - 3\,a + 10\,c + 2\,n)}{2}} ) \\
\mbox{} - \mathrm{cos}({\displaystyle \frac {\pi \,(b - 4\,m + a
 + 2\,c + 2\,n)}{2}} ) - \mathrm{cos}({\displaystyle \frac {\pi 
\,(6\,n + 4\,m + 2\,c - b - a)}{2}} ) \\
\mbox{} + \mathrm{cos}({\displaystyle \frac {\pi \,(2\,n + 12\,m
 + a + b - 6\,c)}{2}} ) + \mathrm{cos}({\displaystyle \frac {\pi 
\,(6\,n - 4\,m - a - b + 6\,c)}{2}} )))\mbox{} -  \\
\mathrm{WPMR}(2\,m + 2\,n - 1, \,2\,n - 1, \, - b + 1, \, - 2\,c
 + a + 2\,m + 1, \,n + {\displaystyle \frac {a}{2}}  - 
{\displaystyle \frac {b}{2}} ) \\
\Gamma ( - 2\,c + a + 2\,m + 1)\,\Gamma ( - 2\,m + 2\,c)\,\Gamma 
(n + {\displaystyle \frac {a}{2}}  - {\displaystyle \frac {b}{2}
} )\,\mathrm{\%1}\,\Gamma (b - a)( \\
 - \mathrm{cos}({\displaystyle \frac {\pi \,(2\,n - 3\,b - 2\,c
 + 3\,a + 4\,m)}{2}} ) - \mathrm{cos}({\displaystyle \frac {\pi 
\,(6\,n - b + 6\,c - 3\,a - 4\,m)}{2}} ) \\
\mbox{} + \mathrm{cos}({\displaystyle \frac {\pi \,(2\,n - 3\,b
 + 6\,c - a - 4\,m)}{2}} ) + \mathrm{cos}({\displaystyle \frac {
\pi \,(6\,n - b - 2\,c + a + 4\,m)}{2}} ) \\
\mbox{} + \mathrm{cos}({\displaystyle \frac {\pi \,(10\,c - 5\,a
 - 12\,m + b + 2\,n)}{2}} ) + \mathrm{cos}({\displaystyle \frac {
\pi \,( - 2\,c + b - a + 2\,n + 4\,m)}{2}} ) \\
\mbox{} - \mathrm{cos}({\displaystyle \frac {\pi \,(2\,c - a - 4
\,m + b + 2\,n)}{2}} ) - \mathrm{cos}({\displaystyle \frac {\pi 
\,( - 10\,c + 3\,a + 12\,m + b + 2\,n)}{2}} )) \\
 \left/ {\vrule height0.80em width0em depth0.80em} \right. \! 
 \! (\Gamma (c - 2\,m)\,\Gamma (1 - c - {\displaystyle \frac {b}{
2}}  + {\displaystyle \frac {a}{2}}  + n + 2\,m)\,\Gamma ( - c + 
{\displaystyle \frac {b}{2}}  + {\displaystyle \frac {a}{2}}  - n
 + 1)\,\Gamma (c)\,\Gamma (b)( \\
\mathrm{cos}({\displaystyle \frac {\pi \,( - 3\,b - 4\,m - 3\,a
 + 6\,c + 2\,n)}{2}} ) - \mathrm{cos}({\displaystyle \frac {\pi 
\,( - 3\,b - 12\,m - 3\,a + 10\,c + 2\,n)}{2}} ) \\
\mbox{} - \mathrm{cos}({\displaystyle \frac {\pi \,(b - 4\,m + a
 + 2\,c + 2\,n)}{2}} ) - \mathrm{cos}({\displaystyle \frac {\pi 
\,(6\,n + 4\,m + 2\,c - b - a)}{2}} ) \\
\mbox{} + \mathrm{cos}({\displaystyle \frac {\pi \,(2\,n + 12\,m
 + a + b - 6\,c)}{2}} ) + \mathrm{cos}({\displaystyle \frac {\pi 
\,(6\,n - 4\,m - a - b + 6\,c)}{2}} )))\mbox{} - {\displaystyle 
\frac {1}{2}}  \\
\mathrm{WPMR}(2\,m + 2\,n - 1, \,2\,m, \,b, \,2\,c - a - 2\,m, \,
c)\,\Gamma ( - {\displaystyle \frac {a}{2}}  + {\displaystyle 
\frac {b}{2}}  + n)\,\Gamma (c - {\displaystyle \frac {b}{2}}  - 
{\displaystyle \frac {a}{2}}  + n) \\
\mathrm{\%1}(\mathrm{sin}(\pi \,( - 3\,c + a + 4\,m + b)) + 
\mathrm{sin}(\pi \,( - c + a - b)) - \mathrm{sin}(\pi \,(b - c + 
a)) \\
\mbox{} - \mathrm{sin}(\pi \,( - 3\,c + a + 4\,m - b))) \left/ 
{\vrule height0.80em width0em depth0.80em} \right. \!  \! (\pi \,
\Gamma (c + {\displaystyle \frac {b}{2}}  - {\displaystyle 
\frac {a}{2}}  + n) \\
(\mathrm{cos}({\displaystyle \frac {\pi \,( - 6\,c + b + a + 2\,n
 + 8\,m)}{2}} ) - \mathrm{cos}({\displaystyle \frac {\pi \,(2\,c
 + b + a + 2\,n)}{2}} )))\mbox{} + \pi  \\
\mathrm{WPMR}(2\,m + 2\,n - 1, \,2\,m, \, - a + 1, \,1 + b - 2\,c
 + 2\,m, \, - c + 1 + 2\,m) \\
\Gamma ( - 2\,m + 2\,c)^{2}\,\Gamma ( - {\displaystyle \frac {a}{
2}}  + {\displaystyle \frac {b}{2}}  + n)\,\Gamma (c - 
{\displaystyle \frac {b}{2}}  - {\displaystyle \frac {a}{2}}  + n
)\,\mathrm{\%1} \left/ {\vrule height0.80em width0em depth0.80em}
 \right. \!  \! ( \\
\mathrm{sin}({\displaystyle \frac {\pi \,( - 2\,c + b + a + 2\,n
 + 4\,m)}{2}} )\,( - 1 + 2\,c - 2\,m)\,\Gamma (2\,c - 2\,m - b)\,
\Gamma (c - 2\,m) \\
\Gamma (2\,c - a - 2\,m)\,\Gamma ( - c + {\displaystyle \frac {b
}{2}}  - {\displaystyle \frac {a}{2}}  + n + 1 + 2\,m)\,\Gamma (c
)\,\Gamma (b)\,\Gamma (a)) \\
\mathrm{\%1} := \Gamma ({\displaystyle \frac {a}{2}}  + 
{\displaystyle \frac {b}{2}}  - n + 1) }
}
\end{maplelatex}

\begin{maplelatex}
\mapleinline{inert}{2d}{408, ":   ",
-1/4*WPMR(2*m+2*n-1,2*m,-c+2*m+2*b,-a+2*b+2*m,2*m+b)*GAMMA(-1/2*a+b-1/
2*c-n+m+1)*GAMMA(1/2*c+1/2*a-n-m+1)*GAMMA(-1/2*a+b-1/2*c+n+m)*GAMMA(1/
2*c-1/2*a+n+m)*GAMMA(-1/2*c+1/2*a+n+m)*(sin(Pi*(-c+a-b))+sin(Pi*(a-3*b
-4*m+c))-sin(Pi*(b-c+a))-sin(Pi*(a-5*b-4*m+c)))/sin(2*Pi*(b+m))/Pi^2/G
AMMA(-1/2*c-1/2*a+n+3*m+2*b)+1/2*sin(1/2*Pi*(-a-c+2*n+2*m))*WPMR(2*m+2
*n-1,2*m,-a+1,-c+1,-b+1)*(-a-c+2*n+2*m)*GAMMA(-1/2*c+1/2*a+n+m)*GAMMA(
1/2*c-1/2*a+n+m)*GAMMA(-1/2*a+b-1/2*c+n+m)*GAMMA(-1/2*a+b-1/2*c-n+m+1)
*GAMMA(2*m+2*b)^2*GAMMA(1/2*c+1/2*a-n-m)^2/Pi/(2*b-1+2*m)/GAMMA(2*m+b)
/GAMMA(-a+2*b+2*m)/GAMMA(-c+2*m+2*b)/GAMMA(c)/GAMMA(b)/GAMMA(a)+WPMR(2
*m+2*n-1,2*n-1,-a+1,-c+1,-1/2*a+b-1/2*c+n+m)*Pi*GAMMA(-1/2*c+1/2*a+n+m
)*GAMMA(1/2*c+1/2*a-n-m)^2*GAMMA(2*m+2*b)*GAMMA(a-2*b+c-2*m)*(-sin(1/2
*Pi*(-a+c+6*n+2*m-2*b))+sin(1/2*Pi*(a+3*c+2*n-2*m-2*b))+sin(1/2*Pi*(a-
c+2*n-2*m-6*b))+sin(1/2*Pi*(-3*a-c+2*n+6*m+6*b)))*(-a-c+2*n+2*m)/GAMMA
(1/2*a-b+1/2*c+n-m)/GAMMA(-1/2*c+1/2*a-n-m+1)/GAMMA(1/2*a-b+1/2*c-n-m+
1)/GAMMA(2*m+b)/GAMMA(c)/GAMMA(b)/GAMMA(a)/(cos(1/2*Pi*(c-3*a+2*n+6*m+
2*b))-cos(1/2*Pi*(6*b-3*a+c+2*n+6*m))-cos(1/2*Pi*(2*b+a+c+2*n-2*m))-co
s(1/2*Pi*(3*c-a+6*n+2*m-2*b))+cos(1/2*Pi*(-6*b+a+c+2*n-2*m))+cos(1/2*P
i*(2*b-a+3*c+6*n+2*m)))+WPMR(2*m+2*n-1,2*n-1,c-2*m-2*b+1,a-2*b-2*m+1,1
/2*a-b+1/2*c+n-m)*GAMMA(2*m+2*b)*GAMMA(1/2*c+1/2*a-n-m+1)*GAMMA(-1/2*c
+1/2*a+n+m)*GAMMA(a-2*b-2*m+1)*GAMMA(-a+2*b-c+2*m)*(-cos(1/2*Pi*(2*n-8
*b+3*c+3*a-6*m))-cos(1/2*Pi*(6*n+4*b+c-3*a+6*m))+cos(1/2*Pi*(2*n+3*c-a
+2*m))+cos(1/2*Pi*(6*n-4*b+c+a-2*m))+cos(1/2*Pi*(-5*a+12*b+10*m-c+2*n)
)+cos(1/2*Pi*(-a-c+2*n+2*m))-cos(1/2*Pi*(-a+4*b+2*m-c+2*n))-cos(1/2*Pi
*(3*a-8*b-6*m-c+2*n)))/GAMMA(2*m+b)/GAMMA(1+1/2*c+1/2*a+n-m-2*b)/GAMMA
(-1/2*c+1/2*a-n-m+1)/GAMMA(-c+2*m+2*b)/GAMMA(b)/(-cos(1/2*Pi*(c-3*a+2*
n+6*m+2*b))+cos(1/2*Pi*(6*b-3*a+c+2*n+6*m))+cos(1/2*Pi*(2*b+a+c+2*n-2*
m))+cos(1/2*Pi*(3*c-a+6*n+2*m-2*b))-cos(1/2*Pi*(-6*b+a+c+2*n-2*m))-cos
(1/2*Pi*(2*b-a+3*c+6*n+2*m)));}{%
\maplemultiline{
408\mbox{:~~~}  - {\displaystyle \frac {1}{4}} \mathrm{
WPMR}(2\,m + 2\,n - 1, \,2\,m, \, - c + 2\,m + 2\,b, \, - a + 2\,
b + 2\,m, \,2\,m + b) \\
\Gamma ( - {\displaystyle \frac {a}{2}}  + b - {\displaystyle 
\frac {c}{2}}  - n + m + 1)\,\Gamma ({\displaystyle \frac {c}{2}
}  + {\displaystyle \frac {a}{2}}  - n - m + 1)\,\Gamma ( - 
{\displaystyle \frac {a}{2}}  + b - {\displaystyle \frac {c}{2}} 
 + n + m) \\
\Gamma ({\displaystyle \frac {c}{2}}  - {\displaystyle \frac {a}{
2}}  + n + m)\,\mathrm{\%1}(\mathrm{sin}(\pi \,( - c + a - b)) + 
\mathrm{sin}(\pi \,(a - 3\,b - 4\,m + c)) \\
\mbox{} - \mathrm{sin}(\pi \,(b - c + a)) - \mathrm{sin}(\pi \,(a
 - 5\,b - 4\,m + c))) \left/ {\vrule 
height0.80em width0em depth0.80em} \right. \!  \! (\mathrm{sin}(2
\,\pi \,(b + m))\,\pi ^{2} \\
\Gamma ( - {\displaystyle \frac {c}{2}}  - {\displaystyle \frac {
a}{2}}  + n + 3\,m + 2\,b))\mbox{} + {\displaystyle \frac {1}{2}
} \mathrm{sin}({\displaystyle \frac {\pi \,( - a - c + 2\,n + 2\,
m)}{2}} ) \\
\mathrm{WPMR}(2\,m + 2\,n - 1, \,2\,m, \, - a + 1, \, - c + 1, \,
 - b + 1)\,( - a - c + 2\,n + 2\,m)\,\mathrm{\%1} \\
\Gamma ({\displaystyle \frac {c}{2}}  - {\displaystyle \frac {a}{
2}}  + n + m)\,\Gamma ( - {\displaystyle \frac {a}{2}}  + b - 
{\displaystyle \frac {c}{2}}  + n + m)\,\Gamma ( - 
{\displaystyle \frac {a}{2}}  + b - {\displaystyle \frac {c}{2}} 
 - n + m + 1) \\
\Gamma (2\,m + 2\,b)^{2}\,\Gamma ({\displaystyle \frac {c}{2}} 
 + {\displaystyle \frac {a}{2}}  - n - m)^{2}/(\pi \,(2\,b - 1 + 
2\,m)\,\Gamma (2\,m + b)\,\Gamma ( - a + 2\,b + 2\,m) \\
\Gamma ( - c + 2\,m + 2\,b)\,\Gamma (c)\,\Gamma (b)\,\Gamma (a))
\mbox{} +  \\
\mathrm{WPMR}(2\,m + 2\,n - 1, \,2\,n - 1, \, - a + 1, \, - c + 1
, \, - {\displaystyle \frac {a}{2}}  + b - {\displaystyle \frac {
c}{2}}  + n + m)\,\pi \,\mathrm{\%1} \\
\Gamma ({\displaystyle \frac {c}{2}}  + {\displaystyle \frac {a}{
2}}  - n - m)^{2}\,\Gamma (2\,m + 2\,b)\,\Gamma (a - 2\,b + c - 2
\,m)( \\
 - \mathrm{sin}({\displaystyle \frac {\pi \,( - a + c + 6\,n + 2
\,m - 2\,b)}{2}} ) + \mathrm{sin}({\displaystyle \frac {\pi \,(a
 + 3\,c + 2\,n - 2\,m - 2\,b)}{2}} ) \\
\mbox{} + \mathrm{sin}({\displaystyle \frac {\pi \,(a - c + 2\,n
 - 2\,m - 6\,b)}{2}} ) + \mathrm{sin}({\displaystyle \frac {\pi 
\,( - 3\,a - c + 2\,n + 6\,m + 6\,b)}{2}} )) \\
( - a - c + 2\,n + 2\,m) \left/ {\vrule 
height0.80em width0em depth0.80em} \right. \!  \! (\Gamma (
{\displaystyle \frac {a}{2}}  - b + {\displaystyle \frac {c}{2}} 
 + n - m)\,\Gamma ( - {\displaystyle \frac {c}{2}}  + 
{\displaystyle \frac {a}{2}}  - n - m + 1) \\
\Gamma ({\displaystyle \frac {a}{2}}  - b + {\displaystyle 
\frac {c}{2}}  - n - m + 1)\,\Gamma (2\,m + b)\,\Gamma (c)\,
\Gamma (b)\,\Gamma (a)( \\
\mathrm{cos}({\displaystyle \frac {\pi \,(c - 3\,a + 2\,n + 6\,m
 + 2\,b)}{2}} ) - \mathrm{cos}({\displaystyle \frac {\pi \,(6\,b
 - 3\,a + c + 2\,n + 6\,m)}{2}} ) \\
\mbox{} - \mathrm{cos}({\displaystyle \frac {\pi \,(2\,b + a + c
 + 2\,n - 2\,m)}{2}} ) - \mathrm{cos}({\displaystyle \frac {\pi 
\,(3\,c - a + 6\,n + 2\,m - 2\,b)}{2}} ) \\
\mbox{} + \mathrm{cos}({\displaystyle \frac {\pi \,( - 6\,b + a
 + c + 2\,n - 2\,m)}{2}} ) + \mathrm{cos}({\displaystyle \frac {
\pi \,(2\,b - a + 3\,c + 6\,n + 2\,m)}{2}} )))\mbox{} +  \\
\mathrm{WPMR}(2\,m + 2\,n - 1, \,2\,n - 1, \,c - 2\,m - 2\,b + 1
, \,a - 2\,b - 2\,m + 1,  \\
{\displaystyle \frac {a}{2}}  - b + {\displaystyle \frac {c}{2}} 
 + n - m)\Gamma (2\,m + 2\,b)\,\Gamma ({\displaystyle \frac {c}{2
}}  + {\displaystyle \frac {a}{2}}  - n - m + 1)\,\mathrm{\%1}\,
\Gamma (a - 2\,b - 2\,m + 1) \\
\Gamma ( - a + 2\,b - c + 2\,m)( - \mathrm{cos}({\displaystyle 
\frac {\pi \,(2\,n - 8\,b + 3\,c + 3\,a - 6\,m)}{2}} ) \\
\mbox{} - \mathrm{cos}({\displaystyle \frac {\pi \,(6\,n + 4\,b
 + c - 3\,a + 6\,m)}{2}} ) + \mathrm{cos}({\displaystyle \frac {
\pi \,(2\,n + 3\,c - a + 2\,m)}{2}} ) \\
\mbox{} + \mathrm{cos}({\displaystyle \frac {\pi \,(6\,n - 4\,b
 + c + a - 2\,m)}{2}} ) + \mathrm{cos}({\displaystyle \frac {\pi 
\,( - 5\,a + 12\,b + 10\,m - c + 2\,n)}{2}} ) \\
\mbox{} + \mathrm{cos}({\displaystyle \frac {\pi \,( - a - c + 2
\,n + 2\,m)}{2}} ) - \mathrm{cos}({\displaystyle \frac {\pi \,(
 - a + 4\,b + 2\,m - c + 2\,n)}{2}} ) \\
\mbox{} - \mathrm{cos}({\displaystyle \frac {\pi \,(3\,a - 8\,b
 - 6\,m - c + 2\,n)}{2}} )) \left/ {\vrule 
height0.80em width0em depth0.80em} \right. \!  \! (\Gamma (2\,m
 + b) \\
\Gamma (1 + {\displaystyle \frac {c}{2}}  + {\displaystyle 
\frac {a}{2}}  + n - m - 2\,b)\,\Gamma ( - {\displaystyle \frac {
c}{2}}  + {\displaystyle \frac {a}{2}}  - n - m + 1)\,\Gamma ( - 
c + 2\,m + 2\,b)\,\Gamma (b)( \\
 - \mathrm{cos}({\displaystyle \frac {\pi \,(c - 3\,a + 2\,n + 6
\,m + 2\,b)}{2}} ) + \mathrm{cos}({\displaystyle \frac {\pi \,(6
\,b - 3\,a + c + 2\,n + 6\,m)}{2}} ) \\
\mbox{} + \mathrm{cos}({\displaystyle \frac {\pi \,(2\,b + a + c
 + 2\,n - 2\,m)}{2}} ) + \mathrm{cos}({\displaystyle \frac {\pi 
\,(3\,c - a + 6\,n + 2\,m - 2\,b)}{2}} ) \\
\mbox{} - \mathrm{cos}({\displaystyle \frac {\pi \,( - 6\,b + a
 + c + 2\,n - 2\,m)}{2}} ) - \mathrm{cos}({\displaystyle \frac {
\pi \,(2\,b - a + 3\,c + 6\,n + 2\,m)}{2}} ))) \\
\mathrm{\%1} := \Gamma ( - {\displaystyle \frac {c}{2}}  + 
{\displaystyle \frac {a}{2}}  + n + m) }
}
\end{maplelatex}

\begin{maplelatex}
\mapleinline{inert}{2d}{409, ":   ",
GAMMA(-a+c+2*n+2*m-1)*GAMMA(2*a-c-2*n-2*m+2)*GAMMA(1+a-c)*GAMMA(-b+1)/
GAMMA(-2*n-b+2+2*a-c-2*m)*GAMMA(-1+2*a-2*m)/GAMMA(b-c+2*a-2*m)/GAMMA(2
-a-b+2*m)/GAMMA(c+b-2+2*n)/GAMMA(c-b)/GAMMA(a)*GAMMA(c)*WPMR(2*m+2*n-1
,2*m,-b+3-2*n-c,-b+c+2*m-2*a+1,1-a+2*m)+(cos(Pi*(b+2*n-3*a+c+4*m))-cos
(Pi*(b+2*n+a-c)))*WPMR(2*m+2*n-1,2*n-1,-b+3-2*n-c,-b+c+2*m-2*a+1,-b+1)
*GAMMA(c)*GAMMA(1+a-c)*GAMMA(2*b+2*n-2)*GAMMA(2*a-c-2*n-2*m+2)*GAMMA(-
a+c+2*n+2*m-1)/GAMMA(2-a-b+2*m)/GAMMA(b-c+2*a-2*m)/GAMMA(b)/GAMMA(a)/G
AMMA(c+b-2+2*n)/(-cos(Pi*(-c+a))+cos(Pi*(3*a-4*m-c)))+1/2*WPMR(2*m+2*n
-1,2*m,c-b,-2*n-b+2+2*a-c-2*m,a)*GAMMA(2*a-c-2*n-2*m+2)*GAMMA(-a+c+2*n
+2*m-1)*GAMMA(-2*a+2*m+1)*GAMMA(c)*GAMMA(1+a-c)*GAMMA(-b+1)*(cos(2*Pi*
(n-a+c+m))-cos(2*Pi*(b+n-a+m)))/Pi^2/GAMMA(1-a+2*m)/GAMMA(1+a-b)+1/2*W
PMR(2*m+2*n-1,2*n-1,1+b-c,b+2*n-1-2*a+c+2*m,b+2*n-1)*GAMMA(-2*b+2-2*n)
*GAMMA(b+2*n-1)*GAMMA(2*a-c-2*n-2*m+2)*GAMMA(-a+c+2*n+2*m-1)*GAMMA(c)*
GAMMA(1+a-c)*(sin(Pi*(2*b+2*n-3*a+c+4*m))-sin(Pi*(2*n-3*a+c+4*m))-sin(
Pi*(2*b+2*n+a-c))+sin(Pi*(a-c+2*n)))/Pi/GAMMA(-a+b+2*n+2*m)/GAMMA(-2*n
-b+2+2*a-c-2*m)/GAMMA(c-b)/GAMMA(a)/(-cos(Pi*(-c+a))+cos(Pi*(3*a-4*m-c
)));}{%
\maplemultiline{
409\mbox{:~~~} \mathrm{\%1}\,\mathrm{\%2}\,\Gamma (1 + a
 - c)\,\Gamma ( - b + 1)\,\Gamma ( - 1 + 2\,a - 2\,m)\,\Gamma (c)
 \\
\mathrm{WPMR}(2\,m + 2\,n - 1, \,2\,m, \, - b + 3 - 2\,n - c, \,
 - b + c + 2\,m - 2\,a + 1, \,1 - a + 2\,m)/( \\
\Gamma ( - 2\,n - b + 2 + 2\,a - c - 2\,m)\,\Gamma (b - c + 2\,a
 - 2\,m)\,\Gamma (2 - a - b + 2\,m) \\
\Gamma (c + b - 2 + 2\,n)\,\Gamma (c - b)\,\Gamma (a))\mbox{} + 
 \\
(\mathrm{cos}(\pi \,(b + 2\,n - 3\,a + c + 4\,m)) - \mathrm{cos}(
\pi \,(b + 2\,n + a - c))) \\
\mathrm{WPMR}(2\,m + 2\,n - 1, \,2\,n - 1, \, - b + 3 - 2\,n - c
, \, - b + c + 2\,m - 2\,a + 1, \, - b + 1) \\
\Gamma (c)\,\Gamma (1 + a - c)\,\Gamma (2\,b + 2\,n - 2)\,
\mathrm{\%2}\,\mathrm{\%1}/(\Gamma (2 - a - b + 2\,m)\,\Gamma (b
 - c + 2\,a - 2\,m) \\
\Gamma (b)\,\Gamma (a)\,\Gamma (c + b - 2 + 2\,n)\,( - \mathrm{
cos}(\pi \,( - c + a)) + \mathrm{cos}(\pi \,(3\,a - 4\,m - c))))
\mbox{} + {\displaystyle \frac {1}{2}}  \\
\mathrm{WPMR}(2\,m + 2\,n - 1, \,2\,m, \,c - b, \, - 2\,n - b + 2
 + 2\,a - c - 2\,m, \,a)\,\mathrm{\%2}\,\mathrm{\%1} \\
\Gamma ( - 2\,a + 2\,m + 1)\,\Gamma (c)\,\Gamma (1 + a - c)\,
\Gamma ( - b + 1) \\
(\mathrm{cos}(2\,\pi \,(n - a + c + m)) - \mathrm{cos}(2\,\pi \,(
b + n - a + m)))/(\pi ^{2}\,\Gamma (1 - a + 2\,m) \\
\Gamma (1 + a - b))\mbox{} + {\displaystyle \frac {1}{2}}  \\
\mathrm{WPMR}(2\,m + 2\,n - 1, \,2\,n - 1, \,1 + b - c, \,b + 2\,
n - 1 - 2\,a + c + 2\,m, \,b + 2\,n - 1) \\
\Gamma ( - 2\,b + 2 - 2\,n)\,\Gamma (b + 2\,n - 1)\,\mathrm{\%2}
\,\mathrm{\%1}\,\Gamma (c)\,\Gamma (1 + a - c)( \\
\mathrm{sin}(\pi \,(2\,b + 2\,n - 3\,a + c + 4\,m)) - \mathrm{sin
}(\pi \,(2\,n - 3\,a + c + 4\,m)) \\
\mbox{} - \mathrm{sin}(\pi \,(2\,b + 2\,n + a - c)) + \mathrm{sin
}(\pi \,(a - c + 2\,n)))/(\pi \,\Gamma ( - a + b + 2\,n + 2\,m)
 \\
\Gamma ( - 2\,n - b + 2 + 2\,a - c - 2\,m)\,\Gamma (c - b)\,
\Gamma (a) \\
( - \mathrm{cos}(\pi \,( - c + a)) + \mathrm{cos}(\pi \,(3\,a - 4
\,m - c)))) \\
\mathrm{\%1} := \Gamma ( - a + c + 2\,n + 2\,m - 1) \\
\mathrm{\%2} := \Gamma (2\,a - c - 2\,n - 2\,m + 2) }
}
\end{maplelatex}

\mapleinline{inert}{2d}{410, ":   ",
-WPMR(2*m+2*n-1,2*n-1,-b+3-c-2*n-2*m,1+b-c,1+a-c)*GAMMA(2*n-2*a-2+2*c)
*GAMMA(-c+2*a+2-2*n)*GAMMA(b-1+2*n+2*m)*GAMMA(c)^2*(sin(Pi*(2*c-a+2*n-
b))-sin(Pi*(-a-b+2*n))-sin(Pi*(2*c-3*a+2*n+b))+sin(Pi*(-3*a+b+2*n)))/G
AMMA(b-2+c+2*n+2*m)/GAMMA(c-b)/GAMMA(c-a)/GAMMA(b)/GAMMA(a)/(sin(2*Pi*
(a-b))+sin(2*Pi*b)-sin(2*Pi*a))/(c-1)+GAMMA(2*a-2*n+1-c)^2*WPMR(2*m+2*
n-1,2*n-1,-b+1+c-2*a-2*m,c-2*a-1+2*n+b,-a-1+c+2*n)*GAMMA(2*a+2-2*c-2*n
)*GAMMA(b-1+2*n+2*m)*GAMMA(c)*(2*cos(Pi*(-c+2*a-b))+cos(Pi*(3*c-4*a+4*
n-b))-2*cos(Pi*(b-c))-cos(Pi*(c-2*a+4*n-b))-cos(Pi*(3*c-6*a+4*n+b))+co
s(Pi*(c-4*a+4*n+b)))*(-2*a+2*n-1+c)/GAMMA(b-c+2*a+2*m)/GAMMA(-c+2*a+2-
2*n-b)/GAMMA(-c+2+a-2*n)/GAMMA(b)/GAMMA(a)/(-cos(Pi*(c-a+2*n-2*b))+cos
(Pi*(c-3*a+2*n))+cos(Pi*(a+c-2*b+2*n))+cos(Pi*(2*b+2*n-a+c))-cos(Pi*(2
*b+2*n-3*a+c))-cos(Pi*(2*n+a+c)))-1/4*WPMR(2*m+2*n-1,2*m,b-c+2*a+2*m,-
c+2*a+2-2*n-b,a+2*m)*GAMMA(-c+2*a+2-2*n)*GAMMA(b-1+2*n+2*m)*GAMMA(a+2*
m)*GAMMA(1-2*a-2*m)*GAMMA(c)*GAMMA(1+a-c)*GAMMA(-b+1)*(sin(Pi*(-3*a+2*
c-2*m+2*n))-sin(Pi*(-5*a+2*c-2*m+2*n))-sin(Pi*(a+2*b+2*m+2*n))+sin(Pi*
(-a+2*b+2*m+2*n)))/Pi^3/GAMMA(1-c+2*a+2*m)-WPMR(2*m+2*n-1,2*m,-b+3-c-2
*n-2*m,1+b-c,-a+1)*sin(Pi*c)/Pi/(c-1)/GAMMA(b-c+2*a+2*m)*GAMMA(-1+2*a+
2*m)/GAMMA(-c+2*a+2-2*n-b)/GAMMA(b-2+c+2*n+2*m)*GAMMA(-c+2*a+2-2*n)*GA
MMA(b-1+2*n+2*m)/GAMMA(c-b)*GAMMA(c)^2*GAMMA(1+a-c)*GAMMA(-b+1)/GAMMA(
a);}{%
\maplemultiline{
410\mbox{:~~~}  - \mathrm{WPMR}(2\,m + 2\,n - 1, \,2\,n
 - 1, \, - b + 3 - c - 2\,n - 2\,m, \,1 + b - c, \,1 + a - c) \\
\Gamma (2\,n - 2\,a - 2 + 2\,c)\,\Gamma ( - c + 2\,a + 2 - 2\,n)
\,\mathrm{\%1}\,\Gamma (c)^{2}(\mathrm{sin}(\pi \,(2\,c - a + 2\,
n - b)) \\
\mbox{} - \mathrm{sin}(\pi \,( - a - b + 2\,n)) - \mathrm{sin}(
\pi \,(2\,c - 3\,a + 2\,n + b)) + \mathrm{sin}(\pi \,( - 3\,a + b
 + 2\,n)))/( \\
\Gamma (b - 2 + c + 2\,n + 2\,m)\,\Gamma (c - b)\,\Gamma (c - a)
\,\Gamma (b)\,\Gamma (a) \\
(\mathrm{sin}(2\,\pi \,(a - b)) + \mathrm{sin}(2\,\pi \,b) - 
\mathrm{sin}(2\,\pi \,a))\,(c - 1))\mbox{} + \Gamma (2\,a - 2\,n
 + 1 - c)^{2}\mathrm{WPMR}( \\
2\,m + 2\,n - 1, \,2\,n - 1, \, - b + 1 + c - 2\,a - 2\,m, \,c - 
2\,a - 1 + 2\,n + b,  \\
 - a - 1 + c + 2\,n)\Gamma (2\,a + 2 - 2\,c - 2\,n)\,\mathrm{\%1}
\,\Gamma (c)(2\,\mathrm{cos}(\pi \,( - c + 2\,a - b)) \\
\mbox{} + \mathrm{cos}(\pi \,(3\,c - 4\,a + 4\,n - b)) - 2\,
\mathrm{cos}(\pi \,(b - c)) - \mathrm{cos}(\pi \,(c - 2\,a + 4\,n
 - b)) \\
\mbox{} - \mathrm{cos}(\pi \,(3\,c - 6\,a + 4\,n + b)) + \mathrm{
cos}(\pi \,(c - 4\,a + 4\,n + b)))( - 2\,a + 2\,n - 1 + c)/( \\
\Gamma (b - c + 2\,a + 2\,m)\,\Gamma ( - c + 2\,a + 2 - 2\,n - b)
\,\Gamma ( - c + 2 + a - 2\,n)\,\Gamma (b)\,\Gamma (a)( \\
 - \mathrm{cos}(\pi \,(c - a + 2\,n - 2\,b)) + \mathrm{cos}(\pi 
\,(c - 3\,a + 2\,n)) + \mathrm{cos}(\pi \,(a + c - 2\,b + 2\,n))
 \\
\mbox{} + \mathrm{cos}(\pi \,(2\,b + 2\,n - a + c)) - \mathrm{cos
}(\pi \,(2\,b + 2\,n - 3\,a + c)) - \mathrm{cos}(\pi \,(2\,n + a
 + c)))) \\
\mbox{} - {\displaystyle \frac {1}{4}} \mathrm{WPMR}(2\,m + 2\,n
 - 1, \,2\,m, \,b - c + 2\,a + 2\,m, \, - c + 2\,a + 2 - 2\,n - b
, \,a + 2\,m) \\
\Gamma ( - c + 2\,a + 2 - 2\,n)\,\mathrm{\%1}\,\Gamma (a + 2\,m)
\,\Gamma (1 - 2\,a - 2\,m)\,\Gamma (c)\,\Gamma (1 + a - c)\,
\Gamma ( - b + 1) \\
(\mathrm{sin}(\pi \,( - 3\,a + 2\,c - 2\,m + 2\,n)) - \mathrm{sin
}(\pi \,( - 5\,a + 2\,c - 2\,m + 2\,n)) \\
\mbox{} - \mathrm{sin}(\pi \,(a + 2\,b + 2\,m + 2\,n)) + \mathrm{
sin}(\pi \,( - a + 2\,b + 2\,m + 2\,n)))/(\pi ^{3} \\
\Gamma (1 - c + 2\,a + 2\,m))\mbox{} -  \\
\mathrm{WPMR}(2\,m + 2\,n - 1, \,2\,m, \, - b + 3 - c - 2\,n - 2
\,m, \,1 + b - c, \, - a + 1)\,\mathrm{sin}(\pi \,c) \\
\Gamma ( - 1 + 2\,a + 2\,m)\,\Gamma ( - c + 2\,a + 2 - 2\,n)\,
\mathrm{\%1}\,\Gamma (c)^{2}\,\Gamma (1 + a - c)\,\Gamma ( - b + 
1)/(\pi \,(c - 1) \\
\Gamma (b - c + 2\,a + 2\,m)\,\Gamma ( - c + 2\,a + 2 - 2\,n - b)
\,\Gamma (b - 2 + c + 2\,n + 2\,m)\,\Gamma (c - b) \,
\Gamma (a)) \\
\mathrm{\%1} := \Gamma (b - 1 + 2\,n + 2\,m) }
}

\begin{mapleinput}
\end{mapleinput}

\end{maplegroup}

%% file: AppendixB411to443.tex
\begin{maplegroup}
\mapleinline{inert}{2d}{411, ":   ",
WPMR(2*m+2*n-1,2*n,-a+1,-c+1,-1/2*a+b-1/2*c+n+m)*Pi*GAMMA(a-2*b+c+1-2*
m)*GAMMA(-1/2*c+1/2*a+n+m)*GAMMA(2*b-1+2*m)*GAMMA(1/2*c+1/2*a-n-m)^2*(
sin(1/2*Pi*(-a+c+6*n+2*m-2*b))-sin(1/2*Pi*(a+3*c+2*n-2*m-2*b))-sin(1/2
*Pi*(a-c+2*n-2*m-6*b))-sin(1/2*Pi*(-3*a-c+2*n+6*m+6*b)))*(-a-c+2*n+2*m
)/GAMMA(1+1/2*a+1/2*c-b+n-m)/GAMMA(-1/2*c+1/2*a-n-m+1)/GAMMA(1/2*a-b+1
/2*c-n-m+1)/GAMMA(-1+2*m+b)/GAMMA(c)/GAMMA(b)/GAMMA(a)/(-cos(1/2*Pi*(2
*b+6*m-3*a+c+2*n))+cos(1/2*Pi*(6*b+6*m-3*a+c+2*n))+cos(1/2*Pi*(2*b-2*m
+a+c+2*n))+cos(1/2*Pi*(6*n-2*b+3*c-a+2*m))-cos(1/2*Pi*(2*n-6*b+c+a-2*m
))-cos(1/2*Pi*(6*n+2*b+3*c-a+2*m)))+WPMR(2*m+2*n-1,2*n,a-2*b+2-2*m,c+2
-2*b-2*m,1+1/2*a+1/2*c-b+n-m)*GAMMA(-a+2*b-c-1+2*m)*GAMMA(a-2*b+2-2*m)
*GAMMA(1/2*c+1/2*a-n-m+1)*GAMMA(-1/2*c+1/2*a+n+m)*GAMMA(2*b-1+2*m)*(co
s(1/2*Pi*(2*n-8*b+3*c+3*a-6*m))+cos(1/2*Pi*(6*n+4*b+c-3*a+6*m))-cos(1/
2*Pi*(2*n+3*c-a+2*m))-cos(1/2*Pi*(6*n-4*b+c+a-2*m))-cos(1/2*Pi*(-5*a+1
2*b+10*m-c+2*n))-cos(1/2*Pi*(-a-c+2*n+2*m))+cos(1/2*Pi*(-a+4*b+2*m-c+2
*n))+cos(1/2*Pi*(3*a-8*b-6*m-c+2*n)))/GAMMA(-c+2*b-1+2*m)/GAMMA(2+1/2*
c+1/2*a+n-m-2*b)/GAMMA(-1/2*c+1/2*a-n-m+1)/GAMMA(-1+2*m+b)/GAMMA(b)/(-
cos(1/2*Pi*(2*b+6*m-3*a+c+2*n))+cos(1/2*Pi*(6*b+6*m-3*a+c+2*n))+cos(1/
2*Pi*(2*b-2*m+a+c+2*n))+cos(1/2*Pi*(6*n-2*b+3*c-a+2*m))-cos(1/2*Pi*(2*
n-6*b+c+a-2*m))-cos(1/2*Pi*(6*n+2*b+3*c-a+2*m)))+WPMR(2*m+2*n-1,2*m-1,
-a+1,-c+1,-b+1)*sin(1/2*Pi*(-a-c+2*n+2*m))*GAMMA(-2+2*b+2*m)^2*GAMMA(1
/2*c+1/2*a-n-m)^2*GAMMA(-1/2*a+b-1/2*c+n+m)*GAMMA(-1/2*a-1/2*c+b-n+m)*
GAMMA(1/2*c-1/2*a+n+m)*GAMMA(-1/2*c+1/2*a+n+m)*(-1+b+m)*(-a-c+2*n+2*m)
/Pi/GAMMA(-c+2*b-1+2*m)/GAMMA(-a+2*b-1+2*m)/GAMMA(-1+2*m+b)/GAMMA(c)/G
AMMA(b)/GAMMA(a)+1/4*WPMR(2*m+2*n-1,2*m-1,-a+2*b-1+2*m,-c+2*b-1+2*m,-1
+2*m+b)*GAMMA(1/2*c-1/2*a+n+m)*GAMMA(1/2*c+1/2*a-n-m+1)*GAMMA(-1/2*a+b
-1/2*c+n+m)*GAMMA(-1/2*a-1/2*c+b-n+m)*GAMMA(-1/2*c+1/2*a+n+m)*(sin(Pi*
(-c+a-b))+sin(Pi*(a-3*b-4*m+c))-sin(Pi*(b-c+a))-sin(Pi*(a-5*b-4*m+c)))
/sin(2*Pi*(b+m))/Pi^2/GAMMA(-1/2*c-1/2*a+n+3*m+2*b-1);}{%
\maplemultiline{
411\mbox{:~~~} \mathrm{WPMR}(2\,m + 2\,n - 1, \,2\,n, \,
 - a + 1, \, - c + 1, \, - {\displaystyle \frac {a}{2}}  + b - 
{\displaystyle \frac {c}{2}}  + n + m)\,\pi  \\
\Gamma (a - 2\,b + c + 1 - 2\,m)\,\mathrm{\%1}\,\Gamma (2\,b - 1
 + 2\,m)\,\Gamma ({\displaystyle \frac {c}{2}}  + {\displaystyle 
\frac {a}{2}}  - n - m)^{2}( \\
\mathrm{sin}({\displaystyle \frac {\pi \,( - a + c + 6\,n + 2\,m
 - 2\,b)}{2}} ) - \mathrm{sin}({\displaystyle \frac {\pi \,(a + 3
\,c + 2\,n - 2\,m - 2\,b)}{2}} ) \\
\mbox{} - \mathrm{sin}({\displaystyle \frac {\pi \,(a - c + 2\,n
 - 2\,m - 6\,b)}{2}} ) - \mathrm{sin}({\displaystyle \frac {\pi 
\,( - 3\,a - c + 2\,n + 6\,m + 6\,b)}{2}} )) \\
( - a - c + 2\,n + 2\,m) \left/ {\vrule 
height0.80em width0em depth0.80em} \right. \!  \! (\Gamma (1 + 
{\displaystyle \frac {a}{2}}  + {\displaystyle \frac {c}{2}}  - b
 + n - m)\,\Gamma ( - {\displaystyle \frac {c}{2}}  + 
{\displaystyle \frac {a}{2}}  - n - m + 1) \\
\Gamma ({\displaystyle \frac {a}{2}}  - b + {\displaystyle 
\frac {c}{2}}  - n - m + 1)\,\Gamma ( - 1 + 2\,m + b)\,\Gamma (c)
\,\Gamma (b)\,\Gamma (a)( \\
 - \mathrm{cos}({\displaystyle \frac {\pi \,(2\,b + 6\,m - 3\,a
 + c + 2\,n)}{2}} ) + \mathrm{cos}({\displaystyle \frac {\pi \,(6
\,b + 6\,m - 3\,a + c + 2\,n)}{2}} ) \\
\mbox{} + \mathrm{cos}({\displaystyle \frac {\pi \,(2\,b - 2\,m
 + a + c + 2\,n)}{2}} ) + \mathrm{cos}({\displaystyle \frac {\pi 
\,(6\,n - 2\,b + 3\,c - a + 2\,m)}{2}} ) \\
\mbox{} - \mathrm{cos}({\displaystyle \frac {\pi \,(2\,n - 6\,b
 + c + a - 2\,m)}{2}} ) - \mathrm{cos}({\displaystyle \frac {\pi 
\,(6\,n + 2\,b + 3\,c - a + 2\,m)}{2}} )))\mbox{} +  \\
\mathrm{WPMR}(2\,m + 2\,n - 1, \,2\,n, \,a - 2\,b + 2 - 2\,m, \,c
 + 2 - 2\,b - 2\,m,  \\
1 + {\displaystyle \frac {a}{2}}  + {\displaystyle \frac {c}{2}} 
 - b + n - m)\Gamma ( - a + 2\,b - c - 1 + 2\,m)\,\Gamma (a - 2\,
b + 2 - 2\,m) \\
\Gamma ({\displaystyle \frac {c}{2}}  + {\displaystyle \frac {a}{
2}}  - n - m + 1)\,\mathrm{\%1}\,\Gamma (2\,b - 1 + 2\,m)(
\mathrm{cos}({\displaystyle \frac {\pi \,(2\,n - 8\,b + 3\,c + 3
\,a - 6\,m)}{2}} ) \\
\mbox{} + \mathrm{cos}({\displaystyle \frac {\pi \,(6\,n + 4\,b
 + c - 3\,a + 6\,m)}{2}} ) - \mathrm{cos}({\displaystyle \frac {
\pi \,(2\,n + 3\,c - a + 2\,m)}{2}} ) \\
\mbox{} - \mathrm{cos}({\displaystyle \frac {\pi \,(6\,n - 4\,b
 + c + a - 2\,m)}{2}} ) - \mathrm{cos}({\displaystyle \frac {\pi 
\,( - 5\,a + 12\,b + 10\,m - c + 2\,n)}{2}} ) \\
\mbox{} - \mathrm{cos}({\displaystyle \frac {\pi \,( - a - c + 2
\,n + 2\,m)}{2}} ) + \mathrm{cos}({\displaystyle \frac {\pi \,(
 - a + 4\,b + 2\,m - c + 2\,n)}{2}} ) \\
\mbox{} + \mathrm{cos}({\displaystyle \frac {\pi \,(3\,a - 8\,b
 - 6\,m - c + 2\,n)}{2}} )) \left/ {\vrule 
height0.80em width0em depth0.80em} \right. \!  \! (\Gamma ( - c
 + 2\,b - 1 + 2\,m) \\
\Gamma (2 + {\displaystyle \frac {c}{2}}  + {\displaystyle 
\frac {a}{2}}  + n - m - 2\,b)\,\Gamma ( - {\displaystyle \frac {
c}{2}}  + {\displaystyle \frac {a}{2}}  - n - m + 1)\,\Gamma ( - 
1 + 2\,m + b)\,\Gamma (b)( \\
 - \mathrm{cos}({\displaystyle \frac {\pi \,(2\,b + 6\,m - 3\,a
 + c + 2\,n)}{2}} ) + \mathrm{cos}({\displaystyle \frac {\pi \,(6
\,b + 6\,m - 3\,a + c + 2\,n)}{2}} ) \\
\mbox{} + \mathrm{cos}({\displaystyle \frac {\pi \,(2\,b - 2\,m
 + a + c + 2\,n)}{2}} ) + \mathrm{cos}({\displaystyle \frac {\pi 
\,(6\,n - 2\,b + 3\,c - a + 2\,m)}{2}} ) \\
\mbox{} - \mathrm{cos}({\displaystyle \frac {\pi \,(2\,n - 6\,b
 + c + a - 2\,m)}{2}} ) - \mathrm{cos}({\displaystyle \frac {\pi 
\,(6\,n + 2\,b + 3\,c - a + 2\,m)}{2}} )))\mbox{} +  \\
\mathrm{WPMR}(2\,m + 2\,n - 1, \,2\,m - 1, \, - a + 1, \, - c + 1
, \, - b + 1)\,\mathrm{sin}({\displaystyle \frac {\pi \,( - a - c
 + 2\,n + 2\,m)}{2}} ) \\
\Gamma ( - 2 + 2\,b + 2\,m)^{2}\,\Gamma ({\displaystyle \frac {c
}{2}}  + {\displaystyle \frac {a}{2}}  - n - m)^{2}\,\Gamma ( - 
{\displaystyle \frac {a}{2}}  + b - {\displaystyle \frac {c}{2}} 
 + n + m) \\
\Gamma ( - {\displaystyle \frac {a}{2}}  - {\displaystyle \frac {
c}{2}}  + b - n + m)\,\Gamma ({\displaystyle \frac {c}{2}}  - 
{\displaystyle \frac {a}{2}}  + n + m)\,\mathrm{\%1}\,( - 1 + b
 + m)\,( - a - c + 2\,n + 2\,m)/(\pi  \\
\Gamma ( - c + 2\,b - 1 + 2\,m)\,\Gamma ( - a + 2\,b - 1 + 2\,m)
\,\Gamma ( - 1 + 2\,m + b)\,\Gamma (c)\,\Gamma (b)\,\Gamma (a))
\mbox{} + {\displaystyle \frac {1}{4}}  \\
\mathrm{WPMR}(2\,m + 2\,n - 1, \,2\,m - 1, \, - a + 2\,b - 1 + 2
\,m, \, - c + 2\,b - 1 + 2\,m,  \\
 - 1 + 2\,m + b)\Gamma ({\displaystyle \frac {c}{2}}  - 
{\displaystyle \frac {a}{2}}  + n + m)\,\Gamma ({\displaystyle 
\frac {c}{2}}  + {\displaystyle \frac {a}{2}}  - n - m + 1)\,
\Gamma ( - {\displaystyle \frac {a}{2}}  + b - {\displaystyle 
\frac {c}{2}}  + n + m) \\
\Gamma ( - {\displaystyle \frac {a}{2}}  - {\displaystyle \frac {
c}{2}}  + b - n + m)\,\mathrm{\%1}(\mathrm{sin}(\pi \,( - c + a
 - b)) + \mathrm{sin}(\pi \,(a - 3\,b - 4\,m + c)) \\
\mbox{} - \mathrm{sin}(\pi \,(b - c + a)) - \mathrm{sin}(\pi \,(a
 - 5\,b - 4\,m + c))) \left/ {\vrule 
height0.80em width0em depth0.80em} \right. \!  \! (\mathrm{sin}(2
\,\pi \,(b + m))\,\pi ^{2} \\
\Gamma ( - {\displaystyle \frac {c}{2}}  - {\displaystyle \frac {
a}{2}}  + n + 3\,m + 2\,b - 1)) \\
\mathrm{\%1} := \Gamma ( - {\displaystyle \frac {c}{2}}  + 
{\displaystyle \frac {a}{2}}  + n + m) }
}

\mapleresult
\begin{maplelatex}
\mapleinline{inert}{2d}{412, ":   ",
(cos(Pi*(c+b-3*a+2*n+4*m))-cos(Pi*(b+2*n+a-c)))*WPMR(2*m+2*n-1,2*n,-c+
2-b-2*n,-b+c+2*m-2*a,-b+1)*GAMMA(c)*GAMMA(1+a-c)*GAMMA(-1+2*b+2*n)*GAM
MA(2*a-c-2*n-2*m+2)*GAMMA(-a+c+2*n+2*m-1)/GAMMA(-b+1-a+2*m)/GAMMA(b-c+
2*a-2*m+1)/GAMMA(b)/GAMMA(c+b-1+2*n)/GAMMA(a)/(cos(Pi*(-c+a))-cos(Pi*(
3*a-4*m-c)))+1/2*WPMR(2*m+2*n-1,2*n,b+2*n-1-2*a+c+2*m,1+b-c,b+2*n)*GAM
MA(-2*b+1-2*n)*GAMMA(b+2*n)*GAMMA(2*a-c-2*n-2*m+2)*GAMMA(-a+c+2*n+2*m-
1)*GAMMA(c)*GAMMA(1+a-c)*(sin(Pi*(2*b+c-3*a+2*n+4*m))-sin(Pi*(c-3*a+2*
n+4*m))-sin(Pi*(2*b+2*n+a-c))+sin(Pi*(a-c+2*n)))/Pi/GAMMA(-a+b+2*n+2*m
)/GAMMA(-2*n-b+2+2*a-c-2*m)/GAMMA(c-b)/GAMMA(a)/(cos(Pi*(-c+a))-cos(Pi
*(3*a-4*m-c)))+GAMMA(-a+c+2*n+2*m-1)*GAMMA(2*a-c-2*n-2*m+2)*GAMMA(1+a-
c)*GAMMA(-b+1)/GAMMA(-2*n-b+2+2*a-c-2*m)*GAMMA(-2*m+2*a)/GAMMA(b-c+2*a
-2*m+1)/GAMMA(c-b)/GAMMA(-b+1-a+2*m)/GAMMA(c+b-1+2*n)/GAMMA(a)*GAMMA(c
)*WPMR(2*m+2*n-1,2*m-1,-c+2-b-2*n,-b+c+2*m-2*a,-a+2*m)+1/2*WPMR(2*m+2*
n-1,2*m-1,-2*n-b+2+2*a-c-2*m,c-b,a)*GAMMA(2*a-c-2*n-2*m+2)*GAMMA(-a+c+
2*n+2*m-1)*GAMMA(-2*a+2*m)*GAMMA(c)*GAMMA(1+a-c)*GAMMA(-b+1)*(cos(2*Pi
*(c-a+n+m))-cos(2*Pi*(b-a+n+m)))/Pi^2/GAMMA(-a+2*m)/GAMMA(1+a-b);}{%
\maplemultiline{
412\mbox{:~~~} (\mathrm{cos}(\pi \,(c + b - 3\,a + 2\,n
 + 4\,m)) - \mathrm{cos}(\pi \,(b + 2\,n + a - c))) \\
\mathrm{WPMR}(2\,m + 2\,n - 1, \,2\,n, \, - c + 2 - b - 2\,n, \,
 - b + c + 2\,m - 2\,a, \, - b + 1)\,\Gamma (c) \\
\Gamma (1 + a - c)\,\Gamma ( - 1 + 2\,b + 2\,n)\,\mathrm{\%2}\,
\mathrm{\%1}/(\Gamma ( - b + 1 - a + 2\,m) \\
\Gamma (b - c + 2\,a - 2\,m + 1)\,\Gamma (b)\,\Gamma (c + b - 1
 + 2\,n)\,\Gamma (a) \\
(\mathrm{cos}(\pi \,( - c + a)) - \mathrm{cos}(\pi \,(3\,a - 4\,m
 - c))))\mbox{} + {\displaystyle \frac {1}{2}}  \\
\mathrm{WPMR}(2\,m + 2\,n - 1, \,2\,n, \,b + 2\,n - 1 - 2\,a + c
 + 2\,m, \,1 + b - c, \,b + 2\,n) \\
\Gamma ( - 2\,b + 1 - 2\,n)\,\Gamma (b + 2\,n)\,\mathrm{\%2}\,
\mathrm{\%1}\,\Gamma (c)\,\Gamma (1 + a - c)( \\
\mathrm{sin}(\pi \,(2\,b + c - 3\,a + 2\,n + 4\,m)) - \mathrm{sin
}(\pi \,(c - 3\,a + 2\,n + 4\,m)) \\
\mbox{} - \mathrm{sin}(\pi \,(2\,b + 2\,n + a - c)) + \mathrm{sin
}(\pi \,(a - c + 2\,n)))/(\pi \,\Gamma ( - a + b + 2\,n + 2\,m)
 \\
\Gamma ( - 2\,n - b + 2 + 2\,a - c - 2\,m)\,\Gamma (c - b)\,
\Gamma (a) \\
(\mathrm{cos}(\pi \,( - c + a)) - \mathrm{cos}(\pi \,(3\,a - 4\,m
 - c))))\mbox{} + \mathrm{\%1}\,\mathrm{\%2}\,\Gamma (1 + a - c)
\,\Gamma ( - b + 1) \\
\Gamma ( - 2\,m + 2\,a)\,\Gamma (c) \\
\mathrm{WPMR}(2\,m + 2\,n - 1, \,2\,m - 1, \, - c + 2 - b - 2\,n
, \, - b + c + 2\,m - 2\,a, \, - a + 2\,m)/( \\
\Gamma ( - 2\,n - b + 2 + 2\,a - c - 2\,m)\,\Gamma (b - c + 2\,a
 - 2\,m + 1)\,\Gamma (c - b) \\
\Gamma ( - b + 1 - a + 2\,m)\,\Gamma (c + b - 1 + 2\,n)\,\Gamma (
a))\mbox{} + {\displaystyle \frac {1}{2}}  \\
\mathrm{WPMR}(2\,m + 2\,n - 1, \,2\,m - 1, \, - 2\,n - b + 2 + 2
\,a - c - 2\,m, \,c - b, \,a)\,\mathrm{\%2}\,\mathrm{\%1} \\
\Gamma ( - 2\,a + 2\,m)\,\Gamma (c)\,\Gamma (1 + a - c)\,\Gamma (
 - b + 1) \\
(\mathrm{cos}(2\,\pi \,(c - a + n + m)) - \mathrm{cos}(2\,\pi \,(
b - a + n + m)))/(\pi ^{2}\,\Gamma ( - a + 2\,m)\,\Gamma (1 + a
 - b)) \\
\mathrm{\%1} := \Gamma ( - a + c + 2\,n + 2\,m - 1) \\
\mathrm{\%2} := \Gamma (2\,a - c - 2\,n - 2\,m + 2) }
}
\end{maplelatex}

\begin{maplelatex}
\mapleinline{inert}{2d}{413, ":   ",
(GAMMA(b-1+2*n+2*m)*GAMMA(2*a-2*n+1-c)*GAMMA(b-a)*GAMMA(-a+1)*GAMMA(1+
a-b)*GAMMA(-1+2*c-2*a+2*n)/GAMMA(c-b)/GAMMA(-c+1+a-2*n)/GAMMA(c-a)/GAM
MA(b)/GAMMA(-a+c+2*n)/GAMMA(2-c)/GAMMA(b-2+c+2*n+2*m)*GAMMA(c)-GAMMA(b
-1+2*n+2*m)*GAMMA(2*a-2*n+1-c)*GAMMA(-b+1)*GAMMA(1+a-b)*GAMMA(b-a)*GAM
MA(-1+2*c-2*a+2*n)/GAMMA(2*a-b-c+1-2*n)/GAMMA(c-b)/GAMMA(c-a)/GAMMA(2-
c)/GAMMA(b-2+c+2*n+2*m)/GAMMA(-2*a+b+c+2*n)*GAMMA(c)/GAMMA(a))*WPMR(2*
m+2*n-1,2*n,-b+3-c-2*n-2*m,1+b-c,1+a-c)+(GAMMA(b-1+2*n+2*m)*GAMMA(2*a-
2*n+1-c)*GAMMA(b-a)/GAMMA(2*a-b-c+1-2*n)*GAMMA(-a+1)*GAMMA(1+a-b)*GAMM
A(-2*c+1+2*a-2*n)/GAMMA(-c+1+a-2*n)/GAMMA(c-a)/GAMMA(b)/GAMMA(1+a-c)/G
AMMA(-1+b-c+2*a+2*m)/GAMMA(1-2*a+c+2*n)*GAMMA(c)-GAMMA(b-1+2*n+2*m)*GA
MMA(2*a-2*n+1-c)*GAMMA(-b+1)*GAMMA(-a+c+2*n)*GAMMA(1+a-b)*GAMMA(b-a)*G
AMMA(-2*c+1+2*a-2*n)/GAMMA(2*a-b-c+1-2*n)^2/GAMMA(c-a)/GAMMA(1+a-c)/GA
MMA(-1+b-c+2*a+2*m)/GAMMA(1-2*a+c+2*n)/GAMMA(-2*a+b+c+2*n)*GAMMA(c)/GA
MMA(a))*WPMR(2*m+2*n-1,2*n,-2*a+b+c+2*n,-b+2+c-2*a-2*m,-a+c+2*n)+GAMMA
(b-1+2*n+2*m)*GAMMA(2*a-2*n+1-c)*GAMMA(-b+1)*GAMMA(1+a-c)/GAMMA(2*a-b-
c+1-2*n)*GAMMA(2*m-2+2*a)/GAMMA(a)/GAMMA(c-b)/GAMMA(-1+b-c+2*a+2*m)/GA
MMA(2-c)/GAMMA(b-2+c+2*n+2*m)*GAMMA(c)*WPMR(2*m+2*n-1,2*m-1,-b+3-c-2*n
-2*m,1+b-c,-a+1)+GAMMA(b-1+2*n+2*m)*GAMMA(2*a-2*n+1-c)*GAMMA(-b+1)*GAM
MA(1+a-c)/GAMMA(2*a-b-c+1-2*n)*GAMMA(2-2*a-2*m)*WPMR(2*m+2*n-1,2*m-1,2
*a-b-c+1-2*n,-1+b-c+2*a+2*m,-1+a+2*m)/GAMMA(a)/GAMMA(-c+2*a+2*m)/GAMMA
(-1+b-c+2*a+2*m)/GAMMA(-a+1)/GAMMA(-2*a+b+c+2*n)/GAMMA(-b+2+c-2*a-2*m)
*GAMMA(-1+a+2*m)*GAMMA(c);}{%
\maplemultiline{
413\mbox{:~~~} (\mathrm{\%4}\,\mathrm{\%3}\,\Gamma (b - 
a)\,\Gamma ( - a + 1)\,\Gamma (1 + a - b)\,\Gamma ( - 1 + 2\,c - 
2\,a + 2\,n)\,\Gamma (c)/(\Gamma (c - b) \\
\Gamma ( - c + 1 + a - 2\,n)\,\Gamma (c - a)\,\Gamma (b)\,\Gamma 
( - a + c + 2\,n)\,\Gamma (2 - c)\,\Gamma (b - 2 + c + 2\,n + 2\,
m) \\
)\mbox{} - {\displaystyle \frac {\mathrm{\%4}\,\mathrm{\%3}\,
\Gamma ( - b + 1)\,\Gamma (1 + a - b)\,\Gamma (b - a)\,\Gamma (
 - 1 + 2\,c - 2\,a + 2\,n)\,\Gamma (c)}{\Gamma (\mathrm{\%2})\,
\Gamma (c - b)\,\Gamma (c - a)\,\Gamma (2 - c)\,\Gamma (b - 2 + c
 + 2\,n + 2\,m)\,\Gamma ( - 2\,a + b + c + 2\,n)\,\Gamma (a)}} 
 \\
)\mathrm{WPMR}(2\,m + 2\,n - 1, \,2\,n, \, - b + 3 - c - 2\,n - 2
\,m, \,1 + b - c, \,1 + a - c)\mbox{} + ( \\
{\displaystyle \frac {\mathrm{\%4}\,\mathrm{\%3}\,\Gamma (b - a)
\,\Gamma ( - a + 1)\,\Gamma (1 + a - b)\,\Gamma ( - 2\,c + 1 + 2
\,a - 2\,n)\,\Gamma (c)}{\Gamma (\mathrm{\%2})\,\Gamma ( - c + 1
 + a - 2\,n)\,\Gamma (c - a)\,\Gamma (b)\,\Gamma (1 + a - c)\,
\Gamma (\mathrm{\%1})\,\Gamma (1 - 2\,a + c + 2\,n)}}  \\
\mbox{} - {\displaystyle \frac {\mathrm{\%4}\,\mathrm{\%3}\,
\Gamma ( - b + 1)\,\Gamma ( - a + c + 2\,n)\,\Gamma (1 + a - b)\,
\Gamma (b - a)\,\Gamma ( - 2\,c + 1 + 2\,a - 2\,n)\,\Gamma (c)}{
\Gamma (\mathrm{\%2})^{2}\,\Gamma (c - a)\,\Gamma (1 + a - c)\,
\Gamma (\mathrm{\%1})\,\Gamma (1 - 2\,a + c + 2\,n)\,\Gamma ( - 2
\,a + b + c + 2\,n)\,\Gamma (a)}}  \\
)\mathrm{WPMR}(2\,m + 2\,n - 1, \,2\,n, \, - 2\,a + b + c + 2\,n
, \, - b + 2 + c - 2\,a - 2\,m,  \\
 - a + c + 2\,n)\mbox{} + \mathrm{\%4}\,\mathrm{\%3}\,\Gamma ( - 
b + 1)\,\Gamma (1 + a - c)\,\Gamma (2\,m - 2 + 2\,a)\,\Gamma (c)
 \\
\mathrm{WPMR}(2\,m + 2\,n - 1, \,2\,m - 1, \, - b + 3 - c - 2\,n
 - 2\,m, \,1 + b - c, \, - a + 1)/(\Gamma (\mathrm{\%2}) \\
\Gamma (a)\,\Gamma (c - b)\,\Gamma (\mathrm{\%1})\,\Gamma (2 - c)
\,\Gamma (b - 2 + c + 2\,n + 2\,m))\mbox{} + \mathrm{\%4}\,
\mathrm{\%3}\,\Gamma ( - b + 1) \\
\Gamma (1 + a - c)\,\Gamma (2 - 2\,a - 2\,m) \\
\mathrm{WPMR}(2\,m + 2\,n - 1, \,2\,m - 1, \,\mathrm{\%2}, \,
\mathrm{\%1}, \, - 1 + a + 2\,m)\,\Gamma ( - 1 + a + 2\,m)\,
\Gamma (c)/( \\
\Gamma (\mathrm{\%2})\,\Gamma (a)\,\Gamma ( - c + 2\,a + 2\,m)\,
\Gamma (\mathrm{\%1})\,\Gamma ( - a + 1)\,\Gamma ( - 2\,a + b + c
 + 2\,n) \\
\Gamma ( - b + 2 + c - 2\,a - 2\,m)) \\
\mathrm{\%1} :=  - 1 + b - c + 2\,a + 2\,m \\
\mathrm{\%2} := 2\,a - b - c + 1 - 2\,n \\
\mathrm{\%3} := \Gamma (2\,a - 2\,n + 1 - c) \\
\mathrm{\%4} := \Gamma (b - 1 + 2\,n + 2\,m) }
}
\end{maplelatex}

\begin{maplelatex}
\mapleinline{inert}{2d}{414, ":   ",
(-GAMMA(-1/6+n-a)*GAMMA(a+2/3)*GAMMA(1/6+n-a)*GAMMA(1/2+n-a)*GAMMA(3*a
-2*n+1/2)*V1(-1/6+n-a,n)/GAMMA(-1/2-3*a+2*n)/GAMMA(2*a+1/3-n)/GAMMA(2*
a+2/3-n)/GAMMA(2*a+1-n)/Pi^(1/2)+GAMMA(-1/6+n-a)*GAMMA(a+2/3)*GAMMA(1/
6+n-a)*GAMMA(1/2+n-a)*V1(a,n)/GAMMA(-1/2-3*a+2*n)/Pi/GAMMA(3*a+1-n))*G
AMMA(2*a+5/6-n)/GAMMA(1/2+n);}{%
\maplemultiline{
414\mbox{:~~~}  \left( {\vrule 
height1.67em width0em depth1.67em} \right. \!  \!  - 
{\displaystyle \frac {\Gamma ( - {\displaystyle \frac {1}{6}}  + 
n - a)\,\Gamma (a + {\displaystyle \frac {2}{3}} )\,\Gamma (
{\displaystyle \frac {1}{6}}  + n - a)\,\Gamma ({\displaystyle 
\frac {1}{2}}  + n - a)\,\Gamma (3\,a - 2\,n + {\displaystyle 
\frac {1}{2}} )\,\mathrm{V1}( - {\displaystyle \frac {1}{6}}  + n
 - a, \,n)}{\Gamma ( - {\displaystyle \frac {1}{2}}  - 3\,a + 2\,
n)\,\Gamma (2\,a + {\displaystyle \frac {1}{3}}  - n)\,\Gamma (2
\,a + {\displaystyle \frac {2}{3}}  - n)\,\Gamma (2\,a + 1 - n)\,
\sqrt{\pi }}}  \\
\mbox{} + {\displaystyle \frac {\Gamma ( - {\displaystyle \frac {
1}{6}}  + n - a)\,\Gamma (a + {\displaystyle \frac {2}{3}} )\,
\Gamma ({\displaystyle \frac {1}{6}}  + n - a)\,\Gamma (
{\displaystyle \frac {1}{2}}  + n - a)\,\mathrm{V1}(a, \,n)}{
\Gamma ( - {\displaystyle \frac {1}{2}}  - 3\,a + 2\,n)\,\pi \,
\Gamma (3\,a + 1 - n)}}  \! \! \left. {\vrule 
height1.67em width0em depth1.67em} \right) \Gamma (2\,a + 
{\displaystyle \frac {5}{6}}  - n) \left/ {\vrule 
height0.80em width0em depth0.80em} \right. \!  \!
\Gamma ({\displaystyle \frac {1}{2}}  + n) }
}
\end{maplelatex}

\begin{maplelatex}
\mapleinline{inert}{2d}{415, ":   ",
(-GAMMA(a+1/3)*GAMMA(1/6+n-a)*GAMMA(-1/6+n-a)*GAMMA(1/2+n-a)*GAMMA(7/6
+2*a-n)*GAMMA(3*a-2*n+1/2)*V1(-1/6+n-a,n)/GAMMA(-1/2-3*a+2*n)/GAMMA(2*
a+2/3-n)/GAMMA(2*a+1/3-n)/GAMMA(2*a+1-n)/Pi^(1/2)+GAMMA(a+1/3)*GAMMA(1
/6+n-a)*GAMMA(-1/6+n-a)*GAMMA(1/2+n-a)*GAMMA(7/6+2*a-n)*V1(a,n)/GAMMA(
-1/2-3*a+2*n)/Pi/GAMMA(3*a+1-n))/GAMMA(1/2+n);}{%
\maplemultiline{
415\mbox{:~~~}  \left( {\vrule 
height1.67em width0em depth1.67em} \right. \!  \!  - \Gamma (a + 
{\displaystyle \frac {1}{3}} )\,\Gamma ({\displaystyle \frac {1}{
6}}  + n - a)\,\Gamma ( - {\displaystyle \frac {1}{6}}  + n - a)
\,\Gamma ({\displaystyle \frac {1}{2}}  + n - a)\,\Gamma (
{\displaystyle \frac {7}{6}}  + 2\,a - n) \\
\Gamma (3\,a - 2\,n + {\displaystyle \frac {1}{2}} )\,\mathrm{V1}
( - {\displaystyle \frac {1}{6}}  + n - a, \,n) \left/ {\vrule 
height0.80em width0em depth0.80em} \right. \!  \! (\Gamma ( - 
{\displaystyle \frac {1}{2}}  - 3\,a + 2\,n)\,\Gamma (2\,a + 
{\displaystyle \frac {2}{3}}  - n) \\
\Gamma (2\,a + {\displaystyle \frac {1}{3}}  - n)\,\Gamma (2\,a
 + 1 - n)\,\sqrt{\pi }) \\
\mbox{} + {\displaystyle \frac {\Gamma (a + {\displaystyle 
\frac {1}{3}} )\,\Gamma ({\displaystyle \frac {1}{6}}  + n - a)\,
\Gamma ( - {\displaystyle \frac {1}{6}}  + n - a)\,\Gamma (
{\displaystyle \frac {1}{2}}  + n - a)\,\Gamma ({\displaystyle 
\frac {7}{6}}  + 2\,a - n)\,\mathrm{V1}(a, \,n)}{\Gamma ( - 
{\displaystyle \frac {1}{2}}  - 3\,a + 2\,n)\,\pi \,\Gamma (3\,a
 + 1 - n)}}  \! \! \left. {\vrule 
height1.67em width0em depth1.67em} \right)  \left/ {\vrule 
height0.80em width0em depth0.80em} \right. \!  \! 
\Gamma ({\displaystyle \frac {1}{2}}  + n) }
}
\end{maplelatex}

\begin{maplelatex}
\mapleinline{inert}{2d}{416, ":   ",
(-2/9*GAMMA(1/2+n-a)*GAMMA(-1/6+n-a)*GAMMA(a+1/3)*3^(1/2)/GAMMA(2/3)*G
AMMA(a+2/3-n)*V1(-1/6+n-a,n)*GAMMA(1/3+n-a)/GAMMA(-a-2/3+n)/GAMMA(2*a+
2/3-n)/GAMMA(2*a+1-n)/GAMMA(-3*a+2*n+1/2)+2/9*GAMMA(1/2+n-a)*GAMMA(-1/
6+n-a)*GAMMA(a+1/3)*Pi^(1/2)*3^(1/2)/GAMMA(2/3)*V1(a,n)/GAMMA(-a-2/3+n
)/GAMMA(2/3+n-2*a)/GAMMA(3*a+1-n))*GAMMA(a)/GAMMA(1/2+n);}{%
\maplemultiline{
416\mbox{:~~~}  \left( {\vrule 
height1.67em width0em depth1.67em} \right. \!  \!  - 
{\displaystyle \frac {2}{9}} \,{\displaystyle \frac {\Gamma (
{\displaystyle \frac {1}{2}}  + n - a)\,\Gamma ( - 
{\displaystyle \frac {1}{6}}  + n - a)\,\Gamma (a + 
{\displaystyle \frac {1}{3}} )\,\sqrt{3}\,\Gamma (a + 
{\displaystyle \frac {2}{3}}  - n)\,\mathrm{V1}( - 
{\displaystyle \frac {1}{6}}  + n - a, \,n)\,\Gamma (
{\displaystyle \frac {1}{3}}  + n - a)}{\Gamma ({\displaystyle 
\frac {2}{3}} )\,\Gamma ( - a - {\displaystyle \frac {2}{3}}  + n
)\,\Gamma (2\,a + {\displaystyle \frac {2}{3}}  - n)\,\Gamma (2\,
a + 1 - n)\,\Gamma ( - 3\,a + 2\,n + {\displaystyle \frac {1}{2}
} )}}  \\
\mbox{} + {\displaystyle \frac {2}{9}} \,{\displaystyle \frac {
\Gamma ({\displaystyle \frac {1}{2}}  + n - a)\,\Gamma ( - 
{\displaystyle \frac {1}{6}}  + n - a)\,\Gamma (a + 
{\displaystyle \frac {1}{3}} )\,\sqrt{\pi }\,\sqrt{3}\,\mathrm{V1
}(a, \,n)}{\Gamma ({\displaystyle \frac {2}{3}} )\,\Gamma ( - a
 - {\displaystyle \frac {2}{3}}  + n)\,\Gamma ({\displaystyle 
\frac {2}{3}}  + n - 2\,a)\,\Gamma (3\,a + 1 - n)}}  \! 
\! \left. {\vrule height1.67em width0em depth1.67em} \right) 
\Gamma (a) \left/ {\vrule height0.80em width0em depth0.80em}
 \right. \!  \! \Gamma ({\displaystyle \frac {1}{2}}  + n) }
}
\end{maplelatex}

\begin{maplelatex}
\mapleinline{inert}{2d}{417, ":   ",
(-GAMMA(1/2+n-a)*GAMMA(-1/6+n-a)*GAMMA(1/6+n-a)*GAMMA(2*a+3/2-n)*GAMMA
(3*a-2*n+1/2)*V1(-1/6+n-a,n)/GAMMA(-1/2-3*a+2*n)/GAMMA(2*a+1-n)/GAMMA(
2*a+1/3-n)/GAMMA(2*a+2/3-n)/Pi^(1/2)+GAMMA(1/2+n-a)*GAMMA(-1/6+n-a)*GA
MMA(1/6+n-a)*GAMMA(2*a+3/2-n)*V1(a,n)/GAMMA(-1/2-3*a+2*n)/Pi/GAMMA(3*a
+1-n))*GAMMA(a)/GAMMA(1/2+n);}{%
\maplemultiline{
417\mbox{:~~~}  \left( {\vrule 
height1.67em width0em depth1.67em} \right. \!  \!  - 
{\displaystyle \frac {\Gamma ({\displaystyle \frac {1}{2}}  + n
 - a)\,\Gamma ( - {\displaystyle \frac {1}{6}}  + n - a)\,\Gamma 
({\displaystyle \frac {1}{6}}  + n - a)\,\Gamma (2\,a + 
{\displaystyle \frac {3}{2}}  - n)\,\Gamma (3\,a - 2\,n + 
{\displaystyle \frac {1}{2}} )\,\mathrm{V1}( - {\displaystyle 
\frac {1}{6}}  + n - a, \,n)}{\Gamma ( - {\displaystyle \frac {1
}{2}}  - 3\,a + 2\,n)\,\Gamma (2\,a + 1 - n)\,\Gamma (2\,a + 
{\displaystyle \frac {1}{3}}  - n)\,\Gamma (2\,a + 
{\displaystyle \frac {2}{3}}  - n)\,\sqrt{\pi }}}  \\
\mbox{} + {\displaystyle \frac {\Gamma ({\displaystyle \frac {1}{
2}}  + n - a)\,\Gamma ( - {\displaystyle \frac {1}{6}}  + n - a)
\,\Gamma ({\displaystyle \frac {1}{6}}  + n - a)\,\Gamma (2\,a + 
{\displaystyle \frac {3}{2}}  - n)\,\mathrm{V1}(a, \,n)}{\Gamma (
 - {\displaystyle \frac {1}{2}}  - 3\,a + 2\,n)\,\pi \,\Gamma (3
\,a + 1 - n)}}  \! \! \left. {\vrule 
height1.67em width0em depth1.67em} \right) \Gamma (a) \left/ 
{\vrule height0.80em width0em depth0.80em} \right. \!  \!  
\Gamma ({\displaystyle \frac {1}{2}}  + n) }
}
\end{maplelatex}

\begin{maplelatex}
\mapleinline{inert}{2d}{418, ":   ",
(-2/3*GAMMA(a+1/3)*GAMMA(1/6+n-a)*GAMMA(-1/6+n-a)*GAMMA(2/3)*GAMMA(a+2
/3-n)*V1(-1/6+n-a,n)*GAMMA(1/3+n-a)/GAMMA(-a-2/3+n)/GAMMA(2*a+1-n)/Pi/
GAMMA(2*a+2/3-n)/GAMMA(-3*a+2*n+1/2)+2/3*GAMMA(a+1/3)*GAMMA(1/6+n-a)*G
AMMA(-1/6+n-a)*GAMMA(2/3)*V1(a,n)/GAMMA(-a-2/3+n)/GAMMA(2/3+n-2*a)/Pi^
(1/2)/GAMMA(3*a+1-n))*GAMMA(a)/GAMMA(1/2+n);}{%
\maplemultiline{
418\mbox{:~~~}  \left( {\vrule 
height1.67em width0em depth1.67em} \right. \!  \!  - 
{\displaystyle \frac {2}{3}} \,{\displaystyle \frac {\Gamma (a + 
{\displaystyle \frac {1}{3}} )\,\Gamma ({\displaystyle \frac {1}{
6}}  + n - a)\,\Gamma ( - {\displaystyle \frac {1}{6}}  + n - a)
\,\Gamma ({\displaystyle \frac {2}{3}} )\,\Gamma (a + 
{\displaystyle \frac {2}{3}}  - n)\,\mathrm{V1}( - 
{\displaystyle \frac {1}{6}}  + n - a, \,n)\,\Gamma (
{\displaystyle \frac {1}{3}}  + n - a)}{\Gamma ( - a - 
{\displaystyle \frac {2}{3}}  + n)\,\Gamma (2\,a + 1 - n)\,\pi \,
\Gamma (2\,a + {\displaystyle \frac {2}{3}}  - n)\,\Gamma ( - 3\,
a + 2\,n + {\displaystyle \frac {1}{2}} )}}  \\
\mbox{} + {\displaystyle \frac {2}{3}} \,{\displaystyle \frac {
\Gamma (a + {\displaystyle \frac {1}{3}} )\,\Gamma (
{\displaystyle \frac {1}{6}}  + n - a)\,\Gamma ( - 
{\displaystyle \frac {1}{6}}  + n - a)\,\Gamma ({\displaystyle 
\frac {2}{3}} )\,\mathrm{V1}(a, \,n)}{\Gamma ( - a - 
{\displaystyle \frac {2}{3}}  + n)\,\Gamma ({\displaystyle 
\frac {2}{3}}  + n - 2\,a)\,\sqrt{\pi }\,\Gamma (3\,a + 1 - n)}} 
 \! \! \left. {\vrule height1.67em width0em depth1.67em} \right) 
\Gamma (a) \left/ {\vrule height0.80em width0em depth0.80em}
 \right. \!  \! \Gamma ({\displaystyle \frac {1}{2}}  + n) }
}
\end{maplelatex}

\begin{maplelatex}
\mapleinline{inert}{2d}{419, ":   ",
(-2/9*GAMMA(-1/6+n-a)*GAMMA(a+2/3)*GAMMA(1/6+n-a)*3^(1/2)/GAMMA(2/3)*G
AMMA(a+1/3-n)*V1(-1/6+n-a,n)*GAMMA(2/3+n-a)/GAMMA(-a-1/3+n)/GAMMA(2*a+
1-n)/GAMMA(2*a+1/3-n)/GAMMA(-3*a+2*n+1/2)+2/9*GAMMA(-1/6+n-a)*GAMMA(a+
2/3)*GAMMA(1/6+n-a)*Pi^(1/2)*3^(1/2)/GAMMA(2/3)*V1(a,n)/GAMMA(-a-1/3+n
)/GAMMA(1/3+n-2*a)/GAMMA(3*a+1-n))*GAMMA(a)/GAMMA(1/2+n);}{%
\maplemultiline{
419\mbox{:~~~}  \left( {\vrule 
height1.67em width0em depth1.67em} \right. \!  \!  - 
{\displaystyle \frac {2}{9}} \,{\displaystyle \frac {\Gamma ( - 
{\displaystyle \frac {1}{6}}  + n - a)\,\Gamma (a + 
{\displaystyle \frac {2}{3}} )\,\Gamma ({\displaystyle \frac {1}{
6}}  + n - a)\,\sqrt{3}\,\Gamma (a + {\displaystyle \frac {1}{3}
}  - n)\,\mathrm{V1}( - {\displaystyle \frac {1}{6}}  + n - a, \,
n)\,\Gamma ({\displaystyle \frac {2}{3}}  + n - a)}{\Gamma (
{\displaystyle \frac {2}{3}} )\,\Gamma ( - a - {\displaystyle 
\frac {1}{3}}  + n)\,\Gamma (2\,a + 1 - n)\,\Gamma (2\,a + 
{\displaystyle \frac {1}{3}}  - n)\,\Gamma ( - 3\,a + 2\,n + 
{\displaystyle \frac {1}{2}} )}}  \\
\mbox{} + {\displaystyle \frac {2}{9}} \,{\displaystyle \frac {
\Gamma ( - {\displaystyle \frac {1}{6}}  + n - a)\,\Gamma (a + 
{\displaystyle \frac {2}{3}} )\,\Gamma ({\displaystyle \frac {1}{
6}}  + n - a)\,\sqrt{\pi }\,\sqrt{3}\,\mathrm{V1}(a, \,n)}{\Gamma
 ({\displaystyle \frac {2}{3}} )\,\Gamma ( - a - {\displaystyle 
\frac {1}{3}}  + n)\,\Gamma ({\displaystyle \frac {1}{3}}  + n - 
2\,a)\,\Gamma (3\,a + 1 - n)}}  \! \! \left. {\vrule 
height1.67em width0em depth1.67em} \right) \Gamma (a) \left/ 
{\vrule height0.80em width0em depth0.80em} \right. \!  \! \Gamma 
({\displaystyle \frac {1}{2}}  + n) }
}
\end{maplelatex}

\begin{maplelatex}
\mapleinline{inert}{2d}{420, ":   ",
(GAMMA(-1/6+n-a)*GAMMA(1/6+n-a)*GAMMA(1/2+n-a)*V1(a,n)/GAMMA(-1/2-3*a+
2*n)/Pi*GAMMA(a-n+1/2)*GAMMA(2*a+3/2-n)/GAMMA(3*a-2*n+3/2)-GAMMA(-1/6+
n-a)*GAMMA(1/6+n-a)*GAMMA(1/2+n-a)*GAMMA(3*a-2*n+1/2)*V1(-1/6+n-a,n)*G
AMMA(3*a+1-n)/GAMMA(-1/2-3*a+2*n)/GAMMA(2*a+1/3-n)/GAMMA(2*a+2/3-n)/GA
MMA(2*a+1-n)/Pi^(1/2)*GAMMA(a-n+1/2)*GAMMA(2*a+3/2-n)/GAMMA(3*a-2*n+3/
2))/GAMMA(1/2+n);}{%
\maplemultiline{
420\mbox{:~~~}  \left( {\vrule 
height1.67em width0em depth1.67em} \right. \!  \! {\displaystyle 
\frac {\Gamma ( - {\displaystyle \frac {1}{6}}  + n - a)\,\Gamma 
({\displaystyle \frac {1}{6}}  + n - a)\,\Gamma ({\displaystyle 
\frac {1}{2}}  + n - a)\,\mathrm{V1}(a, \,n)\,\Gamma (a - n + 
{\displaystyle \frac {1}{2}} )\,\Gamma (2\,a + {\displaystyle 
\frac {3}{2}}  - n)}{\Gamma ( - {\displaystyle \frac {1}{2}}  - 3
\,a + 2\,n)\,\pi \,\Gamma (3\,a - 2\,n + {\displaystyle \frac {3
}{2}} )}}  -  \\
\Gamma ( - {\displaystyle \frac {1}{6}}  + n - a)\,\Gamma (
{\displaystyle \frac {1}{6}}  + n - a)\,\Gamma ({\displaystyle 
\frac {1}{2}}  + n - a)\,\Gamma (3\,a - 2\,n + {\displaystyle 
\frac {1}{2}} )\,\mathrm{V1}( - {\displaystyle \frac {1}{6}}  + n
 - a, \,n) \\
\Gamma (3\,a + 1 - n)\,\Gamma (a - n + {\displaystyle \frac {1}{2
}} )\,\Gamma (2\,a + {\displaystyle \frac {3}{2}}  - n) \left/ 
{\vrule height0.80em width0em depth0.80em} \right. \!  \! (\Gamma
 ( - {\displaystyle \frac {1}{2}}  - 3\,a + 2\,n)\,\Gamma (2\,a
 + {\displaystyle \frac {1}{3}}  - n) \\
\Gamma (2\,a + {\displaystyle \frac {2}{3}}  - n)\,\Gamma (2\,a
 + 1 - n)\,\sqrt{\pi }\,\Gamma (3\,a - 2\,n + {\displaystyle 
\frac {3}{2}} )) \! \! \left. {\vrule 
height1.67em width0em depth1.67em} \right)  \left/ {\vrule 
height0.80em width0em depth0.80em} \right. \!  \! \Gamma (
{\displaystyle \frac {1}{2}}  + n) }
}
\end{maplelatex}

\begin{maplelatex}
\mapleinline{inert}{2d}{421, ":   ",
(GAMMA(-1/6+n-a)*GAMMA(1/6+n-a)*GAMMA(1/2+n-a)*V1(a,n)/GAMMA(-1/2-3*a+
2*n)/Pi*GAMMA(5/6+a-n)*GAMMA(7/6+2*a-n)/GAMMA(3*a-2*n+3/2)-GAMMA(-1/6+
n-a)*GAMMA(1/6+n-a)*GAMMA(1/2+n-a)*GAMMA(3*a-2*n+1/2)*V1(-1/6+n-a,n)*G
AMMA(3*a+1-n)/GAMMA(-1/2-3*a+2*n)/GAMMA(2*a+1/3-n)/GAMMA(2*a+2/3-n)/GA
MMA(2*a+1-n)/Pi^(1/2)*GAMMA(5/6+a-n)*GAMMA(7/6+2*a-n)/GAMMA(3*a-2*n+3/
2))/GAMMA(1/2+n);}{%
\maplemultiline{
421\mbox{:~~~}  \left( {\vrule 
height1.67em width0em depth1.67em} \right. \!  \! {\displaystyle 
\frac {\Gamma ( - {\displaystyle \frac {1}{6}}  + n - a)\,\Gamma 
({\displaystyle \frac {1}{6}}  + n - a)\,\Gamma ({\displaystyle 
\frac {1}{2}}  + n - a)\,\mathrm{V1}(a, \,n)\,\Gamma (
{\displaystyle \frac {5}{6}}  + a - n)\,\Gamma ({\displaystyle 
\frac {7}{6}}  + 2\,a - n)}{\Gamma ( - {\displaystyle \frac {1}{2
}}  - 3\,a + 2\,n)\,\pi \,\Gamma (3\,a - 2\,n + {\displaystyle 
\frac {3}{2}} )}}  -  \\
\Gamma ( - {\displaystyle \frac {1}{6}}  + n - a)\,\Gamma (
{\displaystyle \frac {1}{6}}  + n - a)\,\Gamma ({\displaystyle 
\frac {1}{2}}  + n - a)\,\Gamma (3\,a - 2\,n + {\displaystyle 
\frac {1}{2}} )\,\mathrm{V1}( - {\displaystyle \frac {1}{6}}  + n
 - a, \,n) \\
\Gamma (3\,a + 1 - n)\,\Gamma ({\displaystyle \frac {5}{6}}  + a
 - n)\,\Gamma ({\displaystyle \frac {7}{6}}  + 2\,a - n) \left/ 
{\vrule height0.80em width0em depth0.80em} \right. \!  \! (\Gamma
 ( - {\displaystyle \frac {1}{2}}  - 3\,a + 2\,n)\,\Gamma (2\,a
 + {\displaystyle \frac {1}{3}}  - n) \\
\Gamma (2\,a + {\displaystyle \frac {2}{3}}  - n)\,\Gamma (2\,a
 + 1 - n)\,\sqrt{\pi }\,\Gamma (3\,a - 2\,n + {\displaystyle 
\frac {3}{2}} )) \! \! \left. {\vrule 
height1.67em width0em depth1.67em} \right)  \left/ {\vrule 
height0.80em width0em depth0.80em} \right. \!  \! \Gamma (
{\displaystyle \frac {1}{2}}  + n) }
}
\end{maplelatex}

\begin{maplelatex}
\mapleinline{inert}{2d}{422, ":   ",
(-GAMMA(-1/6+n-a)*GAMMA(1/6+n-a)*GAMMA(1/2+n-a)*GAMMA(3*a-2*n+1/2)*V1(
-1/6+n-a,n)/GAMMA(-1/2-3*a+2*n)/GAMMA(2*a+1/3-n)/GAMMA(2*a+2/3-n)/GAMM
A(2*a+1-n)/Pi^(1/2)/GAMMA(3*a-2*n+3/2)*GAMMA(7/6+a-n)*GAMMA(3*a+1-n)+G
AMMA(-1/6+n-a)*GAMMA(1/6+n-a)*GAMMA(1/2+n-a)*V1(a,n)/GAMMA(-1/2-3*a+2*
n)/Pi/GAMMA(3*a-2*n+3/2)*GAMMA(7/6+a-n))*GAMMA(2*a+5/6-n)/GAMMA(1/2+n)
;}{%
\maplemultiline{
422\mbox{:~~~}  \left( {\vrule 
height1.67em width0em depth1.67em} \right. \!  \!  - \Gamma ( - 
{\displaystyle \frac {1}{6}}  + n - a)\,\Gamma ({\displaystyle 
\frac {1}{6}}  + n - a)\,\Gamma ({\displaystyle \frac {1}{2}}  + 
n - a)\,\Gamma (3\,a - 2\,n + {\displaystyle \frac {1}{2}} )\,
\mathrm{V1}( - {\displaystyle \frac {1}{6}}  + n - a, \,n) \\
\Gamma ({\displaystyle \frac {7}{6}}  + a - n)\,\Gamma (3\,a + 1
 - n) \left/ {\vrule height0.80em width0em depth0.80em}
 \right. \!  \! (\Gamma ( - {\displaystyle \frac {1}{2}}  - 3\,a
 + 2\,n)\,\Gamma (2\,a + {\displaystyle \frac {1}{3}}  - n)\,
\Gamma (2\,a + {\displaystyle \frac {2}{3}}  - n) \\
\Gamma (2\,a + 1 - n)\,\sqrt{\pi }\,\Gamma (3\,a - 2\,n + 
{\displaystyle \frac {3}{2}} )) \\
\mbox{} + {\displaystyle \frac {\Gamma ( - {\displaystyle \frac {
1}{6}}  + n - a)\,\Gamma ({\displaystyle \frac {1}{6}}  + n - a)
\,\Gamma ({\displaystyle \frac {1}{2}}  + n - a)\,\mathrm{V1}(a, 
\,n)\,\Gamma ({\displaystyle \frac {7}{6}}  + a - n)}{\Gamma ( - 
{\displaystyle \frac {1}{2}}  - 3\,a + 2\,n)\,\pi \,\Gamma (3\,a
 - 2\,n + {\displaystyle \frac {3}{2}} )}}  \! \! \left. {\vrule 
height1.67em width0em depth1.67em} \right) \Gamma (2\,a + 
{\displaystyle \frac {5}{6}}  - n) 
 \left/ {\vrule height0.80em width0em depth0.80em} \right. \! 
 \! \Gamma ({\displaystyle \frac {1}{2}}  + n) }
}
\end{maplelatex}

\begin{maplelatex}
\mapleinline{inert}{2d}{423, ":   ",
(-GAMMA(1/2+n-a)*GAMMA(-1/6+n-a)*GAMMA(1/6+n-a)*GAMMA(3*a-2*n+1/2)*V1(
-1/6+n-a,n)/GAMMA(-1/2-3*a+2*n)/GAMMA(2*a+1-n)/GAMMA(2*a+1/3-n)/GAMMA(
2*a+2/3-n)*GAMMA(3*a+1-n)+GAMMA(1/2+n-a)*GAMMA(-1/6+n-a)*GAMMA(1/6+n-a
)*V1(a,n)/GAMMA(-1/2-3*a+2*n)/Pi^(1/2))/GAMMA(1/2+n);}{%
\maplemultiline{
423\mbox{:~~~}  \left( {\vrule 
height1.67em width0em depth1.67em} \right. \!  \!  - 
{\displaystyle \frac {\Gamma ({\displaystyle \frac {1}{2}}  + n
 - a)\,\Gamma ( - {\displaystyle \frac {1}{6}}  + n - a)\,\Gamma 
({\displaystyle \frac {1}{6}}  + n - a)\,\Gamma (3\,a - 2\,n + 
{\displaystyle \frac {1}{2}} )\,\mathrm{V1}( - {\displaystyle 
\frac {1}{6}}  + n - a, \,n)\,\Gamma (3\,a + 1 - n)}{\Gamma ( - 
{\displaystyle \frac {1}{2}}  - 3\,a + 2\,n)\,\Gamma (2\,a + 1 - 
n)\,\Gamma (2\,a + {\displaystyle \frac {1}{3}}  - n)\,\Gamma (2
\,a + {\displaystyle \frac {2}{3}}  - n)}}  \\
\mbox{} + {\displaystyle \frac {\Gamma ({\displaystyle \frac {1}{
2}}  + n - a)\,\Gamma ( - {\displaystyle \frac {1}{6}}  + n - a)
\,\Gamma ({\displaystyle \frac {1}{6}}  + n - a)\,\mathrm{V1}(a, 
\,n)}{\Gamma ( - {\displaystyle \frac {1}{2}}  - 3\,a + 2\,n)\,
\sqrt{\pi }}}  \! \! \left. {\vrule 
height1.67em width0em depth1.67em} \right)  \left/ {\vrule 
height0.80em width0em depth0.80em} \right. \!  \! \Gamma (
{\displaystyle \frac {1}{2}}  + n) }
}
\end{maplelatex}

\begin{maplelatex}
\mapleinline{inert}{2d}{424, ":   ",
(GAMMA(-1/6+n-a)*GAMMA(1/6+n-a)*GAMMA(1/2+n-a)*V1(a,n)/GAMMA(-1/2-3*a+
2*n)/Pi/GAMMA(3*a+1-n)*GAMMA(2*a+1/3-n)*GAMMA(2*a+3/2-n)*GAMMA(7/6+2*a
-n)/GAMMA(3*a-2*n+3/2)-GAMMA(-1/6+n-a)*GAMMA(1/6+n-a)*GAMMA(1/2+n-a)*G
AMMA(3*a-2*n+1/2)*V1(-1/6+n-a,n)/GAMMA(-1/2-3*a+2*n)/GAMMA(2*a+2/3-n)/
GAMMA(2*a+1-n)/Pi^(1/2)*GAMMA(2*a+3/2-n)*GAMMA(7/6+2*a-n)/GAMMA(3*a-2*
n+3/2))/GAMMA(1/2+n);}{%
\maplemultiline{
424\mbox{:~~~} (\Gamma ( - {\displaystyle \frac {1}{6}} 
 + n - a)\,\Gamma ({\displaystyle \frac {1}{6}}  + n - a)\,\Gamma
 ({\displaystyle \frac {1}{2}}  + n - a)\,\mathrm{V1}(a, \,n)\,
\Gamma (2\,a + {\displaystyle \frac {1}{3}}  - n)\,\Gamma (2\,a
 + {\displaystyle \frac {3}{2}}  - n) \\
\Gamma ({\displaystyle \frac {7}{6}}  + 2\,a - n) \left/ {\vrule 
height0.80em width0em depth0.80em} \right. \!  \! (\Gamma ( - 
{\displaystyle \frac {1}{2}}  - 3\,a + 2\,n)\,\pi \,\Gamma (3\,a
 + 1 - n)\,\Gamma (3\,a - 2\,n + {\displaystyle \frac {3}{2}} ))
\mbox{} -  \\
\Gamma ( - {\displaystyle \frac {1}{6}}  + n - a)\,\Gamma (
{\displaystyle \frac {1}{6}}  + n - a)\,\Gamma ({\displaystyle 
\frac {1}{2}}  + n - a)\,\Gamma (3\,a - 2\,n + {\displaystyle 
\frac {1}{2}} )\,\mathrm{V1}( - {\displaystyle \frac {1}{6}}  + n
 - a, \,n) \\
\Gamma (2\,a + {\displaystyle \frac {3}{2}}  - n)\,\Gamma (
{\displaystyle \frac {7}{6}}  + 2\,a - n) \left/ {\vrule 
height0.80em width0em depth0.80em} \right. \!  \! (\Gamma ( - 
{\displaystyle \frac {1}{2}}  - 3\,a + 2\,n)\,\Gamma (2\,a + 
{\displaystyle \frac {2}{3}}  - n) \\
\Gamma (2\,a + 1 - n)\,\sqrt{\pi }\,\Gamma (3\,a - 2\,n + 
{\displaystyle \frac {3}{2}} ))) \left/ {\vrule 
height0.80em width0em depth0.80em} \right. \!  \! \Gamma (
{\displaystyle \frac {1}{2}}  + n) }
}
\end{maplelatex}

\begin{maplelatex}
\mapleinline{inert}{2d}{425, ":   ",
(GAMMA(-1/6+n-a)*GAMMA(1/6+n-a)*GAMMA(1/2+n-a)*V1(a,n)/GAMMA(-1/2-3*a+
2*n)/Pi/GAMMA(3*a+1-n)*GAMMA(2*a+2/3-n)*GAMMA(2*a+3/2-n)/GAMMA(3*a-2*n
+3/2)-GAMMA(-1/6+n-a)*GAMMA(1/6+n-a)*GAMMA(1/2+n-a)*GAMMA(3*a-2*n+1/2)
*V1(-1/6+n-a,n)/GAMMA(-1/2-3*a+2*n)/GAMMA(2*a+1/3-n)/GAMMA(2*a+1-n)/Pi
^(1/2)*GAMMA(2*a+3/2-n)/GAMMA(3*a-2*n+3/2))*GAMMA(2*a+5/6-n)/GAMMA(1/2
+n);}{%
\maplemultiline{
425\mbox{:~~~}  \left( {\vrule 
height1.67em width0em depth1.67em} \right. \!  \! {\displaystyle 
\frac {\Gamma ( - {\displaystyle \frac {1}{6}}  + n - a)\,\Gamma 
({\displaystyle \frac {1}{6}}  + n - a)\,\Gamma ({\displaystyle 
\frac {1}{2}}  + n - a)\,\mathrm{V1}(a, \,n)\,\Gamma (2\,a + 
{\displaystyle \frac {2}{3}}  - n)\,\Gamma (2\,a + 
{\displaystyle \frac {3}{2}}  - n)}{\Gamma ( - {\displaystyle 
\frac {1}{2}}  - 3\,a + 2\,n)\,\pi \,\Gamma (3\,a + 1 - n)\,
\Gamma (3\,a - 2\,n + {\displaystyle \frac {3}{2}} )}}  \\
\mbox{} - {\displaystyle \frac {\Gamma ( - {\displaystyle \frac {
1}{6}}  + n - a)\,\Gamma ({\displaystyle \frac {1}{6}}  + n - a)
\,\Gamma ({\displaystyle \frac {1}{2}}  + n - a)\,\Gamma (3\,a - 
2\,n + {\displaystyle \frac {1}{2}} )\,\mathrm{V1}( - 
{\displaystyle \frac {1}{6}}  + n - a, \,n)\,\Gamma (2\,a + 
{\displaystyle \frac {3}{2}}  - n)}{\Gamma ( - {\displaystyle 
\frac {1}{2}}  - 3\,a + 2\,n)\,\Gamma (2\,a + {\displaystyle 
\frac {1}{3}}  - n)\,\Gamma (2\,a + 1 - n)\,\sqrt{\pi }\,\Gamma (
3\,a - 2\,n + {\displaystyle \frac {3}{2}} )}}  
 \! \! \left. {\vrule height1.67em width0em depth1.67em} \right) 
\Gamma (2\,a + {\displaystyle \frac {5}{6}}  - n) \left/ {\vrule 
height0.80em width0em depth0.80em} \right. \!  \! \Gamma (
{\displaystyle \frac {1}{2}}  + n) }
}
\end{maplelatex}

\begin{maplelatex}
\mapleinline{inert}{2d}{426, ":   ",
(GAMMA(-1/6+n-a)*GAMMA(1/6+n-a)*GAMMA(1/2+n-a)*V1(a,n)/GAMMA(-1/2-3*a+
2*n)/Pi/GAMMA(3*a+1-n)*GAMMA(2*a+1-n)*GAMMA(7/6+2*a-n)/GAMMA(3*a-2*n+3
/2)-GAMMA(-1/6+n-a)*GAMMA(1/6+n-a)*GAMMA(1/2+n-a)*GAMMA(3*a-2*n+1/2)*V
1(-1/6+n-a,n)/GAMMA(-1/2-3*a+2*n)/GAMMA(2*a+1/3-n)/GAMMA(2*a+2/3-n)/Pi
^(1/2)*GAMMA(7/6+2*a-n)/GAMMA(3*a-2*n+3/2))*GAMMA(2*a+5/6-n)/GAMMA(1/2
+n);}{%
\maplemultiline{
426\mbox{:~~~}  \left( {\vrule 
height1.67em width0em depth1.67em} \right. \!  \! {\displaystyle 
\frac {\Gamma ( - {\displaystyle \frac {1}{6}}  + n - a)\,\Gamma 
({\displaystyle \frac {1}{6}}  + n - a)\,\Gamma ({\displaystyle 
\frac {1}{2}}  + n - a)\,\mathrm{V1}(a, \,n)\,\Gamma (2\,a + 1 - 
n)\,\Gamma ({\displaystyle \frac {7}{6}}  + 2\,a - n)}{\Gamma (
 - {\displaystyle \frac {1}{2}}  - 3\,a + 2\,n)\,\pi \,\Gamma (3
\,a + 1 - n)\,\Gamma (3\,a - 2\,n + {\displaystyle \frac {3}{2}} 
)}}  \\
\mbox{} - {\displaystyle \frac {\Gamma ( - {\displaystyle \frac {
1}{6}}  + n - a)\,\Gamma ({\displaystyle \frac {1}{6}}  + n - a)
\,\Gamma ({\displaystyle \frac {1}{2}}  + n - a)\,\Gamma (3\,a - 
2\,n + {\displaystyle \frac {1}{2}} )\,\mathrm{V1}( - 
{\displaystyle \frac {1}{6}}  + n - a, \,n)\,\Gamma (
{\displaystyle \frac {7}{6}}  + 2\,a - n)}{\Gamma ( - 
{\displaystyle \frac {1}{2}}  - 3\,a + 2\,n)\,\Gamma (2\,a + 
{\displaystyle \frac {1}{3}}  - n)\,\Gamma (2\,a + 
{\displaystyle \frac {2}{3}}  - n)\,\sqrt{\pi }\,\Gamma (3\,a - 2
\,n + {\displaystyle \frac {3}{2}} )}}  \! \! \left. {\vrule height1.67em width0em depth1.67em} \right) \\ 
\Gamma (2\,a + {\displaystyle \frac {5}{6}}  - n) \left/ {\vrule 
height0.80em width0em depth0.80em} \right. \!  \! \Gamma (
{\displaystyle \frac {1}{2}}  + n) }
}
\end{maplelatex}

\begin{maplelatex}
\mapleinline{inert}{2d}{427, ":   ",
(-2/3*GAMMA(-1/6+n-a)*GAMMA(a+1/3)/Pi/GAMMA(2/3)*GAMMA(a+2/3-n)*V1(-1/
6+n-a,n)*GAMMA(1/3+n-a)/GAMMA(-a-2/3+n)/GAMMA(2*a+1-n)/GAMMA(-3*a+2*n+
1/2)*GAMMA(7/6-2*a+n)/GAMMA(5/3+a-n)+2/3*GAMMA(-1/6+n-a)*GAMMA(a+1/3)/
Pi^(1/2)/GAMMA(2/3)*V1(a,n)/GAMMA(-a-2/3+n)/GAMMA(2/3+n-2*a)*GAMMA(2*a
+2/3-n)*GAMMA(7/6-2*a+n)/GAMMA(5/3+a-n)/GAMMA(3*a+1-n))*GAMMA(a)/GAMMA
(1/2+n)*GAMMA(2/3)^2;}{%
\maplemultiline{
427\mbox{:~~~}  \left( {\vrule 
height1.67em width0em depth1.67em} \right. \!  \!  - 
{\displaystyle \frac {2}{3}}  
{\displaystyle \frac {\Gamma ( - {\displaystyle \frac {1}{6}}  + 
n - a)\,\Gamma (a + {\displaystyle \frac {1}{3}} )\,\Gamma (a + 
{\displaystyle \frac {2}{3}}  - n)\,\mathrm{V1}( - 
{\displaystyle \frac {1}{6}}  + n - a, \,n)\,\Gamma (
{\displaystyle \frac {1}{3}}  + n - a)\,\Gamma ({\displaystyle 
\frac {7}{6}}  - 2\,a + n)}{\pi \,\Gamma ({\displaystyle \frac {2
}{3}} )\,\Gamma ( - a - {\displaystyle \frac {2}{3}}  + n)\,
\Gamma (2\,a + 1 - n)\,\Gamma ( - 3\,a + 2\,n + {\displaystyle 
\frac {1}{2}} )\,\Gamma ({\displaystyle \frac {5}{3}}  + a - n)}
}  \\
\mbox{} + {\displaystyle \frac {2}{3}} \,{\displaystyle \frac {
\Gamma ( - {\displaystyle \frac {1}{6}}  + n - a)\,\Gamma (a + 
{\displaystyle \frac {1}{3}} )\,\mathrm{V1}(a, \,n)\,\Gamma (2\,a
 + {\displaystyle \frac {2}{3}}  - n)\,\Gamma ({\displaystyle 
\frac {7}{6}}  - 2\,a + n)}{\sqrt{\pi }\,\Gamma ({\displaystyle 
\frac {2}{3}} )\,\Gamma ( - a - {\displaystyle \frac {2}{3}}  + n
)\,\Gamma ({\displaystyle \frac {2}{3}}  + n - 2\,a)\,\Gamma (
{\displaystyle \frac {5}{3}}  + a - n)\,\Gamma (3\,a + 1 - n)}} 
 \! \! \left. {\vrule height1.67em width0em depth1.67em} \right) 
\Gamma (a)\,\Gamma ({\displaystyle \frac {2}{3}} )^{2} 
 \left/ {\vrule height0.80em width0em depth0.80em} \right. \! 
 \! \Gamma ({\displaystyle \frac {1}{2}}  + n) }
}
\end{maplelatex}

\begin{maplelatex}
\mapleinline{inert}{2d}{428, ":   ",
(-GAMMA(-1/6+n-a)*GAMMA(a+1/3)/Pi/GAMMA(2/3)*GAMMA(a+2/3-n)*V1(-1/6+n-
a,n)*GAMMA(1/3+n-a)/GAMMA(-a-2/3+n)/GAMMA(2*a+2/3-n)/GAMMA(2*a+1-n)/GA
MMA(-3*a+2*n+1/2)*GAMMA(7/6-2*a+n)*GAMMA(7/6+a)/GAMMA(5/3+a-n)*GAMMA(7
/6+a-n)+GAMMA(-1/6+n-a)*GAMMA(a+1/3)/Pi^(1/2)/GAMMA(2/3)*V1(a,n)/GAMMA
(-a-2/3+n)/GAMMA(2/3+n-2*a)*GAMMA(7/6-2*a+n)*GAMMA(7/6+a)/GAMMA(5/3+a-
n)*GAMMA(7/6+a-n)/GAMMA(3*a+1-n))*GAMMA(a)/GAMMA(1/2+n)*GAMMA(2/3);}{%
\maplemultiline{
428\mbox{:~~~}  \left( {\vrule 
height1.67em width0em depth1.67em} \right. \!  \!  - \Gamma ( - 
{\displaystyle \frac {1}{6}}  + n - a)\,\Gamma (a + 
{\displaystyle \frac {1}{3}} )\,\Gamma (a + {\displaystyle 
\frac {2}{3}}  - n)\,\mathrm{V1}( - {\displaystyle \frac {1}{6}} 
 + n - a, \,n)\,\Gamma ({\displaystyle \frac {1}{3}}  + n - a)
 \\
\Gamma ({\displaystyle \frac {7}{6}}  - 2\,a + n)\,\Gamma (
{\displaystyle \frac {7}{6}}  + a)\,\Gamma ({\displaystyle 
\frac {7}{6}}  + a - n) \left/ {\vrule 
height0.80em width0em depth0.80em} \right. \!  \! (\pi \,\Gamma (
{\displaystyle \frac {2}{3}} )\,\Gamma ( - a - {\displaystyle 
\frac {2}{3}}  + n)\,\Gamma (2\,a + {\displaystyle \frac {2}{3}} 
 - n) \\
\Gamma (2\,a + 1 - n)\,\Gamma ( - 3\,a + 2\,n + {\displaystyle 
\frac {1}{2}} )\,\Gamma ({\displaystyle \frac {5}{3}}  + a - n))
 \\
\mbox{} + {\displaystyle \frac {\Gamma ( - {\displaystyle \frac {
1}{6}}  + n - a)\,\Gamma (a + {\displaystyle \frac {1}{3}} )\,
\mathrm{V1}(a, \,n)\,\Gamma ({\displaystyle \frac {7}{6}}  - 2\,a
 + n)\,\Gamma ({\displaystyle \frac {7}{6}}  + a)\,\Gamma (
{\displaystyle \frac {7}{6}}  + a - n)}{\sqrt{\pi }\,\Gamma (
{\displaystyle \frac {2}{3}} )\,\Gamma ( - a - {\displaystyle 
\frac {2}{3}}  + n)\,\Gamma ({\displaystyle \frac {2}{3}}  + n - 
2\,a)\,\Gamma ({\displaystyle \frac {5}{3}}  + a - n)\,\Gamma (3
\,a + 1 - n)}}  \! \! \left. {\vrule 
height1.67em width0em depth1.67em} \right) \Gamma (a)\,\Gamma (
{\displaystyle \frac {2}{3}} )  \left/ {\vrule height0.80em width0em depth0.80em} \right. \! 
 \! \Gamma ({\displaystyle \frac {1}{2}}  + n) }
}
\end{maplelatex}

\begin{maplelatex}
\mapleinline{inert}{2d}{429, ":   ",
(-GAMMA(-1/6+n-a)*GAMMA(a+1/3)/Pi/GAMMA(2/3)*GAMMA(a+2/3-n)*V1(-1/6+n-
a,n)*GAMMA(1/3+n-a)/GAMMA(-a-2/3+n)/GAMMA(2*a+2/3-n)/GAMMA(2*a+1-n)/GA
MMA(-3*a+2*n+1/2)*GAMMA(7/6-2*a+n)*GAMMA(a+2/3)+GAMMA(-1/6+n-a)*GAMMA(
a+1/3)/Pi^(1/2)/GAMMA(2/3)*V1(a,n)/GAMMA(-a-2/3+n)/GAMMA(2/3+n-2*a)*GA
MMA(7/6-2*a+n)*GAMMA(a+2/3)/GAMMA(3*a+1-n))*GAMMA(a)/GAMMA(1/2+n)*GAMM
A(2/3);}{%
\maplemultiline{
429\mbox{:~~~}  \left( {\vrule 
height1.67em width0em depth1.67em} \right. \!  \!  - \Gamma ( - 
{\displaystyle \frac {1}{6}}  + n - a)\,\Gamma (a + 
{\displaystyle \frac {1}{3}} )\,\Gamma (a + {\displaystyle 
\frac {2}{3}}  - n)\,\mathrm{V1}( - {\displaystyle \frac {1}{6}} 
 + n - a, \,n)\,\Gamma ({\displaystyle \frac {1}{3}}  + n - a)
 \\
\Gamma ({\displaystyle \frac {7}{6}}  - 2\,a + n)\,\Gamma (a + 
{\displaystyle \frac {2}{3}} ) \left/ {\vrule 
height0.80em width0em depth0.80em} \right. \!  \! (\pi \,\Gamma (
{\displaystyle \frac {2}{3}} )\,\Gamma ( - a - {\displaystyle 
\frac {2}{3}}  + n)\,\Gamma (2\,a + {\displaystyle \frac {2}{3}} 
 - n)\,\Gamma (2\,a + 1 - n) \\
\Gamma ( - 3\,a + 2\,n + {\displaystyle \frac {1}{2}} ))\mbox{}
 + {\displaystyle \frac {\Gamma ( - {\displaystyle \frac {1}{6}} 
 + n - a)\,\Gamma (a + {\displaystyle \frac {1}{3}} )\,\mathrm{V1
}(a, \,n)\,\Gamma ({\displaystyle \frac {7}{6}}  - 2\,a + n)\,
\Gamma (a + {\displaystyle \frac {2}{3}} )}{\sqrt{\pi }\,\Gamma (
{\displaystyle \frac {2}{3}} )\,\Gamma ( - a - {\displaystyle 
\frac {2}{3}}  + n)\,\Gamma ({\displaystyle \frac {2}{3}}  + n - 
2\,a)\,\Gamma (3\,a + 1 - n)}}  \! \! \left. {\vrule 
height1.67em width0em depth1.67em} \right) \Gamma (a) 
\Gamma ({\displaystyle \frac {2}{3}} ) \left/ {\vrule 
height0.80em width0em depth0.80em} \right. \!  \! \Gamma (
{\displaystyle \frac {1}{2}}  + n) }
}
\end{maplelatex}

\begin{maplelatex}
\mapleinline{inert}{2d}{430, ":   ",
(2/9*GAMMA(-1/6+n-a)*GAMMA(a+1/3)*Pi^(1/2)*3^(1/2)/GAMMA(2/3)*V1(a,n)/
GAMMA(-a-2/3+n)/GAMMA(3*a+1-n)/GAMMA(2/3+n-2*a)*GAMMA(2*a+1-n)*GAMMA(7
/6-2*a+n)/GAMMA(5/3+a-n)-2/9*GAMMA(-1/6+n-a)*GAMMA(a+1/3)*3^(1/2)/GAMM
A(2/3)*GAMMA(a+2/3-n)*V1(-1/6+n-a,n)*GAMMA(1/3+n-a)/GAMMA(-a-2/3+n)/GA
MMA(2*a+2/3-n)/GAMMA(-3*a+2*n+1/2)*GAMMA(7/6-2*a+n)/GAMMA(5/3+a-n))*GA
MMA(a)/GAMMA(1/2+n);}{%
\maplemultiline{
430\mbox{:~~~}  \left( {\vrule 
height1.67em width0em depth1.67em} \right. \!  \! {\displaystyle 
\frac {2}{9}} \,{\displaystyle \frac {\Gamma ( - {\displaystyle 
\frac {1}{6}}  + n - a)\,\Gamma (a + {\displaystyle \frac {1}{3}
} )\,\sqrt{\pi }\,\sqrt{3}\,\mathrm{V1}(a, \,n)\,\Gamma (2\,a + 1
 - n)\,\Gamma ({\displaystyle \frac {7}{6}}  - 2\,a + n)}{\Gamma 
({\displaystyle \frac {2}{3}} )\,\Gamma ( - a - {\displaystyle 
\frac {2}{3}}  + n)\,\Gamma (3\,a + 1 - n)\,\Gamma (
{\displaystyle \frac {2}{3}}  + n - 2\,a)\,\Gamma (
{\displaystyle \frac {5}{3}}  + a - n)}} \\ - {\displaystyle 
\frac {2}{9}} \,{\displaystyle \frac {\Gamma ( - {\displaystyle 
\frac {1}{6}}  + n - a)\,\Gamma (a + {\displaystyle \frac {1}{3}
} )\,\sqrt{3}\,\Gamma (a + {\displaystyle \frac {2}{3}}  - n)\,
\mathrm{V1}( - {\displaystyle \frac {1}{6}}  + n - a, \,n)\,
\Gamma ({\displaystyle \frac {1}{3}}  + n - a)\,\Gamma (
{\displaystyle \frac {7}{6}}  - 2\,a + n)}{\Gamma (
{\displaystyle \frac {2}{3}} )\,\Gamma ( - a - {\displaystyle 
\frac {2}{3}}  + n)\,\Gamma (2\,a + {\displaystyle \frac {2}{3}} 
 - n)\,\Gamma ( - 3\,a + 2\,n + {\displaystyle \frac {1}{2}} )\,
\Gamma ({\displaystyle \frac {5}{3}}  + a - n)}}  
 \! \! \left. {\vrule height1.67em width0em depth1.67em} \right) 
\Gamma (a) \left/ {\vrule height0.80em width0em depth0.80em}
 \right. \!  \! \Gamma ({\displaystyle \frac {1}{2}}  + n) }
}
\end{maplelatex}

\begin{maplelatex}
\mapleinline{inert}{2d}{431, ":   ",
(-2/3*GAMMA(-1/6+n-a)*GAMMA(a+1/3)/Pi/GAMMA(2/3)*GAMMA(a+2/3-n)*V1(-1/
6+n-a,n)*GAMMA(1/3+n-a)/GAMMA(-a-2/3+n)/GAMMA(2*a+2/3-n)/GAMMA(2*a+1-n
)/GAMMA(-3*a+2*n+1/2)*GAMMA(2/3-a)*GAMMA(7/6+a)/GAMMA(5/3+a-n)+2/3*GAM
MA(-1/6+n-a)*GAMMA(a+1/3)/Pi^(1/2)/GAMMA(2/3)*V1(a,n)/GAMMA(-a-2/3+n)/
GAMMA(2/3+n-2*a)*GAMMA(2/3-a)*GAMMA(7/6+a)/GAMMA(5/3+a-n)/GAMMA(3*a+1-
n))*GAMMA(a)/GAMMA(1/2+n)*GAMMA(2/3)^2;}{%
\maplemultiline{
431\mbox{:~~~}  \left( {\vrule 
height1.67em width0em depth1.67em} \right. \!  \!  - 
{\displaystyle \frac {2}{3}} \,{\displaystyle \frac {\Gamma ( - 
{\displaystyle \frac {1}{6}}  + n - a)\,\Gamma (a + 
{\displaystyle \frac {1}{3}} )\,\Gamma (a + {\displaystyle 
\frac {2}{3}}  - n)\,\mathrm{V1}( - {\displaystyle \frac {1}{6}} 
 + n - a, \,n)\,\Gamma ({\displaystyle \frac {1}{3}}  + n - a)\,
\Gamma ({\displaystyle \frac {2}{3}}  - a)\,\Gamma (
{\displaystyle \frac {7}{6}}  + a)}{\pi \,\Gamma ({\displaystyle 
\frac {2}{3}} )\,\Gamma ( - a - {\displaystyle \frac {2}{3}}  + n
)\,\Gamma (2\,a + {\displaystyle \frac {2}{3}}  - n)\,\Gamma (2\,
a + 1 - n)\,\Gamma ( - 3\,a + 2\,n + {\displaystyle \frac {1}{2}
} )\,\Gamma ({\displaystyle \frac {5}{3}}  + a - n)}}  \\
\mbox{} + {\displaystyle \frac {2}{3}} \,{\displaystyle \frac {
\Gamma ( - {\displaystyle \frac {1}{6}}  + n - a)\,\Gamma (a + 
{\displaystyle \frac {1}{3}} )\,\mathrm{V1}(a, \,n)\,\Gamma (
{\displaystyle \frac {2}{3}}  - a)\,\Gamma ({\displaystyle 
\frac {7}{6}}  + a)}{\sqrt{\pi }\,\Gamma ({\displaystyle \frac {2
}{3}} )\,\Gamma ( - a - {\displaystyle \frac {2}{3}}  + n)\,
\Gamma ({\displaystyle \frac {2}{3}}  + n - 2\,a)\,\Gamma (
{\displaystyle \frac {5}{3}}  + a - n)\,\Gamma (3\,a + 1 - n)}} 
 \! \! \left. {\vrule height1.67em width0em depth1.67em} \right) 
\Gamma (a)\,\Gamma ({\displaystyle \frac {2}{3}} )^{2}
 \left/ {\vrule height0.80em width0em depth0.80em} \right. \! 
 \! \Gamma ({\displaystyle \frac {1}{2}}  + n) }
}
\end{maplelatex}

\begin{maplelatex}
\mapleinline{inert}{2d}{432, ":   ",
(-4/27*GAMMA(-1/6+n-a)*GAMMA(a+1/3)*Pi^(1/2)*3^(1/2)/GAMMA(2/3)*GAMMA(
a+2/3-n)*V1(-1/6+n-a,n)*GAMMA(1/3+n-a)/GAMMA(-a-2/3+n)/GAMMA(2*a+2/3-n
)/GAMMA(2*a+1-n)/GAMMA(-3*a+2*n+1/2)/GAMMA(5/3+a-n)+4/27*GAMMA(-1/6+n-
a)*GAMMA(a+1/3)*Pi*3^(1/2)/GAMMA(2/3)*V1(a,n)/GAMMA(-a-2/3+n)/GAMMA(2/
3+n-2*a)/GAMMA(5/3+a-n)/GAMMA(3*a+1-n))*GAMMA(a)/GAMMA(1/2+n)*GAMMA(2/
3);}{%
\maplemultiline{
432\mbox{:~~~}  \left( {\vrule 
height1.67em width0em depth1.67em} \right. \!  \!  -  
{\displaystyle \frac {4}{27}} \,{\displaystyle \frac {\Gamma ( - 
{\displaystyle \frac {1}{6}}  + n - a)\,\Gamma (a + 
{\displaystyle \frac {1}{3}} )\,\sqrt{\pi }\,\sqrt{3}\,\Gamma (a
 + {\displaystyle \frac {2}{3}}  - n)\,\mathrm{V1}( - 
{\displaystyle \frac {1}{6}}  + n - a, \,n)\,\Gamma (
{\displaystyle \frac {1}{3}}  + n - a)}{\Gamma ({\displaystyle 
\frac {2}{3}} )\,\Gamma ( - a - {\displaystyle \frac {2}{3}}  + n
)\,\Gamma (2\,a + {\displaystyle \frac {2}{3}}  - n)\,\Gamma (2\,
a + 1 - n)\,\Gamma ( - 3\,a + 2\,n + {\displaystyle \frac {1}{2}
} )\,\Gamma ({\displaystyle \frac {5}{3}}  + a - n)}}  \\
\mbox{} + {\displaystyle \frac {4}{27}} \,{\displaystyle \frac {
\Gamma ( - {\displaystyle \frac {1}{6}}  + n - a)\,\Gamma (a + 
{\displaystyle \frac {1}{3}} )\,\pi \,\sqrt{3}\,\mathrm{V1}(a, \,
n)}{\Gamma ({\displaystyle \frac {2}{3}} )\,\Gamma ( - a - 
{\displaystyle \frac {2}{3}}  + n)\,\Gamma ({\displaystyle 
\frac {2}{3}}  + n - 2\,a)\,\Gamma ({\displaystyle \frac {5}{3}} 
 + a - n)\,\Gamma (3\,a + 1 - n)}}  \! \! \left. {\vrule 
height1.67em width0em depth1.67em} \right) \Gamma (a)\,\Gamma (
{\displaystyle \frac {2}{3}} ) 
 \left/ {\vrule height0.80em width0em depth0.80em} \right. \! 
 \! \Gamma ({\displaystyle \frac {1}{2}}  + n) }
}
\end{maplelatex}

\begin{maplelatex}
\mapleinline{inert}{2d}{433, ":   ",
(-GAMMA(-1/6+n-a)*GAMMA(a+1/3)/Pi/GAMMA(2/3)*GAMMA(a+2/3-n)*V1(-1/6+n-
a,n)*GAMMA(1/3+n-a)/GAMMA(-a-2/3+n)/GAMMA(2*a+2/3-n)/GAMMA(2*a+1-n)/GA
MMA(-3*a+2*n+1/2)*GAMMA(2/3+n-2*a)*GAMMA(7/6+a)+GAMMA(-1/6+n-a)*GAMMA(
a+1/3)/Pi^(1/2)/GAMMA(2/3)*V1(a,n)/GAMMA(-a-2/3+n)*GAMMA(7/6+a)/GAMMA(
3*a+1-n))*GAMMA(a)/GAMMA(1/2+n)*GAMMA(2/3);}{%
\maplemultiline{
433\mbox{:~~~}  \left( {\vrule 
height1.67em width0em depth1.67em} \right. \!  \!  - \Gamma ( - 
{\displaystyle \frac {1}{6}}  + n - a)\,\Gamma (a + 
{\displaystyle \frac {1}{3}} )\,\Gamma (a + {\displaystyle 
\frac {2}{3}}  - n)\,\mathrm{V1}( - {\displaystyle \frac {1}{6}} 
 + n - a, \,n)\,\Gamma ({\displaystyle \frac {1}{3}}  + n - a)
 \\
\Gamma ({\displaystyle \frac {2}{3}}  + n - 2\,a)\,\Gamma (
{\displaystyle \frac {7}{6}}  + a) \left/ {\vrule 
height0.80em width0em depth0.80em} \right. \!  \! (\pi \,\Gamma (
{\displaystyle \frac {2}{3}} )\,\Gamma ( - a - {\displaystyle 
\frac {2}{3}}  + n)\,\Gamma (2\,a + {\displaystyle \frac {2}{3}} 
 - n)\,\Gamma (2\,a + 1 - n) \\
\Gamma ( - 3\,a + 2\,n + {\displaystyle \frac {1}{2}} ))\mbox{}
 + {\displaystyle \frac {\Gamma ( - {\displaystyle \frac {1}{6}} 
 + n - a)\,\Gamma (a + {\displaystyle \frac {1}{3}} )\,\mathrm{V1
}(a, \,n)\,\Gamma ({\displaystyle \frac {7}{6}}  + a)}{\sqrt{\pi 
}\,\Gamma ({\displaystyle \frac {2}{3}} )\,\Gamma ( - a - 
{\displaystyle \frac {2}{3}}  + n)\,\Gamma (3\,a + 1 - n)}}  \! 
\! \left. {\vrule height1.67em width0em depth1.67em} \right) 
\Gamma (a)\,\Gamma ({\displaystyle \frac {2}{3}} ) \left/ 
{\vrule height0.80em width0em depth0.80em} \right. \!  \! 
\Gamma ({\displaystyle \frac {1}{2}}  + n) }
}
\end{maplelatex}

\begin{maplelatex}
\mapleinline{inert}{2d}{434, ":   ",
(2/9*GAMMA(-1/6+n-a)*GAMMA(a+1/3)*Pi^(1/2)*3^(1/2)/GAMMA(2/3)*V1(a,n)/
GAMMA(-a-2/3+n)/GAMMA(3*a+1-n)/GAMMA(2/3+n-2*a)*GAMMA(7/6+a)/GAMMA(5/3
+a-n)-2/9*GAMMA(-1/6+n-a)*GAMMA(a+1/3)*3^(1/2)/GAMMA(2/3)*GAMMA(a+2/3-
n)*V1(-1/6+n-a,n)*GAMMA(1/3+n-a)/GAMMA(-a-2/3+n)/GAMMA(2*a+2/3-n)/GAMM
A(2*a+1-n)/GAMMA(-3*a+2*n+1/2)*GAMMA(7/6+a)/GAMMA(5/3+a-n))*GAMMA(-a+1
)*GAMMA(a)/GAMMA(1/2+n);}{%
\maplemultiline{
434\mbox{:~~~}  \left( {\vrule 
height1.67em width0em depth1.67em} \right. \!  \! {\displaystyle 
\frac {2}{9}} \,{\displaystyle \frac {\Gamma ( - {\displaystyle 
\frac {1}{6}}  + n - a)\,\Gamma (a + {\displaystyle \frac {1}{3}
} )\,\sqrt{\pi }\,\sqrt{3}\,\mathrm{V1}(a, \,n)\,\Gamma (
{\displaystyle \frac {7}{6}}  + a)}{\Gamma ({\displaystyle 
\frac {2}{3}} )\,\Gamma ( - a - {\displaystyle \frac {2}{3}}  + n
)\,\Gamma (3\,a + 1 - n)\,\Gamma ({\displaystyle \frac {2}{3}} 
 + n - 2\,a)\,\Gamma ({\displaystyle \frac {5}{3}}  + a - n)}} 
 \\ - 
{\displaystyle \frac {2}{9}} \,{\displaystyle \frac {\Gamma ( - 
{\displaystyle \frac {1}{6}}  + n - a)\,\Gamma (a + 
{\displaystyle \frac {1}{3}} )\,\sqrt{3}\,\Gamma (a + 
{\displaystyle \frac {2}{3}}  - n)\,\mathrm{V1}( - 
{\displaystyle \frac {1}{6}}  + n - a, \,n)\,\Gamma (
{\displaystyle \frac {1}{3}}  + n - a)\,\Gamma ({\displaystyle 
\frac {7}{6}}  + a)}{\Gamma ({\displaystyle \frac {2}{3}} )\,
\Gamma ( - a - {\displaystyle \frac {2}{3}}  + n)\,\Gamma (2\,a
 + {\displaystyle \frac {2}{3}}  - n)\,\Gamma (2\,a + 1 - n)\,
\Gamma ( - 3\,a + 2\,n + {\displaystyle \frac {1}{2}} )\,\Gamma (
{\displaystyle \frac {5}{3}}  + a - n)}}  \! \! \left. {\vrule 
height1.67em width0em depth1.67em} \right)  \\
\Gamma ( - a + 1)\,\Gamma (a) \left/ {\vrule 
height0.80em width0em depth0.80em} \right. \!  \! \Gamma (
{\displaystyle \frac {1}{2}}  + n) }
}
\end{maplelatex}

\begin{maplelatex}
\mapleinline{inert}{2d}{435, ":   ",
(-GAMMA(a+2/3)*GAMMA(1/6+n-a)/Pi/GAMMA(2/3)*GAMMA(a+1/3-n)*V1(-1/6+n-a
,n)*GAMMA(2/3+n-a)/GAMMA(-a-1/3+n)/GAMMA(2*a+1-n)/GAMMA(2*a+1/3-n)/GAM
MA(-3*a+2*n+1/2)*GAMMA(5/6+a-n)*GAMMA(5/6+n-2*a)/GAMMA(4/3+a-n)*GAMMA(
a+5/6)+GAMMA(a+2/3)*GAMMA(1/6+n-a)/Pi^(1/2)/GAMMA(2/3)*V1(a,n)/GAMMA(-
a-1/3+n)/GAMMA(1/3+n-2*a)*GAMMA(5/6+a-n)*GAMMA(5/6+n-2*a)/GAMMA(4/3+a-
n)*GAMMA(a+5/6)/GAMMA(3*a+1-n))*GAMMA(a)/GAMMA(1/2+n)*GAMMA(2/3);}{%
\maplemultiline{
435\mbox{:~~~}  \left( {\vrule 
height1.67em width0em depth1.67em} \right. \!  \!  - \Gamma (a + 
{\displaystyle \frac {2}{3}} )\,\Gamma ({\displaystyle \frac {1}{
6}}  + n - a)\,\Gamma (a + {\displaystyle \frac {1}{3}}  - n)\,
\mathrm{V1}( - {\displaystyle \frac {1}{6}}  + n - a, \,n)\,
\Gamma ({\displaystyle \frac {2}{3}}  + n - a)\,\Gamma (
{\displaystyle \frac {5}{6}}  + a - n) \\
\Gamma ({\displaystyle \frac {5}{6}}  + n - 2\,a)\,\Gamma (a + 
{\displaystyle \frac {5}{6}} ) \left/ {\vrule 
height0.80em width0em depth0.80em} \right. \!  \! (\pi \,\Gamma (
{\displaystyle \frac {2}{3}} )\,\Gamma ( - a - {\displaystyle 
\frac {1}{3}}  + n)\,\Gamma (2\,a + 1 - n)\,\Gamma (2\,a + 
{\displaystyle \frac {1}{3}}  - n) \\
\Gamma ( - 3\,a + 2\,n + {\displaystyle \frac {1}{2}} )\,\Gamma (
{\displaystyle \frac {4}{3}}  + a - n)) \\
\mbox{} + {\displaystyle \frac {\Gamma (a + {\displaystyle 
\frac {2}{3}} )\,\Gamma ({\displaystyle \frac {1}{6}}  + n - a)\,
\mathrm{V1}(a, \,n)\,\Gamma ({\displaystyle \frac {5}{6}}  + a - 
n)\,\Gamma ({\displaystyle \frac {5}{6}}  + n - 2\,a)\,\Gamma (a
 + {\displaystyle \frac {5}{6}} )}{\sqrt{\pi }\,\Gamma (
{\displaystyle \frac {2}{3}} )\,\Gamma ( - a - {\displaystyle 
\frac {1}{3}}  + n)\,\Gamma ({\displaystyle \frac {1}{3}}  + n - 
2\,a)\,\Gamma ({\displaystyle \frac {4}{3}}  + a - n)\,\Gamma (3
\,a + 1 - n)}}  \! \! \left. {\vrule 
height1.67em width0em depth1.67em} \right) \Gamma (a)\,\Gamma (
{\displaystyle \frac {2}{3}} )
 \left/ {\vrule height0.80em width0em depth0.80em} \right. \! 
 \! \Gamma ({\displaystyle \frac {1}{2}}  + n) }
}
\end{maplelatex}

\begin{maplelatex}
\mapleinline{inert}{2d}{436, ":   ",
(-GAMMA(a+2/3)*GAMMA(1/6+n-a)/Pi/GAMMA(2/3)*GAMMA(a+1/3-n)*V1(-1/6+n-a
,n)*GAMMA(2/3+n-a)/GAMMA(-a-1/3+n)/GAMMA(2*a+1/3-n)/GAMMA(-3*a+2*n+1/2
)*GAMMA(5/6+n-2*a)/GAMMA(4/3+a-n)+GAMMA(a+2/3)*GAMMA(1/6+n-a)/Pi^(1/2)
/GAMMA(2/3)*V1(a,n)/GAMMA(-a-1/3+n)/GAMMA(1/3+n-2*a)*GAMMA(2*a+1-n)*GA
MMA(5/6+n-2*a)/GAMMA(4/3+a-n)/GAMMA(3*a+1-n))*GAMMA(a)/GAMMA(1/2+n)*GA
MMA(2/3)^2;}{%
\maplemultiline{
436\mbox{:~~~}  \left( {\vrule 
height1.67em width0em depth1.67em} \right. \!  \!  
 - {\displaystyle \frac {\Gamma (a + {\displaystyle \frac {2}{3}
} )\,\Gamma ({\displaystyle \frac {1}{6}}  + n - a)\,\Gamma (a + 
{\displaystyle \frac {1}{3}}  - n)\,\mathrm{V1}( - 
{\displaystyle \frac {1}{6}}  + n - a, \,n)\,\Gamma (
{\displaystyle \frac {2}{3}}  + n - a)\,\Gamma ({\displaystyle 
\frac {5}{6}}  + n - 2\,a)}{\pi \,\Gamma ({\displaystyle \frac {2
}{3}} )\,\Gamma ( - a - {\displaystyle \frac {1}{3}}  + n)\,
\Gamma (2\,a + {\displaystyle \frac {1}{3}}  - n)\,\Gamma ( - 3\,
a + 2\,n + {\displaystyle \frac {1}{2}} )\,\Gamma (
{\displaystyle \frac {4}{3}}  + a - n)}}  \\
\mbox{} + {\displaystyle \frac {\Gamma (a + {\displaystyle 
\frac {2}{3}} )\,\Gamma ({\displaystyle \frac {1}{6}}  + n - a)\,
\mathrm{V1}(a, \,n)\,\Gamma (2\,a + 1 - n)\,\Gamma (
{\displaystyle \frac {5}{6}}  + n - 2\,a)}{\sqrt{\pi }\,\Gamma (
{\displaystyle \frac {2}{3}} )\,\Gamma ( - a - {\displaystyle 
\frac {1}{3}}  + n)\,\Gamma ({\displaystyle \frac {1}{3}}  + n - 
2\,a)\,\Gamma ({\displaystyle \frac {4}{3}}  + a - n)\,\Gamma (3
\,a + 1 - n)}}  \! \! \left. {\vrule 
height1.67em width0em depth1.67em} \right) \Gamma (a)\,\Gamma (
{\displaystyle \frac {2}{3}} )^{2}
 \left/ {\vrule height0.80em width0em depth0.80em} \right. \! 
 \! \Gamma ({\displaystyle \frac {1}{2}}  + n) }
}
\end{maplelatex}

\begin{maplelatex}
\mapleinline{inert}{2d}{437, ":   ",
(-GAMMA(a+2/3)*GAMMA(1/6+n-a)/Pi/GAMMA(2/3)*GAMMA(a+1/3-n)*V1(-1/6+n-a
,n)*GAMMA(2/3+n-a)/GAMMA(-a-1/3+n)/GAMMA(2*a+1-n)/GAMMA(2*a+1/3-n)/GAM
MA(-3*a+2*n+1/2)*GAMMA(5/6+n-2*a)*GAMMA(a+1/3)+GAMMA(a+2/3)*GAMMA(1/6+
n-a)/Pi^(1/2)/GAMMA(2/3)*V1(a,n)/GAMMA(-a-1/3+n)/GAMMA(1/3+n-2*a)*GAMM
A(5/6+n-2*a)*GAMMA(a+1/3)/GAMMA(3*a+1-n))*GAMMA(a)/GAMMA(1/2+n)*GAMMA(
2/3);}{%
\maplemultiline{
437\mbox{:~~~}  \left( {\vrule 
height1.67em width0em depth1.67em} \right. \!  \!  - 
{\displaystyle \frac {\Gamma (a + {\displaystyle \frac {2}{3}} )
\,\Gamma ({\displaystyle \frac {1}{6}}  + n - a)\,\Gamma (a + 
{\displaystyle \frac {1}{3}}  - n)\,\mathrm{V1}( - 
{\displaystyle \frac {1}{6}}  + n - a, \,n)\,\Gamma (
{\displaystyle \frac {2}{3}}  + n - a)\,\Gamma ({\displaystyle 
\frac {5}{6}}  + n - 2\,a)\,\Gamma (a + {\displaystyle \frac {1}{
3}} )}{\pi \,\Gamma ({\displaystyle \frac {2}{3}} )\,\Gamma ( - a
 - {\displaystyle \frac {1}{3}}  + n)\,\Gamma (2\,a + 1 - n)\,
\Gamma (2\,a + {\displaystyle \frac {1}{3}}  - n)\,\Gamma ( - 3\,
a + 2\,n + {\displaystyle \frac {1}{2}} )}}  \\
\mbox{} + {\displaystyle \frac {\Gamma (a + {\displaystyle 
\frac {2}{3}} )\,\Gamma ({\displaystyle \frac {1}{6}}  + n - a)\,
\mathrm{V1}(a, \,n)\,\Gamma ({\displaystyle \frac {5}{6}}  + n - 
2\,a)\,\Gamma (a + {\displaystyle \frac {1}{3}} )}{\sqrt{\pi }\,
\Gamma ({\displaystyle \frac {2}{3}} )\,\Gamma ( - a - 
{\displaystyle \frac {1}{3}}  + n)\,\Gamma ({\displaystyle 
\frac {1}{3}}  + n - 2\,a)\,\Gamma (3\,a + 1 - n)}}  \! 
\! \left. {\vrule height1.67em width0em depth1.67em} \right) 
\Gamma (a)\,\Gamma ({\displaystyle \frac {2}{3}} ) \left/ 
{\vrule height0.80em width0em depth0.80em} \right. \!  \! \Gamma 
({\displaystyle \frac {1}{2}}  + n) }
}
\end{maplelatex}

\begin{maplelatex}
\mapleinline{inert}{2d}{438, ":   ",
(2/9*GAMMA(a+2/3)*GAMMA(1/6+n-a)*Pi^(1/2)*3^(1/2)/GAMMA(2/3)*V1(a,n)/G
AMMA(-a-1/3+n)/GAMMA(3*a+1-n)/GAMMA(1/3+n-2*a)*GAMMA(2*a+1/3-n)*GAMMA(
5/6+n-2*a)/GAMMA(4/3+a-n)-2/9*GAMMA(a+2/3)*GAMMA(1/6+n-a)*3^(1/2)/GAMM
A(2/3)*GAMMA(a+1/3-n)*V1(-1/6+n-a,n)*GAMMA(2/3+n-a)/GAMMA(-a-1/3+n)/GA
MMA(2*a+1-n)/GAMMA(-3*a+2*n+1/2)*GAMMA(5/6+n-2*a)/GAMMA(4/3+a-n))*GAMM
A(a)/GAMMA(1/2+n);}{%
\maplemultiline{
438\mbox{:~~~}  \left( {\vrule 
height1.67em width0em depth1.67em} \right. \!  \! {\displaystyle 
\frac {2}{9}} \,{\displaystyle \frac {\Gamma (a + {\displaystyle 
\frac {2}{3}} )\,\Gamma ({\displaystyle \frac {1}{6}}  + n - a)\,
\sqrt{\pi }\,\sqrt{3}\,\mathrm{V1}(a, \,n)\,\Gamma (2\,a + 
{\displaystyle \frac {1}{3}}  - n)\,\Gamma ({\displaystyle 
\frac {5}{6}}  + n - 2\,a)}{\Gamma ({\displaystyle \frac {2}{3}} 
)\,\Gamma ( - a - {\displaystyle \frac {1}{3}}  + n)\,\Gamma (3\,
a + 1 - n)\,\Gamma ({\displaystyle \frac {1}{3}}  + n - 2\,a)\,
\Gamma ({\displaystyle \frac {4}{3}}  + a - n)}} \\ - 
{\displaystyle \frac {2}{9}} \,{\displaystyle \frac {\Gamma (a + 
{\displaystyle \frac {2}{3}} )\,\Gamma ({\displaystyle \frac {1}{
6}}  + n - a)\,\sqrt{3}\,\Gamma (a + {\displaystyle \frac {1}{3}
}  - n)\,\mathrm{V1}( - {\displaystyle \frac {1}{6}}  + n - a, \,
n)\,\Gamma ({\displaystyle \frac {2}{3}}  + n - a)\,\Gamma (
{\displaystyle \frac {5}{6}}  + n - 2\,a)}{\Gamma (
{\displaystyle \frac {2}{3}} )\,\Gamma ( - a - {\displaystyle 
\frac {1}{3}}  + n)\,\Gamma (2\,a + 1 - n)\,\Gamma ( - 3\,a + 2\,
n + {\displaystyle \frac {1}{2}} )\,\Gamma ({\displaystyle 
\frac {4}{3}}  + a - n)}}  
 \! \! \left. {\vrule height1.67em width0em depth1.67em} \right) 
\Gamma (a) \left/ {\vrule height0.80em width0em depth0.80em}
 \right. \!  \! \Gamma ({\displaystyle \frac {1}{2}}  + n) }
}
\end{maplelatex}

\begin{maplelatex}
\mapleinline{inert}{2d}{439, ":   ",
(-GAMMA(a+2/3)*GAMMA(1/6+n-a)/Pi/GAMMA(2/3)*GAMMA(a+1/3-n)*V1(-1/6+n-a
,n)*GAMMA(2/3+n-a)/GAMMA(-a-1/3+n)/GAMMA(2*a+1-n)/GAMMA(2*a+1/3-n)/GAM
MA(-3*a+2*n+1/2)/GAMMA(4/3+a-n)*GAMMA(a+5/6)+GAMMA(a+2/3)*GAMMA(1/6+n-
a)/Pi^(1/2)/GAMMA(2/3)*V1(a,n)/GAMMA(-a-1/3+n)/GAMMA(1/3+n-2*a)/GAMMA(
4/3+a-n)*GAMMA(a+5/6)/GAMMA(3*a+1-n))*GAMMA(-a+1)*GAMMA(a)/GAMMA(1/2+n
)*GAMMA(2/3)^2;}{%
\maplemultiline{
439\mbox{:~~~}  \left( {\vrule 
height1.67em width0em depth1.67em} \right. \!  \!  -  
{\displaystyle \frac {\Gamma (a + {\displaystyle \frac {2}{3}} )
\,\Gamma ({\displaystyle \frac {1}{6}}  + n - a)\,\Gamma (a + 
{\displaystyle \frac {1}{3}}  - n)\,\mathrm{V1}( - 
{\displaystyle \frac {1}{6}}  + n - a, \,n)\,\Gamma (
{\displaystyle \frac {2}{3}}  + n - a)\,\Gamma (a + 
{\displaystyle \frac {5}{6}} )}{\pi \,\Gamma ({\displaystyle 
\frac {2}{3}} )\,\Gamma ( - a - {\displaystyle \frac {1}{3}}  + n
)\,\Gamma (2\,a + 1 - n)\,\Gamma (2\,a + {\displaystyle \frac {1
}{3}}  - n)\,\Gamma ( - 3\,a + 2\,n + {\displaystyle \frac {1}{2}
} )\,\Gamma ({\displaystyle \frac {4}{3}}  + a - n)}}  \\
\mbox{} + {\displaystyle \frac {\Gamma (a + {\displaystyle 
\frac {2}{3}} )\,\Gamma ({\displaystyle \frac {1}{6}}  + n - a)\,
\mathrm{V1}(a, \,n)\,\Gamma (a + {\displaystyle \frac {5}{6}} )}{
\sqrt{\pi }\,\Gamma ({\displaystyle \frac {2}{3}} )\,\Gamma ( - a
 - {\displaystyle \frac {1}{3}}  + n)\,\Gamma ({\displaystyle 
\frac {1}{3}}  + n - 2\,a)\,\Gamma ({\displaystyle \frac {4}{3}} 
 + a - n)\,\Gamma (3\,a + 1 - n)}}  \! \! \left. {\vrule 
height1.67em width0em depth1.67em} \right) \Gamma ( - a + 1) 
\Gamma (a)\,\Gamma ({\displaystyle \frac {2}{3}} )^{2} \left/ 
{\vrule height0.80em width0em depth0.80em} \right. \!  \! \Gamma 
({\displaystyle \frac {1}{2}}  + n) }
}
\end{maplelatex}

\begin{maplelatex}
\mapleinline{inert}{2d}{440, ":   ",
(-GAMMA(a+2/3)*GAMMA(1/6+n-a)/Pi/GAMMA(2/3)*GAMMA(a+1/3-n)*V1(-1/6+n-a
,n)*GAMMA(2/3+n-a)/GAMMA(-a-1/3+n)/GAMMA(2*a+1-n)/GAMMA(2*a+1/3-n)/GAM
MA(-3*a+2*n+1/2)*GAMMA(1/3+n-2*a)*GAMMA(a+5/6)+GAMMA(a+2/3)*GAMMA(1/6+
n-a)/Pi^(1/2)/GAMMA(2/3)*V1(a,n)/GAMMA(-a-1/3+n)*GAMMA(a+5/6)/GAMMA(3*
a+1-n))*GAMMA(a)/GAMMA(1/2+n)*GAMMA(2/3);}{%
\maplemultiline{
440\mbox{:~~~}  \left( {\vrule 
height1.67em width0em depth1.67em} \right. \!  \!  - 
{\displaystyle \frac {\Gamma (a + {\displaystyle \frac {2}{3}} )
\,\Gamma ({\displaystyle \frac {1}{6}}  + n - a)\,\Gamma (a + 
{\displaystyle \frac {1}{3}}  - n)\,\mathrm{V1}( - 
{\displaystyle \frac {1}{6}}  + n - a, \,n)\,\Gamma (
{\displaystyle \frac {2}{3}}  + n - a)\,\Gamma ({\displaystyle 
\frac {1}{3}}  + n - 2\,a)\,\Gamma (a + {\displaystyle \frac {5}{
6}} )}{\pi \,\Gamma ({\displaystyle \frac {2}{3}} )\,\Gamma ( - a
 - {\displaystyle \frac {1}{3}}  + n)\,\Gamma (2\,a + 1 - n)\,
\Gamma (2\,a + {\displaystyle \frac {1}{3}}  - n)\,\Gamma ( - 3\,
a + 2\,n + {\displaystyle \frac {1}{2}} )}}  \\
\mbox{} + {\displaystyle \frac {\Gamma (a + {\displaystyle 
\frac {2}{3}} )\,\Gamma ({\displaystyle \frac {1}{6}}  + n - a)\,
\mathrm{V1}(a, \,n)\,\Gamma (a + {\displaystyle \frac {5}{6}} )}{
\sqrt{\pi }\,\Gamma ({\displaystyle \frac {2}{3}} )\,\Gamma ( - a
 - {\displaystyle \frac {1}{3}}  + n)\,\Gamma (3\,a + 1 - n)}} 
 \! \! \left. {\vrule height1.67em width0em depth1.67em} \right) 
\Gamma (a)\,\Gamma ({\displaystyle \frac {2}{3}} ) \left/ 
{\vrule height0.80em width0em depth0.80em} \right. \!  \! \Gamma 
({\displaystyle \frac {1}{2}}  + n) }
}
\end{maplelatex}

\begin{maplelatex}
\mapleinline{inert}{2d}{441, ":   ",
(-2/9*GAMMA(a+2/3)*GAMMA(1/6+n-a)*3^(1/2)/GAMMA(2/3)*GAMMA(a+1/3-n)*V1
(-1/6+n-a,n)*GAMMA(2/3+n-a)/GAMMA(-a-1/3+n)/GAMMA(2*a+1-n)/GAMMA(2*a+1
/3-n)/GAMMA(-3*a+2*n+1/2)/GAMMA(4/3+a-n)*GAMMA(a+5/6)*GAMMA(-a+1/3)+2/
9*GAMMA(a+2/3)*GAMMA(1/6+n-a)*Pi^(1/2)*3^(1/2)/GAMMA(2/3)*V1(a,n)/GAMM
A(-a-1/3+n)/GAMMA(1/3+n-2*a)/GAMMA(4/3+a-n)*GAMMA(a+5/6)*GAMMA(-a+1/3)
/GAMMA(3*a+1-n))*GAMMA(a)/GAMMA(1/2+n);}{%
\maplemultiline{
441\mbox{:~~~}  \left( {\vrule 
height1.67em width0em depth1.67em} \right. \!  \!  - 
{\displaystyle \frac {2}{9}} \,{\displaystyle \frac {\Gamma (a + 
{\displaystyle \frac {2}{3}} )\,\Gamma ({\displaystyle \frac {1}{
6}}  + n - a)\,\sqrt{3}\,\Gamma (a + {\displaystyle \frac {1}{3}
}  - n)\,\mathrm{V1}( - {\displaystyle \frac {1}{6}}  + n - a, \,
n)\,\Gamma ({\displaystyle \frac {2}{3}}  + n - a)\,\Gamma (a + 
{\displaystyle \frac {5}{6}} )\,\Gamma ( - a + {\displaystyle 
\frac {1}{3}} )}{\Gamma ({\displaystyle \frac {2}{3}} )\,\Gamma (
 - a - {\displaystyle \frac {1}{3}}  + n)\,\Gamma (2\,a + 1 - n)
\,\Gamma (2\,a + {\displaystyle \frac {1}{3}}  - n)\,\Gamma ( - 3
\,a + 2\,n + {\displaystyle \frac {1}{2}} )\,\Gamma (
{\displaystyle \frac {4}{3}}  + a - n)}}  \\
\mbox{} + {\displaystyle \frac {2}{9}} \,{\displaystyle \frac {
\Gamma (a + {\displaystyle \frac {2}{3}} )\,\Gamma (
{\displaystyle \frac {1}{6}}  + n - a)\,\sqrt{\pi }\,\sqrt{3}\,
\mathrm{V1}(a, \,n)\,\Gamma (a + {\displaystyle \frac {5}{6}} )\,
\Gamma ( - a + {\displaystyle \frac {1}{3}} )}{\Gamma (
{\displaystyle \frac {2}{3}} )\,\Gamma ( - a - {\displaystyle 
\frac {1}{3}}  + n)\,\Gamma ({\displaystyle \frac {1}{3}}  + n - 
2\,a)\,\Gamma ({\displaystyle \frac {4}{3}}  + a - n)\,\Gamma (3
\,a + 1 - n)}}  \! \! \left. {\vrule 
height1.67em width0em depth1.67em} \right) \Gamma (a) \left/ 
{\vrule height0.80em width0em depth0.80em} \right. \!  \!  
\Gamma ({\displaystyle \frac {1}{2}}  + n) }
}
\end{maplelatex}

\begin{maplelatex}
\mapleinline{inert}{2d}{442, ":   ",
(-GAMMA(a+2/3)*GAMMA(1/6+n-a)/Pi/GAMMA(2/3)*GAMMA(a+1/3-n)*V1(-1/6+n-a
,n)*GAMMA(2/3+n-a)/GAMMA(-a-1/3+n)/GAMMA(2*a+1-n)/GAMMA(2*a+1/3-n)/GAM
MA(-3*a+2*n+1/2)*GAMMA(1/2+n-a)+GAMMA(a+2/3)*GAMMA(1/6+n-a)/Pi^(1/2)/G
AMMA(2/3)*V1(a,n)/GAMMA(-a-1/3+n)/GAMMA(1/3+n-2*a)*GAMMA(1/2+n-a)/GAMM
A(3*a+1-n))*GAMMA(a)/GAMMA(1/2+n)*GAMMA(2/3)^2;}{%
\maplemultiline{
442\mbox{:~~~}  \left( {\vrule 
height1.67em width0em depth1.67em} \right. \!  \!  
 - {\displaystyle \frac {\Gamma (a + {\displaystyle \frac {2}{3}
} )\,\Gamma ({\displaystyle \frac {1}{6}}  + n - a)\,\Gamma (a + 
{\displaystyle \frac {1}{3}}  - n)\,\mathrm{V1}( - 
{\displaystyle \frac {1}{6}}  + n - a, \,n)\,\Gamma (
{\displaystyle \frac {2}{3}}  + n - a)\,\Gamma ({\displaystyle 
\frac {1}{2}}  + n - a)}{\pi \,\Gamma ({\displaystyle \frac {2}{3
}} )\,\Gamma ( - a - {\displaystyle \frac {1}{3}}  + n)\,\Gamma (
2\,a + 1 - n)\,\Gamma (2\,a + {\displaystyle \frac {1}{3}}  - n)
\,\Gamma ( - 3\,a + 2\,n + {\displaystyle \frac {1}{2}} )}}  \\
\mbox{} + {\displaystyle \frac {\Gamma (a + {\displaystyle 
\frac {2}{3}} )\,\Gamma ({\displaystyle \frac {1}{6}}  + n - a)\,
\mathrm{V1}(a, \,n)\,\Gamma ({\displaystyle \frac {1}{2}}  + n - 
a)}{\sqrt{\pi }\,\Gamma ({\displaystyle \frac {2}{3}} )\,\Gamma (
 - a - {\displaystyle \frac {1}{3}}  + n)\,\Gamma (
{\displaystyle \frac {1}{3}}  + n - 2\,a)\,\Gamma (3\,a + 1 - n)}
}  \! \! \left. {\vrule height1.67em width0em depth1.67em}
 \right) \Gamma (a)\,\Gamma ({\displaystyle \frac {2}{3}} )^{2}
 \left/ {\vrule height0.80em width0em depth0.80em} \right. \! 
 \! \Gamma ({\displaystyle \frac {1}{2}}  + n) }
}
\end{maplelatex}

\begin{maplelatex}
\mapleinline{inert}{2d}{443, ":   ",
(-2/9*GAMMA(a+2/3)*GAMMA(1/6+n-a)*Pi^(1/2)*3^(1/2)/GAMMA(2/3)*GAMMA(a+
1/3-n)*V1(-1/6+n-a,n)*GAMMA(2/3+n-a)/GAMMA(-a-1/3+n)/GAMMA(2*a+1-n)/GA
MMA(2*a+1/3-n)/GAMMA(-3*a+2*n+1/2)/GAMMA(4/3+a-n)+2/9*GAMMA(a+2/3)*GAM
MA(1/6+n-a)*Pi*3^(1/2)/GAMMA(2/3)*V1(a,n)/GAMMA(-a-1/3+n)/GAMMA(1/3+n-
2*a)/GAMMA(4/3+a-n)/GAMMA(3*a+1-n))*GAMMA(a)/GAMMA(1/2+n)*GAMMA(2/3);}
{%
\maplemultiline{
443\mbox{:~~~}  \left( {\vrule 
height1.67em width0em depth1.67em} \right. \!  \!  
 - {\displaystyle \frac {2}{9}} \,{\displaystyle \frac {\Gamma (a
 + {\displaystyle \frac {2}{3}} )\,\Gamma ({\displaystyle \frac {
1}{6}}  + n - a)\,\sqrt{\pi }\,\sqrt{3}\,\Gamma (a + 
{\displaystyle \frac {1}{3}}  - n)\,\mathrm{V1}( - 
{\displaystyle \frac {1}{6}}  + n - a, \,n)\,\Gamma (
{\displaystyle \frac {2}{3}}  + n - a)}{\Gamma ({\displaystyle 
\frac {2}{3}} )\,\Gamma ( - a - {\displaystyle \frac {1}{3}}  + n
)\,\Gamma (2\,a + 1 - n)\,\Gamma (2\,a + {\displaystyle \frac {1
}{3}}  - n)\,\Gamma ( - 3\,a + 2\,n + {\displaystyle \frac {1}{2}
} )\,\Gamma ({\displaystyle \frac {4}{3}}  + a - n)}}  \\
\mbox{} + {\displaystyle \frac {2}{9}} \,{\displaystyle \frac {
\Gamma (a + {\displaystyle \frac {2}{3}} )\,\Gamma (
{\displaystyle \frac {1}{6}}  + n - a)\,\pi \,\sqrt{3}\,\mathrm{
V1}(a, \,n)}{\Gamma ({\displaystyle \frac {2}{3}} )\,\Gamma ( - a
 - {\displaystyle \frac {1}{3}}  + n)\,\Gamma ({\displaystyle 
\frac {1}{3}}  + n - 2\,a)\,\Gamma ({\displaystyle \frac {4}{3}} 
 + a - n)\,\Gamma (3\,a + 1 - n)}}  \! \! \left. {\vrule 
height1.67em width0em depth1.67em} \right) \Gamma (a)\,\Gamma (
{\displaystyle \frac {2}{3}} ) 
 \left/ {\vrule height0.80em width0em depth0.80em} \right. \! 
 \! \Gamma ({\displaystyle \frac {1}{2}}  + n) }
}
\end{maplelatex}

\end{maplegroup}

%% file: AppendixB444to469.tex
\begin{maplegroup}
\mapleinline{inert}{2d}{444, ":   ",
(2/3*V1(a,n)*(2*sin(1/6*Pi*(-24*a+12*n-1))+1)*3^(1/2)*Pi*GAMMA(2*a+2/3
-n)*GAMMA(2*a+1-n)*GAMMA(a+1/3)*GAMMA(a+2/3)/GAMMA(1/6+a)/GAMMA(3*a+1-
n)/GAMMA(5/6+a-n)/GAMMA(a-n+1/2)/(cos(Pi*a)+cos(Pi*(-3*a+2*n))+sin(1/6
*Pi*(6*a-1))+sin(1/6*Pi*(-6*a+12*n-1)))+2/3*1/Pi^(1/2)*3^(1/2)*V1(-1/6
+n-a,n)*cos(Pi*(-3*a+2*n))/GAMMA(1/6+a)*GAMMA(a+1/3)/GAMMA(2*a+1/3-n)*
GAMMA(3*a-2*n+1/2)*GAMMA(1/2+n-a)*GAMMA(1/6+n-a)*GAMMA(a+2/3))/GAMMA(1
/2+n)/GAMMA(2/3);}{%
\maplemultiline{
444\mbox{:~~~} ({\displaystyle \frac {2}{3}} \mathrm{V1}
(a, \,n)\,(2\,\mathrm{sin}({\displaystyle \frac {\pi \,( - 24\,a
 + 12\,n - 1)}{6}} ) + 1)\,\sqrt{3}\,\pi \,\Gamma (2\,a + 
{\displaystyle \frac {2}{3}}  - n) \\
\Gamma (2\,a + 1 - n)\,\Gamma (a + {\displaystyle \frac {1}{3}} )
\,\Gamma (a + {\displaystyle \frac {2}{3}} ) \left/ {\vrule 
height0.80em width0em depth0.80em} \right. \!  \! (\Gamma (
{\displaystyle \frac {1}{6}}  + a)\,\Gamma (3\,a + 1 - n)\,\Gamma
 ({\displaystyle \frac {5}{6}}  + a - n) 
\Gamma (a - n + {\displaystyle \frac {1}{2}} ) \\
(\mathrm{cos}(\pi \,a) + \mathrm{cos}(\pi \,( - 3\,a + 2\,n)) + 
\mathrm{sin}({\displaystyle \frac {\pi \,(6\,a - 1)}{6}} ) + 
\mathrm{sin}({\displaystyle \frac {\pi \,( - 6\,a + 12\,n - 1)}{6
}} )))\mbox{} + {\displaystyle \frac {2}{3}}  \\
\sqrt{3}\,\mathrm{V1}( - {\displaystyle \frac {1}{6}}  + n - a, 
\,n)\,\mathrm{cos}(\pi \,( - 3\,a + 2\,n))\,\Gamma (a + 
{\displaystyle \frac {1}{3}} )\,\Gamma (3\,a - 2\,n + 
{\displaystyle \frac {1}{2}} )\,\Gamma ({\displaystyle \frac {1}{
2}}  + n - a) \\
\Gamma ({\displaystyle \frac {1}{6}}  + n - a)\,\Gamma (a + 
{\displaystyle \frac {2}{3}} ) \left/ {\vrule 
height0.80em width0em depth0.80em} \right. \!  \! (\sqrt{\pi }\,
\Gamma ({\displaystyle \frac {1}{6}}  + a)\,\Gamma (2\,a + 
{\displaystyle \frac {1}{3}}  - n))) \left/ {\vrule 
height0.80em width0em depth0.80em} \right. \!  \! (\Gamma (
{\displaystyle \frac {1}{2}}  + n)\,\Gamma ({\displaystyle 
\frac {2}{3}} )) }
}

\mapleresult
\begin{maplelatex}
\mapleinline{inert}{2d}{445, ":   ",
(2*(-sin(1/3*Pi*(-9*a+1+6*n))-sin(Pi*a)+3^(1/2)*cos(Pi*(-3*a+2*n)))*V1
(a,n)*GAMMA(1/2+n-a)*GAMMA(2*a+2/3-n)*GAMMA(a+2/3)*GAMMA(5/6-a)/Pi/GAM
MA(3*a+1-n)/GAMMA(n-2*a)/GAMMA(5/6+a-n)/(2*cos(2*Pi*(n-a))-1)+V1(-1/6+
n-a,n)*GAMMA(3*a-2*n+1/2)*GAMMA(a+2/3)*GAMMA(5/6-a)*GAMMA(1/2+n-a)*(co
s(Pi*(-3*a+n))+cos(3*Pi*(n-a))-sin(1/6*Pi*(-6*a+6*n+1))+sin(1/6*Pi*(-3
0*a+18*n-1))-3^(1/2)*sin(Pi*(-5*a+3*n))+3^(1/2)*sin(Pi*(n-a)))/Pi^(3/2
)/GAMMA(2*a+1/3-n)/GAMMA(5/6+a-n)/(2*cos(2*Pi*(n-a))-1))*GAMMA(a)/GAMM
A(1/2+n)*GAMMA(2/3);}{%
\maplemultiline{
445\mbox{:~~~} (2\,( - \mathrm{sin}({\displaystyle 
\frac {\pi \,( - 9\,a + 1 + 6\,n)}{3}} ) - \mathrm{sin}(\pi \,a)
 + \sqrt{3}\,\mathrm{cos}(\pi \,( - 3\,a + 2\,n)))\,\mathrm{V1}(a
, \,n) \\
\Gamma ({\displaystyle \frac {1}{2}}  + n - a)\,\Gamma (2\,a + 
{\displaystyle \frac {2}{3}}  - n)\,\Gamma (a + {\displaystyle 
\frac {2}{3}} )\,\Gamma ({\displaystyle \frac {5}{6}}  - a)
 \left/ {\vrule height0.80em width0em depth0.80em} \right. \! 
 \! (\pi \,\Gamma (3\,a + 1 - n)\,\Gamma (n - 2\,a) \\
\Gamma ({\displaystyle \frac {5}{6}}  + a - n)\,(2\,\mathrm{cos}(
2\,\pi \,(n - a)) - 1))\mbox{} + \mathrm{V1}( - {\displaystyle 
\frac {1}{6}}  + n - a, \,n)\,\Gamma (3\,a - 2\,n + 
{\displaystyle \frac {1}{2}} )\,\Gamma (a + {\displaystyle 
\frac {2}{3}} ) \\
\Gamma ({\displaystyle \frac {5}{6}}  - a)\,\Gamma (
{\displaystyle \frac {1}{2}}  + n - a)(\mathrm{cos}(\pi \,( - 3\,
a + n)) + \mathrm{cos}(3\,\pi \,(n - a)) \\
\mbox{} - \mathrm{sin}({\displaystyle \frac {\pi \,( - 6\,a + 6\,
n + 1)}{6}} ) + \mathrm{sin}({\displaystyle \frac {\pi \,( - 30\,
a + 18\,n - 1)}{6}} ) - \sqrt{3}\,\mathrm{sin}(\pi \,( - 5\,a + 3
\,n)) \\
\mbox{} + \sqrt{3}\,\mathrm{sin}(\pi \,(n - a))) \left/ {\vrule 
height0.80em width0em depth0.80em} \right. \!  \! (\pi ^{(3/2)}\,
\Gamma (2\,a + {\displaystyle \frac {1}{3}}  - n)\,\Gamma (
{\displaystyle \frac {5}{6}}  + a - n) \\
(2\,\mathrm{cos}(2\,\pi \,(n - a)) - 1)))\Gamma (a)\,\Gamma (
{\displaystyle \frac {2}{3}} ) \left/ {\vrule 
height0.80em width0em depth0.80em} \right. \!  \! \Gamma (
{\displaystyle \frac {1}{2}}  + n) }
}
\end{maplelatex}

\begin{maplelatex}
\mapleinline{inert}{2d}{446, ":   ",
(2*sin(Pi*(n-2*a))*V1(a,n)*GAMMA(a+2/3)*GAMMA(2*a+1-n)*GAMMA(1/2+n-a)*
GAMMA(-1/6+n-a)/GAMMA(1/6+a)/GAMMA(5/6+a-n)/GAMMA(3*a+1-n)/GAMMA(1/3+n
-2*a)/(cos(Pi*(n-2*a))+sin(1/6*Pi*(-1+6*n)))+V1(-1/6+n-a,n)/Pi^(1/2)/G
AMMA(-3*a+2*n+1/2)/GAMMA(2*a+1/3-n)*GAMMA(1/2+n-a)*GAMMA(1/6+n-a)*GAMM
A(-1/6+n-a)/GAMMA(1/6+a)*GAMMA(a+2/3))*GAMMA(2*a+5/6-n)/GAMMA(1/2+n);}
{%
\maplemultiline{
446\mbox{:~~~}  \left( {\vrule 
height1.67em width0em depth1.67em} \right. \!  \! {\displaystyle 
\frac {2\,\mathrm{sin}(\pi \,(n - 2\,a))\,\mathrm{V1}(a, \,n)\,
\Gamma (a + {\displaystyle \frac {2}{3}} )\,\Gamma (2\,a + 1 - n)
\,\Gamma ({\displaystyle \frac {1}{2}}  + n - a)\,\Gamma ( - 
{\displaystyle \frac {1}{6}}  + n - a)}{\Gamma ({\displaystyle 
\frac {1}{6}}  + a)\,\Gamma ({\displaystyle \frac {5}{6}}  + a - 
n)\,\Gamma (3\,a + 1 - n)\,\Gamma ({\displaystyle \frac {1}{3}} 
 + n - 2\,a)\,(\mathrm{cos}(\pi \,(n - 2\,a)) + \mathrm{sin}(
{\displaystyle \frac {\pi \,( - 1 + 6\,n)}{6}} ))}}  \\
\mbox{} + {\displaystyle \frac {\mathrm{V1}( - {\displaystyle 
\frac {1}{6}}  + n - a, \,n)\,\Gamma ({\displaystyle \frac {1}{2}
}  + n - a)\,\Gamma ({\displaystyle \frac {1}{6}}  + n - a)\,
\Gamma ( - {\displaystyle \frac {1}{6}}  + n - a)\,\Gamma (a + 
{\displaystyle \frac {2}{3}} )}{\sqrt{\pi }\,\Gamma ( - 3\,a + 2
\,n + {\displaystyle \frac {1}{2}} )\,\Gamma (2\,a + 
{\displaystyle \frac {1}{3}}  - n)\,\Gamma ({\displaystyle 
\frac {1}{6}}  + a)}}  \! \! \left. {\vrule 
height1.67em width0em depth1.67em} \right)  
\Gamma (2\,a + {\displaystyle \frac {5}{6}}  - n) \left/ {\vrule 
height0.80em width0em depth0.80em} \right. \!  \! \Gamma (
{\displaystyle \frac {1}{2}}  + n) }
}
\end{maplelatex}

\begin{maplelatex}
\mapleinline{inert}{2d}{447, ":   ",
(-V1(a,n)*(2*cos(Pi*(-3*a+n))+cos(3*Pi*(n-a)))*GAMMA(4/3+a-n)*GAMMA(a+
1/3)*GAMMA(1/2+n-a)*GAMMA(a-n+1)*GAMMA(a+2/3)/Pi^(3/2)/GAMMA(5/6+a-n)/
GAMMA(3*a+1-n)/(2*cos(2*Pi*(n-a))-1)+2*V1(-1/6+n-a,n)*GAMMA(4/3+a-n)*G
AMMA(3*a-2*n+1/2)*GAMMA(a+1/3)*GAMMA(a-n+1)*GAMMA(1/2+n-a)*GAMMA(a+2/3
)*(cos(Pi*(-3*a+n))+cos(3*Pi*(n-a)))/Pi/GAMMA(2*a+2/3-n)/GAMMA(5/6+a-n
)/GAMMA(2*a+1-n)/GAMMA(2*a+1/3-n)/(2*cos(2*Pi*(n-a))-1))*GAMMA(a)/GAMM
A(1/2+n);}{%
\maplemultiline{
447\mbox{:~~~} ( - \mathrm{V1}(a, \,n)\,(2\,\mathrm{cos}
(\pi \,( - 3\,a + n)) + \mathrm{cos}(3\,\pi \,(n - a)))\,\Gamma (
{\displaystyle \frac {4}{3}}  + a - n)\,\Gamma (a + 
{\displaystyle \frac {1}{3}} ) \\
\Gamma ({\displaystyle \frac {1}{2}}  + n - a)\,\Gamma (a - n + 1
)\,\Gamma (a + {\displaystyle \frac {2}{3}} ) \left/ {\vrule 
height0.80em width0em depth0.80em} \right. \!  \! (\pi ^{(3/2)}\,
\Gamma ({\displaystyle \frac {5}{6}}  + a - n)\,\Gamma (3\,a + 1
 - n) \\
(2\,\mathrm{cos}(2\,\pi \,(n - a)) - 1))\mbox{} + 2\,\mathrm{V1}(
 - {\displaystyle \frac {1}{6}}  + n - a, \,n)\,\Gamma (
{\displaystyle \frac {4}{3}}  + a - n)\,\Gamma (3\,a - 2\,n + 
{\displaystyle \frac {1}{2}} ) \\
\Gamma (a + {\displaystyle \frac {1}{3}} )\,\Gamma (a - n + 1)\,
\Gamma ({\displaystyle \frac {1}{2}}  + n - a)\,\Gamma (a + 
{\displaystyle \frac {2}{3}} )\,(\mathrm{cos}(\pi \,( - 3\,a + n)
) + \mathrm{cos}(3\,\pi \,(n - a))) \\
 \left/ {\vrule height0.80em width0em depth0.80em} \right. \! 
 \! (\pi \,\Gamma (2\,a + {\displaystyle \frac {2}{3}}  - n)\,
\Gamma ({\displaystyle \frac {5}{6}}  + a - n)\,\Gamma (2\,a + 1
 - n)\,\Gamma (2\,a + {\displaystyle \frac {1}{3}}  - n) \\
(2\,\mathrm{cos}(2\,\pi \,(n - a)) - 1)))\Gamma (a) \left/ 
{\vrule height0.80em width0em depth0.80em} \right. \!  \! \Gamma 
({\displaystyle \frac {1}{2}}  + n) }
}
\end{maplelatex}

\begin{maplelatex}
\mapleinline{inert}{2d}{448, ":   ",
(-4/3*V1(a,n)*(3^(1/2)*sin(1/3*Pi*(-9*a+1+6*n))+3^(1/2)*sin(Pi*a)-3*co
s(Pi*(-3*a+2*n)))*Pi^(1/2)*GAMMA(1/2+n-a)*GAMMA(2*a+2/3-n)*GAMMA(a+2/3
)/GAMMA(3*a+1-n)/GAMMA(n-2*a)/GAMMA(5/6+a-n)/(2*cos(2*Pi*(n-a))-1)+2/3
*V1(-1/6+n-a,n)*GAMMA(3*a-2*n+1/2)*GAMMA(1/2+n-a)*GAMMA(a+2/3)*(3^(1/2
)*cos(Pi*(-3*a+n))+3^(1/2)*cos(3*Pi*(n-a))-3^(1/2)*sin(1/6*Pi*(-6*a+6*
n+1))+3^(1/2)*sin(1/6*Pi*(-30*a+18*n-1))-3*sin(Pi*(-5*a+3*n))+3*sin(Pi
*(n-a)))/GAMMA(5/6+a-n)/GAMMA(2*a+1/3-n)/(2*cos(2*Pi*(n-a))-1))/GAMMA(
1/2+n);}{%
\maplemultiline{
448\mbox{:~~~} ( - {\displaystyle \frac {4}{3}} \mathrm{
V1}(a, \,n)\,(\sqrt{3}\,\mathrm{sin}({\displaystyle \frac {\pi \,
( - 9\,a + 1 + 6\,n)}{3}} ) + \sqrt{3}\,\mathrm{sin}(\pi \,a) - 3
\,\mathrm{cos}(\pi \,( - 3\,a + 2\,n))) \\
\sqrt{\pi }\,\Gamma ({\displaystyle \frac {1}{2}}  + n - a)\,
\Gamma (2\,a + {\displaystyle \frac {2}{3}}  - n)\,\Gamma (a + 
{\displaystyle \frac {2}{3}} ) \left/ {\vrule 
height0.80em width0em depth0.80em} \right. \!  \! (\Gamma (3\,a
 + 1 - n)\,\Gamma (n - 2\,a) \\
\Gamma ({\displaystyle \frac {5}{6}}  + a - n)\,(2\,\mathrm{cos}(
2\,\pi \,(n - a)) - 1))\mbox{} + {\displaystyle \frac {2}{3}} 
\mathrm{V1}( - {\displaystyle \frac {1}{6}}  + n - a, \,n)\,
\Gamma (3\,a - 2\,n + {\displaystyle \frac {1}{2}} ) \\
\Gamma ({\displaystyle \frac {1}{2}}  + n - a)\,\Gamma (a + 
{\displaystyle \frac {2}{3}} )(\sqrt{3}\,\mathrm{cos}(\pi \,( - 3
\,a + n)) + \sqrt{3}\,\mathrm{cos}(3\,\pi \,(n - a)) \\
\mbox{} - \sqrt{3}\,\mathrm{sin}({\displaystyle \frac {\pi \,( - 
6\,a + 6\,n + 1)}{6}} ) + \sqrt{3}\,\mathrm{sin}({\displaystyle 
\frac {\pi \,( - 30\,a + 18\,n - 1)}{6}} ) \\
\mbox{} - 3\,\mathrm{sin}(\pi \,( - 5\,a + 3\,n)) + 3\,\mathrm{
sin}(\pi \,(n - a))) \left/ {\vrule 
height0.80em width0em depth0.80em} \right. \!  \! (\Gamma (
{\displaystyle \frac {5}{6}}  + a - n)\,\Gamma (2\,a + 
{\displaystyle \frac {1}{3}}  - n) \\
(2\,\mathrm{cos}(2\,\pi \,(n - a)) - 1))) \left/ {\vrule 
height0.80em width0em depth0.80em} \right. \!  \! \Gamma (
{\displaystyle \frac {1}{2}}  + n) }
}
\end{maplelatex}

\begin{maplelatex}
\mapleinline{inert}{2d}{449, ":   ",
(-1/2*(2*cos(Pi*(-3*a+n))+cos(3*Pi*(n-a)))*V1(a,n)*GAMMA(4/3+a-n)*GAMM
A(a+1/3)*GAMMA(2*a+2/3-n)*GAMMA(-1/6+n-a)*GAMMA(1/2+n-a)*GAMMA(a+2/3)/
sin(1/6*Pi*(-6*a+6*n+1))/((-1)^n)/Pi^(3/2)/(2*n-1)/GAMMA(-1/2+n)^2/GAM
MA(3*a+1-n)/GAMMA(5/6+a-n)+V1(-1/6+n-a,n)*GAMMA(4/3+a-n)*GAMMA(a+1/3)*
GAMMA(a+2/3)*GAMMA(3*a-2*n+1/2)*GAMMA(-1/6+n-a)*GAMMA(1/2+n-a)*(cos(Pi
*(-3*a+n))+cos(3*Pi*(n-a)))/sin(1/6*Pi*(-6*a+6*n+1))/Pi/((-1)^n)/(2*n-
1)/GAMMA(-1/2+n)^2/GAMMA(2*a+1/3-n)/GAMMA(2*a+1-n)/GAMMA(5/6+a-n))*GAM
MA(a);}{%
\maplemultiline{
449\mbox{:~~~} ( - {\displaystyle \frac {1}{2}} (2\,
\mathrm{cos}(\pi \,( - 3\,a + n)) + \mathrm{cos}(3\,\pi \,(n - a)
))\,\mathrm{V1}(a, \,n)\,\Gamma ({\displaystyle \frac {4}{3}}  + 
a - n)\,\Gamma (a + {\displaystyle \frac {1}{3}} ) \\
\Gamma (2\,a + {\displaystyle \frac {2}{3}}  - n)\,\Gamma ( - 
{\displaystyle \frac {1}{6}}  + n - a)\,\Gamma ({\displaystyle 
\frac {1}{2}}  + n - a)\,\Gamma (a + {\displaystyle \frac {2}{3}
} ) \left/ {\vrule height0.80em width0em depth0.80em} \right. \! 
 \! (\mathrm{sin}({\displaystyle \frac {\pi \,( - 6\,a + 6\,n + 1
)}{6}} ) \\
(-1)^{n}\,\pi ^{(3/2)}\,(2\,n - 1)\,\Gamma ( - {\displaystyle 
\frac {1}{2}}  + n)^{2}\,\Gamma (3\,a + 1 - n)\,\Gamma (
{\displaystyle \frac {5}{6}}  + a - n))\mbox{} + \mathrm{V1}( - 
{\displaystyle \frac {1}{6}}  + n - a, \,n) \\
\Gamma ({\displaystyle \frac {4}{3}}  + a - n)\,\Gamma (a + 
{\displaystyle \frac {1}{3}} )\,\Gamma (a + {\displaystyle 
\frac {2}{3}} )\,\Gamma (3\,a - 2\,n + {\displaystyle \frac {1}{2
}} )\,\Gamma ( - {\displaystyle \frac {1}{6}}  + n - a)\,\Gamma (
{\displaystyle \frac {1}{2}}  + n - a) \\
(\mathrm{cos}(\pi \,( - 3\,a + n)) + \mathrm{cos}(3\,\pi \,(n - a
))) \left/ {\vrule height0.80em width0em depth0.80em} \right. \! 
 \! (\mathrm{sin}({\displaystyle \frac {\pi \,( - 6\,a + 6\,n + 1
)}{6}} )\,\pi \,(-1)^{n} \\
(2\,n - 1)\,\Gamma ( - {\displaystyle \frac {1}{2}}  + n)^{2}\,
\Gamma (2\,a + {\displaystyle \frac {1}{3}}  - n)\,\Gamma (2\,a
 + 1 - n)\,\Gamma ({\displaystyle \frac {5}{6}}  + a - n)))\Gamma
 (a) }
}
\end{maplelatex}

\begin{maplelatex}
\mapleinline{inert}{2d}{450, ":   ",
(-1/8*(2*cos(Pi*(-3*a+n))+cos(3*Pi*(n-a)))*(2*n-1)*GAMMA(a+1/3)*GAMMA(
a+2/3)*GAMMA(-1/6+n-a)*GAMMA(1/6+n-a)*GAMMA(a-n+1)*GAMMA(2*a+1-n)*V1(a
,n)/cos(Pi*(n-a))/((-1)^n)/Pi^(3/2)/GAMMA(a-n+1/2)/GAMMA(3*a+1-n)+1/4*
V1(-1/6+n-a,n)*GAMMA(a-n+1)*GAMMA(a+1/3)*GAMMA(3*a-2*n+1/2)*GAMMA(a+2/
3)*GAMMA(1/6+n-a)*GAMMA(-1/6+n-a)*(cos(Pi*(-3*a+n))+cos(3*Pi*(n-a)))*(
2*n-1)/cos(Pi*(n-a))/Pi/((-1)^n)/GAMMA(2*a+2/3-n)/GAMMA(a-n+1/2)/GAMMA
(2*a+1/3-n))*GAMMA(a)/GAMMA(1/2+n)^2;}{%
\maplemultiline{
450\mbox{:~~~} ( - {\displaystyle \frac {1}{8}} (2\,
\mathrm{cos}(\pi \,( - 3\,a + n)) + \mathrm{cos}(3\,\pi \,(n - a)
))\,(2\,n - 1)\,\Gamma (a + {\displaystyle \frac {1}{3}} )\,
\Gamma (a + {\displaystyle \frac {2}{3}} ) \\
\Gamma ( - {\displaystyle \frac {1}{6}}  + n - a)\,\Gamma (
{\displaystyle \frac {1}{6}}  + n - a)\,\Gamma (a - n + 1)\,
\Gamma (2\,a + 1 - n)\,\mathrm{V1}(a, \,n) \left/ {\vrule 
height0.80em width0em depth0.80em} \right. \!  \! ( \\
\mathrm{cos}(\pi \,(n - a))\,(-1)^{n}\,\pi ^{(3/2)}\,\Gamma (a - 
n + {\displaystyle \frac {1}{2}} )\,\Gamma (3\,a + 1 - n))\mbox{}
 + {\displaystyle \frac {1}{4}} \mathrm{V1}( - {\displaystyle 
\frac {1}{6}}  + n - a, \,n) \\
\Gamma (a - n + 1)\,\Gamma (a + {\displaystyle \frac {1}{3}} )\,
\Gamma (3\,a - 2\,n + {\displaystyle \frac {1}{2}} )\,\Gamma (a
 + {\displaystyle \frac {2}{3}} )\,\Gamma ({\displaystyle \frac {
1}{6}}  + n - a)\,\Gamma ( - {\displaystyle \frac {1}{6}}  + n - 
a) \\
(\mathrm{cos}(\pi \,( - 3\,a + n)) + \mathrm{cos}(3\,\pi \,(n - a
)))\,(2\,n - 1) \left/ {\vrule height0.80em width0em depth0.80em}
 \right. \!  \! (\mathrm{cos}(\pi \,(n - a))\,\pi \,(-1)^{n} \\
\Gamma (2\,a + {\displaystyle \frac {2}{3}}  - n)\,\Gamma (a - n
 + {\displaystyle \frac {1}{2}} )\,\Gamma (2\,a + {\displaystyle 
\frac {1}{3}}  - n)))\Gamma (a) \left/ {\vrule 
height0.84em width0em depth0.84em} \right. \!  \! \Gamma (
{\displaystyle \frac {1}{2}}  + n)^{2} }
}
\end{maplelatex}

\begin{maplelatex}
\mapleinline{inert}{2d}{451, ":   ",
(1/2*V1(a,n)*(2*cos(Pi*(-3*a+n))+cos(3*Pi*(n-a)))*GAMMA(4/3+a-n)*GAMMA
(3*a-2*n+3/2)*GAMMA(1/6+n-a)*GAMMA(-1/6+n-a)*GAMMA(a+2/3)/Pi^(3/2)/GAM
MA(3*a+1-n)/GAMMA(a-n+1/2)/(sin(1/3*Pi*(1+3*n))-sin(1/3*Pi*(-6*a-1+3*n
)))-V1(-1/6+n-a,n)*GAMMA(4/3+a-n)*GAMMA(2/3+n-2*a)*GAMMA(3*a-2*n+3/2)^
2*GAMMA(1/6+n-a)*GAMMA(-1/6+n-a)*GAMMA(a+2/3)*(sin(1/3*Pi*(3*a-1))+sin
(1/3*Pi*(-15*a+12*n-1))-sin(1/3*Pi*(-3*a+6*n+1))+sin(1/3*Pi*(-15*a+6*n
-1)))/Pi^2/GAMMA(2*a+2/3-n)/GAMMA(a-n+1/2)/GAMMA(2*a+1-n)/(sin(1/3*Pi*
(1+3*n))-sin(1/3*Pi*(-6*a-1+3*n)))/(-6*a+4*n-1))*GAMMA(a)/GAMMA(1/2+n)
;}{%
\maplemultiline{
451\mbox{:~~~} ({\displaystyle \frac {1}{2}} \mathrm{V1}
(a, \,n)\,(2\,\mathrm{cos}(\pi \,( - 3\,a + n)) + \mathrm{cos}(3
\,\pi \,(n - a)))\,\Gamma ({\displaystyle \frac {4}{3}}  + a - n)
 \\
\Gamma (3\,a - 2\,n + {\displaystyle \frac {3}{2}} )\,\Gamma (
{\displaystyle \frac {1}{6}}  + n - a)\,\Gamma ( - 
{\displaystyle \frac {1}{6}}  + n - a)\,\Gamma (a + 
{\displaystyle \frac {2}{3}} ) \left/ {\vrule 
height0.80em width0em depth0.80em} \right. \!  \! (\pi ^{(3/2)}\,
\Gamma (3\,a + 1 - n) \\
\Gamma (a - n + {\displaystyle \frac {1}{2}} )\,(\mathrm{sin}(
{\displaystyle \frac {\pi \,(1 + 3\,n)}{3}} ) - \mathrm{sin}(
{\displaystyle \frac {\pi \,( - 6\,a - 1 + 3\,n)}{3}} )))\mbox{}
 - \mathrm{V1}( - {\displaystyle \frac {1}{6}}  + n - a, \,n) \\
\Gamma ({\displaystyle \frac {4}{3}}  + a - n)\,\Gamma (
{\displaystyle \frac {2}{3}}  + n - 2\,a)\,\Gamma (3\,a - 2\,n + 
{\displaystyle \frac {3}{2}} )^{2}\,\Gamma ({\displaystyle 
\frac {1}{6}}  + n - a)\,\Gamma ( - {\displaystyle \frac {1}{6}} 
 + n - a)\,\Gamma (a + {\displaystyle \frac {2}{3}} )( \\
\mathrm{sin}({\displaystyle \frac {\pi \,(3\,a - 1)}{3}} ) + 
\mathrm{sin}({\displaystyle \frac {\pi \,( - 15\,a + 12\,n - 1)}{
3}} ) - \mathrm{sin}({\displaystyle \frac {\pi \,( - 3\,a + 6\,n
 + 1)}{3}} ) \\
\mbox{} + \mathrm{sin}({\displaystyle \frac {\pi \,( - 15\,a + 6
\,n - 1)}{3}} )) \left/ {\vrule height0.80em width0em depth0.80em
} \right. \!  \! (\pi ^{2}\,\Gamma (2\,a + {\displaystyle \frac {
2}{3}}  - n)\,\Gamma (a - n + {\displaystyle \frac {1}{2}} )\,
\Gamma (2\,a + 1 - n) \\
(\mathrm{sin}({\displaystyle \frac {\pi \,(1 + 3\,n)}{3}} ) - 
\mathrm{sin}({\displaystyle \frac {\pi \,( - 6\,a - 1 + 3\,n)}{3}
} ))\,( - 6\,a + 4\,n - 1)))\Gamma (a) \left/ {\vrule 
height0.80em width0em depth0.80em} \right. \!  \! \Gamma (
{\displaystyle \frac {1}{2}}  + n) }
}
\end{maplelatex}

\begin{maplelatex}
\mapleinline{inert}{2d}{452, ":   ",
(-1/8*V1(a,n)*GAMMA(a+1/3)*GAMMA(3*a-2*n+3/2)*GAMMA(1/6+n-a)*GAMMA(1/2
+n-a)*GAMMA(a+2/3)*(2*cos(Pi*(-3*a+2*n))+2*sin(1/6*Pi*(6*a+1))+(-1)^(-
n)*sin(1/6*Pi*(-1+6*n-6*a))+(-1)^(-n)*cos(Pi*(-3*a+n))-(-1)^(-n)*sin(1
/6*Pi*(1+6*a+6*n)))*(2*n-1)/Pi/sin(1/6*Pi*(-1+6*n-6*a))/GAMMA(7/6+a-n)
/GAMMA(3*a+1-n)-cos(Pi*(-3*a+2*n))*(2*n-1)*V1(-1/6+n-a,n)/(-6*a+4*n-1)
/Pi^(3/2)/GAMMA(2*a+1/3-n)*GAMMA(-1/6+n-a)*GAMMA(3*a-2*n+3/2)^2*GAMMA(
1/6+n-a)*GAMMA(1/2+n-a)/GAMMA(2*a+1-n)/GAMMA(2*a+2/3-n)*GAMMA(a+1/3)*G
AMMA(a+2/3))*GAMMA(a)/GAMMA(1/2+n)^2;}{%
\maplemultiline{
452\mbox{:~~~} ( - {\displaystyle \frac {1}{8}} \mathrm{
V1}(a, \,n)\,\Gamma (a + {\displaystyle \frac {1}{3}} )\,\Gamma (
3\,a - 2\,n + {\displaystyle \frac {3}{2}} )\,\Gamma (
{\displaystyle \frac {1}{6}}  + n - a)\,\Gamma ({\displaystyle 
\frac {1}{2}}  + n - a)\,\Gamma (a + {\displaystyle \frac {2}{3}
} )( \\
2\,\mathrm{cos}(\pi \,( - 3\,a + 2\,n)) + 2\,\mathrm{sin}(
{\displaystyle \frac {\pi \,(6\,a + 1)}{6}} ) + (-1)^{( - n)}\,
\mathrm{sin}({\displaystyle \frac {\pi \,( - 1 + 6\,n - 6\,a)}{6}
} ) \\
\mbox{} + (-1)^{( - n)}\,\mathrm{cos}(\pi \,( - 3\,a + n)) - (-1)
^{( - n)}\,\mathrm{sin}({\displaystyle \frac {\pi \,(1 + 6\,a + 6
\,n)}{6}} ))(2\,n - 1) \left/ {\vrule 
height0.80em width0em depth0.80em} \right. \!  \! (\pi  \\
\mathrm{sin}({\displaystyle \frac {\pi \,( - 1 + 6\,n - 6\,a)}{6}
} )\,\Gamma ({\displaystyle \frac {7}{6}}  + a - n)\,\Gamma (3\,a
 + 1 - n))\mbox{} - \mathrm{cos}(\pi \,( - 3\,a + 2\,n))\,(2\,n
 - 1) \\
\mathrm{V1}( - {\displaystyle \frac {1}{6}}  + n - a, \,n)\,
\Gamma ( - {\displaystyle \frac {1}{6}}  + n - a)\,\Gamma (3\,a
 - 2\,n + {\displaystyle \frac {3}{2}} )^{2}\,\Gamma (
{\displaystyle \frac {1}{6}}  + n - a)\,\Gamma ({\displaystyle 
\frac {1}{2}}  + n - a) \\
\Gamma (a + {\displaystyle \frac {1}{3}} )\,\Gamma (a + 
{\displaystyle \frac {2}{3}} ) \left/ {\vrule 
height0.84em width0em depth0.84em} \right. \!  \! (( - 6\,a + 4\,
n - 1)\,\pi ^{(3/2)}\,\Gamma (2\,a + {\displaystyle \frac {1}{3}
}  - n)\,\Gamma (2\,a + 1 - n) \\
\Gamma (2\,a + {\displaystyle \frac {2}{3}}  - n)))\Gamma (a)
 \left/ {\vrule height0.87em width0em depth0.87em} \right. \! 
 \! \Gamma ({\displaystyle \frac {1}{2}}  + n)^{2} }
}
\end{maplelatex}

\begin{maplelatex}
\mapleinline{inert}{2d}{453, ":   ",
(-1/2*V1(a,n)*(2*cos(Pi*(-3*a+n))+cos(3*Pi*(n-a)))*GAMMA(a-n+1)*GAMMA(
3*a-2*n+3/2)*GAMMA(1/2+n-a)*GAMMA(-1/6+n-a)*GAMMA(a+1/3)*GAMMA(a+2/3)/
Pi^(3/2)/GAMMA(5/6+a-n)/GAMMA(3*a+1-n)/(sin(1/3*Pi*(-6*a-1+3*n))-sin(1
/3*Pi*(3*n-1)))-1/4*V1(-1/6+n-a,n)*GAMMA(3*a-2*n+1/2)^2*GAMMA(2/3+n-2*
a)*GAMMA(1/2+n-a)*GAMMA(-1/6+n-a)*GAMMA(a-n+1)*GAMMA(a+1/3)*GAMMA(a+2/
3)*(sin(1/3*Pi*(-3*a+6*n+1))-sin(1/3*Pi*(-15*a+6*n-1))-sin(1/3*Pi*(3*a
-1))-sin(1/3*Pi*(-15*a+12*n-1)))*(-6*a+4*n-1)/Pi^2/GAMMA(2*a+1-n)/GAMM
A(5/6+a-n)/GAMMA(2*a+2/3-n)/(sin(1/3*Pi*(-6*a-1+3*n))-sin(1/3*Pi*(3*n-
1))))/GAMMA(1/2+n);}{%
\maplemultiline{
453\mbox{:~~~} ( - {\displaystyle \frac {1}{2}} \mathrm{
V1}(a, \,n)\,(2\,\mathrm{cos}(\pi \,( - 3\,a + n)) + \mathrm{cos}
(3\,\pi \,(n - a)))\,\Gamma (a - n + 1) \\
\Gamma (3\,a - 2\,n + {\displaystyle \frac {3}{2}} )\,\Gamma (
{\displaystyle \frac {1}{2}}  + n - a)\,\Gamma ( - 
{\displaystyle \frac {1}{6}}  + n - a)\,\Gamma (a + 
{\displaystyle \frac {1}{3}} )\,\Gamma (a + {\displaystyle 
\frac {2}{3}} ) \left/ {\vrule height0.80em width0em depth0.80em}
 \right. \!  \! \\ (\pi ^{(3/2)} 
\Gamma ({\displaystyle \frac {5}{6}}  + a - n)\,\Gamma (3\,a + 1
 - n)\,(\mathrm{sin}({\displaystyle \frac {\pi \,( - 6\,a - 1 + 3
\,n)}{3}} ) - \mathrm{sin}({\displaystyle \frac {\pi \,(3\,n - 1)
}{3}} )))\mbox{} - {\displaystyle \frac {1}{4}}  \\
\mathrm{V1}( - {\displaystyle \frac {1}{6}}  + n - a, \,n)\,
\Gamma (3\,a - 2\,n + {\displaystyle \frac {1}{2}} )^{2}\,\Gamma 
({\displaystyle \frac {2}{3}}  + n - 2\,a)\,\Gamma (
{\displaystyle \frac {1}{2}}  + n - a)\,\Gamma ( - 
{\displaystyle \frac {1}{6}}  + n - a) \\
\Gamma (a - n + 1)\,\Gamma (a + {\displaystyle \frac {1}{3}} )\,
\Gamma (a + {\displaystyle \frac {2}{3}} )(\mathrm{sin}(
{\displaystyle \frac {\pi \,( - 3\,a + 6\,n + 1)}{3}} ) - 
\mathrm{sin}({\displaystyle \frac {\pi \,( - 15\,a + 6\,n - 1)}{3
}} ) \\
\mbox{} - \mathrm{sin}({\displaystyle \frac {\pi \,(3\,a - 1)}{3}
} ) - \mathrm{sin}({\displaystyle \frac {\pi \,( - 15\,a + 12\,n
 - 1)}{3}} ))( - 6\,a + 4\,n - 1) \left/ {\vrule 
height0.80em width0em depth0.80em} \right. \!  \! (\pi ^{2} \\
\Gamma (2\,a + 1 - n)\,\Gamma ({\displaystyle \frac {5}{6}}  + a
 - n)\,\Gamma (2\,a + {\displaystyle \frac {2}{3}}  - n) \\
(\mathrm{sin}({\displaystyle \frac {\pi \,( - 6\,a - 1 + 3\,n)}{3
}} ) - \mathrm{sin}({\displaystyle \frac {\pi \,(3\,n - 1)}{3}} )
))) \left/ {\vrule height0.80em width0em depth0.80em} \right. \! 
 \! \Gamma ({\displaystyle \frac {1}{2}}  + n) }
}
\end{maplelatex}

\begin{maplelatex}
\mapleinline{inert}{2d}{454, ":   ",
(-(2*cos(Pi*(-3*a+n))+cos(3*Pi*(n-a)))*V1(a,n)*GAMMA(a-n+1)*GAMMA(5/3+
a-n)*GAMMA(1/2+n-a)*GAMMA(a+1/3)*GAMMA(a+2/3)/Pi^(3/2)/GAMMA(7/6+a-n)/
GAMMA(3*a+1-n)/(2*cos(2*Pi*(n-a))-1)+2*V1(-1/6+n-a,n)*(-sin(1/3*Pi*(6*
a-1))+sin(1/3*Pi*(-6*a+6*n+1)))*GAMMA(a-n+1)*GAMMA(1/3+n-2*a)*GAMMA(5/
3+a-n)*GAMMA(1/2+n-a)*GAMMA(a+1/3)*GAMMA(a+2/3)/Pi/GAMMA(2*a+1-n)/GAMM
A(-3*a+2*n+1/2)/GAMMA(7/6+a-n)/GAMMA(2*a+1/3-n)/(2*cos(2*Pi*(n-a))-1))
*GAMMA(a)/GAMMA(1/2+n);}{%
\maplemultiline{
454\mbox{:~~~} ( - (2\,\mathrm{cos}(\pi \,( - 3\,a + n))
 + \mathrm{cos}(3\,\pi \,(n - a)))\,\mathrm{V1}(a, \,n)\,\Gamma (
a - n + 1)\,\Gamma ({\displaystyle \frac {5}{3}}  + a - n) \\
\Gamma ({\displaystyle \frac {1}{2}}  + n - a)\,\Gamma (a + 
{\displaystyle \frac {1}{3}} )\,\Gamma (a + {\displaystyle 
\frac {2}{3}} ) \left/ {\vrule height0.80em width0em depth0.80em}
 \right. \!  \! (\pi ^{(3/2)}\,\Gamma ({\displaystyle \frac {7}{6
}}  + a - n)\,\Gamma (3\,a + 1 - n) \\
(2\,\mathrm{cos}(2\,\pi \,(n - a)) - 1))\mbox{} + 2\,\mathrm{V1}(
 - {\displaystyle \frac {1}{6}}  + n - a, \,n) \\
( - \mathrm{sin}({\displaystyle \frac {\pi \,(6\,a - 1)}{3}} ) + 
\mathrm{sin}({\displaystyle \frac {\pi \,( - 6\,a + 6\,n + 1)}{3}
} ))\,\Gamma (a - n + 1)\,\Gamma ({\displaystyle \frac {1}{3}} 
 + n - 2\,a) \\
\Gamma ({\displaystyle \frac {5}{3}}  + a - n)\,\Gamma (
{\displaystyle \frac {1}{2}}  + n - a)\,\Gamma (a + 
{\displaystyle \frac {1}{3}} )\,\Gamma (a + {\displaystyle 
\frac {2}{3}} ) \left/ {\vrule height0.80em width0em depth0.80em}
 \right. \!  \! (\pi \,\Gamma (2\,a + 1 - n) \\
\Gamma ( - 3\,a + 2\,n + {\displaystyle \frac {1}{2}} )\,\Gamma (
{\displaystyle \frac {7}{6}}  + a - n)\,\Gamma (2\,a + 
{\displaystyle \frac {1}{3}}  - n)\,(2\,\mathrm{cos}(2\,\pi \,(n
 - a)) - 1)))\Gamma (a) \left/ {\vrule 
height0.80em width0em depth0.80em} \right. \!  \! 
\Gamma ({\displaystyle \frac {1}{2}}  + n) }
}
\end{maplelatex}

\begin{maplelatex}
\mapleinline{inert}{2d}{455, ":   ",
(-1/3*(-3^(1/2)*cos(Pi*(-3*a+n))+3^(1/2)*sin(1/6*Pi*(-1+6*n-6*a))-3*si
n(1/3*Pi*(-3*a+3*n-1)))*V1(a,n)*GAMMA(2*a+1-n)*GAMMA(1/2+n-a)*GAMMA(1/
6+n-a)*GAMMA(a+2/3)/Pi/sin(Pi*a)/GAMMA(3*a+1-n)+1/6*(3+2*3^(1/2)*sin(1
/3*Pi*(6*a-1)))*V1(-1/6+n-a,n)*GAMMA(1/3+n-2*a)*GAMMA(1/2+n-a)*GAMMA(1
/6+n-a)*GAMMA(a+2/3)/sin(Pi*a)/Pi^(1/2)/GAMMA(-3*a+2*n+1/2)/GAMMA(2*a+
1/3-n))*GAMMA(2*a+5/6-n)/GAMMA(1/2+n)/GAMMA(2/3);}{%
\maplemultiline{
455\mbox{:~~~} ( - {\displaystyle \frac {1}{3}} 
( - \sqrt{3}\,\mathrm{cos}(\pi \,( - 3\,a + n)) + \sqrt{3}\,
\mathrm{sin}({\displaystyle \frac {\pi \,( - 1 + 6\,n - 6\,a)}{6}
} ) - 3\,\mathrm{sin}({\displaystyle \frac {\pi \,( - 3\,a + 3\,n
 - 1)}{3}} )) \\
\mathrm{V1}(a, \,n)\,\Gamma (2\,a + 1 - n)\,\Gamma (
{\displaystyle \frac {1}{2}}  + n - a)\,\Gamma ({\displaystyle 
\frac {1}{6}}  + n - a)\,\Gamma (a + {\displaystyle \frac {2}{3}
} )/(\pi \,\mathrm{sin}(\pi \,a) \\
\Gamma (3\,a + 1 - n))\mbox{} + {\displaystyle \frac {1}{6}} (3
 + 2\,\sqrt{3}\,\mathrm{sin}({\displaystyle \frac {\pi \,(6\,a - 
1)}{3}} ))\,\mathrm{V1}( - {\displaystyle \frac {1}{6}}  + n - a
, \,n)\,\Gamma ({\displaystyle \frac {1}{3}}  + n - 2\,a) \\
\Gamma ({\displaystyle \frac {1}{2}}  + n - a)\,\Gamma (
{\displaystyle \frac {1}{6}}  + n - a)\,\Gamma (a + 
{\displaystyle \frac {2}{3}} ) \left/ {\vrule 
height0.80em width0em depth0.80em} \right. \!  \! (\mathrm{sin}(
\pi \,a)\,\sqrt{\pi }\,\Gamma ( - 3\,a + 2\,n + {\displaystyle 
\frac {1}{2}} ) \\
\Gamma (2\,a + {\displaystyle \frac {1}{3}}  - n)))\Gamma (2\,a
 + {\displaystyle \frac {5}{6}}  - n) \left/ {\vrule 
height0.80em width0em depth0.80em} \right. \!  \! (\Gamma (
{\displaystyle \frac {1}{2}}  + n)\,\Gamma ({\displaystyle 
\frac {2}{3}} )) }
}
\end{maplelatex}

\begin{maplelatex}
\mapleinline{inert}{2d}{456, ":   ",
((2*sin(1/6*Pi*(-24*a+12*n-1))+1)*V1(a,n)*GAMMA(2*a+1-n)*GAMMA(1/3+n-a
)*GAMMA(2*a+2/3-n)*GAMMA(a+2/3)/Pi/GAMMA(3*a+1-n)/GAMMA(5/6+a-n)/(2*si
n(1/6*Pi*(-12*a+12*n+1))+1)-1/2*V1(-1/6+n-a,n)*GAMMA(3*a-2*n+1/2)*GAMM
A(1/3+n-a)*GAMMA(1/6+n-a)*GAMMA(a+2/3)*(-sin(1/6*Pi*(-12*a+12*n+1))+si
n(1/6*Pi*(-24*a+12*n-1)))/cos(Pi*(n-a))/Pi^(3/2)/GAMMA(2*a+1/3-n))/GAM
MA(1/2+n)*GAMMA(2/3);}{%
\maplemultiline{
456\mbox{:~~~} ((2\,\mathrm{sin}({\displaystyle \frac {
\pi \,( - 24\,a + 12\,n - 1)}{6}} ) + 1)\,\mathrm{V1}(a, \,n)\,
\Gamma (2\,a + 1 - n)\,\Gamma ({\displaystyle \frac {1}{3}}  + n
 - a) \\
\Gamma (2\,a + {\displaystyle \frac {2}{3}}  - n)\,\Gamma (a + 
{\displaystyle \frac {2}{3}} ) \left/ {\vrule 
height0.80em width0em depth0.80em} \right. \!  \! (\pi \,\Gamma (
3\,a + 1 - n)\,\Gamma ({\displaystyle \frac {5}{6}}  + a - n) \\
(2\,\mathrm{sin}({\displaystyle \frac {\pi \,( - 12\,a + 12\,n + 
1)}{6}} ) + 1))\mbox{} - {\displaystyle \frac {1}{2}} \mathrm{V1}
( - {\displaystyle \frac {1}{6}}  + n - a, \,n)\,\Gamma (3\,a - 2
\,n + {\displaystyle \frac {1}{2}} ) \\
\Gamma ({\displaystyle \frac {1}{3}}  + n - a)\,\Gamma (
{\displaystyle \frac {1}{6}}  + n - a)\,\Gamma (a + 
{\displaystyle \frac {2}{3}} ) \\
( - \mathrm{sin}({\displaystyle \frac {\pi \,( - 12\,a + 12\,n + 
1)}{6}} ) + \mathrm{sin}({\displaystyle \frac {\pi \,( - 24\,a + 
12\,n - 1)}{6}} )) \left/ {\vrule 
height0.80em width0em depth0.80em} \right. \!  \! (\mathrm{cos}(
\pi \,(n - a)) \\
\pi ^{(3/2)}\,\Gamma (2\,a + {\displaystyle \frac {1}{3}}  - n)))
\Gamma ({\displaystyle \frac {2}{3}} ) \left/ {\vrule 
height0.80em width0em depth0.80em} \right. \!  \! \Gamma (
{\displaystyle \frac {1}{2}}  + n) }
}
\end{maplelatex}

\begin{maplelatex}
\mapleinline{inert}{2d}{457, ":   ",
((2*sin(1/6*Pi*(-24*a+12*n-1))+1)*V1(a,n)*GAMMA(2*a+1-n)*GAMMA(1/3+n-a
)*GAMMA(2*a+2/3-n)/GAMMA(3*a+1-n)/GAMMA(5/6+a-n)/GAMMA(a-n+1/2)/(sin(1
/3*Pi*(3*a+1))+sin(Pi*(-3*a+2*n))+sin(1/3*Pi*(-3*a+6*n+1))-sin(Pi*a))-
1/2*V1(-1/6+n-a,n)*GAMMA(1/3+n-a)*GAMMA(3*a-2*n+1/2)*GAMMA(1/2+n-a)*GA
MMA(1/6+n-a)*(sin(1/6*Pi*(-12*a+12*n+1))-sin(1/6*Pi*(-24*a+12*n-1)))/s
in(1/3*Pi*(3*a-1))/Pi^(3/2)/GAMMA(2*a+1/3-n))*GAMMA(2*a+5/6-n)/GAMMA(1
/2+n);}{%
\maplemultiline{
457\mbox{:~~~} ((2\,\mathrm{sin}({\displaystyle \frac {
\pi \,( - 24\,a + 12\,n - 1)}{6}} ) + 1)\,\mathrm{V1}(a, \,n)\,
\Gamma (2\,a + 1 - n)\,\Gamma ({\displaystyle \frac {1}{3}}  + n
 - a) \\
\Gamma (2\,a + {\displaystyle \frac {2}{3}}  - n) \left/ {\vrule 
height0.80em width0em depth0.80em} \right. \!  \! (\Gamma (3\,a
 + 1 - n)\,\Gamma ({\displaystyle \frac {5}{6}}  + a - n)\,\Gamma
 (a - n + {\displaystyle \frac {1}{2}} ) \\
(\mathrm{sin}({\displaystyle \frac {\pi \,(3\,a + 1)}{3}} ) + 
\mathrm{sin}(\pi \,( - 3\,a + 2\,n)) + \mathrm{sin}(
{\displaystyle \frac {\pi \,( - 3\,a + 6\,n + 1)}{3}} ) - 
\mathrm{sin}(\pi \,a)))\mbox{} - {\displaystyle \frac {1}{2}} 
 \\
\mathrm{V1}( - {\displaystyle \frac {1}{6}}  + n - a, \,n)\,
\Gamma ({\displaystyle \frac {1}{3}}  + n - a)\,\Gamma (3\,a - 2
\,n + {\displaystyle \frac {1}{2}} )\,\Gamma ({\displaystyle 
\frac {1}{2}}  + n - a)\,\Gamma ({\displaystyle \frac {1}{6}}  + 
n - a) \\
(\mathrm{sin}({\displaystyle \frac {\pi \,( - 12\,a + 12\,n + 1)
}{6}} ) - \mathrm{sin}({\displaystyle \frac {\pi \,( - 24\,a + 12
\,n - 1)}{6}} )) \left/ {\vrule height0.80em width0em depth0.80em
} \right. \!  \! (\mathrm{sin}({\displaystyle \frac {\pi \,(3\,a
 - 1)}{3}} ) \\
\pi ^{(3/2)}\,\Gamma (2\,a + {\displaystyle \frac {1}{3}}  - n)))
\Gamma (2\,a + {\displaystyle \frac {5}{6}}  - n) \left/ {\vrule 
height0.80em width0em depth0.80em} \right. \!  \! \Gamma (
{\displaystyle \frac {1}{2}}  + n) }
}
\end{maplelatex}

\begin{maplelatex}
\mapleinline{inert}{2d}{458, ":   ",
((2*sin(1/6*Pi*(-24*a+12*n-1))+1)*V1(a,n)*GAMMA(2*a+1-n)*GAMMA(1/3+n-a
)*GAMMA(2*a+1/3-n)*GAMMA(2*a+2/3-n)*GAMMA(a+2/3)/GAMMA(1/6+a)/GAMMA(3*
a+1-n)/GAMMA(5/6+a-n)/GAMMA(a-n+1/2)/(cos(Pi*a)+cos(Pi*(-3*a+2*n))+sin
(1/6*Pi*(6*a-1))+sin(1/6*Pi*(-6*a+12*n-1)))+1/Pi^(3/2)*V1(-1/6+n-a,n)*
cos(Pi*(-3*a+2*n))/GAMMA(1/6+a)*GAMMA(3*a-2*n+1/2)*GAMMA(1/2+n-a)*GAMM
A(1/3+n-a)*GAMMA(1/6+n-a)*GAMMA(a+2/3))/GAMMA(1/2+n);}{%
\maplemultiline{
458\mbox{:~~~} ((2\,\mathrm{sin}({\displaystyle \frac {
\pi \,( - 24\,a + 12\,n - 1)}{6}} ) + 1)\,\mathrm{V1}(a, \,n)\,
\Gamma (2\,a + 1 - n)\,\Gamma ({\displaystyle \frac {1}{3}}  + n
 - a) \\
\Gamma (2\,a + {\displaystyle \frac {1}{3}}  - n)\,\Gamma (2\,a
 + {\displaystyle \frac {2}{3}}  - n)\,\Gamma (a + 
{\displaystyle \frac {2}{3}} ) \left/ {\vrule 
height0.80em width0em depth0.80em} \right. \!  \! (\Gamma (
{\displaystyle \frac {1}{6}}  + a)\,\Gamma (3\,a + 1 - n)\,\Gamma
 ({\displaystyle \frac {5}{6}}  + a - n) \\
\Gamma (a - n + {\displaystyle \frac {1}{2}} ) 
(\mathrm{cos}(\pi \,a) + \mathrm{cos}(\pi \,( - 3\,a + 2\,n)) + 
\mathrm{sin}({\displaystyle \frac {\pi \,(6\,a - 1)}{6}} ) + 
\mathrm{sin}({\displaystyle \frac {\pi \,( - 6\,a + 12\,n - 1)}{6
}} )))\mbox{} +  \\
\mathrm{V1}( - {\displaystyle \frac {1}{6}}  + n - a, \,n)\,
\mathrm{cos}(\pi \,( - 3\,a + 2\,n))\,\Gamma (3\,a - 2\,n + 
{\displaystyle \frac {1}{2}} )\,\Gamma ({\displaystyle \frac {1}{
2}}  + n - a)\,\Gamma ({\displaystyle \frac {1}{3}}  + n - a) \\
\Gamma ({\displaystyle \frac {1}{6}}  + n - a)\,\Gamma (a + 
{\displaystyle \frac {2}{3}} ) \left/ {\vrule 
height0.80em width0em depth0.80em} \right. \!  \! (\pi ^{(3/2)}\,
\Gamma ({\displaystyle \frac {1}{6}}  + a))) \left/ {\vrule 
height0.80em width0em depth0.80em} \right. \!  \! \Gamma (
{\displaystyle \frac {1}{2}}  + n) }
}
\end{maplelatex}

\begin{maplelatex}
\mapleinline{inert}{2d}{459, ":   ",
(2/3*3^(1/2)*(2*sin(1/6*Pi*(-24*a+12*n-1))+1)*V1(a,n)*GAMMA(2*a+1-n)*G
AMMA(1/3+n-a)*GAMMA(2*a+2/3-n)*GAMMA(a+2/3)/GAMMA(3*a+1-n)/GAMMA(a-n+1
/2)/(2*sin(1/6*Pi*(-12*a+12*n+1))+1)-1/3*V1(-1/6+n-a,n)*3^(1/2)*GAMMA(
3*a-2*n+1/2)*GAMMA(1/2+n-a)*GAMMA(1/3+n-a)*GAMMA(a+2/3)*(-sin(1/6*Pi*(
-12*a+12*n+1))+sin(1/6*Pi*(-24*a+12*n-1)))/sin(1/6*Pi*(-6*a+6*n+1))/Pi
^(1/2)/GAMMA(2*a+1/3-n))/GAMMA(1/2+n)/GAMMA(2/3);}{%
\maplemultiline{
459\mbox{:~~~} ({\displaystyle \frac {2}{3}} \sqrt{3}\,(
2\,\mathrm{sin}({\displaystyle \frac {\pi \,( - 24\,a + 12\,n - 1
)}{6}} ) + 1)\,\mathrm{V1}(a, \,n)\,\Gamma (2\,a + 1 - n)\,\Gamma
 ({\displaystyle \frac {1}{3}}  + n - a) \\
\Gamma (2\,a + {\displaystyle \frac {2}{3}}  - n)\,\Gamma (a + 
{\displaystyle \frac {2}{3}} ) \left/ {\vrule 
height0.80em width0em depth0.80em} \right. \!  \! (\Gamma (3\,a
 + 1 - n)\,\Gamma (a - n + {\displaystyle \frac {1}{2}} ) \\
(2\,\mathrm{sin}({\displaystyle \frac {\pi \,( - 12\,a + 12\,n + 
1)}{6}} ) + 1))\mbox{} - {\displaystyle \frac {1}{3}} \mathrm{V1}
( - {\displaystyle \frac {1}{6}}  + n - a, \,n)\,\sqrt{3}\,\Gamma
 (3\,a - 2\,n + {\displaystyle \frac {1}{2}} ) \\
\Gamma ({\displaystyle \frac {1}{2}}  + n - a)\,\Gamma (
{\displaystyle \frac {1}{3}}  + n - a)\,\Gamma (a + 
{\displaystyle \frac {2}{3}} ) \\
( - \mathrm{sin}({\displaystyle \frac {\pi \,( - 12\,a + 12\,n + 
1)}{6}} ) + \mathrm{sin}({\displaystyle \frac {\pi \,( - 24\,a + 
12\,n - 1)}{6}} )) \left/ {\vrule 
height0.80em width0em depth0.80em} \right. \!  \! \\ (
\mathrm{sin}({\displaystyle \frac {\pi \,( - 6\,a + 6\,n + 1)}{6}
} )\,\sqrt{\pi }\,\Gamma (2\,a + {\displaystyle \frac {1}{3}}  - 
n))) \left/ {\vrule height0.80em width0em depth0.80em}
 \right. \!  \! (\Gamma ({\displaystyle \frac {1}{2}}  + n)\,
\Gamma ({\displaystyle \frac {2}{3}} )) }
}
\end{maplelatex}

\begin{maplelatex}
\mapleinline{inert}{2d}{460, ":   ",
(V1(a,n)*(2*sin(1/6*Pi*(-24*a+12*n-1))+1)*GAMMA(2*a+2/3-n)*GAMMA(a+2/3
)/GAMMA(3*a+1-n)/GAMMA(5/6+a-n)/GAMMA(a-n+1/2)/(sin(1/3*Pi*(1-6*a+3*n)
)+sin(1/3*Pi*(-1-12*a+9*n))+sin(1/3*Pi*(-6*a-1+3*n))-sin(1/3*Pi*(3*n-1
)))-1/2*V1(-1/6+n-a,n)*GAMMA(3*a-2*n+1/2)*GAMMA(n-2*a)*GAMMA(1/2+n-a)*
GAMMA(1/6+n-a)*GAMMA(a+2/3)*(-sin(1/6*Pi*(-12*a+12*n+1))+sin(1/6*Pi*(-
24*a+12*n-1)))/Pi^(5/2)/GAMMA(2*a+1/3-n))*GAMMA(2*a+5/6-n)/GAMMA(1/2+n
)*GAMMA(2/3);}{%
\maplemultiline{
460\mbox{:~~~} (\mathrm{V1}(a, \,n)\,(2\,\mathrm{sin}(
{\displaystyle \frac {\pi \,( - 24\,a + 12\,n - 1)}{6}} ) + 1)\,
\Gamma (2\,a + {\displaystyle \frac {2}{3}}  - n)\,\Gamma (a + 
{\displaystyle \frac {2}{3}} ) \left/ {\vrule 
height0.80em width0em depth0.80em} \right. \!  \! ( \\
\Gamma (3\,a + 1 - n)\,\Gamma ({\displaystyle \frac {5}{6}}  + a
 - n)\,\Gamma (a - n + {\displaystyle \frac {1}{2}} )(\mathrm{sin
}({\displaystyle \frac {\pi \,(1 - 6\,a + 3\,n)}{3}} ) \\
\mbox{} + \mathrm{sin}({\displaystyle \frac {\pi \,( - 1 - 12\,a
 + 9\,n)}{3}} ) + \mathrm{sin}({\displaystyle \frac {\pi \,( - 6
\,a - 1 + 3\,n)}{3}} ) - \mathrm{sin}({\displaystyle \frac {\pi 
\,(3\,n - 1)}{3}} )))\mbox{} - {\displaystyle \frac {1}{2}}  \\
\mathrm{V1}( - {\displaystyle \frac {1}{6}}  + n - a, \,n)\,
\Gamma (3\,a - 2\,n + {\displaystyle \frac {1}{2}} )\,\Gamma (n
 - 2\,a)\,\Gamma ({\displaystyle \frac {1}{2}}  + n - a)\,\Gamma 
({\displaystyle \frac {1}{6}}  + n - a)\,\Gamma (a + 
{\displaystyle \frac {2}{3}} ) \\
( - \mathrm{sin}({\displaystyle \frac {\pi \,( - 12\,a + 12\,n + 
1)}{6}} ) + \mathrm{sin}({\displaystyle \frac {\pi \,( - 24\,a + 
12\,n - 1)}{6}} )) \left/ {\vrule 
height0.80em width0em depth0.80em} \right. \!  \! (\pi ^{(5/2)}
 \\
\Gamma (2\,a + {\displaystyle \frac {1}{3}}  - n)))\Gamma (2\,a
 + {\displaystyle \frac {5}{6}}  - n)\,\Gamma ({\displaystyle 
\frac {2}{3}} ) \left/ {\vrule height0.80em width0em depth0.80em}
 \right. \!  \! \Gamma ({\displaystyle \frac {1}{2}}  + n) }
}
\end{maplelatex}

\begin{maplelatex}
\mapleinline{inert}{2d}{461, ":   ",
-1/2*(2*cos(Pi*(-3*a+n))+cos(3*Pi*(n-a)))*V1(a,n)*GAMMA(5/3+a-n)*GAMMA
(2*a+1/3-n)*GAMMA(-1/6+n-a)*GAMMA(1/2+n-a)*GAMMA(a+1/3)*GAMMA(a+2/3)*G
AMMA(a)/(2*n-1)/((-1)^n)/sin(1/6*Pi*(-6*a+6*n+1))/Pi^(3/2)/GAMMA(-1/2+
n)^2/GAMMA(5/6+a-n)/GAMMA(3*a+1-n)+1/4*V1(-1/6+n-a,n)*GAMMA(3*a-2*n+1/
2)*GAMMA(5/3+a-n)*GAMMA(-1/6+n-a)*GAMMA(1/2+n-a)*GAMMA(a+1/3)*GAMMA(a+
2/3)*GAMMA(a)*(cos(Pi*(-3*a+n))+cos(3*Pi*(n-a)))*(2*n-1)/sin(1/6*Pi*(-
6*a+6*n+1))/((-1)^n)/Pi/GAMMA(1/2+n)^2/GAMMA(2*a+2/3-n)/GAMMA(5/6+a-n)
/GAMMA(2*a+1-n);}{%
\maplemultiline{
461\mbox{:~~~}  - {\displaystyle \frac {1}{2}} (2\,
\mathrm{cos}(\pi \,( - 3\,a + n)) + \mathrm{cos}(3\,\pi \,(n - a)
))\,\mathrm{V1}(a, \,n)\,\Gamma ({\displaystyle \frac {5}{3}}  + 
a - n)\,\Gamma (2\,a + {\displaystyle \frac {1}{3}}  - n) \\
\Gamma ( - {\displaystyle \frac {1}{6}}  + n - a)\,\Gamma (
{\displaystyle \frac {1}{2}}  + n - a)\,\Gamma (a + 
{\displaystyle \frac {1}{3}} )\,\Gamma (a + {\displaystyle 
\frac {2}{3}} )\,\Gamma (a) \left/ {\vrule 
height0.80em width0em depth0.80em} \right. \!  \! ((2\,n - 1)\,(
-1)^{n} \\
\mathrm{sin}({\displaystyle \frac {\pi \,( - 6\,a + 6\,n + 1)}{6}
} )\,\pi ^{(3/2)}\,\Gamma ( - {\displaystyle \frac {1}{2}}  + n)
^{2}\,\Gamma ({\displaystyle \frac {5}{6}}  + a - n)\,\Gamma (3\,
a + 1 - n))\mbox{} + {\displaystyle \frac {1}{4}}  \\
\mathrm{V1}( - {\displaystyle \frac {1}{6}}  + n - a, \,n)\,
\Gamma (3\,a - 2\,n + {\displaystyle \frac {1}{2}} )\,\Gamma (
{\displaystyle \frac {5}{3}}  + a - n)\,\Gamma ( - 
{\displaystyle \frac {1}{6}}  + n - a)\,\Gamma ({\displaystyle 
\frac {1}{2}}  + n - a) \\
\Gamma (a + {\displaystyle \frac {1}{3}} )\,\Gamma (a + 
{\displaystyle \frac {2}{3}} )\,\Gamma (a)\,(\mathrm{cos}(\pi \,(
 - 3\,a + n)) + \mathrm{cos}(3\,\pi \,(n - a)))\,(2\,n - 1)
 \left/ {\vrule height0.80em width0em depth0.80em} \right. \! 
 \! \\ ( 
\mathrm{sin}({\displaystyle \frac {\pi \,( - 6\,a + 6\,n + 1)}{6}
} )\,(-1)^{n}\,\pi \,\Gamma ({\displaystyle \frac {1}{2}}  + n)^{
2}\,\Gamma (2\,a + {\displaystyle \frac {2}{3}}  - n)\,\Gamma (
{\displaystyle \frac {5}{6}}  + a - n)\,\Gamma (2\,a + 1 - n) )
 }
}
\end{maplelatex}

\begin{maplelatex}
\mapleinline{inert}{2d}{462, ":   ",
(1/2*(2*cos(Pi*(-3*a+n))+cos(3*Pi*(n-a)))*V1(a,n)*GAMMA(3*a-2*n+3/2)*G
AMMA(5/3+a-n)*GAMMA(a+1/3)*GAMMA(-1/6+n-a)*GAMMA(1/2+n-a)/Pi^(3/2)/GAM
MA(5/6+a-n)/GAMMA(3*a+1-n)/(sin(1/3*Pi*(1+3*n))+sin(Pi*(n-2*a)))+1/2*V
1(-1/6+n-a,n)*GAMMA(5/3+a-n)*GAMMA(3*a-2*n+1/2)^2*GAMMA(-1/6+n-a)*GAMM
A(1/2+n-a)*GAMMA(a+1/3)*(cos(Pi*(-3*a+n))+cos(3*Pi*(n-a)))*(-6*a+4*n-1
)/Pi/GAMMA(2*a+1-n)/GAMMA(2*a+2/3-n)/GAMMA(5/6+a-n)/GAMMA(2*a+1/3-n)/(
sin(1/3*Pi*(1+3*n))+sin(Pi*(n-2*a))))*GAMMA(a)/GAMMA(1/2+n);}{%
\maplemultiline{
462\mbox{:~~~} ({\displaystyle \frac {1}{2}} (2\,
\mathrm{cos}(\pi \,( - 3\,a + n)) + \mathrm{cos}(3\,\pi \,(n - a)
))\,\mathrm{V1}(a, \,n)\,\Gamma (3\,a - 2\,n + {\displaystyle 
\frac {3}{2}} ) \\
\Gamma ({\displaystyle \frac {5}{3}}  + a - n)\,\Gamma (a + 
{\displaystyle \frac {1}{3}} )\,\Gamma ( - {\displaystyle \frac {
1}{6}}  + n - a)\,\Gamma ({\displaystyle \frac {1}{2}}  + n - a)
 \left/ {\vrule height0.80em width0em depth0.80em} \right. \! 
 \! (\pi ^{(3/2)}\,\Gamma ({\displaystyle \frac {5}{6}}  + a - n)
 \\
\Gamma (3\,a + 1 - n)\,(\mathrm{sin}({\displaystyle \frac {\pi \,
(1 + 3\,n)}{3}} ) + \mathrm{sin}(\pi \,(n - 2\,a))))\mbox{} + 
{\displaystyle \frac {1}{2}} \mathrm{V1}( - {\displaystyle 
\frac {1}{6}}  + n - a, \,n) \\
\Gamma ({\displaystyle \frac {5}{3}}  + a - n)\,\Gamma (3\,a - 2
\,n + {\displaystyle \frac {1}{2}} )^{2}\,\Gamma ( - 
{\displaystyle \frac {1}{6}}  + n - a)\,\Gamma ({\displaystyle 
\frac {1}{2}}  + n - a)\,\Gamma (a + {\displaystyle \frac {1}{3}
} ) \\
(\mathrm{cos}(\pi \,( - 3\,a + n)) + \mathrm{cos}(3\,\pi \,(n - a
)))\,( - 6\,a + 4\,n - 1) \left/ {\vrule 
height0.80em width0em depth0.80em} \right. \!  \! (\pi \,\Gamma (
2\,a + 1 - n) \\
\Gamma (2\,a + {\displaystyle \frac {2}{3}}  - n)\,\Gamma (
{\displaystyle \frac {5}{6}}  + a - n)\,\Gamma (2\,a + 
{\displaystyle \frac {1}{3}}  - n)\,(\mathrm{sin}({\displaystyle 
\frac {\pi \,(1 + 3\,n)}{3}} ) + \mathrm{sin}(\pi \,(n - 2\,a))))
) \\
\Gamma (a) \left/ {\vrule height0.80em width0em depth0.80em}
 \right. \!  \! \Gamma ({\displaystyle \frac {1}{2}}  + n) }
}
\end{maplelatex}

\begin{maplelatex}
\mapleinline{inert}{2d}{463, ":   ",
(V1(a,n)*(2*cos(Pi*(-3*a+n))+cos(3*Pi*(n-a)))*GAMMA(4/3+a-n)*GAMMA(a+1
/3)*GAMMA(a+2/3)*GAMMA(5/3+a-n)*GAMMA(-1/6+n-a)/Pi^(3/2)/GAMMA(5/6+a-n
)/GAMMA(3*a+1-n)/(2*sin(1/6*Pi*(-12*a+12*n+1))+1)-2*V1(-1/6+n-a,n)*GAM
MA(4/3+a-n)*GAMMA(3*a-2*n+1/2)*GAMMA(5/3+a-n)*GAMMA(-1/6+n-a)*GAMMA(a+
1/3)*GAMMA(a+2/3)*(cos(Pi*(-3*a+n))+cos(3*Pi*(n-a)))/Pi/GAMMA(2*a+2/3-
n)/GAMMA(5/6+a-n)/GAMMA(2*a+1-n)/GAMMA(2*a+1/3-n)/(2*sin(1/6*Pi*(-12*a
+12*n+1))+1))*GAMMA(a)/GAMMA(1/2+n);}{%
\maplemultiline{
463\mbox{:~~~} (\mathrm{V1}(a, \,n)\,(2\,\mathrm{cos}(
\pi \,( - 3\,a + n)) + \mathrm{cos}(3\,\pi \,(n - a)))\,\Gamma (
{\displaystyle \frac {4}{3}}  + a - n)\,\Gamma (a + 
{\displaystyle \frac {1}{3}} ) \\
\Gamma (a + {\displaystyle \frac {2}{3}} )\,\Gamma (
{\displaystyle \frac {5}{3}}  + a - n)\,\Gamma ( - 
{\displaystyle \frac {1}{6}}  + n - a) \left/ {\vrule 
height0.80em width0em depth0.80em} \right. \!  \! (\pi ^{(3/2)}\,
\Gamma ({\displaystyle \frac {5}{6}}  + a - n)\,\Gamma (3\,a + 1
 - n) \\
(2\,\mathrm{sin}({\displaystyle \frac {\pi \,( - 12\,a + 12\,n + 
1)}{6}} ) + 1))\mbox{} - 2\,\mathrm{V1}( - {\displaystyle \frac {
1}{6}}  + n - a, \,n)\,\Gamma ({\displaystyle \frac {4}{3}}  + a
 - n) \\
\Gamma (3\,a - 2\,n + {\displaystyle \frac {1}{2}} )\,\Gamma (
{\displaystyle \frac {5}{3}}  + a - n)\,\Gamma ( - 
{\displaystyle \frac {1}{6}}  + n - a)\,\Gamma (a + 
{\displaystyle \frac {1}{3}} )\,\Gamma (a + {\displaystyle 
\frac {2}{3}} ) \\
(\mathrm{cos}(\pi \,( - 3\,a + n)) + \mathrm{cos}(3\,\pi \,(n - a
))) \left/ {\vrule height0.80em width0em depth0.80em} \right. \! 
 \! (\pi \,\Gamma (2\,a + {\displaystyle \frac {2}{3}}  - n)\,
\Gamma ({\displaystyle \frac {5}{6}}  + a - n) \\
\Gamma (2\,a + 1 - n)\,\Gamma (2\,a + {\displaystyle \frac {1}{3}
}  - n)\,(2\,\mathrm{sin}({\displaystyle \frac {\pi \,( - 12\,a
 + 12\,n + 1)}{6}} ) + 1)))\Gamma (a) \left/ {\vrule 
height0.80em width0em depth0.80em} \right. \!  \! 
\Gamma ({\displaystyle \frac {1}{2}}  + n) }
}
\end{maplelatex}

\begin{maplelatex}
\mapleinline{inert}{2d}{469, ":   ",
sin(Pi*(a+a^2+a*c+a*b+b*c)/a)+GAMMA(-(a*c-a+a*b+b*c+a^2)/a)+GAMMA((a+a
*c+a*b+b*c)/a)+GAMMA(1+b+a+c)+GAMMA(-c*(a+b)/a)+GAMMA(-b*(a+c)/a)+a+1/
GAMMA(c)+1/GAMMA(b)+1/GAMMA(-a+1)+1/sin(Pi*a)+1/GAMMA(2+a)+1/GAMMA(-(b
*c+a*c-a)/a)+1/GAMMA(-(a*b-a+b*c)/a);}{%
\maplemultiline{
469\mbox{:~~~} \mathrm{sin}({\displaystyle \frac {\pi \,
(a + a^{2} + a\,c + a\,b + b\,c)}{a}} ) + \Gamma ( - 
{\displaystyle \frac {a\,c - a + a\,b + b\,c + a^{2}}{a}} ) \\
\mbox{} + \Gamma ({\displaystyle \frac {a + a\,c + a\,b + b\,c}{a
}} ) + \Gamma (1 + b + a + c) + \Gamma ( - {\displaystyle \frac {
c\,(a + b)}{a}} ) + \Gamma ( - {\displaystyle \frac {b\,(a + c)}{
a}} ) + a \\
\mbox{} + {\displaystyle \frac {1}{\Gamma (c)}}  + 
{\displaystyle \frac {1}{\Gamma (b)}}  + {\displaystyle \frac {1
}{\Gamma ( - a + 1)}}  + {\displaystyle \frac {1}{\mathrm{sin}(
\pi \,a)}}  + {\displaystyle \frac {1}{\Gamma (2 + a)}}  + 
{\displaystyle \frac {1}{\Gamma ( - {\displaystyle \frac {b\,c + 
a\,c - a}{a}} )}} 
\mbox{} + {\displaystyle \frac {1}{\Gamma ( - {\displaystyle 
\frac {a\,b - a + b\,c}{a}} )}}  }
}
\end{maplelatex}

\end{maplegroup}

%% file: AppendixC.tex
\begin{maplegroup}

\mapleresult
\begin{maplelatex}
\mapleinline{inert}{2d}{`  1: `, `Watsons Theorem from Chu, Math. Comp., Theorem 5, m
symbolic`;}{%
\[
\mathit{\ \ 1:\ }  \,\mbox{Watsons\ Theorem\ from\ Chu,\ Math.
\ Comp.,\ Theorem\ 5,\ m\ symbolic}
\]
}
\end{maplelatex}

\begin{maplelatex}
\mapleinline{inert}{2d}{`  2: `, `Watsons Theorem from Chu, Math. Comp., Theorem 5, m
even`;}{%
\[
\mathit{\ \ 2:\ }  \,\mbox{Watsons\ Theorem\ from\ Chu,\ Math.
\ Comp.,\ Theorem\ 5,\ m\ even}
\]
}
\end{maplelatex}

\begin{maplelatex}
\mapleinline{inert}{2d}{`  3: `, `Watsons Theorem from Chu, Math. Comp., Theorem 5, m
odd`;}{%
\[
\mathit{\ \ 3:\ }  \,\mbox{Watsons\ Theorem\ from\ Chu,\ Math.
\ Comp.,\ Theorem\ 5,\ m\ odd}
\]
}
\end{maplelatex}

\begin{maplelatex}
\mapleinline{inert}{2d}{`  4: `, `Watsons Theorem from Chu, Math. Comp., Theorem 6, m
symbolic`;}{%
\[
\mathit{\ \ 4:\ }  \,\mbox{Watsons\ Theorem\ from\ Chu,\ Math.
\ Comp.,\ Theorem\ 6,\ m\ symbolic}
\]
}
\end{maplelatex}

\begin{maplelatex}
\mapleinline{inert}{2d}{`  5: `, `Watsons Theorem from Chu, Math. Comp., Theorem 6, m
even`;}{%
\[
\mathit{\ \ 5:\ }  \,\mbox{Watsons\ Theorem\ from\ Chu,\ Math.
\ Comp.,\ Theorem\ 6,\ m\ even}
\]
}
\end{maplelatex}

\begin{maplelatex}
\mapleinline{inert}{2d}{`  6: `, `Watsons Theorem from Chu, Math. Comp., Theorem 6, m
odd`;}{%
\[
\mathit{\ \ 6:\ }  \,\mbox{Watsons\ Theorem\ from\ Chu,\ Math.
\ Comp.,\ Theorem\ 6,\ m\ odd}
\]
}
\end{maplelatex}

\begin{maplelatex}
\mapleinline{inert}{2d}{`  7: `, `Watsons Theorem from Chu, Math. Comp., Theorem 7, m
symbolic`;}{%
\[
\mathit{\ \ 7:\ }  \,\mbox{Watsons\ Theorem\ from\ Chu,\ Math.
\ Comp.,\ Theorem\ 7,\ m\ symbolic}
\]
}
\end{maplelatex}

\begin{maplelatex}
\mapleinline{inert}{2d}{`  8: `, `Watsons Theorem from Chu, Math. Comp., Theorem 7, m
even`;}{%
\[
\mathit{\ \ 8:\ }  \,\mbox{Watsons\ Theorem\ from\ Chu,\ Math.
\ Comp.,\ Theorem\ 7,\ m\ even}
\]
}
\end{maplelatex}

\begin{maplelatex}
\mapleinline{inert}{2d}{`  9: `, `Watsons Theorem from Chu, Math. Comp., Theorem 7, m
odd`;}{%
\[
\mathit{\ \ 9:\ }  \,\mbox{Watsons\ Theorem\ from\ Chu,\ Math.
\ Comp.,\ Theorem\ 7,\ m\ odd}
\]
}
\end{maplelatex}

\begin{maplelatex}
\mapleinline{inert}{2d}{`  10: `, `Watsons Theorem from Chu, Math. Comp., Theorem 8, m
symbolic`;}{%
\[
\mathit{\ \ 10:\ }  \,\mbox{Watsons\ Theorem\ from\ Chu,\ Math.
\ Comp.,\ Theorem\ 8,\ m\ symbolic}
\]
}
\end{maplelatex}

\begin{maplelatex}
\mapleinline{inert}{2d}{`  11: `, `Watsons Theorem from Chu, Math. Comp., Theorem 8, m
even`;}{%
\[
\mathit{\ \ 11:\ }  \,\mbox{Watsons\ Theorem\ from\ Chu,\ Math.
\ Comp.,\ Theorem\ 8,\ m\ even}
\]
}
\end{maplelatex}

\begin{maplelatex}
\mapleinline{inert}{2d}{`  12: `, `Watsons Theorem from Chu, Math. Comp., Theorem 8, m
odd`;}{%
\[
\mathit{\ \ 12:\ }  \,\mbox{Watsons\ Theorem\ from\ Chu,\ Math.
\ Comp.,\ Theorem\ 8,\ m\ odd}
\]
}
\end{maplelatex}

\begin{maplelatex}
\mapleinline{inert}{2d}{`  13: `, `Whipples Theorem from Chu, Math. Comp.; m,n>=0 `;}{%
\[
\mathit{\ \ 13:\ }  \,\mbox{Whipples\ Theorem\ from\ Chu,\ 
Math.\ Comp.;\ ${m,n>=0}$\ }
\]
}
\end{maplelatex}

\begin{maplelatex}
\mapleinline{inert}{2d}{`  14: `, `Whipples Theorem from Chu, Math. Comp.; m>0,n<0 `;}{%
\[
\mathit{\ \ 14:\ }  \,\mbox{Whipples\ Theorem\ from\ Chu,\ 
Math.\ Comp.;\ ${m>0,n<0}$\ }
\]
}
\end{maplelatex}

\begin{maplelatex}
\mapleinline{inert}{2d}{`  15: `, `Whipples Theorem from Chu, Math. Comp.; m<0,n>=0 `;}{%
\[
\mathit{\ \ 15:\ }  \,\mbox{Whipples\ Theorem\ from\ Chu,\ 
Math.\ Comp.;\ ${m<0,n>=0}$\ }
\]
}
\end{maplelatex}

\begin{maplelatex}
\mapleinline{inert}{2d}{`  16: `, `Whipples Theorem from Chu, Math. Comp.; m<0,n<0 `;}{%
\[
\mathit{\ \ 16:\ }  \,\mbox{Whipples\ Theorem\ from\ Chu,\ 
Math.\ Comp.;\ ${m<0,n<0}$\ }
\]
}
\end{maplelatex}

\begin{maplelatex}
\mapleinline{inert}{2d}{`  17: `, `Dixon's Theorem from Chu, Eq. 10, Math. Comp.; m,n>=0
`;}{%
\[
\mathit{\ \ 17:\ }  \,\mbox{Dixon$^{\prime }s$ \ Theorem\ from\ 
Chu,\ Eq.\ 10,\ Math.\ Comp.;\ ${m,n>=0}$ \ }
\]
}
\end{maplelatex}

\begin{maplelatex}
\mapleinline{inert}{2d}{`  18: `, `Dixons Theorem from Chu, Eq. 10, Math. Comp.; m>0,n<0
`;}{%
\[
\mathit{\ \ 18:\ }  \,\mbox{Dixons\ Theorem\ from\ Chu,\ Eq.\ 
10,\ Math.\ Comp.;\ ${m>0,n<0}$ \ }
\]
}
\end{maplelatex}

\begin{maplelatex}
\mapleinline{inert}{2d}{`  19: `, `Dixons Theorem from Chu, Eq. 10, Math. Comp.; m,n<0 `;}{%
\[
\mathit{\ \ 19:\ }  \,\mbox{Dixons\ Theorem\ from\ Chu,\ Eq.\ 
10,\ Math.\ Comp.;\ ${m,n<0}$ \ }
\]
}
\end{maplelatex}

\begin{maplelatex}
\mapleinline{inert}{2d}{`  25: `, "This is Lemma 2.2";}{%
\[
\mathit{\ \ 25:\ }  \,\mbox{``This~is~Lemma~2.2''}
\]
}
\end{maplelatex}

\begin{maplelatex}
\mapleinline{inert}{2d}{`  26: `, "This is Lemma 2.2";}{%
\[
\mathit{\ \ 26:\ }  \,\mbox{``This~is~Lemma~2.2''}
\]
}
\end{maplelatex}

\begin{maplelatex}
\mapleinline{inert}{2d}{`  27: `, "This is Lemma 2.2";}{%
\[
\mathit{\ \ 27:\ }  \,\mbox{``This~is~Lemma~2.2''}
\]
}
\end{maplelatex}

\begin{maplelatex}
\mapleinline{inert}{2d}{`  28: `, "This is Lemma 2.2";}{%
\[
\mathit{\ \ 28:\ }  \,\mbox{``This~is~Lemma~2.2''}
\]
}
\end{maplelatex}

\begin{maplelatex}
\mapleinline{inert}{2d}{`  29: `, "This is Lemma 2.2";}{%
\[
\mathit{\ \ 29:\ }  \,\mbox{``This~is~Lemma~2.2''}
\]
}
\end{maplelatex}

\begin{maplelatex}
\mapleinline{inert}{2d}{`  30: `, "Prudnikov 7.4.4.23";}{%
\[
\mathit{\ \ 30:\ }  \,\mbox{``Prudnikov~7.4.4.23''}
\]
}
\end{maplelatex}

\begin{maplelatex}
\mapleinline{inert}{2d}{`  31: `, "Prudnikov 7.4.4.23";}{%
\[
\mathit{\ \ 31:\ }  \,\mbox{``Prudnikov~7.4.4.23''}
\]
}
\end{maplelatex}

\begin{maplelatex}
\mapleinline{inert}{2d}{`  33: `, "Prudnikov 7.4.4.23";}{%
\[
\mathit{\ \ 33:\ }  \,\mbox{``Prudnikov~7.4.4.23''}
\]
}
\end{maplelatex}

\begin{maplelatex}
\mapleinline{inert}{2d}{`  34: `, "Prudnikov 7.4.4.23";}{%
\[
\mathit{\ \ 34:\ }  \,\mbox{``Prudnikov~7.4.4.23''}
\]
}
\end{maplelatex}

\begin{maplelatex}
\mapleinline{inert}{2d}{`  35: `, "Prudnikov 7.4.4.23";}{%
\[
\mathit{\ \ 35:\ }  \,\mbox{``Prudnikov~7.4.4.23''}
\]
}
\end{maplelatex}

\begin{maplelatex}
\mapleinline{inert}{2d}{`  36: `, "Prudnikov 7.4.4.23";}{%
\[
\mathit{\ \ 36:\ }  \,\mbox{``Prudnikov~7.4.4.23''}
\]
}
\end{maplelatex}

\begin{maplelatex}
\mapleinline{inert}{2d}{`  37: `, "Prudnikov 7.4.4.32";}{%
\[
\mathit{\ \ 37:\ }  \,\mbox{``Prudnikov~7.4.4.32''}
\]
}
\end{maplelatex}

\begin{maplelatex}
\mapleinline{inert}{2d}{`  38: `, "Prudnikov 7.4.4.32";}{%
\[
\mathit{\ \ 38:\ }  \,\mbox{``Prudnikov~7.4.4.32''}
\]
}
\end{maplelatex}

\begin{maplelatex}
\mapleinline{inert}{2d}{`  39: `, "Prudnikov 7.4.4.32";}{%
\[
\mathit{\ \ 39:\ }  \,\mbox{``Prudnikov~7.4.4.32''}
\]
}
\end{maplelatex}

\begin{maplelatex}
\mapleinline{inert}{2d}{`  63: `, "Prudnikov 7.4.4.32, n odd";}{%
\[
\mathit{\ \ 63:\ }  \,\mbox{``Prudnikov~7.4.4.32,~n~odd''}
\]
}
\end{maplelatex}

\begin{maplelatex}
\mapleinline{inert}{2d}{`  64: `, "Prudnikov 7.4.4.32, n odd";}{%
\[
\mathit{\ \ 64:\ }  \,\mbox{``Prudnikov~7.4.4.32,~n~odd''}
\]
}
\end{maplelatex}

\begin{maplelatex}
\mapleinline{inert}{2d}{`  65: `, "Prudnikov 7.4.4.32, n odd";}{%
\[
\mathit{\ \ 65:\ }  \,\mbox{``Prudnikov~7.4.4.32,~n~odd''}
\]
}
\end{maplelatex}

\begin{maplelatex}
\mapleinline{inert}{2d}{`  80: `, "Lewanowicz, J. Comp & Appl. Math. 86,375(1997) Generalized
Watson; Eq.(2.15), n odd, m odd";}{%
\maplemultiline{
\mathit{\ \ 80:\ }  \mbox{``Lewanowicz,~J.~Comp~\&~Appl.~Math.~
86,375(1997)~Generalized~Watson;~  } \\
\mbox{Eq.(2.15),~n~odd,~m~odd''} }
}
\end{maplelatex}

\begin{maplelatex}
\mapleinline{inert}{2d}{`  81: `, "Lewanowicz, J. Comp & Appl. Math. 86,375(1997) Generalized
Watson; Eq.(2.15), m even, n odd";}{%
\maplemultiline{
\mathit{\ \ 81:\ }  \mbox{``Lewanowicz,~J.~Comp~\&~Appl.~Math.~
86,375(1997)~Generalized~Watson;~ } \\
\mbox{Eq.(2.15),~m~even,~n~odd''} }
}
\end{maplelatex}

\begin{maplelatex}
\mapleinline{inert}{2d}{`  82: `, "Lewanowicz, J. Comp & Appl. Math. 86,375(1997) Generalized
Watson; Eq.(2.15), m odd, n even";}{%
\maplemultiline{
\mathit{\ \ 82:\ }  \mbox{``Lewanowicz,~J.~Comp~\&~Appl.~Math.~
86,375(1997)~Generalized~Watson;~ } \\
\mbox{Eq.(2.15),~m~odd,~n~even''} }
}
\end{maplelatex}

\begin{maplelatex}
\mapleinline{inert}{2d}{`  84: `, "Gessel and Stanton, SIAM J. Math. Anal.,13,295(1982)
Eq.(5.16)";}{%
\[
\mathit{\ \ 84:\ }  \,\mbox{``Gessel~and~Stanton,~SIAM~J.~Math.~
Anal.,13,295(1982)~Eq.(5.16)''}
\]
}
\end{maplelatex}

\begin{maplelatex}
\mapleinline{inert}{2d}{`  85: `, "Gessel and Stanton, SIAM J. Math. Anal.,13,295(1982)
Eq.(5.16)";}{%
\[
\mathit{\ \ 85:\ }  \,\mbox{``Gessel~and~Stanton,~SIAM~J.~Math.~
Anal.,13,295(1982)~Eq.(5.16)''}
\]
}
\end{maplelatex}

\begin{maplelatex}
\mapleinline{inert}{2d}{`  86: `, "Gessel and Stanton, SIAM J. Math. Anal.,13,295(1982)
Eq.(5.16)";}{%
\[
\mathit{\ \ 86:\ }  \,\mbox{``Gessel~and~Stanton,~SIAM~J.~Math.~
Anal.,13,295(1982)~Eq.(5.16)''}
\]
}
\end{maplelatex}

\begin{maplelatex}
\mapleinline{inert}{2d}{`  87: `, "Gessel and Stanton, SIAM J. Math. Anal.,13,295(1982)
Eq.(5.16)";}{%
\[
\mathit{\ \ 87:\ }  \,\mbox{``Gessel~and~Stanton,~SIAM~J.~Math.~
Anal.,13,295(1982)~Eq.(5.16)''}
\]
}
\end{maplelatex}

\begin{maplelatex}
\mapleinline{inert}{2d}{`  88: `, "Gessel and Stanton, SIAM J. Math. Anal.,13,295(1982)
Eq.(5.16)";}{%
\[
\mathit{\ \ 88:\ }  \,\mbox{``Gessel~and~Stanton,~SIAM~J.~Math.~
Anal.,13,295(1982)~Eq.(5.16)''}
\]
}
\end{maplelatex}

\begin{maplelatex}
\mapleinline{inert}{2d}{`  89: `, "Gessel and Stanton, SIAM J. Math. Anal.,13,295(1982)
Eq.(5.16)";}{%
\[
\mathit{\ \ 89:\ }  \,\mbox{``Gessel~and~Stanton,~SIAM~J.~Math.~
Anal.,13,295(1982)~Eq.(5.16)''}
\]
}
\end{maplelatex}

\begin{maplelatex}
\mapleinline{inert}{2d}{`  91: `, "CHU Wenshang";}{%
\[
\mathit{\ \ 91:\ }  \,\mbox{``CHU~Wenshang''}
\]
}
\end{maplelatex}

\begin{maplelatex}
\mapleinline{inert}{2d}{`  92: `, "CHU Wenshang";}{%
\[
\mathit{\ \ 92:\ }  \,\mbox{``CHU~Wenshang''}
\]
}
\end{maplelatex}

\begin{maplelatex}
\mapleinline{inert}{2d}{`  104: `, "Maier Eq. (7.1) with L=n";}{%
\[
\mathit{\ \ 104:\ }  \,\mbox{``Maier Eq (7.1)~with~L=n''}
\]
}
\end{maplelatex}

\begin{maplelatex}
\mapleinline{inert}{2d}{`  105: `, "Maier Eq. (7.1) with L=n";}{%
\[
\mathit{\ \ 105:\ }  \,\mbox{``Maier Eq. (7.1)~with~L=n''}
\]
}
\end{maplelatex}

\begin{maplelatex}
\mapleinline{inert}{2d}{`  106: `, "Maier Eq. (7.1) with L=n";}{%
\[
\mathit{\ \ 106:\ }  \,\mbox{``Maier Eq. (7.1)~with~L=n''}
\]
}
\end{maplelatex}

\begin{maplelatex}
\mapleinline{inert}{2d}{`  107: `, "Maier Eq. (7.1) with L=n";}{%
\[
\mathit{\ \ 107:\ }  \,\mbox{``Maier Eq. (7.1)~with~L=n''}
\]
}
\end{maplelatex}

\begin{maplelatex}
\mapleinline{inert}{2d}{`  108: `, "Maier Eq. (7.1) with L=n";}{%
\[
\mathit{\ \ 108:\ }  \,\mbox{``Maier Eq. (7.1)~with~L=n''}
\]
}
\end{maplelatex}

\begin{maplelatex}
\mapleinline{inert}{2d}{`  109: `, "Maier Eq. (7.1) with L=n";}{%
\[
\mathit{\ \ 109:\ }  \,\mbox{``Maier Eq. (7.1)~with~L=n''}
\]
}
\end{maplelatex}

\begin{maplelatex}
\mapleinline{inert}{2d}{`  110: `, "Maier Eq. (7.1) with L=n";}{%
\[
\mathit{\ \ 110:\ }  \,\mbox{``Maier Eq. (7.1)~with~L=n''}
\]
}
\end{maplelatex}

\begin{maplelatex}
\mapleinline{inert}{2d}{`  111: `, "Maier Eq. (7.1) with L=n";}{%
\[
\mathit{\ \ 111:\ }  \,\mbox{``Maier Eq. (7.1)~with~L=n''}
\]
}
\end{maplelatex}

\begin{maplelatex}
\mapleinline{inert}{2d}{`  112: `, "Maier Eq. (7.1) with L=n";}{%
\[
\mathit{\ \ 112:\ }  \,\mbox{``Maier Eq. (7.1)~with~L=n''}
\]
}
\end{maplelatex}

\begin{maplelatex}
\mapleinline{inert}{2d}{`  113: `, "Maier Eq. (7.1) with L=n";}{%
\[
\mathit{\ \ 113:\ }  \,\mbox{``Maier Eq. (7.1)~with~L=n''}
\]
}
\end{maplelatex}

\begin{maplelatex}
\mapleinline{inert}{2d}{`  114: `, "Maier Thm.(7.3) with L=n variation 2";}{%
\[
\mathit{\ \ 114:\ }  \,\mbox{``Maier~Thm.(7.3)~with~L=n~variation
~2''}
\]
}
\end{maplelatex}

\begin{maplelatex}
\mapleinline{inert}{2d}{`  115: `, "Maier Thm.(7.3) with L=n variation 2";}{%
\[
\mathit{\ \ 115:\ }  \,\mbox{``Maier~Thm.(7.3)~with~L=n~variation
~2''}
\]
}
\end{maplelatex}

\begin{maplelatex}
\mapleinline{inert}{2d}{`  116: `, "Maier Thm.(7.3) with L=n variation 2";}{%
\[
\mathit{\ \ 116:\ }  \,\mbox{``Maier~Thm.(7.3)~with~L=n~variation
~2''}
\]
}
\end{maplelatex}

\begin{maplelatex}
\mapleinline{inert}{2d}{`  117: `, "Maier Thm.(7.3) with L=n variation 2";}{%
\[
\mathit{\ \ 117:\ }  \,\mbox{``Maier~Thm.(7.3)~with~L=n~variation
~2''}
\]
}
\end{maplelatex}

\begin{maplelatex}
\mapleinline{inert}{2d}{`  118: `, "Maier Thm.(7.3) with L=n variation 2";}{%
\[
\mathit{\ \ 118:\ }  \,\mbox{``Maier~Thm.(7.3)~with~L=n~variation
~2''}
\]
}
\end{maplelatex}

\begin{maplelatex}
\mapleinline{inert}{2d}{`  119: `, "Maier Thm.(7.3) with L=n variation 2";}{%
\[
\mathit{\ \ 119:\ }  \,\mbox{``Maier~Thm.(7.3)~with~L=n~variation
~2''}
\]
}
\end{maplelatex}

\begin{maplelatex}
\mapleinline{inert}{2d}{`  120: `, "Maier Thm.(7.3) with L=n variation 2";}{%
\[
\mathit{\ \ 120:\ }  \,\mbox{``Maier~Thm.(7.3)~with~L=n~variation
~2''}
\]
}
\end{maplelatex}

\begin{maplelatex}
\mapleinline{inert}{2d}{`  121: `, "Maier Thm.(7.3) with L=n variation 2";}{%
\[
\mathit{\ \ 121:\ }  \,\mbox{``Maier~Thm.(7.3)~with~L=n~variation
~2''}
\]
}
\end{maplelatex}

\begin{maplelatex}
\mapleinline{inert}{2d}{`  122: `, "Maier Thm.(7.3) with L=n variation 2";}{%
\[
\mathit{\ \ 122:\ }  \,\mbox{``Maier~Thm.(7.3)~with~L=n~variation
~2''}
\]
}
\end{maplelatex}

\begin{maplelatex}
\mapleinline{inert}{2d}{`  123: `, "Maier Thm.(7.3) with L=n variation 2";}{%
\[
\mathit{\ \ 123:\ }  \,\mbox{``Maier~Thm.(7.3)~with~L=n~variation
~2''}
\]
}
\end{maplelatex}

\begin{maplelatex}
\mapleinline{inert}{2d}{`  124: `, "Maier Thm.(7.3) with L=n variation 1";}{%
\[
\mathit{\ \ 124:\ }  \,\mbox{``Maier~Thm.(7.3)~with~L=n~variation
~1''}
\]
}
\end{maplelatex}

\begin{maplelatex}
\mapleinline{inert}{2d}{`  125: `, "Maier Thm.(7.3) with L=n variation 1";}{%
\[
\mathit{\ \ 125:\ }  \,\mbox{``Maier~Thm.(7.3)~with~L=n~variation
~1''}
\]
}
\end{maplelatex}

\begin{maplelatex}
\mapleinline{inert}{2d}{`  126: `, "Maier Thm.(7.3) with L=n variation 1";}{%
\[
\mathit{\ \ 126:\ }  \,\mbox{``Maier~Thm.(7.3)~with~L=n~variation
~1''}
\]
}
\end{maplelatex}

\begin{maplelatex}
\mapleinline{inert}{2d}{`  127: `, "Maier Thm.(7.3) with L=n variation 1";}{%
\[
\mathit{\ \ 127:\ }  \,\mbox{``Maier~Thm.(7.3)~with~L=n~variation
~1''}
\]
}
\end{maplelatex}

\begin{maplelatex}
\mapleinline{inert}{2d}{`  128: `, "Maier Thm.(7.3) with L=n variation 1";}{%
\[
\mathit{\ \ 128:\ }  \,\mbox{``Maier~Thm.(7.3)~with~L=n~variation
~1''}
\]
}
\end{maplelatex}

\begin{maplelatex}
\mapleinline{inert}{2d}{`  129: `, "Maier Thm.(7.3) with L=n variation 1";}{%
\[
\mathit{\ \ 129:\ }  \,\mbox{``Maier~Thm.(7.3)~with~L=n~variation
~1''}
\]
}
\end{maplelatex}

\begin{maplelatex}
\mapleinline{inert}{2d}{`  130: `, "Maier Thm.(7.3) with L=n variation 1";}{%
\[
\mathit{\ \ 130:\ }  \,\mbox{``Maier~Thm.(7.3)~with~L=n~variation
~1''}
\]
}
\end{maplelatex}

\begin{maplelatex}
\mapleinline{inert}{2d}{`  131: `, "Maier Thm.(7.3) with L=n variation 1";}{%
\[
\mathit{\ \ 131:\ }  \,\mbox{``Maier~Thm.(7.3)~with~L=n~variation
~1''}
\]
}
\end{maplelatex}

\begin{maplelatex}
\mapleinline{inert}{2d}{`  132: `, "Maier Thm.(7.3) with L=n variation 1";}{%
\[
\mathit{\ \ 132:\ }  \,\mbox{``Maier~Thm.(7.3)~with~L=n~variation
~1''}
\]
}
\end{maplelatex}

\begin{maplelatex}
\mapleinline{inert}{2d}{`  133: `, "Maier Thm.(7.3) with L=n variation 1";}{%
\[
\mathit{\ \ 133:\ }  \,\mbox{``Maier~Thm.(7.3)~with~L=n~variation
~1''}
\]
}
\end{maplelatex}

\begin{maplelatex}
\mapleinline{inert}{2d}{`  134: `, "Prudnikov 7.4.4.32 + Miller J.Phys.A., Eq. 1.1";}{%
\[
\mathit{\ \ 134:\ }  \,\mbox{``Prudnikov~7.4.4.32~+~Miller~
J.Phys.A.,~Eq.~1.1''}
\]
}
\end{maplelatex}

\begin{maplelatex}
\mapleinline{inert}{2d}{`  139: `, "Prudnikov 7.4.4.32, n odd + Miller J.Phys.A., Eq.
1.1";}{%
\[
\mathit{\ \ 139:\ }  \,\mbox{``Prudnikov~7.4.4.32,~n~odd~+~Miller
~J.Phys.A.,~Eq.~1.1''}
\]
}
\end{maplelatex}

\begin{maplelatex}
\mapleinline{inert}{2d}{`  148: `, "Prudnikov 7.4.4.32 + Miller J.Phys.A., Eq. 1.1+T:1";}{%
\[
\mathit{\ \ 148:\ }  \,\mbox{``Prudnikov~7.4.4.32~+~Miller~
J.Phys.A.,~Eq.~1.1+T:1''}
\]
}
\end{maplelatex}

\begin{maplelatex}
\mapleinline{inert}{2d}{`  149: `, "Prudnikov 7.4.4.32 + Miller J.Phys.A., Eq. 1.1+T:3";}{%
\[
\mathit{\ \ 149:\ }  \,\mbox{``Prudnikov~7.4.4.32~+~Miller~
J.Phys.A.,~Eq.~1.1+T:3''}
\]
}
\end{maplelatex}

\begin{maplelatex}
\mapleinline{inert}{2d}{`  150: `, "Prudnikov 7.4.4.32 + Miller J.Phys.A., Eq. 1.1+T:5";}{%
\[
\mathit{\ \ 150:\ }  \,\mbox{``Prudnikov~7.4.4.32~+~Miller~
J.Phys.A.,~Eq.~1.1+T:5''}
\]
}
\end{maplelatex}

\begin{maplelatex}
\mapleinline{inert}{2d}{`  151: `, "Prudnikov 7.4.4.32 + Miller J.Phys.A., Eq. 1.1+T:6";}{%
\[
\mathit{\ \ 151:\ }  \,\mbox{``Prudnikov~7.4.4.32~+~Miller~
J.Phys.A.,~Eq.~1.1+T:6''}
\]
}
\end{maplelatex}

\begin{maplelatex}
\mapleinline{inert}{2d}{`  152: `, "Prudnikov 7.4.4.32 + Miller J.Phys.A., Eq. 1.1+T:7";}{%
\[
\mathit{\ \ 152:\ }  \,\mbox{``Prudnikov~7.4.4.32~+~Miller~
J.Phys.A.,~Eq.~1.1+T:7''}
\]
}
\end{maplelatex}

\begin{maplelatex}
\mapleinline{inert}{2d}{`  160: `, "Prudnikov 7.4.4.32, n odd + Miller J.Phys.A., Eq.
1.1+T:1";}{%
\[
\mathit{\ \ 160:\ }  \,\mbox{``Prudnikov~7.4.4.32,~n~odd~+~Miller
~J.Phys.A.,~Eq.~1.1+T:1''}
\]
}
\end{maplelatex}

\begin{maplelatex}
\mapleinline{inert}{2d}{`  161: `, "Prudnikov 7.4.4.32, n odd + Miller J.Phys.A., Eq.
1.1+T:2";}{%
\[
\mathit{\ \ 161:\ }  \,\mbox{``Prudnikov~7.4.4.32,~n~odd~+~Miller
~J.Phys.A.,~Eq.~1.1+T:2''}
\]
}
\end{maplelatex}

\begin{maplelatex}
\mapleinline{inert}{2d}{`  162: `, "Prudnikov 7.4.4.32, n odd + Miller J.Phys.A., Eq.
1.1+T:3";}{%
\[
\mathit{\ \ 162:\ }  \,\mbox{``Prudnikov~7.4.4.32,~n~odd~+~Miller
~J.Phys.A.,~Eq.~1.1+T:3''}
\]
}
\end{maplelatex}

\begin{maplelatex}
\mapleinline{inert}{2d}{`  163: `, "Prudnikov 7.4.4.32, n odd + Miller J.Phys.A., Eq.
1.1+T:4";}{%
\[
\mathit{\ \ 163:\ }  \,\mbox{``Prudnikov~7.4.4.32,~n~odd~+~Miller
~J.Phys.A.,~Eq.~1.1+T:4''}
\]
}
\end{maplelatex}

\begin{maplelatex}
\mapleinline{inert}{2d}{`  164: `, "Prudnikov 7.4.4.32, n odd + Miller J.Phys.A., Eq.
1.1+T:5";}{%
\[
\mathit{\ \ 164:\ }  \,\mbox{``Prudnikov~7.4.4.32,~n~odd~+~Miller
~J.Phys.A.,~Eq.~1.1+T:5''}
\]
}
\end{maplelatex}

\begin{maplelatex}
\mapleinline{inert}{2d}{`  165: `, "Prudnikov 7.4.4.32, n odd + Miller J.Phys.A., Eq.
1.1+T:6";}{%
\[
\mathit{\ \ 165:\ }  \,\mbox{``Prudnikov~7.4.4.32,~n~odd~+~Miller
~J.Phys.A.,~Eq.~1.1+T:6''}
\]
}
\end{maplelatex}

\begin{maplelatex}
\mapleinline{inert}{2d}{`  166: `, "Prudnikov 7.4.4.32, n odd + Miller J.Phys.A., Eq.
1.1+T:7";}{%
\[
\mathit{\ \ 166:\ }  \,\mbox{``Prudnikov~7.4.4.32,~n~odd~+~Miller
~J.Phys.A.,~Eq.~1.1+T:7''}
\]
}
\end{maplelatex}

\begin{maplelatex}
\mapleinline{inert}{2d}{`  167: `, "Prudnikov 7.4.4.32, n odd + Miller J.Phys.A., Eq.
1.1+T:8";}{%
\[
\mathit{\ \ 167:\ }  \,\mbox{``Prudnikov~7.4.4.32,~n~odd~+~Miller
~J.Phys.A.,~Eq.~1.1+T:8''}
\]
}
\end{maplelatex}

\begin{maplelatex}
\mapleinline{inert}{2d}{`  168: `, "Prudnikov 7.4.4.32, n odd + Miller J.Phys.A., Eq.
1.1+T:9";}{%
\[
\mathit{\ \ 168:\ }  \,\mbox{``Prudnikov~7.4.4.32,~n~odd~+~Miller
~J.Phys.A.,~Eq.~1.1+T:9''}
\]
}
\end{maplelatex}

\begin{maplelatex}
\mapleinline{inert}{2d}{`  194: `, "Prudnikov 7.4.4.32 + Miller J.Phys.A., Eq. 1.1";}{%
\[
\mathit{\ \ 194:\ }  \,\mbox{``Prudnikov~7.4.4.32~+~Miller~
J.Phys.A.,~Eq.~1.1''}
\]
}
\end{maplelatex}

\begin{maplelatex}
\mapleinline{inert}{2d}{`  195: `, "Prudnikov 7.4.4.32, n odd + Miller J.Phys.A., Eq.
1.1";}{%
\[
\mathit{\ \ 195:\ }  \,\mbox{``Prudnikov~7.4.4.32,~n~odd~+~Miller
~J.Phys.A.,~Eq.~1.1''}
\]
}
\end{maplelatex}

\begin{maplelatex}
\mapleinline{inert}{2d}{`  196: `, "Prudnikov 7.4.4.32 + Miller J.Phys.A., Eq. 1.1+T:3 +
Miller J.Phys.A., Eq. 1.1";}{%
\[
\mathit{\ \ 196:\ }  \,\mbox{``Prudnikov~7.4.4.32~+~Miller~
J.Phys.A.,~Eq.~1.1+T:3~+~Miller~J.Phys.A.,~Eq.~1.1''}
\]
}
\end{maplelatex}

\begin{maplelatex}
\mapleinline{inert}{2d}{`  197: `, "Prudnikov 7.4.4.32 + Miller J.Phys.A., Eq. 1.1+T:7 +
Miller J.Phys.A., Eq. 1.1";}{%
\[
\mathit{\ \ 197:\ }  \,\mbox{``Prudnikov~7.4.4.32~+~Miller~
J.Phys.A.,~Eq.~1.1+T:7~+~Miller~J.Phys.A.,~Eq.~1.1''}
\]
}
\end{maplelatex}

\begin{maplelatex}
\mapleinline{inert}{2d}{`  198: `, "Prudnikov 7.4.4.32, n odd + Miller J.Phys.A., Eq. 1.1+T:3
+ Miller J.Phys.A., Eq. 1.1";}{%
\maplemultiline{
\mathit{\ \ 198:\ }  \mbox{``Prudnikov~7.4.4.32,~n~odd~+~Miller~
J.Phys.A.,~Eq.~1.1+T:3~+~Miller~J.Phys.A.,} 
\mbox{~Eq.~1.1''} }
}
\end{maplelatex}

\begin{maplelatex}
\mapleinline{inert}{2d}{`  199: `, "Prudnikov 7.4.4.32, n odd + Miller J.Phys.A., Eq. 1.1+T:4
+ Miller J.Phys.A., Eq. 1.1";}{%
\maplemultiline{
\mathit{\ \ 199:\ }  \mbox{``Prudnikov~7.4.4.32,~n~odd~+~Miller~
J.Phys.A.,~Eq.~1.1+T:4~+~Miller~J.Phys.A.,} 
\mbox{~Eq.~1.1''} }
}
\end{maplelatex}

\begin{maplelatex}
\mapleinline{inert}{2d}{`  200: `, "Prudnikov 7.4.4.32, n odd + Miller J.Phys.A., Eq. 1.1+T:6
+ Miller J.Phys.A., Eq. 1.1";}{%
\maplemultiline{
\mathit{\ \ 200:\ }  \mbox{``Prudnikov~7.4.4.32,~n~odd~+~Miller~
J.Phys.A.,~Eq.~1.1+T:6~+~Miller~J.Phys.A.,} 
\mbox{~Eq.~1.1''} }
}
\end{maplelatex}

\begin{maplelatex}
\mapleinline{inert}{2d}{`  201: `, "Prudnikov 7.4.4.32, n odd + Miller J.Phys.A., Eq. 1.1+T:8
+ Miller J.Phys.A., Eq. 1.1";}{%
\maplemultiline{
\mathit{\ \ 201:\ }  \mbox{``Prudnikov~7.4.4.32,~n~odd~+~Miller~
J.Phys.A.,~Eq.~1.1+T:8~+~Miller~J.Phys.A.,} \mbox{~Eq.~1.1''} }
}
\end{maplelatex}

\begin{maplelatex}
\mapleinline{inert}{2d}{`  210: `, "Prudnikov 7.4.4.32 + Miller J.Phys.A., Eq. 1.1+T:1";}{%
\[
\mathit{\ \ 210:\ }  \,\mbox{``Prudnikov~7.4.4.32~+~Miller~
J.Phys.A.,~Eq.~1.1+T:1''}
\]
}
\end{maplelatex}

\begin{maplelatex}
\mapleinline{inert}{2d}{`  211: `, "Prudnikov 7.4.4.32 + Miller J.Phys.A., Eq. 1.1+T:2";}{%
\[
\mathit{\ \ 211:\ }  \,\mbox{``Prudnikov~7.4.4.32~+~Miller~
J.Phys.A.,~Eq.~1.1+T:2''}
\]
}
\end{maplelatex}

\begin{maplelatex}
\mapleinline{inert}{2d}{`  212: `, "Prudnikov 7.4.4.32 + Miller J.Phys.A., Eq. 1.1+T:3";}{%
\[
\mathit{\ \ 212:\ }  \,\mbox{``Prudnikov~7.4.4.32~+~Miller~
J.Phys.A.,~Eq.~1.1+T:3''}
\]
}
\end{maplelatex}

\begin{maplelatex}
\mapleinline{inert}{2d}{`  213: `, "Prudnikov 7.4.4.32 + Miller J.Phys.A., Eq. 1.1+T:4";}{%
\[
\mathit{\ \ 213:\ }  \,\mbox{``Prudnikov~7.4.4.32~+~Miller~
J.Phys.A.,~Eq.~1.1+T:4''}
\]
}
\end{maplelatex}

\begin{maplelatex}
\mapleinline{inert}{2d}{`  214: `, "Prudnikov 7.4.4.32 + Miller J.Phys.A., Eq. 1.1+T:5";}{%
\[
\mathit{\ \ 214:\ }  \,\mbox{``Prudnikov~7.4.4.32~+~Miller~
J.Phys.A.,~Eq.~1.1+T:5''}
\]
}
\end{maplelatex}

\begin{maplelatex}
\mapleinline{inert}{2d}{`  215: `, "Prudnikov 7.4.4.32 + Miller J.Phys.A., Eq. 1.1+T:6";}{%
\[
\mathit{\ \ 215:\ }  \,\mbox{``Prudnikov~7.4.4.32~+~Miller~
J.Phys.A.,~Eq.~1.1+T:6''}
\]
}
\end{maplelatex}

\begin{maplelatex}
\mapleinline{inert}{2d}{`  216: `, "Prudnikov 7.4.4.32 + Miller J.Phys.A., Eq. 1.1+T:7";}{%
\[
\mathit{\ \ 216:\ }  \,\mbox{``Prudnikov~7.4.4.32~+~Miller~
J.Phys.A.,~Eq.~1.1+T:7''}
\]
}
\end{maplelatex}

\begin{maplelatex}
\mapleinline{inert}{2d}{`  217: `, "Prudnikov 7.4.4.32 + Miller J.Phys.A., Eq. 1.1+T:8";}{%
\[
\mathit{\ \ 217:\ }  \,\mbox{``Prudnikov~7.4.4.32~+~Miller~
J.Phys.A.,~Eq.~1.1+T:8''}
\]
}
\end{maplelatex}

\begin{maplelatex}
\mapleinline{inert}{2d}{`  218: `, "Prudnikov 7.4.4.32 + Miller J.Phys.A., Eq. 1.1+T:9";}{%
\[
\mathit{\ \ 218:\ }  \,\mbox{``Prudnikov~7.4.4.32~+~Miller~
J.Phys.A.,~Eq.~1.1+T:9''}
\]
}
\end{maplelatex}

\begin{maplelatex}
\mapleinline{inert}{2d}{`  219: `, "Prudnikov 7.4.4.32, n odd + Miller J.Phys.A., Eq.
1.1+T:1";}{%
\[
\mathit{\ \ 219:\ }  \,\mbox{``Prudnikov~7.4.4.32,~n~odd~+~Miller
~J.Phys.A.,~Eq.~1.1+T:1''}
\]
}
\end{maplelatex}

\begin{maplelatex}
\mapleinline{inert}{2d}{`  220: `, "Prudnikov 7.4.4.32, n odd + Miller J.Phys.A., Eq.
1.1+T:2";}{%
\[
\mathit{\ \ 220:\ }  \,\mbox{``Prudnikov~7.4.4.32,~n~odd~+~Miller
~J.Phys.A.,~Eq.~1.1+T:2''}
\]
}
\end{maplelatex}

\begin{maplelatex}
\mapleinline{inert}{2d}{`  221: `, "Prudnikov 7.4.4.32, n odd + Miller J.Phys.A., Eq.
1.1+T:6";}{%
\[
\mathit{\ \ 221:\ }  \,\mbox{``Prudnikov~7.4.4.32,~n~odd~+~Miller
~J.Phys.A.,~Eq.~1.1+T:6''}
\]
}
\end{maplelatex}

\begin{maplelatex}
\mapleinline{inert}{2d}{`  222: `, "Prudnikov 7.4.4.32, n odd + Miller J.Phys.A., Eq.
1.1+T:7";}{%
\[
\mathit{\ \ 222:\ }  \,\mbox{``Prudnikov~7.4.4.32,~n~odd~+~Miller
~J.Phys.A.,~Eq.~1.1+T:7''}
\]
}
\end{maplelatex}

\begin{maplelatex}
\mapleinline{inert}{2d}{`  223: `, "Prudnikov 7.4.4.32, n odd + Miller J.Phys.A., Eq.
1.1+T:9";}{%
\[
\mathit{\ \ 223:\ }  \,\mbox{``Prudnikov~7.4.4.32,~n~odd~+~Miller
~J.Phys.A.,~Eq.~1.1+T:9''}
\]
}
\end{maplelatex}

\begin{maplelatex}
\mapleinline{inert}{2d}{`  224: `, "Prudnikov 7.4.4.32 + Miller J.Phys.A., Eq. 1.1+T:3 +
Miller J.Phys.A., Eq. 1.1+T:1";}{%
\maplemultiline{
\mathit{\ \ 224:\ }  \mbox{``Prudnikov~7.4.4.32~+~Miller~
J.Phys.A.,~Eq.~1.1+T:3~+~Miller~J.Phys.A.,~}
\mbox{Eq.~1.1+T:1''} }
}
\end{maplelatex}

\begin{maplelatex}
\mapleinline{inert}{2d}{`  225: `, "Prudnikov 7.4.4.32 + Miller J.Phys.A., Eq. 1.1+T:3 +
Miller J.Phys.A., Eq. 1.1+T:2";}{%
\maplemultiline{
\mathit{\ \ 225:\ }  \mbox{``Prudnikov~7.4.4.32~+~Miller~
J.Phys.A.,~Eq.~1.1+T:3~+~Miller~J.Phys.A.,~}
\mbox{Eq.~1.1+T:2''} }
}
\end{maplelatex}

\begin{maplelatex}
\mapleinline{inert}{2d}{`  226: `, "Prudnikov 7.4.4.32 + Miller J.Phys.A., Eq. 1.1+T:3 +
Miller J.Phys.A., Eq. 1.1+T:3";}{%
\maplemultiline{
\mathit{\ \ 226:\ }  \mbox{``Prudnikov~7.4.4.32~+~Miller~
J.Phys.A.,~Eq.~1.1+T:3~+~Miller~J.Phys.A.,~ }
\mbox{Eq.~1.1+T:3''} }
}
\end{maplelatex}

\begin{maplelatex}
\mapleinline{inert}{2d}{`  227: `, "Prudnikov 7.4.4.32 + Miller J.Phys.A., Eq. 1.1+T:3 +
Miller J.Phys.A., Eq. 1.1+T:4";}{%
\maplemultiline{
\mathit{\ \ 227:\ }  \mbox{``Prudnikov~7.4.4.32~+~Miller~
J.Phys.A.,~Eq.~1.1+T:3~+~Miller~J.Phys.A.,~}
\mbox{Eq.~1.1+T:4''} }
}
\end{maplelatex}

\begin{maplelatex}
\mapleinline{inert}{2d}{`  228: `, "Prudnikov 7.4.4.32 + Miller J.Phys.A., Eq. 1.1+T:3 +
Miller J.Phys.A., Eq. 1.1+T:5";}{%
\maplemultiline{
\mathit{\ \ 228:\ }  \mbox{``Prudnikov~7.4.4.32~+~Miller~
J.Phys.A.,~Eq.~1.1+T:3~+~Miller~J.Phys.A.,~}
\mbox{Eq.~1.1+T:5''} }
}
\end{maplelatex}

\begin{maplelatex}
\mapleinline{inert}{2d}{`  229: `, "Prudnikov 7.4.4.32 + Miller J.Phys.A., Eq. 1.1+T:3 +
Miller J.Phys.A., Eq. 1.1+T:7";}{%
\maplemultiline{
\mathit{\ \ 229:\ }  \mbox{``Prudnikov~7.4.4.32~+~Miller~
J.Phys.A.,~Eq.~1.1+T:3~+~Miller~J.Phys.A.,~}
\mbox{Eq.~1.1+T:7''} }
}
\end{maplelatex}

\begin{maplelatex}
\mapleinline{inert}{2d}{`  230: `, "Prudnikov 7.4.4.32 + Miller J.Phys.A., Eq. 1.1+T:3 +
Miller J.Phys.A., Eq. 1.1+T:8";}{%
\maplemultiline{
\mathit{\ \ 230:\ }  \mbox{``Prudnikov~7.4.4.32~+~Miller~
J.Phys.A.,~Eq.~1.1+T:3~+~Miller~J.Phys.A.,~}
\mbox{Eq.~1.1+T:8''} }
}
\end{maplelatex}

\begin{maplelatex}
\mapleinline{inert}{2d}{`  231: `, "Prudnikov 7.4.4.32 + Miller J.Phys.A., Eq. 1.1+T:3 +
Miller J.Phys.A., Eq. 1.1+T:9";}{%
\maplemultiline{
\mathit{\ \ 231:\ }  \mbox{``Prudnikov~7.4.4.32~+~Miller~
J.Phys.A.,~Eq.~1.1+T:3~+~Miller~J.Phys.A.,~}
\mbox{Eq.~1.1+T:9''} }
}
\end{maplelatex}

\begin{maplelatex}
\mapleinline{inert}{2d}{`  232: `, "Prudnikov 7.4.4.32, n odd + Miller J.Phys.A., Eq. 1.1+T:3
+ Miller J.Phys.A., Eq. 1.1+T:1";}{%
\maplemultiline{
\mathit{\ \ 232:\ }  \mbox{``Prudnikov~7.4.4.32,~n~odd~+~Miller~
J.Phys.A.,~Eq.~1.1+T:3~+~Miller~J.Phys}
\mbox{.A.,~Eq.~1.1+T:1''} }
}
\end{maplelatex}

\begin{maplelatex}
\mapleinline{inert}{2d}{`  233: `, "Prudnikov 7.4.4.32, n odd + Miller J.Phys.A., Eq. 1.1+T:3
+ Miller J.Phys.A., Eq. 1.1+T:2";}{%
\maplemultiline{
\mathit{\ \ 233:\ }  \mbox{``Prudnikov~7.4.4.32,~n~odd~+~Miller~
J.Phys.A.,~Eq.~1.1+T:3~+~Miller~J.Phys}
\mbox{.A.,~Eq.~1.1+T:2''} }
}
\end{maplelatex}

\begin{maplelatex}
\mapleinline{inert}{2d}{`  234: `, "Prudnikov 7.4.4.32, n odd + Miller J.Phys.A., Eq. 1.1+T:3
+ Miller J.Phys.A., Eq. 1.1+T:4";}{%
\maplemultiline{
\mathit{\ \ 234:\ }  \mbox{``Prudnikov~7.4.4.32,~n~odd~+~Miller~
J.Phys.A.,~Eq.~1.1+T:3~+~Miller~J.Phys}
\mbox{.A.,~Eq.~1.1+T:4''} }
}
\end{maplelatex}

\begin{maplelatex}
\mapleinline{inert}{2d}{`  235: `, "Prudnikov 7.4.4.32, n odd + Miller J.Phys.A., Eq. 1.1+T:3
+ Miller J.Phys.A., Eq. 1.1+T:5";}{%
\maplemultiline{
\mathit{\ \ 235:\ }  \mbox{``Prudnikov~7.4.4.32,~n~odd~+~Miller~
J.Phys.A.,~Eq.~1.1+T:3~+~Miller~J.Phys}
\mbox{.A.,~Eq.~1.1+T:5''} }
}
\end{maplelatex}

\begin{maplelatex}
\mapleinline{inert}{2d}{`  236: `, "Prudnikov 7.4.4.32, n odd + Miller J.Phys.A., Eq. 1.1+T:3
+ Miller J.Phys.A., Eq. 1.1+T:6";}{%
\maplemultiline{
\mathit{\ \ 236:\ }  \mbox{``Prudnikov~7.4.4.32,~n~odd~+~Miller~
J.Phys.A.,~Eq.~1.1+T:3~+~Miller~J.Phys}
\mbox{.A.,~Eq.~1.1+T:6''} }
}
\end{maplelatex}

\begin{maplelatex}
\mapleinline{inert}{2d}{`  237: `, "Prudnikov 7.4.4.32, n odd + Miller J.Phys.A., Eq. 1.1+T:3
+ Miller J.Phys.A., Eq. 1.1+T:7";}{%
\maplemultiline{
\mathit{\ \ 237:\ }  \mbox{``Prudnikov~7.4.4.32,~n~odd~+~Miller~
J.Phys.A.,~Eq.~1.1+T:3~+~Miller~J.Phys}
\mbox{.A.,~Eq.~1.1+T:7''} }
}
\end{maplelatex}

\begin{maplelatex}
\mapleinline{inert}{2d}{`  238: `, "Prudnikov 7.4.4.32 + Miller J.Phys.A., Eq. 1.1+T:6 +
Miller J.Phys.A., Eq. 1.1";}{%
\[
\mathit{\ \ 238:\ }  \,\mbox{``Prudnikov~7.4.4.32~+~Miller~
J.Phys.A.,~Eq.~1.1+T:6~+~Miller~J.Phys.A.,~Eq.~1.1''}
\]
}
\end{maplelatex}

\begin{maplelatex}
\mapleinline{inert}{2d}{`  239: `, "Prudnikov 7.4.4.32, n odd + Miller J.Phys.A., Eq. 1.1+T:5
+ Miller J.Phys.A., Eq. 1.1";}{%
\maplemultiline{
\mathit{\ \ 239:\ }  \mbox{``Prudnikov~7.4.4.32,~n~odd~+~Miller~
J.Phys.A.,~Eq.~1.1+T:5~+~Miller~J.Phys}
\mbox{.A.,~Eq.~1.1''} }
}
\end{maplelatex}

\begin{maplelatex}
\mapleinline{inert}{2d}{`  240: `, "Prudnikov 7.4.4.32, n odd + Miller J.Phys.A., Eq. 1.1+T:7
+ Miller J.Phys.A., Eq. 1.1";}{%
\maplemultiline{
\mathit{\ \ 240:\ }  \mbox{``Prudnikov~7.4.4.32,~n~odd~+~Miller~
J.Phys.A.,~Eq.~1.1+T:7~+~Miller~J.Phys}
\mbox{.A.,~Eq.~1.1''} }
}
\end{maplelatex}

\begin{maplelatex}
\mapleinline{inert}{2d}{`  241: `, "Prudnikov 7.4.4.32 + Miller J.Phys.A., Eq. 1.1+T:2 +
Miller J.Phys.A., Eq. 1.1";}{%
\[
\mathit{\ \ 241:\ }  \,\mbox{``Prudnikov~7.4.4.32~+~Miller~
J.Phys.A.,~Eq.~1.1+T:2~+~Miller~J.Phys.A.,~Eq.~1.1''}
\]
}
\end{maplelatex}

\begin{maplelatex}
\mapleinline{inert}{2d}{`  242: `, "Prudnikov 7.4.4.32 + Miller J.Phys.A., Eq. 1.1+T:4 +
Miller J.Phys.A., Eq. 1.1";}{%
\[
\mathit{\ \ 242:\ }  \,\mbox{``Prudnikov~7.4.4.32~+~Miller~
J.Phys.A.,~Eq.~1.1+T:4~+~Miller~J.Phys.A.,~Eq.~1.1''}
\]
}
\end{maplelatex}

\begin{maplelatex}
\mapleinline{inert}{2d}{`  243: `, "Prudnikov 7.4.4.32, n odd + Miller J.Phys.A., Eq. 1.1+T:6
+ Miller J.Phys.A., Eq. 1.1";}{%
\maplemultiline{
\mathit{\ \ 243:\ }  \mbox{``Prudnikov~7.4.4.32,~n~odd~+~Miller~
J.Phys.A.,~Eq.~1.1+T:6~+~Miller~J.Phys}
\mbox{.A.,~Eq.~1.1''} }
}
\end{maplelatex}

\begin{maplelatex}
\mapleinline{inert}{2d}{`  244: `, "Gessel and Stanton, SIAM J. Math. Anal.,13,295(1982)
Eq.(5.16) + Miller, J.Phys.A.,Eq. 1.2 (e=c)";}{%
\maplemultiline{
\mathit{\ \ 244:\ }  \mbox{``Gessel~and~Stanton,~SIAM~J.~Math.~
Anal.,13,295(1982)~Eq.(5.16)~+~Mille}
\mbox{r,~J.Phys.A.,Eq.~1.2~(e=c)''} }
}
\end{maplelatex}

\begin{maplelatex}
\mapleinline{inert}{2d}{`  245: `, "Gessel and Stanton, SIAM J. Math. Anal.,13,295(1982)
Eq.(5.16) + Miller, J.Phys.A.,Eq. 1.2 (e=c)";}{%
\maplemultiline{
\mathit{\ \ 245:\ }  \mbox{``Gessel~and~Stanton,~SIAM~J.~Math.~
Anal.,13,295(1982)~Eq.(5.16)~+~Mille}
\mbox{r,~J.Phys.A.,Eq.~1.2~(e=c)''} }
}
\end{maplelatex}

\begin{maplelatex}
\mapleinline{inert}{2d}{`  246: `, "Gessel and Stanton, SIAM J. Math. Anal.,13,295(1982)
Eq.(5.16) + Miller, J.Phys.A.,Eq. 1.2 (e=c)+T:1";}{%
\maplemultiline{
\mathit{\ \ 246:\ }  \mbox{``Gessel~and~Stanton,~SIAM~J.~Math.~
Anal.,13,295(1982)~Eq.(5.16)~+~Mille}
\mbox{r,~J.Phys.A.,Eq.~1.2~(e=c)+T:1''} }
}
\end{maplelatex}

\begin{maplelatex}
\mapleinline{inert}{2d}{`  247: `, "Gessel and Stanton, SIAM J. Math. Anal.,13,295(1982)
Eq.(5.16) + Miller, J.Phys.A.,Eq. 1.2 (e=c)+T:4";}{%
\maplemultiline{
\mathit{\ \ 247:\ }  \mbox{``Gessel~and~Stanton,~SIAM~J.~Math.~
Anal.,13,295(1982)~Eq.(5.16)~+~Mille}
\mbox{r,~J.Phys.A.,Eq.~1.2~(e=c)+T:4''} }
}
\end{maplelatex}

\begin{maplelatex}
\mapleinline{inert}{2d}{`  248: `, "Gessel and Stanton, SIAM J. Math. Anal.,13,295(1982)
Eq.(5.16) + Miller, J.Phys.A.,Eq. 1.2 (e=c)+T:6";}{%
\maplemultiline{
\mathit{\ \ 248:\ }  \mbox{``Gessel~and~Stanton,~SIAM~J.~Math.~
Anal.,13,295(1982)~Eq.(5.16)~+~Mille}
\mbox{r,~J.Phys.A.,Eq.~1.2~(e=c)+T:6''} }
}
\end{maplelatex}

\begin{maplelatex}
\mapleinline{inert}{2d}{`  253: `, "Maier Eq. (7.1) with L=n+ Miller, J.Phys.A.,Eq. 1.2 (g =
d)";}{%
\[
\mathit{\ \ 253:\ }  \,\mbox{``Maier Eq. (7.1)~with~L=n+~Miller,~
J.Phys.A.,Eq.~1.2~(g~=~d)''}
\]
}
\end{maplelatex}

\begin{maplelatex}
\mapleinline{inert}{2d}{`  254: `, "Maier Eq. (7.1) with L=n+ Miller, J.Phys.A.,Eq. 1.2 (g =
d)";}{%
\[
\mathit{\ \ 254:\ }  \,\mbox{``Maier Eq. (7.1)~with~L=n+~Miller,~
J.Phys.A.,Eq.~1.2~(g~=~d)''}
\]
}
\end{maplelatex}

\begin{maplelatex}
\mapleinline{inert}{2d}{`  255: `, "Maier Eq. (7.1) with L=n+ Miller, J.Phys.A.,Eq. 1.2 (g =
d)";}{%
\[
\mathit{\ \ 255:\ }  \,\mbox{``Maier Eq. (7.1)~with~L=n+~Miller,~
J.Phys.A.,Eq.~1.2~(g~=~d)''}
\]
}
\end{maplelatex}

\begin{maplelatex}
\mapleinline{inert}{2d}{`  256: `, "Maier Eq. (7.1) with L=n+ Miller, J.Phys.A.,Eq. 1.2 (g =
d)";}{%
\[
\mathit{\ \ 256:\ }  \,\mbox{``Maier Eq. (7.1)~with~L=n+~Miller,~
J.Phys.A.,Eq.~1.2~(g~=~d)''}
\]
}
\end{maplelatex}

\begin{maplelatex}
\mapleinline{inert}{2d}{`  257: `, "Maier Eq. (7.1) with L=n+ Miller, J.Phys.A.,Eq. 1.2 (g =
d)";}{%
\[
\mathit{\ \ 257:\ }  \,\mbox{``Maier Eq. (7.1)~with~L=n+~Miller,~
J.Phys.A.,Eq.~1.2~(g~=~d)''}
\]
}
\end{maplelatex}

\begin{maplelatex}
\mapleinline{inert}{2d}{`  258: `, "Maier Eq. (7.1) with L=n+ Miller, J.Phys.A.,Eq. 1.2 (g =
d)";}{%
\[
\mathit{\ \ 258:\ }  \,\mbox{``Maier Eq. (7.1)~with~L=n+~Miller,~
J.Phys.A.,Eq.~1.2~(g~=~d)''}
\]
}
\end{maplelatex}

\begin{maplelatex}
\mapleinline{inert}{2d}{`  259: `, "Minton-Karlsson, J.Math. Phys., 12,2(1971), p370,
Eq.(8)";}{%
\[
\mathit{\ \ 259:\ }  \,\mbox{``Minton-Karlsson,~J.Math.~Phys.,~
12,2(1971),~p370,~Eq.(8)''}
\]
}
\end{maplelatex}

\begin{maplelatex}
\mapleinline{inert}{2d}{`  260: `, "standard expression for top negative parameter";}{%
\[
\mathit{\ \ 260:\ }  \,\mbox{``standard~expression~for~top~
negative~parameter''}
\]
}
\end{maplelatex}

\begin{maplelatex}
\mapleinline{inert}{2d}{`  267: `, "Maier Eq. (7.1) with L=n+ Miller, J.Phys.A.,Eq. 1.2 (g =
d)+T:1";}{%
\[
\mathit{\ \ 267:\ }  \,\mbox{``Maier Eq. (7.1)~with~L=n+~Miller,~
J.Phys.A.,Eq.~1.2~(g~=~d)+T:1''}
\]
}
\end{maplelatex}

\begin{maplelatex}
\mapleinline{inert}{2d}{`  268: `, "Maier Eq. (7.1) with L=n+ Miller, J.Phys.A.,Eq. 1.2 (g =
d)+T:2";}{%
\[
\mathit{\ \ 268:\ }  \,\mbox{``Maier Eq. (7.1)~with~L=n+~Miller,~
J.Phys.A.,Eq.~1.2~(g~=~d)+T:2''}
\]
}
\end{maplelatex}

\begin{maplelatex}
\mapleinline{inert}{2d}{`  269: `, "Maier Eq. (7.1) with L=n+ Miller, J.Phys.A.,Eq. 1.2 (g =
d)+T:3";}{%
\[
\mathit{\ \ 269:\ }  \,\mbox{``Maier Eq. (7.1)~with~L=n+~Miller,~
J.Phys.A.,Eq.~1.2~(g~=~d)+T:3''}
\]
}
\end{maplelatex}

\begin{maplelatex}
\mapleinline{inert}{2d}{`  270: `, "Maier Eq. (7.1) with L=n+ Miller, J.Phys.A.,Eq. 1.2 (g =
d)+T:4";}{%
\[
\mathit{\ \ 270:\ }  \,\mbox{``Maier Eq. (7.1)~with~L=n+~Miller,~
J.Phys.A.,Eq.~1.2~(g~=~d)+T:4''}
\]
}
\end{maplelatex}

\begin{maplelatex}
\mapleinline{inert}{2d}{`  271: `, "Maier Eq. (7.1) with L=n+ Miller, J.Phys.A.,Eq. 1.2 (g =
d)+T:5";}{%
\[
\mathit{\ \ 271:\ }  \,\mbox{``Maier Eq. (7.1)~with~L=n+~Miller,~
J.Phys.A.,Eq.~1.2~(g~=~d)+T:5''}
\]
}
\end{maplelatex}

\begin{maplelatex}
\mapleinline{inert}{2d}{`  272: `, "Maier Eq. (7.1) with L=n+ Miller, J.Phys.A.,Eq. 1.2 (g =
d)+T:6";}{%
\[
\mathit{\ \ 272:\ }  \,\mbox{``Maier Eq. (7.1)~with~L=n+~Miller,~
J.Phys.A.,Eq.~1.2~(g~=~d)+T:6''}
\]
}
\end{maplelatex}

\begin{maplelatex}
\mapleinline{inert}{2d}{`  273: `, "Maier Eq. (7.1) with L=n+ Miller, J.Phys.A.,Eq. 1.2 (g =
d)+T:7";}{%
\[
\mathit{\ \ 273:\ }  \,\mbox{``Maier Eq. (7.1)~with~L=n+~Miller,~
J.Phys.A.,Eq.~1.2~(g~=~d)+T:7''}
\]
}
\end{maplelatex}

\begin{maplelatex}
\mapleinline{inert}{2d}{`  274: `, "Maier Eq. (7.1) with L=n+ Miller, J.Phys.A.,Eq. 1.2 (g =
d)+T:8";}{%
\[
\mathit{\ \ 274:\ }  \,\mbox{``Maier Eq. (7.1)~with~L=n+~Miller,~
J.Phys.A.,Eq.~1.2~(g~=~d)+T:8''}
\]
}
\end{maplelatex}

\begin{maplelatex}
\mapleinline{inert}{2d}{`  275: `, "Maier Eq. (7.1) with L=n+ Miller, J.Phys.A.,Eq. 1.2 (g =
d)+T:9";}{%
\[
\mathit{\ \ 275:\ }  \,\mbox{``Maier Eq. (7.1)~with~L=n+~Miller,~
J.Phys.A.,Eq.~1.2~(g~=~d)+T:9''}
\]
}
\end{maplelatex}

\begin{maplelatex}
\mapleinline{inert}{2d}{`  276: `, "Maier Eq. (7.1) with L=n+ Miller, J.Phys.A.,Eq. 1.2 (g =
d)+T:1";}{%
\[
\mathit{\ \ 276:\ }  \,\mbox{``Maier Eq. (7.1)~with~L=n+~Miller,~
J.Phys.A.,Eq.~1.2~(g~=~d)+T:1''}
\]
}
\end{maplelatex}

\begin{maplelatex}
\mapleinline{inert}{2d}{`  277: `, "Maier Eq. (7.1) with L=n+ Miller, J.Phys.A.,Eq. 1.2 (g =
d)+T:2";}{%
\[
\mathit{\ \ 277:\ }  \,\mbox{``Maier Eq. (7.1)~with~L=n+~Miller,~
J.Phys.A.,Eq.~1.2~(g~=~d)+T:2''}
\]
}
\end{maplelatex}

\begin{maplelatex}
\mapleinline{inert}{2d}{`  278: `, "Maier Eq. (7.1) with L=n+ Miller, J.Phys.A.,Eq. 1.2 (g =
d)+T:3";}{%
\[
\mathit{\ \ 278:\ }  \,\mbox{``Maier Eq. (7.1)~with~L=n+~Miller,~
J.Phys.A.,Eq.~1.2~(g~=~d)+T:3''}
\]
}
\end{maplelatex}

\begin{maplelatex}
\mapleinline{inert}{2d}{`  279: `, "Maier Eq. (7.1) with L=n+ Miller, J.Phys.A.,Eq. 1.2 (g =
d)+T:5";}{%
\[
\mathit{\ \ 279:\ }  \,\mbox{``Maier Eq. (7.1)~with~L=n+~Miller,~
J.Phys.A.,Eq.~1.2~(g~=~d)+T:5''}
\]
}
\end{maplelatex}

\begin{maplelatex}
\mapleinline{inert}{2d}{`  280: `, "Maier Eq. (7.1) with L=n+ Miller, J.Phys.A.,Eq. 1.2 (g =
d)+T:6";}{%
\[
\mathit{\ \ 280:\ }  \,\mbox{``Maier Eq. (7.1)~with~L=n+~Miller,~
J.Phys.A.,Eq.~1.2~(g~=~d)+T:6''}
\]
}
\end{maplelatex}

\begin{maplelatex}
\mapleinline{inert}{2d}{`  281: `, "Maier Eq. (7.1) with L=n+ Miller, J.Phys.A.,Eq. 1.2 (g =
d)+T:7";}{%
\[
\mathit{\ \ 281:\ }  \,\mbox{``Maier Eq. (7.1)~with~L=n+~Miller,~
J.Phys.A.,Eq.~1.2~(g~=~d)+T:7''}
\]
}
\end{maplelatex}

\begin{maplelatex}
\mapleinline{inert}{2d}{`  282: `, "Maier Eq. (7.1) with L=n+ Miller, J.Phys.A.,Eq. 1.2 (g =
d)+T:8";}{%
\[
\mathit{\ \ 282:\ }  \,\mbox{``Maier Eq. (7.1)~with~L=n+~Miller,~
J.Phys.A.,Eq.~1.2~(g~=~d)+T:8''}
\]
}
\end{maplelatex}

\begin{maplelatex}
\mapleinline{inert}{2d}{`  283: `, "Maier Eq. (7.1) with L=n+ Miller, J.Phys.A.,Eq. 1.2 (g =
d)+T:9";}{%
\[
\mathit{\ \ 283:\ }  \,\mbox{``Maier Eq. (7.1)~with~L=n+~Miller,~
J.Phys.A.,Eq.~1.2~(g~=~d)+T:9''}
\]
}
\end{maplelatex}

\begin{maplelatex}
\mapleinline{inert}{2d}{`  284: `, "Maier Eq. (7.1) with L=n+ Miller, J.Phys.A.,Eq. 1.2 (g =
d)+T:1";}{%
\[
\mathit{\ \ 284:\ }  \,\mbox{``Maier Eq. (7.1)~with~L=n+~Miller,~
J.Phys.A.,Eq.~1.2~(g~=~d)+T:1''}
\]
}
\end{maplelatex}

\begin{maplelatex}
\mapleinline{inert}{2d}{`  285: `, "Maier Eq. (7.1) with L=n+ Miller, J.Phys.A.,Eq. 1.2 (g =
d)+T:2";}{%
\[
\mathit{\ \ 285:\ }  \,\mbox{``Maier Eq. (7.1)~with~L=n+~Miller,~
J.Phys.A.,Eq.~1.2~(g~=~d)+T:2''}
\]
}
\end{maplelatex}

\begin{maplelatex}
\mapleinline{inert}{2d}{`  286: `, "Maier Eq. (7.1) with L=n+ Miller, J.Phys.A.,Eq. 1.2 (g =
d)+T:4";}{%
\[
\mathit{\ \ 286:\ }  \,\mbox{``Maier Eq. (7.1)~with~L=n+~Miller,~
J.Phys.A.,Eq.~1.2~(g~=~d)+T:4''}
\]
}
\end{maplelatex}

\begin{maplelatex}
\mapleinline{inert}{2d}{`  287: `, "Maier Eq. (7.1) with L=n+ Miller, J.Phys.A.,Eq. 1.2 (g =
d)+T:5";}{%
\[
\mathit{\ \ 287:\ }  \,\mbox{``Maier Eq. (7.1)~with~L=n+~Miller,~
J.Phys.A.,Eq.~1.2~(g~=~d)+T:5''}
\]
}
\end{maplelatex}

\begin{maplelatex}
\mapleinline{inert}{2d}{`  288: `, "Maier Eq. (7.1) with L=n+ Miller, J.Phys.A.,Eq. 1.2 (g =
d)+T:7";}{%
\[
\mathit{\ \ 288:\ }  \,\mbox{``Maier Eq. (7.1)~with~L=n+~Miller,~
J.Phys.A.,Eq.~1.2~(g~=~d)+T:7''}
\]
}
\end{maplelatex}

\begin{maplelatex}
\mapleinline{inert}{2d}{`  289: `, "Maier Eq. (7.1) with L=n+ Miller, J.Phys.A.,Eq. 1.2 (g =
d)+T:8";}{%
\[
\mathit{\ \ 289:\ }  \,\mbox{``Maier Eq. (7.1)~with~L=n+~Miller,~
J.Phys.A.,Eq.~1.2~(g~=~d)+T:8''}
\]
}
\end{maplelatex}

\begin{maplelatex}
\mapleinline{inert}{2d}{`  290: `, "Maier Eq. (7.1) with L=n+ Miller, J.Phys.A.,Eq. 1.2 (g =
d)+T:9";}{%
\[
\mathit{\ \ 290:\ }  \,\mbox{``Maier Eq. (7.1)~with~L=n+~Miller,~
J.Phys.A.,Eq.~1.2~(g~=~d)+T:9''}
\]
}
\end{maplelatex}

\begin{maplelatex}
\mapleinline{inert}{2d}{`  293: `, "Maier Eq. (7.1) with L=n+ Miller, J.Phys.A.,Eq. 1.2 (g = d)+
Miller, J.Phys.A.,Eq. 1.2 (g = d)";}{%
\maplemultiline{
\mathit{\ \ 293:\ }  \mbox{``Maier Eq. (7.1)~with~L=n+~Miller,~
J.Phys.A.,Eq.~1.2~(g~=~d)+~Miller,~J.Phys.}
\mbox{A.,Eq.~1.2~(g~=~d)''} }
}
\end{maplelatex}

\begin{maplelatex}
\mapleinline{inert}{2d}{`  294: `, "Maier Eq. (7.1) with L=n+ Miller, J.Phys.A.,Eq. 1.2 (g = d)+
Miller, J.Phys.A.,Eq. 1.2 (g = d)";}{%
\maplemultiline{
\mathit{\ \ 294:\ }  \mbox{``Maier Eq. (7.1)~with~L=n+~Miller,~
J.Phys.A.,Eq.~1.2~(g~=~d)+~Miller,~J.Phys.}
\mbox{A.,Eq.~1.2~(g~=~d)''} }
}
\end{maplelatex}

\begin{maplelatex}
\mapleinline{inert}{2d}{`  295: `, "Maier Eq. (7.1) with L=n+ Miller, J.Phys.A.,Eq. 1.2 (g = d)+
Miller, J.Phys.A.,Eq. 1.2 (g = d)";}{%
\maplemultiline{
\mathit{\ \ 295:\ }  \mbox{``Maier Eq. (7.1)~with~L=n+~Miller,~
J.Phys.A.,Eq.~1.2~(g~=~d)+~Miller,~J.Phys.}
\mbox{A.,Eq.~1.2~(g~=~d)''} }
}
\end{maplelatex}

\begin{maplelatex}
\mapleinline{inert}{2d}{`  298: `, "Maier Eq. (7.1) with L=n+ Miller, J.Phys.A.,Eq. 1.2 (g =
d)+T:2+ Miller, J.Phys.A.,Eq. 1.2 (g = d)";}{%
\maplemultiline{
\mathit{\ \ 298:\ }  \mbox{``Maier Eq. (7.1)~with~L=n+~Miller,~
J.Phys.A.,Eq.~1.2~(g~=~d)+T:2+~Miller,~J.P}
\mbox{hys.A.,Eq.~1.2~(g~=~d)''} }
}
\end{maplelatex}

\begin{maplelatex}
\mapleinline{inert}{2d}{`  299: `, "Maier Eq. (7.1) with L=n+ Miller, J.Phys.A.,Eq. 1.2 (g =
d)+T:3+ Miller, J.Phys.A.,Eq. 1.2 (g = d)";}{%
\maplemultiline{
\mathit{\ \ 299:\ }  \mbox{``Maier Eq. (7.1)~with~L=n+~Miller,~
J.Phys.A.,Eq.~1.2~(g~=~d)+T:3+~Miller,~J.P}
\mbox{hys.A.,Eq.~1.2~(g~=~d)''} }
}
\end{maplelatex}

\begin{maplelatex}
\mapleinline{inert}{2d}{`  300: `, "Maier Eq. (7.1) with L=n+ Miller, J.Phys.A.,Eq. 1.2 (g =
d)+T:4+ Miller, J.Phys.A.,Eq. 1.2 (g = d)";}{%
\maplemultiline{
\mathit{\ \ 300:\ }  \mbox{``Maier Eq. (7.1)~with~L=n+~Miller,~
J.Phys.A.,Eq.~1.2~(g~=~d)+T:4+~Miller,~J.P}
\mbox{hys.A.,Eq.~1.2~(g~=~d)''} }
}
\end{maplelatex}

\begin{maplelatex}
\mapleinline{inert}{2d}{`  301: `, "Maier Eq. (7.1) with L=n+ Miller, J.Phys.A.,Eq. 1.2 (g =
d)+T:8+ Miller, J.Phys.A.,Eq. 1.2 (g = d)";}{%
\maplemultiline{
\mathit{\ \ 301:\ }  \mbox{``Maier Eq. (7.1)~with~L=n+~Miller,~
J.Phys.A.,Eq.~1.2~(g~=~d)+T:8+~Miller,~J.P}
\mbox{hys.A.,Eq.~1.2~(g~=~d)''} }
}
\end{maplelatex}

\begin{maplelatex}
\mapleinline{inert}{2d}{`  302: `, "Maier Eq. (7.1) with L=n+ Miller, J.Phys.A.,Eq. 1.2 (g =
d)+T:9+ Miller, J.Phys.A.,Eq. 1.2 (g = d)";}{%
\maplemultiline{
\mathit{\ \ 302:\ }  \mbox{``Maier Eq. (7.1)~with~L=n+~Miller,~
J.Phys.A.,Eq.~1.2~(g~=~d)+T:9+~Miller,~J.P}
\mbox{hys.A.,Eq.~1.2~(g~=~d)''} }
}
\end{maplelatex}

\begin{maplelatex}
\mapleinline{inert}{2d}{`  303: `, "Maier Eq. (7.1) with L=n+ Miller, J.Phys.A.,Eq. 1.2 (g =
d)+T:2+ Miller, J.Phys.A.,Eq. 1.2 (g = d)";}{%
\maplemultiline{
\mathit{\ \ 303:\ }  \mbox{``Maier Eq. (7.1)~with~L=n+~Miller,~
J.Phys.A.,Eq.~1.2~(g~=~d)+T:2+~Miller,~J.P}
\mbox{hys.A.,Eq.~1.2~(g~=~d)''} }
}
\end{maplelatex}

\begin{maplelatex}
\mapleinline{inert}{2d}{`  304: `, "Maier Eq. (7.1) with L=n+ Miller, J.Phys.A.,Eq. 1.2 (g =
d)+T:8+ Miller, J.Phys.A.,Eq. 1.2 (g = d)";}{%
\maplemultiline{
\mathit{\ \ 304:\ }  \mbox{``Maier Eq. (7.1)~with~L=n+~Miller,~
J.Phys.A.,Eq.~1.2~(g~=~d)+T:8+~Miller,~J.P}
\mbox{hys.A.,Eq.~1.2~(g~=~d)''} }
}
\end{maplelatex}

\begin{maplelatex}
\mapleinline{inert}{2d}{`  305: `, "Maier Eq. (7.1) with L=n+ Miller, J.Phys.A.,Eq. 1.2 (g = d)+
Miller, J.Phys.A.,Eq. 1.2 (g = d)+T:1";}{%
\maplemultiline{
\mathit{\ \ 305:\ }  \mbox{``Maier Eq. (7.1)~with~L=n+~Miller,~
J.Phys.A.,Eq.~1.2~(g~=~d)+~Miller,~J.Phys.}
\mbox{A.,Eq.~1.2~(g~=~d)+T:1''} }
}
\end{maplelatex}

\begin{maplelatex}
\mapleinline{inert}{2d}{`  306: `, "Maier Eq. (7.1) with L=n+ Miller, J.Phys.A.,Eq. 1.2 (g = d)+
Miller, J.Phys.A.,Eq. 1.2 (g = d)+T:2";}{%
\maplemultiline{
\mathit{\ \ 306:\ }  \mbox{``Maier Eq. (7.1)~with~L=n+~Miller,~
J.Phys.A.,Eq.~1.2~(g~=~d)+~Miller,~J.Phys.}
\mbox{A.,Eq.~1.2~(g~=~d)+T:2''} }
}
\end{maplelatex}

\begin{maplelatex}
\mapleinline{inert}{2d}{`  307: `, "Maier Eq. (7.1) with L=n+ Miller, J.Phys.A.,Eq. 1.2 (g = d)+
Miller, J.Phys.A.,Eq. 1.2 (g = d)+T:3";}{%
\maplemultiline{
\mathit{\ \ 307:\ }  \mbox{``Maier Eq. (7.1)~with~L=n+~Miller,~
J.Phys.A.,Eq.~1.2~(g~=~d)+~Miller,~J.Phys.}
\mbox{A.,Eq.~1.2~(g~=~d)+T:3''} }
}
\end{maplelatex}

\begin{maplelatex}
\mapleinline{inert}{2d}{`  308: `, "Maier Eq. (7.1) with L=n+ Miller, J.Phys.A.,Eq. 1.2 (g = d)+
Miller, J.Phys.A.,Eq. 1.2 (g = d)+T:4";}{%
\maplemultiline{
\mathit{\ \ 308:\ }  \mbox{``Maier Eq. (7.1)~with~L=n+~Miller,~
J.Phys.A.,Eq.~1.2~(g~=~d)+~Miller,~J.Phys.}
\mbox{A.,Eq.~1.2~(g~=~d)+T:4''} }
}
\end{maplelatex}

\begin{maplelatex}
\mapleinline{inert}{2d}{`  309: `, "Maier Eq. (7.1) with L=n+ Miller, J.Phys.A.,Eq. 1.2 (g = d)+
Miller, J.Phys.A.,Eq. 1.2 (g = d)+T:5";}{%
\maplemultiline{
\mathit{\ \ 309:\ }  \mbox{``Maier Eq. (7.1)~with~L=n+~Miller,~
J.Phys.A.,Eq.~1.2~(g~=~d)+~Miller,~J.Phys.}
\mbox{A.,Eq.~1.2~(g~=~d)+T:5''} }
}
\end{maplelatex}

\begin{maplelatex}
\mapleinline{inert}{2d}{`  310: `, "Maier Eq. (7.1) with L=n+ Miller, J.Phys.A.,Eq. 1.2 (g = d)+
Miller, J.Phys.A.,Eq. 1.2 (g = d)+T:7";}{%
\maplemultiline{
\mathit{\ \ 310:\ }  \mbox{``Maier Eq. (7.1)~with~L=n+~Miller,~
J.Phys.A.,Eq.~1.2~(g~=~d)+~Miller,~J.Phys.}
\mbox{A.,Eq.~1.2~(g~=~d)+T:7''} }
}
\end{maplelatex}

\begin{maplelatex}
\mapleinline{inert}{2d}{`  311: `, "Maier Eq. (7.1) with L=n+ Miller, J.Phys.A.,Eq. 1.2 (g = d)+
Miller, J.Phys.A.,Eq. 1.2 (g = d)+T:9";}{%
\maplemultiline{
\mathit{\ \ 311:\ }  \mbox{``Maier Eq. (7.1)~with~L=n+~Miller,~
J.Phys.A.,Eq.~1.2~(g~=~d)+~Miller,~J.Phys.}
\mbox{A.,Eq.~1.2~(g~=~d)+T:9''} }
}
\end{maplelatex}

\begin{maplelatex}
\mapleinline{inert}{2d}{`  312: `, "Maier Eq. (7.1) with L=n+ Miller, J.Phys.A.,Eq. 1.2 (g = d)+
Miller, J.Phys.A.,Eq. 1.2 (g = d)+T:1";}{%
\maplemultiline{
\mathit{\ \ 312:\ }  \mbox{``Maier Eq. (7.1)~with~L=n+~Miller,~
J.Phys.A.,Eq.~1.2~(g~=~d)+~Miller,~J.Phys.}
\mbox{A.,Eq.~1.2~(g~=~d)+T:1''} }
}
\end{maplelatex}

\begin{maplelatex}
\mapleinline{inert}{2d}{`  313: `, "Maier Eq. (7.1) with L=n+ Miller, J.Phys.A.,Eq. 1.2 (g = d)+
Miller, J.Phys.A.,Eq. 1.2 (g = d)+T:2";}{%
\maplemultiline{
\mathit{\ \ 313:\ }  \mbox{``Maier Eq. (7.1)~with~L=n+~Miller,~
J.Phys.A.,Eq.~1.2~(g~=~d)+~Miller,~J.Phys.}
\mbox{A.,Eq.~1.2~(g~=~d)+T:2''} }
}
\end{maplelatex}

\begin{maplelatex}
\mapleinline{inert}{2d}{`  314: `, "Maier Eq. (7.1) with L=n+ Miller, J.Phys.A.,Eq. 1.2 (g = d)+
Miller, J.Phys.A.,Eq. 1.2 (g = d)+T:4";}{%
\maplemultiline{
\mathit{\ \ 314:\ }  \mbox{``Maier Eq. (7.1)~with~L=n+~Miller,~
J.Phys.A.,Eq.~1.2~(g~=~d)+~Miller,~J.Phys.}
\mbox{A.,Eq.~1.2~(g~=~d)+T:4''} }
}
\end{maplelatex}

\begin{maplelatex}
\mapleinline{inert}{2d}{`  315: `, "Maier Eq. (7.1) with L=n+ Miller, J.Phys.A.,Eq. 1.2 (g = d)+
Miller, J.Phys.A.,Eq. 1.2 (g = d)+T:5";}{%
\maplemultiline{
\mathit{\ \ 315:\ }  \mbox{``Maier Eq. (7.1)~with~L=n+~Miller,~
J.Phys.A.,Eq.~1.2~(g~=~d)+~Miller,~J.Phys.}
\mbox{A.,Eq.~1.2~(g~=~d)+T:5''} }
}
\end{maplelatex}

\begin{maplelatex}
\mapleinline{inert}{2d}{`  316: `, "Maier Eq. (7.1) with L=n+ Miller, J.Phys.A.,Eq. 1.2 (g = d)+
Miller, J.Phys.A.,Eq. 1.2 (g = d)+T:6";}{%
\maplemultiline{
\mathit{\ \ 316:\ }  \mbox{``Maier Eq. (7.1)~with~L=n+~Miller,~
J.Phys.A.,Eq.~1.2~(g~=~d)+~Miller,~J.Phys.}
\mbox{A.,Eq.~1.2~(g~=~d)+T:6''} }
}
\end{maplelatex}

\begin{maplelatex}
\mapleinline{inert}{2d}{`  317: `, "Maier Eq. (7.1) with L=n+ Miller, J.Phys.A.,Eq. 1.2 (g = d)+
Miller, J.Phys.A.,Eq. 1.2 (g = d)+T:7";}{%
\maplemultiline{
\mathit{\ \ 317:\ }  \mbox{``Maier Eq. (7.1)~with~L=n+~Miller,~
J.Phys.A.,Eq.~1.2~(g~=~d)+~Miller,~J.Phys.}
\mbox{A.,Eq.~1.2~(g~=~d)+T:7''} }
}
\end{maplelatex}

\begin{maplelatex}
\mapleinline{inert}{2d}{`  318: `, "Maier Eq. (7.1) with L=n+ Miller, J.Phys.A.,Eq. 1.2 (g = d)+
Miller, J.Phys.A.,Eq. 1.2 (g = d)+T:9";}{%
\maplemultiline{
\mathit{\ \ 318:\ }  \mbox{``Maier Eq. (7.1)~with~L=n+~Miller,~
J.Phys.A.,Eq.~1.2~(g~=~d)+~Miller,~J.Phys.}
\mbox{A.,Eq.~1.2~(g~=~d)+T:9''} }
}
\end{maplelatex}

\begin{maplelatex}
\mapleinline{inert}{2d}{`  319: `, "Maier Eq. (7.1) with L=n+ Miller, J.Phys.A.,Eq. 1.2 (g =
d)+T:2+ Miller, J.Phys.A.,Eq. 1.2 (g = d)+T:1";}{%
\maplemultiline{
\mathit{\ \ 319:\ }  \mbox{``Maier Eq. (7.1)~with~L=n+~Miller,~
J.Phys.A.,Eq.~1.2~(g~=~d)+T:2+~Miller,~J.P}
\mbox{hys.A.,Eq.~1.2~(g~=~d)+T:1''} }
}
\end{maplelatex}

\begin{maplelatex}
\mapleinline{inert}{2d}{`  320: `, "Maier Eq. (7.1) with L=n+ Miller, J.Phys.A.,Eq. 1.2 (g =
d)+T:2+ Miller, J.Phys.A.,Eq. 1.2 (g = d)+T:2";}{%
\maplemultiline{
\mathit{\ \ 320:\ }  \mbox{``Maier Eq. (7.1)~with~L=n+~Miller,~
J.Phys.A.,Eq.~1.2~(g~=~d)+T:2+~Miller,~J.P}
\mbox{hys.A.,Eq.~1.2~(g~=~d)+T:2''} }
}
\end{maplelatex}

\begin{maplelatex}
\mapleinline{inert}{2d}{`  321: `, "Maier Eq. (7.1) with L=n+ Miller, J.Phys.A.,Eq. 1.2 (g =
d)+T:2+ Miller, J.Phys.A.,Eq. 1.2 (g = d)+T:3";}{%
\maplemultiline{
\mathit{\ \ 321:\ }  \mbox{``Maier Eq. (7.1)~with~L=n+~Miller,~
J.Phys.A.,Eq.~1.2~(g~=~d)+T:2+~Miller,~J.P}
\mbox{hys.A.,Eq.~1.2~(g~=~d)+T:3''} }
}
\end{maplelatex}

\begin{maplelatex}
\mapleinline{inert}{2d}{`  322: `, "Maier Eq. (7.1) with L=n+ Miller, J.Phys.A.,Eq. 1.2 (g =
d)+T:2+ Miller, J.Phys.A.,Eq. 1.2 (g = d)+T:4";}{%
\maplemultiline{
\mathit{\ \ 322:\ }  \mbox{``Maier Eq. (7.1)~with~L=n+~Miller,~
J.Phys.A.,Eq.~1.2~(g~=~d)+T:2+~Miller,~J.P}
\mbox{hys.A.,Eq.~1.2~(g~=~d)+T:4''} }
}
\end{maplelatex}

\begin{maplelatex}
\mapleinline{inert}{2d}{`  323: `, "Maier Eq. (7.1) with L=n+ Miller, J.Phys.A.,Eq. 1.2 (g =
d)+T:2+ Miller, J.Phys.A.,Eq. 1.2 (g = d)+T:5";}{%
\maplemultiline{
\mathit{\ \ 323:\ }  \mbox{``Maier Eq. (7.1)~with~L=n+~Miller,~
J.Phys.A.,Eq.~1.2~(g~=~d)+T:2+~Miller,~J.P}
\mbox{hys.A.,Eq.~1.2~(g~=~d)+T:5''} }
}
\end{maplelatex}

\begin{maplelatex}
\mapleinline{inert}{2d}{`  324: `, "Maier Eq. (7.1) with L=n+ Miller, J.Phys.A.,Eq. 1.2 (g =
d)+T:2+ Miller, J.Phys.A.,Eq. 1.2 (g = d)+T:6";}{%
\maplemultiline{
\mathit{\ \ 324:\ }  \mbox{``Maier Eq. (7.1)~with~L=n+~Miller,~
J.Phys.A.,Eq.~1.2~(g~=~d)+T:2+~Miller,~J.P}
\mbox{hys.A.,Eq.~1.2~(g~=~d)+T:6''} }
}
\end{maplelatex}

\begin{maplelatex}
\mapleinline{inert}{2d}{`  325: `, "Maier Eq. (7.1) with L=n+ Miller, J.Phys.A.,Eq. 1.2 (g =
d)+T:2+ Miller, J.Phys.A.,Eq. 1.2 (g = d)+T:7";}{%
\maplemultiline{
\mathit{\ \ 325:\ }  \mbox{``Maier Eq. (7.1)~with~L=n+~Miller,~
J.Phys.A.,Eq.~1.2~(g~=~d)+T:2+~Miller,~J.P}
\mbox{hys.A.,Eq.~1.2~(g~=~d)+T:7''} }
}
\end{maplelatex}

\begin{maplelatex}
\mapleinline{inert}{2d}{`  326: `, "Maier Eq. (7.1) with L=n+ Miller, J.Phys.A.,Eq. 1.2 (g =
d)+T:2+ Miller, J.Phys.A.,Eq. 1.2 (g = d)+T:8";}{%
\maplemultiline{
\mathit{\ \ 326:\ }  \mbox{``Maier Eq. (7.1)~with~L=n+~Miller,~
J.Phys.A.,Eq.~1.2~(g~=~d)+T:2+~Miller,~J.P}
\mbox{hys.A.,Eq.~1.2~(g~=~d)+T:8''} }
}
\end{maplelatex}

\begin{maplelatex}
\mapleinline{inert}{2d}{`  327: `, "Maier Eq. (7.1) with L=n+ Miller, J.Phys.A.,Eq. 1.2 (g =
d)+T:2+ Miller, J.Phys.A.,Eq. 1.2 (g = d)+T:9";}{%
\maplemultiline{
\mathit{\ \ 327:\ }  \mbox{``Maier Eq. (7.1)~with~L=n+~Miller,~
J.Phys.A.,Eq.~1.2~(g~=~d)+T:2+~Miller,~J.P}
\mbox{hys.A.,Eq.~1.2~(g~=~d)+T:9''} }
}
\end{maplelatex}

\begin{maplelatex}
\mapleinline{inert}{2d}{`  328: `, "Maier Eq. (7.1) with L=n+ Miller, J.Phys.A.,Eq. 1.2 (g =
d)+T:3+ Miller, J.Phys.A.,Eq. 1.2 (g = d)+T:1";}{%
\maplemultiline{
\mathit{\ \ 328:\ }  \mbox{``Maier Eq. (7.1)~with~L=n+~Miller,~
J.Phys.A.,Eq.~1.2~(g~=~d)+T:3+~Miller,~J.P}
\mbox{hys.A.,Eq.~1.2~(g~=~d)+T:1''} }
}
\end{maplelatex}

\begin{maplelatex}
\mapleinline{inert}{2d}{`  329: `, "Maier Eq. (7.1) with L=n+ Miller, J.Phys.A.,Eq. 1.2 (g =
d)+T:3+ Miller, J.Phys.A.,Eq. 1.2 (g = d)+T:2";}{%
\maplemultiline{
\mathit{\ \ 329:\ }  \mbox{``Maier Eq. (7.1)~with~L=n+~Miller,~
J.Phys.A.,Eq.~1.2~(g~=~d)+T:3+~Miller,~J.P}
\mbox{hys.A.,Eq.~1.2~(g~=~d)+T:2''} }
}
\end{maplelatex}

\begin{maplelatex}
\mapleinline{inert}{2d}{`  330: `, "Maier Eq. (7.1) with L=n+ Miller, J.Phys.A.,Eq. 1.2 (g =
d)+T:3+ Miller, J.Phys.A.,Eq. 1.2 (g = d)+T:3";}{%
\maplemultiline{
\mathit{\ \ 330:\ }  \mbox{``Maier Eq. (7.1)~with~L=n+~Miller,~
J.Phys.A.,Eq.~1.2~(g~=~d)+T:3+~Miller,~J.P}
\mbox{hys.A.,Eq.~1.2~(g~=~d)+T:3''} }
}
\end{maplelatex}

\begin{maplelatex}
\mapleinline{inert}{2d}{`  331: `, "Maier Eq. (7.1) with L=n+ Miller, J.Phys.A.,Eq. 1.2 (g =
d)+T:3+ Miller, J.Phys.A.,Eq. 1.2 (g = d)+T:4";}{%
\maplemultiline{
\mathit{\ \ 331:\ }  \mbox{``Maier Eq. (7.1)~with~L=n+~Miller,~
J.Phys.A.,Eq.~1.2~(g~=~d)+T:3+~Miller,~J.P}
\mbox{hys.A.,Eq.~1.2~(g~=~d)+T:4''} }
}
\end{maplelatex}

\begin{maplelatex}
\mapleinline{inert}{2d}{`  332: `, "Maier Eq. (7.1) with L=n+ Miller, J.Phys.A.,Eq. 1.2 (g =
d)+T:3+ Miller, J.Phys.A.,Eq. 1.2 (g = d)+T:5";}{%
\maplemultiline{
\mathit{\ \ 332:\ }  \mbox{``Maier Eq. (7.1)~with~L=n+~Miller,~
J.Phys.A.,Eq.~1.2~(g~=~d)+T:3+~Miller,~J.P}
\mbox{hys.A.,Eq.~1.2~(g~=~d)+T:5''} }
}
\end{maplelatex}

\begin{maplelatex}
\mapleinline{inert}{2d}{`  333: `, "Maier Eq. (7.1) with L=n+ Miller, J.Phys.A.,Eq. 1.2 (g =
d)+T:3+ Miller, J.Phys.A.,Eq. 1.2 (g = d)+T:6";}{%
\maplemultiline{
\mathit{\ \ 333:\ }  \mbox{``Maier Eq. (7.1)~with~L=n+~Miller,~
J.Phys.A.,Eq.~1.2~(g~=~d)+T:3+~Miller,~J.P}
\mbox{hys.A.,Eq.~1.2~(g~=~d)+T:6''} }
}
\end{maplelatex}

\begin{maplelatex}
\mapleinline{inert}{2d}{`  334: `, "Maier Eq. (7.1) with L=n+ Miller, J.Phys.A.,Eq. 1.2 (g =
d)+T:3+ Miller, J.Phys.A.,Eq. 1.2 (g = d)+T:7";}{%
\maplemultiline{
\mathit{\ \ 334:\ }  \mbox{``Maier Eq. (7.1)~with~L=n+~Miller,~
J.Phys.A.,Eq.~1.2~(g~=~d)+T:3+~Miller,~J.P}
\mbox{hys.A.,Eq.~1.2~(g~=~d)+T:7''} }
}
\end{maplelatex}

\begin{maplelatex}
\mapleinline{inert}{2d}{`  335: `, "Maier Eq. (7.1) with L=n+ Miller, J.Phys.A.,Eq. 1.2 (g =
d)+T:3+ Miller, J.Phys.A.,Eq. 1.2 (g = d)+T:9";}{%
\maplemultiline{
\mathit{\ \ 335:\ }  \mbox{``Maier Eq. (7.1)~with~L=n+~Miller,~
J.Phys.A.,Eq.~1.2~(g~=~d)+T:3+~Miller,~J.P}
\mbox{hys.A.,Eq.~1.2~(g~=~d)+T:9''} }
}
\end{maplelatex}

\begin{maplelatex}
\mapleinline{inert}{2d}{`  336: `, "Maier Eq. (7.1) with L=n+ Miller, J.Phys.A.,Eq. 1.2 (g =
d)+T:2+ Miller, J.Phys.A.,Eq. 1.2 (g = d)+T:1";}{%
\maplemultiline{
\mathit{\ \ 336:\ }  \mbox{``Maier Eq. (7.1)~with~L=n+~Miller,~
J.Phys.A.,Eq.~1.2~(g~=~d)+T:2+~Miller,~J.P}
\mbox{hys.A.,Eq.~1.2~(g~=~d)+T:1''} }
}
\end{maplelatex}

\begin{maplelatex}
\mapleinline{inert}{2d}{`  337: `, "Maier Eq. (7.1) with L=n+ Miller, J.Phys.A.,Eq. 1.2 (g =
d)+T:2+ Miller, J.Phys.A.,Eq. 1.2 (g = d)+T:2";}{%
\maplemultiline{
\mathit{\ \ 337:\ }  \mbox{``Maier Eq. (7.1)~with~L=n+~Miller,~
J.Phys.A.,Eq.~1.2~(g~=~d)+T:2+~Miller,~J.P}
\mbox{hys.A.,Eq.~1.2~(g~=~d)+T:2''} }
}
\end{maplelatex}

\begin{maplelatex}
\mapleinline{inert}{2d}{`  338: `, "Maier Eq. (7.1) with L=n+ Miller, J.Phys.A.,Eq. 1.2 (g =
d)+T:2+ Miller, J.Phys.A.,Eq. 1.2 (g = d)+T:3";}{%
\maplemultiline{
\mathit{\ \ 338:\ }  \mbox{``Maier Eq. (7.1)~with~L=n+~Miller,~
J.Phys.A.,Eq.~1.2~(g~=~d)+T:2+~Miller,~J.P}
\mbox{hys.A.,Eq.~1.2~(g~=~d)+T:3''} }
}
\end{maplelatex}

\begin{maplelatex}
\mapleinline{inert}{2d}{`  339: `, "Maier Eq. (7.1) with L=n+ Miller, J.Phys.A.,Eq. 1.2 (g =
d)+T:2+ Miller, J.Phys.A.,Eq. 1.2 (g = d)+T:4";}{%
\maplemultiline{
\mathit{\ \ 339:\ }  \mbox{``Maier Eq. (7.1)~with~L=n+~Miller,~
J.Phys.A.,Eq.~1.2~(g~=~d)+T:2+~Miller,~J.P}
\mbox{hys.A.,Eq.~1.2~(g~=~d)+T:4''} }
}
\end{maplelatex}

\begin{maplelatex}
\mapleinline{inert}{2d}{`  340: `, "Maier Eq. (7.1) with L=n+ Miller, J.Phys.A.,Eq. 1.2 (g =
d)+T:2+ Miller, J.Phys.A.,Eq. 1.2 (g = d)+T:5";}{%
\maplemultiline{
\mathit{\ \ 340:\ }  \mbox{``Maier Eq. (7.1)~with~L=n+~Miller,~
J.Phys.A.,Eq.~1.2~(g~=~d)+T:2+~Miller,~J.P}
\mbox{hys.A.,Eq.~1.2~(g~=~d)+T:5''} }
}
\end{maplelatex}

\begin{maplelatex}
\mapleinline{inert}{2d}{`  341: `, "Maier Eq. (7.1) with L=n+ Miller, J.Phys.A.,Eq. 1.2 (g =
d)+T:2+ Miller, J.Phys.A.,Eq. 1.2 (g = d)+T:6";}{%
\maplemultiline{
\mathit{\ \ 341:\ }  \mbox{``Maier Eq. (7.1)~with~L=n+~Miller,~
J.Phys.A.,Eq.~1.2~(g~=~d)+T:2+~Miller,~J.P}
\mbox{hys.A.,Eq.~1.2~(g~=~d)+T:6''} }
}
\end{maplelatex}

\begin{maplelatex}
\mapleinline{inert}{2d}{`  342: `, "Maier Eq. (7.1) with L=n+ Miller, J.Phys.A.,Eq. 1.2 (g =
d)+T:2+ Miller, J.Phys.A.,Eq. 1.2 (g = d)+T:7";}{%
\maplemultiline{
\mathit{\ \ 342:\ }  \mbox{``Maier Eq. (7.1)~with~L=n+~Miller,~
J.Phys.A.,Eq.~1.2~(g~=~d)+T:2+~Miller,~J.P}
\mbox{hys.A.,Eq.~1.2~(g~=~d)+T:7''} }
}
\end{maplelatex}

\begin{maplelatex}
\mapleinline{inert}{2d}{`  343: `, "Maier Eq. (7.1) with L=n+ Miller, J.Phys.A.,Eq. 1.2 (g =
d)+T:2+ Miller, J.Phys.A.,Eq. 1.2 (g = d)+T:8";}{%
\maplemultiline{
\mathit{\ \ 343:\ }  \mbox{``Maier Eq. (7.1)~with~L=n+~Miller,~
J.Phys.A.,Eq.~1.2~(g~=~d)+T:2+~Miller,~J.P}
\mbox{hys.A.,Eq.~1.2~(g~=~d)+T:8''} }
}
\end{maplelatex}

\begin{maplelatex}
\mapleinline{inert}{2d}{`  344: `, "Maier Eq. (7.1) with L=n+ Miller, J.Phys.A.,Eq. 1.2 (g =
d)+T:2+ Miller, J.Phys.A.,Eq. 1.2 (g = d)+T:9";}{%
\maplemultiline{
\mathit{\ \ 344:\ }  \mbox{``Maier Eq. (7.1)~with~L=n+~Miller,~
J.Phys.A.,Eq.~1.2~(g~=~d)+T:2+~Miller,~J.Phys.A.,} \\
\mbox{Eq.~1.2~(g~=~d)+T:9''} }
}
\end{maplelatex}

\begin{maplelatex}
\mapleinline{inert}{2d}{`  345: `, "Maier Eq. (7.1) with L=n+ Miller, J.Phys.A.,Eq. 1.2 (g = d)+
Miller, J.Phys.A.,Eq. 1.2 (g = d)+ Miller, J.Phys.A.,Eq. 1.2 (g =
d)";}{%
\maplemultiline{
\mathit{\ \ 345:\ }  \mbox{``Maier Eq. (7.1)~with~L=n+~Miller,~
J.Phys.A.,Eq.~1.2~(g~=~d)+~Miller,~J.Phys.A.,} \\
\mbox{Eq.~1.2~(g~=~d)+~Miller,~J.Phys.A.,Eq.~1.2~(g~=~d)''} }
}
\end{maplelatex}

\begin{maplelatex}
\mapleinline{inert}{2d}{`  346: `, "Maier Eq. (7.1) with L=n+ Miller, J.Phys.A.,Eq. 1.2 (g = d)+
Miller, J.Phys.A.,Eq. 1.2 (g = d)+ Miller, J.Phys.A.,Eq. 1.2 (g =
d)";}{%
\maplemultiline{
\mathit{\ \ 346:\ }  \mbox{``Maier Eq. (7.1)~with~L=n+~Miller,~
J.Phys.A.,Eq.~1.2~(g~=~d)+~Miller,~J.Phys.A.,} \\
\mbox{Eq.~1.2~(g~=~d)+~Miller,~J.Phys.A.,Eq.~1.2~(g~=~d)''} }
}
\end{maplelatex}

\begin{maplelatex}
\mapleinline{inert}{2d}{`  347: `, "Maier Eq. (7.1) with L=n+ Miller, J.Phys.A.,Eq. 1.2 (g = d)+
Miller, J.Phys.A.,Eq. 1.2 (g = d)+T:1+ Miller, J.Phys.A.,Eq. 1.2 (g =
d)";}{%
\maplemultiline{
\mathit{\ \ 347:\ }  \mbox{``Maier Eq. (7.1)~with~L=n+~Miller,~
J.Phys.A.,Eq.~1.2~(g~=~d)+~Miller,~J.Phys.A.,} \\
\mbox{Eq.~1.2~(g~=~d)+T:1+~Miller,~J.Phys.A.,Eq.~1.2~(g~=~d)''
} }
}
\end{maplelatex}

\begin{maplelatex}
\mapleinline{inert}{2d}{`  348: `, "Maier Eq. (7.1) with L=n+ Miller, J.Phys.A.,Eq. 1.2 (g = d)+
Miller, J.Phys.A.,Eq. 1.2 (g = d)+T:3+ Miller, J.Phys.A.,Eq. 1.2 (g =
d)";}{%
\maplemultiline{
\mathit{\ \ 348:\ }  \mbox{``Maier Eq. (7.1)~with~L=n+~Miller,~
J.Phys.A.,Eq.~1.2~(g~=~d)+~Miller,~J.Phys.A., } \\
\mbox{Eq.~1.2~(g~=~d)+T:3+~Miller,~J.Phys.A.,Eq.~1.2~(g~=~d)''
} }
}
\end{maplelatex}

\begin{maplelatex}
\mapleinline{inert}{2d}{`  349: `, "Maier Eq. (7.1) with L=n+ Miller, J.Phys.A.,Eq. 1.2 (g = d)+
Miller, J.Phys.A.,Eq. 1.2 (g = d)+T:4+ Miller, J.Phys.A.,Eq. 1.2 (g =
d)";}{%
\maplemultiline{
\mathit{\ \ 349:\ }  \mbox{``Maier Eq. (7.1)~with~L=n+~Miller,~
J.Phys.A.,Eq.~1.2~(g~=~d)+~Miller,~J.Phys.A., } \\
\mbox{Eq.~1.2~(g~=~d)+T:4+~Miller,~J.Phys.A.,Eq.~1.2~(g~=~d)''
} }
}
\end{maplelatex}

\begin{maplelatex}
\mapleinline{inert}{2d}{`  350: `, "Maier Eq. (7.1) with L=n+ Miller, J.Phys.A.,Eq. 1.2 (g = d)+
Miller, J.Phys.A.,Eq. 1.2 (g = d)+T:5+ Miller, J.Phys.A.,Eq. 1.2 (g =
d)";}{%
\maplemultiline{
\mathit{\ \ 350:\ }  \mbox{``Maier Eq. (7.1)~with~L=n+~Miller,~
J.Phys.A.,Eq.~1.2~(g~=~d)+~Miller,~J.Phys.A., } \\
\mbox{Eq.~1.2~(g~=~d)+T:5+~Miller,~J.Phys.A.,Eq.~1.2~(g~=~d)''
} }
}
\end{maplelatex}

\begin{maplelatex}
\mapleinline{inert}{2d}{`  351: `, "Maier Eq. (7.1) with L=n+ Miller, J.Phys.A.,Eq. 1.2 (g = d)+
Miller, J.Phys.A.,Eq. 1.2 (g = d)+T:9+ Miller, J.Phys.A.,Eq. 1.2 (g =
d)";}{%
\maplemultiline{
\mathit{\ \ 351:\ }  \mbox{``Maier Eq. (7.1)~with~L=n+~Miller,~
J.Phys.A.,Eq.~1.2~(g~=~d)+~Miller,~J.Phys.A., } \\
\mbox{Eq.~1.2~(g~=~d)+T:9+~Miller,~J.Phys.A.,Eq.~1.2~(g~=~d)''
} }
}
\end{maplelatex}

\begin{maplelatex}
\mapleinline{inert}{2d}{`  352: `, "Maier Eq. (7.1) with L=n+ Miller, J.Phys.A.,Eq. 1.2 (g = d)+
Miller, J.Phys.A.,Eq. 1.2 (g = d)+T:2+ Miller, J.Phys.A.,Eq. 1.2 (g =
d)";}{%
\maplemultiline{
\mathit{\ \ 352:\ }  \mbox{``Maier Eq. (7.1)~with~L=n+~Miller,~
J.Phys.A.,Eq.~1.2~(g~=~d)+~Miller,~J.Phys.A.,} \\
\mbox{Eq.~1.2~(g~=~d)+T:2+~Miller,~J.Phys.A.,Eq.~1.2~(g~=~d)''
} }
}
\end{maplelatex}

\begin{maplelatex}
\mapleinline{inert}{2d}{`  353: `, "Maier Eq. (7.1) with L=n+ Miller, J.Phys.A.,Eq. 1.2 (g = d)+
Miller, J.Phys.A.,Eq. 1.2 (g = d)+T:9+ Miller, J.Phys.A.,Eq. 1.2 (g =
d)";}{%
\maplemultiline{
\mathit{\ \ 353:\ }  \mbox{``Maier Eq. (7.1)~with~L=n+~Miller,~
J.Phys.A.,Eq.~1.2~(g~=~d)+~Miller,~J.Phys.A., } \\
\mbox{Eq.~1.2~(g~=~d)+T:9+~Miller,~J.Phys.A.,Eq.~1.2~(g~=~d)''
} }
}
\end{maplelatex}

\begin{maplelatex}
\mapleinline{inert}{2d}{`  354: `, "Maier Eq. (7.1) with L=n+ Miller, J.Phys.A.,Eq. 1.2 (g =
d)+T:2+ Miller, J.Phys.A.,Eq. 1.2 (g = d)+T:2+ Miller, J.Phys.A.,Eq.
1.2 (g = d)";}{%
\maplemultiline{
\mathit{\ \ 354:\ }  \mbox{``Maier Eq. (7.1)~with~L=n+~Miller,~
J.Phys.A.,Eq.~1.2~(g~=~d)+T:2+~Miller,~J.Phys.A., } \\
\mbox{Eq.~1.2~(g~=~d)+T:2+~Miller,~J.Phys.A.,Eq.~1.2~(g~=~
d)''} }
}
\end{maplelatex}

\begin{maplelatex}
\mapleinline{inert}{2d}{`  355: `, "Maier Eq. (7.1) with L=n+ Miller, J.Phys.A.,Eq. 1.2 (g =
d)+T:3+ Miller, J.Phys.A.,Eq. 1.2 (g = d)+T:1+ Miller, J.Phys.A.,Eq.
1.2 (g = d)";}{%
\maplemultiline{
\mathit{\ \ 355:\ }  \mbox{``Maier Eq. (7.1)~with~L=n+~Miller,~
J.Phys.A.,Eq.~1.2~(g~=~d)+T:3+~Miller,~J.Phys.A., } \\
\mbox{Eq.~1.2~(g~=~d)+T:1+~Miller,~J.Phys.A.,Eq.~1.2~(g~=~
d)''} }
}
\end{maplelatex}

\begin{maplelatex}
\mapleinline{inert}{2d}{`  356: `, "Maier Eq. (7.1) with L=n+ Miller, J.Phys.A.,Eq. 1.2 (g =
d)+T:3+ Miller, J.Phys.A.,Eq. 1.2 (g = d)+T:2+ Miller, J.Phys.A.,Eq.
1.2 (g = d)";}{%
\maplemultiline{
\mathit{\ \ 356:\ }  \mbox{``Maier Eq. (7.1)~with~L=n+~Miller,~
J.Phys.A.,Eq.~1.2~(g~=~d)+T:3+~Miller,~J.Phys.A., } \\
\mbox{Eq.~1.2~(g~=~d)+T:2+~Miller,~J.Phys.A.,Eq.~1.2~(g~=~
d)''} }
}
\end{maplelatex}

\begin{maplelatex}
\mapleinline{inert}{2d}{`  357: `, "Maier Eq. (7.1) with L=n+ Miller, J.Phys.A.,Eq. 1.2 (g =
d)+T:3+ Miller, J.Phys.A.,Eq. 1.2 (g = d)+T:5+ Miller, J.Phys.A.,Eq.
1.2 (g = d)";}{%
\maplemultiline{
\mathit{\ \ 357:\ }  \mbox{``Maier Eq. (7.1)~with~L=n+~Miller,~
J.Phys.A.,Eq.~1.2~(g~=~d)+T:3+~Miller,~J.Phys.A., } \\
\mbox{Eq.~1.2~(g~=~d)+T:5+~Miller,~J.Phys.A.,Eq.~1.2~(g~=~
d)''} }
}
\end{maplelatex}

\begin{maplelatex}
\mapleinline{inert}{2d}{`  358: `, "Maier Eq. (7.1) with L=n+ Miller, J.Phys.A.,Eq. 1.2 (g = d)+
Miller, J.Phys.A.,Eq. 1.2 (g = d)+ Miller, J.Phys.A.,Eq. 1.2 (g =
d)+T:2";}{%
\maplemultiline{
\mathit{\ \ 358:\ }  \mbox{``Maier Eq. (7.1)~with~L=n+~Miller,~
J.Phys.A.,Eq.~1.2~(g~=~d)+~Miller,~J.Phys.A., }  \\
\mbox{Eq.~1.2~(g~=~d)+~Miller,~J.Phys.A.,Eq.~1.2~(g~=~d)+T:2''
} }
}
\end{maplelatex}

\begin{maplelatex}
\mapleinline{inert}{2d}{`  359: `, "Maier Eq. (7.1) with L=n+ Miller, J.Phys.A.,Eq. 1.2 (g = d)+
Miller, J.Phys.A.,Eq. 1.2 (g = d)+ Miller, J.Phys.A.,Eq. 1.2 (g =
d)+T:4";}{%
\maplemultiline{
\mathit{\ \ 359:\ }  \mbox{``Maier Eq. (7.1)~with~L=n+~Miller,~
J.Phys.A.,Eq.~1.2~(g~=~d)+~Miller,~J.Phys.A., }  \\
\mbox{Eq.~1.2~(g~=~d)+~Miller,~J.Phys.A.,Eq.~1.2~(g~=~d)+T:4''
} }
}
\end{maplelatex}

\begin{maplelatex}
\mapleinline{inert}{2d}{`  360: `, "Maier Eq. (7.1) with L=n+ Miller, J.Phys.A.,Eq. 1.2 (g = d)+
Miller, J.Phys.A.,Eq. 1.2 (g = d)+ Miller, J.Phys.A.,Eq. 1.2 (g =
d)+T:9";}{%
\maplemultiline{
\mathit{\ \ 360:\ }  \mbox{``Maier Eq. (7.1)~with~L=n+~Miller,~
J.Phys.A.,Eq.~1.2~(g~=~d)+~Miller,~J.Phys.A., }  \\
\mbox{Eq.~1.2~(g~=~d)+~Miller,~J.Phys.A.,Eq.~1.2~(g~=~d)+T:9''
} }
}
\end{maplelatex}

\begin{maplelatex}
\mapleinline{inert}{2d}{`  361: `, "Maier Eq. (7.1) with L=n+ Miller, J.Phys.A.,Eq. 1.2 (g = d)+
Miller, J.Phys.A.,Eq. 1.2 (g = d)+T:1+ Miller, J.Phys.A.,Eq. 1.2 (g =
d)+T:1";}{%
\maplemultiline{
\mathit{\ \ 361:\ }  \mbox{``Maier Eq. (7.1)~with~L=n+~Miller,~
J.Phys.A.,Eq.~1.2~(g~=~d)+~Miller,~J.Phys.A., }  \\
\mbox{Eq.~1.2~(g~=~d)+T:1+~Miller,~J.Phys.A.,Eq.~1.2~(g~=~
d)+T:1''} }
}
\end{maplelatex}

\begin{maplelatex}
\mapleinline{inert}{2d}{`  362: `, "Maier Eq. (7.1) with L=n+ Miller, J.Phys.A.,Eq. 1.2 (g = d)+
Miller, J.Phys.A.,Eq. 1.2 (g = d)+T:1+ Miller, J.Phys.A.,Eq. 1.2 (g =
d)+T:4";}{%
\maplemultiline{
\mathit{\ \ 362:\ }  \mbox{``Maier Eq. (7.1)~with~L=n+~Miller,~
J.Phys.A.,Eq.~1.2~(g~=~d)+~Miller,~J.Phys.A., }  \\
\mbox{Eq.~1.2~(g~=~d)+T:1+~Miller,~J.Phys.A.,Eq.~1.2~(g~=~
d)+T:4''} }
}
\end{maplelatex}

\begin{maplelatex}
\mapleinline{inert}{2d}{`  363: `, "Maier Eq. (7.1) with L=n+ Miller, J.Phys.A.,Eq. 1.2 (g = d)+
Miller, J.Phys.A.,Eq. 1.2 (g = d)+T:1+ Miller, J.Phys.A.,Eq. 1.2 (g =
d)+T:7";}{%
\maplemultiline{
\mathit{\ \ 363:\ }  \mbox{``Maier Eq. (7.1)~with~L=n+~Miller,~
J.Phys.A.,Eq.~1.2~(g~=~d)+~Miller,~J.Phys.A., }  \\
\mbox{Eq.~1.2~(g~=~d)+T:1+~Miller,~J.Phys.A.,Eq.~1.2~(g~=~
d)+T:7''} }
}
\end{maplelatex}

\begin{maplelatex}
\mapleinline{inert}{2d}{`  364: `, "Maier Eq. (7.1) with L=n+ Miller, J.Phys.A.,Eq. 1.2 (g = d)+
Miller, J.Phys.A.,Eq. 1.2 (g = d)+T:1+ Miller, J.Phys.A.,Eq. 1.2 (g =
d)+T:8";}{%
\maplemultiline{
\mathit{\ \ 364:\ }  \mbox{``Maier Eq. (7.1)~with~L=n+~Miller,~
J.Phys.A.,Eq.~1.2~(g~=~d)+~Miller,~J.Phys.A., }  \\
\mbox{Eq.~1.2~(g~=~d)+T:1+~Miller,~J.Phys.A.,Eq.~1.2~(g~=~
d)+T:8''} }
}
\end{maplelatex}

\begin{maplelatex}
\mapleinline{inert}{2d}{`  366: `, "Lewanowicz, J. Comp & Appl. Math. 86,375(1997)
Generalized Watson; Eq.(2.15), n odd, m odd+ Miller, J.Phys.A.,Eq. 1.2
(e = a)";}{%
\maplemultiline{
\mathit{\ \ 366:\ }  \mbox{``Lewanowicz,~J.~Comp~\&~Appl.~Math.~
86,375(1997)~Generalized~Watson; } \\
\mbox{Eq.(2.15),~n~odd,~m~odd+~Miller,~J.Phys.A.,Eq.~1.2~(e~=~
a)''} }
}
\end{maplelatex}

\begin{maplelatex}
\mapleinline{inert}{2d}{`  367: `, "Lewanowicz, J. Comp & Appl. Math. 86,375(1997)
Generalized Watson; Eq.(2.15), m even, n odd+ Miller, J.Phys.A.,Eq.
1.2 (e = a)";}{%
\maplemultiline{
\mathit{\ \ 367:\ }  \mbox{``Lewanowicz,~J.~Comp~\&~Appl.~Math.~
86,375(1997)~Generalized~Watson;~ } \\
\mbox{Eq.(2.15),~m~even,~n~odd+~Miller,~J.Phys.A.,Eq.~1.2~(e~=~
a)''} }
}
\end{maplelatex}

\begin{maplelatex}
\mapleinline{inert}{2d}{`  368: `, "Lewanowicz, J. Comp & Appl. Math. 86,375(1997)
Generalized Watson; Eq.(2.15), m odd, n even+ Miller, J.Phys.A.,Eq.
1.2 (e = a)";}{%
\maplemultiline{
\mathit{\ \ 368:\ }  \mbox{``Lewanowicz,~J.~Comp~\&~Appl.~Math.~
86,375(1997)~Generalized~Watson;~ } \\
\mbox{Eq.(2.15),~m~odd,~n~even+~Miller,~J.Phys.A.,Eq.~1.2~(e~=~
a)''} }
}
\end{maplelatex}

\begin{maplelatex}
\mapleinline{inert}{2d}{`  375: `, "Lewanowicz, J. Comp & Appl. Math. 86,375(1997)
Generalized Watson; Eq.(2.15), n odd, m odd+ Miller, J.Phys.A.,Eq. 1.2
(e = a)+T:1";}{%
\maplemultiline{
\mathit{\ \ 375:\ }  \mbox{``Lewanowicz,~J.~Comp~\&~Appl.~Math.~
86,375(1997)~Generalized~Watson;~ } \\
\mbox{Eq.(2.15),~n~odd,~m~odd+~Miller,~J.Phys.A.,Eq.~1.2~(e~=~
a)+T:1''} }
}
\end{maplelatex}

\begin{maplelatex}
\mapleinline{inert}{2d}{`  376: `, "Lewanowicz, J. Comp & Appl. Math. 86,375(1997)
Generalized Watson; Eq.(2.15), n odd, m odd+ Miller, J.Phys.A.,Eq. 1.2
(e = a)+T:3";}{%
\maplemultiline{
\mathit{\ \ 376:\ }  \mbox{``Lewanowicz,~J.~Comp~\&~Appl.~Math.~
86,375(1997)~Generalized~Watson;~ } \\
\mbox{Eq.(2.15),~n~odd,~m~odd+~Miller,~J.Phys.A.,Eq.~1.2~(e~=~
a)+T:3''} }
}
\end{maplelatex}

\begin{maplelatex}
\mapleinline{inert}{2d}{`  377: `, "Lewanowicz, J. Comp & Appl. Math. 86,375(1997)
Generalized Watson; Eq.(2.15), n odd, m odd+ Miller, J.Phys.A.,Eq. 1.2
(e = a)+T:8";}{%
\maplemultiline{
\mathit{\ \ 377:\ }  \mbox{``Lewanowicz,~J.~Comp~\&~Appl.~Math.~
86,375(1997)~Generalized~Watson;~ } \\
\mbox{Eq.(2.15),~n~odd,~m~odd+~Miller,~J.Phys.A.,Eq.~1.2~(e~=~
a)+T:8''} }
}
\end{maplelatex}

\begin{maplelatex}
\mapleinline{inert}{2d}{`  378: `, "Lewanowicz, J. Comp & Appl. Math. 86,375(1997)
Generalized Watson; Eq.(2.15), m even, n odd+ Miller, J.Phys.A.,Eq.
1.2 (e = a)+T:1";}{%
\maplemultiline{
\mathit{\ \ 378:\ }  \mbox{``Lewanowicz,~J.~Comp~\&~Appl.~Math.~
86,375(1997)~Generalized~Watson;~ } \\
\mbox{Eq.(2.15),~m~even,~n~odd+~Miller,~J.Phys.A.,Eq.~1.2~(e~=~
a)+T:1''} }
}
\end{maplelatex}

\begin{maplelatex}
\mapleinline{inert}{2d}{`  379: `, "Lewanowicz, J. Comp & Appl. Math. 86,375(1997)
Generalized Watson; Eq.(2.15), m even, n odd+ Miller, J.Phys.A.,Eq.
1.2 (e = a)+T:3";}{%
\maplemultiline{
\mathit{\ \ 379:\ }  \mbox{``Lewanowicz,~J.~Comp~\&~Appl.~Math.~
86,375(1997)~Generalized~Watson;~ } \\
\mbox{Eq.(2.15),~m~even,~n~odd+~Miller,~J.Phys.A.,Eq.~1.2~(e~=~
a)+T:3''} }
}
\end{maplelatex}

\begin{maplelatex}
\mapleinline{inert}{2d}{`  380: `, "Lewanowicz, J. Comp & Appl. Math. 86,375(1997)
Generalized Watson; Eq.(2.15), m even, n odd+ Miller, J.Phys.A.,Eq.
1.2 (e = a)+T:8";}{%
\maplemultiline{
\mathit{\ \ 380:\ }  \mbox{``Lewanowicz,~J.~Comp~\&~Appl.~Math.~
86,375(1997)~Generalized~Watson;~ } \\
\mbox{Eq.(2.15),~m~even,~n~odd+~Miller,~J.Phys.A.,Eq.~1.2~(e~=~
a)+T:8''} }
}
\end{maplelatex}

\begin{maplelatex}
\mapleinline{inert}{2d}{`  381: `, "Lewanowicz, J. Comp & Appl. Math. 86,375(1997)
Generalized Watson; Eq.(2.15), m odd, n even+ Miller, J.Phys.A.,Eq.
1.2 (e = a)+T:1";}{%
\maplemultiline{
\mathit{\ \ 381:\ }  \mbox{``Lewanowicz,~J.~Comp~\&~Appl.~Math.~
86,375(1997)~Generalized~Watson;~ } \\
\mbox{Eq.(2.15),~m~odd,~n~even+~Miller,~J.Phys.A.,Eq.~1.2~(e~=~
a)+T:1''} }
}
\end{maplelatex}

\begin{maplelatex}
\mapleinline{inert}{2d}{`  382: `, "Lewanowicz, J. Comp & Appl. Math. 86,375(1997)
Generalized Watson; Eq.(2.15), m odd, n even+ Miller, J.Phys.A.,Eq.
1.2 (e = a)+T:3";}{%
\maplemultiline{
\mathit{\ \ 382:\ }  \mbox{``Lewanowicz,~J.~Comp~\&~Appl.~Math.~
86,375(1997)~Generalized~Watson;~ } \\
\mbox{Eq.(2.15),~m~odd,~n~even+~Miller,~J.Phys.A.,Eq.~1.2~(e~=~
a)+T:3''} }
}
\end{maplelatex}

\begin{maplelatex}
\mapleinline{inert}{2d}{`  383: `, "Lewanowicz, J. Comp & Appl. Math. 86,375(1997)
Generalized Watson; Eq.(2.15), m odd, n even+ Miller, J.Phys.A.,Eq.
1.2 (e = a)+T:8";}{%
\maplemultiline{
\mathit{\ \ 383:\ }  \mbox{``Lewanowicz,~J.~Comp~\&~Appl.~Math.~
86,375(1997)~Generalized~Watson;~ } \\
\mbox{Eq.(2.15),~m~odd,~n~even+~Miller,~J.Phys.A.,Eq.~1.2~(e~=~
a)+T:8''} }
}
\end{maplelatex}

\begin{maplelatex}
\mapleinline{inert}{2d}{`  384: `, "Lewanowicz, J. Comp & Appl. Math. 86,375(1997)
Generalized Watson; Eq.(2.15), n odd, m odd+ Miller, J.Phys.A.,Eq. 1.2
(e = a)+ Miller, J.Phys.A.,Eq. 1.2 (e = a)";}{%
\maplemultiline{
\mathit{\ \ 384:\ }  \mbox{``Lewanowicz,~J.~Comp~\&~Appl.~Math.~
86,375(1997)~Generalized~Watson;~ } \\
\mbox{Eq.(2.15),~n~odd,~m~odd+~Miller,~J.Phys.A.,Eq.~1.2~(e~=~a)+
~Miller,~J.Phys.A.,} \\
\mbox{Eq.~1.2~(e~=~a)''} }
}
\end{maplelatex}

\begin{maplelatex}
\mapleinline{inert}{2d}{`  385: `, "Lewanowicz, J. Comp & Appl. Math. 86,375(1997)
Generalized Watson; Eq.(2.15), m even, n odd+ Miller, J.Phys.A.,Eq.
1.2 (e = a)+ Miller, J.Phys.A.,Eq. 1.2 (e = a)";}{%
\maplemultiline{
\mathit{\ \ 385:\ }  \mbox{``Lewanowicz,~J.~Comp~\&~Appl.~Math.~
86,375(1997)~Generalized~Watson;~ } \\
\mbox{Eq.(2.15),~m~even,~n~odd+~Miller,~J.Phys.A.,Eq.~1.2~(e~=~
a)+~Miller,~J.Phys.A., } \\
\mbox{Eq.~1.2~(e~=~a)''} }
}
\end{maplelatex}

\begin{maplelatex}
\mapleinline{inert}{2d}{`  386: `, "Lewanowicz, J. Comp & Appl. Math. 86,375(1997)
Generalized Watson; Eq.(2.15), m odd, n even+ Miller, J.Phys.A.,Eq.
1.2 (e = a)+ Miller, J.Phys.A.,Eq. 1.2 (e = a)";}{%
\maplemultiline{
\mathit{\ \ 386:\ }  \mbox{``Lewanowicz,~J.~Comp~\&~Appl.~Math.~
86,375(1997)~Generalized~Watson;~ } \\
\mbox{Eq.(2.15),~m~odd,~n~even+~Miller,~J.Phys.A.,Eq.~1.2~(e~=~
a)+~Miller, \\ ~J.Phys.A., }\\
\mbox{Eq.~1.2~(e~=~a)''} }
}
\end{maplelatex}

\begin{maplelatex}
\mapleinline{inert}{2d}{`  387: `, "Lewanowicz, J. Comp & Appl. Math. 86,375(1997)
Generalized Watson; Eq.(2.15), n odd, m odd+ Miller, J.Phys.A.,Eq. 1.2
(e = a)+T:3+ Miller, J.Phys.A.,Eq. 1.2 (e = a)";}{%
\maplemultiline{
\mathit{\ \ 387:\ }  \mbox{``Lewanowicz,~J.~Comp~\&~Appl.~Math.~
86,375(1997)~Generalized~Watson;~ } \\
\mbox{Eq.(2.15),~n~odd,~m~odd+~Miller,~J.Phys.A.,Eq.~1.2~(e~=~
a)+T:3+~Miller,~ } \\
\mbox{J.Phys.A.,Eq.~1.2~(e~=~a)''} }
}
\end{maplelatex}

\begin{maplelatex}
\mapleinline{inert}{2d}{`  388: `, "Lewanowicz, J. Comp & Appl. Math. 86,375(1997)
Generalized Watson; Eq.(2.15), n odd, m odd+ Miller, J.Phys.A.,Eq. 1.2
(e = a)+T:8+ Miller, J.Phys.A.,Eq. 1.2 (e = a)";}{%
\maplemultiline{
\mathit{\ \ 388:\ }  \mbox{``Lewanowicz,~J.~Comp~\&~Appl.~Math.~
86,375(1997)~Generalized~Watson;~ } \\
\mbox{Eq.(2.15),~n~odd,~m~odd+~Miller,~J.Phys.A.,Eq.~1.2~(e~=~
a)+T:8+~Miller,~ } \\
\mbox{J.Phys.A.,Eq.~1.2~(e~=~a)''} }
}
\end{maplelatex}

\begin{maplelatex}
\mapleinline{inert}{2d}{`  389: `, "Lewanowicz, J. Comp & Appl. Math. 86,375(1997)
Generalized Watson; Eq.(2.15), m even, n odd+ Miller, J.Phys.A.,Eq.
1.2 (e = a)+T:3+ Miller, J.Phys.A.,Eq. 1.2 (e = a)";}{%
\maplemultiline{
\mathit{\ \ 389:\ }  \mbox{``Lewanowicz,~J.~Comp~\&~Appl.~Math.~
86,375(1997)~Generalized~Watson;~ } \\
\mbox{Eq.(2.15),~m~even,~n~odd+~Miller,~J.Phys.A.,Eq.~1.2~(e~=~
a)+T:3+~Miller,~ } \\
\mbox{J.Phys.A.,Eq.~1.2~(e~=~a)''} }
}
\end{maplelatex}

\begin{maplelatex}
\mapleinline{inert}{2d}{`  390: `, "Lewanowicz, J. Comp & Appl. Math. 86,375(1997)
Generalized Watson; Eq.(2.15), m even, n odd+ Miller, J.Phys.A.,Eq.
1.2 (e = a)+T:8+ Miller, J.Phys.A.,Eq. 1.2 (e = a)";}{%
\maplemultiline{
\mathit{\ \ 390:\ }  \mbox{``Lewanowicz,~J.~Comp~\&~Appl.~Math.~
86,375(1997)~Generalized~Watson;~ } \\
\mbox{Eq.(2.15),~m~even,~n~odd+~Miller,~J.Phys.A.,Eq.~1.2~(e~=~
a)+T:8+~Miller,~ } \\
\mbox{J.Phys.A.,Eq.~1.2~(e~=~a)''} }
}
\end{maplelatex}

\begin{maplelatex}
\mapleinline{inert}{2d}{`  391: `, "Lewanowicz, J. Comp & Appl. Math. 86,375(1997)
Generalized Watson; Eq.(2.15), m odd, n even+ Miller, J.Phys.A.,Eq.
1.2 (e = a)+T:3+ Miller, J.Phys.A.,Eq. 1.2 (e = a)";}{%
\maplemultiline{
\mathit{\ \ 391:\ }  \mbox{``Lewanowicz,~J.~Comp~\&~Appl.~Math.~
86,375(1997)~Generalized~Watson;~ } \\
\mbox{Eq.(2.15),~m~odd,~n~even+~Miller,~J.Phys.A.,Eq.~1.2~(e~=~
a)+T:3+~Miller,~ } \\
\mbox{J.Phys.A.,Eq.~1.2~(e~=~a)''} }
}
\end{maplelatex}

\begin{maplelatex}
\mapleinline{inert}{2d}{`  392: `, "Lewanowicz, J. Comp & Appl. Math. 86,375(1997)
Generalized Watson; Eq.(2.15), m odd, n even+ Miller, J.Phys.A.,Eq.
1.2 (e = a)+T:8+ Miller, J.Phys.A.,Eq. 1.2 (e = a)";}{%
\maplemultiline{
\mathit{\ \ 392:\ }  \mbox{``Lewanowicz,~J.~Comp~\&~Appl.~Math.~
86,375(1997)~Generalized~Watson;~ } \\
\mbox{Eq.(2.15),~m~odd,~n~even+~Miller,~J.Phys.A.,Eq.~1.2~(e~=~
a)+T:8+~Miller,~ } \\
\mbox{J.Phys.A.,Eq.~1.2~(e~=~a)''} }
}
\end{maplelatex}

\begin{maplelatex}
\mapleinline{inert}{2d}{`  393: `, "Lewanowicz, J. Comp & Appl. Math. 86,375(1997)
Generalized Watson; Eq.(2.15), n odd, m odd+ Miller, J.Phys.A.,Eq. 1.2
(e = a)+ Miller, J.Phys.A.,Eq. 1.2 (e = a)+T:1";}{%
\maplemultiline{
\mathit{\ \ 393:\ }  \mbox{``Lewanowicz,~J.~Comp~\&~Appl.~Math.~
86,375(1997)~Generalized~Watson;~ } \\
\mbox{Eq.(2.15),~n~odd,~m~odd+~Miller,~J.Phys.A.,Eq.~1.2~(e~=~a)+
~Miller,~J.Phys.A.,} \\
\mbox{Eq.~1.2~(e~=~a)+T:1''} }
}
\end{maplelatex}

\begin{maplelatex}
\mapleinline{inert}{2d}{`  394: `, "Lewanowicz, J. Comp & Appl. Math. 86,375(1997)
Generalized Watson; Eq.(2.15), n odd, m odd+ Miller, J.Phys.A.,Eq. 1.2
(e = a)+ Miller, J.Phys.A.,Eq. 1.2 (e = a)+T:4";}{%
\maplemultiline{
\mathit{\ \ 394:\ }  \mbox{``Lewanowicz,~J.~Comp~\&~Appl.~Math.~
86,375(1997)~Generalized~Watson;~ } \\
\mbox{Eq.(2.15),~n~odd,~m~odd+~Miller,~J.Phys.A.,Eq.~1.2~(e~=~a)+
~Miller,~J.Phys.A., } \\
\mbox{Eq.~1.2~(e~=~a)+T:4''} }
}
\end{maplelatex}

\begin{maplelatex}
\mapleinline{inert}{2d}{`  395: `, "Lewanowicz, J. Comp & Appl. Math. 86,375(1997)
Generalized Watson; Eq.(2.15), m even, n odd+ Miller, J.Phys.A.,Eq.
1.2 (e = a)+ Miller, J.Phys.A.,Eq. 1.2 (e = a)+T:1";}{%
\maplemultiline{
\mathit{\ \ 395:\ }  \mbox{``Lewanowicz,~J.~Comp~\&~Appl.~Math.~
86,375(1997)~Generalized~Watson;~ } \\
\mbox{Eq.(2.15),~m~even,~n~odd+~Miller,~J.Phys.A.,Eq.~1.2~(e~=~
a)+~Miller,~J.Phys.A., } \\
\mbox{Eq.~1.2~(e~=~a)+T:1''} }
}
\end{maplelatex}

\begin{maplelatex}
\mapleinline{inert}{2d}{`  396: `, "Lewanowicz, J. Comp & Appl. Math. 86,375(1997)
Generalized Watson; Eq.(2.15), m even, n odd+ Miller, J.Phys.A.,Eq.
1.2 (e = a)+ Miller, J.Phys.A.,Eq. 1.2 (e = a)+T:4";}{%
\maplemultiline{
\mathit{\ \ 396:\ }  \mbox{``Lewanowicz,~J.~Comp~\&~Appl.~Math.~
86,375(1997)~Generalized~Watson;~ } \\
\mbox{Eq.(2.15),~m~even,~n~odd+~Miller,~J.Phys.A.,Eq.~1.2~(e~=~
a)+~Miller,~J.Phys.A., } \\
\mbox{Eq.~1.2~(e~=~a)+T:4''} }
}
\end{maplelatex}

\begin{maplelatex}
\mapleinline{inert}{2d}{`  397: `, "Lewanowicz, J. Comp & Appl. Math. 86,375(1997)
Generalized Watson; Eq.(2.15), m odd, n even+ Miller, J.Phys.A.,Eq.
1.2 (e = a)+ Miller, J.Phys.A.,Eq. 1.2 (e = a)+T:1";}{%
\maplemultiline{
\mathit{\ \ 397:\ }  \mbox{``Lewanowicz,~J.~Comp~\&~Appl.~Math.~
86,375(1997)~Generalized~Watson;~ } \\
\mbox{Eq.(2.15),~m~odd,~n~even+~Miller,~J.Phys.A.,Eq.~1.2~(e~=~
a)+~Miller,~J.Phys.A., } \\
\mbox{Eq.~1.2~(e~=~a)+T:1''} }
}
\end{maplelatex}

\begin{maplelatex}
\mapleinline{inert}{2d}{`  398: `, "Lewanowicz, J. Comp & Appl. Math. 86,375(1997)
Generalized Watson; Eq.(2.15), m odd, n even+ Miller, J.Phys.A.,Eq.
1.2 (e = a)+ Miller, J.Phys.A.,Eq. 1.2 (e = a)+T:4";}{%
\maplemultiline{
\mathit{\ \ 398:\ }  \mbox{``Lewanowicz,~J.~Comp~\&~Appl.~Math.~
86,375(1997)~Generalized~Watson;~ } \\
\mbox{Eq.(2.15),~m~odd,~n~even+~Miller,~J.Phys.A.,Eq.~1.2~(e~=~
a)+~Miller,~J.Phys.A.,} \\
\mbox{Eq.~1.2~(e~=~a)+T:4''} }
}
\end{maplelatex}

\begin{maplelatex}
\mapleinline{inert}{2d}{`  399: `, "Lewanowicz, J. Comp & Appl. Math. 86,375(1997)
Generalized Watson; Eq.(2.15), n odd, m odd+ Miller, J.Phys.A.,Eq. 1.2
(e = a)+T:3+ Miller, J.Phys.A.,Eq. 1.2 (e = a)+ Miller, J.Phys.A.,Eq.
1.2 (e = a)";}{%
\maplemultiline{
\mathit{\ \ 399:\ }  \mbox{``Lewanowicz,~J.~Comp~\&~Appl.~Math.~
86,375(1997)~Generalized~Watson;~ } \\
\mbox{Eq.(2.15),~n~odd,~m~odd+~Miller,~J.Phys.A.,Eq.~1.2~(e~=~
a)+T:3+~Miller,~ } \\
\mbox{J.Phys.A.,Eq.~1.2~(e~=~a)+~Miller,~J.Phys.A.,Eq.~1.2~(e~=~
a)''} }
}
\end{maplelatex}

\begin{maplelatex}
\mapleinline{inert}{2d}{`  400: `, "Lewanowicz, J. Comp & Appl. Math. 86,375(1997)
Generalized Watson; Eq.(2.15), m even, n odd+ Miller, J.Phys.A.,Eq.
1.2 (e = a)+T:3+ Miller, J.Phys.A.,Eq. 1.2 (e = a)+ Miller,
J.Phys.A.,Eq. 1.2 (e = a)";}{%
\maplemultiline{
\mathit{\ \ 400:\ }  \mbox{``Lewanowicz,~J.~Comp~\&~Appl.~Math.~
86,375(1997)~Generalized~Watson;~ } \\
\mbox{Eq.(2.15),~m~even,~n~odd+~Miller,~J.Phys.A.,Eq.~1.2~(e~=~
a)+T:3+~Miller,~}  \\
\mbox{J.Phys.A.,Eq.~1.2~(e~=~a)+~Miller,~J.Phys.A.,Eq.~1.2~(e~=~
a)''} }
}
\end{maplelatex}

\begin{maplelatex}
\mapleinline{inert}{2d}{`  401: `, "Lewanowicz, J. Comp & Appl. Math. 86,375(1997)
Generalized Watson; Eq.(2.15), m odd, n even+ Miller, J.Phys.A.,Eq.
1.2 (e = a)+T:3+ Miller, J.Phys.A.,Eq. 1.2 (e = a)+ Miller,
J.Phys.A.,Eq. 1.2 (e = a)";}{%
\maplemultiline{
\mathit{\ \ 401:\ }  \mbox{``Lewanowicz,~J.~Comp~\&~Appl.~Math.~
86,375(1997)~Generalized~Watson;~ } \\
\mbox{Eq.(2.15),~m~odd,~n~even+~Miller,~J.Phys.A.,Eq.~1.2~(e~=~
a)+T:3+~Miller,~ }  \\
\mbox{J.Phys.A.,Eq.~1.2~(e~=~a)+~Miller,~J.Phys.A.,Eq.~1.2~(e~=~
a)''} }
}
\end{maplelatex}

\begin{maplelatex}
\mapleinline{inert}{2d}{`  402: `, "Lewanowicz, J. Comp & Appl. Math. 86,375(1997)
Generalized Watson; Eq.(2.15), m odd, n even+ Miller, J.Phys.A.,Eq.
1.2 (e = a)+ Miller, J.Phys.A.,Eq. 1.2 (e = a)+T:1+ Miller,
J.Phys.A.,Eq. 1.2 (e = a)";}{%
\maplemultiline{
\mathit{\ \ 402:\ }  \mbox{``Lewanowicz,~J.~Comp~\&~Appl.~Math.~
86,375(1997)~Generalized~Watson;~ } \\
\mbox{Eq.(2.15),~m~odd,~n~even+~Miller,~J.Phys.A.,Eq.~1.2~(e~=~
a)+~Miller,~J.Phys.A., } \\
\mbox{Eq.~1.2~(e~=~a)+T:1+~Miller,~J.Phys.A.,Eq.~1.2~(e~=~
a)''} }
}
\end{maplelatex}

\begin{maplelatex}
\mapleinline{inert}{2d}{`  403: `, "Lewanowicz, J. Comp & Appl. Math. 86,375(1997)
Generalized Watson; Eq.(2.15), n odd, m odd+ Miller, J.Phys.A.,Eq. 1.2
(e = a)+T:3+ Miller, J.Phys.A.,Eq. 1.2 (e = a)+ Miller, J.Phys.A.,Eq.
1.2 (e = a)+T:1";}{%
\maplemultiline{
\mathit{\ \ 403:\ }  \mbox{``Lewanowicz,~J.~Comp~\&~Appl.~Math.~
86,375(1997)~Generalized~Watson;~ } \\
\mbox{Eq.(2.15),~n~odd,~m~odd+~Miller,~J.Phys.A.,Eq.~1.2~(e~=~
a)+T:3+~Miller,~ }  \\
\mbox{J.Phys.A.,Eq.~1.2~(e~=~a)+~Miller,~J.Phys.A.,Eq.~1.2~(e~=~
a)+T:1''} }
}
\end{maplelatex}

\begin{maplelatex}
\mapleinline{inert}{2d}{`  404: `, "Lewanowicz, J. Comp & Appl. Math. 86,375(1997)
Generalized Watson; Eq.(2.15), n odd, m odd+ Miller, J.Phys.A.,Eq. 1.2
(e = a)+T:3+ Miller, J.Phys.A.,Eq. 1.2 (e = a)+ Miller, J.Phys.A.,Eq.
1.2 (e = a)+T:2";}{%
\maplemultiline{
\mathit{\ \ 404:\ }  \mbox{``Lewanowicz,~J.~Comp~\&~Appl.~Math.~
86,375(1997)~Generalized~Watson;~ } \\
\mbox{Eq.(2.15),~n~odd,~m~odd+~Miller,~J.Phys.A.,Eq.~1.2~(e~=~
a)+T:3+~Miller,~ }  \\
\mbox{J.Phys.A.,Eq.~1.2~(e~=~a)+~Miller,~J.Phys.A.,Eq.~1.2~(e~=~
a)+T:2''} }
}
\end{maplelatex}

\begin{maplelatex}
\mapleinline{inert}{2d}{`  405: `, "Lewanowicz, J. Comp & Appl. Math. 86,375(1997)
Generalized Watson; Eq.(2.15), n odd, m odd+ Miller, J.Phys.A.,Eq. 1.2
(e = a)+T:3+ Miller, J.Phys.A.,Eq. 1.2 (e = a)+ Miller, J.Phys.A.,Eq.
1.2 (e = a)+T:7";}{%
\maplemultiline{
\mathit{\ \ 405:\ }  \mbox{``Lewanowicz,~J.~Comp~\&~Appl.~Math.~
86,375(1997)~Generalized~Watson;~ } \\
\mbox{Eq.(2.15),~n~odd,~m~odd+~Miller,~J.Phys.A.,Eq.~1.2~(e~=~
a)+T:3+~Miller,~ }  \\
\mbox{J.Phys.A.,Eq.~1.2~(e~=~a)+~Miller,~J.Phys.A.,Eq.~1.2~(e~=~
a)+T:7''} }
}
\end{maplelatex}

\begin{maplelatex}
\mapleinline{inert}{2d}{`  406: `, "Lewanowicz, J. Comp & Appl. Math. 86,375(1997)
Generalized Watson; Eq.(2.15), n odd, m odd+ Miller, J.Phys.A.,Eq. 1.2
(e = a)+T:3+ Miller, J.Phys.A.,Eq. 1.2 (e = a)+ Miller, J.Phys.A.,Eq.
1.2 (e = a)+T:8";}{%
\maplemultiline{
\mathit{\ \ 406:\ }  \mbox{``Lewanowicz,~J.~Comp~\&~Appl.~Math.~
86,375(1997)~Generalized~Watson;~ } \\
\mbox{Eq.(2.15),~n~odd,~m~odd+~Miller,~J.Phys.A.,Eq.~1.2~(e~=~
a)+T:3+~Miller,~ }  \\
\mbox{J.Phys.A.,Eq.~1.2~(e~=~a)+~Miller,~J.Phys.A.,Eq.~1.2~(e~=~
a)+T:8''} }
}
\end{maplelatex}

\begin{maplelatex}
\mapleinline{inert}{2d}{`  407: `, "Lewanowicz, J. Comp & Appl. Math. 86,375(1997)
Generalized Watson; Eq.(2.15), m even, n odd+ Miller, J.Phys.A.,Eq.
1.2 (e = a)+T:3+ Miller, J.Phys.A.,Eq. 1.2 (e = a)+ Miller,
J.Phys.A.,Eq. 1.2 (e = a)+T:1";}{%
\maplemultiline{
\mathit{\ \ 407:\ }  \mbox{``Lewanowicz,~J.~Comp~\&~Appl.~Math.~
86,375(1997)~Generalized~Watson;~ } \\
\mbox{Eq.(2.15),~m~even,~n~odd+~Miller,~J.Phys.A.,Eq.~1.2~(e~=~
a)+T:3+~Miller,~ }  \\
\mbox{J.Phys.A.,Eq.~1.2~(e~=~a)+~Miller,~J.Phys.A.,Eq.~1.2~(e~=~
a)+T:1''} }
}
\end{maplelatex}

\begin{maplelatex}
\mapleinline{inert}{2d}{`  408: `, "Lewanowicz, J. Comp & Appl. Math. 86,375(1997)
Generalized Watson; Eq.(2.15), m even, n odd+ Miller, J.Phys.A.,Eq.
1.2 (e = a)+T:3+ Miller, J.Phys.A.,Eq. 1.2 (e = a)+ Miller,
J.Phys.A.,Eq. 1.2 (e = a)+T:2";}{%
\maplemultiline{
\mathit{\ \ 408:\ }  \mbox{``Lewanowicz,~J.~Comp~\&~Appl.~Math.~
86,375(1997)~Generalized~Watson;~ } \\
\mbox{Eq.(2.15),~m~even,~n~odd+~Miller,~J.Phys.A.,Eq.~1.2~(e~=~
a)+T:3+~Miller,~ }  \\
\mbox{J.Phys.A.,Eq.~1.2~(e~=~a)+~Miller,~J.Phys.A.,Eq.~1.2~(e~=~
a)+T:2''} }
}
\end{maplelatex}

\begin{maplelatex}
\mapleinline{inert}{2d}{`  409: `, "Lewanowicz, J. Comp & Appl. Math. 86,375(1997)
Generalized Watson; Eq.(2.15), m even, n odd+ Miller, J.Phys.A.,Eq.
1.2 (e = a)+T:3+ Miller, J.Phys.A.,Eq. 1.2 (e = a)+ Miller,
J.Phys.A.,Eq. 1.2 (e = a)+T:7";}{%
\maplemultiline{
\mathit{\ \ 409:\ }  \mbox{``Lewanowicz,~J.~Comp~\&~Appl.~Math.~
86,375(1997)~Generalized~Watson;~ } \\
\mbox{Eq.(2.15),~m~even,~n~odd+~Miller,~J.Phys.A.,Eq.~1.2~(e~=~
a)+T:3+~Miller,~ }  \\
\mbox{J.Phys.A.,Eq.~1.2~(e~=~a)+~Miller,~J.Phys.A.,Eq.~1.2~(e~=~
a)+T:7''} }
}
\end{maplelatex}

\begin{maplelatex}
\mapleinline{inert}{2d}{`  410: `, "Lewanowicz, J. Comp & Appl. Math. 86,375(1997)
Generalized Watson; Eq.(2.15), m even, n odd+ Miller, J.Phys.A.,Eq.
1.2 (e = a)+T:3+ Miller, J.Phys.A.,Eq. 1.2 (e = a)+ Miller,
J.Phys.A.,Eq. 1.2 (e = a)+T:8";}{%
\maplemultiline{
\mathit{\ \ 410:\ }  \mbox{``Lewanowicz,~J.~Comp~\&~Appl.~Math.~
86,375(1997)~Generalized~Watson;~ } \\
\mbox{Eq.(2.15),~m~even,~n~odd+~Miller,~J.Phys.A.,Eq.~1.2~(e~=~
a)+T:3+~Miller,~ }  \\
\mbox{J.Phys.A.,Eq.~1.2~(e~=~a)+~Miller,~J.Phys.A.,Eq.~1.2~(e~=~
a)+T:8''} }
}
\end{maplelatex}

\begin{maplelatex}
\mapleinline{inert}{2d}{`  411: `, "Lewanowicz, J. Comp & Appl. Math. 86,375(1997)
Generalized Watson; Eq.(2.15), m odd, n even+ Miller, J.Phys.A.,Eq.
1.2 (e = a)+T:3+ Miller, J.Phys.A.,Eq. 1.2 (e = a)+ Miller,
J.Phys.A.,Eq. 1.2 (e = a)+T:2";}{%
\maplemultiline{
\mathit{\ \ 411:\ }  \mbox{``Lewanowicz,~J.~Comp~\&~Appl.~Math.~
86,375(1997)~Generalized~Watson;~ } \\
\mbox{Eq.(2.15),~m~odd,~n~even+~Miller,~J.Phys.A.,Eq.~1.2~(e~=~
a)+T:3+~Miller,~ }  \\
\mbox{J.Phys.A.,Eq.~1.2~(e~=~a)+~Miller,~J.Phys.A.,Eq.~1.2~(e~=~
a)+T:2''} }
}
\end{maplelatex}

\begin{maplelatex}
\mapleinline{inert}{2d}{`  412: `, "Lewanowicz, J. Comp & Appl. Math. 86,375(1997)
Generalized Watson; Eq.(2.15), m odd, n even+ Miller, J.Phys.A.,Eq.
1.2 (e = a)+T:3+ Miller, J.Phys.A.,Eq. 1.2 (e = a)+ Miller,
J.Phys.A.,Eq. 1.2 (e = a)+T:7";}{%
\maplemultiline{
\mathit{\ \ 412:\ }  \mbox{``Lewanowicz,~J.~Comp~\&~Appl.~Math.~
86,375(1997)~Generalized~Watson;~ } \\
\mbox{Eq.(2.15),~m~odd,~n~even+~Miller,~J.Phys.A.,Eq.~1.2~(e~=~
a)+T:3+~Miller,~ }  \\
\mbox{J.Phys.A.,Eq.~1.2~(e~=~a)+~Miller,~J.Phys.A.,Eq.~1.2~(e~=~
a)+T:7''} }
}
\end{maplelatex}

\begin{maplelatex}
\mapleinline{inert}{2d}{`  413: `, "Lewanowicz, J. Comp & Appl. Math. 86,375(1997)
Generalized Watson; Eq.(2.15), m odd, n even+ Miller, J.Phys.A.,Eq.
1.2 (e = a)+T:3+ Miller, J.Phys.A.,Eq. 1.2 (e = a)+ Miller,
J.Phys.A.,Eq. 1.2 (e = a)+T:8";}{%
\maplemultiline{
\mathit{\ \ 413:\ }  \mbox{``Lewanowicz,~J.~Comp~\&~Appl.~Math.~
86,375(1997)~Generalized~Watson;~ } \\
\mbox{Eq.(2.15),~m~odd,~n~even+~Miller,~J.Phys.A.,Eq.~1.2~(e~=~
a)+T:3+~Miller,~ }  \\
\mbox{J.Phys.A.,Eq.~1.2~(e~=~a)+~Miller,~J.Phys.A.,Eq.~1.2~(e~=~
a)+T:8''} }
}
\end{maplelatex}

\begin{maplelatex}
\mapleinline{inert}{2d}{`  414: `, "Maier Thm.(7.3) with L=n variation 1+ Miller,
J.Phys.A.,Eq. 1.2 (a = g)";}{%
\[
\mathit{\ \ 414:\ }  \,\mbox{``Maier~Thm.(7.3)~with~L=n~variation
~1+~Miller,~J.Phys.A.,Eq.~1.2~(a~=~g)''}
\]
}
\end{maplelatex}

\begin{maplelatex}
\mapleinline{inert}{2d}{`  415: `, "Maier Thm.(7.3) with L=n variation 1+ Miller,
J.Phys.A.,Eq. 1.2 (a = g)";}{%
\[
\mathit{\ \ 415:\ }  \,\mbox{``Maier~Thm.(7.3)~with~L=n~variation
~1+~Miller,~J.Phys.A.,Eq.~1.2~(a~=~g)''}
\]
}
\end{maplelatex}

\begin{maplelatex}
\mapleinline{inert}{2d}{`  416: `, "Maier Thm.(7.3) with L=n variation 1+ Miller,
J.Phys.A.,Eq. 1.2 (a = g)";}{%
\[
\mathit{\ \ 416:\ }  \,\mbox{``Maier~Thm.(7.3)~with~L=n~variation
~1+~Miller,~J.Phys.A.,Eq.~1.2~(a~=~g)''}
\]
}
\end{maplelatex}

\begin{maplelatex}
\mapleinline{inert}{2d}{`  417: `, "Maier Thm.(7.3) with L=n variation 1+ Miller,
J.Phys.A.,Eq. 1.2 (a = g)";}{%
\[
\mathit{\ \ 417:\ }  \,\mbox{``Maier~Thm.(7.3)~with~L=n~variation
~1+~Miller,~J.Phys.A.,Eq.~1.2~(a~=~g)''}
\]
}
\end{maplelatex}

\begin{maplelatex}
\mapleinline{inert}{2d}{`  418: `, "Maier Thm.(7.3) with L=n variation 1+ Miller,
J.Phys.A.,Eq. 1.2 (a = g)";}{%
\[
\mathit{\ \ 418:\ }  \,\mbox{``Maier~Thm.(7.3)~with~L=n~variation
~1+~Miller,~J.Phys.A.,Eq.~1.2~(a~=~g)''}
\]
}
\end{maplelatex}

\begin{maplelatex}
\mapleinline{inert}{2d}{`  419: `, "Maier Thm.(7.3) with L=n variation 1+ Miller,
J.Phys.A.,Eq. 1.2 (a = g)";}{%
\[
\mathit{\ \ 419:\ }  \,\mbox{``Maier~Thm.(7.3)~with~L=n~variation
~1+~Miller,~J.Phys.A.,Eq.~1.2~(a~=~g)''}
\]
}
\end{maplelatex}

\begin{maplelatex}
\mapleinline{inert}{2d}{`  420: `, "Maier Thm.(7.3) with L=n variation 1+ Miller,
J.Phys.A.,Eq. 1.2 (a = g)+T:1";}{%
\[
\mathit{\ \ 420:\ }  \,\mbox{``Maier~Thm.(7.3)~with~L=n~variation
~1+~Miller,~J.Phys.A.,Eq.~1.2~(a~=~g)+T:1''}
\]
}
\end{maplelatex}

\begin{maplelatex}
\mapleinline{inert}{2d}{`  421: `, "Maier Thm.(7.3) with L=n variation 1+ Miller,
J.Phys.A.,Eq. 1.2 (a = g)+T:2";}{%
\[
\mathit{\ \ 421:\ }  \,\mbox{``Maier~Thm.(7.3)~with~L=n~variation
~1+~Miller,~J.Phys.A.,Eq.~1.2~(a~=~g)+T:2''}
\]
}
\end{maplelatex}

\begin{maplelatex}
\mapleinline{inert}{2d}{`  422: `, "Maier Thm.(7.3) with L=n variation 1+ Miller,
J.Phys.A.,Eq. 1.2 (a = g)+T:3";}{%
\[
\mathit{\ \ 422:\ }  \,\mbox{``Maier~Thm.(7.3)~with~L=n~variation
~1+~Miller,~J.Phys.A.,Eq.~1.2~(a~=~g)+T:3''}
\]
}
\end{maplelatex}

\begin{maplelatex}
\mapleinline{inert}{2d}{`  423: `, "Maier Thm.(7.3) with L=n variation 1+ Miller,
J.Phys.A.,Eq. 1.2 (a = g)+T:4";}{%
\[
\mathit{\ \ 423:\ }  \,\mbox{``Maier~Thm.(7.3)~with~L=n~variation
~1+~Miller,~J.Phys.A.,Eq.~1.2~(a~=~g)+T:4''}
\]
}
\end{maplelatex}

\begin{maplelatex}
\mapleinline{inert}{2d}{`  424: `, "Maier Thm.(7.3) with L=n variation 1+ Miller,
J.Phys.A.,Eq. 1.2 (a = g)+T:5";}{%
\[
\mathit{\ \ 424:\ }  \,\mbox{``Maier~Thm.(7.3)~with~L=n~variation
~1+~Miller,~J.Phys.A.,Eq.~1.2~(a~=~g)+T:5''}
\]
}
\end{maplelatex}

\begin{maplelatex}
\mapleinline{inert}{2d}{`  425: `, "Maier Thm.(7.3) with L=n variation 1+ Miller,
J.Phys.A.,Eq. 1.2 (a = g)+T:6";}{%
\[
\mathit{\ \ 425:\ }  \,\mbox{``Maier~Thm.(7.3)~with~L=n~variation
~1+~Miller,~J.Phys.A.,Eq.~1.2~(a~=~g)+T:6''}
\]
}
\end{maplelatex}

\begin{maplelatex}
\mapleinline{inert}{2d}{`  426: `, "Maier Thm.(7.3) with L=n variation 1+ Miller,
J.Phys.A.,Eq. 1.2 (a = g)+T:8";}{%
\[
\mathit{\ \ 426:\ }  \,\mbox{``Maier~Thm.(7.3)~with~L=n~variation
~1+~Miller,~J.Phys.A.,Eq.~1.2~(a~=~g)+T:8''}
\]
}
\end{maplelatex}

\begin{maplelatex}
\mapleinline{inert}{2d}{`  427: `, "Maier Thm.(7.3) with L=n variation 1+ Miller,
J.Phys.A.,Eq. 1.2 (a = g)+T:1";}{%
\[
\mathit{\ \ 427:\ }  \,\mbox{``Maier~Thm.(7.3)~with~L=n~variation
~1+~Miller,~J.Phys.A.,Eq.~1.2~(a~=~g)+T:1''}
\]
}
\end{maplelatex}

\begin{maplelatex}
\mapleinline{inert}{2d}{`  428: `, "Maier Thm.(7.3) with L=n variation 1+ Miller,
J.Phys.A.,Eq. 1.2 (a = g)+T:2";}{%
\[
\mathit{\ \ 428:\ }  \,\mbox{``Maier~Thm.(7.3)~with~L=n~variation
~1+~Miller,~J.Phys.A.,Eq.~1.2~(a~=~g)+T:2''}
\]
}
\end{maplelatex}

\begin{maplelatex}
\mapleinline{inert}{2d}{`  429: `, "Maier Thm.(7.3) with L=n variation 1+ Miller,
J.Phys.A.,Eq. 1.2 (a = g)+T:3";}{%
\[
\mathit{\ \ 429:\ }  \,\mbox{``Maier~Thm.(7.3)~with~L=n~variation
~1+~Miller,~J.Phys.A.,Eq.~1.2~(a~=~g)+T:3''}
\]
}
\end{maplelatex}

\begin{maplelatex}
\mapleinline{inert}{2d}{`  430: `, "Maier Thm.(7.3) with L=n variation 1+ Miller,
J.Phys.A.,Eq. 1.2 (a = g)+T:4";}{%
\[
\mathit{\ \ 430:\ }  \,\mbox{``Maier~Thm.(7.3)~with~L=n~variation
~1+~Miller,~J.Phys.A.,Eq.~1.2~(a~=~g)+T:4''}
\]
}
\end{maplelatex}

\begin{maplelatex}
\mapleinline{inert}{2d}{`  431: `, "Maier Thm.(7.3) with L=n variation 1+ Miller,
J.Phys.A.,Eq. 1.2 (a = g)+T:5";}{%
\[
\mathit{\ \ 431:\ }  \,\mbox{``Maier~Thm.(7.3)~with~L=n~variation
~1+~Miller,~J.Phys.A.,Eq.~1.2~(a~=~g)+T:5''}
\]
}
\end{maplelatex}

\begin{maplelatex}
\mapleinline{inert}{2d}{`  432: `, "Maier Thm.(7.3) with L=n variation 1+ Miller,
J.Phys.A.,Eq. 1.2 (a = g)+T:7";}{%
\[
\mathit{\ \ 432:\ }  \,\mbox{``Maier~Thm.(7.3)~with~L=n~variation
~1+~Miller,~J.Phys.A.,Eq.~1.2~(a~=~g)+T:7''}
\]
}
\end{maplelatex}

\begin{maplelatex}
\mapleinline{inert}{2d}{`  433: `, "Maier Thm.(7.3) with L=n variation 1+ Miller,
J.Phys.A.,Eq. 1.2 (a = g)+T:8";}{%
\[
\mathit{\ \ 433:\ }  \,\mbox{``Maier~Thm.(7.3)~with~L=n~variation
~1+~Miller,~J.Phys.A.,Eq.~1.2~(a~=~g)+T:8''}
\]
}
\end{maplelatex}

\begin{maplelatex}
\mapleinline{inert}{2d}{`  434: `, "Maier Thm.(7.3) with L=n variation 1+ Miller,
J.Phys.A.,Eq. 1.2 (a = g)+T:9";}{%
\[
\mathit{\ \ 434:\ }  \,\mbox{``Maier~Thm.(7.3)~with~L=n~variation
~1+~Miller,~J.Phys.A.,Eq.~1.2~(a~=~g)+T:9''}
\]
}
\end{maplelatex}

\begin{maplelatex}
\mapleinline{inert}{2d}{`  435: `, "Maier Thm.(7.3) with L=n variation 1+ Miller,
J.Phys.A.,Eq. 1.2 (a = g)+T:1";}{%
\[
\mathit{\ \ 435:\ }  \,\mbox{``Maier~Thm.(7.3)~with~L=n~variation
~1+~Miller,~J.Phys.A.,Eq.~1.2~(a~=~g)+T:1''}
\]
}
\end{maplelatex}

\begin{maplelatex}
\mapleinline{inert}{2d}{`  436: `, "Maier Thm.(7.3) with L=n variation 1+ Miller,
J.Phys.A.,Eq. 1.2 (a = g)+T:2";}{%
\[
\mathit{\ \ 436:\ }  \,\mbox{``Maier~Thm.(7.3)~with~L=n~variation
~1+~Miller,~J.Phys.A.,Eq.~1.2~(a~=~g)+T:2''}
\]
}
\end{maplelatex}

\begin{maplelatex}
\mapleinline{inert}{2d}{`  437: `, "Maier Thm.(7.3) with L=n variation 1+ Miller,
J.Phys.A.,Eq. 1.2 (a = g)+T:3";}{%
\[
\mathit{\ \ 437:\ }  \,\mbox{``Maier~Thm.(7.3)~with~L=n~variation
~1+~Miller,~J.Phys.A.,Eq.~1.2~(a~=~g)+T:3''}
\]
}
\end{maplelatex}

\begin{maplelatex}
\mapleinline{inert}{2d}{`  438: `, "Maier Thm.(7.3) with L=n variation 1+ Miller,
J.Phys.A.,Eq. 1.2 (a = g)+T:4";}{%
\[
\mathit{\ \ 438:\ }  \,\mbox{``Maier~Thm.(7.3)~with~L=n~variation
~1+~Miller,~J.Phys.A.,Eq.~1.2~(a~=~g)+T:4''}
\]
}
\end{maplelatex}

\begin{maplelatex}
\mapleinline{inert}{2d}{`  439: `, "Maier Thm.(7.3) with L=n variation 1+ Miller,
J.Phys.A.,Eq. 1.2 (a = g)+T:5";}{%
\[
\mathit{\ \ 439:\ }  \,\mbox{``Maier~Thm.(7.3)~with~L=n~variation
~1+~Miller,~J.Phys.A.,Eq.~1.2~(a~=~g)+T:5''}
\]
}
\end{maplelatex}

\begin{maplelatex}
\mapleinline{inert}{2d}{`  440: `, "Maier Thm.(7.3) with L=n variation 1+ Miller,
J.Phys.A.,Eq. 1.2 (a = g)+T:6";}{%
\[
\mathit{\ \ 440:\ }  \,\mbox{``Maier~Thm.(7.3)~with~L=n~variation
~1+~Miller,~J.Phys.A.,Eq.~1.2~(a~=~g)+T:6''}
\]
}
\end{maplelatex}

\begin{maplelatex}
\mapleinline{inert}{2d}{`  441: `, "Maier Thm.(7.3) with L=n variation 1+ Miller,
J.Phys.A.,Eq. 1.2 (a = g)+T:7";}{%
\[
\mathit{\ \ 441:\ }  \,\mbox{``Maier~Thm.(7.3)~with~L=n~variation
~1+~Miller,~J.Phys.A.,Eq.~1.2~(a~=~g)+T:7''}
\]
}
\end{maplelatex}

\begin{maplelatex}
\mapleinline{inert}{2d}{`  442: `, "Maier Thm.(7.3) with L=n variation 1+ Miller,
J.Phys.A.,Eq. 1.2 (a = g)+T:8";}{%
\[
\mathit{\ \ 442:\ }  \,\mbox{``Maier~Thm.(7.3)~with~L=n~variation
~1+~Miller,~J.Phys.A.,Eq.~1.2~(a~=~g)+T:8''}
\]
}
\end{maplelatex}

\begin{maplelatex}
\mapleinline{inert}{2d}{`  443: `, "Maier Thm.(7.3) with L=n variation 1+ Miller,
J.Phys.A.,Eq. 1.2 (a = g)+T:9";}{%
\[
\mathit{\ \ 443:\ }  \,\mbox{``Maier~Thm.(7.3)~with~L=n~variation
~1+~Miller,~J.Phys.A.,Eq.~1.2~(a~=~g)+T:9''}
\]
}
\end{maplelatex}

\begin{maplelatex}
\mapleinline{inert}{2d}{`  444: `, "Maier Thm.(7.3) with L=n variation 1+ Miller,
J.Phys.A.,Eq. 1.2 (a = g)+ Miller, J.Phys.A.,Eq. 1.2 (a = g)";}{%
\maplemultiline{
\mathit{\ \ 444:\ }  \mbox{``Maier~Thm.(7.3)~with~L=n~variation~
1+~Miller,~J.Phys.A.,Eq.~1.2~(a~=~g)+~ } \\
\mbox{Miller,~J.Phys.A.,Eq.~1.2~(a~=~g)''} }
}
\end{maplelatex}

\begin{maplelatex}
\mapleinline{inert}{2d}{`  445: `, "Maier Thm.(7.3) with L=n variation 1+ Miller,
J.Phys.A.,Eq. 1.2 (a = g)+ Miller, J.Phys.A.,Eq. 1.2 (a = g)";}{%
\maplemultiline{
\mathit{\ \ 445:\ }  \mbox{``Maier~Thm.(7.3)~with~L=n~variation~
1+~Miller,~J.Phys.A.,Eq.~1.2~(a~=~g)+~ } \\
\mbox{Miller,~J.Phys.A.,Eq.~1.2~(a~=~g)''} }
}
\end{maplelatex}

\begin{maplelatex}
\mapleinline{inert}{2d}{`  446: `, "Maier Thm.(7.3) with L=n variation 1+ Miller,
J.Phys.A.,Eq. 1.2 (a = g)+ Miller, J.Phys.A.,Eq. 1.2 (a = g)";}{%
\maplemultiline{
\mathit{\ \ 446:\ }  \mbox{``Maier~Thm.(7.3)~with~L=n~variation~
1+~Miller,~J.Phys.A.,Eq.~1.2~(a~=~g)+~ } \\
\mbox{Miller,~J.Phys.A.,Eq.~1.2~(a~=~g)''} }
}
\end{maplelatex}

\begin{maplelatex}
\mapleinline{inert}{2d}{`  447: `, "Maier Thm.(7.3) with L=n variation 1+ Miller,
J.Phys.A.,Eq. 1.2 (a = g)+T:3+ Miller, J.Phys.A.,Eq. 1.2 (a = g)";}{%
\maplemultiline{
\mathit{\ \ 447:\ }  \mbox{``Maier~Thm.(7.3)~with~L=n~variation~
1+~Miller,~J.Phys.A.,Eq.~1.2~(a~=~g)+} \\
\mbox{T:3+~Miller,~J.Phys.A.,Eq.~1.2~(a~=~g)''} }
}
\end{maplelatex}

\begin{maplelatex}
\mapleinline{inert}{2d}{`  448: `, "Maier Thm.(7.3) with L=n variation 1+ Miller,
J.Phys.A.,Eq. 1.2 (a = g)+T:4+ Miller, J.Phys.A.,Eq. 1.2 (a = g)";}{%
\maplemultiline{
\mathit{\ \ 448:\ }  \mbox{``Maier~Thm.(7.3)~with~L=n~variation~
1+~Miller,~J.Phys.A.,Eq.~1.2~(a~=~g)+} \\
\mbox{T:4+~Miller,~J.Phys.A.,Eq.~1.2~(a~=~g)''} }
}
\end{maplelatex}

\begin{maplelatex}
\mapleinline{inert}{2d}{`  449: `, "Maier Thm.(7.3) with L=n variation 1+ Miller,
J.Phys.A.,Eq. 1.2 (a = g)+T:6+ Miller, J.Phys.A.,Eq. 1.2 (a = g)";}{%
\maplemultiline{
\mathit{\ \ 449:\ }  \mbox{``Maier~Thm.(7.3)~with~L=n~variation~
1+~Miller,~J.Phys.A.,Eq.~1.2~(a~=~g)+} \\
\mbox{T:6+~Miller,~J.Phys.A.,Eq.~1.2~(a~=~g)''} }
}
\end{maplelatex}

\begin{maplelatex}
\mapleinline{inert}{2d}{`  450: `, "Maier Thm.(7.3) with L=n variation 1+ Miller,
J.Phys.A.,Eq. 1.2 (a = g)+T:8+ Miller, J.Phys.A.,Eq. 1.2 (a = g)";}{%
\maplemultiline{
\mathit{\ \ 450:\ }  \mbox{``Maier~Thm.(7.3)~with~L=n~variation~
1+~Miller,~J.Phys.A.,Eq.~1.2~(a~=~g)+} \\
\mbox{T:8+~Miller,~J.Phys.A.,Eq.~1.2~(a~=~g)''} }
}
\end{maplelatex}

\begin{maplelatex}
\mapleinline{inert}{2d}{`  451: `, "Maier Thm.(7.3) with L=n variation 1+ Miller,
J.Phys.A.,Eq. 1.2 (a = g)+T:5+ Miller, J.Phys.A.,Eq. 1.2 (a = g)";}{%
\maplemultiline{
\mathit{\ \ 451:\ }  \mbox{``Maier~Thm.(7.3)~with~L=n~variation~
1+~Miller,~J.Phys.A.,Eq.~1.2~(a~=~g)+} \\
\mbox{T:5+~Miller,~J.Phys.A.,Eq.~1.2~(a~=~g)''} }
}
\end{maplelatex}

\begin{maplelatex}
\mapleinline{inert}{2d}{`  452: `, "Maier Thm.(7.3) with L=n variation 1+ Miller,
J.Phys.A.,Eq. 1.2 (a = g)+T:7+ Miller, J.Phys.A.,Eq. 1.2 (a = g)";}{%
\maplemultiline{
\mathit{\ \ 452:\ }  \mbox{``Maier~Thm.(7.3)~with~L=n~variation~
1+~Miller,~J.Phys.A.,Eq.~1.2~(a~=~g)+} \\
\mbox{T:7+~Miller,~J.Phys.A.,Eq.~1.2~(a~=~g)''} }
}
\end{maplelatex}

\begin{maplelatex}
\mapleinline{inert}{2d}{`  453: `, "Maier Thm.(7.3) with L=n variation 1+ Miller,
J.Phys.A.,Eq. 1.2 (a = g)+T:9+ Miller, J.Phys.A.,Eq. 1.2 (a = g)";}{%
\maplemultiline{
\mathit{\ \ 453:\ }  \mbox{``Maier~Thm.(7.3)~with~L=n~variation~
1+~Miller,~J.Phys.A.,Eq.~1.2~(a~=~g)+} \\
\mbox{T:9+~Miller,~J.Phys.A.,Eq.~1.2~(a~=~g)''} }
}
\end{maplelatex}

\begin{maplelatex}
\mapleinline{inert}{2d}{`  454: `, "Maier Thm.(7.3) with L=n variation 1+ Miller,
J.Phys.A.,Eq. 1.2 (a = g)+T:1+ Miller, J.Phys.A.,Eq. 1.2 (a = g)";}{%
\maplemultiline{
\mathit{\ \ 454:\ }  \mbox{``Maier~Thm.(7.3)~with~L=n~variation~
1+~Miller,~J.Phys.A.,Eq.~1.2~(a~=~g)+} \\
\mbox{T:1+~Miller,~J.Phys.A.,Eq.~1.2~(a~=~g)''} }
}
\end{maplelatex}

\begin{maplelatex}
\mapleinline{inert}{2d}{`  455: `, "Maier Thm.(7.3) with L=n variation 1+ Miller,
J.Phys.A.,Eq. 1.2 (a = g)+T:6+ Miller, J.Phys.A.,Eq. 1.2 (a = g)";}{%
\maplemultiline{
\mathit{\ \ 455:\ }  \mbox{``Maier~Thm.(7.3)~with~L=n~variation~
1+~Miller,~J.Phys.A.,Eq.~1.2~(a~=~g)+} \\
\mbox{T:6+~Miller,~J.Phys.A.,Eq.~1.2~(a~=~g)''} }
}
\end{maplelatex}

\begin{maplelatex}
\mapleinline{inert}{2d}{`  456: `, "Maier Thm.(7.3) with L=n variation 1+ Miller,
J.Phys.A.,Eq. 1.2 (a = g)+ Miller, J.Phys.A.,Eq. 1.2 (a = g)+T:1";}{%
\maplemultiline{
\mathit{\ \ 456:\ }  \mbox{``Maier~Thm.(7.3)~with~L=n~variation~
1+~Miller,~J.Phys.A.,Eq.~1.2~(a~=~g)+~ } \\
\mbox{Miller,~J.Phys.A.,Eq.~1.2~(a~=~g)+T:1''} }
}
\end{maplelatex}

\begin{maplelatex}
\mapleinline{inert}{2d}{`  457: `, "Maier Thm.(7.3) with L=n variation 1+ Miller,
J.Phys.A.,Eq. 1.2 (a = g)+ Miller, J.Phys.A.,Eq. 1.2 (a = g)+T:2";}{%
\maplemultiline{
\mathit{\ \ 457:\ }  \mbox{``Maier~Thm.(7.3)~with~L=n~variation~
1+~Miller,~J.Phys.A.,Eq.~1.2~(a~=~g)+~ } \\
\mbox{Miller,~J.Phys.A.,Eq.~1.2~(a~=~g)+T:2''} }
}
\end{maplelatex}

\begin{maplelatex}
\mapleinline{inert}{2d}{`  458: `, "Maier Thm.(7.3) with L=n variation 1+ Miller,
J.Phys.A.,Eq. 1.2 (a = g)+ Miller, J.Phys.A.,Eq. 1.2 (a = g)+T:3";}{%
\maplemultiline{
\mathit{\ \ 458:\ }  \mbox{``Maier~Thm.(7.3)~with~L=n~variation~
1+~Miller,~J.Phys.A.,Eq.~1.2~(a~=~g)+~ } \\
\mbox{Miller,~J.Phys.A.,Eq.~1.2~(a~=~g)+T:3''} }
}
\end{maplelatex}

\begin{maplelatex}
\mapleinline{inert}{2d}{`  459: `, "Maier Thm.(7.3) with L=n variation 1+ Miller,
J.Phys.A.,Eq. 1.2 (a = g)+ Miller, J.Phys.A.,Eq. 1.2 (a = g)+T:4";}{%
\maplemultiline{
\mathit{\ \ 459:\ }  \mbox{``Maier~Thm.(7.3)~with~L=n~variation~
1+~Miller,~J.Phys.A.,Eq.~1.2~(a~=~g)+~ } \\
\mbox{Miller,~J.Phys.A.,Eq.~1.2~(a~=~g)+T:4''} }
}
\end{maplelatex}

\begin{maplelatex}
\mapleinline{inert}{2d}{`  460: `, "Maier Thm.(7.3) with L=n variation 1+ Miller,
J.Phys.A.,Eq. 1.2 (a = g)+ Miller, J.Phys.A.,Eq. 1.2 (a = g)+T:5";}{%
\maplemultiline{
\mathit{\ \ 460:\ }  \mbox{``Maier~Thm.(7.3)~with~L=n~variation~
1+~Miller,~J.Phys.A.,Eq.~1.2~(a~=~g)+~ } \\
\mbox{Miller,~J.Phys.A.,Eq.~1.2~(a~=~g)+T:5''} }
}
\end{maplelatex}

\begin{maplelatex}
\mapleinline{inert}{2d}{`  461: `, "Maier Thm.(7.3) with L=n variation 1+ Miller,
J.Phys.A.,Eq. 1.2 (a = g)+T:3+ Miller, J.Phys.A.,Eq. 1.2 (a =
g)+T:1";}{%
\maplemultiline{
\mathit{\ \ 461:\ }  \mbox{``Maier~Thm.(7.3)~with~L=n~variation~
1+~Miller,~J.Phys.A.,Eq.~1.2~(a~=~g)+} \\
\mbox{T:3+~Miller,~J.Phys.A.,Eq.~1.2~(a~=~g)+T:1''} }
}
\end{maplelatex}

\begin{maplelatex}
\mapleinline{inert}{2d}{`  462: `, "Maier Thm.(7.3) with L=n variation 1+ Miller,
J.Phys.A.,Eq. 1.2 (a = g)+T:3+ Miller, J.Phys.A.,Eq. 1.2 (a =
g)+T:2";}{%
\maplemultiline{
\mathit{\ \ 462:\ }  \mbox{``Maier~Thm.(7.3)~with~L=n~variation~
1+~Miller,~J.Phys.A.,Eq.~1.2~(a~=~g)+} \\
\mbox{T:3+~Miller,~J.Phys.A.,Eq.~1.2~(a~=~g)+T:2''} }
}
\end{maplelatex}

\begin{maplelatex}
\mapleinline{inert}{2d}{`  463: `, "Maier Thm.(7.3) with L=n variation 1+ Miller,
J.Phys.A.,Eq. 1.2 (a = g)+T:3+ Miller, J.Phys.A.,Eq. 1.2 (a =
g)+T:3";}{%
\maplemultiline{
\mathit{\ \ 463:\ }  \mbox{``Maier~Thm.(7.3)~with~L=n~variation~
1+~Miller,~J.Phys.A.,Eq.~1.2~(a~=~g)+} \\
\mbox{T:3+~Miller,~J.Phys.A.,Eq.~1.2~(a~=~g)+T:3''} }
}
\end{maplelatex}

\begin{maplelatex}
\mapleinline{inert}{2d}{`  469: `, "CHU Wenshang+ Miller, J.Phys.A.,Eq. 1.2 (Luke Eq.
23)";}{%
\[
\mathit{\ \ 469:\ }  \,\mbox{``CHU~Wenshang+~Miller,~
J.Phys.A.,Eq.~1.2~(Luke~Eq.~23)''}
\]
}
\end{maplelatex}

\end{maplegroup}
\begin{maplegroup}
\begin{mapleinput}
\end{mapleinput}

\end{maplegroup}

%% file: AppendixD.tex
\begin{maplegroup}
\mapleresult

\begin{maplelatex}
\mapleinline{inert}{2d}{`  20: `, _3F_2([a, b, c],[c+m, b-n],1);}{%
\[
\mathit{\ \ 20:\ } \,\mathrm{_3F_2}([a, \,b, \,c], \,[c + m, \,b
 - n], \,1)
\]
}
\end{maplelatex}

\begin{maplelatex}
\mapleinline{inert}{2d}{`  21: `, _3F_2([a, m, b],[c, -n+m],1);}{%
\[
\mathit{\ \ 21:\ } \,\mathrm{_3F_2}([a, \,m, \,b], \,[c, \, - n + 
m], \,1)
\]
}
\end{maplelatex}

\begin{maplelatex}
\mapleinline{inert}{2d}{`  23: `, _3F_2([a, b, -n],[c, -n+m],1);}{%
\[
\mathit{\ \ 23:\ } \,\mathrm{_3F_2}([a, \,b, \, - n], \,[c, \, - n
 + m], \,1)
\]
}
\end{maplelatex}

\begin{maplelatex}
\mapleinline{inert}{2d}{`  24: `, _3F_2([a, -n, b],[c, a-c+b-n+m],1);}{%
\[
\mathit{\ \ 24:\ } \,\mathrm{_3F_2}([a, \, - n, \,b], \,[c, \,a - 
c + b - n + m], \,1)
\]
}
\end{maplelatex}

\begin{maplelatex}
\mapleinline{inert}{2d}{`  40: `, _3F_2([1-4*a+2*n, 1/2+n-a, n-a],[1-3*a+2*n,
-3*a+2*n+3/2],1);}{%
\[
\mathit{\ \ 40:\ } \,\mathrm{_3F_2}([1 - 4\,a + 2\,n, \,
{\displaystyle \frac {1}{2}}  + n - a, \,n - a], \,[1 - 3\,a + 2
\,n, \, - 3\,a + 2\,n + {\displaystyle \frac {3}{2}} ], \,1)
\]
}
\end{maplelatex}

\begin{maplelatex}
\mapleinline{inert}{2d}{`  41: `, _3F_2([1-4*a+2*n, 1-2*a+n, 1/2-2*a+n],[1-3*a+2*n,
2-4*a+2*n],1);}{%
\[
\mathit{\ \ 41:\ } \,\mathrm{_3F_2}([1 - 4\,a + 2\,n, \,1 - 2\,a
 + n, \,{\displaystyle \frac {1}{2}}  - 2\,a + n], \,[1 - 3\,a + 
2\,n, \,2 - 4\,a + 2\,n], \,1)
\]
}
\end{maplelatex}

\begin{maplelatex}
\mapleinline{inert}{2d}{`  42: `, _3F_2([1/2+n-a, 1-2*a+n, a],[1-3*a+2*n, 3/2-a+n],1);}{%
\[
\mathit{\ \ 42:\ } \,\mathrm{_3F_2}([{\displaystyle \frac {1}{2}} 
 + n - a, \,1 - 2\,a + n, \,a], \,[1 - 3\,a + 2\,n, \,
{\displaystyle \frac {3}{2}}  - a + n], \,1)
\]
}
\end{maplelatex}

\begin{maplelatex}
\mapleinline{inert}{2d}{`  43: `, _3F_2([n-a, 1/2-2*a+n, a],[1-3*a+2*n, n+1-a],1);}{%
\[
\mathit{\ \ 43:\ } \,\mathrm{_3F_2}([n - a, \,{\displaystyle 
\frac {1}{2}}  - 2\,a + n, \,a], \,[1 - 3\,a + 2\,n, \,n + 1 - a]
, \,1)
\]
}
\end{maplelatex}

\begin{maplelatex}
\mapleinline{inert}{2d}{`  44: `, _3F_2([1/2+n-a, 3/2-2*a+n, a+1/2],[-3*a+2*n+3/2,
3/2-a+n],1);}{%
\[
\mathit{\ \ 44:\ } \,\mathrm{_3F_2}([{\displaystyle \frac {1}{2}} 
 + n - a, \,{\displaystyle \frac {3}{2}}  - 2\,a + n, \,a + 
{\displaystyle \frac {1}{2}} ], \,[ - 3\,a + 2\,n + 
{\displaystyle \frac {3}{2}} , \,{\displaystyle \frac {3}{2}}  - 
a + n], \,1)
\]
}
\end{maplelatex}

\begin{maplelatex}
\mapleinline{inert}{2d}{`  45: `, _3F_2([n-a, 1-2*a+n, a+1/2],[-3*a+2*n+3/2, n+1-a],1);}{%
\[
\mathit{\ \ 45:\ } \,\mathrm{_3F_2}([n - a, \,1 - 2\,a + n, \,a + 
{\displaystyle \frac {1}{2}} ], \,[ - 3\,a + 2\,n + 
{\displaystyle \frac {3}{2}} , \,n + 1 - a], \,1)
\]
}
\end{maplelatex}

\begin{maplelatex}
\mapleinline{inert}{2d}{`  46: `, _3F_2([1/2-2*a+n, 1-2*a+n, 1],[2-4*a+2*n, n+1-a],1);}{%
\[
\mathit{\ \ 46:\ } \,\mathrm{_3F_2}([{\displaystyle \frac {1}{2}} 
 - 2\,a + n, \,1 - 2\,a + n, \,1], \,[2 - 4\,a + 2\,n, \,n + 1 - 
a], \,1)
\]
}
\end{maplelatex}

\begin{maplelatex}
\mapleinline{inert}{2d}{`  47: `, _3F_2([a, a+1/2, 1],[3/2-a+n, n+1-a],1);}{%
\[
\mathit{\ \ 47:\ } \,\mathrm{_3F_2}([a, \,a + {\displaystyle 
\frac {1}{2}} , \,1], \,[{\displaystyle \frac {3}{2}}  - a + n, 
\,n + 1 - a], \,1)
\]
}
\end{maplelatex}

\begin{maplelatex}
\mapleinline{inert}{2d}{`  48: `, _3F_2([a, b, c],[e, c-1],1);}{%
\[
\mathit{\ \ 48:\ } \,\mathrm{_3F_2}([a, \,b, \,c], \,[e, \,c - 1]
, \,1)
\]
}
\end{maplelatex}

\begin{maplelatex}
\mapleinline{inert}{2d}{`  50: `, _3F_2([a, -1, b],[c, e],1);}{%
\[
\mathit{\ \ 50:\ } \,\mathrm{_3F_2}([a, \,-1, \,b], \,[c, \,e], \,
1)
\]
}
\end{maplelatex}

\begin{maplelatex}
\mapleinline{inert}{2d}{`  51: `, _3F_2([a, b, c],[c+2, 2+b],1);}{%
\[
\mathit{\ \ 51:\ } \,\mathrm{_3F_2}([a, \,b, \,c], \,[c + 2, \,2
 + b], \,1)
\]
}
\end{maplelatex}

\begin{maplelatex}
\mapleinline{inert}{2d}{`  52: `, _3F_2([a, 2, b],[c, 4],1);}{%
\[
\mathit{\ \ 52:\ } \,\mathrm{_3F_2}([a, \,2, \,b], \,[c, \,4], \,1
)
\]
}
\end{maplelatex}

\begin{maplelatex}
\mapleinline{inert}{2d}{`  53: `, _3F_2([a, 2, b],[c, 4+a-c+b],1);}{%
\[
\mathit{\ \ 53:\ } \,\mathrm{_3F_2}([a, \,2, \,b], \,[c, \,4 + a
 - c + b], \,1)
\]
}
\end{maplelatex}

\begin{maplelatex}
\mapleinline{inert}{2d}{`  54: `, _3F_2([-n, a, b],[b-m, a-n+m],1);}{%
\[
\mathit{\ \ 54:\ } \,\mathrm{_3F_2}([ - n, \,a, \,b], \,[b - m, \,
a - n + m], \,1)
\]
}
\end{maplelatex}

\begin{maplelatex}
\mapleinline{inert}{2d}{`  55: `, _3F_2([-L-m+n, a, b],[b-m, a-L],1);}{%
\[
\mathit{\ \ 55:\ } \,\mathrm{_3F_2}([ - L - m + n, \,a, \,b], \,[b
 - m, \,a - L], \,1)
\]
}
\end{maplelatex}

\begin{maplelatex}
\mapleinline{inert}{2d}{`  56: `, _3F_2([a, b, b+L-n],[b-m, b-n],1);}{%
\[
\mathit{\ \ 56:\ } \,\mathrm{_3F_2}([a, \,b, \,b + L - n], \,[b - 
m, \,b - n], \,1)
\]
}
\end{maplelatex}

\begin{maplelatex}
\mapleinline{inert}{2d}{`  57: `, _3F_2([-n, a, b],[b-m, a-L],1);}{%
\[
\mathit{\ \ 57:\ } \,\mathrm{_3F_2}([ - n, \,a, \,b], \,[b - m, \,
a - L], \,1)
\]
}
\end{maplelatex}

\begin{maplelatex}
\mapleinline{inert}{2d}{`  61: `, _3F_2([-n, a, b],[b-m, a-2-n+m],1);}{%
\[
\mathit{\ \ 61:\ } \,\mathrm{_3F_2}([ - n, \,a, \,b], \,[b - m, \,
a - 2 - n + m], \,1)
\]
}
\end{maplelatex}

\begin{maplelatex}
\mapleinline{inert}{2d}{`  62: `, _3F_2([a, b, -n*(n-b+a)/(b-n)],[-(a*n+b*n-b+n-b^2)/(b-n),
n+1+a],1);}{%
\[
\mathit{\ \ 62:\ } \,\mathrm{_3F_2}([a, \,b, \, - {\displaystyle 
\frac {n\,(n - b + a)}{b - n}} ], \,[ - {\displaystyle \frac {a\,
n + b\,n - b + n - b^{2}}{b - n}} , \,n + 1 + a], \,1)
\]
}
\end{maplelatex}

\begin{maplelatex}
\mapleinline{inert}{2d}{`  66: `, _3F_2([3-4*a-2*n, -a-n+1, -a+3/2-n],[3-3*a-2*n,
7/2-3*a-2*n],1);}{%
\[
\mathit{\ \ 66:\ } \,\mathrm{_3F_2}([3 - 4\,a - 2\,n, \, - a - n
 + 1, \, - a + {\displaystyle \frac {3}{2}}  - n], \,[3 - 3\,a - 
2\,n, \,{\displaystyle \frac {7}{2}}  - 3\,a - 2\,n], \,1)
\]
}
\end{maplelatex}

\begin{maplelatex}
\mapleinline{inert}{2d}{`  67: `, _3F_2([-a-n+1, 3/2-2*a-n, a],[3-3*a-2*n, 2-a-n],1);}{%
\[
\mathit{\ \ 67:\ } \,\mathrm{_3F_2}([ - a - n + 1, \,
{\displaystyle \frac {3}{2}}  - 2\,a - n, \,a], \,[3 - 3\,a - 2\,
n, \,2 - a - n], \,1)
\]
}
\end{maplelatex}

\begin{maplelatex}
\mapleinline{inert}{2d}{`  68: `, _3F_2([-a+3/2-n, -2*a-n+2, a],[3-3*a-2*n, 5/2-a-n],1);}{%
\[
\mathit{\ \ 68:\ } \,\mathrm{_3F_2}([ - a + {\displaystyle \frac {
3}{2}}  - n, \, - 2\,a - n + 2, \,a], \,[3 - 3\,a - 2\,n, \,
{\displaystyle \frac {5}{2}}  - a - n], \,1)
\]
}
\end{maplelatex}

\begin{maplelatex}
\mapleinline{inert}{2d}{`  69: `, _3F_2([3-4*a-2*n, -2*a-n+2, 5/2-2*a-n],[7/2-3*a-2*n,
4-4*a-2*n],1);}{%
\[
\mathit{\ \ 69:\ } \,\mathrm{_3F_2}([3 - 4\,a - 2\,n, \, - 2\,a - 
n + 2, \,{\displaystyle \frac {5}{2}}  - 2\,a - n], \,[
{\displaystyle \frac {7}{2}}  - 3\,a - 2\,n, \,4 - 4\,a - 2\,n], 
\,1)
\]
}
\end{maplelatex}

\begin{maplelatex}
\mapleinline{inert}{2d}{`  70: `, _3F_2([-a-n+1, -2*a-n+2, a+1/2],[7/2-3*a-2*n, 2-a-n],1);}{%
\[
\mathit{\ \ 70:\ } \,\mathrm{_3F_2}([ - a - n + 1, \, - 2\,a - n
 + 2, \,a + {\displaystyle \frac {1}{2}} ], \,[{\displaystyle 
\frac {7}{2}}  - 3\,a - 2\,n, \,2 - a - n], \,1)
\]
}
\end{maplelatex}

\begin{maplelatex}
\mapleinline{inert}{2d}{`  71: `, _3F_2([-a+3/2-n, 5/2-2*a-n, a+1/2],[7/2-3*a-2*n,
5/2-a-n],1);}{%
\[
\mathit{\ \ 71:\ } \,\mathrm{_3F_2}([ - a + {\displaystyle \frac {
3}{2}}  - n, \,{\displaystyle \frac {5}{2}}  - 2\,a - n, \,a + 
{\displaystyle \frac {1}{2}} ], \,[{\displaystyle \frac {7}{2}} 
 - 3\,a - 2\,n, \,{\displaystyle \frac {5}{2}}  - a - n], \,1)
\]
}
\end{maplelatex}

\begin{maplelatex}
\mapleinline{inert}{2d}{`  72: `, _3F_2([-2*a-n+2, 5/2-2*a-n, 1],[4-4*a-2*n, 5/2-a-n],1);}{%
\[
\mathit{\ \ 72:\ } \,\mathrm{_3F_2}([ - 2\,a - n + 2, \,
{\displaystyle \frac {5}{2}}  - 2\,a - n, \,1], \,[4 - 4\,a - 2\,
n, \,{\displaystyle \frac {5}{2}}  - a - n], \,1)
\]
}
\end{maplelatex}

\begin{maplelatex}
\mapleinline{inert}{2d}{`  73: `, _3F_2([a, a+1/2, 1],[2-a-n, 5/2-a-n],1);}{%
\[
\mathit{\ \ 73:\ } \,\mathrm{_3F_2}([a, \,a + {\displaystyle 
\frac {1}{2}} , \,1], \,[2 - a - n, \,{\displaystyle \frac {5}{2}
}  - a - n], \,1)
\]
}
\end{maplelatex}

\begin{maplelatex}
\mapleinline{inert}{2d}{`  74: `, _3F_2([-4*a+2*n, 1/2-2*a+n, n-2*a],[-3*a+2*n,
1-4*a+2*n],1);}{%
\[
\mathit{\ \ 74:\ } \,\mathrm{_3F_2}([ - 4\,a + 2\,n, \,
{\displaystyle \frac {1}{2}}  - 2\,a + n, \,n - 2\,a], \,[ - 3\,a
 + 2\,n, \,1 - 4\,a + 2\,n], \,1)
\]
}
\end{maplelatex}

\begin{maplelatex}
\mapleinline{inert}{2d}{`  75: `, _3F_2([-a+n-1/2, n-2*a, a],[-3*a+2*n, 1/2+n-a],1);}{%
\[
\mathit{\ \ 75:\ } \,\mathrm{_3F_2}([ - a + n - {\displaystyle 
\frac {1}{2}} , \,n - 2\,a, \,a], \,[ - 3\,a + 2\,n, \,
{\displaystyle \frac {1}{2}}  + n - a], \,1)
\]
}
\end{maplelatex}

\begin{maplelatex}
\mapleinline{inert}{2d}{`  76: `, _3F_2([-4*a+2*n, 1-2*a+n, 1/2-2*a+n],[-3*a+2*n+1/2,
1-4*a+2*n],1);}{%
\[
\mathit{\ \ 76:\ } \,\mathrm{_3F_2}([ - 4\,a + 2\,n, \,1 - 2\,a + 
n, \,{\displaystyle \frac {1}{2}}  - 2\,a + n], \,[ - 3\,a + 2\,n
 + {\displaystyle \frac {1}{2}} , \,1 - 4\,a + 2\,n], \,1)
\]
}
\end{maplelatex}

\begin{maplelatex}
\mapleinline{inert}{2d}{`  77: `, _3F_2([n-a, 1-2*a+n, a+1/2],[-3*a+2*n+1/2, n+1-a],1);}{%
\[
\mathit{\ \ 77:\ } \,\mathrm{_3F_2}([n - a, \,1 - 2\,a + n, \,a + 
{\displaystyle \frac {1}{2}} ], \,[ - 3\,a + 2\,n + 
{\displaystyle \frac {1}{2}} , \,n + 1 - a], \,1)
\]
}
\end{maplelatex}

\begin{maplelatex}
\mapleinline{inert}{2d}{`  78: `, _3F_2([1/2-2*a+n, 1-2*a+n, 1],[1-4*a+2*n, n+1-a],1);}{%
\[
\mathit{\ \ 78:\ } \,\mathrm{_3F_2}([{\displaystyle \frac {1}{2}} 
 - 2\,a + n, \,1 - 2\,a + n, \,1], \,[1 - 4\,a + 2\,n, \,n + 1 - 
a], \,1)
\]
}
\end{maplelatex}

\begin{maplelatex}
\mapleinline{inert}{2d}{`  79: `, _3F_2([n-2*a, 1/2-2*a+n, 1],[1-4*a+2*n, 1/2+n-a],1);}{%
\[
\mathit{\ \ 79:\ } \,\mathrm{_3F_2}([n - 2\,a, \,{\displaystyle 
\frac {1}{2}}  - 2\,a + n, \,1], \,[1 - 4\,a + 2\,n, \,
{\displaystyle \frac {1}{2}}  + n - a], \,1)
\]
}
\end{maplelatex}

\begin{maplelatex}
\mapleinline{inert}{2d}{`  90: `, _3F_2([3*n-3/2, -a+n-1/2, 2*a],[2*a+1/2+n, 3*n-a-1/2],1);}{%
\[
\mathit{\ \ 90:\ } \,\mathrm{_3F_2}([3\,n - {\displaystyle \frac {
3}{2}} , \, - a + n - {\displaystyle \frac {1}{2}} , \,2\,a], \,[
2\,a + {\displaystyle \frac {1}{2}}  + n, \,3\,n - a - 
{\displaystyle \frac {1}{2}} ], \,1)
\]
}
\end{maplelatex}

\begin{maplelatex}
\mapleinline{inert}{2d}{`  93: `, _3F_2([a, b,
c],[1/2*c+1+1/2*b+1/2*a+1/2*(a^2-2*a*c-2*a*b+c^2-2*b*c+b^2)^(1/2),
1/2*c+1+1/2*b+1/2*a-1/2*(a^2-2*a*c-2*a*b+c^2-2*b*c+b^2)^(1/2)],1);}{%
\maplemultiline{
\mathit{\ \ 93:\ }, \mathrm{_3F_2}([a, \,b, \,c], [{\displaystyle 
\frac {c}{2}}  + 1 + {\displaystyle \frac {b}{2}}  + 
{\displaystyle \frac {a}{2}}  + {\displaystyle \frac {\sqrt{a^{2}
 - 2\,a\,c - 2\,a\,b + c^{2} - 2\,b\,c + b^{2}}}{2}} ,  \\
{\displaystyle \frac {c}{2}}  + 1 + {\displaystyle \frac {b}{2}} 
 + {\displaystyle \frac {a}{2}}  - {\displaystyle \frac {\sqrt{a
^{2} - 2\,a\,c - 2\,a\,b + c^{2} - 2\,b\,c + b^{2}}}{2}} ], \,1)
 }
}
\end{maplelatex}

\begin{maplelatex}
\mapleinline{inert}{2d}{`  94: `, _3F_2([3/2, a, a-1/2],[2*a-1/6, 1/6+2*a],1);}{%
\[
\mathit{\ \ 94:\ } \,\mathrm{_3F_2}([{\displaystyle \frac {3}{2}} 
, \,a, \,a - {\displaystyle \frac {1}{2}} ], \,[2\,a - 
{\displaystyle \frac {1}{6}} , \,{\displaystyle \frac {1}{6}}  + 
2\,a], \,1)
\]
}
\end{maplelatex}

\begin{maplelatex}
\mapleinline{inert}{2d}{`  95: `, _3F_2([3/2, a, a-1/2],[2*a-5/6, 2*a-1/6],1);}{%
\[
\mathit{\ \ 95:\ } \,\mathrm{_3F_2}([{\displaystyle \frac {3}{2}} 
, \,a, \,a - {\displaystyle \frac {1}{2}} ], \,[2\,a - 
{\displaystyle \frac {5}{6}} , \,2\,a - {\displaystyle \frac {1}{
6}} ], \,1)
\]
}
\end{maplelatex}

\begin{maplelatex}
\mapleinline{inert}{2d}{`  96: `, _3F_2([a, a+1/3, 2*a-5/3],[2*a-1/6, 3*a-1],1);}{%
\[
\mathit{\ \ 96:\ } \,\mathrm{_3F_2}([a, \,a + {\displaystyle 
\frac {1}{3}} , \,2\,a - {\displaystyle \frac {5}{3}} ], \,[2\,a
 - {\displaystyle \frac {1}{6}} , \,3\,a - 1], \,1)
\]
}
\end{maplelatex}

\begin{maplelatex}
\mapleinline{inert}{2d}{`  97: `, _3F_2([a, a+1/3, 2*a-2/3],[2*a+5/6, 3*a],1);}{%
\[
\mathit{\ \ 97:\ } \,\mathrm{_3F_2}([a, \,a + {\displaystyle 
\frac {1}{3}} , \,2\,a - {\displaystyle \frac {2}{3}} ], \,[2\,a
 + {\displaystyle \frac {5}{6}} , \,3\,a], \,1)
\]
}
\end{maplelatex}

\begin{maplelatex}
\mapleinline{inert}{2d}{`  98: `, _3F_2([3/2, a, a-1/2],[-7/6+2*a, 2*a-5/6],1);}{%
\[
\mathit{\ \ 98:\ } \,\mathrm{_3F_2}([{\displaystyle \frac {3}{2}} 
, \,a, \,a - {\displaystyle \frac {1}{2}} ], \,[ - 
{\displaystyle \frac {7}{6}}  + 2\,a, \,2\,a - {\displaystyle 
\frac {5}{6}} ], \,1)
\]
}
\end{maplelatex}

\begin{maplelatex}
\mapleinline{inert}{2d}{`  99: `, _3F_2([a, a+2/3, 2*a-4/3],[1/6+2*a, 3*a-1],1);}{%
\[
\mathit{\ \ 99:\ } \,\mathrm{_3F_2}([a, \,a + {\displaystyle 
\frac {2}{3}} , \,2\,a - {\displaystyle \frac {4}{3}} ], \,[
{\displaystyle \frac {1}{6}}  + 2\,a, \,3\,a - 1], \,1)
\]
}
\end{maplelatex}

\begin{maplelatex}
\mapleinline{inert}{2d}{`  100: `, _3F_2([a, a+2/3, 2*a-1/3],[7/6+2*a, 3*a],1);}{%
\[
\mathit{\ \ 100:\ } \,\mathrm{_3F_2}([a, \,a + {\displaystyle 
\frac {2}{3}} , \,2\,a - {\displaystyle \frac {1}{3}} ], \,[
{\displaystyle \frac {7}{6}}  + 2\,a, \,3\,a], \,1)
\]
}
\end{maplelatex}

\begin{maplelatex}
\mapleinline{inert}{2d}{`  101: `, _3F_2([a, a+1/3, 2*a-5/3],[2*a-1/6, 3*a-2],1);}{%
\[
\mathit{\ \ 101:\ } \,\mathrm{_3F_2}([a, \,a + {\displaystyle 
\frac {1}{3}} , \,2\,a - {\displaystyle \frac {5}{3}} ], \,[2\,a
 - {\displaystyle \frac {1}{6}} , \,3\,a - 2], \,1)
\]
}
\end{maplelatex}

\begin{maplelatex}
\mapleinline{inert}{2d}{`  102: `, _3F_2([a, a+1/3, 2*a-2/3],[2*a+5/6, 3*a-1],1);}{%
\[
\mathit{\ \ 102:\ } \,\mathrm{_3F_2}([a, \,a + {\displaystyle 
\frac {1}{3}} , \,2\,a - {\displaystyle \frac {2}{3}} ], \,[2\,a
 + {\displaystyle \frac {5}{6}} , \,3\,a - 1], \,1)
\]
}
\end{maplelatex}

\begin{maplelatex}
\mapleinline{inert}{2d}{`  103: `, _3F_2([a, a+1/3, a+2/3],[3/2*a+3/2, 1+3/2*a],1);}{%
\[
\mathit{\ \ 103:\ } \,\mathrm{_3F_2}([a, \,a + {\displaystyle 
\frac {1}{3}} , \,a + {\displaystyle \frac {2}{3}} ], \,[
{\displaystyle \frac {3\,a}{2}}  + {\displaystyle \frac {3}{2}} 
, \,1 + {\displaystyle \frac {3\,a}{2}} ], \,1)
\]
}
\end{maplelatex}

\begin{maplelatex}
\mapleinline{inert}{2d}{`  135: `, _3F_2([a, 4*a-2*n, 2*a-1/2-n],[3*a-n, 1/2+2*a-n],1);}{%
\[
\mathit{\ \ 135:\ } \,\mathrm{_3F_2}([a, \,4\,a - 2\,n, \,2\,a - 
{\displaystyle \frac {1}{2}}  - n], \,[3\,a - n, \,
{\displaystyle \frac {1}{2}}  + 2\,a - n], \,1)
\]
}
\end{maplelatex}

\begin{maplelatex}
\mapleinline{inert}{2d}{`  136: `, _3F_2([a, 2*a-n, 4*a-2*n],[2*a+1-n, 3*a+1/2-n],1);}{%
\[
\mathit{\ \ \dagger 136:\ } \,\mathrm{_3F_2}([a, \,2\,a - n, \,4\,a - 2\,n
], \,[2\,a + 1 - n, \,3\,a + {\displaystyle \frac {1}{2}}  - n], 
\,1)
\]
}
\end{maplelatex}

\begin{maplelatex}
\mapleinline{inert}{2d}{`  137: `, _3F_2([2*a-n, 4*a-2*n, a+1/2],[3*a-n, 2*a+1-n],1);}{%
\[
\mathit{\ \ \dagger 137:\ } \,\mathrm{_3F_2}([2\,a - n, \,4\,a - 2\,n, \,a
 + {\displaystyle \frac {1}{2}} ], \,[3\,a - n, \,2\,a + 1 - n], 
\,1)
\]
}
\end{maplelatex}

\begin{maplelatex}
\mapleinline{inert}{2d}{`  138: `, _3F_2([1, 4*a-2*n, a-n+1],[2*a+1-n, 2*a+3/2-n],1);}{%
\[
\mathit{\ \ 138:\ } \,\mathrm{_3F_2}([1, \,4\,a - 2\,n, \,a - n + 
1], \,[2\,a + 1 - n, \,2\,a + {\displaystyle \frac {3}{2}}  - n]
, \,1)
\]
}
\end{maplelatex}

\begin{maplelatex}
\mapleinline{inert}{2d}{`  140: `, _3F_2([a, -2+4*a+2*n, -1+n+2*a],[2*a+n, 3*a-1/2+n],1);}{%
\[
\mathit{\ \ \dagger 140:\ } \,\mathrm{_3F_2}([a, \, - 2 + 4\,a + 2\,n, \,
 - 1 + n + 2\,a], \,[2\,a + n, \,3\,a - {\displaystyle \frac {1}{
2}}  + n], \,1)
\]
}
\end{maplelatex}

\begin{maplelatex}
\mapleinline{inert}{2d}{`  141: `, _3F_2([a+1/2, -2+4*a+2*n, 2*a-1/2+n],[2*a+1/2+n,
3*a-1/2+n],1);}{%
\[
\mathit{\ \ 141:\ } \,\mathrm{_3F_2}([a + {\displaystyle \frac {1
}{2}} , \, - 2 + 4\,a + 2\,n, \,2\,a - {\displaystyle \frac {1}{2
}}  + n], \,[2\,a + {\displaystyle \frac {1}{2}}  + n, \,3\,a - 
{\displaystyle \frac {1}{2}}  + n], \,1)
\]
}
\end{maplelatex}

\begin{maplelatex}
\mapleinline{inert}{2d}{`  142: `, _3F_2([a+1/2, -2+4*a+2*n, -1+n+2*a],[2*a+n, 3*a+n-1],1);}{%
\[
\mathit{\ \ \dagger 142:\ } \,\mathrm{_3F_2}([a + {\displaystyle \frac {1
}{2}} , \, - 2 + 4\,a + 2\,n, \, - 1 + n + 2\,a], \,[2\,a + n, \,
3\,a + n - 1], \,1)
\]
}
\end{maplelatex}

\begin{maplelatex}
\mapleinline{inert}{2d}{`  143: `, _3F_2([1, -2+4*a+2*n, n-1/2+a],[2*a+n, 2*a-1/2+n],1);}{%
\[
\mathit{\ \ 143:\ } \,\mathrm{_3F_2}([1, \, - 2 + 4\,a + 2\,n, \,n
 - {\displaystyle \frac {1}{2}}  + a], \,[2\,a + n, \,2\,a - 
{\displaystyle \frac {1}{2}}  + n], \,1)
\]
}
\end{maplelatex}

\begin{maplelatex}
\mapleinline{inert}{2d}{`  144: `, _3F_2([1, a-n+1, 1+4*a-2*n],[2*a+1-n, 2*a+3/2-n],1);}{%
\[
\mathit{\ \ 144:\ } \,\mathrm{_3F_2}([1, \,a - n + 1, \,1 + 4\,a
 - 2\,n], \,[2\,a + 1 - n, \,2\,a + {\displaystyle \frac {3}{2}} 
 - n], \,1)
\]
}
\end{maplelatex}

\begin{maplelatex}
\mapleinline{inert}{2d}{`  145: `, _3F_2([1, 1+4*a-2*n, 3/2-n+a],[2-n+2*a, 2*a+3/2-n],1);}{%
\[
\mathit{\ \ 145:\ } \,\mathrm{_3F_2}([1, \,1 + 4\,a - 2\,n, \,
{\displaystyle \frac {3}{2}}  - n + a], \,[2 - n + 2\,a, \,2\,a
 + {\displaystyle \frac {3}{2}}  - n], \,1)
\]
}
\end{maplelatex}

\begin{maplelatex}
\mapleinline{inert}{2d}{`  146: `, _3F_2([n-a, 2*a-n, 2*a-1/2-n],[1/2, 3*a-n],1);}{%
\[
\mathit{\ \ \dagger 146:\ } \,\mathrm{_3F_2}([n - a, \,2\,a - n, \,2
\,a - {\displaystyle \frac {1}{2}}  - n], \,[{\displaystyle 
\frac {1}{2}} , \,3\,a - n], \,1)
\]
}
\end{maplelatex}

\begin{maplelatex}
\mapleinline{inert}{2d}{`  147: `, _3F_2([-1+n+2*a, -a+3/2-n, 2*a-1/2+n],[3/2,
3*a-1/2+n],1);}{%
\[
\mathit{\ \ \dagger 147:\ } \,\mathrm{_3F_2}([ - 1 + n + 2\,a, \, - 
a + {\displaystyle \frac {3}{2}}  - n, \,2\,a - {\displaystyle 
\frac {1}{2}}  + n], \,[{\displaystyle \frac {3}{2}} , \,3\,a - 
{\displaystyle \frac {1}{2}}  + n], \,1)
\]
}
\end{maplelatex}

\begin{maplelatex}
\mapleinline{inert}{2d}{`  153: `, _3F_2([1-2*a+n, 1, a+1/2],[n+1-a, 2*a+1-n],1);}{%
\[
\mathit{\ \ \dagger 153:\ } \,\mathrm{_3F_2}([1 - 2\,a + n, \,1, \,a + 
{\displaystyle \frac {1}{2}} ], \,[n + 1 - a, \,2\,a + 1 - n], \,
1)
\]
}
\end{maplelatex}

\begin{maplelatex}
\mapleinline{inert}{2d}{`  154: `, _3F_2([1-2*a+n, 1/2-2*a+n, n-a],[n+1-a, 1/2],1);}{%
\[
\mathit{\ \ \dagger 154:\ } \,\mathrm{_3F_2}([1 - 2\,a + n, \,
{\displaystyle \frac {1}{2}}  - 2\,a + n, \,n - a], \,[n + 1 - a
, \,{\displaystyle \frac {1}{2}} ], \,1)
\]
}
\end{maplelatex}

\begin{maplelatex}
\mapleinline{inert}{2d}{`  155: `, _3F_2([1, 1/2-2*a+n, a],[n+1-a, 1/2+2*a-n],1);}{%
\[
\mathit{\ \ 155:\ } \,\mathrm{_3F_2}([1, \,{\displaystyle \frac {1
}{2}}  - 2\,a + n, \,a], \,[n + 1 - a, \,{\displaystyle \frac {1
}{2}}  + 2\,a - n], \,1)
\]
}
\end{maplelatex}

\begin{maplelatex}
\mapleinline{inert}{2d}{`  156: `, _3F_2([a+1/2, n-a, a],[n+1-a, 3*a-n],1);}{%
\[
\mathit{\ \ 156:\ } \,\mathrm{_3F_2}([a + {\displaystyle \frac {1
}{2}} , \,n - a, \,a], \,[n + 1 - a, \,3\,a - n], \,1)
\]
}
\end{maplelatex}

\begin{maplelatex}
\mapleinline{inert}{2d}{`  157: `, _3F_2([1/2-2*a+n, a-n+1/2, 2*a-1/2-n],[1/2,
1/2+2*a-n],1);}{%
\[
\mathit{\ \ 157:\ } \,\mathrm{_3F_2}([{\displaystyle \frac {1}{2}
}  - 2\,a + n, \,a - n + {\displaystyle \frac {1}{2}} , \,2\,a - 
{\displaystyle \frac {1}{2}}  - n], \,[{\displaystyle \frac {1}{2
}} , \,{\displaystyle \frac {1}{2}}  + 2\,a - n], \,1)
\]
}
\end{maplelatex}

\begin{maplelatex}
\mapleinline{inert}{2d}{`  158: `, _3F_2([1-2*a+n, 1, a],[3/2-a+n, 2*a+1-n],1);}{%
\[
\mathit{\ \ \dagger 158:\ } \,\mathrm{_3F_2}([1 - 2\,a + n, \,1, \,a], \,[
{\displaystyle \frac {3}{2}}  - a + n, \,2\,a + 1 - n], \,1)
\]
}
\end{maplelatex}

\begin{maplelatex}
\mapleinline{inert}{2d}{`  159: `, _3F_2([1-2*a+n, a-n+1, 2*a-n],[3/2, 2*a+1-n],1);}{%
\[
\mathit{\ \ 159:\ } \,\mathrm{_3F_2}([1 - 2\,a + n, \,a - n + 1, 
\,2\,a - n], \,[{\displaystyle \frac {3}{2}} , \,2\,a + 1 - n], 
\,1)
\]
}
\end{maplelatex}

\begin{maplelatex}
\mapleinline{inert}{2d}{`  169: `, _3F_2([5/2-2*a-n, a+1/2, 1],[5/2-a-n, 2*a+1/2+n],1);}{%
\[
\mathit{\ \ 169:\ } \,\mathrm{_3F_2}([{\displaystyle \frac {5}{2}
}  - 2\,a - n, \,a + {\displaystyle \frac {1}{2}} , \,1], \,[
{\displaystyle \frac {5}{2}}  - a - n, \,2\,a + {\displaystyle 
\frac {1}{2}}  + n], \,1)
\]
}
\end{maplelatex}

\begin{maplelatex}
\mapleinline{inert}{2d}{`  170: `, _3F_2([5/2-2*a-n, -a+3/2-n, -2*a-n+2],[5/2-a-n, 3/2],1);}{%
\[
\mathit{\ \ \dagger 170:\ } \,\mathrm{_3F_2}([{\displaystyle \frac {5}{2}
}  - 2\,a - n, \, - a + {\displaystyle \frac {3}{2}}  - n, \, - 2
\,a - n + 2], \,[{\displaystyle \frac {5}{2}}  - a - n, \,
{\displaystyle \frac {3}{2}} ], \,1)
\]
}
\end{maplelatex}

\begin{maplelatex}
\mapleinline{inert}{2d}{`  171: `, _3F_2([a+1/2, -a+3/2-n, a],[5/2-a-n, 3*a-1/2+n],1);}{%
\[
\mathit{\ \ 171:\ } \,\mathrm{_3F_2}([a + {\displaystyle \frac {1
}{2}} , \, - a + {\displaystyle \frac {3}{2}}  - n, \,a], \,[
{\displaystyle \frac {5}{2}}  - a - n, \,3\,a - {\displaystyle 
\frac {1}{2}}  + n], \,1)
\]
}
\end{maplelatex}

\begin{maplelatex}
\mapleinline{inert}{2d}{`  172: `, _3F_2([1, -2*a-n+2, a],[5/2-a-n, 2*a+n],1);}{%
\[
\mathit{\ \ \dagger 172:\ } \,\mathrm{_3F_2}([1, \, - 2\,a - n + 2, \,a], 
\,[{\displaystyle \frac {5}{2}}  - a - n, \,2\,a + n], \,1)
\]
}
\end{maplelatex}

\begin{maplelatex}
\mapleinline{inert}{2d}{`  173: `, _3F_2([5/2-2*a-n, 2*a-1/2+n, a+n],[2*a+1/2+n, 3/2],1);}{%
\[
\mathit{\ \ 173:\ } \,\mathrm{_3F_2}([{\displaystyle \frac {5}{2}
}  - 2\,a - n, \,2\,a - {\displaystyle \frac {1}{2}}  + n, \,a + 
n], \,[2\,a + {\displaystyle \frac {1}{2}}  + n, \,
{\displaystyle \frac {3}{2}} ], \,1)
\]
}
\end{maplelatex}

\begin{maplelatex}
\mapleinline{inert}{2d}{`  174: `, _3F_2([1, -2*a-n+2, a+1/2],[2-a-n, 2*a+n],1);}{%
\[
\mathit{\ \ \dagger 174:\ } \,\mathrm{_3F_2}([1, \, - 2\,a - n + 2, \,a + 
{\displaystyle \frac {1}{2}} ], \,[2 - a - n, \,2\,a + n], \,1)
\]
}
\end{maplelatex}

\begin{maplelatex}
\mapleinline{inert}{2d}{`  175: `, _3F_2([-2*a-n+2, n-1/2+a, -1+n+2*a],[1/2, 2*a+n],1);}{%
\[
\mathit{\ \ 175:\ } \,\mathrm{_3F_2}([ - 2\,a - n + 2, \,n - 
{\displaystyle \frac {1}{2}}  + a, \, - 1 + n + 2\,a], \,[
{\displaystyle \frac {1}{2}} , \,2\,a + n], \,1)
\]
}
\end{maplelatex}

\begin{maplelatex}
\mapleinline{inert}{2d}{`  176: `, _3F_2([-a+n-1/2, 1/2-2*a+n, n-2*a],[1/2+n-a, 1/2],1);}{%
\[
\mathit{\ \ \dagger 176:\ } \,\mathrm{_3F_2}([ - a + n - {\displaystyle 
\frac {1}{2}} , \,{\displaystyle \frac {1}{2}}  - 2\,a + n, \,n
 - 2\,a], \,[{\displaystyle \frac {1}{2}}  + n - a, \,
{\displaystyle \frac {1}{2}} ], \,1)
\]
}
\end{maplelatex}

\begin{maplelatex}
\mapleinline{inert}{2d}{`  177: `, _3F_2([-a+n-1/2, a+1/2, a],[1/2+n-a, 3*a+1/2-n],1);}{%
\[
\mathit{\ \ 177:\ } \,\mathrm{_3F_2}([ - a + n - {\displaystyle 
\frac {1}{2}} , \,a + {\displaystyle \frac {1}{2}} , \,a], \,[
{\displaystyle \frac {1}{2}}  + n - a, \,3\,a + {\displaystyle 
\frac {1}{2}}  - n], \,1)
\]
}
\end{maplelatex}

\begin{maplelatex}
\mapleinline{inert}{2d}{`  178: `, _3F_2([1/2-2*a+n, a+1/2, 1],[1/2+n-a, 2*a+3/2-n],1);}{%
\[
\mathit{\ \ \dagger 178:\ } \,\mathrm{_3F_2}([{\displaystyle \frac {1}{2}
}  - 2\,a + n, \,a + {\displaystyle \frac {1}{2}} , \,1], \,[
{\displaystyle \frac {1}{2}}  + n - a, \,2\,a + {\displaystyle 
\frac {3}{2}}  - n], \,1)
\]
}
\end{maplelatex}

\begin{maplelatex}
\mapleinline{inert}{2d}{`  179: `, _3F_2([n-2*a, a, 1],[1/2+n-a, 2*a+1-n],1);}{%
\[
\mathit{\ \ 179:\ } \,\mathrm{_3F_2}([n - 2\,a, \,a, \,1], \,[
{\displaystyle \frac {1}{2}}  + n - a, \,2\,a + 1 - n], \,1)
\]
}
\end{maplelatex}

\begin{maplelatex}
\mapleinline{inert}{2d}{`  180: `, _3F_2([-a+n-1/2, 1/2+2*a-n, 2*a-n],[1/2, 3*a+1/2-n],1);}{%
\[
\mathit{\ \ \dagger 180:\ } \,\mathrm{_3F_2}([ - a + n - {\displaystyle 
\frac {1}{2}} , \,{\displaystyle \frac {1}{2}}  + 2\,a - n, \,2\,
a - n], \,[{\displaystyle \frac {1}{2}} , \,3\,a + 
{\displaystyle \frac {1}{2}}  - n], \,1)
\]
}
\end{maplelatex}

\begin{maplelatex}
\mapleinline{inert}{2d}{`  181: `, _3F_2([1/2-2*a+n, 1/2+2*a-n, a-n+1],[1/2, 2*a+3/2-n],1);}{%
\[
\mathit{\ \ 181:\ } \,\mathrm{_3F_2}([{\displaystyle \frac {1}{2}
}  - 2\,a + n, \,{\displaystyle \frac {1}{2}}  + 2\,a - n, \,a - 
n + 1], \,[{\displaystyle \frac {1}{2}} , \,2\,a + 
{\displaystyle \frac {3}{2}}  - n], \,1)
\]
}
\end{maplelatex}

\begin{maplelatex}
\mapleinline{inert}{2d}{`  182: `, _3F_2([n-2*a, 2*a-n, a-n+1],[1/2, 2*a+1-n],1);}{%
\[
\mathit{\ \ 182:\ } \,\mathrm{_3F_2}([n - 2\,a, \,2\,a - n, \,a - 
n + 1], \,[{\displaystyle \frac {1}{2}} , \,2\,a + 1 - n], \,1)
\]
}
\end{maplelatex}

\begin{maplelatex}
\mapleinline{inert}{2d}{`  183: `, _3F_2([a+1/2, 1/2+2*a-n, 1+4*a-2*n],[3*a+1/2-n,
2*a+3/2-n],1);}{%
\[
\mathit{\ \ \dagger 183:\ } \,\mathrm{_3F_2}([a + {\displaystyle \frac {1
}{2}} , \,{\displaystyle \frac {1}{2}}  + 2\,a - n, \,1 + 4\,a - 
2\,n], \,[3\,a + {\displaystyle \frac {1}{2}}  - n, \,2\,a + 
{\displaystyle \frac {3}{2}}  - n], \,1)
\]
}
\end{maplelatex}

\begin{maplelatex}
\mapleinline{inert}{2d}{`  184: `, _3F_2([a, 2*a-n, 1+4*a-2*n],[3*a+1/2-n, 2*a+1-n],1);}{%
\[
\mathit{\ \ 184:\ } \,\mathrm{_3F_2}([a, \,2\,a - n, \,1 + 4\,a - 
2\,n], \,[3\,a + {\displaystyle \frac {1}{2}}  - n, \,2\,a + 1 - 
n], \,1)
\]
}
\end{maplelatex}

\begin{maplelatex}
\mapleinline{inert}{2d}{`  185: `, _3F_2([n-a, a, a+1/2],[n+1-a, 3*a+1-n],1);}{%
\[
\mathit{\ \ 185:\ } \,\mathrm{_3F_2}([n - a, \,a, \,a + 
{\displaystyle \frac {1}{2}} ], \,[n + 1 - a, \,3\,a + 1 - n], \,
1)
\]
}
\end{maplelatex}

\begin{maplelatex}
\mapleinline{inert}{2d}{`  186: `, _3F_2([n-a, 1/2-2*a+n, 1-2*a+n],[n+1-a, 3/2],1);}{%
\[
\mathit{\ \ \dagger 186:\ } \,\mathrm{_3F_2}([n - a, \,{\displaystyle 
\frac {1}{2}}  - 2\,a + n, \,1 - 2\,a + n], \,[n + 1 - a, \,
{\displaystyle \frac {3}{2}} ], \,1)
\]
}
\end{maplelatex}

\begin{maplelatex}
\mapleinline{inert}{2d}{`  187: `, _3F_2([a, 1/2-2*a+n, 1],[n+1-a, 2*a+3/2-n],1);}{%
\[
\mathit{\ \ \dagger 187:\ } \,\mathrm{_3F_2}([a, \,{\displaystyle \frac {1
}{2}}  - 2\,a + n, \,1], \,[n + 1 - a, \,2\,a + {\displaystyle 
\frac {3}{2}}  - n], \,1)
\]
}
\end{maplelatex}

\begin{maplelatex}
\mapleinline{inert}{2d}{`  188: `, _3F_2([a+1/2, 1-2*a+n, 1],[n+1-a, 2-n+2*a],1);}{%
\[
\mathit{\ \ 188:\ } \,\mathrm{_3F_2}([a + {\displaystyle \frac {1
}{2}} , \,1 - 2\,a + n, \,1], \,[n + 1 - a, \,2 - n + 2\,a], \,1)
\]
}
\end{maplelatex}

\begin{maplelatex}
\mapleinline{inert}{2d}{`  189: `, _3F_2([n-a, 1/2+2*a-n, 2*a+1-n],[3*a+1-n, 3/2],1);}{%
\[
\mathit{\ \ \dagger 189:\ } \,\mathrm{_3F_2}([n - a, \,{\displaystyle 
\frac {1}{2}}  + 2\,a - n, \,2\,a + 1 - n], \,[3\,a + 1 - n, \,
{\displaystyle \frac {3}{2}} ], \,1)
\]
}
\end{maplelatex}

\begin{maplelatex}
\mapleinline{inert}{2d}{`  190: `, _3F_2([a, 1/2+2*a-n, 1+4*a-2*n],[3*a+1-n, 2*a+3/2-n],1);}{%
\[
\mathit{\ \ \dagger 190:\ } \,\mathrm{_3F_2}([a, \,{\displaystyle \frac {1
}{2}}  + 2\,a - n, \,1 + 4\,a - 2\,n], \,[3\,a + 1 - n, \,2\,a + 
{\displaystyle \frac {3}{2}}  - n], \,1)
\]
}
\end{maplelatex}

\begin{maplelatex}
\mapleinline{inert}{2d}{`  191: `, _3F_2([a+1/2, 2*a+1-n, 1+4*a-2*n],[3*a+1-n, 2-n+2*a],1);}{%
\[
\mathit{\ \ 191:\ } \,\mathrm{_3F_2}([a + {\displaystyle \frac {1
}{2}} , \,2\,a + 1 - n, \,1 + 4\,a - 2\,n], \,[3\,a + 1 - n, \,2
 - n + 2\,a], \,1)
\]
}
\end{maplelatex}

\begin{maplelatex}
\mapleinline{inert}{2d}{`  192: `, _3F_2([1/2-2*a+n, 1/2+2*a-n, 3/2-n+a],[3/2,
2*a+3/2-n],1);}{%
\[
\mathit{\ \ 192:\ } \,\mathrm{_3F_2}([{\displaystyle \frac {1}{2}
}  - 2\,a + n, \,{\displaystyle \frac {1}{2}}  + 2\,a - n, \,
{\displaystyle \frac {3}{2}}  - n + a], \,[{\displaystyle \frac {
3}{2}} , \,2\,a + {\displaystyle \frac {3}{2}}  - n], \,1)
\]
}
\end{maplelatex}

\begin{maplelatex}
\mapleinline{inert}{2d}{`  193: `, _3F_2([1-2*a+n, 2*a+1-n, 3/2-n+a],[3/2, 2-n+2*a],1);}{%
\[
\mathit{\ \ 193:\ } \,\mathrm{_3F_2}([1 - 2\,a + n, \,2\,a + 1 - n
, \,{\displaystyle \frac {3}{2}}  - n + a], \,[{\displaystyle 
\frac {3}{2}} , \,2 - n + 2\,a], \,1)
\]
}
\end{maplelatex}

\begin{maplelatex}
\mapleinline{inert}{2d}{`  202: `, _3F_2([1, 1/2-2*a+n, 3/2-n+a],[-a+3/2, 2*a+3/2-n],1);}{%
\[
\mathit{\ \ \dagger 202:\ } \,\mathrm{_3F_2}([1, \,{\displaystyle \frac {1
}{2}}  - 2\,a + n, \,{\displaystyle \frac {3}{2}}  - n + a], \,[
 - a + {\displaystyle \frac {3}{2}} , \,2\,a + {\displaystyle 
\frac {3}{2}}  - n], \,1)
\]
}
\end{maplelatex}

\begin{maplelatex}
\mapleinline{inert}{2d}{`  203: `, _3F_2([1, 1-2*a+n, 3/2-n+a],[-a+2, 2-n+2*a],1);}{%
\[
\mathit{\ \ 203:\ } \,\mathrm{_3F_2}([1, \,1 - 2\,a + n, \,
{\displaystyle \frac {3}{2}}  - n + a], \,[ - a + 2, \,2 - n + 2
\,a], \,1)
\]
}
\end{maplelatex}

\begin{maplelatex}
\mapleinline{inert}{2d}{`  204: `, _3F_2([a-n+1, 1/2-a, 3/2-n+a],[-a+3/2, 3*a-2*n+3/2],1);}{%
\[
\mathit{\ \ 204:\ } \,\mathrm{_3F_2}([a - n + 1, \,{\displaystyle 
\frac {1}{2}}  - a, \,{\displaystyle \frac {3}{2}}  - n + a], \,[
 - a + {\displaystyle \frac {3}{2}} , \,3\,a - 2\,n + 
{\displaystyle \frac {3}{2}} ], \,1)
\]
}
\end{maplelatex}

\begin{maplelatex}
\mapleinline{inert}{2d}{`  205: `, _3F_2([-a+1, a-n+1, 3/2-n+a],[-a+2, 3*a-2*n+2],1);}{%
\[
\mathit{\ \ 205:\ } \,\mathrm{_3F_2}([ - a + 1, \,a - n + 1, \,
{\displaystyle \frac {3}{2}}  - n + a], \,[ - a + 2, \,3\,a - 2\,
n + 2], \,1)
\]
}
\end{maplelatex}

\begin{maplelatex}
\mapleinline{inert}{2d}{`  206: `, _3F_2([1+4*a-2*n, 1/2+2*a-n, 3/2-n+a],[2*a+3/2-n,
3*a-2*n+3/2],1);}{%
\[
\mathit{\ \ \dagger 206:\ } \,\mathrm{_3F_2}([1 + 4\,a - 2\,n, \,
{\displaystyle \frac {1}{2}}  + 2\,a - n, \,{\displaystyle 
\frac {3}{2}}  - n + a], \,[2\,a + {\displaystyle \frac {3}{2}} 
 - n, \,3\,a - 2\,n + {\displaystyle \frac {3}{2}} ], \,1)
\]
}
\end{maplelatex}

\begin{maplelatex}
\mapleinline{inert}{2d}{`  207: `, _3F_2([1, a-n+1, 1/2-2*a+n],[-a+2, 2*a+3/2-n],1);}{%
\[
\mathit{\ \ \dagger 207:\ } \,\mathrm{_3F_2}([1, \,a - n + 1, \,
{\displaystyle \frac {1}{2}}  - 2\,a + n], \,[ - a + 2, \,2\,a + 
{\displaystyle \frac {3}{2}}  - n], \,1)
\]
}
\end{maplelatex}

\begin{maplelatex}
\mapleinline{inert}{2d}{`  208: `, _3F_2([1, n-2*a, a-n+1],[2*a+1-n, -a+3/2],1);}{%
\[
\mathit{\ \ 208:\ } \,\mathrm{_3F_2}([1, \,n - 2\,a, \,a - n + 1]
, \,[2\,a + 1 - n, \, - a + {\displaystyle \frac {3}{2}} ], \,1)
\]
}
\end{maplelatex}

\begin{maplelatex}
\mapleinline{inert}{2d}{`  209: `, _3F_2([a-n+1, 1+4*a-2*n, 1/2+2*a-n],[3*a-2*n+2,
2*a+3/2-n],1);}{%
\[
\mathit{\ \ \dagger 209:\ } \,\mathrm{_3F_2}([a - n + 1, \,1 + 4\,a - 2\,n
, \,{\displaystyle \frac {1}{2}}  + 2\,a - n], \,[3\,a - 2\,n + 2
, \,2\,a + {\displaystyle \frac {3}{2}}  - n], \,1)
\]
}
\end{maplelatex}

\begin{maplelatex}
\mapleinline{inert}{2d}{`  249: `, _3F_2([b, 2*b-a-m, b-a+c-n],[1+b-c, b-a+1],1);}{%
\[
\mathit{\ \ 249:\ } \,\mathrm{_3F_2}([b, \,2\,b - a - m, \,b - a
 + c - n], \,[1 + b - c, \,b - a + 1], \,1)
\]
}
\end{maplelatex}

\begin{maplelatex}
\mapleinline{inert}{2d}{`  250: `, _3F_2([b, 2*b-a+m, b-a+c+n],[1+b-c, b-a+1],1);}{%
\[
\mathit{\ \ 250:\ } \,\mathrm{_3F_2}([b, \,2\,b - a + m, \,b - a
 + c + n], \,[1 + b - c, \,b - a + 1], \,1)
\]
}
\end{maplelatex}

\begin{maplelatex}
\mapleinline{inert}{2d}{`  252: `, _3F_2([-3*a, 1-2*n, -a+n-1/2],[-a-2*n+2, 1/2-3*a+n],1);}{%
\[
\mathit{\ \ 252:\ } \,\mathrm{_3F_2}([ - 3\,a, \,1 - 2\,n, \, - a
 + n - {\displaystyle \frac {1}{2}} ], \,[ - a - 2\,n + 2, \,
{\displaystyle \frac {1}{2}}  - 3\,a + n], \,1)
\]
}
\end{maplelatex}

\begin{maplelatex}
\mapleinline{inert}{2d}{`  261: `, _3F_2([a, b, 1-b+m],[c, 1+m-c+2*a-n],1);}{%
\[
\mathit{\ \ 261:\ } \,\mathrm{_3F_2}([a, \,b, \,1 - b + m], \,[c, 
\,1 + m - c + 2\,a - n], \,1)
\]
}
\end{maplelatex}

\begin{maplelatex}
\mapleinline{inert}{2d}{`  262: `, _3F_2([a, b, c],[2*c+n, 1/2+1/2*a+1/2*b+1/2*m-1/2*n],1);}{%
\[
\mathit{\ \ 262:\ } \,\mathrm{_3F_2}([a, \,b, \,c], \,[2\,c + n, 
\,{\displaystyle \frac {1}{2}}  + {\displaystyle \frac {a}{2}} 
 + {\displaystyle \frac {b}{2}}  + {\displaystyle \frac {m}{2}} 
 - {\displaystyle \frac {n}{2}} ], \,1)
\]
}
\end{maplelatex}

\begin{maplelatex}
\mapleinline{inert}{2d}{`  263: `, _3F_2([a, b, -a+m+1-n],[c, 1-c+2*b+m],1);}{%
\[
\mathit{\ \ 263:\ } \,\mathrm{_3F_2}([a, \,b, \, - a + m + 1 - n]
, \,[c, \,1 - c + 2\,b + m], \,1)
\]
}
\end{maplelatex}

\begin{maplelatex}
\mapleinline{inert}{2d}{`  264: `, _3F_2([a, b, -b-m+1],[c, -m+1-c+2*a+n],1);}{%
\[
\mathit{\ \ 264:\ } \,\mathrm{_3F_2}([a, \,b, \, - b - m + 1], \,[
c, \, - m + 1 - c + 2\,a + n], \,1)
\]
}
\end{maplelatex}

\begin{maplelatex}
\mapleinline{inert}{2d}{`  265: `, _3F_2([a, b, c],[2*c-n, 1/2+1/2*a+1/2*b-1/2*m+1/2*n],1);}{%
\[
\mathit{\ \ 265:\ } \,\mathrm{_3F_2}([a, \,b, \,c], \,[2\,c - n, 
\,{\displaystyle \frac {1}{2}}  + {\displaystyle \frac {a}{2}} 
 + {\displaystyle \frac {b}{2}}  - {\displaystyle \frac {m}{2}} 
 + {\displaystyle \frac {n}{2}} ], \,1)
\]
}
\end{maplelatex}

\begin{maplelatex}
\mapleinline{inert}{2d}{`  266: `, _3F_2([a, b, 1-a+n-m],[c, 1-c+2*b-m],1);}{%
\[
\mathit{\ \ 266:\ } \,\mathrm{_3F_2}([a, \,b, \,1 - a + n - m], \,
[c, \,1 - c + 2\,b - m], \,1)
\]
}
\end{maplelatex}

\begin{maplelatex}
\mapleinline{inert}{2d}{`  296: `, _3F_2([a, 1+a-c, a+c-2*b-m],[1+a-b, 2*a-m+n],1);}{%
\[
\mathit{\ \ 296:\ } \,\mathrm{_3F_2}([a, \,1 + a - c, \,a + c - 2
\,b - m], \,[1 + a - b, \,2\,a - m + n], \,1)
\]
}
\end{maplelatex}

\begin{maplelatex}
\mapleinline{inert}{2d}{`  297: `, _3F_2([a, 1+a-c, a+c-2*b+m],[1+a-b, 2*a-n+m],1);}{%
\[
\mathit{\ \ 297:\ } \,\mathrm{_3F_2}([a, \,1 + a - c, \,a + c - 2
\,b + m], \,[1 + a - b, \,2\,a - n + m], \,1)
\]
}
\end{maplelatex}

\end{maplegroup}

%% file: AppendixE.tex
\begin{maplegroup}
\mapleresult
20: Subset of case 259.  Use the mapping  \{a = a, b = b-n, c = c, e = c+m, n = n\} \\ 21: Subset of case 259. This case requires n$<$m Use the mapping  \{a = a, b = -n+m, c = b, e = c, n = n\} \\ 23: Subset of case 260. This case requires n$<$m Use the mapping  \{a = a, b = b, c = c, e = -n+m, n = n\} \\ 24: Subset of case 260.  Use the mapping  \{a = a, b = b, c = c, e = a-c+b-n+m, n = n\} \\ 40: Subset of case 249.  Use the mapping  \{a = 2*a-n, b = 1/2+n-a, c = 1/2+2*a-n, m = n, n = 1\} \\ 41: Subset of case 262.  Use the mapping  \{a = 1-4*a+2*n, b = 1-2*a+n, c = 1/2-2*a+n, m = n, n = 1\} \\ 42: Subset of case 249.  Use the mapping  \{a = 0, b = 1/2+n-a, c = 1/2+2*a-n, m = n, n = 1\} \\ 43: Subset of case 263.  Use the mapping  \{a = n-a, b = 1/2-2*a+n, c = 1-3*a+2*n, m = n, n = 1\} \\ 44: Subset of case 261.  Use the mapping  \{a = 3/2-2*a+n, b = 1/2+n-a, c = -3*a+2*n+3/2, m = n, n = 1\} \\ 45: Subset of case 249.  Use the mapping  \{a = 4*a-2*n, b = a+1/2, c = 1/2+2*a-n, m = n, n = 1\} \\ 46: Subset of case 262.  Use the mapping  \{a = 1-2*a+n, b = 1, c = 1/2-2*a+n, m = n, n = 1\} \\ 47: Subset of case 249.  Use the mapping  \{a = 2*a-n, b = a+1/2, c = 1/2+2*a-n, m = n, n = 1\} \\ 48: Subset of case 259.  Use the mapping  \{a = a, b = c-1, c = b, e = e, n = 1\} \\ 50: Subset of case 260.  Use the mapping  \{a = a, b = b, c = c, e = e, n = 1\} \\ 51: Subset of case 246.  Use the mapping  \{a = -b+a-1, b = -1+a-c, c = -b+1, n = 2\} \\ 52: Subset of case 244.  Use the mapping  \{a = 3-c, b = 3-b, c = 3-a, n = 2\} \\ 53: Subset of case 247.  Use the mapping  \{a = a, b = b, c = c, n = 2\} \\ 54: Subset of case 259. Prudnikov 7.4.4.(91) to (93) have a sign error Use the mapping  \{a = a, b = b, c = b-m, e = a-n+m, n = n\} \\ 55: Subset of case 259.  Use the mapping  \{a = -L-m+n, b = a, c = b, e = b-m, n = n\} \\ 56: Subset of case 259. Requires n $\geq$ L+m+1. Evaluate numerically with care. Use the mapping  \{a = a, b = b-m, c = b+L-n, e = b-n, n = m\} \\ 57: Subset of case 259. Requires n$\geq$L+m+1. Use the mapping  \{a = a, b = b, c = b-m, e = a-L, n = n\} \\ 61: Subset of case 259. Prudnikov 7.4.4.(91) to (93) have a sign error. Use the mapping  \{a = -n, b = a, c = b, e = b-m, n = 2\} \\ 62: Subset of case 469.  Use the mapping  \{a = -a*b/(b-n), b = b, c = a\} \\ 66: Subset of case 250.  Use the mapping  \{a = -1+n+2*a, b = -a-n+1, c = 2*a-3/2+n, m = n, n = 1\} \\ 67: Subset of case 264.  Use the mapping  \{a = 3/2-2*a-n, b = -a-n+1, c = 3-3*a-2*n, m = n, n = 1\} \\ 68: Subset of case 250.  Use the mapping  \{a = -2+4*a+2*n, b = a, c = 2*a-3/2+n, m = n, n = 1\} \\ 69: Subset of case 265.  Use the mapping  \{a = 3-4*a-2*n, b = -2*a-n+2, c = 5/2-2*a-n, m = n, n = 1\} \\ 70: Subset of case 250.  Use the mapping  \{a = 0, b = -a-n+1, c = 2*a-3/2+n, m = n, n = 1\} \\ 71: Subset of case 266.  Use the mapping  \{a = -a+3/2-n, b = 5/2-2*a-n, c = 7/2-3*a-2*n, m = n, n = 1\} \\ 72: Subset of case 265.  Use the mapping  \{a = -2*a-n+2, b = 1, c = 5/2-2*a-n, m = n, n = 1\} \\ 73: Subset of case 250.  Use the mapping  \{a = -1+n+2*a, b = a, c = 2*a-3/2+n, m = n, n = 1\} \\ 74: Subset of case 241.  Use the mapping  \{a = 1/2-2*a+n, b = 3/2-n+a, n = n+1\} \\ 75: Subset of case 238.  Use the mapping  \{a = n-2*a, b = a-n+1, n = n+1\} \\ 76: Subset of case 241.  Use the mapping  \{a = 1/2-2*a+n, b = a-n+1, n = n+1\} \\ 77: Subset of case 238.  Use the mapping  \{a = 1-2*a+n, b = 3/2-n+a, n = n+1\} \\ 78: Subset of case 242.  Use the mapping  \{a = a-n+1, b = 1+4*a-2*n, n = n+1\} \\ 79: Subset of case 242.  Use the mapping  \{a = 3/2-n+a, b = 1+4*a-2*n, n = n+1\} \\ 90: Subset of case 469.  Use the mapping  \{a = 3/2+3*a-3*n, b = -a+n-1/2, c = 3*n-3/2\} \\ 93: Subset of case 469.  Use the mapping  \{a = -1/2*b-1/2*a+1/2*c+1/2*(a\symbol{94}2-2*a*c-2*a*b+c\symbol{94}2-2*b*c+b\symbol{94}2)\symbol{94}(1/2), b = b, c = a\} \\ 94: Subset of case 350.  Use the mapping  \{a = 7/6-a, n = 1\} \\ 95: Subset of case 357.  Use the mapping  \{a = -a+3/2, n = 1\} \\ 96: Subset of case 363.  Use the mapping  \{a = 7/6-a, n = 1\} \\ 97: Subset of case 361.  Use the mapping  \{a = 2/3-a, n = 1\} \\ 98: Subset of case 364.  Use the mapping  \{a = 11/6-a, n = 1\} \\ 99: Subset of case 347.  Use the mapping  \{a = 7/6-a, n = 1\} \\ 100: Subset of case 352.  Use the mapping  \{a = 2/3-a, n = 1\} \\ 101: Subset of case 355.  Use the mapping  \{a = -a+3/2, n = 1\} \\ 102: Subset of case 356.  Use the mapping  \{a = -a+1, n = 1\} \\ 103: Subset of case 362.  Use the mapping  \{a = 1/3-1/2*a, n = 1\} \\ 135: Subset of case 250.  Use the mapping  \{a = 0, b = 2*a-1/2-n, c = 1/2-a, m = 1, n = n\} \\ 136: Subset of case 296.  Use the mapping  \{a = a, b = 1/2-2*a+n, c = n+1-a, m = n, n = 1\} \\ 137: Subset of case 265.  Use the mapping  \{a = 2*a-n, b = 4*a-2*n, c = a+1/2, m = 1, n = n\} \\ 138: Subset of case 296.  Use the mapping  \{a = a-n+1, b = 1/2-a, c = a-n+1, m = 1, n = n\} \\ 140: Subset of case 262.  Use the mapping  \{a = -2+4*a+2*n, b = -1+n+2*a, c = a, m = 1, n = n\} \\ 141: Subset of case 249.  Use the mapping  \{a = 0, b = 2*a-1/2+n, c = -a+1, m = 1, n = n\} \\ 142: Subset of case 297.  Use the mapping  \{a = a+1/2, b = 5/2-2*a-n, c = 7/2-3*a-2*n, m = n, n = 1\} \\ 143: Subset of case 297.  Use the mapping  \{a = n-1/2+a, b = -a+1, c = n-1/2+a, m = 1, n = n\} \\ 144: Subset of case 213.  Use the mapping  \{a = a-n+1, b = 1+4*a-2*n, n = n+1\} \\ 145: Subset of case 213.  Use the mapping  \{a = 3/2-n+a, b = 1+4*a-2*n, n = n+1\} \\ 146: Subset of case 250.  Use the mapping  \{a = 2*a-n, b = 2*a-1/2-n, c = 1/2-a, m = 1, n = n\} \\ 147: Subset of case 249.  Use the mapping  \{a = -1+n+2*a, b = 2*a-1/2+n, c = -a+1, m = 1, n = n\} \\ 153: Subset of case 265.  Use the mapping  \{a = 1-2*a+n, b = 1, c = a+1/2, m = 1, n = n\} \\ 154: Subset of case 250.  Use the mapping  \{a = 1-2*a+n, b = 1/2-2*a+n, c = 1/2-a, m = 1, n = n\} \\ 155: Subset of case 250.  Use the mapping  \{a = 1-4*a+2*n, b = 1/2-2*a+n, c = 1/2-a, m = 1, n = n\} \\ 156: Subset of case 266.  Use the mapping  \{a = n-a, b = a+1/2, c = n+1-a, m = 1, n = n\} \\ 157: Subset of case 264.  Use the mapping  \{a = a-n+1/2, b = 1/2-2*a+n, c = 1/2, m = 1, n = n\} \\ 158: Subset of case 296.  Use the mapping  \{a = a, b = 2*a-1/2-n, c = 3*a-n, m = n, n = 1\} \\ 159: Subset of case 296.  Use the mapping  \{a = a-n+1, b = a-n+1/2, c = 3*a-2*n+1, m = 1, n = n\} \\ 169: Subset of case 249.  Use the mapping  \{a = 3-4*a-2*n, b = 5/2-2*a-n, c = -a+1, m = 1, n = n\} \\ 170: Subset of case 249.  Use the mapping  \{a = -2*a-n+2, b = 5/2-2*a-n, c = -a+1, m = 1, n = n\} \\ 171: Subset of case 263.  Use the mapping  \{a = a+1/2, b = a, c = 5/2-a-n, m = 1, n = n\} \\ 172: Subset of case 262.  Use the mapping  \{a = 1, b = -2*a-n+2, c = a, m = 1, n = n\} \\ 173: Subset of case 261.  Use the mapping  \{a = a+n, b = 5/2-2*a-n, c = 2*a+1/2+n, m = 1, n = n\} \\ 174: Subset of case 297.  Use the mapping  \{a = a+1/2, b = 2*a-1/2+n, c = a+1/2, m = n, n = 1\} \\ 175: Subset of case 297.  Use the mapping  \{a = n-1/2+a, b = a+n, c = -3/2+3*a+2*n, m = 1, n = n\} \\ 176: Subset of case 214.  Use the mapping  \{a = 1/2-2*a+n, b = -a+n-1/2, n = n+1\} \\ 177: Subset of case 215.  Use the mapping  \{a = a+1/2, b = a, n = n+1\} \\ 178: Subset of case 210.  Use the mapping  \{a = 1/2-2*a+n, b = a+1/2, n = n+1\} \\ 179: Subset of case 216.  Use the mapping  \{a = n-2*a, b = a, n = n+1\} \\ 180: Subset of case 217.  Use the mapping  \{a = 1/2+2*a-n, b = -a+n-1/2, n = n+1\} \\ 181: Subset of case 211.  Use the mapping  \{a = 1/2-2*a+n, b = a-n+1, n = n+1\} \\ 182: Subset of case 218.  Use the mapping  \{a = a-n+1, b = n-2*a, n = n+1\} \\ 183: Subset of case 212.  Use the mapping  \{a = a+1/2, b = 1/2+2*a-n, n = n+1\} \\ 184: Subset of case 194.  Use the mapping  \{a = 1/2-2*a+n, b = n-a, n = n+1\} \\ 185: Subset of case 215.  Use the mapping  \{a = a, b = a+1/2, n = n+1\} \\ 186: Subset of case 214.  Use the mapping  \{a = 1/2-2*a+n, b = n-a, n = n+1\} \\ 187: Subset of case 210.  Use the mapping  \{a = 1/2-2*a+n, b = a, n = n+1\} \\ 188: Subset of case 216.  Use the mapping  \{a = 1-2*a+n, b = a+1/2, n = n+1\} \\ 189: Subset of case 217.  Use the mapping  \{a = 1/2+2*a-n, b = n-a, n = n+1\} \\ 190: Subset of case 212.  Use the mapping  \{a = a, b = 1/2+2*a-n, n = n+1\} \\ 191: Subset of case 194.  Use the mapping  \{a = 1/2-2*a+n, b = -a+n-1/2, n = n+1\} \\ 192: Subset of case 211.  Use the mapping  \{a = 1/2-2*a+n, b = 3/2-n+a, n = n+1\} \\ 193: Subset of case 218.  Use the mapping  \{a = 3/2-n+a, b = 1-2*a+n, n = n+1\} \\ 202: Subset of case 227.  Use the mapping  \{a = 3/2-n+a, b = 1/2-2*a+n, n = n+1\} \\ 203: Subset of case 226.  Use the mapping  \{a = 3/2-n+a, b = 1-2*a+n, n = n+1\} \\ 204: Subset of case 224.  Use the mapping  \{a = 1/2-a, b = a-n+1, n = n+1\} \\ 205: Subset of case 224.  Use the mapping  \{a = -a+1, b = 3/2-n+a, n = n+1\} \\ 206: Subset of case 229.  Use the mapping  \{a = 3/2-n+a, b = 1/2+2*a-n, n = n+1\} \\ 207: Subset of case 227.  Use the mapping  \{a = a-n+1, b = 1/2-2*a+n, n = n+1\} \\ 208: Subset of case 226.  Use the mapping  \{a = a-n+1, b = n-2*a, n = n+1\} \\ 209: Subset of case 229.  Use the mapping  \{a = a-n+1, b = 1/2+2*a-n, n = n+1\} \\ 249: Enlargens Dixon/Whipple/Watson case 18. Use the mapping  \{a = 2*b-a-m, b = b, c = b-a+c+n, m = m, n = n+m\} \\ 250: Enlargens Dixon/Whipple/Watson case 18. Use the mapping  \{a = 2*c-a+n, b = c, c = c-a+b-m, m = n, n = n+m\} \\ 252: Subset of case 469.  Use the mapping  \{a = 2*a, b = 1-2*n, c = -3*a\} \\ 261: Enlargens Dixon/Whipple/Watson case 14. Use the mapping  \{a = b, b = a, c = c, m = m, n = n+m\} \\ 262: Enlargens Dixon/Whipple/Watson case 1. Use the mapping  \{a = a, b = b, c = c, m = n+m, n = n\} \\ 263: Enlargens Dixon/Whipple/Watson case 15. Use the mapping  \{a = a, b = b, c = c, m = n, n = n+m\} \\ 264: Enlargens Dixon/Whipple/Watson case 15. Use the mapping  \{a = b, b = a, c = c, m = m, n = n+m\} \\ 265: Enlargens Dixon/Whipple/Watson case 10. Use the mapping  \{a = a, b = b, c = c, m = n+m, n = n\} \\ 266: Enlargens Dixon/Whipple/Watson case 14. Use the mapping  \{a = a, b = b, c = c, m = n, n = n+m\} \\ 296: Enlargens Dixon/Whipple/Watson case 1. Use the mapping  \{a = c, b = 1/2+c-1/2*a-1/2*b-1/2*m, c = c+1-a, m = m, n = n+m\} \\ 297: Enlargens Dixon/Whipple/Watson case 10. Use the mapping  \{a = c, b = 1/2+c-1/2*a-1/2*b+1/2*m, c = c+1-a, m = m, n = n+m\} \\ 
\end{maplegroup}

%% file: AppendixF.tex
\begin{maplegroup}
\mapleresult
\begin{maplelatex}
\mapleinline{inert}{2d}{` 20: `,
GAMMA(a+n-m+1)*GAMMA(-c+a-m+1)*Sum((-1)^k*GAMMA(b-a-n+k)*GAMMA(1+b-c-m
+k)/GAMMA(m-k)/GAMMA(b-n-c-m+1+k)/GAMMA(b-m+1+k)/GAMMA(k+1),k = 0 ..
m-1)/GAMMA(n+1-c+a-m)/GAMMA(b-a-n)*GAMMA(-a-n+m)*GAMMA(c+m)*GAMMA(b-n)
/GAMMA(a)/GAMMA(c-a+m-n);}{%
\maplemultiline{
\mathit{\ 20:\ }   \Gamma (a + n - m + 1)\,\Gamma ( - c + a - m + 
1) \\
 \left(  \! {\displaystyle \sum _{k=0}^{m - 1}} \,{\displaystyle 
\frac {(-1)^{k}\,\Gamma (b - a - n + k)\,\Gamma (1 + b - c - m + 
k)}{\Gamma (m - k)\,\Gamma (b - n - c - m + 1 + k)\,\Gamma (b - m
 + 1 + k)\,\Gamma (k + 1)}}  \!  \right)  \\
\Gamma ( - a - n + m)\,\Gamma (c + m)\,\Gamma (b - n)/(\Gamma (n
 + 1 - c + a - m)\,\Gamma (b - a - n)\,\Gamma (a) 
\Gamma (c - a + m - n)) }
}
\end{maplelatex}

\begin{maplelatex}
\mapleinline{inert}{2d}{` 21: `,
GAMMA(a+n-m+1)*GAMMA(1+a-c)*Sum((-1)^k*GAMMA(-b+c-a-n+k)*GAMMA(1-b+k)/
GAMMA(m-k)/GAMMA(-b-n+1+k)/GAMMA(-b+c-m+1+k)/GAMMA(k+1),k = 0 ..
m-1)/GAMMA(1+a-c+n)*GAMMA(c)*GAMMA(-n+m)/GAMMA(a)/GAMMA(c-a-n);}{%
\maplemultiline{
\mathit{\ 21:\ }   \Gamma (a + n - m + 1)\,\Gamma (1 + a - c) \\
 \left(  \! {\displaystyle \sum _{k=0}^{m - 1}} \,{\displaystyle 
\frac {(-1)^{k}\,\Gamma ( - b + c - a - n + k)\,\Gamma (1 - b + k
)}{\Gamma (m - k)\,\Gamma ( - b - n + 1 + k)\,\Gamma ( - b + c - 
m + 1 + k)\,\Gamma (k + 1)}}  \!  \right) \,\Gamma (c)\,\Gamma (
 - n + m) \\
/(\Gamma (1 + a - c + n)\,\Gamma (a)\,\Gamma (c - a - n)) } 
}
\end{maplelatex}
\begin{flushright}
n$<$m
\end{flushright} 

\begin{maplelatex}
\mapleinline{inert}{2d}{` 23: `,
GAMMA(a+n-m+1)*GAMMA(b-c+a-m+1)*Sum((-1)^k*GAMMA(c-a+k)*GAMMA(1+b-m+n+
k)/GAMMA(m-k)/GAMMA(b-m+1+k)/GAMMA(c+n-m+1+k)/GAMMA(k+1),k = 0 ..
m-1)/GAMMA(n+1+b-c+a-m)/GAMMA(c-a)*GAMMA(-b+c-a+m)*GAMMA(c)*GAMMA(-n+m
)/GAMMA(a)/GAMMA(-b+c-n+m-a);}{%
\maplemultiline{
\mathit{\ 23:\ }   \Gamma (a + n - m + 1)\,\Gamma (b - c + a - m
 + 1) \\
 \left(  \! {\displaystyle \sum _{k=0}^{m - 1}} \,{\displaystyle 
\frac {(-1)^{k}\,\Gamma (c - a + k)\,\Gamma (1 + b - m + n + k)}{
\Gamma (m - k)\,\Gamma (b - m + 1 + k)\,\Gamma (c + n - m + 1 + k
)\,\Gamma (k + 1)}}  \!  \right)  \\
\Gamma ( - b + c - a + m)\,\Gamma (c)\,\Gamma ( - n + m)/(\Gamma 
(n + 1 + b - c + a - m)\,\Gamma (c - a)\,\Gamma (a) 
\Gamma ( - b + c - n + m - a))  } 
}
\end{maplelatex}
\begin{flushright}
n$<$m
\end{flushright} 

\begin{maplelatex}
\mapleinline{inert}{2d}{` 24: `,
GAMMA(c-b+n-m+1)*GAMMA(-a+c-b-m+1)*Sum((-1)^k*GAMMA(b+k)*GAMMA(1-a+c-m
+n+k)/GAMMA(m-k)/GAMMA(-a+c-m+1+k)/GAMMA(c+n-m+1+k)/GAMMA(k+1),k = 0
.. m-1)/GAMMA(n+1-a+c-b-m)/GAMMA(b)*GAMMA(m)*GAMMA(c)/GAMMA(c-b);}{%
\maplemultiline{
\mathit{\ 24:\ }   \Gamma (c - b + n - m + 1)\,\Gamma ( - a + c - 
b - m + 1) \\
 \left(  \! {\displaystyle \sum _{k=0}^{m - 1}} \,{\displaystyle 
\frac {(-1)^{k}\,\Gamma (b + k)\,\Gamma (1 - a + c - m + n + k)}{
\Gamma (m - k)\,\Gamma ( - a + c - m + 1 + k)\,\Gamma (c + n - m
 + 1 + k)\,\Gamma (k + 1)}}  \!  \right) \,\Gamma (m)\,\Gamma (c)/ \\
(\Gamma (n + 1 - a + c - b - m)\,\Gamma (b)\,\Gamma (c - b)) }
}
\end{maplelatex}

\begin{maplelatex}
\mapleinline{inert}{2d}{` 40: `,
GAMMA(2-2*a+2*n)*2^(2*n-1-4*a)/GAMMA(1-4*a+2*n)/GAMMA(2-4*a+2*n)*GAMMA
(1-3*a+2*n)*GAMMA(-3*a+2*n+3/2)/GAMMA(3/2-a+n)/GAMMA(n+1-a)*GAMMA(1-2*
a)-1/2*GAMMA(2-2*a+2*n)*Sum((-1)^k*GAMMA(2*a-2*n+k)*GAMMA(1-2*a)/GAMMA
(2*a-2*n)/GAMMA(1-2*a+k),k = 1 ..
2*n-1)/GAMMA(2-4*a+2*n)/GAMMA(1-2*a+2*n)*GAMMA(1-3*a+2*n)*GAMMA(-3*a+2
*n+3/2)/GAMMA(3/2-a+n)/GAMMA(n+1-a);}{%
\maplemultiline{
\mathit{\ 40:\ }   {\displaystyle \frac {\Gamma (2 - 2\,a + 2\,n)
\,2^{(2\,n - 1 - 4\,a)}\,\Gamma (1 - 3\,a + 2\,n)\,\Gamma ( - 3\,
a + 2\,n + {\displaystyle \frac {3}{2}} )\,\Gamma (1 - 2\,a)}{
\Gamma (1 - 4\,a + 2\,n)\,\Gamma (2 - 4\,a + 2\,n)\,\Gamma (
{\displaystyle \frac {3}{2}}  - a + n)\,\Gamma (n + 1 - a)}}  \\ - 
{\displaystyle \frac {1}{2}}  
\Gamma (2 - 2\,a + 2\,n)\, \left(  \! {\displaystyle \sum _{k=1}
^{2\,n - 1}} \,{\displaystyle \frac {(-1)^{k}\,\Gamma (2\,a - 2\,
n + k)\,\Gamma (1 - 2\,a)}{\Gamma (2\,a - 2\,n)\,\Gamma (1 - 2\,a
 + k)}}  \!  \right) \,\Gamma (1 - 3\,a + 2\,n) \\
\Gamma ( - 3\,a + 2\,n + {\displaystyle \frac {3}{2}} ) \left/ 
{\vrule height0.80em width0em depth0.80em} \right. \!  \! (\Gamma
 (2 - 4\,a + 2\,n)\,\Gamma (1 - 2\,a + 2\,n)\,\Gamma (
{\displaystyle \frac {3}{2}}  - a + n) 
\Gamma (n + 1 - a)) }
}
\end{maplelatex}

\begin{maplelatex}
\mapleinline{inert}{2d}{` 41: `,
-1/2*GAMMA(2-2*a+2*n)/GAMMA(1-2*a+2*n)*GAMMA(1-3*a+2*n)/GAMMA(3/2-a+n)
/GAMMA(n+1-a)*GAMMA(a+1/2)*Sum((-1)^k*GAMMA(2*a-2*n+k)*GAMMA(1-2*a)/GA
MMA(2*a-2*n)/GAMMA(1-2*a+k),k = 1 ..
2*n-1)+GAMMA(2-2*a+2*n)*2^(2*n-1-4*a)/GAMMA(1-4*a+2*n)*GAMMA(1-3*a+2*n
)/GAMMA(3/2-a+n)/GAMMA(n+1-a)*GAMMA(1-2*a)*GAMMA(a+1/2);}{%
\maplemultiline{
\mathit{\ 41:\ }    - {\displaystyle \frac {1}{2}} \,
{\displaystyle \frac {\Gamma (2 - 2\,a + 2\,n)\,\Gamma (1 - 3\,a
 + 2\,n)\,\Gamma (a + {\displaystyle \frac {1}{2}} )\, \left( 
 \! {\displaystyle \sum _{k=1}^{2\,n - 1}} \,{\displaystyle 
\frac {(-1)^{k}\,\Gamma (2\,a - 2\,n + k)\,\Gamma (1 - 2\,a)}{
\Gamma (2\,a - 2\,n)\,\Gamma (1 - 2\,a + k)}}  \!  \right) }{
\Gamma (1 - 2\,a + 2\,n)\,\Gamma ({\displaystyle \frac {3}{2}} 
 - a + n)\,\Gamma (n + 1 - a)}}  \\
\mbox{} + {\displaystyle \frac {\Gamma (2 - 2\,a + 2\,n)\,2^{(2\,
n - 1 - 4\,a)}\,\Gamma (1 - 3\,a + 2\,n)\,\Gamma (1 - 2\,a)\,
\Gamma (a + {\displaystyle \frac {1}{2}} )}{\Gamma (1 - 4\,a + 2
\,n)\,\Gamma ({\displaystyle \frac {3}{2}}  - a + n)\,\Gamma (n
 + 1 - a)}}  }
}
\end{maplelatex}

\begin{maplelatex}
\mapleinline{inert}{2d}{` 42: `,
-1/2*GAMMA(2-2*a+2*n)/GAMMA(2-4*a+2*n)/GAMMA(1-2*a+2*n)*GAMMA(1-3*a+2*
n)/GAMMA(n+1-a)*GAMMA(1-2*a+n)*Sum((-1)^k*GAMMA(2*a-2*n+k)*GAMMA(1-2*a
)/GAMMA(2*a-2*n)/GAMMA(1-2*a+k),k = 1 ..
2*n-1)+GAMMA(2-2*a+2*n)*2^(2*n-1-4*a)/GAMMA(1-4*a+2*n)/GAMMA(2-4*a+2*n
)*GAMMA(1-3*a+2*n)/GAMMA(n+1-a)*GAMMA(1-2*a)*GAMMA(1-2*a+n);}{%
\maplemultiline{
\mathit{\ 42:\ }    - {\displaystyle \frac {1}{2}} \,
{\displaystyle \frac {\Gamma (2 - 2\,a + 2\,n)\,\Gamma (1 - 3\,a
 + 2\,n)\,\Gamma (1 - 2\,a + n)\, \left(  \! {\displaystyle \sum 
_{k=1}^{2\,n - 1}} \,{\displaystyle \frac {(-1)^{k}\,\Gamma (2\,a
 - 2\,n + k)\,\Gamma (1 - 2\,a)}{\Gamma (2\,a - 2\,n)\,\Gamma (1
 - 2\,a + k)}}  \!  \right) }{\Gamma (2 - 4\,a + 2\,n)\,\Gamma (1
 - 2\,a + 2\,n)\,\Gamma (n + 1 - a)}}  \\
\mbox{} + {\displaystyle \frac {\Gamma (2 - 2\,a + 2\,n)\,2^{(2\,
n - 1 - 4\,a)}\,\Gamma (1 - 3\,a + 2\,n)\,\Gamma (1 - 2\,a)\,
\Gamma (1 - 2\,a + n)}{\Gamma (1 - 4\,a + 2\,n)\,\Gamma (2 - 4\,a
 + 2\,n)\,\Gamma (n + 1 - a)}}  }
}
\end{maplelatex}

\begin{maplelatex}
\mapleinline{inert}{2d}{` 43: `,
-1/2*GAMMA(2-2*a+2*n)/GAMMA(2-4*a+2*n)/GAMMA(1-2*a+2*n)*GAMMA(1-3*a+2*
n)/GAMMA(3/2-a+n)*GAMMA(3/2-2*a+n)*Sum((-1)^k*GAMMA(2*a-2*n+k)*GAMMA(1
-2*a)/GAMMA(2*a-2*n)/GAMMA(1-2*a+k),k = 1 ..
2*n-1)+GAMMA(2-2*a+2*n)*2^(2*n-1-4*a)/GAMMA(1-4*a+2*n)/GAMMA(2-4*a+2*n
)*GAMMA(1-3*a+2*n)/GAMMA(3/2-a+n)*GAMMA(1-2*a)*GAMMA(3/2-2*a+n);}{%
\maplemultiline{
\mathit{\ 43:\ }    - {\displaystyle \frac {1}{2}} \,
{\displaystyle \frac {\Gamma (2 - 2\,a + 2\,n)\,\Gamma (1 - 3\,a
 + 2\,n)\,\Gamma ({\displaystyle \frac {3}{2}}  - 2\,a + n)\,
 \left(  \! {\displaystyle \sum _{k=1}^{2\,n - 1}} \,
{\displaystyle \frac {(-1)^{k}\,\Gamma (2\,a - 2\,n + k)\,\Gamma 
(1 - 2\,a)}{\Gamma (2\,a - 2\,n)\,\Gamma (1 - 2\,a + k)}}  \! 
 \right) }{\Gamma (2 - 4\,a + 2\,n)\,\Gamma (1 - 2\,a + 2\,n)\,
\Gamma ({\displaystyle \frac {3}{2}}  - a + n)}}  \\
\mbox{} + {\displaystyle \frac {\Gamma (2 - 2\,a + 2\,n)\,2^{(2\,
n - 1 - 4\,a)}\,\Gamma (1 - 3\,a + 2\,n)\,\Gamma (1 - 2\,a)\,
\Gamma ({\displaystyle \frac {3}{2}}  - 2\,a + n)}{\Gamma (1 - 4
\,a + 2\,n)\,\Gamma (2 - 4\,a + 2\,n)\,\Gamma ({\displaystyle 
\frac {3}{2}}  - a + n)}}  }
}
\end{maplelatex}

\begin{maplelatex}
\mapleinline{inert}{2d}{` 44: `,
GAMMA(2-2*a+2*n)*2^(2*n-1-4*a)/GAMMA(1-4*a+2*n)/GAMMA(2-4*a+2*n)*GAMMA
(1/2-2*a+n)*GAMMA(-3*a+2*n+3/2)/GAMMA(n+1-a)*GAMMA(1-2*a)-1/2*GAMMA(2-
2*a+2*n)*Sum((-1)^k*GAMMA(2*a-2*n+k)*GAMMA(1-2*a)/GAMMA(2*a-2*n)/GAMMA
(1-2*a+k),k = 1 ..
2*n-1)/GAMMA(2-4*a+2*n)/GAMMA(1-2*a+2*n)*GAMMA(1/2-2*a+n)*GAMMA(-3*a+2
*n+3/2)/GAMMA(n+1-a);}{%
\maplemultiline{
\mathit{\ 44:\ }   {\displaystyle \frac {\Gamma (2 - 2\,a + 2\,n)
\,2^{(2\,n - 1 - 4\,a)}\,\Gamma ({\displaystyle \frac {1}{2}}  - 
2\,a + n)\,\Gamma ( - 3\,a + 2\,n + {\displaystyle \frac {3}{2}} 
)\,\Gamma (1 - 2\,a)}{\Gamma (1 - 4\,a + 2\,n)\,\Gamma (2 - 4\,a
 + 2\,n)\,\Gamma (n + 1 - a)}}  \\ - {\displaystyle \frac {1}{2}} 
\Gamma (2 - 2\,a + 2\,n)\, \left(  \! {\displaystyle \sum _{k=1}
^{2\,n - 1}} \,{\displaystyle \frac {(-1)^{k}\,\Gamma (2\,a - 2\,
n + k)\,\Gamma (1 - 2\,a)}{\Gamma (2\,a - 2\,n)\,\Gamma (1 - 2\,a
 + k)}}  \!  \right) \,\Gamma ({\displaystyle \frac {1}{2}}  - 2
\,a + n) \\
\Gamma ( - 3\,a + 2\,n + {\displaystyle \frac {3}{2}} )/(\Gamma (
2 - 4\,a + 2\,n)\,\Gamma (1 - 2\,a + 2\,n)\,\Gamma (n + 1 - a))
 }
}
\end{maplelatex}

\begin{maplelatex}
\mapleinline{inert}{2d}{` 45: `,
-1/2*GAMMA(2-2*a+2*n)/GAMMA(2-4*a+2*n)/GAMMA(1-2*a+2*n)*GAMMA(-3*a+2*n
+3/2)/GAMMA(3/2-a+n)*GAMMA(1-2*a+n)*Sum((-1)^k*GAMMA(2*a-2*n+k)*GAMMA(
1-2*a)/GAMMA(2*a-2*n)/GAMMA(1-2*a+k),k = 1 ..
2*n-1)+GAMMA(2-2*a+2*n)*2^(2*n-1-4*a)/GAMMA(1-4*a+2*n)/GAMMA(2-4*a+2*n
)*GAMMA(-3*a+2*n+3/2)/GAMMA(3/2-a+n)*GAMMA(1-2*a)*GAMMA(1-2*a+n);}{%
\maplemultiline{
\mathit{\ 45:\ }    - {\displaystyle \frac {1}{2}} \Gamma (2 - 2\,
a + 2\,n)\,\Gamma ( - 3\,a + 2\,n + {\displaystyle \frac {3}{2}} 
)\,\Gamma (1 - 2\,a + n) \\
 \left(  \! {\displaystyle \sum _{k=1}^{2\,n - 1}} \,
{\displaystyle \frac {(-1)^{k}\,\Gamma (2\,a - 2\,n + k)\,\Gamma 
(1 - 2\,a)}{\Gamma (2\,a - 2\,n)\,\Gamma (1 - 2\,a + k)}}  \! 
 \right)  \left/ {\vrule height0.80em width0em depth0.80em}
 \right. \!  \! (\Gamma (2 - 4\,a + 2\,n)\,\Gamma (1 - 2\,a + 2\,
n) 
\Gamma ({\displaystyle \frac {3}{2}}  - a + n)) \\
\mbox{} + {\displaystyle \frac {\Gamma (2 - 2\,a + 2\,n)\,2^{(2\,
n - 1 - 4\,a)}\,\Gamma ( - 3\,a + 2\,n + {\displaystyle \frac {3
}{2}} )\,\Gamma (1 - 2\,a)\,\Gamma (1 - 2\,a + n)}{\Gamma (1 - 4
\,a + 2\,n)\,\Gamma (2 - 4\,a + 2\,n)\,\Gamma ({\displaystyle 
\frac {3}{2}}  - a + n)}}  }
}
\end{maplelatex}

\begin{maplelatex}
\mapleinline{inert}{2d}{` 46: `,
-1/2*GAMMA(2-2*a+2*n)/GAMMA(1-2*a+2*n)/GAMMA(3/2-a+n)*GAMMA(1/2+n-a)*S
um((-1)^k*GAMMA(2*a-2*n+k)*GAMMA(1-2*a)/GAMMA(2*a-2*n)/GAMMA(1-2*a+k),
k = 1 ..
2*n-1)+GAMMA(2-2*a+2*n)*2^(2*n-1-4*a)/GAMMA(1-4*a+2*n)/GAMMA(3/2-a+n)*
GAMMA(1-2*a)*GAMMA(1/2+n-a);}{%
\maplemultiline{
\mathit{\ 46:\ }    - {\displaystyle \frac {1}{2}} \,
{\displaystyle \frac {\Gamma (2 - 2\,a + 2\,n)\,\Gamma (
{\displaystyle \frac {1}{2}}  + n - a)\, \left(  \! 
{\displaystyle \sum _{k=1}^{2\,n - 1}} \,{\displaystyle \frac {(
-1)^{k}\,\Gamma (2\,a - 2\,n + k)\,\Gamma (1 - 2\,a)}{\Gamma (2\,
a - 2\,n)\,\Gamma (1 - 2\,a + k)}}  \!  \right) }{\Gamma (1 - 2\,
a + 2\,n)\,\Gamma ({\displaystyle \frac {3}{2}}  - a + n)}}  \\
\mbox{} + {\displaystyle \frac {\Gamma (2 - 2\,a + 2\,n)\,2^{(2\,
n - 1 - 4\,a)}\,\Gamma (1 - 2\,a)\,\Gamma ({\displaystyle \frac {
1}{2}}  + n - a)}{\Gamma (1 - 4\,a + 2\,n)\,\Gamma (
{\displaystyle \frac {3}{2}}  - a + n)}}  }
}
\end{maplelatex}

\begin{maplelatex}
\mapleinline{inert}{2d}{` 47: `,
GAMMA(2-2*a+2*n)*2^(2*n-1-4*a)/GAMMA(2-4*a+2*n)*GAMMA(1-2*a)-1/2*GAMMA
(2-2*a+2*n)*GAMMA(1-4*a+2*n)*Sum((-1)^k*GAMMA(2*a-2*n+k)*GAMMA(1-2*a)/
GAMMA(2*a-2*n)/GAMMA(1-2*a+k),k = 1 ..
2*n-1)/GAMMA(2-4*a+2*n)/GAMMA(1-2*a+2*n);}{%
\maplemultiline{
\mathit{\ 47:\ }   {\displaystyle \frac {\Gamma (2 - 2\,a + 2\,n)
\,2^{(2\,n - 1 - 4\,a)}\,\Gamma (1 - 2\,a)}{\Gamma (2 - 4\,a + 2
\,n)}}  \\
\mbox{} - {\displaystyle \frac {1}{2}} \,{\displaystyle \frac {
\Gamma (2 - 2\,a + 2\,n)\,\Gamma (1 - 4\,a + 2\,n)\, \left(  \! 
{\displaystyle \sum _{k=1}^{2\,n - 1}} \,{\displaystyle \frac {(
-1)^{k}\,\Gamma (2\,a - 2\,n + k)\,\Gamma (1 - 2\,a)}{\Gamma (2\,
a - 2\,n)\,\Gamma (1 - 2\,a + k)}}  \!  \right) }{\Gamma (2 - 4\,
a + 2\,n)\,\Gamma (1 - 2\,a + 2\,n)}}  }
}
\end{maplelatex}

\begin{maplelatex}
\mapleinline{inert}{2d}{` 48: `,
GAMMA(1+a)/GAMMA(e-b)/GAMMA(c)*(1+(-a+c-1)*(e-a-1-b)/(e-a-1)/a)*GAMMA(
e-a-1-b)*GAMMA(e)*GAMMA(c-1)/GAMMA(a)/GAMMA(e-a-1);}{%
\[
\mathit{\ 48:\ }   \,{\displaystyle \frac {\Gamma (1 + a)\,(1 + 
{\displaystyle \frac {( - a + c - 1)\,(e - a - 1 - b)}{(e - a - 1
)\,a}} )\,\Gamma (e - a - 1 - b)\,\Gamma (e)\,\Gamma (c - 1)}{
\Gamma (e - b)\,\Gamma (c)\,\Gamma (a)\,\Gamma (e - a - 1)}} 
\]
}
\end{maplelatex}

\begin{maplelatex}
\mapleinline{inert}{2d}{` 50: `,
GAMMA(1+a)/GAMMA(e+1)/GAMMA(c+1)*(1+(c-a)*(e-a)/(c-b+e-a)/a)*GAMMA(c-b
+1+e-a)*GAMMA(c)*GAMMA(e)/GAMMA(a)/GAMMA(c-b+e-a);}{%
\[
\mathit{\ 50:\ }   \,{\displaystyle \frac {\Gamma (1 + a)\,(1 + 
{\displaystyle \frac {(c - a)\,(e - a)}{(c - b + e - a)\,a}} )\,
\Gamma (c - b + 1 + e - a)\,\Gamma (c)\,\Gamma (e)}{\Gamma (e + 1
)\,\Gamma (c + 1)\,\Gamma (a)\,\Gamma (c - b + e - a)}} 
\]
}
\end{maplelatex}

\begin{maplelatex}
\mapleinline{inert}{2d}{` 51: `,
-(c-1+a+a*b-a*c-3*b)*Pi*(c+1)*c/sin(Pi*a)/(b-1-c)/(1+b-c)/(b-c)/GAMMA(
b+2-a)*GAMMA(2+b)/GAMMA(a)-Pi*b*(b+1)*(-b+1+a*b-a-a*c+3*c)/sin(Pi*a)/(
b-1-c)/(1+b-c)/(b-c)/GAMMA(2+c-a)*GAMMA(c+2)/GAMMA(a);}{%
\maplemultiline{
\mathit{\ 51:\ }    - {\displaystyle \frac {(c - 1 + a + a\,b - a
\,c - 3\,b)\,\pi \,(c + 1)\,c\,\Gamma (2 + b)}{\mathrm{sin}(\pi 
\,a)\,(b - 1 - c)\,(1 + b - c)\,(b - c)\,\Gamma (b + 2 - a)\,
\Gamma (a)}}  \\
\mbox{} - {\displaystyle \frac {\pi \,b\,(b + 1)\,( - b + 1 + a\,
b - a - a\,c + 3\,c)\,\Gamma (c + 2)}{\mathrm{sin}(\pi \,a)\,(b
 - 1 - c)\,(1 + b - c)\,(b - c)\,\Gamma (2 + c - a)\,\Gamma (a)}
}  }
}
\end{maplelatex}

\begin{maplelatex}
\mapleinline{inert}{2d}{` 52: `,
6*(5-b-2*c+a*b-a)/(b-2)/(-2+a)/(a-1)/(a-3)/(b-3)/(b-1)*GAMMA(-b+2+c-a)
/GAMMA(c-a)*GAMMA(c)/GAMMA(c-b)+6*(c-2)*(c-1)*(3+2*c+a*b-3*a-3*b)/(b-2
)/(-2+a)/(a-1)/(a-3)/(b-3)/(b-1);}{%
\maplemultiline{
\mathit{\ 52:\ }   {\displaystyle \frac {6\,(5 - b - 2\,c + a\,b
 - a)\,\Gamma ( - b + 2 + c - a)\,\Gamma (c)}{(b - 2)\,( - 2 + a)
\,(a - 1)\,(a - 3)\,(b - 3)\,(b - 1)\,\Gamma (c - a)\,\Gamma (c
 - b)}}  \\
\mbox{} + {\displaystyle \frac {6\,(c - 2)\,(c - 1)\,(3 + 2\,c + 
a\,b - 3\,a - 3\,b)}{(b - 2)\,( - 2 + a)\,(a - 1)\,(a - 3)\,(b - 
3)\,(b - 1)}}  }
}
\end{maplelatex}

\begin{maplelatex}
\mapleinline{inert}{2d}{` 53: `,
(a+5-4*c-a*c+c^2+b+a*b-b*c)/(2+a-c)/(-c+b+2)/(1+b-c)/(-c+3+b)/(a+3-c)/
(1+a-c)/GAMMA(b)/GAMMA(a)*GAMMA(c)*GAMMA(4+a-c+b)+(a+2-c+b)*(a+3-c+b)*
(c-2)*(c-1)*(3-4*c-a*c+c^2+3*a+a*b+3*b-b*c)/(2+a-c)/(-c+b+2)/(1+b-c)/(
-c+3+b)/(a+3-c)/(1+a-c);}{%
\maplemultiline{
\mathit{\ 53:\ }   \,{\displaystyle \frac {(a + 5 - 4\,c - a\,c + 
c^{2} + b + a\,b - b\,c)\,\Gamma (c)\,\Gamma (4 + a - c + b)}{(2
 + a - c)\,( - c + b + 2)\,(1 + b - c)\,( - c + 3 + b)\,(a + 3 - 
c)\,(1 + a - c)\,\Gamma (b)\,\Gamma (a)}} \\ + {\displaystyle 
\frac {(a + 2 - c + b)\,(a + 3 - c + b)\,(c - 2)\,(c - 1)\,(3 - 4
\,c - a\,c + c^{2} + 3\,a + a\,b + 3\,b - b\,c)}{(2 + a - c)\,(
 - c + b + 2)\,(1 + b - c)\,( - c + 3 + b)\,(a + 3 - c)\,(1 + a
 - c)}} 
  }
}
\end{maplelatex}

\begin{maplelatex}
\mapleinline{inert}{2d}{` 54: `,
GAMMA(n+1)/GAMMA(1-a+n-m)*GAMMA(-a+1)/GAMMA(1-b+m)*GAMMA(-b+1);}{%
\[
\mathit{\ 54:\ }   \,{\displaystyle \frac {\Gamma (n + 1)\,\Gamma 
( - a + 1)\,\Gamma ( - b + 1)}{\Gamma (1 - a + n - m)\,\Gamma (1
 - b + m)}} 
\]
}
\end{maplelatex}

\begin{maplelatex}
\mapleinline{inert}{2d}{` 55: `, 0;}{%
\[
\mathit{\ 55:\ }   \,0
\]
}
\end{maplelatex}

\begin{maplelatex}
\mapleinline{inert}{2d}{` 56: `, 0;}{%
\[
\mathit{\ 56:\ }   \,0
\]
}
\end{maplelatex}
\begin{flushright}
$n\geq L+m+1$
\end{flushright} 

\begin{maplelatex}
\mapleinline{inert}{2d}{` 57: `, 0;}{%
\[
\mathit{\ 57:\ }   \,0
\]
}
\end{maplelatex}
\begin{flushright}
$n\geq L+m+1$
\end{flushright} 

\begin{maplelatex}
\mapleinline{inert}{2d}{` 61: `,
(1/2*m*(m-1)*(b+n-m)*(b-m+n+1)+1/2*(2+n-m)*(-m+1+n)*(-2+a)*(a-1)+(2+n-
m)*m*(b+n-m)*(a-1))*GAMMA(n+1)*GAMMA(-b+1)*GAMMA(-a+1)/GAMMA(1-b+m)/GA
MMA(3-a+n-m);}{%
\maplemultiline{
\mathit{\ 61:\ }   ({\displaystyle \frac {m\,(m - 1)\,(b + n - m)
\,(b - m + n + 1)}{2}}  + {\displaystyle \frac {(2 + n - m)\,( - 
m + 1 + n)\,( - 2 + a)\,(a - 1)}{2}}  \\
\mbox{} + (2 + n - m)\,m\,(b + n - m)\,(a - 1))\Gamma (n + 1)\,
\Gamma ( - b + 1)\,\Gamma ( - a + 1)/(\Gamma (1 - b + m) 
\Gamma (3 - a + n - m)) }
}
\end{maplelatex}

\begin{maplelatex}
\mapleinline{inert}{2d}{` 62: `,
-GAMMA((a*n+b*n-b+n-b^2)/(n-b))*GAMMA(-a+1)*(a*n-b*n+n^2-b+n)*GAMMA(-b
)/GAMMA((a*n-b*n+n^2-2*b+2*n)/(n-b))/GAMMA(-n+b+1)/a/GAMMA(-a-n)/GAMMA
(n+1-b);}{%
\[
\mathit{\ 62:\ }   \, - {\displaystyle \frac {\Gamma (
{\displaystyle \frac {a\,n + b\,n - b + n - b^{2}}{n - b}} )\,
\Gamma ( - a + 1)\,(a\,n - b\,n + n^{2} - b + n)\,\Gamma ( - b)}{
\Gamma ({\displaystyle \frac {a\,n - b\,n + n^{2} - 2\,b + 2\,n}{
n - b}} )\,\Gamma ( - n + b + 1)\,a\,\Gamma ( - a - n)\,\Gamma (n
 + 1 - b)}} 
\]
}
\end{maplelatex}

\begin{maplelatex}
\mapleinline{inert}{2d}{` 66: `,
1/2*1/GAMMA(3-2*a-2*n)*GAMMA(3-3*a-2*n)*GAMMA(7/2-3*a-2*n)/GAMMA(2-a-n
)/GAMMA(5/2-a-n)*GAMMA(4-2*a-2*n)/GAMMA(4-4*a-2*n)*Sum((-1)^k*GAMMA(2*
a+k)*GAMMA(3-2*a-2*n)/GAMMA(2*a)/GAMMA(-2*a-2*n+3+k),k = 0 ..
2*n-2)+2^(-4*a-2*n+1)*GAMMA(1-2*a)*GAMMA(3-3*a-2*n)*GAMMA(7/2-3*a-2*n)
/GAMMA(2-a-n)/GAMMA(5/2-a-n)/GAMMA(3-4*a-2*n)*GAMMA(4-2*a-2*n)/GAMMA(4
-4*a-2*n);}{%
\maplemultiline{
\mathit{\ 66:\ }   {\displaystyle \frac {1}{2}} \Gamma (3 - 3\,a
 - 2\,n)\,\Gamma ({\displaystyle \frac {7}{2}}  - 3\,a - 2\,n)\,
\Gamma (4 - 2\,a - 2\,n) \\
 \left(  \! {\displaystyle \sum _{k=0}^{2\,n - 2}} \,
{\displaystyle \frac {(-1)^{k}\,\Gamma (2\,a + k)\,\Gamma (3 - 2
\,a - 2\,n)}{\Gamma (2\,a)\,\Gamma ( - 2\,a - 2\,n + 3 + k)}} 
 \!  \right)  \left/ {\vrule height0.80em width0em depth0.80em}
 \right. \!  \! (\Gamma (3 - 2\,a - 2\,n)\,\Gamma (2 - a - n) \\
\Gamma ({\displaystyle \frac {5}{2}}  - a - n)\,\Gamma (4 - 4\,a
 - 2\,n)) \\
\mbox{} + {\displaystyle \frac {2^{( - 4\,a - 2\,n + 1)}\,\Gamma 
(1 - 2\,a)\,\Gamma (3 - 3\,a - 2\,n)\,\Gamma ({\displaystyle 
\frac {7}{2}}  - 3\,a - 2\,n)\,\Gamma (4 - 2\,a - 2\,n)}{\Gamma (
2 - a - n)\,\Gamma ({\displaystyle \frac {5}{2}}  - a - n)\,
\Gamma (3 - 4\,a - 2\,n)\,\Gamma (4 - 4\,a - 2\,n)}}  }
}
\end{maplelatex}

\begin{maplelatex}
\mapleinline{inert}{2d}{` 67: `,
1/2*1/GAMMA(3-2*a-2*n)*GAMMA(3-3*a-2*n)/GAMMA(5/2-a-n)*GAMMA(4-2*a-2*n
)/GAMMA(4-4*a-2*n)*GAMMA(5/2-2*a-n)*Sum((-1)^k*GAMMA(2*a+k)*GAMMA(3-2*
a-2*n)/GAMMA(2*a)/GAMMA(-2*a-2*n+3+k),k = 0 ..
2*n-2)+2^(-4*a-2*n+1)*GAMMA(1-2*a)*GAMMA(3-3*a-2*n)/GAMMA(5/2-a-n)/GAM
MA(3-4*a-2*n)*GAMMA(4-2*a-2*n)/GAMMA(4-4*a-2*n)*GAMMA(5/2-2*a-n);}{%
\maplemultiline{
\mathit{\ 67:\ }   {\displaystyle \frac {1}{2}} \,{\displaystyle 
\frac {\Gamma (3 - 3\,a - 2\,n)\,\Gamma (4 - 2\,a - 2\,n)\,\Gamma
 ({\displaystyle \frac {5}{2}}  - 2\,a - n)\, \left(  \! 
{\displaystyle \sum _{k=0}^{2\,n - 2}} \,{\displaystyle \frac {(
-1)^{k}\,\Gamma (2\,a + k)\,\Gamma (3 - 2\,a - 2\,n)}{\Gamma (2\,
a)\,\Gamma ( - 2\,a - 2\,n + 3 + k)}}  \!  \right) }{\Gamma (3 - 
2\,a - 2\,n)\,\Gamma ({\displaystyle \frac {5}{2}}  - a - n)\,
\Gamma (4 - 4\,a - 2\,n)}}  \\
\mbox{} + {\displaystyle \frac {2^{( - 4\,a - 2\,n + 1)}\,\Gamma 
(1 - 2\,a)\,\Gamma (3 - 3\,a - 2\,n)\,\Gamma (4 - 2\,a - 2\,n)\,
\Gamma ({\displaystyle \frac {5}{2}}  - 2\,a - n)}{\Gamma (
{\displaystyle \frac {5}{2}}  - a - n)\,\Gamma (3 - 4\,a - 2\,n)
\,\Gamma (4 - 4\,a - 2\,n)}}  }
}
\end{maplelatex}

\begin{maplelatex}
\mapleinline{inert}{2d}{` 68: `,
1/2*1/GAMMA(3-2*a-2*n)*GAMMA(3-3*a-2*n)/GAMMA(2-a-n)*GAMMA(4-2*a-2*n)/
GAMMA(4-4*a-2*n)*GAMMA(-2*a-n+2)*Sum((-1)^k*GAMMA(2*a+k)*GAMMA(3-2*a-2
*n)/GAMMA(2*a)/GAMMA(-2*a-2*n+3+k),k = 0 ..
2*n-2)+2^(-4*a-2*n+1)*GAMMA(1-2*a)*GAMMA(3-3*a-2*n)/GAMMA(2-a-n)/GAMMA
(3-4*a-2*n)*GAMMA(4-2*a-2*n)/GAMMA(4-4*a-2*n)*GAMMA(-2*a-n+2);}{%
\maplemultiline{
\mathit{\ 68:\ }   {\displaystyle \frac {1}{2}} \Gamma (3 - 3\,a
 - 2\,n)\,\Gamma (4 - 2\,a - 2\,n)\,\Gamma ( - 2\,a - n + 2) \\
 \left(  \! {\displaystyle \sum _{k=0}^{2\,n - 2}} \,
{\displaystyle \frac {(-1)^{k}\,\Gamma (2\,a + k)\,\Gamma (3 - 2
\,a - 2\,n)}{\Gamma (2\,a)\,\Gamma ( - 2\,a - 2\,n + 3 + k)}} 
 \!  \right) /(\Gamma (3 - 2\,a - 2\,n)\,\Gamma (2 - a - n)
\Gamma (4 - 4\,a - 2\,n)) \\
\mbox{} + {\displaystyle \frac {2^{( - 4\,a - 2\,n + 1)}\,\Gamma 
(1 - 2\,a)\,\Gamma (3 - 3\,a - 2\,n)\,\Gamma (4 - 2\,a - 2\,n)\,
\Gamma ( - 2\,a - n + 2)}{\Gamma (2 - a - n)\,\Gamma (3 - 4\,a - 
2\,n)\,\Gamma (4 - 4\,a - 2\,n)}}  }
}
\end{maplelatex}

\begin{maplelatex}
\mapleinline{inert}{2d}{` 69: `,
1/2*1/GAMMA(3-2*a-2*n)*GAMMA(7/2-3*a-2*n)/GAMMA(2-a-n)/GAMMA(5/2-a-n)*
GAMMA(4-2*a-2*n)*GAMMA(a)*Sum((-1)^k*GAMMA(2*a+k)*GAMMA(3-2*a-2*n)/GAM
MA(2*a)/GAMMA(-2*a-2*n+3+k),k = 0 ..
2*n-2)+2^(-4*a-2*n+1)*GAMMA(1-2*a)*GAMMA(7/2-3*a-2*n)/GAMMA(2-a-n)/GAM
MA(5/2-a-n)/GAMMA(3-4*a-2*n)*GAMMA(4-2*a-2*n)*GAMMA(a);}{%
\maplemultiline{
\mathit{\ 69:\ }   {\displaystyle \frac {1}{2}} \,{\displaystyle 
\frac {\Gamma ({\displaystyle \frac {7}{2}}  - 3\,a - 2\,n)\,
\Gamma (4 - 2\,a - 2\,n)\,\Gamma (a)\, \left(  \! {\displaystyle 
\sum _{k=0}^{2\,n - 2}} \,{\displaystyle \frac {(-1)^{k}\,\Gamma 
(2\,a + k)\,\Gamma (3 - 2\,a - 2\,n)}{\Gamma (2\,a)\,\Gamma ( - 2
\,a - 2\,n + 3 + k)}}  \!  \right) }{\Gamma (3 - 2\,a - 2\,n)\,
\Gamma (2 - a - n)\,\Gamma ({\displaystyle \frac {5}{2}}  - a - n
)}}  \\
\mbox{} + {\displaystyle \frac {2^{( - 4\,a - 2\,n + 1)}\,\Gamma 
(1 - 2\,a)\,\Gamma ({\displaystyle \frac {7}{2}}  - 3\,a - 2\,n)
\,\Gamma (4 - 2\,a - 2\,n)\,\Gamma (a)}{\Gamma (2 - a - n)\,
\Gamma ({\displaystyle \frac {5}{2}}  - a - n)\,\Gamma (3 - 4\,a
 - 2\,n)}}  }
}
\end{maplelatex}

\begin{maplelatex}
\mapleinline{inert}{2d}{` 70: `,
1/2*1/GAMMA(3-2*a-2*n)*GAMMA(7/2-3*a-2*n)/GAMMA(5/2-a-n)*GAMMA(4-2*a-2
*n)/GAMMA(4-4*a-2*n)*GAMMA(-2*a-n+2)*Sum((-1)^k*GAMMA(2*a+k)*GAMMA(3-2
*a-2*n)/GAMMA(2*a)/GAMMA(-2*a-2*n+3+k),k = 0 ..
2*n-2)+2^(-4*a-2*n+1)*GAMMA(1-2*a)*GAMMA(7/2-3*a-2*n)/GAMMA(5/2-a-n)/G
AMMA(3-4*a-2*n)*GAMMA(4-2*a-2*n)/GAMMA(4-4*a-2*n)*GAMMA(-2*a-n+2);}{%
\maplemultiline{
\mathit{\ 70:\ }   {\displaystyle \frac {1}{2}} \Gamma (
{\displaystyle \frac {7}{2}}  - 3\,a - 2\,n)\,\Gamma (4 - 2\,a - 
2\,n)\,\Gamma ( - 2\,a - n + 2) \\
 \left(  \! {\displaystyle \sum _{k=0}^{2\,n - 2}} \,
{\displaystyle \frac {(-1)^{k}\,\Gamma (2\,a + k)\,\Gamma (3 - 2
\,a - 2\,n)}{\Gamma (2\,a)\,\Gamma ( - 2\,a - 2\,n + 3 + k)}} 
 \!  \right)  \left/ {\vrule height0.80em width0em depth0.80em}
 \right. \!  \! (\Gamma (3 - 2\,a - 2\,n)\,\Gamma (
{\displaystyle \frac {5}{2}}  - a - n) \\
\Gamma (4 - 4\,a - 2\,n)) \\
\mbox{} + {\displaystyle \frac {2^{( - 4\,a - 2\,n + 1)}\,\Gamma 
(1 - 2\,a)\,\Gamma ({\displaystyle \frac {7}{2}}  - 3\,a - 2\,n)
\,\Gamma (4 - 2\,a - 2\,n)\,\Gamma ( - 2\,a - n + 2)}{\Gamma (
{\displaystyle \frac {5}{2}}  - a - n)\,\Gamma (3 - 4\,a - 2\,n)
\,\Gamma (4 - 4\,a - 2\,n)}}  }
}
\end{maplelatex}

\begin{maplelatex}
\mapleinline{inert}{2d}{` 71: `,
1/2*1/GAMMA(3-2*a-2*n)*GAMMA(7/2-3*a-2*n)/GAMMA(2-a-n)*GAMMA(4-2*a-2*n
)/GAMMA(4-4*a-2*n)*GAMMA(3/2-2*a-n)*Sum((-1)^k*GAMMA(2*a+k)*GAMMA(3-2*
a-2*n)/GAMMA(2*a)/GAMMA(-2*a-2*n+3+k),k = 0 ..
2*n-2)+2^(-4*a-2*n+1)*GAMMA(1-2*a)*GAMMA(7/2-3*a-2*n)/GAMMA(2-a-n)/GAM
MA(3-4*a-2*n)*GAMMA(4-2*a-2*n)/GAMMA(4-4*a-2*n)*GAMMA(3/2-2*a-n);}{%
\maplemultiline{
\mathit{\ 71:\ }   {\displaystyle \frac {1}{2}} \,{\displaystyle 
\frac {\Gamma ({\displaystyle \frac {7}{2}}  - 3\,a - 2\,n)\,
\Gamma (4 - 2\,a - 2\,n)\,\Gamma ({\displaystyle \frac {3}{2}} 
 - 2\,a - n)\, \left(  \! {\displaystyle \sum _{k=0}^{2\,n - 2}} 
\,{\displaystyle \frac {(-1)^{k}\,\Gamma (2\,a + k)\,\Gamma (3 - 
2\,a - 2\,n)}{\Gamma (2\,a)\,\Gamma ( - 2\,a - 2\,n + 3 + k)}} 
 \!  \right) }{\Gamma (3 - 2\,a - 2\,n)\,\Gamma (2 - a - n)\,
\Gamma (4 - 4\,a - 2\,n)}}  \\
\mbox{} + {\displaystyle \frac {2^{( - 4\,a - 2\,n + 1)}\,\Gamma 
(1 - 2\,a)\,\Gamma ({\displaystyle \frac {7}{2}}  - 3\,a - 2\,n)
\,\Gamma (4 - 2\,a - 2\,n)\,\Gamma ({\displaystyle \frac {3}{2}} 
 - 2\,a - n)}{\Gamma (2 - a - n)\,\Gamma (3 - 4\,a - 2\,n)\,
\Gamma (4 - 4\,a - 2\,n)}}  }
}
\end{maplelatex}

\begin{maplelatex}
\mapleinline{inert}{2d}{` 72: `,
1/2*1/GAMMA(3-2*a-2*n)*GAMMA(-a-n+1)/GAMMA(2-a-n)*GAMMA(4-2*a-2*n)*Sum
((-1)^k*GAMMA(2*a+k)*GAMMA(3-2*a-2*n)/GAMMA(2*a)/GAMMA(-2*a-2*n+3+k),k
= 0 ..
2*n-2)+2^(-4*a-2*n+1)*GAMMA(1-2*a)*GAMMA(-a-n+1)/GAMMA(2-a-n)/GAMMA(3-
4*a-2*n)*GAMMA(4-2*a-2*n);}{%
\maplemultiline{
\mathit{\ 72:\ }   {\displaystyle \frac {1}{2}} \,{\displaystyle 
\frac {\Gamma ( - a - n + 1)\,\Gamma (4 - 2\,a - 2\,n)\, \left( 
 \! {\displaystyle \sum _{k=0}^{2\,n - 2}} \,{\displaystyle 
\frac {(-1)^{k}\,\Gamma (2\,a + k)\,\Gamma (3 - 2\,a - 2\,n)}{
\Gamma (2\,a)\,\Gamma ( - 2\,a - 2\,n + 3 + k)}}  \!  \right) }{
\Gamma (3 - 2\,a - 2\,n)\,\Gamma (2 - a - n)}}  \\
\mbox{} + {\displaystyle \frac {2^{( - 4\,a - 2\,n + 1)}\,\Gamma 
(1 - 2\,a)\,\Gamma ( - a - n + 1)\,\Gamma (4 - 2\,a - 2\,n)}{
\Gamma (2 - a - n)\,\Gamma (3 - 4\,a - 2\,n)}}  }
}
\end{maplelatex}

\begin{maplelatex}
\mapleinline{inert}{2d}{` 73: `,
1/2*1/GAMMA(3-2*a-2*n)*GAMMA(3-4*a-2*n)*GAMMA(4-2*a-2*n)/GAMMA(4-4*a-2
*n)*Sum((-1)^k*GAMMA(2*a+k)*GAMMA(3-2*a-2*n)/GAMMA(2*a)/GAMMA(-2*a-2*n
+3+k),k = 0 ..
2*n-2)+2^(-4*a-2*n+1)*GAMMA(1-2*a)*GAMMA(4-2*a-2*n)/GAMMA(4-4*a-2*n);}
{%
\maplemultiline{
\mathit{\ 73:\ }   {\displaystyle \frac {1}{2}} \,{\displaystyle 
\frac {\Gamma (3 - 4\,a - 2\,n)\,\Gamma (4 - 2\,a - 2\,n)\,
 \left(  \! {\displaystyle \sum _{k=0}^{2\,n - 2}} \,
{\displaystyle \frac {(-1)^{k}\,\Gamma (2\,a + k)\,\Gamma (3 - 2
\,a - 2\,n)}{\Gamma (2\,a)\,\Gamma ( - 2\,a - 2\,n + 3 + k)}} 
 \!  \right) }{\Gamma (3 - 2\,a - 2\,n)\,\Gamma (4 - 4\,a - 2\,n)
}}  \\
\mbox{} + {\displaystyle \frac {2^{( - 4\,a - 2\,n + 1)}\,\Gamma 
(1 - 2\,a)\,\Gamma (4 - 2\,a - 2\,n)}{\Gamma (4 - 4\,a - 2\,n)}} 
 }
}
\end{maplelatex}

\begin{maplelatex}
\mapleinline{inert}{2d}{` 74: `,
-1/2*GAMMA(1-2*a+2*n)/GAMMA(2*n-2*a)*GAMMA(-3*a+2*n)/GAMMA(n+1-a)/GAMM
A(1/2+n-a)*GAMMA(a+1/2)*Sum((-1)^k*GAMMA(2*a-2*n+1+k)*GAMMA(1-2*a)/GAM
MA(2*a-2*n+1)/GAMMA(1-2*a+k),k = 1 ..
2*n-2)+GAMMA(1-2*a+2*n)*2^(2*n-2-4*a)/GAMMA(-4*a+2*n)*GAMMA(-3*a+2*n)/
GAMMA(n+1-a)/GAMMA(1/2+n-a)*GAMMA(1-2*a)*GAMMA(a+1/2);}{%
\maplemultiline{
\mathit{\ 74:\ }    - {\displaystyle \frac {1}{2}} \,
{\displaystyle \frac {\Gamma (1 - 2\,a + 2\,n)\,\Gamma ( - 3\,a
 + 2\,n)\,\Gamma (a + {\displaystyle \frac {1}{2}} )\, \left( 
 \! {\displaystyle \sum _{k=1}^{2\,n - 2}} \,{\displaystyle 
\frac {(-1)^{k}\,\Gamma (2\,a - 2\,n + 1 + k)\,\Gamma (1 - 2\,a)
}{\Gamma (2\,a - 2\,n + 1)\,\Gamma (1 - 2\,a + k)}}  \!  \right) 
}{\Gamma (2\,n - 2\,a)\,\Gamma (n + 1 - a)\,\Gamma (
{\displaystyle \frac {1}{2}}  + n - a)}}  \\
\mbox{} + {\displaystyle \frac {\Gamma (1 - 2\,a + 2\,n)\,2^{(2\,
n - 2 - 4\,a)}\,\Gamma ( - 3\,a + 2\,n)\,\Gamma (1 - 2\,a)\,
\Gamma (a + {\displaystyle \frac {1}{2}} )}{\Gamma ( - 4\,a + 2\,
n)\,\Gamma (n + 1 - a)\,\Gamma ({\displaystyle \frac {1}{2}}  + n
 - a)}}  }
}
\end{maplelatex}

\begin{maplelatex}
\mapleinline{inert}{2d}{` 75: `,
-1/2*GAMMA(1-2*a+2*n)/GAMMA(1-4*a+2*n)/GAMMA(2*n-2*a)*GAMMA(-3*a+2*n)/
GAMMA(n+1-a)*GAMMA(1-2*a+n)*Sum((-1)^k*GAMMA(2*a-2*n+1+k)*GAMMA(1-2*a)
/GAMMA(2*a-2*n+1)/GAMMA(1-2*a+k),k = 1 ..
2*n-2)+GAMMA(1-2*a+2*n)*2^(2*n-2-4*a)/GAMMA(-4*a+2*n)/GAMMA(1-4*a+2*n)
*GAMMA(-3*a+2*n)/GAMMA(n+1-a)*GAMMA(1-2*a)*GAMMA(1-2*a+n);}{%
\maplemultiline{
\mathit{\ 75:\ }    - {\displaystyle \frac {1}{2}} \Gamma (1 - 2\,
a + 2\,n)\,\Gamma ( - 3\,a + 2\,n)\,\Gamma (1 - 2\,a + n) \\
 \left(  \! {\displaystyle \sum _{k=1}^{2\,n - 2}} \,
{\displaystyle \frac {(-1)^{k}\,\Gamma (2\,a - 2\,n + 1 + k)\,
\Gamma (1 - 2\,a)}{\Gamma (2\,a - 2\,n + 1)\,\Gamma (1 - 2\,a + k
)}}  \!  \right) /(\Gamma (1 - 4\,a + 2\,n)\,\Gamma (2\,n - 2\,a)
 \\
\Gamma (n + 1 - a)) \\
\mbox{} + {\displaystyle \frac {\Gamma (1 - 2\,a + 2\,n)\,2^{(2\,
n - 2 - 4\,a)}\,\Gamma ( - 3\,a + 2\,n)\,\Gamma (1 - 2\,a)\,
\Gamma (1 - 2\,a + n)}{\Gamma ( - 4\,a + 2\,n)\,\Gamma (1 - 4\,a
 + 2\,n)\,\Gamma (n + 1 - a)}}  }
}
\end{maplelatex}

\begin{maplelatex}
\mapleinline{inert}{2d}{` 76: `,
-1/2*GAMMA(1-2*a+2*n)/GAMMA(2*n-2*a)*GAMMA(-3*a+2*n+1/2)/GAMMA(n+1-a)/
GAMMA(1/2+n-a)*GAMMA(a)*Sum((-1)^k*GAMMA(2*a-2*n+1+k)*GAMMA(1-2*a)/GAM
MA(2*a-2*n+1)/GAMMA(1-2*a+k),k = 1 ..
2*n-2)+GAMMA(1-2*a+2*n)*2^(2*n-2-4*a)/GAMMA(-4*a+2*n)*GAMMA(-3*a+2*n+1
/2)/GAMMA(n+1-a)/GAMMA(1/2+n-a)*GAMMA(1-2*a)*GAMMA(a);}{%
\maplemultiline{
\mathit{\ 76:\ }    - {\displaystyle \frac {1}{2}} \,
{\displaystyle \frac {\Gamma (1 - 2\,a + 2\,n)\,\Gamma ( - 3\,a
 + 2\,n + {\displaystyle \frac {1}{2}} )\,\Gamma (a)\, \left( 
 \! {\displaystyle \sum _{k=1}^{2\,n - 2}} \,{\displaystyle 
\frac {(-1)^{k}\,\Gamma (2\,a - 2\,n + 1 + k)\,\Gamma (1 - 2\,a)
}{\Gamma (2\,a - 2\,n + 1)\,\Gamma (1 - 2\,a + k)}}  \!  \right) 
}{\Gamma (2\,n - 2\,a)\,\Gamma (n + 1 - a)\,\Gamma (
{\displaystyle \frac {1}{2}}  + n - a)}}  \\
\mbox{} + {\displaystyle \frac {\Gamma (1 - 2\,a + 2\,n)\,2^{(2\,
n - 2 - 4\,a)}\,\Gamma ( - 3\,a + 2\,n + {\displaystyle \frac {1
}{2}} )\,\Gamma (1 - 2\,a)\,\Gamma (a)}{\Gamma ( - 4\,a + 2\,n)\,
\Gamma (n + 1 - a)\,\Gamma ({\displaystyle \frac {1}{2}}  + n - a
)}}  }
}
\end{maplelatex}

\begin{maplelatex}
\mapleinline{inert}{2d}{` 77: `,
GAMMA(1-2*a+2*n)*2^(2*n-2-4*a)/GAMMA(-4*a+2*n)/GAMMA(1-4*a+2*n)*GAMMA(
n-2*a)*GAMMA(-3*a+2*n+1/2)/GAMMA(1/2+n-a)*GAMMA(1-2*a)-1/2*GAMMA(1-2*a
+2*n)*Sum((-1)^k*GAMMA(2*a-2*n+1+k)*GAMMA(1-2*a)/GAMMA(2*a-2*n+1)/GAMM
A(1-2*a+k),k = 1 ..
2*n-2)/GAMMA(1-4*a+2*n)/GAMMA(2*n-2*a)*GAMMA(n-2*a)*GAMMA(-3*a+2*n+1/2
)/GAMMA(1/2+n-a);}{%
\maplemultiline{
\mathit{\ 77:\ }   {\displaystyle \frac {\Gamma (1 - 2\,a + 2\,n)
\,2^{(2\,n - 2 - 4\,a)}\,\Gamma (n - 2\,a)\,\Gamma ( - 3\,a + 2\,
n + {\displaystyle \frac {1}{2}} )\,\Gamma (1 - 2\,a)}{\Gamma (
 - 4\,a + 2\,n)\,\Gamma (1 - 4\,a + 2\,n)\,\Gamma (
{\displaystyle \frac {1}{2}}  + n - a)}} \\ - {\displaystyle 
\frac {1}{2}}  
\Gamma (1 - 2\,a + 2\,n)\, \left(  \! {\displaystyle \sum _{k=1}
^{2\,n - 2}} \,{\displaystyle \frac {(-1)^{k}\,\Gamma (2\,a - 2\,
n + 1 + k)\,\Gamma (1 - 2\,a)}{\Gamma (2\,a - 2\,n + 1)\,\Gamma (
1 - 2\,a + k)}}  \!  \right) \,\Gamma (n - 2\,a) \\
\Gamma ( - 3\,a + 2\,n + {\displaystyle \frac {1}{2}} ) \left/ 
{\vrule height0.80em width0em depth0.80em} \right. \!  \! (\Gamma
 (1 - 4\,a + 2\,n)\,\Gamma (2\,n - 2\,a)\,\Gamma ({\displaystyle 
\frac {1}{2}}  + n - a)) }
}
\end{maplelatex}

\begin{maplelatex}
\mapleinline{inert}{2d}{` 78: `,
GAMMA(1-2*a+2*n)*2^(2*n-2-4*a)/GAMMA(-4*a+2*n)*GAMMA(-a+n-1/2)/GAMMA(1
/2+n-a)*GAMMA(1-2*a)-1/2*GAMMA(1-2*a+2*n)*Sum((-1)^k*GAMMA(2*a-2*n+1+k
)*GAMMA(1-2*a)/GAMMA(2*a-2*n+1)/GAMMA(1-2*a+k),k = 1 ..
2*n-2)/GAMMA(2*n-2*a)*GAMMA(-a+n-1/2)/GAMMA(1/2+n-a);}{%
\maplemultiline{
\mathit{\ 78:\ }   {\displaystyle \frac {\Gamma (1 - 2\,a + 2\,n)
\,2^{(2\,n - 2 - 4\,a)}\,\Gamma ( - a + n - {\displaystyle 
\frac {1}{2}} )\,\Gamma (1 - 2\,a)}{\Gamma ( - 4\,a + 2\,n)\,
\Gamma ({\displaystyle \frac {1}{2}}  + n - a)}}  \\
\mbox{} - {\displaystyle \frac {1}{2}} \,{\displaystyle \frac {
\Gamma (1 - 2\,a + 2\,n)\, \left(  \! {\displaystyle \sum _{k=1}
^{2\,n - 2}} \,{\displaystyle \frac {(-1)^{k}\,\Gamma (2\,a - 2\,
n + 1 + k)\,\Gamma (1 - 2\,a)}{\Gamma (2\,a - 2\,n + 1)\,\Gamma (
1 - 2\,a + k)}}  \!  \right) \,\Gamma ( - a + n - {\displaystyle 
\frac {1}{2}} )}{\Gamma (2\,n - 2\,a)\,\Gamma ({\displaystyle 
\frac {1}{2}}  + n - a)}}  }
}
\end{maplelatex}

\begin{maplelatex}
\mapleinline{inert}{2d}{` 79: `,
-1/2*GAMMA(1-2*a+2*n)/GAMMA(2*n-2*a)/GAMMA(n+1-a)*GAMMA(n-a)*Sum((-1)^
k*GAMMA(2*a-2*n+1+k)*GAMMA(1-2*a)/GAMMA(2*a-2*n+1)/GAMMA(1-2*a+k),k =
1 ..
2*n-2)+GAMMA(1-2*a+2*n)*2^(2*n-2-4*a)/GAMMA(-4*a+2*n)/GAMMA(n+1-a)*GAM
MA(1-2*a)*GAMMA(n-a);}{%
\maplemultiline{
\mathit{\ 79:\ }    - {\displaystyle \frac {1}{2}} \,
{\displaystyle \frac {\Gamma (1 - 2\,a + 2\,n)\,\Gamma (n - a)\,
 \left(  \! {\displaystyle \sum _{k=1}^{2\,n - 2}} \,
{\displaystyle \frac {(-1)^{k}\,\Gamma (2\,a - 2\,n + 1 + k)\,
\Gamma (1 - 2\,a)}{\Gamma (2\,a - 2\,n + 1)\,\Gamma (1 - 2\,a + k
)}}  \!  \right) }{\Gamma (2\,n - 2\,a)\,\Gamma (n + 1 - a)}} 
 \\
\mbox{} + {\displaystyle \frac {\Gamma (1 - 2\,a + 2\,n)\,2^{(2\,
n - 2 - 4\,a)}\,\Gamma (1 - 2\,a)\,\Gamma (n - a)}{\Gamma ( - 4\,
a + 2\,n)\,\Gamma (n + 1 - a)}}  }
}
\end{maplelatex}

\begin{maplelatex}
\mapleinline{inert}{2d}{` 90: `,
(2/3*a+1/3)/GAMMA(1/2+n)/GAMMA(3/2-3*n+a)*GAMMA(-3*n+5/2)*GAMMA(a-n+1/
2)/GAMMA(2*a+2)/GAMMA(3/2-n)*GAMMA(2*a+1/2+n);}{%
\[
\mathit{\ 90:\ }   \,{\displaystyle \frac {({\displaystyle \frac {
2\,a}{3}}  + {\displaystyle \frac {1}{3}} )\,\Gamma ( - 3\,n + 
{\displaystyle \frac {5}{2}} )\,\Gamma (a - n + {\displaystyle 
\frac {1}{2}} )\,\Gamma (2\,a + {\displaystyle \frac {1}{2}}  + n
)}{\Gamma ({\displaystyle \frac {1}{2}}  + n)\,\Gamma (
{\displaystyle \frac {3}{2}}  - 3\,n + a)\,\Gamma (2\,a + 2)\,
\Gamma ({\displaystyle \frac {3}{2}}  - n)}} 
\]
}
\end{maplelatex}

\begin{maplelatex}
\mapleinline{inert}{2d}{` 93: `,
1/GAMMA(c+1)/GAMMA(b+1)/GAMMA(1+a)*GAMMA(1/2*c+1+1/2*b+1/2*a-1/2*(a^2-
2*a*c-2*a*b+c^2-2*b*c+b^2)^(1/2))*GAMMA(1/2*c+1+1/2*b+1/2*a+1/2*(a^2-2
*a*c-2*a*b+c^2-2*b*c+b^2)^(1/2));}{%
\maplemultiline{
\mathit{\ 93:\ }   \Gamma ({\displaystyle \frac {c}{2}}  + 1 + 
{\displaystyle \frac {b}{2}}  + {\displaystyle \frac {a}{2}}  - 
{\displaystyle \frac {\sqrt{a^{2} - 2\,a\,c - 2\,a\,b + c^{2} - 2
\,b\,c + b^{2}}}{2}} ) \\
\Gamma ({\displaystyle \frac {c}{2}}  + 1 + {\displaystyle 
\frac {b}{2}}  + {\displaystyle \frac {a}{2}}  + {\displaystyle 
\frac {\sqrt{a^{2} - 2\,a\,c - 2\,a\,b + c^{2} - 2\,b\,c + b^{2}}
}{2}} )/(\Gamma (c + 1)\,\Gamma (b + 1)
\Gamma (1 + a)) }
}
\end{maplelatex}

\begin{maplelatex}
\mapleinline{inert}{2d}{` 94: `,
2*3^(6*a-5)*2^(-6*a+5)/(2*a-2/3)*GAMMA(2*a-1)*GAMMA(2*a-1/6)*GAMMA(1/6
+2*a)/Pi^(1/2)/GAMMA(3*a-1)/GAMMA(-3/2+3*a);}{%
\[
\mathit{\ 94:\ }   \,{\displaystyle \frac {2\,3^{(6\,a - 5)}\,2^{(
 - 6\,a + 5)}\,\Gamma (2\,a - 1)\,\Gamma (2\,a - {\displaystyle 
\frac {1}{6}} )\,\Gamma ({\displaystyle \frac {1}{6}}  + 2\,a)}{(
2\,a - {\displaystyle \frac {2}{3}} )\,\sqrt{\pi }\,\Gamma (3\,a
 - 1)\,\Gamma ( - {\displaystyle \frac {3}{2}}  + 3\,a)}} 
\]
}
\end{maplelatex}

\begin{maplelatex}
\mapleinline{inert}{2d}{` 95: `,
2*3^(6*a-7)*2^(-6*a+7)/(2*a-4/3)*GAMMA(2*a-2)*GAMMA(2*a-5/6)*GAMMA(2*a
-1/6)/Pi^(1/2)/GAMMA(3*a-2)/GAMMA(-5/2+3*a);}{%
\[
\mathit{\ 95:\ }   \,{\displaystyle \frac {2\,3^{(6\,a - 7)}\,2^{(
 - 6\,a + 7)}\,\Gamma (2\,a - 2)\,\Gamma (2\,a - {\displaystyle 
\frac {5}{6}} )\,\Gamma (2\,a - {\displaystyle \frac {1}{6}} )}{(
2\,a - {\displaystyle \frac {4}{3}} )\,\sqrt{\pi }\,\Gamma (3\,a
 - 2)\,\Gamma ( - {\displaystyle \frac {5}{2}}  + 3\,a)}} 
\]
}
\end{maplelatex}

\begin{maplelatex}
\mapleinline{inert}{2d}{` 96: `,
2*3^(6*a-5)*2^(-6*a+5)/(2*a-2/3)*GAMMA(1/6+a)*GAMMA(2*a-1/6)/Pi^(1/2)/
GAMMA(-3/2+3*a);}{%
\[
\mathit{\ 96:\ }   \,{\displaystyle \frac {2\,3^{(6\,a - 5)}\,2^{(
 - 6\,a + 5)}\,\Gamma ({\displaystyle \frac {1}{6}}  + a)\,\Gamma
 (2\,a - {\displaystyle \frac {1}{6}} )}{(2\,a - {\displaystyle 
\frac {2}{3}} )\,\sqrt{\pi }\,\Gamma ( - {\displaystyle \frac {3
}{2}}  + 3\,a)}} 
\]
}
\end{maplelatex}

\begin{maplelatex}
\mapleinline{inert}{2d}{` 97: `,
2*3^(6*a-2)*2^(-6*a+2)/(2*a+1/3)*GAMMA(7/6+a)*GAMMA(2*a+5/6)/Pi^(1/2)/
GAMMA(1/2+3*a);}{%
\[
\mathit{\ 97:\ }   \,{\displaystyle \frac {2\,3^{(6\,a - 2)}\,2^{(
 - 6\,a + 2)}\,\Gamma ({\displaystyle \frac {7}{6}}  + a)\,\Gamma
 (2\,a + {\displaystyle \frac {5}{6}} )}{(2\,a + {\displaystyle 
\frac {1}{3}} )\,\sqrt{\pi }\,\Gamma ({\displaystyle \frac {1}{2}
}  + 3\,a)}} 
\]
}
\end{maplelatex}

\begin{maplelatex}
\mapleinline{inert}{2d}{` 98: `,
2*3^(6*a-9)*2^(-6*a+9)/(2*a-2)*GAMMA(2*a-3)*GAMMA(-7/6+2*a)*GAMMA(2*a-
5/6)/Pi^(1/2)/GAMMA(3*a-3)/GAMMA(-7/2+3*a);}{%
\[
\mathit{\ 98:\ }   \,{\displaystyle \frac {2\,3^{(6\,a - 9)}\,2^{(
 - 6\,a + 9)}\,\Gamma (2\,a - 3)\,\Gamma ( - {\displaystyle 
\frac {7}{6}}  + 2\,a)\,\Gamma (2\,a - {\displaystyle \frac {5}{6
}} )}{(2\,a - 2)\,\sqrt{\pi }\,\Gamma (3\,a - 3)\,\Gamma ( - 
{\displaystyle \frac {7}{2}}  + 3\,a)}} 
\]
}
\end{maplelatex}

\begin{maplelatex}
\mapleinline{inert}{2d}{` 99: `,
2*3^(6*a-5)*2^(-6*a+5)/(2*a-2/3)*GAMMA(a-1/6)*GAMMA(1/6+2*a)/Pi^(1/2)/
GAMMA(-3/2+3*a);}{%
\[
\mathit{\ 99:\ }   \,{\displaystyle \frac {2\,3^{(6\,a - 5)}\,2^{(
 - 6\,a + 5)}\,\Gamma (a - {\displaystyle \frac {1}{6}} )\,\Gamma
 ({\displaystyle \frac {1}{6}}  + 2\,a)}{(2\,a - {\displaystyle 
\frac {2}{3}} )\,\sqrt{\pi }\,\Gamma ( - {\displaystyle \frac {3
}{2}}  + 3\,a)}} 
\]
}
\end{maplelatex}

\begin{maplelatex}
\mapleinline{inert}{2d}{` 100: `,
2*3^(6*a-2)*2^(-6*a+2)/(2*a+1/3)*GAMMA(a+5/6)*GAMMA(7/6+2*a)/Pi^(1/2)/
GAMMA(1/2+3*a);}{%
\[
\mathit{\ 100:\ }   \,{\displaystyle \frac {2\,3^{(6\,a - 2)}\,2^{
( - 6\,a + 2)}\,\Gamma (a + {\displaystyle \frac {5}{6}} )\,
\Gamma ({\displaystyle \frac {7}{6}}  + 2\,a)}{(2\,a + 
{\displaystyle \frac {1}{3}} )\,\sqrt{\pi }\,\Gamma (
{\displaystyle \frac {1}{2}}  + 3\,a)}} 
\]
}
\end{maplelatex}

\begin{maplelatex}
\mapleinline{inert}{2d}{` 101: `,
2*3^(6*a-7)*2^(-6*a+7)/(2*a-4/3)*GAMMA(a-5/6)*GAMMA(2*a-1/6)/Pi^(1/2)/
GAMMA(-5/2+3*a);}{%
\[
\mathit{\ 101:\ }   \,{\displaystyle \frac {2\,3^{(6\,a - 7)}\,2^{
( - 6\,a + 7)}\,\Gamma (a - {\displaystyle \frac {5}{6}} )\,
\Gamma (2\,a - {\displaystyle \frac {1}{6}} )}{(2\,a - 
{\displaystyle \frac {4}{3}} )\,\sqrt{\pi }\,\Gamma ( - 
{\displaystyle \frac {5}{2}}  + 3\,a)}} 
\]
}
\end{maplelatex}

\begin{maplelatex}
\mapleinline{inert}{2d}{` 102: `,
2*3^(6*a-4)*2^(-6*a+4)/(2*a-1/3)*GAMMA(1/6+a)*GAMMA(2*a+5/6)/Pi^(1/2)/
GAMMA(3*a-1/2);}{%
\[
\mathit{\ 102:\ }   \,{\displaystyle \frac {2\,3^{(6\,a - 4)}\,2^{
( - 6\,a + 4)}\,\Gamma ({\displaystyle \frac {1}{6}}  + a)\,
\Gamma (2\,a + {\displaystyle \frac {5}{6}} )}{(2\,a - 
{\displaystyle \frac {1}{3}} )\,\sqrt{\pi }\,\Gamma (3\,a - 
{\displaystyle \frac {1}{2}} )}} 
\]
}
\end{maplelatex}

\begin{maplelatex}
\mapleinline{inert}{2d}{` 103: `, 3^(3*a)*2^(-3*a)/(1+a);}{%
\[
\mathit{\ 103:\ }   \,{\displaystyle \frac {3^{(3\,a)}\,2^{( - 3\,
a)}}{1 + a}} 
\]
}
\end{maplelatex}

\begin{maplelatex}
\mapleinline{inert}{2d}{` 135: `,
1/4*(sin(Pi*a)+sin(3*Pi*a))*sin(2*Pi*a)*GAMMA(1-2*a+n)^2*GAMMA(1/2+n-a
)*2^(2*n-2*a)*GAMMA(1-2*a)/sin(3*Pi*a)/GAMMA(1-3*a+n)/cos(Pi*a)/Pi^(3/
2)*Sum((-1)^k*GAMMA(2*a-2*n+1+k)/GAMMA(2-2*a+k),k = 0 ..
2*n-2)-1/4*Pi*2^(2*n-2*a)*(-sin(2*Pi*(n-2*a))+2*sin(2*Pi*a))*GAMMA(1-2
*a+n)*GAMMA(1/2+n-a)/sin(Pi*(-3*a+n))/sin(2*Pi*a)/cos(Pi*(n-2*a))/GAMM
A(1-3*a+n)/GAMMA(1/2-2*a+n)/GAMMA(2*a);}{%
\maplemultiline{
\mathit{\ 135:\ }   {\displaystyle \frac {1}{4}} (\mathrm{sin}(\pi
 \,a) + \mathrm{sin}(3\,\pi \,a))\,\mathrm{sin}(2\,\pi \,a)\,
\Gamma (1 - 2\,a + n)^{2}\,\Gamma ({\displaystyle \frac {1}{2}} 
 + n - a)\,2^{(2\,n - 2\,a)} \\
\Gamma (1 - 2\,a)\, \left(  \! {\displaystyle \sum _{k=0}^{2\,n
 - 2}} \,{\displaystyle \frac {(-1)^{k}\,\Gamma (2\,a - 2\,n + 1
 + k)}{\Gamma (2 - 2\,a + k)}}  \!  \right)  \left/ {\vrule 
height0.63em width0em depth0.63em} \right. \!  \! (\mathrm{sin}(3
\,\pi \,a)\,\Gamma (1 - 3\,a + n) 
\mathrm{cos}(\pi \,a)\,\pi ^{(3/2)})\mbox{} \\ -  
{\displaystyle \frac {1}{4}} \,{\displaystyle \frac {\pi \,2^{(2
\,n - 2\,a)}\,( - \mathrm{sin}(2\,\pi \,(n - 2\,a)) + 2\,\mathrm{
sin}(2\,\pi \,a))\,\Gamma (1 - 2\,a + n)\,\Gamma ({\displaystyle 
\frac {1}{2}}  + n - a)}{\mathrm{sin}(\pi \,( - 3\,a + n))\,
\mathrm{sin}(2\,\pi \,a)\,\mathrm{cos}(\pi \,(n - 2\,a))\,\Gamma 
(1 - 3\,a + n)\,\Gamma ({\displaystyle \frac {1}{2}}  - 2\,a + n)
\,\Gamma (2\,a)}}  }
}
\end{maplelatex}

\begin{maplelatex}
\mapleinline{inert}{2d}{` 136: `,
1/4*(cos(3*Pi*a)+cos(Pi*a))*sin(2*Pi*a)*GAMMA(1/2-2*a+n)^2*2^(2*n-2*a)
*GAMMA(n-a)*GAMMA(1-2*a)*(n-2*a)/cos(3*Pi*a)/GAMMA(1/2-3*a+n)/cos(Pi*a
)/Pi^(3/2)*Sum((-1)^k*GAMMA(2*a-2*n+1+k)/GAMMA(2-2*a+k),k = 0 ..
2*n-2)-1/4*Pi*2^(2*n-2*a)*(-sin(4*Pi*a)+2*sin(2*Pi*a))*GAMMA(n-a)*GAMM
A(1/2-2*a+n)/cos(3*Pi*a)/sin(2*Pi*a)^2/GAMMA(1/2-3*a+n)/GAMMA(2*a)/GAM
MA(n-2*a);}{%
\maplemultiline{
\mathit{\ 136:\ }   {\displaystyle \frac {1}{4}} (\mathrm{cos}(3\,
\pi \,a) + \mathrm{cos}(\pi \,a))\,\mathrm{sin}(2\,\pi \,a)\,
\Gamma ({\displaystyle \frac {1}{2}}  - 2\,a + n)^{2}\,2^{(2\,n
 - 2\,a)}\,\Gamma (n - a)\,\Gamma (1 - 2\,a) \\
(n - 2\,a)\, \left(  \! {\displaystyle \sum _{k=0}^{2\,n - 2}} \,
{\displaystyle \frac {(-1)^{k}\,\Gamma (2\,a - 2\,n + 1 + k)}{
\Gamma (2 - 2\,a + k)}}  \!  \right)  \left/ {\vrule 
height0.84em width0em depth0.84em} \right. \!  \! (\mathrm{cos}(3
\,\pi \,a)\,\Gamma ({\displaystyle \frac {1}{2}}  - 3\,a + n) 
\mathrm{cos}(\pi \,a)\,\pi ^{(3/2)}) \\
\mbox{} - {\displaystyle \frac {1}{4}} \,{\displaystyle \frac {
\pi \,2^{(2\,n - 2\,a)}\,( - \mathrm{sin}(4\,\pi \,a) + 2\,
\mathrm{sin}(2\,\pi \,a))\,\Gamma (n - a)\,\Gamma (
{\displaystyle \frac {1}{2}}  - 2\,a + n)}{\mathrm{cos}(3\,\pi \,
a)\,\mathrm{sin}(2\,\pi \,a)^{2}\,\Gamma ({\displaystyle \frac {1
}{2}}  - 3\,a + n)\,\Gamma (2\,a)\,\Gamma (n - 2\,a)}}  }
}
\end{maplelatex}

\begin{maplelatex}
\mapleinline{inert}{2d}{` 137: `,
1/2*(sin(Pi*a)-sin(3*Pi*a))*cos(Pi*a)*GAMMA(1/2-2*a+n)^2*GAMMA(1/2+n-a
)*2^(2*n-2*a)*GAMMA(1-2*a)*(n-2*a)/sin(3*Pi*a)/GAMMA(1-3*a+n)/Pi^(3/2)
*Sum((-1)^k*GAMMA(2*a-2*n+1+k)/GAMMA(2-2*a+k),k = 0 ..
2*n-2)-1/8*Pi*2^(2*n-2*a)*(cos(3*Pi*a)+cos(5*Pi*a)-2*cos(Pi*a))*GAMMA(
1/2-2*a+n)*GAMMA(1/2+n-a)/sin(3*Pi*a)/(-1+cos(2*Pi*a))/cos(Pi*a)/sin(2
*Pi*a)/GAMMA(1-3*a+n)/GAMMA(2*a)/GAMMA(n-2*a);}{%
\maplemultiline{
\mathit{\ 137:\ }   {\displaystyle \frac {1}{2}} (\mathrm{sin}(\pi
 \,a) - \mathrm{sin}(3\,\pi \,a))\,\mathrm{cos}(\pi \,a)\,\Gamma 
({\displaystyle \frac {1}{2}}  - 2\,a + n)^{2}\,\Gamma (
{\displaystyle \frac {1}{2}}  + n - a)\,2^{(2\,n - 2\,a)}\,\Gamma
 (1 - 2\,a) \\
(n - 2\,a)\, \left(  \! {\displaystyle \sum _{k=0}^{2\,n - 2}} \,
{\displaystyle \frac {(-1)^{k}\,\Gamma (2\,a - 2\,n + 1 + k)}{
\Gamma (2 - 2\,a + k)}}  \!  \right)  \left/ {\vrule 
height0.63em width0em depth0.63em} \right. \!  \! (\mathrm{sin}(3
\,\pi \,a)\,\Gamma (1 - 3\,a + n)\,\pi ^{(3/2)}) \\
\mbox{} -  
{\displaystyle \frac {1}{8}} \,{\displaystyle \frac {\pi \,2^{(2
\,n - 2\,a)}\,(\mathrm{cos}(3\,\pi \,a) + \mathrm{cos}(5\,\pi \,a
) - 2\,\mathrm{cos}(\pi \,a))\,\Gamma ({\displaystyle \frac {1}{2
}}  - 2\,a + n)\,\Gamma ({\displaystyle \frac {1}{2}}  + n - a)}{
\mathrm{sin}(3\,\pi \,a)\,( - 1 + \mathrm{cos}(2\,\pi \,a))\,
\mathrm{cos}(\pi \,a)\,\mathrm{sin}(2\,\pi \,a)\,\Gamma (1 - 3\,a
 + n)\,\Gamma (2\,a)\,\Gamma (n - 2\,a)}}  }
}
\end{maplelatex}

\begin{maplelatex}
\mapleinline{inert}{2d}{` 138: `,
1/2*GAMMA(1/2-2*a+n)^2*(n-2*a)/GAMMA(3/2-2*a+n)/GAMMA(-2*a-1/2+n)/(n-a
)*Sum((-1)^k*GAMMA(2*a-2*n+k)*GAMMA(1-2*a)/GAMMA(2*a-2*n)/GAMMA(1-2*a+
k),k = 1 ..
2*n-1)-1/8*Pi*16^(-a)*4^n*(sin(4*Pi*a)-4*sin(Pi*a)*cos(Pi*a))*GAMMA(2*
n-2*a)/(1-4*a+2*n)/sin(Pi*a)/cos(Pi*a)/cos(2*Pi*a)/sin(2*Pi*a)/GAMMA(2
*n-1-4*a)/GAMMA(2*a);}{%
\maplemultiline{
\mathit{\ 138:\ }   {\displaystyle \frac {1}{2}} \,{\displaystyle 
\frac {\Gamma ({\displaystyle \frac {1}{2}}  - 2\,a + n)^{2}\,(n
 - 2\,a)\, \left(  \! {\displaystyle \sum _{k=1}^{2\,n - 1}} \,
{\displaystyle \frac {(-1)^{k}\,\Gamma (2\,a - 2\,n + k)\,\Gamma 
(1 - 2\,a)}{\Gamma (2\,a - 2\,n)\,\Gamma (1 - 2\,a + k)}}  \! 
 \right) }{\Gamma ({\displaystyle \frac {3}{2}}  - 2\,a + n)\,
\Gamma ( - 2\,a - {\displaystyle \frac {1}{2}}  + n)\,(n - a)}} 
 \\
\mbox{} - {\displaystyle \frac {1}{8}} \,{\displaystyle \frac {
\pi \,16^{( - a)}\,4^{n}\,(\mathrm{sin}(4\,\pi \,a) - 4\,\mathrm{
sin}(\pi \,a)\,\mathrm{cos}(\pi \,a))\,\Gamma (2\,n - 2\,a)}{(1
 - 4\,a + 2\,n)\,\mathrm{sin}(\pi \,a)\,\mathrm{cos}(\pi \,a)\,
\mathrm{cos}(2\,\pi \,a)\,\mathrm{sin}(2\,\pi \,a)\,\Gamma (2\,n
 - 1 - 4\,a)\,\Gamma (2\,a)}}  }
}
\end{maplelatex}

\begin{maplelatex}
\mapleinline{inert}{2d}{` 140: `,
-2*2^(-2*a-2*n)/Pi*(1/2*sin(4*Pi*a)-sin(2*Pi*a))/sin(2*Pi*a)/GAMMA(2*a
)*GAMMA(-a-n+1)*GAMMA(3/2-2*a-n)*GAMMA(2*a+n)*GAMMA(3*a-1/2+n)+2^(-2*a
+1-2*n)*Sum((-1)^k*GAMMA(2*a+k)/GAMMA(-2*a-2*n+3+k),k = 0 ..
2*n-2)*(-1+n+2*a)/Pi^(3/2)*((-1)^n*sin(2*Pi*a)*sin(Pi*a)+(-1)^n*cos(3*
Pi*a))/cos(Pi*a)*GAMMA(3*a-1/2+n)/GAMMA(2*a)*GAMMA(3/2-2*a-n)^2*GAMMA(
-a-n+1);}{%
\maplemultiline{
\mathit{\ 140:\ }    - {\displaystyle \frac {2\,2^{( - 2\,a - 2\,n
)}\,({\displaystyle \frac {1}{2}} \,\mathrm{sin}(4\,\pi \,a) - 
\mathrm{sin}(2\,\pi \,a))\,\Gamma ( - a - n + 1)\,\Gamma (
{\displaystyle \frac {3}{2}}  - 2\,a - n)\,\Gamma (2\,a + n)\,
\Gamma (3\,a - {\displaystyle \frac {1}{2}}  + n)}{\pi \,\mathrm{
sin}(2\,\pi \,a)\,\Gamma (2\,a)}}  \\
\mbox{} + 2^{( - 2\,a + 1 - 2\,n)}\, \left(  \! {\displaystyle 
\sum _{k=0}^{2\,n - 2}} \,{\displaystyle \frac {(-1)^{k}\,\Gamma 
(2\,a + k)}{\Gamma ( - 2\,a - 2\,n + 3 + k)}}  \!  \right) \,( - 
1 + n + 2\,a) 
((-1)^{n}\,\mathrm{sin}(2\,\pi \,a)\,\mathrm{sin}(\pi \,a) + (-1)
^{n}\,\mathrm{cos}(3\,\pi \,a)) \\ \,\Gamma (3\,a - {\displaystyle 
\frac {1}{2}}  + n)\,\Gamma ({\displaystyle \frac {3}{2}}  - 2\,a
 - n)^{2}\,\Gamma ( - a - n + 1) 
 \left/ {\vrule height0.51em width0em depth0.51em} \right. \! 
 \! (\pi ^{(3/2)}\,\mathrm{cos}(\pi \,a)\,\Gamma (2\,a)) }
}
\end{maplelatex}

\begin{maplelatex}
\mapleinline{inert}{2d}{` 141: `,
2*(-1+4*a+2*n)/Pi*(1/4*cos(5*Pi*a)-3/4*cos(3*Pi*a)+1/2*cos(Pi*a))*2^(-
2*a-2*n)/(-1+cos(2*Pi*a))/sin(2*Pi*a)/GAMMA(2*a)*sin(Pi*a)*GAMMA(3*a-1
/2+n)*GAMMA(2*a-3/2+n)*GAMMA(-2*a-n+2)*GAMMA(-a-n+1)-2^(-2*a+1-2*n)*Su
m((-1)^k*GAMMA(2*a+k)/GAMMA(-2*a-2*n+3+k),k = 0 ..
2*n-2)*(-1+4*a+2*n)*GAMMA(3*a-1/2+n)*GAMMA(-a-n+1)*GAMMA(-2*a-n+2)^2*(
-1)^n*(cos(2*Pi*a)*cos(Pi*a)-cos(3*Pi*a))/(4*a-3+2*n)/Pi^(3/2)/GAMMA(2
*a)/sin(Pi*a);}{%
\maplemultiline{
\mathit{\ 141:\ }   2\,( - 1 + 4\,a + 2\,n)\,({\displaystyle 
\frac {1}{4}} \,\mathrm{cos}(5\,\pi \,a) - {\displaystyle \frac {
3}{4}} \,\mathrm{cos}(3\,\pi \,a) + {\displaystyle \frac {1}{2}} 
\,\mathrm{cos}(\pi \,a))\,2^{( - 2\,a - 2\,n)}\,\mathrm{sin}(\pi 
\,a) \\
\Gamma (3\,a - {\displaystyle \frac {1}{2}}  + n)\,\Gamma (2\,a
 - {\displaystyle \frac {3}{2}}  + n)\,\Gamma ( - 2\,a - n + 2)\,
\Gamma ( - a - n + 1)/(\pi \,( - 1 + \mathrm{cos}(2\,\pi \,a))
 \\
\mathrm{sin}(2\,\pi \,a)\,\Gamma (2\,a))\mbox{} - 2^{( - 2\,a + 1
 - 2\,n)}\, \left(  \! {\displaystyle \sum _{k=0}^{2\,n - 2}} \,
{\displaystyle \frac {(-1)^{k}\,\Gamma (2\,a + k)}{\Gamma ( - 2\,
a - 2\,n + 3 + k)}}  \!  \right) \,( - 1 + 4\,a + 2\,n) \\
\Gamma (3\,a - {\displaystyle \frac {1}{2}}  + n)\,\Gamma ( - a
 - n + 1)\,\Gamma ( - 2\,a - n + 2)^{2}\,(-1)^{n} \\
(\mathrm{cos}(2\,\pi \,a)\,\mathrm{cos}(\pi \,a) - \mathrm{cos}(3
\,\pi \,a)) \left/ {\vrule height0.63em width0em depth0.63em}
 \right. \!  \! ((4\,a - 3 + 2\,n)\,\pi ^{(3/2)}\,\Gamma (2\,a)\,
\mathrm{sin}(\pi \,a)) }
}
\end{maplelatex}

\begin{maplelatex}
\mapleinline{inert}{2d}{` 142: `,
4*2^(-2*a-2*n)/Pi*cos(Pi*a)^2/GAMMA(2*a)*GAMMA(2*a+n)*GAMMA(3*a+n-1)*G
AMMA(3/2-2*a-n)*GAMMA(-a+3/2-n)-2^(-2*a+1-2*n)*Sum((-1)^k*GAMMA(2*a+k)
/GAMMA(-2*a-2*n+3+k),k = 0 ..
2*n-2)*(-1+n+2*a)*(-1)^n*(2*cos(Pi*a)^2-1)*GAMMA(3*a+n-1)*GAMMA(3/2-2*
a-n)^2*GAMMA(-a+3/2-n)/Pi^(3/2)/GAMMA(2*a);}{%
\maplemultiline{
\mathit{\ 142:\ }   {\displaystyle \frac {4\,2^{( - 2\,a - 2\,n)}
\,\mathrm{cos}(\pi \,a)^{2}\,\Gamma (2\,a + n)\,\Gamma (3\,a + n
 - 1)\,\Gamma ({\displaystyle \frac {3}{2}}  - 2\,a - n)\,\Gamma 
( - a + {\displaystyle \frac {3}{2}}  - n)}{\pi \,\Gamma (2\,a)}
} \\ -  
2^{( - 2\,a + 1 - 2\,n)}\, \left(  \! {\displaystyle \sum _{k=0}
^{2\,n - 2}} \,{\displaystyle \frac {(-1)^{k}\,\Gamma (2\,a + k)
}{\Gamma ( - 2\,a - 2\,n + 3 + k)}}  \!  \right) \,( - 1 + n + 2
\,a)\,(-1)^{n}\,(2\,\mathrm{cos}(\pi \,a)^{2} - 1) \\
\Gamma (3\,a + n - 1)\,\Gamma ({\displaystyle \frac {3}{2}}  - 2
\,a - n)^{2}\,\Gamma ( - a + {\displaystyle \frac {3}{2}}  - n)
 \left/ {\vrule height0.51em width0em depth0.51em} \right. \! 
 \! (\pi ^{(3/2)}\,\Gamma (2\,a)) }
}
\end{maplelatex}

\begin{maplelatex}
\mapleinline{inert}{2d}{` 143: `,
8*16^(-a)*4^(-n)*(-1/8*cos(8*Pi*a)-1/4*cos(6*Pi*a)+1/4*cos(2*Pi*a)+1/8
)/sin(4*Pi*a)/sin(2*Pi*a)/GAMMA(2*a)*GAMMA(-1+4*a+2*n)*GAMMA(-2*a-2*n+
2)+Sum((-1)^k*GAMMA(2*a+k)/GAMMA(-2*a-2*n+3+k),k = 0 ..
2*n-2)*(-1+n+2*a)*GAMMA(-2*a-2*n+2)/GAMMA(2*a);}{%
\maplemultiline{
\mathit{\ 143:\ }   8\,16^{( - a)}\,4^{( - n)}\,( - 
{\displaystyle \frac {1}{8}} \,\mathrm{cos}(8\,\pi \,a) - 
{\displaystyle \frac {1}{4}} \,\mathrm{cos}(6\,\pi \,a) + 
{\displaystyle \frac {1}{4}} \,\mathrm{cos}(2\,\pi \,a) + 
{\displaystyle \frac {1}{8}} )\,\Gamma ( - 1 + 4\,a + 2\,n) \\
\Gamma ( - 2\,a - 2\,n + 2)/(\mathrm{sin}(4\,\pi \,a)\,\mathrm{
sin}(2\,\pi \,a)\,\Gamma (2\,a)) \\
\mbox{} + {\displaystyle \frac { \left(  \! {\displaystyle \sum 
_{k=0}^{2\,n - 2}} \,{\displaystyle \frac {(-1)^{k}\,\Gamma (2\,a
 + k)}{\Gamma ( - 2\,a - 2\,n + 3 + k)}}  \!  \right) \,( - 1 + n
 + 2\,a)\,\Gamma ( - 2\,a - 2\,n + 2)}{\Gamma (2\,a)}}  }
}
\end{maplelatex}

\begin{maplelatex}
\mapleinline{inert}{2d}{` 144: `,
-1/8*1/sin(2*Pi*a)^2/cos(2*Pi*a)*Pi*16^(-a)*2^(2*n)*(-2*sin(2*Pi*a)+si
n(4*Pi*a))/GAMMA(2*n-1-4*a)/GAMMA(2*a)*GAMMA(2*n-2*a-1)-1/2*Sum((-1)^k
*GAMMA(2*a-2*n+1+k)/GAMMA(1-2*a+k),k = 1 ..
2*n-2)*sin(2*Pi*a)*(2*n-1-4*a)/Pi*GAMMA(2*n-2*a-1)*GAMMA(1-2*a);}{%
\maplemultiline{
\mathit{\ 144:\ }   \, - {\displaystyle \frac {1}{8}} \,
{\displaystyle \frac {\pi \,16^{( - a)}\,2^{(2\,n)}\,( - 2\,
\mathrm{sin}(2\,\pi \,a) + \mathrm{sin}(4\,\pi \,a))\,\Gamma (2\,
n - 2\,a - 1)}{\mathrm{sin}(2\,\pi \,a)^{2}\,\mathrm{cos}(2\,\pi 
\,a)\,\Gamma (2\,n - 1 - 4\,a)\,\Gamma (2\,a)}} \\ - 
{\displaystyle \frac {1}{2}} \,{\displaystyle \frac { \left(  \! 
{\displaystyle \sum _{k=1}^{2\,n - 2}} \,{\displaystyle \frac {(
-1)^{k}\,\Gamma (2\,a - 2\,n + 1 + k)}{\Gamma (1 - 2\,a + k)}} 
 \!  \right) \,\mathrm{sin}(2\,\pi \,a)\,(2\,n - 1 - 4\,a)\,
\Gamma (2\,n - 2\,a - 1)\,\Gamma (1 - 2\,a)}{\pi }} }
}
\end{maplelatex}

\begin{maplelatex}
\mapleinline{inert}{2d}{` 145: `,
-Pi/(n-2*a)*16^(-1-a)*2^(2*n)*(2*sin(2*Pi*a)+sin(4*Pi*a))/sin(2*Pi*a)^
2/cos(2*Pi*a)/GAMMA(2*n-2-4*a)/GAMMA(2*a)*GAMMA(2*n-2*a-1)-1/2*Sum((-1
)^k*GAMMA(2*a-2*n+1+k)/GAMMA(1-2*a+k),k = 1 ..
2*n-2)*sin(2*Pi*a)*(-1+n-2*a)*(2*n-1-4*a)/Pi/(n-2*a)*GAMMA(2*n-2*a-1)*
GAMMA(1-2*a);}{%
\maplemultiline{
\mathit{\ 145:\ }    - {\displaystyle \frac {\pi \,16^{( - 1 - a)}
\,2^{(2\,n)}\,(2\,\mathrm{sin}(2\,\pi \,a) + \mathrm{sin}(4\,\pi 
\,a))\,\Gamma (2\,n - 2\,a - 1)}{(n - 2\,a)\,\mathrm{sin}(2\,\pi 
\,a)^{2}\,\mathrm{cos}(2\,\pi \,a)\,\Gamma (2\,n - 2 - 4\,a)\,
\Gamma (2\,a)}} \\ - {\displaystyle \frac {1}{2}}  
 \left(  \! {\displaystyle \sum _{k=1}^{2\,n - 2}} \,
{\displaystyle \frac {(-1)^{k}\,\Gamma (2\,a - 2\,n + 1 + k)}{
\Gamma (1 - 2\,a + k)}}  \!  \right) \,\mathrm{sin}(2\,\pi \,a)\,
( - 1 + n - 2\,a)\,(2\,n - 1 - 4\,a) \\
\Gamma (2\,n - 2\,a - 1)\,\Gamma (1 - 2\,a)/(\pi \,(n - 2\,a)) }
}
\end{maplelatex}

\begin{maplelatex}
\mapleinline{inert}{2d}{` 146: `, -1/2*Sum((-1)^k*GAMMA(2*a-2*n+1+k)/GAMMA(2-2*a+k),k = 0 ..
2*n-2)*4^a/Pi^(3/2)/GAMMA(1-3*a+n)*GAMMA(1/2+n-a)*GAMMA(1-2*a)*GAMMA(1
-4*a+2*n)*sin(4*Pi*a)*cos(Pi*a)*sin(2*Pi*a)/sin(3*Pi*a)/cos(Pi*(-3*a+n
))+1/2*2^(2*n-2*a)/Pi^(1/2)*sin(4*Pi*a)*cos(Pi*a)/sin(3*Pi*a)/cos(Pi*(
-3*a+n))/GAMMA(1-3*a+n)*GAMMA(1/2+n-a)*GAMMA(1-2*a)+sin(Pi*(-a+3/2))*(
1-4*a+2*n)*DixonPP(n,1,2*a-n,1/2+n-a,1/2+2*a-n)/sin(Pi*a)/GAMMA(1-3*a+
n)*GAMMA(1/2+n-a)*GAMMA(1/2-3*a+n)/GAMMA(n-a);}{%
\maplemultiline{
\mathit{\ 146:\ }    - {\displaystyle \frac {1}{2}}  \left(  \! 
{\displaystyle \sum _{k=0}^{2\,n - 2}} \,{\displaystyle \frac {(
-1)^{k}\,\Gamma (2\,a - 2\,n + 1 + k)}{\Gamma (2 - 2\,a + k)}} 
 \!  \right) \,4^{a}\,\Gamma ({\displaystyle \frac {1}{2}}  + n
 - a)\,\Gamma (1 - 2\,a)\,\Gamma (1 - 4\,a + 2\,n) \\
\mathrm{sin}(4\,\pi \,a)\,\mathrm{cos}(\pi \,a)\,\mathrm{sin}(2\,
\pi \,a) \left/ {\vrule height0.63em width0em depth0.63em}
 \right. \!  \! (\pi ^{(3/2)}\,\Gamma (1 - 3\,a + n)\,\mathrm{sin
}(3\,\pi \,a) 
\mathrm{cos}(\pi \,( - 3\,a + n)))\mbox{} \\ + {\displaystyle 
\frac {1}{2}} \,{\displaystyle \frac {2^{(2\,n - 2\,a)}\,\mathrm{
sin}(4\,\pi \,a)\,\mathrm{cos}(\pi \,a)\,\Gamma ({\displaystyle 
\frac {1}{2}}  + n - a)\,\Gamma (1 - 2\,a)}{\sqrt{\pi }\,\mathrm{
sin}(3\,\pi \,a)\,\mathrm{cos}(\pi \,( - 3\,a + n))\,\Gamma (1 - 
3\,a + n)}} \\ +  
\mathrm{sin}(\pi \,( - a + {\displaystyle \frac {3}{2}} ))\,(1 - 
4\,a + 2\,n)\,\mathrm{DixonPP}(n, \,1, \,2\,a - n, \,
{\displaystyle \frac {1}{2}}  + n - a, \,{\displaystyle \frac {1
}{2}}  + 2\,a - n) \\
\Gamma ({\displaystyle \frac {1}{2}}  + n - a)\,\Gamma (
{\displaystyle \frac {1}{2}}  - 3\,a + n) \left/ {\vrule 
height0.37em width0em depth0.37em} \right. \!  \! (\mathrm{sin}(
\pi \,a)\,\Gamma (1 - 3\,a + n)\,\Gamma (n - a)) }
}
\end{maplelatex}

\begin{maplelatex}
\mapleinline{inert}{2d}{` 147: `, 1/2*4^a*Sum((-1)^k*GAMMA(2*a+k)/GAMMA(-2*a-2*n+3+k),k = 0
..
2*n-2)/(4*a-3+2*n)/Pi^(3/2)*GAMMA(3*a-1/2+n)*GAMMA(3-4*a-2*n)*GAMMA(-a
-n+1)/GAMMA(2*a)*sin(4*Pi*a)*sin(Pi*a)*cos(Pi*a)/sin(Pi*(3*a+n))/cos(P
i*(a+n))-2*2^(-2*a-2*n)/Pi^(3/2)/(4*a-3+2*n)*sin(4*Pi*a)*sin(Pi*a)*cos
(Pi*a)/sin(Pi*(3*a+n))/cos(Pi*(a+n))*GAMMA(3*a-1/2+n)*GAMMA(-a-n+1)*GA
MMA(1-2*a)+sin(Pi*a)*sin(Pi*(-3/2+3*a))*DixonMM(n,1,-1+n+2*a,-a-n+1,2*
a-3/2+n)/(4*a-3+2*n)/sin(3*Pi*a)/sin(Pi*(a-3/2))/GAMMA(3*a+n-1)*GAMMA(
3*a-1/2+n)/GAMMA(-a+3/2-n)*GAMMA(-a-n+1);}{%
\maplemultiline{
\mathit{\ 147:\ }   {\displaystyle \frac {1}{2}} 4^{a}\, \left( 
 \! {\displaystyle \sum _{k=0}^{2\,n - 2}} \,{\displaystyle 
\frac {(-1)^{k}\,\Gamma (2\,a + k)}{\Gamma ( - 2\,a - 2\,n + 3 + 
k)}}  \!  \right) \,\Gamma (3\,a - {\displaystyle \frac {1}{2}} 
 + n)\,\Gamma (3 - 4\,a - 2\,n)\,\Gamma ( - a - n + 1) \\
\mathrm{sin}(4\,\pi \,a)\,\mathrm{sin}(\pi \,a)\,\mathrm{cos}(\pi
 \,a) \left/ {\vrule height0.63em width0em depth0.63em}
 \right. \!  \! ((4\,a - 3 + 2\,n)\,\pi ^{(3/2)}\,\Gamma (2\,a)\,
\mathrm{sin}(\pi \,(3\,a + n)) 
\mathrm{cos}(\pi \,(a + n)))\mbox{}  \\ -
{\displaystyle \frac {2\,2^{( - 2\,a - 2\,n)}\,\mathrm{sin}(4\,
\pi \,a)\,\mathrm{sin}(\pi \,a)\,\mathrm{cos}(\pi \,a)\,\Gamma (3
\,a - {\displaystyle \frac {1}{2}}  + n)\,\Gamma ( - a - n + 1)\,
\Gamma (1 - 2\,a)}{\pi ^{(3/2)}\,(4\,a - 3 + 2\,n)\,\mathrm{sin}(
\pi \,(3\,a + n))\,\mathrm{cos}(\pi \,(a + n))}}  \\
\mbox{} + \mathrm{sin}(\pi \,a)\,\mathrm{sin}(\pi \,( - 
{\displaystyle \frac {3}{2}}  + 3\,a))\,\mathrm{DixonMM}(n, \,1, 
\, - 1 + n + 2\,a, \, - a - n + 1, \,2\,a - {\displaystyle 
\frac {3}{2}}  + n) \\
\Gamma (3\,a - {\displaystyle \frac {1}{2}}  + n)\,\Gamma ( - a
 - n + 1) \left/ {\vrule height0.80em width0em depth0.80em}
 \right. \!  \! ((4\,a - 3 + 2\,n)\,\mathrm{sin}(3\,\pi \,a)\,
\mathrm{sin}(\pi \,(a - {\displaystyle \frac {3}{2}} )) 
\Gamma (3\,a + n - 1)\,\Gamma ( - a + {\displaystyle \frac {3}{2}
}  - n)) }
}
\end{maplelatex}

\begin{maplelatex}
\mapleinline{inert}{2d}{` 153: `, 1/2*(n-2*a)*Sum((-1)^k*GAMMA(2*a-2*n+1+k)/GAMMA(2-2*a+k),k
= 0 ..
2*n-2)/(1-4*a+2*n)/Pi*GAMMA(1-2*a+2*n)*GAMMA(1-2*a)*(sin(Pi*a)+sin(3*P
i*a))/cos(Pi*a)-1/4*Pi^(1/2)/GAMMA(3/2-2*a+n)*GAMMA(1-2*a+2*n)/GAMMA(2
*a)*GAMMA(2*a+1-n)*(-1)^n*(2*sin(2*Pi*a)+sin(4*Pi*a))/cos(2*Pi*a)/sin(
2*Pi*a);}{%
\maplemultiline{
\mathit{\ 153:\ }   {\displaystyle \frac {1}{2}} (n - 2\,a)\,
 \left(  \! {\displaystyle \sum _{k=0}^{2\,n - 2}} \,
{\displaystyle \frac {(-1)^{k}\,\Gamma (2\,a - 2\,n + 1 + k)}{
\Gamma (2 - 2\,a + k)}}  \!  \right) \,\Gamma (1 - 2\,a + 2\,n)\,
\Gamma (1 - 2\,a) \\
(\mathrm{sin}(\pi \,a) + \mathrm{sin}(3\,\pi \,a))/((1 - 4\,a + 2
\,n)\,\pi \,\mathrm{cos}(\pi \,a)) \\
\mbox{} - {\displaystyle \frac {1}{4}} \,{\displaystyle \frac {
\sqrt{\pi }\,\Gamma (1 - 2\,a + 2\,n)\,\Gamma (2\,a + 1 - n)\,(-1
)^{n}\,(2\,\mathrm{sin}(2\,\pi \,a) + \mathrm{sin}(4\,\pi \,a))}{
\Gamma ({\displaystyle \frac {3}{2}}  - 2\,a + n)\,\Gamma (2\,a)
\,\mathrm{cos}(2\,\pi \,a)\,\mathrm{sin}(2\,\pi \,a)}}  }
}
\end{maplelatex}

\begin{maplelatex}
\mapleinline{inert}{2d}{` 154: `, 1/8*Sum((-1)^k*GAMMA(2*a-2*n+1+k)/GAMMA(2-2*a+k),k = 0 ..
2*n-2)*16^a*4^(-n)/Pi*GAMMA(1-2*a+2*n)*GAMMA(1-2*a)*(sin(Pi*a)+sin(3*P
i*a))/cos(Pi*a)+1/4*GAMMA(1-2*a+2*n)/GAMMA(2*a)*GAMMA(4*a-2*n)*(2*sin(
2*Pi*a)+sin(4*Pi*a))/sin(2*Pi*a);}{%
\maplemultiline{
\mathit{\ 154:\ }   {\displaystyle \frac {1}{8}}  \left(  \! 
{\displaystyle \sum _{k=0}^{2\,n - 2}} \,{\displaystyle \frac {(
-1)^{k}\,\Gamma (2\,a - 2\,n + 1 + k)}{\Gamma (2 - 2\,a + k)}} 
 \!  \right) \,16^{a}\,4^{( - n)}\,\Gamma (1 - 2\,a + 2\,n)\,
\Gamma (1 - 2\,a) \\
(\mathrm{sin}(\pi \,a) + \mathrm{sin}(3\,\pi \,a))/(\pi \,
\mathrm{cos}(\pi \,a)) \\
\mbox{} + {\displaystyle \frac {1}{4}} \,{\displaystyle \frac {
\Gamma (1 - 2\,a + 2\,n)\,\Gamma (4\,a - 2\,n)\,(2\,\mathrm{sin}(
2\,\pi \,a) + \mathrm{sin}(4\,\pi \,a))}{\Gamma (2\,a)\,\mathrm{
sin}(2\,\pi \,a)}}  }
}
\end{maplelatex}

\begin{maplelatex}
\mapleinline{inert}{2d}{` 155: `, 1/4*Sum((-1)^k*GAMMA(2*a-2*n+1+k)/GAMMA(2-2*a+k),k = 0 ..
2*n-2)/Pi*GAMMA(1-2*a+2*n)*GAMMA(1-2*a)*(sin(Pi*a)+sin(3*Pi*a))/cos(Pi
*a)+1/4*Pi^(1/2)*GAMMA(1-2*a+2*n)/GAMMA(2*a)/GAMMA(1/2-2*a+n)*GAMMA(2*
a-n)*(-1)^n*(2*sin(2*Pi*a)+sin(4*Pi*a))/cos(2*Pi*a)/sin(2*Pi*a);}{%
\maplemultiline{
\mathit{\ 155:\ }   {\displaystyle \frac {1}{4}} \,{\displaystyle 
\frac { \left(  \! {\displaystyle \sum _{k=0}^{2\,n - 2}} \,
{\displaystyle \frac {(-1)^{k}\,\Gamma (2\,a - 2\,n + 1 + k)}{
\Gamma (2 - 2\,a + k)}}  \!  \right) \,\Gamma (1 - 2\,a + 2\,n)\,
\Gamma (1 - 2\,a)\,(\mathrm{sin}(\pi \,a) + \mathrm{sin}(3\,\pi 
\,a))}{\pi \,\mathrm{cos}(\pi \,a)}}  \\
\mbox{} + {\displaystyle \frac {1}{4}} \,{\displaystyle \frac {
\sqrt{\pi }\,\Gamma (1 - 2\,a + 2\,n)\,\Gamma (2\,a - n)\,(-1)^{n
}\,(2\,\mathrm{sin}(2\,\pi \,a) + \mathrm{sin}(4\,\pi \,a))}{
\Gamma (2\,a)\,\Gamma ({\displaystyle \frac {1}{2}}  - 2\,a + n)
\,\mathrm{cos}(2\,\pi \,a)\,\mathrm{sin}(2\,\pi \,a)}}  }
}
\end{maplelatex}

\begin{maplelatex}
\mapleinline{inert}{2d}{` 156: `,
2^(-2*a-2+2*n)*Sum((-1)^k*GAMMA(2*a-2*n+1+k)/GAMMA(2-2*a+k),k = 0 ..
2*n-2)/Pi^(1/2)/GAMMA(1-3*a+n)/GAMMA(1/2+2*a-n)*GAMMA(1-2*a)*GAMMA(n+1
-a)*GAMMA(1-2*a+n)*sin(2*Pi*a)*(sin(Pi*a)+sin(3*Pi*a))/cos(Pi*(n-a))/c
os(Pi*a)/sin(3*Pi*a)+2^(-2*a-2+2*n)*Pi/GAMMA(1-3*a+n)/GAMMA(2*a)*GAMMA
(n+1-a)*(2*sin(2*Pi*a)+sin(4*Pi*a))/cos(Pi*a)/sin(2*Pi*a)/sin(3*Pi*a);
}{%
\maplemultiline{
\mathit{\ 156:\ }   2^{( - 2\,a - 2 + 2\,n)}\, \left(  \! 
{\displaystyle \sum _{k=0}^{2\,n - 2}} \,{\displaystyle \frac {(
-1)^{k}\,\Gamma (2\,a - 2\,n + 1 + k)}{\Gamma (2 - 2\,a + k)}} 
 \!  \right) \,\Gamma (1 - 2\,a)\,\Gamma (n + 1 - a)\,\Gamma (1
 - 2\,a + n) \\
\mathrm{sin}(2\,\pi \,a)\,(\mathrm{sin}(\pi \,a) + \mathrm{sin}(3
\,\pi \,a)) \left/ {\vrule height0.80em width0em depth0.80em}
 \right. \!  \! (\sqrt{\pi }\,\Gamma (1 - 3\,a + n)\,\Gamma (
{\displaystyle \frac {1}{2}}  + 2\,a - n)\,\mathrm{cos}(\pi \,(n
 - a)) \\
\mathrm{cos}(\pi \,a)\,\mathrm{sin}(3\,\pi \,a))\mbox{} + 
{\displaystyle \frac {2^{( - 2\,a - 2 + 2\,n)}\,\pi \,\Gamma (n
 + 1 - a)\,(2\,\mathrm{sin}(2\,\pi \,a) + \mathrm{sin}(4\,\pi \,a
))}{\Gamma (1 - 3\,a + n)\,\Gamma (2\,a)\,\mathrm{cos}(\pi \,a)\,
\mathrm{sin}(2\,\pi \,a)\,\mathrm{sin}(3\,\pi \,a)}}  }
}
\end{maplelatex}

\begin{maplelatex}
\mapleinline{inert}{2d}{` 157: `,
2^(-2*a-2+2*n)*Sum((-1)^k*GAMMA(2*a-2*n+1+k)/GAMMA(2-2*a+k),k = 0 ..
2*n-2)/Pi^2*GAMMA(a+1/2)*GAMMA(1/2+n-a)*GAMMA(1-2*a)*GAMMA(1-2*a+n)*(-
1)^n*sin(2*Pi*a)*(sin(Pi*a)+sin(3*Pi*a))/cos(Pi*a)+1/2*1/GAMMA(a)/GAMM
A(1/2-2*a+n)/cos(2*Pi*a)/sin(2*Pi*a)*(2*sin(2*Pi*a)+sin(4*Pi*a))*4^n*1
6^(-a)*Pi*GAMMA(1/2+n-a)*(-1)^n;}{%
\maplemultiline{
\mathit{\ 157:\ }   2^{( - 2\,a - 2 + 2\,n)}\, \left(  \! 
{\displaystyle \sum _{k=0}^{2\,n - 2}} \,{\displaystyle \frac {(
-1)^{k}\,\Gamma (2\,a - 2\,n + 1 + k)}{\Gamma (2 - 2\,a + k)}} 
 \!  \right) \,\Gamma (a + {\displaystyle \frac {1}{2}} )\,\Gamma
 ({\displaystyle \frac {1}{2}}  + n - a)\,\Gamma (1 - 2\,a)\,
\Gamma (1 - 2\,a + n) \\
(-1)^{n}\,\mathrm{sin}(2\,\pi \,a)\,(\mathrm{sin}(\pi \,a) + 
\mathrm{sin}(3\,\pi \,a)) \left/ {\vrule 
height0.43em width0em depth0.43em} \right. \!  \! (\pi ^{2}\,
\mathrm{cos}(\pi \,a)) 
\mbox{} + {\displaystyle \frac {1}{2}} \,{\displaystyle \frac {(2
\,\mathrm{sin}(2\,\pi \,a) + \mathrm{sin}(4\,\pi \,a))\,4^{n}\,16
^{( - a)}\,\pi \,\Gamma ({\displaystyle \frac {1}{2}}  + n - a)\,
(-1)^{n}}{\Gamma (a)\,\Gamma ({\displaystyle \frac {1}{2}}  - 2\,
a + n)\,\mathrm{cos}(2\,\pi \,a)\,\mathrm{sin}(2\,\pi \,a)}}  }
}
\end{maplelatex}

\begin{maplelatex}
\mapleinline{inert}{2d}{` 158: `, 1/2*Sum((-1)^k*GAMMA(2*a-2*n+1+k)/GAMMA(2-2*a+k),k = 0 ..
2*n-2)*(n-2*a)*(1-2*a+2*n)/Pi/(1-4*a+2*n)*GAMMA(2*n-2*a)*GAMMA(1-2*a)*
cos(Pi*(-3*a+n))*sin(2*Pi*a)*(cos(3*Pi*a)+cos(Pi*a))/cos(Pi*(n-2*a))/c
os(Pi*a)/cos(3*Pi*a)-2^(2*n-2)*Pi*16^(-a)*(1-2*a+2*n)/(1-4*a+2*n)*GAMM
A(2*n-2*a)/GAMMA(-4*a+2*n)/GAMMA(2*a)*(-sin(4*Pi*a)+2*sin(2*Pi*a))/cos
(2*Pi*a)/sin(2*Pi*a)^2;}{%
\maplemultiline{
\mathit{\ 158:\ }   {\displaystyle \frac {1}{2}}  \left(  \! 
{\displaystyle \sum _{k=0}^{2\,n - 2}} \,{\displaystyle \frac {(
-1)^{k}\,\Gamma (2\,a - 2\,n + 1 + k)}{\Gamma (2 - 2\,a + k)}} 
 \!  \right) \,(n - 2\,a)\,(1 - 2\,a + 2\,n)\,\Gamma (2\,n - 2\,a
) \\
\Gamma (1 - 2\,a)\,\mathrm{cos}(\pi \,( - 3\,a + n))\,\mathrm{sin
}(2\,\pi \,a)\,(\mathrm{cos}(3\,\pi \,a) + \mathrm{cos}(\pi \,a))
/(\pi \,(1 - 4\,a + 2\,n) \\
\mathrm{cos}(\pi \,(n - 2\,a))\,\mathrm{cos}(\pi \,a)\,\mathrm{
cos}(3\,\pi \,a)) \\
\mbox{} - {\displaystyle \frac {2^{(2\,n - 2)}\,\pi \,16^{( - a)}
\,(1 - 2\,a + 2\,n)\,\Gamma (2\,n - 2\,a)\,( - \mathrm{sin}(4\,
\pi \,a) + 2\,\mathrm{sin}(2\,\pi \,a))}{(1 - 4\,a + 2\,n)\,
\Gamma ( - 4\,a + 2\,n)\,\Gamma (2\,a)\,\mathrm{cos}(2\,\pi \,a)
\,\mathrm{sin}(2\,\pi \,a)^{2}}}  }
}
\end{maplelatex}

\begin{maplelatex}
\mapleinline{inert}{2d}{` 159: `,
2^(-2*a-2+2*n)/(1-4*a+2*n)*Sum((-1)^k*GAMMA(2*a-2*n+1+k)/GAMMA(2-2*a+k
),k = 0 ..
2*n-2)*(n-2*a)/Pi^2*GAMMA(1/2-2*a+n)*GAMMA(a+1/2)*GAMMA(1-2*a)*GAMMA(n
-a)*(-1)^n*sin(2*Pi*a)*(cos(3*Pi*a)+cos(Pi*a))/cos(Pi*a)-1/2*16^(-a)*4
^n*Pi/(1-4*a+2*n)/GAMMA(a)/GAMMA(n-2*a)*GAMMA(n-a)*(-1)^n*(-sin(4*Pi*a
)+2*sin(2*Pi*a))/sin(2*Pi*a)^2;}{%
\maplemultiline{
\mathit{\ 159:\ }   2^{( - 2\,a - 2 + 2\,n)}\, \left(  \! 
{\displaystyle \sum _{k=0}^{2\,n - 2}} \,{\displaystyle \frac {(
-1)^{k}\,\Gamma (2\,a - 2\,n + 1 + k)}{\Gamma (2 - 2\,a + k)}} 
 \!  \right) \,(n - 2\,a)\,\Gamma ({\displaystyle \frac {1}{2}} 
 - 2\,a + n)\,\Gamma (a + {\displaystyle \frac {1}{2}} )\,\Gamma 
(1 - 2\,a) \\
\Gamma (n - a)\,(-1)^{n}\,\mathrm{sin}(2\,\pi \,a)\,(\mathrm{cos}
(3\,\pi \,a) + \mathrm{cos}(\pi \,a)) \left/ {\vrule 
height0.43em width0em depth0.43em} \right. \!  \! ((1 - 4\,a + 2
\,n)\,\pi ^{2}\,\mathrm{cos}(\pi \,a)) \\
\mbox{} - {\displaystyle \frac {1}{2}} \,{\displaystyle \frac {16
^{( - a)}\,4^{n}\,\pi \,\Gamma (n - a)\,(-1)^{n}\,( - \mathrm{sin
}(4\,\pi \,a) + 2\,\mathrm{sin}(2\,\pi \,a))}{(1 - 4\,a + 2\,n)\,
\Gamma (a)\,\Gamma (n - 2\,a)\,\mathrm{sin}(2\,\pi \,a)^{2}}}  }
}
\end{maplelatex}

\begin{maplelatex}
\mapleinline{inert}{2d}{` 169: `, -1/2*Sum((-1)^k*GAMMA(2*a+k)/GAMMA(-2*a-2*n+3+k),k = 0 ..
2*n-2)*(-1+4*a+2*n)*(2*a-3+2*n)/(4*a-3+2*n)*GAMMA(-2*a-2*n+2)/GAMMA(2*
a)*(sin(2*Pi*a)*sin(Pi*a)+cos(3*Pi*a))/cos(2*Pi*a)/cos(Pi*a)+2*16^(-a)
*4^(-n)*(-1+4*a+2*n)*(2*a-3+2*n)*(-sin(4*Pi*a)+2*sin(2*Pi*a))/sin(2*Pi
*a)*GAMMA(-2*a-2*n+2)/GAMMA(2*a)*GAMMA(4*a-3+2*n);}{%
\maplemultiline{
\mathit{\ 169:\ }    - {\displaystyle \frac {1}{2}}  \left(  \! 
{\displaystyle \sum _{k=0}^{2\,n - 2}} \,{\displaystyle \frac {(
-1)^{k}\,\Gamma (2\,a + k)}{\Gamma ( - 2\,a - 2\,n + 3 + k)}} 
 \!  \right) \,( - 1 + 4\,a + 2\,n)\,(2\,a - 3 + 2\,n)\,\Gamma (
 - 2\,a - 2\,n + 2) \\
(\mathrm{sin}(2\,\pi \,a)\,\mathrm{sin}(\pi \,a) + \mathrm{cos}(3
\,\pi \,a))/((4\,a - 3 + 2\,n)\,\Gamma (2\,a)\,\mathrm{cos}(2\,
\pi \,a)\,\mathrm{cos}(\pi \,a))\mbox{} + 2 \\
16^{( - a)}\,4^{( - n)}\,( - 1 + 4\,a + 2\,n)\,(2\,a - 3 + 2\,n)
\,( - \mathrm{sin}(4\,\pi \,a) + 2\,\mathrm{sin}(2\,\pi \,a)) \\
\Gamma ( - 2\,a - 2\,n + 2)\,\Gamma (4\,a - 3 + 2\,n)/(\mathrm{
sin}(2\,\pi \,a)\,\Gamma (2\,a)) }
}
\end{maplelatex}

\begin{maplelatex}
\mapleinline{inert}{2d}{` 170: `, -2^(4*a-4)*Sum((-1)^k*GAMMA(2*a+k)/GAMMA(-2*a-2*n+3+k),k =
0 ..
2*n-2)*4^n*(2*a-3+2*n)/(4*a-3+2*n)*GAMMA(-2*a-2*n+2)/GAMMA(2*a)*(sin(2
*Pi*a)*sin(Pi*a)+cos(3*Pi*a))/cos(2*Pi*a)/cos(Pi*a)+1/4*(2*a-3+2*n)*(-
sin(4*Pi*a)+2*sin(2*Pi*a))/sin(2*Pi*a)*GAMMA(4*a-3+2*n)*GAMMA(-2*a-2*n
+2)/GAMMA(2*a);}{%
\maplemultiline{
\mathit{\ 170:\ }    - 2^{(4\,a - 4)}\, \left(  \! {\displaystyle 
\sum _{k=0}^{2\,n - 2}} \,{\displaystyle \frac {(-1)^{k}\,\Gamma 
(2\,a + k)}{\Gamma ( - 2\,a - 2\,n + 3 + k)}}  \!  \right) \,4^{n
}\,(2\,a - 3 + 2\,n)\,\Gamma ( - 2\,a - 2\,n + 2) \\
(\mathrm{sin}(2\,\pi \,a)\,\mathrm{sin}(\pi \,a) + \mathrm{cos}(3
\,\pi \,a))/((4\,a - 3 + 2\,n)\,\Gamma (2\,a)\,\mathrm{cos}(2\,
\pi \,a)\,\mathrm{cos}(\pi \,a)) \\
\mbox{} + {\displaystyle \frac {1}{4}} \,{\displaystyle \frac {(2
\,a - 3 + 2\,n)\,( - \mathrm{sin}(4\,\pi \,a) + 2\,\mathrm{sin}(2
\,\pi \,a))\,\Gamma (4\,a - 3 + 2\,n)\,\Gamma ( - 2\,a - 2\,n + 2
)}{\mathrm{sin}(2\,\pi \,a)\,\Gamma (2\,a)}}  }
}
\end{maplelatex}

\begin{maplelatex}
\mapleinline{inert}{2d}{` 171: `,
-2^(-2*a-2*n+2)*Sum((-1)^k*GAMMA(2*a+k)/GAMMA(-2*a-2*n+3+k),k = 0 ..
2*n-2)/(4*a-3+2*n)/Pi^(1/2)/GAMMA(-1+n+2*a)*GAMMA(3*a-1/2+n)*GAMMA(3/2
-2*a-n)/GAMMA(2*a)*GAMMA(5/2-a-n)*(sin(2*Pi*a)*sin(Pi*a)+cos(3*Pi*a))/
sin(Pi*a)/cos(Pi*a)-2^(-2*a+1-2*n)/(4*a-3+2*n)*GAMMA(3*a-1/2+n)/GAMMA(
2*a)*GAMMA(5/2-a-n)*(-sin(4*Pi*a)+2*sin(2*Pi*a))/sin(Pi*(a+n))/sin(2*P
i*a);}{%
\maplemultiline{
\mathit{\ 171:\ }    - 2^{( - 2\,a - 2\,n + 2)}\, \left(  \! 
{\displaystyle \sum _{k=0}^{2\,n - 2}} \,{\displaystyle \frac {(
-1)^{k}\,\Gamma (2\,a + k)}{\Gamma ( - 2\,a - 2\,n + 3 + k)}} 
 \!  \right) \,\Gamma (3\,a - {\displaystyle \frac {1}{2}}  + n)
\,\Gamma ({\displaystyle \frac {3}{2}}  - 2\,a - n)\,\Gamma (
{\displaystyle \frac {5}{2}}  - a - n) \\
(\mathrm{sin}(2\,\pi \,a)\,\mathrm{sin}(\pi \,a) + \mathrm{cos}(3
\,\pi \,a)) \left/ {\vrule height0.41em width0em depth0.41em}
 \right. \!  \! ((4\,a - 3 + 2\,n)\,\sqrt{\pi }\,\Gamma ( - 1 + n
 + 2\,a)\,\Gamma (2\,a)\,\mathrm{sin}(\pi \,a) \\
\mathrm{cos}(\pi \,a))\mbox{} - {\displaystyle \frac {2^{( - 2\,a
 + 1 - 2\,n)}\,\Gamma (3\,a - {\displaystyle \frac {1}{2}}  + n)
\,\Gamma ({\displaystyle \frac {5}{2}}  - a - n)\,( - \mathrm{sin
}(4\,\pi \,a) + 2\,\mathrm{sin}(2\,\pi \,a))}{(4\,a - 3 + 2\,n)\,
\Gamma (2\,a)\,\mathrm{sin}(\pi \,(a + n))\,\mathrm{sin}(2\,\pi 
\,a)}}  }
}
\end{maplelatex}

\begin{maplelatex}
\mapleinline{inert}{2d}{` 172: `,
-1/(4*a-3+2*n)/GAMMA(2*a)/cos(Pi*a)/cos(2*Pi*a)*Sum((-1)^k*GAMMA(2*a+k
)/GAMMA(-2*a-2*n+3+k),k = 0 ..
2*n-2)*(-1+n+2*a)*(2*a-3+2*n)*GAMMA(-2*a-2*n+2)*(sin(2*Pi*a)*sin(Pi*a)
+cos(3*Pi*a))+2*(-sin(4*Pi*a)+2*sin(2*Pi*a))*16^(-a)*4^(-n)*(2*a-3+2*n
)/sin(2*Pi*a)/(4*a-3+2*n)*GAMMA(-1+4*a+2*n)*GAMMA(-2*a-2*n+2)/GAMMA(2*
a);}{%
\maplemultiline{
\mathit{\ 172:\ }    -  \left(  \! {\displaystyle \sum _{k=0}^{2\,
n - 2}} \,{\displaystyle \frac {(-1)^{k}\,\Gamma (2\,a + k)}{
\Gamma ( - 2\,a - 2\,n + 3 + k)}}  \!  \right) \,( - 1 + n + 2\,a
)\,(2\,a - 3 + 2\,n)\,\Gamma ( - 2\,a - 2\,n + 2) \\
(\mathrm{sin}(2\,\pi \,a)\,\mathrm{sin}(\pi \,a) + \mathrm{cos}(3
\,\pi \,a))/((4\,a - 3 + 2\,n)\,\Gamma (2\,a)\,\mathrm{cos}(\pi 
\,a)\,\mathrm{cos}(2\,\pi \,a))\mbox{} \\ + 2 \,
( - \mathrm{sin}(4\,\pi \,a) + 2\,\mathrm{sin}(2\,\pi \,a))\,16^{
( - a)}\,4^{( - n)}\,(2\,a - 3 + 2\,n)\,\Gamma ( - 1 + 4\,a + 2\,
n) \\
\Gamma ( - 2\,a - 2\,n + 2)/(\mathrm{sin}(2\,\pi \,a)\,(4\,a - 3
 + 2\,n)\,\Gamma (2\,a)) }
}
\end{maplelatex}

\begin{maplelatex}
\mapleinline{inert}{2d}{` 173: `,
2*16^(-a)*(-1+4*a+2*n)*Sum((-1)^k*GAMMA(2*a+k)/GAMMA(-2*a-2*n+3+k),k =
0 ..
2*n-2)/(4*a-3+2*n)*Pi^(1/2)/GAMMA(-1+n+2*a)/GAMMA(a+1/2)*GAMMA(-a-n+1)
*(sin(2*Pi*a)*sin(Pi*a)+cos(3*Pi*a))/cos(Pi*a)/cos(2*Pi*a)*4^(-n)-2*16
^(-a)*4^(-n)*(-sin(4*Pi*a)+2*sin(2*Pi*a))/sin(2*Pi*a)/(4*a-3+2*n)/GAMM
A(a+1/2)*GAMMA(-a-n+1)*GAMMA(2*a+1/2+n);}{%
\maplemultiline{
\mathit{\ 173:\ }   2\,16^{( - a)}\,( - 1 + 4\,a + 2\,n)\, \left( 
 \! {\displaystyle \sum _{k=0}^{2\,n - 2}} \,{\displaystyle 
\frac {(-1)^{k}\,\Gamma (2\,a + k)}{\Gamma ( - 2\,a - 2\,n + 3 + 
k)}}  \!  \right) \,\sqrt{\pi }\,\Gamma ( - a - n + 1) \\
(\mathrm{sin}(2\,\pi \,a)\,\mathrm{sin}(\pi \,a) + \mathrm{cos}(3
\,\pi \,a))\,4^{( - n)} \left/ {\vrule 
height0.80em width0em depth0.80em} \right. \!  \! ((4\,a - 3 + 2
\,n)\,\Gamma ( - 1 + n + 2\,a)\,\Gamma (a + {\displaystyle 
\frac {1}{2}} ) \,
\mathrm{cos}(\pi \,a)\,\mathrm{cos}(2\,\pi \,a)) \\
\mbox{} - {\displaystyle \frac {2\,16^{( - a)}\,4^{( - n)}\,( - 
\mathrm{sin}(4\,\pi \,a) + 2\,\mathrm{sin}(2\,\pi \,a))\,\Gamma (
 - a - n + 1)\,\Gamma (2\,a + {\displaystyle \frac {1}{2}}  + n)
}{\mathrm{sin}(2\,\pi \,a)\,(4\,a - 3 + 2\,n)\,\Gamma (a + 
{\displaystyle \frac {1}{2}} )}}  }
}
\end{maplelatex}

\begin{maplelatex}
\mapleinline{inert}{2d}{` 174: `, -(-1+n+2*a)*Sum((-1)^k*GAMMA(2*a+k)/GAMMA(-2*a-2*n+3+k),k =
0 ..
2*n-2)/(4*a-3+2*n)*GAMMA(3-2*a-2*n)/GAMMA(2*a)*(sin(2*Pi*a)*cos(Pi*a)-
sin(3*Pi*a))/cos(2*Pi*a)/sin(Pi*a)+2^(2-2*n)*16^(-a)*(cos(3*Pi*a)+cos(
5*Pi*a)-2*cos(Pi*a))*sin(Pi*a)/(-1+cos(2*Pi*a))/sin(2*Pi*a)*(-1+n+2*a)
*GAMMA(4*a-3+2*n)*GAMMA(3-2*a-2*n)/GAMMA(2*a);}{%
\maplemultiline{
\mathit{\ 174:\ }    - ( - 1 + n + 2\,a)\, \left(  \! 
{\displaystyle \sum _{k=0}^{2\,n - 2}} \,{\displaystyle \frac {(
-1)^{k}\,\Gamma (2\,a + k)}{\Gamma ( - 2\,a - 2\,n + 3 + k)}} 
 \!  \right) \,\Gamma (3 - 2\,a - 2\,n) \\
(\mathrm{sin}(2\,\pi \,a)\,\mathrm{cos}(\pi \,a) - \mathrm{sin}(3
\,\pi \,a))/((4\,a - 3 + 2\,n)\,\Gamma (2\,a)\,\mathrm{cos}(2\,
\pi \,a)\,\mathrm{sin}(\pi \,a))\mbox{} +  \\
2^{(2 - 2\,n)}\,16^{( - a)}\,(\mathrm{cos}(3\,\pi \,a) + \mathrm{
cos}(5\,\pi \,a) - 2\,\mathrm{cos}(\pi \,a))\,\mathrm{sin}(\pi \,
a)\,( - 1 + n + 2\,a) \\
\Gamma (4\,a - 3 + 2\,n)\,\Gamma (3 - 2\,a - 2\,n)/(( - 1 + 
\mathrm{cos}(2\,\pi \,a))\,\mathrm{sin}(2\,\pi \,a)\,\Gamma (2\,a
)) }
}
\end{maplelatex}

\begin{maplelatex}
\mapleinline{inert}{2d}{` 175: `,
(-1)^n*2^(2-2*n)*16^(-a)*(-1+n+2*a)*(sin(2*Pi*a)*cos(Pi*a)-sin(3*Pi*a)
)*Sum((-1)^k*GAMMA(2*a+k)/GAMMA(-2*a-2*n+3+k),k = 0 ..
2*n-2)/Pi^(1/2)/sin(Pi*a)/GAMMA(a+1/2)*GAMMA(3/2-2*a-n)*GAMMA(-a+3/2-n
)+2*16^(-a)*4^(-n)*(cos(3*Pi*a)+cos(5*Pi*a)-2*cos(Pi*a))*sin(Pi*a)/(-1
+cos(2*Pi*a))/sin(2*Pi*a)/GAMMA(a+1/2)*GAMMA(-a+3/2-n)*GAMMA(2*a+n);}{
\maplemultiline{
\mathit{\ 175:\ }   (-1)^{n}\,2^{(2 - 2\,n)}\,16^{( - a)}\,( - 1
 + n + 2\,a)\,(\mathrm{sin}(2\,\pi \,a)\,\mathrm{cos}(\pi \,a) - 
\mathrm{sin}(3\,\pi \,a)) \\
 \left(  \! {\displaystyle \sum _{k=0}^{2\,n - 2}} \,
{\displaystyle \frac {(-1)^{k}\,\Gamma (2\,a + k)}{\Gamma ( - 2\,
a - 2\,n + 3 + k)}}  \!  \right) \,\Gamma ({\displaystyle \frac {
3}{2}}  - 2\,a - n)\,\Gamma ( - a + {\displaystyle \frac {3}{2}} 
 - n) \left/ {\vrule height0.80em width0em depth0.80em}
 \right. \!  \! (\sqrt{\pi }\,\mathrm{sin}(\pi \,a) \, 
\Gamma (a + {\displaystyle \frac {1}{2}} ))\mbox{} \\ + 
{\displaystyle \frac {2\,16^{( - a)}\,4^{( - n)}\,(\mathrm{cos}(3
\,\pi \,a) + \mathrm{cos}(5\,\pi \,a) - 2\,\mathrm{cos}(\pi \,a))
\,\mathrm{sin}(\pi \,a)\,\Gamma ( - a + {\displaystyle \frac {3}{
2}}  - n)\,\Gamma (2\,a + n)}{( - 1 + \mathrm{cos}(2\,\pi \,a))\,
\mathrm{sin}(2\,\pi \,a)\,\Gamma (a + {\displaystyle \frac {1}{2}
} )}}  }
}
\end{maplelatex}

\begin{maplelatex}
\mapleinline{inert}{2d}{` 176: `, 1/2*Sum((-1)^k*GAMMA(2*a-2*n+1+k)/GAMMA(1-2*a+k),k = 1 ..
2*n-2)*sin(2*Pi*a)*16^a*4^(-n)/Pi*GAMMA(2*n-2*a)*GAMMA(1-2*a)-1/4*Pi*G
AMMA(2*n-2*a)/GAMMA(2*a)/GAMMA(-4*a+2*n)*sin(Pi*a)/cos(Pi*a)/cos(2*Pi*
a);}{%
\maplemultiline{
\mathit{\ 176:\ }   {\displaystyle \frac {1}{2}} \,{\displaystyle 
\frac { \left(  \! {\displaystyle \sum _{k=1}^{2\,n - 2}} \,
{\displaystyle \frac {(-1)^{k}\,\Gamma (2\,a - 2\,n + 1 + k)}{
\Gamma (1 - 2\,a + k)}}  \!  \right) \,\mathrm{sin}(2\,\pi \,a)\,
16^{a}\,4^{( - n)}\,\Gamma (2\,n - 2\,a)\,\Gamma (1 - 2\,a)}{\pi 
}}  \\
\mbox{} - {\displaystyle \frac {1}{4}} \,{\displaystyle \frac {
\pi \,\Gamma (2\,n - 2\,a)\,\mathrm{sin}(\pi \,a)}{\Gamma (2\,a)
\,\Gamma ( - 4\,a + 2\,n)\,\mathrm{cos}(\pi \,a)\,\mathrm{cos}(2
\,\pi \,a)}}  }
}
\end{maplelatex}

\begin{maplelatex}
\mapleinline{inert}{2d}{` 177: `, 2^(2*a-2)*Sum((-1)^k*GAMMA(2*a-2*n+1+k)/GAMMA(1-2*a+k),k =
1 ..
2*n-2)/Pi/GAMMA(1+4*a-2*n)*GAMMA(1/2+n-a)*GAMMA(1-2*a)*GAMMA(3*a+1/2-n
)*sin(2*Pi*a)/sin(Pi*(n-a))-1/4*2^(2*n-2*a-1)*GAMMA(1/2+n-a)/GAMMA(2*a
)*GAMMA(3*a+1/2-n)*sin(4*Pi*a)/cos(Pi*a)/((-1)^n)/cos(2*Pi*a);}{%
\maplemultiline{
\mathit{\ 177:\ }   {\displaystyle \frac {2^{(2\,a - 2)}\, \left( 
 \! {\displaystyle \sum _{k=1}^{2\,n - 2}} \,{\displaystyle 
\frac {(-1)^{k}\,\Gamma (2\,a - 2\,n + 1 + k)}{\Gamma (1 - 2\,a
 + k)}}  \!  \right) \,\Gamma ({\displaystyle \frac {1}{2}}  + n
 - a)\,\Gamma (1 - 2\,a)\,\Gamma (3\,a + {\displaystyle \frac {1
}{2}}  - n)\,\mathrm{sin}(2\,\pi \,a)}{\pi \,\Gamma (1 + 4\,a - 2
\,n)\,\mathrm{sin}(\pi \,(n - a))}}  \\
\mbox{} - {\displaystyle \frac {1}{4}} \,{\displaystyle \frac {2
^{(2\,n - 2\,a - 1)}\,\Gamma ({\displaystyle \frac {1}{2}}  + n
 - a)\,\Gamma (3\,a + {\displaystyle \frac {1}{2}}  - n)\,
\mathrm{sin}(4\,\pi \,a)}{\Gamma (2\,a)\,\mathrm{cos}(\pi \,a)\,(
-1)^{n}\,\mathrm{cos}(2\,\pi \,a)}}  }
}
\end{maplelatex}

\begin{maplelatex}
\mapleinline{inert}{2d}{` 178: `, 1/4*Sum((-1)^k*GAMMA(2*a-2*n+1+k)/GAMMA(1-2*a+k),k = 1 ..
2*n-2)*sin(2*Pi*a)*(2*n-1-4*a)/Pi/(n-2*a)*GAMMA(2*n-2*a)*GAMMA(1-2*a)-
1/8*Pi/(n-2*a)*GAMMA(2*n-2*a)/GAMMA(2*n-1-4*a)/GAMMA(2*a)*16^(-a)*4^n*
sin(Pi*a)/cos(Pi*a)/cos(2*Pi*a);}{%
\maplemultiline{
\mathit{\ 178:\ }   {\displaystyle \frac {1}{4}} \,{\displaystyle 
\frac { \left(  \! {\displaystyle \sum _{k=1}^{2\,n - 2}} \,
{\displaystyle \frac {(-1)^{k}\,\Gamma (2\,a - 2\,n + 1 + k)}{
\Gamma (1 - 2\,a + k)}}  \!  \right) \,\mathrm{sin}(2\,\pi \,a)\,
(2\,n - 1 - 4\,a)\,\Gamma (2\,n - 2\,a)\,\Gamma (1 - 2\,a)}{\pi 
\,(n - 2\,a)}}  \\
\mbox{} - {\displaystyle \frac {1}{8}} \,{\displaystyle \frac {
\pi \,\Gamma (2\,n - 2\,a)\,16^{( - a)}\,4^{n}\,\mathrm{sin}(\pi 
\,a)}{(n - 2\,a)\,\Gamma (2\,n - 1 - 4\,a)\,\Gamma (2\,a)\,
\mathrm{cos}(\pi \,a)\,\mathrm{cos}(2\,\pi \,a)}}  }
}
\end{maplelatex}

\begin{maplelatex}
\mapleinline{inert}{2d}{` 179: `,
1/2*sin(2*Pi*a)*Sum((-1)^k*GAMMA(2*a-2*n+1+k)/GAMMA(1-2*a+k),k = 1 ..
2*n-2)/Pi*GAMMA(2*n-2*a)*GAMMA(1-2*a)-1/4*Pi*GAMMA(2*n-2*a)/GAMMA(-4*a
+2*n)/GAMMA(2*a)*16^(-a)*4^n*sin(Pi*a)/cos(Pi*a)/cos(2*Pi*a);}{%
\maplemultiline{
\mathit{\ 179:\ }   {\displaystyle \frac {1}{2}} \,{\displaystyle 
\frac {\mathrm{sin}(2\,\pi \,a)\, \left(  \! {\displaystyle \sum 
_{k=1}^{2\,n - 2}} \,{\displaystyle \frac {(-1)^{k}\,\Gamma (2\,a
 - 2\,n + 1 + k)}{\Gamma (1 - 2\,a + k)}}  \!  \right) \,\Gamma (
2\,n - 2\,a)\,\Gamma (1 - 2\,a)}{\pi }}  \\
\mbox{} - {\displaystyle \frac {1}{4}} \,{\displaystyle \frac {
\pi \,\Gamma (2\,n - 2\,a)\,16^{( - a)}\,4^{n}\,\mathrm{sin}(\pi 
\,a)}{\Gamma ( - 4\,a + 2\,n)\,\Gamma (2\,a)\,\mathrm{cos}(\pi \,
a)\,\mathrm{cos}(2\,\pi \,a)}}  }
}
\end{maplelatex}

\begin{maplelatex}
\mapleinline{inert}{2d}{` 180: `,
2^(2*a-2)*sin(2*Pi*a)*Sum((-1)^k*GAMMA(2*a-2*n+1+k)/GAMMA(1-2*a+k),k =
1 ..
2*n-2)/Pi^(3/2)*GAMMA(1-2*a)*GAMMA(n-a)*GAMMA(3*a+1/2-n)/GAMMA(1+4*a-2
*n)+1/4*2^(2*n-2*a-1)*GAMMA(n-a)*GAMMA(3*a+1/2-n)/Pi^(1/2)/GAMMA(2*a)*
sin(4*Pi*a)*sin(Pi*a)/cos(Pi*a)/cos(2*Pi*a);}{%
\maplemultiline{
\mathit{\ 180:\ }   {\displaystyle \frac {2^{(2\,a - 2)}\,\mathrm{
sin}(2\,\pi \,a)\, \left(  \! {\displaystyle \sum _{k=1}^{2\,n - 
2}} \,{\displaystyle \frac {(-1)^{k}\,\Gamma (2\,a - 2\,n + 1 + k
)}{\Gamma (1 - 2\,a + k)}}  \!  \right) \,\Gamma (1 - 2\,a)\,
\Gamma (n - a)\,\Gamma (3\,a + {\displaystyle \frac {1}{2}}  - n)
}{\pi ^{(3/2)}\,\Gamma (1 + 4\,a - 2\,n)}}  \\
\mbox{} + {\displaystyle \frac {1}{4}} \,{\displaystyle \frac {2
^{(2\,n - 2\,a - 1)}\,\Gamma (n - a)\,\Gamma (3\,a + 
{\displaystyle \frac {1}{2}}  - n)\,\mathrm{sin}(4\,\pi \,a)\,
\mathrm{sin}(\pi \,a)}{\sqrt{\pi }\,\Gamma (2\,a)\,\mathrm{cos}(
\pi \,a)\,\mathrm{cos}(2\,\pi \,a)}}  }
}
\end{maplelatex}

\begin{maplelatex}
\mapleinline{inert}{2d}{` 181: `, -1/8*Sum((-1)^k*GAMMA(2*a-2*n+1+k)/GAMMA(1-2*a+k),k = 1 ..
2*n-2)*sin(2*Pi*a)*(2*n-1-4*a)/Pi*2^(2*n-2*a)*GAMMA(n-a)/GAMMA(2*a+1-n
)*GAMMA(1-2*a)*GAMMA(a)+1/8*Pi^(3/2)*GAMMA(n-a)*(16^(-a))^2*16^n*sin(P
i*a)/cos(Pi*a)/GAMMA(a+1/2)/GAMMA(2*n-1-4*a)/GAMMA(2*a+1-n)/cos(2*Pi*a
);}{%
\maplemultiline{
\mathit{\ 181:\ }    - {\displaystyle \frac {1}{8}}  \left(  \! 
{\displaystyle \sum _{k=1}^{2\,n - 2}} \,{\displaystyle \frac {(
-1)^{k}\,\Gamma (2\,a - 2\,n + 1 + k)}{\Gamma (1 - 2\,a + k)}} 
 \!  \right) \,\mathrm{sin}(2\,\pi \,a)\,(2\,n - 1 - 4\,a)\,2^{(2
\,n - 2\,a)}\,\Gamma (n - a) \\
\Gamma (1 - 2\,a)\,\Gamma (a)/(\pi \,\Gamma (2\,a + 1 - n)) \\
\mbox{} + {\displaystyle \frac {1}{8}} \,{\displaystyle \frac {
\pi ^{(3/2)}\,\Gamma (n - a)\,(16^{( - a)})^{2}\,16^{n}\,\mathrm{
sin}(\pi \,a)}{\mathrm{cos}(\pi \,a)\,\Gamma (a + {\displaystyle 
\frac {1}{2}} )\,\Gamma (2\,n - 1 - 4\,a)\,\Gamma (2\,a + 1 - n)
\,\mathrm{cos}(2\,\pi \,a)}}  }
}
\end{maplelatex}

\begin{maplelatex}
\mapleinline{inert}{2d}{` 182: `,
2^(-2*a-2+2*n)*Sum((-1)^k*GAMMA(2*a-2*n+1+k)/GAMMA(1-2*a+k),k = 1 ..
2*n-2)*sin(2*Pi*a)/Pi/GAMMA(1/2+2*a-n)*GAMMA(n-a)*GAMMA(1-2*a)*GAMMA(a
+1/2)+1/8*Pi^(3/2)*GAMMA(n-a)*256^(-a)*16^n*sin(Pi*a)/cos(Pi*a)/GAMMA(
a)/GAMMA(2*n-1-4*a)/GAMMA(2*a+3/2-n)/cos(2*Pi*a);}{%
\maplemultiline{
\mathit{\ 182:\ }   {\displaystyle \frac {2^{( - 2\,a - 2 + 2\,n)}
\, \left(  \! {\displaystyle \sum _{k=1}^{2\,n - 2}} \,
{\displaystyle \frac {(-1)^{k}\,\Gamma (2\,a - 2\,n + 1 + k)}{
\Gamma (1 - 2\,a + k)}}  \!  \right) \,\mathrm{sin}(2\,\pi \,a)\,
\Gamma (n - a)\,\Gamma (1 - 2\,a)\,\Gamma (a + {\displaystyle 
\frac {1}{2}} )}{\pi \,\Gamma ({\displaystyle \frac {1}{2}}  + 2
\,a - n)}}  \\
\mbox{} + {\displaystyle \frac {1}{8}} \,{\displaystyle \frac {
\pi ^{(3/2)}\,\Gamma (n - a)\,256^{( - a)}\,16^{n}\,\mathrm{sin}(
\pi \,a)}{\mathrm{cos}(\pi \,a)\,\Gamma (a)\,\Gamma (2\,n - 1 - 4
\,a)\,\Gamma (2\,a + {\displaystyle \frac {3}{2}}  - n)\,\mathrm{
cos}(2\,\pi \,a)}}  }
}
\end{maplelatex}

\begin{maplelatex}
\mapleinline{inert}{2d}{` 183: `, 1/8*Sum((-1)^k*GAMMA(2*a-2*n+1+k)/GAMMA(1-2*a+k),k = 1 ..
2*n-2)*(2*n-1-4*a)/Pi^(3/2)*2^(2*n-2*a)*GAMMA(n-a)*GAMMA(n-2*a)^2*sin(
2*Pi*a)/GAMMA(2*a)*GAMMA(3*a+1/2-n)*(-1)^n-1/8*GAMMA(n-a)*GAMMA(n-2*a)
^2*GAMMA(3*a+1/2-n)*64^(-a)*16^n/Pi^(1/2)/GAMMA(2*n-1-4*a)/GAMMA(2*a)*
(-1)^n*sin(Pi*a)^2/cos(2*Pi*a);}{%
\maplemultiline{
\mathit{\ 183:\ }   {\displaystyle \frac {1}{8}}  \left(  \! 
{\displaystyle \sum _{k=1}^{2\,n - 2}} \,{\displaystyle \frac {(
-1)^{k}\,\Gamma (2\,a - 2\,n + 1 + k)}{\Gamma (1 - 2\,a + k)}} 
 \!  \right) \,(2\,n - 1 - 4\,a)\,2^{(2\,n - 2\,a)}\,\Gamma (n - 
a)\,\Gamma (n - 2\,a)^{2} 
\mathrm{sin}(2\,\pi \,a) \\ \,\Gamma (3\,a + {\displaystyle \frac {1
}{2}}  - n)\,(-1)^{n} \left/ {\vrule 
height0.63em width0em depth0.63em} \right. \!  \! (\pi ^{(3/2)}\,
\Gamma (2\,a)) 
\mbox{} - {\displaystyle \frac {1}{8}} \,{\displaystyle \frac {
\Gamma (n - a)\,\Gamma (n - 2\,a)^{2}\,\Gamma (3\,a + 
{\displaystyle \frac {1}{2}}  - n)\,64^{( - a)}\,16^{n}\,(-1)^{n}
\,\mathrm{sin}(\pi \,a)^{2}}{\sqrt{\pi }\,\Gamma (2\,n - 1 - 4\,a
)\,\Gamma (2\,a)\,\mathrm{cos}(2\,\pi \,a)}}  }
}
\end{maplelatex}

\begin{maplelatex}
\mapleinline{inert}{2d}{` 184: `,
2^(-2*a-2+2*n)*Sum((-1)^k*GAMMA(2*a-2*n+1+k)/GAMMA(1-2*a+k),k = 1 ..
2*n-2)/Pi^(3/2)*GAMMA(n-a)/GAMMA(2*a)*GAMMA(1/2-2*a+n)^2*GAMMA(3*a+1/2
-n)*cos(Pi*(n-2*a))-1/8*GAMMA(n-a)*GAMMA(1/2-2*a+n)^2*GAMMA(3*a+1/2-n)
*64^(-a)*16^n/Pi^(1/2)/GAMMA(-4*a+2*n)/GAMMA(2*a)*(-1)^n*sin(Pi*a)/cos
(Pi*a);}{%
\maplemultiline{
\mathit{\ 184:\ }   2^{( - 2\,a - 2 + 2\,n)}\, \left(  \! 
{\displaystyle \sum _{k=1}^{2\,n - 2}} \,{\displaystyle \frac {(
-1)^{k}\,\Gamma (2\,a - 2\,n + 1 + k)}{\Gamma (1 - 2\,a + k)}} 
 \!  \right) \,\Gamma (n - a)\,\Gamma ({\displaystyle \frac {1}{2
}}  - 2\,a + n)^{2}\,\Gamma (3\,a + {\displaystyle \frac {1}{2}} 
 - n) \\
\mathrm{cos}(\pi \,(n - 2\,a)) \left/ {\vrule 
height0.51em width0em depth0.51em} \right. \!  \! (\pi ^{(3/2)}\,
\Gamma (2\,a)) \\
\mbox{} - {\displaystyle \frac {1}{8}} \,{\displaystyle \frac {
\Gamma (n - a)\,\Gamma ({\displaystyle \frac {1}{2}}  - 2\,a + n)
^{2}\,\Gamma (3\,a + {\displaystyle \frac {1}{2}}  - n)\,64^{( - 
a)}\,16^{n}\,(-1)^{n}\,\mathrm{sin}(\pi \,a)}{\sqrt{\pi }\,\Gamma
 ( - 4\,a + 2\,n)\,\Gamma (2\,a)\,\mathrm{cos}(\pi \,a)}}  }
}
\end{maplelatex}

\begin{maplelatex}
\mapleinline{inert}{2d}{` 185: `, 2^(2*a-2)*Sum((-1)^k*GAMMA(2*a-2*n+1+k)/GAMMA(1-2*a+k),k =
1 ..
2*n-2)/(n-2*a)*GAMMA(n+1-a)/GAMMA(1+4*a-2*n)*GAMMA(3*a+1-n)/GAMMA(2*a)
/((-1)^n)/cos(Pi*a)-1/4*2^(2*n-2*a-1)/(n-2*a)*GAMMA(3*a+1-n)/GAMMA(2*a
)*GAMMA(n+1-a)*sin(4*Pi*a)/((-1)^n)/cos(2*Pi*a)/sin(Pi*a);}{%
\maplemultiline{
\mathit{\ 185:\ }   {\displaystyle \frac {2^{(2\,a - 2)}\, \left( 
 \! {\displaystyle \sum _{k=1}^{2\,n - 2}} \,{\displaystyle 
\frac {(-1)^{k}\,\Gamma (2\,a - 2\,n + 1 + k)}{\Gamma (1 - 2\,a
 + k)}}  \!  \right) \,\Gamma (n + 1 - a)\,\Gamma (3\,a + 1 - n)
}{(n - 2\,a)\,\Gamma (1 + 4\,a - 2\,n)\,\Gamma (2\,a)\,(-1)^{n}\,
\mathrm{cos}(\pi \,a)}}  \\
\mbox{} - {\displaystyle \frac {1}{4}} \,{\displaystyle \frac {2
^{(2\,n - 2\,a - 1)}\,\Gamma (3\,a + 1 - n)\,\Gamma (n + 1 - a)\,
\mathrm{sin}(4\,\pi \,a)}{(n - 2\,a)\,\Gamma (2\,a)\,(-1)^{n}\,
\mathrm{cos}(2\,\pi \,a)\,\mathrm{sin}(\pi \,a)}}  }
}
\end{maplelatex}

\begin{maplelatex}
\mapleinline{inert}{2d}{` 186: `, -1/2*Sum((-1)^k*GAMMA(2*a-2*n+1+k)/GAMMA(1-2*a+k),k = 1 ..
2*n-2)*(n-a)*16^a*4^(-n)/(n-2*a)/GAMMA(2*a)*GAMMA(2*n-2*a-1)-1/2*(n-a)
*Pi/GAMMA(2*a)/GAMMA(1-4*a+2*n)*GAMMA(2*n-2*a-1)/cos(2*Pi*a)/sin(Pi*a)
*cos(Pi*a);}{%
\maplemultiline{
\mathit{\ 186:\ }    - {\displaystyle \frac {1}{2}} \,
{\displaystyle \frac { \left(  \! {\displaystyle \sum _{k=1}^{2\,
n - 2}} \,{\displaystyle \frac {(-1)^{k}\,\Gamma (2\,a - 2\,n + 1
 + k)}{\Gamma (1 - 2\,a + k)}}  \!  \right) \,(n - a)\,16^{a}\,4
^{( - n)}\,\Gamma (2\,n - 2\,a - 1)}{(n - 2\,a)\,\Gamma (2\,a)}} 
 \\
\mbox{} - {\displaystyle \frac {1}{2}} \,{\displaystyle \frac {(n
 - a)\,\pi \,\Gamma (2\,n - 2\,a - 1)\,\mathrm{cos}(\pi \,a)}{
\Gamma (2\,a)\,\Gamma (1 - 4\,a + 2\,n)\,\mathrm{cos}(2\,\pi \,a)
\,\mathrm{sin}(\pi \,a)}}  }
}
\end{maplelatex}

\begin{maplelatex}
\mapleinline{inert}{2d}{` 187: `, 1/2*Sum((-1)^k*GAMMA(2*a-2*n+1+k)/GAMMA(1-2*a+k),k = 1 ..
2*n-2)*(n-a)*(2*n-1-4*a)/(n-2*a)/GAMMA(2*a)*GAMMA(2*n-2*a-1)+1/4*(n-a)
*Pi/(n-2*a)/GAMMA(2*n-1-4*a)/GAMMA(2*a)*GAMMA(2*n-2*a-1)/cos(2*Pi*a)/s
in(Pi*a)*cos(Pi*a)*16^(-a)*4^n;}{%
\maplemultiline{
\mathit{\ 187:\ }   {\displaystyle \frac {1}{2}} \,{\displaystyle 
\frac { \left(  \! {\displaystyle \sum _{k=1}^{2\,n - 2}} \,
{\displaystyle \frac {(-1)^{k}\,\Gamma (2\,a - 2\,n + 1 + k)}{
\Gamma (1 - 2\,a + k)}}  \!  \right) \,(n - a)\,(2\,n - 1 - 4\,a)
\,\Gamma (2\,n - 2\,a - 1)}{(n - 2\,a)\,\Gamma (2\,a)}}  \\
\mbox{} + {\displaystyle \frac {1}{4}} \,{\displaystyle \frac {(n
 - a)\,\pi \,\Gamma (2\,n - 2\,a - 1)\,\mathrm{cos}(\pi \,a)\,16
^{( - a)}\,4^{n}}{(n - 2\,a)\,\Gamma (2\,n - 1 - 4\,a)\,\Gamma (2
\,a)\,\mathrm{cos}(2\,\pi \,a)\,\mathrm{sin}(\pi \,a)}}  }
}
\end{maplelatex}

\begin{maplelatex}
\mapleinline{inert}{2d}{` 188: `, Sum((-1)^k*GAMMA(2*a-2*n+1+k)/GAMMA(1-2*a+k),k = 1 ..
2*n-2)*(n-a)*(-1+n-2*a)/(n-2*a)/GAMMA(2*a)*GAMMA(2*n-2*a-1)+1/4*Pi*(n-
a)/(n-2*a)/(2*n-1-4*a)/GAMMA(2*a)*GAMMA(2*n-2*a-1)/GAMMA(2*n-2-4*a)/si
n(Pi*a)*cos(Pi*a)/cos(2*Pi*a)*16^(-a)*4^n;}{%
\maplemultiline{
\mathit{\ 188:\ }   {\displaystyle \frac { \left(  \! 
{\displaystyle \sum _{k=1}^{2\,n - 2}} \,{\displaystyle \frac {(
-1)^{k}\,\Gamma (2\,a - 2\,n + 1 + k)}{\Gamma (1 - 2\,a + k)}} 
 \!  \right) \,(n - a)\,( - 1 + n - 2\,a)\,\Gamma (2\,n - 2\,a - 
1)}{(n - 2\,a)\,\Gamma (2\,a)}}  \\
\mbox{} + {\displaystyle \frac {1}{4}} \,{\displaystyle \frac {
\pi \,(n - a)\,\Gamma (2\,n - 2\,a - 1)\,\mathrm{cos}(\pi \,a)\,
16^{( - a)}\,4^{n}}{(n - 2\,a)\,(2\,n - 1 - 4\,a)\,\Gamma (2\,a)
\,\Gamma (2\,n - 2 - 4\,a)\,\mathrm{sin}(\pi \,a)\,\mathrm{cos}(2
\,\pi \,a)}}  }
}
\end{maplelatex}

\begin{maplelatex}
\mapleinline{inert}{2d}{` 189: `, -1/8*Sum((-1)^k*GAMMA(2*a-2*n+1+k)/GAMMA(1-2*a+k),k = 1 ..
2*n-2)/Pi^(1/2)/(n-2*a)*4^a*GAMMA(-a+n-1/2)/GAMMA(1+4*a-2*n)/GAMMA(2*a
)*GAMMA(3*a+1-n)+1/4*2^(-2*a-2+2*n)/(n-2*a)/Pi^(1/2)*GAMMA(3*a+1-n)/GA
MMA(2*a)*GAMMA(-a+n-1/2)*sin(4*Pi*a)/sin(Pi*a)*cos(Pi*a)/cos(2*Pi*a);}
{%
\maplemultiline{
\mathit{\ 189:\ }    - {\displaystyle \frac {1}{8}} \,
{\displaystyle \frac { \left(  \! {\displaystyle \sum _{k=1}^{2\,
n - 2}} \,{\displaystyle \frac {(-1)^{k}\,\Gamma (2\,a - 2\,n + 1
 + k)}{\Gamma (1 - 2\,a + k)}}  \!  \right) \,4^{a}\,\Gamma ( - a
 + n - {\displaystyle \frac {1}{2}} )\,\Gamma (3\,a + 1 - n)}{
\sqrt{\pi }\,(n - 2\,a)\,\Gamma (1 + 4\,a - 2\,n)\,\Gamma (2\,a)}
}  \\
\mbox{} + {\displaystyle \frac {1}{4}} \,{\displaystyle \frac {2
^{( - 2\,a - 2 + 2\,n)}\,\Gamma (3\,a + 1 - n)\,\Gamma ( - a + n
 - {\displaystyle \frac {1}{2}} )\,\mathrm{sin}(4\,\pi \,a)\,
\mathrm{cos}(\pi \,a)}{(n - 2\,a)\,\sqrt{\pi }\,\Gamma (2\,a)\,
\mathrm{sin}(\pi \,a)\,\mathrm{cos}(2\,\pi \,a)}}  }
}
\end{maplelatex}

\begin{maplelatex}
\mapleinline{inert}{2d}{` 190: `,
-1/8*(2*n-1-4*a)*Sum((-1)^k*GAMMA(2*a-2*n+1+k)/GAMMA(1-2*a+k),k = 1 ..
2*n-2)/Pi^(3/2)*2^(2*n-2*a)*GAMMA(n-2*a)^2*GAMMA(-a+n-1/2)*GAMMA(3*a+1
-n)*sin(2*Pi*a)/GAMMA(2*a)*(-1)^n-1/8*16^n*64^(-a)*GAMMA(n-2*a)^2*GAMM
A(-a+n-1/2)*GAMMA(3*a+1-n)/Pi^(1/2)/GAMMA(2*n-1-4*a)/GAMMA(2*a)*(-1)^n
*cos(Pi*a)^2/cos(2*Pi*a);}{%
\maplemultiline{
\mathit{\ 190:\ }    - {\displaystyle \frac {1}{8}} (2\,n - 1 - 4
\,a)\, \left(  \! {\displaystyle \sum _{k=1}^{2\,n - 2}} \,
{\displaystyle \frac {(-1)^{k}\,\Gamma (2\,a - 2\,n + 1 + k)}{
\Gamma (1 - 2\,a + k)}}  \!  \right) \,2^{(2\,n - 2\,a)}\,\Gamma 
(n - 2\,a)^{2} \\
\Gamma ( - a + n - {\displaystyle \frac {1}{2}} )\,\Gamma (3\,a
 + 1 - n)\,\mathrm{sin}(2\,\pi \,a)\,(-1)^{n} \left/ {\vrule 
height0.63em width0em depth0.63em} \right. \!  \! (\pi ^{(3/2)}\,
\Gamma (2\,a)) \\
\mbox{} - {\displaystyle \frac {1}{8}} \,{\displaystyle \frac {16
^{n}\,64^{( - a)}\,\Gamma (n - 2\,a)^{2}\,\Gamma ( - a + n - 
{\displaystyle \frac {1}{2}} )\,\Gamma (3\,a + 1 - n)\,(-1)^{n}\,
\mathrm{cos}(\pi \,a)^{2}}{\sqrt{\pi }\,\Gamma (2\,n - 1 - 4\,a)
\,\Gamma (2\,a)\,\mathrm{cos}(2\,\pi \,a)}}  }
}
\end{maplelatex}

\begin{maplelatex}
\mapleinline{inert}{2d}{` 191: `,
2^(-2*a-2+2*n)*(-1+n-2*a)*Sum((-1)^k*GAMMA(2*a-2*n+1+k)/GAMMA(1-2*a+k)
,k = 1 ..
2*n-2)/(n-2*a)/Pi^(3/2)*GAMMA(-a+n-1/2)*GAMMA(3*a+1-n)/GAMMA(2*a)*GAMM
A(1/2-2*a+n)^2*cos(Pi*(n-2*a))+1/4*2^(4*n-4)*64^(-a)*(2*n-1-4*a)/(n-2*
a)/Pi^(1/2)*GAMMA(3*a+1-n)*GAMMA(-2*a-1/2+n)^2/GAMMA(2*a)/GAMMA(2*n-2-
4*a)*GAMMA(-a+n-1/2)*(-1)^n/sin(Pi*a)*cos(Pi*a);}{%
\maplemultiline{
\mathit{\ 191:\ }   2^{( - 2\,a - 2 + 2\,n)}\,( - 1 + n - 2\,a)\,
 \left(  \! {\displaystyle \sum _{k=1}^{2\,n - 2}} \,
{\displaystyle \frac {(-1)^{k}\,\Gamma (2\,a - 2\,n + 1 + k)}{
\Gamma (1 - 2\,a + k)}}  \!  \right) \,\Gamma ( - a + n - 
{\displaystyle \frac {1}{2}} )\,\Gamma (3\,a + 1 - n) \\
\Gamma ({\displaystyle \frac {1}{2}}  - 2\,a + n)^{2}\,\mathrm{
cos}(\pi \,(n - 2\,a)) \left/ {\vrule 
height0.51em width0em depth0.51em} \right. \!  \! ((n - 2\,a)\,
\pi ^{(3/2)}\,\Gamma (2\,a)) \\
\mbox{} + {\displaystyle \frac {1}{4}} \,{\displaystyle \frac {2
^{(4\,n - 4)}\,64^{( - a)}\,(2\,n - 1 - 4\,a)\,\Gamma (3\,a + 1
 - n)\,\Gamma ( - 2\,a - {\displaystyle \frac {1}{2}}  + n)^{2}\,
\Gamma ( - a + n - {\displaystyle \frac {1}{2}} )\,(-1)^{n}\,
\mathrm{cos}(\pi \,a)}{(n - 2\,a)\,\sqrt{\pi }\,\Gamma (2\,a)\,
\Gamma (2\,n - 2 - 4\,a)\,\mathrm{sin}(\pi \,a)}}  }
}
\end{maplelatex}

\begin{maplelatex}
\mapleinline{inert}{2d}{` 192: `, 1/16*Sum((-1)^k*GAMMA(2*a-2*n+1+k)/GAMMA(1-2*a+k),k = 1 ..
2*n-2)*(2*n-1-4*a)/(n-2*a)*2^(2*n-2*a)/GAMMA(2*a+1-n)*GAMMA(-a+n-1/2)/
GAMMA(2*a)*GAMMA(a+1/2)-1/32*16^(n-2*a)*Pi^(3/2)/(n-2*a)/GAMMA(a)*GAMM
A(-a+n-1/2)/GAMMA(2*n-2-4*a)/GAMMA(2-n+2*a)/cos(2*Pi*a)/sin(Pi*a)*cos(
Pi*a);}{%
\maplemultiline{
\mathit{\ 192:\ }   {\displaystyle \frac {1}{16}} \,
{\displaystyle \frac { \left(  \! {\displaystyle \sum _{k=1}^{2\,
n - 2}} \,{\displaystyle \frac {(-1)^{k}\,\Gamma (2\,a - 2\,n + 1
 + k)}{\Gamma (1 - 2\,a + k)}}  \!  \right) \,(2\,n - 1 - 4\,a)\,
2^{(2\,n - 2\,a)}\,\Gamma ( - a + n - {\displaystyle \frac {1}{2}
} )\,\Gamma (a + {\displaystyle \frac {1}{2}} )}{(n - 2\,a)\,
\Gamma (2\,a + 1 - n)\,\Gamma (2\,a)}}  \\
\mbox{} - {\displaystyle \frac {1}{32}} \,{\displaystyle \frac {
16^{(n - 2\,a)}\,\pi ^{(3/2)}\,\Gamma ( - a + n - {\displaystyle 
\frac {1}{2}} )\,\mathrm{cos}(\pi \,a)}{(n - 2\,a)\,\Gamma (a)\,
\Gamma (2\,n - 2 - 4\,a)\,\Gamma (2 - n + 2\,a)\,\mathrm{cos}(2\,
\pi \,a)\,\mathrm{sin}(\pi \,a)}}  }
}
\end{maplelatex}

\begin{maplelatex}
\mapleinline{inert}{2d}{` 193: `, 1/8*Sum((-1)^k*GAMMA(2*a-2*n+1+k)/GAMMA(1-2*a+k),k = 1 ..
2*n-2)*(-1+n-2*a)/(n-2*a)*2^(2*n-2*a)/GAMMA(1/2+2*a-n)*GAMMA(-a+n-1/2)
/GAMMA(2*a)*GAMMA(a)-1/32*Pi^(3/2)/(n-2*a)*16^(n-2*a)*GAMMA(-a+n-1/2)/
GAMMA(a+1/2)/GAMMA(2*n-2-4*a)/GAMMA(2*a+3/2-n)/cos(2*Pi*a)/sin(Pi*a)*c
os(Pi*a);}{%
\maplemultiline{
\mathit{\ 193:\ }   {\displaystyle \frac {1}{8}} \,{\displaystyle 
\frac { \left(  \! {\displaystyle \sum _{k=1}^{2\,n - 2}} \,
{\displaystyle \frac {(-1)^{k}\,\Gamma (2\,a - 2\,n + 1 + k)}{
\Gamma (1 - 2\,a + k)}}  \!  \right) \,( - 1 + n - 2\,a)\,2^{(2\,
n - 2\,a)}\,\Gamma ( - a + n - {\displaystyle \frac {1}{2}} )\,
\Gamma (a)}{(n - 2\,a)\,\Gamma ({\displaystyle \frac {1}{2}}  + 2
\,a - n)\,\Gamma (2\,a)}}  \\
\mbox{} - {\displaystyle \frac {1}{32}} \,{\displaystyle \frac {
\pi ^{(3/2)}\,16^{(n - 2\,a)}\,\Gamma ( - a + n - {\displaystyle 
\frac {1}{2}} )\,\mathrm{cos}(\pi \,a)}{(n - 2\,a)\,\Gamma (a + 
{\displaystyle \frac {1}{2}} )\,\Gamma (2\,n - 2 - 4\,a)\,\Gamma 
(2\,a + {\displaystyle \frac {3}{2}}  - n)\,\mathrm{cos}(2\,\pi 
\,a)\,\mathrm{sin}(\pi \,a)}}  }
}
\end{maplelatex}

\begin{maplelatex}
\mapleinline{inert}{2d}{` 202: `,
1/2*2^(2*a-1)*(2*n-1-4*a)*Sum((-1)^k*GAMMA(2*a-2*n+1+k)/GAMMA(1-2*a+k)
,k = 1 ..
2*n-2)/(2^(2*a))/(n-2*a)*GAMMA(2*n-2*a-1)/GAMMA(2*a-1)-1/64*(cos(Pi*a)
+1)*(-1+cos(Pi*a))*(2*cos(Pi*a)^2-1)*(2*cos(Pi*a)-1)*(2*cos(Pi*a)+1)*4
^n*2^(4-4*a)*cos(Pi*a)*GAMMA(2*n-2*a-1)*Pi/sin(3*Pi*a)/(cos(2*Pi*a)+1)
/cos(2*Pi*a)^2/GAMMA(2*n-1-4*a)/GAMMA(2*a-1)/(n-2*a);}{%
\maplemultiline{
\mathit{\ 202:\ }   {\displaystyle \frac {1}{2}} \,{\displaystyle 
\frac {2^{(2\,a - 1)}\,(2\,n - 1 - 4\,a)\, \left(  \! 
{\displaystyle \sum _{k=1}^{2\,n - 2}} \,{\displaystyle \frac {(
-1)^{k}\,\Gamma (2\,a - 2\,n + 1 + k)}{\Gamma (1 - 2\,a + k)}} 
 \!  \right) \,\Gamma (2\,n - 2\,a - 1)}{2^{(2\,a)}\,(n - 2\,a)\,
\Gamma (2\,a - 1)}} \\ - {\displaystyle \frac {1}{64}} \,
(\mathrm{cos}(\pi \,a) + 1)\,( - 1 + \mathrm{cos}(\pi \,a))\,(2\,
\mathrm{cos}(\pi \,a)^{2} - 1)\,(2\,\mathrm{cos}(\pi \,a) - 1)\,(
2\,\mathrm{cos}(\pi \,a) + 1)\,4^{n}\,2^{(4 - 4\,a)} \\
\mathrm{cos}(\pi \,a)\,\Gamma (2\,n - 2\,a - 1)\,\pi  \left/ 
{\vrule height0.44em width0em depth0.44em} \right. \!  \! (
\mathrm{sin}(3\,\pi \,a)\,(\mathrm{cos}(2\,\pi \,a) + 1)\,
\mathrm{cos}(2\,\pi \,a)^{2}\,\Gamma (2\,n - 1 - 4\,a) 
\Gamma (2\,a - 1)\,(n - 2\,a)) }
}
\end{maplelatex}

\begin{maplelatex}
\mapleinline{inert}{2d}{` 203: `, 1/2*Sum((-1)^k*GAMMA(2*a-2*n+1+k)/GAMMA(1-2*a+k),k = 1 ..
2*n-2)*(-1+n-2*a)/(n-2*a)/(2*a-1)/GAMMA(2*a-2)*GAMMA(2*n-2*a-1)-1/2*(2
*cos(Pi*a)^2-1)*(2*cos(Pi*a)-1)*(2*cos(Pi*a)+1)*2^(-4*a)*4^n*cos(Pi*a)
*GAMMA(2*n-2*a-1)*Pi*(-1+n-2*a)/sin(3*Pi*a)/cos(2*Pi*a)^2/GAMMA(2*a-2)
/GAMMA(1-4*a+2*n)/(2*a-1);}{%
\maplemultiline{
\mathit{\ 203:\ }   {\displaystyle \frac {1}{2}} \,{\displaystyle 
\frac { \left(  \! {\displaystyle \sum _{k=1}^{2\,n - 2}} \,
{\displaystyle \frac {(-1)^{k}\,\Gamma (2\,a - 2\,n + 1 + k)}{
\Gamma (1 - 2\,a + k)}}  \!  \right) \,( - 1 + n - 2\,a)\,\Gamma 
(2\,n - 2\,a - 1)}{(n - 2\,a)\,(2\,a - 1)\,\Gamma (2\,a - 2)}} 
 \\ - {\displaystyle \frac {1}{2}} \,  
(2\,\mathrm{cos}(\pi \,a)^{2} - 1)\,(2\,\mathrm{cos}(\pi \,a) - 1
)\,(2\,\mathrm{cos}(\pi \,a) + 1)\,2^{( - 4\,a)}\,4^{n}\,\mathrm{
cos}(\pi \,a) \\
\Gamma (2\,n - 2\,a - 1)\,\pi \,( - 1 + n - 2\,a) \left/ {\vrule 
height0.44em width0em depth0.44em} \right. \!  \! (\mathrm{sin}(3
\,\pi \,a)\,\mathrm{cos}(2\,\pi \,a)^{2}\,\Gamma (2\,a - 2) 
\Gamma (1 - 4\,a + 2\,n)\,(2\,a - 1)) }
}
\end{maplelatex}

\begin{maplelatex}
\mapleinline{inert}{2d}{` 204: `,
2^(2-2*n+2*a)*sin(1/2*Pi*(-3+4*a))*sin(2*Pi*a)*Sum((-1)^k*GAMMA(2*a-2*
n+1+k)/GAMMA(1-2*a+k),k = 1 ..
2*n-2)/sin(1/2*Pi*(1+6*a))*GAMMA(2*n-2*a-1)/GAMMA(a-1/2)*GAMMA(-4*a+2*
n)/GAMMA(-1/2-3*a+2*n)-2^(-2*a+2)*Pi*GAMMA(2*n-2*a-1)*(-1+cos(Pi*a))*(
cos(Pi*a)+1)/GAMMA(a-1/2)/GAMMA(-1/2-3*a+2*n)/cos(3*Pi*a);}{%
\maplemultiline{
\mathit{\ 204:\ }   2^{(2 - 2\,n + 2\,a)}\,\mathrm{sin}(
{\displaystyle \frac {\pi \,( - 3 + 4\,a)}{2}} )\,\mathrm{sin}(2
\,\pi \,a)\, \left(  \! {\displaystyle \sum _{k=1}^{2\,n - 2}} \,
{\displaystyle \frac {(-1)^{k}\,\Gamma (2\,a - 2\,n + 1 + k)}{
\Gamma (1 - 2\,a + k)}}  \!  \right) \,\Gamma (2\,n - 2\,a - 1)
 \\
\Gamma ( - 4\,a + 2\,n) \left/ {\vrule 
height0.80em width0em depth0.80em} \right. \!  \! (\mathrm{sin}(
{\displaystyle \frac {\pi \,(1 + 6\,a)}{2}} )\,\Gamma (a - 
{\displaystyle \frac {1}{2}} )\,\Gamma ( - {\displaystyle \frac {
1}{2}}  - 3\,a + 2\,n)) \\
\mbox{} - {\displaystyle \frac {2^{( - 2\,a + 2)}\,\pi \,\Gamma (
2\,n - 2\,a - 1)\,( - 1 + \mathrm{cos}(\pi \,a))\,(\mathrm{cos}(
\pi \,a) + 1)}{\Gamma (a - {\displaystyle \frac {1}{2}} )\,\Gamma
 ( - {\displaystyle \frac {1}{2}}  - 3\,a + 2\,n)\,\mathrm{cos}(3
\,\pi \,a)}}  }
}
\end{maplelatex}

\begin{maplelatex}
\mapleinline{inert}{2d}{` 205: `,
8*sin(Pi*(1/2-2*a))*sin(2*Pi*a)*Sum((-1)^k*GAMMA(2*a-2*n+1+k)/GAMMA(1-
2*a+k),k = 1 ..
2*n-2)/(n-2*a)/sin(3*Pi*a)*2^(2*a-2*n-1)/GAMMA(2*n-3*a-1)*GAMMA(2*n-2*
a-1)/GAMMA(a-1)*GAMMA(-4*a+2*n)-2^(-2*a+2)*Pi*GAMMA(2*n-2*a-1)*cos(Pi*
a)^2/(n-2*a)/GAMMA(a-1)/GAMMA(2*n-3*a-1)/sin(3*Pi*a);}{%
\maplemultiline{
\mathit{\ 205:\ }   8\,\mathrm{sin}(\pi \,({\displaystyle \frac {1
}{2}}  - 2\,a))\,\mathrm{sin}(2\,\pi \,a)\, \left(  \! 
{\displaystyle \sum _{k=1}^{2\,n - 2}} \,{\displaystyle \frac {(
-1)^{k}\,\Gamma (2\,a - 2\,n + 1 + k)}{\Gamma (1 - 2\,a + k)}} 
 \!  \right) \,2^{(2\,a - 2\,n - 1)} \\
\Gamma (2\,n - 2\,a - 1)\,\Gamma ( - 4\,a + 2\,n)/((n - 2\,a)\,
\mathrm{sin}(3\,\pi \,a)\,\Gamma (2\,n - 3\,a - 1)\,\Gamma (a - 1
)) \\
\mbox{} - {\displaystyle \frac {2^{( - 2\,a + 2)}\,\pi \,\Gamma (
2\,n - 2\,a - 1)\,\mathrm{cos}(\pi \,a)^{2}}{(n - 2\,a)\,\Gamma (
a - 1)\,\Gamma (2\,n - 3\,a - 1)\,\mathrm{sin}(3\,\pi \,a)}}  }
}
\end{maplelatex}

\begin{maplelatex}
\mapleinline{inert}{2d}{` 206: `, -Sum((-1)^k*GAMMA(2*a-2*n+1+k)/GAMMA(1-2*a+k),k = 1 ..
2*n-2)*sin(1/2*Pi*(3+2*a))*sin(2*Pi*a)/Pi^(1/2)*(-1)^n*(2*n-1-4*a)*4^(
-a)/sin(1/2*Pi*(1+6*a))*GAMMA(2*n-2*a-1)/GAMMA(-1/2-3*a+2*n)/GAMMA(a)*
GAMMA(n-2*a)^2+2*Pi*GAMMA(2*n-2*a-1)*GAMMA(n-2*a)*cos(Pi*a)*(-1)^n*(-2
^(-2*a+2)*cos(Pi*a)^4+6*4^(-a)*cos(Pi*a)^2-2*4^(-a))/GAMMA(a)/GAMMA(-2
*a-1/2+n)/GAMMA(-1/2-3*a+2*n)/cos(3*Pi*a)/cos(2*Pi*a)^2;}{%
\maplemultiline{
\mathit{\ 206:\ }    -  \left(  \! {\displaystyle \sum _{k=1}^{2\,
n - 2}} \,{\displaystyle \frac {(-1)^{k}\,\Gamma (2\,a - 2\,n + 1
 + k)}{\Gamma (1 - 2\,a + k)}}  \!  \right) \,\mathrm{sin}(
{\displaystyle \frac {\pi \,(3 + 2\,a)}{2}} )\,\mathrm{sin}(2\,
\pi \,a)\,(-1)^{n}\,(2\,n - 1 - 4\,a) \\
4^{( - a)}\,\Gamma (2\,n - 2\,a - 1)\,\Gamma (n - 2\,a)^{2}
 \left/ {\vrule height0.80em width0em depth0.80em} \right. \! 
 \! (\sqrt{\pi }\,\mathrm{sin}({\displaystyle \frac {\pi \,(1 + 6
\,a)}{2}} )\,\Gamma ( - {\displaystyle \frac {1}{2}}  - 3\,a + 2
\,n) 
\Gamma (a)) \\ \mbox{} + 2\,\pi \,\Gamma (2\,n - 2\,a - 1)\,\Gamma (n
 - 2\,a)\,\mathrm{cos}(\pi \,a)\,(-1)^{n} 
( - 2^{( - 2\,a + 2)}\,\mathrm{cos}(\pi \,a)^{4} + 6\,4^{( - a)}
\,\mathrm{cos}(\pi \,a)^{2} - 2\,4^{( - a)}) \left/ {\vrule 
height0.56em width0em depth0.56em} \right. \! \\ \! (\Gamma (a)\,
\Gamma ( - 2\,a - {\displaystyle \frac {1}{2}}  + n) 
\Gamma ( - {\displaystyle \frac {1}{2}}  - 3\,a + 2\,n)\,\mathrm{
cos}(3\,\pi \,a)\,\mathrm{cos}(2\,\pi \,a)^{2}) }
}
\end{maplelatex}

\begin{maplelatex}
\mapleinline{inert}{2d}{` 207: `,
1/4*(2*n-1-4*a)*Sum((-1)^k*GAMMA(2*a-2*n+1+k)/GAMMA(1-2*a+k),k = 1 ..
2*n-2)/(2*a-1)/(n-2*a)/GAMMA(2*a-2)*GAMMA(2*n-2*a-1)-1/4*(2*cos(Pi*a)^
2-1)*(4*cos(Pi*a)^2-3)*2^(-4*a)*4^n*cos(Pi*a)^3*GAMMA(2*n-2*a-1)*Pi/co
s(3*Pi*a)/sin(2*Pi*a)/cos(2*Pi*a)^2/GAMMA(2*a-2)/GAMMA(2*n-1-4*a)/(2*a
-1)/(n-2*a);}{%
\maplemultiline{
\mathit{\ 207:\ }   {\displaystyle \frac {1}{4}} \,{\displaystyle 
\frac {(2\,n - 1 - 4\,a)\, \left(  \! {\displaystyle \sum _{k=1}
^{2\,n - 2}} \,{\displaystyle \frac {(-1)^{k}\,\Gamma (2\,a - 2\,
n + 1 + k)}{\Gamma (1 - 2\,a + k)}}  \!  \right) \,\Gamma (2\,n
 - 2\,a - 1)}{(2\,a - 1)\,(n - 2\,a)\,\Gamma (2\,a - 2)}}  -  \\
{\displaystyle \frac {1}{4}} \,{\displaystyle \frac {(2\,\mathrm{
cos}(\pi \,a)^{2} - 1)\,(4\,\mathrm{cos}(\pi \,a)^{2} - 3)\,2^{(
 - 4\,a)}\,4^{n}\,\mathrm{cos}(\pi \,a)^{3}\,\Gamma (2\,n - 2\,a
 - 1)\,\pi }{\mathrm{cos}(3\,\pi \,a)\,\mathrm{sin}(2\,\pi \,a)\,
\mathrm{cos}(2\,\pi \,a)^{2}\,\Gamma (2\,a - 2)\,\Gamma (2\,n - 1
 - 4\,a)\,(2\,a - 1)\,(n - 2\,a)}}  }
}
\end{maplelatex}

\begin{maplelatex}
\mapleinline{inert}{2d}{` 208: `, 1/2*Sum((-1)^k*GAMMA(2*a-2*n+1+k)/GAMMA(1-2*a+k),k = 1 ..
2*n-2)*GAMMA(2*n-2*a-1)/GAMMA(2*a-1)+1/4*(2*cos(Pi*a)^2-1)*(4*cos(Pi*a
)^2-3)*2^(-4*a)*4^n*sin(Pi*a)*GAMMA(2*n-2*a-1)*Pi/cos(3*Pi*a)/cos(2*Pi
*a)^2/GAMMA(-4*a+2*n)/GAMMA(2*a-1);}{%
\maplemultiline{
\mathit{\ 208:\ }   {\displaystyle \frac {1}{2}} \,{\displaystyle 
\frac { \left(  \! {\displaystyle \sum _{k=1}^{2\,n - 2}} \,
{\displaystyle \frac {(-1)^{k}\,\Gamma (2\,a - 2\,n + 1 + k)}{
\Gamma (1 - 2\,a + k)}}  \!  \right) \,\Gamma (2\,n - 2\,a - 1)}{
\Gamma (2\,a - 1)}}  \\
\mbox{} + {\displaystyle \frac {1}{4}} \,{\displaystyle \frac {(2
\,\mathrm{cos}(\pi \,a)^{2} - 1)\,(4\,\mathrm{cos}(\pi \,a)^{2}
 - 3)\,2^{( - 4\,a)}\,4^{n}\,\mathrm{sin}(\pi \,a)\,\Gamma (2\,n
 - 2\,a - 1)\,\pi }{\mathrm{cos}(3\,\pi \,a)\,\mathrm{cos}(2\,\pi
 \,a)^{2}\,\Gamma ( - 4\,a + 2\,n)\,\Gamma (2\,a - 1)}}  }
}
\end{maplelatex}

\begin{maplelatex}
\mapleinline{inert}{2d}{` 209: `,
-(-1)^n*sin(2*Pi*a)*sin(Pi*(1+a))*(2*n-1-4*a)*Sum((-1)^k*GAMMA(2*a-2*n
+1+k)/GAMMA(1-2*a+k),k = 1 ..
2*n-2)*GAMMA(n-2*a)^2*GAMMA(2*n-2*a-1)/sin(3*Pi*a)*4^(-a)/Pi^(1/2)/GAM
MA(2*n-3*a-1)/GAMMA(a+1/2)+4*Pi*GAMMA(2*n-2*a-1)*GAMMA(n-2*a)*cos(Pi*a
)^3*(-1)^n*(-6*4^(-a)*cos(Pi*a)^2+2*4^(-a)+2^(-2*a+2)*cos(Pi*a)^4)/GAM
MA(a+1/2)/GAMMA(-2*a-1/2+n)/GAMMA(2*n-3*a-1)/sin(2*Pi*a)/sin(3*Pi*a)/c
os(2*Pi*a)^2;}{%
\maplemultiline{
\mathit{\ 209:\ }    - (-1)^{n}\,\mathrm{sin}(2\,\pi \,a)\,
\mathrm{sin}(\pi \,(1 + a))\,(2\,n - 1 - 4\,a)\, \left(  \! 
{\displaystyle \sum _{k=1}^{2\,n - 2}} \,{\displaystyle \frac {(
-1)^{k}\,\Gamma (2\,a - 2\,n + 1 + k)}{\Gamma (1 - 2\,a + k)}} 
 \!  \right)  \\
\Gamma (n - 2\,a)^{2}\,\Gamma (2\,n - 2\,a - 1)\,4^{( - a)}
 \left/ {\vrule height0.80em width0em depth0.80em} \right. \! 
 \! (\mathrm{sin}(3\,\pi \,a)\,\sqrt{\pi }\,\Gamma (2\,n - 3\,a
 - 1)\,\Gamma (a + {\displaystyle \frac {1}{2}} )) \\
\mbox{} + 4\,\pi \,\Gamma (2\,n - 2\,a - 1)\,\Gamma (n - 2\,a)\,
\mathrm{cos}(\pi \,a)^{3}\,(-1)^{n} 
( - 6\,4^{( - a)}\,\mathrm{cos}(\pi \,a)^{2} + 2\,4^{( - a)} + 2
^{( - 2\,a + 2)}\,\mathrm{cos}(\pi \,a)^{4}) \left/ {\vrule 
height0.56em width0em depth0.56em} \right. \\ \!  \! (\Gamma (a + 
{\displaystyle \frac {1}{2}} )\,\Gamma ( - 2\,a - {\displaystyle 
\frac {1}{2}}  + n) \,
\Gamma (2\,n - 3\,a - 1)\,\mathrm{sin}(2\,\pi \,a)\,\mathrm{sin}(
3\,\pi \,a)\,\mathrm{cos}(2\,\pi \,a)^{2}) }
}
\end{maplelatex}

\begin{maplelatex}
\mapleinline{inert}{2d}{` 249: `,
-Gamma(a)*Gamma(1+a-c+n-b)*Gamma(1+a-2*b+m)*Gamma(1+b-c)*(Gamma(b-a)*D
ixonMM(m,n,a,b-m,c-n)*Gamma(1+a-c+n)*Gamma(1+a-b+m)*Gamma(-c+1)-DixonP
P(m,n,a,b,c)*Gamma(b)*Gamma(n+1-c)*Gamma(1-b+m)*Gamma(1+a-c))/Gamma(a-
b)/Gamma(1+a-c+n)/Gamma(1+a-b+m)/Gamma(-c+1)/Gamma(b)/Gamma(n+1-c)/Gam
ma(1-b+m)/Gamma(1+a-c);}{%
\maplemultiline{
\mathit{\ 249:\ }    - \Gamma (a)\,\Gamma (1 + a - c + n - b)\,
\Gamma (1 + a - 2\,b + m)\,\Gamma (1 + b - c)(\Gamma (b - a) \\
\mathrm{DixonMM}(m, \,n, \,a, \,b - m, \,c - n)\,\Gamma (1 + a - 
c + n)\,\Gamma (1 + a - b + m)\,\Gamma ( - c + 1) \\
\mbox{} - \mathrm{DixonPP}(m, \,n, \,a, \,b, \,c)\,\Gamma (b)\,
\Gamma (n + 1 - c)\,\Gamma (1 - b + m)\,\Gamma (1 + a - c))/(
\Gamma (a - b) \\
\Gamma (1 + a - c + n)\,\Gamma (1 + a - b + m)\,\Gamma ( - c + 1)
\,\Gamma (b)\,\Gamma (n + 1 - c)\,\Gamma (1 - b + m) 
\Gamma (1 + a - c)) }
}
\end{maplelatex}

\begin{maplelatex}
\mapleinline{inert}{2d}{` 250: `,
-Gamma(a)*Gamma(1+a-c-n-b)*Gamma(1+a-2*b-m)*Gamma(1+b-c)*(Gamma(b-a)*D
ixonPP(m,n,a,b+m,c+n)*Gamma(1+a-c-n)*Gamma(1+a-b-m)*Gamma(-c+1)-DixonM
M(m,n,a,b,c)*Gamma(b)*Gamma(-c-n+1)*Gamma(-b-m+1)*Gamma(1+a-c))/Gamma(
a-b)/Gamma(1+a-c-n)/Gamma(1+a-b-m)/Gamma(-c+1)/Gamma(b)/Gamma(-c-n+1)/
Gamma(-b-m+1)/Gamma(1+a-c);}{%
\maplemultiline{
\mathit{\ 250:\ }    - \Gamma (a)\,\Gamma (1 + a - c - n - b)\,
\Gamma (1 + a - 2\,b - m)\,\Gamma (1 + b - c)(\Gamma (b - a) \\
\mathrm{DixonPP}(m, \,n, \,a, \,b + m, \,c + n)\,\Gamma (1 + a - 
c - n)\,\Gamma (1 + a - b - m)\,\Gamma ( - c + 1) \\
\mbox{} - \mathrm{DixonMM}(m, \,n, \,a, \,b, \,c)\,\Gamma (b)\,
\Gamma ( - c - n + 1)\,\Gamma ( - b - m + 1)\,\Gamma (1 + a - c))
/ \\ (
\Gamma (a - b)\,\Gamma (1 + a - c - n)\,\Gamma (1 + a - b - m)\,
\Gamma ( - c + 1)\,\Gamma (b)\,\Gamma ( - c - n + 1) 
\Gamma ( - b - m + 1)\,\Gamma (1 + a - c)) }
}
\end{maplelatex}

\begin{maplelatex}
\mapleinline{inert}{2d}{` 252: `, 0;}{%
\[
\mathit{\ 252:\ }   \,0
\]
}
\end{maplelatex}

\begin{maplelatex}
\mapleinline{inert}{2d}{` 261: `,
DixonPP(m,n,c-b+1-2*a+n,1-b+m,n+1-a)/GAMMA(c-2*a+n-m)/GAMMA(-b+c-a+n+1
)/GAMMA(c+1-2*a+n)*GAMMA(c)/GAMMA(b-1+c-m)*GAMMA(c+b-2*a+n-m)*GAMMA(c-
a+n-m)*GAMMA(c-b+1-2*a+n)-DixonMM(m,n,c-b+1-2*a+n,-b+1,-a+1)*GAMMA(a-n
)/GAMMA(c-2*a+n-m)/GAMMA(a)/GAMMA(b)/GAMMA(-b+c-a+1)*GAMMA(c)/GAMMA(b-
1+c-m)*GAMMA(m-c+2*a-n)*GAMMA(c+b-2*a+n-m)*GAMMA(c-a+n-m)*GAMMA(c-b+1-
2*a+n)/GAMMA(1-b+m);}{%
\maplemultiline{
\mathit{\ 261:\ }   \mathrm{DixonPP}(m, \,n, \,c - b + 1 - 2\,a + 
n, \,1 - b + m, \,n + 1 - a)\,\Gamma (c) \\
\Gamma (c + b - 2\,a + n - m)\,\Gamma (c - a + n - m)\,\Gamma (c
 - b + 1 - 2\,a + n)/(\Gamma (c - 2\,a + n - m) \\
\Gamma ( - b + c - a + n + 1)\,\Gamma (c + 1 - 2\,a + n)\,\Gamma 
(b - 1 + c - m))\mbox{}  \\ -
\mathrm{DixonMM}(m, \,n, \,c - b + 1 - 2\,a + n, \, - b + 1, \,
 - a + 1)\,\Gamma (a - n)\,\Gamma (c) \\
\Gamma (m - c + 2\,a - n)\,\Gamma (c + b - 2\,a + n - m)\,\Gamma 
(c - a + n - m)\,\Gamma (c - b + 1 - 2\,a + n)/ \\ (
\Gamma (c - 2\,a + n - m)\,\Gamma (a)\,\Gamma (b)\,\Gamma ( - b
 + c - a + 1)\,\Gamma (b - 1 + c - m)\,\Gamma (1 - b + m)) }
}
\end{maplelatex}

\begin{maplelatex}
\mapleinline{inert}{2d}{` 262: `,
GAMMA(-b+1)*GAMMA(-1/2*b+1/2+1/2*n+1/2*a-1/2*m)*GAMMA(1+a-2*c-n)*Dixon
PP(m,n,-b+1,-1/2*a+c-1/2*b+1/2+1/2*m+1/2*n,-c+1)*GAMMA(-1/2*a+c-1/2*b+
1/2+1/2*m+1/2*n)/GAMMA(-1/2*b+1/2+1/2*a-c-1/2*m-1/2*n)/GAMMA(c+1-b+n)/
GAMMA(3/2-1/2*b+1/2*a-c+1/2*m-1/2*n)/GAMMA(c)*GAMMA(2*c+n)*GAMMA(1/2+1
/2*a+1/2*b+1/2*m-1/2*n)/GAMMA(a)/GAMMA(-1/2*a+c+1/2*b+1/2+1/2*m+1/2*n)
-GAMMA(-b+1)*GAMMA(-1/2*b+1/2+1/2*n+1/2*a-1/2*m)*GAMMA(1+a-2*c-n)*GAMM
A(-1/2*a+c+1/2*b-1/2+1/2*m+1/2*n)*DixonMM(m,n,-b+1,-1/2*a+c-1/2*b+1/2-
1/2*m+1/2*n,-c-n+1)/GAMMA(-1/2*b+1/2+1/2*a-c-1/2*m-1/2*n)/GAMMA(1/2+1/
2*a-c+1/2*b+1/2*m-1/2*n)/GAMMA(c+1-b)*GAMMA(2*c+n)*GAMMA(1/2+1/2*a+1/2
*b+1/2*m-1/2*n)/GAMMA(a)/GAMMA(-1/2*a+c+1/2*b+1/2+1/2*m+1/2*n)/GAMMA(c
+n);}{%
\maplemultiline{
\mathit{\ 262:\ }   \Gamma ( - b + 1)\,\Gamma ( - {\displaystyle 
\frac {b}{2}}  + {\displaystyle \frac {1}{2}}  + {\displaystyle 
\frac {n}{2}}  + {\displaystyle \frac {a}{2}}  - {\displaystyle 
\frac {m}{2}} )\,\Gamma (1 + a - 2\,c - n) \\
\mathrm{DixonPP}(m, \,n, \, - b + 1, \, - {\displaystyle \frac {a
}{2}}  + c - {\displaystyle \frac {b}{2}}  + {\displaystyle 
\frac {1}{2}}  + {\displaystyle \frac {m}{2}}  + {\displaystyle 
\frac {n}{2}} , \, - c + 1) \\
\Gamma ( - {\displaystyle \frac {a}{2}}  + c - {\displaystyle 
\frac {b}{2}}  + {\displaystyle \frac {1}{2}}  + {\displaystyle 
\frac {m}{2}}  + {\displaystyle \frac {n}{2}} )\,\Gamma (2\,c + n
)\,\Gamma ({\displaystyle \frac {1}{2}}  + {\displaystyle \frac {
a}{2}}  + {\displaystyle \frac {b}{2}}  + {\displaystyle \frac {m
}{2}}  - {\displaystyle \frac {n}{2}} ) \left/ {\vrule 
height0.80em width0em depth0.80em} \right. \!  \!  \\ (
\Gamma ( - {\displaystyle \frac {b}{2}}  + {\displaystyle \frac {
1}{2}}  + {\displaystyle \frac {a}{2}}  - c - {\displaystyle 
\frac {m}{2}}  - {\displaystyle \frac {n}{2}} )\,\Gamma (c + 1 - 
b + n)\,\Gamma ({\displaystyle \frac {3}{2}}  - {\displaystyle 
\frac {b}{2}}  + {\displaystyle \frac {a}{2}}  - c + 
{\displaystyle \frac {m}{2}}  - {\displaystyle \frac {n}{2}} )\,
\Gamma (c)\,\Gamma (a) \\
\Gamma ( - {\displaystyle \frac {a}{2}}  + c + {\displaystyle 
\frac {b}{2}}  + {\displaystyle \frac {1}{2}}  + {\displaystyle 
\frac {m}{2}}  + {\displaystyle \frac {n}{2}} ))\mbox{} - \Gamma 
( - b + 1)\,\Gamma ( - {\displaystyle \frac {b}{2}}  + 
{\displaystyle \frac {1}{2}}  + {\displaystyle \frac {n}{2}}  + 
{\displaystyle \frac {a}{2}}  - {\displaystyle \frac {m}{2}} )
 \\
\Gamma (1 + a - 2\,c - n)\,\Gamma ( - {\displaystyle \frac {a}{2}
}  + c + {\displaystyle \frac {b}{2}}  - {\displaystyle \frac {1
}{2}}  + {\displaystyle \frac {m}{2}}  + {\displaystyle \frac {n
}{2}} ) \\
\mathrm{DixonMM}(m, \,n, \, - b + 1, \, - {\displaystyle \frac {a
}{2}}  + c - {\displaystyle \frac {b}{2}}  + {\displaystyle 
\frac {1}{2}}  - {\displaystyle \frac {m}{2}}  + {\displaystyle 
\frac {n}{2}} , \, - c - n + 1)\,\Gamma (2\,c + n) \\
\Gamma ({\displaystyle \frac {1}{2}}  + {\displaystyle \frac {a}{
2}}  + {\displaystyle \frac {b}{2}}  + {\displaystyle \frac {m}{2
}}  - {\displaystyle \frac {n}{2}} ) \left/ {\vrule 
height0.80em width0em depth0.80em} \right. \!  \! (\Gamma ( - 
{\displaystyle \frac {b}{2}}  + {\displaystyle \frac {1}{2}}  + 
{\displaystyle \frac {a}{2}}  - c - {\displaystyle \frac {m}{2}} 
 - {\displaystyle \frac {n}{2}} ) \\
\Gamma ({\displaystyle \frac {1}{2}}  + {\displaystyle \frac {a}{
2}}  - c + {\displaystyle \frac {b}{2}}  + {\displaystyle \frac {
m}{2}}  - {\displaystyle \frac {n}{2}} )\,\Gamma (c + 1 - b)\,
\Gamma (a)\,\Gamma ( - {\displaystyle \frac {a}{2}}  + c + 
{\displaystyle \frac {b}{2}}  + {\displaystyle \frac {1}{2}}  + 
{\displaystyle \frac {m}{2}}  + {\displaystyle \frac {n}{2}} )\,
\Gamma (c + n) 
) }
}
\end{maplelatex}

\begin{maplelatex}
\mapleinline{inert}{2d}{` 263: `,
-GAMMA(1+a-c)*GAMMA(n+a-m)*GAMMA(a+c-2*b-m)*GAMMA(1-c+2*b+m)*DixonMM(m
,n,1+a-c,1+b-c,-b+1-n)/GAMMA(-b+a-m)/GAMMA(-c+1+b+m)/GAMMA(c-b)/GAMMA(
-c+b+1+a)*GAMMA(c)/GAMMA(c+a-m-1+n)/GAMMA(-a+1+m+b)*GAMMA(-a+m+b)+GAMM
A(1+a-c)*GAMMA(n+a-m)*GAMMA(a+c-2*b-m)*GAMMA(1-c+2*b+m)*DixonPP(m,n,1+
a-c,-c+1+b+m,-b+1)/GAMMA(-b+a-m)/GAMMA(-c+b+a+n+1)/GAMMA(1+a-b)/GAMMA(
b)*GAMMA(c)/GAMMA(c+a-m-1+n)*GAMMA(b+n)/GAMMA(-a+1+m+b);}{%
\maplemultiline{
\mathit{\ 263:\ }    - \Gamma (1 + a - c)\,\Gamma (n + a - m)\,
\Gamma (a + c - 2\,b - m)\,\Gamma (1 - c + 2\,b + m) \\
\mathrm{DixonMM}(m, \,n, \,1 + a - c, \,1 + b - c, \, - b + 1 - n
)\,\Gamma (c)\,\Gamma ( - a + m + b)/ \\ (
\Gamma ( - b + a - m)\,\Gamma ( - c + 1 + b + m)\,\Gamma (c - b)
\,\Gamma ( - c + b + 1 + a)\,\Gamma (c + a - m - 1 + n) \\
\Gamma ( - a + 1 + m + b))\mbox{} + \Gamma (1 + a - c)\,\Gamma (n
 + a - m)\,\Gamma (a + c - 2\,b - m) \\
\Gamma (1 - c + 2\,b + m)\,\mathrm{DixonPP}(m, \,n, \,1 + a - c, 
\, - c + 1 + b + m, \, - b + 1)\,\Gamma (c)\,\Gamma (b + n)/ \\
(\Gamma ( - b + a - m)\,\Gamma ( - c + b + a + n + 1)\,\Gamma (1
 + a - b)\,\Gamma (b)\,\Gamma (c + a - m - 1 + n) 
\Gamma ( - a + 1 + m + b)) }
}
\end{maplelatex}

\begin{maplelatex}
\mapleinline{inert}{2d}{` 264: `,
GAMMA(c-a+m-n)*DixonMM(m,n,c-b+1-2*a-n,-b-m+1,-a-n+1)/GAMMA(c-2*a-n+m)
/GAMMA(-a-b+1+c-n)/GAMMA(c+1-2*a-n)*GAMMA(c)/GAMMA(b+c+m-1)*GAMMA(c+b-
2*a-n+m)*GAMMA(c-b+1-2*a-n)-GAMMA(c-a+m-n)*DixonPP(m,n,c-b+1-2*a-n,-b+
1,-a+1)*GAMMA(a+n)/GAMMA(c-2*a-n+m)/GAMMA(a)/GAMMA(b)/GAMMA(-b+c-a+1)*
GAMMA(c)/GAMMA(b+c+m-1)*GAMMA(-m-c+2*a+n)*GAMMA(c+b-2*a-n+m)*GAMMA(c-b
+1-2*a-n)/GAMMA(-b-m+1);}{%
\maplemultiline{
\mathit{\ 264:\ }   \Gamma (c - a + m - n)\,\mathrm{DixonMM}(m, \,
n, \,c - b + 1 - 2\,a - n, \, - b - m + 1, \, - a - n + 1)\,
\Gamma (c) \\
\Gamma (c + b - 2\,a - n + m)\,\Gamma (c - b + 1 - 2\,a - n)/(
\Gamma (c - 2\,a - n + m) \\
\Gamma ( - a - b + 1 + c - n)\,\Gamma (c + 1 - 2\,a - n)\,\Gamma 
(b + c + m - 1))\mbox{} - \Gamma (c - a + m - n) \\
\mathrm{DixonPP}(m, \,n, \,c - b + 1 - 2\,a - n, \, - b + 1, \,
 - a + 1)\,\Gamma (a + n)\,\Gamma (c) \\
\Gamma ( - m - c + 2\,a + n)\,\Gamma (c + b - 2\,a - n + m)\,
\Gamma (c - b + 1 - 2\,a - n)/ \\ (
\Gamma (c - 2\,a - n + m)\,\Gamma (a)\,\Gamma (b)\,\Gamma ( - b
 + c - a + 1)\,\Gamma (b + c + m - 1)\,\Gamma ( - b - m + 1)) }
}
\end{maplelatex}

\begin{maplelatex}
\mapleinline{inert}{2d}{` 265: `,
(GAMMA(-b+1)*GAMMA(1+a-2*c+n)*DixonMM(m,n,-b+1,-1/2*a+c-1/2*b+1/2-1/2*
m-1/2*n,-c+1)/GAMMA(c-b+1-n)/GAMMA(3/2-1/2*b+1/2*a-c-1/2*m+1/2*n)/GAMM
A(c)*GAMMA(2*c-n)*GAMMA(1/2+1/2*a+1/2*b-1/2*m+1/2*n)/GAMMA(a)/GAMMA(-1
/2*a+c+1/2*b+1/2-1/2*m-1/2*n)*GAMMA(-1/2*a+c-1/2*b+1/2-1/2*m-1/2*n)-GA
MMA(-b+1)*GAMMA(1+a-2*c+n)*GAMMA(-1/2*a+c+1/2*b-1/2-1/2*m-1/2*n)*Dixon
PP(m,n,-b+1,-1/2*a+c-1/2*b+1/2+1/2*m-1/2*n,n+1-c)/GAMMA(c-n)/GAMMA(1/2
*a-c+1/2*b+1/2-1/2*m+1/2*n)/GAMMA(c+1-b)*GAMMA(2*c-n)*GAMMA(1/2+1/2*a+
1/2*b-1/2*m+1/2*n)/GAMMA(a)/GAMMA(-1/2*a+c+1/2*b+1/2-1/2*m-1/2*n))/GAM
MA(-1/2*b+1/2+1/2*a-c+1/2*m+1/2*n)*GAMMA(1/2-1/2*b-1/2*n+1/2*a+1/2*m);
}{%
\maplemultiline{
\mathit{\ 265:\ }   (\Gamma ( - b + 1)\,\Gamma (1 + a - 2\,c + n)
\mathrm{DixonMM}(m, \,n, \, - b + 1, \, - {\displaystyle \frac {a
}{2}}  + c - {\displaystyle \frac {b}{2}}  + {\displaystyle 
\frac {1}{2}}  - {\displaystyle \frac {m}{2}}  - {\displaystyle 
\frac {n}{2}} , \, - c + 1)\, \\ \Gamma (2\,c - n) \,
\Gamma ({\displaystyle \frac {1}{2}}  + {\displaystyle \frac {a}{
2}}  + {\displaystyle \frac {b}{2}}  - {\displaystyle \frac {m}{2
}}  + {\displaystyle \frac {n}{2}} )\,\Gamma ( - {\displaystyle 
\frac {a}{2}}  + c - {\displaystyle \frac {b}{2}}  + 
{\displaystyle \frac {1}{2}}  - {\displaystyle \frac {m}{2}}  - 
{\displaystyle \frac {n}{2}} ) \left/ {\vrule 
height0.80em width0em depth0.80em} \right. \!  \! (\Gamma (c - b
 + 1 - n) \\
\Gamma ({\displaystyle \frac {3}{2}}  - {\displaystyle \frac {b}{
2}}  + {\displaystyle \frac {a}{2}}  - c - {\displaystyle \frac {
m}{2}}  + {\displaystyle \frac {n}{2}} )\,\Gamma (c)\,\Gamma (a)
\,\Gamma ( - {\displaystyle \frac {a}{2}}  + c + {\displaystyle 
\frac {b}{2}}  + {\displaystyle \frac {1}{2}}  - {\displaystyle 
\frac {m}{2}}  - {\displaystyle \frac {n}{2}} ))\mbox{} - \Gamma 
( - b + 1) \,
\Gamma (1 + a - 2\,c + n)\, \\ \Gamma ( - {\displaystyle \frac {a}{2}
}  + c + {\displaystyle \frac {b}{2}}  - {\displaystyle \frac {1
}{2}}  - {\displaystyle \frac {m}{2}}  - {\displaystyle \frac {n
}{2}} ) 
\mathrm{DixonPP}(m, \,n, \, - b + 1, \, - {\displaystyle \frac {a
}{2}}  + c - {\displaystyle \frac {b}{2}}  + {\displaystyle 
\frac {1}{2}}  + {\displaystyle \frac {m}{2}}  - {\displaystyle 
\frac {n}{2}} , \,n + 1 - c)\, \\ \Gamma (2\,c - n) \,
\Gamma ({\displaystyle \frac {1}{2}}  + {\displaystyle \frac {a}{
2}}  + {\displaystyle \frac {b}{2}}  - {\displaystyle \frac {m}{2
}}  + {\displaystyle \frac {n}{2}} ) \left/ {\vrule 
height0.80em width0em depth0.80em} \right. \!  \! (\Gamma (c - n)
\,\Gamma ({\displaystyle \frac {a}{2}}  - c + {\displaystyle 
\frac {b}{2}}  + {\displaystyle \frac {1}{2}}  - {\displaystyle 
\frac {m}{2}}  + {\displaystyle \frac {n}{2}} )\,\Gamma (c + 1 - 
b) \\
\Gamma (a)\,\Gamma ( - {\displaystyle \frac {a}{2}}  + c + 
{\displaystyle \frac {b}{2}}  + {\displaystyle \frac {1}{2}}  - 
{\displaystyle \frac {m}{2}}  - {\displaystyle \frac {n}{2}} )))
\Gamma ({\displaystyle \frac {1}{2}}  - {\displaystyle \frac {b}{
2}}  - {\displaystyle \frac {n}{2}}  + {\displaystyle \frac {a}{2
}}  + {\displaystyle \frac {m}{2}} ) \left/ {\vrule 
height0.80em width0em depth0.80em} \right. \!  \!  
\Gamma ( - {\displaystyle \frac {b}{2}}  + {\displaystyle \frac {
1}{2}}  + {\displaystyle \frac {a}{2}}  - c + {\displaystyle 
\frac {m}{2}}  + {\displaystyle \frac {n}{2}} ) }
}
\end{maplelatex}

\begin{maplelatex}
\mapleinline{inert}{2d}{` 266: `,
(-GAMMA(1+a-c)*GAMMA(a+c-2*b+m)*GAMMA(1-c+2*b-m)*GAMMA(-a-m+b)*DixonPP
(m,n,1+a-c,1+b-c,n+1-b)/GAMMA(-c+1+b-m)/GAMMA(c-b)/GAMMA(-c+b+1+a)*GAM
MA(c)/GAMMA(c+a+m-1-n)/GAMMA(-a-m+1+b)+GAMMA(1+a-c)*GAMMA(a+c-2*b+m)*G
AMMA(1-c+2*b-m)*DixonMM(m,n,1+a-c,-c+1+b-m,-b+1)/GAMMA(1+a-b)/GAMMA(b)
*GAMMA(c)/GAMMA(c+a+m-1-n)*GAMMA(b-n)/GAMMA(-a-m+1+b)/GAMMA(1+a+b-c-n)
)*GAMMA(a-n+m)/GAMMA(-b+a+m);}{%
\maplemultiline{
\mathit{\ 266:\ }   ( - \Gamma (1 + a - c)\,\Gamma (a + c - 2\,b
 + m)\,\Gamma (1 - c + 2\,b - m)\,\Gamma ( - a - m + b) \\
\mathrm{DixonPP}(m, \,n, \,1 + a - c, \,1 + b - c, \,n + 1 - b)\,
\Gamma (c)/(\Gamma ( - c + 1 + b - m)\,\Gamma (c - b) \\
\Gamma ( - c + b + 1 + a)\,\Gamma (c + a + m - 1 - n)\,\Gamma (
 - a - m + 1 + b))\mbox{} \\ + \, \Gamma (1 + a - c) 
\Gamma (a + c - 2\,b + m)\,\Gamma (1 - c + 2\,b - m) \\
\mathrm{DixonMM}(m, \,n, \,1 + a - c, \, - c + 1 + b - m, \, - b
 + 1)\,\Gamma (c)\,\Gamma (b - n)/(\Gamma (1 + a - b) \\
\Gamma (b)\,\Gamma (c + a + m - 1 - n)\,\Gamma ( - a - m + 1 + b)
\,\Gamma (1 + a + b - c - n)))\Gamma (a - n + m)/ \\
\Gamma ( - b + a + m) }
}
\end{maplelatex}

\begin{maplelatex}
\mapleinline{inert}{2d}{` 296: `,
(-GAMMA(c-a)*GAMMA(2*a-m+n)*GAMMA(1+a-c)*GAMMA(a+c-2*b-m)*DixonMM(m,n,
1+a-c,1+b-c,-b+1-n)*GAMMA(b)/GAMMA(c-b)/GAMMA(-c+1+b+m)/GAMMA(-b+a-m)/
GAMMA(-c+b+1+a)/GAMMA(c+a-m-1+n)/GAMMA(-a+1+m+b)/GAMMA(b-a)*GAMMA(-a+m
+b)+GAMMA(c-a)*GAMMA(2*a-m+n)*GAMMA(1+a-c)*GAMMA(a+c-2*b-m)*DixonPP(m,
n,1+a-c,-c+1+b+m,-b+1)/GAMMA(-b+a-m)/GAMMA(-c+b+a+n+1)/GAMMA(1+a-b)/GA
MMA(c+a-m-1+n)/GAMMA(-a+1+m+b)*GAMMA(b+n)/GAMMA(b-a))*GAMMA(1-c+2*b+m-
a)-GAMMA(c-a)*GAMMA(2*a-m+n)*GAMMA(a-b)*GAMMA(-2*b-a+1+c-n)*GAMMA(-a+m
+1-n)*GAMMA(-a+2-c-n+m)*DixonMM(m,n,-2*b-a+1+c-n,-b+c-m,-b+1-n)*GAMMA(
b)/GAMMA(a)/GAMMA(c-b)/GAMMA(-c+1+b+m)/GAMMA(b+a-m+n)/GAMMA(-b-a+1-n+m
)/GAMMA(-a-b+1+c-n)/GAMMA(-b-a+2-n)*GAMMA(b-a+1)/GAMMA(b-a)+GAMMA(c-a)
*GAMMA(2*a-m+n)*GAMMA(a-b)*GAMMA(-2*b-a+1+c-n)*GAMMA(-a+m+1-n)*GAMMA(-
a+2-c-n+m)*GAMMA(b-a+1)*GAMMA(-m+b+a-1+n)*DixonPP(m,n,-2*b-a+1+c-n,c-b
,-b+1)*GAMMA(b+n)/GAMMA(b-a)/GAMMA(a)/GAMMA(c-b)/GAMMA(-c+1+b+m)/GAMMA
(b+a-m+n)/GAMMA(-b-a+1-n+m)/GAMMA(-b+c-m)/GAMMA(-b+c-a+1)/GAMMA(1+b-c)
;}{%
\maplemultiline{
\mathit{\ 296:\ }   ( - \Gamma (c - a)\,\Gamma (2\,a - m + n)\,
\Gamma (1 + a - c)\,\Gamma (a + c - 2\,b - m) \\
\mathrm{DixonMM}(m, \,n, \,1 + a - c, \,1 + b - c, \, - b + 1 - n
)\,\Gamma (b)\,\Gamma ( - a + m + b)/(\Gamma (c - b) \\
\Gamma ( - c + 1 + b + m)\,\Gamma ( - b + a - m)\,\Gamma ( - c + 
b + 1 + a)\,\Gamma (c + a - m - 1 + n) \\
\Gamma ( - a + 1 + m + b)\,\Gamma (b - a))\mbox{} + \Gamma (c - a
)\,\Gamma (2\,a - m + n)\,\Gamma (1 + a - c) \\
\Gamma (a + c - 2\,b - m)\,\mathrm{DixonPP}(m, \,n, \,1 + a - c, 
\, - c + 1 + b + m, \, - b + 1)\,\Gamma (b + n)/( \\
\Gamma ( - b + a - m)\,\Gamma ( - c + b + a + n + 1)\,\Gamma (1
 + a - b)\,\Gamma (c + a - m - 1 + n) \\
\Gamma ( - a + 1 + m + b)\,\Gamma (b - a)))\Gamma (1 - c + 2\,b
 + m - a)\mbox{} - \Gamma (c - a)\,\Gamma (2\,a - m + n) \\
\Gamma (a - b)\,\Gamma ( - 2\,b - a + 1 + c - n)\,\Gamma ( - a + 
m + 1 - n)\,\Gamma ( - a + 2 - c - n + m) \\
\mathrm{DixonMM}(m, \,n, \, - 2\,b - a + 1 + c - n, \, - b + c - 
m, \, - b + 1 - n)\,\Gamma (b)\,\Gamma (b - a + 1)/( \\
\Gamma (a)\,\Gamma (c - b)\,\Gamma ( - c + 1 + b + m)\,\Gamma (b
 + a - m + n)\,\Gamma ( - b - a + 1 - n + m) \\
\Gamma ( - a - b + 1 + c - n)\,\Gamma ( - b - a + 2 - n)\,\Gamma 
(b - a))\mbox{} + \Gamma (c - a)\,\Gamma (2\,a - m + n) \\
\Gamma (a - b)\,\Gamma ( - 2\,b - a + 1 + c - n)\,\Gamma ( - a + 
m + 1 - n)\,\Gamma ( - a + 2 - c - n + m) \\
\Gamma (b - a + 1)\,\Gamma ( - m + b + a - 1 + n) \\
\mathrm{DixonPP}(m, \,n, \, - 2\,b - a + 1 + c - n, \,c - b, \,
 - b + 1)\,\Gamma (b + n)/(\Gamma (b - a)\,\Gamma (a) \\
\Gamma (c - b)\,\Gamma ( - c + 1 + b + m)\,\Gamma (b + a - m + n)
\,\Gamma ( - b - a + 1 - n + m)\,\Gamma ( - b + c - m) \\
\Gamma ( - b + c - a + 1)\,\Gamma (1 + b - c)) }
}
\end{maplelatex}
\emptyline
\begin{maplelatex}
\mapleinline{inert}{2d}{` 297: `,
(-GAMMA(c-a)*GAMMA(2*a-n+m)*GAMMA(1+a-c)*GAMMA(a+c-2*b+m)*GAMMA(-a-m+b
)*DixonPP(m,n,1+a-c,1+b-c,n+1-b)*GAMMA(b)/GAMMA(c-b)/GAMMA(-c+1+b-m)/G
AMMA(-c+b+1+a)/GAMMA(c+a+m-1-n)/GAMMA(-a-m+1+b)+GAMMA(c-a)*GAMMA(2*a-n
+m)*GAMMA(1+a-c)*GAMMA(a+c-2*b+m)*DixonMM(m,n,1+a-c,-c+1+b-m,-b+1)*GAM
MA(b-n)/GAMMA(1+a-b)/GAMMA(c+a+m-1-n)/GAMMA(-a-m+1+b)/GAMMA(1+a+b-c-n)
)/GAMMA(b-a)/GAMMA(-b+a+m)*GAMMA(1-c+2*b-m-a)-GAMMA(c-a)*GAMMA(2*a-n+m
)*GAMMA(a-b)*GAMMA(-2*b-a+1+c+n)*GAMMA(1-a+n-m)*GAMMA(-a+2-c+n-m)*GAMM
A(b-a+1)*DixonPP(m,n,-2*b-a+1+c+n,c+m-b,n+1-b)*GAMMA(b)/GAMMA(b-a)/GAM
MA(a)/GAMMA(c-b)/GAMMA(-c+1+b-m)/GAMMA(b+a-n+m)/GAMMA(-b-a+1+n-m)/GAMM
A(-b+c-a+n+1)/GAMMA(-b-a+2+n)+GAMMA(c-a)*GAMMA(2*a-n+m)*GAMMA(a-b)*GAM
MA(-2*b-a+1+c+n)*GAMMA(1-a+n-m)*GAMMA(-a+2-c+n-m)*GAMMA(b-a+1)*GAMMA(m
+b+a-1-n)*DixonMM(m,n,-2*b-a+1+c+n,c-b,-b+1)*GAMMA(b-n)/GAMMA(b-a)/GAM
MA(a)/GAMMA(c-b)/GAMMA(-c+1+b-m)/GAMMA(b+a-n+m)/GAMMA(-b-a+1+n-m)/GAMM
A(c+m-b)/GAMMA(-b+c-a+1)/GAMMA(1+b-c);}{%
\maplemultiline{
\mathit{\ 297:\ }   ( - \Gamma (c - a)\,\Gamma (2\,a - n + m)\,
\Gamma (1 + a - c)\,\Gamma (a + c - 2\,b + m)\,\Gamma ( - a - m
 + b) \\
\mathrm{DixonPP}(m, \,n, \,1 + a - c, \,1 + b - c, \,n + 1 - b)\,
\Gamma (b)/(\Gamma (c - b)\,\Gamma ( - c + 1 + b - m) \\
\Gamma ( - c + b + 1 + a)\,\Gamma (c + a + m - 1 - n)\,\Gamma (
 - a - m + 1 + b))\mbox{} + \Gamma (c - a) \\
\Gamma (2\,a - n + m)\,\Gamma (1 + a - c)\,\Gamma (a + c - 2\,b
 + m) \\
\mathrm{DixonMM}(m, \,n, \,1 + a - c, \, - c + 1 + b - m, \, - b
 + 1)\,\Gamma (b - n)/(\Gamma (1 + a - b) \\
\Gamma (c + a + m - 1 - n)\,\Gamma ( - a - m + 1 + b)\,\Gamma (1
 + a + b - c - n))) \\
\Gamma (1 - c + 2\,b - m - a)/(\Gamma (b - a)\,\Gamma ( - b + a
 + m))\mbox{} - \Gamma (c - a)\,\Gamma (2\,a - n + m)\,\Gamma (a
 - b) \\
\Gamma ( - 2\,b - a + 1 + c + n)\,\Gamma (1 - a + n - m)\,\Gamma 
( - a + 2 - c + n - m)\,\Gamma (b - a + 1) \\
\mathrm{DixonPP}(m, \,n, \, - 2\,b - a + 1 + c + n, \,c + m - b, 
\,n + 1 - b)\,\Gamma (b)/(\Gamma (b - a)\,\Gamma (a) \\
\Gamma (c - b)\,\Gamma ( - c + 1 + b - m)\,\Gamma (b + a - n + m)
\,\Gamma ( - b - a + 1 + n - m) \\
\Gamma ( - b + c - a + n + 1)\,\Gamma ( - b - a + 2 + n))\mbox{}
 + \Gamma (c - a)\,\Gamma (2\,a - n + m)\,\Gamma (a - b) \\
\Gamma ( - 2\,b - a + 1 + c + n)\,\Gamma (1 - a + n - m)\,\Gamma 
( - a + 2 - c + n - m)\,\Gamma (b - a + 1) \\
\Gamma (m + b + a - 1 - n)\,\mathrm{DixonMM}(m, \,n, \, - 2\,b - 
a + 1 + c + n, \,c - b, \, - b + 1)\,\Gamma (b - n)/ \\
(\Gamma (b - a)\,\Gamma (a)\,\Gamma (c - b)\,\Gamma ( - c + 1 + b
 - m)\,\Gamma (b + a - n + m)\,\Gamma ( - b - a + 1 + n - m) \\
\Gamma (c + m - b)\,\Gamma ( - b + c - a + 1)\,\Gamma (1 + b - c)
) }
}
\end{maplelatex}

\end{maplegroup}
\begin{maplegroup}
\emptyline
\end{maplegroup}
\begin{maplegroup}
\begin{mapleinput}
\mapleinline{active}{1d}{

}{%
}
\end{mapleinput}

\end{maplegroup}

%% file: cut_paste.tex
\pagestyle{empty}
\begin{maplegroup}
\mapleresult
C( 1 ):=F32([a, b, c] , [2*c+n, 1/2+1/2*a+1/2*b+1/2*m] ,1);
C( 2 ):=F32([a, b, c] , [2*c+n, 1/2+1/2*a+1/2*b+m] ,1);C( 3 ):=F32([a, b, c] , [2*c+n, 1+1/2*a+1/2*b+m] ,1);C( 4 ):=F32([a, b, c] , [2*c-n, 1/2+1/2*a+1/2*b+1/2*m] ,1);C( 5 ):=F32([a, b, c] , [2*c-n, 1/2+1/2*a+1/2*b+m] ,1);C( 6 ):=F32([a, b, c] , [2*c-n, 1+1/2*a+1/2*b+m] ,1);C( 7 ):=F32([a, b, c] , [2*c+n, 1/2+1/2*a+1/2*b-1/2*m] ,1);C( 8 ):=F32([a, b, c] , [2*c+n, 1/2+1/2*a+1/2*b-m] ,1);C( 9 ):=F32([a, b, c] , [2*c+n, 1/2*a+1/2*b-m] ,1);C( 10 ):=F32([a, b, c] , [2*c-n, 1/2+1/2*a+1/2*b-1/2*m] ,1);C( 11 ):=F32([a, b, c] , [2*c-n, 1/2+1/2*a+1/2*b-m] ,1);C( 12 ):=F32([a, b, c] , [2*c-n, 1/2*a+1/2*b-m] ,1);C( 13 ):=F32([a, b, 1-a+m] , [c, 1+2*b-c+n] ,1);C( 14 ):=F32([a, b, 1-a+m] , [c, 1+2*b-c-n] ,1);C( 15 ):=F32([a, b, 1-a-m] , [c, 1+2*b-c+n] ,1);C( 16 ):=F32([a, b, 1-a-m] , [c, 1+2*b-c-n] ,1);C( 17 ):=F32([a, b, c] , [1+a-c+n, 1+a-b+m] ,1);C( 18 ):=F32([a, b, c] , [1+a-c-n, 1+a-b+m] ,1);C( 19 ):=F32([a, b, c] , [1+a-c-n, 1+a-b-m] ,1);C( 20 ):=F32([a, b, c] , [c+m, b-n] ,1);C( 21 ):=F32([a, m, b] , [c, -n+m] ,1);C( 23 ):=F32([a, b, -n] , [c, -n+m] ,1);C( 24 ):=F32([a, -n, b] , [c, a-c+b-n+m] ,1);C( 25 ):=F32([a, 1, b] , [n+1, c] ,1);C( 26 ):=F32([a, b, n] , [n+1, c] ,1);C( 27 ):=F32([a, b, c] , [c+1, b+n] ,1);C( 28 ):=F32([1, a, b] , [c, -c+b+a+n+1] ,1);C( 29 ):=F32([n, a, b] , [c, -c+b+a+n+1] ,1);C( 30 ):=F32([a, b, c] , [b+1, -n+b+1] ,1);C( 31 ):=F32([a, b, 1] , [c, -n+b+1] ,1);C( 33 ):=F32([a, b, 1-n] , [c, b+1] ,1);C( 34 ):=F32([a, b, c] , [a-n+1, c+b+n] ,1);C( 35 ):=F32([a, a-n, b] , [a-n+1, c] ,1);C( 36 ):=F32([1, 1-n, a] , [b, c] ,1);C( 37 ):=F32([a, b, n+1/2*a-b-3/2] , [1/2*a+b+1/2, n+a-b-1] ,1);C( 38 ):=F32([a, b, n-b-1] , [n-b-2+2*a, b+1] ,1);C( 39 ):=F32([1, a, b] , [n-b, n-a] ,1);C( 40 ):=F32([1-4*a+2*n, 1/2+n-a, n-a] , [1-3*a+2*n, -3*a+2*n+3/2] ,1);C( 41 ):=F32([1-4*a+2*n, 1-2*a+n, 1/2-2*a+n] , [1-3*a+2*n, 2-4*a+2*n] ,1);C( 42 ):=F32([1/2+n-a, 1-2*a+n, a] , [1-3*a+2*n, 3/2-a+n] ,1);C( 43 ):=F32([n-a, 1/2-2*a+n, a] , [1-3*a+2*n, n+1-a] ,1);C( 44 ):=F32([1/2+n-a, 3/2-2*a+n, a+1/2] , [-3*a+2*n+3/2, 3/2-a+n] ,1);C( 45 ):=F32([n-a, 1-2*a+n, a+1/2] , [-3*a+2*n+3/2, n+1-a] ,1);C( 46 ):=F32([1/2-2*a+n, 1-2*a+n, 1] , [2-4*a+2*n, n+1-a] ,1);C( 47 ):=F32([a, a+1/2, 1] , [3/2-a+n, n+1-a] ,1);C( 48 ):=F32([a, b, c] , [e, c-1] ,1);C( 50 ):=F32([a, -1, b] , [c, e] ,1);C( 51 ):=F32([a, b, c] , [c+2, 2+b] ,1);C( 52 ):=F32([a, 2, b] , [c, 4] ,1);C( 53 ):=F32([a, 2, b] , [c, 4+a-c+b] ,1);C( 54 ):=F32([-n, a, b] , [b-m, a-n+m] ,1);C( 55 ):=F32([-L-m+n, a, b] , [b-m, a-L] ,1);C( 56 ):=F32([a, b, b+L-n] , [b-m, b-n] ,1);C( 57 ):=F32([-n, a, b] , [b-m, a-L] ,1);C( 61 ):=F32([-n, a, b] , [b-m, a-2-n+m] ,1);C( 62 ):=F32([a, b, -n*(n-b+a)/(b-n)] , [-(a*n+b*n-b+n-b\symbol{94}2)/(b-n), n+1+a] ,1);C( 63 ):=F32([a, b, 1/2*a+1/2] , [1+a, n-1+1/2*a+1/2*b] ,1);C( 64 ):=F32([a, b, 1] , [2*b, n-1/2+1/2*a] ,1);C( 65 ):=F32([a, n-3/2+1/2*a, b] , [2*n-2+a-b, n-1/2+1/2*a] ,1);C( 66 ):=F32([3-4*a-2*n, -a-n+1, -a+3/2-n] , [3-3*a-2*n, 7/2-3*a-2*n] ,1);C( 67 ):=F32([-a-n+1, 3/2-2*a-n, a] , [3-3*a-2*n, 2-a-n] ,1);C( 68 ):=F32([-a+3/2-n, -2*a-n+2, a] , [3-3*a-2*n, 5/2-a-n] ,1);C( 69 ):=F32([3-4*a-2*n, -2*a-n+2, 5/2-2*a-n] , [7/2-3*a-2*n, 4-4*a-2*n] ,1);C( 70 ):=F32([-a-n+1, -2*a-n+2, a+1/2] , [7/2-3*a-2*n, 2-a-n] ,1);C( 71 ):=F32([-a+3/2-n, 5/2-2*a-n, a+1/2] , [7/2-3*a-2*n, 5/2-a-n] ,1);C( 72 ):=F32([-2*a-n+2, 5/2-2*a-n, 1] , [4-4*a-2*n, 5/2-a-n] ,1);C( 73 ):=F32([a, a+1/2, 1] , [2-a-n, 5/2-a-n] ,1);C( 74 ):=F32([-4*a+2*n, 1/2-2*a+n, n-2*a] , [-3*a+2*n, 1-4*a+2*n] ,1);C( 75 ):=F32([-a+n-1/2, n-2*a, a] , [-3*a+2*n, 1/2+n-a] ,1);C( 76 ):=F32([-4*a+2*n, 1-2*a+n, 1/2-2*a+n] , [-3*a+2*n+1/2, 1-4*a+2*n] ,1);C( 77 ):=F32([n-a, 1-2*a+n, a+1/2] , [-3*a+2*n+1/2, n+1-a] ,1);C( 78 ):=F32([1/2-2*a+n, 1-2*a+n, 1] , [1-4*a+2*n, n+1-a] ,1);C( 79 ):=F32([n-2*a, 1/2-2*a+n, 1] , [1-4*a+2*n, 1/2+n-a] ,1);C( 80 ):=F32([a, b, c] , [1/2*c+1/2*b+n+m-1/2, 2*a-2*n+1] ,1);C( 81 ):=F32([a, b, c] , [1/2*c+1/2*b+n+m, 2*a-2*n+1] ,1);C( 82 ):=F32([a, b, c] , [1/2*c+1/2*b+n+m, 2*a-2*n] ,1);C( 84 ):=F32([a, b, c] , [2+a, a-n+1] ,1);C( 85 ):=F32([a, b, 2] , [c, a-n+1] ,1);C( 86 ):=F32([a, b, c] , [n+1+c+b, a-n+1] ,1);C( 87 ):=F32([a, a-n-1, b] , [c, a-n+1] ,1);C( 88 ):=F32([a, b, 1-n] , [c, 2+b] ,1);C( 89 ):=F32([2, a, 1-n] , [b, c] ,1);C( 90 ):=F32([3*n-3/2, -a+n-1/2, 2*a] , [2*a+1/2+n, 3*n-a-1/2] ,1);C( 91 ):=F32([a, b, c] , [2+a, (a*c+a*b-a-b*c+b+c-1)/a] ,1);C( 92 ):=F32([2, a, b] , [c, (-3+a*b-b-a+2*c)/(c-2)] ,1);C( 93 ):=F32([a, b, c] , [1/2*c+1+1/2*b+1/2*a+1/2*U, 1/2*c+1+1/2*b+1/2*a-1/2*U] ,1);C( 94 ):=F32([3/2, a, a-1/2] , [2*a-1/6, 1/6+2*a] ,1);C( 95 ):=F32([3/2, a, a-1/2] , [2*a-5/6, 2*a-1/6] ,1);C( 96 ):=F32([a, a+1/3, 2*a-5/3] , [2*a-1/6, 3*a-1] ,1);C( 97 ):=F32([a, a+1/3, 2*a-2/3] , [2*a+5/6, 3*a] ,1);C( 98 ):=F32([3/2, a, a-1/2] , [-7/6+2*a, 2*a-5/6] ,1);C( 99 ):=F32([a, a+2/3, 2*a-4/3] , [1/6+2*a, 3*a-1] ,1);C( 100 ):=F32([a, a+2/3, 2*a-1/3] , [7/6+2*a, 3*a] ,1);C( 101 ):=F32([a, a+1/3, 2*a-5/3] , [2*a-1/6, 3*a-2] ,1);C( 102 ):=F32([a, a+1/3, 2*a-2/3] , [2*a+5/6, 3*a-1] ,1);C( 103 ):=F32([a, a+1/3, a+2/3] , [3/2*a+3/2, 1+3/2*a] ,1);C( 104 ):=F32([2*a, a+1/3-1/2*n, a-1/3-1/2*n] , [1/2+2*a-n, 3*a-1/2*n] ,1);C( 105 ):=F32([2*a, a-1/2*n+5/6, a+1/6-1/2*n] , [1/2+2*a-n, 3*a+1/2-1/2*n] ,1);C( 106 ):=F32([a+1/3-1/2*n, a-1/2*n+5/6, -n+1/2] , [1/2+2*a-n, 2*a+5/6-n] ,1);C( 107 ):=F32([a-1/3-1/2*n, a+1/6-1/2*n, -n+1/2] , [1/2+2*a-n, 2*a+1/6-n] ,1);C( 108 ):=F32([2*a, 2*a+1/3, 2*a-1/3] , [3*a-1/2*n, 3*a+1/2-1/2*n] ,1);C( 109 ):=F32([a+1/3-1/2*n, 2*a+1/3, a-1/2*n] , [3*a-1/2*n, 2*a+5/6-n] ,1);C( 110 ):=F32([a-1/3-1/2*n, 2*a-1/3, a-1/2*n] , [3*a-1/2*n, 2*a+1/6-n] ,1);C( 111 ):=F32([a-1/2*n+5/6, 2*a+1/3, 1/2+a-1/2*n] , [3*a+1/2-1/2*n, 2*a+5/6-n] ,1);C( 112 ):=F32([a+1/6-1/2*n, 2*a-1/3, 1/2+a-1/2*n] , [3*a+1/2-1/2*n, 2*a+1/6-n] ,1);C( 113 ):=F32([-n+1/2, a-1/2*n, 1/2+a-1/2*n] , [2*a+5/6-n, 2*a+1/6-n] ,1);C( 114 ):=F32([-1/2+n, 2*a+n-2/3, -1/6-a] , [n-1/2+a, -1/6+a+n] ,1);C( 115 ):=F32([-1/2+n, 2*a-1/3+n, -a+1/6] , [n-1/2+a, 1/6+a+n] ,1);C( 116 ):=F32([2*a+n-2/3, 2*a-1/3+n, a] , [n-1/2+a, 3*a+n] ,1);C( 117 ):=F32([-1/6-a, -a+1/6, a] , [n-1/2+a, 1/2] ,1);C( 118 ):=F32([-1/2+n, 2*a+n, 1/2-a] , [-1/6+a+n, 1/6+a+n] ,1);C( 119 ):=F32([2*a+n-2/3, 2*a+n, a+1/3] , [-1/6+a+n, 3*a+n] ,1);C( 120 ):=F32([-1/6-a, 1/2-a, a+1/3] , [-1/6+a+n, 1/2] ,1);C( 121 ):=F32([2*a-1/3+n, 2*a+n, a+2/3] , [1/6+a+n, 3*a+n] ,1);C( 122 ):=F32([-a+1/6, 1/2-a, a+2/3] , [1/6+a+n, 1/2] ,1);C( 123 ):=F32([a, a+1/3, a+2/3] , [3*a+n, 1/2] ,1);C( 124 ):=F32([1/2, 2*a+1/3-n, -1/6+n-a] , [a+1/2, a+5/6] ,1);C( 125 ):=F32([1/2, 2*a+2/3-n, 1/6+n-a] , [a+1/2, 7/6+a] ,1);C( 126 ):=F32([2*a+1/3-n, 2*a+2/3-n, a] , [a+1/2, 3*a+1-n] ,1);C( 127 ):=F32([-1/6+n-a, 1/6+n-a, a] , [a+1/2, 1/2+n] ,1);C( 128 ):=F32([1/2, 2*a+1-n, 1/2+n-a] , [a+5/6, 7/6+a] ,1);C( 129 ):=F32([2*a+1/3-n, 2*a+1-n, a+1/3] , [a+5/6, 3*a+1-n] ,1);C( 130 ):=F32([-1/6+n-a, 1/2+n-a, a+1/3] , [a+5/6, 1/2+n] ,1);C( 131 ):=F32([2*a+2/3-n, 2*a+1-n, a+2/3] , [7/6+a, 3*a+1-n] ,1);C( 132 ):=F32([1/6+n-a, 1/2+n-a, a+2/3] , [7/6+a, 1/2+n] ,1);C( 133 ):=F32([a, a+1/3, a+2/3] , [3*a+1-n, 1/2+n] ,1);C( 134 ):=F32([1, 2-n+a, b-n+2] , [-a+2, -b+2] ,1);C( 135 ):=F32([a, 4*a-2*n, 2*a-1/2-n] , [3*a-n, 1/2+2*a-n] ,1);C( 136 ):=F32([a, 2*a-n, 4*a-2*n] , [2*a+1-n, 3*a+1/2-n] ,1);C( 137 ):=F32([2*a-n, 4*a-2*n, a+1/2] , [3*a-n, 2*a+1-n] ,1);C( 138 ):=F32([1, 4*a-2*n, a-n+1] , [2*a+1-n, 2*a+3/2-n] ,1);C( 139 ):=F32([1, 2-2*b, 5/2-n-1/2*a] , [-a+2, -b+2] ,1);C( 140 ):=F32([a, -2+4*a+2*n, -1+n+2*a] , [2*a+n, 3*a-1/2+n] ,1);C( 141 ):=F32([a+1/2, -2+4*a+2*n, 2*a-1/2+n] , [2*a+1/2+n, 3*a-1/2+n] ,1);C( 142 ):=F32([a+1/2, -2+4*a+2*n, -1+n+2*a] , [2*a+n, 3*a+n-1] ,1);C( 143 ):=F32([1, -2+4*a+2*n, n-1/2+a] , [2*a+n, 2*a-1/2+n] ,1);C( 144 ):=F32([1, a-n+1, 1+4*a-2*n] , [2*a+1-n, 2*a+3/2-n] ,1);C( 145 ):=F32([1, 1+4*a-2*n, 3/2-n+a] , [2-n+2*a, 2*a+3/2-n] ,1);C( 146 ):=F32([2*a-n, n-a, 2*a-1/2-n] , [1/2, 3*a-n] ,1);C( 147 ):=F32([-1+n+2*a, 2*a-1/2+n, -a+3/2-n] , [3/2, 3*a-1/2+n] ,1);C( 148 ):=F32([a, b, 1/2*a+1/2] , [1+a, 1/2*a+3/2+1/2*b-1/2*n] ,1);C( 149 ):=F32([a, b, 1] , [2*b, 1/2*a-1/2*n+2] ,1);C( 150 ):=F32([a, b, -b+5/2+1/2*a-n] , [1/2*a+b+1/2, 3+a-n-b] ,1);C( 151 ):=F32([a, 1/2*a-1/2*n+1, b] , [3+a-n-b, 1/2*a-1/2*n+2] ,1);C( 152 ):=F32([a, b, -b+3-n] , [-n-b+2+2*a, b+1] ,1);C( 153 ):=F32([1-2*a+n, 1, a+1/2] , [n+1-a, 2*a+1-n] ,1);C( 154 ):=F32([1-2*a+n, 1/2-2*a+n, n-a] , [n+1-a, 1/2] ,1);C( 155 ):=F32([1, 1/2-2*a+n, a] , [n+1-a, 1/2+2*a-n] ,1);C( 156 ):=F32([a+1/2, n-a, a] , [n+1-a, 3*a-n] ,1);C( 157 ):=F32([1/2-2*a+n, a-n+1/2, 2*a-1/2-n] , [1/2, 1/2+2*a-n] ,1);C( 158 ):=F32([1-2*a+n, 1, a] , [3/2-a+n, 2*a+1-n] ,1);C( 159 ):=F32([1-2*a+n, a-n+1, 2*a-n] , [3/2, 2*a+1-n] ,1);C( 160 ):=F32([a, 2*n-2-a, b] , [1+a, 2*b-a] ,1);C( 161 ):=F32([a, b, 2*a-2*n+3] , [1+a, 2*a-2*n+4-b] ,1);C( 162 ):=F32([a, b, 1] , [2*n-1-a, -b+2] ,1);C( 163 ):=F32([a, b, 1] , [n-1/2+1/2*b, 2*a-2*n+3] ,1);C( 164 ):=F32([a, b, 2*a-2*n+2*b+2] , [2*b+a, 2*a-2*n+b+3] ,1);C( 165 ):=F32([a, b, 2*b] , [-a+2*n-2+2*b, b+1] ,1);C( 166 ):=F32([a, 2-2*n+2*a, b] , [a+1/2*b, 2*a-2*n+3] ,1);C( 167 ):=F32([a, -a+1, b] , [a+2*b-4+2*n, -a+2] ,1);C( 168 ):=F32([a, b, 2-n+1/2*b] , [1/2*a+1/2+1/2*b, b+1] ,1);C( 169 ):=F32([5/2-2*a-n, a+1/2, 1] , [5/2-a-n, 2*a+1/2+n] ,1);C( 170 ):=F32([5/2-2*a-n, -a+3/2-n, -2*a-n+2] , [5/2-a-n, 3/2] ,1);C( 171 ):=F32([a+1/2, -a+3/2-n, a] , [5/2-a-n, 3*a-1/2+n] ,1);C( 172 ):=F32([1, -2*a-n+2, a] , [5/2-a-n, 2*a+n] ,1);C( 173 ):=F32([5/2-2*a-n, 2*a-1/2+n, a+n] , [2*a+1/2+n, 3/2] ,1);C( 174 ):=F32([1, -2*a-n+2, a+1/2] , [2-a-n, 2*a+n] ,1);C( 175 ):=F32([-2*a-n+2, n-1/2+a, -1+n+2*a] , [1/2, 2*a+n] ,1);C( 176 ):=F32([-a+n-1/2, 1/2-2*a+n, n-2*a] , [1/2+n-a, 1/2] ,1);C( 177 ):=F32([-a+n-1/2, a+1/2, a] , [1/2+n-a, 3*a+1/2-n] ,1);C( 178 ):=F32([1/2-2*a+n, a+1/2, 1] , [1/2+n-a, 2*a+3/2-n] ,1);C( 179 ):=F32([n-2*a, a, 1] , [1/2+n-a, 2*a+1-n] ,1);C( 180 ):=F32([-a+n-1/2, 1/2+2*a-n, 2*a-n] , [1/2, 3*a+1/2-n] ,1);C( 181 ):=F32([1/2-2*a+n, 1/2+2*a-n, a-n+1] , [1/2, 2*a+3/2-n] ,1);C( 182 ):=F32([n-2*a, 2*a-n, a-n+1] , [1/2, 2*a+1-n] ,1);C( 183 ):=F32([a+1/2, 1/2+2*a-n, 1+4*a-2*n] , [3*a+1/2-n, 2*a+3/2-n] ,1);C( 184 ):=F32([a, 2*a-n, 1+4*a-2*n] , [3*a+1/2-n, 2*a+1-n] ,1);C( 185 ):=F32([n-a, a, a+1/2] , [n+1-a, 3*a+1-n] ,1);C( 186 ):=F32([n-a, 1/2-2*a+n, 1-2*a+n] , [n+1-a, 3/2] ,1);C( 187 ):=F32([a, 1/2-2*a+n, 1] , [n+1-a, 2*a+3/2-n] ,1);C( 188 ):=F32([a+1/2, 1-2*a+n, 1] , [n+1-a, 2-n+2*a] ,1);C( 189 ):=F32([n-a, 1/2+2*a-n, 2*a+1-n] , [3*a+1-n, 3/2] ,1);C( 190 ):=F32([a, 1/2+2*a-n, 1+4*a-2*n] , [3*a+1-n, 2*a+3/2-n] ,1);C( 191 ):=F32([a+1/2, 2*a+1-n, 1+4*a-2*n] , [3*a+1-n, 2-n+2*a] ,1);C( 192 ):=F32([1/2-2*a+n, 1/2+2*a-n, 3/2-n+a] , [3/2, 2*a+3/2-n] ,1);C( 193 ):=F32([1-2*a+n, 2*a+1-n, 3/2-n+a] , [3/2, 2-n+2*a] ,1);C( 194 ):=F32([-2*a+2, n-b-1, n-2*b-1] , [n-2*b, n-a-b] ,1);C( 195 ):=F32([a, 3-2*n+b, 1/2*a-n+3/2] , [1+a-b, 5/2+1/2*a-n] ,1);C( 196 ):=F32([1, 2-2*b, -1/2*a+1/2*n] , [-b+2, -a+2] ,1);C( 197 ):=F32([-2*a+2, -b+3-n, -2*b+3-n] , [-n-2*b+4, -n-a+4-b] ,1);C( 198 ):=F32([1, b, 3-2*n+a] , [-b+2, -a+2] ,1);C( 199 ):=F32([1, 2*n-2*a-1, 5/2-n-1/2*b] , [-b+2, -a+2] ,1);C( 200 ):=F32([b, 2*b, 3-2*n+a] , [b+1, 2*b+1-a] ,1);C( 201 ):=F32([b, -1+a+b, 5-a-b-2*n] , [a+b, b-a+1] ,1);C( 202 ):=F32([1, 1/2-2*a+n, 3/2-n+a] , [-a+3/2, 2*a+3/2-n] ,1);C( 203 ):=F32([1, 1-2*a+n, 3/2-n+a] , [-a+2, 2-n+2*a] ,1);C( 204 ):=F32([a-n+1, 1/2-a, 3/2-n+a] , [-a+3/2, 3*a-2*n+3/2] ,1);C( 205 ):=F32([-a+1, a-n+1, 3/2-n+a] , [-a+2, 3*a-2*n+2] ,1);C( 206 ):=F32([1+4*a-2*n, 1/2+2*a-n, 3/2-n+a] , [2*a+3/2-n, 3*a-2*n+3/2] ,1);C( 207 ):=F32([1, a-n+1, 1/2-2*a+n] , [-a+2, 2*a+3/2-n] ,1);C( 208 ):=F32([1, n-2*a, a-n+1] , [2*a+1-n, -a+3/2] ,1);C( 209 ):=F32([a-n+1, 1+4*a-2*n, 1/2+2*a-n] , [3*a-2*n+2, 2*a+3/2-n] ,1);C( 210 ):=F32([a, b, 1] , [-a+2, n-b] ,1);C( 211 ):=F32([a, -a+1, b] , [-a+2, n+a+2*b-3] ,1);C( 212 ):=F32([a, b, 2*b] , [b+1, n+2*b-1-a] ,1);C( 213 ):=F32([1, a, b] , [1/2*b+1, n+2*a-2] ,1);C( 214 ):=F32([a, b, -n+2*b+2] , [b+1, -n-a+2*b+3] ,1);C( 215 ):=F32([a, n-a-1, b] , [n-a, a+2*b-n+1] ,1);C( 216 ):=F32([1, a, b] , [1/2*n+1/2*a, -n+2*b+2] ,1);C( 217 ):=F32([a, b, -n+2*b+2*a+1] , [-n+a+2*b+2, 2*a+b] ,1);C( 218 ):=F32([a, b, 2*a+n-3] , [1/2*n+1/2*b+a-1, n+2*a-2] ,1);C( 219 ):=F32([a, 1, b] , [2*b, 5/2+1/2*a-n] ,1);C( 220 ):=F32([a, b, 2*b-1] , [2*b, b+1/2*a-n+3/2] ,1);C( 221 ):=F32([1, a, b] , [5-b-2*n, 5-a-2*n] ,1);C( 222 ):=F32([a, b, -b+4-2*n] , [b+1, 3+2*a-2*n-b] ,1);C( 223 ):=F32([a, b, 1/2*a-2*n+7/2-b] , [1/2*a+b+1/2, -b+4+a-2*n] ,1);C( 224 ):=F32([a, b, -n-a+3] , [1+a, 2*b-a] ,1);C( 225 ):=F32([a, n+2*a-2, b] , [1+a, -1-b+2*a+n] ,1);C( 226 ):=F32([a, b, 1] , [1/2*b-1/2*n+2, n+2*a-2] ,1);C( 227 ):=F32([a, b, 1] , [-n-a+4, -b+2] ,1);C( 228 ):=F32([a, b, -1/2*n-a+3/2+1/2*b] , [b+3-a-n, a+1/2*b+1/2*n-1/2] ,1);C( 229 ):=F32([a, b, 2*b] , [-n-a+2*b+3, b+1] ,1);C( 230 ):=F32([a, b, 1/2*b-1/2+1/2*n] , [1/2*a+1/2+1/2*b, b+1] ,1);C( 231 ):=F32([a, -a+1, b] , [a+2*b-n+1, -a+2] ,1);C( 232 ):=F32([a, b, -1+1/2*a+n] , [1+a, 1/2*a+1/2+1/2*b] ,1);C( 233 ):=F32([a, 2-n+1/2*a, b] , [1+a, 1/2*a+1/2*b+2-n] ,1);C( 234 ):=F32([a, b, 1] , [2*a-2*n+3, 1/2*b+1] ,1);C( 235 ):=F32([a, b, 2-n+1/2*a-b] , [1/2*a+n+b-1, -b+4+a-2*n] ,1);C( 236 ):=F32([a, 1/2*a-n+3/2, b] , [1+a-b, 5/2+1/2*a-n] ,1);C( 237 ):=F32([a, b, -b+1] , [3+2*a-2*n-b, b+1] ,1);C( 238 ):=F32([a, n+b-2, 1/2*n+1/2*a-1] , [1+a-b, 1/2*n+1/2*a] ,1);C( 239 ):=F32([a, -2*b+1, 2*n-2-a-b] , [1+a-b, -2*b+2*n-a-1] ,1);C( 240 ):=F32([a, 2*n-2-a, -1/2*b+1] , [2*n-1-a, 1+a-b] ,1);C( 241 ):=F32([a, 2*a-1, -n-2*b+4] , [2*a, 1+a-b] ,1);C( 242 ):=F32([1, -1/2*b+1, -n-2*a+4] , [-a+2, -b+2] ,1);C( 243 ):=F32([1, 2*n-3+b, a-3+2*n] , [-a+2, -b+2] ,1);C( 244 ):=F32([n, b-a+n, c-a+n] , [2+n, n+1-a] ,1);C( 245 ):=F32([2, 3-c, 3-b] , [3-a, 2+n] ,1);C( 246 ):=F32([a-c+3-b-n, -c+1, 2+a-c] , [3-b+a-c, 3-c] ,1);C( 247 ):=F32([a, b, n] , [c, n+2-c+b+a] ,1);C( 248 ):=F32([a, b, 2] , [c, n+2-c+b+a] ,1);C( 249 ):=F32([b, 2*b-a-m, b-a+c-n] , [1+b-c, b-a+1] ,1);C( 250 ):=F32([b, 2*b-a+m, b-a+c+n] , [1+b-c, b-a+1] ,1);C( 252 ):=F32([-3*a, 1-2*n, -a+n-1/2] , [-a-2*n+2, 1/2-3*a+n] ,1);C( 253 ):=F32([4/3-2*a, a-1/2*n+5/6, 4/3-a+1/2*n] , [5/3, -a-1/2*n+11/6] ,1);C( 254 ):=F32([1-a+1/2*n, a-1/2*n+5/6, 4/3-a+1/2*n] , [3/2, 4/3+a+1/2*n] ,1);C( 255 ):=F32([1-a+1/2*n, 2/3-a+1/2*n, a+1/6-1/2*n] , [3/2, 2/3+a+1/2*n] ,1);C( 256 ):=F32([4/3-2*a, 1-a+1/2*n, a-1/2*n+5/6] , [4/3, 3/2-a-1/2*n] ,1);C( 257 ):=F32([2/3-2*a, 1-a+1/2*n, a+1/6-1/2*n] , [2/3, 3/2-a-1/2*n] ,1);C( 258 ):=F32([2/3-2*a, -n+1/2, 4/3-2*a] , [1-1/2*n-a, 3/2-a-1/2*n] ,1);C( 259 ):=F32([a, c, b+n] , [b, e] ,1);C( 260 ):=F32([a, b, -n] , [c, e] ,1);C( 261 ):=F32([a, b, 1-b+m] , [c, 1+m-c+2*a-n] ,1);C( 262 ):=F32([a, b, c] , [2*c+n, 1/2+1/2*a+1/2*b+1/2*m-1/2*n] ,1);C( 263 ):=F32([a, b, -a+m+1-n] , [c, 1-c+2*b+m] ,1);C( 264 ):=F32([a, b, -b-m+1] , [c, -m+1-c+2*a+n] ,1);C( 265 ):=F32([a, b, c] , [2*c-n, 1/2+1/2*a+1/2*b-1/2*m+1/2*n] ,1);C( 266 ):=F32([a, b, 1-a+n-m] , [c, 1-c+2*b-m] ,1);C( 267 ):=F32([a-1/2*n, -n+1/2, a+1/3-1/2*n] , [4/3-a-1/2*n, 2*a+5/6-n] ,1);C( 268 ):=F32([a-1/2*n, 1-2*a, -a+5/6+1/2*n] , [4/3-a-1/2*n, 4/3] ,1);C( 269 ):=F32([-n+1/2, 1-2*a, 4/3-2*a] , [4/3-a-1/2*n, -a-1/2*n+11/6] ,1);C( 270 ):=F32([a+1/3-1/2*n, -a+5/6+1/2*n, 4/3-2*a] , [4/3-a-1/2*n, 5/3] ,1);C( 271 ):=F32([a-1/2*n, 1/2+a-1/2*n, 2*a+1/3] , [2*a+5/6-n, 4/3] ,1);C( 272 ):=F32([-n+1/2, 1/2+a-1/2*n, a-1/2*n+5/6] , [2*a+5/6-n, -a-1/2*n+11/6] ,1);C( 273 ):=F32([a+1/3-1/2*n, 2*a+1/3, a-1/2*n+5/6] , [2*a+5/6-n, 5/3] ,1);C( 274 ):=F32([1-2*a, 1/2+a-1/2*n, 4/3-a+1/2*n] , [4/3, -a-1/2*n+11/6] ,1);C( 275 ):=F32([-a+5/6+1/2*n, 2*a+1/3, 4/3-a+1/2*n] , [4/3, 5/3] ,1);C( 276 ):=F32([2*a-1/3, 2*a, a+1/6-1/2*n] , [2/3+a+1/2*n, 3*a+1/2-1/2*n] ,1);C( 277 ):=F32([2*a-1/3, 1/2+n, 2/3-a+1/2*n] , [2/3+a+1/2*n, 1+a+1/2*n] ,1);C( 278 ):=F32([2*a, 1/2+n, 1-a+1/2*n] , [2/3+a+1/2*n, 4/3+a+1/2*n] ,1);C( 279 ):=F32([2*a-1/3, 2*a+1/3, 1/2+a-1/2*n] , [3*a+1/2-1/2*n, 1+a+1/2*n] ,1);C( 280 ):=F32([2*a, 2*a+1/3, a-1/2*n+5/6] , [3*a+1/2-1/2*n, 4/3+a+1/2*n] ,1);C( 281 ):=F32([a+1/6-1/2*n, 1/2+a-1/2*n, a-1/2*n+5/6] , [3*a+1/2-1/2*n, 3/2] ,1);C( 282 ):=F32([1/2+n, 2*a+1/3, 4/3-a+1/2*n] , [1+a+1/2*n, 4/3+a+1/2*n] ,1);C( 283 ):=F32([2/3-a+1/2*n, 1/2+a-1/2*n, 4/3-a+1/2*n] , [1+a+1/2*n, 3/2] ,1);C( 284 ):=F32([a-1/3-1/2*n, 2/3-2*a, 1/2-a+1/2*n] , [1-1/2*n-a, 2/3] ,1);C( 285 ):=F32([a-1/3-1/2*n, -n+1/2, a+1/3-1/2*n] , [1-1/2*n-a, 1/2+2*a-n] ,1);C( 286 ):=F32([1/2-a+1/2*n, a+1/3-1/2*n, 4/3-2*a] , [1-1/2*n-a, 4/3] ,1);C( 287 ):=F32([a-1/3-1/2*n, a+1/6-1/2*n, 2*a] , [2/3, 1/2+2*a-n] ,1);C( 288 ):=F32([1/2-a+1/2*n, 2*a, 1-a+1/2*n] , [2/3, 4/3] ,1);C( 289 ):=F32([-n+1/2, a+1/6-1/2*n, a-1/2*n+5/6] , [1/2+2*a-n, 3/2-a-1/2*n] ,1);C( 290 ):=F32([a+1/3-1/2*n, 2*a, a-1/2*n+5/6] , [1/2+2*a-n, 4/3] ,1);C( 293 ):=F32([4/3-2*a, 2/3-2*a, 1/2-a+1/2*n] , [1-1/2*n-a, 3/2-3*a+1/2*n] ,1);C( 294 ):=F32([2/3-2*a, 1-a+1/2*n, 1/2-a+1/2*n] , [2/3, 7/6-2*a+n] ,1);C( 295 ):=F32([2/3-2*a, 2/3-a+1/2*n, 1/6-a+1/2*n] , [1/3, 7/6-2*a+n] ,1);C( 296 ):=F32([a, 1+a-c, a+c-2*b-m] , [1+a-b, 2*a-m+n] ,1);C( 297 ):=F32([a, 1+a-c, a+c-2*b+m] , [1+a-b, 2*a-n+m] ,1);C( 298 ):=F32([1-2*a, 2/3-2*a, 2/3-a+1/2*n] , [2-3*a+1/2*n, 7/6-a-1/2*n] ,1);C( 299 ):=F32([1-2*a, 2/3-a+1/2*n, 1/6-a+1/2*n] , [2/3, 3/2-2*a+n] ,1);C( 300 ):=F32([1/2+n, -a+5/6+1/2*n, 1/6-a+1/2*n] , [3/2-2*a+n, a+1/2+1/2*n] ,1);C( 301 ):=F32([1-2*a, 2/3-2*a, 1/6-a+1/2*n] , [3/2-3*a+1/2*n, 2/3-a-1/2*n] ,1);C( 302 ):=F32([-a+5/6+1/2*n, 1/2-a+1/2*n, 1/6-a+1/2*n] , [1/2, 3/2-3*a+1/2*n] ,1);C( 303 ):=F32([1/2+n, -a+5/6+1/2*n, 1/2-a+1/2*n] , [11/6+n-2*a, 5/6+a+1/2*n] ,1);C( 304 ):=F32([1/2+n, 1/2-a+1/2*n, 1/6-a+1/2*n] , [7/6-2*a+n, 1/6+a+1/2*n] ,1);C( 305 ):=F32([a-1/2*n, 1-2*a, -n+1/2] , [4/3-a-1/2*n, 2/3-a-1/2*n] ,1);C( 306 ):=F32([a-1/2*n, -a+5/6+1/2*n, a+1/3-1/2*n] , [4/3-a-1/2*n, 1/2] ,1);C( 307 ):=F32([1-2*a, -a+5/6+1/2*n, 4/3-2*a] , [4/3-a-1/2*n, 3/2-3*a+1/2*n] ,1);C( 308 ):=F32([-n+1/2, a+1/3-1/2*n, 4/3-2*a] , [4/3-a-1/2*n, 1-1/2*n-a] ,1);C( 309 ):=F32([a-1/2*n, 1/6-a+1/2*n, a-1/3-1/2*n] , [2/3-a-1/2*n, 1/2] ,1);C( 310 ):=F32([-n+1/2, a-1/3-1/2*n, 2/3-2*a] , [2/3-a-1/2*n, 1-1/2*n-a] ,1);C( 311 ):=F32([a+1/3-1/2*n, a-1/3-1/2*n, 1/2-a+1/2*n] , [1/2, 1-1/2*n-a] ,1);C( 312 ):=F32([2*a-1/3, 2/3-a+1/2*n, a+1/6-1/2*n] , [1/3, 2/3+a+1/2*n] ,1);C( 313 ):=F32([2*a-1/3, 1/6-a+1/2*n, a-1/3-1/2*n] , [1/3, 1/6+a+1/2*n] ,1);C( 314 ):=F32([a+1/6-1/2*n, a-1/3-1/2*n, 2/3-2*a] , [1/3, 2/3] ,1);C( 315 ):=F32([2*a-1/3, 1/2+n, 2*a] , [2/3+a+1/2*n, 1/6+a+1/2*n] ,1);C( 316 ):=F32([2/3-a+1/2*n, 1/2+n, 1-a+1/2*n] , [2/3+a+1/2*n, 7/6-2*a+n] ,1);C( 317 ):=F32([a+1/6-1/2*n, 2*a, 1-a+1/2*n] , [2/3+a+1/2*n, 2/3] ,1);C( 318 ):=F32([a-1/3-1/2*n, 2*a, 1/2-a+1/2*n] , [1/6+a+1/2*n, 2/3] ,1);C( 319 ):=F32([a-1/2*n+5/6, -n+1/2, 4/3-2*a] , [-a-1/2*n+11/6, 3/2-a-1/2*n] ,1);C( 320 ):=F32([a-1/2*n+5/6, 1/2+a-1/2*n, 4/3-a+1/2*n] , [-a-1/2*n+11/6, 3/2] ,1);C( 321 ):=F32([-n+1/2, 1/2+a-1/2*n, 1-2*a] , [-a-1/2*n+11/6, 7/6-a-1/2*n] ,1);C( 322 ):=F32([4/3-2*a, 4/3-a+1/2*n, 1-2*a] , [-a-1/2*n+11/6, 2-3*a+1/2*n] ,1);C( 323 ):=F32([a-1/2*n+5/6, a+1/6-1/2*n, 1-a+1/2*n] , [3/2-a-1/2*n, 3/2] ,1);C( 324 ):=F32([-n+1/2, a+1/6-1/2*n, 2/3-2*a] , [3/2-a-1/2*n, 7/6-a-1/2*n] ,1);C( 325 ):=F32([4/3-2*a, 1-a+1/2*n, 2/3-2*a] , [3/2-a-1/2*n, 2-3*a+1/2*n] ,1);C( 326 ):=F32([1/2+a-1/2*n, a+1/6-1/2*n, 2/3-a+1/2*n] , [3/2, 7/6-a-1/2*n] ,1);C( 327 ):=F32([4/3-a+1/2*n, 1-a+1/2*n, 2/3-a+1/2*n] , [3/2, 2-3*a+1/2*n] ,1);C( 328 ):=F32([2*a+1/3, 4/3-a+1/2*n, 1/2+a-1/2*n] , [4/3, 1+a+1/2*n] ,1);C( 329 ):=F32([2*a+1/3, -a+5/6+1/2*n, a-1/2*n] , [4/3, a+1/2+1/2*n] ,1);C( 330 ):=F32([4/3-a+1/2*n, -a+5/6+1/2*n, 1-2*a] , [4/3, 3/2-2*a+n] ,1);C( 331 ):=F32([1/2+a-1/2*n, a-1/2*n, 1-2*a] , [4/3, 2/3] ,1);C( 332 ):=F32([2*a+1/3, 1/2+n, 2*a-1/3] , [1+a+1/2*n, a+1/2+1/2*n] ,1);C( 333 ):=F32([4/3-a+1/2*n, 1/2+n, 2/3-a+1/2*n] , [1+a+1/2*n, 3/2-2*a+n] ,1);C( 334 ):=F32([1/2+a-1/2*n, 2*a-1/3, 2/3-a+1/2*n] , [1+a+1/2*n, 2/3] ,1);C( 335 ):=F32([a-1/2*n, 2*a-1/3, 1/6-a+1/2*n] , [a+1/2+1/2*n, 2/3] ,1);C( 336 ):=F32([a-1/2*n+5/6, 2*a+1/3, 4/3-a+1/2*n] , [4/3+a+1/2*n, 5/3] ,1);C( 337 ):=F32([a-1/2*n+5/6, 2*a, 1-a+1/2*n] , [4/3+a+1/2*n, 4/3] ,1);C( 338 ):=F32([2*a+1/3, 2*a, 1/2+n] , [4/3+a+1/2*n, 5/6+a+1/2*n] ,1);C( 339 ):=F32([4/3-a+1/2*n, 1-a+1/2*n, 1/2+n] , [4/3+a+1/2*n, 11/6+n-2*a] ,1);C( 340 ):=F32([a-1/2*n+5/6, a+1/3-1/2*n, 4/3-2*a] , [5/3, 4/3] ,1);C( 341 ):=F32([2*a+1/3, a+1/3-1/2*n, -a+5/6+1/2*n] , [5/3, 5/6+a+1/2*n] ,1);C( 342 ):=F32([4/3-a+1/2*n, 4/3-2*a, -a+5/6+1/2*n] , [5/3, 11/6+n-2*a] ,1);C( 343 ):=F32([2*a, a+1/3-1/2*n, 1/2-a+1/2*n] , [4/3, 5/6+a+1/2*n] ,1);C( 344 ):=F32([1-a+1/2*n, 4/3-2*a, 1/2-a+1/2*n] , [4/3, 11/6+n-2*a] ,1);C( 345 ):=F32([2/3-2*a, a+1/6-1/2*n, 2/3-a+1/2*n] , [1/3, 7/6-a-1/2*n] ,1);C( 346 ):=F32([1-2*a, 2/3-2*a, -n+1/2] , [2/3-a-1/2*n, 7/6-a-1/2*n] ,1);C( 347 ):=F32([1-2*a, 2/3-a+1/2*n, 4/3-a+1/2*n] , [2-3*a+1/2*n, 3/2-2*a+n] ,1);C( 348 ):=F32([1-2*a, 2/3-a+1/2*n, 1/2+a-1/2*n] , [2/3, 7/6-a-1/2*n] ,1);C( 349 ):=F32([-n+1/2, a+1/6-1/2*n, 1/2+a-1/2*n] , [2*a+1/6-n, 7/6-a-1/2*n] ,1);C( 350 ):=F32([1/2+n, 1/6-a+1/2*n, 2/3-a+1/2*n] , [7/6-2*a+n, 3/2-2*a+n] ,1);C( 351 ):=F32([2*a-1/3, a-1/3-1/2*n, a+1/6-1/2*n] , [1/3, 2*a+1/6-n] ,1);C( 352 ):=F32([1-2*a, 1/6-a+1/2*n, -a+5/6+1/2*n] , [3/2-2*a+n, 3/2-3*a+1/2*n] ,1);C( 353 ):=F32([a-1/2*n, -n+1/2, a-1/3-1/2*n] , [2*a+1/6-n, 2/3-a-1/2*n] ,1);C( 354 ):=F32([a-1/2*n, 2*a-1/3, 1/2+a-1/2*n] , [2/3, 2*a+1/6-n] ,1);C( 355 ):=F32([4/3-2*a, 1-a+1/2*n, 4/3-a+1/2*n] , [2-3*a+1/2*n, 11/6+n-2*a] ,1);C( 356 ):=F32([4/3-2*a, 1/2-a+1/2*n, -a+5/6+1/2*n] , [11/6+n-2*a, 3/2-3*a+1/2*n] ,1);C( 357 ):=F32([1/2+n, 1-a+1/2*n, 1/2-a+1/2*n] , [7/6-2*a+n, 11/6+n-2*a] ,1);C( 358 ):=F32([a-1/2*n, 1-2*a, 1/6-a+1/2*n] , [2/3-a-1/2*n, 2/3] ,1);C( 359 ):=F32([a-1/3-1/2*n, 1/6-a+1/2*n, 2/3-2*a] , [2/3-a-1/2*n, 1/3] ,1);C( 360 ):=F32([1/6-a+1/2*n, 2*a-1/3, 2/3-a+1/2*n] , [2/3, 1/3] ,1);C( 361 ):=F32([1/2-a+1/2*n, 1/6-a+1/2*n, 2/3-2*a] , [3/2-3*a+1/2*n, 7/6-2*a+n] ,1);C( 362 ):=F32([2/3-2*a, 4/3-2*a, 1-2*a] , [3/2-3*a+1/2*n, 2-3*a+1/2*n] ,1);C( 363 ):=F32([2/3-2*a, 1-a+1/2*n, 2/3-a+1/2*n] , [7/6-2*a+n, 2-3*a+1/2*n] ,1);C( 364 ):=F32([-a+5/6+1/2*n, 1/2+n, 4/3-a+1/2*n] , [11/6+n-2*a, 3/2-2*a+n] ,1);C( 366 ):=F32([a, -a+2*n, 3/2+a-1/2*c-1/2*b-n-m] , [1+a-c, 1+a-b] ,1);C( 367 ):=F32([a, -a+2*n, 1+a-1/2*c-1/2*b-n-m] , [1+a-c, 1+a-b] ,1);C( 368 ):=F32([a, 1-a+2*n, 1+a-1/2*c-1/2*b-n-m] , [1+a-c, 1+a-b] ,1);C( 375 ):=F32([a, b, -1+2*n-b+2*m] , [c, -c+2*a+2*n] ,1);C( 376 ):=F32([a, b, c] , [-a-1+b+2*n+2*m, -c+b+2*n] ,1);C( 377 ):=F32([a, b, c] , [2*c+2*m-1, 1/2*a+1/2*b+n] ,1);C( 378 ):=F32([a, b, 2*n-b+2*m] , [c, -c+2*a+2*n] ,1);C( 379 ):=F32([a, b, c] , [-a+b+2*n+2*m, -c+b+2*n] ,1);C( 380 ):=F32([a, b, c] , [2*c+2*m, 1/2*a+1/2*b+n] ,1);C( 381 ):=F32([a, b, 2*n-b+2*m] , [c, -c+2*a+2*n+1] ,1);C( 382 ):=F32([a, b, c] , [-a+b+2*n+2*m, 1-c+b+2*n] ,1);C( 383 ):=F32([a, b, c] , [2*c+2*m-1, 1/2+1/2*a+n+1/2*b] ,1);C( 384 ):=F32([b, -2*a+2*n+b, 3/2-1/2*c+1/2*b-n-m] , [1+b-c, b-a+1] ,1);C( 385 ):=F32([b, -2*a+2*n+b, 1-1/2*c+1/2*b-n-m] , [1+b-c, b-a+1] ,1);C( 386 ):=F32([b, 1+b-2*a+2*n, 1-1/2*c+1/2*b-n-m] , [1+b-c, b-a+1] ,1);C( 387 ):=F32([c-2*n+1, 1+2*c-b-2*n, a+1-b+c-2*n] , [-a+2*m+c, 2+c-b-2*n] ,1);C( 388 ):=F32([1/2*a-1/2*b-n+1, -1/2*a-n+1+1/2*b, 1+c-1/2*a-1/2*b-n] , [2-1/2*a-1/2*b-n, 2*c+2*m-1/2*a-1/2*b-n] ,1);C( 389 ):=F32([c-2*n+1, 1+2*c-b-2*n, a+1-b+c-2*n] , [1-a+2*m+c, 2+c-b-2*n] ,1);C( 390 ):=F32([1/2*a-1/2*b-n+1, -1/2*a-n+1+1/2*b, 1+c-1/2*a-1/2*b-n] , [2-1/2*a-1/2*b-n, 1+2*c+2*m-1/2*a-1/2*b-n] ,1);C( 391 ):=F32([c-2*n, 2*c-b-2*n, a-b+c-2*n] , [-a+2*m+c, 1+c-b-2*n] ,1);C( 392 ):=F32([1/2-1/2*a-n+1/2*b, 1/2+1/2*a-n-1/2*b, 1/2+c-1/2*a-n-1/2*b] , [3/2-1/2*a-n-1/2*b, 2*c+2*m-1/2-1/2*a-n-1/2*b] ,1);C( 393 ):=F32([a, -a+2*m, b] , [c, 2-c+2*b-2*n] ,1);C( 394 ):=F32([a, b, c] , [1/2*c+m+1/2*b, 2+2*a-2*n-2*m] ,1);C( 395 ):=F32([a, 1-a+2*m, b] , [c, 2-c+2*b-2*n] ,1);C( 396 ):=F32([a, b, c] , [1/2*b+1/2+1/2*c+m, 1+2*a-2*n-2*m] ,1);C( 397 ):=F32([a, -a+2*m, b] , [c, 1-c+2*b-2*n] ,1);C( 398 ):=F32([a, b, c] , [1/2*c+m+1/2*b, 1+2*a-2*n-2*m] ,1);C( 399 ):=F32([b, c-2*n+1, 2-2*n+a-2*m] , [1+b-c, b-a+1] ,1);C( 400 ):=F32([b, c-2*n+1, 1-2*n+a-2*m] , [1+b-c, b-a+1] ,1);C( 401 ):=F32([b, c-2*n, 1-2*n+a-2*m] , [1+b-c, b-a+1] ,1);C( 402 ):=F32([a, 1+a-c, a+c-2*b+2*n] , [1+a-b, 1+2*a-2*m] ,1);C( 403 ):=F32([a, b, c] , [1+2*c-2*m, 1/2*a+1/2*b-n+1] ,1);C( 404 ):=F32([a, b, c] , [2*b-1+2*m, 1/2*a+1/2*c+3/2-n-m] ,1);C( 405 ):=F32([a, b, -b+2-2*n] , [c, 3-c+2*a-2*n-2*m] ,1);C( 406 ):=F32([a, b, 3-2*n-b-2*m] , [c, -c+2*a+2-2*n] ,1);C( 407 ):=F32([a, b, c] , [-2*m+2*c, 1/2*a+1/2*b-n+1] ,1);C( 408 ):=F32([a, b, c] , [2*m+2*b, 1/2*c+1/2*a-n-m+1] ,1);C( 409 ):=F32([a, b, -b+2-2*n] , [c, 2*a-c-2*n-2*m+2] ,1);C( 410 ):=F32([a, b, 2-2*n-b-2*m] , [c, -c+2*a+2-2*n] ,1);C( 411 ):=F32([a, b, c] , [2*b-1+2*m, 1/2*c+1/2*a-n-m+1] ,1);C( 412 ):=F32([a, b, -b+1-2*n] , [c, 2*a-c-2*n-2*m+2] ,1);C( 413 ):=F32([a, b, 2-2*n-b-2*m] , [c, 2*a-2*n+1-c] ,1);C( 414 ):=F32([2*a+1/3-n, 5/6+a-n, a-n+1/2] , [3*a-2*n+3/2, 2*a+5/6-n] ,1);C( 415 ):=F32([7/6+a-n, 2*a+2/3-n, a-n+1/2] , [3*a-2*n+3/2, 7/6+2*a-n] ,1);C( 416 ):=F32([2/3-a, 7/6+a-n, 2*a+2/3-n] , [4/3, 5/3+a-n] ,1);C( 417 ):=F32([2*a+1-n, 7/6+a-n, 5/6+a-n] , [3*a-2*n+3/2, 2*a+3/2-n] ,1);C( 418 ):=F32([-a+1, 2*a+1-n, 7/6+a-n] , [5/3, 5/3+a-n] ,1);C( 419 ):=F32([-a+1, 2*a+1-n, 5/6+a-n] , [4/3, 4/3+a-n] ,1);C( 420 ):=F32([a+2/3, a+1/3, 2*a+1-n] , [3*a+1-n, 2*a+3/2-n] ,1);C( 421 ):=F32([a+2/3, a, 2*a+2/3-n] , [3*a+1-n, 7/6+2*a-n] ,1);C( 422 ):=F32([a+1/3, a, 2*a+1/3-n] , [3*a+1-n, 2*a+5/6-n] ,1);C( 423 ):=F32([2*a+1-n, 2*a+2/3-n, 2*a+1/3-n] , [3*a+1-n, 3*a-2*n+3/2] ,1);C( 424 ):=F32([a+2/3, 1/2, 7/6+a-n] , [2*a+3/2-n, 7/6+2*a-n] ,1);C( 425 ):=F32([a+1/3, 1/2, 5/6+a-n] , [2*a+3/2-n, 2*a+5/6-n] ,1);C( 426 ):=F32([a, 1/2, a-n+1/2] , [7/6+2*a-n, 2*a+5/6-n] ,1);C( 427 ):=F32([1/2+n-a, -a+1, 2/3+n-2*a] , [7/6-2*a+n, 5/3] ,1);C( 428 ):=F32([1/2+n-a, 1/2, 1/6+n-a] , [7/6-2*a+n, 7/6+a] ,1);C( 429 ):=F32([-a+1, 1/2, 2/3-a] , [7/6-2*a+n, 5/3+a-n] ,1);C( 430 ):=F32([2/3+n-2*a, 1/6+n-a, 2/3-a] , [7/6-2*a+n, 4/3] ,1);C( 431 ):=F32([1/2+n-a, 2*a+1-n, a+2/3] , [5/3, 7/6+a] ,1);C( 432 ):=F32([2/3+n-2*a, a+2/3, 7/6+a-n] , [5/3, 4/3] ,1);C( 433 ):=F32([1/2, 2*a+1-n, 2*a+2/3-n] , [7/6+a, 5/3+a-n] ,1);C( 434 ):=F32([1/6+n-a, a+2/3, 2*a+2/3-n] , [7/6+a, 4/3] ,1);C( 435 ):=F32([-1/6+n-a, 1/2, 1/2+n-a] , [5/6+n-2*a, a+5/6] ,1);C( 436 ):=F32([-1/6+n-a, -a+1/3, 1/3+n-2*a] , [5/6+n-2*a, 2/3] ,1);C( 437 ):=F32([1/2, -a+1/3, -a+1] , [5/6+n-2*a, 4/3+a-n] ,1);C( 438 ):=F32([1/2+n-a, 1/3+n-2*a, -a+1] , [5/6+n-2*a, 4/3] ,1);C( 439 ):=F32([-1/6+n-a, 2*a+1/3-n, a+1/3] , [a+5/6, 2/3] ,1);C( 440 ):=F32([1/2, 2*a+1/3-n, 2*a+1-n] , [a+5/6, 4/3+a-n] ,1);C( 441 ):=F32([1/2+n-a, a+1/3, 2*a+1-n] , [a+5/6, 4/3] ,1);C( 442 ):=F32([-a+1/3, 2*a+1/3-n, 5/6+a-n] , [2/3, 4/3+a-n] ,1);C( 443 ):=F32([1/3+n-2*a, a+1/3, 5/6+a-n] , [2/3, 4/3] ,1);C( 444 ):=F32([n-2*a, -a+1/3, a-n+1/2] , [1/3, 5/6-a] ,1);C( 445 ):=F32([-a+1/3, 5/6+a-n, 1/3+n-2*a] , [2/3, 5/6-a] ,1);C( 446 ):=F32([1/2, 5/6+a-n, a-n+1/2] , [5/6-a, 2*a+5/6-n] ,1);C( 447 ):=F32([1/2, -a+1/3, 2*a+1/3-n] , [a-n+1, 4/3+a-n] ,1);C( 448 ):=F32([-a+1/3, 2*a+1/3-n, -1/6+n-a] , [1/3, 2/3] ,1);C( 449 ):=F32([-a+1, -a+1/3, 5/6+a-n] , [3/2-n, 4/3+a-n] ,1);C( 450 ):=F32([2/3-a, -a+1/3, a-n+1/2] , [a-n+1, 3/2-n] ,1);C( 451 ):=F32([2*a+1-n, 5/6+a-n, 2*a+1/3-n] , [3*a-2*n+3/2, 4/3+a-n] ,1);C( 452 ):=F32([7/6+a-n, 5/6+a-n, a-n+1/2] , [3/2-n, 3*a-2*n+3/2] ,1);C( 453 ):=F32([2*a+2/3-n, a-n+1/2, 2*a+1/3-n] , [a-n+1, 3*a-2*n+3/2] ,1);C( 454 ):=F32([1/2, 2/3-a, 2*a+2/3-n] , [a-n+1, 5/3+a-n] ,1);C( 455 ):=F32([a, a-n+1/2, 2*a+1/3-n] , [1/3, 2*a+5/6-n] ,1);C( 456 ):=F32([a+1/3, 1/3+n-2*a, -1/6+n-a] , [1/3+n-a, 2/3] ,1);C( 457 ):=F32([a+1/3, 1/2, a] , [1/3+n-a, 2*a+5/6-n] ,1);C( 458 ):=F32([1/3+n-2*a, 1/2, n-2*a] , [1/3+n-a, 5/6-a] ,1);C( 459 ):=F32([-1/6+n-a, a, n-2*a] , [1/3+n-a, 1/3] ,1);C( 460 ):=F32([a+1/3, 5/6+a-n, 2*a+1/3-n] , [2/3, 2*a+5/6-n] ,1);C( 461 ):=F32([7/6+a-n, -a+1, 2/3-a] , [5/3+a-n, 3/2-n] ,1);C( 462 ):=F32([7/6+a-n, 2*a+1-n, 2*a+2/3-n] , [5/3+a-n, 3*a-2*n+3/2] ,1);C( 463 ):=F32([-a+1, 2*a+1-n, 1/2] , [5/3+a-n, 4/3+a-n] ,1);C( 469 ):=F32([c, b, (a+c)*(a+b)/a] , [1+b+a+c, (a+a*c+a*b+b*c)/a] ,1);U:= (a**2-2*a*c-2*a*b+c**2-2*b*c+b**2)**(1/2) ;\end{maplegroup}